\documentclass{extarticle}

\usepackage[english]{babel}
\usepackage[T1]{fontenc}
\usepackage[utf8]{inputenc}
\usepackage{amsfonts}
\usepackage{amsmath}
\usepackage{dsfont}
\usepackage{hyperref}
\usepackage[thmmarks,amsmath]{ntheorem}
\usepackage{amssymb}
\usepackage{geometry}
\usepackage{xcolor}
\usepackage{tikz}
\usepackage{subcaption}
\usepackage{subfiles}
\usepackage{alphalph}
\usepackage{etoolbox}

\AtBeginDocument{%
  \AtBeginEnvironment{subequations}{%
  }
}

\makeatletter
\newcommand{\pushright}[1]{\ifmeasuring@#1\else\omit\hfill$\displaystyle{#1}$\fi\ignorespaces}
\newcommand{\pushleft}[1]{\ifmeasuring@#1\else\omit$\displaystyle{#1}$\hfill\fi\ignorespaces}
\makeatother

\numberwithin{equation}{section}

{
\newtheorem{Def}{Definition}[section]
}

{
\newtheorem{Lem}[Def]{\textit{Lemma}}
}

{
\theorembodyfont{}
\newtheorem{Rem}[Def]{\textit{Remark}}
}

{
\newtheorem{Cor}[Def]{\textit{Corollary}}
}

\newtheorem{Prop}[Def]{\textit{Proposition}}

\newtheorem{Theo}[Def]{Theorem}

{
\theoremstyle{break}
\theoremsymbol{$\square$}
\theorembodyfont{}
\newtheorem*{Dem}{Proof}
}

\allowdisplaybreaks

\title{Scattering of the 3D Zakharov-Kuznetsov equation}

\author{Philippe Anjolras}

\begin{document}

\maketitle

\begin{abstract}
We consider the Zakharov-Kuznetsov equation in space dimension $3$: 
\begin{align*}
\left\{ \begin{array}{l}
\partial_t u + \partial_x \Delta u + \partial_x \frac{u^2}{2} = 0 \\
u(t = 0) = u_0 
\end{array} \right. 
\end{align*}
where $u : (t, x, y) \in \mathbb{R} \times \mathbb{R} \times \mathbb{R}^2 \mapsto u(t, x, y) \in \mathbb{R}$, and $\Delta = \partial_x^2 + \Delta_y$ is the full Laplacian. We show that, for any $u_0$ satisfying 
\begin{align*}
\Vert (1 + x^2 + |y|^2) u_0 \Vert_{H^1} \ll 1 
\end{align*}
then the global solution exhibits scattering in $H^1$. This is done using the method of space-time resonances, and more precisely the partial symmetries approach \cite{GPWpartialsymmEuler} in order to treat the anisotropy. We introduce well suited anisotropic weighted norms, prove dispersive decay estimates adapted to these norms and an a priori estimate allowing to close by a bootstrap argument. 
\end{abstract}

\tableofcontents

\section{Introduction} 

We consider the Zakharov-Kuznetsov equation in dimension $3$: 
\begin{equation} \left\{ \begin{array}{l}
\partial_t u + \partial_x \Delta u + \partial_x \frac{u^2}{2} = 0 \\
u(t = 0, x) = u_0(x) 
\end{array} \right. \label{equ-ZK} \end{equation}
where $u : (t, x, y) \in \mathbb{R} \times \mathbb{R} \times \mathbb{R}^2 \mapsto u(t, x, y) \in \mathbb{R}$, and $\Delta = \partial_x^2 + \Delta_y^2$ is the full Laplacian. \eqref{equ-ZK} was introduced by Zakharov and Kuznetsov in \cite{ZKoriginal} to model acoustic waves in a plasma, and rigourously derived from the Euler-Poisson system by Lannes, Linares and Saut \cite{LLS-derivation}, see also Guo and Pu \cite{GuoPuDerivationKdV} for the 1D case and Pu \cite{PuDerivationZK} for the 3D case. As a generalisation of the Korteweg-de-Vries equation (which would correspond to $d = 1$), it preserves the $L^2$ norm and the energy 
\begin{align*}
\frac{1}{2} \int_{\mathbb{R}^3} |\nabla_{x, y} u(t, x, y)|^2 ~ dx dy ~ - \frac{1}{6} \int_{\mathbb{R}^3} u(t, x, y)^3 ~ dx dy 
\end{align*}
The equation is known to be locally (and therefore globally) well-posed in $L^2$ by a result of Herr and Kinoshita \cite{HerrKinoshitaLWP}. 

The goal of this paper is to obtain a stability result of the zero solution and scattering for small data, which was an open problem about the long-time behaviour of solutions of \eqref{equ-ZK}. Our main theorem is the following: 

\begin{Theo} There exist $\varepsilon_0, C > 0$ such that, for any $0 < \varepsilon \leq \varepsilon_0$, if 
\[ \Vert (1 + x^2 + |y|^2) u_0 \Vert_{H^1} \leq \varepsilon \]
then the global solution $u(t)$ of \eqref{equ-ZK} satisfies 
\[ \sup_{t > 0} \left( \Vert u(t) \Vert_{H^1} + t^{\frac{5}{6}} \langle t \rangle^{\frac{1}{5}} \Vert \partial_x u(t) \Vert_{L^{\infty}} \right) \leq C \varepsilon \]
Furthermore, there exist $f_{\infty} \in H^1$ such that $e^{t \partial_x \Delta} u(t) \underset{t \to \infty}{\overset{H^1}{\longrightarrow}} f_{\infty}$. \label{theo-principal}
\end{Theo}

\begin{Rem} \begin{itemize} 
\item Obviously, by reversibility of the equation, the same applies for negative times $t \to -\infty$ (with a different limit profile $f_{-\infty}$). 
\item If the initial data $u_0$ also satisfies $\Vert u_0 \Vert_{H^n} \ll 1$ for some integer $n \geq 1$, then $\Vert u(t) \Vert_{H^n}$ remains uniformly bounded in time and $f_{\infty} \in H^n$, so that scattering occurs in $H^s$ for any $s < n$. (In fact, $\Vert u_0 \Vert_{H^n}$ does not need to be small, as long as $\Vert (1 + x^2 + |y|^2) u_0 \Vert_{H^1}$ is sufficiently small.) 
\item This scattering result does \emph{not} hold in the usual weighted norms, i.e. we do not prove that $(x, y) \left( e^{t \partial_x \Delta} u(t) \right)$ converges as $t \to \infty$. In fact, we will not propagate the weighted norm above, or rather we will propagate only specific anisotropic weights. 
\end{itemize}
\end{Rem}

\paragraph{Previous works} Equation \eqref{equ-ZK} has been extensively studied over the last decades, in various space dimensions, as well as its modified or generalised versions: 
\begin{align*}
\partial_t u + \partial_x \Delta u + \partial_x u^{k+1} = 0 
\end{align*}
for $k \geq 1$, where the case $k = 2$ of the cubic nonlinearity is called the modified Zakharov-Kuznetsov equation and other values than $k = 1, 2$ are called the generalized Zakharov-Kuznetsov equation. The generalized Zakharov-Kuznetsov equation with nonlinear parameter $k$ enjoys the scaling invariance
\begin{align*}
u_{\lambda}(t, x, y) &= \lambda^{\frac{2}{k}} u(\lambda^3 t, \lambda x, \lambda y)
\end{align*}
so that, in space dimension $d$, its critical Sobolev regularity is $s_c = \frac{d}{2}-\frac{2}{k}$, for \eqref{equ-ZK} this is $s_c = -\frac{1}{2}$. In dimension $d = 3$, the first local well-posedness result is due to Linares and Saut \cite{LinaresSautLWP} in $H^s$, $s > \frac{9}{8}$, which was then improved to $s > 1$ by Ribaud and Vento \cite{RibaudVentoLWP3d}. Molinet and Pilot \cite{MolinetPilodLWP} proved global well-posedness for $s > 1$, before Herr and Kinoshita \cite{HerrKinoshitaLWP} recently proved local well-posedness in any $d \geq 3$ for $s > \frac{d-4}{2}$, which is optimal up to the endpoint. In particular, \eqref{equ-ZK} is globally well-posed in $L^2$ as the $L^2$ norm is preserved. In dimension $d = 2$, the first result is due to Faminskii \cite{FaminskiiLWP} of global well-posedness in $H^1$, Grünrock and Herr \cite{GrunrockHerrLWP} and Molinet and Pilod \cite{MolinetPilodLWP} proved local well-posedness for $s > \frac{1}{2}$, and Kinoshita \cite{KinoshitaGWP} proved local well-posedness for $s > -\frac{1}{4}$, which is optimal up to the endpoint. For other dimensions and the generalized equation, we refer to \cite{LinaresPastormZK2DCauchy, LinaresPastorLWP2DgZK, FarahLP2DgZKCauchy, RibaudVentogZKCauchy, GrunrockLWPmZK3D, GrunrockLWPgZK23D, KinoshitaGWPmZK, HerrKinLWPZKhighd, LinaresRamosLWPgZK}. 

For the long-time behaviour of solutions, the soliton was proven to be orbitally stable by de Bouard \cite{deBouardOrbstab} and to be asymptotically stable recently by Farah, Holmer, Roudenko and Yang in \cite{FHRYStability3DsolitaryZK} for the three-dimensional case, building on a result by C{\^o}te, Mu{\~n}oz, Pilod and Simpson for the two-dimensional case \cite{CoteStab2DZKsoliton}: to avoid the radiative term, the authors restrict the stability result to a finite $(x, y)$ area in the direction of propagation of linear waves, noting that the solitary waves travel in the opposite direction. In \cite{PilodValetStabmultisolitonZK}, Pilod and Valet extended the result to a finite sum of (well-ordered) solitons. Scattering results on the other hand were only obtained either for the generalized equation \cite{FarahLP2DgZKCauchy}, or for the modified equation, in dimension $1$ that is for mKdV \cite{HayashiNaumkinmKdV1, HayashiNaumkinmKdV2, HarropGriffLongtimemKdV, GermainPusateriRousset_mKDV} or in dimension $2$ by the author \cite{AnjolrasmZK} and by Correia and Kinoshita \cite{CorreiaKinoshitaGWPmZK}. Very recently, Segata obtained scattering results in two \cite{SegataZK2D} and three dimensions \cite{SegataZK3D} by constructing solutions that scatters to a given asymptotic profile, that is by fixing $f_{\infty}$ with the notations of Theorem \ref{theo-principal} instead of fixing the initial data $u_0$. However, the functional space in which this scattering holds, although not directly comparable to the norm used in Theorem \ref{theo-principal}, requires more regularity and is singular in the sense that it does not contain all Schwartz functions (nor all smooth, compactly supported functions). 

The method of proof is based on the space-time resonances method introduced by Germain, Masmoudi and Shatah in \cite{Germain_3DSchrod, Germain_2DSchrod, GMSWaterwaves} for the quadratic Schrödinger equation and gravity water waves, and simultaneously by Gustafson, Nakanishi and Tsai in \cite{NakanishiResonances, NakanishiResonances2} for the Gross-Pitaevskii equation. The space-time resonance method has been widely used in the recent years to prove global regularity and scattering of nonlinear PDEs, e.g. water waves \cite{GMSWaterwaves,IonescuPusateri_waterwaves,GMS_waterwavescap}, the Klein-Gordon equation \cite{IonescuPausader_KG}, the wave equation \cite{ShatahPusateri,AnjolrasABI}, the Euler-Maxwell system \cite{Germain_EulerMaxwell, GuoIonescuPausader_plasma,GuoIonescuPausader_EulerMaxwell,DengIonescuPausader_EulerMaxwell}, the Euler-Poisson system \cite{GuoPausader_EulerPoisson,GuoIonescuPausader_plasma}, the modified Korteweg-de-Vries equation \cite{GermainPusateriRousset_mKDV,StewartCmZKasymp}. Similar results concerning global regularity and scattering have also been obtained in \cite{DelortAlazard_gravitywaterwaves2015,HunterIfrimTataru_ww,IfrimTataru_ww}. There were also recent advances for cubic dispersive equations by Ifrim and Tataru concerning the long-time behaviour of solutions in dimension $1$, see for instance \cite{IfrimTataru1Dconj}. 

Due to the anisotropy of Equation \eqref{equ-ZK}, we will use the extension of the space-time resonances method introduced by Guo, Pausader and Widmayer \cite{GPWpartialsymmEuler} for the stability of rotating axisymmetric 3D Euler flows, named the partial symmetries method. This approach has been used by Ren and Tian \cite{RenTianEulerCoriolis3D} to extend the previous result to non-axisymmetric flow, by Jurja and Widmayer \cite{JurjaWidBoussinesq} for the 2D Boussinesq equation and by Ko \cite{KoRotatingNS} for the Navier-Stokes-Coriolis equation (uniformly in the viscosity parameter). We will explain the method more in details in Subsection \ref{subsection-outline} once the notations have been set. 

\paragraph{Structure of the paper} The next sections are devoted to the proof of Theorem \ref{theo-principal}. In Section \ref{section-setting-notations}, we introduce important notations and state the main a priori estimate, from which the result follows. At the end of this section, once the notations have been introduced, we also give an overview of the rest of the proof and describe the method. In Section \ref{section-harmonicanalysis}, we recall useful theorems of linear and multilinear harmonic analysis. Section  \ref{section-spacetimeresanal} is devoted to the analysis of space-time resonances in the general and in the anisotropic case. In Sections \ref{section-Linfdispersive} and \ref{section-L4disp}, we prove $L^{\infty}$ and $L^4$ (linear) dispersive estimates. We state and prove the key decomposition lemma in Section \ref{section-decomplemma}, which we then use in Section \ref{section-structurepoids} to give an expression for weighted terms. Section \ref{section-estimees-simples-fug}, \ref{section-estimeepoidssimple}, \ref{section-estimees-quadratiques-horsbb}, \ref{section-estimee-quad-bb}, \ref{section-estimee-mauvaise-b} and \ref{section-continuite} are then devoted to the proof of the a priori argument. 

\section{Setting of the bootstrap argument} \label{section-setting-notations} 

\subsection{Notations} 

We denote by $\overline{\xi} = (\xi_0, \xi) \in \mathbb{R} \times \mathbb{R}^2$ the Fourier variable associated to $(x, y)$. We will also decompose $\xi = (\xi_1, \xi_2)$ and, by analogy, $(x, y) = (x_0, x_1, x_2)$. We will denote by $\partial_x, \nabla_y$ the directional partial derivatives (or gradient), and by $\nabla = \nabla_{x, y}$ the full gradient. Also, we denote by $D = i \nabla$. 

Let us set $\omega(D) = - i \partial_x \Delta$ the symbol associated to the linear evolution of \eqref{equ-ZK}, so that the linear evolution corresponds to the Fourier multiplier $e^{-i t \omega\left( D \right)}$, and  
\begin{equation} \omega\left( \overline{\xi} \right) = \xi_0^3 + \xi_0 |\xi|^2 \label{equ-def-partielin} \end{equation}
We compute that 
\[ \mbox{det} \nabla^2 \omega(\overline{\xi}) = \mbox{det} \begin{pmatrix} 6 \xi_0 & 2 \xi^T \\ 2 \xi & 2 \xi_0 I_2 \end{pmatrix} = 8 \xi_0 (3 \xi_0^2 - |\xi|^2) \]
In particular, the linear evolution suffers from some dispersion degeneracy even away from $\overline{\xi} = 0$. This justifies the introduction of the two dual cones
\begin{equation} \mathcal{C} := \{ (x, y) \in \mathbb{R}^3, x^2 = 3 |y|^2 \}, \quad \widehat{\mathcal{C}} := \{ (\xi_0, \xi) \in \mathbb{R}^3, 3 \xi_0^2 = |\xi|^2 \} \label{equ-def-cones} \end{equation}
and of the dual plane and line
\begin{equation} \mathcal{L} := \{ (x, y) \in \mathbb{R}^3, y = 0 \}, \quad \widehat{\mathcal{P}} := \{ (\xi_0, \xi) \in \mathbb{R}^3, \xi_0 = 0 \} \label{equ-def-plan} \end{equation}
We will also need to consider the line 
\begin{equation} \widehat{\mathcal{L}} := \{ (\xi_0, \xi) \in \mathbb{R}^3, \xi = 0 \} \label{equ-def-ligne} \end{equation}
These sets are dual in the sense that 
\begin{align*}
\overline{\xi} \in \widehat{\mathcal{C}} ~ \iff ~ \nabla \omega(\overline{\xi}) \in \mathcal{C} 
\end{align*}
and 
\begin{align*}
\overline{\xi} \in \widehat{\mathcal{P}} ~ \Longrightarrow ~ \nabla \omega(\overline{\xi}) \in \mathcal{L} ~ \iff ~ \overline{\xi} \in \widehat{\mathcal{P}} \cup \widehat{\mathcal{L}} 
\end{align*}

Let us define the vector fields
\begin{equation} \begin{aligned}
\widehat{X}_a(\overline{\xi}) &:= \begin{pmatrix} \dfrac{\xi_0}{|\overline{\xi}|}, & \dfrac{\xi}{|\overline{\xi}|} \end{pmatrix} = \frac{\overline{\xi}}{|\overline{\xi}|} \\
\widehat{X}_b(\overline{\xi}) &:= \begin{pmatrix} \dfrac{|\xi|}{|\overline{\xi}|}, & -\dfrac{\xi_0 \xi}{|\overline{\xi}| |\xi|} \end{pmatrix} \\
\widehat{X}_c(\overline{\xi}) &:= \begin{pmatrix} 0,& \dfrac{J \xi}{|\xi|} \end{pmatrix}
\end{aligned} \label{equ-def-champsvecteurs}
\end{equation}
Here above, $J : \mathbb{R}^2 \to \mathbb{R}^2$ is the linear $2 \times 2$ matrix of rotation by the angle $\frac{\pi}{2}$. 
The vector field $\widehat{X}_a$ is the radial vector field and corresponds to the scaling invariance of the equation. The vector field $\widehat{X}_c$ corresponds to the axial symmetry in $y$. $(\widehat{X}_a, \widehat{X}_b, \widehat{X}_c)$ form an orthonormal basis at any point in $\mathbb{R}^3 \setminus \widehat{\mathcal{L}}$. 

\begin{figure}
\begin{center}
\scalebox{1.5}{
\begin{tikzpicture}[line join=round]
\filldraw[color=black!80,fill=red!50](0,0)--(-.082,-1.027)--(.566,-1.037)--(0,0)--cycle;
\filldraw[color=black!80,fill=red!50](0,0)--(-.724,-1.044)--(-.082,-1.027)--(0,0)--cycle;
\filldraw[color=black!80,fill=red!50](0,0)--(.566,-1.037)--(1.178,-1.076)--(0,0)--cycle;
\filldraw[color=black!80,fill=red!50](0,0)--(-1.321,-1.089)--(-.724,-1.044)--(0,0)--cycle;
\filldraw[color=black!80,fill=red!50](0,0)--(1.178,-1.076)--(1.716,-1.139)--(0,0)--cycle;
\filldraw[color=black!80,fill=red!50](0,0)--(-1.835,-1.159)--(-1.321,-1.089)--(0,0)--cycle;
\filldraw[color=black!80,fill=red!50](0,0)--(1.716,-1.139)--(2.146,-1.224)--(0,0)--cycle;
\filldraw[color=black!80,fill=red!50](0,0)--(-2.234,-1.248)--(-1.835,-1.159)--(0,0)--cycle;
\filldraw[color=black!80,fill=red!50](0,0)--(2.146,-1.224)--(2.442,-1.325)--(0,0)--cycle;
\filldraw[color=black!80,fill=red!50](0,0)--(-2.492,-1.351)--(-2.234,-1.248)--(0,0)--cycle;
\filldraw[color=black!80,fill=red!50](0,0)--(2.442,-1.325)--(2.583,-1.435)--(0,0)--cycle;
\filldraw[color=black!80,fill=red!50](0,0)--(-2.594,-1.463)--(-2.492,-1.351)--(0,0)--cycle;
\draw[arrows=-](0,-2.954)--(0,0);
\filldraw[color=black!80,fill=red!50](0,0)--(2.583,-1.435)--(2.563,-1.548)--(0,0)--cycle;
\draw[color=green!50,line width=2pt](0,-2.659)--(0,0);
\filldraw[color=black!80,fill=red!50](0,0)--(-2.533,-1.576)--(-2.594,-1.463)--(0,0)--cycle;
\filldraw[color=black!80,fill=red!50](0,0)--(2.563,-1.548)--(2.382,-1.656)--(0,0)--cycle;
\filldraw[color=black!80,fill=red!50](0,0)--(2.382,-1.656)--(2.05,-1.753)--(0,0)--cycle;
\filldraw[color=black!80,fill=red!50](0,0)--(2.05,-1.753)--(1.59,-1.833)--(0,0)--cycle;
\filldraw[color=black!80,fill=red!50](0,0)--(1.59,-1.833)--(1.031,-1.891)--(0,0)--cycle;
\filldraw[color=black!80,fill=red!50](0,0)--(1.031,-1.891)--(.406,-1.922)--(0,0)--cycle;
\filldraw[color=black!80,fill=red!50](0,0)--(.406,-1.922)--(-.244,-1.926)--(0,0)--cycle;
\filldraw[color=black!80,fill=red!50](0,0)--(-.244,-1.926)--(-.879,-1.901)--(0,0)--cycle;
\filldraw[color=black!80,fill=red!50](0,0)--(-.879,-1.901)--(-1.459,-1.85)--(0,0)--cycle;
\filldraw[color=black!80,fill=red!50](0,0)--(-1.459,-1.85)--(-1.947,-1.775)--(0,0)--cycle;
\filldraw[color=black!80,fill=red!50](0,0)--(-1.947,-1.775)--(-2.312,-1.682)--(0,0)--cycle;
\filldraw[color=black!80,fill=red!50](0,0)--(-2.312,-1.682)--(-2.533,-1.576)--(0,0)--cycle;
\filldraw[fill=blue!50](0,.737)--(-4.243,0)--(0,-.737)--(4.243,0)--cycle;
\draw[arrows=->](3.182,.553)--(-3.182,-.553);
\draw[arrows=->](-3.182,.553)--(3.182,-.553);
\filldraw[color=black!80,fill=red!50](0,0)--(-.406,1.922)--(.244,1.926)--(0,0)--cycle;
\filldraw[color=black!80,fill=red!50](0,0)--(.244,1.926)--(.879,1.901)--(0,0)--cycle;
\filldraw[color=black!80,fill=red!50](0,0)--(-1.031,1.891)--(-.406,1.922)--(0,0)--cycle;
\filldraw[color=black!80,fill=red!50](0,0)--(.879,1.901)--(1.459,1.85)--(0,0)--cycle;
\filldraw[color=black!80,fill=red!50](0,0)--(-1.59,1.833)--(-1.031,1.891)--(0,0)--cycle;
\filldraw[color=black!80,fill=red!50](0,0)--(1.459,1.85)--(1.947,1.775)--(0,0)--cycle;
\filldraw[color=black!80,fill=red!50](0,0)--(-2.05,1.753)--(-1.59,1.833)--(0,0)--cycle;
\filldraw[color=black!80,fill=red!50](0,0)--(1.947,1.775)--(2.312,1.682)--(0,0)--cycle;
\filldraw[color=black!80,fill=red!50](0,0)--(-2.382,1.656)--(-2.05,1.753)--(0,0)--cycle;
\draw[dotted,line width=0.5pt](1.326,2.078)--(1.482,2.139)--(1.545,2.206)--(1.511,2.273)--(1.382,2.336)--(1.166,2.392)--(.877,2.437)--(.533,2.468)--(.155,2.483)--(-.232,2.481)--(-.605,2.463)--(-.94,2.429)--(-1.216,2.382)--(-1.415,2.324)--(-1.526,2.26)--(-1.54,2.192)--(-1.458,2.126)--(-1.284,2.066)--(-1.03,2.016)--(-.711,1.977)--(-.347,1.954)--(.039,1.947)--(.422,1.958)--(.779,1.984)--(1.086,2.025)--(1.326,2.078);
\filldraw[color=black!80,fill=red!50](0,0)--(2.312,1.682)--(2.533,1.576)--(0,0)--cycle;
\filldraw[color=black!80,fill=red!50](0,0)--(-2.563,1.548)--(-2.382,1.656)--(0,0)--cycle;
\filldraw[color=black!80,fill=red!50](0,0)--(2.533,1.576)--(2.594,1.463)--(0,0)--cycle;
\filldraw[color=black!80,fill=red!50](0,0)--(-2.583,1.435)--(-2.563,1.548)--(0,0)--cycle;
\draw[color=green!50,line width=2pt](0,0)--(0,2.659);
\draw[dotted,line width=0.5pt](0,0)--(1.989,3.117);
\filldraw[color=black!80,fill=red!50](0,0)--(1.835,1.159)--(1.321,1.089)--(0,0)--cycle;
\filldraw[color=black!80,fill=red!50](0,0)--(1.321,1.089)--(.724,1.044)--(0,0)--cycle;
\filldraw[color=black!80,fill=red!50](0,0)--(.724,1.044)--(.082,1.027)--(0,0)--cycle;
\draw[arrows=->](0,0)--(0,2.954);
\filldraw[color=black!80,fill=red!50](0,0)--(.082,1.027)--(-.566,1.037)--(0,0)--cycle;
\filldraw[color=black!80,fill=red!50](0,0)--(-.566,1.037)--(-1.178,1.076)--(0,0)--cycle;
\filldraw[color=black!80,fill=red!50](0,0)--(-1.178,1.076)--(-1.716,1.139)--(0,0)--cycle;
\filldraw[color=black!80,fill=red!50](0,0)--(-1.716,1.139)--(-2.146,1.224)--(0,0)--cycle;
\filldraw[color=black!80,fill=red!50](0,0)--(-2.442,1.325)--(-2.583,1.435)--(0,0)--cycle;
\filldraw[color=black!80,fill=red!50](0,0)--(-2.146,1.224)--(-2.442,1.325)--(0,0)--cycle;
\filldraw[color=black!80,fill=red!50](0,0)--(2.594,1.463)--(2.492,1.351)--(0,0)--cycle;
\filldraw[color=black!80,fill=red!50](0,0)--(2.492,1.351)--(2.234,1.248)--(0,0)--cycle;
\filldraw[color=black!80,fill=red!50](0,0)--(2.234,1.248)--(1.835,1.159)--(0,0)--cycle;
\draw[arrows=->,color=magenta](1.326,2.078)--(.94,2.422);
\filldraw[line width=1pt,color=black](1.326,2.078)--(1.326,2.078)--(1.326,2.078)--(1.326,2.078)--cycle;
\filldraw[line width=1pt,color=black](1.326,2.078)--(1.326,2.078)--(1.326,2.078)--(1.326,2.078)--cycle;
\filldraw[line width=1pt,color=black](1.326,2.078)--(1.325,2.078)--(1.326,2.078)--(1.326,2.078)--cycle;
\filldraw[line width=1pt,color=black](1.326,2.078)--(1.326,2.078)--(1.327,2.078)--(1.326,2.078)--cycle;
\filldraw[line width=1pt,color=black](1.326,2.078)--(1.325,2.078)--(1.325,2.078)--(1.326,2.078)--cycle;
\filldraw[line width=1pt,color=black](1.326,2.078)--(1.327,2.078)--(1.327,2.078)--(1.326,2.078)--cycle;
\filldraw[line width=1pt,color=black](1.326,2.078)--(1.325,2.078)--(1.325,2.078)--(1.326,2.078)--cycle;
\filldraw[line width=1pt,color=black](1.326,2.078)--(1.327,2.078)--(1.327,2.078)--(1.326,2.078)--cycle;
\filldraw[line width=1pt,color=black](1.326,2.078)--(1.324,2.078)--(1.325,2.078)--(1.326,2.078)--cycle;
\filldraw[line width=1pt,color=black](1.326,2.078)--(1.327,2.078)--(1.327,2.078)--(1.326,2.078)--cycle;
\filldraw[line width=1pt,color=black](1.326,2.078)--(1.324,2.078)--(1.324,2.078)--(1.326,2.078)--cycle;
\filldraw[line width=1pt,color=black](1.326,2.078)--(1.327,2.078)--(1.327,2.078)--(1.326,2.078)--cycle;
\filldraw[line width=1pt,color=black](1.326,2.078)--(1.324,2.078)--(1.324,2.078)--(1.326,2.078)--cycle;
\filldraw[line width=1pt,color=black](1.326,2.078)--(1.327,2.078)--(1.327,2.078)--(1.326,2.078)--cycle;
\filldraw[line width=1pt,color=black](1.326,2.078)--(1.327,2.078)--(1.327,2.078)--(1.326,2.078)--cycle;
\filldraw[line width=1pt,color=black](1.326,2.078)--(1.326,2.077)--(1.325,2.077)--(1.326,2.078)--cycle;
\filldraw[line width=1pt,color=black](1.326,2.078)--(1.327,2.077)--(1.327,2.077)--(1.326,2.078)--cycle;
\filldraw[line width=1pt,color=black](1.326,2.078)--(1.326,2.077)--(1.326,2.077)--(1.326,2.078)--cycle;
\filldraw[line width=1pt,color=black](1.326,2.078)--(1.327,2.077)--(1.326,2.077)--(1.326,2.078)--cycle;
\draw[arrows=->,color=purple](1.326,2.078)--(1.591,2.493);
\filldraw[line width=1pt,color=black](1.326,2.078)--(1.327,2.078)--(1.327,2.077)--(1.326,2.078)--cycle;
\filldraw[line width=1pt,color=black](1.326,2.078)--(1.324,2.078)--(1.324,2.078)--(1.326,2.078)--cycle;
\filldraw[line width=1pt,color=black](1.326,2.078)--(1.325,2.078)--(1.324,2.078)--(1.326,2.078)--cycle;
\filldraw[line width=1pt,color=black](1.326,2.078)--(1.325,2.077)--(1.325,2.078)--(1.326,2.078)--cycle;
\filldraw[line width=1pt,color=black](1.326,2.078)--(1.325,2.077)--(1.325,2.077)--(1.326,2.078)--cycle;
\filldraw[line width=1pt,color=black](1.326,2.078)--(1.325,2.077)--(1.325,2.077)--(1.326,2.078)--cycle;
\draw[arrows=->,color=orange](1.326,2.078)--(1.045,1.996);

\node[] at (.424,2.881) {$\xi_0$};
\node[] at (3.182,-.109) {$\xi_1$};
\node[] at (-3.182,-.109) {$\xi_2$};
\node[] at (3.182,2.032) {\textcolor{red}{$\widehat{\mathcal{C}}$}};
\node[] at (4.243,.739) {\textcolor{blue}{$\widehat{\mathcal{P}}$}};
\node[] at (-.424,3.028) {\textcolor{green}{$\widehat{\mathcal{L}}$}};
\node[] at (2.,2.5) {\textcolor{purple}{$\widehat{X}_a$}};
\node[] at (.8,2.2) {\textcolor{orange}{$\widehat{X}_c$}};
\node[] at (1.,2.75) {\textcolor{magenta}{$\widehat{X}_b$}};
\end{tikzpicture}
}
\caption{The geometric areas and the vector fields} \label{figure3dZK}
\end{center}
\end{figure}
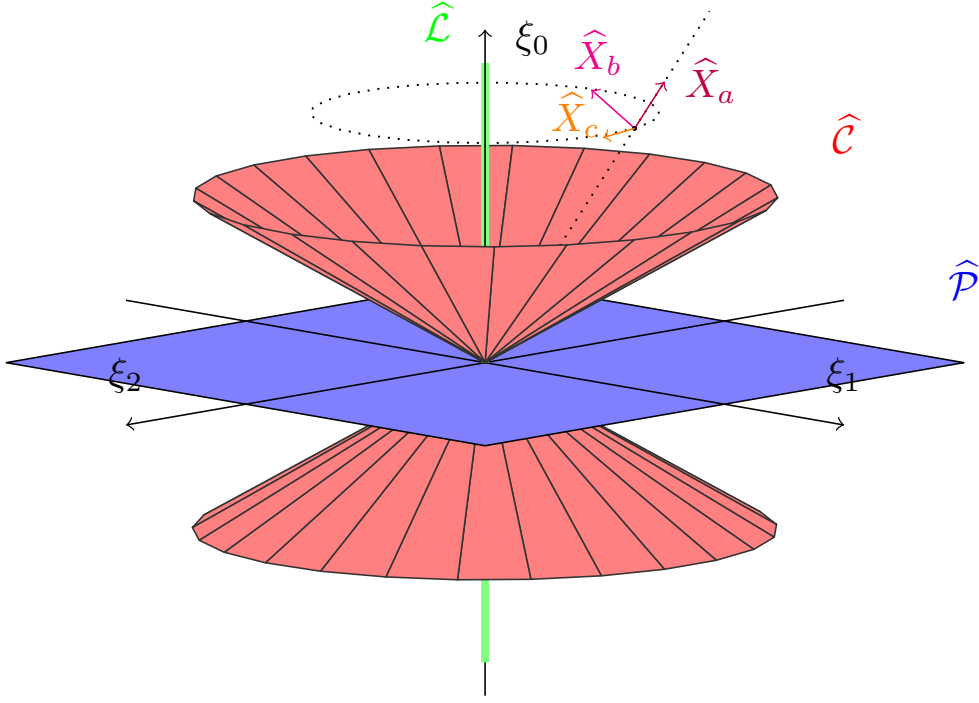

We represented the geometric areas and the vector fields on Figure \ref{figure3dZK}. 

We denote by $X_{\alpha}$ the operator defined as
\begin{align*}
X_{\alpha} F = \mathcal{F}^{-1} \left( \widehat{X}_{\alpha}(\overline{\xi}) \cdot \nabla_{\overline{\xi}} \widehat{F}(\overline{\xi}) \right) 
\end{align*}
for $\alpha = a, b, c$. 

Let $m_{\widehat{\mathcal{C}}}, m_{\widehat{\mathcal{P}}}, m_{\widehat{\mathcal{L}}}, m_{\widehat{\mathcal{R}}}$ be homogeneous functions of degree $0$, belonging to $C^{\infty}(\mathbb{R}^3 \setminus \{ 0 \}, [0, 1])$, localizing on (sufficiently small) angular neighborhoods of $\widehat{\mathcal{C}}, \widehat{\mathcal{P}}, \widehat{\mathcal{L}}$ respectively for the first three, and such that
\[ 1 = m_{\widehat{\mathcal{C}}}(\overline{\xi}) + m_{\widehat{\mathcal{P}}}(\overline{\xi}) + m_{\widehat{\mathcal{L}}}(\overline{\xi}) + m_{\widehat{\mathcal{R}}}(\overline{\xi}), \quad \forall ~ \overline{\xi} \in \mathbb{R}^3 \setminus \{ 0 \} \]

We will also denote by $m_a(\overline{\xi}) := 1$, $m_c(\overline{\xi}) := \frac{|\xi|}{|\overline{\xi}|}$ and
\[ m_b(\overline{\xi}) := \frac{|\xi| (3 \xi_0^2 - |\xi|^2)}{|\overline{\xi}|^3} \]
These are (homogeneous) Hörmander-Mikhlin symbols, with $m_c$ vanishing linearly on $\widehat{\mathcal{L}} \cap S^2$ and $m_b$ vanishing linearly on $\left( \widehat{\mathcal{L}} \cup \widehat{\mathcal{C}} \right) \cap S^2$. 

Let us set
\begin{equation} f(t) = e^{i t \omega(D)} u(t) \label{equ-def-fonctionprofil} \end{equation}
$f(t)$ is called the profile function associated to $u(t)$. Then $f$ satisfies the equation 
\begin{equation} \partial_t f(t) + \frac{1}{2} e^{i t \omega(D)} \partial_x \left( e^{-i t \omega(D)} f(t) \right)^2 = 0 \label{equ-ZK-fonctionprofil} \end{equation}

In particular, 
\[ \widehat{f}\left(t, \overline{\xi}\right) = \widehat{f}\left(0, \overline{\xi}\right) - \frac{1}{2} \int_0^t \int e^{i s \varphi\left(\overline{\xi}, \overline{\eta}\right)} i \xi_0 \widehat{f}\left(s, \overline{\eta}\right) \widehat{f}\left(s, \overline{\xi}-\overline{\eta}\right) d\overline{\eta} ds \]
by Duhamel's formula (up to a universal multiplicative constant due to the choice of the Fourier transform, which we will never take into account), and where
\begin{equation} \varphi\left(\overline{\xi}, \overline{\eta}\right) = \omega\left(\overline{\xi}\right) - \omega\left(\overline{\eta}\right) - \omega\left(\overline{\xi}-\overline{\eta}\right) = \xi_0^3 + \xi_0 |\xi|^2 - \eta_0^3 - \eta_0 |\eta|^2 - (\xi_0 - \eta_0)^3 - (\xi_0 - \eta_0) |\xi - \eta|^2 \label{equ-def-varphi2} \end{equation}
is the (quadratic) interaction phase. 

In the case of a $n$-linear interaction, we will consider products under the form 
\[ \int e^{i t \varphi_n\left( \overline{\eta}^0, ..., \overline{\eta}^{n-1} \right)} m\left( \overline{\eta}^0, ..., \overline{\eta}^{n-1} \right) \widehat{f}\left( t, \overline{\eta}^1 \right) ... \widehat{f}\left( t, \overline{\eta}^{n-1} \right) \widehat{f}\left( t, \overline{\eta}^n \right) d\overline{\eta}^1 ... d\overline{\eta}^{n-1} \]
where $\overline{\eta}^0 = \overline{\xi}$ and $\overline{\eta}^n = \overline{\xi} - \sum_{i = 1}^{n-1} \overline{\eta}^i$, for some symbol $m$, while 
\begin{equation} \varphi_n\left( \overline{\eta}^0, ..., \overline{\eta}^{n-1} \right) = \omega\left( \overline{\eta}^0 \right) - \sum_{i = 1}^n \omega\left( \overline{\eta}^i \right) \label{equ-def-varphin} \end{equation}
is the $n$-interaction phase. 

In the case $n = 2$, we will use the convention $\left( \overline{\eta}^0, \overline{\eta}^1, \overline{\eta}^2 \right) = \left( \overline{\xi}, \overline{\eta}, \overline{\sigma} \right)$; in the case $n = 3$, we will use the convention $\left( \overline{\eta}^0, \overline{\eta}^1, \overline{\eta}^2, \overline{\eta}^3 \right) = \left( \overline{\xi}, \overline{\eta}, \overline{\sigma}, \overline{\rho} \right)$. 

We will write $A \lesssim B$ if there exists a universal constant $C > 0$ such that $A \leq C B$. We will write $A \simeq B$ if $A \lesssim B$ and $B \lesssim A$. 

\subsection{A priori estimate}

\begin{Prop}[a priori estimate] There exist $C, \delta > 0$ as follows. Let $u_0$ be such that 
\[ \Vert u_0 \Vert_{X_0} := \Vert (1 + x^2 + |y|^2) u_0 \Vert_{H^1} < \infty \]
and let $u$ be the global solution of \eqref{equ-ZK} with this initial data. Denote by $f$ the profile function associated to $u$ through \eqref{equ-def-fonctionprofil}. 

For any $\alpha \in \{ a, b, c \}$, there exist a decomposition 
\begin{align*}
m_{\alpha}(D) X_{\alpha} f(t) = h_{\alpha}(t) + g_{\alpha}(t) 
\end{align*}
and for $\alpha \in \{ b, c \}$ : 
\begin{align*}
X_{\alpha} f(t) = h_{\alpha *}(t) + g_{\alpha *}(t) 
\end{align*}
such that, if we denote by 
\begin{align*}
\Vert u \Vert_{X, T, \delta} := \sup_{t \in [0, T]} \Bigl( &\Vert u(t) \Vert_{H^1} + \sum_{\alpha = a, b, c} \Vert h_{\alpha}(t) \Vert_{H^1} ~ + \sum_{(\alpha, \beta) \in \{ a, b, c \}^2 \setminus \{ (b, b) \}} \Vert \nabla m_{\alpha}(D) X_{\alpha} h_{\beta}(t) \Vert_{L^2} \\
&\quad + \langle t \rangle^{-\frac{1}{24}-301\delta} \Vert \partial_x m_{\widehat{\mathcal{R}}}(D) m_b(D) X_b h_b(t) \Vert_{L^2} ~ + \Vert \partial_x (1 - m_{\widehat{\mathcal{R}}}(D)) m_b(D) X_b h_b(t) \Vert_{L^2}  \\
&+ \sum_{\alpha = b, c} \Big( \langle t \rangle^{-\frac{5}{48}-202\delta} \Vert m_{\widehat{\mathcal{L}}}(D) h_{\alpha *}(t) \Vert_{L^2} 
+ \langle t \rangle^{-202 \delta} \Vert m_{\widehat{\mathcal{L}}}(D) |\nabla|^{\frac{1}{2}} \langle \nabla \rangle^{\frac{1}{2}} h_{\alpha *}(t) \Vert_{L^2} \Big) \\
&+ \langle t \rangle^{-\frac{1}{2}-101\delta} \Vert m_{\widehat{\mathcal{C}}}(D) h_{b*}(t) \Vert_{L^2}
+ \langle t \rangle^{-\frac{1}{6}-101\delta} \Vert m_{\widehat{\mathcal{C}}}(D) \nabla h_{b*}(t) \Vert_{L^2} \\
&+ \langle t \rangle^{\frac{1}{12}-100\delta} \sum_{\alpha = a, b, c} \Vert g_{\alpha}(t) \Vert_{H^1} + \langle t \rangle^{\frac{1}{6}-100\delta} \sum_{\alpha = a, b, c} \Vert |\nabla|^{\frac{1}{2}} g_{\alpha}(t) \Vert_{H^1} \\
&+ \langle t \rangle^{-\frac{5}{48}-202\delta} \Vert m_{\widehat{\mathcal{L}}}(D) g_{\alpha *}(t) \Vert_{L^2} 
+ \langle t \rangle^{-202 \delta} \Vert m_{\widehat{\mathcal{L}}}(D) |\nabla|^{\frac{1}{2}} \langle \nabla \rangle^{\frac{1}{2}} g_{\alpha *}(t) \Vert_{L^2} \\
&+ \langle t \rangle^{-\frac{1}{2}-101\delta} \Vert m_{\widehat{\mathcal{C}}}(D) g_{b*}(t) \Vert_{L^2} + \langle t \rangle^{-\frac{1}{6}-101\delta} \Vert m_{\widehat{\mathcal{C}}}(D) \nabla g_{b*}(t) \Vert_{L^2} \Bigl)
\end{align*}

Then: 
\begin{itemize}
\item (local existence) the function $T \mapsto \Vert u \Vert_{X, T, \delta}$ is defined on an interval $[0, T_{max})$ for some $T_{max} > 0$ and is continuous;
\item (a priori estimate) for any $T < T_{max}$, assuming $\Vert u \Vert_{X, T, \delta} \ll 1$, then 
\[ \Vert u \Vert_{X, T, \delta} ~ \leq C \left( \Vert u \Vert_{X_0} + \Vert u \Vert_{X, T, \delta}^{3/2} \right) \]
In particular, $\Vert u \Vert_{X, T, \delta} \lesssim \Vert u \Vert_{X_0}$ and we may choose $T = \infty$. 
\end{itemize}
\label{prop-estimeeapriori-complete}
\end{Prop}

The Cauchy theory ensures that a solution $u$ with data $u_0 \in H^1$ is defined globally in $C(\mathbb{R}, H^1)$. 

$\delta$ is chosen small enough with respect to $1$: for instance, we may fix $\delta = 2^{-100}$. 

In particular, if $\Vert u_0 \Vert_{X_0} \leq \varepsilon$ for some small enough $\varepsilon > 0$, then since $\Vert u \Vert_{X, 0, \delta} \leq C \Vert u_0 \Vert_{X_0}$, by continuity and local existence we may choose $T$ maximal such that $\Vert u \Vert_{X, T, \delta} \leq 2 C \Vert u_0 \Vert_{X_0}$. If this $T$ is not $+\infty$, then the a priori estimate ensure that 
\begin{align*}
\Vert u \Vert_{X, T, \delta} &\leq C \left( \varepsilon + \left( 2 C \varepsilon \right)^{\frac{3}{2}} \right) \\
&\leq C \varepsilon \left( 1 + \left( 2 C \right)^{\frac{3}{2}} \varepsilon^{\frac{1}{2}} \right) 
\end{align*}
$C > 0$ is fixed and universal, so choosing $\varepsilon > 0$ small enough, we get 
\begin{align*}
\Vert u \Vert_{X, T, \delta} &\leq \frac{3}{2} C \varepsilon
\end{align*}
which is a contradiction with the continuity of the norm with respect to $T$ and the choice of $T$ maximal. 

We deduce that, for small enough initial data, the global solution satisfies $\Vert u \Vert_{X, \infty, \delta} \leq 2 C \Vert u_0 \Vert_{X_0}$. This means that Proposition \ref{prop-estimeeapriori-complete} implies part of the result of Theorem \ref{theo-principal} already. To get scattering, we will prove the dispersive estimate 
\begin{align*}
\Vert \partial_x u(t) \Vert_{L^{\infty}} \lesssim t^{-\frac{5}{6}} \langle t \rangle^{-\frac{1}{4}+100\delta} \Vert u \Vert_{X, T, \delta} 
\end{align*}
so that, for small $\delta$, $\partial_x u \in L^1_t L^{\infty}_{x, y}$, and therefore 
\begin{align*}
\Vert \partial_t f(t) \Vert_{L^2} &= \Vert u(t) \partial_x u(t) \Vert_{L^2} \leq \Vert u(t) \Vert_{L^2} \Vert \partial_x u(t) \Vert_{L^{\infty}} \lesssim t^{-\frac{5}{6}} \langle t \rangle^{-\frac{1}{5}} \Vert u \Vert_{X, T, \delta} 
\end{align*}
and therefore $\partial_t f \in L^1_t L^2_{x, y}$, hence scattering in $L^2$. In fact, we will show in Proposition \ref{propositionestimeesstandardscattering} that $\partial_t f \in L^1_t H^1_{x, y}$ as well, hence scattering in $H^1$. 

In what follows, we drop the dependency in $\delta, T$ of the $X$-norm, since we may assume $\delta > 0$ to be fixed once and for all (small enough), while as shown before $T$ can be pushed to $+\infty$. Note that the universal constants will \emph{not} depend on $T$ in any case, but will certainly depend on $\delta$. 

\begin{Rem} Note that, $u$ being real-valued, so is $f$. Therefore, $|\widehat{f}(t, \overline{\xi})|^2 = \widehat{f}(t, \overline{\xi}) \widehat{f}(t, -\overline{\xi})$. 

Furthermore, from the choice of $X_{\alpha}$ and $m_{\alpha}$, the operators $|\nabla| X_a, |\nabla| m_c(D) X_c$ and $\partial_x m_b(D) X_b$ map real-valued functions to real-valued functions as well. 
\end{Rem}

We will use repeatedly the following consequences of the previous definition of the $X$-norm: 

\begin{Lem} The following norms are controlled by the $X$-norm: 
\begin{subequations}
\begin{align}
\Vert |\nabla|^{-1} f(t) \Vert_{L^2} &\lesssim \Vert u \Vert_X \label{lem-estimees-directes-Hardy} \\
\Vert m_{\alpha}(D) X_{\alpha} f(t) \Vert_{H^1} &\lesssim \Vert u \Vert_X \label{lem-estimees-directes-Xa} \\
\Vert m_{\widehat{\mathcal{R}}}(D) (x, y) f(t) \Vert_{H^1} &\lesssim \Vert u \Vert_X \label{lem-estimees-directes-XR} \\
\Vert m_{\widehat{\mathcal{P}}}(D) (x, y) f(t) \Vert_{H^1} &\lesssim \Vert u \Vert_X \label{lem-estimees-directes-XP}
\end{align} 
\end{subequations} \label{lem-estimees-directes-uX}
\end{Lem}

\begin{Dem}
\eqref{lem-estimees-directes-Hardy} follows from Hardy's inequality and Parseval's identity: 
\begin{align*}
\Vert |\nabla|^{-1} f(t) \Vert_{L^2} &= \Vert |\overline{\xi}|^{-1} \widehat{f}(t) \Vert_{L^2} \\
&\lesssim \Vert \partial_r \widehat{f}(t) \Vert_{L^2}
\end{align*}
where $\partial_r$ denotes the radial derivative. But precisely $\partial_r = \widehat{X}_a(\overline{\xi}) \cdot \nabla_{\overline{\xi}}$, hence 
\begin{align*}
\Vert |\nabla|^{-1} f(t) \Vert_{L^2} &\lesssim \Vert X_a f(t) \Vert_{L^2}
\end{align*}
and the bound by the $X$-norm will follow from the second inequality. 

For any $\alpha$, we may decompose by definition of $h_{\alpha}$ and $g_{\alpha}$: 
\begin{align*}
\Vert m_{\alpha}(D) X_{\alpha} f(t) \Vert_{H^1} &\leq \Vert h_{\alpha}(t) \Vert_{H^1} + \Vert g_{\alpha}(t) \Vert_{H^1} \lesssim \Vert u \Vert_X 
\end{align*}
since $g_{\alpha}(t)$ satisfies better bounds than $h_{\alpha}$. 

Finally, if we localise near $\widehat{\mathcal{P}}$ or $\widehat{\mathcal{R}}$, then $m_{\alpha}(\overline{\xi}) \simeq 1$ for every $\alpha$, so that 
\begin{align*}
\Vert m_{\widehat{\mathcal{P}}}(D) X_{\alpha} f(t) \Vert_{H^1} &\lesssim \Vert m_{\widehat{\mathcal{P}}}(D) m_{\alpha}(D) X_{\alpha} f(t) \Vert_{H^1} \\
&\lesssim \Vert m_{\alpha}(D) X_{\alpha} f(t) \Vert_{H^1} \lesssim \Vert u \Vert_X 
\end{align*}
and similarly near $\widehat{\mathcal{R}}$. But since $\left( \widehat{X}_{\alpha}(\overline{\xi}) \right)_{\alpha = a, b, c}$ is an orthonormal basis for every $\overline{\xi} \notin \widehat{\mathcal{L}}$, we deduce 
\begin{align*}
\Vert m_{\widehat{\mathcal{P}}}(D) (x, y) f(t) \Vert_{H^1} &\lesssim \sum_{\alpha = a, b, c} \Vert m_{\widehat{\mathcal{P}}}(D) X_{\alpha} f(t) \Vert_{L^2} \\
&\lesssim \Vert u \Vert_X 
\end{align*}
as desired. 
\end{Dem}

\subsection{Outline of the proof} \label{subsection-outline}

We now explain the main steps of the proof, and the differences and challenges with respect to existing literature. 

The strategy is based on the space-time resonances method \cite{Germain_3DSchrod, Germain_2DSchrod, GMSWaterwaves} which can be explained as follows. Recall that the profile function $f$ satisfies 
\begin{align}
\widehat{f}(t, \overline{\xi}) &= \widehat{f}(0, \overline{\xi}) + \int_0^t \int e^{i s \varphi(\overline{\xi}, \overline{\eta})} i \xi_0 \widehat{f}(s, \overline{\eta}) \widehat{f}(s, \overline{\sigma}) ~ d\overline{\eta} ds \label{outline-Duhamelform} 
\end{align}
up to an universal constant. 
In order to prove scattering, it is sufficient to prove that $f(t)$ converges as $t \to \infty$ as we saw, or in particular to bound the quadratic term on the right here above, uniformly in $t$. The space-time resonances method is based on the stationary phase heuristics and, using the integrals in both time and frequency, it states that the nonlinear dynamics should be driven by the frequency interactions such that the phase is stationary in both $s$ and $\overline{\eta}$, that is on 
\begin{align*}
\left\{ (\overline{\xi}, \overline{\eta}) \in \mathbb{R}^6, ~ \varphi(\overline{\xi}, \overline{\eta}) = 0, ~ \nabla_{\overline{\eta}} \varphi(\overline{\xi}, \overline{\eta}) = 0 \right\} 
\end{align*}
Points where $\varphi$ vanishes are called time-resonances, and points where $\nabla_{\overline{\eta}} \varphi$ vanishes are called space-resonances. Now the goal will be to identify some additional structure in the nonlinear term, using typically the cancellations of $\xi_0$, in order to compensate the interactions on this resonant set. 

In our case, since the linear $L^{\infty}$ dispersive decay will be better than $t^{-1}$ (and thus integrable), the unweighted bound is easy to get, but we will use the space-time resonances approach in order to bound weighted estimates, which are used to close the $L^{\infty}$ bound. The typical approach is to try to bound $\Vert (x, y)^2 f(t) \Vert_{L^2}$ so that, applying a standard dispersive estimate, one would have 
\begin{align*}
\Vert \partial_x u(t) \Vert_{L^{\infty}} &= \Vert \partial_x e^{it \omega(D)} f(t) \Vert_{L^{\infty}} \lesssim t^{-a} \Vert f(t) \Vert_{L^1} \\
&\lesssim t^{-a} \Vert \langle x, y \rangle^2 f(t) \Vert_{L^2} 
\end{align*}
for some $a > 1$. Since $(x, y)$ weights correspond to a derivation in Fourier, one can see that after applying $\nabla_{\overline{\xi}}$ to \eqref{outline-Duhamelform}, the dominant term corresponds to the case where the derivative hits the exponential and adds a factor $s$: 
\begin{align}
\int_0^t \int e^{i s \varphi} i s \xi_0 \nabla_{\overline{\xi}} \varphi ~ \widehat{f}(s, \overline{\eta}) \widehat{f}(s, \overline{\sigma}) ~ d\overline{\eta} ds \label{outlinetermedominants} 
\end{align}
Now, ideally, we would like to have a relation of the form 
\begin{align}
\xi_0 \nabla_{\overline{\xi}} \varphi &= \mu_1 \varphi + \mu_2 \nabla_{\overline{\eta}} \varphi \label{relationidealespacetimeres} 
\end{align}
in order to apply integrations by parts either in $\overline{\eta}$ to absorb the factor $s$, or in $s$ to obtain a cubic term (or a boundary in time term), which are all easier to control. This would be a null structure in the spirit of Klainerman \cite{Klainermannullform}. 

The main difficulty is that \eqref{relationidealespacetimeres} does not hold. In fact, in the case where $\overline{\eta}$ is close to $\widehat{\mathcal{C}}$ and $\overline{\xi}$ is close to $0$, both $\varphi$ and $\nabla_{\overline{\eta}} \varphi$ vanish at order at least $2$ (in a certain direction), while $\xi_0$ vanishes at order $1$ and $\nabla_{\overline{\xi}} \varphi$ does not vanish. There is therefore some cancellation in the symbol, but it is not enough to apply the space-time resonance method. In \cite{AnjolrasmZK}, the author studied the modified Zakhrov-Kuznetsov equation in dimension $2$ where a similar problem occurs, but it was possible to overcome this difficulty by making use of the linear dispersive decay alone. This is not possible here as the linear decay will typically be of order $t^{-\frac{7}{6}}$ at most, and only $t^{-\frac{13}{12}}$ near singular geometric areas (like the cone). Although the linear decay is weaker in dimension $2$, the fact that the equation is modified (that is, the nonlinear term is cubic instead of quadratic) allows for a stronger decay of the nonlinear term. If the dispersive decay is high enough as in the 2d case, it is possible to compensate this defect by separating frequencies near the resonant point and away from it: away from it, the usual space-time resonance method applies, and near it we can use the cancellation of the left-hand side of \eqref{relationidealespacetimeres} (or the equivalent) without further integrations by parts, and then optimize the (time dependant) threshold between the near and away area. Applying this method allows to compensate part of the additional $s$ factor in the dominant part \eqref{outlinetermedominants}, and the rest is absorbed using the linear dispersive decay. 

We thus follow a different strategy by observing that, in \eqref{relationidealespacetimeres}, not only are the point where the relation fails well localised, but also $\nabla_{\overline{\xi}} \varphi$ does not vanish, so it brings no additional cancellation. Instead of trying to get a bound on the whole gradient (so on an isotropic weight), we therefore try to bound only specific directions, hence the introduction above of the vector fields $\widehat{X}_{\alpha}$. We then have to rethink entirely relation \ref{relationidealespacetimeres}, as we cannot anymore keep $\nabla_{\overline{\eta}} \varphi$ as a whole to apply integrations by parts but want to consider only specific $\overline{\eta}$ directions of derivations which will be adapted to \emph{both} $\overline{\eta}$ and $\overline{\sigma}$. Also, the proof of the linear dispersive estimate has to be done from scratch as well, as we want to avoid ``weights'' (in the generalised sense of derivatives in Fourier space) in certain directions, which means that we will in fact prove no $L^1$ estimate on $f$, so we will redo the proof of the dispersive estimate by using fine stationary phase analysis. Another difficulty in the $L^{\infty}$ bound comes from the fact that we need to separate $m_{\alpha} X_{\alpha} f = h_{\alpha} + g_{\alpha}$, where $g_{\alpha}$ is the boundary term in time coming from the normal form transform (that is some integration by parts in time in \eqref{outlinetermedominants}), and which needs to be estimated differently: we cannot put another weight on it, since another integration by parts in time is not allowed, but since there is no time integral in $g$, it automatically has some time decay. 

As can be easily seen, coming back to the physical space, the vector field $X_a$ can be rewritten as
\begin{align*}
|\nabla| X_a = x \partial_x + y \cdot \nabla_y + 3 
\end{align*}
which (up to the factor $3$ coming from the commutation of $(x, y)$ and $\nabla$) is the (space) scaling vector field of the equation. Indeed, from scaling invariance we have that the following (physical) vector field commutes with the equation: 
\begin{align}
S = x \partial_x + y \cdot \nabla_y + 3 t \partial_t \label{equscalingvfZK3D}
\end{align}
The absence of $3 t \partial_t$ in $|\nabla| X_a$ means that an integration by parts in time will be used in the estimate of \eqref{outlinetermedominants} in this case. Likewise, 
\begin{align*}
|\nabla| m_c(D) X_c = Jy \cdot \nabla_y 
\end{align*}
which is exactly the vector field associated with the rotation invariance around the $x$-axis of the equation. On the other hand, $X_b$ is not associated to any invariance. 

This connects this strategy to the method of partial symmetries of Guo, Pausader and Widmayer \cite{GPWpartialsymmEuler}, where a similar approach was used in order to treat the anisotropy of the Euler-Coriolis equation. Let us however note the following novelties: 
\begin{itemize}
\item The number of singular geometric areas ($3$) and their geometries is new in the literature, as in previous works one would not have the cone $\widehat{\mathcal{C}}$. This leads to many case exhaustions depending on the interactions. We need to separate the interactions depending on their angular localisations, using the symbols $m_A$, $A \in \{ \widehat{\mathcal{C}}, \widehat{\mathcal{L}}, \widehat{\mathcal{P}}, \widehat{\mathcal{R}} \}$, which localise angularly close to each of the three geometric areas, respectively away enough from them, but also according to the relative size of the frequencies, since we want to be able to control where the (physical) derivatives fall. 
\item From the viewpoint of harmonic analysis, the fact that the cone $\widehat{\mathcal{C}}$ is not a vector subspace also complicates the situation. At the linear level, it is not even clear that the localisation symbol 
\begin{align*}
\psi_{j, k}^{\widehat{\mathcal{C}}}(\overline{\xi}) = \psi\left( 2^{-j} \overline{\xi} \right) \psi\left( 2^{-j-k} (\sqrt{3} |\xi_0| - |\xi|) \right) 
\end{align*}
where $\psi$ is the standard Littlewood-Paley localisation function on an annulus, satisfies good $L^p \to L^p$ estimates (except for $p = 2$), while it is always the case for 
\begin{align*}
\psi_{j, k}^{\widehat{\mathcal{P}}}(\overline{\xi}) = \psi\left( 2^{-j} \overline{\xi} \right) \psi\left( 2^{-j-k} \xi_0 \right) 
\end{align*}
or for the analogous localisation symbol near $\widehat{\mathcal{L}}$, simply by applying the Hörmander-Mikhlin theorem. 

The situation is even more dramatic for multilinear estimate: if we consider the interaction of a frequency $\overline{\eta}$ with $\overline{\sigma} = \overline{\xi}-\overline{\eta}$ near the plane, we may want to localise according to a symbol having the form 
\begin{align*}
\mu_1(\overline{\xi}, \overline{\eta}) \mu_2(\xi_0, \eta_0) 
\end{align*}
where $\mu_1, \mu_2$, for instance, are homogeneous of order $0$. This allows in particular to apply a first localisation on the relative sizes of frequencies, and a second on the relative sizes of their $0$-components only. The Coifman-Meyer theorem applies only for $\mu_1$, but extensions due to Muscalu, Pipher, Tao and Thiele \cite{MPTTBiparameter2004, MPTT2006} allow to generalise it. By analogy, we could try to localise as well the distance to the cone of $\overline{\eta}, \overline{\sigma}$ in a finer way than the global size of frequencies: but there exist no result in this direction, and there are even good reasons to believe such results cannot hold, since the distance of the cone cannot be treated as another coordinate, because it does not satisfy the properties allowing to apply a Bony decomposition. 

This difficulty is overcome by considering directly $\Vert \psi_{j, k}^{\widehat{\mathcal{C}}}(D) u(t) \Vert_{L^{\infty}}$ in the dispersive estimate, and trying to prove an estimate summable in $j, k$, which allows to get good bounds even after the application of linear symbols and, in a certain extend, multi-linear symbols as well. 
\item As seen above, the scaling vector field $S$ contains time derivatives, which is related to the fact that the linear operator of Zakharov-Kuznetsov is of order $3$. In the case of Euler-Coriolis or Boussinesq, the linear operator is of order $0$: in particular, one can link directly $|\nabla| X_a$ to $S$, and no integration by parts in time (ie no normal form) is necessary in \eqref{outlinetermedominants}. This is not the case for us, and it is a potential problem, as after one integration by parts in time, one has to deal with boundary terms in time, without time integral, on which we cannot apply the same strategy as new integrations by parts in time are not possible. On the other hand, since such terms are not integrated in time, they automatically decay in time. 

In such a case, in the usual literature using space-time resonances, one tries to apply the normal form directly on $f(t)$, remove the boundary term in time and then treat the rest the same way for both linear and quadratic weights. But this is not possible here, as we use crucially the null structure not only of $\xi_0$, but of $\xi_0 \widehat{X}_{\alpha} \cdot \nabla_{\overline{\xi}} \varphi$ in order to be able to perform the normal form. This leads to the decomposition $X_{\alpha} f(t) = h_{\alpha}(t) + g_{\alpha}(t)$, and we can only apply a second weight on $h_{\alpha}(t)$. In comparison, in previous works where $\omega$ is of order $0$, $g_{\alpha} \equiv 0$ so one can apply recursively $N$ weights without difficulty. 

Now, for the $L^{\infty}$ dispersive estimate, in order to apply the stationary phase method and obtain a better decay than $t^{-1}$, one needs to apply at least two integrations by parts in frequency, and thus consider two weights. We therefore need to use the following strategy: apply a first integration by parts, separate into $h_{\alpha}+g_{\alpha}$, apply a second integration by parts only on the $h$-term, and treat the other one using its time decay. This means that the linear $L^{\infty}$ estimate depends in fact on the nonlinear structure. 
\item Connected to the previous point, we manage to keep the number of weights and (physical) derivatives low, namely $2$ and $1$, which is likely to be optimal considering only integer exponents. In previous works, the number of weights was typically $N \gg 1$, which is also linked to the fact that there is no difficulty propagating them recursively. Since the scaling vector field $S$ adds one (physical) derivative with each weight, the number of derivatives is therefore high as well. 

Note that it is possible in our case as well to propagate $a$ and $c$-weights of order $N$ up to removing the boundary terms in time at each step. 

Also, we avoid the use of Besov-type estimates in most of the argument. 
\item When estimating the quadratic weighted terms, we try to recover terms of the form 
\begin{align*}
e^{i t \omega(D)} \left( \left( e^{-i t \omega(D)} \nabla m_{\alpha}(D) X_{\alpha} m_{\beta}(D) X_{\beta} f(t) \right) \left( e^{-i t \omega(D)} \partial_x f(t) \right) \right) 
\end{align*}
up to a multilinear symbol, and we can apply a $L^2 \times L^{\infty}$ estimate to conclude. However, we cannot avoid to have terms of the form 
\begin{align*}
e^{i t \omega(D)} \left( \left( e^{-i t \omega(D)} \nabla m_{\alpha}(D) X_{\alpha} f(t) \right) \left( e^{-i t \omega(D)} \nabla m_{\beta}(D) X_{\beta} f(t) \right) \right) 
\end{align*}
in general. 

The usual strategy, also applied for the wave equation above, is to obtain a bound in $L^4$ on each of the weighted interacting terms, with a decay obtained by Riesz-Thorin interpolation theorem between the dispersive decay $L^1 \to L^{\infty}$ and the unitary nature $L^2 \to L^2$. But, as we explained before, we cannot apply the standard $L^1 \to L^{\infty}$ in our case and we must therefore also redo this $L^4$ estimate that does not follow directly from the $L^{\infty}$ dispersive estimate. 

Worse: we do not reach (here at least) a sufficient decay on weighted estimates in $L^4$. It is however possible to exploit again the space-time resonances by trying to get a bound on these weighted terms in $L^4_{x, y}$ but also integrated in time. The advantage of this new integral is that we can apply again an integration by parts in time, which allows to reach the desired decay in $L^4$, up to the time integration. 
\item Surprisingly, there are cubic resonances. The choice of $m_b, m_c$ allows, so to speak, to absorb all the quadratic resonances, but it is not the case at the cubic level. But cubic interaction terms appear as soon as we apply a normal form, so after the first weight, and need to be treated differently when going to the quadratic weights. 

These resonances are only partially linked to the geometric areas identified at the linear level: they appear in the interaction of three frequencies where one only is situated on a singular geometric area (typically the cone $\widehat{\mathcal{C}}$) while the other two are away enough from $\widehat{\mathcal{C}} \cup \widehat{\mathcal{L}} \cup \widehat{\mathcal{P}}$, as well as the output frequency. More precisely, these frequencies are localised along a specific cone that does not appear at the linear level. 

In particular, some quadratic weighted estimates are degenerate and add a time growth. It should be possible to compensate them by adding a new cancellation in $m_b$ (since only the $b$ weight leads to these resonances); however, we can also deal with it using the following observation. In a simple case where we would use the standard dispersive estimate and Hölder's inequality, we could estimate: 
\begin{align*}
\Vert f(t) \Vert_{L^1(\mathbb{R}^3)} \lesssim \Vert \langle x \rangle^{\frac{1}{2}+} \langle y \rangle^{1+} f(t) \Vert_{L^2(\mathbb{R}^3)} 
\end{align*}
by differentiating the two directions $x, y$, which is strictly better than a bound by the isotropic weight $\Vert \langle x, y \rangle^{\frac{3}{2}+} f(t) \Vert_{L^2}$ if the two directions are not equivalent. Therefore, here, we can arrange so that, in the $L^{\infty}$ dispersive estimate, we alway avoid the double $b$ weight, and only consider crossed weights $a-b$ or $c-b$. 
\item A crucial point is that the weights $a, c$, even after applying a normal form, make no cubic term appear (which would have resonances). There is indeed a ``magical'' cancellation that allows to show that, for these two weights, the different cubic terms obtained by integrations by parts cancel each other. 

This cancellation is weaker for the $b$ weight, but also needs to be used in order to compensate certain critical interactions at the cubic level, namely when three frequencies interact, one being very large with respect to the other two, allowing so to speak to distribute the derivatives between the different interacting factors and avoid that they all fall on the same. 
\end{itemize}

We now explain each step of the proof more precisely. 

In Section \ref{section-harmonicanalysis}, we recall useful harmonic analysis theory, both in the linear and in the multilinear setting, which we will use repeatedly during the argument. In Section \ref{section-spacetimeresanal}, we compute the space-time resonances in both the quadratic and the cubic case: even if the equation is quadratic, note that, when applying an integration by parts in time in \eqref{outlinetermedominants} or similar terms, one can obtain a cubic term when developing $\partial_s \widehat{f}$, so that a cubic analysis is necessary when applying another weight after. We also introduce conical coordinates which are suitable for later computations, and provide some new vector fields, typically depending on $\overline{\eta}$ and $\overline{\sigma}$ and differentiating in good directions for both of them. 

In Section \ref{section-Linfdispersive}, we prove a linear $L^{\infty}$ decay estimate. This is done by separating the frequency space in several areas: low frequencies, where a KdV-type singularity occurs, which can be easily bounded using Sobolev's embedding; the remainder area $\widehat{\mathcal{R}}$ and the three singular geometric areas. We also treat short times differently in order to reach integrability up to $t = 0^{+}$. 

For the $L^{\infty}$ dispersive estimate, we write 
\begin{align*}
u(t, x, y) &= e^{-i t \omega(D)} f(t, x, y) = \int e^{i t \Phi(\overline{\xi}, x/t, y/t)} \widehat{f}(t, \overline{\xi}) ~ d\overline{\xi}
\end{align*}
where the phase is
\begin{align*}
\Phi\left( \overline{\xi}, \frac{x}{t}, \frac{y}{t} \right) &= \xi_0^3 + \xi_0 |\xi|^2 - \frac{\xi_0 x}{t} - \frac{\xi \cdot y}{t} 
\end{align*}
The low frequency estimate is easily obtained by Sobolev's embedding. On the remainder area, on the support of $m_{\widehat{\mathcal{R}}}$, $m_b \simeq m_c \simeq 1$ uniformly and $\nabla^2 \omega$ is never degenerate, which allows to apply estimates in a relatively easy way and to reach the decay 
\begin{align*}
\Vert m_{\widehat{\mathcal{R}}}(D) \nabla u(t) \Vert_{L^{\infty}} &\lesssim t^{-\frac{7}{6}+101\delta} \Vert u \Vert_X \\
\Vert m_{\widehat{\mathcal{R}}}(D) |\nabla|^{\frac{3}{2}} u(t) \Vert_{L^{\infty}} &\lesssim t^{-\frac{7}{6}+100\delta} \Vert u \Vert_X 
\end{align*}
The main difference with respect to the usual stationary phase estimate is that we want to avoid a double $b$ weight. 

Degeneracies appear however on other areas. Near the cone $\widehat{\mathcal{C}}$, $m_c \simeq 1$ but $m_b$ vanishes, and degeneracies may occur if $(x, y)$ is localised near the cone $\mathcal{C}$, or near the line $\mathcal{L}$ (other cases being similar to the one already treated when $\overline{\xi}$ is localised by $m_{\widehat{\mathcal{R}}}$). If $(x, y)$ is localised near $\mathcal{C}$, then $\widehat{X}_a \cdot \nabla_{\overline{\xi}} \Phi$ and $\widehat{X}_c \cdot \nabla_{\overline{\xi}} \Phi$ can vanish at order $1$, but $\widehat{X}_b \cdot \nabla_{\overline{\xi}} \Phi$ can vanish at order $2$ (this is the degenarcy of $\nabla^2 \omega$ on the cone), which is aggravated by the presence of the symbol $m_b$ which itself vanish linearly on the cone. It is therefore necessary to apply a frequency decomposition depending on the size of each of these derivatives of the phase and estimate the corresponding volumes and regularity. On the other hand, if $(x, y)$ is localised near $\mathcal{L}$, then since 
\begin{align*}
\widehat{X}_c(\overline{\xi}) \cdot \nabla_{\overline{\xi}} \Phi &= - \frac{J \xi \cdot y}{t |\xi|} 
\end{align*}
this direction can become very degenerate if $y = 0$. On the other hand, $\widehat{X}_a \cdot \nabla_{\overline{\xi}} \Phi$ can only vanish at order $1$ and, in this case, $\widehat{X}_b \cdot \nabla_{\overline{\xi}} \Phi$ does not vanish, but $m_b$ adds a cancellation that, in some extend, makes the phase stationary, and we must again proceed carefully. We finally obtain a total decay of the form 
\begin{align*}
\Vert m_{\widehat{\mathcal{C}}}(D) \nabla u(t) \Vert_{L^{\infty}} &\lesssim t^{-\frac{13}{12}+100\delta} \Vert u \Vert_X \\
\Vert m_{\widehat{\mathcal{C}}}(D) |\nabla|^{\frac{3}{2}} u(t) \Vert_{L^{\infty}} &\lesssim t^{-\frac{7}{6}+100\delta} \Vert u \Vert_X
\end{align*}

Near $\widehat{\mathcal{L}}$, degeneracies only happen due to the symbols $m_b, m_c$ (because $\nabla^2 \omega$ is not degenerate here) and only if $(x, y)$ is near $\mathcal{L}$. In particular, we have the same $y = 0$ degeneracy as on the cone, and even if $\widehat{X}_a \cdot \nabla_{\overline{\xi}} \Phi$ can only vanish at order $1$, $m_b \widehat{X}_b \cdot \nabla_{\overline{\xi}} \Phi$ vanish at order $2$. Again, we reach a decay of the form 
\begin{align*}
\Vert m_{\widehat{\mathcal{L}}}(D) \nabla u(t) \Vert_{L^{\infty}} &\lesssim t^{-\frac{13}{12}+100\delta} \Vert u \Vert_X \\
\Vert m_{\widehat{\mathcal{L}}}(D) |\nabla|^{\frac{3}{2}} u(t) \Vert_{L^{\infty}} &\lesssim t^{-\frac{7}{6}+100\delta} \Vert u \Vert_X 
\end{align*}

Finally, near the plane $\widehat{\mathcal{P}}$, $m_b \simeq m_c \simeq 1$, but $\nabla^2 \omega$ is heavily degenerate: $\nabla \omega$ is constant on any circle $\xi_0 = 0$, $|\xi| = cte$. In order to compensate this degeneracy, we use the presence of the $\partial_x$ derivative in the nonlinearity, and show a good decay if we have at least some $\partial_x$ derivative in the dispersive estimate. We then obtain 
\begin{align*}
\Vert m_{\widehat{\mathcal{P}}}(D) |\partial_x|^{\frac{1}{3}} |\nabla|^{\frac{2}{3}} u(t) \Vert_{L^{\infty}} &\lesssim t^{-\frac{7}{6}+101\delta} \Vert u \Vert_X \\
\Vert m_{\widehat{\mathcal{P}}}(D) |\partial_x|^{\frac{1}{3}} |\nabla|^{\frac{7}{6}} u(t) \Vert_{L^{\infty}} &\lesssim t^{-\frac{7}{6}+100\delta} \Vert u \Vert_X \\
\Vert m_{\widehat{\mathcal{P}}}(D) |\nabla|^{\frac{1}{2}} \langle \nabla \rangle u(t) \Vert_{L^{\infty}} &\lesssim t^{-1+100\delta} \Vert u \Vert_X
\end{align*}

The $L^{\infty}$ estimates are delicate notably because we try to get an optimal decay, this optimal decay being important for the optimal growth of degenerate weighted terms $h_{\alpha *}(t)$, and these degenerate weighted norms being themselves useful in the $L^{\infty}$ estimates in certain regimes. These estimates are performed localising the size of the frequency as some $2^j$, and the distance to one of the singular geometric areas as some $2^{j+k}$, $k \leq -10$, in a summable way in $(j, k)$. This means we can automatically absorb Hörmaner-Mikhlin (and even more degenerate) symbols in the $L^{\infty}$ estimate. These bounds are also useful for some specific delicate nonlinear interactions. 

In Section \ref{section-L4disp}, we provide $L^4$ bounds for the $h_{\alpha}(t)$. Indeed, for the quadratic weight, we will typically need to bound products of the form 
\begin{align*}
\mathcal{N}\left( e^{-it \omega(D)} |\nabla| m_{\alpha}(D) X_{\alpha} f(t), ~ e^{-i t \omega(D)} |\nabla| m_{\beta}(D) X_{\beta} f(t) \right) 
\end{align*}
in $L^2$ for some $\alpha, \beta$, because we could never force two weights to fall on the same factor, and thus $L^4$ bounds are necessary. Unlike in the usual space-time resonance method, such $L^4$ estimates do not follow immediately from interpolation between the unitary $L^2$ estimate and the $L^{\infty}$ bound. 

The weighted $L^4$ estimates can be obtained in a simpler way than the $L^{\infty}$ ones because we only need an asymptotic decay better than $t^{-\frac{1}{2}}$ to conclude, and not the optimal decay, these bounds not being used in the singular estimates. Moreover, as long as we are localised far enough from $\widehat{\mathcal{C}}$, we can simplify the argument by estimating simply the kernel $K_t$, $\widehat{K_t} = e^{-i t \omega(\overline{\xi})}$, in some anisotropic Lebesgue spaces, noting that, near $\widehat{\mathcal{L}}$ for instance, the $x$ weight corresponds to the vector field $\partial_{\xi_0}$, which is a linear combination with smooth symbols coefficients of the $m_{\alpha} \widehat{X}_{\alpha} \cdot \nabla_{\overline{\xi}}$, and similarly near $\widehat{\mathcal{P}}$. 

We must however use again a stationary phase estimate and decomposition as in the $L^{\infty}$ case when we are near the cone, introducing also a pseudo-differential localisation lemma well-suited for our purpose: before, it was possible to fix $(x, y)$ since we wanted to compute a supremum in $(x, y)$, so that all the symbols localising the derivatives of the phase $\Phi$, depending on both $\overline{\xi}$ and $x, y$, could be considered as depending only on $\overline{\xi}$ in the estimates. In $L^4$ norm, this simplification is no longer possible. 

Moreover, the pointwise in time estimate obtained this way does not reach in general a better decay than $t^{-\frac{1}{2}}$. We can however improve them by integrating in time: 
\begin{align*}
\left( \int_0^t s^{1-q} \langle s \rangle^r \Vert e^{-i t \omega(D)} \nabla h_{\alpha}(s) \Vert_{L^4}^4 ~ ds \right)^{\frac{1}{4}} &\lesssim \Vert u \Vert_X 
\end{align*}
for $r, q > 0$ some small (explicit) parameters such that $r > q$. Since the integral is bounded time, the goal is to obtain a maximal exponent on the $L^4$ norm: the exponent $4$ is a priori optimal since it allows to rewrite as 
\begin{align*}
\left( \int_0^t s^{1-q} \langle s \rangle^r \int \left| e^{-i t \omega(D)} \nabla h_{\alpha}(s, x, y) \right|^4 ~ dx dy ds \right)^{\frac{1}{4}}
\end{align*}
and then switch to the Fourier side in order to apply space-time resonances techniques. 

Section \ref{section-decomplemma} focuses on the key decomposition lemma of this paper. The idea is to write an interaction as
\begin{align*}
I[F_1, F_2](t) &= \int_0^t e^{is \omega(D)} \left[ \left( e^{-i s \omega(D)} F_1(s) \right) \left( e^{-i s \omega(D)} F_2(s) \right) \right] ~ ds 
\end{align*}
or possibly adding some (linear or multilinear) symbols, and then trying to decompose 
\begin{align*}
|\nabla| m_{\alpha}(D) X_{\alpha} I[F_1, F_2](t), \quad \partial_x m_b(D) X_b I[F_1, F_2](t) 
\end{align*}
for $\alpha = a, c$ as ``good interaction terms'' for which nonlinear estimates can be proven in subsequent sections. We prove the lemma in this abstract form in order to be able to apply it twice when considering the quadratic weights. For $\alpha = a, c$, the lemma is very easy to prove since it relies on the invariance of the equation. However, for the $b$-weight, we need to decompose the interaction $I[F_1, F_2](t)$ depending on the Fourier localisations of $F_1, F_2$. Their directions in Fourier, depending on which geometric area they are close to, is of particular importance, but also their relative sizes as we want to be able to distribute the derivatives in a suitable way. 

In Section \ref{section-structurepoids}, we apply the key decomposition lemma in order to derive expressions for the $h_{\alpha}, g_{\alpha}$ and $m_{\beta} X_{\beta} h_{\alpha}$, as by definition 
\begin{align*}
f(t) = f(0) + \partial_x I[f, f](t)
\end{align*}
An important step at this point is to prove some cancellation of the cubic terms coming from integrations by parts in time. It turns out that $h_a$ and $h_c$ are purely quadratic and contain no cubic term, while $h_b$ contains a cubic term, which still enjoys a cancellation for the High-Low-Low interactions. Also, when applying a weight $a$ or $c$ on a cubic term (coming from $h_b$), the invariance of the equation allows to obtain a nice decomposition as well. 

In Section \ref{section-estimees-simples-fug}, we prove useful linear bounds on $u$, and the a priori estimates for $f, g_{\alpha}$. In Section \ref{section-estimeepoidssimple}, we prove the a priori estimates for $h_{\alpha}$. These estimates are now pretty straightforward as the structure and dispersive decay has already been obtained. Likewise, in Section \ref{section-estimees-quadratiques-horsbb}, we prove the a priori estimate for quadratic weights, except the double $b$-weight. 

Section \ref{section-estimee-quad-bb} deals with the quadratic $b$-weight estimate. While adding a $b$ weight on the quadratic part of $h_b$ can be tackled applying again the key decomposition lemma, the most delicate part comes from the cubic term. In fact, there are some new resonances arising when we apply a $b$ weight in this cubic interaction, which explain why $\Vert \partial_x m_b(D) X_b h_b(t) \Vert_{L^2}$ has some time growth. Even when there is no time growth, there can still be delicate interactions to bound. Again, we separate cubic interactions depending on the directions of the interacting frequencies in a similar fashion as what was done for the proof of the decomposition lemma. 

Finally, in Section \ref{section-estimee-mauvaise-b}, we prove estimates for simple weights $X_b, X_c$ without adding the symbols $m_b, m_c$, and thus getting some time growth, these bounds being important for the $L^{\infty}$ dispersive estimate. The strategy is again to decompose the interactions depending on the directions and sizes of the interacting frequencies, but are somewhat simpler since we now have access to more directions to apply integrations by parts in $\overline{\eta}$. It turns out that the time growth is driven by very specific types of interactions, for which we follow again a stationary phase approach to obtain optimal time growth. This closes the a priori estimate. 

Section \ref{section-continuite} explains how to show the continuity of the introduced weighted norms $\Vert \cdot \Vert_{X, T}$ with respect to $T$. This follows from the already done estimates, which are enough to show that boundedness of this norm implies its continuity on the same time interval. Then, we may approximate the initial solution in a high regularity space, where local continuity (and hence boundedness) is automatic, then pass to the limit to obtain boundedness at the limit, and hence extend the solution for all times. 

\section{Linear and multilinear harmonic analysis} \label{section-harmonicanalysis} 

\begin{Theo}[Hörmander-Mikhlin] Let $d \in \mathbb{N}^{*}$. Let $m : \mathbb{R}^d \to \mathbb{C}$ be such that 
\[ \Vert m \Vert_{HM} := \sup_{\substack{\alpha \in \mathbb{N}^d, \\ |\alpha| \leq N}} \sup_{\xi \in \mathbb{R}^d} \left[ \left| \partial^{\alpha} m(\xi) \right| |\xi|^{|\alpha|} \right] ~ < \infty \]
where $N$ is a large enough integer. 

Then the linear operator $f \mapsto m(D) f$ is continuous from $L^p(\mathbb{R}^d)$ to $L^p(\mathbb{R}^d)$ for every $1 < p < \infty$, and its operator norm is bounded by
\[ \Vert m(D) \Vert_{L^p \to L^p} \leq C(d, p) \Vert m \Vert_{HM} \]
for some universal constant $C(d, p)$. \label{theoremHormanderMikhlin}
\end{Theo}

\begin{Cor} Let $1 \leq n \leq 3$, $m : \mathbb{R}^n \to \mathbb{C}$ with
\[ \Vert m \Vert_{HM}' := \sup_{\substack{\alpha \in \mathbb{N}^n, \\ |\alpha| \leq N}} \sup_{\xi \in \mathbb{R}^n} \left[ \left| \partial^{\alpha} m(\xi) \right| |\xi|^{|\alpha|} \right] ~ < \infty \]
Let us extend $m$ into a function $\mathbb{R}^3 \to \mathbb{C}$ invariant in the $3-n$ last coordinates. Then $f \mapsto m(D) f$ is continuous as in Theorem \ref{theoremHormanderMikhlin}. \label{corolHormanderMikhlin} 
\end{Cor}

\begin{Theo}[Coifman-Meyer \cite{CoifmanMeyer}] Let $m : \mathbb{R}^3 \times \mathbb{R}^3 \to \mathbb{C}$ and define $\Vert m \Vert_{CM}$ as the $\Vert \cdot \Vert_{HM}$-norm of $m$ seen as a Hörmander-Mikhlin multiplier on functions $\mathbb{R}^6 \to \mathbb{C}$. Define the bilinear operator 
\begin{align*}
T_m : (f, g) \mapsto \mathcal{F}_{\overline{\xi}}^{-1} \left( \int_{\mathbb{R}^3} m(\overline{\xi}, \overline{\eta}) \widehat{f}(\overline{\eta}) \widehat{g}(\overline{\xi}-\overline{\eta}) ~ d\overline{\eta} \right) 
\end{align*}
Assume $\Vert m \Vert_{CM} < \infty$. Then $T_m$ is continuous from $L^p \times L^q$ to $L^r$ for every $(p, q, r) \in (1, \infty]^2 \times [1, \infty)$ satisfying $\frac{1}{p}+\frac{1}{q}=\frac{1}{r}$. Furthermore, its operator norm is bounded by 
\[ \Vert T_m \Vert_{L^p \times L^q \to L^r} \leq C(p, q, r) \Vert m \Vert_{CM} \]
for some universal constant $C(p, q, r)$. 

The same holds in the multilinear case. 
\end{Theo} 

For the theorems above, we refer to \cite{MuscaluHAbook1, MuscaluHAbook2}. 

Near the plane, we will also need the following extension: 

\begin{Theo}[Muscalu-Pipher-Tao-Thiele \cite{MPTT2006}] Let $m : \mathbb{R}^6 \to \mathbb{C}$ be a bounded function, such that 
\[ \Vert m \Vert_{MPTT} := \sup_{\substack{\alpha \in \mathbb{N}^2, \\ \beta \in \left( \mathbb{N}^2 \right)^2, \\ |\alpha| + |\beta| \leq N}} ~ \sup_{\substack{(\overline{\xi}, \overline{\eta}) \in \mathbb{R}^3 \times \mathbb{R}^3}} |(\xi_0, \eta_0)|^{|\alpha|} |(\xi, \eta)|^{|\beta|} |\partial_{\xi_0}^{\alpha_1} \nabla_{\xi}^{\beta_1} \partial_{\eta_0}^{\alpha_2} \nabla_{\eta}^{\beta_2} m(\xi, \eta)| ~~ < \infty \]
where $N \in \mathbb{N}$ is a large enough (universal) number. Then the bilinear operator 
\[ T_m : ~ (f, g) \in \mathcal{S}(\mathbb{R}^3)^2 ~ \mapsto ~ \mathcal{F}^{-1} \left[ \int m(\overline{\xi}, \overline{\eta}) \widehat{f}(\overline{\eta}) \widehat{g}(\overline{\sigma}) ~ d\overline{\eta} \right] \in \mathcal{S}'(\mathbb{R}^2) \]
extends to a continuous bilinear operator $L^p(\mathbb{R}^3) \times L^q(\mathbb{R}^3) \to L^r(\mathbb{R}^3)$ for any $1 < p, q \leq \infty$, $0 < r < \infty$ such that $\frac{1}{r} = \frac{1}{p} + \frac{1}{q}$, with 
\[ \Vert T_m \Vert_{L^p \times L^q \to L^r} \leq C(p, q, r) \Vert m \Vert_{MPTT} \]
for some universal constant $C(p, q, r)$ depending only on $(p, q, r)$. 

This extends to the multilinear case. 
\end{Theo}

\begin{Rem} \cite{MPTT2006} is more general than the case of space dimension $3$. \end{Rem}

\begin{Theo}[Muscalu-Tao-Thiele \cite{MTT2002}] Let $m : \left( \mathbb{R}^3 \right)^n \to \mathbb{C}$ with $n \geq 3$ be a symbol $m = m(\overline{\eta}^1, ..., \overline{\eta}^n)$ and $\Gamma \subset \mathbb{R}^{3n}$ is a $3$-dimensional non-degenerate vector subspace, in the sense that $\Gamma$ is a graph in $\overline{\eta}^i$ for every $i$, and in $\overline{\xi} := \overline{\eta}^1 + ... + \overline{\eta}^n$. Let us assume that $m$ satisfies 
\begin{align*}
\Vert m \Vert_{MTT} := \sup_{\alpha \in \mathbb{N}^{3n}} \sup_{\overline{\eta}^1, ..., \overline{\eta}^n \in \mathbb{R}^3} \left| \partial^{\alpha} m(\overline{\eta}^1, ..., \overline{\eta}^n) \right| \mbox{dist}\left( (\overline{\eta}^1, ..., \overline{\eta}^n), \Gamma \right)^{|\alpha|} ~~ < \infty
\end{align*}
where $N$ is a large enough integer. 

Then the $n$-linear operator 
\begin{align*}
T_m : (f_1, ..., f_n) \mapsto \mathcal{F}_{\overline{\xi}}^{-1} \left( \int_{\mathbb{R}^{3n}} \delta(\overline{\xi} - \overline{\eta}^1 - ... - \overline{\eta}^n) m(\overline{\eta}^1, ..., \overline{\eta}^n) \prod_{j = 1}^n \widehat{f}_j(\overline{\eta}^j) ~ d\overline{\eta}^1 ... d\overline{\eta}^n \right) 
\end{align*}
is continuous from $L^{p_1} \times ... \times L^{p_n}$ to $L^q$ for every $(p_1, ..., p_n, q) \in (1, \infty]^n \times [1, \infty)$ satis-fying $\frac{1}{p_1}+...+\frac{1}{p_n} = \frac{1}{q}$. Furthermore, its operator norm is bounded by 
\begin{align*}
\Vert T_m \Vert_{L^{p_1} \times ... \times L^p_n \to L^q} \leq C(p_1, ..., p_n, q) \Vert m \Vert_{MTT}
\end{align*}
for some universal constant $C(p, q, r)$. \label{theoMTTmultilinsingGamma}
\end{Theo}

\begin{Rem} Theorem \ref{theoMTTmultilinsingGamma} was proven in \cite{MTT2002} in the one-dimensional setting but with a more general geometry (namely, any non-degenerate $\Gamma$ of dimension $k < \frac{n}{2}$ was allowed). The proof in the multi-dimensional case can be easily obtained from \cite{MTT2002} in our case, and a more general extension was proven in \cite{DPT2010} by , Demeter, Pramanik and Thiele, also allowing for the dimension of $\Gamma$ not to be a multiple of the space dimension. We also refer to \cite{BeneaMuscaluHelicoidal} and later works of the same authors for more advanced results regarding singular operators (for instance, combining a standard Coifman-Meyer operator in some space directions and a singular operator along a subspace in others). 
\end{Rem}

Note that we may combine several bilinear (or multilinear) operators with linear operators to get symbols not satisfying $\Vert m \Vert_{CM} < \infty$, but that can be written as a linear combination of products of bounded operators. For instance, a symbol of the form 
\begin{align*}
m(\overline{\xi}) \mu_1(\overline{\xi}, \overline{\eta}) + m(\overline{\eta}) \mu_2(\overline{\xi}, \overline{\eta})
\end{align*}
for $m$ of type Hörmander-Mikhlin, $\mu_1, \mu_2$ of type Coifman-Meyer, satisfies any Hölder bound $L^p \times L^q \to L^r$ for $p, q, r \in (1, \infty)$ such that $\frac{1}{p} + \frac{1}{q} = \frac{1}{r}$, but can typically not be expressed a a pure product of Hörmander-Mikhlin and Coifman-Meyer symbols. 

By considering now a series of such symbols, each expressed as a pure product, one can easily see that this class is not closed. In particular, we now present a framework to work with general symbols obtained as infinite series of pure products, namely by considering them under the form $f \circ m$ for $f$ having a power series development and $m$ a pure product symbol. 

\begin{Lem} Let $\Vert \cdot \Vert$ be one of the norms $\Vert \cdot \Vert_{HM}, \Vert \cdot \Vert_{HM}', \Vert \cdot \Vert_{CM}, \Vert \cdot \Vert_{MPTT}, \Vert \cdot \Vert_{MTT}$, and $m$ a corresponding symbol, and let $N$ be the integer appearing in the definition of the norm. Let $k \in \mathbb{N}^{*}$. Then
\[ \Vert m^k \Vert \leq \Vert m \Vert^{\min(k, N)} k^N \Vert m \Vert_{L^{\infty}}^{\max(0, k-N)} \]
\label{lemsymboles-puissancek}
\end{Lem}

\begin{Dem}
It is a direct consequence of the Leibniz rule. 
\end{Dem}

\begin{Cor} Let $m$ be a bilinear symbol that can be written as a product of Hörman-der-Mikhlin, Coifman-Meyer, Muscalu-Pipher-Tao-Thiele or Muscalu-Tao-Thiele. Let $f : \mathbb{C} \to \mathbb{C}$ be a function having a power series expansion near $0$, with radius of convergence $R > 0$. Assume $\Vert m \Vert_{L^{\infty}} < R$. Then $f \circ m$ corresponds to a bilinear operator continuous from $L^p \times L^q$ to $L^r$ for any $1 < p, q, r < \infty$ such that $\frac{1}{p} + \frac{1}{q} = \frac{1}{r}$. \label{corsymboles-dse}
\end{Cor}

\begin{Dem}
Let us write
\[ f(z) = \sum_{k \geq 0} a_k z^k \]
for $|z| < R$, so that, pointwise, 
\[ f \circ m = \sum_{k \geq 0} a_k m^k \]
Every symbol $a_k m^k$ corresponds to a continuous bilinear operator. Furthermore, by Lemma \ref{lemsymboles-puissancek}, the associated operator norm is dominated for $k \geq N$ by 
\begin{align*}
|a_k| \Vert m \Vert_{L^{\infty}}^{k-N} k^N C(p, q, r, m, N)
\end{align*} 
for some constant $C$ independant of $k$. In particular, this series of operator norms is convergent, so we deduce that $f \circ m$ also corresponds to a continuous bilinear operator. 
\end{Dem}

\begin{Def} In what follows, we simply call ``symbol'' a $m$ satisfying, by Corollary \ref{corsymboles-dse} or one of the previous theorems, sufficient conditions for the associated operator to be continuous. 
\end{Def}

\section{Spacetime resonances} \label{section-spacetimeresanal} 
\subsection{Quadratic case}

Recall that
\[ \omega\left(\overline{\xi}\right) = \xi_0^3 + \xi_0 |\xi|^2 \]
and 
\[ \varphi\left(\overline{\xi}, \overline{\eta}\right) = \omega\left(\overline{\xi}\right) -  \omega\left(\overline{\eta}\right) - \omega\left(\overline{\sigma}\right) \]
where we used the quadratic interaction convention $\overline{\sigma} = \overline{\xi} - \overline{\eta}$. 

\begin{Lem} Let $(\overline{\xi}, \overline{\eta}) \in S^5$ be such that
\[ \varphi\left(\overline{\xi}, \overline{\eta}\right) = \partial_{\eta_0} \varphi\left(\overline{\xi}, \overline{\eta}\right) = 0, \quad \nabla_{\eta} \varphi\left(\overline{\xi}, \overline{\eta}\right) = 0 \]
Then $\xi_0 = 0$ and one of the following hold:
\begin{itemize}
\item either $\eta_0 = 0$ and $|\eta| = |\sigma|$; 
\item or $\xi = 0$. 
\end{itemize} \label{lem-calculespacesresonantssimple}
\end{Lem}

\begin{Dem}
On the one hand, 
\[ \left\{ \begin{array}{l}
\partial_{\eta_0} \varphi = 0 \\
\nabla_{\eta} \varphi = 0 
\end{array} \right. ~ \iff ~ \left\{ \begin{array}{l}
3 \eta_0^2 + |\eta|^2 = 3 \sigma_0^2 + |\sigma|^2 \\
2 \eta_0 \eta = 2 \sigma_0 \sigma
\end{array} \right. ~ \Longrightarrow ~ \left\{ \begin{array}{l}
(\sqrt{3} |\eta_0| + |\eta|)^2 = (\sqrt{3} |\sigma_0| + |\sigma|)^2 \\
(\sqrt{3} |\eta_0| - |\eta|)^2 = (\sqrt{3} |\sigma_0| - |\sigma|)^2 
\end{array} \right. \]
In particular, $\sqrt{3} |\eta_0| + |\eta| = \sqrt{3} |\sigma_0| + |\sigma|$, and $\sqrt{3} |\eta_0| - |\eta| = \epsilon_1 \left( \sqrt{3} |\sigma_0| - |\sigma| \right)$ for some $\epsilon_1 \in \{ -1, 1 \}$. 

\paragraph{Case $\eta_0 = 0$} Assume first $\eta_0 = 0$. Then $2 \sigma_0 \sigma = 0$, so $\sigma_0 = 0$ or $\sigma = 0$. In the first case, we deduce $\xi_0 = 0$ and $|\eta| = |\sigma|$. 

In the second case, we deduce that $|\eta| = \sqrt{3} |\sigma_0| = \sqrt{3} |\xi_0|$, and $\xi = \eta$. Therefore, 
\[ \varphi\left(\overline{\xi}, \overline{\eta}\right) = \xi_0^3 + \xi_0 |\xi|^2 - \eta_0^3 - \eta_0 |\eta|^2 - \sigma_0^3 - \sigma_0 |\sigma|^2 = \xi_0 |\xi|^2 \]
The cancelation of $\varphi$ implies $\xi_0 = 0$ or $\xi = 0$, which is a contradiction since this would mean $\left(\overline{\xi}, \overline{\eta}\right) = 0 \notin S^5$. 

\paragraph{Case $\eta_0 \neq 0$} Assume now $\eta_0 \neq 0$. By symmetry, we may also assume $\sigma_0 \neq 0$. Hence, 
\[ \nabla_{\eta} \varphi\left(\overline{\xi}, \overline{\eta}\right) = 0 ~ \Longrightarrow ~ \eta = \frac{\sigma_0}{\eta_0} \sigma \]
Assume first $\epsilon_1 = 1$, so that
\[ |\sigma| = |\eta|, \quad |\sigma_0| = |\eta_0| \]
Let us set $\epsilon_2 = \frac{\sigma_0}{\eta_0} \in \{ -1 , 1 \}$, so that $\eta = \epsilon_2 \sigma$ and $\eta_0 = \epsilon_2 \sigma$. 

If $\epsilon_2 = -1$, then $\xi_0 = 0$, $\xi = 0$. 

If $\epsilon_2 = 1$, then $\overline{\xi} = 2 \overline{\eta} = 2 \overline{\sigma}$ and 
\[ \varphi(\overline{\xi}, \overline{\eta}) = \xi_0^3 - \eta_0^3 - \sigma_0^3 + \xi_0 |\xi|^2 - \eta_0 |\eta|^2 - \sigma_0 |\sigma|^2 
= 6 \eta_0 (\eta_0^2 + |\eta|^2) \]
We assumed $\eta_0 \neq 0$ so this is a contradiction. 

Assume now $\epsilon_1 = -1$, so that
\[ |\sigma| = \sqrt{3} |\eta_0|, \quad |\eta| = \sqrt{3} |\sigma_0| \]
We may now compute 
\begin{align*}
\varphi(\overline{\xi}, \overline{\eta}) 
&= \xi_0^3 - \eta_0^3 - \sigma_0^3 + \xi_0 |\xi|^2 - \eta_0 |\eta|^2 - \sigma_0 |\sigma|^2 
= 3 \xi_0 \eta_0 \sigma_0 + \eta_0 |\sigma|^2 + \sigma_0 |\eta|^2 + 2 \xi_0 \eta \cdot \sigma \\
&= 3 \xi_0 \eta_0 \sigma_0 + 3 \eta_0^3 + 3 \sigma_0^3 + 2 \sigma_0 |\sigma|^2 + 2 \eta_0 |\eta|^2 
= 3 \left( \eta_0^3 + \sigma_0^3 + 3 \eta_0^2 \sigma_0 + 3 \eta_0 \sigma_0^2 \right) 
= 3 \xi_0^3 
\end{align*}
Therefore, $\xi_0 = 0$, and we deduce that $|\sigma| = |\eta| = \sqrt{3} |\eta_0| = \sqrt{3} |\sigma_0|$. We recover the case $\epsilon_1 = 1$. 
\end{Dem}

\begin{Lem} Let $a, b_0, c_0 \in \mathbb{R}$ and $b \in \mathbb{R}^3$. Let us set
\[ M = \begin{pmatrix} 3 b_0 & b^T & \frac{3a + 3 c_0 - 6 b_0}{2} & - b^T \\
b & (a + b_0 + c_0) I_2 & - b & - \frac{2 b_0 + a + 3 c_0}{2} I_2  \\
\frac{3a + 3 c_0 - 6 b_0}{2} & - b^T & - (3 a + 6 c_0) & 0 \\
- b & - \frac{2 b_0 + a + 3 c_0}{2} I_2 & 0 & (a + 2 c_0) I_2 \end{pmatrix} \]
Let $\mathfrak{m}$ be the quadratic form associated to $M$. Then 
\[ a \varphi\left(\overline{\xi}, \overline{\eta}\right) + (b_0 \xi_0 + c_0 \eta_0) \partial_{\eta_0} \varphi\left(\overline{\xi}, \overline{\eta}\right) + \left( \frac{1}{2} (a + c_0) (\xi - 2 \eta) + \xi_0 b \right) \cdot \nabla_{\eta} \varphi\left(\overline{\xi}, \overline{\eta}\right) = \xi_0 \mathfrak{m}\left(\overline{\xi}, \overline{\eta}\right) \]
\end{Lem}

\begin{Dem}
This is a computation. 
\end{Dem}

\begin{Lem} Let $(\overline{\xi}, \overline{\eta}) \in S^5$ be such that $\overline{\xi} \neq 0$ or $3 \eta_0^2 \neq |\eta|^2$. Then, there exist a choice of $(a, b_0, c_0, b) \in \mathbb{R}^6$ such that $\mathfrak{m}$ does not vanish in the neighborhood of this point. 
\label{lem-non-res-loin0Cone}
\end{Lem}

\begin{Dem}
Assume first $\overline{\xi} \neq 0$. Then we choose $b = 0$, $b_0 = 1$, $c_0 = -2$, $a = 4$, so that 
\begin{align*} 
\mathfrak{m}\left(\overline{\xi}, \overline{\eta}\right) &= \Big( 3 \xi_0^2
+ \xi_0 \eta_0 (12 - 6 -6)
+ |\xi|^2 (4 + 1 -2)
+ \xi \cdot \eta (-4-2 +6) \\
&\quad + \eta_0^2 (-12 + 12)
+ |\eta|^2 (4 -4) \Big) \\
&= 3 (\xi_0^2 + |\xi|^2) 
\end{align*}
which does not vanish near the considered point. 

Then, if $\overline{\xi} = 0$ but $3 \eta_0^2 \neq |\eta|^2$, at first order we have: 
\[ \mathfrak{m}\left(\overline{\xi}, \overline{\eta}\right) = O\left(\overline{\xi}\right) + (a + 2 c_0) (|\eta|^2 - 3 \eta_0^2) \]
and the hypothesis $3 \eta_0^2 \neq |\eta|^2$ is therefore sufficient for $\mathfrak{m}$ not to vanish nearby, choosing for instance $c_0 = 1$, $a = 0$. 
\end{Dem}

We conclude this section by the following computation (similar to the one-dimensional resonance computation for Korteweg-de-Vries): 

\begin{Lem} Let $\xi_0, \eta_0, \sigma_0 \in \mathbb{R}$ be such that $\xi_0 = \eta_0 + \sigma_0$. Then
\begin{align*}
\xi_0^3 - \eta_0^3 - \sigma_0^3 &= 3 \xi_0 \eta_0 \sigma_0 
\end{align*} 
More generally, if $a_i, b_i, c_i \in \mathbb{R}$ are such that $a_i = b_i + c_i$ for $i = 1, 2, 3$, then 
\begin{align*}
\prod_{i = 1}^3 a_i ~ - \prod_{i = 1}^3 b_i ~ - \prod_{i = 1}^3 c_i = a_1 b_2 c_3 + b_1 c_2 a_3 + c_1 a_2 b_3
\end{align*}
\label{lemcalculvarphiKdV1D} 
\end{Lem}

\begin{Dem}
This is an elementary computation. 
\end{Dem}

\subsection{Cubic case}

We now consider
\begin{align*}
\varphi_3(\overline{\xi}, \overline{\eta}, \overline{\sigma}) = \omega(\overline{\xi}) - \omega(\overline{\eta}) - \omega(\overline{\sigma}) - \omega(\overline{\rho})
\end{align*}
where we follow the cubic interaction convention $\overline{\rho} = \overline{\xi} - \overline{\eta} - \overline{\sigma}$. 

\begin{Lem} Let $(\overline{\xi}, \overline{\eta}, \overline{\sigma}) \in S^8$ be such that 
\[ \varphi_3\left(\overline{\xi}, \overline{\eta}, \overline{\sigma}\right) = 0, \quad \nabla_{\overline{\eta}, \overline{\sigma}} \varphi_3\left(\overline{\xi}, \overline{\eta}, \overline{\sigma}\right) = 0 \]
Then $\xi_0 \nabla_{\overline{\xi}} \varphi_3 = 0$ and one of the following hold: 
\begin{itemize}
\item either $\overline{\xi} = \overline{\rho}, \overline{\eta} = - \overline{\sigma}$
\item or $\overline{\xi} = \overline{\eta}, \overline{\sigma} = - \overline{\rho}$
\item or $\overline{\xi} = \overline{\sigma}, \overline{\eta} = - \overline{\rho}$
\item or $\xi_0 = 0 = \eta_0 = \sigma_0 = \rho_0$. 
\end{itemize} 
Furthermore, in any case, $|\overline{\eta}| \simeq |\overline{\sigma}| \simeq |\overline{\rho}|$. \label{lemcalculresonancescubiquesgen} 
\end{Lem}

\begin{Dem}
First, $\partial_{\eta_0} \varphi_3 = 3 \eta_0^2 + |\eta|^2 - 3 \rho_0^2 - |\rho|^2$, so for this quantity to vanish it is needed to have $|\overline{\eta}| \simeq |\overline{\rho}|$. Likewise, we get $|\overline{\sigma}| \simeq |\overline{\rho}|$. In particular, here, since we chose the frequencies on the unit sphere, $|\overline{\eta}| \simeq |\overline{\sigma}| \simeq |\overline{\rho}| \simeq 1$. 

Since $\nabla_{\overline{\eta}} \varphi_3 = 0$, then 
\begin{align*}
\left\{ \begin{array}{l}
3 \eta_0^2 + |\eta|^2 = 3 \rho_0^2 + |\rho|^2 \\
2 \eta_0 \eta = 2 \rho_0 \rho
\end{array} \right. ~ \Longrightarrow ~ \left\{ \begin{array}{l}
\sqrt{3} |\eta_0| + |\eta| = \sqrt{3} |\rho_0| + |\rho| \\
\left| \sqrt{3} |\eta_0| - |\eta| \right| = \left| \sqrt{3} |\rho_0| - |\rho| \right| 
\end{array} \right.
\end{align*}
We may distinguish two cases: either $\sqrt{3} |\eta_0| = \sqrt{3} |\rho_0|$ and $|\eta| = |\rho|$, either $\sqrt{3} |\eta_0| = |\rho|$ and $\sqrt{3} |\rho_0| = |\eta|$. We have the same distinction replacing $\overline{\eta}$ by $\overline{\sigma}$. Up to some permutation of $\overline{\eta}, \overline{\sigma}, \overline{\rho}$, there are only two cases to consider. 

\paragraph{1.} If $|\eta_0| = |\sigma_0| = |\rho_0|$ and $|\eta| = |\sigma| = |\rho|$, coming back to $\eta_0 \eta = \rho_0 \rho = \sigma_0 \sigma$, we deduce that $\overline{\eta} = \epsilon_1 \overline{\sigma} = \epsilon_2 \overline{\rho}$ for some $\epsilon_1, \epsilon_2 \in \{ -1, 1 \}$, or that $\eta_0 = \sigma_0 = \rho_0 = 0$. 

\paragraph{1.1.} If $\overline{\eta} = \epsilon_1 \overline{\sigma} = \epsilon_2 \overline{\rho}$, up to a permutation of $\overline{\eta}, \overline{\sigma}, \overline{\rho}$, we may assume $\epsilon_1 = -1$ and $\epsilon_2 = 1$, or $\epsilon_1 = \epsilon_2 = 1$. In the first case, we have indeed $\overline{\xi} = \overline{\rho}$ and $\overline{\eta} = - \overline{\sigma}$, and therefore $\nabla_{\overline{\xi}} \varphi_3 = 0$. In the second case, an easy computation shows that $\varphi$ does not vanish. 

\paragraph{1.2.} If $\eta_0 = \sigma_0 = \rho_0 = 0$, then $\xi_0 = 0$. 

\paragraph{2.} If $\sqrt{3} |\eta_0| = \sqrt{3} |\sigma_0| = |\rho|$ and $|\eta| = |\sigma| = \sqrt{3} |\rho_0|$, then as before we also have $\overline{\eta} = \pm \overline{\sigma}$ or $\eta_0 = \sigma_0 = 0$. If $\overline{\eta} = - \overline{\sigma}$, then $\overline{\xi} = \overline{\rho}$ and therefore $\nabla_{\overline{\xi}} \varphi_3 = 0$. 

Let us then assume that $\overline{\eta} = \overline{\sigma}$. The equation $\eta_0 \eta = \rho_0 \rho$ ensures that $\rho$ is aligned with $\eta$ : if $\eta \neq 0$, we have also $\rho_0 \neq 0$. Up to a rotation, we may assume $\eta_2 = \rho_2 = \sigma_2 = 0$. Then we have $(\sqrt{3} \rho_0, \rho_1) = \epsilon_2 (\eta_1, \sqrt{3} \eta_0)$ for some $\epsilon_2 \in \{ -1, 1 \}$. We now develop: 
\begin{align*}
\varphi_3(\overline{\xi}, \overline{\eta}, \overline{\sigma}) &= \xi_0^3 + \xi_0 |\xi|^2 - \eta_0^3 - \eta_0 |\eta|^2 - \sigma_0^3 - \sigma_0 |\sigma|^2 - \rho_0^3 - \rho_0 |\rho|^2 \\
&= \left( \rho_0 + 2 \eta_0 \right)^3 + \left( \rho_0 + 2 \eta_0 \right) \left( \rho_1 + 2 \eta_1 \right)^2 - 2 \eta_0^3 - 2 \eta_0 \eta_1^2 - \rho_0^3 - \rho_0 \rho_1^2 \\
&= \left( \rho_0 + 2 \eta_0 \right)^3 + 3 \left( \rho_0 + 2 \eta_0 \right) \left( \eta_0 + 2 \rho_0 \right)^2 - 2 \eta_0^3 - 6 \eta_0 \rho_0^2 - \rho_0^3 - 3 \rho_0 \eta_0^2 \\
&= \rho_0^3 + 6 \eta_0 \rho_0^2 + 12 \eta_0^2 \rho_0 + 8 \eta_0^3 + 12 \rho_0^3 + 36 \eta_0 \rho_0^2 + 27 \eta_0^2 \rho_0 + 6 \eta_0^3 \\
&- 2 \eta_0^3 - 6 \eta_0 \rho_0^2 - \rho_0^3 - 3 \rho_0 \eta_0^2 \\
&= 12 \rho_0^3 + 36 \eta_0 \rho_0^2 + 36 \eta_0^2 \rho_0 + 12 \eta_0^3 \\
&= 12 (\rho_0 + \eta_0)^3
\end{align*}
For this quantity to vanish, it is needed that $\rho_0 + \eta_0 = 0$, so that $|\eta_0| = |\sigma_0| = |\rho_0|$ and $|\eta| = |\sigma| = |\rho|$, so we came back to case 1. 

Finally, if we assume that $\eta_0 = \sigma_0 = 0$, we also have $\rho = 0$, $\xi_0 = \rho_0$ and thus 
\begin{align*}
\varphi_3\left( \overline{\xi}, \overline{\eta}, \overline{\sigma} \right) &= \xi_0^3 + \xi_0 |\xi|^2 - \rho_0^3 \\
&= \xi_0 |\xi|^2 
\end{align*}
This implies $\xi_0 = 0$ or $\xi = 0 = \eta + \sigma$. If $\xi_0 = 0$, then $\rho_0 = 0$ and we came back to case 1. If $\xi = 0$, then $\overline{\eta} = - \overline{\sigma}$, $\overline{\xi} = \overline{\rho}$ and $\nabla_{\overline{\xi}} \varphi_3 = 0$. 
\end{Dem}

We end this subsection by a cubic equivalent of Lemma \ref{lemcalculvarphiKdV1D}: 

\begin{Lem} Let $a, b, c, d \in \mathbb{R}$ be such that $a = b - c + d$. Then 
\begin{align*}
a^3 - b^3 + c^3 - d^3 &= 3 (c-d)(c-b)(b+d)
\end{align*}
\label{lemcalculvarphiSchrodcub1D} 
\end{Lem}

\begin{Dem}
This is an elementary computation. 
\end{Dem}

\subsection{Conical coordinates} 

For $n$ the interaction order and $\overline{\eta}^i, \overline{\eta}^j$ two vectors, we now define coordinates adapted to these vectors: 

\begin{Def} Let $\overline{\eta}^i, \overline{\eta}^j \in \mathbb{R}^3$. 

Assume that $\overline{\eta}^j \notin \widehat{\mathcal{L}} \cup \widehat{\mathcal{P}}$. Then we define
\begin{align*}
\left( \overline{\eta}^i \right)_a^{\overline{\eta}^j} &:= \sqrt{3} \frac{\eta^j_0}{|\eta^j_0|} \eta^i_0 + \frac{\eta^j}{|\eta^j|} \cdot \eta^i \\
\left( \overline{\eta}^i \right)_b^{\overline{\eta}^j} &:= \sqrt{3} \frac{\eta^j_0}{|\eta^j_0|} \eta^i_0 - \frac{\eta^j}{|\eta^j|} \cdot \eta^i 
\end{align*}
We may extend this definition for $\left( \overline{\eta}^i \right)_{\alpha}^{\overline{\eta}^i}$ even if $\overline{\eta}^i \in \widehat{\mathcal{L}} \cup \widehat{\mathcal{P}}$. 

Assume $\eta^i, \eta^j \neq 0$. We also define the angular coordinates 
\begin{align*}
\xi_t^{\overline{\eta}^i \overline{\eta}^j} &:= \frac{J \eta^i \cdot \eta^j}{|\eta^i| |\eta^j|} \\
\theta^{\overline{\eta}^i \overline{\eta}^j} &:= \frac{\eta^i \cdot \eta^j}{|\eta^i| |\eta^j|} 
\end{align*}

Finally, assume $\eta^i_0 \neq 0, \eta^j_0 \neq 0$. We define the sign
\begin{align*}
\epsilon^{\overline{\eta}^i \overline{\eta}^j} &:= \frac{\eta^i_0 \eta^j_0}{|\eta^i_0| |\eta^j_0|} 
\end{align*}
\end{Def}

Note that, unlike $\left( \overline{\eta}^i \right)^{\overline{\eta}^j}_a$ and $\left( \overline{\eta}^i \right)^{\overline{\eta}^j}_b$ which are of order $1$ in $\overline{\eta}^i$ and of order $0$ in $\overline{\eta}^j$, the coordinates $\xi_t^{\overline{\eta}^i \overline{\eta}^j}$, $\theta^{\overline{\eta}^i \overline{\eta}^j}$ and $\epsilon^{\overline{\eta}^i \overline{\eta}^j}$ are of order $0$ in both $\overline{\eta}^i$ and $\overline{\eta}^j$, and can be seen as angular coordinates (hence the choice of the subscript $t$, as ``tangent''). 

\begin{Lem} We have for any $\overline{\eta}^j$: 
\begin{align*}
6 \sqrt{3} \frac{\eta_0^j}{|\eta_0^j|} \omega(\overline{\eta}^j) &= \left( \left( \overline{\eta}^j \right)_a^{\overline{\eta}^j} \right)^3 + \left( \left( \overline{\eta}^j \right)_b^{\overline{\eta}^j} \right)^3 
\end{align*}
In particular, for any choice of $\overline{\eta}^0, ..., \overline{\eta}^{n-1}$ and any $j \geq 1$: 
\begin{align*}
6 \sqrt{3} \varphi_n &= \frac{\eta_0^0}{|\eta_0|} \left( \left( \left( \overline{\eta}^0 \right)_a^{\overline{\eta}^0} \right)^3 + \left( \left( \overline{\eta}^0 \right)_b^{\overline{\eta}^0} \right)^3 \right) 
- \sum_{i = 1}^n \frac{\eta_0^i}{|\eta_0^i|} \left( \left( \left( \overline{\eta}^i \right)_a^{\overline{\eta}^i} \right)^3 + \left( \left( \overline{\eta}^i \right)_b^{\overline{\eta}^i} \right)^3 \right) \\
2 \partial_{\eta_0^j} \varphi_n &= \left( \left( \overline{\eta}^n \right)_a^{\overline{\eta}^n} \right)^2 + \left( \left( \overline{\eta}^n \right)_b^{\overline{\eta}^n} \right)^2 - \left( \left( \overline{\eta}^j \right)_a^{\overline{\eta}^j} \right)^2 - \left( \left( \overline{\eta}^j \right)_b^{\overline{\eta}^j} \right)^2 
\end{align*}
If $\eta^j \neq 0$: 
\begin{align*}
\frac{J \eta^j}{|\eta^j|} \cdot \nabla_{\eta^j} \varphi_n &= 2 \eta^n_0 |\eta^n| \xi_t^{\overline{\eta}^j \overline{\eta}^n} \\
2 \sqrt{3} \frac{\eta_0^j}{|\eta_0^j|} \frac{\eta^j}{|\eta^j|} \cdot \nabla_{\eta^j} \varphi_n &= \left( \left( \overline{\eta}^n \right)_a^{\overline{\eta}^j} \right)^2 - \left( \left( \overline{\eta}^n \right)_b^{\overline{\eta}^j} \right)^2 - \left( \left( \overline{\eta}^j \right)_a^{\overline{\eta}^j} \right)^2 + \left( \left( \overline{\eta}^j \right)_b^{\overline{\eta}^j} \right)^2 \\
2 \partial_{\eta_0^j} \varphi_n &= \left( \left( \overline{\eta}^n \right)_a^{\overline{\eta}^j} \right)^2 + \left( \left( \overline{\eta}^n \right)_b^{\overline{\eta}^j} \right)^2 + |\eta^n|^2 \left( \xi_t^{\overline{\eta}^j \overline{\eta}^n} \right)^2 - \left( \left( \overline{\eta}^j \right)_a^{\overline{\eta}^j} \right)^2 - \left( \left( \overline{\eta}^j \right)_b^{\overline{\eta}^j} \right)^2
\end{align*} 
Similar identities are valid for $j = 0$ up to the sign. 
\label{lemcalculssimplescoordonneesconiquevarphi} 
\end{Lem}

\begin{Dem}
These are elementary computations. 
\end{Dem}

\begin{Lem} Let $\overline{\eta}^i \notin \widehat{\mathcal{P}} \cup \widehat{\mathcal{L}}$, $\overline{\eta}^n$ be arbitrary. Then, 
\begin{align*}
\widehat{X}_a(\overline{\eta}^i) \cdot \nabla_{\overline{\eta}^i} \varphi_n &= \frac{\eta^i_0}{|\overline{\eta}^i|} \left( \left( \left( \overline{\eta}^n \right)_a^{\overline{\eta}^i} \right)^2 - \left( \left( \overline{\eta}^i \right)_a^{\overline{\eta}^i} \right)^2 \right) 
+ \frac{\eta^i_0}{|\overline{\eta}^i|} |\eta^n|^2 \left( \xi_t^{\overline{\eta}^i \overline{\eta}^n} \right)^2 \\
&\quad - \frac{\left( \overline{\eta}^i \right)_b^{\overline{\eta}^i}}{|\overline{\eta}^i|} \frac{\eta^i}{|\eta^i|} \cdot \nabla_{\eta^i} \varphi_n \\
\widehat{X}_c(\overline{\eta}^i) \cdot \nabla_{\overline{\eta}^i} \varphi_n &= 2 \eta^n_0 |\eta^n| \xi_t^{\overline{\eta}^i \overline{\eta}^n} \\
\left( \frac{\eta_0^i}{|\eta_0^i|} \partial_{\eta_0^i} - \sqrt{3} \frac{\eta^i}{|\eta^i|} \cdot \nabla_{\eta^i} \right) \varphi_n &= \frac{\eta_0^i}{|\eta_0^i|} \left( \left( \left( \overline{\eta}^n \right)_b^{\overline{\eta}^i} \right)^2 - \left( \left( \overline{\eta}^i \right)_b^{\overline{\eta}^i} \right)^2 \right) + \frac{1}{2} \frac{\eta_0^i}{|\eta_0^i|} |\eta^n|^2 \left( \xi_t^{\overline{\eta}^i \overline{\eta}^n} \right)^2 \\
6 \sqrt{3} \frac{\eta_0^i}{|\eta_0^i|} \varphi_n &= \left( \left( \overline{\eta}^0 \right)_a^{\overline{\eta}^i} \right)^3 + \left( \left( \overline{\eta}^0 \right)_b^{\overline{\eta}^i} \right)^3 - \sum_{j = 1}^n \left( \left( \left( \overline{\eta}^j \right)_a^{\overline{\eta}^i} \right)^3 + \left( \left( \overline{\eta}^j \right)_b^{\overline{\eta}^i} \right)^3 \right) \\
&\quad + 6 \sqrt{3} \frac{\eta_0^i}{|\eta_0^i|} \left( \eta_0^0 |\eta^0|^2 \left( \xi_t^{\overline{\eta}^i \overline{\eta}^0} \right)^2 - \sum_{j = 1}^n \eta_0^j |\eta^j|^2 \left( \xi_t^{\overline{\eta}^i \overline{\eta}^j} \right)^2 \right) 
\end{align*}
\label{lemcalculsconecoordonneesconiquesvarphi}
\end{Lem}

\begin{Dem}
These are elementary computations, following from Lemma \ref{lemcalculssimplescoordonneesconiquevarphi}. 
\end{Dem}

\subsection{Computations on the vector fields} 

\begin{Lem} For any choice of $\overline{\eta}^i, \overline{\eta}^j \notin \widehat{\mathcal{L}}$ and every $\alpha = a, b, c$ we have
\[ \widehat{X}_{\alpha}(\overline{\eta}^i) = \sum_{\beta = a, b, c} P^{\beta}_{\alpha}\left(\overline{\eta}^i, \overline{\eta}^j\right) \widehat{X}_{\beta}(\overline{\eta}^j) \]
where
\begin{equation}
P\left(\overline{\xi}, \overline{\eta}\right) = \left( \begin{array}{ccc} \dfrac{\xi_0 \eta_0 + \xi \cdot \eta}{|\overline{\xi}| |\overline{\eta}|} & \dfrac{\xi_0 |\eta|^2 - \eta_0 \xi \cdot \eta}{|\overline{\xi}| |\overline{\eta}| |\eta|} & \dfrac{\xi \cdot J \eta}{|\overline{\xi}| |\eta|} \\[5mm]
\dfrac{\eta_0 |\xi|^2 - \xi_0 \xi \cdot \eta}{|\overline{\xi}| |\overline{\eta}| |\xi|} & \dfrac{|\xi|^2 |\eta|^2 + \xi_0 \eta_0 \xi \cdot \eta}{|\overline{\xi}| |\overline{\eta}| |\xi| |\eta|} & - \dfrac{\xi_0 \xi \cdot J \eta}{|\overline{\xi}| |\xi| |\eta|} \\[5mm]
\dfrac{J \xi \cdot \eta}{|\xi| |\eta|} & - \dfrac{\eta_0 J \xi \cdot \eta}{|\overline{\eta}| |\xi| |\eta|} & \dfrac{\xi \cdot \eta}{|\xi| |\eta|} \end{array} \right) ~ = \begin{pmatrix} P_a^a & P_a^b & P_a^c \\ P_b^a & P_b^b & P_b^c \\ P_c^a & P_c^b & P_c^c \end{pmatrix} 
\end{equation}
\label{lem-projection-chpsvect}
\end{Lem}

\begin{Dem}
The vector fields $\left( \widehat{X}_{\alpha} \right)_{\alpha = a, b, c}$ form an orthonormal basis of $\mathbb{R}^3$, therefore the projection is obtained uniquely by computing scalar products between these vector fields. 
\end{Dem}

In the neighborhood of the cone $\widehat{\mathcal{C}}$, we will also use the following ``straightened'' vector fields: 
\begin{equation}
\begin{aligned}
\widehat{X}_a'(\overline{\xi}) &:= \begin{pmatrix} \frac{\xi_0}{2|\xi_0|} \\ \frac{\sqrt{3} \xi}{2 |\xi|} \end{pmatrix} \\
\widehat{X}_b'(\overline{\xi}) &:= \begin{pmatrix} \frac{\sqrt{3} \xi_0}{2 |\xi_0|} \\ -\frac{\xi}{2 |\xi|} \end{pmatrix} \\
\widehat{X}_c'(\overline{\xi}) &:= \begin{pmatrix} 0 \\ \frac{J \xi}{|\xi|} \end{pmatrix}
\end{aligned} \label{equdef-chpsvectanc-voiscone}
\end{equation}
Unlike the previous ones, they only depend on the sign of $\xi_0$ and the direction of $\xi$ (and not the direction of $\overline{\xi}$ for instance). They will become handy to compute the algebra of resonances. The following lemma shows that, near the cone, the use of this new basis is essentially equivalent to the old one: 

\begin{Lem} Let $\alpha \neq \beta \in \{ a, b, c \}$, and $\overline{\xi}$ be near the cone. Then 
\[ \left| \widehat{X}_{\alpha}(\overline{\xi}) \cdot \widehat{X}_{\beta}'(\overline{\xi}) \right| \lesssim \frac{\left| \sqrt{3} |\xi_0| - |\xi| \right|}{|\overline{\xi}|} \]
Furthermore, if $\alpha = a$ or $c$, 
\[ \left| \widehat{X}_{\alpha}(\overline{\xi}) \cdot \widehat{X}_{\alpha}'(\overline{\xi}) ~ - 1 \right| \lesssim \frac{\left| \sqrt{3} |\xi_0| - |\xi| \right|}{|\overline{\xi}|} \] 
and if $\alpha = b$: 
\[ \left| \widehat{X}_b(\overline{\xi}) \cdot \widehat{X}_b'(\overline{\xi}) - \frac{\xi_0}{|\xi_0|} \right| \lesssim \frac{\left| \sqrt{3} |\xi_0| - |\xi| \right|}{|\overline{\xi}|} \]\label{lemlienentrebaseancbasenouv} 
\end{Lem}

\begin{Dem}
If $\alpha = c$ or $\beta = c$, the result is obvious since $\widehat{X}_c = \widehat{X}_c'$ is orthogonal to $\widehat{X}_a, \widehat{X}_a', \widehat{X}_b, \widehat{X}_b'$. Then, 
\begin{align*}
\widehat{X}_a(\overline{\xi}) \cdot \widehat{X}_b'(\overline{\xi}) &= \frac{\sqrt{3} \xi_0^2}{2 |\overline{\xi}| |\xi_0|} - \frac{|\xi|}{2 |\overline{\xi}|} \\
&= \frac{\sqrt{3} |\xi_0| - |\xi|}{|\overline{\xi}|}
\end{align*}
Then, since the matrix $(\widehat{X}_{\alpha}(\overline{\xi}) \cdot \widehat{X}_{\beta}(\overline{\xi}))_{\alpha, \beta}$ is orthogonal, 
\begin{align*}
(\widehat{X}_a'(\overline{\xi}) \cdot \widehat{X}_b(\overline{\xi}))^2 &= 1 - (\widehat{X}_a'(\overline{\xi}) \cdot \widehat{X}_a(\overline{\xi})^2 - (\widehat{X}_a'(\overline{\xi}) \cdot \widehat{X}_c(\overline{\xi}))^2 \\
&= 1 - \left( 1 - (\widehat{X}_b'(\overline{\xi}) \cdot \widehat{X}_a(\overline{\xi}))^2 - (\widehat{X}_c'(\overline{\xi}) \cdot \widehat{X}_a(\overline{\xi}))^2 \right) - (\widehat{X}_a'(\overline{\xi}) \cdot \widehat{X}_c(\overline{\xi}))^2 \\
&= (\widehat{X}_b'(\overline{\xi}) \cdot \widehat{X}_a(\overline{\xi}))^2 
\end{align*}
Using orthogonality again, as we just showed that every non-diagonal terms are bounded by a constant times $\frac{\left| \sqrt{3} |\xi_0| - |\xi| \right|}{|\overline{\xi}|}$, we deduce that the diagonal is close to $\pm 1$ up to an error of size $\frac{\left| \sqrt{3} |\xi_0| - |\xi| \right|}{|\overline{\xi}|}$. It is then enough to identify the signs to obtain the second part of the lemma. Yet the $0$-coordinate of $\widehat{X}_a(\overline{\xi})$ and of $\widehat{X}_a'(\overline{\xi})$, respectively $\frac{\xi_0}{|\overline{\xi}|}$ and $\frac{\xi_0}{2 |\xi_0|}$, have the same sign, while the sign of the $0$-coordinate of $\widehat{X}_b(\overline{\xi})$ and of $\widehat{X}_b'(\overline{\xi})$, respectively $\frac{|\xi|}{|\overline{\xi}|}$ and $\frac{\sqrt{3} \xi_0}{2 |\xi_0|}$, differ by a factor $\frac{\xi_0}{|\xi_0|}$. 
\end{Dem}

Similarly, we have projection formulas on this new basis: 

\begin{Lem} For every choice of $\overline{\eta}^i, \overline{\eta}^j$ and every $\alpha = a, b, c$ we have 
\[ \widehat{X}_{\alpha}'(\overline{\eta}^i) = \sum_{\beta = a, b, c} \widetilde{P}^{\beta}_{\alpha}\left(\overline{\eta}^i, \overline{\eta}^j\right) \widehat{X}_{\beta}'(\overline{\eta}^j) \]
where 
\begin{equation}
\widetilde{P}\left(\overline{\xi}, \overline{\eta}\right) = \left( \begin{array}{ccc} \dfrac{1}{4} \epsilon^{\overline{\xi} \overline{\eta}} + \dfrac{3}{4} \theta^{\overline{\xi} \overline{\eta}} & \dfrac{\sqrt{3}}{4} \epsilon^{\overline{\xi} \overline{\eta}} - \dfrac{\sqrt{3}}{4} \theta^{\overline{\xi} \overline{\eta}} & \dfrac{\sqrt{3}}{2} \xi_t^{\overline{\eta} \overline{\xi}} \\[5mm]
\dfrac{\sqrt{3}}{4} \epsilon^{\overline{\xi} \overline{\eta}} - \dfrac{\sqrt{3}}{4} \theta^{\overline{\xi} \overline{\eta}} & \dfrac{3}{4} \epsilon^{\overline{\xi} \overline{\eta}} + \dfrac{1}{4} \theta^{\overline{\xi} \overline{\eta}} & \dfrac{1}{2} \xi_t^{\overline{\xi} \overline{\eta}} \\[5mm]
\dfrac{\sqrt{3}}{2} \xi_t^{\overline{\xi} \overline{\eta}} & - \dfrac{1}{2} \xi_t^{\overline{\eta} \overline{\xi}} & \theta^{\overline{\xi} \overline{\eta}} \end{array} \right) ~ = \begin{pmatrix} \widetilde{P}_a^{a} & \widetilde{P}_a^{b} & \widetilde{P}_a^{c} \\ \widetilde{P}_b^{a} & \widetilde{P}_b^{b} & \widetilde{P}_b^{c} \\ \widetilde{P}_c^{a} & \widetilde{P}_c^{b} & \widetilde{P}_c^{c} \end{pmatrix} 
\end{equation}
\label{lem-projection-chpsvectanc}
\end{Lem}

In the case where two Fourier variables are simultaneously on the cone (or nearby), the definition of the $X$-norm shows we have a good control over $\Vert m_b X_b f \Vert_{L^2}$ but not on $\Vert X_b f \Vert_{L^2}$ (which has a growth): in particular, later, as we apply the space-time resonances strategy and integrations by parts, we will want to avoid the $b$-direction for both interacting frequencies, or equivalently to apply integrations by parts along vector fields which are in the plane generated by $\widehat{X}_a$ and $\widehat{X}_c$, for both interacting frequencies. These two planes have no reason to coincide, but there exist at least one (non-zero) vector belonging to the intersection of these planes. 

We will thus now define a vector field depending on two Fourier variables such that it can be decomposed as a sum of $m_{\alpha} X_{\alpha}$ for both of the variables. 

\begin{Def} Let $i \in \{ 1, ..., n-1 \}$, and assume that $\overline{\eta}^i$ and $\overline{\eta}^n$ are close to the cone. 

\textbf{1.} Away from $\epsilon^{\overline{\eta}^i \overline{\eta}^n} \theta^{\overline{\eta}^i \overline{\eta}^n} = 1$, we define the modified vector field: 
\begin{equation}
\widehat{Y}(\overline{\eta}^i, \overline{\eta}^n) = \widetilde{P}_c^{b}(\overline{\eta}^i, \overline{\eta}^n) \widehat{X}_a'(\overline{\eta}^i) - \widetilde{P}_a^{b}(\overline{\eta}^i, \overline{\eta}^n) \widehat{X}_c'(\overline{\eta}^i) \label{equdef-champmodifie1} 
\end{equation}

\textbf{2.} Near $\epsilon^{\overline{\eta}^i \overline{\eta}^n} \theta^{\overline{\eta}^i \overline{\eta}^n} = 1$, we define it as: 
\begin{equation} \widehat{Y}(\overline{\eta}^i, \overline{\eta}^n) = \widehat{X}_a'(\overline{\eta}^i) - \frac{\widetilde{P}_a^{b}\left( \overline{\eta}^i, \overline{\eta}^n \right)}{1 - \widetilde{P}_b^{b}\left( \overline{\eta}^i, \overline{\eta}^n \right)^2} \left( \widehat{X}_b'(\overline{\eta}^n) - \widetilde{P}_b^{b}\left( \overline{\eta}^i, \overline{\eta}^n \right) \widehat{X}_b'(\overline{\eta}^i) \right) \label{equdef-champmodifie2} \end{equation}
\label{def-champmodifies}
\end{Def} 

The usefulness of this choice is given by the following lemma:

\begin{Lem} The vector field $\widehat{Y}(\overline{\eta}^i, \overline{\eta}^n)$ belongs to the intersection
\[ \mbox{Vect}(\widehat{X}_a'(\overline{\eta}^i), \widehat{X}_c'(\overline{\eta}^i)) \cap \mbox{Vect}(\widehat{X}_a'(\overline{\eta}^n), \widehat{X}_c'(\overline{\eta}^n)) \] \label{lemchampmodifieYpropfond}
\end{Lem}

\begin{Dem}
Near $\epsilon^{\overline{\eta}^i \overline{\eta}^n} \theta^{\overline{\eta}^i \overline{\eta}^n} = 1$, we have that 
\begin{align*}
\widehat{Y}(\overline{\eta}^i, \overline{\eta}^n) &= \frac{\widetilde{P}_c^{b}(\overline{\eta}^i, \overline{\eta}^n)^2}{1 - \widetilde{P}_b^{b}(\overline{\eta}^i, \overline{\eta}^n)^2} \widehat{X}_a'(\overline{\eta}^i) - \frac{\widetilde{P}_a^{b}(\overline{\eta}^i, \overline{\eta}^n) \widetilde{P}_c^{b}(\overline{\eta}^i, \overline{\eta}^n)}{1 - \widetilde{P}_b^{b}(\overline{\eta}^i, \overline{\eta}^n)^2} \widehat{X}_c'(\overline{\eta}^i)
\end{align*}
which corresponds to the vector field defined away from $\epsilon^{\overline{\eta}^i \overline{\eta}^n} \theta^{\overline{\eta}^i \overline{\eta}^n} = 1$, multiplied by $\frac{\widetilde{P}_c^{b}(\overline{\eta}^i, \overline{\eta}^n)}{1 - \widetilde{P}_b^{b}(\overline{\eta}^i, \overline{\eta}^n)^2}$. 

Considering now the expression in the case away from $\epsilon^{\overline{\eta}^i \overline{\eta}^n} \theta^{\overline{\eta}^i \overline{\eta}^n} = 1$, it is clear that $\widehat{Y}(\overline{\eta}^i, \overline{\eta}^n)$ is a linear combination of $\widehat{X}_a'(\overline{\eta}^i), \widehat{X}_c'(\overline{\eta}^i)$. On the other hand, 
\begin{align*}
\widehat{Y}(\overline{\eta}^i, \overline{\eta}^n) \cdot \widehat{X}_b(\overline{\eta}^n) 
&= \widetilde{P}_c^b(\overline{\eta}^i, \overline{\eta}^n) \widetilde{P}_a^b(\overline{\eta}^i, \overline{\eta}^n) - \widetilde{P}_a^b(\overline{\eta}^i, \overline{\eta}^n) \widetilde{P}_c^b(\overline{\eta}^i, \overline{\eta}^n) = 0
\end{align*}
therefore, the basis being orthonormal, $\widehat{Y}(\overline{\eta}^i, \overline{\eta}^n)$ is also a linear combination of $\widehat{X}_a'(\overline{\eta}^n)$ and $\widehat{X}_c'(\overline{\eta}^n)$. 
\end{Dem}

\begin{Rem} 
Even if $\widetilde{P}_b^b(\overline{\eta}^i, \overline{\eta}^n) = \frac{3 \epsilon^{\overline{\eta}^i \overline{\eta}^n} + \theta^{\overline{\eta}^i \overline{\eta}^n}}{4}$ is close to $\pm 1$ when $\epsilon^{\overline{\eta}^i \overline{\eta}^n} \theta^{\overline{\eta}^i \overline{\eta}^n}$ is close to $1$, and thus $\frac{1}{1 - \widetilde{P}_b^b(\overline{\eta}^i, \overline{\eta}^n)^2}$ is singular, we have that
\begin{align*} 
\frac{\widetilde{P}_a^b(\overline{\eta}^i, \overline{\eta}^n)}{1 - \widetilde{P}_b^b(\overline{\eta}^i, \overline{\eta}^n)^2} &= \frac{\sqrt{3}}{4} \frac{\epsilon^{\overline{\eta}^i \overline{\eta}^n} - \theta^{\overline{\eta}^i \overline{\eta}^n}}{1 - \left( \frac{3 \epsilon^{\overline{\eta}^i \overline{\eta}^n} + \theta^{\overline{\eta}^i \overline{\eta}^n}}{4} \right)^2} \\
&= 4 \sqrt{3} \frac{\epsilon^{\overline{\eta}^i \overline{\eta}^n} - \theta^{\overline{\eta}^i \overline{\eta}^n}}{\left( 4 \epsilon^{\overline{\eta}^i \overline{\eta}^n} - 3 \epsilon^{\overline{\eta}^i \overline{\eta}^n} - \theta^{\overline{\eta}^i \overline{\eta}^n} \right) \left( 4 \epsilon^{\overline{\eta}^i \overline{\eta}^n} + 3 \epsilon^{\overline{\eta}^i \overline{\eta}^n} + \theta^{\overline{\eta}^i \overline{\eta}^n} \right)} \\
&= \frac{4 \sqrt{3}}{7 \epsilon^{\overline{\eta}^i \overline{\eta}^n} + \theta^{\overline{\eta}^i \overline{\eta}^n}} 
\end{align*}
which is not singular. 

Furthermore, the distinction above is useful in order to keep normalized vector fields. Indeed, if $\epsilon^{\overline{\eta}^i \overline{\eta}^n} \theta^{\overline{\eta}^i \overline{\eta}^n}$ is close to $1$, then  $\widetilde{P}_a^b(\overline{\eta}^i, \overline{\eta}^n) = \frac{\sqrt{3}}{4} (\epsilon^{\overline{\eta}^i \overline{\eta}^n} - \theta^{\overline{\eta}^i \overline{\eta}^n})$ et $\widetilde{P}_c^b(\overline{\eta}^i, \overline{\eta}^n)$ are close to $0$, so that the vector field defined by \eqref{equdef-champmodifie1} tends to $0$. Conversely, when $\epsilon^{\overline{\eta}^i \overline{\eta}^n} \theta^{\overline{\eta}^i \overline{\eta}^n}$ is close to $-1$, we have that 
\begin{align*} &\widehat{X}_a'(\overline{\eta}^i) - \frac{\widetilde{P}_a^b\left( \overline{\eta}^i, \overline{\eta}^n \right)}{1 - \widetilde{P}_b^b\left( \overline{\eta}^i, \overline{\eta}^n \right)^2} \left( \widehat{X}_b'(\overline{\eta}^n) - P_b^b\left( \overline{\eta}^i, \overline{\eta}^n \right) \widehat{X}_b'(\overline{\eta}^i) \right) \\
&\quad = \frac{\widetilde{P}_c^b(\overline{\eta}^i, \overline{\eta}^n)^2}{1 - \widetilde{P}_b^b(\overline{\eta}^i, \overline{\eta}^n)^2} \widehat{X}_a'(\overline{\eta}^i) + \frac{\widetilde{P}_a^b(\overline{\eta}^i, \overline{\eta}^n) \widetilde{P}_b^c(\overline{\eta}^i, \overline{\eta}^n)}{1 - \widetilde{P}_b^b(\overline{\eta}^i, \overline{\eta}^n)^2} \widehat{X}_c'(\overline{\eta}^i) \\
&\quad = \frac{\frac{1}{4} (1 - \left( \theta^{\overline{\eta}^i \overline{\eta}^n} \right)^2)}{\frac{1}{16} (\epsilon^{\overline{\eta}^i \overline{\eta}^n} - \theta^{\overline{\eta}^i \overline{\eta}^n}) (7 \epsilon^{\overline{\eta}^i \overline{\eta}^n} + \theta^{\overline{\eta}^i \overline{\eta}^n})} \widehat{X}_a'(\overline{\eta}^i) + \frac{4 \sqrt{3}}{7 \epsilon^{\overline{\eta}^i \overline{\eta}^n} + \theta^{\overline{\eta}^i \overline{\eta}^n}} \widetilde{P}_b^c(\overline{\eta}^i, \overline{\eta}^n) \widehat{X}_c'(\overline{\eta}^i)
\end{align*}
so it also goes to $0$. Separating the cases is unavoidable as the intersection of the two planes spanned by $\widehat{X}_a, \widehat{X}_c$ when one of the Fourier variables is fixed and the other is moving around the cone $\widehat{\mathcal{C}}$ has the topological structure of a Möbius band. 

The case $\epsilon^{\overline{\eta}^i \overline{\eta}^n} \theta^{\overline{\eta}^i \overline{\eta}^n}$ close to $1$ corresponds to the case where $\overline{\eta}^i$ and $\overline{\eta}^n$ are close to colinearity. 
\end{Rem}

\begin{Lem} Let $i \in \{ 1, ..., n-1 \}$, $\overline{\eta}^i$ and $\overline{\eta}^n$ be close to the cone. Denote by $\epsilon_i = \frac{\eta_0^i}{|\eta_0^i|}$. 

\textbf{1.} Away from $\epsilon^{\overline{\eta}^i \overline{\eta}^n} \theta^{\overline{\eta}^i \overline{\eta}^n} = 1$, we have
\[ \widehat{Y}(\overline{\eta}^i, \overline{\eta}^n) \cdot \nabla_{\overline{\eta}^i} \varphi_n = - \frac{1}{4} \epsilon_i \xi_t^{\overline{\eta}^i \overline{\eta}^n} \left( \left( \left( \overline{\eta}^n \right)_a^{\overline{\eta}^n} \right)^2 - \left( \left( \overline{\eta}^i \right)_a^{\overline{\eta}^i} \right)^2 \right) \]

\textbf{2.} Near $\epsilon^{\overline{\eta}^i \overline{\eta}^n} \theta^{\overline{\eta}^i \overline{\eta}^n} = 1$, we have
\[ \widehat{Y}(\overline{\eta}^i, \overline{\eta}^n) \cdot \nabla_{\overline{\eta}^i} \varphi_n = 2 \epsilon_i \frac{1 + \epsilon^{\overline{\eta}^i \overline{\eta}^n} \theta^{\overline{\eta}^i \overline{\eta}^n}}{7 + \epsilon^{\overline{\eta}^i \overline{\eta}^n} \theta^{\overline{\eta}^i \overline{\eta}^n}} \left( \left( \left( \overline{\eta}^n \right)_a^{\overline{\eta}^n} \right)^2 - \left( \left( \overline{\eta}^i \right)_a^{\overline{\eta}^i} \right)^2 \right) \]
\label{lemcalculchampmodifie}
\end{Lem}

\begin{Dem}
These are tedious but elementary computations. 
\end{Dem}

\begin{Rem} A symbol of the form $\frac{1}{7 + \frac{\eta_0^i \eta_0^n}{|\eta_0^i||\eta_0^n|} \frac{\eta^i \cdot \eta^j}{|\eta^i| |\eta^j|}}$ does satisfy every desired Hölder estimates by applying Corollary \ref{corsymboles-dse}. Indeed, every $\frac{\eta^i}{|\eta^i|}$ is a (vector-valued) Hörmander-Mikhlin symbol type, of $L^{\infty}$ norm equal to $1$, and the function $z \mapsto \frac{1}{7 + z}$ admits an expansion in power series, with radius of convergence $R = 7 > 1$. 
\end{Rem}

In presence of a factor of the form $\widehat{X}_b'(\overline{\eta}^n) \cdot \nabla_{\overline{\eta}^i} \widehat{f}(t, \overline{\eta}^n)$, a priori not controlled, it will be needed to apply an integration by parts in $\overline{\eta}^i$ so to recover controlled terms. However, if $\overline{\eta}^i$ is near $\widehat{\mathcal{C}}$ or $\widehat{\mathcal{L}}$, a direct integration by parts gives a factor $\widetilde{P}_b^b(\overline{\eta}^n, \overline{\eta}^i) \widehat{X}_b'(\overline{\eta}^i) \cdot \nabla_{\overline{\eta}^i} \widehat{f}(t, \overline{\eta}^i)$, also not controlled. This motivates the introduction of the following vector field: 

\begin{Def} Let $\overline{\eta}^i, \overline{\eta}^n$ be near the cone. We defined the modified $b$-vector field as: 
\begin{equation} \widehat{X}_{b-\widehat{\mathcal{C}}}(\overline{\eta}^i, \overline{\eta}^n) = \widehat{X}_b'(\overline{\eta}^n) - \frac{\widetilde{P}_b^b(\overline{\eta}^n, \overline{\eta}^i)}{\widetilde{P}_a^b(\overline{\eta}^i, \overline{\eta}^n)^2 + \widetilde{P}_c^b(\overline{\eta}^i, \overline{\eta}^n)^2} \left( \widetilde{P}_a^b(\overline{\eta}^n, \overline{\eta}^i) \widehat{X}_a'(\overline{\eta}^n) + \widetilde{P}_c^b(\overline{\eta}^n, \overline{\eta}^i) \widehat{X}_c'(\overline{\eta}^n) \right) \label{equdef-champmodifiebCC} \end{equation}
\label{def-champmodifiebCC}
\end{Def}

\begin{Lem} The vector $\widehat{X}_{b-\widehat{\mathcal{C}}}(\overline{\eta}^i, \overline{\eta}^n)$ belongs to
\[ \left( \widehat{X}_b'(\overline{\eta}^n) + \mbox{Vect}(\widehat{X}_a'(\overline{\eta}^n), \widehat{X}_c'(\overline{\eta}^n)) \right) \cap \mbox{Vect}(\widehat{X}_a'(\overline{\eta}^i), \widehat{X}_c'(\overline{\eta}^i)) \]
\label{lemchampmodifiebCC-propfond}
\end{Lem}

\begin{Rem} The vector field $\widehat{X}_{b-\widehat{\mathcal{C}}}$ may be singular if $P_a^b, P_c^b$ vanish (in which case $P_b^b$ is close to $\pm 1$). With previous notations, this only happens if $\epsilon^{\overline{\eta}^i \overline{\eta}^n} \theta^{\overline{\eta}^i \overline{\eta}^n}$ is close to $1$. 
\end{Rem}

\begin{Lem} We have that
\begin{align*} 
&\widehat{X}_{b-\widehat{\mathcal{C}}}(\overline{\eta}^i, \overline{\eta}^n) \cdot \nabla_{\overline{\eta}^i} \varphi_n \\
&\quad = \epsilon_i \frac{2 \sqrt{3}}{7 \epsilon^{\overline{\eta}^i \overline{\eta}^n} + \theta^{\overline{\eta}^i \overline{\eta}^n}} \left( \left( \left( \overline{\eta}^n \right)_b^{\overline{\eta}^n} \right)^2 - \left( \left( \overline{\eta}^i \right)_a^{\overline{\eta}^i} \right)^2 - 2 \sqrt{3} \frac{1 + \epsilon^{\overline{\eta}^i \overline{\eta}^n} \theta^{\overline{\eta}^i \overline{\eta}^n}}{3} |\eta_0^n| |\eta^n| \right) 
\end{align*}
using again the notations from Lemma \ref{lemcalculchampmodifie}. \label{lemcalculchampmodifiebCC} 
\end{Lem}

\begin{Dem}
Again, this is a computation. 
\end{Dem}

This concludes the case of cone-cone interactions. The case of line-line interactions is simpler because $(1, 0, 0)$ always belong to $\mbox{Vect}(\widehat{X}_a, m_b \widehat{X}_b, m_c \widehat{X}_c)$, for any choice of Fourier variable close to $\widehat{\mathcal{L}}$. Finally, in the case of a cone-line interaction, if $\overline{\eta}^n$ is close to $\widehat{\mathcal{C}}$ and $\overline{\eta}^i$ close to $\widehat{\mathcal{L}}$ for instance, then it is easy to see that $\mbox{Vect}(\widehat{X}_a(\overline{\eta}^n), \widehat{X}_c(\overline{\eta}^n)) \cap \mbox{Vect}(\widehat{X}_a(\overline{\eta}^i)) = \{ (0, 0, 0) \}$. 

However, we may compute a modified $b$-vector field as before. 

\begin{Def} Let $\overline{\eta}^i$ be near $\widehat{\mathcal{L}}$ and $\overline{\eta}^n$ be near $\widehat{\mathcal{C}}$. We define the modified $b$-vector field: 
\begin{equation} \widehat{X}_{b-\widehat{\mathcal{L}}}(\overline{\eta}^i, \overline{\eta}^n) = \widehat{X}_b'(\overline{\eta}^n) + \frac{1}{\sqrt{3}} \widehat{X}_a'(\overline{\eta}^n) \label{equdefchampmodifiebLC} \end{equation}
\label{defchampmodifiebLC}
\end{Def}

\begin{Lem} The vector $\widehat{X}_{b-\widehat{\mathcal{L}}}(\overline{\eta}^i, \overline{\eta}^n)$ belongs to: 
\[ \left( \widehat{X}_b'(\overline{\eta}^n) + \mbox{Vect}(\widehat{X}_a'(\overline{\eta}^n), \widehat{X}_c'(\overline{\eta}^n)) \right) \cap \mbox{Vect}\left( (1, 0, 0) \right) \] \label{lemchampmodifiebLCpropfond}
\end{Lem}

\begin{Lem} We have that
\[ \widehat{X}_{b-\widehat{\mathcal{L}}}(\overline{\eta}^i, \overline{\eta}^n) \cdot \nabla_{\overline{\eta}^i} \varphi_n = \frac{2 \sqrt{3}}{3} \frac{\eta_0^n}{|\eta_0^n|} \left( 3 (\eta_0^n)^2 + |\eta^n|^2 - 3 (\eta_0^i)^2 - |\eta^i|^2 \right) \]
\label{lemcalculchampmodifiebLC} 
\end{Lem}

\begin{Dem}
This is an elementary computation. 
\end{Dem}

\begin{Rem} 
If $\overline{\eta}^i$ is near $\widehat{\mathcal{P}}$, we may express $\nabla_{\eta^i}$ as a linear combination with (smooth) symbol coefficients of the $\widehat{X}_a(\overline{\eta}^i) \cdot \nabla_{\overline{\eta}^i}, \widehat{X}_c(\overline{\eta}^i) \cdot \nabla_{\overline{\eta}^i}, \frac{\eta^i_0}{|\overline{\eta}^i|} \widehat{X}_b(\overline{\eta}^i) \cdot \nabla_{\overline{\eta}^i}$. 
\end{Rem}

\section{Dispersive estimate} \label{section-Linfdispersive} 

The aim of this section is to prove a $L^{\infty}$ dispersive estimate for the linear part of Equation \eqref{equ-ZK}, but in a way that suits the choice of the weighted norm. In particular, it is not enough to have an estimate of the form
\[ \Vert e^{i t \omega(D)} \Vert_{L^1 \to L^{\infty}} \leq t^{-\gamma} \]
for some $\gamma > 1$ because $\Vert u \Vert_X$ provides no control over $\Vert f(t) \Vert_{L^1}$ due to the lack of control over $X_b f, X_b X_b f$, which corresponds heuristically to a loss of integrability in a (frequency-depending) direction. Furthermore, even in presence of the symbol $m_b(D)$ to compensate this loss, there is a time growth of $\partial_x (m_b(D) X_b)^2 f(t)$ in $L^2$ in general. 

However, to prove a standard $L^1 \to L^{\infty}$ dispersive estimate, one typically applies stationary phase arguments and thus integrations by parts. The strategy for our proof is therefore to do such estimates in a finer way, writing $e^{-i t \omega(D)} f(t)$ as an integral by Fourier's inversion formula and trying to apply integrations by parts only using the directions $m_{\alpha} \widehat{X}_{\alpha}$, separating different cases depending on the direction and sizes of the position $(x, y)$ and the frequency $\overline{\xi}$. 

In all of this section, we will denote by $(\chi, \psi)$ Littlewood-Paley localisation functions, with $\psi$ localizing on an annulus and $\chi$ on a ball. For simplicity, we write them the same way, whatever the dimension of the domain, setting for an arbitrary vector $\vec v \in \mathbb{R}^d$, $d \geq 1$, 
\begin{align*}
\chi(\vec v) = \chi(|\vec v|), \quad \psi(\vec v) = \psi(|\vec v|)
\end{align*}

The main result of this section is the following: 

\begin{Prop} Assume that $\delta > 0$ in the norm $\Vert \cdot \Vert_X$ is small enough. Then there exists $C > 0$ such that, for any $t > 0$, 
\begin{align*}
\Vert e^{-i t \omega(D)} \partial_x f(t) \Vert_{L^{\infty}} &\leq C t^{-\frac{5}{6}} \langle t \rangle^{-\frac{1}{4}+100\delta} \Vert u \Vert_X \\
\Vert e^{-i t \omega(D)} \partial_x |\nabla|^{\frac{1}{2}} f(t) \Vert_{L^{\infty}} &\leq C t^{-\frac{5}{6}} \langle t \rangle^{-\frac{1}{3}+100\delta} \Vert u \Vert_X \\
\Vert e^{-i t \omega(D)} \nabla f(t) \Vert_{L^{\infty}} &\leq C t^{-\frac{5}{6}} \langle t \rangle^{-\frac{1}{6}+100\delta} \Vert u \Vert_X \\
\Vert e^{-i t \omega(D)} f(t) \Vert_{L^{\infty}} &\leq C t^{-\frac{5}{6}} \langle t \rangle^{\frac{1}{6}+100\delta} \Vert u \Vert_X 
\end{align*} \label{prop-estimeedispersivetotale}
\end{Prop} 

Recall that the $\Vert \cdot \Vert_X$ depends on $T, \delta$ a priori, but since $\delta > 0$ can be considered as a fixed and very small parameter, and $T$ will be pushed to $+\infty$, we do not write the dependance explicitely. In particular, the universal constant $C$ above does not depend on $T$ but will depend on $\delta$. 

In order to apply stationary phase estimates, we will consider the following phase: 
\begin{align*}
\Phi = \Phi\left( \overline{\xi}, \frac{x}{t}, \frac{y}{t} \right) &= \xi_0^3 + \xi_0 |\xi|^2 - \frac{\xi_0 x}{t} - \frac{\xi \cdot y}{t} 
\end{align*}
so that, for $F$ such that $\widehat{F} \in L^1$, 
\begin{align*}
e^{-i t \omega(D)} F(x, y) &= \int e^{i t \Phi} \widehat{F}(\overline{\xi}) ~ d\overline{\xi}
\end{align*}
up to an universal multiplicative constant, which we will omit. 

For fixed $(x, y, t)$, we will denote by 
\begin{align*}
\overline{\xi_a} &:= \sqrt{\frac{|x|+|y|}{t}}
\end{align*}
the characteristic size of resonant frequencies. If frequencies are localized of size $2^j$, we will use the notation 
\begin{align*}
\mathfrak{t} &:= t 2^{3j}
\end{align*}
which is the expected scaling parameter of the argument. 

We will also use the following localisation functions: first, for $j \in \mathbb{Z}$, 
\begin{align*}
\psi_j(\overline{\xi}) &:= \psi\left( 2^{-j} \overline{\xi} \right) 
\end{align*}
$\psi_j$ localises on frequencies of size $2^j$. For $k \leq -10$, we also define 
\begin{align*}
\psi_{j, k}^{\widehat{\mathcal{C}}}(\overline{\xi}) &:= \psi_j(\overline{\xi}) \psi\left( 2^{-k} \frac{\sqrt{3} |\xi_0| - |\xi|}{|\overline{\xi}|} \right) \\
\psi_{j, k}^{\widehat{\mathcal{L}}}(\overline{\xi}) &:= \psi_j(\overline{\xi}) \psi\left( 2^{-k} \frac{\xi}{|\overline{\xi}|} \right) \\
\psi_{j, k}^{\widehat{\mathcal{P}}}(\overline{\xi}) &:= \psi_j(\overline{\xi}) \psi\left( 2^{-k} \frac{\xi_0}{|\overline{\xi}|} \right) 
\end{align*}
that allow to localise also the distance to one of the geometric areas. 

Along the argument, we will not use the precise form of $\psi$, other than by assuming that $\psi$ is regular enough, and estimates will implicitely depend on bounds on $\psi$ and its derivatives, up to a finite (universal) order. In particular, it will always be possible to replace $\psi$ in the definition of the $\psi_{j, k}$ by another function, possibly depending on $j$ (as long as the estimates are uniform in $j$), so one can easily deduce the validity of Remark \ref{remarqueestimeedispersiveavecsymbole} (and even with more singular symbols, notably near $\widehat{\mathcal{C}}$). 

Note that, for any fixed $j$, $\psi_j(D) f$ is the Fourier transform of a $L^1$ function, so that we can always apply Fourier's inversion formula on localized version of $f$ (and we may add derivatives and other symbols). 

To prove the estimates, we will also need the following bounds on $g$, for any $\alpha$: 
\begin{align*}
\Vert \psi_{j, k}^{\widehat{\mathcal{C}}}(D) g_b(t) \Vert_{L^2} &\lesssim 2^{k-\frac{j}{2}} \langle 2^j \rangle^{-1} t^{-\frac{1}{6}+100\delta} \Vert u \Vert_X^2 + t^{-\frac{1}{2}+201\delta} 2^{-3j} \Vert u \Vert_X^2 \\
\Vert \psi_j(D) m_{\widehat{\mathcal{C}}}(D) X_c g_{\alpha}(t) \Vert_{L^2} &\lesssim \left( t^{-\frac{1}{12}+\frac{1}{32}+50\delta} + t^{-\frac{1}{8}} 2^{-\frac{7j}{6}+1000\delta j} \right) \Vert u \Vert_X^2
\end{align*}
Note that these estimates depend on an implicit constant, that will itself depend on the implicit constant in the dispersive estimate that we are proving in this section. However, since all these estimates are at least in $\Vert u \Vert_X^2$, the bootstrap argument allows to absorb the implicit constant by one of the factors $\Vert u \Vert_X$. 

The estimates on $g$ will be proven in Sections \ref{section-estimees-simples-fug}, \ref{section-estimees-quadratiques-horsbb} and \ref{section-estimee-quad-bb}. 

\subsection{Low frequencies} 

\begin{Lem} Let $j \in \mathbb{Z}$, $t > 0$. Then
\begin{align*}
\Vert e^{i t \omega(D)} \psi_j(D) f(t) \Vert_{L^{\infty}} &\lesssim 2^{\frac{5j}{2}} \langle 2^j \rangle^{-2} \Vert u \Vert_X 
\end{align*}
\label{lemestdisph-bassesfreq}
\end{Lem}

\begin{Dem}
Let $(x, y) \in \mathbb{R}^3$ be fixed. Then
\begin{align*} 
\left| \left[ e^{i t \omega(D)} \psi_j(D) f(t) \right](x, y) \right| &= \left| \int e^{i t \Phi} \psi_j(\overline{\xi}) \widehat{f}(t, \xi_0, \xi) ~ d\xi_0 d\xi \right| \\
&\lesssim \int |\psi_j(\overline{\xi})| 2^j \langle 2^j \rangle^{-2} \langle \overline{\xi} \rangle^2 |\overline{\xi}|^{-1} |\widehat{f}(t, \xi_0, \xi)| ~ d\xi_0 d\xi \\
&\lesssim 2^j \langle 2^j \rangle^{-2} \Vert \psi_j \Vert_{L^2(\mathbb{R}^3)} \Vert \langle \overline{\xi} \rangle^2 |\overline{\xi}|^{-1} \widehat{f}(t) \Vert_{L^2} \\
&\lesssim 2^{\frac{5j}{2}} \langle 2^j \rangle^{-2} \Vert |\nabla|^{-1} f(t) \Vert_{L^2} \lesssim 2^{\frac{5j}{2}} \langle 2^j \rangle^{-2} \Vert u \Vert_X
\end{align*}
by Parseval's identity and Lemma \ref{lem-estimees-directes-uX}. 
\end{Dem}

\subsection{Remainder area}

\begin{Lem} Let $t > 1$, $j \in \mathbb{Z}$ be such that $2^j \gg t^{-\frac{1}{3}}$. Then 
\begin{align*}
\Vert e^{i t \omega(D)} m_{\widehat{\mathcal{R}}}(D) \psi_j(D) f(t) \Vert_{L^{\infty}} &\lesssim t^{-\frac{7}{6}} \langle t \rangle^{100 \delta} 2^{-j} \langle 2^j \rangle^{-1} \Vert u \Vert_X 
\end{align*} \label{lem-estdisph-zonereste}
\end{Lem}

\begin{Dem}
Let us fix $(x, y) \in \mathbb{R}^3, t > 0$. Then by Fourier's inversion formula it is enough to bound: 
\begin{align*}
e^{i t \omega(D)} m_{\widehat{\mathcal{R}}}(D) \psi_j(D) f(t, x, y) &= \int e^{i t \Phi} m_{\widehat{\mathcal{R}}}(\overline{\xi}) \psi_j(\overline{\xi}) \widehat{f}(t, \overline{\xi}) ~ d\overline{\xi}
\end{align*}
up to a multiplicative constant. 

We have that 
\begin{align*}
\widehat{X}_a(\overline{\xi}) \cdot \nabla_{\overline{\xi}} \Phi &= \frac{\xi_0}{|\overline{\xi}|} \left( 3 \xi_0^2 + |\xi|^2 \right) + \frac{\xi}{|\overline{\xi}|} \cdot 2 \xi_0 \xi - \frac{\xi_0 x + \xi \cdot y}{|\overline{\xi}| t} \\
&= 3 \xi_0 |\overline{\xi}| - \frac{(x, y)}{t} \cdot \frac{\overline{\xi}}{|\overline{\xi}|} \\
\widehat{X}_c(\overline{\xi}) \cdot \nabla_{\overline{\xi}} \Phi &= - \frac{J \xi \cdot y}{t |\xi|} \\
\widehat{X}_b(\overline{\xi}) \cdot \nabla_{\overline{\xi}} \Phi &= \frac{|\xi|}{|\overline{\xi}|} \left( 3 \xi_0^2 + |\xi|^2 \right) - \frac{\xi_0 \xi}{|\overline{\xi}| |\xi|} \cdot 2 \xi_0 \xi - \frac{|\xi|^2 x - \xi_0 \xi \cdot y}{|\overline{\xi}| |\xi| t} \\
&= |\xi| |\overline{\xi}| - \frac{|\xi|^2 x - \xi_0 \xi \cdot y}{|\overline{\xi}| |\xi| t} 
\end{align*}

\paragraph{High frequencies} Assume that $\overline{\xi_a} \ll 2^j$. Then, we have that
\begin{align*}
\widehat{X}_a(\overline{\xi}) \cdot \nabla_{\overline{\xi}} \Phi &\simeq 2^{2j} 
\end{align*}
Therefore, we can apply an integration by parts: 
\begin{align*}
&\int e^{i t \Phi} m_{\widehat{\mathcal{R}}}(\overline{\xi}) \psi_j(\overline{\xi}) \widehat{f}(t, \overline{\xi}) ~ d\overline{\xi} \\
&\quad = t^{-1} \int e^{i t \Phi} \nabla_{\overline{\xi}} \cdot \left( \frac{\widehat{X}_a(\overline{\xi})}{\widehat{X}_a(\overline{\xi}) \cdot \nabla_{\overline{\xi}} \Phi} m_{\widehat{\mathcal{R}}}(\overline{\xi}) \psi_j(\overline{\xi}) \widehat{f}(t, \overline{\xi}) \right) ~ d\overline{\xi} \\
&\quad = t^{-1} \int e^{i t \Phi} 2^{-3j} \psi_j(\overline{\xi}) \widehat{f}(t, \overline{\xi}) ~ d\overline{\xi} \\
&\quad \quad + t^{-1} \int e^{i t \Phi} 2^{-2j} \psi_j(\overline{\xi}) \widehat{h_a}(t, \overline{\xi}) ~ d\overline{\xi} \\
&\quad \quad + t^{-1} \int e^{i t \Phi} 2^{-2j} \psi_j(\overline{\xi}) \widehat{g_a}(t, \overline{\xi}) ~ d\overline{\xi} 
\end{align*}
where $\psi_j$ may vary from line to line as long as it keeps similar properties (of regularity and localisation). On the two first terms, we apply a second integration by part: 
\begin{subequations}
\begin{align}
&\int e^{i t \Phi} m_{\widehat{\mathcal{R}}}(\overline{\xi}) \psi_j(\overline{\xi}) \widehat{f}(t, \overline{\xi}) ~ d\overline{\xi} \notag \\
&\quad = t^{-2} 2^{-6j}\int e^{i t \Phi} \psi_j(\overline{\xi}) \widehat{f}(t, \overline{\xi}) ~ d\overline{\xi} \label{estdispRHF-1} \\
&\quad \quad t^{-2} 2^{-5j} \int e^{i t \Phi} \psi_j(\overline{\xi}) \widehat{X}_a(\overline{\xi}) \cdot \nabla_{\overline{\xi}} \widehat{f}(t, \overline{\xi}) ~ d\overline{\xi} \label{estdispRHF-2} \\
&\quad \quad + t^{-2} 2^{-5j} \int e^{i t \Phi} \psi_j(\overline{\xi}) \widehat{h_a}(t, \overline{\xi}) ~ d\overline{\xi} \label{estdispRHF-3} \\
&\quad \quad + t^{-2} 2^{-4j} \int e^{i t \Phi} \psi_j(\overline{\xi}) \widehat{X}_a(\overline{\xi}) \cdot \nabla_{\overline{\xi}} \widehat{h_a}(t, \overline{\xi}) ~ d\overline{\xi} \label{estdispRHF-4} \\
&\quad \quad + t^{-1} 2^{-2j} \int e^{i t \Phi} \psi_j(\overline{\xi}) \widehat{g_a}(t, \overline{\xi}) ~ d\overline{\xi} \label{estdispRHF-5} 
\end{align}
\end{subequations}
Finally, we estimate using Hölder's inequality: 
\begin{align*}
\eqref{estdispRHF-1} &\lesssim t^{-2} 2^{-5j} \Vert \psi_j \Vert_{L^2} \Vert |\nabla|^{-1} f \Vert_{L^2} \\
&\lesssim t^{-2} 2^{-\frac{7j}{2}} \Vert u \Vert_X \\
&\lesssim t^{-\frac{7}{6}} 2^{-j} \langle 2^j \rangle^{-2} \Vert u \Vert_X \\
\eqref{estdispRHF-2} &\lesssim t^{-2} 2^{-5j} \Vert \psi_j \Vert_{L^2} \Vert X_a f \Vert_{L^2} \\
&\lesssim t^{-\frac{7}{6}} 2^{-j} \langle 2^j \rangle^{-2} \Vert u \Vert_X \\
\eqref{estdispRHF-4} &\lesssim t^{-2} 2^{-5j} \Vert \psi_j \Vert_{L^2} \Vert \nabla X_a h_a \Vert_{L^2} \\
&\lesssim t^{-\frac{7}{6}} 2^{-j} \langle 2^j \rangle^{-2} \Vert u \Vert_X \\
\eqref{estdispRHF-5} &\lesssim t^{-1} 2^{-\frac{5j}{2}} \langle 2^j \rangle^{-1} \Vert \psi_j \Vert_{L^2} \Vert \langle \overline{\xi} \rangle |\overline{\xi}|^{\frac{1}{2}} \widehat{g}_a(t) \Vert_{L^2} \\
&\lesssim t^{-\frac{7}{6}+100\delta} 2^{-j} \langle 2^j \rangle^{-1} \Vert u \Vert_X 
\end{align*}
\eqref{estdispRHF-3} is similar to \eqref{estdispRHF-2}. 

\paragraph{Low frequencies} Assume that $2^j \simeq \overline{\xi_a}$. In particular, since $2^j \gg t^{-\frac{1}{3}}$, we also have $\overline{\xi_a} \gg t^{-\frac{1}{3}}$ and $\mathfrak{t} \gg 1$. 

Let us define, for $l_a, l_b, l_c \in \mathbb{Z}$: 
\begin{align*}
\Psi_{(l_a, l_b, l_c)}^{a-b}\left( \overline{\xi}, \frac{x}{t}, \frac{y}{t} \right) &:= \psi\left( 2^{-l_a-2j} \widehat{X}_a(\overline{\xi}) \cdot \nabla_{\overline{\xi}} \Phi \right) \psi\left( 2^{-l_b-2j} m_b(\overline{\xi}) \widehat{X}_b(\overline{\xi}) \cdot \nabla_{\overline{\xi}} \Phi \right) \\
&\quad \chi\left( 2^{-l_c-2j} m_c(\overline{\xi}) \widehat{X}_c(\overline{\xi}) \cdot \nabla_{\overline{\xi}} \Phi \right) \\
\Psi_{(l_a, l_b, l_c)}^{c-b}\left( \overline{\xi}, \frac{x}{t}, \frac{y}{t} \right) &:= \chi\left( 2^{-l_a-2j} \widehat{X}_a(\overline{\xi}) \cdot \nabla_{\overline{\xi}} \Phi \right) \psi\left( 2^{-l_b-2j} m_b(\overline{\xi}) \widehat{X}_b(\overline{\xi}) \cdot \nabla_{\overline{\xi}} \Phi \right) \\
&\quad \psi\left( 2^{-l_c-2j} m_c(\overline{\xi}) \widehat{X}_c(\overline{\xi}) \cdot \nabla_{\overline{\xi}} \Phi \right) \\
\Psi_{(l_a, l_b, l_c)}^{a-int}\left( \overline{\xi}, \frac{x}{t}, \frac{y}{t} \right) &:= \psi\left( 2^{-l_a-2j} \widehat{X}_a(\overline{\xi}) \cdot \nabla_{\overline{\xi}} \Phi \right) \chi\left( 2^{-l_b-2j} m_b(\overline{\xi}) \widehat{X}_b(\overline{\xi}) \cdot \nabla_{\overline{\xi}} \Phi \right) \\
&\quad \chi\left( 2^{-l_c-2j} m_c(\overline{\xi}) \widehat{X}_c(\overline{\xi}) \cdot \nabla_{\overline{\xi}} \Phi \right) \stepcounter{equation}\tag{\theequation}\label{equdefPsigenerique} \\
\Psi_{(l_a, l_b, l_c)}^{c-int}\left( \overline{\xi}, \frac{x}{t}, \frac{y}{t} \right) &:= \chi\left( 2^{-l_a-2j} \widehat{X}_a(\overline{\xi}) \cdot \nabla_{\overline{\xi}} \Phi \right) \chi\left( 2^{-l_b-2j} m_b(\overline{\xi}) \widehat{X}_b(\overline{\xi}) \cdot \nabla_{\overline{\xi}} \Phi \right) \\
&\quad \psi\left( 2^{-l_c-2j} m_c(\overline{\xi}) \widehat{X}_c(\overline{\xi}) \cdot \nabla_{\overline{\xi}} \Phi \right) \\
\Psi_{(l_a, l_b, l_c)}^b\left( \overline{\xi}, \frac{x}{t}, \frac{y}{t} \right) &:= \chi\left( 2^{-l_a-2j} \widehat{X}_a(\overline{\xi}) \cdot \nabla_{\overline{\xi}} \Phi \right) \psi\left( 2^{-l_b-2j} m_b(\overline{\xi}) \widehat{X}_b(\overline{\xi}) \cdot \nabla_{\overline{\xi}} \Phi \right) \\
&\quad \chi\left( 2^{-l_c-2j} m_c(\overline{\xi}) \widehat{X}_c(\overline{\xi}) \cdot \nabla_{\overline{\xi}} \Phi \right) \\
\Psi_{(l_a, l_b, l_c)}^{int}\left( \overline{\xi}, \frac{x}{t}, \frac{y}{t} \right) &:= \chi\left( 2^{-l_a-2j} \widehat{X}_a(\overline{\xi}) \cdot \nabla_{\overline{\xi}} \Phi \right) \chi\left( 2^{-l_b-2j} m_b(\overline{\xi}) \widehat{X}_b(\overline{\xi}) \cdot \nabla_{\overline{\xi}} \Phi \right) \\
&\quad \chi\left( 2^{-l_c-2j} m_c(\overline{\xi}) \widehat{X}_c(\overline{\xi}) \cdot \nabla_{\overline{\xi}} \Phi \right) 
\end{align*}
We may then decompose
\begin{subequations} \label{estdispsupR} 
\begin{align}
&\int e^{i t \Phi} m_{\widehat{\mathcal{R}}}(\overline{\xi}) \psi_j(\overline{\xi}) \widehat{f}(t, \overline{\xi}) ~ d\overline{\xi} \notag \\
&\quad = \sum_{l_a = l_0}^{200} \sum_{l_b = l_a}^{200} \int e^{i t \Phi} \Psi_{(l_a, l_b, l_a)}^{a-b} m_{\widehat{\mathcal{R}}}(\overline{\xi}) \psi_j(\overline{\xi}) \widehat{f}(t, \overline{\xi}) ~ d\overline{\xi} \label{estdispsupR-ab} \\
&\quad \quad + \sum_{l_a = l_0}^{200} \sum_{l_b = l_a}^{200} \int e^{i t \Phi} \Psi_{(l_a, l_b, l_a)}^{c-b} m_{\widehat{\mathcal{R}}}(\overline{\xi}) \psi_j(\overline{\xi}) \widehat{f}(t, \overline{\xi}) ~ d\overline{\xi} \label{estdispsupR-cb} \\
&\quad \quad + \sum_{l_a = l_0}^{200} \int e^{i t \Phi} \Psi_{(l_a, l_a, l_a)}^{a-int} m_{\widehat{\mathcal{R}}}(\overline{\xi}) \psi_j(\overline{\xi}) \widehat{f}(t, \overline{\xi}) ~ d\overline{\xi} \label{estdispsupR-aint} \\
&\quad \quad + \sum_{l_a = l_0}^{200} \int e^{i t \Phi} \Psi_{(l_a, l_a, l_a)}^{c-int} m_{\widehat{\mathcal{R}}}(\overline{\xi}) \psi_j(\overline{\xi}) \widehat{f}(t, \overline{\xi}) ~ d\overline{\xi} \label{estdispsupR-cint} \\
&\quad \quad + \sum_{l_b = l_0}^{200} \int e^{i t \Phi} \Psi_{(l_0, l_b, l_0)}^b m_{\widehat{\mathcal{R}}}(\overline{\xi}) \psi_j(\overline{\xi}) \widehat{f}(t, \overline{\xi}) ~ d\overline{\xi} \label{estdispsupR-bint} \\
&\quad \quad + \int e^{i t \Phi} \Psi_{(l_0, l_0, l_0)}^{int} m_{\widehat{\mathcal{R}}}(\overline{\xi}) \psi_j(\overline{\xi}) \widehat{f}(t, \overline{\xi}) ~ d\overline{\xi} \label{estdispsupR-int} 
\end{align}
\end{subequations}
for $l_0$ a negative integer satisfying $2^{l_0} \simeq \mathfrak{t}^{-\frac{1}{2}}$. 

The point of this decomposition is the following: on the support of $\Psi_{(l_a, l_b, l_a)}^{a-b}$ or $\Psi_{(l_a, l_b, l_c)}^{c-b}$, $\widehat{X}_b \cdot \nabla \Phi$ is the largest, but we will only apply one integration by parts in this direction to avoid losses from the $\Vert \cdot \Vert_X$-norm, so we add a finer decomposition ($l_a \leq l_b$) in the directions $\widehat{X}_a, \widehat{X}_c$; on the support of $\Psi_{(l_a, l_b, l_a)}^{a-int}$ or $\Psi_{(l_a, l_b, l_a)}^{c-int}$, $\widehat{X}_b \cdot \nabla \Phi$ is smaller than $\widehat{X}_a \cdot \nabla \Phi$ or $\widehat{X}_c \cdot \nabla \Phi$ so it is enough to apply one or two integrations by parts in the dominant direction; finally, we stop the decomposition at scale $2^{l_0}$, under which volume estimates (i.e. Sobolev's embeddings) will be more interesting than integrations by parts. 

Let us estimate, for $l_a = l_c \leq l_b \leq 200$, the volume of the support of $\Psi_{(l_a, l_b, l_a)}^{*}$ for $*$ one of the possible exponents, that is $\Vert \Psi_{(l_a, l_b, l_a)}^{*} \Vert_{L^1_{\overline{\xi}}}$: on the support of $m_{\widehat{\mathcal{R}}}$, 
\begin{align*}
\left\{ \begin{array}{l}
\left| \widehat{X}_a(\overline{\xi}) \cdot \nabla_{\overline{\xi}} \Phi \right| \lesssim 2^{l_a+2j} \\[5pt]
\left| m_b(\overline{\xi}) \widehat{X}_b(\overline{\xi}) \cdot \nabla_{\overline{\xi}} \Phi \right| \lesssim 2^{l_b+2j} \\[5pt]
\left| m_c(\overline{\xi}) \widehat{X}_c(\overline{\xi}) \cdot \nabla_{\overline{\xi}} \Phi \right| \lesssim 2^{l_c+2j} 
\end{array} \right. 
~ & \Rightarrow ~ \left\{ 
\begin{array}{l}
\left| 3 \xi_0 |\overline{\xi}| - \widehat{X}_a(\overline{\xi}) \cdot \frac{(x, y)}{t} \right| \lesssim 2^{l_a+2j} \\[5pt]
\left| |\xi| |\overline{\xi}| - \widehat{X}_b(\overline{\xi}) \cdot \frac{(x, y)}{t} \right| \lesssim 2^{l_b+2j} \\[5pt]
\left| \widehat{X}_c(\overline{\xi}) \cdot \frac{(x, y)}{t} \right| \lesssim 2^{l_c+2j} 
\end{array} \right. 
\end{align*}
The application 
\begin{align*}
\Lambda : \overline{\xi} \mapsto \left( 3 \xi_0 |\overline{\xi}| - \widehat{X}_a(\overline{\xi}) \cdot \frac{(x, y)}{t}, 
|\xi| |\overline{\xi}| - \widehat{X}_b(\overline{\xi}) \cdot \frac{(x, y)}{t}, 
 - \widehat{X}_c(\overline{\xi}) \cdot \frac{(x, y)}{t} \right)
\end{align*}
has the following derivative: {\tiny
\begin{gather*}
D \Lambda (\overline{\xi}) = \begin{pmatrix} 3 |\overline{\xi}| + \frac{3 \xi_0^2}{|\overline{\xi}|} - \frac{x}{t|\overline{\xi}|} + \frac{\xi_0}{|\overline{\xi}|^2} \widehat{X}_a(\overline{\xi}) \cdot \frac{(x, y)}{t} 
& \frac{|\xi| \xi_0}{|\overline{\xi}|} + \frac{\xi \cdot y}{|\overline{\xi}| |\xi| t} + \frac{\xi_0}{|\overline{\xi}|^2} \widehat{X}_b(\overline{\xi}) \cdot \frac{(x, y)}{t} 
& 0 \\
3 \xi_0 \frac{\xi}{|\overline{\xi}|} - \frac{y}{t |\overline{\xi}|} + \frac{\xi}{|\overline{\xi}|^2} \widehat{X}_a(\overline{\xi}) \cdot \frac{(x, y)}{t} 
& \frac{\xi}{|\xi|} |\overline{\xi}| + |\xi| \frac{\xi}{|\overline{\xi}|} - \frac{x \xi}{t |\xi| |\overline{\xi}|} + \frac{\xi_0 y}{t |\xi| |\overline{\xi}|} - \frac{\xi_0 \xi \xi \cdot y}{|\xi|^3 |\overline{\xi}|t} + \frac{\xi}{|\overline{\xi}|^2} \widehat{X}_b(\overline{\xi}) \cdot \frac{(x, y)}{t} 
& \frac{J y}{|\xi| t} + \frac{\xi J \xi \cdot y}{t |\xi|^3} \end{pmatrix}
\end{gather*} }
In particular, one can compute that
\begin{align} 
\begin{aligned}
\mbox{det} D \Lambda(\overline{\xi}) 
&= 2 \frac{\xi \cdot y}{|\xi|^2 |\overline{\xi}|^2 t} \left( 3 |\overline{\xi}|^4 - \frac{x}{t} (|\xi|^2 + 3 \xi_0^2) - 2 \frac{\xi \cdot y}{|\xi| t} \xi_0 |\xi| \right) 
\end{aligned} \label{estdispsup-calculdeterminantbase} 
\end{align}

First, if $l_a \geq -200$, we may simply bound the volume by $2^{3j} \simeq 2^{2l_a+l_b+3j}$ (up to universal multiplicative constants). 

If $l_a \leq -200 \leq l_b$, then since $\widehat{X}_c \cdot \nabla \Phi$ is invariant by scaling of $\overline{\xi}$, while $\widehat{X}_a \cdot \nabla \Phi$ is not, we deduce easily that the volume can be bounded by $2^{2l_a+3j} \simeq 2^{2l_a+l_b+3j}$. 

Now, if $l_a \leq l_b \leq -200$, we have that
\begin{align*}
-\frac{x}{t} (|\xi|^2 + 3 \xi_0^2) - 2 \frac{\xi \cdot y}{|\xi| t} \xi_0 |\xi| 
&= - \frac{(x, y)}{t} \cdot \nabla_{\overline{\xi}} L \\
&= \left( \nabla_{\overline{\xi}} \Phi - \nabla_{\overline{\xi}} L \right) \cdot \nabla_{\overline{\xi}} L \\
&= O\left( 2^{l_b} |\overline{\xi_a}|^4 \right) - |\nabla_{\overline{\xi}} L|^2 \\
&= O\left( 2^{l_b} |\overline{\xi_a}|^4 \right) - (3 \xi_0^2 + |\xi|^2)^2 - 4 \xi_0^2 |\xi|^2 \\
&= O\left( 2^{l_b} |\overline{\xi_a}|^4 \right) - 9 \xi_0^4 - |\xi|^2 - 10 \xi_0^2 |\xi|^2 
\end{align*}
Hence, 
\begin{align*}
3 |\overline{\xi}|^4 - \frac{x}{t} (|\xi|^2 + 3 \xi_0^2) - 2 \frac{\xi \cdot y}{|\xi| t} \xi_0 |\xi|
&= -6 \xi_0^4 + 2 |\xi|^4 - 4 \xi_0^2 |\xi|^2 + O\left( 2^{\max(l_a, l_b, l_c)} |\overline{\xi_a}|^4 \right) \\
&= 2 (|\xi|^2 - 3 \xi_0^2) (\xi_0^2 + |\xi|^2) + O\left( 2^{\max(l_a, l_b, l_c)} |\overline{\xi_a}|^4 \right)
\end{align*}
But since $|\overline{\xi}| \simeq |\overline{\xi_a}|$, if $\overline{\xi}$ is in the support of $m_{\widehat{\mathcal{R}}}$, then 
\begin{align*}
\mbox{det} D \Lambda(\overline{\xi}) &\simeq 2 \frac{\xi \cdot y}{t}
\end{align*}
If $|y| \gtrsim |x|$, the only potential singularity occurs when $\xi$ is orthogonal to $y$. In this case, 
\begin{align*}
\widehat{X}_c(\overline{\xi}) \cdot \nabla_{\overline{\xi}} \Phi &= \frac{J \xi \cdot y}{|\xi| t} \simeq |\overline{\xi_a}|^2 
\end{align*}
which contradicts the hypothesis $l_a = l_c \leq -100$. 

On the other hand, if $|y| \ll |x|$, 
\begin{align*}
2^{2j} \gg \left| \widehat{X}_a(\overline{\xi}) \cdot \nabla_{\overline{\xi}} \Phi \right| &= \left| 3 \xi_0 |\overline{\xi}| - \frac{\xi_0}{|\overline{\xi}|} \frac{x}{t} + \frac{\xi \cdot y}{|\overline{\xi}| t} \right| \\
&\geq |\xi_0| \left| 3 |\overline{\xi}| - \frac{x}{t |\overline{\xi}|} \right| - \frac{|y|}{t} \\
2^{2j} \gg \left| \widehat{X}_b(\overline{\xi}) \cdot \nabla_{\overline{\xi}} \Phi \right| &= \left| |\xi| |\overline{\xi}| - \frac{|\xi| x}{|\overline{\xi}| t} + \frac{\xi_0 \xi \cdot y}{|\overline{\xi}| |\xi| t} \right| \\
&\geq |\xi| \left| |\overline{\xi}| - \frac{x}{t |\overline{\xi}|} \right| - \frac{|y|}{t}  
\end{align*}
But $\frac{|y|}{t} \ll 2^{2j}$, $|\xi| \simeq |\xi_0| \simeq 2^j$ and we cannot have at the same time $\left| |\overline{\xi}| - \frac{x}{t} \right| \ll |\overline{\xi_a}|$ and $\left| 3 |\overline{\xi}| - \frac{x}{t} \right| \ll |\overline{\xi_a}|$, so this would contradict $l_a, l_b, l_c \leq -200$. 

Therefore, $D \Lambda$ is never degenerate near resonant points, so the volume can be estimated easily. 

In conclusion, we may always bound
\begin{align*}
\Vert \Psi_{(l_a, l_b, l_a)}^{*} \Vert_{L^1} \lesssim 2^{2l_a + l_b+3j} 
\end{align*}

We now want to apply integrations by parts on the terms of \eqref{estdispsupR}. In the event when the derivation hits $\Psi_{(l_a, l_b, l_a)}^{*}$, we need estimates of $\widehat{X}_{\alpha} \cdot \nabla \Psi^{*}$. To that end, we compute: 
\begin{align*}
\widehat{X}_a(\overline{\xi}) \cdot \nabla_{\overline{\xi}} \left[ \widehat{X}_a(\overline{\xi}) \cdot \nabla_{\overline{\xi}} \Phi \right] &= \widehat{X}_a(\overline{\xi}) \cdot \nabla_{\overline{\xi}} \left[ 3 \xi_0 |\overline{\xi}| - \frac{\xi_0 x}{t |\overline{\xi}|} - \frac{\xi \cdot y}{t |\overline{\xi}|} \right] \\
&= \frac{\overline{\xi}}{|\overline{\xi}|} \cdot \begin{pmatrix} 3 |\overline{\xi}| + \frac{3 \xi_0^2}{|\overline{\xi}|} \\ \frac{3 \xi_0 \xi}{|\overline{\xi}|} \end{pmatrix} \\
&= 6 \xi_0 \\
\widehat{X}_a(\overline{\xi}) \cdot \nabla_{\overline{\xi}} \left[ \widehat{X}_b(\overline{\xi}) \cdot \nabla_{\overline{\xi}} \Phi \right] &= \widehat{X}_a(\overline{\xi}) \cdot \nabla_{\overline{\xi}} \left[ |\xi| |\overline{\xi}| - \frac{|\xi| x}{t |\overline{\xi}|} + \frac{\xi_0 \xi \cdot y}{t |\overline{\xi}| |\xi|} \right] \\
&= \frac{\overline{\xi}}{|\overline{\xi}|} \cdot \begin{pmatrix} \frac{\xi_0 |\xi|}{|\overline{\xi}|} \\ \frac{\xi |\overline{\xi}|}{|\xi|} + \frac{|\xi| \xi}{|\overline{\xi}|} \end{pmatrix} \\
&= 2 |\xi| \\
\widehat{X}_a(\overline{\xi}) \cdot \nabla_{\overline{\xi}} \left[ \widehat{X}_c(\overline{\xi}) \cdot \nabla_{\overline{\xi}} \Phi \right] &= \widehat{X}_a(\overline{\xi}) \cdot \nabla_{\overline{\xi}} \left[ - \frac{J \xi \cdot y}{|\xi| t} \right] \\
&= 0 \\
\widehat{X}_b(\overline{\xi}) \cdot \nabla_{\overline{\xi}} \left[ \widehat{X}_a(\overline{\xi}) \cdot \nabla_{\overline{\xi}} \Phi \right] &= \widehat{X}_b(\overline{\xi}) \cdot \nabla_{\overline{\xi}} \left[ 3 \xi_0 |\overline{\xi}| - \frac{\xi_0 x}{t |\overline{\xi}|} - \frac{\xi \cdot y}{t |\overline{\xi}|} \right] \\
&= \frac{|\xi|}{|\overline{\xi}|} \left( 3 |\overline{\xi}| - \frac{x}{t |\overline{\xi}|} \right) - \frac{\xi_0 \xi}{|\overline{\xi}| |\xi|} \cdot \left( - \frac{y}{t |\overline{\xi}|} \right) \\
&= \frac{1}{|\overline{\xi}|} \left( 3 |\xi| |\overline{\xi}| - \frac{|\xi| x}{t |\overline{\xi}|} + \frac{\xi_0 \xi \cdot y}{t |\overline{\xi}| |\xi|} \right) \\
\widehat{X}_b(\overline{\xi}) \cdot \nabla_{\overline{\xi}} \left[ \widehat{X}_b(\overline{\xi}) \cdot \nabla_{\overline{\xi}} \Phi \right] &= \widehat{X}_b(\overline{\xi}) \cdot \nabla_{\overline{\xi}} \left[ |\xi| |\overline{\xi}| - \frac{|\xi| x}{t |\overline{\xi}|} + \frac{\xi_0 \xi \cdot y}{t |\overline{\xi}| |\xi|} \right] \\
&= \frac{|\xi|}{|\overline{\xi}|} \left( \frac{\xi \cdot y}{t |\overline{\xi}| |\xi|} \right) - \frac{\xi_0 \xi}{|\overline{\xi}| |\xi|} \cdot \left( \frac{\xi |\overline{\xi}|}{|\xi|} - \frac{\xi x}{t |\overline{\xi}| |\xi|} \right) \\
&= - \frac{1}{|\overline{\xi}|} \left( \xi_0 |\overline{\xi}| - \frac{\xi_0 x}{t |\overline{\xi}|} - \frac{\xi \cdot y}{t |\overline{\xi}|} \right) \\
\widehat{X}_b(\overline{\xi}) \cdot \nabla_{\overline{\xi}} \left[ \widehat{X}_c(\overline{\xi}) \cdot \nabla_{\overline{\xi}} \Phi \right] &= \widehat{X}_b(\overline{\xi}) \cdot \left[ - \frac{J \xi \cdot y}{|\xi| t} \right] \\
&= 0 \\
\widehat{X}_c(\overline{\xi}) \cdot \nabla_{\overline{\xi}} \left[ \widehat{X}_a(\overline{\xi}) \cdot \nabla_{\overline{\xi}} \Phi \right] &= \frac{J \xi}{|\xi|} \cdot \nabla_{\xi} \left[ 3 \xi_0 |\overline{\xi}| - \frac{\xi_0 x}{t |\overline{\xi}|} - \frac{\xi \cdot y}{t |\overline{\xi}|} \right] \\
&= - \frac{J \xi \cdot y}{t |\xi| |\overline{\xi}|} \\
&= \frac{1}{|\overline{\xi}|} \widehat{X}_c(\overline{\xi}) \cdot \nabla_{\overline{\xi}} \Phi \\
\widehat{X}_c(\overline{\xi}) \cdot \nabla_{\overline{\xi}} \left[ \widehat{X}_b(\overline{\xi}) \cdot \nabla_{\overline{\xi}} \Phi \right] &= \frac{J \xi}{|\xi|} \cdot \nabla_{\xi} \left[ |\xi| |\overline{\xi}| - \frac{|\xi| x}{t |\overline{\xi}|} + \frac{\xi_0 \xi \cdot y}{t |\overline{\xi}| |\xi|} \right] \\
&= \frac{\xi_0 J \xi \cdot y}{t |\overline{\xi}| |\xi|^2} \\
&= - \frac{\xi_0}{|\overline{\xi}| |\xi|} \widehat{X}_c(\overline{\xi}) \cdot \nabla_{\overline{\xi}} \Phi \\
\widehat{X}_c(\overline{\xi}) \cdot \nabla_{\overline{\xi}} \cdot \left[ \widehat{X}_c(\overline{\xi}) \cdot \nabla_{\overline{\xi}} \Phi \right] &= \frac{J \xi}{|\xi|} \cdot \nabla_{\xi} \left[ - \frac{J \xi \cdot y}{|\xi| t} \right] \\
&= \frac{J \xi \cdot Jy}{|\xi|^2 t} = \frac{\xi \cdot y}{|\xi|^2 t} 
\end{align*}

Therefore, when $l_a \leq l_b$, we have 
\begin{align*}
\left| \widehat{X}_a(\overline{\xi}) \cdot \nabla_{\overline{\xi}} \Psi_{(l_a, l_b, l_a)}^{*} \right| &\lesssim 2^{-l_a-j} \left| \Psi_{(l_a+1, l_b+1, l_a+1)}^{int} \right| \\
\left| \widehat{X}_b(\overline{\xi}) \cdot \nabla_{\overline{\xi}} \Psi_{(l_a, l_b, l_a)}^{*} \right| &\lesssim 2^{-l_a-j} \left| \Psi_{(l_a+1, l_b+1, l_a+1)}^{int} \right| \\
\left| \widehat{X}_c(\overline{\xi}) \cdot \nabla_{\overline{\xi}} \Psi_{(l_a, l_b, l_a)}^{*} \right| &\lesssim 2^{-l_a-j} \left| \Psi_{(l_a+1, l_b+1, l_a+1)}^{int} \right| 
\end{align*}
For $\widehat{X}_b$, if $l_b \leq l_a+100$, the estimate above implies
\begin{align*}
\left| \widehat{X}_b(\overline{\xi}) \cdot \nabla_{\overline{\xi}} \Psi_{(l_a, l_b, l_a)}^{*} \right| &\lesssim 2^{-l_b-j} \left| \Psi_{(l_a+1, l_b+1, l_a+1)}^{int} \right|
\end{align*}
up to changing the implicit constant. If $l_b \geq l_a+100$, we set 
\begin{align*}
\widehat{X}_{b-corr} := \widehat{X}_b - \frac{\widehat{X}_b(\overline{\xi}) \cdot \nabla_{\overline{\xi}} \Phi}{\widehat{X}_a(\overline{\xi}) \cdot \nabla_{\overline{\xi}} \Phi} \widehat{X}_a
\end{align*}
This vector field satisfies
\begin{align*}
\widehat{X}_{b-corr}(\overline{\xi}) \cdot \nabla_{\overline{\xi}} \left[ \widehat{X}_a(\overline{\xi}) \cdot \nabla_{\overline{\xi}} \Phi \right] &= 0 \\
\widehat{X}_{b-corr}(\overline{\xi}) \cdot \nabla_{\overline{\xi}} \left[ \widehat{X}_c(\overline{\xi}) \cdot \nabla_{\overline{\xi}} \Phi \right] &= 0 
\end{align*}
while, on the support of $m_{\widehat{\mathcal{R}}} \Psi_{(l_a, l_b, l_a)}^{a-b}$ or $m_{\widehat{\mathcal{R}}} \Psi_{(l_a, l_b, l_a)}^{c-b}$, 
\begin{align*}
\left| \widehat{X}_{b-corr}(\overline{\xi}) \cdot \nabla_{\overline{\xi}} \Phi \right| &\geq \left| \widehat{X}_b(\overline{\xi}) \cdot \nabla_{\overline{\xi}} \Phi \right| - \left\Vert \frac{\widehat{X}_b(\overline{\xi}) \cdot \nabla_{\overline{\xi}} \Phi}{\widehat{X}_a(\overline{\xi}) \cdot \nabla_{\overline{\xi}} \Phi} \right\Vert_{L^{\infty}} \left| \widehat{X}_a(\overline{\xi}) \cdot \nabla_{\overline{\xi}} \Phi \right| \\
&\gtrsim 2^{l_b+2j}
\end{align*}
since $\left\Vert \frac{\widehat{X}_b(\overline{\xi}) \cdot \nabla_{\overline{\xi}} \Phi}{\widehat{X}_a(\overline{\xi}) \cdot \nabla_{\overline{\xi}} \Phi} \right\Vert_{L^{\infty}} \lesssim 1$ and we supposed here $l_b \geq l_a+100$. We therefore obtain
\begin{align*}
\left| \widehat{X}_{b-corr}(\overline{\xi}) \cdot \nabla_{\overline{\xi}} \Psi_{(l_a, l_b, l_a)}^{*} \right| \lesssim 2^{-l_b-j} \left| \Psi_{(l_a, l_b, l_a)}^{int} \right| 
\end{align*}

In what follows, we write $\widehat{X}_{b-corr}$ for any values of $(l_a, l_b)$, setting $\widehat{X}_{b-corr} = \widehat{X}_b$ if $l_b \leq l_a+100$.  

On \eqref{estdispsupR-ab}, we first apply an integration by parts along $\widehat{X}_{b-corr}$: 
\begin{align*}
\eqref{estdispsupR-ab} &= \sum_{l_b \geq l_a = l_0}^{200} \int e^{i t \Phi} \Psi_{(l_a, l_b, l_a)}^{a-b} m_{\widehat{\mathcal{R}}}(\overline{\xi}) \psi_j(\overline{\xi}) \widehat{f}(t, \overline{\xi}) ~ d\overline{\xi} \\
&= \sum_{l_b \geq l_a = l_0}^{200} \int t^{-1} e^{i t \Phi} \nabla_{\overline{\xi}} \cdot \left( \frac{\widehat{X}_{b-corr}(\overline{\xi})}{\widehat{X}_{b-corr}(\overline{\xi}) \cdot \nabla_{\overline{\xi}} \Phi} \Psi_{(l_a, l_b, l_a)}^{a-b} m_{\widehat{\mathcal{R}}}(\overline{\xi}) \psi_j(\overline{\xi}) \widehat{f}(t, \overline{\xi}) \right) ~ d\overline{\xi} \\
&= \sum_{l_b \geq l_a = l_0}^{200} \int t^{-1} e^{i t \Phi} \nabla_{\overline{\xi}} \cdot \left( \frac{\widehat{X}_{b-corr}(\overline{\xi})}{\widehat{X}_{b-corr}(\overline{\xi}) \cdot \nabla_{\overline{\xi}} \Phi} \Psi_{(l_a, l_b, l_a)}^{a-b} m_{\widehat{\mathcal{R}}}(\overline{\xi}) \psi_j(\overline{\xi}) \right) \widehat{f}(t, \overline{\xi}) ~ d\overline{\xi} \\
&+ \sum_{l_b \geq l_a = l_0}^{200} \int t^{-1} e^{i t \Phi} \frac{1}{\widehat{X}_a(\overline{\xi}) \cdot \nabla_{\overline{\xi}} \Phi} \Psi_{(l_a, l_b, l_a)}^{a-b} m_{\widehat{\mathcal{R}}}(\overline{\xi}) \psi_j(\overline{\xi}) \widehat{X}_{b-corr}(\overline{\xi}) \cdot \nabla_{\overline{\xi}} \widehat{f}(t, \overline{\xi}) ~ d\overline{\xi} 
\end{align*}
We may then decompose $\widehat{X}_{b-corr} \cdot \nabla \widehat{f}(t) = \widehat{h}_{b-corr}(t) + \widehat{g}_{b-corr}(t)$, by projecting $\widehat{X}_{b-corr}$ on the basis $\left( \widehat{X}_a, \widehat{X}_b, \widehat{X}_c \right)$ and by using the decomposition for each of these canonical vector fields, and we apply a second integration by parts except when there is $\widehat{g}$. More precisely, as $\left( \widehat{X}_a, \widehat{X}_b, \widehat{X}_c \right)$ is an orthonormal basis, we can write 
\begin{align*}
\widehat{X}_{b-corr} = \sum_{\alpha = a, b, c} \left( \widehat{X}_{b-corr} \cdot \widehat{X}_{\alpha} \right) \widehat{X}_{\alpha}
\end{align*}
and then decompose linearly 
\begin{align*}
\widehat{h}_{b-corr}(t) &= \sum_{\alpha = a, b, c} \left( \widehat{X}_{b-corr} \cdot \widehat{X}_{\alpha} \right) \widehat{h}_{\alpha}(t) \\
\widehat{g}_{b-corr}(t) &= \sum_{\alpha = a, b, c} \left( \widehat{X}_{b-corr} \cdot \widehat{X}_{\alpha} \right) \widehat{g}_{\alpha}(t)
\end{align*}
Note that the symbols $\left( \widehat{X}_{b-corr} \cdot \widehat{X}_{\alpha} \right)$ are of Hörmander-Mikhlin type. We get
\begin{subequations}
\begin{align}
\eqref{estdispsupR-ab} &= \sum_{l_b \geq l_a = l_0}^{200} t^{-2} 2^{-2l_b-2l_a-6j} \int e^{i t \Phi} \Psi_{(l_a, l_b, l_a)}^{a-b} \psi_j(\overline{\xi}) \widehat{f}(t, \overline{\xi}) d\overline{\xi} \label{estdispsupR-ab-1} \\
&\quad + \sum_{l_b \geq l_a = l_0}^{200} t^{-2} 2^{-2l_b-l_a-5j} \int e^{i t \Phi} \Psi_{(l_a, l_b, l_a)}^{a-b} \psi_j(\overline{\xi}) \widehat{X}_a(\overline{\xi}) \cdot \nabla_{\overline{\xi}} \widehat{f}(t, \overline{\xi}) d\overline{\xi} \label{estdispsupR-ab-2} \\
&\quad + \sum_{l_b \geq l_a = l_0}^{200} t^{-2} 2^{-l_b-2l_a-5j} \int e^{i t \Phi} \Psi_{(l_a, l_b, l_a)}^{a-b} \psi_j(\overline{\xi}) \widehat{h}_{b-corr}(t, \overline{\xi}) d\overline{\xi} \label{estdispsupR-ab-3} \\
&\quad + \sum_{l_b \geq l_a = l_0}^{200} t^{-2} 2^{-l_b-l_a-4j} \int e^{i t \Phi} \Psi_{(l_a, l_b, l_a)}^{a-b} \psi_j(\overline{\xi}) \widehat{X}_a(\overline{\xi}) \cdot \nabla_{\overline{\xi}} \widehat{h}_{b-corr}(t, \overline{\xi}) d\overline{\xi} \label{estdispsupR-ab-4} \\
&\quad + \sum_{l_b \geq l_a = l_0}^{200} t^{-1} 2^{-l_b-2j} \int e^{i t \Phi} \Psi_{(l_a, l_b, l_a)}^{a-b} \psi_j(\overline{\xi}) \widehat{g}_{b-corr}(t, \overline{\xi}) d\overline{\xi} \label{estdispsupR-ab-5} 
\end{align}
\end{subequations} 
where the symbols $\Psi_{(l_a, l_b, l_a)}, \psi_j$ may vary from line to line as long as they keep similar properties. We now estimate: 
\begin{align*}
\eqref{estdispsupR-ab-1} &\lesssim \sum_{l_b \geq l_a = l_0}^{200} t^{-2} 2^{-2l_b-2l_a-6j} \Vert \Psi_{(l_a, l_b, l_a)}^{a-b} \Vert_{L^{\frac{6}{5}}} \langle 2^j \rangle^{-1} \Vert \langle \overline{\xi} \rangle m_{\widehat{\mathcal{R}}} \widehat{f}(t) \Vert_{L^6} \\
&\lesssim \sum_{l_b \geq l_a = l_0}^{200} t^{-2} 2^{-\frac{7l_b}{6}-\frac{l_a}{3}-\frac{7j}{2}} \langle 2^j \rangle^{-1} \Vert \langle \overline{\xi} \rangle m_{\widehat{\mathcal{R}}} \widehat{f}(t) \Vert_{\dot{H}^1} \\
&\lesssim t^{-\frac{7}{6}} \mathfrak{t}^{-\frac{5}{6}} 2^{-\frac{3l_0}{2}-j} \langle 2^j \rangle^{-1} \Vert u \Vert_X \\
&\lesssim t^{-\frac{7}{6}} \mathfrak{t}^{-\frac{1}{12}} 2^{-j} \langle 2^j \rangle^{-1} \Vert u \Vert_X \\
&\lesssim t^{-\frac{7}{6}} 2^{-j} \langle 2^j \rangle^{-1} \Vert u \Vert_X \\
\eqref{estdispsupR-ab-2} &\lesssim \sum_{l_b \geq l_a = l_0}^{200} t^{-2} 2^{-2l_b-l_a-5j} \Vert \Psi_{(l_a, l_b, l_a)}^{a-b} \Vert_{L^2} \langle 2^j \rangle^{-1} \Vert \langle \nabla \rangle X_a f(t) \Vert_{L^2} \\
&\lesssim \sum_{l_b \geq l_a = l_0}^{200} t^{-2} 2^{-\frac{3l_b}{2}-\frac{7j}{2}} \langle 2^j \rangle^{-1} \Vert u \Vert_X \\
&\lesssim t^{-\frac{7}{6}} 2^{-j} \langle 2^j \rangle^{-1} \Vert u \Vert_X \\
\eqref{estdispsupR-ab-4} &\lesssim \sum_{l_b \geq l_a = l_0}^{200} t^{-2} 2^{-l_b-l_a-5j} \Vert \Psi_{(l_a, l_b, l_a)}^{a-b} \Vert_{L^2} \Vert m_{\widehat{\mathcal{R}}}(D) \nabla X_a h_{b-corr}(t) \Vert_{L^2} \\
&\lesssim \sum_{l_b \geq l_a = l_0}^{200} t^{-\frac{7}{6}} \mathfrak{t}^{-\frac{5}{6}} 2^{-\frac{l_b}{2}-j} \Vert u \Vert_X \\
&\lesssim t^{-\frac{7}{6}} \mathfrak{t}^{-\frac{1}{3}} 2^{-j} \Vert u \Vert_X \\
&\lesssim t^{-\frac{7}{6}} 2^{-j} \langle 2^j \rangle^{-1} \Vert u \Vert_X \\
\eqref{estdispsupR-ab-5} &\lesssim \sum_{l_b \geq l_a = l_0}^{200} t^{-1} 2^{-l_b-\frac{5j}{2}} \Vert \Psi_{(l_a, l_b, l_a)}^{a-b} \Vert_{L^2} \langle 2^j \rangle^{-1} \Vert m_{\widehat{\mathcal{R}}}(D) |\nabla|^{\frac{1}{2}} \langle \nabla \rangle g_{b-corr}(t) \Vert_{L^2} \\
&\lesssim \sum_{l_b \geq l_a = l_0}^{200} t^{-\frac{7}{6}+100\delta} 2^{l_a-\frac{l_b}{2}-j} \langle 2^j \rangle^{-1} \Vert u \Vert_X \\
&\lesssim t^{-\frac{7}{6}+100\delta} 2^{-j} \langle 2^j \rangle^{-1} \Vert u \Vert_X 
\end{align*}
\eqref{estdispsupR-ab-3} is similar to \eqref{estdispsupR-ab-2}. 

We can proceed the same way on \eqref{estdispsupR-aint}, replacing the integration by parts along $\widehat{X}_{b-corr}$ by an integration by parts along $\widehat{X}_a$; and likewise \eqref{estdispsupR-cb} and \eqref{estdispsupR-cint} replacing $\widehat{X}_a$ by $\widehat{X}_c$. 

For \eqref{estdispsupR-bint}, we apply a single integration by parts in the direction $\widehat{X}_{b-corr}$: 
\begin{subequations}
\begin{align}
\eqref{estdispsupR-bint} &= \sum_{l_b = l_0}^{200} t^{-1} 2^{-2l_b-3j} \int e^{i t \Phi} \Psi_{(l_0, l_b, l_0)}^b \psi_j(\overline{\xi}) \widehat{f}(t, \overline{\xi}) ~ d\overline{\xi} \label{estdispsupR-bint-1} \\
&\quad + \sum_{l_b = l_0}^{200} t^{-1} 2^{-l_b-2j} \int e^{i t \Phi} \Psi_{(l_0, l_b, l_0)}^b \psi_j(\overline{\xi}) \widehat{X}_{b-corr}(\overline{\xi}) \cdot \nabla_{\overline{\xi}} \widehat{f}(t, \overline{\xi}) ~ d\overline{\xi} \label{estdispsupR-bint-2} 
\end{align}
\end{subequations}
We then estimate: 
\begin{align*}
\eqref{estdispsupR-bint-1} &\lesssim \sum_{l_b = l_0}^{200} t^{-1} 2^{-2l_b-3j} \Vert \Psi_{(l_0, l_b, l_0)} \Vert_{L^{\frac{6}{5}}} \langle 2^j \rangle^{-1} \Vert \langle \overline{\xi} \rangle m_{\widehat{\mathcal{R}}}(\overline{\xi}) \widehat{f}(t) \Vert_{L^6} \\
&\lesssim \sum_{l_b = l_0}^{200} t^{-\frac{7}{6}} \mathfrak{t}^{\frac{1}{6}} 2^{-\frac{7l_b}{6}+\frac{5l_0}{3}-j} \langle 2^j \rangle^{-1} \Vert u \Vert_X \\ 
&\lesssim t^{-\frac{7}{6}} \mathfrak{t}^{-\frac{1}{12}} 2^{-j} \langle 2^j \rangle^{-1} \Vert u \Vert_X \\
&\lesssim t^{-\frac{7}{6}} 2^{-j} \langle 2^j \rangle^{-1} \Vert u \Vert_X \\
\eqref{estdispsupR-bint-2} &\lesssim \sum_{l_b = l_0}^{200} t^{-1} 2^{-l_b-2j} \Vert \Psi_{(l_0, l_b, l_0)} \Vert_{L^2} \langle 2^j \rangle^{-1} \Vert \langle \nabla \rangle m_{\widehat{\mathcal{R}}}(D) X_{b-corr} f(t) \Vert_{L^2} \\
&\lesssim \sum_{l_b = l_0}^{200} t^{-\frac{7}{6}} \mathfrak{t}^{\frac{1}{6}} 2^{-\frac{l_b}{2}+l_0-j} \langle 2^j \rangle^{-1} \Vert u \Vert_X \\
&\lesssim t^{-\frac{7}{6}} \mathfrak{t}^{-\frac{1}{12}} 2^{-j} \langle 2^j \rangle^{-1} \Vert u \Vert_X \\
&\lesssim t^{-\frac{7}{6}} 2^{-j} \langle 2^j \rangle^{-1} \Vert u \Vert_X
\end{align*}

Finally, for \eqref{estdispsupR-int}, we apply directly without any integration by parts: 
\begin{align*}
\eqref{estdispsupR-int} &\lesssim \Vert \Psi_{(l_0, l_0, l_0)}^{int} \Vert_{L^{\frac{6}{5}}} \langle 2^j \rangle^{-1} \Vert \langle \overline{\xi} \rangle m_{\widehat{\mathcal{R}}} \widehat{f} \Vert_{L^6} \\
&\lesssim 2^{\frac{5l_0}{2}+\frac{5j}{2}} \langle 2^j \rangle^{-1} \Vert u \Vert_X \\
&\lesssim t^{-\frac{7}{6}} \mathfrak{t}^{-\frac{1}{12}} 2^{-j} \langle 2^j \rangle^{-1} \Vert u \Vert_X \\
&\lesssim t^{-\frac{7}{6}} 2^{-j} \langle 2^j \rangle^{-1} \Vert u \Vert_X 
\end{align*}

\paragraph{Low frequencies} Assume that $2^j \ll \overline{\xi_a}$. 

Then since $(\widehat{X}_a, \widehat{X}_b, \widehat{X}_c)$ is an orthonormal family everywhere, we deduce that $\widehat{X}_{\alpha} \cdot \nabla \Phi \simeq |\overline{\xi_a}|^2$ for at least one $\alpha$. Up to decomposing the frequency space into three parts, we can assume this $\alpha$ is fixed. If $\alpha = a$ or $\alpha = c$, there is no difficulty applying the same estimates as in the high frequency regime. On the other hand, if $\alpha = b$, we may apply a simpler version of the decomposition used in the resonant case above $2^j \simeq \overline{\xi_a}$, localising only in the directions $\widehat{X}_a$ and $\widehat{X}_c$. We skip the details. 

This concludes the proof of Lemma \ref{lem-estdisph-zonereste}. 
\end{Dem}

\subsection{Neighborhood of the cone}

\begin{Lem} Assume $\delta > 0$ is small enough. There exists $C > 0$ such that, for every $t \geq 1$ and every $j, k \in \mathbb{Z}$ satisfying $2^j \gg t^{-\frac{1}{3}}$, $k \leq -10$, for every solution $u$ of the equation, 
\begin{align*}
\Vert e^{i t\omega(D)} \psi_{j, k}^{\widehat{\mathcal{C}}}(D) f(t) \Vert_{L^{\infty}} &\lesssim t^{-\frac{13}{12}+100\delta} 2^{-j+\delta j+\delta k} \langle 2^j \rangle^{-\frac{1}{2}} \Vert u \Vert_X \\
\Vert e^{it\omega(D)} \psi_{j, k}^{\widehat{\mathcal{C}}}(D) f(t) \Vert_{L^{\infty}} &\lesssim t^{-\frac{7}{6}+100\delta} 2^{-\frac{3j}{2}+\delta j+\delta k} \langle 2^j \rangle^{-\frac{1}{4}} \Vert u \Vert_X 
\end{align*} \label{lemestdispvoiscone} 
\end{Lem}

\begin{Dem}
As before, fixing $(x, y, t)$ we need to bound
\begin{align}
\int e^{i t \Phi} \psi_{j, k}^{\widehat{\mathcal{C}}}(\overline{\xi}) \widehat{f}(t, \overline{\xi}) d\overline{\xi} \label{estdispsupC-termeinit}
\end{align} 

\paragraph{High frequencies} Assume first $2^j \gg |\overline{\xi_a}|$, then on the support of the integral
\begin{align*}
\widehat{X}_a(\overline{\xi}) \cdot \nabla_{\overline{\xi}} \Phi \simeq 2^{2j}
\end{align*}
We can therefore apply the same estimates as in the case $\widehat{\mathcal{R}}$. 

In a similar fashion, if $2^j \lesssim |\overline{\xi_a}|$, as long as
\begin{align*}
\widehat{X}_a(\overline{\xi}) \cdot \nabla_{\overline{\xi}} \Phi \gtrsim 2^{2j} \quad \mbox{ or } \quad \widehat{X}_c(\overline{\xi}) \cdot \nabla_{\overline{\xi}} \Phi \gtrsim 2^{2j}
\end{align*}
we may proceed the same way and get a similar bound. 

We first consider the case $|\overline{\xi_a}| \simeq 2^j$ and it is enough to get a bound on
\begin{align}
\int e^{i t \Phi} \psi_{j, k}^{\widehat{\mathcal{C}}}(\overline{\xi}) \chi\left( 2^{100} \widehat{X}_a(\overline{\xi}) \cdot \nabla_{\overline{\xi}} \Phi \right) \chi\left( 2^{100} \widehat{X}_c(\overline{\xi}) \cdot \nabla_{\overline{\xi}} \Phi \right) \widehat{f}(t, \overline{\xi}) d\overline{\xi} \label{estdispsupC-termeinitbis}
\end{align}

\paragraph{Away from $\mathcal{L}$, close to $\mathcal{C}$} Let us start by assuming that $|y| \gtrsim |x|$, that is $(x, y)$ is away enough from $\mathcal{L}$, and even $\tau := \sqrt{\frac{\left| |x| - \sqrt{3}|y| \right|}{|x|+|y|}} \ll 2^k$. 

We define the corrected vector field 
\begin{align*}
\widehat{X}_{b-corr}(\overline{\xi}) \cdot \nabla_{\overline{\xi}} := \partial_{\xi_0} - \frac{\sqrt{3} \xi_0 \xi}{|\xi_0| |\xi|} \cdot \nabla_{\xi} 
\end{align*}
On the support of $m_{\widehat{\mathcal{C}}}$, $\widehat{X}_{b-corr}(\overline{\xi})$ is well-defined and completes $(\widehat{X}_a(\overline{\xi}), \widehat{X}_c(\overline{\xi}))$ into a basis. Furthermore, 
\begin{align*}
\widehat{X}_{b-corr}(\overline{\xi}) \cdot \nabla_{\overline{\xi}} \Phi 
&= \left( 3 \xi_0^2 + |\xi|^2 - \frac{x}{t} \right) - \frac{\sqrt{3} \xi_0 \xi}{|\xi_0| |\xi|} \cdot \left( 2 \xi_0 \xi - \frac{y}{t} \right) \\
&= 3 \xi_0^2 + |\xi|^2 - 2 \sqrt{3} |\xi_0| |\xi| - \frac{x}{t} + \frac{\sqrt{3} \xi_0 \xi \cdot y}{t |\xi_0| |\xi|} \\
&= \left( \sqrt{3} |\xi_0| - |\xi| \right)^2 - \frac{x}{t} + \frac{\sqrt{3} \xi_0 \xi \cdot y}{t |\xi_0| |\xi|} 
\end{align*}
hence $\widehat{X}_{b-corr}(\overline{\xi}) \cdot \nabla_{\overline{\xi}} \Phi$ is quadratically degenerate but controls exactly the distance to the cone $\widehat{\mathcal{C}}$. 

Moreover, we also have 
\begin{align*}
\widehat{X}_a(\overline{\xi}) \cdot \nabla_{\overline{\xi}} \Phi &= 3 \xi_0 |\overline{\xi}| - \frac{x \xi_0}{t |\overline{\xi}|} - \frac{y \cdot \xi}{t |\overline{\xi}|} \\
\widehat{X}_c(\overline{\xi}) \cdot \nabla_{\overline{\xi}} \Phi &= - \frac{J \xi \cdot y}{t |\xi|} 
\end{align*}
In particular, since $k \leq -10$, $\sqrt{3} |\xi_0| - |\xi| \ll |\overline{\xi}|$, the localisation also ensures that $\xi$ and $y$ are close to alignment. Therefore, the integral in $\overline{\xi}$ can be separated into two disjoint supports: either $x \xi_0$ and $y \cdot \xi$ have the same sign, or they have opposite sign. If the sign is opposite, then 
\begin{align*}
\left| \frac{x \xi_0}{t |\overline{\xi}|} + \frac{y \cdot \xi}{t |\overline{\xi}|} \right| &\lesssim \left| \widehat{X}_c(\overline{\xi}) \cdot \nabla_{\overline{\xi}} \Phi \right| + \left| \sqrt{3} |\xi_0| - |\xi| \right| + \frac{1}{|\overline{\xi}|} \left| \frac{|x| - \sqrt{3} |y|}{t} \right| \\
&\ll |\overline{\xi}|  
\end{align*}
In particular, we cannot have $\widehat{X}_a(\overline{\xi}) \cdot \nabla_{\overline{\xi}} \Phi$ small with respect to $2^{2j}$ in this case, and we deduce that the localisation ensures that $x \xi_0$ and $y \cdot \xi$ have the same sign. Furthermore, we can see that $x > 0$, otherwise the localisation makes the integral identically zero as well. 

We now consider the localisation symbols
\begin{align*}
\Psi_{l_a}^a\left( \overline{\xi}, \frac{x}{t}, \frac{y}{t} \right) &:= \psi\left( 2^{-l_a-2j} \widehat{X}_a(\overline{\xi}) \cdot \nabla_{\overline{\xi}} \Phi \right) \chi\left( 2^{-l_a-2j} \widehat{X}_c(\overline{\xi}) \cdot \nabla_{\overline{\xi}} \Phi \right) \\
\Psi_{l_a}^c\left( \overline{\xi}, \frac{x}{t}, \frac{y}{t} \right) &:= \chi\left( 2^{-l_a-2j} \widehat{X}_a(\overline{\xi}) \cdot \nabla_{\overline{\xi}} \Phi \right) \psi\left( 2^{-l_a-2j} \widehat{X}_c(\overline{\xi}) \cdot \nabla_{\overline{\xi}} \Phi \right) \\
\Psi_{l_a}^{int}\left( \overline{\xi}, \frac{x}{t}, \frac{y}{t} \right) &:= \chi\left( 2^{-l_a-2j} \widehat{X}_a(\overline{\xi}) \cdot \nabla_{\overline{\xi}} \Phi \right) \chi\left( 2^{-l_a-2j} \widehat{X}_c(\overline{\xi}) \cdot \nabla_{\overline{\xi}} \Phi \right)
\end{align*}
Now we separate
\begin{subequations}
\begin{align}
\eqref{estdispsupC-termeinitbis} &= \sum_{l_a = l_0}^{-100} 1_{3k - j_{+} \leq 2l_a} \int e^{i t \Phi} \psi_{j, k}^{\widehat{\mathcal{C}}}(\overline{\xi}) \Psi_{l_a}^a \widehat{f}(t, \overline{\xi}) d\overline{\xi} \label{estdispsupCloinLpetittau-a} \\
&\quad + \sum_{l_a = l_0}^{-100} 1_{2 l_a < 3k - j_{+} \leq l_a} \int e^{i t \Phi} \psi_{j, k}^{\widehat{\mathcal{C}}}(\overline{\xi}) \Psi_{l_a}^a \widehat{f}(t, \overline{\xi}) d\overline{\xi} \label{estdispsupCloinLpetittau-ab1} \\
&\quad + \sum_{l_a = l_0}^{-100} 1_{l_a < 3k - j_{+}} \int e^{i t \Phi} \psi_{j, k}^{\widehat{\mathcal{C}}}(\overline{\xi}) \Psi_{l_a}^a \widehat{f}(t, \overline{\xi}) d\overline{\xi} \label{estdispsupCloinLpetittau-ab2} \\
&\quad + \sum_{l_a = l_0}^{-100} 1_{3k - j_{+} \leq 2l_a} \int e^{i t \Phi} \psi_{j, k}^{\widehat{\mathcal{C}}}(\overline{\xi}) \Psi_{l_a}^c \widehat{f}(t, \overline{\xi}) d\overline{\xi} \label{estdispsupCloinLpetittau-c} \\
&\quad + \sum_{l_a = l_0}^{-100} 1_{2 l_a < 3k - j_{+} \leq l_a} \int e^{i t \Phi} \psi_{j, k}^{\widehat{\mathcal{C}}}(\overline{\xi}) \Psi_{l_a}^c \widehat{f}(t, \overline{\xi}) d\overline{\xi} \label{estdispsupCloinLpetittau-cb1} \\
&\quad + \sum_{l_a = l_0}^{-100} 1_{l_a < 3k - j_{+}} \int e^{i t \Phi} \psi_{j, k}^{\widehat{\mathcal{C}}}(\overline{\xi}) \Psi_{l_a}^c \widehat{f}(t, \overline{\xi}) d\overline{\xi} \label{estdispsupCloinLpetittau-cb2} \\
&\quad + \int e^{i t \Phi} \psi_{j, k}^{\widehat{\mathcal{C}}}(\overline{\xi}) \Psi_{l_0}^{int} \widehat{f}(t, \overline{\xi}) d\overline{\xi} \label{estdispsupCloinLpetittau-int} 
\end{align}
\end{subequations}
for $l_0$ such that $2^{l_0} \simeq \mathfrak{t}^{-\frac{1}{2}+\varrho}$ for some $\varrho > 0$ small enough with respect to $\delta$. Here above, $j_{+}$ is used for the positive part of $j$, i.e. $j_{+} = \max(j, 0)$. 

On the support of the considered integrals, 
\begin{align*}
\widehat{X}_a(\overline{\xi}) \cdot \nabla_{\overline{\xi}} \psi_{j, k}^{\widehat{\mathcal{C}}}(\overline{\xi}) &\lesssim 2^{-j} \\
\widehat{X}_{b-corr}(\overline{\xi}) \cdot \nabla_{\overline{\xi}} \psi_{j, k}^{\widehat{\mathcal{C}}}(\overline{\xi}) &\lesssim 2^{-j-k} \\
\widehat{X}_c(\overline{\xi}) \cdot \nabla_{\overline{\xi}} \psi_{j, k}^{\widehat{\mathcal{C}}}(\overline{\xi}) &= 0 \\
\widehat{X}_a(\overline{\xi}) \cdot \nabla_{\overline{\xi}} \left[ \widehat{X}_a(\overline{\xi}) \cdot \nabla_{\overline{\xi}} \Phi \right] &\simeq 2^j \\
\widehat{X}_{b-corr}(\overline{\xi}) \cdot \nabla_{\overline{\xi}} \left[ \widehat{X}_a(\overline{\xi}) \cdot \nabla_{\overline{\xi}} \Phi \right] &\lesssim 2^{j+k} + 2^{j+l_a} \\
\widehat{X}_c(\overline{\xi}) \cdot \nabla_{\overline{\xi}} \left[ \widehat{X}_a(\overline{\xi}) \cdot \nabla_{\overline{\xi}} \Phi \right] &\lesssim 2^{j+l_a} \\
\widehat{X}_a(\overline{\xi}) \cdot \nabla_{\overline{\xi}} \left[ \widehat{X}_c(\overline{\xi}) \cdot \nabla_{\overline{\xi}} \Phi \right] &= 0 \\
\widehat{X}_{b-corr}(\overline{\xi}) \cdot \nabla_{\overline{\xi}} \left[ \widehat{X}_c(\overline{\xi}) \cdot \nabla_{\overline{\xi}} \Phi \right] &\lesssim 2^{j+l_a} \\
\widehat{X}_c(\overline{\xi}) \cdot \nabla_{\overline{\xi}} \left[ \widehat{X}_c(\overline{\xi}) \cdot \nabla_{\overline{\xi}} \Phi \right] &\lesssim 2^j 
\end{align*}
and similar bounds by below. 
In particular, if $l_a \geq 2k-10$, we have that 
\begin{align*}
m_b(\overline{\xi}) \widehat{X}_{b-corr}(\overline{\xi}) \cdot \nabla_{\overline{\xi}} \left[ \widehat{X}_a(\overline{\xi}) \cdot \nabla_{\overline{\xi}} \Phi \right] &\lesssim 2^{j+l_a} + 2^{j+2k} \lesssim 2^{j+l_a} 
\end{align*}
On the other hand, if $l_a \leq 2k-10$, we may correct again $\widehat{X}_{b-corr}$: 
\begin{align*}
\widehat{X}_{b-corr}'(\overline{\xi}) &:= \widehat{X}_{b-corr}(\overline{\xi}) - \frac{\widehat{X}_{b-corr}(\overline{\xi}) \cdot \nabla_{\overline{\xi}} \Phi}{\widehat{X}_a(\overline{\xi}) \cdot \nabla_{\overline{\xi}} \Phi} \widehat{X}_a(\overline{\xi}) 
\end{align*}
so that 
\begin{align*}
\widehat{X}_{b-corr}'(\overline{\xi}) \cdot \nabla_{\overline{\xi}} \psi_{j, k}^{\widehat{\mathcal{C}}}(\overline{\xi}) &\lesssim 2^{-j-k} \\
\widehat{X}_{b-corr}'(\overline{\xi}) \cdot \nabla_{\overline{\xi}} \left[ \widehat{X}_a(\overline{\xi}) \cdot \nabla_{\overline{\xi}} \Phi \right] &= 0 \\
\widehat{X}_{b-corr}'(\overline{\xi}) \cdot \nabla_{\overline{\xi}} \left[ \widehat{X}_c(\overline{\xi}) \cdot \nabla_{\overline{\xi}} \Phi \right] &\lesssim 2^{j+l_a} 
\end{align*}
Moreover, since $l_a \leq 2k-10$, on the support of $\Psi_{l_a}^{*}$, we have that 
\begin{align*}
&\frac{\widehat{X}_{b-corr}(\overline{\xi}) \cdot \nabla_{\overline{\xi}} \Phi}{\widehat{X}_a(\overline{\xi}) \cdot \nabla_{\overline{\xi}} \Phi} \widehat{X}_a(\overline{\xi}) \cdot \nabla_{\overline{\xi}} \Phi \\
&\quad \lesssim \left( 2^{l_a} + 2^k \right) 2^{l_a+2j} \ll 2^{2k+2j} 
\end{align*}
therefore
\begin{align*}
\widehat{X}_{b-corr}'(\overline{\xi}) \cdot \nabla_{\overline{\xi}} \Phi &= o(2^{2k+2j}) + \widehat{X}_{b-corr}(\overline{\xi}) \cdot \nabla_{\overline{\xi}} \Phi \\
&= o(2^{2k+2j}) + O(2^{l_a+2j}) + 2^{2k} \simeq 2^{2k}  
\end{align*}
To avoid an overburden of notations, we will simply denote by $\widehat{X}_{b-corr}$ the field $\widehat{X}_{b-corr}'$ when $l_a \leq 2k-10$. 

Let us introduce the following local coordinates near a fixed point $\overline{\xi}^Z$: 
\begin{align*}
\Lambda : \overline{\xi} \mapsto \left( |\overline{\xi}|, 2^j \frac{\sqrt{3} |\xi_0| - |\xi|}{|\overline{\xi}|}, 2^j \arcsin \frac{J \xi^Z \cdot \xi}{|\xi^Z| |\xi|} \right) 
\end{align*}
Then $\Lambda$ is locally invertible: indeed, 
\begin{align*}
\widehat{X}_a(\overline{\xi}) \cdot \nabla_{\overline{\xi}} \Lambda &= \left( 1, 0, 0 \right) \\
\widehat{X}_b(\overline{\xi}) \cdot \nabla_{\overline{\xi}} \Lambda &= \frac{2^j}{|\overline{\xi}|} \left( 0, \frac{\sqrt{3} |\xi| \xi_0}{|\overline{\xi}| |\xi_0|} + \frac{\xi_0}{|\overline{\xi}|}, 0 \right) \\
\widehat{X}_c(\overline{\xi}) \cdot \nabla_{\overline{\xi}} \Lambda &= \frac{2^j}{|\overline{\xi}|} \left( 0, 0, \arcsin'\left( \frac{J \xi^Z \cdot \xi}{|\xi^Z| |\xi|} \right) \frac{\xi^Z \cdot \xi}{|\xi^Z| |\xi|} \right) 
\end{align*}
In particular, the determinant of $D \Lambda$ is uniformly of order $1$. Moreover, if we see $\Lambda$ as a change of variables and compute in the coordinates given by $\Lambda$, that we denote by $a, b, c$, we have that 
\begin{align*}
\partial_{\alpha} (F \circ \Lambda^{-1}) &\simeq \left( \widehat{X}_{\alpha} \cdot \nabla F \right) \circ \Lambda^{-1} 
\end{align*}

We will then denote by $L^{p_a}_a L^{p_b}_b L^{p_c}_c$ the norm
\begin{align*}
\Vert F \Vert_{L^{p_a}_a L^{p_b}_b L^{p_c}_c} &= \left( \int \left( \int \left( \int |(F \circ \Lambda^{-1})(\xi_a, \xi_b, \xi_c)|^{p_c} ~ d\xi_c \right)^{\frac{p_b}{p_c}} ~ d\xi_b \right)^{\frac{p_a}{p_b}} ~ d\xi_a \right)^{\frac{1}{p_a}} 
\end{align*}
for $F$ a function compactly supported on a neighborhood of $\overline{\xi}^Z$. Here above, the integration in $\xi_a$ is on a neighborhood of $2^j$, and in $\xi_b, \xi_c$ on a neighborhood of $0$. Analogously, we define $L^{p_b}_b L^{p_a}_a L^{p_c}_c$, or any permutation of the order of integration. In the case $p_a = p_b = p_c = p$, denoting by $L^{p}_{a, b, c}$ the corresponding norm, we have that 
\begin{align*}
\Vert F \Vert_{L^{p}_{a, b, c}} &\simeq \Vert F \Vert_{L^{p}}
\end{align*}
again for $F$ compactly supported near $\overline{\xi}^Z$. It is clear that usual Hölder inequalities are preserved for $L^{p_a}_a L^{p_b}_b L^{p_c}_c$ spaces or their permutations. Finally, we may define Sobolev norms
\begin{align*}
\Vert F \Vert_{\dot{H}^{s_1}_a \dot{H}^{s_2}_b \dot{H}^{s_3}_c} &= \left( \int \int \int \left[ |D_a|^{s_1} |D_b|^{s_2} |D_c|^{s_3} (F \circ \Lambda^{-1}) \right]^2 ~ d\xi_a d\xi_b d\xi_c \right)^{\frac{1}{2}}
\end{align*}
We also have Sobolev's embeddings in these coordinates and standard interpolations, and the correspondence 
\begin{align*}
\Vert F \Vert_{\dot{H}^1_a L^2_{b, c}} \simeq \Vert \widehat{X}_a \cdot \nabla F \Vert_{L^2} 
\end{align*}
and all others exchanging the role of $a, b, c$. 

We may now compute anisotropic norms of $\left( \Psi_{l_a} \psi_{j, k}^{\widehat{\mathcal{C}}} \right) \circ \Lambda^{-1}$: 
\begin{align*}
&\Vert \Psi_{l_a} \psi_{j, k}^{\widehat{\mathcal{C}}} \Vert_{L^{p_c}_c L^{p_b}_b L^{p_a}_a} \\
&\lesssim \left( \int 1_{\xi_c \lesssim 2^{l_a+j}} \left( \int 1_{\xi_b \lesssim 2^{j+k}} \left( \int \psi\left( 2^{-l_a-2j} \left( 3 \xi_0 |\overline{\xi}| - \frac{\xi_0 x}{t |\overline{\xi}|} - \frac{\xi \cdot y}{t |\overline{\xi}|} \right) \right) ~ d\xi_a \right)^{\frac{p_b}{p_a}} ~ d\xi_b \right)^{\frac{p_c}{p_b}} ~ d\xi_c \right)^{\frac{1}{p_c}}
\end{align*}
Considering $\xi_b, \xi_c$ fixed as well as $(x, y)$, the last $\psi$ localisation above forces $|\overline{\xi}|$ to be close to a point of size $2^j$ with precision $2^{l_a+j}$. Therefore, 
\begin{align}
\Vert \Psi_{l_a} \psi_{j, k}^{\widehat{\mathcal{C}}} \Vert_{L^{p_c}_c L^{p_b}_b L^{p_a}_a} &\lesssim 2^{\frac{l_a+j}{p_a} + \frac{k+j}{p_b} + \frac{l_a+j}{p_c}} \label{estimeevolumiquePsilacasCloinLloinC}
\end{align}

On \eqref{estdispsupCloinLpetittau-a}, we apply one or two integrations by parts along $\widehat{X}_a$: 
\begin{subequations}
\begin{align}
\eqref{estdispsupCloinLpetittau-a} &= \sum_{l_a = l_0}^{-100} 1_{2 l_a + j_{+} \geq 3k} t^{-2} 2^{-6j-4l_a} \int e^{i t \Phi} \psi_{j, k}^{\widehat{\mathcal{C}}}(\overline{\xi}) \Psi_{l_a}^a \widehat{f}(t, \overline{\xi}) d\overline{\xi} \label{estdispsupCloinLpetittau-a-1} \\
&\quad + \sum_{l_a = l_0}^{-100} 1_{2 l_a + j_{+} \geq 3k} t^{-2} 2^{-5j-3l_a} \int e^{i t \Phi} \psi_{j, k}^{\widehat{\mathcal{C}}}(\overline{\xi}) \Psi_{l_a}^a \widehat{X}_a(\overline{\xi}) \cdot \nabla_{\overline{\xi}} \widehat{f}(t, \overline{\xi}) d\overline{\xi} \label{estdispsupCloinLpetittau-a-2} \\
&\quad + \sum_{l_a = l_0}^{-100} 1_{2 l_a + j_{+} \geq 3k} t^{-2} 2^{-5j-3l_a} \int e^{i t \Phi} \psi_{j, k}^{\widehat{\mathcal{C}}}(\overline{\xi}) \Psi_{l_a}^a \widehat{h}_a(t, \overline{\xi}) d\overline{\xi} \label{estdispsupCloinLpetittau-a-3} \\
&\quad + \sum_{l_a = l_0}^{-100} 1_{2 l_a + j_{+} \geq 3k} t^{-2} 2^{-4j-2l_a} \int e^{i t \Phi} \psi_{j, k}^{\widehat{\mathcal{C}}}(\overline{\xi}) \Psi_{l_a}^a \widehat{X}_a(\overline{\xi}) \cdot \nabla_{\overline{\xi}} \widehat{h}_a(t, \overline{\xi}) d\overline{\xi} \label{estdispsupCloinLpetittau-a-4} \\
&\quad + \sum_{l_a = l_0}^{-100} 1_{2 l_a + j_{+} \geq 3k} t^{-1} 2^{-2j-l_a} \int e^{i t \Phi} \psi_{j, k}^{\widehat{\mathcal{C}}}(\overline{\xi}) \Psi_{l_a}^a \widehat{g}_a(t, \overline{\xi}) d\overline{\xi} \label{estdispsupCloinLpetittau-a-5}
\end{align}
\end{subequations} 
where symbols may vary from line to line as long as they keep similar properties. 

We now estimate for $\kappa > 0$ a small enough parameter: 
\begin{align*}
\eqref{estdispsupCloinLpetittau-a-1} &\lesssim \sum_{l_a = l_0}^{-100} 1_{2 l_a + j_{+} \geq 3k} t^{-2} 2^{-6j-4l_a} \Vert \Psi_{l_a}^a \psi_{j, k}^{\widehat{\mathcal{C}}} \Vert_{L^{\frac{1}{1-\kappa}}_c L^2_b L^{\frac{1}{1-\kappa}}_a} \langle 2^j \rangle^{-1} \Vert \langle \overline{\xi} \rangle \widehat{f}(t) \Vert_{L^{\frac{1}{\kappa}}_c L^2_b L^{\frac{1}{\kappa}}_a} \\
&\lesssim \sum_{l_a = l_0}^{-100} 1_{2 l_a + j_{+} \geq 3k} t^{-2} 2^{-\frac{7j}{2}-2l_a-2\kappa l_a+\frac{k}{2}-2\kappa j} \langle 2^j \rangle^{-1} \Vert \langle \overline{\xi} \rangle \widehat{f}(t) \Vert_{H^1_{a, c} L^2_b} \\
&\lesssim t^{-2} 2^{-\frac{7j}{2}-\frac{5l_0}{3}-2\delta l_0-2\kappa l_0-\frac{k}{2}+\frac{j_{+}}{6}+3 \delta k+\frac{k}{2}-2\kappa j} \langle 2^j \rangle^{-1} \Vert u \Vert_X \\
&\lesssim t^{-\frac{7}{6}+\delta+\kappa} 2^{-j+3\delta j + \kappa j+3 \delta k} \langle 2^j \rangle^{-\frac{5}{6}} \Vert u \Vert_X \\
\eqref{estdispsupCloinLpetittau-a-2} &\lesssim \sum_{l_a = l_0}^{-100} 1_{2 l_a + j_{+} \geq 3k} t^{-2} 2^{-5j-3l_a} \Vert \Psi_{l_a}^a \psi_{j, k}^{\widehat{\mathcal{C}}} \Vert_{L^2} \langle 2^j \rangle^{-1} \Vert \langle \nabla \rangle X_a f(t) \Vert_{L^2} \\
&\lesssim \sum_{l_a = l_0}^{-100} 1_{2 l_a + j_{+} \geq 3k} t^{-2} 2^{-\frac{7j}{2}-2l_a+\frac{k}{2}} \langle 2^j \rangle^{-1} \Vert u \Vert_X \\
&\lesssim t^{-\frac{7}{6}+\delta} 2^{-j+3\delta j+3 \delta k} \langle 2^j \rangle^{-\frac{5}{6}} \Vert u \Vert_X \\
\eqref{estdispsupCloinLpetittau-a-4} &\lesssim \sum_{l_a = l_0}^{-100} 1_{2 l_a + j_{+} \geq 3k} t^{-2} 2^{-5j-2l_a} \Vert \Psi_{l_a}^a \psi_{j, k}^{\widehat{\mathcal{C}}} \Vert_{L^2} \Vert \nabla X_a h_a(t) \Vert_{L^2} \\
&\lesssim \sum_{l_a = l_0}^{-100} 1_{2 l_a + j_{+} \geq 3k} t^{-2} 2^{-\frac{7j}{2}-l_a+\frac{k}{2}} \Vert u \Vert_X \\
&\lesssim t^{-2} 2^{\beta j - \frac{7j}{2}-l_0+\frac{k}{2}} \Vert u \Vert_X \\
&\lesssim t^{-\frac{3}{2}} 2^{ - 2j+\frac{k}{2}} \Vert u \Vert_X \\
&\lesssim t^{-\frac{7}{6}} 2^{ - j + \frac{k}{2}} \langle 2^j \rangle^{-1} \Vert u \Vert_X \\
\eqref{estdispsupCloinLpetittau-a-5} &\lesssim \sum_{l_a = l_0}^{-100} 1_{2 l_a + j_{+} \geq 3k} t^{-1} 2^{-\frac{5j}{2}-l_a} \Vert \Psi_{l_a}^a \psi_{j, k}^{\widehat{\mathcal{C}}} \Vert_{L^2} \langle 2^j \rangle^{-1} \Vert m_{\widehat{\mathcal{C}}}(D) \langle \nabla \rangle |\nabla|^{\frac{1}{2}} g_a(t) \Vert_{L^2} \\
&\lesssim \sum_{l_a = l_0}^{-100} 1_{2 l_a + j_{+} \geq 3k} t^{-\frac{7}{6}+100\delta} 2^{-j+\frac{k}{2}} \langle 2^j \rangle^{-1} \Vert u \Vert_X \\
&\lesssim t^{-\frac{7}{6}+100\delta} 2^{-j+\frac{k}{3}} \langle 2^j \rangle^{-\frac{5}{6}} \Vert u \Vert_X
\end{align*}
\eqref{estdispsupCloinLpetittau-a-3} is similar to \eqref{estdispsupCloinLpetittau-a-2}. 

We estimate \eqref{estdispsupCloinLpetittau-c} in a similar fashion, replacing $\widehat{X}_a$ by $\widehat{X}_c$. 

Then, \eqref{estdispsupCloinLpetittau-ab1} only exists if $3k > 2 l_0 + j_{+}$, which we now assume. Having fixed $\varrho$ small with respect to $\delta$, we apply at most $n$, large with respect to $\varrho^{-1}$, integrations by parts along $\widehat{X}_a$, as long as the integration by parts does not hit $\widehat{f}(t)$: 
\begin{align*}
\eqref{estdispsupCloinLpetittau-ab1} &= \sum_{l_a = l_0}^{-100} 1_{2 l_a < 3k- j_{+} \leq l_a} t^{-n} 2^{-3nj-2nl_a} \int e^{i t \Phi} \psi_{j, k}^{\widehat{\mathcal{C}}}(\overline{\xi}) \Psi_{l_a}^a \widehat{f}(t, \overline{\xi}) d\overline{\xi} \\
&\quad + \sum_{i = 1}^n \sum_{l_a = l_0}^{-100} 1_{2 l_a < 3k - j_{+} \leq l_a} t^{-i} 2^{-3ij+j-2il_a+l_a} \int e^{i t \Phi} \psi_{j, k}^{\widehat{\mathcal{C}}}(\overline{\xi}) \Psi_{l_a}^a \widehat{h}_a(t, \overline{\xi}) d\overline{\xi} \\
&\quad + \sum_{i = 1}^n \sum_{l_a = l_0}^{-100} 1_{2 l_a < 3k - j_{+} \leq l_a} t^{-i} 2^{-3ij+j-2il_a+l_a} \int e^{i t \Phi} \psi_{j, k}^{\widehat{\mathcal{C}}}(\overline{\xi}) \Psi_{l_a}^a \widehat{g}_a(t, \overline{\xi}) d\overline{\xi}
\end{align*}
Again, symbols may change from line to line. $n$ is chosen big enough with respect to $\varrho^{-1}$, itself big enough with respect to $\delta^{-1}$, but all these parameters can be considered fixed (the universal constants being allowed to depend on $\delta$), so that the total amount of derivations of symbols is fixed and therefore the universal constant also remains universal (or rather, only depending on $\delta$). Then, on the central term in which a $\widehat{h}_a$ appeared, we apply another integration by parts in the direction $\widehat{X}_{b-corr}$: 
\begin{subequations}
\begin{align}
\eqref{estdispsupCloinLpetittau-ab1} &= \sum_{l_a = l_0}^{-100} 1_{2 l_a < 3k- j_{+} \leq l_a} t^{-n} 2^{-3nj-2nl_a} \int e^{i t \Phi} \psi_{j, k}^{\widehat{\mathcal{C}}}(\overline{\xi}) \Psi_{l_a}^a \widehat{f}(t, \overline{\xi}) d\overline{\xi} \label{estdispsupCloinLpetittau-ab1-1} \\
&\quad + \sum_{i = 1}^n \sum_{l_a = l_0}^{-100} 1_{2 l_a < 3k- j_{+} \leq l_a} t^{-i-1} 2^{-3ij+j-2il_a+l_a-3j-3k} \int e^{i t \Phi} \psi_{j, k}^{\widehat{\mathcal{C}}}(\overline{\xi}) \Psi_{l_a}^a \widehat{h}_a(t, \overline{\xi}) d\overline{\xi} \label{estdispsupCloinLpetittau-ab1-2} \\
&\quad \begin{aligned}
&+ \sum_{i = 1}^n \sum_{l_a = l_0}^{-100} 1_{2 l_a < 3k- j_{+} \leq l_a} t^{-i-1} 2^{-3ij+j-2il_a+l_a-3j-3k} \\
&\quad \quad \quad \int e^{i t \Phi} \psi_{j, k}^{\widehat{\mathcal{C}}}(\overline{\xi}) \Psi_{l_a}^a m_b(\overline{\xi}) \widehat{X}_{b-corr}(\overline{\xi}) \cdot \nabla_{\overline{\xi}} \widehat{h}_a(t, \overline{\xi}) d\overline{\xi}
\end{aligned} \label{estdispsupCloinLpetittau-ab1-3} \\
&\quad + \sum_{i = 1}^n \sum_{l_a = l_0}^{-100} 1_{2 l_a < 3k- j_{+} \leq l_a} t^{-i} 2^{-3ij+j-2il_a+l_a} \int e^{i t \Phi} \psi_{j, k}^{\widehat{\mathcal{C}}}(\overline{\xi}) \Psi_{l_a}^a \widehat{g}_a(t, \overline{\xi}) d\overline{\xi} \label{estdispsupCloinLpetittau-ab1-4} 
\end{align}
\end{subequations} 
We now estimate: 
\begin{align*}
\eqref{estdispsupCloinLpetittau-ab1-1} &\lesssim \sum_{l_a = l_0}^{-100} 1_{2 l_a < 3k- j_{+} \leq l_a} t^{-n} 2^{-3nj-2nl_a+j} \Vert \Psi_{l_a}^a \psi_{j, k}^{\widehat{\mathcal{C}}} \Vert_{L^2} \langle 2^j \rangle^{-2} \Vert \langle \nabla \rangle^2 |\nabla|^{-1} f(t) \Vert_{L^2} \\
&\lesssim \sum_{l_a = l_0}^{-100} 1_{2 l_a + j_{+} < 3k} t^{-n} 2^{-3nj-2nl_a+l_a+\frac{k}{2}+\frac{5j}{2}} \langle 2^j \rangle^{-2} \Vert u \Vert_X \\
&\lesssim \mathfrak{t}^{-n} \mathfrak{t}^{n-2n\varrho} 2^{l_0+\frac{k}{2}+\frac{5j}{2}} \langle 2^j \rangle^{-2} \Vert u \Vert_X \\
&\lesssim t^{-\frac{7}{6}} \mathfrak{t}^{\frac{2}{3}-2n\varrho} 2^{-j+\frac{k}{2}} \langle 2^j \rangle^{-2} \Vert u \Vert_X \\
&\lesssim t^{-\frac{7}{6}} 2^{- j +\frac{k}{2}} \langle 2^j \rangle^{-2} \\
\eqref{estdispsupCloinLpetittau-ab1-3} &\lesssim \sum_{i = 1}^n \sum_{l_a = l_0}^{-100} 1_{2 l_a < 3k- j_{+} \leq l_a} t^{-i-1} 2^{-3ij-2il_a+l_a-3j-3k} \Vert \Psi_{l_a}^a \psi_{j, k}^{\widehat{\mathcal{C}}} \Vert_{L^2} \Vert \nabla m_b(D) X_{b-corr} h_a(t) \Vert_{L^2} \\
&\lesssim \sum_{i = 1}^n \sum_{l_a = l_0}^{-100} 1_{2 l_a + j_{+} < 3k} t^{-\frac{7}{6}} \mathfrak{t}^{-i+\frac{1}{6}} 2^{-j-2il_a+2l_a-\frac{5k}{2}} \Vert u \Vert_X \\
&\lesssim t^{-\frac{7}{6}} \mathfrak{t}^{\delta} 2^{-j+3\delta k} \langle 2^j \rangle^{-\frac{5}{6}} \Vert u \Vert_X \\
&\lesssim t^{-\frac{7}{6}+\delta} 2^{-j+3\delta j + 3 \delta k} \langle 2^j \rangle^{-\frac{5}{6}} \Vert u \Vert_X \\
\eqref{estdispsupCloinLpetittau-ab1-4} &\lesssim \sum_{i = 1}^n \sum_{l_a = l_0}^{-100} 1_{2 l_a < 3k- j_{+} \leq l_a} t^{-i} 2^{-3ij+\frac{j}{2}-2il_a+l_a} \Vert \Psi_{l_a}^a \psi_{j, k}^{\widehat{\mathcal{C}}} \Vert_{L^2} \langle 2^j \rangle^{-1} \Vert m_{\widehat{\mathcal{C}}}(D) \langle \nabla \rangle |\nabla|^{\frac{1}{2}} g_a(t) \Vert_{L^2} \\
&\lesssim \sum_{i = 1}^n \sum_{l_a = l_0}^{-100} 1_{2 l_a < 3k- j_{+} \leq l_a} t^{-\frac{7}{6}+100\delta} \mathfrak{t}^{-i+1} 2^{-j-2il_a+2l_a+\frac{k}{2}} \langle 2^j \rangle^{-1} \Vert u \Vert_X \\
&\lesssim t^{-\frac{7}{6}+100\delta} 2^{-j+\frac{k}{3}} \langle 2^j \rangle^{-\frac{5}{6}} \Vert u \Vert_X
\end{align*}
where we used that $n \varrho$ is large enough, and that $\mathfrak{t} \simeq t 2^{3j} \gtrsim 1$. \eqref{estdispsupCloinLpetittau-ab1-2} can be estimated like \eqref{estdispsupCloinLpetittau-ab1-3}. Here above, only for \eqref{estdispsupCloinLpetittau-ab1-4}, we used that the summation stopped at $l_a \geq 3k - j_{+}$ when $i = 1$, for which the summation in $l_a$ adds a factor $|3k-j_{+}|$, which can be absorbed into $2^{\frac{k}{2}} \langle 2^j \rangle^{-1}$. 

For \eqref{estdispsupCloinLpetittau-ab2}, we proceed the same way but reversing the order of the integration by parts and starting by the one along $\widehat{X}_{b-corr}$: 
\begin{align*}
\eqref{estdispsupCloinLpetittau-ab2} &= \sum_{l_a = l_0}^{-100} 1_{l_a < 3k - j_{+}} t^{-1} 2^{-3j-3k} \int e^{i t \Phi} \psi_{j, k}^{\widehat{\mathcal{C}}}(\overline{\xi}) \Psi_{l_a}^a \widehat{f}(t, \overline{\xi}) d\overline{\xi} \\
&+ \sum_{l_a = l_0}^{-100} 1_{l_a < 3k - j_{+}} t^{-1} 2^{-2j-3k} \int e^{i t \Phi} \psi_{j, k}^{\widehat{\mathcal{C}}}(\overline{\xi}) \Psi_{l_a}^a \widehat{h}_{b-corr}(t, \overline{\xi}) d\overline{\xi} \\
&+ \sum_{l_a = l_0}^{-100} 1_{l_a < 3k - j_{+}} t^{-1} 2^{-2j-3k} \int e^{i t \Phi} \psi_{j, k}^{\widehat{\mathcal{C}}}(\overline{\xi}) \Psi_{l_a}^a \widehat{g}_{b-corr}(t, \overline{\xi}) d\overline{\xi}
\end{align*}
where $h_{b-corr}, g_{b-corr}$ correspond to the decomposition of $m_b \widehat{X}_{b-corr} \cdot \nabla \widehat{f}(t)$ obtained by projecting $m_b \widehat{X}_{b-corr}$ on the basis $(\widehat{X}_a$, $\widehat{X}_b$, $\widehat{X}_c)$ and then using the canonical decompositions. Then, we apply an integration by parts along $\widehat{X}_a$: 
\begin{subequations}
\begin{align}
\eqref{estdispsupCloinLpetittau-ab2} &= \sum_{l_a = l_0}^{-100} 1_{l_a < 3k - j_{+}} t^{-2} 2^{-6j-3k-2l_a} \int e^{i t \Phi} \psi_{j, k}^{\widehat{\mathcal{C}}}(\overline{\xi}) \Psi_{l_a}^a \widehat{f}(t, \overline{\xi}) d\overline{\xi} \label{estdispsupCloinLpetittau-ab2-1} \\
&+ \sum_{l_a = l_0}^{-100} 1_{l_a < 3k - j_{+}} t^{-2} 2^{-5j-3k-l_a} \int e^{i t \Phi} \psi_{j, k}^{\widehat{\mathcal{C}}}(\overline{\xi}) \Psi_{l_a}^a \widehat{X}_a(\overline{\xi}) \cdot \nabla_{\overline{\xi}} \widehat{f}(t, \overline{\xi}) d\overline{\xi} \label{estdispsupCloinLpetittau-ab2-2} \\
&+ \sum_{l_a = l_0}^{-100} 1_{l_a < 3k - j_{+}} t^{-2} 2^{-5j-3k-2l_a} \int e^{i t \Phi} \psi_{j, k}^{\widehat{\mathcal{C}}}(\overline{\xi}) \Psi_{l_a}^a \widehat{h}_{b-corr}(t, \overline{\xi}) d\overline{\xi} \label{estdispsupCloinLpetittau-ab2-3} \\
&+ \sum_{l_a = l_0}^{-100} 1_{l_a < 3k - j_{+}} t^{-2} 2^{-4j-3k-l_a} \int e^{i t \Phi} \psi_{j, k}^{\widehat{\mathcal{C}}}(\overline{\xi}) \Psi_{l_a}^a \widehat{X}_a(\overline{\xi}) \cdot \nabla_{\overline{\xi}} \widehat{h}_{b-corr}(t, \overline{\xi}) d\overline{\xi} \label{estdispsupCloinLpetittau-ab2-4} \\
&+ \sum_{l_a = l_0}^{-100} 1_{l_a < 3k - j_{+}} t^{-1} 2^{-2j-3k} \int e^{i t \Phi} \psi_{j, k}^{\widehat{\mathcal{C}}}(\overline{\xi}) \Psi_{l_a}^a \widehat{g}_{b-corr}(t, \overline{\xi}) d\overline{\xi} \label{estdispsupCloinLpetittau-ab2-5} 
\end{align}
\end{subequations}
We now estimate: 
\begin{align*}
\eqref{estdispsupCloinLpetittau-ab2-3} &\lesssim \sum_{l_a = l_0}^{-100} 1_{l_a < 3k} t^{-2} 2^{-6j-3k-2l_a} \Vert \Psi_{l_a}^a \psi_{j, k}^{\widehat{\mathcal{C}}} \Vert_{L^{\frac{1}{1-\kappa}}_c L^2_b L^{\frac{1}{1-\kappa}}_a} \Vert \overline{\xi} \widehat{h}_{b-corr}(t) \Vert_{L^{\frac{1}{\kappa}}_c L^2_b L^{\frac{1}{\kappa}}_a} \\
&\lesssim \sum_{l_a = l_0}^{-100} 1_{l_a < 3k} t^{-\frac{7}{6}+\frac{2\kappa}{3}} \mathfrak{t}^{-\frac{5}{6}-\frac{2\kappa}{3}} 2^{-j-l_a - 2\kappa l_a + \frac{k}{2}} \Vert u \Vert_X \\
&\lesssim t^{-\frac{7}{6}+\frac{2\kappa}{3}} \mathfrak{t}^{-\frac{1}{3}+\frac{\kappa}{3}} 2^{-j + \frac{k}{2}} \Vert u \Vert_X \\
&\lesssim t^{-\frac{7}{6}+\frac{2\kappa}{3}} 2^{-j+\frac{k}{2}} \langle 2^j \rangle^{-\frac{5}{6}} \Vert u \Vert_X \\
\eqref{estdispsupCloinLpetittau-ab2-4} &\lesssim \sum_{l_a = l_0}^{-100} 1_{l_a < 3k} t^{-2} 2^{-5j-3k-l_a} \Vert \Psi_{l_a}^a \psi_{j, k}^{\widehat{\mathcal{C}}} \Vert_{L^2} \Vert \nabla X_a h_{b-corr} \Vert_{L^2} \\
&\lesssim \sum_{l_a = l_0}^{-100} 1_{l_a < 3k} t^{-\frac{7}{6}} \mathfrak{t}^{-\frac{5}{6}} 2^{-j-l_a+\frac{k}{2}} \Vert u \Vert_X \\
&\lesssim t^{-\frac{7}{6}} \mathfrak{t}^{-\frac{1}{3}} 2^{-j+\frac{k}{2}} \Vert u \Vert_X \\
&\lesssim t^{-\frac{7}{6}} 2^{-j+\frac{k}{2}} \langle 2^j \rangle^{-1} \Vert u \Vert_X \\
\eqref{estdispsupCloinLpetittau-ab2-5} &\lesssim \sum_{l_a = l_0}^{-100} 1_{l_a < 3k} t^{-1} 2^{-2j-3k-\frac{j}{2}} \Vert \Psi_{l_a}^a \psi_{j, k}^{\widehat{\mathcal{C}}} \Vert_{L^2} \langle 2^j \rangle^{-1} \Vert m_{\widehat{\mathcal{C}}}(D) \langle \nabla \rangle |\nabla|^{\frac{1}{2}} g_{b-corr}(t) \Vert_{L^2} \\
&\lesssim \sum_{l_a = l_0}^{-100} 1_{l_a < 3k} t^{-\frac{7}{6}+100\delta} 2^{-j-3k+l_a+\frac{k}{2}} \langle 2^j \rangle^{-1} \Vert u \Vert_X \\
&\lesssim t^{-\frac{7}{6}+100\delta} 2^{-j+\frac{k}{2}} \langle 2^j \rangle^{-1} \Vert u \Vert_X 
\end{align*}
for small enough $\kappa$. \eqref{estdispsupCloinLpetittau-ab2-1} is similar to \eqref{estdispsupCloinLpetittau-ab2-3} and \eqref{estdispsupCloinLpetittau-ab2-2} to \eqref{estdispsupCloinLpetittau-ab2-4}. 

We proceed the same way for \eqref{estdispsupCloinLpetittau-cb1} and \eqref{estdispsupCloinLpetittau-cb2} replacing $\widehat{X}_a$ by $\widehat{X}_c$. 

Finally, for the internal term \eqref{estdispsupCloinLpetittau-int}, we exhaust the cases depending on $j, k, t$. 

Assume first that $2^j \lesssim t^{-\frac{3}{44}}$. If $2^k \lesssim t^{-\frac{7}{24}} 2^{-\frac{5j}{6}}$, then we estimate directly
\begin{align*}
\eqref{estdispsupCloinLpetittau-int} &\lesssim \Vert \Psi_{l_0}^{int} \psi_{j, k}^{\widehat{\mathcal{C}}} \Vert_{L^{\frac{1}{1-\kappa}}_c L^2_b L^{\frac{1}{1-\kappa}}_a} \Vert \widehat{f}(t) \Vert_{L^{\frac{1}{\kappa}}_c L^2_b L^{\frac{1}{\kappa}}_a} \\
&\lesssim 2^{2l_0-2\kappa l_0+\frac{k}{2}+\frac{5j}{2}-2\kappa j} \Vert u \Vert_X \\
&\lesssim t^{-\frac{7}{6}+\frac{2\kappa}{3}} \mathfrak{t}^{\frac{1}{6}+2\varrho+\frac{\kappa}{3}} t^{-\frac{7}{48}+\frac{7\delta}{24}} 2^{-\frac{5j}{12}-j+\delta k+\frac{5\delta j}{6}} \Vert u \Vert_X \\
&\lesssim t^{-\frac{7}{6}+\frac{\kappa}{3}+2\varrho+\frac{7\delta}{24}} t^{\frac{1}{48}} 2^{\frac{j}{12}-j+\delta k + \frac{5\delta j}{6}} \Vert u \Vert_X
\end{align*}
In particular, 
\begin{align*}
\eqref{estdispsupCloinLpetittau-int} &\lesssim t^{-\frac{13}{12}} 2^{-j+\delta k} \Vert u \Vert_X
\end{align*}
and 
\begin{align*}
\eqref{estdispsupCloinLpetittau-int} &\lesssim t^{-\frac{7}{6}} 2^{-\frac{3j}{2}+ \delta k} \Vert u \Vert_X
\end{align*}
if $\kappa, \varrho, \delta$ are small enough. If however $2^k \gtrsim t^{-\frac{7}{24}} 2^{-\frac{5j}{6}}$, then we apply an integration by parts along $\widehat{X}_{b-corr}$: 
\begin{subequations}
\begin{align}
\eqref{estdispsupCloinLpetittau-int} &= t^{-1} 2^{-3k-3j} \int e^{i t \Phi} \psi_{j, k}^{\widehat{\mathcal{C}}} \Psi_{l_0}^{int} \widehat{f}(t) d\overline{\xi} \label{estdispsupCloinLpetittau-int-1} \\
&+ t^{-1} 2^{-3k-2j} \int e^{i t \Phi} \psi_{j, k}^{\widehat{\mathcal{C}}} \Psi_{l_0}^{int} \widehat{h}_{b-corr}(t) d\overline{\xi} \label{estdispsupCloinLpetittau-int-2} \\
&+ t^{-1} 2^{-3k-2j} \int e^{i t \Phi} \psi_{j, k}^{\widehat{\mathcal{C}}} \Psi_{l_0}^{int} \widehat{g}_{b-corr}(t) d\overline{\xi} \label{estdispsupCloinLpetittau-int-3}
\end{align}
\end{subequations}
We then estimate
\begin{align*}
\eqref{estdispsupCloinLpetittau-int-2} &\lesssim t^{-1} 2^{-3k-3j} \Vert \Psi_{l_0}^{int} \psi_{j, k}^{\widehat{\mathcal{C}}} \Vert_{L^{\frac{1}{1-\kappa}}_c L^2_b L^{\frac{1}{1-\kappa}}_a} \Vert \overline{\xi} \widehat{h}_{b-corr}(t) \Vert_{L^{\frac{1}{\kappa}}_c L^2_b L^{\frac{1}{\kappa}}_a} \\
&\lesssim t^{-1} 2^{2l_0-2\kappa l_0 - \frac{5k}{2}-\frac{j}{2}-2\kappa j} \Vert u \Vert_X \\
&\lesssim t^{-\frac{7}{6}+\frac{2\kappa}{3}} \mathfrak{t}^{-\frac{5}{6}+2\varrho+\frac{\kappa}{3}} t^{\frac{35}{48}+\frac{7 \delta}{24}} 2^{\frac{25j}{12}+\frac{5 \delta j}{6}} 2^{-j+\delta k} \Vert u \Vert_X \\
&\lesssim t^{-\frac{7}{6}+\kappa+\varrho+\frac{7\delta}{24}} t^{-\frac{5}{48}} 2^{-\frac{5j}{12}+\frac{5 \delta j}{6}} 2^{-j+\delta k} \Vert u \Vert_X \\
\eqref{estdispsupCloinLpetittau-int-3} &\lesssim t^{-1} 2^{-3k-3j} \Vert \Psi_{l_0}^{int} \psi_{j, k}^{\widehat{\mathcal{C}}} \Vert_{L^{\frac{1}{1-\kappa}}_c L^2_{b,a}} \Vert \psi\left( 2^{-j} \overline{\xi} \right) \overline{\xi} \widehat{g}_{b-corr}(t) \Vert_{L^{\frac{1}{\kappa}}_c L^2_{b, a}} \\
&\lesssim t^{-1} 2^{\frac{3l_0}{2}-\kappa l_0-\frac{5k}{2}-j-\kappa j} \Vert \psi\left( 2^{-j} \overline{\xi} \right) \overline{\xi} \widehat{g}_{b-corr}(t) \Vert_{H^1_c L^2_{b, a}}^{\frac{1}{2}} \Vert \psi\left( 2^{-j} \overline{\xi} \right) \overline{\xi} \widehat{g}_{b-corr}(t) \Vert_{L^2}^{\frac{1}{2}} \\
&\lesssim t^{-1} \mathfrak{t}^{-\frac{3}{4}+\frac{3\varrho}{2}+\frac{\kappa}{2}} \left( t^{-\frac{7}{24}} 2^{-\frac{5j}{6}} \right)^{-\frac{5}{2}-\delta} 2^{\delta k-j-\kappa j} \left( t^{-\frac{1}{12}+\frac{1}{32}+50\delta} + t^{-\frac{1}{8}} 2^{-\frac{7j}{6}+1000\delta} \right)^{\frac{1}{2}} \\
&\quad \quad \quad \left( t^{-\frac{1}{6}+100\delta} 2^{\frac{j}{2}} \right)^{\frac{1}{2}} \Vert u \Vert_X \\
&\lesssim t^{-1-\frac{7}{48}+\frac{1}{64}+76\delta+\frac{3\varrho}{2}+\frac{\kappa}{2}} 2^{-\frac{11j}{12}+\frac{5\delta j}{6}+\frac{9\varrho j}{2}+\frac{\kappa j}{2}+\delta k} \Vert u \Vert_X \\
&\quad + t^{-\frac{7}{6}+51\delta+\frac{3\varrho}{2}+\frac{\kappa}{2}} 2^{-\frac{3j}{12}+500\delta j+\frac{5\delta j}{6}+\frac{9\varrho j}{2}+\frac{\kappa j}{2}+\delta k} \Vert u \Vert_X \\
&\lesssim t^{-1-\frac{7}{48}+\frac{1}{64}-\frac{7}{176}+76\delta+\frac{3\varrho}{2}+\frac{\kappa}{2}} 2^{-\frac{3j}{2}+\frac{5\delta j}{6}+\frac{9\varrho j}{2}+\frac{\kappa j}{2}+\delta k} \Vert u \Vert_X \\
&\quad + t^{-\frac{7}{6}+51\delta+\frac{3\varrho}{2}+\frac{\kappa}{2}} 2^{-\frac{3j}{2}+500\delta j+\frac{5\delta j}{6}+\frac{9\varrho j}{2}+\frac{\kappa j}{2}+\delta k} \Vert u \Vert_X \\
&\lesssim t^{-\frac{7}{6}+100\delta} 2^{-\frac{3j}{2}+\delta j+\delta k} \Vert u \Vert_X 
\end{align*}
In particular, 
\begin{align*}
\eqref{estdispsupCloinLpetittau-int-2} &\lesssim t^{-\frac{7}{6}+\kappa+\varrho+\frac{7\delta}{24}} t^{-\frac{5}{48}+\frac{5}{36}} 2^{-j+\delta k} \Vert u \Vert_X \\
&\lesssim t^{-\frac{13}{12}} 2^{-j+\delta k} \Vert u \Vert_X 
\end{align*}
and 
\begin{align*}
\eqref{estdispsupCloinLpetittau-int-2} &\lesssim t^{-\frac{7}{6}} 2^{-\frac{3j}{2}+\delta k} \Vert u \Vert_X 
\end{align*}
for small enough $\kappa, \varrho, \delta$. \eqref{estdispsupCloinLpetittau-int-1} can be estimated as \eqref{estdispsupCloinLpetittau-int-2}. For \eqref{estdispsupCloinLpetittau-int-3}, we already have one of the estimates, and the other is obtained in a simpler way: 
\begin{align*}
\eqref{estdispsupCloinLpetittau-int-3} &\lesssim t^{-1} 2^{-3k-3j} \Vert \Psi_{l_0}^{int} \psi_{j, k}^{\widehat{\mathcal{C}}} \Vert_{L^{\frac{1}{1-\kappa}}_c L^2_{b,a}} \Vert \psi\left( 2^{-j} \overline{\xi} \right) \overline{\xi} \widehat{g}_{b-corr}(t) \Vert_{L^{\frac{1}{\kappa}}_c L^2_{b, a}} \\
&\lesssim t^{-1} 2^{\frac{3l_0}{2}-\kappa l_0-\frac{5k}{2}-j-\kappa j} \Vert \psi\left( 2^{-j} \overline{\xi} \right) \overline{\xi} \widehat{g}_{b-corr}(t) \Vert_{H^1_c L^2_{b, a}}^{\frac{1}{2}} \Vert \psi\left( 2^{-j} \overline{\xi} \right) \overline{\xi} \widehat{g}_{b-corr}(t) \Vert_{L^2}^{\frac{1}{2}} \\
&\lesssim t^{-1} \mathfrak{t}^{-\frac{3}{4}+\frac{3\varrho}{2}+\frac{\kappa}{2}} \left( t^{-\frac{7}{24}} 2^{-\frac{5j}{6}} \right)^{-\frac{5}{2}-\delta} 2^{-j-\kappa j+\delta k} \left( t^{-\frac{1}{6}+100\delta} 2^{\frac{j}{2}} \right)^{\frac{1}{2}} \Vert u \Vert_X \\
&\lesssim t^{-1-\frac{1}{48}-\frac{1}{12}+51\delta+\frac{3\varrho}{2}+\frac{\kappa}{2}} 2^{-\frac{11j}{12}+\frac{9\varrho j}{2}+\frac{\kappa j}{2}+\frac{5\delta j}{6}+\delta k} \Vert u \Vert_X \\
&\lesssim t^{-\frac{13}{12}} 2^{-j+\delta j+\delta k} \Vert u \Vert_X
\end{align*}
as wanted. 

Let now assume $2^j \gtrsim t^{-\frac{3}{44}}$. If first $2^k \lesssim t^{-\frac{9}{35}} 2^{-\frac{4j}{7}} \langle 2^j \rangle^{-\frac{1}{7}}$, then we estimate
\begin{align*}
\eqref{estdispsupCloinLpetittau-int} &\lesssim 2^{-j} \Vert \Psi_{l_0}^{int} \psi_{j, k}^{\widehat{\mathcal{C}}} \Vert_{L^{\frac{1}{1-\kappa}}} \Vert \overline{\xi} \widehat{f}(t) \Vert_{L^{\frac{1}{\kappa}}} \\
&\lesssim 2^{2l_0-2\kappa l_0+k - \kappa k + 2j - 3 \kappa j} \Vert \overline{\xi} \widehat{f}(t) \Vert_{H^1_a H^1_c L^2_b}^{\frac{1}{2}} \Vert \overline{\xi} \widehat{f}(t) \Vert_{H^1_b L^2_{a, b}}^{\frac{1}{2}} \\
&\lesssim \mathfrak{t}^{-1+2\varrho+\kappa} t^{-\frac{9}{35}+\frac{9 (\kappa+\delta)}{35}} 2^{-\frac{4j}{7}+\frac{4\kappa j}{7}+\frac{4\delta j}{7}} \langle 2^j \rangle^{-\frac{1}{7}+\frac{\kappa+\delta}{7}} 2^{\delta k + 2j - 3 \kappa j} t^{\frac{1}{12}} \Vert u \Vert_X \\
&\lesssim t^{-\frac{7}{6}} t^{-\frac{1}{140}+2\varrho+\kappa+\frac{9 (\kappa+\delta)}{35}} 2^{-\frac{4j}{7}} \langle 2^j \rangle^{-\frac{1}{14}} 2^{-j+\delta k} \Vert u \Vert_X
\end{align*}
In particular, 
\begin{align*}
\eqref{estdispsupCloinLpetittau-int} &\lesssim t^{-\frac{7}{6}} t^{-\frac{1}{140}+2\varrho+\kappa+\frac{9 (\kappa+\delta)}{35}+\frac{3}{77}} 2^{-\frac{4j}{7}} \langle 2^j \rangle^{-\frac{1}{14}} 2^{-j+\delta k} \Vert u \Vert_X \\
&\lesssim t^{-\frac{13}{12}} 2^{-j+\delta k} \langle 2^j \rangle^{-\frac{1}{14}} \Vert u \Vert_X 
\end{align*}
and 
\begin{align*}
\eqref{estdispsupCloinLpetittau-int} &\lesssim t^{-\frac{7}{6}} t^{-\frac{1}{140}+2\varrho+\kappa+\frac{9(\kappa+\delta)}{35}+\frac{3}{616}} \langle 2^j \rangle^{-\frac{1}{14}} 2^{-\frac{3j}{2}+\delta k} \Vert u \Vert_X \\
&\lesssim t^{-\frac{7}{6}} 2^{-\frac{3j}{2}+\delta k} \langle 2^j \rangle^{-\frac{1}{14}} \Vert u \Vert_X 
\end{align*}
Finally, if $2^k \gtrsim t^{-\frac{9}{35}} 2^{-\frac{4j}{7}} \langle 2^j \rangle^{-\frac{1}{7}}$, we apply the same integrations by parts as above and we estimate: 
\begin{align*}
\eqref{estdispsupCloinLpetittau-int-2} &\lesssim t^{-1} 2^{-3k-3j} \Vert \Psi_{l_0}^{int} \psi_{j, k}^{\widehat{\mathcal{C}}} \Vert_{L^{\frac{1}{1-\kappa}}_c L^2_b L^{\frac{1}{1-\kappa}}_a} \Vert \overline{\xi} \widehat{h}_{b-corr}(t) \Vert_{L^{\frac{1}{\kappa}}_c L^2_b L^{\frac{1}{\kappa}}_a} \\
&\lesssim t^{-1} 2^{2l_0-2\kappa l_0-\frac{5k}{2}-\frac{j}{2}-2\kappa j} \Vert u \Vert_X \\
&\lesssim t^{-\frac{7}{6}} t^{\frac{1}{6}+\frac{9}{14}+\frac{9\delta}{35}} \mathfrak{t}^{-1+2\varrho+\kappa} 2^{\frac{27j}{14}+\frac{4\delta j}{7}-2\kappa j} \langle 2^j \rangle^{\frac{5}{14}+\frac{\delta}{7}} 2^{-j+\delta k} \Vert u \Vert_X \\
&\lesssim t^{-\frac{7}{6}} t^{-\frac{4}{21}+\delta+2\varrho+\kappa} 2^{-\frac{15j}{14}+\frac{4\delta j}{7}+6\varrho j + \kappa j} \langle 2^j \rangle^{\frac{5}{14}+\frac{\delta}{7}} 2^{-j+\delta k} \Vert u \Vert_X \\
&\lesssim t^{-\frac{7}{6}} t^{-\frac{4}{21}+\frac{45}{616}+\frac{23\delta}{84}+2\varrho+\kappa} \langle 2^j \rangle^{\frac{5}{14}-1+\frac{\delta}{7}} 2^{-j+\delta k} \Vert u \Vert_X \\
&\lesssim t^{-\frac{7}{6}} 2^{-j+\delta k} \langle 2^j \rangle^{-\frac{4}{7}} \Vert u \Vert_X \\
\eqref{estdispsupCloinLpetittau-int-3} &\lesssim t^{-1} 2^{-3k-3j} \Vert \Psi_{l_0}^{int} \psi_{j, k}^{\widehat{\mathcal{C}}} \Vert_{L^{\frac{1}{1-\kappa}}_c L^2_{a, b}} \Vert \psi\left( 2^{-j} \overline{\xi} \right) \overline{\xi} \widehat{g}_{b-corr}(t) \Vert_{L^{\frac{1}{1-\kappa}}_c L^2_{a, b}} \\
&\lesssim t^{-1} 2^{\frac{3l_0}{2}-\kappa l_0-\frac{5k}{2}-j-\kappa j} \left( t^{-\frac{1}{12}+\frac{1}{32}+50\delta} + t^{-\frac{1}{8}} 2^{-\frac{7j}{6}+1000\delta j} \right)^{\frac{1}{2}} \left( t^{-\frac{1}{6}+100\delta} 2^{\frac{j}{2}} \langle 2^j \rangle^{-1} \right)^{\frac{1}{2}} \Vert u \Vert_X \\
&\lesssim t^{-1-\frac{1}{24}-\frac{1}{12}+\frac{1}{64}+75\delta} \mathfrak{t}^{-\frac{3}{4}+\frac{3\varrho}{2}+\frac{\kappa}{2}} \left( t^{-\frac{9}{35}} 2^{-\frac{4j}{7}} \langle 2^j \rangle^{-\frac{1}{7}} \right)^{-\frac{5}{2}-\delta} 2^{-\frac{3j}{4}-\kappa j+\delta k} \langle 2^j \rangle^{-\frac{1}{2}} \Vert u \Vert_X \\
&\quad \quad + t^{-1-\frac{1}{16}-\frac{1}{12}+50\delta} \mathfrak{t}^{-\frac{3}{4}+\frac{3\varrho}{2}+\frac{\kappa}{2}} \left( t^{-\frac{9}{35}} 2^{-\frac{4j}{7}} \langle 2^j \rangle^{-\frac{1}{7}} \right)^{-\frac{5}{2}-\delta} 2^{-\frac{4j}{3}+500\delta j-\kappa j+\delta k} \langle 2^j \rangle^{-\frac{1}{2}} \Vert u \Vert_X \\
&\lesssim t^{-1-\frac{97}{448}+76\delta+\frac{3\varrho}{2}+\frac{\kappa}{2}} 2^{-\frac{11j}{7}+\frac{4\delta j}{7}+\frac{9\varrho j}{2}+\frac{\kappa j}{2}+\delta k} \langle 2^j \rangle^{-\frac{1}{7}+\frac{\delta}{7}} \Vert u \Vert_X \\
&\quad \quad + t^{-1-\frac{85}{336}+51\delta+\frac{3\varrho}{2}+\frac{\kappa}{2}} 2^{-\frac{181j}{84}+500\delta j+\frac{4\delta j}{7}+\frac{9\varrho j}{2}+\frac{\kappa j}{2}+\delta k} \langle 2^j \rangle^{-\frac{1}{7}+\frac{\delta}{7}} \Vert u \Vert_X \\
&\lesssim t^{-\frac{7}{6}} 2^{-\frac{3j}{2}+\delta j+\delta k} \langle 2^j \rangle^{-\frac{1}{7}} \Vert u \Vert_X 
\end{align*}
The other estimate for \eqref{estdispsupCloinLpetittau-int-3} follows using that, here, $t^{-\frac{1}{12}} 2^{-\frac{j}{2}} \lesssim 1$. This concludes the case $\tau \ll 2^k$. 

\paragraph{Away from $\mathcal{L}$, at distance $2^k$ of $\mathcal{C}$} Let us now assume that $|y| \gtrsim |x|$ and $\tau = \sqrt{\frac{\left| |x| - \sqrt{3} |y| \right|}{|x|+|y|}} \simeq 2^k$. 

In this situation, we proceed as before, but we observe that $\widehat{X}_{b-corr} \cdot \nabla_{\overline{\xi}} \Phi$ does not behave anymore like $2^{2k}$ because $\tau^2$ is no longer negligible with respect to $2^{2k}$. In particular, we now add localisations of $\widehat{X}_{b-corr} \cdot \nabla_{\overline{\xi}} \Phi$. 

We apply the following decomposition: 
\begin{subequations}
\begin{align}
\eqref{estdispsupC-termeinitbis} &= \sum_{l_a = l_0}^{-200} 1_{l_a \geq \frac{3k}{2}-\frac{j_{+}}{2}} \int e^{i t \Phi} \Psi_{(l_a, 0, l_a)}^{a-int} \psi_{j, k}^{\widehat{\mathcal{C}}}(\overline{\xi}) \widehat{f}(t, \overline{\xi}) ~ d\overline{\xi} \label{estdispsupCloinLmoyentau-aext} \\
&+ \sum_{l_a = l_0}^{-200} 1_{l_a \geq \frac{3k}{2}-\frac{j_{+}}{2}} \int e^{i t \Phi} \Psi_{(l_a, 0, l_a)}^{c-int} \psi_{j, k}^{\widehat{\mathcal{C}}}(\overline{\xi}) \widehat{f}(t, \overline{\xi}) ~ d\overline{\xi} \label{estdispsupCloinLmoyentau-cext} \\
&+ \sum_{l_a = l_0}^{-200} \sum_{l_b = F(l_a, k, j)}^{3k} 1_{l_a < \frac{3k}{2}-\frac{j_{+}}{2}} \int e^{i t \Phi} \Psi_{(l_a, l_b, l_a)}^{a-b} \psi_{j, k}^{\widehat{\mathcal{C}}}(\overline{\xi}) \widehat{f}(t, \overline{\xi}) ~ d\overline{\xi} \label{estdispsupCloinLmoyentau-ab} \\
&+ \sum_{l_a = l_0}^{-200} 1_{l_a < \frac{3k}{2}-\frac{j_{+}}{2}} \int e^{i t \Phi} \Psi_{(l_a, F(l_a, k, j), l_a)}^{a-int} \psi_{j, k}^{\widehat{\mathcal{C}}}(\overline{\xi}) \widehat{f}(t, \overline{\xi}) ~ d\overline{\xi} \label{estdispsupCloinLmoyentau-aint} \\
&+ \sum_{l_a = l_0}^{-200} \sum_{l_b = F(l_a, k, j)}^{3k} 1_{l_a < \frac{3k}{2}-\frac{j_{+}}{2}} \int e^{i t \Phi} \Psi_{(l_a, l_b, l_a)}^{c-b} \psi_{j, k}^{\widehat{\mathcal{C}}}(\overline{\xi}) \widehat{f}(t, \overline{\xi}) ~ d\overline{\xi} \label{estdispsupCloinLmoyentau-cb} \\
&+ \sum_{l_a = l_0}^{-200} 1_{l_a < \frac{3k}{2}-\frac{j_{+}}{2}} \int e^{i t \Phi} \Psi_{(l_a, F(l_a, k, j), l_a)}^{c-int} \psi_{j, k}^{\widehat{\mathcal{C}}}(\overline{\xi}) \widehat{f}(t, \overline{\xi}) ~ d\overline{\xi} \label{estdispsupCloinLmoyentau-cint} \\
&+ \sum_{l_b = G(l_0, j, k)}^{3k} 1_{l_0 < \frac{3k}{2}-\frac{j_{+}}{2}} \int e^{i t \Phi} \Psi_{(l_0, l_b, l_0)}^{b} \psi_{j, k}^{\widehat{\mathcal{C}}}(\overline{\xi}) \widehat{f}(t, \overline{\xi}) ~ d\overline{\xi} \label{estdispsupCloinLmoyentau-b} \\
&+ \int e^{i t \Phi} \Psi_{(l_0, G(l_0, j, k), l_0)}^{int} \psi_{j, k}^{\widehat{\mathcal{C}}}(\overline{\xi}) \widehat{f}(t, \overline{\xi}) ~ d\overline{\xi} \label{estdispsupCloinLmoyentau-int} 
\end{align}
\end{subequations} 
where the symbols $\Psi_{(l_a, l_b, l_c)}^{*}$, $* \in \{ a-b, a-int, c-b, c-int, b, int \}$ are the ones defined in the proof of Lemma \ref{lem-estdisph-zonereste} by \eqref{equdefPsigenerique}. In fact, we will rather consider symbols where $\widehat{X}_b$ has been replaced by $\widehat{X}_{b-corr}$. We choose $l_0$ such that $2^{l_0} \simeq \mathfrak{t}^{-\frac{1}{2}+\varrho}$ for some small enough $\varrho > 0$. Furthermore, $F(l_a, k, j) = l_a + \frac{3k}{2} + \frac{j_{+}}{2}$, where $j_{+} = \max(0, j)$ is the positive part of $j$. Finally, $G$ is a function whose form we will specify later, satisfying $G(l_0, j, k) = 0$ if $l_0 \geq \frac{3k}{2}-\frac{j_{+}}{2}$. 

Note that, for $* \in \{ a-b, c-b, b \}$, $\Psi_{(l_a, l_b, l_a)}^{*}$ localises $m_b \widehat{X}_{b-corr} \cdot \nabla_{\overline{\xi}} \Phi \simeq 2^{l_b+2j}$, but here we necessarily have 
\begin{align*}
m_b(\overline{\xi}) \widehat{X}_{b-corr}(\overline{\xi}) \cdot \nabla_{\overline{\xi}} \Phi &\simeq 2^k \left( \left( \sqrt{3} |\xi_0| - |\xi| \right)^2 - \widehat{X}_{b-corr}(\overline{\xi}) \cdot \frac{(x, y)}{t} \right) 
\end{align*}
Yet 
\begin{align*}
\widehat{X}_{b-corr}(\overline{\xi}) \cdot \frac{(x, y)}{t} &= \frac{x}{t} - \frac{\sqrt{3} \xi_0 \xi \cdot y}{t |\xi_0| |\xi|} \\
&= \frac{x}{|x|} \frac{|x| - \sqrt{3} |y|}{t} + O\left( \left( \frac{J \xi \cdot y}{|\xi| t} \right)^2 \right) \\
&= \tau^2 2^{2j} + O\left( 2^{2l_a+2j} \right) 
\end{align*}
Therefore, if $l_a < \frac{3k}{2}-\frac{j_{+}}{2} \leq \frac{3k}{2}$, 
\begin{align*}
m_b(\overline{\xi}) \widehat{X}_{b-corr}(\overline{\xi}) \cdot \nabla_{\overline{\xi}} \Phi &\lesssim 2^{3k} 2^{2j} 
\end{align*}
It is thus coherent to only sum up to $l_b = 3k$. 

In a similar way as in the case $\tau \ll 2^k$, we estimate the volumes. First, if $l_b = 0$ et $l_a = l_c \geq \frac{3k}{2}$, then $\Psi_{(l_a, 0, l_a)}$ localises the angle of $\xi$ with a precision $2^{l_a}$, and the size of $\overline{\xi}$ with a precision $2^{l_a+j}$, while $\psi_{j, k}^{\widehat{\mathcal{C}}}$ localises the distance to the cone with a precision $2^{k+j}$, hence
\begin{align*}
\Vert \Psi_{(l_a, 0, l_a)}^{*} \Vert_{L^{p_c}_c L^{p_b}_b L^{p_a}_a} &\lesssim 2^{\frac{l_a}{p_a}+\frac{l_a}{p_c}+\frac{k}{p_b}} |\overline{\xi_a}|^{\frac{1}{p_a}+\frac{1}{p_b}+\frac{1}{p_c}} 
\end{align*}
On the other hand, if $l_a = l_c < \frac{3k}{2}$ and $l_b \leq 3k$, as $m_b \simeq 2^k$, we know that the distance to the cone squared is of size $2^{2k+2j}$ and can be localised with precision $2^{l_b-k+2j} \lesssim 2^{2k+2j}$. Therefore, we compute likewise that 
\begin{align}
\Vert \Psi_{(l_a, l_b, l_a)}^{*} \Vert_{L^{p_c}_c L^{p_b}_b L^{p_a}_a} &\lesssim 2^{\frac{l_a}{p_a}+\frac{l_a}{p_c}+\frac{l_b-2k}{p_b}} |\overline{\xi_a}|^{\frac{1}{p_a}+\frac{1}{p_b}+\frac{1}{p_c}} \label{estvolumiquePsilalbcasCCtausimk}
\end{align}

Finally, we compute directional derivatives of the symbols: 
\begin{align*}
&\widehat{X}_a(\overline{\xi}) \cdot \nabla_{\overline{\xi}} \left[ \widehat{X}_a(\overline{\xi}) \cdot \nabla_{\overline{\xi}} \Phi \right] \simeq 2^j \\
&\widehat{X}_a(\overline{\xi}) \cdot \nabla_{\overline{\xi}} \left[ m_b(\overline{\xi}) \widehat{X}_{b-corr}(\overline{\xi}) \cdot \nabla_{\overline{\xi}} \Phi \right] \lesssim 2^{3k+j} \\
&\widehat{X}_a(\overline{\xi}) \cdot \nabla_{\overline{\xi}} \left[ \widehat{X}_c(\overline{\xi}) \cdot \nabla_{\overline{\xi}} \Phi \right] = 0 \\
&m_b(\overline{\xi}) \widehat{X}_{b-corr}(\overline{\xi}) \cdot \nabla_{\overline{\xi}} \left[ \widehat{X}_a(\overline{\xi}) \cdot \nabla_{\overline{\xi}} \Phi \right] \lesssim 2^{l_a+j} + 2^{3k+j} \\
&m_b(\overline{\xi}) \widehat{X}_{b-corr}(\overline{\xi}) \cdot \nabla_{\overline{\xi}} \left[ m_b(\overline{\xi}) \widehat{X}_{b-corr}(\overline{\xi}) \cdot \nabla_{\overline{\xi}} \Phi \right] \lesssim 2^{3k+j} \\
&m_b(\overline{\xi}) \widehat{X}_{b-corr}(\overline{\xi}) \cdot \nabla_{\overline{\xi}} \left[ \widehat{X}_c(\overline{\xi}) \cdot \nabla_{\overline{\xi}} \Phi \right] \lesssim 2^{l_a+k+j} \\
&\widehat{X}_c(\overline{\xi}) \cdot \nabla_{\overline{\xi}} \left[ \widehat{X}_a(\overline{\xi}) \cdot \nabla_{\overline{\xi}} \Phi \right] \lesssim 2^{l_a+j} \\
&\widehat{X}_c(\overline{\xi}) \cdot \nabla_{\overline{\xi}} \left[ m_b(\overline{\xi}) \widehat{X}_{b-corr}(\overline{\xi}) \cdot \nabla_{\overline{\xi}} \Phi \right] \lesssim 2^{l_a+k+j} \\
&\widehat{X}_c(\overline{\xi}) \cdot \nabla_{\overline{\xi}} \left[ \widehat{X}_c(\overline{\xi}) \cdot \nabla_{\overline{\xi}} \Phi \right] \lesssim 2^j 
\end{align*}
In particular, 
\begin{align*}
\widehat{X}_a(\overline{\xi}) \cdot \nabla_{\overline{\xi}} \Psi_{(l_a, l_b, l_a)} &\lesssim 2^{-l_a-j} + 2^{3k-l_b-j} \\
&\lesssim 2^{-l_a-j} \\
\widehat{X}_c(\overline{\xi}) \cdot \nabla_{\overline{\xi}} \Psi_{(l_a, l_b, l_a)} &\lesssim 2^{-l_a-j} \\
\widehat{X}_a(\overline{\xi}) \cdot \nabla_{\overline{\xi}} \psi_{j, k}^{\widehat{\mathcal{C}}} &\lesssim 2^{-j} \\
\widehat{X}_c(\overline{\xi}) \cdot \nabla_{\overline{\xi}} \psi_{j, k}^{\widehat{\mathcal{C}}} &\lesssim 2^{-j} 
\end{align*}
using that $l_b \geq l_a + \frac{3k}{2}$. Then, if $l_a \geq l_b$ or $l_a \geq 3k$, we have 
\begin{align*}
m_b(\overline{\xi}) \widehat{X}_{b-corr}(\overline{\xi}) \cdot \nabla_{\overline{\xi}} \left( \Psi_{(l_a, l_b, l_a)} \psi_{j, k}^{\widehat{\mathcal{C}}} \right) &\lesssim 2^{3k-l_a-j} + 2^{3k-l_b-j} + 2^{-j} \\
&\lesssim 2^{3k-l_b-j} + 2^{-j} 
\end{align*}
On the other hand, if $l_a \leq l_b \leq 3k$, we may correct $\widehat{X}_{b-corr}$ again so to get 
\begin{align*}
\widehat{X}_{b-corr} \cdot \nabla_{\overline{\xi}} \left[ \widehat{X}_a(\overline{\xi}) \cdot \nabla_{\overline{\xi}} \Phi \right] &= 0 
\end{align*}
It is then easy to check that, since $l_a \leq l_b \leq 3k$, the new vector field $\widehat{X}_{b-corr}$ satisfies the same estimates as the precedent on the support of $\Psi_{(l_a, l_b, l_a)}$. 

On \eqref{estdispsupCloinLmoyentau-aext}, we apply two integrations by parts along $\widehat{X}_a$: 
\begin{subequations}
\begin{align}
\eqref{estdispsupCloinLmoyentau-aext} &= \sum_{l_a = l_0}^{-200} 1_{l_a \geq \frac{3k}{2}-\frac{j_{+}}{2}} t^{-2} 2^{-4l_a-6j} \int e^{i t \Phi} \Psi_{(l_a, 0, l_a)}^{a-int} \psi_{j, k}^{\widehat{\mathcal{C}}}(\overline{\xi}) \widehat{f}(t, \overline{\xi}) ~ d\overline{\xi} \label{estdispsupCloinLmoyentau-aext-1} \\
&\quad + \sum_{l_a = l_0}^{-200} 1_{l_a \geq \frac{3k}{2}-\frac{j_{+}}{2}} t^{-2} 2^{-3l_a-5j} \int e^{i t \Phi} \Psi_{(l_a, 0, l_a)}^{a-int} \psi_{j, k}^{\widehat{\mathcal{C}}}(\overline{\xi}) \widehat{X}_a(\overline{\xi}) \cdot \nabla_{\overline{\xi}} \widehat{f}(t, \overline{\xi}) ~ d\overline{\xi} \label{estdispsupCloinLmoyentau-aext-2} \\
&\quad + \sum_{l_a = l_0}^{-200} 1_{l_a \geq \frac{3k}{2}-\frac{j_{+}}{2}} t^{-2} 2^{-3l_a-5j} \int e^{i t \Phi} \Psi_{(l_a, 0, l_a)}^{a-int} \psi_{j, k}^{\widehat{\mathcal{C}}}(\overline{\xi}) \widehat{h}_a(t, \overline{\xi}) ~ d\overline{\xi} \label{estdispsupCloinLmoyentau-aext-3} \\
&\quad + \sum_{l_a = l_0}^{-200} 1_{l_a \geq \frac{3k}{2}-\frac{j_{+}}{2}} t^{-2} 2^{-2l_a-4j} \int e^{i t \Phi} \Psi_{(l_a, 0, l_a)}^{a-int} \psi_{j, k}^{\widehat{\mathcal{C}}}(\overline{\xi}) \widehat{X}_a(\overline{\xi}) \cdot \nabla_{\overline{\xi}} \widehat{h}_a(t, \overline{\xi}) ~ d\overline{\xi} \label{estdispsupCloinLmoyentau-aext-4} \\
&\quad + \sum_{l_a = l_0}^{-200} 1_{l_a \geq \frac{3k}{2}-\frac{j_{+}}{2}} t^{-1} 2^{-l_a-2j} \int e^{i t \Phi} \Psi_{(l_a, 0, l_a)}^{a-int} \psi_{j, k}^{\widehat{\mathcal{C}}}(\overline{\xi}) \widehat{g}_a(t, \overline{\xi}) ~ d\overline{\xi} \label{estdispsupCloinLmoyentau-aext-5}
\end{align}
\end{subequations}
where as before the symbols may vary from line to line. We estimate: 
\begin{align*}
\eqref{estdispsupCloinLmoyentau-aext-1} &\lesssim \sum_{l_a = l_0}^{-200} 1_{l_a \geq \frac{3k}{2}-\frac{j_{+}}{2}} t^{-2} 2^{-4l_a-6j} \Vert \Psi_{(l_a, 0, l_a)}^{a-int} \psi_{j, k}^{\widehat{\mathcal{C}}} \Vert_{L^{\frac{1}{1-\kappa}}_c L^2_b L^{\frac{1}{1-\kappa}}_a} \langle 2^j \rangle^{-1} \Vert \langle \overline{\xi} \rangle \widehat{f}(t) \Vert_{L^{\frac{1}{\kappa}}_c L^2_b L^{\frac{1}{\kappa}}_a} \\
&\lesssim \sum_{l_a = l_0}^{-200} 1_{l_a \geq \frac{3k}{2}-\frac{j_{+}}{2}} t^{-\frac{7}{6}+\frac{2\kappa}{3}} \mathfrak{t}^{-\frac{5}{6}-\frac{2\kappa}{3}} 2^{-2l_a-2\kappa l_a-j+\frac{k}{2}} \langle 2^j \rangle^{-1} \Vert u \Vert_X \\
&\lesssim \sum_{l_a = l_0}^{-200} t^{-\frac{7}{6}+\frac{2\kappa}{3}} \mathfrak{t}^{-\frac{5}{6}-\frac{2\kappa}{3}} 2^{-\frac{5l_a}{3}-\frac{2\delta l_a}{3}-2\kappa l_a-j+\delta k} \langle 2^j \rangle^{-1+\frac{1}{6}-\frac{\delta}{3}} \Vert u \Vert_X \\
&\lesssim t^{-\frac{7}{6}+\frac{2\kappa}{3}} \mathfrak{t}^{\frac{\delta}{3}+\frac{\kappa}{3}} 2^{-j+\delta k} \langle 2^j \rangle^{-\frac{5}{6}-\frac{\delta}{3}} \Vert u \Vert_X \\
&\lesssim t^{-\frac{7}{6}+\frac{\delta}{3}+\kappa} 2^{-j+\delta k} \langle 2^j \rangle^{-\frac{5}{6}+\frac{2\delta}{3}+\kappa} \Vert u \Vert_X \\
\eqref{estdispsupCloinLmoyentau-aext-2} &\lesssim \sum_{l_a = l_0}^{-200} 1_{l_a \geq \frac{3k}{2}-\frac{j_{+}}{2}} t^{-2} 2^{-3l_a-5j} \Vert \Psi_{(l_a, 0, l_a)}^{a-int} \psi_{j, k}^{\widehat{\mathcal{C}}} \Vert_{L^2} \langle 2^j \rangle^{-1} \Vert \langle \nabla \rangle X_a f(t) \Vert_{L^2} \\
&\lesssim \sum_{l_a = l_0}^{-200} 1_{l_a \geq \frac{3k}{2}-\frac{j_{+}}{2}} t^{-\frac{7}{6}} \mathfrak{t}^{-\frac{5}{6}} 2^{-2l_a+\frac{k}{2}-j} \langle 2^j \rangle^{-1} \Vert u \Vert_X \\
&\lesssim t^{-\frac{7}{6}+\frac{\delta}{3}} 2^{-j+\delta k} \langle 2^j \rangle^{-\frac{5}{6}+\frac{2\delta}{3}} \Vert u \Vert_X \\
\eqref{estdispsupCloinLmoyentau-aext-4} &\lesssim \sum_{l_a = l_0}^{-200} 1_{l_a \geq \frac{3k}{2}-\frac{j_{+}}{2}} t^{-2} 2^{-2l_a-5j} \Vert \Psi_{(l_a, 0, l_a)}^{a-int} \psi_{j, k}^{\widehat{\mathcal{C}}} \Vert_{L^2} \Vert \nabla X_a h_a(t) \Vert_{L^2} \\
&\lesssim \sum_{l_a = l_0}^{-200} 1_{l_a \geq \frac{3k}{2}-\frac{j_{+}}{2}} t^{-\frac{7}{6}} \mathfrak{t}^{-\frac{5}{6}} 2^{-l_a-j+\frac{k}{2}} \Vert u \Vert_X \\
&\lesssim t^{-\frac{7}{6}} \mathfrak{t}^{-\frac{1}{3}} 2^{-j+\frac{k}{2}} \Vert u \Vert_X \\
&\lesssim t^{-\frac{7}{6}} 2^{-j+\frac{k}{2}} \langle 2^j \rangle^{-1} \Vert u \Vert_X \\
\eqref{estdispsupCloinLmoyentau-aext-5} &\lesssim \sum_{l_a = l_0}^{-200} 1_{l_a \geq \frac{3k}{2}-\frac{j_{+}}{2}} t^{-1} 2^{-l_a-\frac{5j}{2}} \Vert \Psi_{(l_a, 0, l_a)}^{a-int} \psi_{j, k}^{\widehat{\mathcal{C}}} \Vert_{L^2} \langle 2^j \rangle^{-1} \Vert m_{\widehat{\mathcal{C}}}(D) \langle \nabla \rangle |\nabla|^{\frac{1}{2}} g_a(t) \Vert_{L^2} \\
&\lesssim \sum_{l_a = l_0}^{-200} 1_{l_a \geq \frac{3k}{2}-\frac{j_{+}}{2}} t^{-\frac{7}{6}+100\delta} 2^{-j+\frac{k}{2}} \langle 2^j \rangle^{-1} \Vert u \Vert_X \\
&\lesssim t^{-\frac{7}{6}+100\delta} 2^{-j+\frac{k}{4}} \langle 2^j \rangle^{-\frac{5}{6}} \Vert u \Vert_X 
\end{align*}
\eqref{estdispsupCloinLmoyentau-aext-3} is similar to \eqref{estdispsupCloinLmoyentau-aext-2}. 

Pour \eqref{estdispsupCloinLmoyentau-ab}, we distinguish $l_a \geq 6k-\frac{j_{+}}{2}$ and $l_a \leq 6k - \frac{j_{+}}{2}$: 
\begin{subequations}
\begin{align}
\eqref{estdispsupCloinLmoyentau-ab} &= \sum_{l_a = l_0}^{-200} \sum_{l_b = F(l_a, k, j)}^{3k} 1_{6k \leq l_a + \frac{j_{+}}{2} < \frac{3k}{2}} \int e^{i t \Phi} \Psi_{(l_a, l_b, l_a)}^{a-b} \psi_{j, k}^{\widehat{\mathcal{C}}}(\overline{\xi}) \widehat{f}(t, \overline{\xi}) ~ d\overline{\xi} \label{estdispsupCloinLmoyentau-ab1} \\
&\quad + \sum_{l_a = l_0}^{-200} \sum_{l_b = F(l_a, k, j)}^{3k} 1_{l_a +\frac{j_{+}}{2} < 6k} \int e^{i t \Phi} \Psi_{(l_a, l_b, l_a)}^{a-b} \psi_{j, k}^{\widehat{\mathcal{C}}}(\overline{\xi}) \widehat{f}(t, \overline{\xi}) ~ d\overline{\xi} \label{estdispsupCloinLmoyentau-ab2} 
\end{align}
\end{subequations}
On \eqref{estdispsupCloinLmoyentau-ab}, we apply at most $n$ (large with respect to $\varrho^{-1}$) integrations by parts along $\widehat{X}_a$, stopping as soon as a derivative hits $\widehat{f}(t)$, and then possibly one along $\widehat{X}_{b-corr}$: 
\begin{subequations}
\begin{align}
\eqref{estdispsupCloinLmoyentau-ab1} &= \sum_{l_a = l_0}^{-200} \sum_{l_b = F(l_a, k, j)}^{3k} 1_{6k \leq l_a + \frac{j_{+}}{2} < \frac{3k}{2}} t^{-n} 2^{-2nl_a-3nj} \int e^{i t \Phi} \Psi_{(l_a, l_b, l_a)}^{a-b} \psi_{j, k}^{\widehat{\mathcal{C}}}(\overline{\xi}) \widehat{f}(t, \overline{\xi}) ~ d\overline{\xi} \label{estdispsupCloinLmoyentau-ab1-1} \\
&\quad \begin{aligned}
+ \sum_{i = 1}^n \sum_{l_a = l_0}^{-200} \sum_{l_b = F(l_a, k, j)}^{3k} 1_{6k \leq l_a + \frac{j_{+}}{2} < \frac{3k}{2}} &t^{-i-1} 2^{-2il_a+l_a-2l_b+3k-3ij-2j} \\
&\int e^{i t \Phi} \Psi_{(l_a, l_b, l_a)}^{a-b} \psi_{j, k}^{\widehat{\mathcal{C}}}(\overline{\xi}) \widehat{h}_a(t, \overline{\xi}) ~ d\overline{\xi} 
\end{aligned}\label{estdispsupCloinLmoyentau-ab1-2} \\
&\quad \begin{aligned}
+ \sum_{i = 1}^n &\sum_{l_a = l_0}^{-200} \sum_{l_b = F(l_a, k, j)}^{3k} 1_{6k \leq l_a + \frac{j_{+}}{2} < \frac{3k}{2}} t^{-i-1} 2^{-2il_a+l_a-l_b-3ij-j} \\
&\int e^{i t \Phi} \Psi_{(l_a, l_b, l_a)}^{a-b} \psi_{j, k}^{\widehat{\mathcal{C}}}(\overline{\xi}) m_b(\overline{\xi}) \widehat{X}_{b-corr}(\overline{\xi}) \cdot \nabla_{\overline{\xi}} \widehat{h}_a(t, \overline{\xi}) ~ d\overline{\xi}
\end{aligned} \label{estdispsupCloinLmoyentau-ab1-3} \\
&\quad \begin{aligned}
+ \sum_{i = 1}^n \sum_{l_a = l_0}^{-200} \sum_{l_b = F(l_a, k, j)}^{3k} 1_{6k \leq l_a + \frac{j_{+}}{2} < \frac{3k}{2}} &t^{-i} 2^{-2il_a+l_a-3ij+j} \\
&\int e^{i t \Phi} \Psi_{(l_a, l_b, l_a)}^{a-b} \psi_{j, k}^{\widehat{\mathcal{C}}}(\overline{\xi}) \widehat{g}_a(t, \overline{\xi}) ~ d\overline{\xi} 
\end{aligned} \label{estdispsupCloinLmoyentau-ab1-4} 
\end{align}
\end{subequations}
We then estimate: 
\begin{align*}
\eqref{estdispsupCloinLmoyentau-ab1-1} &\lesssim \sum_{l_a = l_0}^{-200} \sum_{l_b = F(l_a, k, j)}^{3k} 1_{6k \leq l_a + \frac{j_{+}}{2} < \frac{3k}{2}} \mathfrak{t}^{-n} 2^{-2nl_a+j} \Vert \Psi_{(l_a, l_b, l_a)}^{a-b} \psi_{j, k}^{\widehat{\mathcal{C}}} \Vert_{L^2} \langle 2^j \rangle^{-1} \Vert \langle \nabla \rangle |\nabla|^{-1} f(t) \Vert_{L^2} \\
&\lesssim \sum_{l_a = l_0}^{-200} \sum_{l_b = F(l_a, k, j)}^{3k} 1_{6k \leq l_a + \frac{j_{+}}{2} < \frac{3k}{2}} t^{-\frac{7}{6}} \mathfrak{t}^{-n+\frac{7}{6}} 2^{-2nl_a+l_a+\frac{l_b}{2}-k-j} \langle 2^j \rangle^{-1} \Vert u \Vert_X \\ 
&\lesssim t^{-\frac{7}{6}} \mathfrak{t}^{\frac{7}{6}-2n\varrho} 2^{\frac{k}{2}-k-j} \langle 2^j \rangle^{-1} \Vert u \Vert_X \\
&\lesssim t^{-\frac{7}{6}} 2^{-j+\frac{k}{2}} \langle 2^j \rangle^{-1} \Vert u \Vert_X \\
\eqref{estdispsupCloinLmoyentau-ab1-2} &\lesssim \sum_{i = 1}^n \sum_{l_a = l_0}^{-200} \sum_{l_b = F(l_a, k, j)}^{3k} 1_{6k \leq l_a + \frac{j_{+}}{2} < \frac{3k}{2}} t^{-1} \mathfrak{t}^{-i} 2^{-2il_a+l_a-2l_b+3k-2j} \\
&\pushright{\Vert \Psi_{(l_a, l_b, l_a)}^{a-b} \psi_{j, k}^{\widehat{\mathcal{C}}} \Vert_{L^2} \langle 2^j \rangle^{-1} \Vert \langle \nabla \rangle h_a(t) \Vert_{L^2}} \\
&\lesssim \sum_{i = 1}^n \sum_{l_a = l_0}^{-200} \sum_{l_b = F(l_a, k, j)}^{3k} 1_{6k \leq l_a + \frac{j_{+}}{2} < \frac{3k}{2}} t^{-\frac{7}{6}} \mathfrak{t}^{-i+\frac{1}{6}} 2^{-2il_a+2l_a-\frac{3l_b}{2}+2k-j} \langle 2^j \rangle^{-1} \Vert u \Vert_X \\
&\lesssim \sum_{i = 1}^n \sum_{l_a = l_0}^{-200} 1_{6k \leq l_a + \frac{j_{+}}{2} < \frac{3k}{2}} t^{-\frac{7}{6}} \mathfrak{t}^{-i+\frac{1}{6}} 2^{-2il_a+\frac{l_a}{2}-\frac{k}{4}-j} \langle 2^j \rangle^{-1} \Vert u \Vert_X \\
&\lesssim \sum_{i = 1}^n \sum_{l_a = l_0}^{-200} t^{-\frac{7}{6}} \mathfrak{t}^{-i+\frac{1}{6}} 2^{-2il_a+\frac{l_a}{3}-\frac{2\delta l_a}{3}-j+\delta k} \langle 2^j \rangle^{-\frac{13}{12}-\frac{\delta}{3}} \Vert u \Vert_X \\
&\lesssim t^{-\frac{7}{6}} \mathfrak{t}^{\frac{\delta}{3}} 2^{-j+\delta k} \langle 2^j \rangle^{-\frac{13}{12}-\frac{\delta}{3}} \Vert u \Vert_X \\
&\lesssim t^{-\frac{7}{6}+\frac{\delta}{3}} 2^{-j+\delta k} \langle 2^j \rangle^{-\frac{13}{12}+\frac{2\delta}{3}} \Vert u \Vert_X \\
\eqref{estdispsupCloinLmoyentau-ab1-3} &\lesssim \sum_{i = 1}^n \sum_{l_a = l_0}^{-200} \sum_{l_b = F(l_a, k, j)}^{3k} 1_{6k \leq l_a + \frac{j_{+}}{2} < \frac{3k}{2}} t^{-1} \mathfrak{t}^{-i} 2^{-2il_a+l_a-l_b-2j} \\
&\pushright{\Vert \Psi_{(l_a, l_b, l_a)} \psi_{j, k}^{\widehat{\mathcal{C}}} \Vert_{L^2} \Vert \nabla m_b(D) X_{b-corr} h_a(t) \Vert_{L^2}} \\
&\lesssim \sum_{i = 1}^n \sum_{l_a = l_0}^{-200} \sum_{l_b = F(l_a, k, j)}^{3k} 1_{6k \leq l_a + \frac{j_{+}}{2} < \frac{3k}{2}} t^{-\frac{7}{6}} \mathfrak{t}^{-i+\frac{1}{6}} 2^{-2il_a+2l_a-\frac{l_b}{2}-k-j} \Vert u \Vert_X \\
&\lesssim \sum_{i = 1}^n \sum_{l_a = l_0}^{-200} 1_{6k \leq l_a + \frac{j_{+}}{2} < \frac{3k}{2}} t^{-\frac{7}{6}} \mathfrak{t}^{-i+\frac{1}{6}} 2^{-2il_a+\frac{3l_a}{2}-\frac{7k}{4}-j} \langle 2^j \rangle^{-\frac{1}{4}} \Vert u \Vert_X \\
&\lesssim t^{-\frac{7}{6}} \mathfrak{t}^{\frac{\delta}{3}} 2^{-j+\delta k} \langle 2^j \rangle^{-\frac{5}{6}-\frac{\delta}{3}} \Vert u \Vert_X \\
&\lesssim t^{-\frac{7}{6}+\frac{\delta}{3}} 2^{-j+\delta k} \langle 2^j \rangle^{-\frac{5}{6}+\frac{2\delta}{3}} \Vert u \Vert_X \\
\eqref{estdispsupCloinLmoyentau-ab1-4} &\lesssim \sum_{i = 1}^n \sum_{l_a = l_0}^{-200} \sum_{l_b = F(l_a, k, j)}^{3k} 1_{6k \leq l_a + \frac{j_{+}}{2} < \frac{3k}{2}} \mathfrak{t}^{-i} 2^{-2il_a+l_a+\frac{j}{2}} \\
&\pushright{\Vert \Psi_{(l_a, l_b, l_a)}^{a-b} \psi_{j, k}^{\widehat{\mathcal{C}}} \Vert_{L^2} \langle 2^j \rangle^{-1} \Vert m_{\widehat{\mathcal{C}}}(D) \langle \nabla \rangle |\nabla|^{\frac{1}{2}} g_a(t) \Vert_{L^2}} \\
&\lesssim \sum_{i = 1}^n \sum_{l_a = l_0}^{-200} \sum_{l_b = F(l_a, k, j)}^{3k} 1_{6k \leq l_a + \frac{j_{+}}{2} < \frac{3k}{2}} t^{-\frac{7}{6}+100\delta} \mathfrak{t}^{-i+1} 2^{-2il_a+2l_a+\frac{l_b}{2}-k-j} \langle 2^j \rangle^{-1} \Vert u \Vert_X \\
&\lesssim \sum_{i = 1}^n \sum_{l_a = l_0}^{-200} 1_{3k \leq l_a + \frac{j_{+}}{2} < \frac{3k}{2}} t^{-\frac{7}{6}+100\delta} \mathfrak{t}^{-i+1} 2^{-2il_a+2l_a+\frac{k}{2}-j} \langle 2^j \rangle^{-1} \Vert u \Vert_X \\
&\lesssim t^{-\frac{7}{6}+100\delta} 2^{\frac{k}{3}-j} \langle 2^j \rangle^{-1} \Vert u \Vert_X
\end{align*}
for $n \varrho$ large enough with respect to $1$. On the last line, we used that, if $i \geq 2$, the sum in $l_a$ is geometric, while if $i = 1$, we lose $2^{-\kappa l_a} \lesssim 2^{-6\kappa k}$ for arbitrary $\kappa > 0$ (up to an universal constant). 

Then, for \eqref{estdispsupCloinLmoyentau-ab2}, we apply first an integration by parts along $\widehat{X}_{b-corr}$, then possibly an integration by parts along $\widehat{X}_a$: 
\begin{subequations}
\begin{align}
&\eqref{estdispsupCloinLmoyentau-ab2} = \sum_{l_a = l_0}^{-200} \sum_{l_b = F(l_a, k, j)}^{3k} 1_{l_a +\frac{j_{+}}{2} < 6k} t^{-2} 2^{-2l_b+3k-2l_a-6j} \int e^{i t \Phi} \Psi_{(l_a, l_b, l_a)}^{a-b} \psi_{j, k}^{\widehat{\mathcal{C}}}(\overline{\xi}) \widehat{f}(t, \overline{\xi}) ~ d\overline{\xi} \label{estdispsupCloinLmoyentau-ab2-1} \\
&\quad + \sum_{l_a = l_0}^{-200} \sum_{l_b = F(l_a, k, j)}^{3k} 1_{l_a +\frac{j_{+}}{2} < 6k} t^{-2} 2^{-2l_b+3k-l_a-5j} \int e^{i t \Phi} \Psi_{(l_a, l_b, l_a)}^{a-b} \psi_{j, k}^{\widehat{\mathcal{C}}}(\overline{\xi}) \widehat{X}_a(\overline{\xi}) \cdot \nabla_{\overline{\xi}} \widehat{f}(t, \overline{\xi}) ~ d\overline{\xi} \label{estdispsupCloinLmoyentau-ab2-2} \\
&\quad + \sum_{l_a = l_0}^{-200} \sum_{l_b = F(l_a, k, j)}^{3k} 1_{l_a +\frac{j_{+}}{2} < 6k} t^{-2} 2^{-l_b-2l_a-5j} \int e^{i t \Phi} \Psi_{(l_a, l_b, l_a)}^{a-b} \psi_{j, k}^{\widehat{\mathcal{C}}}(\overline{\xi}) \widehat{h}_{b-corr}(t, \overline{\xi}) ~ d\overline{\xi} \label{estdispsupCloinLmoyentau-ab2-3} \\
&\quad + \sum_{l_a = l_0}^{-200} \sum_{l_b = F(l_a, k, j)}^{3k} 1_{l_a +\frac{j_{+}}{2} < 6k} t^{-2} 2^{-l_b-l_a-4j} \int e^{i t \Phi} \Psi_{(l_a, l_b, l_a)}^{a-b} \psi_{j, k}^{\widehat{\mathcal{C}}}(\overline{\xi}) \widehat{X}_a(\overline{\xi}) \cdot \nabla_{\overline{\xi}} \widehat{h}_{b-corr}(t, \overline{\xi}) ~ d\overline{\xi} \label{estdispsupCloinLmoyentau-ab2-4} \\
&\quad + \sum_{l_a = l_0}^{-200} \sum_{l_b = F(l_a, k, j)}^{3k} 1_{l_a +\frac{j_{+}}{2} < 6k} t^{-1} 2^{-l_b-2j} \int e^{i t \Phi} \Psi_{(l_a, l_b, l_a)}^{a-b} \psi_{j, k}^{\widehat{\mathcal{C}}}(\overline{\xi}) \widehat{g}_{b-corr}(t, \overline{\xi}) ~ d\overline{\xi} \label{estdispsupCloinLmoyentau-ab2-5} 
\end{align}
\end{subequations}
We then estimate: 
\begin{align*}
\eqref{estdispsupCloinLmoyentau-ab2-1} &\lesssim \sum_{l_a = l_0}^{-200} \sum_{l_b = F(l_a, k, j)}^{3k} 1_{l_a +\frac{j_{+}}{2} < 6k} t^{-2} 2^{-2l_b+3k-2l_a-6j} \\
&\pushright{\Vert \Psi_{(l_a, l_b, l_a)}^{a-b} \psi_{j, k}^{\widehat{\mathcal{C}}} \Vert_{L^{\frac{1}{1-\kappa}}_c L^2_b L^{\frac{1}{1-\kappa}}_a} \langle 2^j \rangle^{-1} \Vert \langle \overline{\xi} \rangle \widehat{f}(t) \Vert_{L^{\frac{1}{\kappa}}_c L^2_b L^{\frac{1}{\kappa}}_a}} \\
&\lesssim \sum_{l_a = l_0}^{-200} \sum_{l_b = F(l_a, k, j)}^{3k} 1_{l_a +\frac{j_{+}}{2} < 6k} t^{-\frac{7}{6}+\frac{2\kappa}{3}} \mathfrak{t}^{-\frac{5}{6}-\frac{2\kappa}{3}} 2^{-\frac{3l_b}{2}+2k-2 \kappa l_a-j} \langle 2^j \rangle^{-1} \Vert u \Vert_X \\
&\lesssim \sum_{l_a = l_0}^{-200} 1_{l_a +\frac{j_{+}}{2} < 6k} t^{-\frac{7}{6}+\frac{2\kappa}{3}} \mathfrak{t}^{-\frac{5}{6}-\frac{2\kappa}{3}} 2^{-\frac{3l_a}{2}-\frac{k}{4}-2 \kappa l_a-j} \langle 2^j \rangle^{-1} \Vert u \Vert_X \\
&\lesssim t^{-\frac{7}{6}+\frac{2\kappa}{3}} \mathfrak{t}^{-\frac{1}{16}+\frac{\delta}{6}+\frac{\kappa}{3}} 2^{-j+\delta k} \langle 2^j \rangle^{-1} \Vert u \Vert_X \\
&\lesssim t^{-\frac{7}{6}+\frac{2\kappa}{3}} 2^{-j+\delta k} \langle 2^j \rangle^{-1} \Vert u \Vert_X \\
\eqref{estdispsupCloinLmoyentau-ab2-2} &\lesssim \sum_{l_a = l_0}^{-200} \sum_{l_b = F(l_a, k, j)}^{3k} 1_{l_a +\frac{j_{+}}{2} < 6k} t^{-2} 2^{-2l_b+3k-l_a-5j} \Vert \Psi_{(l_a, l_b, l_a)}^{a-b} \Vert_{L^2} \langle 2^j \rangle^{-1} \Vert \langle \nabla \rangle X_a f(t) \Vert_{L^2} \\
&\lesssim \sum_{l_a = l_0}^{-200} \sum_{l_b = F(l_a, k, j)}^{3k} 1_{l_a +\frac{j_{+}}{2} < 6k} t^{-\frac{7}{6}} \mathfrak{t}^{-\frac{5}{6}} 2^{-\frac{3l_b}{2}+2k-j} \langle 2^j \rangle^{-1} \Vert u \Vert_X \\
&\lesssim t^{-\frac{7}{6}} 2^{-j+\delta k} \langle 2^j \rangle^{-1} \Vert u \Vert_X \\
\eqref{estdispsupCloinLmoyentau-ab2-3} &\lesssim \sum_{l_a = l_0}^{-200} \sum_{l_b = F(l_a, k, j)}^{3k} 1_{l_a +\frac{j_{+}}{2} < 6k} t^{-2} 2^{-l_b-2l_a-6j} \\
&\pushright{\Vert \Psi_{(l_a, l_b, l_a)}^{a-b} \psi_{j, k}^{\widehat{\mathcal{C}}} \Vert_{L^{\frac{1}{1-\kappa}}_c L^2_b L^{\frac{1}{1-\kappa}}_a} \Vert \overline{\xi} \widehat{h}_{b-corr}(t) \Vert_{L^{\frac{1}{\kappa}}_c L^2_b L^{\frac{1}{\kappa}}_a}} \\
&\lesssim \sum_{l_a = l_0}^{-200} \sum_{l_b = F(l_a, k, j)}^{3k} 1_{l_a +\frac{j_{+}}{2} < 6k} t^{-\frac{7}{6}+\frac{2\kappa}{3}} \mathfrak{t}^{-\frac{5}{6}-\frac{2\kappa}{3}} 2^{-\frac{l_b}{2}-k-2\kappa l_a-j} \Vert u \Vert_X \\
&\lesssim \sum_{l_a = l_0}^{-200} 1_{l_a +\frac{j_{+}}{2} < 6k} t^{-\frac{7}{6}+\frac{2\kappa}{3}} \mathfrak{t}^{-\frac{5}{6}-\frac{2\kappa}{3}} 2^{-\frac{l_a}{2}-\frac{7k}{4}-2\kappa l_a-j} \langle 2^j \rangle^{-\frac{1}{4}} \Vert u \Vert_X \\
&\lesssim \sum_{l_a = l_0}^{-200} t^{-\frac{7}{6}+\frac{2\kappa}{3}} \mathfrak{t}^{-\frac{5}{6}-\frac{2\kappa}{3}} 2^{-\frac{l_a}{2}-\frac{7l_a}{24}-\frac{\delta l_a}{6}-2\kappa l_a-j+\delta k} \langle 2^j \rangle^{-\frac{1}{4}-\frac{7}{48}-\frac{\delta}{12}} \Vert u \Vert_X \\
&\lesssim t^{-\frac{7}{6}+\frac{2\kappa}{3}} \mathfrak{t}^{-\frac{7}{16}+\frac{\delta}{12}+\frac{\kappa}{3}} 2^{-j+\delta k} \langle 2^j \rangle^{-\frac{19}{48}-\frac{\delta}{12}} \Vert u \Vert_X \\
&\lesssim t^{-\frac{7}{6}+\frac{2\kappa}{3}} 2^{-j+\delta k} \langle 2^j \rangle^{-1} \Vert u \Vert_X \\
\eqref{estdispsupCloinLmoyentau-ab2-4} &\lesssim \sum_{l_a = l_0}^{-200} \sum_{l_b = F(l_a, k, j)}^{3k} 1_{l_a +\frac{j_{+}}{2} < 6k} t^{-2} 2^{-l_b-l_a-5j} \Vert \Psi_{(l_a, l_b, l_a)}^{a-b} \psi_{j, k}^{\widehat{\mathcal{C}}} \Vert_{L^2} \Vert \nabla X_a h_{b-corr}(t) \Vert_{L^2} \\
&\lesssim \sum_{l_a = l_0}^{-200} \sum_{l_b = F(l_a, k, j)}^{3k} 1_{l_a +\frac{j_{+}}{2} < 6k} t^{-\frac{7}{6}} \mathfrak{t}^{-\frac{5}{6}} 2^{-\frac{l_b}{2}-k-j} \Vert u \Vert_X \\
&\lesssim t^{-\frac{7}{6}} 2^{-j+\delta k} \langle 2^j \rangle^{-1} \Vert u \Vert_X \\
\eqref{estdispsupCloinLmoyentau-ab2-5} &\lesssim \sum_{l_a = l_0}^{-200} \sum_{l_b = F(l_a, k, j)}^{3k} 1_{l_a +\frac{j_{+}}{2} < 6k} t^{-1} 2^{-l_b-\frac{5j}{2}} \\
&\pushright{\Vert \Psi_{(l_a, l_b, l_a)}^{a-b} \psi_{j, k}^{\widehat{\mathcal{C}}} \Vert_{L^2} \langle 2^j \rangle^{-1} \Vert m_{\widehat{\mathcal{C}}}(D) \langle \nabla \rangle |\nabla|^{\frac{1}{2}} g_{b-corr}(t) \Vert_{L^2}} \\
&\lesssim \sum_{l_a = l_0}^{-200} \sum_{l_b = F(l_a, k, j)}^{3k} 1_{l_a +\frac{j_{+}}{2} < 6k} t^{-\frac{7}{6}+100\delta} 2^{l_a-\frac{l_b}{2}-k-j} \langle 2^j \rangle^{-1} \Vert u \Vert_X \\
&\lesssim \sum_{l_a = l_0}^{-200} 1_{l_a +\frac{j_{+}}{2} < 6k} t^{-\frac{7}{6}+100\delta} 2^{\frac{l_a}{2}-\frac{7k}{4}-j} \langle 2^j \rangle^{-1} \Vert u \Vert_X \\
&\lesssim \sum_{l_a = l_0}^{-200} t^{-\frac{7}{6}+100\delta} 2^{\frac{l_a}{2}-\frac{7l_a}{24}-\frac{\delta l_a}{6}-j+\delta k} \langle 2^j \rangle^{-1} \Vert u \Vert_X \\
&\lesssim t^{-\frac{7}{6}+100\delta} 2^{-j+\delta k} \langle 2^j \rangle^{-1} \Vert u \Vert_X 
\end{align*}
if $\delta, \kappa$ are small enough.  

Then, for \eqref{estdispsupCloinLmoyentau-aint}, we apply two integrations by parts along $\widehat{X}_a$: 
\begin{subequations}
\begin{align}
\eqref{estdispsupCloinLmoyentau-aint} &= \sum_{l_a = l_0}^{-200} 1_{l_a < \frac{3k}{2}-\frac{j_{+}}{2}} t^{-2} 2^{-4l_a-6j} \int e^{i t \Phi} \Psi_{(l_a, F(l_a, k, j), l_a)}^{a-int} \psi_{j, k}^{\widehat{\mathcal{C}}}(\overline{\xi}) \widehat{f}(t, \overline{\xi}) ~ d\overline{\xi} \label{estdispsupCloinLmoyentau-aint-1} \\
&\quad + \sum_{l_a = l_0}^{-200} 1_{l_a < \frac{3k}{2}-\frac{j_{+}}{2}} t^{-2} 2^{-3l_a-5j} \int e^{i t \Phi} \Psi_{(l_a, F(l_a, k, j), l_a)}^{a-int} \psi_{j, k}^{\widehat{\mathcal{C}}}(\overline{\xi}) \widehat{X}_a(\overline{\xi}) \cdot \nabla_{\overline{\xi}} \widehat{f}(t, \overline{\xi}) ~ d\overline{\xi} \label{estdispsupCloinLmoyentau-aint-2} \\
&\quad + \sum_{l_a = l_0}^{-200} 1_{l_a < \frac{3k}{2}-\frac{j_{+}}{2}} t^{-2} 2^{-3l_a-5j} \int e^{i t \Phi} \Psi_{(l_a, F(l_a, k, j), l_a)}^{a-int} \psi_{j, k}^{\widehat{\mathcal{C}}}(\overline{\xi}) \widehat{h}_a(t, \overline{\xi}) ~ d\overline{\xi} \label{estdispsupCloinLmoyentau-aint-3} \\
&\quad + \sum_{l_a = l_0}^{-200} 1_{l_a < \frac{3k}{2}-\frac{j_{+}}{2}} t^{-2} 2^{-2l_a-4j} \int e^{i t \Phi} \Psi_{(l_a, F(l_a, k, j), l_a)}^{a-int} \psi_{j, k}^{\widehat{\mathcal{C}}}(\overline{\xi}) \widehat{X}_a(\overline{\xi}) \cdot \nabla_{\overline{\xi}} \widehat{h}_a(t, \overline{\xi}) ~ d\overline{\xi} \label{estdispsupCloinLmoyentau-aint-4} \\
&\quad + \sum_{l_a = l_0}^{-200} 1_{l_a < \frac{3k}{2}-\frac{j_{+}}{2}} t^{-1} 2^{-l_a-2j} \int e^{i t \Phi} \Psi_{(l_a, F(l_a, k, j), l_a)}^{a-int} \psi_{j, k}^{\widehat{\mathcal{C}}}(\overline{\xi}) \widehat{g}_a(t, \overline{\xi}) ~ d\overline{\xi} \label{estdispsupCloinLmoyentau-aint-5} 
\end{align}
\end{subequations}
We then estimate: 
\begin{align*}
\eqref{estdispsupCloinLmoyentau-aint-1} &\lesssim \sum_{l_a = l_0}^{-200} 1_{l_a < \frac{3k}{2}-\frac{j_{+}}{2}} t^{-2} 2^{-4l_a-6j} \Vert \Psi_{(l_a, F(l_a, k, j), l_a)} \Vert_{L^{\frac{1}{1-\kappa}}_c L^2_b L^{\frac{1}{1-\kappa}}_a} \langle 2^j \rangle^{-1} \Vert \langle \overline{\xi} \rangle \widehat{f}(t) \Vert_{L^{\frac{1}{\kappa}}_c L^2_b L^{\frac{1}{\kappa}}_a} \\
&\lesssim \sum_{l_a = l_0}^{-200} 1_{l_a < \frac{3k}{2}-\frac{j_{+}}{2}} t^{-\frac{7}{6}+\frac{2\kappa}{3}} \mathfrak{t}^{-\frac{5}{6}-\frac{2\kappa}{3}} 2^{-\frac{3l_a}{2}-2\kappa l_a-\frac{k}{4}-j} \langle 2^j \rangle^{-\frac{3}{4}} \Vert u \Vert_X \\
&\lesssim t^{-\frac{7}{6}+\frac{2\kappa}{3}} \mathfrak{t}^{\frac{\delta}{3}+\frac{\kappa}{3}} 2^{-j+\delta k} \langle 2^j \rangle^{-\frac{5}{6}-\frac{\delta}{3}} \Vert u \Vert_X \\
&\lesssim t^{-\frac{7}{6}+\frac{\delta}{3}+\kappa} 2^{-j+\delta k} \langle 2^j \rangle^{-\frac{5}{6}+\frac{2\delta}{3}+\kappa} \Vert u \Vert_X \\
\eqref{estdispsupCloinLmoyentau-aint-2} &\lesssim \sum_{l_a = l_0}^{-200} 1_{l_a < \frac{3k}{2}-\frac{j_{+}}{2}} t^{-2} 2^{-3l_a-5j} \Vert \Psi_{(l_a, F(l_a, k, j), l_a)} \psi_{j, k}^{\widehat{\mathcal{C}}} \Vert_{L^2} \langle 2^j \rangle^{-1} \Vert \langle \nabla \rangle X_a f(t) \Vert_{L^2} \\
&\lesssim \sum_{l_a = l_0}^{-200} 1_{l_a < \frac{3k}{2}-\frac{j_{+}}{2}} t^{-\frac{7}{6}} \mathfrak{t}^{-\frac{5}{6}} 2^{-\frac{3l_a}{2}-\frac{k}{4}-j} \langle 2^j \rangle^{-\frac{3}{4}} \Vert u \Vert_X \\
&\lesssim t^{-\frac{7}{6}+\frac{\delta}{3}} 2^{-j+\delta k} \langle 2^j \rangle^{-\frac{5}{6}+\frac{2\delta}{3}} \Vert u \Vert_X \\
\eqref{estdispsupCloinLmoyentau-aint-4} &\lesssim \sum_{l_a = l_0}^{-200} 1_{l_a < \frac{3k}{2}-\frac{j_{+}}{2}} t^{-2} 2^{-2l_a-5j} \Vert \Psi_{(l_a, F(l_a, k, j), l_a)} \psi_{j, k}^{\widehat{\mathcal{C}}} \Vert_{L^2} \Vert \nabla X_a h_a(t) \Vert_{L^2} \\
&\lesssim \sum_{l_a = l_0}^{-200} 1_{l_a < \frac{3k}{2}-\frac{j_{+}}{2}} t^{-\frac{7}{6}} \mathfrak{t}^{-\frac{5}{6}} 2^{-\frac{l_a}{2}-\frac{k}{4}-j} \Vert u \Vert_X \\
&\lesssim t^{-\frac{7}{6}} \mathfrak{t}^{-\frac{1}{2}+\frac{\delta}{3}} 2^{-j+\delta k} \langle 2^j \rangle^{-\frac{1}{12}-\frac{\delta}{3}} \Vert u \Vert_X \\
&\lesssim t^{-\frac{7}{6}} 2^{-j+\delta k} \langle 2^j \rangle^{-1} \Vert u \Vert_X \\
\eqref{estdispsupCloinLmoyentau-aint-5} &\lesssim \sum_{l_a = l_0}^{-200} 1_{l_a < \frac{3k}{2}-\frac{j_{+}}{2}} t^{-1} 2^{-l_a-\frac{5j}{2}} \Vert \Psi_{(l_a, F(l_a, k, j), l_a)} \psi_{j, k}^{\widehat{\mathcal{C}}} \Vert_{L^2} \langle 2^j \rangle^{-1} \Vert m_{\widehat{\mathcal{C}}}(D) \langle \nabla \rangle |\nabla|^{\frac{1}{2}} g_a(t) \Vert_{L^2} \\
&\lesssim \sum_{l_a = l_0}^{-200} 1_{l_a < \frac{3k}{2}-\frac{j_{+}}{2}} t^{-\frac{7}{6}+100\delta} 2^{\frac{l_a}{2}-\frac{k}{4}-j} \langle 2^j \rangle^{-\frac{3}{4}} \Vert u \Vert_X \\
&\lesssim \sum_{l_a = l_0}^{-200} t^{-\frac{7}{6}+100\delta} 2^{\frac{l_a}{3}-\frac{2\delta l_a}{3}-j+\delta k} \langle 2^j \rangle^{-\frac{5}{6}-\frac{\delta}{3}} \Vert u \Vert_X \\
&\lesssim t^{-\frac{7}{6}+100\delta} 2^{-j+\delta k} \langle 2^j \rangle^{-\frac{5}{6}} \Vert u \Vert_X 
\end{align*}
\eqref{estdispsupCloinLmoyentau-aint-3} is similar to \eqref{estdispsupCloinLmoyentau-aint-2}. 

\eqref{estdispsupCloinLmoyentau-cext}, \eqref{estdispsupCloinLmoyentau-cb}, \eqref{estdispsupCloinLmoyentau-cint} can be estimated the same way, replacing $\widehat{X}_a$ by $\widehat{X}_c$. 

It remains to estimate \eqref{estdispsupCloinLmoyentau-b} and \eqref{estdispsupCloinLmoyentau-int}, specifying the form $G$. 

First, if $l_0 \geq \frac{3k}{2}-\frac{j_{+}}{2}$, then we already fixed $G(l_0, j, k) = 0$, and $\eqref{estdispsupCloinLmoyentau-b} = 0$, while
\begin{align*}
\eqref{estdispsupCloinLmoyentau-int} &\lesssim \Vert \Psi_{(l_0, 0, l_0)} \psi_{j, k}^{\widehat{\mathcal{C}}} \Vert_{L^{\frac{1}{1-\kappa}}_c L^2_b L^{\frac{1}{1-\kappa}}_a} \langle 2^j \rangle^{-1} \Vert \langle \overline{\xi} \rangle \widehat{f}(t) \Vert_{L^{\frac{1}{\kappa}}_c L^2_b L^{\frac{1}{\kappa}}_a} \\
&\lesssim 2^{2l_0-2\kappa l_0+\frac{k}{2}+\frac{5j}{2}-2\kappa j} \langle 2^j \rangle^{-1} \Vert u \Vert_X \\
&\lesssim t^{-\frac{7}{6}+\frac{2\kappa}{3}} \mathfrak{t}^{\frac{7}{6}-\frac{2\kappa}{3}} 2^{\frac{7l_0}{3}-\frac{2\delta l_a}{3}-2\kappa l_0-j+\delta k} \langle 2^j \rangle^{-\frac{5}{6}-\frac{\delta}{3}} \Vert u \Vert_X \\
&\lesssim t^{-\frac{7}{6}+\frac{2\kappa}{3}} \mathfrak{t}^{\frac{7\varrho}{3}+\frac{\kappa}{3}+\frac{\delta}{3}} 2^{-j+\delta k} \langle 2^j \rangle^{-\frac{5}{6}-\frac{\delta}{3}} \Vert u \Vert_X \\
&\lesssim t^{-\frac{7}{6}+\frac{7\varrho}{3}+\kappa+\frac{\delta}{3}} 2^{-j+\delta k} \langle 2^j \rangle^{-\frac{5}{6}+\frac{2\delta}{3}+\kappa+7\varrho} \Vert u \Vert_X 
\end{align*}
which is enough if $\delta, \kappa, \varrho$ are small enough with respect to $1$. 

We now assume $l_0 < \frac{3k}{2}-\frac{j_{+}}{2}$. Let us define $G$ by 
\begin{align*}
G(l_0, j, k) := \left\{ \begin{array}{ll}
2k - \frac{7}{24} \log_2(t) - \frac{5j}{6} & \mbox{ if } 2^j \lesssim t^{-\frac{3}{44}} \\
2k - \frac{9}{35} \log_2(t) - \frac{4j}{7} - \frac{j_{+}}{7} &\mbox{ if } 2^j \gtrsim t^{-\frac{3}{44}} 
\end{array} \right. 
\end{align*}
This choice is motivated by the following heuristics: we use the case $\tau \ll 2^k$ but replace $2^k$ by $2^{l_b-2k}$. (Above, we expressed $G$ as a function of $\log_2(t), j, k$, but it is clear that $\log_2(t)$ can be expressed as a function of $l_0, j, \varrho$, and $\varrho$ is fixed here.) 

On \eqref{estdispsupCloinLmoyentau-b}, we apply one integration by parts along $\widehat{X}_{b-corr}$: 
\begin{subequations}
\begin{align}
\eqref{estdispsupCloinLmoyentau-b} &= \sum_{l_b = G(l_0, j, k)}^{3k} 1_{l_0 < \frac{3k}{2}-\frac{j_{+}}{2}} t^{-1} 2^{-2l_b+3k-3j} \int e^{i t \Phi} \Psi_{(l_0, l_b, l_0)}^b \psi_{j, k}^{\widehat{\mathcal{C}}}(\overline{\xi}) \widehat{f}(t, \overline{\xi}) ~ d\overline{\xi} \label{estdispsupCloinLmoyentau-b-1} \\
&\quad + \sum_{l_b = G(l_0, j, k)}^{3k} 1_{l_0 < \frac{3k}{2}-\frac{j_{+}}{2}} t^{-1} 2^{-l_b-2j} \int e^{i t \Phi} \Psi_{(l_0, l_b, l_0)}^b \psi_{j, k}^{\widehat{\mathcal{C}}}(\overline{\xi}) \widehat{h}_{b-corr}(t, \overline{\xi}) ~ d\overline{\xi} \label{estdispsupCloinLmoyentau-b-2} \\
&\quad + \sum_{l_b = G(l_0, j, k)}^{3k} 1_{l_0 < \frac{3k}{2}-\frac{j_{+}}{2}} t^{-1} 2^{-l_b-2j} \int e^{i t \Phi} \Psi_{(l_0, l_b, l_0)}^b \psi_{j, k}^{\widehat{\mathcal{C}}}(\overline{\xi}) \widehat{g}_{b-corr}(t, \overline{\xi}) ~ d\overline{\xi} \label{estdispsupCloinLmoyentau-b-3} 
\end{align}
\end{subequations}
In particular, here, $k \geq \frac{l_b}{3}$. Therefore, 
\begin{align*}
2^{3k} \gtrsim 2^{G(l_0, j, k)} = 2^{2k} 2^{G(l_0, j, k)-2k} 
\end{align*}
Hence, $2^{-k} \lesssim 2^{-(G(l_0, j, k)-2k)}$. 

Assume first $2^j \lesssim t^{-\frac{3}{44}}$. We then estimate: 
\begin{align*}
\eqref{estdispsupCloinLmoyentau-b-1} &\lesssim \sum_{l_b = G(l_0, j, k)}^{3k} 1_{l_0 < \frac{3k}{2}-\frac{j_{+}}{2}} t^{-1} 2^{-2l_b+3k-3j} \Vert \Psi_{(l_0, l_b, l_0)}^b \psi_{j, k}^{\widehat{\mathcal{C}}} \Vert_{L^{\frac{1}{1-\kappa}}_c L^2_b L^{\frac{1}{1-\kappa}}_a} \Vert \widehat{f}(t) \Vert_{L^{\frac{1}{\kappa}}_c L^2_b L^{\frac{1}{\kappa}}_a} \\
&\lesssim \sum_{l_b = G(l_0, j, k)}^{3k} 1_{l_0 < \frac{3k}{2}-\frac{j_{+}}{2}} t^{-\frac{7}{6}+\frac{2\kappa}{3}} \mathfrak{t}^{\frac{1}{6}-\frac{2\kappa}{3}} 2^{-2l_b+3k+2l_0-2\kappa l_0+\frac{l_b-2k}{2}-j} \Vert u \Vert_X \\
\eqref{estdispsupCloinLmoyentau-b-2} &\lesssim \sum_{l_b = G(l_0, j, k)}^{3k} 1_{l_0 < \frac{3k}{2}-\frac{j_{+}}{2}} t^{-1} 2^{-l_b-3j} \Vert \Psi_{(l_0, l_b, l_0)}^b \psi_{j, k}^{\widehat{\mathcal{C}}} \Vert_{L^{\frac{1}{1-\kappa}}_c L^2_b L^{\frac{1}{1-\kappa}}_a} \Vert \overline{\xi} \widehat{h}_{b-corr}(t) \Vert_{L^{\frac{1}{\kappa}}_c L^2_b L^{\frac{1}{\kappa}}_a} \\
&\lesssim \sum_{l_b = G(l_0, j, k)}^{3k} 1_{l_0 < \frac{3k}{2}-\frac{j_{+}}{2}} t^{-1} 2^{-2l_b+3k-3j} \Vert \Psi_{(l_0, l_b, l_0)}^b \psi_{j, k}^{\widehat{\mathcal{C}}} \Vert_{L^{\frac{1}{1-\kappa}}_c L^2_b L^{\frac{1}{1-\kappa}}_a} \Vert u \Vert_X \\
\eqref{estdispsupCloinLmoyentau-b-3} &\lesssim \sum_{l_b = G(l_0, j, k)}^{3k} 1_{l_0 < \frac{3k}{2}-\frac{j_{+}}{2}} t^{-1} 2^{-l_b-3j} \Vert \Psi_{(l_0, l_b, l_0)}^b \psi_{j, k}^{\widehat{\mathcal{C}}} \Vert_{L^{\frac{1}{1-\kappa}}_c L^2_{b, a}} \Vert \psi(2^{-j} \overline{\xi})  \overline{\xi} \widehat{g}_{b-corr}(t) \Vert_{L^{\frac{1}{\kappa}}_c L^2_{b, a}} \\
&\lesssim \sum_{l_b = G(l_0, j, k)}^{3k} 1_{l_0 < \frac{3k}{2}-\frac{j_{+}}{2}} t^{-1} 2^{-l_b+\frac{l_b-2k}{2}+\frac{3l_0}{2}-\kappa l_0-j-\kappa j} \\
&\pushright{\Vert \psi(2^{-j} \overline{\xi}) \overline{\xi} \widehat{g}_{b-corr}(t) \Vert_{H^1_c L^2_{b, a}}^{\frac{1}{2}} \Vert \psi(2^{-j} \overline{\xi}) \overline{\xi} \widehat{g}_{b-corr}(t) \Vert_{L^2}^{\frac{1}{2}}}
\end{align*}
if $\kappa, \varrho, \delta$ are small enough. The end of the estimates is the same as in the case $\tau \ll 2^k$, using that $2^{-l_b+\frac{l_b-2k}{2}} \lesssim 2^{-\frac{5(l_b-2k)}{2}}$ or $2^{-2l_b+3k+\frac{l_b-2k}{2}} \lesssim 2^{-\frac{5(l_b-2k)}{2}}$, which can be summed in $l_b \geq G(l_0, j, k)$, and we recover estimates similar to the previous case. 

Assume now $2^j \gtrsim t^{-\frac{3}{44}}$. We estimate: 
\begin{align*}
\eqref{estdispsupCloinLmoyentau-b-1} &\lesssim \sum_{l_b = G(l_0, j, k)}^{3k} 1_{l_0 < \frac{3k}{2}-\frac{j_{+}}{2}} t^{-1} 2^{-2l_b+3k-3j} \Vert \Psi_{(l_0, l_b, l_0)}^b \psi_{j, k}^{\widehat{\mathcal{C}}} \Vert_{L^{\frac{1}{1-\kappa}}_c L^2_b L^{\frac{1}{1-\kappa}}} \langle 2^j \rangle^{-1} \Vert \langle \overline{\xi} \rangle \widehat{f}(t) \Vert_{L^{\frac{1}{\kappa}}_c L^2_b L^{\frac{1}{\kappa}}_a} \\
&\lesssim \sum_{l_b = G(l_0, j, k)}^{3k} 1_{l_0 < \frac{3k}{2}-\frac{j_{+}}{2}} t^{-1} 2^{-2(l_b-2k)-k+2l_0-2\kappa l_0 + \frac{l_b-2k}{2}-\frac{j}{2}-2\kappa j} \langle 2^j \rangle^{-1} \Vert u \Vert_X \\
\eqref{estdispsupCloinLmoyentau-b-2} &\lesssim \sum_{l_b = G(l_0, j, k)}^{3k} 1_{l_0 < \frac{3k}{2}-\frac{j_{+}}{2}} t^{-1} 2^{-l_b-3j} \Vert \Psi_{(l_0, l_b, l_0)}^b \psi_{j, k}^{\widehat{\mathcal{C}}} \Vert_{L^{\frac{1}{1-\kappa}}_c L^2_b L^{\frac{1}{1-\kappa}}_a} \Vert \overline{\xi} \widehat{h}_{b-corr}(t) \Vert_{L^{\frac{1}{\kappa}}_c L^2_b L^{\frac{1}{\kappa}}_a} \\
&\lesssim \sum_{l_b = G(l_0, j, k)}^{3k} 1_{l_0 < \frac{3k}{2}-\frac{j_{+}}{2}} t^{-1} 2^{-2l_b+3k-3j} \Vert \Psi_{(l_0, l_b, l_0)}^b \psi_{j, k}^{\widehat{\mathcal{C}}} \Vert_{L^{\frac{1}{1-\kappa}}_c L^2_b L^{\frac{1}{1-\kappa}}_a} \Vert u \Vert_X \\
\eqref{estdispsupCloinLmoyentau-b-3} &\lesssim \sum_{l_b = G(l_0, j, k)}^{3k} 1_{l_0 < \frac{3k}{2}-\frac{j_{+}}{2}} t^{-1} 2^{-l_b-3j} \Vert \Psi_{(l_0, l_b, l_0)}^b \psi_{j, k}^{\widehat{\mathcal{C}}} \Vert_{L^{\frac{1}{1-\kappa}}_c L^2_{b, a}} \Vert \psi(2^{-j} \overline{\xi}) \overline{\xi} \widehat{g}_{b-corr}(t) \Vert_{L^{\frac{1}{\kappa}}_c L^2_{b, a}} 
\end{align*}
for small enough $\kappa, \delta, \varrho$. Again we can always recover similar bounds as in the previous case. 

For \eqref{estdispsupCloinLmoyentau-int}, we first estimate if $2^j \lesssim t^{-\frac{3}{44}}$: 
\begin{align*}
\eqref{estdispsupCloinLmoyentau-int} &\lesssim \Vert \Psi_{(l_0, G(l_0, j, k), l_0)} \psi_{j, k}^{\widehat{\mathcal{C}}} \Vert_{L^{\frac{1}{1-\kappa}}_c L^2_b L^{\frac{1}{1-\kappa}}_a} \Vert \widehat{f}(t) \Vert_{L^{\frac{1}{\kappa}}_c L^2_b L^{\frac{1}{\kappa}}_a} \\
&\lesssim 2^{2l_0-2\kappa l_0+\frac{G(l_0, j, k)-2k}{2}+\frac{5j}{2}-2 \kappa j} \Vert u \Vert_X \\
&\lesssim t^{-\frac{7}{6}} t^{\frac{1}{48}+\frac{7\delta}{24}+2\varrho+\kappa} 2^{\frac{j}{12} + \kappa j + 6 \varrho j + \frac{5 \delta j}{6}} 2^{-j+\delta k} \Vert u \Vert_X \\
&\lesssim t^{-\frac{7}{6}} t^{\frac{1}{44}} 2^{\frac{j}{12}} 2^{-j+\delta k} \Vert u \Vert_X 
\end{align*}
that satisfies the same estimates as \eqref{estdispsupCloinLmoyentau-b-3}. 

On the other hand, if $2^j \gtrsim t^{-\frac{3}{44}}$: 
\begin{align*}
\eqref{estdispsupCloinLmoyentau-int} &\lesssim 2^{-j} \Vert \Psi_{(l_0, G(l_0, j, k), l_0)} \psi_{j, k}^{\widehat{\mathcal{C}}} \Vert_{L^{\frac{1}{1-\kappa}}} \Vert \overline{\xi} \widehat{f}(t) \Vert_{L^{\frac{1}{\kappa}}} \\
&\lesssim 2^{2l_0-2\kappa l_0 + G(l_0, j, k) - 2k - \kappa (G(l_0, j, k)- 2k)+2j-3\kappa j} \Vert \overline{\xi} \widehat{f}(t) \Vert_{H^1_c L^2_b H^1_a}^{\frac{1}{2}} \Vert \overline{\xi} \widehat{f}(t) \Vert_{L^2_c H^1_b L^2_a}^{\frac{1}{2}} 
\end{align*}
and again the estimate is similar. 

\paragraph{Away from $\mathcal{L}$, away from $\mathcal{C}$} Let us now assume that $|y| \gtrsim |x|$ and $\tau = \sqrt{\frac{\left| |x| - \sqrt{3} |y| \right|}{|x|+|y|}} \gg 2^k$. 

In this situation, once localised to have $\widehat{X}_a \cdot \nabla_{\overline{\xi}} \Phi$ and $\widehat{X}_c \cdot \nabla_{\overline{\xi}} \Phi$ small enough, we have that $\widehat{X}_{b-corr} \cdot \nabla_{\overline{\xi}} \Phi \simeq \tau^2$. In particular, we can proceed with estimates similar to the case $\tau \ll 2^k$. We skip the details. 

\paragraph{Close to $\mathcal{L}$} Let us finally assume that $|y| \ll |x|$, so that $\frac{|y|}{t} \ll 2^{2j}$. 

In this case, once localised on $\widehat{X}_a \cdot \nabla_{\overline{\xi}} \Phi, \widehat{X}_c \cdot \nabla_{\overline{\xi}} \Phi \ll 2^{2j}$, we have
\begin{align*}
\left| 3 \xi_0 |\overline{\xi}| - \frac{x \xi_0}{t |\overline{\xi}|} \right| \ll 2^{2j} 
\end{align*}
Therefore, 
\begin{align*}
\widehat{X}_{b-corr}(\overline{\xi}) \cdot \nabla_{\overline{\xi}} \Phi &= \left( \sqrt{3} |\xi_0| - |\xi| \right)^2 - \frac{x}{t} + \frac{\sqrt{3} \xi_0 \xi \cdot y}{t |\xi_0| |\xi|} \\
&\simeq - \frac{x}{t} \simeq 2^{2j} 
\end{align*}
This means that the only resonance in $\widehat{X}_{b-corr}$ comes from $m_b$, which behaves like $2^k$. 

Considering again the localisation symbols $\Psi_{l_a}^{*}$, $* \in \{ a, c, int \}$ introduced in the case $\tau \ll 2^k$ above (away from $\mathcal{L}$), we can separate: 
\begin{subequations}
\begin{align}
\eqref{estdispsupC-termeinitbis} &= \sum_{l_a = l_0}^{-200} 1_{l_a \geq (1-2\delta) k} \int e^{i t \Phi} \Psi_{l_a}^a \psi_{j, k}^{\widehat{\mathcal{C}}}(\overline{\xi}) \widehat{f}(t, \overline{\xi}) d\overline{\xi} \label{estdispsupCprocheL-a1} \\
&\quad + \sum_{l_a = l_0}^{-200} 1_{l_a < (1-2\delta) k} \int e^{i t \Phi} \Psi_{l_a}^a \psi_{j, k}^{\widehat{\mathcal{C}}}(\overline{\xi}) \widehat{f}(t, \overline{\xi}) d\overline{\xi} \label{estdispsupCprocheL-a2} \\
&\quad + \sum_{l_a = l_0}^{-200} 1_{l_a \geq (1-2\delta) k} \int e^{i t \Phi} \Psi_{l_a}^ca \psi_{j, k}^{\widehat{\mathcal{C}}}(\overline{\xi}) \widehat{f}(t, \overline{\xi}) d\overline{\xi} \label{estdispsupCprocheL-c1} \\
&\quad + \sum_{l_a = l_0}^{-200} 1_{l_a < (1-2\delta) k} \int e^{i t \Phi} \Psi_{l_a}^a \psi_{j, k}^{\widehat{\mathcal{C}}}(\overline{\xi}) \widehat{f}(t, \overline{\xi}) d\overline{\xi} \label{estdispsupCprocheL-c2} \\
&\quad + \int e^{i t \Phi} \Psi_{l_0}^{int} \psi_{j, k}^{\widehat{\mathcal{C}}}(\overline{\xi}) \widehat{f}(t, \overline{\xi}) d\overline{\xi} \label{estdispsupCprocheL-int} 
\end{align}
\end{subequations} 
where $l_0$ is such that $2^{l_0} \simeq \mathfrak{t}^{-\frac{1}{2}+\varrho}$ for some $\varrho > 0$ small enough. 

Since $|y| \ll |x|$, the localisation $\widehat{X}_c \cdot \nabla \Phi \lesssim 2^{l_a+2j}$ does not give as much information as before. We get the bounds 
\begin{align}
\Vert \Psi_{l_a}^{*} \psi_{j, k}^{\widehat{\mathcal{C}}} \Vert_{L^{p_c}_c L^{p_b}_b L^{p_a}_a} &\lesssim 2^{\frac{l_a}{p_a}+\frac{k}{p_b}} |\overline{\xi_a}|^{\frac{1}{p_a}+\frac{1}{p_b}+\frac{1}{p_c}}  \label{estimeevolumiquePsilacasCvoisL}
\end{align}

On the other hand, we can compute 
\begin{align*}
\widehat{X}_a(\overline{\xi}) \cdot \nabla_{\overline{\xi}} \left[ \widehat{X}_a(\overline{\xi}) \cdot \nabla_{\overline{\xi}} \Phi \right] &\simeq 2^j \\
\widehat{X}_a(\overline{\xi}) \cdot \nabla_{\overline{\xi}} \left[ m_b(\overline{\xi}) \widehat{X}_b(\overline{\xi}) \cdot \nabla_{\overline{\xi}} \Phi \right] &\lesssim 2^{k+j} \\
\widehat{X}_a(\overline{\xi}) \cdot \nabla_{\overline{\xi}} \left[ \widehat{X}_c(\overline{\xi}) \cdot \nabla_{\overline{\xi}} \Phi \right] &= 0 \\
m_b(\overline{\xi}) \widehat{X}_b(\overline{\xi}) \cdot \nabla_{\overline{\xi}} \left[ \widehat{X}_a(\overline{\xi}) \cdot \nabla_{\overline{\xi}} \Phi \right] &\lesssim 2^{k+j} \\
m_b(\overline{\xi}) \widehat{X}_b(\overline{\xi}) \cdot \nabla_{\overline{\xi}} \left[ m_b(\overline{\xi}) \widehat{X}_b(\overline{\xi}) \cdot \nabla_{\overline{\xi}} \Phi \right] &\lesssim 2^{k+j} \\
m_b(\overline{\xi}) \widehat{X}_b(\overline{\xi}) \cdot \nabla_{\overline{\xi}} \left[ \widehat{X}_c(\overline{\xi}) \cdot \nabla_{\overline{\xi}} \Phi \right] &\lesssim 2^{l_a+j} \\
\widehat{X}_c(\overline{\xi}) \cdot \nabla_{\overline{\xi}} \left[ m_b(\overline{\xi}) \widehat{X}_b(\overline{\xi}) \cdot \nabla_{\overline{\xi}} \Phi \right] &\lesssim 2^{j+k} \\
\widehat{X}_c(\overline{\xi}) \cdot \nabla_{\overline{\xi}} \left[ \widehat{X}_c(\overline{\xi}) \cdot \nabla_{\overline{\xi}} \Phi \right] &\lesssim 2^j 
\end{align*}
In particular, 
\begin{align*}
\widehat{X}_a(\overline{\xi}) \cdot \nabla_{\overline{\xi}} \left[ \Psi_{l_a}^{*} \psi_{j, k}^{\widehat{\mathcal{C}}} \right] &\lesssim 2^{-l_a-j} \\
\widehat{X}_c(\overline{\xi}) \cdot \nabla_{\overline{\xi}} \left[ \Psi_{l_a}^{*} \psi_{j, k}^{\widehat{\mathcal{C}}} \right] &\lesssim 2^{-l_a-j} 
\end{align*}
On the other hand, if $k \leq l_a+10$, we also have
\begin{align*}
m_b(\overline{\xi}) \widehat{X}_b(\overline{\xi}) \cdot \nabla_{\overline{\xi}} \left[ \Psi_{l_a}^{*} \psi_{j, k}^{\widehat{\mathcal{C}}} \right] &\lesssim 2^{-j} 
\end{align*}
While if $k \geq l_a+10$, as before, we correct $\widehat{X}_b$ into a $\widehat{X}_{b-corr}$ by adding a part in the direction $\widehat{X}_a$, so that 
\begin{align*}
\widehat{X}_{b-corr}(\overline{\xi}) \cdot \nabla_{\overline{\xi}} \left[ \widehat{X}_a(\overline{\xi}) \cdot \nabla_{\overline{\xi}} \Phi \right] &= 0 
\end{align*}
Hence, 
\begin{align*}
m_b(\overline{\xi}) \widehat{X}_{b-corr}(\overline{\xi}) \cdot \nabla_{\overline{\xi}} \left[ \Psi_{l_a}^{*} \psi_{j, k}^{\widehat{\mathcal{C}}} \right] &\lesssim 2^{-j} 
\end{align*}
as well. Since $l_a \leq -200$, we also have
\begin{align*}
\widehat{X}_{b-corr}(\overline{\xi}) \cdot \nabla_{\overline{\xi}} \Phi &\simeq \widehat{X}_b(\overline{\xi}) \cdot \nabla_{\overline{\xi}} \Phi \simeq 2^{2j} 
\end{align*}

On \eqref{estdispsupCprocheL-a1}, we apply at most $n$ (large with respect to $\varrho^{-1}$) integrations by parts along $\widehat{X}_a$, as long as no derivative hits $\widehat{f}(t)$, and if a derivative hits $\widehat{f}(t)$, we separate into $h_a+g_a$ and apply a last integration by parts on $h_a$: 
\begin{subequations}
\begin{align}
\eqref{estdispsupCprocheL-a1} &= \sum_{l_a = l_0}^{-200} 1_{l_a \geq (1-2\delta) k} t^{-n} 2^{-2nl_a-3nj} \int e^{i t \Phi} \Psi_{l_a}^a \psi_{j, k}^{\widehat{\mathcal{C}}}(\overline{\xi}) \widehat{f}(t, \overline{\xi}) d\overline{\xi} \label{estdispsupCprocheL-a1-1} \\
&\quad + \sum_{i = 1}^n \sum_{l_a = l_0}^{-200} 1_{l_a \geq (1-2\delta) k} t^{-i-1} 2^{-2il_a-l_a-3ij-2j} \int e^{i t \Phi} \Psi_{l_a}^a \psi_{j, k}^{\widehat{\mathcal{C}}}(\overline{\xi}) \widehat{h}_a(t, \overline{\xi}) d\overline{\xi} \label{estdispsupCprocheL-a1-2} \\
&\quad + \sum_{i = 1}^n \sum_{l_a = l_0}^{-200} 1_{l_a \geq (1-2\delta) k} t^{-i-1} 2^{-2il_a-3ij-j} \int e^{i t \Phi} \Psi_{l_a}^a \psi_{j, k}^{\widehat{\mathcal{C}}}(\overline{\xi}) \widehat{X}_a(\overline{\xi}) \cdot \nabla_{\overline{\xi}} \widehat{h}_a(t, \overline{\xi}) d\overline{\xi} \label{estdispsupCprocheL-a1-3} \\
&\quad + \sum_{i = 1}^n \sum_{l_a = l_0}^{-200} 1_{l_a \geq (1-2\delta) k} t^{-i} 2^{-2il_a+l_a-3ij+j} \int e^{i t \Phi} \Psi_{l_a}^a \psi_{j, k}^{\widehat{\mathcal{C}}}(\overline{\xi}) \widehat{g}_a(t, \overline{\xi}) d\overline{\xi} \label{estdispsupCprocheL-a1-4} 
\end{align}
\end{subequations}
We now estimate: 
\begin{align*}
\eqref{estdispsupCprocheL-a1-1} &\lesssim \sum_{l_a = l_0}^{-200} 1_{l_a \geq (1-2\delta) k} \mathfrak{t}^{-n} 2^{-2nl_a+j} \Vert \Psi_{l_a}^a \psi_{j, k}^{\widehat{\mathcal{C}}} \Vert_{L^2} \langle 2^j \rangle^{-1} \Vert \langle \nabla \rangle |\nabla|^{-1} f(t) \Vert_{L^2} \\
&\lesssim \sum_{l_a = l_0}^{-200} 1_{l_a \geq k} t^{-\frac{7}{6}} \mathfrak{t}^{-n+\frac{7}{6}} 2^{-2nl_a+\frac{l_a}{2}-j+\frac{k}{2}} \langle 2^j \rangle^{-1} \Vert u \Vert_X \\
&\lesssim t^{-\frac{7}{6}} \mathfrak{t}^{-2n\varrho+\frac{7}{6}} 2^{-j+\frac{k}{2}} \langle 2^j \rangle^{-1} \Vert u \Vert_X \\
&\lesssim t^{-\frac{7}{6}} 2^{-j+\frac{k}{2}} \langle 2^j \rangle^{-1} \Vert u \Vert_X \\
\eqref{estdispsupCprocheL-a1-2} &\lesssim \sum_{i = 1}^n \sum_{l_a = l_0}^{-200} 1_{l_a \geq (1-2\delta) k} t^{-\frac{7}{6}} \mathfrak{t}^{-i+\frac{1}{6}} 2^{-2il_a-l_a-\frac{7j}{2}} \Vert \Psi_{l_a}^a \psi_{j, k}^{\widehat{\mathcal{C}}} \Vert_{L^2_{c, b} L^{\frac{1}{1-\kappa}}_a} \Vert \psi_j(\overline{\xi}) \overline{\xi} \widehat{h}_a(t) \Vert_{L^2_{c, b} L^{\frac{1}{\kappa}}_a} \\
&\lesssim \sum_{i = 1}^n \sum_{l_a = l_0}^{-200} 1_{l_a \geq (1-2\delta) k} t^{-\frac{7}{6}} \mathfrak{t}^{-i+\frac{1}{6}} 2^{-2il_a-\kappa l_a-\frac{3j}{2}-\kappa j+\frac{k}{2}} \Vert \psi_j(\overline{\xi}) \overline{\xi} \widehat{h}_a(t) \Vert_{L^2_{c, b} H^1_a}^{\frac{1}{2}} \Vert \psi_j(\overline{\xi}) \overline{\xi} \widehat{h}_a(t) \Vert_{L^2}^{\frac{1}{2}} \\
&\lesssim \sum_{i = 1}^n \sum_{l_a = l_0}^{-200} 1_{l_a \geq (1-2\delta) k} t^{-\frac{7}{6}+\frac{\kappa}{3}} \mathfrak{t}^{-i+\frac{1}{6}-\frac{\kappa}{3}} 2^{-2il_a-\kappa l_a-j+\frac{k}{2}} \langle 2^j \rangle^{-\frac{1}{2}} \Vert u \Vert_X \\
&\lesssim t^{-\frac{7}{6}+\frac{\kappa}{3}} \mathfrak{t}^{-\frac{1}{12}+\frac{\kappa}{6}} 2^{-j+2\delta k} \langle 2^j \rangle^{-\frac{1}{2}} \Vert u \Vert_X \\
&\lesssim t^{-\frac{7}{6}+\frac{\kappa}{3}} 2^{-j+\delta k} \langle 2^j \rangle^{-\frac{3}{4}+\frac{\kappa}{2}} \Vert u \Vert_X \\
\eqref{estdispsupCprocheL-a1-3} &\lesssim \sum_{i = 1}^n \sum_{l_a = l_0}^{-200} 1_{l_a \geq (1-2\delta) k} t^{-\frac{7}{6}} \mathfrak{t}^{-i+\frac{1}{6}} 2^{-2il_a-\frac{5j}{2}} \Vert \Psi_{l_a}^a \psi_{j, k}^{\widehat{\mathcal{C}}} \Vert_{L^2} \Vert \nabla X_a h_a(t) \Vert_{L^2} \\
&\lesssim \sum_{i = 1}^n \sum_{l_a = l_0}^{-200} 1_{l_a \geq (1-2\delta) k} t^{-\frac{7}{6}} \mathfrak{t}^{-i+\frac{1}{6}} 2^{-2il_a+\frac{l_a}{2}-j+\frac{k}{2}} \Vert u \Vert_X \\
&\lesssim t^{-\frac{7}{6}} \mathfrak{t}^{-\frac{1}{3}} 2^{-j+2\delta k} \Vert u \Vert_X \\
&\lesssim t^{-\frac{7}{6}} 2^{-j+\delta k} \langle 2^j \rangle^{-1} \Vert u \Vert_X \\
\eqref{estdispsupCprocheL-a1-4} &\lesssim \sum_{i = 1}^n \sum_{l_a = l_0}^{-200} 1_{l_a \geq (1-2\delta)k} t^{-1} \mathfrak{t}^{-i+1} 2^{-2il_a+l_a-\frac{5j}{2}} \Vert \Psi_{l_a}^a \psi_{j, k}^{\widehat{\mathcal{C}}} \Vert_{L^2} \langle 2^j \rangle^{-1} \Vert m_{\widehat{\mathcal{C}}}(D) \langle \nabla \rangle |\nabla|^{\frac{1}{2}} g_a(t) \Vert_{L^2} \\
&\lesssim \sum_{i = 1}^n \sum_{l_a = l_0}^{-200} 1_{l_a \geq (1-2\delta)k} t^{-\frac{7}{6}+100\delta} \mathfrak{t}^{-i+1} 2^{-2il_a+\frac{3l_a}{2}-j+\frac{k}{2}} \langle 2^j \rangle^{-1} \Vert u \Vert_X \\
&\lesssim t^{-\frac{7}{6}+100\delta} 2^{-j+\delta k} \langle 2^j \rangle^{-1} \Vert u \Vert_X
\end{align*}
for $n \varrho$ large enough with respect to $1$. 

On the other hand, for \eqref{estdispsupCprocheL-a2}, we start by an integration by parts along $\widehat{X}_{b-corr}$, then only we apply possibly integrations by parts along $\widehat{X}_a$. Assume first that $2^{|j|} \geq t^{1000\delta}$. Then we decompose: 
\begin{subequations}
\begin{align}
\eqref{estdispsupCprocheL-a2} &= \sum_{l_a = l_0}^{-200} 1_{l_a < (1-2\delta) k} t^{-1} 2^{-k-3j} \int e^{i t \Phi} \Psi_{l_a}^a \psi_{j, k}^{\widehat{\mathcal{C}}}(\overline{\xi}) \widehat{f}(t, \overline{\xi}) d\overline{\xi} \notag \\
&\quad + \sum_{l_a = l_0}^{-200} 1_{l_a < (1-2\delta) k} t^{-1} 2^{-k-2j} \int e^{i t \Phi} \Psi_{l_a}^a \psi_{j, k}^{\widehat{\mathcal{C}}}(\overline{\xi}) \widehat{h}_{b-corr}(t, \overline{\xi}) d\overline{\xi} \notag \\
&\quad + \sum_{l_a = l_0}^{-200} 1_{l_a < (1-2\delta) k} t^{-1} 2^{-k-2j} \int e^{i t \Phi} \Psi_{l_a}^a \psi_{j, k}^{\widehat{\mathcal{C}}}(\overline{\xi}) \widehat{g}_{b-corr}(t, \overline{\xi}) d\overline{\xi} \notag \\
&= \sum_{l_a = l_0}^{-200} 1_{l_a < (1-2\delta) k} t^{-1-n} 2^{-k-3j-2nl_a-3nj} \int e^{i t \Phi} \Psi_{l_a}^a \psi_{j, k}^{\widehat{\mathcal{C}}}(\overline{\xi}) \widehat{f}(t, \overline{\xi}) d\overline{\xi} \label{estdispsupCprocheL-a2-1} \\
&\quad + \sum_{i = 1}^n \sum_{l_a = l_0}^{-200} 1_{l_a < (1-2\delta) k} t^{-1-i} 2^{-k-3j-2il_a+l_a-3ij+j} \int e^{i t \Phi} \Psi_{l_a}^a \psi_{j, k}^{\widehat{\mathcal{C}}}(\overline{\xi}) \widehat{X}_a(\overline{\xi}) \cdot \nabla_{\overline{\xi}} \widehat{f}(t, \overline{\xi}) d\overline{\xi} \label{estdispsupCprocheL-a2-2} \\
&\quad + \sum_{l_a = l_0}^{-200} 1_{l_a < (1-2\delta) k} t^{-1-n} 2^{-k-2j-2nl_a-3nj} \int e^{i t \Phi} \Psi_{l_a}^a \psi_{j, k}^{\widehat{\mathcal{C}}}(\overline{\xi}) \widehat{h}_{b-corr}(t, \overline{\xi}) d\overline{\xi} \label{estdispsupCprocheL-a2-3} \\
&\quad \begin{aligned}
+ \sum_{i = 1}^n \sum_{l_a = l_0}^{-200} 1_{l_a < (1-2\delta) k} t^{-1-i} &2^{-k-2j-2il_a+l_a-3ij+j} \\
&\int e^{i t \Phi} \Psi_{l_a}^a \psi_{j, k}^{\widehat{\mathcal{C}}}(\overline{\xi}) \widehat{X}_a(\overline{\xi}) \cdot \nabla_{\overline{\xi}} \widehat{h}_{b-corr}(t, \overline{\xi}) d\overline{\xi} 
\end{aligned} \label{estdispsupCprocheL-a2-4} \\
&\quad + \sum_{l_a = l_0}^{-200} 1_{l_a < (1-2\delta) k} t^{-1} 2^{-k-2j} \int e^{i t \Phi} \Psi_{l_a}^a \psi_{j, k}^{\widehat{\mathcal{C}}}(\overline{\xi}) \widehat{g}_{b-corr}(t, \overline{\xi}) d\overline{\xi} \label{estdispsupCprocheL-a2-5} 
\end{align}
\end{subequations} 
We then estimate:  
\begin{align*}
\eqref{estdispsupCprocheL-a2-3} &\lesssim \sum_{l_a = l_0}^{-200} 1_{l_a < (1-2\delta) k} t^{-\frac{7}{6}} \mathfrak{t}^{-n+\frac{1}{6}} 2^{-k-\frac{5j}{2}-2nl_a} \Vert \Psi_{l_a}^a \psi_{j, k}^{\widehat{\mathcal{C}}} \Vert_{L^2} \langle 2^j \rangle^{-1} \Vert \langle \nabla \rangle h_{b-corr}(t) \Vert_{L^2} \\
&\lesssim \sum_{l_a = l_0}^{-200} 1_{l_a < (1-2\delta) k} t^{-\frac{7}{6}} \mathfrak{t}^{-n+\frac{1}{6}} 2^{-2nl_a+\frac{l_a}{2}-\frac{k}{2}-j} \langle 2^j \rangle^{-1} \Vert u \Vert_X \\
&\lesssim t^{-\frac{7}{6}} \mathfrak{t}^{-2n\varrho+\frac{1}{2}} 2^{-j+\frac{k}{2}} \langle 2^j \rangle^{-1} \Vert u \Vert_X \\
&\lesssim t^{-\frac{7}{6}} 2^{-j+\frac{k}{2}} \langle 2^j \rangle^{-1} \Vert u \Vert_X \\
\eqref{estdispsupCprocheL-a2-4} &\lesssim \sum_{i = 1}^n \sum_{l_a = l_0}^{-200} 1_{l_a < (1-2\delta) k} t^{-\frac{7}{6}} \mathfrak{t}^{-i+\frac{1}{6}} 2^{-k-2il_a+l_a-\frac{5j}{2}} \Vert \Psi_{l_a}^a \psi_{j, k}^{\widehat{\mathcal{C}}} \Vert_{L^2} \Vert \nabla X_a h_{b-corr}(t) \Vert_{L^2} \\
&\lesssim \sum_{i = 1}^n \sum_{l_a = l_0}^{-200} 1_{l_a < (1-2\delta) k} t^{-\frac{7}{6}} \mathfrak{t}^{-i+\frac{1}{6}} 2^{-\frac{k}{2}-2il_a+\frac{3l_a}{2}-j} \Vert u \Vert_X \\
&\lesssim \sum_{i = 1}^n \sum_{l_a = l_0}^{-200} t^{-\frac{7}{6}} \mathfrak{t}^{-i+\frac{1}{6}} 2^{-\frac{1+2\delta}{1-2\delta} \frac{l_a}{2}-2il_a+\frac{3l_a}{2}-j+\delta k} \Vert u \Vert_X \\
&\lesssim t^{-\frac{7}{6}} \mathfrak{t}^{\frac{1}{6}-\frac{3}{4}+\frac{1+2\delta}{1-2\delta} \frac{1}{4} + \frac{3\varrho}{2}} 2^{-j+\delta k} \Vert u \Vert_X \\
&\lesssim t^{-\frac{7}{6}} \mathfrak{t}^{-\frac{1}{3}+ \frac{3\delta}{2}+\frac{3\varrho}{2}} 2^{-j+\delta k} \Vert u \Vert_X \\
&\lesssim t^{-\frac{7}{6}} 2^{-j+\delta k} \langle 2^j \rangle^{-\frac{5}{6}} \Vert u \Vert_X \\
\eqref{estdispsupCprocheL-a2-5} &\lesssim \sum_{l_a = l_0}^{-200} 1_{l_a < (1-2\delta) k} t^{-1} 2^{-k-\frac{5j}{2}} \Vert \Psi_{l_a}^a \psi_{j, k}^{\widehat{\mathcal{C}}} \Vert_{L^2} \langle 2^j \rangle^{-1} \Vert m_{\widehat{\mathcal{C}}}(D) \langle \nabla \rangle |\nabla|^{\frac{1}{2}} g_{b-corr}(t) \Vert_{L^2} \\
&\lesssim \sum_{l_a = l_0}^{-200} 1_{l_a < (1-2\delta) k} t^{-\frac{7}{6}+100\delta} 2^{\frac{l_a}{2}-\left( \frac{1}{2}+\delta \right)k-j+\delta k} \langle 2^j \rangle^{-1} \Vert u \Vert_X \\
&\lesssim \sum_{l_a = l_0}^{-200} t^{-\frac{7}{6}+100\delta} 2^{\frac{l_a}{2}-\frac{1+2\delta}{1-2\delta} \frac{l_a}{2}-j+\delta k} \langle 2^j \rangle^{-1} \Vert u \Vert_X \\
&\lesssim t^{-\frac{7}{6}+100\delta} \mathfrak{t}^{\frac{3\delta}{2}} 2^{-j+\delta k} \langle 2^j \rangle^{-1} \Vert u \Vert_X
\end{align*}
\eqref{estdispsupCprocheL-a2-1} can be estimated as \eqref{estdispsupCprocheL-a2-3} and \eqref{estdispsupCprocheL-a2-2} as \eqref{estdispsupCprocheL-a2-4}. Above, we used that $\delta, \varrho$ are small enough, $n \varrho$ large enough. Finally, for \eqref{estdispsupCprocheL-a2-5}, we have on the one hand
\begin{align*}
\eqref{estdispsupCprocheL-a2-5} &\lesssim t^{-\frac{7}{6}+100\delta+\frac{3\delta}{2}} 2^{-j+\delta k} \langle 2^j \rangle^{-1+\frac{9\delta}{2}} \Vert u \Vert_X \\
&\lesssim t^{-\frac{13}{12}} 2^{-j+\delta k} \langle 2^j \rangle^{-\frac{5}{6}} \Vert u \Vert_X 
\end{align*}
but also, using $2^{|j|} \gtrsim t^{1000 \delta}$, and hence $2^{\frac{j}{2}} \langle 2^j \rangle^{-\frac{3}{4}} \lesssim t^{-500 \delta}$: 
\begin{align*}
\eqref{estdispsupCprocheL-a2-5} &\lesssim t^{-\frac{7}{6}+100\delta+\frac{3\delta}{2}} 2^{\frac{j}{2}} \langle 2^j \rangle^{-\frac{3}{4}} 2^{-\frac{3j}{2}+\delta k} \langle 2^j \rangle^{-\frac{1}{4}+\frac{9\delta}{2}} \Vert u \Vert_X \\
&\lesssim t^{-\frac{7}{6}} 2^{-\frac{3j}{2}+\delta k} \langle 2^j \rangle^{-\frac{1}{6}} \Vert u \Vert_X 
\end{align*}
for $\delta$ small enough. 

It remains to consider \eqref{estdispsupCprocheL-a2} in the case $2^{|j|} \lesssim t^{1000 \delta}$, so $t^{-1000 \delta} \lesssim 2^j \lesssim t^{1000 \delta}$. Note that the only place where we used $2^{|j|} \gtrsim t^{1000 \delta}$ before was to bound  \eqref{estdispsupCprocheL-a2-5}, so we restrict our attention to this term only now. If we consider the sum over $l_a < (1 + 3 \delta) k$ instead of $l_a < (1 - 2 \delta k)$, we can extend the estimates without any trouble: 
\begin{align*}
\eqref{estdispsupCprocheL-a2-5} &\lesssim \sum_{l_a = l_0}^{-200} 1_{l_a < (1+3\delta) k} t^{-1} 2^{-k-\frac{5j}{2}} \Vert \Psi_{l_a}^a \psi_{j, k}^{\widehat{\mathcal{C}}} \Vert_{L^2} \langle 2^j \rangle^{-1} \Vert m_{\widehat{\mathcal{C}}}(D) \langle \nabla \rangle |\nabla|^{\frac{1}{2}} g_{b-corr}(t) \Vert_{L^2} \\
&\lesssim \sum_{l_a = l_0}^{-200} t^{-\frac{7}{6}+100\delta} 2^{\frac{l_a}{2}-\frac{1+2\delta}{1+3\delta} \frac{l_a}{2}-j+\delta k} \langle 2^j \rangle^{-1} \Vert u \Vert_X \\
&\lesssim t^{-\frac{7}{6}+100\delta} 2^{-j+\delta k} \langle 2^j \rangle^{-1} \Vert u \Vert_X 
\end{align*}
which is enough. We can therefore consider the sum in $l_a$ restricted to $\max(l_0, (1 + 3\delta) k) \leq l_a \leq (1-2\delta) k$, so $j$ close $0$ (up to $\delta \log(t)$) and $l_a$ close to $k$ (up to $\delta k$). For $l_a$ close to $k$, we may choose $X_{b-corr} = X_b$: there may be a small loss on \eqref{estdispsupCprocheL-a2-1} to \eqref{estdispsupCprocheL-a2-4} for $l_a \leq k + 10$, but this remains of order $2^{-l_a+k} \lesssim 2^{3 \delta k} \lesssim \mathfrak{t}^{2 \delta}$ (for $\delta$ small enough) and it is clear that such a factor can be absorbed in the estimates. For \eqref{estdispsupCprocheL-a2-5}, this allows to use the fine $L^2$ estimate on $g_b$: 
\begin{align*}
\Vert \psi_j(D) g_b \Vert_{L^2} &\lesssim t^{-\frac{1}{6}+100\delta} 2^{-\frac{j}{2}+k} \langle 2^j \rangle^{-1} \Vert u \Vert_X + t^{-\frac{1}{2}} 2^{-3j} \Vert u \Vert_X 
\end{align*}
Hence, 
\begin{align*}
\eqref{estdispsupCprocheL-a2-5} &\lesssim \sum_{l_a = l_0}^{-200} 1_{|l_a - k| \leq 3 \delta k} t^{-1} 2^{-k-2j} \Vert \Psi_{l_a}^a \psi_{j, k}^{\widehat{\mathcal{C}}} \Vert_{L^2} \Vert \psi_j(D) g_{b-corr}(t) \Vert_{L^2} \\
&\lesssim \sum_{l_a = l_0}^{-200} 1_{|l_a - k| \leq 3 \delta k} t^{-\frac{7}{6}+100\delta} 2^{\frac{l_a}{2}-\frac{k}{2}-\frac{j}{2}} \left( 2^{-\frac{j}{2}+k} \langle 2^j \rangle^{-1} + t^{-\frac{1}{2}-100\delta} 2^{-3j} \right) \Vert u \Vert_X \\
&\lesssim t^{-\frac{7}{6}+100\delta} 2^{-j} \langle 2^j \rangle^{-1} \Vert u \Vert_X 
+ t^{-\frac{7}{6}-\frac{1}{2}+5000\delta} 2^{-j} \langle 2^j \rangle^{-1} \Vert u \Vert_X \\
&\lesssim t^{-\frac{7}{6}+100\delta} 2^{-j} \langle 2^j \rangle^{-1} \Vert u \Vert_X
\end{align*}
for $\delta$ small enough. 

We estimate in a symmetric way \eqref{estdispsupCprocheL-c1} and \eqref{estdispsupCprocheL-c2} replacing $\widehat{X}_a$ by $\widehat{X}_c$. 

It only remains the internal term \eqref{estdispsupCprocheL-int}. We exhaust the cases depending on $j, k, t$. 

If $2^j \lesssim t^{-\frac{1}{5}}$ and $2^k \lesssim t^{-\frac{7}{6}} 2^{-4j}$, we estimate: 
\begin{align*}
\eqref{estdispsupCprocheL-int} &\lesssim \Vert \Psi_{l_0}^{int} \psi_{j, k}^{\widehat{\mathcal{C}}} \Vert_{L^2_{c, b} L^{\frac{1}{1-\kappa}}_a} \Vert \psi_j \widehat{f}(t) \Vert_{L^2_{c, b} L^{\frac{1}{\kappa}}_a} \\
&\lesssim 2^{l_0-\kappa l_0+\frac{k}{2}+2j-\kappa j} \Vert \psi_j \widehat{f}(t) \Vert_{H^1_a L^2_{b, c}}^{\frac{1}{2}} \Vert \psi_j \widehat{f}(t) \Vert_{L^2}^{\frac{1}{2}} \\
&\lesssim \mathfrak{t}^{-\frac{1}{2}+\varrho+\frac{\kappa}{2}} t^{-\frac{7}{12}+\frac{7\delta}{6}} 2^{\frac{3j}{2}+4\delta j-\kappa j} 2^{-j+\delta k} \Vert u \Vert_X \\
&\lesssim t^{-\frac{13}{12}+\frac{7\delta}{6}+\varrho+\frac{\kappa}{2}} 2^{4\delta j+3\varrho j + \frac{\kappa j}{2}} 2^{-j+\delta k} \Vert u \Vert_X \\
&\lesssim t^{-\frac{13}{12}+2\delta} 2^{-j+\delta k} \Vert u \Vert_X 
\end{align*}
for small enough $\kappa, \varrho$ with respect to $\delta$. 

Then, if $2^j \lesssim t^{-\frac{1}{5}}$ and $2^k \gtrsim t^{-\frac{7}{6}} 2^{-4j}$, we apply an integration by parts along $\widehat{X}_{b-corr}$: 
\begin{subequations}
\begin{align}
\eqref{estdispsupCprocheL-int} &= t^{-1} 2^{-k-3j} \int e^{i t \Phi} \Psi_{l_0}^{int} \psi_{j, k}^{\widehat{\mathcal{C}}}(\overline{\xi}) \widehat{f}(t, \overline{\xi}) ~ d\overline{\xi} \label{estdispsupCprocheL-int-1} \\
&\quad + t^{-1} 2^{-k-2j} \int e^{i t \Phi} \Psi_{l_0}^{int} \psi_{j, k}^{\widehat{\mathcal{C}}}(\overline{\xi}) \widehat{h}_{b-corr}(t, \overline{\xi}) ~ d\overline{\xi} \label{estdispsupCprocheL-int-2} \\
&\quad + t^{-1} 2^{-k-2j} \int e^{i t \Phi} \Psi_{l_0}^{int} \psi_{j, k}^{\widehat{\mathcal{C}}}(\overline{\xi}) \widehat{g}_{b-corr}(t, \overline{\xi}) ~ d\overline{\xi} \label{estdispsupCprocheL-int-3} 
\end{align}
\end{subequations}
We can then estimate: 
\begin{align*}
\eqref{estdispsupCprocheL-int-2} &\lesssim t^{-1} 2^{-k-3j} \Vert \Psi_{l_0}^{int} \psi_{j, k}^{\widehat{\mathcal{C}}} \Vert_{L^2_{c, b} L^{\frac{1}{1-\kappa}}_a} \Vert \psi_j(\overline{\xi}) \overline{\xi} \widehat{h}_{b-corr}(t) \Vert_{L^2_{c, b} L^{\frac{1}{\kappa}}_a} \\
&\lesssim t^{-1} 2^{l_0-\kappa l_0-\frac{k}{2}-j-\kappa j} \Vert \psi_j(\overline{\xi}) \overline{\xi} \widehat{h}_{b-corr}(t) \Vert_{H^1_a L^2_{b, c}}^{\frac{1}{2}} \Vert \psi_j(\overline{\xi}) \overline{\xi} \widehat{h}_{b-corr}(t) \Vert_{L^2}^{\frac{1}{2}} \\
&\lesssim t^{-1+\frac{7}{12}+\frac{7\delta}{6}} \mathfrak{t}^{-\frac{1}{2}+\varrho+\frac{\kappa}{2}} 2^{\frac{3j}{2}+4\delta j-\kappa j+\delta k} \Vert u \Vert_X \\
&\lesssim t^{-\frac{11}{12}+\frac{7\delta}{6}+\varrho+\frac{\kappa}{2}} 2^{j+4\delta j+3\varrho j + \frac{\kappa j}{2}} 2^{-j+\delta k} \Vert u \Vert_X \\
&\lesssim t^{-\frac{13}{12}} 2^{-j+\delta k} \Vert u \Vert_X \\
\eqref{estdispsupCprocheL-int-3} &\lesssim t^{-1} 2^{-k-2j} \Vert \Psi_{l_0}^{int} \psi_{j, k}^{\widehat{\mathcal{C}}} \Vert_{L^2} \Vert g_{b-corr}(t) \Vert_{L^2} \\
&\lesssim t^{-\frac{13}{12}+3\delta} 2^{\frac{l_0}{2}-\frac{k}{2}-\frac{j}{2}} \Vert u \Vert_X \\
&\lesssim t^{-\frac{13}{12}+3\delta} \mathfrak{t}^{-\frac{1}{4}+\frac{\varrho}{2}} t^{\frac{7}{12}+\frac{7\delta}{6}} 2^{\frac{5j}{2}+4\delta j} 2^{-j+\delta k} \Vert u \Vert_X \\
&\lesssim t^{-\frac{13}{12}+3\delta} t^{\frac{1}{3}+\frac{7\delta}{6}+\frac{\varrho}{2}} 2^{\frac{7j}{4}+4\delta j+\frac{3\varrho j}{2}} 2^{-j+\delta k} \Vert u \Vert_X \\
&\lesssim t^{-\frac{13}{12}+3\delta} t^{\frac{1}{3}+\frac{7\delta}{6}+\frac{\varrho}{2}-\frac{7}{20}} 2^{-j+\delta k} \Vert u \Vert_X \\
&\lesssim t^{-\frac{13}{12}+2\delta} 2^{-j+\delta k} \Vert u \Vert_X 
\end{align*}
for $\varrho, \kappa$ small enough with respect to $\delta$ and $\delta$ small enough with respect to $1$. We can estimate \eqref{estdispsupCprocheL-int-1} like \eqref{estdispsupCprocheL-int-2}. 

Consider now the case $t^{-\frac{1}{5}} \lesssim 2^j$ and $2^k \lesssim t^{-\frac{2}{3}+40\delta} 2^{-\frac{3j}{2}}$. Then we can estimate
\begin{align*}
\eqref{estdispsupCprocheL-int} &\lesssim \Vert \Psi_{l_0}^{int} \psi_{j, k}^{\widehat{\mathcal{C}}} \Vert_{L^2_c L^{\frac{1}{1-\kappa}}_{b, a}} \Vert \psi_j \widehat{f}(t) \Vert_{L^2_c L^{\frac{1}{\kappa}}_{b, a}} \\
&\lesssim 2^{l_0-\kappa l_0 + k - \kappa k + \frac{5j}{2} - 2 \kappa j} \Vert \psi_j \widehat{f}(t) \Vert_{H^1_a L^2_{b, c}}^{\frac{1}{2}} \Vert \psi_j \widehat{f}(t) \Vert_{H^1_b L^2_{a, c}}^{\frac{1}{2}} \\
&\lesssim \mathfrak{t}^{-\frac{1}{2}+\varrho+\frac{\kappa}{2}} t^{-\frac{2}{3}+40\delta+\frac{2\kappa}{3}+\frac{2\delta}{3}} 2^{2j-\frac{\kappa j}{2}+\frac{3\delta j}{2}} 2^{-j+\delta k} t^{\frac{1}{12}+51\delta} 2^{-\frac{j}{2}} \langle 2^j \rangle^{-\frac{1}{2}} \Vert u \Vert_X \\
&\lesssim t^{-\frac{13}{12}+92\delta} 2^{-j+\delta k} \langle 2^j \rangle^{-\frac{1}{3}} \Vert u \Vert_X 
\end{align*}
for $\varrho, \kappa$ small enough with respect to $\delta$, small with respect to $1$. 

On the other hand, if $2^k \gtrsim t^{-\frac{2}{3}+40\delta} 2^{-\frac{3j}{2}}$, we use again the previous integrations by parts and estimate: 
\begin{align*}
\eqref{estdispsupCprocheL-int-2} &\lesssim t^{-1} 2^{-k-3j} \Vert \Psi_{l_0}^{int} \psi_{j, k}^{\widehat{\mathcal{C}}} \Vert_{L^2_{c, b} L^{\frac{1}{1-\kappa}}_a} \Vert \psi_j(\overline{\xi}) \overline{\xi} \widehat{h}_{b-corr}(t) \Vert_{L^2_{c, b} L^{\frac{1}{\kappa}}_a} \\
&\lesssim t^{-1} 2^{l_0-\kappa l_0-\frac{k}{2}-j-\kappa j} \Vert \psi_j(\overline{\xi}) \overline{\xi} \widehat{h}_{b-corr}(t) \Vert_{H^1_a L^2_{b, c}}^{\frac{1}{2}} \Vert \psi_j(\overline{\xi}) \overline{\xi} \widehat{h}_{b-corr}(t) \Vert_{L^2}^{\frac{1}{2}} \\
&\lesssim t^{-\frac{7}{6}+\frac{2\delta}{3}+\varrho+\frac{\kappa}{2}} 2^{\frac{j}{4}+\frac{\delta j}{2}+3\varrho j + \frac{\kappa j}{2}} 2^{-j+\delta k} \langle 2^j \rangle^{-\frac{1}{2}} \Vert u \Vert_X \\
&\lesssim t^{-\frac{7}{6}+\delta} 2^{-j+\delta k} \langle 2^j \rangle^{-\frac{1}{2}} \Vert u \Vert_X 
\end{align*}
for small enough $\varrho, \kappa$ with respect to $\delta$, and we may estimate \eqref{estdispsupCprocheL-int-1} as well. For \eqref{estdispsupCprocheL-int-3}, we start by the case $2^j \lesssim t^{-\frac{1}{6}+25\delta}$: 
\begin{align*}
\eqref{estdispsupCprocheL-int-3} &\lesssim t^{-1} 2^{-k-2j} \Vert \Psi_{l_0}^{int} \psi_{j, k}^{\widehat{\mathcal{C}}} \Vert_{L^2} \Vert g_{b-corr}(t) \Vert_{L^2} \\
&\lesssim t^{-\frac{13}{12}+100\delta} 2^{\frac{l_0}{2}-\frac{k}{2}-\frac{j}{2}} \Vert u \Vert_X \\
&\lesssim t^{-\frac{13}{12}+100\delta-\frac{1}{4}+\frac{1}{3}+\frac{\varrho}{2}+\frac{2\delta}{3}-20\delta} 2^{\frac{j}{2}+\frac{3\varrho j}{2}+\frac{3\delta j}{2}} 2^{-j+\delta k} \Vert u \Vert_X \\
&\lesssim t^{-\frac{13}{12}+100\delta-6\delta} 2^{-j+\delta k} \Vert u \Vert_X \\
&\lesssim t^{-\frac{13}{12}+94\delta} 2^{-j+\delta k} \Vert u \Vert_X 
\end{align*}
for $\varrho, \kappa, \delta$ small enough. For $2^j \gtrsim t^{-\frac{1}{6}+25\delta}$, we group \eqref{estdispsupCprocheL-int-3} with \eqref{estdispsupCprocheL-int-2}: 
\begin{align*}
\eqref{estdispsupCprocheL-int-2} + \eqref{estdispsupCprocheL-int-3} &\lesssim t^{-1} 2^{-2j} \Vert \Psi_{l_0}^{int} \psi_{j, k}^{\widehat{\mathcal{C}}} \Vert_{L^2} \Vert \psi_j \widehat{X}_{b-corr} \cdot \nabla \widehat{f}(t) \Vert_{L^2} \\
&\lesssim t^{-\frac{5}{6}+101\delta} 2^{\frac{l_0}{2}+\frac{k}{2}-\frac{3j}{2}} \Vert u \Vert_X 
\end{align*} 
On the other hand, we just saw that 
\begin{align*}
\eqref{estdispsupCprocheL-int-2} &\lesssim t^{-1} 2^{l_0-\kappa l_0-\frac{k}{2}-\frac{j}{2}-\kappa j} \langle 2^j \rangle^{-\frac{1}{2}} \Vert u \Vert_X 
\end{align*}
and we can reuse the estimate in the case $2^j \lesssim t^{-\frac{1}{6}}$ for \eqref{estdispsupCprocheL-int-3} as 
\begin{align*}
\eqref{estdispsupCprocheL-int-3} &\lesssim t^{-\frac{7}{6}+100\delta} 2^{\frac{l_0}{2}-\frac{k}{2}-j} \langle 2^j \rangle^{-1} \Vert u \Vert_X
\end{align*}
Finally, we interpolate: 
\begin{align*}
\eqref{estdispsupCprocheL-int-2} + \eqref{estdispsupCprocheL-int-3} &\lesssim \left( t^{-\frac{5}{6}+101\delta} 2^{\frac{l_0}{2}+\frac{k}{2}-\frac{3j}{2}} \Vert u \Vert_X \right)^{\frac{1}{2}} \\
&\quad \left( t^{-1} 2^{l_0-\kappa l_0-\frac{k}{2}-\frac{j}{2}-\kappa j} \langle 2^j \rangle^{-\frac{1}{2}} \Vert u \Vert_X + t^{-\frac{7}{6}+100\delta} 2^{\frac{l_0}{2}-\frac{k}{2}-j} \langle 2^j \rangle^{-1} \Vert u \Vert_X \right)^{\frac{1}{2}} \\
&\lesssim t^{-\frac{13}{12}+100 \delta} 2^{-\delta k} \left( t^{\delta+\frac{\varrho}{2}} 2^{-\frac{5j}{4}} \right)^{\frac{1}{2}} \left( t^{-\frac{5}{12}+\varrho+\frac{\kappa}{2}} 2^{-j} + t^{-\frac{1}{3}+\frac{\varrho}{2}} 2^{-\frac{3j}{4}} \right)^{\frac{1}{2}} \\
&\pushright{2^{-j+\delta k} \langle 2^j \rangle^{-\frac{1}{5}} \Vert u \Vert_X} \\
&\lesssim t^{-\frac{13}{12}+100 \delta} t^{-\frac{1}{6}+2\delta} 2^{-j-\delta k} 2^{-j+\delta k} \langle 2^j \rangle^{-\frac{1}{5}} \Vert u \Vert_X \\
&\lesssim t^{-\frac{13}{12}+100 \delta} t^{-22\delta} 2^{-j+\delta k} \langle 2^j \rangle^{-\frac{1}{6}} \Vert u \Vert_X \\
&\lesssim t^{-\frac{13}{12}+100\delta} 2^{-j+\delta k} \langle 2^j \rangle^{-\frac{1}{6}} \Vert u \Vert_X 
\end{align*}
which is enough for the first part of the Lemma. 

We now address the inequality with bound $t^{-\frac{7}{6}+2\delta} 2^{-\frac{3j}{2}+\delta k} \Vert u \Vert_X$. First, if $2^j \lesssim t^{-\frac{1}{6}+10\delta}$, we have $1 \lesssim t^{-\frac{1}{12}+5\delta} 2^{-\frac{j}{2}}$ so we can immediately obtain it from the precedent, noting that, at this scale, we always have at least $5\delta$ margin. Then, if $2^j \gtrsim t^{-\frac{1}{6}+10\delta}$, we change the threshold only for $k$. 

If $2^k \lesssim t^{-\frac{3}{4}+40\delta} 2^{-2j} \langle 2^j \rangle^{\frac{1}{4}}$, then we estimate
\begin{align*}
\eqref{estdispsupCprocheL-int} &\lesssim \Vert \Psi_{l_0}^{int} \psi_{j, k}^{\widehat{\mathcal{C}}} \Vert_{L^2_c L^{\frac{1}{1-\kappa}}_{b, a}} \Vert \psi_j \widehat{f}(t) \Vert_{L^2_c L^{\frac{1}{\kappa}}_{b, a}} \\
&\lesssim 2^{l_0-\kappa l_0 + k - \kappa k + \frac{5j}{2} - 2 \kappa j} \Vert \psi_j \widehat{f}(t) \Vert_{H^1_a L^2_{b, c}}^{\frac{1}{2}} \Vert \psi_j \widehat{f}(t) \Vert_{H^1_b L^2_{a, c}}^{\frac{1}{2}} \\
&\lesssim t^{-\frac{7}{6}+\varrho+92 \delta + \frac{5 \kappa}{4}} 2^{2 \delta j + 3 \varrho j + \frac{3\kappa j}{2}} 2^{-\frac{3j}{2} + \delta k} \langle 2^j \rangle^{-\frac{1}{4}} \Vert u \Vert_X \\
&\lesssim t^{-\frac{7}{6}+93 \delta} 2^{-\frac{3j}{2}+\delta j + \delta k} \langle 2^j \rangle^{-\frac{1}{5}} \Vert u \Vert_X 
\end{align*}
If $2^k \gtrsim t^{-\frac{3}{4}+40\delta} 2^{-2j} \langle 2^j \rangle^{\frac{1}{4}}$, we apply the integrations by parts and estimate: 
\begin{align*}
\eqref{estdispsupCprocheL-int-1} &\lesssim t^{-1} 2^{-k-3j} \Vert \Psi_{l_0}^{int} \psi_{j, k}^{\widehat{\mathcal{C}}} \Vert_{L^2_c L^{\frac{1}{1-\kappa}}_{b, a}} \Vert \psi_j(\overline{\xi}) \widehat{f}(t) \Vert_{L^2_c L^{\frac{1}{\kappa}}_{a, b}} \\
&\lesssim t^{-1} 2^{l_0-\kappa l_0-\kappa k-\frac{j}{2}-2\kappa j} \Vert \psi_j(\overline{\xi}) \widehat{f}(t) \Vert_{H^1_a L^2_{b, c}}^{\frac{1}{2}} \Vert \psi_j(\overline{\xi}) \widehat{f}(t) \Vert_{L^2_a H^1_b L^2_c}^{\frac{1}{2}} \\
&\lesssim t^{-\frac{3}{2}+\varrho+\frac{\kappa}{2}+\frac{3\kappa}{4}+\frac{3\delta}{4}} 2^{-2j+3\varrho j+\frac{3\kappa j}{2}+2\delta j+\delta k} t^{\frac{1}{12}+51\delta} 2^{-\frac{j}{2}} \langle 2^j \rangle^{-\frac{1}{2}} \Vert u \Vert_X \\
&\lesssim t^{-\frac{17}{12}+52\delta} 2^{-j+3\varrho j+\frac{3\kappa j}{2}+\delta j} 2^{-\frac{3j}{2}+\delta j+\delta k} \langle 2^j \rangle^{-\frac{1}{2}} \Vert u \Vert_X \\
&\lesssim t^{-\frac{5}{4}+52\delta} 2^{-\frac{3j}{2}+\delta j+\delta k} \langle 2^j \rangle^{-\frac{1}{2}} \Vert u \Vert_X 
\end{align*}
Then, we group \eqref{estdispsupCprocheL-int-2} with \eqref{estdispsupCprocheL-int-3} and obtain by interpolation: 
\begin{align*}
&\eqref{estdispsupCprocheL-int-2} +\eqref{estdispsupCprocheL-int-3} \\
&\lesssim \left( t^{-1} 2^{-3j} \Vert \Psi_{l_0}^{int} \psi_{j, k}^{\widehat{\mathcal{C}}} \Vert_{L^2} \Vert \psi_j(D) \nabla X_{b-corr} f(t) \Vert_{L^2} \right)^{\frac{1}{2}} \\
&\quad \left( t^{-1} 2^{-k-2j} \Vert \Psi_{l_0}^{int} \psi_{j, k}^{\widehat{\mathcal{C}}} \Vert_{L^2_{c, b} L^{\frac{1}{1-\kappa}}_a} \Vert \psi_j \widehat{h}_{b-corr}(t) \Vert_{L^2_{c, b} L^{\frac{1}{\kappa}}_a} + t^{-1} 2^{-k-2j} \Vert \Psi_{l_0}^{int} \psi_{j, k}^{\widehat{\mathcal{C}}} \Vert_{L^2} \Vert \psi_j(D) g_{b-corr}(t) \Vert_{L^2} \right)^{\frac{1}{2}} \\
&\lesssim t^{-1} 2^{\frac{l_0}{2}-\frac{3j}{2}} \left( t^{\frac{1}{6}+101\delta} \right)^{\frac{1}{2}} \left( 2^{\frac{l_0}{2}-\kappa l_0+j-\kappa j} \langle 2^j \rangle^{-\frac{1}{2}} + t^{-\frac{1}{6}+100\delta} 2^{\frac{j}{2}} \langle 2^j \rangle^{-1} \right)^{\frac{1}{2}} \Vert u \Vert_X \\
&\lesssim t^{-\frac{5}{4}+101\delta+\frac{\varrho}{2}} 2^{-\delta k-\frac{3j}{4}+\frac{3\varrho j}{2}} \left( t^{-\frac{1}{12}+\frac{\varrho}{2}-101 \delta+\frac{\kappa}{2}} 2^{\frac{j}{4}+\frac{3\varrho j}{2}+\frac{\kappa j}{2}} + 2^{\frac{j}{2}} \langle 2^j \rangle^{-\frac{1}{2}} \right)^{\frac{1}{2}} 2^{-\frac{3j}{2}+\delta k} \langle 2^j \rangle^{-\frac{1}{4}} \Vert u \Vert_X \\
&\lesssim t^{-\frac{5}{4}+101\delta+\frac{\varrho}{2}} 2^{-\delta k-\frac{j}{2}+\frac{3\varrho j}{2}} 2^{-\frac{3j}{2}+\delta k} \langle 2^j \rangle^{-\frac{1}{4}} \Vert u \Vert_X \\
&\lesssim t^{-\frac{7}{6}+97\delta} 2^{-\frac{3j}{2}+\delta j+\delta k} \langle 2^j \rangle^{-\frac{1}{5}} \Vert u \Vert_X 
\end{align*}
as wanted. 

This concludes the cas $|\overline{\xi_a}| \simeq 2^j$. 

\paragraph{Low frequencies} Finally, we turn to the case $|\overline{\xi_a}| \gg 2^j$. 

As before, we can always localise to have $\widehat{X}_a \cdot \nabla \Phi, \widehat{X}_c \cdot \nabla \Phi \ll 2^{2j}$. But this will force $\widehat{X}_b \cdot \nabla \Phi \simeq |\overline{\xi_a}|^2 \gg 2^{2j}$. 

We can then apply simpler decompositions than in the resonant cases. We skip the details. 
\end{Dem}

\subsection{Neighborhood of the line}

\begin{Lem} Let $t \geq 1$, $j \in \mathbb{Z}$ be such that $2^j \gg t^{-\frac{1}{3}}$, $k \leq -10$. Then
\begin{align*}
\Vert e^{it \omega(D)} \psi_{j, k}^{\widehat{\mathcal{L}}}(D) f(t) \Vert_{L^{\infty}} &\lesssim t^{-\frac{13}{12}+100\delta} 2^{-j+\delta j+\delta k} \langle 2^j \rangle^{-\frac{1}{2}} \Vert u \Vert_X \\
\Vert e^{it \omega(D)} \psi_{j, k}^{\widehat{\mathcal{L}}}(D) f(t) \Vert_{L^{\infty}} &\lesssim t^{-\frac{7}{6}+100\delta} 2^{-\frac{3j}{2}+\delta j + \delta k} \langle 2^j \rangle^{-\frac{1}{6}} \Vert u \Vert_X 
\end{align*}
where the universal constants only depend on $\psi, \delta$ (but not on $t, j, k$ or on the solution $u$). 
\label{lemestdisp-voisL} 
\end{Lem}

\begin{Dem}
Again, we fix $(x, y)$ and try to bound 
\begin{align}
e^{i t \omega(D)} \psi_{j, k}^{\widehat{\mathcal{L}}}(D) f(t, x, y) &= \int e^{i t \Phi} \psi_{j, k}^{\widehat{\mathcal{L}}}(\overline{\xi}) \widehat{f}(t, \overline{\xi}) ~ d\overline{\xi} \label{estdispvoisL-termetot}
\end{align}

We may deal with high frequencies, that is when $|\overline{\xi_a}| \ll 2^j$, in the same way as in the proof of Lemma \ref{lem-estdisph-zonereste}, noting that in this case
\begin{align*}
\left| \widehat{X}_a(\overline{\xi}) \cdot \nabla_{\overline{\xi}} \Phi \right| &= \left| 3 \xi_0 |\overline{\xi}| - \widehat{X}_a(\overline{\xi}) \cdot \frac{(x, y)}{t} \right| \\
&\gtrsim |\overline{\xi}|^2
\end{align*}
We can therefore also localise to have $\widehat{X}_a(\overline{\xi}) \cdot \nabla_{\overline{\xi}} \Phi \ll 2^{2j}$. 

It is thus enough to consider
\begin{align}
&\int e^{i t \Phi} \psi\left( 2^{100} 2^{-2j} \widehat{X}_a(\overline{\xi}) \cdot \nabla_{\overline{\xi}} \Phi \right) \psi_{j, k}^{\widehat{\mathcal{L}}}(\overline{\xi}) \widehat{f}(t, \overline{\xi}) ~ d\overline{\xi} \label{estdispsupL-termeinit}
\end{align}

Let us introduce the localisation symbols 
\begin{align*}
\Psi_{(l_a, l_b, l_c)}^{a-b}\left( \overline{\xi}, \frac{x}{t}, \frac{y}{t} \right) &:= \psi\left( 2^{-l_a-2j} \widehat{X}_a(\overline{\xi}) \cdot \nabla_{\overline{\xi}} \Phi \right) \psi\left( 2^{-l_b-2j} m_b(\overline{\xi}) \widehat{X}_b(\overline{\xi}) \cdot \nabla_{\overline{\xi}} \Phi \right) \\
&\quad \quad \chi\left( 2^{-l_c-2j} m_c(\overline{\xi}) \widehat{X}_c(\overline{\xi}) \cdot \nabla_{\overline{\xi}} \Phi \right) \\
\Psi_{(l_a, l_b, l_c)}^{a-c}\left( \overline{\xi}, \frac{x}{t}, \frac{y}{t} \right) &:= \psi\left( 2^{-l_a-2j} \widehat{X}_a(\overline{\xi}) \cdot \nabla_{\overline{\xi}} \Phi \right) \chi\left( 2^{-l_b-2j} m_b(\overline{\xi}) \widehat{X}_b(\overline{\xi}) \cdot \nabla_{\overline{\xi}} \Phi \right) \\
&\quad \quad \psi\left( 2^{-l_c-2j} m_c(\overline{\xi}) \widehat{X}_c(\overline{\xi}) \cdot \nabla_{\overline{\xi}} \Phi \right) \\
\Psi_{(l_a, l_b, l_c)}^{a-int}\left( \overline{\xi}, \frac{x}{t}, \frac{y}{t} \right) &:= \psi\left( 2^{-l_a-2j} \widehat{X}_a(\overline{\xi}) \cdot \nabla_{\overline{\xi}} \Phi \right) \chi\left( 2^{-l_b-2j} m_b(\overline{\xi}) \widehat{X}_b(\overline{\xi}) \cdot \nabla_{\overline{\xi}} \Phi \right) \\
&\quad \quad \chi\left( 2^{-l_c-2j} m_c(\overline{\xi}) \widehat{X}_c(\overline{\xi}) \cdot \nabla_{\overline{\xi}} \Phi \right) \\
\Psi_{(l_a, l_b, l_c)}^{b}\left( \overline{\xi}, \frac{x}{t}, \frac{y}{t} \right) &:= \chi\left( 2^{-l_a-2j} \widehat{X}_a(\overline{\xi}) \cdot \nabla_{\overline{\xi}} \Phi \right) \psi\left( 2^{-l_b-2j} m_b(\overline{\xi}) \widehat{X}_b(\overline{\xi}) \cdot \nabla_{\overline{\xi}} \Phi \right) \\
&\quad \quad \chi\left( 2^{-l_c-2j} m_c(\overline{\xi}) \widehat{X}_c(\overline{\xi}) \cdot \nabla_{\overline{\xi}} \Phi \right) \\
\Psi_{(l_a, l_b, l_c)}^{c}\left( \overline{\xi}, \frac{x}{t}, \frac{y}{t} \right) &:= \chi\left( 2^{-l_a-2j} \widehat{X}_a(\overline{\xi}) \cdot \nabla_{\overline{\xi}} \Phi \right) \chi\left( 2^{-l_b-2j} m_b(\overline{\xi}) \widehat{X}_b(\overline{\xi}) \cdot \nabla_{\overline{\xi}} \Phi \right) \\
&\quad \quad \psi\left( 2^{-l_c-2j} m_c(\overline{\xi}) \widehat{X}_c(\overline{\xi}) \cdot \nabla_{\overline{\xi}} \Phi \right) \\
\Psi_{(l_a, l_b, l_c)}^{int}\left( \overline{\xi}, \frac{x}{t}, \frac{y}{t} \right) &:= \chi\left( 2^{-l_a-2j} \widehat{X}_a(\overline{\xi}) \cdot \nabla_{\overline{\xi}} \Phi \right) \chi\left( 2^{-l_b-2j} m_b(\overline{\xi}) \widehat{X}_b(\overline{\xi}) \cdot \nabla_{\overline{\xi}} \Phi \right) \\
&\quad \quad \chi\left( 2^{-l_c-2j} m_c(\overline{\xi}) \widehat{X}_c(\overline{\xi}) \cdot \nabla_{\overline{\xi}} \Phi \right) 
\end{align*}

\paragraph{Close to $\mathcal{L}$} Assume first $|y| \ll |x|$. We compute: 
\begin{align*}
\widehat{X}_a(\overline{\xi}) \cdot \nabla_{\overline{\xi}} \Phi &= 3 \xi_0 |\overline{\xi}| - \frac{\xi_0 x}{t |\overline{\xi}|} - \frac{\xi \cdot y}{t |\overline{\xi}|} \\
m_c(\overline{\xi}) \widehat{X}_c(\overline{\xi}) \cdot \nabla_{\overline{\xi}} \Phi &= - \frac{J \xi \cdot y}{t |\overline{\xi}|} \\
m_b(\overline{\xi}) &\simeq \frac{|\xi|}{|\overline{\xi}|} \\
\frac{|\xi|}{|\overline{\xi}|} \widehat{X}_b(\overline{\xi}) \cdot \nabla_{\overline{\xi}} \Phi &= |\xi|^2 - \frac{|\xi|^2 x}{t |\overline{\xi}|^2} + \frac{\xi_0 \xi \cdot y}{t |\overline{\xi}|^2} \\
&= |\xi|^2 \left( 1 - \frac{x}{t |\overline{\xi}|^2} \right) + \frac{\xi_0 \xi \cdot y}{t |\overline{\xi}|^2} 
\end{align*}
Let us set
\begin{align*}
m_0 &:= 1 - \frac{x}{t |\overline{\xi}|^2} \\
\eta &:= \frac{\xi_0 y}{2 t |\overline{\xi}|^2 m_0} 
\end{align*}

We have that
\begin{align*}
3 \xi_0 |\overline{\xi}| &\simeq |\overline{\xi}|^2 \\
\frac{\xi_0 x}{t |\overline{\xi}|} &\simeq |\overline{\xi_a}|^2 \\
\frac{\xi \cdot y}{t |\overline{\xi}|} &\ll |\overline{\xi}|^2 \lesssim |\overline{\xi_a}|^2
\end{align*}
This means that, on the support of $\Psi_{(l_a, l_b, l_c)}^{*}$ (with $l_a \leq -100$), $|\overline{\xi}| \simeq |\overline{\xi_a}|$. More precisely, we must have $\left| 3 \xi_0 |\overline{\xi}| - \frac{\xi_0 x}{t |\overline{\xi}|} \right| \ll |\overline{\xi_a}|^2$ hence $\left| 3 |\overline{\xi}|^2 - \frac{x}{t} \right| \ll |\overline{\xi_a}|^2$. Therefore $m_0 \gtrsim 1$. Then we can rewrite 
\begin{align*}
\frac{|\xi|}{|\overline{\xi}|} \widehat{X}_b(\overline{\xi}) \cdot \nabla_{\overline{\xi}} \Phi 
&= m_0 \left( |\xi|^2 + \frac{\xi_0 \xi \cdot y}{t |\overline{\xi}|^2 m_0} \right) \\
&= m_0 \left( |\xi|^2 + 2 \xi \cdot \eta \right) \\
&= m_0 \left( |\xi + \eta|^2 - |\eta|^2 \right) 
\end{align*}

Furthermore, $|\eta| \simeq \frac{|y|}{t |\overline{\xi_a}|}$. We will denote by $\tau:= \frac{|y|}{|x|} \simeq \frac{|\eta|}{|\overline{\xi_a}|}$. 

Note that, as long as $2^{l_b} \gtrsim \tau^2$, the $b$-localisation of $\Psi_{(l_a, l_b, l_c)}$ can be approximated by
\begin{align*}
|\xi|^2 \lesssim 2^{l_b+2j}
\end{align*}
which describes a disc. In turn, if $2^{l_b} \ll \tau^2$, then the $b$-localisation of $\Psi_{(l_a, l_b, l_c)}$ describes an annulus of center $-\eta$ and radius $|\eta|$. Besides, if $l_c = l_b$ and $2^{l_b} \ll \tau^2$, then the $c$-localisation allows to separate the support of $\Psi_{(l_a, l_b, l_c)}$ in $\xi$ into two connected components, each approximated by a disc of radius $2^{l_b+j} \tau^{-1}$, one centered at $-2\eta$ and the other at $0$. 

On the other hand, the localisation $\psi_{j, k}^{\widehat{\mathcal{L}}}$ ensures $|\xi| \simeq 2^{j+k}$. 

\paragraph{1.} Let us first assume that $2^k \gg \tau$. Denote by $\Psi_{l_a}^a := \Psi_{(l_a, 0, 0)}^{a-int}$. In this case, we decompose
\begin{subequations}
\begin{align}
\eqref{estdispsupL-termeinit} &= \sum_{l_a = l_0}^{-100} 1_{2(1-\delta)k \leq l_a} \int e^{i t \Phi} \Psi_{l_a}^a \psi_{j, k}^{\widehat{\mathcal{L}}}(\overline{\xi}) \widehat{f}(t, \overline{\xi}) ~ d\overline{\xi} \label{estdispsupLpetittaua2} \\
&\quad + \sum_{l_a = l_0}^{-100} 1_{2(1-\delta)k \geq l_a} \int e^{i t \Phi} \Psi_{l_a}^a \psi_{j, k}^{\widehat{\mathcal{L}}}(\overline{\xi}) \widehat{f}(t, \overline{\xi}) ~ d\overline{\xi} \label{estdispsupLpetittauab} \\
&\quad + \int e^{i t \Phi} \Psi_{(l_a, 0, 0)}^{int} \psi_{j, k}^{\widehat{\mathcal{L}}}(\overline{\xi}) \widehat{f}(t, \overline{\xi}) ~ d\overline{\xi} \label{estdispsupLpetittauint} 
\end{align}
\end{subequations} 
for $l_0$ such that $2^{l_0} \simeq \mathfrak{t}^{-\frac{1}{2}+\varrho}$, $\varrho > 0$ small enough with respect $\delta$. 

Just like near the cone, we introduce local coordinates $a, b, c$ through a function $\Lambda$ defined locally as
\begin{align*}
\Lambda(\overline{\xi}) &:= \left( |\overline{\xi}|, 2^j \frac{\xi}{|\overline{\xi}|} \right) =: (\xi_a, \xi_b, \xi_c)
\end{align*}
Hence, 
\begin{align*}
\widehat{X}_a(\overline{\xi}) \cdot \nabla_{\overline{\xi}} \Lambda(\overline{\xi}) &= \left( 1, 0, 0 \right) \\
\widehat{X}_b(\overline{\xi}) \cdot \nabla_{\overline{\xi}} \Lambda(\overline{\xi}) &= \left( 0, -\frac{2^j \xi_0}{|\overline{\xi}|^2} \frac{\xi}{|\xi|} \right) \\
\widehat{X}_c(\overline{\xi}) \cdot \nabla_{\overline{\xi}} \Lambda(\overline{\xi}) &= \left( 0, \frac{2^j \xi_0}{|\overline{\xi}|^2} \frac{J \xi}{|\xi|} \right) 
\end{align*}
and in particular, $\mbox{det} D \Lambda \simeq 1$ uniformly. Moreover, 
\begin{align*}
(\widehat{X}_a \cdot \nabla G) \circ \Lambda^{-1} &= \partial_a (G \circ \Lambda^{-1})
\end{align*}
To avoid any singularity at $\xi = 0$, we will always remain isotropic in the coordinates $b, c$. 

We can then compute the anisotropic volumes of $\Psi_{l_a}^a$: 
\begin{align*}
\Vert \Psi_{l_a}^a \Vert_{L^{p_b}_{b, c} L^{p_a}_a} 
&\lesssim 2^{\frac{l_a}{p_a}+\frac{2k}{p_b}} |\overline{\xi_a}|^{\frac{1}{p_a}+\frac{2}{p_b}} 
\end{align*}

On the other hand, we have for derivatives of symbols that 
\begin{align*}
\widehat{X}_a(\overline{\xi}) \cdot \nabla_{\overline{\xi}} \left[ \widehat{X}_a(\overline{\xi}) \cdot \nabla_{\overline{\xi}} \Phi \right] &= 6 \xi_0 \simeq 2^j \\
m_b(\overline{\xi}) \widehat{X}_b(\overline{\xi}) \cdot \nabla_{\overline{\xi}} \left[ \widehat{X}_a(\overline{\xi}) \cdot \nabla_{\overline{\xi}} \Phi \right] &\lesssim 2^{j+2k} + \tau^2 2^j \lesssim 2^{j+2k} \\
m_c(\overline{\xi}) \widehat{X}_c(\overline{\xi}) \cdot \nabla_{\overline{\xi}} \left[ \widehat{X}_a(\overline{\xi}) \cdot \nabla_{\overline{\xi}} \Phi \right] &\lesssim \tau^2 2^j \ll 2^{j+2k} \\
\widehat{X}_a(\overline{\xi}) \cdot \nabla_{\overline{\xi}} \psi_{j, k}^{\widehat{\mathcal{L}}} &\lesssim 2^{-j} \\
m_b(\overline{\xi}) \widehat{X}_b(\overline{\xi}) \cdot \nabla_{\overline{\xi}} \psi_{j, k}^{\widehat{\mathcal{L}}} &\lesssim 2^{-j} \\
m_c(\overline{\xi}) \widehat{X}_c(\overline{\xi}) \cdot \nabla_{\overline{\xi}} \psi_{j, k}^{\widehat{\mathcal{L}}} &\lesssim 2^{-j} 
\end{align*}
In particular, 
\begin{align*}
\widehat{X}_a(\overline{\xi}) \cdot \nabla_{\overline{\xi}} \left[ \Psi_{l_a}^a \psi_{j, k}^{\widehat{\mathcal{L}}} \right] &\lesssim 2^{-l_a-j} 
\end{align*}
If $2k \leq l_a+10$, we also have
\begin{align*}
m_b(\overline{\xi}) \widehat{X}_b(\overline{\xi}) \cdot \nabla_{\overline{\xi}} \left[ \Psi_{l_a}^a \psi_{j, k}^{\widehat{\mathcal{L}}} \right] &\lesssim 2^{-j} \\
m_c(\overline{\xi}) \widehat{X}_c(\overline{\xi}) \cdot \nabla_{\overline{\xi}} \left[ \Psi_{l_a}^a \psi_{j, k}^{\widehat{\mathcal{L}}} \right] &\lesssim 2^{-j} 
\end{align*}
If $2k \geq l_a+10$, we introduce the corrected fields
\begin{align*}
\widehat{X}_{b-corr}(\overline{\xi}) &:= \widehat{X}_b(\overline{\xi}) - \frac{\widehat{X}_b(\overline{\xi}) \cdot \nabla_{\overline{\xi}} \left[ \widehat{X}_a(\overline{\xi}) \cdot \nabla_{\overline{\xi}} \Phi \right]}{\widehat{X}_a(\overline{\xi}) \cdot \nabla_{\overline{\xi}} \left[ \widehat{X}_a(\overline{\xi}) \cdot \nabla_{\overline{\xi}} \Phi \right]} \widehat{X}_a(\overline{\xi}) \\
\widehat{X}_{c-corr}(\overline{\xi}) &:= \widehat{X}_c(\overline{\xi}) - \frac{\widehat{X}_c(\overline{\xi}) \cdot \nabla_{\overline{\xi}} \left[ \widehat{X}_a(\overline{\xi}) \cdot \nabla_{\overline{\xi}} \Phi \right]}{\widehat{X}_a(\overline{\xi}) \cdot \nabla_{\overline{\xi}} \left[ \widehat{X}_a(\overline{\xi}) \cdot \nabla_{\overline{\xi}} \Phi \right]} \widehat{X}_a(\overline{\xi})
\end{align*}
so that 
\begin{align*}
m_b(\overline{\xi}) \widehat{X}_{b-corr}(\overline{\xi}) \cdot \nabla_{\overline{\xi}} \left[ \Psi_{l_a}^a \psi_{j, k}^{\widehat{\mathcal{L}}} \right] &\lesssim 2^{-j} \\
m_c(\overline{\xi}) \widehat{X}_{c-corr}(\overline{\xi}) \cdot \nabla_{\overline{\xi}} \left[ \Psi_{l_a}^a \psi_{j, k}^{\widehat{\mathcal{L}}} \right] &\lesssim 2^{-j}
\end{align*}
and, on the support of $\Psi_{l_a}^a \psi_{j, k}^{\widehat{\mathcal{L}}}$ we have indeed
\begin{align*}
m_b(\overline{\xi}) \widehat{X}_{b-corr}(\overline{\xi}) \cdot \nabla_{\overline{\xi}} \Phi &= m_b(\overline{\xi}) \widehat{X}_b(\overline{\xi}) \cdot \nabla_{\overline{\xi}} \Phi + O(2^{2k+l_a+2j}) \simeq 2^{2k+2j} 
\end{align*}
as before. 

We now apply integrations by parts. On \eqref{estdispsupLpetittaua2}, we apply at most $n$ (large with respect to $\varrho^{-1}$) integrations by parts along $\widehat{X}_a$, as long as none hits $\widehat{f}(t)$: 
\begin{subequations}
\begin{align}
\eqref{estdispsupLpetittaua2} &= \sum_{l_a = l_0}^{-200} 1_{2(1-\delta)k \leq l_a} t^{-n} 2^{-2nl_a-3nj} \int e^{i t \Phi} \Psi_{l_a}^a \psi_{j, k}^{\widehat{\mathcal{L}}}(\overline{\xi}) \widehat{f}(t, \overline{\xi}) ~ d\overline{\xi} \label{estdispsupLpetittaua2-1} \\
&\quad + \sum_{i = 1}^n \sum_{l_a = l_0}^{-200} 1_{2(1-\delta)k \leq l_a} t^{-i-1} 2^{-2il_a-l_a-3ij-2j} \int e^{i t \Phi} \Psi_{l_a}^a \psi_{j, k}^{\widehat{\mathcal{L}}}(\overline{\xi}) \widehat{h}_a(t, \overline{\xi}) ~ d\overline{\xi} \label{estdispsupLpetittaua2-2} \\
&\quad + \sum_{i = 1}^n \sum_{l_a = l_0}^{-200} 1_{2(1-\delta)k \leq l_a} t^{-i-1} 2^{-2il_a-3ij-j} \int e^{i t \Phi} \Psi_{l_a}^a \psi_{j, k}^{\widehat{\mathcal{L}}}(\overline{\xi}) \widehat{X}_a(\overline{\xi}) \cdot \nabla_{\overline{\xi}} \widehat{h}_a(t, \overline{\xi}) ~ d\overline{\xi} \label{estdispsupLpetittaua2-3} \\
&\quad + \sum_{i = 1}^n \sum_{l_a = l_0}^{-200} 1_{2(1-\delta)k \leq l_a} t^{-i} 2^{-2il_a+l_a-3ij+j} \int e^{i t \Phi} \Psi_{l_a}^a \psi_{j, k}^{\widehat{\mathcal{L}}}(\overline{\xi}) \widehat{g}_a(t, \overline{\xi}) ~ d\overline{\xi} \label{estdispsupLpetittaua2-4} 
\end{align}
\end{subequations} 
We then estimate: 
\begin{align*}
\eqref{estdispsupLpetittaua2-1} &\lesssim \sum_{l_a = l_0}^{-200} 1_{2(1-\delta)k \leq l_a} t^{-n} 2^{-2nl_a-3nj+j} \Vert \Psi_{l_a}^a \psi_{j, k}^{\widehat{\mathcal{L}}} \Vert_{L^2} \langle 2^j \rangle^{-1} \Vert \langle \nabla \rangle |\nabla|^{-1} f(t) \Vert_{L^2} \\
&\lesssim \sum_{l_a = l_0}^{-200} 1_{2k \leq l_a} \mathfrak{t}^{-n} 2^{-2nl_a+\frac{l_a}{2}+k+\frac{5j}{2}} \langle 2^j \rangle^{-1} \Vert u \Vert_X \\
&\lesssim t^{-\frac{7}{6}} \mathfrak{t}^{-n+\frac{7}{6}} 2^{-2nl_0+\frac{l_0}{2}} 2^{-j+k} \langle 2^j \rangle^{-1} \Vert u \Vert_X \\
&\lesssim t^{-\frac{7}{6}} \mathfrak{t}^{\frac{7}{6}-2n\varrho} 2^{-j+k} \langle 2^j \rangle^{-1} \Vert u \Vert_X \\
&\lesssim t^{-\frac{7}{6}} 2^{-j+k} \langle 2^j \rangle^{-1} \Vert u \Vert_X \\
\eqref{estdispsupLpetittaua2-2} &\lesssim \sum_{i = 1}^n \sum_{l_a = l_0}^{-200} 1_{2(1-\delta)k \leq l_a} t^{-i-1} 2^{-2il_a-l_a-3ij-2j} \Vert \Psi_{l_a}^a \psi_{j, k}^{\widehat{\mathcal{L}}} \Vert_{L^2_{b, c} L^{\frac{1}{1-\kappa}}_a} \Vert \psi_j \widehat{h}_a(t) \Vert_{L^2_{b, c} L^{\frac{1}{\kappa}}_a} \\
&\lesssim \sum_{i = 1}^n \sum_{l_a = l_0}^{-200} 1_{2k \leq l_a} t^{-\frac{7}{6}} \mathfrak{t}^{-i+\frac{1}{6}} 2^{-2il_a-\kappa l_a+k-\frac{j}{2}-\kappa j} \Vert \psi_j \widehat{h}_a(t) \Vert_{L^2_{b, c} H^1_a}^{\frac{1}{2}} \Vert \psi_j \widehat{h}_a(t) \Vert_{L^2}^{\frac{1}{2}} \\
&\lesssim t^{-\frac{7}{6}+\frac{\kappa}{3}} \mathfrak{t}^{-\frac{1}{12}+\frac{\kappa}{6}+\frac{\delta}{4}} 2^{-j+\delta k} \langle 2^j \rangle^{-\frac{1}{2}} \Vert u \Vert_X \\
&\lesssim t^{-\frac{7}{6}+\frac{\kappa}{3}} 2^{-j+\delta k} \langle 2^j \rangle^{-\frac{3}{5}} \Vert u \Vert_X \\
\eqref{estdispsupLpetittaua2-3} &\lesssim \sum_{i = 1}^n \sum_{l_a = l_0}^{-200} 1_{2(1-\delta)k \leq l_a} t^{-i-1} 2^{-2il_a-3ij-2j} \Vert \Psi_{l_a}^a \psi_{j, k}^{\widehat{\mathcal{L}}} \Vert_{L^2} \Vert \nabla X_a h_a(t) \Vert_{L^2} \\
&\lesssim \sum_{i = 1}^n \sum_{l_a = l_0}^{-200} 1_{2k \leq l_a} t^{-\frac{7}{6}} \mathfrak{t}^{-i+\frac{1}{6}} 2^{-2il_a+\frac{l_a}{2}+k-j} \Vert u \Vert_X \\
&\lesssim t^{-\frac{7}{6}} \mathfrak{t}^{-\frac{1}{3}+\frac{\delta}{4}} 2^{-j+\delta k} \Vert u \Vert_X \\
&\lesssim t^{-\frac{7}{6}} 2^{-j+\delta k} \langle 2^j \rangle^{-\frac{3}{4}} \Vert u \Vert_X \\
\eqref{estdispsupLpetittaua2-4} &\lesssim \sum_{i = 1}^n \sum_{l_a = l_0}^{-200} 1_{2(1-\delta)k \leq l_a} t^{-i} 2^{-2il_a+l_a-3ij+\frac{j}{2}} \Vert \Psi_{l_a}^a \psi_{j, k}^{\widehat{\mathcal{L}}} \Vert_{L^2} \langle 2^j \rangle^{-1} \Vert m_{\widehat{\mathcal{L}}}(D) \langle \nabla \rangle |\nabla|^{\frac{1}{2}} g_a(t) \Vert_{L^2} \\
&\lesssim \sum_{i = 1}^n \sum_{l_a = l_0}^{-200} 1_{2(1-\delta)k \leq l_a} t^{-\frac{1}{6}+100\delta} \mathfrak{t}^{-i} 2^{-2il_a+\frac{3l_a}{2}+k+2j} \langle 2^j \rangle^{-1} \Vert u \Vert_X \\
&\lesssim t^{-\frac{7}{6}+100\delta} 2^{-j+\delta k} \langle 2^j \rangle^{-1} \Vert u \Vert_X 
\end{align*}
for $n \varrho$ large enough, $\delta$ small enough.  

Then, for \eqref{estdispsupLpetittauab}, we start by an integration by parts along $\widehat{X}_{b-corr}$, and then only apply integrations by parts along $\widehat{X}_a$: 
\begin{subequations}
\begin{align}
\eqref{estdispsupLpetittauab} &= \sum_{l_a = l_0}^{-200} 1_{2(1-\delta) k \geq l_a} t^{-2} 2^{-2l_a-2k-6j} \int e^{i t \Phi} \Psi_{l_a}^a \psi_{j, k}^{\widehat{\mathcal{L}}}(\overline{\xi}) \widehat{f}(t, \overline{\xi}) ~ d\overline{\xi} \label{estdispsupLpetittauab-1} \\
&\quad + \sum_{l_a = l_0}^{-200} 1_{2(1-\delta) k \geq l_a} t^{-2} 2^{-l_a-2k-5j} \int e^{i t \Phi} \Psi_{l_a}^a \psi_{j, k}^{\widehat{\mathcal{L}}}(\overline{\xi}) \widehat{X}_a(\overline{\xi}) \cdot \nabla_{\overline{\xi}} \widehat{f}(t, \overline{\xi}) ~ d\overline{\xi} \label{estdispsupLpetittauab-2} \\
&\quad + \sum_{l_a = l_0}^{-200} 1_{2(1-\delta) k \geq l_a} t^{-2} 2^{-2l_a-2k-5j} \int e^{i t \Phi} \Psi_{l_a}^a \psi_{j, k}^{\widehat{\mathcal{L}}}(\overline{\xi}) \widehat{h}_{b-corr}(t, \overline{\xi}) ~ d\overline{\xi} \label{estdispsupLpetittauab-3} \\
&\quad + \sum_{l_a = l_0}^{-200} 1_{2(1-\delta) k \geq l_a} t^{-2} 2^{-l_a-2k-4j} \int e^{i t \Phi} \Psi_{l_a}^a \psi_{j, k}^{\widehat{\mathcal{L}}}(\overline{\xi}) \widehat{X}_a(\overline{\xi}) \cdot \nabla_{\overline{\xi}} \widehat{h}_{b-corr}(t, \overline{\xi}) ~ d\overline{\xi} \label{estdispsupLpetittauab-4} \\
&\quad + \sum_{l_a = l_0}^{-200} 1_{2(1-\delta) k \geq l_a} t^{-1} 2^{-2k-2j} \int e^{i t \Phi} \Psi_{l_a}^a \psi_{j, k}^{\widehat{\mathcal{L}}}(\overline{\xi}) \widehat{g}_{b-corr}(t, \overline{\xi}) ~ d\overline{\xi} \label{estdispsupLpetittauab-5} 
\end{align}
\end{subequations} 
We then estimate: 
\begin{align*}
\eqref{estdispsupLpetittauab-1} &\lesssim \sum_{l_a = l_0}^{-200} 1_{2(1-\delta) k \geq l_a} t^{-2} 2^{-2l_a-2k-6j} \Vert \Psi_{l_a}^a \psi_{j, k}^{\widehat{\mathcal{L}}} \Vert_{L^2_{b, c} L^{\frac{1}{1-\kappa}}_a} \Vert \psi_j \widehat{f}(t) \Vert_{L^2_{b, c} L^{\frac{1}{\kappa}}_a} \\
&\lesssim \sum_{l_a = l_0}^{-200} 1_{2(1-\delta) k \geq l_a} t^{-\frac{7}{6}+\frac{\kappa}{3}} \mathfrak{t}^{-\frac{5}{6}-\frac{\kappa}{3}} 2^{-l_a-\kappa l_a-k-\frac{3j}{2}} \Vert \psi_j \widehat{f}(t) \Vert_{L^2_{b, c} H^1_a}^{\frac{1}{2}} \Vert \psi_j \widehat{f}(t) \Vert_{L^2}^{\frac{1}{2}} \\
&\lesssim t^{-\frac{7}{6}+\frac{\kappa}{3}} \mathfrak{t}^{-\frac{5}{6}-\frac{\kappa}{3}} 2^{-l_0-\kappa l_0-\frac{1+\delta}{2(1-\delta)}l_0} 2^{-j+\delta k} \langle 2^j \rangle^{-1} \Vert u \Vert_X \\
&\lesssim t^{-\frac{7}{6}+\frac{\kappa}{3}} \mathfrak{t}^{-\frac{1}{12}+\frac{\kappa}{6}+\delta} 2^{-j+\delta k} \langle 2^j \rangle^{-1} \Vert u \Vert_X \\
&\lesssim t^{-\frac{7}{6}+\frac{\kappa}{3}} 2^{-j+\delta k} \langle 2^j \rangle^{-1} \Vert u \Vert_X \\
\eqref{estdispsupLpetittauab-2} &\lesssim \sum_{l_a = l_0}^{-200} 1_{2(1-\delta) k \geq l_a} t^{-2} 2^{-l_a-2k-5j} \Vert \Psi_{l_a}^a \psi_{j, k}^{\widehat{\mathcal{L}}} \Vert_{L^2} \langle 2^j \rangle^{-1} \Vert \langle \nabla \rangle X_a f(t) \Vert_{L^2} \\
&\lesssim \sum_{l_a = l_0}^{-200} 1_{2(1-\delta) k \geq l_a} t^{-\frac{7}{6}} \mathfrak{t}^{-\frac{5}{6}} 2^{-\frac{l_a}{2}-k-j} \langle 2^j \rangle^{-1} \Vert u \Vert_X \\
&\lesssim t^{-\frac{7}{6}} \mathfrak{t}^{-\frac{1}{12}} 2^{-j+\delta k} \langle 2^j \rangle^{-1} \Vert u \Vert_X \\
&\lesssim t^{-\frac{7}{6}} 2^{-j+\delta k} \langle 2^j \rangle^{-1} \Vert u \Vert_X \\
\eqref{estdispsupLpetittauab-3} &\lesssim \sum_{l_a = l_0}^{-200} 1_{2(1-\delta) k \geq l_a} t^{-2} 2^{-2l_a-2k-5j} \Vert \Psi_{l_a}^a \psi_{j, k}^{\widehat{\mathcal{L}}} \Vert_{L^2_{b, c} L^{\frac{1}{1-\kappa}}_a} \Vert \psi_j \widehat{h}_{b-corr}(t) \Vert_{L^2_{b, c} L^{\frac{1}{\kappa}}_a} \\
&\lesssim \sum_{l_a = l_0}^{-200} 1_{2(1-\delta) k \geq l_a} t^{-\frac{7}{6}} \mathfrak{t}^{-\frac{5}{6}} 2^{-l_a-\kappa l_a-k-\frac{j}{2}-\kappa j} \Vert \psi_j \widehat{h}_{b-corr}(t) \Vert_{L^2_{b, c} H^1_a}^{\frac{1}{2}} \Vert \psi_j \widehat{h}_{b-corr}(t) \Vert_{L^2}^{\frac{1}{2}} \\
&\lesssim t^{-\frac{7}{6}+\frac{\kappa}{3}} \mathfrak{t}^{-\frac{1}{12}+\frac{\kappa}{6}+\delta} 2^{-j+\delta k} \langle 2^j \rangle^{-\frac{1}{2}} \Vert u \Vert_X \\
&\lesssim t^{-\frac{7}{6}+\frac{\kappa}{3}} 2^{-j+\delta k} \langle 2^j \rangle^{-\frac{3}{5}} \Vert u \Vert_X \\
\eqref{estdispsupLpetittauab-4} &\lesssim \sum_{l_a = l_0}^{-200} 1_{2(1-\delta) k \geq l_a} t^{-2} 2^{-l_a-2k-5j} \Vert \Psi_{l_a}^a \psi_{j, k}^{\widehat{\mathcal{L}}} \Vert_{L^2} \Vert \nabla X_a h_{b-corr}(t) \Vert_{L^2} \\
&\lesssim \sum_{l_a = l_0}^{-200} 1_{2(1-\delta) k \geq l_a} t^{-\frac{7}{6}} \mathfrak{t}^{-\frac{5}{6}} 2^{-\frac{l_a}{2}-k-j} \Vert u \Vert_X \\
&\lesssim t^{-\frac{7}{6}} \mathfrak{t}^{-\frac{1}{3}+\delta} 2^{-j+\delta k} \Vert u \Vert_X \\
&\lesssim t^{-\frac{7}{6}} 2^{-j+\delta k} \langle 2^j \rangle^{-\frac{3}{4}} \Vert u \Vert_X 
\end{align*}
Finally, for \eqref{estdispsupLpetittauab-5}, we separate into: 
\begin{subequations}
\begin{align}
\eqref{estdispsupLpetittauab-5} &= \sum_{l_a = l_0}^{-200} 1_{2(1-\delta) k \geq l_a \geq 2(1+\delta) k} t^{-1} 2^{-2k-2j} \int e^{i t \Phi} \Psi_{l_a}^a \psi_{j, k}^{\widehat{\mathcal{L}}}(\overline{\xi}) \widehat{g}_{b-corr}(t, \overline{\xi}) ~ d\overline{\xi} \label{estdispsupLpetittauab-5-1} \\
&\quad + \sum_{l_a = l_0}^{-200} 1_{2(1+\delta) k \geq l_a} t^{-1} 2^{-2k-2j} \int e^{i t \Phi} \Psi_{l_a}^a \psi_{j, k}^{\widehat{\mathcal{L}}}(\overline{\xi}) \widehat{g}_{b-corr}(t, \overline{\xi}) ~ d\overline{\xi} \label{estdispsupLpetittauab-5-2} 
\end{align}
\end{subequations}
Then, 
\begin{align*}
\eqref{estdispsupLpetittauab-5-2} &\lesssim \sum_{l_a = l_0}^{-200} 1_{2(1+\delta) k \geq l_a} t^{-1} 2^{-2k-\frac{5j}{2}} \Vert \Psi_{l_a}^a \psi_{j, k}^{\widehat{\mathcal{L}}} \Vert_{L^2} \langle 2^j \rangle^{-1} \Vert m_{\widehat{\mathcal{L}}}(D) \langle \nabla \rangle |\nabla|^{\frac{1}{2}} g_{b-corr}(t) \Vert_{L^2} \\
&\lesssim \sum_{l_a = l_0}^{-200} 1_{2(1+\delta) k \geq l_a} t^{-\frac{7}{6}+100\delta} 2^{-k+\frac{l_a}{2}-j} \langle 2^j \rangle^{-1} \Vert u \Vert_X \\
&\lesssim t^{-\frac{7}{6}+100\delta} 2^{-j+\delta k} \langle 2^j \rangle^{-1} \Vert u \Vert_X
\end{align*}
On the other hand, for \eqref{estdispsupLpetittauab-5-1}, if $2^{|j|} \gtrsim t^{1000\delta}$, we can estimate in the same way: 
\begin{align*}
\eqref{estdispsupLpetittauab-5-1} &\lesssim \sum_{l_a = l_0}^{-200} 1_{2(1-\delta) k \geq l_a \geq 2(1+\delta) k} t^{-1} 2^{-2k-\frac{5j}{2}} \Vert \Psi_{l_a}^a \psi_{j, k}^{\widehat{\mathcal{L}}} \Vert_{L^2} \langle 2^j \rangle^{-1} \Vert m_{\widehat{\mathcal{L}}}(D) \langle \nabla \rangle |\nabla|^{\frac{1}{2}} g_{b-corr}(t) \Vert_{L^2} \\
&\lesssim \sum_{l_a = l_0}^{-200} 1_{2(1-\delta) k \geq l_a \geq 2(1+\delta) k} t^{-\frac{7}{6}+100\delta} 2^{-k+\frac{l_a}{2}-j} \langle 2^j \rangle^{-1} \Vert u \Vert_X \\
&\lesssim t^{-\frac{7}{6}+100\delta} \mathfrak{t}^{\delta} 2^{-j+\delta k} \langle 2^j \rangle^{-1} \Vert u \Vert_X
\end{align*}
so that 
\begin{align*}
\eqref{estdispsupLpetittauab-5-1} &\lesssim t^{-\frac{13}{12}} 2^{-j+\delta k} \langle 2^j \rangle^{-\frac{3}{4}} \Vert u \Vert_X 
\end{align*}
for small enough $\delta$, but also using $2^{|j|} \gtrsim t^{1000 \delta}$: 
\begin{align*}
\eqref{estdispsupLpetittauab-5-1} &\lesssim t^{-\frac{7}{6}} 2^{-\frac{3j}{2}+\delta k} \langle 2^j \rangle^{-\frac{1}{4}} \Vert u \Vert_X 
\end{align*}
as wanted. Finally, if $2^{|j|} \lesssim t^{1000 \delta}$, since we restricted our attention to a sum in $l_a$ such that $2(1-\delta) k \geq l_a \geq 2(1+\delta) k$, we can use $\widehat{X}_b$ everywhere instead of $\widehat{X}_{b-corr}$: this modifies slightly the estimates of \eqref{estdispsupLpetittauab-1} to \eqref{estdispsupLpetittauab-4}, losing possibly a factor $2^{2k-l_a} \lesssim \mathfrak{t}^{\delta}$ that can be absorbed into all the estimates. We then use the fine estimate on $g_b$: 
\begin{align*}
\eqref{estdispsupLpetittauab-5-1} &\lesssim \sum_{l_a = l_0}^{-200} 1_{2(1-\delta) k \geq l_a \geq 2(1+\delta) k} t^{-1} 2^{-2k-2j} \Vert \Psi_{l_a}^a \psi_{j, k}^{\widehat{\mathcal{L}}} \Vert_{L^2} \Vert \psi_{j, k}^{\widehat{\mathcal{L}}}(D) g_b(t) \Vert_{L^2} \\
&\lesssim \sum_{l_a = l_0}^{-200} 1_{2(1-\delta) k \geq l_a \geq 2(1+\delta) k} t^{-1} 2^{\frac{l_a}{2}-k-\frac{j}{2}} \left( t^{-\frac{1}{6}+100\delta} 2^{k-\frac{j}{2}} \langle 2^j \rangle^{-1} + t^{-\frac{1}{2}} 2^{-3j} \right) \Vert u \Vert_X \\
&\lesssim t^{-\frac{7}{6}+100\delta} 2^{-j} \langle 2^j \rangle^{-1} \Vert u \Vert_X + t^{-\frac{3}{2}} \mathfrak{t}^{\delta} 2^{-\frac{7j}{2}+\delta k} \Vert u \Vert_X \\
&\lesssim t^{-\frac{7}{6}+100\delta} 2^{-j} \langle 2^j \rangle^{-1} \Vert u \Vert_X + t^{-\frac{7}{6}} 2^{-j+\delta k} \Vert u \Vert_X 
\end{align*}
We also used $2^{|j|} \lesssim t^{1000 \delta}$ and that $\delta$ is small enough.  

It only remains \eqref{estdispsupLpetittauint}. We exhaust the cases depending on $t, j, k$. 

Assume first $2^j \lesssim t^{-\frac{2}{9}}$. If $2^k \lesssim t^{-\frac{5}{12}+10\delta} 2^{-j}$, then we estimate
\begin{align*}
\eqref{estdispsupLpetittauint} &\lesssim \Vert \Psi_{(l_0, 0, 0)}^{int} \psi_{j, k}^{\widehat{\mathcal{L}}} \Vert_{L^2_{b, c} L^{\frac{1}{1-\kappa}}_a} \Vert \psi_j \widehat{f}(t) \Vert_{L^2_{b, c} L^{\frac{1}{1-\kappa}}_a} \\
&\lesssim 2^{l_0-\kappa l_0+k+2j-\kappa j} \Vert \psi_j \widehat{f}(t) \Vert_{L^2_{b, c} H^1_a}^{\frac{1}{2}} \Vert \psi_j \widehat{f}(t) \Vert_{L^2}^{\frac{1}{2}} \\
&\lesssim t^{-\frac{13}{12}} t^{\frac{1}{6}+11\delta+\varrho+\frac{\kappa}{2}} 2^{j+3\varrho j+\frac{\kappa j}{2}} 2^{-j+\delta j+\delta k} \Vert u \Vert_X \\
&\lesssim t^{-\frac{13}{12}} t^{-\frac{1}{18}+11\delta+\varrho+\frac{\kappa}{2}} 2^{-j+\delta j+\delta k} \Vert u \Vert_X \\
&\lesssim t^{-\frac{13}{12}} 2^{-j+\delta j+\delta k} \Vert u \Vert_X 
\end{align*}
for $\varrho, \kappa, \delta$ small enough. On the other hand, if $2^k \gtrsim t^{-\frac{5}{12}+10\delta} 2^{-j}$, then we apply an integration by parts along $\widehat{X}_{b-corr}$: 
\begin{subequations}
\begin{align}
\eqref{estdispsupLpetittauint} &= t^{-1} 2^{-2k-3j} \int e^{i t \Phi} \Psi_{(l_0, 0, 0)}^{int} \psi_{j, k}^{\widehat{\mathcal{L}}}(\overline{\xi}) \widehat{f}(t, \overline{\xi}) ~ d\overline{\xi} \label{estdispsupLpetittauint-1} \\
&\quad + t^{-1} 2^{-2k-2j} \int e^{i t \Phi} \Psi_{(l_0, 0, 0)}^{int} \psi_{j, k}^{\widehat{\mathcal{L}}}(\overline{\xi}) \widehat{h}_{b-corr}(t, \overline{\xi}) ~ d\overline{\xi} \label{estdispsupLpetittauint-2} \\
&\quad + t^{-1} 2^{-2k-2j} \int e^{i t \Phi} \Psi_{(l_0, 0, 0)}^{int} \psi_{j, k}^{\widehat{\mathcal{L}}}(\overline{\xi}) \widehat{g}_{b-corr}(t, \overline{\xi}) ~ d\overline{\xi} \label{estdispsupLpetittauint-3} 
\end{align}
\end{subequations}
We then estimate: 
\begin{align*}
\eqref{estdispsupLpetittauint-2} &\lesssim t^{-1} 2^{-2k-2j} \Vert \Psi_{(l_0, 0, 0)}^{int} \psi_{j, k}^{\widehat{\mathcal{L}}} \Vert_{L^2_{b, c} L^{\frac{1}{1-\kappa}}_a} \Vert \psi_j \widehat{h}_{b-corr}(t) \Vert_{L^2_{b, c} L^{\frac{1}{\kappa}}_a} \\
&\lesssim t^{-1} 2^{l_0-\kappa l_0-k-\kappa j} \Vert \psi_j \widehat{h}_{b-corr}(t) \Vert_{L^2_{b, c} H^1_a}^{\frac{1}{2}} \Vert \psi_j \widehat{h}_{b-corr}(t) \Vert_{L^2}^{\frac{1}{2}} \\
&\lesssim t^{-\frac{13}{12}} t^{\varrho+\frac{\kappa}{2}-9\delta} 2^{3\varrho j + \frac{\kappa j}{2}} 2^{-j+\delta j+\delta k} \Vert u \Vert_X \\
&\lesssim t^{-\frac{13}{12}} 2^{-j+\delta j+\delta k} \Vert u \Vert_X \\
\eqref{estdispsupLpetittauint-3} &\lesssim t^{-1} 2^{-2k-2j} \Vert \Psi_{(l_0, 0, 0)}^{int} \psi_{j, k}^{\widehat{\mathcal{L}}} \Vert_{L^2} \Vert m_{\widehat{\mathcal{L}}}(D) g_{b-corr}(t) \Vert_{L^2} \\
&\lesssim t^{-\frac{13}{12}+100\delta} 2^{\frac{l_0}{2}-k-\frac{j}{2}} \Vert u \Vert_X \\
&\lesssim t^{-\frac{13}{12}+100\delta} t^{\frac{1}{6}+\frac{\varrho}{2}-9\delta} 2^{\frac{3j}{4}+3\varrho j} 2^{-j+\delta j+\delta k} \Vert u \Vert_X \\
&\lesssim t^{-\frac{13}{12}+100\delta} t^{\frac{\varrho}{2}-9\delta} 2^{-j+\delta j+\delta k} \Vert u \Vert_X \\
&\lesssim t^{-\frac{13}{12}+100\delta} 2^{-j+\delta j+\delta k} \Vert u \Vert_X 
\end{align*}
for $\varrho, \kappa$ small enough with respect to $\delta$. \eqref{estdispsupLpetittauint-1} can be estimated like \eqref{estdispsupLpetittauint-2}. 

Above, we can automatically improve the estimate to get the $t^{-\frac{7}{6}}$ decay up to losing a factor $2^{\frac{j}{2}}$: indeed, $1 \lesssim t^{-\frac{1}{9}} 2^{-\frac{j}{2}} \lesssim t^{-\frac{1}{12}} 2^{-\frac{j}{2}}$. 

Assume now that $t^{-\frac{2}{9}} \lesssim 2^j \lesssim t^{-\frac{5}{24}}$. If $2^k \lesssim t^{-\frac{1}{4}+10\delta} 2^{-\frac{j}{4}}$, then we estimate
\begin{align*}
\eqref{estdispsupLpetittauint} &\lesssim \Vert \Psi_{(l_0, 0, 0)}^{int} \psi_{j, k}^{\widehat{\mathcal{L}}} \Vert_{L^2_{b, c} L^{\frac{1}{1-\kappa}}_a} \Vert \psi_j \widehat{f}(t) \Vert_{L^2_{b, c} L^{\frac{1}{1-\kappa}}_a} \\
&\lesssim 2^{l_0-\kappa l_0+k+2j-\kappa j} \Vert \psi_j \widehat{f}(t) \Vert_{L^2_{b, c} H^1_a}^{\frac{1}{2}} \Vert \psi_j \widehat{f}(t) \Vert_{L^2}^{\frac{1}{2}} \\
&\lesssim t^{-\frac{13}{12}} t^{\frac{1}{3}+\varrho+\frac{\kappa}{2}+10\delta+\frac{\delta}{4}} 2^{\frac{7j}{4}+3\varrho j+\frac{\kappa j}{2}+\frac{\delta j}{4}} 2^{-j+\delta k} \Vert u \Vert_X \\
&\lesssim t^{-\frac{13}{12}} t^{-\frac{1}{32}+\varrho+\frac{\kappa}{2}+10\delta+\frac{\delta}{4}} 2^{-j+\delta k} \Vert u \Vert_X \\
&\lesssim t^{-\frac{13}{12}} 2^{-j+\delta k} \Vert u \Vert_X 
\end{align*}
for small enough $\varrho, \kappa, \delta$. On the other hand, if $2^k \gtrsim t^{-\frac{1}{4}+10\delta} 2^{-\frac{j}{4}}$, then we apply the same integration by parts as before and estimate: 
\begin{align*}
\eqref{estdispsupLpetittauint-2} &\lesssim t^{-1} 2^{-2k-2j} \Vert \Psi_{(l_0, 0, 0)}^{int} \psi_{j, k}^{\widehat{\mathcal{L}}} \Vert_{L^2_{b, c} L^{\frac{1}{1-\kappa}}_a} \Vert \psi_j \widehat{h}_{b-corr}(t) \Vert_{L^2_{b, c} L^{\frac{1}{\kappa}}_a} \\
&\lesssim t^{-1} 2^{l_0-\kappa l_0-k-\kappa j} \Vert \psi_j \widehat{h}_{b-corr}(t) \Vert_{L^2_{b, c} H^1_a}^{\frac{1}{2}} \Vert \psi_j \widehat{h}_{b-corr}(t) \Vert_{L^2}^{\frac{1}{2}} \\
&\lesssim t^{-\frac{13}{12}} t^{-\frac{1}{6}-9\delta+\varrho+\frac{\kappa}{2}} 2^{-\frac{3j}{4}+3\varrho j+\frac{\kappa j}{2}+\frac{\delta j}{4}} 2^{-j+\delta k} \Vert u \Vert_X \\
&\lesssim t^{-\frac{13}{12}} 2^{-j+\delta k} \Vert u \Vert_X \\
\eqref{estdispsupLpetittauint-3} &\lesssim t^{-1} 2^{-2k-2j} \Vert \Psi_{(l_0, 0, 0)}^{int} \psi_{j, k}^{\widehat{\mathcal{L}}} \Vert_{L^2} \Vert m_{\widehat{\mathcal{L}}}(D) g_{b-corr}(t) \Vert_{L^2} \\
&\lesssim t^{-\frac{13}{12}+100\delta} 2^{\frac{l_0}{2}-k-\frac{j}{2}} \Vert u \Vert_X \\
&\lesssim t^{-\frac{13}{12}+100\delta} t^{\frac{\varrho}{2}-10\delta+\frac{\delta}{4}} 2^{\frac{3\varrho j}{2}+\frac{\delta j}{4}} 2^{-j+\delta k} \Vert u \Vert_X \\
&\lesssim t^{-\frac{13}{12}+91\delta} 2^{-j+\delta k} \Vert u \Vert_X 
\end{align*}
\eqref{estdispsupLpetittauint-1} is similar to \eqref{estdispsupLpetittauint-2}. 

Here we have $2^j \lesssim t^{-\frac{5}{24}} \lesssim t^{-\frac{1}{6}}$, so again $1 \lesssim t^{-\frac{1}{12}} 2^{-\frac{j}{2}}$ and we get the other estimate of the Lemma automatically. 

Assume now $t^{-\frac{5}{24}} \lesssim 2^j$. If $2^k \lesssim \mathfrak{t}^{-\frac{1}{2}}$, then we estimate
\begin{align*}
\eqref{estdispsupLpetittauint} &\lesssim \Vert \Psi_{(l_0, 0, 0)}^{int} \psi_{j, k}^{\widehat{\mathcal{L}}} \Vert_{L^{\frac{4}{3}}_{b, c} L^{\frac{1}{1-\kappa}}_a} \Vert \psi_j \widehat{f}(t) \Vert_{L^4_{b, c} L^{\frac{1}{1-\kappa}}_a} \\
&\lesssim 2^{l_0-\kappa l_0+\frac{3k}{2}+\frac{5j}{2}-\kappa j} \Vert \psi_j \widehat{f}(t) \Vert_{L^2_{b, c} H^1_a}^{\frac{1}{2}} \Vert \psi_j \widehat{f}(t) \Vert_{\dot{H}^1_{b, c} L^2_a}^{\frac{1}{2}} \\
&\lesssim t^{-\frac{7}{6}} t^{-\frac{1}{12}+\varrho+\frac{\kappa}{2}+101\delta+\frac{\delta}{2}} 2^{-\frac{j}{2}+3\varrho j + \frac{\kappa j}{2}+\frac{3\delta j}{2}} 2^{-j+\delta k} \langle 2^j \rangle^{-\frac{1}{2}} \Vert u \Vert_X
\end{align*}
In particular, 
\begin{align*}
\eqref{estdispsupLpetittauint} &\lesssim t^{-\frac{7}{6}} 2^{-\frac{3j}{2}+\delta k} \langle 2^j \rangle^{-\frac{1}{3}} \Vert u \Vert_X
\end{align*}
and 
\begin{align*}
\eqref{estdispsupLpetittauint} 
&\lesssim t^{-\frac{7}{6}} t^{\frac{1}{48}+\varrho+\frac{\kappa}{2}+101\delta+\frac{\delta}{2}} 2^{-j+\delta k} \langle 2^j \rangle^{-\frac{1}{2}} \Vert u \Vert_X \\
&\lesssim t^{-\frac{13}{12}} 2^{-j+\delta k} \langle 2^j \rangle^{-\frac{1}{2}} \Vert u \Vert_X 
\end{align*}
as wanted. On the other hand, if $2^k \gtrsim \mathfrak{t}^{-\frac{1}{2}}$, then we apply the same integration by parts as before. We then estimate: 
\begin{align*}
\eqref{estdispsupLpetittauint-1} &\lesssim \mathfrak{t}^{-1} 2^{-2k} \eqref{estdispsupLpetittauint} 
\end{align*}
so that \eqref{estdispsupLpetittauint-1} satisfies the same estimates as \eqref{estdispsupLpetittauint}; then we group
\begin{align*}
\eqref{estdispsupLpetittauint-2} + \eqref{estdispsupLpetittauint-3} &\lesssim t^{-1} 2^{-k-2j} \Vert \Psi_{(l_0, 0, 0)} \psi_{j, k}^{\widehat{\mathcal{L}}} \Vert_{L^2} \Vert \psi_j(D) X_{b-corr} f(t) \Vert_{L^2} \\
&\lesssim t^{-1} 2^{\frac{l_0}{2}-j} \langle 2^j \rangle^{-\frac{1}{2}} \Vert u \Vert_X \\
&\lesssim t^{-\frac{5}{4}+\varrho+\frac{\delta}{2}} 2^{-\frac{3j}{4}+\frac{3\delta j}{2}} 2^{-j+\delta k} \langle 2^j \rangle^{-\frac{1}{2}} \Vert u \Vert_X \\
&\lesssim t^{-\frac{7}{6}} t^{-\frac{1}{32}+\varrho+\frac{\delta}{2}} 2^{-\frac{j}{2}} 2^{-j+\delta k} \langle 2^j \rangle^{-\frac{1}{3}} \Vert u \Vert_X 
\end{align*}
Therefore, 
\begin{align*}
\eqref{estdispsupLpetittauint-2} + \eqref{estdispsupLpetittauint-3} &\lesssim t^{-\frac{7}{6}} 2^{-\frac{3j}{2}+\delta k} \langle 2^j \rangle^{-\frac{1}{3}} \Vert u \Vert_X 
\end{align*}
but also 
\begin{align*}
\eqref{estdispsupLpetittauint-2} + \eqref{estdispsupLpetittauint-3} &\lesssim t^{-\frac{13}{12}} t^{-\frac{1}{32}+\frac{1}{48}+\varrho+\frac{\delta}{2}} 2^{-j+\delta k} \langle 2^j \rangle^{-\frac{1}{3}} \Vert u \Vert_X \\
&\lesssim t^{-\frac{13}{12}} 2^{-j+\delta k} \langle 2^j \rangle^{-\frac{1}{3}} \Vert u \Vert_X 
\end{align*}
which concludes. 

\paragraph{2.} Let us now assume that $\tau \simeq 2^k$. 

In this case, we decompose 
\begin{subequations}
\begin{align}
\eqref{estdispsupL-termeinit} &= \sum_{l_a = l_0}^{-200} 1_{l_a + \frac{j_{+}}{2} \geq k} \int e^{i t \Phi} \Psi_{(l_a, 2k+100, 2k+100)}^{a-int} \psi_{j, k}^{\widehat{\mathcal{L}}}(\overline{\xi}) \widehat{f}(t, \overline{\xi}) ~ d\overline{\xi} \label{estdispsupLmoyentauaext} \\
&\quad + \sum_{l_a = l_0}^{-200} \sum_{l_b = k+l_a}^{2k+100} 1_{l_a + \frac{j_{+}}{2} < k} 1_{(1-2\delta) l_b \leq l_a} \int e^{i t \Phi} \Psi_{(l_a, l_b, l_b)}^{a-b} \psi_{j, k}^{\widehat{\mathcal{L}}}(\overline{\xi}) \widehat{f}(t, \overline{\xi}) ~ d\overline{\xi} \label{estdispsupLmoyentauab} \\
&\quad + \sum_{l_a = l_0}^{-200} \sum_{l_b = k+l_a}^{2k+100} 1_{l_a + \frac{j_{+}}{2} < k} 1_{l_a < (1-2\delta) l_b} \int e^{i t \Phi} \Psi_{(l_a, l_b, l_b)}^{a-b} \psi_{j, k}^{\widehat{\mathcal{L}}}(\overline{\xi}) \widehat{f}(t, \overline{\xi}) ~ d\overline{\xi} \label{estdispsupLmoyentauba} \\
&\quad + \sum_{l_a = l_0}^{-200} \sum_{l_b = k+l_a}^{2k+100} 1_{l_a + \frac{j_{+}}{2} < k} 1_{(1-2\delta) l_b \leq l_a} \int e^{i t \Phi} \Psi_{(l_a, l_b, l_b)}^{a-c} \psi_{j, k}^{\widehat{\mathcal{L}}}(\overline{\xi}) \widehat{f}(t, \overline{\xi}) ~ d\overline{\xi} \label{estdispsupLmoyentauac} \\
&\quad + \sum_{l_a = l_0}^{-200} \sum_{l_b = k+l_a}^{2k+100} 1_{l_a + \frac{j_{+}}{2} < k} 1_{l_a < (1-2\delta) l_b} \int e^{i t \Phi} \Psi_{(l_a, l_b, l_b)}^{a-c} \psi_{j, k}^{\widehat{\mathcal{L}}}(\overline{\xi}) \widehat{f}(t, \overline{\xi}) ~ d\overline{\xi} \label{estdispsupLmoyentauca} \\
&\quad + \sum_{l_a = l_0}^{-200} 1_{l_a + \frac{j_{+}}{2} < k} \int e^{i t \Phi} \Psi_{(l_a, k+l_a, k+l_a)}^{a-int} \psi_{j, k}^{\widehat{\mathcal{L}}}(\overline{\xi}) \widehat{f}(t, \overline{\xi}) ~ d\overline{\xi} \label{estdispsupLmoyentauaint} \\
&\quad + \sum_{l_b = F(l_0, k, j)}^{2k+100} \int e^{i t \Phi} \Psi_{(l_0, l_b, l_b)}^b \psi_{j, k}^{\widehat{\mathcal{L}}}(\overline{\xi}) \widehat{f}(t, \overline{\xi}) ~ d\overline{\xi} \label{estdispsupLmoyentaub} \\
&\quad + \sum_{l_b = F(l_0, k, j)}^{2k+100} \int e^{i t \Phi} \Psi_{(l_0, l_b, l_b)}^c \psi_{j, k}^{\widehat{\mathcal{L}}}(\overline{\xi}) \widehat{f}(t, \overline{\xi}) ~ d\overline{\xi} \label{estdispsupLmoyentauc} \\
&\quad + \int e^{i t \Phi} \Psi_{(l_0, F(l_0, k, j), F(l_0, k, j))}^{int} \psi_{j, k}^{\widehat{\mathcal{L}}}(\overline{\xi}) \widehat{f}(t, \overline{\xi}) ~ d\overline{\xi} \label{estdispsupLmoyentauint} 
\end{align}
\end{subequations}
We stopped the sum in $l_b$ at $l_b \leq 2k+100$ because we always have
\begin{align*}
m_b(\overline{\xi}) \widehat{X}_b(\overline{\xi}) \cdot \nabla_{\overline{\xi}} \Phi \lesssim 2^{2k+2j}, \quad m_b(\overline{\xi}) \widehat{X}_c(\overline{\xi}) \cdot \nabla_{\overline{\xi}} \Phi \lesssim 2^{2k+2j} 
\end{align*}
Above, as before, $l_0$ is such that $2^{l_0} \simeq \mathfrak{t}^{-\frac{1}{2}+\varrho}$ for a small enough $\varrho > 0$ with respect to $\delta$. $F$ is a function whose form we will specify later. 

First we estimate the volume of $\Psi_{(l_a, l_b, l_b)}$. To that end, we apply the change of variables $\xi' = \xi+2\eta$. Then: 
\begin{align*}
\widehat{X}_a(\overline{\xi}) \cdot \nabla_{\overline{\xi}} \Phi &= \frac{\xi_0}{|\overline{\xi}|} \left( 3 |\overline{\xi}|^2 - \frac{x}{t} + \frac{|y|^2}{t^2 |\overline{\xi}|^2 m_0} \right) - \frac{\xi' \cdot y}{t |\overline{\xi}|} \\
m_b(\overline{\xi}) \widehat{X}_b(\overline{\xi}) \cdot \nabla_{\overline{\xi}} \Phi &\simeq |\xi'-\eta|^2 - |\eta|^2 = |\xi'|^2 - 2 \xi' \cdot \eta \\
m_c(\overline{\xi}) \widehat{X}_c(\overline{\xi}) \cdot \nabla_{\overline{\xi}} \Phi &\simeq \xi' \cdot J \eta
\end{align*}
Yet, we have
\begin{align*}
\eta &= \frac{\xi_0 y}{2 t |\overline{\xi}|^2 m_0} \simeq \tau 2^j \simeq 2^{j+k} 
\end{align*}
(where $m_0 = 1 - \frac{x}{t |\overline{\xi}|^2} \simeq 1$ here). In particular, $|\eta| \simeq |\xi|$, $|\eta| \gtrsim |\xi'|$. If $l_b \geq 2k-100$, we have as before
\begin{align*}
\Vert \Psi_{(l_a, l_b, l_b)} \psi_{j, k}^{\widehat{\mathcal{L}}} \Vert_{L^p} &\lesssim 2^{\frac{l_a+2k}{p}+\frac{3j}{p}}
\end{align*}
However, if $l_b \leq 2k-100$, then $\Psi_{(l_a, l_b, l_b)} \psi_{j, k}^{\widehat{\mathcal{L}}}$ localises $\xi'$ close to $0$ with precision $2^{l_b-k+j}$, so that 
\begin{align*}
\Vert \Psi_{(l_a, l_b, l_b)} \psi_{j, k}^{\widehat{\mathcal{L}}} \Vert_{L^p} &\lesssim 2^{\frac{l_a+2l_b-2k}{p}+\frac{3j}{p}}
\end{align*}
Now, using anisotropic coordinates already introduced, we may also estimate the volume in $\xi_a$ while fixing $\xi_b, \xi_c$: 
\begin{align*}
\Vert \Psi_{(l_a, l_b, l_b)} \psi_{j, k}^{\widehat{\mathcal{L}}} \Vert_{L^{\infty}_{b, c} L^1_a} &\lesssim 2^{l_a+j} 
\end{align*}
Therefore, by interpolation, for any $1 \leq p_a \leq p_b \leq \infty$ we have
\begin{align*}
\Vert \Psi_{(l_a, l_b, l_b)} \psi_{j, k}^{\widehat{\mathcal{L}}} \Vert_{L^{p_b}_{b, c} L^{p_a}_a} &\lesssim 2^{\frac{l_a}{p_a} + \frac{\min(2k, 2l_b-2k)}{p_b} + \frac{j}{p_a}+\frac{2j}{p_b}} 
\end{align*}

For the derivatives of symbols, we compute 
\begin{align*}
\widehat{X}_a(\overline{\xi}) \cdot \nabla_{\overline{\xi}} \left[ \widehat{X}_a(\overline{\xi}) \cdot \nabla_{\overline{\xi}} \Phi \right] &\simeq 2^j \\
\widehat{X}_a(\overline{\xi}) \cdot \nabla_{\overline{\xi}} \left[ m_b(\overline{\xi}) \widehat{X}_b(\overline{\xi}) \cdot \nabla_{\overline{\xi}} \Phi \right] 
&= 2 m_b(\overline{\xi}) |\xi| \lesssim 2^{2k+j} \\
\widehat{X}_a(\overline{\xi}) \cdot \nabla_{\overline{\xi}} \left[ m_c(\overline{\xi}) \widehat{X}_c(\overline{\xi}) \cdot \nabla_{\overline{\xi}} \Phi \right] &= 0 \\
m_b(\overline{\xi}) \widehat{X}_b(\overline{\xi}) \cdot \nabla_{\overline{\xi}} \left[ \widehat{X}_a(\overline{\xi}) \cdot \nabla_{\overline{\xi}} \Phi \right] &\lesssim 2^{2k+j} \\
m_b(\overline{\xi}) \widehat{X}_b(\overline{\xi}) \cdot \nabla_{\overline{\xi}} \left[ m_b(\overline{\xi}) \widehat{X}_b(\overline{\xi}) \cdot \nabla_{\overline{\xi}} \Phi \right]
&\lesssim 2^{2k+j} \\
m_b(\overline{\xi}) \widehat{X}_b(\overline{\xi}) \cdot \nabla_{\overline{\xi}} \left[ m_c(\overline{\xi}) \widehat{X}_c(\overline{\xi}) \cdot \nabla_{\overline{\xi}} \Phi \right]
&= m_b(\overline{\xi}) \frac{\xi_0}{|\overline{\xi}|} \frac{\xi \cdot J y}{t |\overline{\xi}| |\xi|} \lesssim 2^{l_b+j} \\
m_c(\overline{\xi}) \widehat{X}_c(\overline{\xi}) \cdot \nabla_{\overline{\xi}} \left[ \widehat{X}_a(\overline{\xi}) \cdot \nabla_{\overline{\xi}} \Phi \right] &= \frac{\xi \cdot Jy}{t |\overline{\xi}|^2} \lesssim 2^{l_b+j} \\
m_c(\overline{\xi}) \widehat{X}_c(\overline{\xi}) \cdot \nabla_{\overline{\xi}} \left[ m_b(\overline{\xi}) \widehat{X}_b(\overline{\xi}) \cdot \nabla_{\overline{\xi}} \Phi \right] &\lesssim 2^{2k+j} \\
m_c(\overline{\xi}) \widehat{X}_c(\overline{\xi}) \cdot \nabla_{\overline{\xi}} \left[ m_c(\overline{\xi}) \widehat{X}_c(\overline{\xi}) \cdot \nabla_{\overline{\xi}} \Phi \right] &\lesssim 2^{2k+j} 
\end{align*}
We then have
\begin{align*}
\widehat{X}_a(\overline{\xi}) \cdot \nabla_{\overline{\xi}} \psi_{j, k}^{\widehat{\mathcal{L}}} &\lesssim 2^{-j} \\
m_b(\overline{\xi}) \nabla_{\overline{\xi}} \psi_{j, k}^{\widehat{\mathcal{L}}} &\lesssim 2^{-j} 
\end{align*}
As before, if $l_b \geq l_a$, we can always correct $\widehat{X}_b, \widehat{X}_c$ in order to get
\begin{align*}
\widehat{X}_{b-corr}(\overline{\xi}) \cdot \nabla_{\overline{\xi}} \left[ \widehat{X}_a(\overline{\xi}) \cdot \nabla_{\overline{\xi}} \Phi \right] &= 0 \\
\widehat{X}_{c-corr}(\overline{\xi}) \cdot \nabla_{\overline{\xi}} \left[ \widehat{X}_a(\overline{\xi}) \cdot \nabla_{\overline{\xi}} \Phi \right] &= 0 
\end{align*}
while keeping similar esimates on $\widehat{X}_{\beta-corr} \cdot \nabla \Phi$, $\beta = b, c$. 

On \eqref{estdispsupLmoyentauaext}, we apply at most $n$ (large with respect to $\varrho^{-1}$) integrations by parts along $\widehat{X}_a$: 
\begin{subequations}
\begin{align}
\eqref{estdispsupLmoyentauaext} &= \sum_{l_a = l_0}^{-200} 1_{l_a + \frac{j_{+}}{2} \geq k} t^{-n} 2^{-2nl_a-3nj} \int e^{i t \Phi} \Psi_{(l_a, 2k+100, 2k+100)}^{a-int} \psi_{j, k}^{\widehat{\mathcal{L}}}(\overline{\xi}) \widehat{f}(t, \overline{\xi}) ~ d\overline{\xi} \label{estdispsupLmoyentauaext-1} \\
&\quad + \sum_{i = 1}^n \sum_{l_a = l_0}^{-200} 1_{l_a + \frac{j_{+}}{2} \geq k} t^{-i-1} 2^{-2il_a-l_a-3ij-2j} \int e^{i t \Phi} \Psi_{(l_a, 2k+100, 2k+100)}^{a-int} \psi_{j, k}^{\widehat{\mathcal{L}}}(\overline{\xi}) \widehat{h}_a(t, \overline{\xi}) ~ d\overline{\xi} \label{estdispsupLmoyentauaext-2} \\
&\quad \begin{aligned}
+ \sum_{i = 1}^n \sum_{l_a = l_0}^{-200} 1_{l_a + \frac{j_{+}}{2} \geq k} &t^{-i-1} 2^{-2il_a-3ij-j} \\
&\int e^{i t \Phi} \Psi_{(l_a, 2k+100, 2k+100)}^{a-int} \psi_{j, k}^{\widehat{\mathcal{L}}}(\overline{\xi}) \widehat{X}_a(\overline{\xi}) \cdot \nabla_{\overline{\xi}} \widehat{h}_a(t, \overline{\xi}) ~ d\overline{\xi} 
\end{aligned} \label{estdispsupLmoyentauaext-3} \\
&\quad + \sum_{i = 1}^n \sum_{l_a = l_0}^{-200} 1_{l_a + \frac{j_{+}}{2} \geq k} t^{-i} 2^{-2il_a+l_a-3ij+j} \int e^{i t \Phi} \Psi_{(l_a, 2k+100, 2k+100)}^{a-int} \psi_{j, k}^{\widehat{\mathcal{L}}}(\overline{\xi}) \widehat{g}_a(t, \overline{\xi}) ~ d\overline{\xi} \label{estdispsupLmoyentauaext-4} 
\end{align}
\end{subequations} 
We then estimate: 
\begin{align*}
\eqref{estdispsupLmoyentauaext-1} &\lesssim \sum_{l_a = l_0}^{-200} 1_{l_a + \frac{j_{+}}{2} \geq k} t^{-\frac{7}{6}} \mathfrak{t}^{-n+\frac{7}{6}} 2^{-2nl_a-\frac{5j}{2}} \Vert \Psi_{(l_a, 2k+100, 2k+100)}^{a-int} \psi_{j, k}^{\widehat{\mathcal{L}}} \Vert_{L^2} \langle 2^j \rangle^{-1} \Vert \langle \nabla \rangle |\nabla|^{-1} f(t) \Vert_{L^2} \\
&\lesssim \sum_{l_a = l_0}^{-200} 1_{l_a + \frac{j_{+}}{2} \geq k} t^{-\frac{7}{6}} \mathfrak{t}^{-n+\frac{7}{6}} 2^{-2nl_a+\frac{l_a}{2}+k-j} \langle 2^j \rangle^{-1} \Vert u \Vert_X \\
&\lesssim t^{-\frac{7}{6}} \mathfrak{t}^{\frac{7}{6}-2n\varrho} 2^{-j+k} \langle 2^j \rangle^{-1} \Vert u \Vert_X \\
&\lesssim t^{-\frac{7}{6}} 2^{-j+k} \langle 2^j \rangle^{-1} \Vert u \Vert_X \\
\eqref{estdispsupLmoyentauaext-2} &\lesssim \sum_{i = 1}^n \sum_{l_a = l_0}^{-200} 1_{l_a + \frac{j_{+}}{2} \geq k} t^{-\frac{7}{6}} \mathfrak{t}^{\frac{1}{6}-i} 2^{-2il_a-l_a-\frac{5j}{2}} \Vert \Psi_{(l_a, 2k+100, 2k+100)}^{a-int} \psi_{j, k}^{\widehat{\mathcal{L}}} \Vert_{L^2_{b, c} L^{\frac{1}{1-\kappa}}_a} \Vert \psi_j \widehat{h}_a(t) \Vert_{L^2_{b, c} L^{\frac{1}{\kappa}}_a} \\
&\lesssim  \sum_{i = 1}^n \sum_{l_a = l_0}^{-200} 1_{l_a + \frac{j_{+}}{2} \geq k} t^{-\frac{7}{6}} \mathfrak{t}^{\frac{1}{6}-i} 2^{-2il_a-\kappa l_a+k-\frac{j}{2}-\kappa j} \Vert \psi_j \widehat{h}_a(t) \Vert_{L^2_{b, c} H^1_a}^{\frac{1}{2}} \Vert \psi_j \widehat{h}_a(t) \Vert_{L^2}^{\frac{1}{2}} \\
&\lesssim t^{-\frac{7}{6}+\frac{\kappa}{3}} \mathfrak{t}^{-\frac{1}{3}+\frac{\kappa}{6}+\frac{\delta}{2}} 2^{-j+\delta k} \Vert u \Vert_X \\
&\lesssim t^{-\frac{7}{6}+\frac{\kappa}{3}} 2^{-j+\delta k} \langle 2^j \rangle^{-\frac{3}{4}} \Vert u \Vert_X \\
\eqref{estdispsupLmoyentauaext-3} &\lesssim \sum_{i = 1}^n \sum_{l_a = l_0}^{-200} 1_{l_a + \frac{j_{+}}{2} \geq k} t^{-\frac{7}{6}} \mathfrak{t}^{\frac{1}{6}-i} 2^{-2il_a-\frac{5j}{2}} \Vert \Psi_{(l_a, 2k+100, 2k+100)} \psi_{j, k}^{\widehat{\mathcal{L}}} \Vert_{L^2} \Vert \nabla X_a h_a(t) \Vert_{L^2} \\
&\lesssim \sum_{i = 1}^n \sum_{l_a = l_0}^{-200} 1_{l_a + \frac{j_{+}}{2} \geq k} t^{-\frac{7}{6}} \mathfrak{t}^{\frac{1}{6}-i} 2^{-2il_a+\frac{l_a}{2}+k-j} \Vert u \Vert_X \\
&\lesssim t^{-\frac{7}{6}} \mathfrak{t}^{-\frac{7}{12}+\frac{\delta}{2}} 2^{-j+\delta k} \langle 2^j \rangle^{\frac{1}{2}} \Vert u \Vert_X \\
&\lesssim t^{-\frac{7}{6}} 2^{-j+\delta k} \langle 2^j \rangle^{-1} \Vert u \Vert_X \\
\eqref{estdispsupLmoyentauaext-4} &\lesssim \sum_{i = 1}^n \sum_{l_a = l_0}^{-200} 1_{l_a + \frac{j_{+}}{2} \geq k} t^{-1} \mathfrak{t}^{1-i} 2^{-2il_a+l_a-\frac{5j}{2}} \\
&\pushright{\Vert \Psi_{(l_a, 2k+100, 2k+100)} \psi_{j, k}^{\widehat{\mathcal{L}}} \Vert_{L^2} \langle 2^j \rangle^{-1} \Vert m_{\widehat{\mathcal{L}}}(D) \langle \nabla \rangle |\nabla|^{\frac{1}{2}} g_a(t) \Vert_{L^2}} \\
&\lesssim \sum_{i = 1}^n \sum_{l_a = l_0}^{-200} 1_{l_a + \frac{j_{+}}{2} \geq k} t^{-\frac{7}{6}+100\delta} \mathfrak{t}^{1-i} 2^{-2il_a+\frac{3l_a}{2}+k-j} \langle 2^j \rangle^{-1} \Vert u \Vert_X \\
&\lesssim t^{-\frac{7}{6}+100\delta} \mathfrak{t}^{-\frac{1}{4}+\frac{\delta}{2}} 2^{-j+\delta k} \langle 2^j \rangle^{-\frac{1}{2}} \Vert u \Vert_X \\
&\lesssim t^{-\frac{7}{6}+100\delta} 2^{-j+\delta k} \langle 2^j \rangle^{-1} \Vert u \Vert_X 
\end{align*}
for $n \varrho$ large enough with respect to $1$, $\delta, \varrho, \kappa$ small enough. 

For \eqref{estdispsupLmoyentauab}, we apply first at most $n$ integrations by parts along $\widehat{X}_a$, then one integration by parts along $\widehat{X}_b$: 
\begin{subequations}
\begin{align}
&\eqref{estdispsupLmoyentauab} = \sum_{l_a = l_0}^{-200} \sum_{l_b = l_a+k}^{2k+100} 1_{(1-2\delta) l_b \leq l_a < k - \frac{j_{+}}{2}} t^{-n} 2^{-2nl_a-3nj} \int e^{i t \Phi} \Psi_{(l_a, l_b, l_b)}^{a-b} \psi_{j, k}^{\widehat{\mathcal{L}}}(\overline{\xi}) \widehat{f}(t, \overline{\xi}) ~ d\overline{\xi} \label{estdispsupLmoyentauab-1} \\
&\quad + \sum_{i = 1}^n \sum_{l_a = l_0}^{-200} \sum_{l_b = l_a+k}^{2k+100} 1_{(1-2\delta) l_b \leq l_a < k - \frac{j_{+}}{2}} t^{-i-1} 2^{-2il_a+l_a-2l_b+2k-3ij-2j} \int e^{i t \Phi} \Psi_{(l_a, l_b, l_b)}^{a-b} \psi_{j, k}^{\widehat{\mathcal{L}}}(\overline{\xi}) \widehat{h}_a(t, \overline{\xi}) ~ d\overline{\xi} \label{estdispsupLmoyentauab-2} \\
&\quad \begin{aligned}
+ \sum_{i = 1}^n \sum_{l_a = l_0}^{-200} \sum_{l_b = l_a+k}^{2k+100} 1_{(1-2\delta) l_b \leq l_a < k - \frac{j_{+}}{2}} &t^{-i-1} 2^{-2il_a+l_a-l_b-3ij-j} \\
&\int e^{i t \Phi} \Psi_{(l_a, l_b, l_b)}^{a-b} \psi_{j, k}^{\widehat{\mathcal{L}}}(\overline{\xi}) m_b(\overline{\xi}) \widehat{X}_b(\overline{\xi}) \cdot \nabla_{\overline{\xi}} \widehat{h}_a(t, \overline{\xi}) ~ d\overline{\xi} 
\end{aligned} \label{estdispsupLmoyentauab-3} \\
&\quad + \sum_{i = 1}^n \sum_{l_a = l_0}^{-200} \sum_{l_b = l_a+k}^{2k+100} 1_{(1-2\delta) l_b \leq l_a < k - \frac{j_{+}}{2}} t^{-i} 2^{-2il_a+l_a-3ij+j} \int e^{i t \Phi} \Psi_{(l_a, l_b, l_b)}^{a-b} \psi_{j, k}^{\widehat{\mathcal{L}}}(\overline{\xi}) \widehat{g}_a(t, \overline{\xi}) ~ d\overline{\xi} \label{estdispsupLmoyentauab-4} 
\end{align}
\end{subequations} 
We then estimate: 
\begin{align*}
\eqref{estdispsupLmoyentauab-1} &\lesssim \sum_{l_a = l_0}^{-200} \sum_{l_b = l_a+k}^{2k+100} 1_{(1-2\delta) l_b \leq l_a < k - \frac{j_{+}}{2}} t^{-\frac{7}{6}} \mathfrak{t}^{\frac{7}{6}-n} 2^{-2nl_a-\frac{5j}{2}} \Vert \Psi_{(l_a, l_b, l_b)}^{a-b} \psi_{j, k}^{\widehat{\mathcal{L}}} \Vert_{L^2} \langle 2^j \rangle^{-1} \Vert \langle \nabla \rangle |\nabla|^{-1} f(t) \Vert_{L^2} \\
&\lesssim \sum_{l_a = l_0}^{-200} \sum_{l_b = l_a+k}^{2k+100} 1_{(1-2\delta) l_b \leq l_a < k - \frac{j_{+}}{2}} t^{-\frac{7}{6}} \mathfrak{t}^{\frac{7}{6}-n} 2^{-2nl_a+\frac{l_a}{2}+l_b-k-j} \langle 2^j \rangle^{-1} \Vert u \Vert_X \\
&\lesssim t^{-\frac{7}{6}} \mathfrak{t}^{\frac{7}{6}-n \varrho} 2^{-j+k} \langle 2^j \rangle^{-1} \Vert u \Vert_X \\
&\lesssim t^{-\frac{7}{6}} 2^{-j+k} \langle 2^j \rangle^{-1} \Vert u \Vert_X \\
\eqref{estdispsupLmoyentauab-2} &\lesssim \sum_{i = 1}^n \sum_{l_a = l_0}^{-200} \sum_{l_b = l_a+k}^{2k+100} 1_{(1-2\delta) l_b \leq l_a < k - \frac{j_{+}}{2}} t^{-\frac{7}{6}} \mathfrak{t}^{\frac{1}{6}-i} 2^{-2il_a+l_a-2l_b+2k-\frac{7j}{2}} \\
&\pushright{\Vert \Psi_{(l_a, l_b, l_b)}^{a-b} \Vert_{L^{\frac{4}{3}}_{b, c} L^{\frac{1}{1-\kappa}}_a} \Vert \overline{\xi} \psi_{j, k}^{\widehat{\mathcal{L}}} \widehat{h}_a(t) \Vert_{L^4_{b, c} L^{\frac{1}{\kappa}}_a}} \\
&\lesssim \sum_{i = 1}^n \sum_{l_a = l_0}^{-200} \sum_{l_b = l_a+k}^{2k+100} 1_{(1-2\delta) l_b \leq l_a < k - \frac{j_{+}}{2}} t^{-\frac{7}{6}} \mathfrak{t}^{\frac{1}{6}-i} 2^{-2il_a+2l_a-\kappa l_a-\frac{l_b}{2}+\frac{k}{2}-j-\kappa j} \\
&\pushright{\Vert \overline{\xi} \psi_{j, k}^{\widehat{\mathcal{L}}} \widehat{h}_a(t) \Vert_{L^2_{b, c} H^1_a}^{\frac{1}{2}} \Vert \overline{\xi} \psi_{j, k}^{\widehat{\mathcal{L}}} \widehat{h}_a(t) \Vert_{\dot{H}^1_{b, c} L^2_a}^{\frac{1}{2}}} \\
&\lesssim t^{-\frac{7}{6}+\frac{\kappa}{3}} \mathfrak{t}^{-\frac{1}{3}+\frac{\kappa}{6}+\frac{\delta}{2}} 2^{-j+\delta k} \langle 2^j \rangle^{\frac{1}{4}} \Vert u \Vert_X \\
&\lesssim t^{-\frac{7}{6}+\frac{\kappa}{3}} 2^{-j+\delta k} \langle 2^j \rangle^{-\frac{3}{5}} \Vert u \Vert_X \\
\eqref{estdispsupLmoyentauab-3} &\lesssim \sum_{i = 1}^n \sum_{l_a = l_0}^{-200} \sum_{l_b = l_a+k}^{2k+100} 1_{(1-2\delta) l_b \leq l_a < k - \frac{j_{+}}{2}} t^{-\frac{7}{6}} \mathfrak{t}^{\frac{1}{6}-i} 2^{-2il_a+l_a-l_b-\frac{5j}{2}} \\
&\pushright{\Vert \Psi_{(l_a, l_b, l_b)}^{a-b} \psi_{j, k}^{\widehat{\mathcal{L}}} \Vert_{L^2} \Vert \nabla m_b(D) X_b h_a(t) \Vert_{L^2}} \\
&\lesssim \sum_{i = 1}^n \sum_{l_a = l_0}^{-200} \sum_{l_b = l_a+k}^{2k+100} 1_{(1-2\delta) l_b \leq l_a < k - \frac{j_{+}}{2}} t^{-\frac{7}{6}} \mathfrak{t}^{\frac{1}{6}-i} 2^{-2il_a+\frac{3l_a}{2}-k-j} \Vert u \Vert_X \\
&\lesssim t^{-\frac{7}{6}} \mathfrak{t}^{-\frac{1}{12}+\frac{\delta}{2}+\frac{\kappa}{2}} 2^{-j+\delta k} \langle 2^j \rangle^{-\frac{1}{2}} \Vert u \Vert_X \\
&\lesssim t^{-\frac{7}{6}} 2^{-j+\delta k} \langle 2^j \rangle^{-\frac{3}{5}} \Vert u \Vert_X \\
\eqref{estdispsupLmoyentauab-4} &\lesssim \sum_{i = 1}^n \sum_{l_a = l_0}^{-200} \sum_{l_b = l_a+k}^{2k+100} 1_{(1-2\delta) l_b \leq l_a < k - \frac{j_{+}}{2}} t^{-1} \mathfrak{t}^{1-i} 2^{-2il_a+l_a-\frac{5j}{2}} \\
&\pushright{\Vert \Psi_{(l_a, l_b, l_b)}^{a-b} \psi_{j, k}^{\widehat{\mathcal{L}}} \Vert_{L^2} \langle 2^j \rangle^{-1} \Vert m_{\widehat{\mathcal{L}}}(D) \langle \nabla \rangle |\nabla|^{\frac{1}{2}} g_a(t) \Vert_{L^2}} \\
&\lesssim \sum_{i = 1}^n \sum_{l_a = l_0}^{-200} \sum_{l_b = l_a+k}^{2k+100} 1_{(1-2\delta) l_b \leq l_a < k - \frac{j_{+}}{2}} t^{-\frac{7}{6}+100\delta} \mathfrak{t}^{1-i} 2^{-2il_a+\frac{3l_a}{2}-j+l_b-k} \langle 2^j \rangle^{-1} \Vert u \Vert_X \\
&\lesssim \sum_{i = 1}^n \sum_{l_a = l_0}^{-200} 1_{l_a < k - \frac{j_{+}}{2}} t^{-\frac{7}{6}+100\delta} \mathfrak{t}^{1-i} 2^{-2il_a+\frac{3l_a}{2}+\frac{(1-\delta)l_a}{2(1-2\delta)}} 2^{-j+\delta k} \langle 2^j \rangle^{-1} \Vert u \Vert_X \\
&\lesssim t^{-\frac{7}{6}+100\delta} 2^{-j+\delta k} \langle 2^j \rangle^{-1} \Vert u \Vert_X 
\end{align*}

For \eqref{estdispsupLmoyentauba}, we reverse the order of integrations by parts: 
\begin{subequations}
\begin{align}
&\eqref{estdispsupLmoyentauba} \notag \\
&= \sum_{l_a = l_0}^{-200} \sum_{l_b = l_a+k}^{2k+100} 1_{l_a + \frac{j_{+}}{2} < k} 1_{l_a < (1-2\delta) l_b} t^{-1-n} 2^{-2l_b+2k-2nl_a-3j-3nj} \int e^{i t \Phi} \Psi_{(l_a, l_b, l_b)}^{a-b} \psi_{j, k}^{\widehat{\mathcal{L}}}(\overline{\xi}) \widehat{f}(t, \overline{\xi}) ~ d\overline{\xi} \label{estdispsupLmoyentauba-1} \\
&\quad \begin{aligned}
+ \sum_{i = 1}^n \sum_{l_a = l_0}^{-200} \sum_{l_b = l_a+k}^{2k+100} 1_{l_a + \frac{j_{+}}{2} < k} 1_{l_a < (1-2\delta) l_b} &t^{-1-i} 2^{-2l_b+2k-2il_a+l_a-2j-3ij} \\
&\int e^{i t \Phi} \Psi_{(l_a, l_b, l_b)}^{a-b} \psi_{j, k}^{\widehat{\mathcal{L}}}(\overline{\xi}) \widehat{X}_a(\overline{\xi}) \cdot \nabla_{\overline{\xi}} \widehat{f}(t, \overline{\xi}) ~ d\overline{\xi} 
\end{aligned} \label{estdispsupLmoyentauba-2} \\
&\quad + \sum_{l_a = l_0}^{-200} \sum_{l_b = l_a+k}^{2k+100} 1_{l_a + \frac{j_{+}}{2} < k} 1_{l_a < (1-2\delta) l_b} t^{-1-n} 2^{-l_b-2nl_a-2j-3nj} \int e^{i t \Phi} \Psi_{(l_a, l_b, l_b)}^{a-b} \psi_{j, k}^{\widehat{\mathcal{L}}}(\overline{\xi}) \widehat{h}_{b-corr}(t, \overline{\xi}) ~ d\overline{\xi} \label{estdispsupLmoyentauba-3} \\
&\quad \begin{aligned}
+ \sum_{i = 1}^n \sum_{l_a = l_0}^{-200} \sum_{l_b = l_a+k}^{2k+100} 1_{l_a + \frac{j_{+}}{2} < k} &1_{l_a < (1-2\delta) l_b} t^{-1-i} 2^{-l_b-2il_a+l_a-j-3ij} \\
&\int e^{i t \Phi} \Psi_{(l_a, l_b, l_b)}^{a-b} \psi_{j, k}^{\widehat{\mathcal{L}}}(\overline{\xi}) \widehat{X}_a(\overline{\xi}) \cdot \nabla_{\overline{\xi}} \widehat{h}_{b-corr}(t, \overline{\xi}) ~ d\overline{\xi}
\end{aligned} \label{estdispsupLmoyentauba-4} \\
&\quad + \sum_{l_a = l_0}^{-200} \sum_{l_b = l_a+k}^{2k+100} 1_{l_a + \frac{j_{+}}{2} < k} 1_{l_a < (1-2\delta) l_b} t^{-1} 2^{-l_b-2j} \int e^{i t \Phi} \Psi_{(l_a, l_b, l_b)}^{a-b} \psi_{j, k}^{\widehat{\mathcal{L}}}(\overline{\xi}) \widehat{g}_{b-corr}(t, \overline{\xi}) ~ d\overline{\xi} \label{estdispsupLmoyentauba-5}  
\end{align}
\end{subequations}
We then estimate: 
\begin{align*}
\eqref{estdispsupLmoyentauba-1} &\lesssim \sum_{l_a = l_0}^{-200} \sum_{l_b = l_a+k}^{2k+100} 1_{l_a + \frac{j_{+}}{2} < k} 1_{l_a < (1-2\delta) l_b} t^{-\frac{7}{6}} \mathfrak{t}^{\frac{1}{6}-n} 2^{-2l_b+2k-2nl_a-\frac{5j}{2}} \\
&\pushright{\Vert \Psi_{(l_a, l_b, l_b)}^{a-b} \psi_{j, k}^{\widehat{\mathcal{L}}} \Vert_{L^2} \langle 2^j \rangle^{-1} \Vert \langle \nabla \rangle |\nabla|^{-1} f(t) \Vert_{L^2}} \\
&\lesssim t^{-\frac{7}{6}} \mathfrak{t}^{\frac{5}{12}-2n\varrho} 2^{-j+k} \langle 2^j \rangle^{-1} \Vert u \Vert_X \\
&\lesssim t^{-\frac{7}{6}} 2^{-j+k} \langle 2^j \rangle^{-1} \Vert u \Vert_X \\
\eqref{estdispsupLmoyentauba-2} &\lesssim \sum_{i = 1}^n \sum_{l_a = l_0}^{-200} \sum_{l_b = l_a+k}^{2k+100} 1_{l_a + \frac{j_{+}}{2} < k} 1_{l_a < (1-2\delta) l_b} t^{-\frac{7}{6}} \mathfrak{t}^{\frac{1}{6}-i} 2^{-2l_b+2k-2il_a+l_a-\frac{5j}{2}} \\
&\pushright{\Vert \Psi_{(l_a, l_b, l_b)}^{a-b} \psi_{j, k}^{\widehat{\mathcal{L}}} \Vert_{L^2} \langle 2^j \rangle^{-1} \Vert \langle \nabla \rangle X_a f(t) \Vert_{L^2}} \\
&\lesssim \sum_{i = 1}^n \sum_{l_a = l_0}^{-200} \sum_{l_b = l_a+k}^{2k+100} 1_{l_a + \frac{j_{+}}{2} < k} 1_{l_a < (1-2\delta) l_b} t^{-\frac{7}{6}} \mathfrak{t}^{\frac{1}{6}-i} 2^{-l_b+k-2il_a+\frac{3l_a}{2}-j} \langle 2^j \rangle^{-1} \Vert u \Vert_X \\
&\lesssim t^{-\frac{7}{6}} \mathfrak{t}^{-\frac{1}{12}+2\delta} 2^{-j+k} \langle 2^j \rangle^{-1} \Vert u \Vert_X \\
&\lesssim t^{-\frac{7}{6}} 2^{-j+k} \langle 2^j \rangle^{-1} \Vert u \Vert_X \\
\eqref{estdispsupLmoyentauba-3} &\lesssim \sum_{l_a = l_0}^{-200} \sum_{l_b = l_a+k}^{2k+100} 1_{l_a + \frac{j_{+}}{2} < k} 1_{l_a < (1-2\delta) l_b} t^{-\frac{7}{6}} \mathfrak{t}^{\frac{1}{6}-n} 2^{-l_b-2nl_a-\frac{5j}{2}} \\
&\pushright{\Vert \Psi_{(l_a, l_b, l_b)}^{a-b} \psi_{j, k}^{\widehat{\mathcal{L}}} \Vert_{L^2} \langle 2^j \rangle^{-1} \Vert \langle \nabla \rangle h_{b-corr}(t) \Vert_{L^2}} \\
&\lesssim \sum_{l_a = l_0}^{-200} \sum_{l_b = l_a+k}^{2k+100} 1_{l_a + \frac{j_{+}}{2} < k} 1_{l_a < (1-2\delta) l_b} t^{-\frac{7}{6}} \mathfrak{t}^{\frac{1}{6}-n} 2^{-2nl_a+\frac{l_a}{2}-k-j} \langle 2^j \rangle^{-1} \Vert u \Vert_X \\
&\lesssim t^{-\frac{7}{6}} \mathfrak{t}^{\frac{2}{3}-2n\varrho} 2^{-j+\delta k} \langle 2^j \rangle^{-1} \Vert u \Vert_X \\
&\lesssim t^{-\frac{7}{6}} 2^{-j+\delta k} \langle 2^j \rangle^{-1} \Vert u \Vert_X \\
\eqref{estdispsupLmoyentauba-4} &\lesssim \sum_{i = 1}^n \sum_{l_a = l_0}^{-200} \sum_{l_b = l_a+k}^{2k+100} 1_{l_a + \frac{j_{+}}{2} < k} 1_{l_a < (1-2\delta) l_b} t^{-\frac{7}{6}} \mathfrak{t}^{\frac{1}{6}-i} 2^{-l_b-2il_a+l_a-\frac{5j}{2}} \\
&\pushright{\Vert \Psi_{(l_a, l_b, l_b)}^{a-b} \psi_{j, k}^{\widehat{\mathcal{L}}} \Vert_{L^2} \Vert \nabla X_a h_{b-corr}(t) \Vert_{L^2}} \\
&\lesssim \sum_{i = 1}^n \sum_{l_a = l_0}^{-200} \sum_{l_b = l_a+k}^{2k+100} 1_{l_a + \frac{j_{+}}{2} < k} 1_{l_a < (1-2\delta) l_b} t^{-\frac{7}{6}} \mathfrak{t}^{\frac{1}{6}-i} 2^{-2il_a+\frac{3l_a}{2}-k-j} \Vert u \Vert_X \\
&\lesssim t^{-\frac{7}{6}} \mathfrak{t}^{-\frac{1}{3}+\frac{\kappa}{2}+\frac{\delta}{4}} 2^{-j+\delta k} \Vert u \Vert_X \\
&\lesssim t^{-\frac{7}{6}} 2^{-j+\delta k} \langle 2^j \rangle^{-\frac{3}{4}} \Vert u \Vert_X \\
\eqref{estdispsupLmoyentauba-5} &\lesssim \sum_{l_a = l_0}^{-200} \sum_{l_b = l_a+k}^{2k+100} 1_{l_a + \frac{j_{+}}{2} < k} 1_{l_a < (1-2\delta) l_b} t^{-1} 2^{-l_b-\frac{5j}{2}} \\
&\pushright{\Vert \Psi_{(l_a, l_b, l_b)}^{a-b} \psi_{j, k}^{\widehat{\mathcal{L}}} \Vert_{L^2} \langle 2^j \rangle^{-1} \Vert \langle \nabla \rangle |\nabla|^{\frac{1}{2}} g_{b-corr}(t) \Vert_{L^2}} \\
&\lesssim \sum_{l_a = l_0}^{-200} \sum_{l_b = l_a+k}^{2k+100} 1_{l_a + \frac{j_{+}}{2} < k} 1_{l_a < (1-2\delta) l_b} t^{-\frac{7}{6}+100\delta} 2^{-k+\frac{l_a}{2}-j} \langle 2^j \rangle^{-1} \Vert u \Vert_X \\
&\lesssim t^{-\frac{7}{6}+100 \delta} \mathfrak{t}^{\delta} 2^{-j+\delta k} \langle 2^j \rangle^{-1} \Vert u \Vert_X 
\end{align*}
All these estimates are enough, except the one for \eqref{estdispsupLmoyentauba-5} on which we lost a factor $t^{\delta}$. As before, we may easily absorb such a factor as soon as $2^{|j|} \gtrsim t^{1000\delta}$. On the other hand, in the case $2^{|j|} \lesssim t^{1000 \delta}$, as before, we can also restrict the sum to $(1+4\delta) l_b \leq l_a < (1-2\delta) l_b$, and in this place use $\widehat{X}_b$ instead of $\widehat{X}_{b-corr}$ in order to use the fine estimate on $g_b$ to improve the estimate of \eqref{estdispsupLmoyentauba-5}. 

We can estimate \eqref{estdispsupLmoyentauac} and \eqref{estdispsupLmoyentauca} the same way, replacing $\widehat{X}_b$ by $\widehat{X}_c$. Note that, when we use $\widehat{X}_c$, $g_c \equiv 0$ so we have nothing to estimate for \eqref{estdispsupLmoyentauba-5}. 

Then, on \eqref{estdispsupLmoyentauaint}, we apply again only integrations by parts along $\widehat{X}_a$: 
\begin{subequations}
\begin{align}
\eqref{estdispsupLmoyentauaint} &= \sum_{l_a = l_0}^{-200} 1_{l_a < k} t^{-n} 2^{-2nl_a-3nj} \int e^{i t \Phi} \Psi_{(l_a, l_a+k, l_a+k)}^{a-int} \psi_{j, k}^{\widehat{\mathcal{L}}}(\overline{\xi}) \widehat{f}(t, \overline{\xi}) ~ d\overline{\xi} \label{estdispsupLmoyentauaint-1} \\
&\quad + \sum_{i = 1}^n \sum_{l_a = l_0}^{-200} 1_{l_a < k} t^{-i-1} 2^{-2il_a-l_a-3ij-2j} \int e^{i t \Phi} \Psi_{(l_a, l_a+k, l_a+k)}^{a-int} \psi_{j, k}^{\widehat{\mathcal{L}}}(\overline{\xi}) \widehat{h}_a(t, \overline{\xi}) ~ d\overline{\xi} \label{estdispsupLmoyentauaint-2} \\
&\quad + \sum_{i = 1}^n \sum_{l_a = l_0}^{-200} 1_{l_a < k} t^{-i-1} 2^{-2il_a-3ij-j} \int e^{i t \Phi} \Psi_{(l_a, l_a+k, l_a+k)}^{a-int} \psi_{j, k}^{\widehat{\mathcal{L}}}(\overline{\xi}) \widehat{X}_a(\overline{\xi}) \cdot \nabla_{\overline{\xi}} \widehat{h}_a(t, \overline{\xi}) ~ d\overline{\xi} \label{estdispsupLmoyentauaint-3} \\
&\quad + \sum_{i = 1}^n \sum_{l_a = l_0}^{-200} 1_{l_a < k} t^{-i} 2^{-2il_a+l_a-3ij+j} \int e^{i t \Phi} \Psi_{(l_a, l_a+k, l_a+k)}^{a-int} \psi_{j, k}^{\widehat{\mathcal{L}}}(\overline{\xi}) \widehat{g}_a(t, \overline{\xi}) ~ d\overline{\xi} \label{estdispsupLmoyentauaint-4} 
\end{align}
\end{subequations} 
We then estimate: 
\begin{align*}
\eqref{estdispsupLmoyentauaint-1} &\lesssim \sum_{l_a = l_0}^{-200} 1_{l_a < k} t^{-\frac{7}{6}} \mathfrak{t}^{\frac{7}{6}-n} 2^{-2nl_a-\frac{5j}{2}} \Vert \Psi_{(l_a, l_a+k, l_a+k)}^{a-int} \psi_{j, k}^{\widehat{\mathcal{L}}} \Vert_{L^2} \langle 2^j \rangle^{-1} \Vert \langle \nabla \rangle |\nabla|^{-1} f(t) \Vert_{L^2} \\
&\lesssim \sum_{l_a = l_0}^{-200} 1_{l_a < k} t^{-\frac{7}{6}} \mathfrak{t}^{\frac{7}{6}-n} 2^{-2nl_a+\frac{3l_a}{2}-j} \langle 2^j \rangle^{-1} \Vert u \Vert_X \\
&\lesssim t^{-\frac{7}{6}} \mathfrak{t}^{\frac{7}{6}-2n\varrho} 2^{-j+k} \langle 2^j \rangle^{-1} \Vert u \Vert_X \\
&\lesssim t^{-\frac{7}{6}} 2^{-j+k} \langle 2^j \rangle^{-1} \Vert u \Vert_X \\
\eqref{estdispsupLmoyentauaint-2} &\lesssim \sum_{i = 1}^n \sum_{l_a = l_0}^{-200} 1_{l_a < k} t^{-\frac{7}{6}} \mathfrak{t}^{\frac{1}{6}-i} 2^{-2il_a-l_a-\frac{5j}{2}} \Vert \Psi_{(l_a, l_a+k, l_a+k)}^{a-int} \psi_{j, k}^{\widehat{\mathcal{L}}} \Vert_{L^2_{b, c} L^{\frac{1}{1-\kappa}}_a} \Vert \psi_j \widehat{h}_a(t) \Vert_{L^2_{b, c} L^{\frac{1}{\kappa}}_a} \\
&\lesssim \sum_{i = 1}^n \sum_{l_a = l_0}^{-200} 1_{l_a < k} t^{-\frac{7}{6}} \mathfrak{t}^{\frac{1}{6}-i} 2^{-2il_a-\kappa l_a+l_a-\frac{j}{2}-\kappa j} \Vert \psi_j \widehat{h}_a(t) \Vert_{L^2_{b, c} H^1_a}^{\frac{1}{2}} \Vert \psi_j \widehat{h}_a(t) \Vert_{L^2}^{\frac{1}{2}} \\
&\lesssim t^{-\frac{7}{6}+\frac{\kappa}{3}} \mathfrak{t}^{-\frac{1}{3}+\frac{\kappa}{6}+\frac{\delta}{2}} 2^{-j+\delta k} \langle 2^j \rangle^{-\frac{1}{2}} \Vert u \Vert_X \\
&\lesssim t^{-\frac{7}{6}+\frac{\kappa}{3}} 2^{-j+\delta k} \langle 2^j \rangle^{-1} \Vert u \Vert_X \\
\eqref{estdispsupLmoyentauaint-3} &\lesssim \sum_{i = 1}^n \sum_{l_a = l_0}^{-200} 1_{l_a < k} t^{-\frac{7}{6}} \mathfrak{t}^{\frac{1}{6}-i} 2^{-2il_a-\frac{5j}{2}} \Vert \Psi_{(l_a, l_a+k, l_a+k)}^{a-int} \psi_{j, k}^{\widehat{\mathcal{L}}} \Vert_{L^2} \Vert \nabla X_a h_a(t) \Vert_{L^2} \\
&\lesssim \sum_{i = 1}^n \sum_{l_a = l_0}^{-200} 1_{l_a < k} t^{-\frac{7}{6}} \mathfrak{t}^{\frac{1}{6}-i} 2^{-2il_a+\frac{3l_a}{2}-j} \Vert u \Vert_X \\
&\lesssim t^{-\frac{7}{6}} \mathfrak{t}^{-\frac{7}{12}+\frac{\delta}{2}} 2^{-j+\delta k} \Vert u \Vert_X \\
&\lesssim t^{-\frac{7}{6}} 2^{-j+\delta k} \langle 2^j \rangle^{-1} \Vert u \Vert_X \\
\eqref{estdispsupLmoyentauaint-4} &\lesssim \sum_{i = 1}^n \sum_{l_a = l_0}^{-200} 1_{l_a < k} t^{-1} \mathfrak{t}^{-i+1} 2^{-2il_a+l_a-\frac{5j}{2}} \Vert \Psi_{(l_a, l_a+k, l_a+k)}^{a-int} \psi_{j, k}^{\widehat{\mathcal{L}}} \Vert_{L^2} \langle 2^j \rangle^{-1} \Vert m_{\widehat{\mathcal{L}}}(D) \langle \nabla \rangle |\nabla|^{\frac{1}{2}} g_a(t) \Vert_{L^2} \\
&\lesssim \sum_{i = 1}^n \sum_{l_a = l_0}^{-200} 1_{l_a < k} t^{-\frac{7}{6}+100\delta} \mathfrak{t}^{-i+1} 2^{-2il_a+\frac{5l_a}{2}-j} \langle 2^j \rangle^{-1} \Vert u \Vert_X \\
&\lesssim t^{-\frac{7}{6}+100\delta} 2^{-j+\frac{k}{4}} \langle 2^j \rangle^{-1} \Vert u \Vert_X
\end{align*}
as wanted. 

For \eqref{estdispsupLmoyentaub}, \eqref{estdispsupLmoyentauc} and \eqref{estdispsupLmoyentauint}, we exhaust the cases depending on $j, k, t$. 

First, if $2^j \lesssim t^{-\frac{2}{9}}$ and $2^k \lesssim t^{-\frac{5}{12}+10\delta} 2^{-j}$, then we may choose $F(l_0, j, k) = 2k+100$ and only \eqref{estdispsupLmoyentauint} remains. We can then apply exactly the same estimates as in the case $\tau \ll 2^k$. We proceed the samy way if $t^{-\frac{2}{9}} \lesssim 2^j \lesssim t^{-\frac{5}{24}}$ and $2^k \lesssim t^{-\frac{1}{4}+10\delta} 2^{-\frac{j}{4}}$, or if $t^{-\frac{5}{24}} \lesssim 2^j$ and $2^k \lesssim \mathfrak{t}^{-\frac{1}{2}}$. We now assume none of these conditions are satisfied. 

Assume first $2^j \lesssim t^{-\frac{1}{6}}$. In this case, we choose $F$ such that 
\begin{align*}
2^{F(l_0, k, j)} \simeq t^{-\frac{5}{6}+20\delta} 2^{-2j} 
\end{align*}
It is easy to check that $F(l_0, k, j) \leq 2k+10$ under our hypothesis. 

On \eqref{estdispsupLmoyentaub}, we apply one integration by parts along $\widehat{X}_{b-corr}$: 
\begin{subequations}
\begin{align}
\eqref{estdispsupLmoyentaub} &= \sum_{l_b = F(l_0, k, j)}^{2k+100} t^{-1} 2^{-2l_b+2k-3j} \int e^{i t \Phi} \Psi_{(l_0, l_b, l_b)}^b \psi_{j, k}^{\widehat{\mathcal{L}}}(\overline{\xi}) \widehat{f}(t, \overline{\xi}) ~ d\overline{\xi} \label{estdispsupLmoyentaub-1} \\
&\quad + \sum_{l_b = F(l_0, k, j)}^{2k+100} t^{-1} 2^{-l_b-2j} \int e^{i t \Phi} \Psi_{(l_0, l_b, l_b)}^b \psi_{j, k}^{\widehat{\mathcal{L}}}(\overline{\xi}) \widehat{h}_{b-corr}(t, \overline{\xi}) ~ d\overline{\xi} \label{estdispsupLmoyentaub-2} \\
&\quad + \sum_{l_b = F(l_0, k, j)}^{2k+100} t^{-1} 2^{-l_b-2j} \int e^{i t \Phi} \Psi_{(l_0, l_b, l_b)}^b \psi_{j, k}^{\widehat{\mathcal{L}}}(\overline{\xi}) \widehat{g}_{b-corr}(t, \overline{\xi}) ~ d\overline{\xi} \label{estdispsupLmoyentaub-3} 
\end{align}
\end{subequations} 
We then estimate: 
\begin{align*}
\eqref{estdispsupLmoyentaub-1} &\lesssim \sum_{l_b = F(l_0, k, j)}^{2k+100} t^{-1} 2^{-2l_b+2k-3j} \Vert \Psi_{(l_0, l_b, l_b)}^b \psi_{j, k}^{\widehat{\mathcal{L}}} \Vert_{L^{\frac{4}{3}}_{b, c} L^{\frac{1}{1-\kappa}}_a} \langle 2^j \rangle^{-1} \Vert \langle \overline{\xi} \rangle \widehat{f}(t) \Vert_{L^4_{b, c} L^{\frac{1}{\kappa}}_a} \\
&\lesssim \sum_{l_b = F(l_0, k, j)}^{2k+100} t^{-1} 2^{l_0-\kappa l_0-\frac{l_b}{2}-\kappa j-\frac{j}{2}} \langle 2^j \rangle^{-1} \Vert u \Vert_X \\
&\lesssim t^{-\frac{13}{12}+\varrho+\frac{\kappa}{2}-10\delta} 2^{3\varrho j+\frac{\kappa j}{2}} 2^{-j+\delta j+\delta k} \langle 2^j \rangle^{-1} \Vert u \Vert_X \\
&\lesssim t^{-\frac{13}{12}} 2^{-j+\delta j+\delta k} \langle 2^j \rangle^{-1} \Vert u \Vert_X 
\end{align*}
Then, for \eqref{estdispsupLmoyentaub-2} and \eqref{estdispsupLmoyentaub-3}, we can reuse the estimate used in the case $\tau \ll 2^k$ up to losing a factor $|F(l_0, k, j)| \lesssim t^{\kappa}$ on the sum in $l_b$: indeed, the volume estimate $L^2_{b, c}$ of $\Psi_{(l_0, l_b, l_b)}^b \psi_{j, k}^{\widehat{\mathcal{L}}}$ gives precisely a factor $2^{l_b-k}$ that compensates $2^{-l_b}$ coming from the integration by parts. We then have exactly the term estimated in the case $\tau \ll 2^k$, summed $|F(l_0, k, j)|$ times. We skip the details. 

\eqref{estdispsupLmoyentauc} can be estimated the same way, replacing $\widehat{X}_b$ by $\widehat{X}_c$. 

For \eqref{estdispsupLmoyentauint}, we estimate directly: 
\begin{align*}
\eqref{estdispsupLmoyentauint} &\lesssim \Vert \Psi_{(l_0, F(l_0, k, j), F(l_0, k, j))} \psi_{j, k}^{\widehat{\mathcal{L}}} \Vert_{L^2_{b, c} L^{\frac{1}{1-\kappa}}_a} \Vert \psi_j \widehat{f}(t) \Vert_{L^2_{b, c} L^{\frac{1}{\kappa}}_a} \\
&\lesssim 2^{l_0-\kappa l_0+l_b-k+2j-\kappa j} \Vert \psi_j \widehat{f}(t) \Vert_{L^2_{b, c} H^1_a}^{\frac{1}{2}} \Vert \psi_j \widehat{f}(t) \Vert_{L^2}^{\frac{1}{2}} \\
&\lesssim t^{-\frac{13}{12}} t^{\frac{1}{6}+\varrho+\frac{\kappa}{2}+11\delta} 2^{j+3\varrho j+\frac{\kappa j}{2}} 2^{-j+\delta j+\delta k} \Vert u \Vert_X \\
&\lesssim t^{-\frac{13}{12}+12\delta} 2^{-j+\delta j+\delta k} \Vert u \Vert_X 
\end{align*}

As we assumed $2^j \lesssim t^{-\frac{1}{6}}$, we can use $1 \lesssim t^{-\frac{1}{12}} 2^{-\frac{j}{2}}$ to get automatically the other estimate of the Lemma. 

Assume now $2^j \gtrsim t^{-\frac{1}{6}}$. In this case, we choose $F$ such that
\begin{align*}
2^{F(l_0, k, j)} &\simeq 2^k \mathfrak{t}^{-\frac{2}{3}} 2^{-\frac{j}{2}} 
\end{align*}
Under our hypothesis, we have here that $2^k \gtrsim \mathfrak{t}^{-\frac{1}{2}}$, hence
\begin{align*}
2^{F(l_0, k, j)} &\lesssim 2^{2k} \mathfrak{t}^{-\frac{1}{6}} 2^{-\frac{j}{2}} \\
&\lesssim 2^{2k} t^{-\frac{1}{6}} 2^{-j} \\
&\lesssim 2^{2k} 
\end{align*}
which ensures we may choose $F$ such that $F(l_0, k, j) \leq 2k+100$. 

Again, on \eqref{estdispsupLmoyentaub}, we apply one integration by parts, and we can estimate \eqref{estdispsupLmoyentaub-2} and \eqref{estdispsupLmoyentaub-3} as in the case $\tau \ll 2^k$, up to the sum. It only remains
\begin{align*}
\eqref{estdispsupLmoyentaub-1} &\lesssim \sum_{l_b = F(l_0, k, j)}^{2k+100} t^{-1} 2^{-2l_b+2k-3j} \Vert \Psi_{(l_0, l_b, l_b)}^b \psi_{j, k}^{\widehat{\mathcal{L}}} \Vert_{L^{\frac{4}{3}}_{b, c} L^{\frac{1}{1-\kappa}}_a} \Vert \psi_j \widehat{f}(t) \Vert_{L^4_{b, c} L^{\frac{1}{\kappa}}_a} \\
&\lesssim \sum_{l_b = F(l_0, k, j)}^{2k+100} t^{-1} 2^{l_0-\kappa l_0-\frac{l_b}{2}+\frac{k}{2}-\frac{j}{2}-\kappa j} \Vert \psi_j \widehat{f}(t) \Vert_{L^2_{b, c} H^1_a}^{\frac{1}{2}} \Vert \psi_j \widehat{f}(t) \Vert_{\dot{H}^1_{b, c} L^2_a}^{\frac{1}{2}} \\
&\lesssim t^{-\frac{7}{6}} \mathfrak{t}^{\varrho+\frac{\kappa}{2}+\frac{\delta}{2}} 2^{-\kappa j} 2^{-j+\delta k} \langle 2^j \rangle^{-\frac{3}{4}} \Vert u \Vert_X \Vert \\
&\lesssim t^{-\frac{7}{6}+\varrho+\frac{\kappa}{2}+\frac{\delta}{2}} 2^{\frac{\kappa j}{2}+3\varrho j+\frac{\delta j}{2}} 2^{-j+\delta j+\delta k} \langle 2^j \rangle^{-\frac{3}{4}} \Vert u \Vert_X \Vert \\
&\lesssim t^{-\frac{7}{6}+\delta} 2^{-j+\delta j+\delta k} \langle 2^j \rangle^{-\frac{3}{5}} \Vert u \Vert_X 
\end{align*}
which is enough. \eqref{estdispsupLmoyentauc} can be estimated in a symmetric way. On the other hand, we can estimate the internal term as: 
\begin{align*}
\eqref{estdispsupLmoyentauint} &\lesssim \Vert \Psi_{(l_0, F(l_0, k, j), F(l_0, k, j))}^{int} \psi_{k, j} \Vert_{L^{\frac{4}{3}}_{b, c} L^{\frac{1}{1-\kappa}}_a} \Vert \psi_j \widehat{f}(t) \Vert_{L^4_{b, c} L^{\frac{1}{\kappa}}_a} \\
&\lesssim 2^{l_0-\kappa l_0+\frac{3F(l_0, k, j)}{2}-\frac{3k}{2}+\frac{5j}{2}-\kappa j} 2^{-\frac{j}{4}} \langle 2^j \rangle^{-\frac{3}{4}} \Vert u \Vert_X \\
&\lesssim t^{-\frac{7}{6}} \mathfrak{t}^{-\frac{1}{3}+\frac{\delta}{2}+\varrho+\frac{\kappa}{2}} 2^{-j-\kappa j} 2^{-j+\delta k} \langle 2^j \rangle^{-\frac{3}{4}} \Vert u \Vert_X \\
&\lesssim t^{-\frac{7}{6}+\delta} t^{-\frac{1}{3}} 2^{-2j} 2^{-j+\delta k} \langle 2^j \rangle^{-\frac{3}{5}} \Vert u \Vert_X \\
&\lesssim t^{-\frac{7}{6}+\delta} 2^{-j+\delta k} \langle 2^j \rangle^{-\frac{3}{5}} \Vert u \Vert_X
\end{align*}
This concludes the case $\tau \simeq 2^k$. 

\paragraph{3.} Finally, let us assume that $\tau \gg 2^k$. Here, we don't need to assume that $|y| \ll |x|$. 

In this case, we may reuse the localisation symbols used for the case $\tau \ll 2^k$ and note that 
\begin{align*}
\nabla_{\xi} \Phi &= 2 \xi_0 \xi - \frac{y}{t} \\
&\simeq \tau 2^{2j} \gg 2^{2j+k} 
\end{align*}
In particular, we may apply simpler estimates than in the case $\tau \ll 2^k$. We skip the details. 

\paragraph{Low frequencies} In the case $\overline{\xi_a} \gg 2^j$, we can again restrict our attention to the case $\widehat{X}_a \cdot \nabla \Phi \ll 2^{2j}$, but this forces $|y| \gtrsim |x|$. 

Then, we have that $\nabla_{\xi} \Phi \simeq |\overline{\xi_a}|^2$, and we can apply simpler estimates than in the resonant case. We skip the details. 
\end{Dem}

\subsection{Neighborhood of the plane}

\begin{Lem} Let $t \geq 1$, $j, k \in \mathbb{Z}$ be such that $2^j \gg t^{-\frac{1}{3}}$ and $k \leq -10$. Then
\begin{align*}
\Vert e^{i t \omega(D)} \psi_{j, k}^{\widehat{\mathcal{P}}}(D) f(t) \Vert_{L^{\infty}} &\lesssim t^{-1+\delta} 2^{-\frac{j}{2}+\delta j+\delta k} \langle 2^j \rangle^{-\frac{6}{5}} \Vert u \Vert_X \\
\Vert e^{i t \omega(D)} \psi_{j, k}^{\widehat{\mathcal{P}}}(D) f(t) \Vert_{L^{\infty}} &\lesssim t^{-\frac{7}{6}+100\delta} 2^{-j-\frac{k}{3}+\delta j+\delta k} \langle 2^j \rangle^{-\frac{3}{5}} \Vert u \Vert_X 
\end{align*}
\label{lemestdispvoisP} 
\end{Lem}

\begin{Dem}
As before, fixing $(x, y), t$ we need to estimate
\begin{align}
e^{tL(D)} \psi_{j, k}^{\widehat{\mathcal{P}}}(D) f(t, x, y) &= \int e^{i t \Phi} \psi_{j, k}^{\widehat{\mathcal{P}}}(\overline{\xi}) \widehat{f}(t, \overline{\xi}) ~ d\overline{\xi} \label{estdispvoisP-termetot}
\end{align}

Recall that, in the neighborhood of $\widehat{\mathcal{P}}$, $m_b \simeq m_c \simeq 1$. We will denote by $\tau = \frac{|y|}{|x|+|y|}$. 

We have that
\begin{align*}
\widehat{X}_c(\overline{\xi}) \cdot \nabla_{\overline{\xi}} \Phi &= \frac{\xi \cdot Jy}{t |\xi|} \\
\widehat{X}_a(\overline{\xi}) \cdot \nabla_{\overline{\xi}} \Phi &= 3 \xi_0 |\overline{\xi}| - \frac{\xi_0 x}{t |\overline{\xi}|} - \frac{\xi \cdot y}{t |\overline{\xi}|} \\
\widehat{X}_b(\overline{\xi}) \cdot \nabla_{\overline{\xi}} \Phi &= |\xi| |\overline{\xi}| - \frac{|\xi| x}{t |\overline{\xi}|} - \frac{\xi_0 \xi \cdot y}{t |\overline{\xi}| |\xi|} 
\end{align*}

\paragraph{Low $x$-frequencies} Assume first that $2^k \lesssim \mathfrak{t}^{-1}$. Then we may directly estimate for a small enough $\kappa > 0$: 
\begin{align*}
\eqref{estdispvoisP-termetot} &\lesssim \Vert \psi_{j, k}^{\widehat{\mathcal{P}}} \Vert_{L^{\frac{1}{1-\kappa}}_{\xi_0} L^2_{\xi}} \Vert \psi_j \widehat{f}(t) \Vert_{L^{\frac{1}{\kappa}}_{\xi_0} L^2_{\xi}} \\
&\lesssim 2^{k-\kappa k+2j-\kappa j} \Vert \psi_j \widehat{f}(t) \Vert_{H^1_{\xi_0} L^2_{\xi}}^{\frac{1}{2}} \Vert \psi_j \widehat{f}(t) \Vert_{L^2}^{\frac{1}{2}} \\
&\lesssim \mathfrak{t}^{-1+\kappa+\delta} 2^{\frac{5j}{2}-\kappa j + \delta k} \langle 2^j \rangle^{-\frac{3}{2}} \Vert u \Vert_X \\
&\lesssim t^{-1+2\delta} 2^{-\frac{j}{2}+\delta j+\delta k} \langle 2^j \rangle^{-\frac{5}{4}} \Vert u \Vert_X 
\end{align*}
The other estimate is automatic noting that, here, $1 \lesssim t^{-\frac{1}{6}} 2^{-\frac{k}{6}-\frac{j}{2}}$. 

In what follows, we always assume $2^k \gtrsim \mathfrak{t}^{-1}$. 

\paragraph{High frequencies} Let us assume that $2^j \gg |\overline{\xi_a}|$. Then, 
\begin{align*}
\widehat{X}_b(\overline{\xi}) \cdot \nabla_{\overline{\xi}} \Phi &\simeq 2^{2j} 
\end{align*}
On the other hand, we may consider three subcases in $(x, y)$. If $2^{k+2j} \gg \tau |\overline{\xi_a}|^2$, then 
\begin{align*}
\widehat{X}_a(\overline{\xi}) \cdot \nabla_{\overline{\xi}} \Phi &\simeq 2^{k+2j} 
\end{align*}
If $2^{k+2j} \ll \tau |\overline{\xi_a}|^2$, then 
\begin{align*}
\left| \widehat{X}_a(\overline{\xi}) \cdot \nabla_{\overline{\xi}} \Phi \right| + \left| \widehat{X}_c(\overline{\xi}) \cdot \nabla_{\overline{\xi}} \Phi \right| &\simeq \tau |\overline{\xi_a}|^2
\end{align*}
Finally, if $2^{k+2j} \simeq \tau |\overline{\xi_a}|^2$, there may be a resonance in $\widehat{X}_a, \widehat{X}_c$. 

In any case, 
\begin{align*}
\widehat{X}_a(\overline{\xi}) \cdot \nabla_{\overline{\xi}} \psi_{j, k}^{\widehat{\mathcal{P}}}(\overline{\xi}) &\lesssim 2^{-j} \\
\widehat{X}_b(\overline{\xi}) \cdot \nabla_{\overline{\xi}} \psi_{j, k}^{\widehat{\mathcal{P}}}(\overline{\xi}) &\lesssim 2^{-j-k} \\
\widehat{X}_c(\overline{\xi}) \cdot \nabla_{\overline{\xi}} \psi_{j, k}^{\widehat{\mathcal{P}}}(\overline{\xi}) &\lesssim 2^{-j} 
\end{align*}

\paragraph{1.} Assume first that $2^j \gg |\overline{\xi_a}|$ and $2^{k+2j} \gg \tau |\overline{\xi_a}|^2$, with $2^k \gtrsim \mathfrak{t}^{-1}$. We can then apply two integration by parts, one along $\widehat{X}_b$, then again along $\widehat{X}_b$ if the first hits the symbol, or along $\widehat{X}_a$ if it hits $\widehat{f}$:  
\begin{subequations}
\begin{align}
\eqref{estdispvoisP-termetot} &= t^{-2} 2^{-6j-2k} \int e^{i t \Phi} \psi_{j, k}^{\widehat{\mathcal{P}}}(\overline{\xi}) \widehat{f}(t, \overline{\xi}) ~ d\overline{\xi} \label{estdispsupPHF1-1} \\
&\quad + t^{-2} 2^{-5j-k} \int e^{i t \Phi} \psi_{j, k}^{\widehat{\mathcal{P}}}(\overline{\xi}) \widehat{X}_b(\overline{\xi}) \cdot \nabla_{\overline{\xi}} \widehat{f}(t, \overline{\xi}) ~ d\overline{\xi} \label{estdispsupPHF1-2} \\
&\quad + t^{-2} 2^{-5j-k} \int e^{i t \Phi} \psi_{j, k}^{\widehat{\mathcal{P}}}(\overline{\xi}) \widehat{h}_b(t, \overline{\xi}) ~ d\overline{\xi} \label{estdispsupPHF1-3} \\
&\quad + t^{-2} 2^{-4j-k} \int e^{i t \Phi} \psi_{j, k}^{\widehat{\mathcal{P}}}(\overline{\xi}) \widehat{X}_a(\overline{\xi}) \cdot \nabla_{\overline{\xi}} \widehat{h}_b(t, \overline{\xi}) ~ d\overline{\xi} \label{estdispsupPHF1-4} \\
&\quad + t^{-1} 2^{-2j} \int e^{i t \Phi} \psi_{j, k}^{\widehat{\mathcal{P}}}(\overline{\xi}) \widehat{g}_b(t, \overline{\xi}) ~ d\overline{\xi} \label{estdispsupPHF1-5} 
\end{align}
\end{subequations} 
where, as before, the symbols can change from line to line as long as they keep similar properties. We then estimate: 
\begin{align*}
\eqref{estdispsupPHF1-1} &\lesssim t^{-2} 2^{-6j-2k} \Vert \psi_{j, k}^{\widehat{\mathcal{P}}} \Vert_{L^1_{\xi_0} L^2_{\xi}} \Vert \psi_j \widehat{f}(t) \Vert_{L^{\infty}_{\xi_0} L^2_{\xi}} \\
&\lesssim t^{-2} 2^{-k-\frac{7j}{2}} \langle 2^j \rangle^{-\frac{3}{2}} \Vert u \Vert_X \\
\eqref{estdispsupPHF1-2} &\lesssim t^{-2} 2^{-5j-k} \Vert \psi_{j, k}^{\widehat{\mathcal{P}}} \Vert_{L^2} \langle 2^j \rangle^{-1} \Vert \langle \nabla \rangle X_b f(t) \Vert_{L^2} \\
&\lesssim t^{-2} 2^{-\frac{k}{2}-\frac{7j}{2}} \langle 2^j \rangle^{-1} \Vert u \Vert_X \\
&\lesssim t^{-\frac{7}{6}+\delta} 2^{-j+\delta j+\delta k} \langle 2^j \rangle^{-1} \Vert u \Vert_X \\
\eqref{estdispsupPHF1-5} &\lesssim t^{-1} 2^{-\frac{5j}{2}} \Vert \psi_{j, k}^{\widehat{\mathcal{P}}} \Vert_{L^2} \langle 2^j \rangle^{-1} \Vert \langle \nabla \rangle |\nabla|^{\frac{1}{2}} g_b(t) \Vert_{L^2} \\
&\lesssim t^{-\frac{7}{6}+100\delta} 2^{-j+\frac{k}{2}} \langle 2^j \rangle^{-1} \Vert u \Vert_X 
\end{align*}
We can estimate \eqref{estdispsupPHF1-3} and \eqref{estdispsupPHF1-4} like \eqref{estdispsupPHF1-2}. We then have indeed
\begin{align*}
\eqref{estdispsupPHF1-1} &\lesssim t^{-1+\delta} 2^{-\frac{j}{2}+\delta j+\delta k} \langle 2^j \rangle^{-\frac{6}{5}} \Vert u \Vert_X 
\end{align*}
and 
\begin{align*}
\eqref{estdispsupPHF1-1} &\lesssim t^{-2} 2^{-\frac{5k}{6}-\frac{5j}{2}} 2^{-j-\frac{k}{6}} \langle 2^j \rangle^{-\frac{3}{2}} \Vert u \Vert_X \\
&\lesssim t^{-\frac{7}{6}} 2^{-j-\frac{k}{6}} \langle 2^j \rangle^{-\frac{3}{2}} \Vert u \Vert_X 
\end{align*}
as wanted. 

\paragraph{2.} Let us now assume that $2^j \gg \overline{\xi_a}$, $2^{k+2j} \ll \tau |\overline{\xi_a}|^2$, $2^k \gtrsim \mathfrak{t}^{-1}$. We can apply the same integrations by parts and the estimates are simpler, we skip the details. 

\paragraph{3.} Let us finally assume that $2^j \gg \overline{\xi_a}$, $2^{k+2j} \simeq \tau |\overline{\xi_a}|^2$, $2^k \gtrsim \mathfrak{t}^{-1}$. 

We start by an integration by parts along $\widehat{X}_b$: 
\begin{subequations}
\begin{align}
\eqref{estdispvoisP-termetot} &= t^{-1} 2^{-3j-k} \int e^{i t \Phi} \psi_{j, k}^{\widehat{\mathcal{P}}}(\overline{\xi}) \widehat{f}(t, \overline{\xi}) ~ d\overline{\xi} \label{estdispsupPHF4-1} \\
&\quad + t^{-1} 2^{-2j} \int e^{i t \Phi} \psi_{j, k}^{\widehat{\mathcal{P}}}(\overline{\xi}) \widehat{h}_b(t, \overline{\xi}) ~ d\overline{\xi} \label{estdispsupPHF4-2} \\
&\quad + t^{-1} 2^{-2j} \int e^{i t \Phi} \psi_{j, k}^{\widehat{\mathcal{P}}}(\overline{\xi}) \widehat{g}_b(t, \overline{\xi}) ~ d\overline{\xi} \label{estdispsupPHF4-3} 
\end{align}
\end{subequations} 
We can apply a second integration by parts along $\widehat{X}_b$ on \eqref{estdispsupPHF4-1} and estimate just like in the case $2^{k+2j} \gg \tau |\overline{\xi_a}|^2$. We can also estimate \eqref{estdispsupPHF4-3} as before. 

Then, to apply integrations by parts on \eqref{estdispsupPHF4-2}, we introduce the localisation symbols 
\begin{align*}
\Psi_{l_a}^a\left( \overline{\xi}, \frac{x}{t}, \frac{y}{t} \right) &:= \psi\left( 2^{-l_a-2j} \widehat{X}_a(\overline{\xi}) \cdot \nabla_{\overline{\xi}} \Phi \right) \chi\left( 2^{-l_a-2j} \widehat{X}_c(\overline{\xi}) \cdot \nabla_{\overline{\xi}} \Phi \right) \\
\Psi_{l_a}^c\left( \overline{\xi}, \frac{x}{t}, \frac{y}{t} \right) &:= \chi\left( 2^{-l_a-2j} \widehat{X}_a(\overline{\xi}) \cdot \nabla_{\overline{\xi}} \Phi \right) \psi\left( 2^{-l_a-2j} \widehat{X}_c(\overline{\xi}) \cdot \nabla_{\overline{\xi}} \Phi \right) \\
\Psi_{l_a}^{int}\left( \overline{\xi}, \frac{x}{t}, \frac{y}{t} \right) &:= \chi\left( 2^{-l_a-2j} \widehat{X}_a(\overline{\xi}) \cdot \nabla_{\overline{\xi}} \Phi \right) \chi\left( 2^{-l_a-2j} \widehat{X}_c(\overline{\xi}) \cdot \nabla_{\overline{\xi}} \Phi \right) 
\end{align*}
Note that, here, we automatically have 
\begin{align*}
\widehat{X}_a(\overline{\xi}) \cdot \nabla_{\overline{\xi}} \Phi, ~ \widehat{X}_c(\overline{\xi}) \cdot \nabla_{\overline{\xi}} \Phi &\lesssim 2^{k+2j} 
\end{align*}

We may then decompose 
\begin{subequations}
\begin{align}
\eqref{estdispsupPHF4-2} &= \sum_{l_a = l_0}^{k+10} t^{-1} 2^{-2j} \int e^{i t \Phi} \Psi_{l_a}^a \psi_{j, k}^{\widehat{\mathcal{P}}}(\overline{\xi}) \widehat{h}_b(t, \overline{\xi}) ~ d\overline{\xi} \label{estdispsupPHF4-2-1} \\
&\quad + \sum_{l_a = l_0}^{k+10} t^{-1} 2^{-2j} \int e^{i t \Phi} \Psi_{l_a}^c \psi_{j, k}^{\widehat{\mathcal{P}}}(\overline{\xi}) \widehat{h}_b(t, \overline{\xi}) ~ d\overline{\xi} \label{estdispsupPHF4-2-2} \\
&\quad + t^{-1} 2^{-2j} \int e^{i t \Phi} \Psi_{l_0}^{int} \psi_{j, k}^{\widehat{\mathcal{P}}}(\overline{\xi}) \widehat{h}_b(t, \overline{\xi}) ~ d\overline{\xi} \label{estdispsupPHF4-2-3} 
\end{align}
\end{subequations} 
for $l_0$ such that $2^{l_0} \simeq \mathfrak{t}^{-1}$: in particular, since $2^k \gtrsim \mathfrak{t}^{-1}$, we have $l_0 \leq k+10$. 

Let us introduce local coordinates near a point $\overline{\xi}^Z$ by 
\begin{align*}
\Lambda : \overline{\xi} \mapsto \left( |\overline{\xi}|, |\overline{\xi}^Z| \frac{\xi_0}{|\overline{\xi}|}, |\overline{\xi}^Z| \frac{\xi \cdot J \xi^Z}{|\xi| |\xi^Z|} \right) 
\end{align*}
Then, 
\begin{align*}
\widehat{X}_a(\overline{\xi}) \cdot \nabla_{\overline{\xi}} \Lambda(\overline{\xi}) &= \left( 1, 0, 0 \right) \\
\widehat{X}_b(\overline{\xi}) \cdot \nabla_{\overline{\xi}} \Lambda(\overline{\xi}) &= \left( 0, \frac{|\overline{\xi}^Z| |\xi|}{|\overline{\xi}|^2}, 0 \right) \\
\widehat{X}_c(\overline{\xi}) \cdot \nabla_{\overline{\xi}} \Lambda(\overline{\xi}) &= \left( 0, 0, \frac{|\overline{\xi}^Z| \xi \cdot \xi^Z}{|\xi|^2 |\xi^Z|} \right) 
\end{align*}
so that, in a neighborhood of $\overline{\xi}^Z$, $\Lambda$ is invertible and 
\begin{align*}
\partial_{\alpha} \left( G \circ \Lambda^{-1} \right) \simeq \left( \widehat{X}_{\alpha} \cdot \nabla G \right) \circ \Lambda^{-1} 
\end{align*}
We denote by $(\xi_a, \xi_b, \xi_c)$ the coordinates defined by $\Lambda$. 

Without loss of generality, since we consider an essential supremum in $y$, we may assume $y \neq 0$ and choose $\overline{\xi}^Z$ the reference point such that $|\overline{\xi}^Z| \simeq 2^j$, $\frac{\xi^Z}{|\xi^Z|} = \frac{y}{|y|}$. 

Let us compute volume estimates in the coordinates introduced by $\Lambda$: $\Psi_{l_a}^{*}$, where $* \in \{ a, c, int \}$, localises on the one hand 
\begin{align*}
\widehat{X}_c(\overline{\xi}) \cdot \nabla_{\overline{\xi}} \Phi &\simeq \tau |\overline{\xi_a}|^2 \frac{\xi_c}{2^j} \lesssim 2^{l_a+2j}
\end{align*}
and therefore localises $\xi_c$ close to $0$ with a precision $2^{l_a-k+j}$. Then, 
\begin{align*}
\widehat{X}_a(\overline{\xi}) \cdot \nabla_{\overline{\xi}} \Phi &= 3 \xi_0 |\overline{\xi}| \left( 1 - \frac{x}{3 t |\overline{\xi}|^2} \right) - \frac{\xi \cdot y}{t |\overline{\xi}|} 
\end{align*}
As before, the first term is of size $2^{2j+k}$, while the second is of size $\tau |\overline{\xi_a}|^2$. Therefore, if $l_a \leq k+10$, we deduce that, $\xi_a, \xi_c$ being fixed, this localises $\xi_b$ (or $\xi_0$) of size $2^{k+j}$ with precision $2^{l_a+j}$. Finally, we deduce that 
\begin{align*}
\Vert \Psi_{l_a}^{*} \psi_{j, k}^{\widehat{\mathcal{P}}} \Vert_{L^{p_a}_a L^{p_c}_c L^{p_b}_b} &\lesssim 2^{\frac{l_a-k}{p_c}+\frac{l_a}{p_b}+\frac{j}{p_a}+\frac{j}{p_b}+\frac{j}{p_c}}
\end{align*}
Note that, unlike previous dispersive estimates, near the plane, it is $\widehat{X}_a$ that localises the local coordinate $\xi_b$ (and, as we will see later, $\widehat{X}_b$ localises $\xi_a$). 

The symbol derivatives are computed through: 
\begin{align*}
\widehat{X}_a(\overline{\xi}) \cdot \nabla_{\overline{\xi}} \left[ \widehat{X}_a(\overline{\xi}) \cdot \nabla_{\overline{\xi}} \Phi \right] &\simeq \xi_0 \simeq 2^{j+k} \\
\widehat{X}_a(\overline{\xi}) \cdot \nabla_{\overline{\xi}} \left[ \widehat{X}_c(\overline{\xi}) \cdot \nabla_{\overline{\xi}} \Phi \right] &= 0 \\
\widehat{X}_c(\overline{\xi}) \cdot \nabla_{\overline{\xi}} \left[ \widehat{X}_a(\overline{\xi}) \cdot \nabla_{\overline{\xi}} \Phi \right] &= - \frac{J \xi \cdot y}{t |\overline{\xi}| |\xi|} \\
&\lesssim 2^{l_a+j} \\
\widehat{X}_c(\overline{\xi}) \cdot \nabla_{\overline{\xi}} \left[ \widehat{X}_c(\overline{\xi}) \cdot \nabla_{\overline{\xi}} \Phi \right] &= \frac{\xi \cdot y}{t |\xi|^2} \\
&\lesssim \tau |\overline{\xi_a}|^2 2^{-j} \simeq 2^{k+j} 
\end{align*}
As before, $\widehat{X}_a \cdot \nabla \psi_{j, k}^{\widehat{\mathcal{P}}} \lesssim 2^{-j}$ and $\widehat{X}_c \cdot \nabla \psi_{j, k}^{\widehat{\mathcal{P}}} = 0$. We therefore have
\begin{align*}
\widehat{X}_a(\overline{\xi}) \cdot \nabla_{\overline{\xi}} \left[ \Psi_{l_a}^{*} \psi_{j, k}^{\widehat{\mathcal{P}}} \right] &\lesssim 2^{-l_a+k-j} \\
\widehat{X}_c(\overline{\xi}) \cdot \nabla_{\overline{\xi}} \left[ \Psi_{l_a}^{*} \psi_{j, k}^{\widehat{\mathcal{P}}} \right] &\lesssim 2^{-l_a+k-j} 
\end{align*}

On \eqref{estdispsupPHF4-2-1}, we apply an integration by parts along $\widehat{X}_a$: 
\begin{subequations}
\begin{align}
\eqref{estdispsupPHF4-2-1} &= \sum_{l_a = l_0}^{k+10} t^{-2} 2^{-2l_a+k-5j} \int e^{i t \Phi} \Psi_{l_a}^a \psi_{j, k}^{\widehat{\mathcal{P}}}(\overline{\xi}) \widehat{h}_b(t, \overline{\xi}) ~ d\overline{\xi} \label{estdispsupPHF4-2-1-1} \\
&\quad + \sum_{l_a = l_0}^{k+10} t^{-2} 2^{-l_a-4j} \int e^{i t \Phi} \Psi_{l_a}^a \psi_{j, k}^{\widehat{\mathcal{P}}}(\overline{\xi}) \widehat{X}_a(\overline{\xi}) \cdot \nabla_{\overline{\xi}} \widehat{h}_b(t, \overline{\xi}) ~ d\overline{\xi} \label{estdispsupPHF4-2-1-2} 
\end{align}
\end{subequations} 
We then estimate: 
\begin{align*}
\eqref{estdispsupPHF4-2-1-1} &\lesssim \sum_{l_a = l_0}^{k+10} t^{-2} 2^{-2l_a+k-5j} \Vert \Psi_{l_a}^a \psi_{j, k}^{\widehat{\mathcal{P}}} \Vert_{L^2_a L^{\frac{1}{1-\delta}}_c L^2_b} \Vert \psi_j \widehat{h}_b(t) \Vert_{L^2_a L^{\frac{1}{\delta}}_c L^2_b} \\
&\lesssim \sum_{l_a = l_0}^{k+10} t^{-2} 2^{-\frac{l_a}{2}-\delta l_a-3j-\delta j+\delta k} \Vert \psi_j \widehat{h}_b(t) \Vert_{L^2_a H^1_c L^2_b}^{\frac{1}{2}} \Vert \psi_j \widehat{h}_b(t) \Vert_{L^2}^{\frac{1}{2}} \\
&\lesssim t^{-\frac{3}{2}+\delta} 2^{-2j+2\delta j+\delta k} \langle 2^j \rangle^{-\frac{1}{2}} \Vert u \Vert_X \\
&\lesssim t^{-\frac{7}{6}+\frac{2\delta}{3}} 2^{-j+\delta j+\delta k} \langle 2^j \rangle^{-1} \Vert u \Vert_X \\
\eqref{estdispsupPHF4-2-1-2} &\lesssim \sum_{l_a = l_0}^{k+10} t^{-2} 2^{-l_a-5j} \Vert \Psi_{l_a}^a \psi_{j, k}^{\widehat{\mathcal{P}}} \Vert_{L^2} \Vert \nabla X_a h_b(t) \Vert_{L^2} \\
&\lesssim \sum_{l_a = l_0}^{k+10} t^{-2} 2^{-\frac{k}{2}-\frac{7j}{2}} \Vert u \Vert_X \\
&\lesssim t^{-\frac{3}{2}+2\delta} 2^{-j+5\delta j} 2^{-j+\delta j+\delta k} \Vert u \Vert_X \\
&\lesssim t^{-\frac{7}{6}+\frac{\delta}{3}} 2^{-j+\delta j+\delta k} \langle 2^j \rangle^{-\frac{2}{3}} \Vert u \Vert_X 
\end{align*}
We also have, using again $2^j \gtrsim t^{-\frac{1}{3}}$ that 
\begin{align*}
\eqref{estdispsupPHF4-2-1-1} &\lesssim t^{-1+\delta} 2^{-\frac{j}{2}+\delta j+\delta k} \langle 2^j \rangle^{-2} \Vert u \Vert_X \\
\eqref{estdispsupPHF4-2-1-2} &\lesssim t^{-1+\frac{\delta}{3}} 2^{-\frac{j}{2}+\delta j+\delta k} \langle 2^j \rangle^{-\frac{2}{3}} \Vert u \Vert_X 
\end{align*}
as wanted. 

We can estimate \eqref{estdispsupPHF4-2-2} the same way, replacing $\widehat{X}_a$ by $\widehat{X}_c$. 

Finally, for \eqref{estdispsupPHF4-2-3}, we estimate directly: 
\begin{align*}
\eqref{estdispsupPHF4-2-3} &\lesssim t^{-1} 2^{-2j} \Vert \Psi_{l_0}^{int} \psi_{j, k}^{\widehat{\mathcal{P}}} \Vert_{L^2_a L^{\frac{1}{1-\delta}}_c L^2_b} \Vert \psi_j \widehat{h}_b(t) \Vert_{L^2_a L^{\frac{1}{\delta}}_c L^2_b} \\
&\lesssim t^{-1} 2^{\frac{3l_0}{2}-\delta l_0-k-\delta j} \Vert \psi_j \widehat{h}_b(t) \Vert_{L^2_a H^1_c L^2_b}^{\frac{1}{2}} \Vert \psi_j \widehat{h}_b(t) \Vert_{L^2}^{\frac{1}{2}} \\
&\lesssim t^{-\frac{3}{2}+2\delta} 2^{-j+4\delta j} 2^{-j+\delta j+\delta k} \langle 2^j \rangle^{-\frac{1}{2}} \Vert u \Vert_X \\
&\lesssim t^{-\frac{7}{6}+\frac{2\delta}{3}} 2^{-j+\delta j+\delta k} \langle 2^j \rangle^{-\frac{3}{4}} \Vert u \Vert_X 
\end{align*}
Using again $2^j \gtrsim t^{-\frac{1}{3}}$, we can obtain the other estimate. 

\paragraph{Resonant frequencies} Let us now assume that $\overline{\xi_a} \simeq 2^j$. 

As on previous geometric areas, it is easy to localise $\overline{\xi}$ such that
\begin{align*}
\widehat{X}_a \cdot \nabla \Phi, ~ \widehat{X}_c \cdot \nabla \Phi ~ \ll 2^{2j} 
\end{align*}
In particular, for this localisation not to be empty, this forces $(x, y)$ to be close to $\mathcal{L}$. Moreover, as long as $\widehat{X}_b \cdot \nabla \Phi \simeq 2^{2j}$, we can apply similar estimates to the ones used for high frequencies, and therefore we can also localise on $\widehat{X}_b \cdot \nabla \Phi \ll 2^{2j}$. 

\paragraph{1.} We first assume that $\overline{\xi_a} \simeq 2^j$ and $\tau = \frac{|y|}{|x|+|y|} \ll 2^k$. We can then use the localisation symbols 
\begin{align*}
\Psi_{l_b}^b\left( \overline{\xi}, \frac{x}{t}, \frac{y}{t} \right) &:= \psi\left( 2^{-l_b-2j} \widehat{X}_b(\overline{\xi}) \cdot \nabla_{\overline{\xi}} \Phi \right) \\
\Psi_{l_b}^{int}\left( \overline{\xi}, \frac{x}{t}, \frac{y}{t} \right) &:= \chi\left( 2^{-l_b-2j} \widehat{X}_b(\overline{\xi}) \cdot \nabla_{\overline{\xi}} \Phi \right)
\end{align*}
to decompose
\begin{subequations}
\begin{align}
\eqref{estdispvoisP-termetot} &= \sum_{l_b = l_0}^{-200} 1_{l_b \geq k} \int e^{i t \Phi} \Psi_{l_b}^b \psi_{j, k}^{\widehat{\mathcal{P}}}(\overline{\xi}) \widehat{f}(t, \overline{\xi}) ~ d\overline{\xi} \label{estdispsupPMFpetittau-b} \\
&\quad + \sum_{l_b = l_0}^{-200} 1_{l_b < k} \int e^{i t \Phi} \Psi_{l_b}^b \psi_{j, k}^{\widehat{\mathcal{P}}}(\overline{\xi}) \widehat{f}(t, \overline{\xi}) ~ d\overline{\xi} \label{estdispsupPMFpetittau-a} \\
&\quad + \int e^{i t \Phi} \Psi_{l_0}^{int} \psi_{j, k}^{\widehat{\mathcal{P}}}(\overline{\xi}) \widehat{f}(t, \overline{\xi}) ~ d\overline{\xi} \label{estdispsupPMFpetittau-int} 
\end{align}
\end{subequations}
where $l_0$ is such that $l_0 \geq 2k$ and will be specified later (it is not the same $l_0$ as for previous geometric areas). 

The localisation $\Psi_{l_b}^b \psi_{j, k}^{\widehat{\mathcal{P}}}$ ensures that 
\begin{align*}
\widehat{X}_b(\overline{\xi}) \cdot \nabla_{\overline{\xi}} \Phi
&\simeq 2^{l_b+2j} 
\end{align*}
Yet, we have that
\begin{align*}
\widehat{X}_b(\overline{\xi}) \cdot \nabla_{\overline{\xi}} \Phi &= |\xi| |\overline{\xi}| \left( 1 - \frac{x}{t |\overline{\xi}|^2} + \frac{\xi_0 \xi \cdot y}{t |\overline{\xi}|^2 |\xi|^2} \right)
\end{align*}
with $|\xi| |\overline{\xi}| \simeq 2^{2j}$, $\frac{x}{t |\overline{\xi}|^2} \gg \frac{\xi_0 \xi \cdot y}{t |\overline{\xi}|^2 |\xi|}$, and $l_b \leq -200$, so that $1 - \frac{x}{t |\overline{\xi}|^2} \ll 1$. In particular, 
\begin{align*}
\widehat{X}_a(\overline{\xi}) \cdot \nabla_{\overline{\xi}} \Phi &= 3 \xi_0 |\overline{\xi}| \left( 1 - \frac{x}{3 t |\overline{\xi}|^2} - \frac{\xi \cdot y}{3 t |\overline{\xi}|^2 \xi_0} \right) 
\end{align*}
Therefore, $1 - \frac{x}{3 t |\overline{\xi}|^2} \simeq 1$ while $\frac{\xi \cdot y}{t |\overline{\xi}|^2 \xi_0} \simeq \tau 2^{-k} \ll 1$ by hypothesis. We deduce that 
\begin{align*}
\widehat{X}_a(\overline{\xi}) \cdot \nabla_{\overline{\xi}} \Phi &\simeq 2^{k+2j} 
\end{align*}

Recall the local coordinates $\xi_a, \xi_b, \xi_c$ introduced in the high frequency estimate. It is clear that $\psi_{j, k}^{\widehat{\mathcal{P}}}$ localises $\xi_b \simeq 2^{j+k}$ with precision $2^{j+k}$. $\xi_c, \xi_b$ being fixed, $\Psi_{l_b}^{*}$ localises $\xi_a \simeq 2^j$ with precision $2^{j+l_b}$. Therefore, 
\begin{align*}
\Vert \Psi_{l_b}^{*} \psi_{j, k}^{\widehat{\mathcal{P}}} \Vert_{L^{p_c}_c L^{p_b}_b L^{p_a}_a} &\lesssim 2^{\frac{k}{p_b}+\frac{l_b}{p_a}+\frac{j}{p_a}+\frac{j}{p_b}+\frac{j}{p_c}}
\end{align*}

We then compute the symbol derivatives: 
\begin{align*}
\widehat{X}_a(\overline{\xi}) \cdot \nabla_{\overline{\xi}} \psi_{j, k}^{\widehat{\mathcal{P}}}(\overline{\xi}) &\lesssim 2^{-j} \\
\widehat{X}_b(\overline{\xi}) \cdot \nabla_{\overline{\xi}} \psi_{j, k}^{\widehat{\mathcal{P}}}(\overline{\xi}) &\lesssim 2^{-j-k} \\
\widehat{X}_a(\overline{\xi}) \cdot \nabla_{\overline{\xi}} \left[ \widehat{X}_b(\overline{\xi}) \cdot \nabla_{\overline{\xi}} \Phi \right] &= \frac{\overline{\xi}}{|\overline{\xi}|} \cdot \nabla_{\overline{\xi}} \left[ |\xi| |\overline{\xi}| \right] \\
&\simeq 2^j \\
\widehat{X}_b(\overline{\xi}) \cdot \nabla_{\overline{\xi}} \left[ \widehat{X}_b(\overline{\xi}) \cdot \nabla_{\overline{\xi}} \Phi \right] &= \left( \frac{|\xi|}{|\overline{\xi}|} \partial_{\xi_0} - \frac{\xi_0 \xi}{|\overline{\xi}| |\xi|} \cdot \nabla_{\xi} \right) \left[ |\xi| |\overline{\xi}| - \frac{|\xi| x}{t |\overline{\xi}|} + \frac{\xi_0 \xi \cdot y}{t |\overline{\xi}| |\xi|} \right] \\
&= - \xi_0 + \frac{\xi_0 x}{t |\overline{\xi}|^2} + \frac{\xi \cdot y}{t |\overline{\xi}|^2} \\
&= - \xi_0 \left( 1 - \frac{x}{t |\overline{\xi}|^2} - \frac{\xi \cdot y}{t |\overline{\xi}|^2 \xi_0} \right) \\
&\lesssim 2^{j+k} + \tau 2^j \lesssim 2^{j+k} 
\end{align*}
Therefore, 
\begin{align*}
\widehat{X}_a(\overline{\xi}) \cdot \nabla_{\overline{\xi}} \left[ \Psi_{l_b}^{*} \psi_{j, k}^{\widehat{\mathcal{P}}} \right] &\lesssim 2^{-l_b-j} \\
\widehat{X}_b(\overline{\xi}) \cdot \nabla_{\overline{\xi}} \left[ \Psi_{l_b}^{*} \psi_{j, k}^{\widehat{\mathcal{P}}} \right] &\lesssim 2^{-l_b+k-j} + 2^{-k-j} \lesssim 2^{-k-j} 
\end{align*}
since here $l_b \geq 2k$. 

We exhaust the cases depending on $k$ and $\mathfrak{t}$. 

\paragraph{1.1.} Let us first consider the case $2^k \gtrsim \mathfrak{t}^{-\frac{1}{2}}$. We then choose $l_0 = 2k(1 - \varrho)$ for some small $\varrho > 0$. On \eqref{estdispsupPMFpetittau-b}, we apply an integration by parts along $\widehat{X}_b$, then one along $\widehat{X}_a$: 
\begin{subequations}
\begin{align}
\eqref{estdispsupPMFpetittau-b} &= \sum_{l_b = k}^{-200} t^{-2} 2^{-2l_b-2k-6j} \int e^{i t \Phi} \Psi_{l_b}^b \psi_{j, k}^{\widehat{\mathcal{P}}}(\overline{\xi}) \widehat{f}(t, \overline{\xi}) ~ d\overline{\xi} \label{estdispsupPMFpetittau1-b-1} \\
&\quad + \sum_{l_b = k}^{-200} t^{-2} 2^{-l_b-2k-5j} \int e^{i t \Phi} \Psi_{l_b}^b \psi_{j, k}^{\widehat{\mathcal{P}}}(\overline{\xi}) \widehat{X}_a(\overline{\xi}) \cdot \nabla_{\overline{\xi}} \widehat{f}(t, \overline{\xi}) ~ d\overline{\xi} \label{estdispsupPMFpetittau1-b-2} \\
&\quad + \sum_{l_b = k}^{-200} t^{-2} 2^{-2l_b-k-5j} \int e^{i t \Phi} \Psi_{l_b}^b \psi_{j, k}^{\widehat{\mathcal{P}}}(\overline{\xi}) \widehat{h}_b(t, \overline{\xi}) ~ d\overline{\xi} \label{estdispsupPMFpetittau1-b-3} \\
&\quad + \sum_{l_b = k}^{-200} t^{-2} 2^{-l_b-k-4j} \int e^{i t \Phi} \Psi_{l_b}^b \psi_{j, k}^{\widehat{\mathcal{P}}}(\overline{\xi}) \widehat{X}_a(\overline{\xi}) \cdot \nabla_{\overline{\xi}} \widehat{h}_b(t, \overline{\xi}) ~ d\overline{\xi} \label{estdispsupPMFpetittau1-b-4} \\
&\quad + \sum_{l_b = k}^{-200} t^{-1} 2^{-l_b-2j} \int e^{i t \Phi} \Psi_{l_b}^b \psi_{j, k}^{\widehat{\mathcal{P}}}(\overline{\xi}) \widehat{g}_b(t, \overline{\xi}) ~ d\overline{\xi} \label{estdispsupPMFpetittau1-b-5} 
\end{align}
\end{subequations} 
We then estimate: 
\begin{align*}
\eqref{estdispsupPMFpetittau1-b-1} &\lesssim \sum_{l_b = k}^{-200} t^{-2} 2^{-2l_b-2k-6j} \Vert \Psi_{l_b}^b \psi_{j, k}^{\widehat{\mathcal{P}}} \Vert_{L^2_c L^{\frac{1}{1-\kappa}}_{a, b}} \langle 2^j \rangle^{-1} \Vert \langle \overline{\xi} \rangle \widehat{f}(t) \Vert_{L^2_c L^{\frac{1}{\kappa}}_{a, b}} \\
&\lesssim \sum_{l_b = k}^{-200} t^{-2} 2^{-l_b+\kappa l_b-k+\kappa k-\frac{7j}{2}-2\kappa j} \langle 2^j \rangle^{-1} \Vert \langle \overline{\xi} \rangle \widehat{f}(t) \Vert_{L^2_c H^1_{a, b}} \\
&\lesssim t^{-\frac{7}{6}} \mathfrak{t}^{-\frac{5}{6}} 2^{-2k-2\kappa k-j-2\kappa j} \langle 2^j \rangle^{-1} \Vert u \Vert_X \\
&\lesssim t^{-\frac{7}{6}+\frac{2\kappa+\delta}{3}} 2^{-j-\frac{k}{3}+\delta j +\delta k} \langle 2^j \rangle^{-1} \Vert u \Vert_X \\
\eqref{estdispsupPMFpetittau1-b-4} &\lesssim \sum_{l_b = k}^{-200} t^{-2} 2^{-l_b-k-5j} \Vert \Psi_{l_b}^b \psi_{j, k}^{\widehat{\mathcal{P}}} \Vert_{L^2} \Vert \nabla X_a h_b(t) \Vert_{L^2} \\
&\lesssim \sum_{l_b = k}^{-200} t^{-\frac{7}{6}} \mathfrak{t}^{-\frac{5}{6}} 2^{-\frac{l_b}{2}-\frac{k}{2}-j} \Vert u \Vert_X \\
&\lesssim t^{-\frac{7}{6}+\frac{\delta}{3}} \mathfrak{t}^{-\frac{1}{2}} 2^{-j+\delta k+\delta j} \Vert u \Vert_X \\
&\lesssim t^{-\frac{7}{6}+\frac{\delta}{3}} 2^{-j+\delta k+\delta j} \langle 2^j \rangle^{-\frac{3}{2}} \Vert u \Vert_X \\
\eqref{estdispsupPMFpetittau1-b-5} &\lesssim \sum_{l_b = k}^{-200} t^{-2} 2^{-l_b-\frac{5j}{2}} \Vert \Psi_{l_b}^b \psi_{j, k}^{\widehat{\mathcal{P}}} \Vert_{L^2} \langle 2^j \rangle^{-1} \Vert \langle \nabla \rangle |\nabla|^{\frac{1}{2}} g_b(t) \Vert_{L^2} \\
&\lesssim \sum_{l_b = k}^{-200} t^{-\frac{7}{6}+100\delta} 2^{-\frac{l_b}{2}+\frac{k}{2}-j} \langle 2^j \rangle^{-1} \Vert u \Vert_X \\
&\lesssim t^{-\frac{7}{6}+100\delta} 2^{-j} \langle 2^j \rangle^{-1} \Vert u \Vert_X 
\end{align*}
\eqref{estdispsupPMFpetittau1-b-2}, \eqref{estdispsupPMFpetittau1-b-3} can be estimated like \eqref{estdispsupPMFpetittau1-b-1}, in a simpler way. We can obtain the second estimate automatically noting that, here, $2^{-\frac{j}{2}-\frac{k}{3}} \lesssim t^{\frac{1}{6}}$. 

Then, for \eqref{estdispsupPMFpetittau-a}, we start by at most $n$ (large with respect to $\varrho$) integrations by parts along $\widehat{X}_a$, as long as none touch $\widehat{f}$, and then one integration by parts along $\widehat{X}_b$: 
\begin{subequations}
\begin{align}
\eqref{estdispsupPMFpetittau-a} &= \sum_{l_b = l_0}^k t^{-n} 2^{-n(l_b+k)-3nj} \int e^{i t \Phi} \Psi_{l_b}^b \psi_{j, k}^{\widehat{\mathcal{P}}}(\overline{\xi}) \widehat{f}(t, \overline{\xi}) ~ d\overline{\xi} \label{estdispsupPMFpetittau1-a-1} \\
&\quad + \sum_{i = 1}^n \sum_{l_b = l_0}^k t^{-i-1} 2^{-i l_b-(i+1)k-3ij-2j} \int e^{i t \Phi} \Psi_{l_b}^b \psi_{j, k}^{\widehat{\mathcal{P}}}(\overline{\xi}) \widehat{h}_a(t, \overline{\xi}) ~ d\overline{\xi} \label{estdispsupPMFpetittau1-a-3} \\
&\quad + \sum_{i = 1}^n \sum_{l_b = l_0}^k t^{-i-1} 2^{-i l_b- i k-3ij-j} \int e^{i t \Phi} \Psi_{l_b}^b \psi_{j, k}^{\widehat{\mathcal{P}}}(\overline{\xi}) \widehat{X}_b(\overline{\xi}) \cdot \nabla_{\overline{\xi}} \widehat{h}_a(t, \overline{\xi}) ~ d\overline{\xi} \label{estdispsupPMFpetittau1-a-4} \\
&\quad + \sum_{i = 1}^n \sum_{l_b = l_0}^k t^{-i} 2^{-(i-1)l_b-ik-3ij+j} \int e^{i t \Phi} \Psi_{l_b}^b \psi_{j, k}^{\widehat{\mathcal{P}}}(\overline{\xi}) \widehat{g}_a(t, \overline{\xi}) ~ d\overline{\xi} \label{estdispsupPMFpetittau1-a-5}
\end{align}
\end{subequations} 
We then estimate: 
\begin{align*}
\eqref{estdispsupPMFpetittau1-a-1} &\lesssim \sum_{l_b = l_0}^k t^{-n} 2^{-n(l_b+k)-3nj} \Vert \Psi_{l_b}^b \psi_{j, k}^{\widehat{\mathcal{P}}} \Vert_{L^2_c L^{\frac{1}{1-\kappa}}_{b, a}} \langle 2^j \rangle^{-1} \Vert \langle \overline{\xi} \rangle \widehat{f}(t) \Vert_{L^2_c L^{\frac{1}{\kappa}}_{b, a}} \\
&\lesssim \sum_{l_b = l_0}^k \mathfrak{t}^{-n} 2^{-n(l_b+k)+l_b-\kappa l_b+k-\kappa k+\frac{5j}{2}-2\kappa j} \langle 2^j \rangle^{-1} \Vert \langle \overline{\xi} \rangle \widehat{f}(t) \Vert_{L^2_c H^1_{b,a}} \\
&\lesssim t^{-\frac{7}{6}} \mathfrak{t}^{-n+\frac{7}{6}} 2^{-(3-2\varrho)nk+3k-3\kappa k-j-2\kappa j} \langle 2^j \rangle^{-1} \Vert u \Vert_X \\
&\lesssim t^{-\frac{7}{6}+\frac{2\kappa+\delta}{3}} 2^{-j+\delta k+\delta j} \langle 2^j \rangle^{-\frac{3}{4}} \Vert u \Vert_X \\
\eqref{estdispsupPMFpetittau1-a-3} &\lesssim \sum_{i = 1}^n \sum_{l_b = l_0}^k t^{-i-1} 2^{-il_b-(i+1)k-3ij-2j} \Vert \Psi_{l_b}^b \psi_{j, k}^{\widehat{\mathcal{P}}} \Vert_{L^2} \langle 2^j \rangle^{-1} \Vert \langle \nabla \rangle h_a(t) \Vert_{L^2} \\
&\lesssim \sum_{i = 1}^n \sum_{l_b = l_0}^k \mathfrak{t}^{-i-1} 2^{-il_b+\frac{l_b}{2}-ik-\frac{k}{2}+\frac{5j}{2}} \langle 2^j \rangle^{-1} \Vert u \Vert_X \\
&\lesssim \mathfrak{t}^{-2} 2^{-\frac{5k}{2}+\frac{5j}{2}} \langle 2^j \rangle^{-1} \Vert u \Vert_X \\
&\lesssim t^{-\frac{7}{6}+\frac{2\delta}{3}} 2^{-j+\delta j + \delta k} \langle 2^j \rangle^{-1} \Vert u \Vert_X \\
\eqref{estdispsupPMFpetittau1-a-4} &\lesssim \sum_{i = 1}^n \sum_{l_b = l_0}^k t^{-i-1} 2^{-i l_b-i k-3ij-j} \Vert \Psi_{l_b}^b \psi_{j, k}^{\widehat{\mathcal{P}}} \Vert_{L^2} \Vert \nabla X_b h_a(t) \Vert_{L^2} \\
&\lesssim t^{-\frac{7}{6}} \mathfrak{t}^{-\frac{1}{3}} 2^{-j} \Vert u \Vert_X \\
&\lesssim t^{-\frac{7}{6}} 2^{-j+\delta j+\delta k} \langle 2^j \rangle^{-\frac{5}{6}} \Vert u \Vert_X \\
\eqref{estdispsupPMFpetittau1-a-5} &\lesssim \sum_{i = 1}^n \sum_{l_b = l_0}^k t^{-i} 2^{-(i-1) l_b - i k-3ij+\frac{j}{2}} \Vert \Psi_{l_b}^b \psi_{j, k}^{\widehat{\mathcal{P}}} \Vert_{L^2} \langle 2^j \rangle^{-1} \Vert \langle \nabla \rangle |\nabla|^{\frac{1}{2}} g_a(t) \Vert_{L^2} \\
&\lesssim \sum_{i = 1}^n \sum_{l_b = l_0}^k \mathfrak{t}^{-i} t^{-\frac{1}{6}+100\delta} 2^{-il_b + \frac{3l_b}{2} - ik + \frac{k}{2}+2j} \langle 2^j \rangle^{-1} \Vert u \Vert_X \\
&\lesssim \mathfrak{t}^{-1} t^{-\frac{1}{6}+100\delta} 2^{2j} \langle 2^j \rangle^{-1} \Vert u \Vert_X 
+ \mathfrak{t}^{-2} t^{-\frac{1}{6}+100\delta} 2^{-\frac{3k}{2}+2j} \langle 2^j \rangle^{-1} \Vert u \Vert_X \\
&\lesssim t^{-\frac{7}{6}+100\delta} 2^{-j} \langle 2^j \rangle^{-1} \Vert u \Vert_X 
\end{align*}
Again, we can obtain the other estimate using $2^{-\frac{j}{2}} \lesssim t^{\frac{1}{6}}$. 

Finally, on the internal term, we apply integrations by parts along $\widehat{X}_a$, at most $n$, as long as none hits $\widehat{f}(t)$: 
\begin{align*}
\eqref{estdispsupPMFpetittau-int} &= t^{-n} 2^{-nl_0-nk-3nj} \int e^{i t \Phi} \Psi_{l_0}^{int} \psi_{j, k}^{\widehat{\mathcal{P}}}(\overline{\xi}) \widehat{f}(t, \overline{\xi}) ~ d\overline{\xi} \\
&\quad + \sum_{i = 1}^n t^{-i-1} 2^{-il_0-(i+1)k-3ij-2j} \int e^{i t \Phi} \Psi_{l_0}^{int} \psi_{j, k}^{\widehat{\mathcal{P}}}(\overline{\xi}) \widehat{h}_a(t, \overline{\xi}) ~ d\overline{\xi} \\
&\quad + \sum_{i = 1}^n t^{-i-1} 2^{-(i-1)l_0-(i+1)k-3ij-j} \int e^{i t \Phi} \Psi_{2k}^{int} \psi_{j, k}^{\widehat{\mathcal{P}}}(\overline{\xi}) \widehat{X}_a(\overline{\xi}) \cdot \nabla_{\overline{\xi}} \widehat{h}_a(t, \overline{\xi}) ~ d\overline{\xi} \\
&\quad + \sum_{i = 1}^n t^{-i} 2^{-il_0+l_0-ik-3ij+j} \int e^{i t \Phi} \Psi_{2k}^{int} \psi_{j, k}^{\widehat{\mathcal{P}}}(\overline{\xi}) \widehat{g}_a(t, \overline{\xi}) ~ d\overline{\xi}
\end{align*}
All these terms are simpler than those obtained after integrations by parts on \eqref{estdispsupPMFpetittau-a}, and can be estimated in a similar way. 

\paragraph{1.2.} Let us now consider the case $\mathfrak{t}^{-1} \lesssim 2^k \lesssim \mathfrak{t}^{-\frac{1}{2}}$. We then choose $l_0$ such that $2^{l_0} \simeq \mathfrak{t}^{-1+\varrho} 2^{-k}$, $\varrho > 0$ small enough with respect to $\delta$. In particular, $l_0 \geq k$. 

On \eqref{estdispsupPMFpetittau-b}, we apply at most $n$ (large with respect to $\varrho^{-1}$) integrations by parts along $\widehat{X}_b$: 
\begin{subequations}
\begin{align}
\eqref{estdispsupPMFpetittau-b} &= \sum_{l_b = l_0}^{-200} 1_{l_b \geq k} t^{-n} 2^{-nl_b-nk-3nj} \int e^{i t \Phi} \Psi_{l_b}^b \psi_{j, k}^{\widehat{\mathcal{P}}}(\overline{\xi}) \widehat{f}(t, \overline{\xi}) ~ d\overline{\xi} \label{estdispsupPMFpetittau2-b-1} \\
&\quad + \sum_{i = 1}^n \sum_{l_b = l_0}^{-200} 1_{l_b \geq k} t^{-i} 2^{-il_b-ik+k-3ij+j} \int e^{i t \Phi} \Psi_{l_b}^b \psi_{j, k}^{\widehat{\mathcal{P}}}(\overline{\xi}) \widehat{h}_b(t, \overline{\xi}) ~ d\overline{\xi} \label{estdispsupPMFpetittau2-b-2} \\
&\quad + \sum_{i = 1}^n \sum_{l_b = l_0}^{-200} 1_{l_b \geq k} t^{-1} 2^{-il_b-ik+k-3ij+j} \int e^{i t \Phi} \Psi_{l_b}^b \psi_{j, k}^{\widehat{\mathcal{P}}}(\overline{\xi}) \widehat{g}_b(t, \overline{\xi}) ~ d\overline{\xi} \label{estdispsupPMFpetittau2-b-3} 
\end{align}
\end{subequations} 
We then estimate: 
\begin{align*}
\eqref{estdispsupPMFpetittau2-b-1} &\lesssim \sum_{l_b = l_0}^{-200} 1_{l_b \geq k} t^{-n} 2^{-nl_b-nk-3nj} \Vert \Psi_{l_b}^b \psi_{j, k}^{\widehat{\mathcal{P}}} \Vert_{L^2_c L^{\frac{1}{1-\kappa}}_{b, a}} \langle 2^j \rangle^{-1} \Vert \langle \overline{\xi} \rangle \widehat{f}(t) \Vert_{L^2_c L^{\frac{1}{\kappa}}_{b, a}} \\
&\lesssim \sum_{l_b = l_0}^{-200} 1_{l_b \geq k} \mathfrak{t}^{-n} 2^{-nl_b-nk+l_b+k-\kappa l_b- \kappa k+\frac{5j}{2}-2\kappa j} \langle 2^j \rangle^{-1} \Vert u \Vert_X \\
&\lesssim t^{-\frac{7}{6}} \mathfrak{t}^{-n\varrho} 2^{-j+\frac{k}{2}} \langle 2^j \rangle^{-1} \Vert u \Vert_X \\
\eqref{estdispsupPMFpetittau2-b-2} &\lesssim \sum_{i = 1}^n \sum_{l_b = l_0}^{-200} 1_{l_b \geq k} t^{-i} 2^{-il_b-ik+k-3ij+j} \Vert \Psi_{l_b}^b \psi_{j, k}^{\widehat{\mathcal{P}}} \Vert_{L^2_{c, b} L^{\frac{1}{1-\kappa}}_a} \Vert \psi_j \widehat{h}_b(t) \Vert_{L^2_{c, b} L^{\frac{1}{\kappa}}_a} \\
&\lesssim \sum_{i = 1}^n \sum_{l_b = l_0}^{-200} 1_{l_b \geq k} \mathfrak{t}^{-i} 2^{-il_b+l_b-\kappa l_b+\frac{3k}{2}-ik+3j-\kappa j} \Vert \psi_j \widehat{h}_b(t) \Vert_{L^2}^{\frac{1}{2}} \Vert \psi_j \widehat{h}_b(t) \Vert_{L^2_{c, b} H^1_a}^{\frac{1}{2}} \\
&\lesssim \mathfrak{t}^{-1} 2^{\frac{k}{2}-\kappa k+\frac{5j}{2}-\kappa j} \langle 2^j \rangle^{-\frac{1}{2}} \Vert u \Vert_X \\
&\lesssim t^{-\frac{7}{6}+\frac{\kappa}{3}+\frac{\delta+\kappa}{3}} 2^{-j+\delta j+\delta k} \langle 2^j \rangle^{-\frac{1}{2}} \Vert u \Vert_X \\
\eqref{estdispsupPMFpetittau2-b-3} &\lesssim \sum_{i = 1}^n \sum_{l_b = l_0}^{-200} 1_{l_b \geq k} t^{-i} 2^{-il_b-ik+k-3ij+\frac{j}{2}} \Vert \Psi_{l_b}^b \psi_{j, k}^{\widehat{\mathcal{P}}} \Vert_{L^2} \langle 2^j \rangle^{-1} \Vert \langle \nabla \rangle |\nabla|^{\frac{1}{2}} g_b(t) \Vert_{L^2} \\
&\lesssim t^{-\frac{7}{6}+100\delta} 2^{-j} \langle 2^j \rangle^{-1} \Vert u \Vert_X 
\end{align*}
We obtain the other estimate the same way as before. 

Since $l_0 \geq k$, \eqref{estdispsupPMFpetittau-a} does not exists here. 

Finally, for the internal term, 
\begin{align*}
\eqref{estdispsupPMFpetittau-int} &\lesssim \Vert \Psi_{l_0}^{int} \psi_{j, k}^{\widehat{\mathcal{P}}} \Vert_{L^2_c L^{\frac{1}{1-\kappa}}_{b, a}} \langle 2^j \rangle^{-1} \Vert \langle \overline{\xi} \rangle \widehat{f}(t) \Vert_{L^2_c L^{\frac{1}{\kappa}}_{b, a}} \\
&\lesssim 2^{(1-\kappa)(l_0+k)+\frac{5j}{2}-2\kappa j} \langle 2^j \rangle^{-1} \Vert u \Vert_X \\
&\lesssim t^{-1+\kappa+\delta} 2^{2\delta j+\kappa j} 2^{-\frac{j}{2}+\delta j+\delta k} \langle 2^j \rangle^{-1} \Vert u \Vert_X \\
&\lesssim t^{-1+2\delta} 2^{-\frac{j}{2}+\delta j+\delta k} \langle 2^j \rangle^{-1+3\delta} \Vert u \Vert_X
\end{align*}
and we can obtain the other estimate using that $2^{\frac{k}{3}+\frac{j}{2}} \lesssim t^{-\frac{1}{6}}$ here. 

\paragraph{2.} Let us now assume that $\overline{\xi_a} \simeq 2^j$ and $\tau \simeq 2^k$. 

\paragraph{2.1.} Assume first that $2^k \lesssim \mathfrak{t}^{-\frac{1}{2}}$. In this case, we use the localisation symbols $\Psi_{l_b}^b, \Psi_{l_b}^{int}$ already introduced and we decompose: 
\begin{subequations}
\begin{align}
\eqref{estdispvoisP-termetot} &= \sum_{l_b = k}^{-200} \int e^{i t \Phi} \Psi_{l_b}^b \psi_{j, k}^{\widehat{\mathcal{P}}}(\overline{\xi}) \widehat{f}(t, \overline{\xi}) ~ d\overline{\xi} \label{estdispsupPMFmoyentau1-b} \\
&\quad + \int e^{i t \Phi} \Psi_{k}^{int} \psi_{j, k}^{\widehat{\mathcal{P}}}(\overline{\xi}) \widehat{f}(t, \overline{\xi}) ~ d\overline{\xi} \label{estdispsupPMFmoyentau1-int} 
\end{align}
\end{subequations} 

We have again the volume estimates
\begin{align*}
\Vert \Psi_{l_b}^{*} \psi_{j, k}^{\widehat{\mathcal{P}}} \Vert_{L^{p_c}_c L^{p_b}_b L^{p_a}_a} &\lesssim 2^{\frac{k}{p_b}+\frac{l_b}{p_a}+\frac{j}{p_a}+\frac{j}{p_b}+\frac{j}{p_c}} 
\end{align*}
and 
\begin{align*}
\widehat{X}_b(\overline{\xi}) \cdot \nabla_{\overline{\xi}} \left[ \widehat{X}_b(\overline{\xi}) \cdot \nabla_{\overline{\xi}} \Phi \right] &\lesssim 2^{k+j} 
\end{align*}
Therefore, 
\begin{align*}
\widehat{X}_b \cdot \nabla \left[ \Psi_{l_b}^b \psi_{j, k}^{\widehat{\mathcal{P}}} \right] &\lesssim 2^{-j-k} + 2^{-l_b+k-j} \lesssim 2^{-j-k} 
\end{align*}
using that $l_b \geq k$. 

On \eqref{estdispsupPMFmoyentau1-b}, we apply an integration by parts: 
\begin{subequations} 
\begin{align}
\eqref{estdispsupPMFmoyentau1-b} &= \sum_{l_b = k}^{-200} t^{-1} 2^{-l_b-k-3j} \int e^{i t \Phi} \Psi_{l_b}^b \psi_{j, k}^{\widehat{\mathcal{P}}}(\overline{\xi}) \widehat{f}(t, \overline{\xi}) ~ d\overline{\xi} \label{estdispsupPMFmoyentau1-b-1} \\
&\quad + \sum_{l_b = k}^{-200} t^{-1} 2^{-l_b-2j} \int e^{i t \Phi} \Psi_{l_b}^b \psi_{j, k}^{\widehat{\mathcal{P}}}(\overline{\xi}) \widehat{h}_b(t, \overline{\xi}) ~ d\overline{\xi} \label{estdispsupPMFmoyentau1-b-2} \\
&\quad + \sum_{l_b = k}^{-200} t^{-1} 2^{-l_b-2j} \int e^{i t \Phi} \Psi_{l_b}^b \psi_{j, k}^{\widehat{\mathcal{P}}}(\overline{\xi}) \widehat{g}_b(t, \overline{\xi}) ~ d\overline{\xi} \label{estdispsupPMFmoyentau1-b-3} 
\end{align}
\end{subequations} 
We then estimate: 
\begin{align*}
\eqref{estdispsupPMFmoyentau1-b-1} &\lesssim \sum_{l_b = k}^{-200} t^{-1} 2^{-l_b-k-3j} \Vert \Psi_{l_b}^b \psi_{j, k}^{\widehat{\mathcal{P}}} \Vert_{L^2_c L^{\frac{1}{1-\kappa}}_{a, b}} \langle 2^j \rangle^{-1} \Vert \langle \overline{\xi} \rangle \widehat{f}(t) \Vert_{L^2_c L^{\frac{1}{\kappa}}_{a, b}} \\
&\lesssim \sum_{l_b = k}^{-200} t^{-1} 2^{-\kappa l_b-\kappa k-\frac{j}{2}} \langle 2^j \rangle^{-1} \Vert u \Vert_X \\
&\lesssim t^{-1+\kappa+\frac{\delta}{2}} 2^{3\kappa j+\frac{\delta j}{2}} 2^{-\frac{j}{2}+\delta j+\delta k} \langle 2^j \rangle^{-1} \Vert u \Vert_X \\
&\lesssim t^{-1+\delta} 2^{-\frac{j}{2}+\delta j+\delta k} \langle 2^j \rangle^{-\frac{5}{6}} \Vert u \Vert_X \\
\eqref{estdispsupPMFmoyentau1-b-2} &\lesssim \sum_{l_b = k}^{-200} t^{-1} 2^{-l_b-2j} \Vert \Psi_{l_b}^b \psi_{j, k}^{\widehat{\mathcal{P}}} \Vert_{L^2_{c, b} L^{\frac{1}{1-\kappa}}_a} \Vert \psi_j \widehat{h}_b(t) \Vert_{L^2_{c, b} L^{\frac{1}{\kappa}}_a} \\
&\lesssim \sum_{l_b = k}^{-200} t^{-1} 2^{\frac{k}{2}-\kappa l_b-\kappa j} \Vert \psi_j \widehat{h}_b(t) \Vert_{L^2_{c, b} H^1_a}^{\frac{1}{2}} \Vert \psi_j \widehat{h}_b(t) \Vert_{L^2}^{\frac{1}{2}} \\
&\lesssim t^{-1} 2^{\frac{k}{2}-\kappa k-\frac{j}{2}-\kappa j} \langle 2^j \rangle^{-1} u \Vert_X \\
&\lesssim t^{-\frac{7}{6}+\frac{\kappa}{3}+\frac{\delta}{3}} \mathfrak{t}^{-\frac{1}{12}+\frac{\kappa}{6}+\frac{\delta}{6}} 2^{-j+\delta j+\delta k} \langle 2^j \rangle^{-1} u \Vert_X \\
&\lesssim t^{-\frac{7}{6}+\delta} 2^{-j+\delta j+\delta k} \langle 2^j \rangle^{-\frac{6}{5}} \Vert u \Vert_X \\
\eqref{estdispsupPMFmoyentau1-b-3} &\lesssim \sum_{l_b = k}^{-200} t^{-1} 2^{-l_b-\frac{5j}{2}} \Vert \Psi_{l_b}^b \psi_{j, k}^{\widehat{\mathcal{P}}} \Vert_{L^2} \langle 2^j \rangle^{-1} \Vert \langle \nabla \rangle |\nabla|^{\frac{1}{2}} g_b(t) \Vert_{L^2} \\
&\lesssim \sum_{l_b = k}^{-200} t^{-\frac{7}{6}+100\delta} 2^{-\frac{l_b}{2}+\frac{k}{2}-j} \langle 2^j \rangle^{-1} \Vert u \Vert_X \\
&\lesssim t^{-\frac{7}{6}+100\delta} 2^{-j} \langle 2^j \rangle^{-1} \Vert u \Vert_X 
\end{align*}
if $\kappa, \delta$ are small enough. We obtain the other estimate using $1 \lesssim t^{-\frac{1}{6}} 2^{-\frac{k}{3}-\frac{j}{2}}$ et $2^{-\frac{j}{2}} \lesssim t^{\frac{1}{6}}$, and for \eqref{estdispsupPMFmoyentau1-b-3}, we have also
\begin{align*}
\eqref{estdispsupPMFmoyentau1-b-3} &\lesssim t^{-1} 2^{-j+\delta k+\delta j} \langle 2^j \rangle^{-\frac{5}{6}} \Vert u \Vert_X
\end{align*}

For the internal term, 
\begin{align*}
\eqref{estdispsupPMFmoyentau1-int} &\lesssim \Vert \Psi_{k}^{int} \psi_{j, k}^{\widehat{\mathcal{P}}} \Vert_{L^2_c L^{\frac{1}{1-\kappa}}_{a, b}} \langle 2^j \rangle^{-1} \Vert \langle \overline{\xi} \rangle \widehat{f}(t) \Vert_{L^2_c L^{\frac{1}{\kappa}}_{a, b}} \\
&\lesssim 2^{2k-2\kappa k+\frac{5j}{2}-2\kappa j} \langle 2^j \rangle^{-1} \Vert u \Vert_X \\
&\lesssim t^{-1+\kappa+\frac{\delta}{3}} 2^{-\frac{j}{2}+\delta j+\delta k} \langle 2^j \rangle^{-\frac{5}{6}} \Vert u \Vert_X 
\end{align*}
and we get the other estimate using $1 \lesssim t^{-\frac{1}{6}} 2^{-\frac{k}{3}-\frac{j}{2}}$. 

\paragraph{2.2.} Assume now that $2^k \gtrsim \mathfrak{t}^{-\frac{1}{2}}$. 

We now introduce localisation symbols $\Psi_{(l_a, l_b, l_c)}^{*}$ for $* \in \{ a-b, c-b, a-int, c-int, b, int \}$ as for previous geometric areas, and then we decompose: 
\begin{subequations}
\begin{align}
\eqref{estdispvoisP-termetot} &= \sum_{l_a = l_0}^{-200} \sum_{l_b = l_a}^{-200} \int e^{i t \Phi} \Psi_{(l_a, l_b, l_a)}^{a-b} \psi_{j, k}^{\widehat{\mathcal{P}}}(\overline{\xi}) \widehat{f}(t, \overline{\xi}) ~ d\overline{\xi} \label{estdispsupPMFmoyentau-ba} \\
&\quad + \sum_{l_a = l_0}^{-200} \sum_{l_b = l_a}^{-200} \int e^{i t \Phi} \Psi_{(l_a, l_b, l_a)}^{c-b} \psi_{j, k}^{\widehat{\mathcal{P}}}(\overline{\xi}) \widehat{f}(t, \overline{\xi}) ~ d\overline{\xi} \label{estdispsupPMFmoyentau-bc} \\
&\quad + \sum_{l_a = l_0}^{-200} \int e^{i t \Phi} \Psi_{(l_a, l_a, l_a)}^{a-int} \psi_{j, k}^{\widehat{\mathcal{P}}}(\overline{\xi}) \widehat{f}(t, \overline{\xi}) ~ d\overline{\xi} \label{estdispsupPMFmoyentau-a2} \\
&\quad + \sum_{l_a = l_0}^{-200} \int e^{i t \Phi} \Psi_{(l_a, l_a, l_a)}^{c-int} \psi_{j, k}^{\widehat{\mathcal{P}}}(\overline{\xi}) \widehat{f}(t, \overline{\xi}) ~ d\overline{\xi} \label{estdispsupPMFmoyentau-c2} \\
&\quad + \sum_{l_b = l_0+k+\frac{j_{+}}{3}}^{-200} \int e^{i t \Phi} \Psi_{(l_0, l_b, l_0)}^{b} \psi_{j, k}^{\widehat{\mathcal{P}}}(\overline{\xi}) \widehat{f}(t, \overline{\xi}) ~ d\overline{\xi} \label{estdispsupPMFmoyentau-b} \\
&\quad + \int e^{i t \Phi} \Psi_{(l_0, l_0+k+\frac{j_{+}}{3}, l_0)}^{int} \psi_{j, k}^{\widehat{\mathcal{P}}}(\overline{\xi}) \widehat{f}(t, \overline{\xi}) ~ d\overline{\xi} \label{estdispsupPMFmoyentau-int} 
\end{align}
\end{subequations}
for $l_0$ such that $2^{l_0} \simeq \mathfrak{t}^{-\frac{1}{2}}$. In particular, $2^{l_0} \lesssim 2^k$. Here above, as before, $j_{+}$ stands for the positive part of $j$. 

Note that
\begin{align*}
\widehat{X}_c(\overline{\xi}) \cdot \nabla_{\overline{\xi}} \Phi &= \frac{\xi \cdot Jy}{t |\xi|} \lesssim 2^{k+2j} \\
\widehat{X}_a(\overline{\xi}) \cdot \nabla_{\overline{\xi}} \Phi &= 3 \xi_0 |\overline{\xi}| \left( 1 - \frac{x}{3t |\overline{\xi}|^2} \right) - \frac{\xi \cdot y}{t |\overline{\xi}|} \lesssim 2^{k+2j} 
\end{align*}
so that the sums in $l_a$ are empty as long as $l_a \geq k+10$. 

Again, in local coordinates $(\xi_a, \xi_b, \xi_c)$, the symbol $\Psi_{(l_a, l_b, l_c)}^{*}$ localises $\xi_c$ close to $0$ with precision $2^{l_a-k}$, then $\xi_c$ being fixed it localises $\xi_b$ close to $2^{k+j}$ with precision $2^{l_a+j}$, and finally $\xi_c, \xi_b$ being fixed it localises $\xi_a$ close to $2^j$ with precision $2^{l_b+j}$. Hence, 
\begin{align*}
\Vert \Psi_{(l_a, l_b, l_c)}^{*} \psi_{j, k}^{\widehat{\mathcal{P}}} \Vert_{L^{p_c}_c L^{p_b}_b L^{p_a}_a} &\lesssim 2^{\frac{l_c-k}{p_c}+\frac{l_a}{p_b}+\frac{l_b}{p_a}+\frac{j}{p_a}+\frac{j}{p_b}+\frac{j}{p_c}}
\end{align*}

We can compute the symbol derivatives: 
\begin{align*}
\widehat{X}_a(\overline{\xi}) \cdot \nabla_{\overline{\xi}} \left[ \widehat{X}_a(\overline{\xi}) \cdot \nabla_{\overline{\xi}} \Phi \right] &\simeq \xi_0 \lesssim 2^{j+k} \\
\widehat{X}_a(\overline{\xi}) \cdot \nabla_{\overline{\xi}} \left[ \widehat{X}_b(\overline{\xi}) \cdot \nabla_{\overline{\xi}} \Phi \right] &\simeq 2^j \\
\widehat{X}_a(\overline{\xi}) \cdot \nabla_{\overline{\xi}} \left[ \widehat{X}_c(\overline{\xi}) \cdot \nabla_{\overline{\xi}} \Phi \right] &= 0 \\
\widehat{X}_b(\overline{\xi}) \cdot \nabla_{\overline{\xi}} \left[ \widehat{X}_a(\overline{\xi}) \cdot \nabla_{\overline{\xi}} \Phi \right] &\simeq 2^j \\
\widehat{X}_b(\overline{\xi}) \cdot \nabla_{\overline{\xi}} \left[ \widehat{X}_b(\overline{\xi}) \cdot \nabla_{\overline{\xi}} \Phi \right] &\lesssim 2^{k+j} \\
\widehat{X}_b(\overline{\xi}) \cdot \nabla_{\overline{\xi}} \left[ \widehat{X}_c(\overline{\xi}) \cdot \nabla_{\overline{\xi}} \Phi \right] &= 0 \\
\widehat{X}_c(\overline{\xi}) \cdot \nabla_{\overline{\xi}} \left[ \widehat{X}_a(\overline{\xi}) \cdot \nabla_{\overline{\xi}} \Phi \right] 
&\lesssim 2^{l_a+j} \\
\widehat{X}_c(\overline{\xi}) \cdot \nabla_{\overline{\xi}} \left[ \widehat{X}_b(\overline{\xi}) \cdot \nabla_{\overline{\xi}} \Phi \right] 
&\lesssim 2^{l_a+k+j} \\
\widehat{X}_c(\overline{\xi}) \cdot \nabla_{\overline{\xi}} \left[ \widehat{X}_c(\overline{\xi}) \cdot \nabla_{\overline{\xi}} \Phi \right] 
&\lesssim 2^{j+k} 
\end{align*}
and
\begin{align*}
\widehat{X}_a \cdot \nabla \psi_{j, k}^{\widehat{\mathcal{P}}} &\lesssim 2^{-j} \\
\widehat{X}_b \cdot \nabla \psi_{j, k}^{\widehat{\mathcal{P}}} &\lesssim 2^{-j-k} \\
\widehat{X}_c \cdot \nabla \psi_{j, k}^{\widehat{\mathcal{P}}} &= 0
\end{align*}
Therefore, 
\begin{align*}
\widehat{X}_a(\overline{\xi}) \cdot \nabla_{\overline{\xi}} \left[ \Psi_{(l_a, l_b, l_c)} \psi_{j, k}^{\widehat{\mathcal{P}}} \right] &\lesssim 2^{-j} \left( 2^{-l_a+k} + 2^{-l_b} \right) \\
\widehat{X}_b(\overline{\xi}) \cdot \nabla_{\overline{\xi}} \left[ \Psi_{(l_a, l_b, l_c)} \psi_{j, k}^{\widehat{\mathcal{P}}} \right] &\lesssim 2^{-j} \left( 2^{-l_a} + 2^{-l_b+k} \right) \\
\widehat{X}_c(\overline{\xi}) \cdot \nabla_{\overline{\xi}} \left[ \Psi_{(l_a, l_b, l_c)} \psi_{j, k}^{\widehat{\mathcal{P}}} \right] &\lesssim 2^{-j} \left( 2^{l_a+k-l_b} + 2^{-l_a+k} \right) \lesssim 2^{-j} \left( 2^{-l_a+k} + 2^{-l_b} \right) 
\end{align*}

On \eqref{estdispsupPMFmoyentau-ba}, we apply an integration by parts along $\widehat{X}_b$ and then one along $\widehat{X}_a$: 
\begin{subequations}
\begin{align}
\eqref{estdispsupPMFmoyentau-ba} &= \sum_{l_a = l_0}^{k+10} \sum_{l_b = l_a}^{-200} t^{-2} 2^{-2l_a-l_b-6j} \left( 2^{-l_a+k} + 2^{-l_b} \right) \int e^{i t \Phi} \Psi_{(l_a, l_b, l_a)}^{a-b} \psi_{j, k}^{\widehat{\mathcal{P}}}(\overline{\xi}) \widehat{f}(t, \overline{\xi}) ~ d\overline{\xi} \label{estdispsupPMFmoyentau-ba-1} \\
&\quad + \sum_{l_a = l_0}^{k+10} \sum_{l_b = l_a}^{-200} t^{-2} 2^{-2l_a-l_b-5j} \int e^{i t \Phi} \Psi_{(l_a, l_b, l_a)}^{a-b} \psi_{j, k}^{\widehat{\mathcal{P}}}(\overline{\xi}) \widehat{X}_a(\overline{\xi}) \cdot \nabla_{\overline{\xi}} \widehat{f}(t, \overline{\xi}) ~ d\overline{\xi} \label{estdispsupPMFmoyentau-ba-2} \\
&\quad + \sum_{l_a = l_0}^{k+10} \sum_{l_b = l_a}^{-200} t^{-2} 2^{-l_a-l_b-5j} \left( 2^{-l_a+k} + 2^{-l_b} \right) \int e^{i t \Phi} \Psi_{(l_a, l_b, l_a)}^{a-b} \psi_{j, k}^{\widehat{\mathcal{P}}}(\overline{\xi}) \widehat{h}_b(t, \overline{\xi}) ~ d\overline{\xi} \label{estdispsupPMFmoyentau-ba-3} \\
&\quad + \sum_{l_a = l_0}^{k+10} \sum_{l_b = l_a}^{-200} t^{-2} 2^{-l_a-l_b-4j} \int e^{i t \Phi} \Psi_{(l_a, l_b, l_a)}^{a-b} \psi_{j, k}^{\widehat{\mathcal{P}}}(\overline{\xi}) \widehat{X}_a(\overline{\xi}) \cdot \nabla_{\overline{\xi}} \widehat{h}_b(t, \overline{\xi}) ~ d\overline{\xi} \label{estdispsupPMFmoyentau-ba-4} \\
&\quad + \sum_{l_a = l_0}^{k+10} \sum_{l_b = l_a}^{-200} t^{-1} 2^{-l_b-2j} \int e^{i t \Phi} \Psi_{(l_a, l_b, l_a)}^{a-b} \psi_{j, k}^{\widehat{\mathcal{P}}}(\overline{\xi}) \widehat{g}_b(t, \overline{\xi}) ~ d\overline{\xi} \label{estdispsupPMFmoyentau-ba-5} 
\end{align}
\end{subequations} 
We then estimate: 
\begin{align*}
\eqref{estdispsupPMFmoyentau-ba-1} &\lesssim \sum_{l_a = l_0}^{k+10} \sum_{l_b = l_a}^{-200} t^{-2} 2^{-2l_a-l_b-6j} \left( 2^{-l_a+k} + 2^{-l_b} \right) \Vert \Psi_{(l_a, l_b, l_a)}^{a-b} \psi_{j, k}^{\widehat{\mathcal{P}}} \Vert_{L^2_c L^{\frac{1}{1-\delta}}_{b, a}} \langle 2^j \rangle^{-1} \Vert \langle \overline{\xi} \rangle \widehat{f}(t) \Vert_{L^2_c L^{\frac{1}{\delta}}_{b, a}} \\
&\lesssim \sum_{l_a = l_0}^{k+10} \sum_{l_b = l_a}^{-200} t^{-\frac{7}{6}} \mathfrak{t}^{-\frac{5}{6}} 2^{-\frac{l_a}{2}-\delta l_a-\delta l_b-\frac{k}{2}-j-2\delta j} \left( 2^{-l_a+k} + 2^{-l_b} \right) \langle 2^j \rangle^{-1} \Vert u \Vert_X \\
&\lesssim \sum_{l_a = l_0}^{k+10} t^{-\frac{7}{6}+\delta} \mathfrak{t}^{-\frac{5}{6}-\delta} 2^{-\frac{3l_a}{2}-2 \delta l_a-\delta k} 2^{-j-\frac{k}{2}+\delta j+\delta k} \langle 2^j \rangle^{-1} \Vert u \Vert_X \\
&\lesssim t^{-\frac{7}{6}+\delta} \mathfrak{t}^{-\frac{1}{12}+\frac{\delta}{2}} 2^{-j-\frac{k}{2}+\delta j+\delta k} \langle 2^j \rangle^{-1} \Vert u \Vert_X \\
&\lesssim t^{-\frac{7}{6}+\delta+\frac{\delta}{2}} 2^{-j-\frac{k}{3}+\delta j+\delta k} \langle 2^j \rangle^{-\frac{5}{6}} \Vert u \Vert_X \\
\eqref{estdispsupPMFmoyentau-ba-2} &\lesssim \sum_{l_a = l_0}^{k+10} \sum_{l_b = l_a}^{-200} t^{-2} 2^{-2l_a-l_b-5j} \Vert \Psi_{(l_a, l_b, l_a)}^{a-b} \psi_{j, k}^{\widehat{\mathcal{P}}} \Vert_{L^2} \langle 2^j \rangle^{-1} \Vert \langle \nabla \rangle X_a f(t) \Vert_{L^2} \\
&\lesssim \sum_{l_a = l_0}^{k+10} \sum_{l_b = l_a}^{-200} t^{-\frac{7}{6}} \mathfrak{t}^{-\frac{5}{6}} 2^{-l_a-\frac{l_b}{2}-\frac{k}{2}-j} \langle 2^j \rangle^{-1} \Vert u \Vert_X \\
&\lesssim t^{-\frac{7}{6}+\frac{\delta}{3}} \mathfrak{t}^{-\frac{1}{12}+\frac{\delta}{6}} 2^{-j-\frac{k}{2}+\delta j+\delta k} \langle 2^j \rangle^{-1} \Vert u \Vert_X \\
&\lesssim t^{-\frac{7}{6}+\frac{\delta}{3}} 2^{-j-\frac{k}{2}+\delta j+\delta k} \langle 2^j \rangle^{-1} \Vert u \Vert_X \\
\eqref{estdispsupPMFmoyentau-ba-4} &\lesssim \sum_{l_a = l_0}^{k+10} \sum_{l_b = l_a}^{-200} t^{-2} 2^{-l_a-l_b-5j} \Vert \Psi_{(l_a, l_b, l_a)}^{a-b} \Vert_{L^2} \Vert \nabla X_a h_b(t) \Vert_{L^2} \\
&\lesssim \sum_{l_a = l_0}^{k+10} \sum_{l_b = l_a}^{-200} t^{-\frac{7}{6}} \mathfrak{t}^{-\frac{5}{6}} 2^{-\frac{k}{2}-\frac{l_b}{2}-j} \Vert u \Vert_X \\
&\lesssim t^{-\frac{7}{6}} \mathfrak{t}^{-\frac{7}{12}+\frac{\delta}{2}} 2^{-\delta j} 2^{-j-\frac{k}{2}+\delta j+\delta k} \Vert u \Vert_X \\
&\lesssim t^{-\frac{7}{6}+\frac{\delta}{3}} 2^{-j-\frac{k}{2}+\delta j+\delta k} \langle 2^j \rangle^{-1} \Vert u \Vert_X \\
\eqref{estdispsupPMFmoyentau-ba-5} &\lesssim \sum_{l_a = l_0}^{k+10} \sum_{l_b = l_a}^{-200} t^{-1} 2^{-l_b-\frac{5j}{2}} \Vert \Psi_{(l_a, l_b, l_a)}^{a-b} \psi_{j, k}^{\widehat{\mathcal{P}}} \Vert_{L^2} \langle 2^j \rangle^{-1} \Vert \langle \nabla \rangle |\nabla|^{\frac{1}{2}} g_b(t) \Vert_{L^2} \\
&\lesssim \sum_{l_a = l_0}^{k+10} \sum_{l_b = l_a}^{-200} t^{-\frac{7}{6}+100\delta} 2^{l_a-\frac{l_b}{2}-\frac{k}{2}-j} \langle 2^j \rangle^{-1} \Vert u \Vert_X \\
&\lesssim t^{-\frac{7}{6}+100\delta} 2^{-j} \langle 2^j \rangle^{-1} \Vert u \Vert_X 
\end{align*}
\eqref{estdispsupPMFmoyentau-ba-3} can be estimated like \eqref{estdispsupPMFmoyentau-ba-2}. We obtain the other estimate using $1 \lesssim t^{\frac{1}{6}} 2^{\frac{k}{3}+\frac{j}{2}}$. 

For \eqref{estdispsupPMFmoyentau-a2}, we apply two integrations by parts along $\widehat{X}_a$: 
\begin{subequations}
\begin{align}
\eqref{estdispsupPMFmoyentau-a2} &= \sum_{l_a = l_0}^{k+10} t^{-2} 2^{-4l_a-6j} \int e^{i t \Phi} \Psi_{(l_a, l_a, l_a)}^{a-int} \psi_{j, k}^{\widehat{\mathcal{P}}}(\overline{\xi}) \widehat{f}(t, \overline{\xi}) ~ d\overline{\xi} \label{estdispsupPMFmoyentau-a2-1} \\
&\quad + \sum_{l_a = l_0}^{k+10} t^{-2} 2^{-3l_a-5j} \int e^{i t \Phi} \Psi_{(l_a, l_a, l_a)}^{a-int} \psi_{j, k}^{\widehat{\mathcal{P}}}(\overline{\xi}) \widehat{X}_a(\overline{\xi}) \cdot \nabla_{\overline{\xi}} \widehat{f}(t, \overline{\xi}) ~ d\overline{\xi} \label{estdispsupPMFmoyentau-a2-2} \\
&\quad + \sum_{l_a = l_0}^{k+10} t^{-2} 2^{-3l_a-5j} \int e^{i t \Phi} \Psi_{(l_a, l_a, l_a)}^{a-int} \psi_{j, k}^{\widehat{\mathcal{P}}}(\overline{\xi}) \widehat{h}_a(t, \overline{\xi}) ~ d\overline{\xi} \label{estdispsupPMFmoyentau-a2-3} \\
&\quad + \sum_{l_a = l_0}^{k+10} t^{-2} 2^{-2l_a-4j} \int e^{i t \Phi} \Psi_{(l_a, l_a, l_a)}^{a-int} \psi_{j, k}^{\widehat{\mathcal{P}}}(\overline{\xi}) \widehat{X}_a(\overline{\xi}) \cdot \nabla_{\overline{\xi}} \widehat{h}_a(t, \overline{\xi}) ~ d\overline{\xi} \label{estdispsupPMFmoyentau-a2-4} \\
&\quad + \sum_{l_a = l_0}^{k+10} t^{-1} 2^{-l_a-2j} \int e^{i t \Phi} \Psi_{(l_a, l_a, l_a)}^{a-int} \psi_{j, k}^{\widehat{\mathcal{P}}}(\overline{\xi}) \widehat{g}_b(t, \overline{\xi}) ~ d\overline{\xi} \label{estdispsupPMFmoyentau-a2-5}
\end{align}
\end{subequations}
We then estimate:  
\begin{align*}
\eqref{estdispsupPMFmoyentau-a2-1} &\lesssim \sum_{l_a = l_0}^{k+10} t^{-2} 2^{-4l_a-6j} \Vert \Psi_{(l_a, l_a, l_a)}^{a-int} \psi_{j, k}^{\widehat{\mathcal{P}}} \Vert_{L^2_c L^{\frac{1}{1-\kappa}}_{b, a}} \langle 2^j \rangle^{-1} \Vert \langle \overline{\xi} \rangle \widehat{f}(t) \Vert_{L^2_c L^{\frac{1}{\kappa}}_{b, a}} \\
&\lesssim \sum_{l_a = l_0}^{k+10} t^{-\frac{7}{6}} \mathfrak{t}^{-\frac{5}{6}} 2^{-\frac{3l_a}{2}-2\kappa l_a-\frac{k}{2}-j-2\kappa j} \langle 2^j \rangle^{-1} \Vert u \Vert_X \\
&\lesssim t^{-\frac{7}{6}+\frac{\delta}{3}+\frac{2\kappa}{3}} \mathfrak{t}^{-\frac{1}{12}+\frac{\kappa}{3}+\frac{\delta}{6}} 2^{-j-\frac{k}{2}+\delta j+\delta k} \langle 2^j \rangle^{-1} \Vert u \Vert_X \\
&\lesssim t^{-\frac{7}{6}+\delta} 2^{-j-\frac{k}{3}+\delta j+\delta k} \langle 2^j \rangle^{-\frac{5}{6}} \Vert u \Vert_X \\
\eqref{estdispsupPMFmoyentau-a2-2} &\lesssim \sum_{l_a = l_0}^{k+10} t^{-2} 2^{-3l_a-5j} \Vert \Psi_{(l_a, l_a, l_a)}^{a-int} \psi_{j, k}^{\widehat{\mathcal{P}}} \Vert_{L^2} \langle 2^j \rangle^{-1} \Vert \langle \nabla \rangle X_a f(t) \Vert_{L^2} \\
&\lesssim \sum_{l_a = l_0}^{k+10} t^{-\frac{7}{6}} \mathfrak{t}^{-\frac{5}{6}} 2^{-\frac{3l_a}{2}-\frac{k}{2}-j} \langle 2^j \rangle^{-1} \Vert u \Vert_X \\
&\lesssim t^{-\frac{7}{6}+\frac{\delta}{2}} 2^{-j-\frac{k}{3}+\delta j+\delta k} \langle 2^j \rangle^{-\frac{5}{6}} \Vert u \Vert_X \\
\eqref{estdispsupPMFmoyentau-a2-4} &\lesssim \sum_{l_a = l_0}^{k+10} t^{-2} 2^{-2l_a-5j} \Vert \Psi_{(l_a, l_a, l_a)}^{a-int} \psi_{j, k}^{\widehat{\mathcal{P}}} \Vert_{L^2} \Vert \nabla X_a h_a(t) \Vert_{L^2} \\
&\lesssim \sum_{l_a = l_0}^{k+10} t^{-\frac{7}{6}} \mathfrak{t}^{-\frac{5}{6}} 2^{-\frac{l_a}{2}-\frac{k}{2}-j} \Vert u \Vert_X \\
&\lesssim t^{-\frac{7}{6}+\frac{\delta}{3}} \mathfrak{t}^{-\frac{7}{12}+\frac{\delta}{6}} 2^{-j-\frac{k}{2}+\delta j+\delta k} \Vert u \Vert_X \\
&\lesssim t^{-\frac{7}{6}+\frac{\delta}{3}} 2^{-j-\frac{k}{3}+\delta j+\delta k} \langle 2^j \rangle^{-1} \Vert u \Vert_X \\
\eqref{estdispsupPMFmoyentau-a2-5} &\lesssim \sum_{l_a = l_0}^{k+10} t^{-1} 2^{-l_a-\frac{5j}{2}} \Vert \Psi_{(l_a, l_a, l_a)}^{a-int} \psi_{j, k}^{\widehat{\mathcal{P}}} \Vert_{L^2} \langle 2^j \rangle^{-1} \Vert \langle \nabla \rangle |\nabla|^{\frac{1}{2}} g_b(t) \Vert_{L^2} \\
&\lesssim \sum_{l_a = l_0}^{k+10} t^{-\frac{7}{6}+100\delta} 2^{\frac{l_a}{2}-\frac{k}{2}-j} \langle 2^j \rangle^{-1} \Vert u \Vert_X \\
&\lesssim t^{-\frac{7}{6}+100\delta} 2^{-j} \langle 2^j \rangle^{-1} \Vert u \Vert_X
\end{align*}
\eqref{estdispsupPMFmoyentau-a2-3} can be estimated like \eqref{estdispsupPMFmoyentau-a2-2}, and the other estimate is obtained as before. 

We can estimate \eqref{estdispsupPMFmoyentau-bc} and \eqref{estdispsupPMFmoyentau-c2} the same way replacing $\widehat{X}_a$ by $\widehat{X}_c$. 

On \eqref{estdispsupPMFmoyentau-b}, we apply an integration by parts along $\widehat{X}_b$: 
\begin{subequations}
\begin{align}
\eqref{estdispsupPMFmoyentau-b} &= \sum_{l_b = l_0+k+\frac{j_{+}}{3}}^{-200} t^{-1} 2^{-l_0-l_b-3j} \int e^{i t \Phi} \Psi_{(l_0, l_b, l_0)}^b \psi_{j, k}^{\widehat{\mathcal{P}}}(\overline{\xi}) \widehat{f}(t, \overline{\xi}) ~ d\overline{\xi} \label{estdispsupPMFmoyentau-b-1} \\
&\quad + \sum_{l_b = l_0+k+\frac{j_{+}}{3}}^{-200} t^{-1} 2^{-l_b-2j} \int e^{i t \Phi} \Psi_{(l_0, l_b, l_0)}^b \psi_{j, k}^{\widehat{\mathcal{P}}}(\overline{\xi}) \widehat{h}_b(t, \overline{\xi}) ~ d\overline{\xi} \label{estdispsupPMFmoyentau-b-2} \\
&\quad + \sum_{l_b = l_0+k+\frac{j_{+}}{3}}^{-200} t^{-1} 2^{-l_b-2j} \int e^{i t \Phi} \Psi_{(l_0, l_b, l_0)}^b \psi_{j, k}^{\widehat{\mathcal{P}}}(\overline{\xi}) \widehat{g}_b(t, \overline{\xi}) ~ d\overline{\xi} \label{estdispsupPMFmoyentau-b-3}
\end{align}
\end{subequations} 
We then estimate: 
\begin{align*}
\eqref{estdispsupPMFmoyentau-b-1} &\lesssim \sum_{l_b = l_0+k+\frac{j_{+}}{3}}^{-200} t^{-1} 2^{-l_0-l_b-3j} \Vert \Psi_{(l_0, l_b, l_0)}^b \psi_{j, k}^{\widehat{\mathcal{P}}} \Vert_{L^2_c L^{\frac{1}{1-\kappa}}_{b, a}} \langle 2^j \rangle^{-1} \Vert \langle \overline{\xi} \rangle \widehat{f}(t) \Vert_{L^2_c L^{\frac{1}{\kappa}}_{b, a}} \\
&\lesssim \sum_{l_b = l_0+k}^{-200} t^{-1} 2^{\frac{l_0}{2}-\kappa l_0-\frac{k}{2}-\kappa l_b-\frac{j}{2}-2\kappa j} \langle 2^j \rangle^{-1} \Vert u \Vert_X \\
&\lesssim t^{-\frac{7}{6}} \mathfrak{t}^{\frac{1}{6}} 2^{\frac{l_0}{2}-2\kappa l_0-\kappa k-\delta k-\delta j-2\kappa j} 2^{-j-\frac{k}{2}+\delta j+\delta k} \langle 2^j \rangle^{-1} \Vert u \Vert_X \\
&\lesssim t^{-\frac{7}{6}+\frac{\delta}{3}+\frac{2\kappa}{3}} \mathfrak{t}^{\frac{5\kappa}{6}+\frac{\delta}{6}} 2^{-2\kappa j} 2^{-j-\frac{k}{3}+\delta j+\delta k} \langle 2^j \rangle^{-1} \Vert u \Vert_X \\
&\lesssim t^{-\frac{7}{6}+\delta} 2^{-2\kappa j} 2^{-j-\frac{k}{3}+\delta j+\delta k} \langle 2^j \rangle^{-\frac{5}{6}} \Vert u \Vert_X \\
\eqref{estdispsupPMFmoyentau-b-2} &\lesssim \sum_{l_b = l_0+k+\frac{j_{+}}{3}}^{-200} t^{-1} 2^{-l_b-2j} \Vert \Psi_{(l_0, l_b, l_0)}^b \psi_{j, k}^{\widehat{\mathcal{P}}} \Vert_{L^2_{c, b} L^{\frac{1}{1-\kappa}}_a} \Vert \psi_j \widehat{h}_b(t) \Vert_{L^2_{c, b} L^{\frac{1}{\kappa}}_a} \\
&\lesssim \sum_{l_b = l_0+k+\frac{j_{+}}{3}}^{-200} t^{-\frac{7}{6}} \mathfrak{t}^{\frac{1}{6}} 2^{\frac{3l_0}{2}-\frac{k}{2}-\frac{l_b}{2}-\kappa l_b-\frac{j}{2}-\kappa j} \Vert \psi_j \widehat{h}_b(t) \Vert_{L^2_{c, b} H^1_a}^{\frac{1}{2}} \Vert \psi_j \widehat{h}_b(t) \Vert_{L^2}^{\frac{1}{2}} \\
&\lesssim t^{-\frac{7}{6}} \mathfrak{t}^{\frac{1}{6}} 2^{l_0-k-\kappa l_b-j-\kappa j} \langle 2^j \rangle^{-\frac{5}{6}} \Vert u \Vert_X \\
&\lesssim t^{-\frac{7}{6}+\frac{\delta}{3}+\frac{\kappa}{3}} \mathfrak{t}^{\frac{\delta}{6}+\frac{2\kappa}{3}} 2^{-j-\frac{k}{3}+\delta j+\delta k} \langle 2^j \rangle^{-\frac{5}{6}} \Vert u \Vert_X \\
&\lesssim t^{-\frac{7}{6}+\delta} 2^{-j-\frac{k}{3}+\delta j+\delta k} \langle 2^j \rangle^{-\frac{3}{4}} \Vert u \Vert_X \\
\eqref{estdispsupPMFmoyentau-b-3} &\lesssim \sum_{l_b = l_0+k+\frac{j_{+}}{3}}^{-200} t^{-1} 2^{-l_b-\frac{5j}{2}} \Vert \Psi_{(l_0, l_b, l_0)}^b \psi_{j, k}^{\widehat{\mathcal{P}}} \Vert_{L^2} \langle 2^j \rangle^{-1} \Vert \langle \nabla \rangle |\nabla|^{\frac{1}{2}} g_b(t) \Vert_{L^2} \\
&\lesssim \sum_{l_b = l_0+k}^{-200} t^{-\frac{7}{6}+100\delta} 2^{l_0-\frac{k}{2}-\frac{l_b}{2}-j} \langle 2^j \rangle^{-1} \Vert u \Vert_X \\
&\lesssim t^{-\frac{7}{6}+100\delta} 2^{-k+\frac{l_0}{2}-j} \langle 2^j \rangle^{-1} \Vert u \Vert_X \\
&\lesssim t^{-\frac{7}{6}+100\delta} 2^{-j-\frac{k}{4}} \langle 2^j \rangle^{-1} \Vert u \Vert_X 
\end{align*}
We can obtain the other estimate using $t^{-\frac{1}{6}} 2^{-\frac{k}{3}-\frac{j}{2}} \lesssim 1$. 

Finally, for the internal term: 
\begin{align*}
\eqref{estdispsupPMFmoyentau-int} &\lesssim \Vert \Psi_{(l_0, l_0+k+\frac{j_{+}}{3}, l_0)}^{int} \psi_{j, k}^{\widehat{\mathcal{P}}} \Vert_{L^2_c L^{\frac{1}{1-\kappa}}_{b, a}} \langle 2^j \rangle^{-1} \Vert \langle \overline{\xi} \rangle \widehat{f}(t) \Vert_{L^2_c L^{\frac{1}{\kappa}}_{b, a}} \\
&\lesssim 2^{\frac{5l_0+k}{2}-2\kappa l_0-\kappa k+\frac{5j}{2}-2\kappa j} \langle 2^j \rangle^{-\frac{2}{3}} \Vert u \Vert_X \\
&\lesssim t^{-\frac{7}{6}} \mathfrak{t}^{-\frac{1}{12}+\kappa} 2^{\frac{k}{2}-\kappa k-j-2\kappa j} \langle 2^j \rangle^{-\frac{2}{3}} \Vert u \Vert_X \\
&\lesssim t^{-\frac{7}{6}+\frac{2\kappa}{3}+\frac{\delta}{3}} 2^{-j+\delta k+\delta j} \langle 2^j \rangle^{-\frac{3}{4}} \Vert u \Vert_X 
\end{align*}
and the other estimate follows the same way. 

\paragraph{3.} Let us assume that $\overline{\xi_a} \simeq 2^j$ and $\tau \gg 2^k$. We then have
\begin{align*}
\left| \widehat{X}_a \cdot \nabla \Phi \right| + \left| \widehat{X}_c \cdot \nabla \Phi \right| \simeq 2^{2j} \tau 
\end{align*}
It is therefore possible to apply estimates simpler than in the case $\tau \ll 2^k$. We skip the details. 

\paragraph{Low frequencies} Let us assume that $\overline{\xi_a} \gg 2^j$. As before, we can localise to have
\begin{align*}
\widehat{X}_a \cdot \nabla \Phi, ~ \widehat{X}_c \cdot \nabla \Phi ~ \ll |\overline{\xi_a}|^2 
\end{align*}
and hence $(x, y)$ close to $\mathcal{L}$. In this case, 
\begin{align*}
\widehat{X}_b(\overline{\xi}) \cdot \nabla_{\overline{\xi}} \Phi &= \frac{|\xi|}{|\overline{\xi}|} \left( 3 \xi_0^2 + |\xi|^2 - \frac{x}{t} \right) - \frac{\xi_0 \xi}{|\overline{\xi}| |\xi|} \cdot \left( 2 \xi_0 \xi - \frac{y}{t} \right) \\
&\simeq |\overline{\xi_a}|^2 \\
\widehat{X}_c(\overline{\xi}) \cdot \nabla_{\overline{\xi}} \Phi &\simeq \tau |\overline{\xi_a}|^2 \xi_c \\
\widehat{X}_a(\overline{\xi}) \cdot \nabla_{\overline{\xi}} \Phi &= 3 \xi_0 |\overline{\xi}| \left( 1 - \frac{x}{3 t |\overline{\xi}|^2} \right) - \frac{\xi \cdot y}{t |\overline{\xi}|^2} 
\end{align*}
We recover a similar situation to the high frequency case, and the estimates are simpler. We skip the details. 
\end{Dem}

\subsection{Short time estimate} 

\begin{Lem} Let $t \in (0, 1]$ and $j, k \in \mathbb{Z}$ be such that $k \leq 10, 2^j \gg t^{-\frac{1}{3}}$. Then 
\begin{align*}
\Vert e^{i t \omega(D)} m_{\widehat{\mathcal{R}}}(D) \psi_j(D) f(t) \Vert_{L^{\infty}} &\lesssim t^{-\frac{9}{10}} 2^{-\frac{13j}{10}+\delta j} \Vert u \Vert_X \\
\Vert e^{i t \omega(D)} \psi_{j, k}^{\widehat{\mathcal{L}}}(D) f(t) \Vert_{L^{\infty}} &\lesssim t^{-\frac{9}{10}} 2^{-\frac{13j}{10}+\delta j+\frac{2k}{5}} \Vert u \Vert_X \\
\Vert e^{i t \omega(D)} \psi_{j, k}^{\widehat{\mathcal{C}}}(D) f(t) \Vert_{L^{\infty}} &\lesssim t^{-1+\delta} 2^{-\frac{3j}{2}+3\delta j+\frac{k}{2}} \Vert u \Vert_X \\
\Vert e^{i t \omega(D)} \psi_{j, k}^{\widehat{\mathcal{P}}}(D) f(t) \Vert_{L^{\infty}} &\lesssim t^{-\frac{9}{10}} 2^{-\frac{13j}{10}+\delta j-\frac{k}{10}} \Vert u \Vert_X
\end{align*}
\label{lemestdisptempscourt} 
\end{Lem}

\begin{Dem}
We start with the estimate for $\widehat{\mathcal{L}}$, the estimate for $\widehat{\mathcal{R}}$ being similar. We write: 
\begin{align*}
e^{i t \omega(D)} \psi_{j, k}^{\widehat{\mathcal{L}}}(D) f(t) &= K_{t, j} \ast \left( \psi_{j, k}^{\widehat{\mathcal{L}}}(D) f(t) \right) 
\end{align*}
where the convolution is only in $(x, y)$, and where
\begin{align*}
\widehat{K}_{t, j}(\overline{\xi}) &= e^{i t \omega(\overline{\xi})} m_{\widehat{\mathcal{L}}}(\overline{\xi}) \widetilde{\psi}_j(\overline{\xi}) 
\end{align*}
for $\widetilde{\psi}_j$ a function with similar localisation properties as $\psi_j$, such that $\psi_j \widetilde{\psi}_j = \psi_j$. To simplify notations, we will drop the tilda. In particular, since $\widehat{K}_{t, j}$ is a Schwartz function, we may write
\begin{align*}
K_{t, j}(x, y) &= \int e^{i t \Phi\left( \overline{\xi}, \frac{x}{t}, \frac{y}{t} \right)} m_{\widehat{\mathcal{L}}}(\overline{\xi}) \psi_j(\overline{\xi}) ~ d\overline{\xi} \\
&= \int e^{i t \left( \xi_0 \left( \xi_0^2 + |\xi|^2 \right) - \xi_0 \frac{x}{t} - \xi \cdot \frac{y}{t} \right)} m_{\widehat{\mathcal{L}}}(\overline{\xi}) \psi\left( 2^{-j} \overline{\xi} \right) ~ d\xi_0 d\xi \\
&= \int e^{i t 2^{3j} \left( \xi_0 \left( \xi_0^2 + |\xi|^2 \right) - \xi_0 \frac{2^{-2j} x}{t} - \xi \cdot \frac{2^{-2j} y}{t} \right)} m_{\widehat{\mathcal{L}}}(\overline{\xi}) \psi\left( \overline{\xi} \right) ~ 2^{3j} d\xi_0 d\xi 
\end{align*}
by change of variables. Yet here, $t 2^{3j} \gg 1$ and the phase
\begin{align*}
\Phi\left( \overline{\xi}, \frac{2^{-2j} x}{t}, \frac{2^{-2j} y}{t} \right) &= \xi_0 \left( \xi_0^2 + |\xi|^2 \right) - \xi_0 \frac{2^{-2j} x}{t} - \xi \cdot \frac{2^{-2j} y}{t}
\end{align*}
satisfies 
\begin{align*}
\nabla^2_{\overline{\xi}} \Phi &= \begin{pmatrix} 6 \xi_0 & 2 \xi^T \\ 2 \xi & 2 \xi_0 I_2 \end{pmatrix} 
\end{align*}
which is uniformly invertible on the support of $m_{\widehat{\mathcal{L}}}$ (and of $m_{\widehat{\mathcal{R}}}$), being homogeneous and with determinant 
\begin{align*}
\mbox{det} \nabla^2_{\overline{\xi}} \Phi &= 8 \xi_0 (3 \xi_0^2 - |\xi|^2) 
\end{align*}
We can therefore apply the standard stationary phase estimate and get that 
\begin{align*}
\Vert K_{t, j} \Vert_{L^{\infty}} &\lesssim \left( t 2^{3j} \right)^{-\frac{3}{2}} 2^{3j} = t^{-\frac{3}{2}} 2^{-\frac{3j}{2}} 
\end{align*}
where the implicit constants are uniform in $t, j$. 

In particular, by Young's convolution inequality, 
\begin{align*} \Vert e^{i t \omega(D)} m_{\widehat{\mathcal{L}}}(D) \psi_j(D) \Vert_{L^1 \to L^{\infty}} \lesssim t^{-\frac{3}{2}} 2^{-\frac{3j}{2}} \end{align*}
Clearly, we also have $\Vert e^{i t \omega(D)} m_{\widehat{\mathcal{L}}}(D) \psi_j(D) \Vert_{L^2 \to L^2} \lesssim 1$, so that by Riesz-Thorin theorem, 
\begin{align*}
\Vert e^{i t \omega(D)} m_{\widehat{\mathcal{L}}}(D) \psi_j(D) \Vert_{L^{p'} \to L^p} \lesssim t^{-\frac{3}{2}+\frac{3}{p}} 2^{-\frac{3j}{2}+\frac{3j}{p}}
\end{align*}
for every $p \in [2, \infty]$. 

Finally, we apply successively Sobolev's embedding, the above inequality and Hölder's inequality: 
\begin{align*}
\Vert e^{i t \omega(D)} \psi_{j, k}^{\widehat{\mathcal{L}}}(D) f(t) \Vert_{L^{\infty}} &\lesssim 2^{-\frac{2j}{5}+\delta j+\frac{2k}{5}} \Vert e^{i t \omega(D)} \psi_{j, k}^{\widehat{\mathcal{L}}}(D) \nabla f(t) \Vert_{L^5} \\
&\lesssim t^{-\frac{9}{10}} 2^{-\frac{13j}{10}+\delta j+\frac{2k}{5}} \Vert \nabla f(t) \Vert_{L^{\frac{5}{4}}} \\
&\lesssim t^{-\frac{9}{10}} 2^{-\frac{13j}{10}+\delta j+\frac{2k}{5}} \Vert \langle x, y \rangle \nabla f(t) \Vert_{L^2} \\
&\lesssim t^{-\frac{9}{10}} 2^{-\frac{13j}{10}+\delta j+\frac{2k}{5}} \Vert u \Vert_X
\end{align*}
where we used that $t \in (0, 1]$ is bounded at the last line, and for $\delta > 0$ small enough. 

We can proceed the same way, without $k$, on $\widehat{\mathcal{R}}$. 

\paragraph{Neighborhood of the cone} Let us now prove the estimate in the neighborhood of the cone. 

First, if we assume that $\overline{\xi_a} \ll 2^j$, $\overline{\xi_a} \gg 2^j$ or $(x, y)$ sufficiently away from the physical cone $\mathcal{C}$, then we may proceed as in the neighborhood of $\widehat{\mathcal{L}}$: indeed, in these case, the phase cannot be stationnary and degenerate at the same time. 

Therefore, we assume that $\overline{\xi_a} \simeq 2^j$, and $(x, y)$ in the neighborhood of $\mathcal{C}$. We can then reuse the decomposition used for $t \geq 1$, but without the $b$ vector field. Let us set 
\begin{align*}
\Psi_{l_a}^a\left( \overline{\xi}, \frac{x}{t}, \frac{y}{t} \right) &= \psi\left( 2^{-2j-l_a} \widehat{X}_a(\overline{\xi}) \cdot \nabla_{\overline{\xi}} \Phi \right) \chi\left( 2^{-2j-l_a} \widehat{X}_c(\overline{\xi}) \cdot \nabla_{\overline{\xi}} \Phi \right) \\
\Psi_{l_a}^c\left( \overline{\xi}, \frac{x}{t}, \frac{y}{t} \right) &= \chi\left( 2^{-2j-l_a} \widehat{X}_a(\overline{\xi}) \cdot \nabla_{\overline{\xi}} \Phi \right) \psi\left( 2^{-2j-l_a} \widehat{X}_c(\overline{\xi}) \cdot \nabla_{\overline{\xi}} \Phi \right) \\
\Psi_{l_a}^{int}\left( \overline{\xi}, \frac{x}{t}, \frac{y}{t} \right) &= \chi\left( 2^{-2j-l_a} \widehat{X}_a(\overline{\xi}) \cdot \nabla_{\overline{\xi}} \Phi \right) \chi\left( 2^{-2j-l_a} \widehat{X}_c(\overline{\xi}) \cdot \nabla_{\overline{\xi}} \Phi \right) 
\end{align*}
and then we decompose: 
\begin{subequations}
\begin{align}
\int e^{i t \Phi} \psi_{j, k}^{\widehat{\mathcal{C}}}(\overline{\xi}) \widehat{f}(t, \overline{\xi}) ~ d\overline{\xi} &= \sum_{l_a = l_0}^{200} \int e^{i t \Phi} \Psi_{l_a}^a \psi_{j, k}^{\widehat{\mathcal{C}}}(\overline{\xi}) \widehat{f}(t, \overline{\xi}) ~ d\overline{\xi} \label{estdispsupCtempscourt-a} \\
&\quad + \sum_{l_a = l_0}^{200} \int e^{i t \Phi} \Psi_{l_a}^c \psi_{j, k}^{\widehat{\mathcal{C}}}(\overline{\xi}) \widehat{f}(t, \overline{\xi}) ~ d\overline{\xi} \label{estdispsupCtempscourt-c} \\
&\quad + \int e^{i t \Phi} \Psi_{l_0}^{int} \psi_{j, k}^{\widehat{\mathcal{C}}}(\overline{\xi}) \widehat{f}(t, \overline{\xi}) ~ d\overline{\xi} \label{estdispsupCtempscourt-int} 
\end{align}
\end{subequations}
for $l_0$ such that $2^{l_0} \simeq \mathfrak{t}^{-\frac{1}{2}}$. 

We have the following anitropic estimates on $\Psi_{l_a}^{*} \psi_{j, k}$: 
\begin{align*}
\Vert \Psi_{l_a}^{*} \psi_{j,k} \Vert_{L^{p_c}_c L^{p_b}_b L^{p_a}_a} &\lesssim 2^{\frac{l_a}{p_c}+\frac{k}{p_b}+\frac{l_a}{p_a}+\frac{j}{p_a}+\frac{j}{p_b}+\frac{j}{p_c}}
\end{align*}
and the symbol estimates 
\begin{align*}
\widehat{X}_{\alpha} \cdot \nabla \left[ \Psi_{l_a}^{*} \psi_{j, k} \right] &\lesssim 2^{-l_a-j} 
\end{align*}
for $\alpha = a, c$. 

On \eqref{estdispsupCtempscourt-a}, we apply only one integration by parts along $\widehat{X}_a$: 
\begin{subequations}
\begin{align}
\eqref{estdispsupCtempscourt-a} &= \sum_{l_a = l_0}^{200} t^{-1} 2^{-2l_a-3j} \int e^{i t \Phi} \Psi_{l_a}^a \psi_{j, k}^{\widehat{\mathcal{C}}}(\overline{\xi}) \widehat{f}(t, \overline{\xi}) ~ d\overline{\xi} \label{estdispsupCtempscourt-a-1} \\
&\quad + \sum_{l_a = l_0}^{200} t^{-1} 2^{-l_a-2j} \int e^{i t \Phi} \Psi_{l_a}^a \psi_{j, k}^{\widehat{\mathcal{C}}}(\overline{\xi}) \widehat{X}_a(\overline{\xi}) \cdot \nabla_{\overline{\xi}} \widehat{f}(t, \overline{\xi}) ~ d\overline{\xi} \label{estdispsupCtempscourt-a-2} 
\end{align}
\end{subequations}
where, as usual, the symbols may change from line to line as long as they keep similar properties. We then estimate: 
\begin{align*}
\eqref{estdispsupCtempscourt-a-1} &\lesssim \sum_{l_a = l_0}^{200} t^{-1} 2^{-2l_a-3j} \Vert \Psi_{l_a}^a \psi_{j, k} \Vert_{L^{\frac{1}{1-\kappa}}_c L^2_b L^{\frac{1}{1-\kappa}}_a} 2^{-j} \Vert \overline{\xi} \widehat{f}(t) \Vert_{L^{\frac{1}{\kappa}}_c L^2_b L^{\frac{1}{\kappa}}_a} \\
&\lesssim \sum_{l_a = l_0}^{200} t^{-1} 2^{-2\kappa l_a-\frac{3j}{2}-2\kappa j+\frac{k}{2}} \Vert u \Vert_X \\
&\lesssim t^{-1+\kappa} 2^{-\frac{3j}{2}+\kappa j+\frac{k}{2}} \Vert u \Vert_X \\
\eqref{estdispsupCtempscourt-a-2} &\lesssim \sum_{l_a = l_0}^{200} t^{-1} 2^{-l_a-3j} \Vert \Psi_{l_a}^a \psi_{j, k} \Vert_{L^2} \Vert X_a f(t) \Vert_{L^2} \\
&\lesssim \sum_{l_a = l_0}^{200} t^{-1} 2^{-\frac{3j}{2}+\frac{k}{2}} \Vert u \Vert_X \\
&\lesssim t^{-1+\kappa} 2^{-\frac{3j}{2}+3\kappa j+\frac{k}{2}} \Vert u \Vert_X 
\end{align*}

We proceed the same way on \eqref{estdispsupCtempscourt-c} replacing $\widehat{X}_a$ by $\widehat{X}_c$. 

Finally, for the internal term: 
\begin{align*}
\eqref{estdispsupCtempscourt-int} &\lesssim \Vert \Psi_{l_0}^{int} \psi_{j, k} \Vert_{L^{\frac{1}{1-\kappa}}_c L^2_b L^{\frac{1}{1-\kappa}}_a} 2^{-j} \Vert \overline{\xi} \widehat{f}(t) \Vert_{L^{\frac{1}{\kappa}}_c L^2_b L^{\frac{1}{\kappa}}_a} \\
&\lesssim 2^{2l_0-2\kappa l_0+\frac{3j}{2}-2\kappa j+\frac{k}{2}} \Vert u \Vert_X \\
&\lesssim t^{-1+\kappa} 2^{-\frac{3j}{2}+\kappa j+\frac{k}{2}} \Vert u \Vert_X
\end{align*}
which concludes. 

\paragraph{Neighborhood of the plane} Again, we write
\begin{align*}
e^{i t \omega(D)} \psi_{j, k}^{\widehat{\mathcal{P}}}(D) f(t) &= K_{t, j, k} \ast \left( \psi_{j, k}^{\widehat{\mathcal{P}}}(D) f(t) \right) 
\end{align*}
where this time $K$ also depends on $k$: 
\begin{align*}
K_{t, j, k}(x, y) &= \int e^{i t \Phi\left( \overline{\xi}, \frac{x}{t}, \frac{y}{t} \right)} \psi_{j, k}^{\widehat{\mathcal{P}}}(\overline{\xi}) ~ d\overline{\xi}
\end{align*}
Again, as near the line, one should replace $\psi_{j, k}^{\widehat{\mathcal{P}}}$ by a function having a slightly larger support, but this does not matter for the estimates. We then rescale: 
\begin{align*}
K_{t, j, k}(x, y) &= \int e^{i t \left( \xi_0 \left( \xi_0^2 + |\xi|^2 \right) - \xi_0 \frac{x}{t} - \xi \cdot \frac{y}{t} \right)} \psi\left( 2^{-j} \overline{\xi} \right) \psi\left( 2^{-j-k} \xi_0 \right) ~ d\xi_0 d\xi \\
&= \int e^{i 2^{3j+k} t \left( \xi_0 \left( 2^{2k} \xi_0^2 + |\xi|^2 \right) - \xi_0 \frac{2^{-2j} x}{t} - \xi \cdot \frac{2^{-2j-k} y}{t} \right)} \psi\left( 2^k \xi_0, \xi \right) \psi\left( \xi_0 \right) ~ 2^{3j+k} d\xi_0 d\xi
\end{align*}
Now, the new phase: 
\begin{align*}
\widetilde{\Phi}\left( \overline{\xi}, \frac{2^{-2j} x}{t}, \frac{2^{-2j-k} y}{t} \right) &= \xi_0 \left( 2^{2k} \xi_0^2 + |\xi|^2 \right) - \xi_0 \frac{2^{-2j} x}{t} - \xi \cdot \frac{2^{-2j-k} y}{t}
\end{align*}
satisfies 
\begin{align*}
\nabla^2_{\overline{\xi}} \widetilde{\Phi} &= \begin{pmatrix}
6 \times 2^{2k} \xi_0 & 2 \xi^T \\ 2 \xi & 2 \xi_0 I_2 
\end{pmatrix}
\end{align*}
and is therefore non-degenerate on the considered support, since 
\begin{align*}
\mbox{det} \nabla^2 \widetilde{\Phi} &= 8 \xi_0 \left( 3 \times 2^{2k} \xi_0^2 - |\xi|^2 \right) 
\end{align*}
We then apply the standard stationary phase estimate and get
\begin{align*}
\Vert K_{t, j, k} \Vert_{L^{\infty}} &\lesssim \left( 2^{3j+k} t \right)^{-\frac{3}{2}} 2^{3j+k} = t^{-\frac{3}{2}} 2^{-\frac{3j}{2}-\frac{k}{2}}
\end{align*}
Using the same computations as near $\widehat{\mathcal{L}}$ we finally get
\begin{align*}
\Vert e^{-i t\omega(D)} \psi_{j, k}^{\widehat{\mathcal{P}}}(D) f(t) \Vert_{L^{\infty}} &\lesssim 2^{-\frac{2j}{5}+\delta j+\frac{k}{5}} \Vert K_{t, j, k} \ast \left( \psi_{j, k}^{\widehat{\mathcal{P}}}(D) \nabla f(t) \right) \Vert_{L^5} \\
&\lesssim t^{-\frac{9}{10}} 2^{-\frac{13j}{10}+\delta j-\frac{k}{10}} \Vert \nabla f(t) \Vert_{L^{\frac{5}{4}}} \\
&\lesssim t^{-\frac{9}{10}} 2^{-\frac{13j}{10}+\delta j-\frac{k}{10}} \Vert u \Vert_X
\end{align*}
as expected. 
\end{Dem}

\subsection{Dispersive estimate} 

\begin{Cor} Let $t > 0$ and $m$ be any Hörmander-Mikhlin symbol. Then, 
\begin{align*}
\Vert e^{i t \omega(D)} m(D) \partial_x f(t) \Vert_{L^{\infty}} &\lesssim t^{-\frac{5}{6}} \langle t \rangle^{-\frac{1}{4}+100\delta} \Vert m \Vert_{HM} \Vert u \Vert_X \\
\Vert e^{i t \omega(D)} m(D) \partial_x |\nabla|^{\frac{1}{2}} f(t) \Vert_{L^{\infty}} &\lesssim t^{-\frac{7}{6}} \langle t \rangle^{100\delta} \Vert m \Vert_{HM} \Vert u \Vert_X \\
\Vert e^{i t \omega(D)} m(D) \nabla f(t) \Vert_{L^{\infty}} &\lesssim t^{-\frac{5}{6}} \langle t \rangle^{-\frac{1}{6}+100\delta} \Vert m \Vert_{HM} \Vert u \Vert_X \\
\Vert e^{i t \omega(D)} m(D) f(t) \Vert_{L^{\infty}} &\lesssim t^{-\frac{5}{6}} \langle t \rangle^{\frac{1}{12}+100\delta} \Vert u \Vert_X 
\end{align*}
\label{corestdisptotale} 
\end{Cor}

\begin{Dem}
We decompose according to the geometric areas and the size of the frequency: 
\begin{align*}
&\Vert e^{i t \omega(D)} m(D) \partial_x f(t) \Vert_{L^{\infty}} \lesssim \sum_{\substack{j \in \mathbb{Z}, \\ j \leq j_0}} 2^j \Vert e^{-i t \omega(D)} \psi_j(D) m(D) f(t) \Vert_{L^{\infty}} \\
&\quad + 
\sum_{\substack{j \in \mathbb{Z}, \\ j \geq j_0}} 2^j \Vert e^{i t \omega(D)} m_{\widehat{\mathcal{R}}}(D) \psi_j(D) m(D) f(t) \Vert_{L^{\infty}} \quad + \sum_{\substack{j, k \in \mathbb{Z}, \\ j \geq j_0, \\ k \leq -10}} 2^j \Vert e^{i t \omega(D)} \psi_{j, k}^{\widehat{\mathcal{C}}}(D) m(D) f(t) \Vert_{L^{\infty}} \\
&\quad + \sum_{\substack{j, k \in \mathbb{Z}, \\ j \geq j_0, \\ k \leq -10}} 2^j \Vert e^{i t \omega(D)} \psi_{j, k}^{\widehat{\mathcal{L}}}(D) m(D) f(t) \Vert_{L^{\infty}} \quad + \sum_{\substack{j, k \in \mathbb{Z}, \\ j \geq j_0, \\ k \leq -10}} 2^j \Vert e^{i t \omega(D)} \psi_{j, k}^{\widehat{\mathcal{P}}}(D) m(D) f(t) \Vert_{L^{\infty}}
\end{align*}
For every term, $\psi_j m = \psi_j \widetilde{m}$ for $\widetilde{m}$ a bounded regular function (which corresponds to $m$ localised on frequencies of size $2^j$). In particular, the Fourier multiplier $\widetilde{m}$ acts boundedly on $L^{\infty}$, with an operator norm controlled by $\Vert m \Vert_{HM}$. 

If $t \geq 1$, we separate frequencies according to $2^{j_0} \simeq t^{-\frac{1}{3}} \leq 1$. For the low-frequency part, we apply Lemma \ref{lemestdisph-bassesfreq}: 
\begin{align*}
\sum_{\substack{j \in \mathbb{Z}, \\ j \leq j_0}} 2^j \Vert e^{-i t \omega(D)} \psi_j(D) m(D) f(t) \Vert_{L^{\infty}}
&\lesssim \Vert m \Vert_{HM} \sum_{\substack{j \in \mathbb{Z}, \\ j \leq j_0}} 2^{\frac{7j}{2}} \Vert u \Vert_X \\
&\lesssim \Vert m \Vert_{HM} 2^{\frac{7j_0}{2}} \Vert u \Vert_X \\
&\lesssim t^{-\frac{7}{6}} \Vert m \Vert_{HM} \Vert u \Vert_X
\end{align*}
Then, it is enough to apply Lemmas \ref{lem-estdisph-zonereste}, \ref{lemestdispvoiscone}, \ref{lemestdisp-voisL}, \ref{lemestdispvoisP} and used the given summability in $k$ to get: 
\begin{align*}
&\sum_{\substack{j \in \mathbb{Z}, \\ j \geq j_0}} 2^j \Vert e^{i t \omega(D)} m_{\widehat{\mathcal{R}}}(D) \psi_j(D) m(D) f(t) \Vert_{L^{\infty}} \quad + \sum_{\substack{j, k \in \mathbb{Z}, \\ j \geq j_0, \\ k \leq -10}} 2^j \Vert e^{i t \omega(D)} \psi_{j, k}^{\widehat{\mathcal{C}}}(D) m(D) f(t) \Vert_{L^{\infty}} \\
&\quad + \sum_{\substack{j, k \in \mathbb{Z}, \\ j \geq j_0, \\ k \leq -10}} 2^j \Vert e^{i t \omega(D)} \psi_{j, k}^{\widehat{\mathcal{L}}}(D) m(D) f(t) \Vert_{L^{\infty}} \quad + \sum_{\substack{j, k \in \mathbb{Z}, \\ j \geq j_0, \\ k \leq -10}} 2^j \Vert e^{i t \omega(D)} \psi_{j, k}^{\widehat{\mathcal{P}}}(D) m(D) f(t) \Vert_{L^{\infty}} \\
&\quad \lesssim \sum_{\substack{j \in \mathbb{Z}, \\ j \geq j_0}} t^{-\frac{13}{12}+100\delta} 2^{\delta j} \langle 2^j \rangle^{-\frac{1}{10}} \Vert m \Vert_{HM} \Vert u \Vert_X \\
&\quad \lesssim t^{-\frac{13}{12}+100\delta} \Vert m \Vert_{HM} \Vert u \Vert_X
\end{align*}

If $t \leq 1$, we separate frequencies according to $2^{j_0} \simeq t^{-\frac{5}{9}} \geq 1$ and we apply Lemmas \ref{lemestdisph-bassesfreq} and \ref{lemestdisptempscourt}: for low frequencies, 
\begin{align*}
\sum_{\substack{j \in \mathbb{Z}, \\ j \leq j_0}} 2^j \Vert e^{-i t \omega(D)} \psi_j(D) m(D) f(t) \Vert_{L^{\infty}} &\lesssim \sum_{\substack{j \in \mathbb{Z}, \\ j \leq j_0}} 2^{\frac{7j}{2}} \langle 2^j \rangle^{-2} \Vert u \Vert_X \\
&\lesssim 2^{\frac{3j_0}{2}} \Vert u \Vert_X \\
&\lesssim t^{-\frac{5}{6}} \Vert u \Vert_X 
\end{align*}
then for high frequencies, 
\begin{align*}
&\sum_{\substack{j \in \mathbb{Z}, \\ j \geq j_0}} 2^j \Vert e^{i t \omega(D)} m_{\widehat{\mathcal{R}}}(D) \psi_j(D) m(D) f(t) \Vert_{L^{\infty}} \quad + \sum_{\substack{j, k \in \mathbb{Z}, \\ j \geq j_0, \\ k \leq -10}} 2^j \Vert e^{i t \omega(D)} \psi_{j, k}^{\widehat{\mathcal{C}}}(D) m(D) f(t) \Vert_{L^{\infty}} \\
&\quad + \sum_{\substack{j, k \in \mathbb{Z}, \\ j \geq j_0, \\ k \leq -10}} 2^j \Vert e^{i t \omega(D)} \psi_{j, k}^{\widehat{\mathcal{L}}}(D) m(D) f(t) \Vert_{L^{\infty}} 
\quad + \sum_{\substack{j, k \in \mathbb{Z}, \\ j \geq j_0, \\ k \leq -10}} 2^j \Vert e^{i t \omega(D)} \psi_{j, k}^{\widehat{\mathcal{P}}}(D) m(D) f(t) \Vert_{L^{\infty}} \\
&\quad \lesssim \sum_{\substack{j \in \mathbb{Z}, \\ j \geq j_0}} \left( t^{-\frac{9}{10}} 2^{-\frac{3j}{10}+\delta j} \Vert u \Vert_X + t^{-1+\delta} 2^{-\frac{j}{2}+3\delta j} \Vert u \Vert_X \right) \\
&\quad \lesssim \left( t^{-\frac{9}{10}+\frac{1}{6}-\frac{5\delta}{9}} + t^{-1+\frac{5}{18}-\frac{2\delta}{3}} \right) \Vert u \Vert_X \\
&\quad \lesssim t^{-\frac{5}{6}} \Vert u \Vert_X 
\end{align*}
for small enough $\delta$. 

We proceed in a similar way for the other estimates depending on the number of applied derivatives. 
\end{Dem}

\begin{Rem} We can also allow $m$ to be more singular: notably, we can allow $m$ to satisfy only
\begin{align*}
\nabla_{\xi_0}^{\alpha} \nabla_{\xi}^{\beta} m(\xi_0, \xi) &\lesssim |\xi_0|^{-|\alpha|} |\xi|^{-|\beta|} 
\end{align*}
\end{Rem}

\begin{Rem} Recall that in Corollary \ref{corsymboles-dse}, we obtained the boundedness of bilinear (or multilinear) operators from $L^p \times L^q$ to $L^r$ for $\frac{1}{r} = \frac{1}{p}+\frac{1}{q}$ and $1 < p, q, r < \infty$, with the endpoints excluded due to the failure of the Hörmander-Mikhlin theorem in $L^1$ and $L^{\infty}$. However, in view of the dispersive estimates at hand, we can now allow for $1 < p, q \leq \infty$ (not $r$) up to leaving such symbols in the $L^{\infty}$ bounds, which will then be absorbed by the dispersive estimate. 

To simplify notations, in the estimates in subsequent sections, we will not write explicitely these linear symbols and simply write $\Vert u \Vert_{L^{\infty}}$ (or with derivatives), being accepted that such quantities will always be bounded either applying the dispersive estimate or a Sobolev embedding to get back to some $L^p$ for finite $p$. \label{remarqueestimeedispersiveavecsymbole} 
\end{Rem}

\section{\texorpdfstring{$L^4$}{L4} estimate} \label{section-L4disp} 

The goal of this section is to prove a $L^4$ bound on $e^{-it \omega(D)} h_{\alpha}(t)$. 

If we did not need to separate $m_{\alpha}(D) X_{\alpha} f = g_{\alpha} + h_{\alpha}$ and that we had access to weights in every direction, the estimate would follow from an interpolation argument with the standard $L^1 \to L^{\infty}$ dispersive estimate: indeed, 
\[ \Vert e^{-it \omega(D)} (|x| + |y|) \partial_x f(t) \Vert_{L^4} \lesssim t^{-\frac{1+c}{2}} \Vert (|x| + |y|) \partial_x f(t) \Vert_{L^{4/3}} \]
using that $e^{it\omega(D)}$ is unitary on $L^2$, for some fixed $c > 0$. 

On the other hand, it will prove useful to dispose of (almost) optimal dispersive $L^{\infty}$ estimate, reaching the $t^{-\frac{7}{6}}$ decay, but we will only need a $t^{-\frac{1+c}{2}}$ decay in $L^4$, for some $c > 0$. Indeed, in the proof of the a priori estimate, the weighted $L^4$ estimate is only used for bounds that are uniform in time, and therefore long-time integrability is enough, while the $L^{\infty}$ estimate is also used for the time-dependant bounds (either with growth or decay), so that its precise decay matters. It turns out that, for most of the geometric areas, it will be enough to apply the usual stationary phase estimates, except near the cone $\widehat{\mathcal{C}}$, using that near $\widehat{\mathcal{L}}$ and $\widehat{\mathcal{P}}$, $X_a, X_b, X_c$ are close to canonical weights $x, y$. 

However, near the cone $\widehat{\mathcal{C}}$, we will need to reuse the heuristics of the $L^{\infty}$ estimate and adapt them to the $L^4$ case. In particular, the localisation symbols of the form 
\[ \chi\left( \hbar^{-1} ((x, y) - m(\overline{\xi})) \right) \]
for a function $m$ (that corresponds to derivatives of $\omega$) and $\hbar$ a small localisation parameter, are harder to estimate since we cannot fix $(x, y)$ and see it as a purely frequential symbol. We will therefore also prove a pseudodifferential localisation lemma when such symbols are present, showing they act boundedly $L^2 \to L^2$ (under some hypothesis). We start by the subsection about the pseudo-differential Lemma, and then prove the $L^4$ estimates in following subsections, one by geometric area. 

\subsection{Pseudodifferential localisation lemmas}

The following lemma is not essential to the proof but is simpler than the following near the cone. Here, we consider functions defined on $\mathbb{R}^d$ (not necessarily $\mathbb{R}^3$) and denote by $z \in \mathbb{R}^d$, $\zeta \in \mathbb{R}^d$ the dual variables. Note that the next Lemma could be used to adapt the proof of the $L^{\infty}$ dispersive estimate to the $L^4$ estimate to optimize the decay rate. 

\begin{Lem} There exists $C > 0$ and an integer $k$ such that, for every function $\chi: \mathbb{R}^d \to \mathbb{R}$ of class Schwartz, every positive definite matrix $\hbar$ of size $d \times d$ and every function $\gamma : \mathbb{R}^d \to \mathbb{R}^d$, if we define the pseudodifferential operator: 
\begin{align*}
\Xi(\chi, \hbar, \gamma) ~ : ~ F \in \mathcal{S}(\mathbb{R}^d) ~ \mapsto ~ \left( z \mapsto \int e^{i z \cdot \zeta} \chi\left( \hbar^{-1} (z - \gamma(\zeta)) \right) \widehat{F}(\zeta) ~ d\zeta \right) \in \mathcal{S}'(\mathbb{R}^d) 
\end{align*}
then it extends into a continuous $L^2(\mathbb{R}^d) \to L^2(\mathbb{R}^d)$ operator with 
\[ \Vert \Xi(\chi, \hbar, \gamma) \Vert_{L^2 \to L^2} ~ \leq ~ C \Vert \chi \Vert_{W^{k, 1}} \] \label{lemlocpseudodiff} 
\end{Lem}

\begin{Dem}
We apply Parseval's identity: 
\begin{align*}
\Vert \Xi(\chi, \hbar, \gamma) F \Vert_{L^2} &\lesssim \left\Vert \int \int e^{i z \cdot \zeta} e^{-i z \cdot \upsilon} \chi(\hbar^{-1} (z-\gamma(\zeta))) \widehat{F}(\zeta) ~ d\zeta dz \right\Vert_{L^2_{\upsilon}} \\
&\lesssim \left\Vert \int \widehat{F}(\zeta) \int e^{i z \cdot (\zeta-\upsilon)} \chi(\hbar^{-1} (z-\gamma(\zeta))) ~ dz d\zeta \right\Vert_{L^2_{\upsilon}} \\
&\lesssim \left\Vert \int \widehat{F}(\zeta) \int e^{i (\hbar z + \gamma(\zeta)) \cdot (\zeta-\upsilon)} \chi(z) ~ \mbox{det}(\hbar) dz d\zeta \right\Vert_{L^2_{\upsilon}} \\
&\lesssim \left\Vert \int \widehat{F}(\zeta) e^{i \gamma(\zeta) \cdot (\zeta-\upsilon)} \widehat{\chi}(\hbar (\zeta-\upsilon)) \mbox{det}(\hbar) d\zeta \right\Vert_{L^2_{\eta}} \\
&\lesssim \left\Vert \int |\widehat{F}(\zeta)| |\widehat{\chi}(\hbar (\zeta-\upsilon))| \mbox{det}(\hbar) d\zeta \right\Vert_{L^2_{\upsilon}} \\
&\lesssim \Vert F \Vert_{L^2} \Vert \widehat{\chi}(\hbar \cdot) \mbox{det}(\hbar) \Vert_{L^1} \\
&\lesssim \Vert F \Vert_{L^2} \Vert \widehat{\chi} \Vert_{L^1}
\end{align*}
where the constants are universal. But $\Vert \widehat{\chi} \Vert_{L^1} \lesssim \Vert \chi \Vert_{W^{k, 1}}$ for some integer $k = k(d)$ only depending on the dimension $d$, and the constants are universal again. 
\end{Dem}

We now come back to dimension $3$ and to the notation $(x, y) \in \mathbb{R}^3$, $\overline{\xi} = (\xi_0, \xi)$ the Fourier dual variable. 

Recall that the modified vector fields $\widehat{X}_a', \widehat{X}_b', \widehat{X}_c'$ are defined near the cone $\widehat{\mathcal{C}}$ by 
\begin{align*}
\widehat{X}_a'(\overline{\xi}) &= \begin{pmatrix} \frac{\xi_0}{2 |\xi_0|} \\ \frac{\sqrt{3} \xi}{2 |\xi|} \end{pmatrix} \\
\widehat{X}_b'(\overline{\xi}) &= \begin{pmatrix} \frac{\sqrt{3} \xi_0}{2 |\xi_0|} \\ -\frac{\xi}{2 |\xi|} \end{pmatrix} \\
\widehat{X}_c'(\overline{\xi}) &= \begin{pmatrix} 0 \\ \frac{J \xi}{|\xi|} \end{pmatrix} 
\end{align*}

\begin{Lem} There exist $C > 0$ such that, for any functions $\chi_a, \chi_b, \chi_c : \mathbb{R} \to \mathbb{R}$ of class Schwartz, any $\hbar_a, \hbar_b, \hbar_c > 0$ satisfying $\hbar_c \lesssim \min(\hbar_a, \hbar_b)$, the pseudo-differential operator: {\footnotesize
\begin{align*}
\Xi(\chi, \hbar) ~ : ~ F ~ \mapsto ~ \left( (x, y) \mapsto \int e^{i (x \xi_0 + y \cdot \xi)} \prod_{\alpha = a, b, c} \chi_{\alpha}\left( \hbar_{\alpha}^{-1} \widehat{X}_{\alpha}'(\overline{\xi}) \cdot \begin{pmatrix} x - 3 \xi_0^2 - |\xi|^2 \\ y - 2 \xi_0 \xi \end{pmatrix} \right) ~ ~ m_{\widehat{\mathcal{C}}}(\overline{\xi}) \widehat{F}(\overline{\xi}) ~ d\overline{\xi} \right)
\end{align*}}
defined as $\mathcal{S} \to \mathcal{S}'$, extends into a continuous $L^2 \to L^2$ operator with
\[ \Vert \Xi(\chi, \hbar) \Vert_{L^2 \to L^2} \leq C \sup_{\alpha = a, b, c} \Vert \chi_{\alpha} \Vert_{W^{1, 1}} \] \label{lemlocpseudonontrivial} 
\end{Lem}

\begin{Dem}
In this proof, we will denote by
\begin{align*}
\xi_a &= \frac{1}{2} \left( |\xi_0| + \sqrt{3} |\xi| \right) = \widehat{X}_a'(\overline{\xi}) \cdot \overline{\xi} \\
\xi_b &= \frac{1}{2} \left( \sqrt{3} |\xi_0| - |\xi| \right) = \widehat{X}_b'(\overline{\xi}) \cdot \overline{\xi}
\end{align*}
We always have $\widehat{X}_c'(\overline{\xi}) \cdot \overline{\xi} = 0$. 

We start with the same computation as in the proof of Lemma \ref{lemlocpseudodiff}: 
\begin{align*}
&\Vert \Xi(\chi, \hbar) F \Vert_{L^2} \\
&\lesssim \left\Vert \int \int e^{i (x \xi_0 + y \cdot \xi)} e^{-i (x \eta_0 + y \cdot \eta)} \prod_{\alpha = a, b, c} \chi_{\alpha}\left(\hbar_{\alpha}^{-1} \widehat{X}_{\alpha}(\overline{\xi}) \cdot \begin{pmatrix} x - 3 \xi_0^2 - |\xi|^2 \\ y - 2 \xi_0 \xi \end{pmatrix} \right) m_{\widehat{\mathcal{C}}}(\overline{\xi}) \widehat{F}(\overline{\xi}) ~ d\overline{\xi} dx dy \right\Vert_{L^2_{\overline{\eta}}} \\
&\lesssim \left\Vert \int m_{\widehat{\mathcal{C}}}(\overline{\xi}) \widehat{F}(\overline{\xi}) \int e^{i (x, y) \cdot (\overline{\xi}-\overline{\eta})} \prod_{\alpha = a, b, c} \chi_{\alpha}\left( \hbar_{\alpha}^{-1} \left( \widehat{X}_{\alpha}(\overline{\xi}) \cdot (x, y) + \mu_{\alpha}(\overline{\xi}) \right) \right) ~ dx dy ~ d\overline{\xi} \right\Vert_{L^2_{\overline{\eta}}}
\end{align*}
for some symbols $\mu_a, \mu_b, \mu_c$. Note that, $\overline{\xi}$ being fixed near the cone, $\left( \widehat{X}_a'(\overline{\xi}), \widehat{X}_b'(\overline{\xi}), \widehat{X}_c'(\overline{\xi}) \right)$ forms an orthonormal basis, and hence the matrix 
\[ M(\overline{\xi}) := \left( \widehat{X}_a'(\overline{\xi}), ~ \widehat{X}_b'(\overline{\xi}), ~ \widehat{X}_c'(\overline{\xi}) \right)^T \]
is an orthogonal matrix as well. We may see $M$ as the matrix allowing to switch between standard $(x, y_1, y_2)$ coordinates and $(x_a, x_b, x_c)$ coordinates (depending on $\overline{\xi}$). Therefore, applying the change of coordinates determined by $M$ and a translation by $\mu$, 
\begin{align*}
&\int e^{i (x, y) \cdot (\overline{\xi}-\overline{\eta})} \prod_{\alpha = a, b, c} \chi_{\alpha}\left( \hbar_{\alpha}^{-1} \left( \left( M(\overline{\xi}) (x, y) \right)_{\alpha} - \mu_{\alpha}(\overline{\xi}) \right) \right) ~ dx dy \\
&= \int e^{i M^{-1}(\overline{\xi}) (x_{\alpha}+\mu_{\alpha}(\overline{\xi}))_{\alpha = a, b, c} \cdot (\overline{\xi}-\overline{\eta})} \prod_{\alpha = a, b, c} \chi_{\alpha}\left( \hbar_{\alpha}^{-1} x_{\alpha} \right) ~ dx_a dx_b dx_c \\
&= \int e^{i (\hbar_{\alpha} x_{\alpha}+\mu_{\alpha}(\overline{\xi}))_{\alpha = a, b, c} \cdot M(\overline{\xi}) (\overline{\xi}-\overline{\eta})} \prod_{\alpha = a, b, c} \chi_{\alpha}\left( x_{\alpha} \right) ~ \hbar_a \hbar_b \hbar_c dx_a dx_b dx_c \\
&= e^{i (\mu_{\alpha}(\overline{\xi})_{\alpha = a, b, c} \cdot M(\overline{\xi}) (\overline{\xi}-\overline{\eta})} 
\prod_{\alpha = a, b, c} \widehat{\chi}_{\alpha}\left( \hbar_{\alpha} \left( M(\overline{\xi}) (\overline{\xi}-\overline{\eta}) \right)_{\alpha} \right) \\
&=: e^{i m(\overline{\xi}, \overline{\eta})} \widetilde{\chi}(\overline{\xi}, \overline{\eta}, \hbar) 
\end{align*}

We now show the estimate by duality. 
\begin{align*}
&\left\Vert \int m_{\widehat{\mathcal{C}}}(\overline{\xi}) \widehat{F}(\overline{\xi}) e^{i m(\overline{\xi}, \overline{\eta})} \widetilde{\chi}(\overline{\xi}, \overline{\eta}, \hbar) ~ d\overline{\xi} \right\Vert_{L^2_{\overline{\eta}}} \\
&= \sup_{G \in L^2, \Vert G \Vert_{L^2} \leq 1} \int \int m_{\widehat{\mathcal{C}}}(\overline{\xi}) \widehat{F}(\overline{\xi}) e^{i m(\overline{\xi}, \overline{\eta})} \widetilde{\chi}(\overline{\xi}, \overline{\eta}, \hbar) G(\overline{\eta}) ~ d\overline{\xi} d\overline{\eta} \\
&\lesssim \left( \int \int m_{\widehat{\mathcal{C}}}(\overline{\xi}) |\widehat{F}(\overline{\xi})|^2 |\widetilde{\chi}(\overline{\xi}, \overline{\eta}, \hbar)| ~ d\overline{\xi} d\overline{\eta} \right)^{\frac{1}{2}} \sup_{G \in L^2, \Vert G \Vert_{L^2} \leq 1} \left( \int \int m_{\widehat{\mathcal{C}}}(\overline{\xi}) |G(\overline{\eta})|^2 |\widetilde{\chi}(\overline{\xi}, \overline{\eta}, \hbar)| ~ d\overline{\xi} d\overline{\eta} \right)^{\frac{1}{2}}
\end{align*}
by Cauchy-Schwarz inequality. 

The first term can be directly estimated. Indeed, since $M(\overline{\xi})$ is orthogonal at any point, we may apply the change of variables $(\eta_a, \eta_b, \eta_c) = M(\overline{\xi}) (\overline{\xi}-\overline{\eta})$ and get 
\begin{align*}
&\left( \int \int m_{\widehat{\mathcal{C}}}(\overline{\xi}) |\widehat{F}(\overline{\xi})|^2 |\widetilde{\chi}(\overline{\xi}, \overline{\eta}, \hbar)| ~ d\overline{\xi} d\overline{\eta} \right)^{\frac{1}{2}} \\
&\lesssim \left( \int m_{\widehat{\mathcal{C}}}(\overline{\xi}) |\widehat{F}(\overline{\xi})|^2 \prod_{\alpha = a, b, c} \left( \int |\widehat{\chi_{\alpha}}(\hbar_{\alpha} \eta_{\alpha})| \hbar_{\alpha} d\eta_{\alpha} \right) ~ d\overline{\xi} \right)^{\frac{1}{2}} \\
&\lesssim \Vert F \Vert_{L^2} \prod_{\alpha = a, b, c} \Vert \widehat{\chi_{\alpha}} \Vert_{L^1} \\
&\lesssim \Vert F \Vert_{L^2} \prod_{\alpha = a, b, c} \Vert \chi_{\alpha} \Vert_{W^{1, 1}}
\end{align*}

It only remains to estimate the second term uniformly in $G$. We rewrite it as: 
\begin{align*}
\left( \int |G(\overline{\eta})|^2 \int m_{\widehat{\mathcal{C}}}(\overline{\xi}) |\widetilde{\chi}(\overline{\xi}, \overline{\eta}, \hbar)| ~ d\overline{\xi} ~ d\overline{\eta} \right)^{\frac{1}{2}} 
\end{align*}
In the internal integral in $\overline{\xi}$, given the localisation by $m_{\widehat{\mathcal{C}}}$, we may rewrite everything with the coordinates 
\begin{align*}
\xi_a = \frac{1}{2} (|\xi_0| + \sqrt{3} |\xi|), \quad \xi_b = \frac{1}{2} (\sqrt{3} |\xi_0| - |\xi|), 
\end{align*}
and $\theta$ the angle between $\xi$ and $\eta$. This change of variables brings to an integral in $\xi_a \in \mathbb{R}^{+}$, $\xi_b \in (-2^{-10} \xi_a, 2^{-10} \xi_a)$, $\theta \in [-\pi, \pi]$, and up to separating the integral in two pieces, we only deal with the case $\xi_0 > 0$ (the other being symmetric). Then, for fixed $\overline{\eta}$ (such that $\eta \neq 0$): 
\begin{align*}
&\int m_{\widehat{\mathcal{C}}}(\overline{\xi}) |\widetilde{\chi}(\overline{\xi}, \overline{\eta}, \hbar)| ~ d\overline{\xi} \\
&\simeq \int_0^{\infty} \int_{-2^{-10} \xi_a}^{2^{-10} \xi_a} \int_{-\pi}^{\pi} \left|\widehat{\chi_a} \left( \hbar_a \left( \xi_a - \frac{\eta_0 + \sqrt{3} |\eta| \cos(\theta)}{2} \right) \right) \right| \left| \widehat{\chi_b}\left(\hbar_b \left( \xi_b - \frac{\sqrt{3} \eta_0 - |\eta| \cos(\theta)}{2} \right) \right) \right| \\
&\pushright{ \left| \widehat{\chi_c}\left( \hbar_c |\eta| \sin(\theta) \right) \right| \hbar_a \hbar_b \hbar_c  ~ \frac{\sqrt{3} \xi_a - \xi_b}{2} d\theta d\xi_b d\xi_a} \\
&\lesssim \hbar_a \hbar_b \hbar_c \int_{-\pi}^{\pi} \left| \widehat{\chi_c}\left( \hbar_c |\eta| \sin(\theta) \right) \right| \int_0^{\infty} \left|\widehat{\chi_a} \left( \hbar_a \left( \xi_a - \frac{\eta_0 + |\eta| \cos(\theta)}{2} \right) \right) \right| \\
&\pushright{ \int_{-2^{-10} \xi_a}^{2^{-10} \xi_a} \left| \widehat{\chi_b}\left(\hbar_b \left( \xi_b - \frac{\sqrt{3} \eta_0 - |\eta| \cos(\theta)}{2} \right) \right) \right| \xi_a ~ d\xi_b d\xi_a ~ d\theta }
\end{align*}

First, we consider the case $\{ |\overline{\eta}| \lesssim |\eta| \}$. Then we have (uniformly in $\xi_a, \theta, \overline{\eta}$): 
\begin{align*}
&\hbar_b \int_{-2^{-10} \xi_a}^{2^{-10} \xi_a} \left| \widehat{\chi_b}\left(\hbar_b \left( \xi_b - \frac{\sqrt{3} \eta_0 - |\eta| \cos(\theta)}{2} \right) \right) \right| ~ d\xi_b ~ \\
&= \int_{-2^{-10} \hbar_b \xi_a}^{2^{-10} \hbar_b \xi_a} \left| \widehat{\chi_b}\left( \xi_b - \hbar_b \frac{\sqrt{3} \eta_0 - |\eta| \cos(\theta)}{2} \right) \right| ~ d\xi_b \lesssim 1
\end{align*}
and then
\begin{align*}
\hbar_a \int_0^{\infty} \left|\widehat{\chi_a} \left( \hbar_a^{-1} \left( \xi_a - \frac{\eta_0 + \sqrt{3} |\eta| \cos(\theta)}{2} \right) \right) \right| \xi_a ~ d\xi_a ~ &\lesssim \left( \hbar_a^{-1} + \left| \frac{\eta_0}{2} + |\eta| \cos(\theta) \right|  \right) \\
&\lesssim \hbar_a^{-1} + |\overline{\eta}| 
\end{align*}
Finally, using $\hbar_c \leq \hbar_a$ and $|\overline{\eta}| \lesssim |\eta|$, 
\begin{align*}
&\hbar_c \int_{-\pi}^{\pi} \left| \widehat{\chi_c}\left( \hbar_c |\eta| \sin(\theta) \right) \right| \left( \hbar_a^{-1} + |\overline{\eta}| \right) ~ d\theta \\
&\quad \lesssim \hbar_c \hbar_a^{-1} \int_{-\pi}^{\pi} 1 ~ d\theta \quad + \int_{-\pi}^{\pi} \left| \widehat{\chi_c}\left( \hbar_c |\eta| \sin(\theta) \right) \right| \hbar_c |\eta| d\theta \\
&\quad \lesssim 1
\end{align*}
Therefore, 
\begin{align*}
&\left( \int 1_{|\eta| \gtrsim |\overline{\eta}|} |G(\overline{\eta})|^2 \int m_{\widehat{\mathcal{C}}}(\overline{\xi}) |\widetilde{\chi}(\overline{\xi}, \overline{\eta}, \hbar)| ~ d\overline{\xi} ~ d\overline{\eta} \right)^{\frac{1}{2}} \\
&\quad \lesssim \left( \int 1_{|\eta| \gtrsim |\overline{\eta}|} |G(\overline{\eta})|^2 ~ d\overline{\eta} \right)^{\frac{1}{2}} \\
&\quad \lesssim \Vert G \Vert_{L^2} \leq 1
\end{align*}

We now turn to the case $\{ |\overline{\eta}| \lesssim \hbar_c^{-1} \}$, we can reuse the previous computation and obtain again 
\begin{align*}
&\hbar_c \int_{-\pi}^{\pi} \left| \widehat{\chi_c}\left( \hbar_c |\eta| \sin(\theta) \right) \right| \left( \hbar_a^{-1} + |\overline{\eta}| \right) ~ d\theta \\
&\quad \lesssim \hbar_c \hbar_a^{-1} \int_{-\pi}^{\pi} 1 ~ d\theta \quad + \int_{-\pi}^{\pi} \left| \widehat{\chi_c}\left( \hbar_c |\eta| \sin(\theta) \right) \right| d\theta \\
&\quad \lesssim 1
\end{align*}

Finally, in the case $|\eta| \ll |\overline{\eta}|$, i.e. $|\overline{\eta}| \simeq |\eta_0|$, and $|\eta_0| \gg \hbar_c^{-1}$, we compute more finely: 
\begin{align*}
&\hbar_b \int_{-2^{-10} \xi_a}^{2^{-10} \xi_a} \left| \widehat{\chi_b}\left(\hbar_b \left( \xi_b - \frac{\sqrt{3} \eta_0 - |\eta| \cos(\theta)}{2} \right) \right) \right| ~ d\xi_b ~ \\
&= \int_{-2^{-10} \hbar_b \xi_a}^{2^{-10} \hbar_b \xi_a} \left| \widehat{\chi_b}\left( \xi_b - \hbar_b \frac{\sqrt{3} \eta_0 - |\eta| \cos(\theta)}{2} \right) \right| ~ d\xi_b \\
&\lesssim \frac{1}{\left( 1 + \hbar_b \mbox{dist}\left( [0, \xi_a], 2^{10} \left| \frac{\sqrt{3} \eta_0 - |\eta| \cos(\theta)}{2} \right| \right) \right)^{100}}
\end{align*}
Then, setting $\eta_a = \eta_a(\theta) = \frac{\eta_0 + \sqrt{3} |\eta| \cos(\theta)}{2}$ and $\eta_b = \eta_b(\theta) = \frac{\sqrt{3} \eta_0 - |\eta| \cos(\theta)}{2}$, the integral in $\xi_a$ becomes: 
\begin{align*}
&\hbar_a \int_0^{\infty} \left|\widehat{\chi_a} \left( \hbar_a \left( \xi_a - \eta_a \right) \right) \right| \frac{\xi_a}{\left( 1 + \hbar_b \mbox{dist}\left( [0, \xi_a], 2^{10} \left| \eta_b \right| \right) \right)^{100}} ~ d\xi_a ~ \\
&\quad \lesssim \int_{\{ |\xi_a - \eta_a| \lesssim |\overline{\eta}| \}} \hbar_a \left|\widehat{\chi_a} \left( \hbar_a \left( \xi_a - \eta_a \right) \right) \right| \frac{\xi_a}{\left( 1 + \hbar_b \mbox{dist}\left( [0, \xi_a], 2^{10} |\eta_b| \right) \right)^{100}} ~ d\xi_a \\
&\quad \quad \quad + \int_{ \{ |\xi_a - \eta_a| \gtrsim |\overline{\eta}| \} } \hbar_a \xi_a \left|\widehat{\chi_a} \left( \hbar_a \left( \xi_a - \eta_a \right) \right) \right| ~ d\xi_a 
\end{align*}
But here $|\overline{\eta}| \simeq |\eta_a| \simeq |\eta_b| \gg \hbar_c^{-1} \gtrsim \hbar_b^{-1}$ so on $\{ |\xi_a - \eta_a| \lesssim |\overline{\eta}| \}$, we have
\[ \frac{\xi_a}{\left( 1 + \hbar_b \mbox{dist}\left( [0, \xi_a], 2^{10} |\eta_b| \right) \right)^{100}} \lesssim \hbar_b^{-100} |\overline{\eta}|^{-99} \]
Hence: 
\begin{align*}
&\int_{\{ |\xi_a - \eta_a| \lesssim |\overline{\eta}| \}} \hbar_a \left|\widehat{\chi_a} \left( \hbar_a \left( \xi_a - \eta_a \right) \right) \right| \frac{\xi_a}{\left( 1 + \hbar_b \mbox{dist}\left( [0, \xi_a], 2^{10} |\eta_b| \right) \right)^{100}} ~ d\xi_a \\
&\quad \lesssim \hbar_b^{-100} |\overline{\eta}|^{-99} \\
&\quad \lesssim \hbar_a \hbar_b^{-100} \hbar_c^{99} \lesssim \hbar_b^{-1}
\end{align*}
using $|\overline{\eta}| \gtrsim \hbar_c \gtrsim \hbar_b$. On the other hand, 
\begin{align*}
&\int_{ \{ |\xi_a - \eta_a| \gtrsim |\overline{\eta}| \} } \hbar_a \xi_a \left|\widehat{\chi_a} \left( \hbar_a \left( \xi_a - \eta_a \right) \right) \right| ~ d\xi_a \\
&\quad = \hbar_a^{-1} \int_{ \{ |\xi_a - \hbar_a \eta_a| \gtrsim \hbar_a |\overline{\eta}| \} } \xi_a  \left|\widehat{\chi_a} \left( \xi_a - \hbar_a \eta_a \right) \right| ~ d\xi_a \\
&\quad = \hbar_a^{-1} \int_{ \{ |\xi_a| \gtrsim \hbar_a |\overline{\eta}| \} } \left( \xi_a + \hbar_a \eta_a \right)  \left|\widehat{\chi_a} \left( \xi_a \right) \right| ~ d\xi_a \\
&\quad \lesssim \hbar_a^{-1} 
\end{align*}
Finally, we integrate in $\xi_c$: 
\begin{align*}
&\hbar_c \int_{-\pi}^{\pi} (\hbar_a^{-1} + \hbar_b^{-1}) \left| \widehat{\chi_c}\left( \hbar_c |\eta| \sin(\theta) \right) \right| ~ d\theta \\
&\quad \lesssim \hbar_c \left( \hbar_a^{-1} + \hbar_b^{-1} \right) \lesssim 1
\end{align*}
and we conclude the same way. 
\end{Dem}

\subsection{Low frequencies}

\begin{Lem} Let $t > 0$, $\alpha \in \{ a, b, c \}$, and $j \in \mathbb{Z}$. Then
\begin{align*}
\Vert e^{-i t \omega(D)} \psi_j(D) h_{\alpha}(t) \Vert_{L^4} &\lesssim 2^{\frac{3j}{4}} \langle 2^j \rangle^{-1} \Vert u \Vert_X
\end{align*} \label{lemestdispL4BF} 
\end{Lem}

\begin{Dem}
We apply Sobolev's embedding: 
\begin{align*}
\Vert e^{-i t \omega(D)} \psi_j(D) h_{\alpha}(t) \Vert_{L^4} &\lesssim 2^{\frac{3j}{4}} \Vert \psi_j(D) h_{\alpha}(t) \Vert_{L^2} \\
&\lesssim 2^{\frac{3j}{4}} \langle 2^j \rangle^{-1} \Vert u \Vert_X
\end{align*}
as wanted. 
\end{Dem}

\subsection{Remainder area}

\begin{Lem} Let $t > 0$ and $\alpha \in \{ a, b, c \}$. Then
\begin{align*}
\Vert e^{-i t \omega(D)} m_{\widehat{\mathcal{R}}}(D) \nabla h_{\alpha}(t) \Vert_{L^4} &\lesssim t^{-\frac{5}{12}} \langle t \rangle^{-\frac{1}{12}-\frac{1}{84}} \Vert u \Vert_X
\end{align*}
\end{Lem}

\begin{Dem}
We will show first that
\begin{align*}
\Vert e^{-i t \omega(D)} \psi_j(D) m_{\widehat{\mathcal{R}}}(D) h_{\alpha}(t) \Vert_{L^4} &\lesssim t^{-\frac{3}{4}} \langle t \rangle^{\frac{1}{24}} 2^{-\frac{13j}{8}} \langle 2^j \rangle^{-\frac{1}{8}} \Vert u \Vert_X
\end{align*}
for $j \in \mathbb{Z}$ such that $2^j \gg t^{-\frac{1}{3}}$. 

As we saw in the proof of Lemma \ref{lemestdisptempscourt}, if we set
\begin{align*}
\widehat{K}_{t, j}(\overline{\xi}) &= e^{-i t \omega(\overline{\xi})} \widetilde{\psi}_j(\overline{\xi}) m_{\widehat{\mathcal{R}}}(\overline{\xi}) 
\end{align*}
for $\widetilde{\psi}$ a function similar to $\psi$ with larger support, 
then $\Vert K_{t, j} \Vert_{L^{\infty}} \lesssim t^{-\frac{3}{2}} 2^{-\frac{3j}{2}}$ uniformly in $t, j$, and
\begin{align*}
e^{-i t \omega(D)} \psi_j(D) m_{\widehat{\mathcal{R}}}(D) h_{\alpha}(t) = K_{t, j} \ast \left( \psi_j(D) m_{\widehat{\mathcal{R}}}(D) h_{\alpha}(t) \right) 
\end{align*}
Therefore, the linear operator $F \mapsto K_{t, j} \ast F$ is continuous from $L^1$ to $L^{\infty}$ by Young's convolution inequality, with an operator norm controlled by $t^{-\frac{3}{2}} 2^{-\frac{3j}{2}}$, and it is also clearly bounded as $L^2 \to L^2$ operator by $1$. Hence, by Riesz-Thorin theorem, 
\begin{align*}
&\Vert e^{-i t \omega(D)} \psi_j(D) m_{\widehat{\mathcal{R}}}(D) h_{\alpha}(t) \Vert_{L^4} \lesssim \left( t^{-\frac{3}{2}} 2^{-\frac{3j}{2}} \right)^{\frac{1}{2}} \Vert \psi_j(D) m_{\widehat{\mathcal{R}}}(D) h_{\alpha}(t) \Vert_{L^{\frac{4}{3}}} \\
&\lesssim t^{-\frac{3}{4}} 2^{-\frac{3j}{4}} \Vert \langle x, y \rangle^{\frac{7}{8}} \psi_j(D) m_{\widehat{\mathcal{R}}}(D) h_{\alpha}(t) \Vert_{L^2} \\
&\lesssim t^{-\frac{3}{4}} 2^{-\frac{3j}{4}} \left( \Vert \psi_j(D) h_{\alpha}(t) \Vert_{L^2} + \Vert \psi_j(D) h_{\alpha}(t) \Vert_{L^2}^{\frac{1}{8}} \Vert (x, y) \psi_j(D) m_{\widehat{\mathcal{R}}}(D) h_{\alpha}(t) \Vert_{L^2}^{\frac{7}{8}} \right) \\
&\lesssim t^{-\frac{3}{4}} 2^{-\frac{3j}{4}} \Biggl( 2^{-\frac{7j}{8}} \Vert |\nabla|^{\frac{7}{8}} \psi_j(D) h_{\alpha}(t) \Vert_{L^2} + 2^{-\frac{7j}{8}} \Vert \psi_j(D) h_{\alpha}(t) \Vert_{L^2} \\
&\quad \quad + 2^{-\frac{7j}{8}} 
\Vert \psi_j(D) h_{\alpha}(t) \Vert_{L^2}^{\frac{1}{8}} \Vert \nabla \psi_j(D) m_{\widehat{\mathcal{R}}}(D) (x, y) h_{\alpha}(t) \Vert_{L^2}^{\frac{7}{8}} \Biggl) \\
&\lesssim t^{-\frac{3}{4}} 2^{-\frac{13j}{8}} \langle 2^j \rangle^{-\frac{1}{8}} \left( \Vert h_{\alpha}(t) \Vert_{H^1} + 
\Vert h_{\alpha}(t) \Vert_{H^1}^{\frac{1}{8}} \langle t \rangle^{\frac{7}{8}\frac{1}{24}+O(\delta)} \Vert u \Vert_X^{\frac{7}{8}} \right) \\
&\lesssim t^{-\frac{3}{4}} \langle t \rangle^{\frac{1}{24}} 2^{-\frac{13j}{8}} \langle 2^j \rangle^{-\frac{1}{8}} \Vert u \Vert_X 
\end{align*}
as wanted. 

To obtain the estimate on $\Vert e^{-i t \omega(D)} m_{\widehat{\mathcal{R}}}(D) \nabla h_{\alpha}(t) \Vert_{L^4}$, we decompose and apply either the above inequality or the result of Lemma \ref{lemestdispL4BF} for low frequencies. First, if $t \geq 1$, 
\begin{align*}
&\Vert e^{-i t \omega(D)} m_{\widehat{\mathcal{R}}}(D) \nabla h_{\alpha}(t) \Vert_{L^4} \\
&\lesssim \sum_{\substack{j \in \mathbb{Z}, \\ 2^j \lesssim t^{-\frac{13}{42}}}} \Vert e^{-i t \omega(D)} m_{\widehat{\mathcal{R}}}(D) \psi_j(D) \nabla h_{\alpha}(t) \Vert_{L^4} \quad + \sum_{\substack{j \in \mathbb{Z}, \\ 2^j \gtrsim t^{-\frac{13}{42}}}} \Vert e^{-i t \omega(D)} m_{\widehat{\mathcal{R}}}(D) \psi_j(D) \nabla h_{\alpha}(t) \Vert_{L^4} \\
&\lesssim \sum_{\substack{j \in \mathbb{Z}, \\ 2^j \lesssim t^{-\frac{13}{42}}}} 2^{\frac{7j}{4}} \Vert u \Vert_X \quad + \sum_{\substack{j \in \mathbb{Z}, \\ 2^j \gg t^{-\frac{13}{42}}}} t^{-\frac{17}{24}} 2^{-\frac{5j}{8}} \Vert u \Vert_X \\
&\lesssim t^{-\frac{13}{24}} \Vert u \Vert_X \quad + t^{-\frac{173}{336}} \Vert u \Vert_X \\
&\lesssim t^{-\frac{1}{2}-\frac{1}{84}} \Vert u \Vert_X 
\end{align*}
On the other hand, if $t \leq 1$: 
\begin{align*}
&\Vert e^{-i t \omega(D)} m_{\widehat{\mathcal{R}}}(D) \nabla h_{\alpha}(t) \Vert_{L^4} \\
&\lesssim \sum_{\substack{j \in \mathbb{Z}, \\ 2^j \lesssim t^{-\frac{1}{2}}}} \Vert e^{-i t \omega(D)} m_{\widehat{\mathcal{R}}}(D) \psi_j(D) \nabla h_{\alpha}(t) \Vert_{L^4} \quad + \sum_{\substack{j \in \mathbb{Z}, \\ 2^j \gtrsim t^{-\frac{1}{2}}}} \Vert e^{-i t \omega(D)} m_{\widehat{\mathcal{R}}}(D) \psi_j(D) \nabla h_{\alpha}(t) \Vert_{L^4} \\
&\lesssim \sum_{\substack{j \in \mathbb{Z}, \\ 2^j \lesssim t^{-\frac{1}{2}}}} 2^{\frac{7j}{4}} \langle 2^j \rangle^{-1} \Vert u \Vert_X \quad + \sum_{\substack{j \in \mathbb{Z}, \\ 2^j \gg t^{-\frac{1}{2}}}} t^{-\frac{3}{4}} 2^{-\frac{3j}{4}} \Vert u \Vert_X \\
&\lesssim t^{-\frac{3}{8}} \Vert u \Vert_X 
\end{align*}
which concludes since $\frac{3}{8} \leq \frac{5}{12}$. 
\end{Dem}

\subsection{Neighborhood of the line}

\begin{Lem} Let $t > 0$. Then
\begin{align*}
\Vert e^{-i t \omega(D)} m_{\widehat{\mathcal{L}}}(D) \nabla (h_a(t), h_c(t)) \Vert_{L^4} &\lesssim t^{-\frac{5}{12}} \langle t \rangle^{-\frac{1}{12}+O(\delta)} \Vert u \Vert_X \\
\Vert e^{-i t \omega(D)} m_{\widehat{\mathcal{L}}}(D) \nabla h_b(t) \Vert_{L^4} &\lesssim t^{-\frac{5}{12}} \langle t \rangle^{-\frac{1}{12}-\frac{1}{42}} \Vert u \Vert_X
\end{align*}
\end{Lem}

\begin{Dem}
We consider 
\begin{align*}
\widehat{K}_{t, j, k}(\overline{\xi}) &:= e^{-i t \omega(\overline{\xi})} \psi_{j, k}^{\widehat{\mathcal{L}}}(\overline{\xi}) 
\end{align*}
Let us also denote 
\begin{align*}
\widehat{K}_{t, j}(\overline{\xi}) &= e^{-i t \omega(\overline{\xi})} \psi_j(\overline{\xi}) m_{\widehat{\mathcal{L}}}(\overline{\xi}) 
\end{align*}

We already saw in the proof of Lemma \ref{lemestdisptempscourt} that $\Vert K_{t, j} \Vert_{L^{\infty}} \lesssim t^{-\frac{3}{2}} 2^{-\frac{3j}{2}}$, and therefore by Riesz-Thorin interpolation inequality, on the one hand, 
\begin{align*}
&\Vert e^{-i t \omega(D)} \psi_{j, k}^{\widehat{\mathcal{L}}}(D) h_{\alpha}(t) \Vert_{L^4} \lesssim t^{-\frac{3}{4}} 2^{-\frac{3j}{4}} \Vert \psi_{j, k}^{\widehat{\mathcal{L}}}(D) h_{\alpha}(t) \Vert_{L^{\frac{4}{3}}} \\
&\lesssim t^{-\frac{3}{4}} 2^{-\frac{3j}{4}} \Vert \langle x \rangle^{\frac{1}{4}+\delta} \langle y \rangle^{\frac{1}{2}+\delta} \psi_{j, k}^{\widehat{\mathcal{L}}}(D) h_{\alpha}(t) \Vert_{L^2} \\
&\lesssim t^{-\frac{3}{4}} 2^{-\frac{3j}{4}} \Vert \psi_{j, k}^{\widehat{\mathcal{L}}}(D) h_{\alpha}(t) \Vert_{L^2}^{\frac{1}{4}-2\delta} \Vert \langle x \rangle \psi_{j, k}^{\widehat{\mathcal{L}}}(D) h_{\alpha}(t) \Vert_{L^2}^{\frac{1}{4}+\delta} \Vert \langle y \rangle \psi_{j, k}^{\widehat{\mathcal{L}}}(D) h_{\alpha}(t) \Vert_{L^2}^{\frac{1}{2}+\delta} \\
&\lesssim t^{-\frac{3}{4}} 2^{-\frac{3j}{4}} \Vert \psi_j(D) h_{\alpha}(t) \Vert_{L^2}^{\frac{1}{4}-2\delta} \left( \left( 1 + 2^{-j} \right) \Vert \psi_j(D) h_{\alpha}(t) \Vert_{L^2} + \Vert \psi_j(D) x h_{\alpha}(t) \Vert_{L^2} \right)^{\frac{1}{4}+\delta} \\
&\quad \quad \left( \left( 1 + 2^{-j-k} \right) \Vert \psi_j(D) h_{\alpha}(t) \Vert_{L^2} + \Vert \psi_{j, k}(D) y h_{\alpha}(t) \Vert_{L^2} \right)^{\frac{1}{2}+\delta} \\
&\lesssim t^{-\frac{3}{4}} 2^{-\frac{3j}{4}} \langle 2^j \rangle^{-\frac{1}{4}+2\delta} \Vert \langle \nabla \rangle \psi_j(D) h_{\alpha}(t) \Vert_{L^2}^{\frac{1}{4}-2\delta} \left( 2^{-j} \Vert u \Vert_X \right)^{\frac{1}{4}+\delta} \left( 2^{-j-k} \langle t \rangle^{O(\delta)} \Vert u \Vert_X \right)^{\frac{1}{2}+\delta} \\
&\lesssim t^{-\frac{3}{4}} 2^{-\frac{3j}{2}-2\delta j-\frac{k}{2}-\delta k} \langle 2^j \rangle^{-\frac{1}{4}+2\delta} \Vert u \Vert_X
\end{align*}

On the other hand, we can show for $K_{t, j, k}$ an anisotropic estimate in $(x, y)$. We have that 
\begin{align*}
\Vert K_{t, j, k} \Vert_{L^{\infty}_x L^2_y} &= \left\Vert \left\Vert \int e^{-i t \omega(\overline{\xi}) - i x \xi_0 - i y \cdot \xi} \psi_{j, k}^{\widehat{\mathcal{L}}}(\overline{\xi}) ~ d\xi_0 d\xi \right\Vert_{L^2_y} \right\Vert_{L^{\infty}_x} \\
&= \left\Vert \left\Vert \int e^{-i t \omega(\overline{\xi}) - i x \xi_0} \psi_{j, k}^{\widehat{\mathcal{L}}}(\overline{\xi}) ~ d\xi_0 \right\Vert_{L^2_{\xi}} \right\Vert_{L^{\infty}_x} \\
&\lesssim \left\Vert \int e^{-i t \omega(\overline{\xi}) - i x \xi_0} \psi_{j, k}^{\widehat{\mathcal{L}}}(\overline{\xi}) ~ d\xi_0 \right\Vert_{L^{\infty}_{x, \xi}} ~ \left( 2^{2j+2k} \right)^{\frac{1}{2}} \\
&\lesssim \left\Vert \int e^{-i t 2^{3j} \xi_0^3 - i t x \xi_0} \psi_{0, 0}^{\widehat{\mathcal{L}}}(\overline{\xi}) ~ 2^j d\xi_0 \right\Vert_{L^{\infty}_{x, \xi}} ~ 2^{j+k} 
\end{align*}
by Parseval's identity and by change of variable. But the phase $\xi_0 \mapsto - \xi_0^3 - x \xi_0$ is non-degenerate, uniformly in $\xi_0$ in the neighborhood of $\{ -1, 1 \}$, and $\psi_{0, 0}^{\widehat{\mathcal{L}}}$ is a smooth function, so by the van der Corput lemma we get
\begin{align*}
\Vert K_{t, j, k} \Vert_{L^{\infty}_x L^2_y} &\lesssim 2^{2j+k} \left( t 2^{3j} \right)^{-\frac{1}{2}} = t^{-\frac{1}{2}} 2^{\frac{j}{2}+k} 
\end{align*}
Applying again Young's convolution inequality, the linear operator $F \mapsto K_{t, j, k} \ast F$ is bounded from $L^1_x L^2_y$ to $L^{\infty}_{x, y}$ and we control its operator norm by $t^{-\frac{1}{2}} 2^{\frac{j}{2}+k}$. On the other hand, this operator is also bounded from $L^2$ to $L^2$, controlled by $1$. By Riesz-Thorin theorem, we deduce a $L^{4/3}_x L^2_y \to L^4_{x, y}$ bound, given by: 
\begin{align*}
&\Vert e^{-i t \omega(D)} \psi_{j, k}^{\widehat{\mathcal{L}}}(D) h_{\alpha}(t) \Vert_{L^4} \lesssim t^{-\frac{1}{4}} 2^{\frac{j}{4}+\frac{k}{2}} \Vert \psi_{j, k}^{\widehat{\mathcal{L}}}(D) m_{\widehat{\mathcal{L}}}(D) h_{\alpha}(t) \Vert_{L^{\frac{4}{3}}_x L^2_y} \\
&\lesssim t^{-\frac{1}{4}} 2^{\frac{j}{4}+\frac{k}{2}} \Vert \langle x \rangle^{\frac{1}{4}+\delta} \psi_{j,k}^{\widehat{\mathcal{L}}}(D) m_{\widehat{\mathcal{L}}}(D) h_{\alpha}(t) \Vert_{L^2} \\
&\lesssim t^{-\frac{1}{4}} 2^{\frac{j}{4}+\frac{k}{2}} \Vert \psi_{j, k}^{\widehat{\mathcal{L}}}(D) h_{\alpha}(t) \Vert_{L^2}^{\frac{3}{4}-\delta} \left( 2^{-j} \Vert \psi_j(D) h_{\alpha}(t) \Vert_{L^2} + 2^{-j} \Vert \nabla \psi_j(D) m_{\widehat{\mathcal{L}}}(D) x h_{\alpha}(t) \Vert_{L^2} \right)^{\frac{1}{4}+\delta} \\
&\lesssim t^{-\frac{1}{4}} 2^{-\delta j+\frac{k}{2}} \Vert \psi_{j, k}^{\widehat{\mathcal{L}}}(D) h_{\alpha}(t) \Vert_{L^2}^{\frac{3}{4}-\delta} \Vert u \Vert_X^{\frac{1}{4}+\delta} 
\end{align*}

We now put together the previous estimates, and Lemma \ref{lemestdispL4BF} for low frequencies. First, if $t \leq 1$, $\alpha$ arbitrary, 
\begin{align*}
&\Vert e^{-i t \omega(D)} m_{\widehat{\mathcal{L}}}(D) \nabla h_{\alpha}(t) \Vert_{L^4} \\
&\lesssim \sum_{\substack{j \in \mathbb{Z}, \\ 2^j \lesssim t^{-\frac{1}{3}}}} \Vert e^{-i t \omega(D)} \psi_j(D) \nabla h_{\alpha}(t) \Vert_{L^4} \quad + \sum_{\substack{j \in \mathbb{Z}, \\ 2^j \gg t^{-\frac{1}{3}}}} \sum_{\substack{k \in \mathbb{Z}, \\ t^{-\frac{1}{2}} 2^{-j} \lesssim 2^k, k \leq -10}} \Vert e^{-i t \omega(D)} \psi_{j, k}^{\widehat{\mathcal{L}}}(D) \nabla h_{\alpha}(t) \Vert_{L^4} \\
&\quad \quad + \sum_{\substack{j \in \mathbb{Z}, \\ 2^j \gg t^{-\frac{1}{3}}}} \sum_{\substack{k \in \mathbb{Z}, \\ 2^k \lesssim t^{-\frac{1}{2}} 2^{-j}, k \leq -10}} \Vert e^{-i t \omega(D)} \psi_{j, k}^{\widehat{\mathcal{L}}}(D) \nabla h_{\alpha}(t) \Vert_{L^4} \\
&\lesssim \sum_{\substack{j \in \mathbb{Z}, \\ 2^j \lesssim t^{-\frac{1}{3}}}} 2^{\frac{7j}{4}} \langle 2^j \rangle^{-1} \Vert u \Vert_X \quad + \sum_{\substack{j \in \mathbb{Z}, \\ 2^j \gg t^{-\frac{1}{3}}}} \sum_{\substack{k \in \mathbb{Z}, \\ t^{-\frac{1}{2}} 2^{-j} \lesssim 2^k, k \leq -10}} t^{-\frac{3}{4}} 2^{-\frac{3j}{4}-\frac{k}{2}-\delta k} \Vert u \Vert_X \\
&\quad \quad + \sum_{\substack{j \in \mathbb{Z}, \\ 2^j \gg t^{-\frac{1}{3}}}} \sum_{\substack{k \in \mathbb{Z}, \\ t^{-\frac{1}{2}} 2^{-j} \gtrsim 2^k, k \leq -10}} t^{-\frac{1}{4}} 2^{j-\delta j+\frac{k}{2}} \Vert \psi_j(D) h_{\alpha}(t) \Vert_{L^2}^{\frac{3}{4}-\delta} \Vert u \Vert_X^{\frac{1}{4}+\delta} \\
&\lesssim t^{-\frac{1}{4}} \Vert u \Vert_X \quad + \sum_{\substack{j \in \mathbb{Z}, \\ 2^j \gg t^{-\frac{1}{3}}}} t^{-\frac{1}{2}+\frac{\delta}{2}} 2^{-\frac{j}{4}+\delta j} \Vert u \Vert_X \\
&\quad \quad + \sum_{\substack{j \in \mathbb{Z}, \\ 2^j \gg t^{-\frac{1}{3}}}} t^{-\frac{1}{2}} 2^{-\frac{j}{4}} \Vert \nabla \psi_j(D) h_{\alpha}(t) \Vert_{L^2}^{\frac{3}{4}-\delta} \Vert u \Vert_X^{\frac{1}{4}+\delta} \\
&\lesssim t^{-\frac{5}{12}} \Vert u \Vert_X 
\end{align*}
as wanted. 

On the other hand, if $t \geq 1$, we differentiate the cases $\alpha \in \{ a, c \}$ and $\alpha = b$. 

First, if $\alpha \in \{ a, c \}$, we decompose: 
\begin{align*}
&\Vert e^{-i t \omega(D)} m_{\widehat{\mathcal{L}}}(D) \nabla h_{\alpha}(t) \Vert_{L^4} \\
&\lesssim \sum_{\substack{j \in \mathbb{Z}, \\ 2^j \lesssim t^{-\frac{1}{3}}}} \Vert e^{-i t \omega(D)} \psi_j(D) \nabla h_{\alpha}(t) \Vert_{L^4} \quad + \sum_{\substack{j \in \mathbb{Z}, \\ 2^j \gtrsim t^{-\frac{1}{3}}}} \sum_{\substack{k \in \mathbb{Z}, \\ 2^k \gtrsim t^{-\frac{1}{2}} 2^{-\frac{3j}{2}} \langle 2^j \rangle^{\frac{1}{2}}, \\
k \leq -10}} \Vert e^{-i t \omega(D)} \psi_{j, k}^{\widehat{\mathcal{L}}}(D) \nabla h_{\alpha}(t) \Vert_{L^4} \\
&\quad \quad + \sum_{\substack{j \in \mathbb{Z}, \\ 2^j \gtrsim t^{-\frac{1}{3}}}} \sum_{\substack{k \in \mathbb{Z}, \\ 2^k \lesssim t^{-\frac{1}{2}} 2^{-\frac{3j}{2}} \langle 2^j \rangle^{\frac{1}{2}}, \\
k \leq -10}} \Vert e^{-i t \omega(D)} \psi_{j, k}^{\widehat{\mathcal{L}}}(D) \nabla h_{\alpha}(t) \Vert_{L^4} \\
&\lesssim \sum_{\substack{j \in \mathbb{Z}, \\ 2^j \lesssim t^{-\frac{1}{3}}}} 2^{\frac{7j}{4}} \Vert u \Vert_X \quad + \sum_{\substack{j \in \mathbb{Z}, \\ 2^j \gtrsim t^{-\frac{1}{3}}}} \sum_{\substack{k \in \mathbb{Z}, \\ 2^k \gtrsim t^{-\frac{1}{2}} 2^{-\frac{3j}{2}} \langle 2^j \rangle^{\frac{1}{2}}, \\
k \leq -10}} t^{-\frac{3}{4}} 2^{-\frac{j}{2}-2\delta j-\frac{k}{2}-\delta k} \langle 2^j \rangle^{-\frac{1}{4}+2\delta} \Vert u \Vert_X \\
&\quad \quad + \sum_{\substack{j \in \mathbb{Z}, \\ 2^j \gtrsim t^{-\frac{1}{3}}}} \sum_{\substack{k \in \mathbb{Z}, \\ 2^k \lesssim t^{-\frac{1}{2}} 2^{-\frac{3j}{2}} \langle 2^j \rangle^{\frac{1}{2}}, \\
k \leq -10}} t^{-\frac{1}{4}} 2^{j-\delta j+\frac{k}{2}} \Vert \psi_{j, k}^{\widehat{\mathcal{L}}}(D) h_{\alpha}(t) \Vert_{L^2}^{\frac{3}{4}-\delta} \Vert u \Vert_X^{\frac{1}{4}+\delta} \\
&\lesssim t^{-\frac{7}{12}} \Vert u \Vert_X \quad + \sum_{\substack{j \in \mathbb{Z}, \\ 2^j \gtrsim t^{-\frac{1}{3}}}} t^{-\frac{1}{2}+\frac{\delta}{2}} 2^{\frac{j}{4}-\frac{\delta j}{2}} \langle 2^j \rangle^{-\frac{1}{2}+\frac{3\delta}{2}} \Vert u \Vert_X \\
&\quad \quad + \sum_{\substack{j \in \mathbb{Z}, \\ 2^j \gtrsim t^{-\frac{1}{3}}}} t^{-\frac{1}{2}} 2^{\frac{j}{4}-\delta j} \langle 2^j \rangle^{-\frac{1}{2}+\delta} \Vert \langle \nabla \rangle \psi_{j, k}^{\widehat{\mathcal{L}}}(D) h_{\alpha}(t) \Vert_{L^2}^{\frac{3}{4}-\delta} \Vert u \Vert_X^{\frac{1}{4}+\delta} \\
&\lesssim t^{-\frac{1}{2}+\frac{\delta}{2}} \Vert u \Vert_X 
\end{align*}
as wanted. 

For $\alpha = b$, we may reuse the previous estimates used for $\alpha = a, c$ if $2^{|j|} \gtrsim t^{\frac{1}{36}}$: indeed, the penultimate line becomes
\begin{align*}
\Vert e^{-i t \omega(D)} m_{\widehat{\mathcal{L}}}(D) \chi\left( t^{\frac{1}{36}} \left| \log_2(|\nabla|) \right|^{-1} \right) \nabla h_b(t) \Vert_{L^4} 
&\lesssim t^{-\frac{7}{12}} \Vert u \Vert_X \quad + \sum_{\substack{j \in \mathbb{Z}, \\ 2^j \gtrsim t^{-\frac{1}{3}}, \\
2^{|j|} \gtrsim t^{\frac{1}{36}}}} t^{-\frac{1}{2}} 2^{\frac{j}{4}-\frac{\delta j}{2}} \langle 2^j \rangle^{-\frac{1}{2}+\frac{3\delta}{2}} \Vert u \Vert_X \\
&\quad \quad + \sum_{\substack{j \in \mathbb{Z}, \\ 2^j \gtrsim t^{-\frac{1}{3}}, \\ 2^{|j|} \gtrsim t^{\frac{1}{36}}}} t^{-\frac{1}{2}} 2^{\frac{j}{4}-\delta j} \langle 2^j \rangle^{-\frac{1}{2}+\delta} \Vert u \Vert_X \\
&\lesssim t^{-\frac{1}{2}-\frac{1}{144}+O(\delta)} \Vert u \Vert_X 
\end{align*}
which is enough. 

It only remains to estimate the contribution of frequencies of size $2^j$ such that $2^{|j|} \lesssim t^{\frac{1}{36}}$. We use that
\begin{align*}
\Vert \psi_{j, k}^{\widehat{\mathcal{L}}}(D) h_b(t) \Vert_{L^2} &\lesssim \Vert \psi_{j, k}^{\widehat{\mathcal{L}}}(D) g_b(t) \Vert_{L^2} + \Vert \psi_{j, k}^{\widehat{\mathcal{L}}}(D) m_b(D) X_b f(t) \Vert_{L^2} \\
&\lesssim \langle 2^j \rangle^{-1} t^{\frac{1}{6}} \langle t \rangle^{-\frac{1}{4}+100\delta} \Vert u \Vert_X 
+ 2^{k-\frac{j}{2}} \langle 2^j \rangle^{-\frac{1}{2}} \Vert m_{\widehat{\mathcal{L}}}(D) \langle \nabla \rangle^{\frac{1}{2}} |\nabla|^{\frac{1}{2}} X_b f(t) \Vert_{L^2} \\
&\lesssim \langle t \rangle^{O(\delta)} \langle 2^j \rangle^{-\frac{1}{2}} \left( \langle t \rangle^{-\frac{1}{12}}
+ 2^{k-\frac{j}{2}} \right) \Vert u \Vert_X
\end{align*}
In particular, we either win a factor $k$, or a factor $\langle t \rangle^{-\frac{1}{12}}$. We can then decompose again: 
\begin{align*}
&\Vert e^{-i t \omega(D)} m_{\widehat{\mathcal{L}}}(D) \chi\left( t^{-\frac{1}{36}} \left| \log_2(|\nabla|) \right| \right) \nabla h_{\alpha}(t) \Vert_{L^4} \\
&\quad \lesssim \sum_{\substack{j \in \mathbb{Z}, \\ 2^j \gtrsim t^{-\frac{1}{36}}, \\ 2^j \lesssim t^{\frac{1}{36}}}} \sum_{\substack{k \in \mathbb{Z}, \\ 2^k \gtrsim t^{-\frac{1}{2}+\frac{1}{16}} 2^{-\frac{3j}{2}} \langle 2^j \rangle^{\frac{1}{2}}, \\ k \leq -10}} \Vert e^{-i t \omega(D)} \psi_{j, k}^{\widehat{\mathcal{L}}}(D) \nabla h_b(t) \Vert_{L^4} \\
&\quad + \sum_{\substack{j \in \mathbb{Z}, \\ 2^j \gtrsim t^{-\frac{1}{36}}, \\ 2^j \lesssim t^{\frac{1}{36}}}} \sum_{\substack{k \in \mathbb{Z}, \\ 2^k \lesssim t^{-\frac{1}{2}+\frac{1}{16}} 2^{-\frac{3j}{2}} \langle 2^j \rangle^{\frac{1}{2}}, \\ k \leq -10}} \Vert e^{-i t \omega(D)} \psi_{j, k}^{\widehat{\mathcal{L}}}(D) \nabla h_b(t) \Vert_{L^4} \\
&\quad \lesssim \sum_{\substack{j \in \mathbb{Z}, \\ 2^j \gtrsim t^{-\frac{1}{36}}, \\ 2^j \lesssim t^{\frac{1}{36}}}} \sum_{\substack{k \in \mathbb{Z}, \\ 2^k \gtrsim t^{-\frac{1}{2}+\frac{1}{16}} 2^{-\frac{3j}{2}} \langle 2^j \rangle^{\frac{1}{2}}, \\ k \leq -10}} t^{-\frac{3}{4}+O(\delta)} 2^{-\frac{j}{2}-2\delta j-\frac{k}{2}-\delta k} \langle 2^j \rangle^{-\frac{1}{4}+2\delta} \Vert u \Vert_X \\
&\quad \quad + \sum_{\substack{j \in \mathbb{Z}, \\ 2^j \gtrsim t^{-\frac{1}{36}}, \\ 2^j \lesssim t^{\frac{1}{36}}}} \sum_{\substack{k \in \mathbb{Z}, \\ 2^k \lesssim t^{-\frac{1}{2}+\frac{1}{16}} 2^{-\frac{3j}{2}} \langle 2^j \rangle^{\frac{1}{2}}, \\ k \leq -10}} t^{-\frac{1}{4}+O(\delta)} 2^{j-\delta j+\frac{k}{2}} \Vert \psi_{j, k}^{\widehat{\mathcal{L}}}(D) h_{\alpha}(t) \Vert_{L^2}^{\frac{3}{4}-\delta} \Vert u \Vert_X^{\frac{1}{4}+\delta} \\
&\quad \lesssim \sum_{\substack{j \in \mathbb{Z}, \\ 2^j \gtrsim t^{-\frac{1}{36}}, \\ 2^j \lesssim t^{\frac{1}{36}}}} \sum_{\substack{k \in \mathbb{Z}, \\ 2^k \gtrsim t^{-\frac{1}{2}+\frac{1}{16}} 2^{-\frac{3j}{2}} \langle 2^j \rangle^{\frac{1}{2}}, \\ k \leq -10}} 
t^{-\frac{3}{4}+O(\delta)} 2^{-\frac{j}{2}-\frac{k}{2}-\delta k} \langle 2^j \rangle^{-\frac{1}{4}} \Vert u \Vert_X \\
&\quad \quad + \sum_{\substack{j \in \mathbb{Z}, \\ 2^j \gtrsim t^{-\frac{1}{36}}, \\ 2^j \lesssim t^{\frac{1}{36}}}} \sum_{\substack{k \in \mathbb{Z}, \\ 2^k \lesssim t^{-\frac{1}{2}+\frac{1}{16}} 2^{-\frac{3j}{2}} \langle 2^j \rangle^{\frac{1}{2}}, \\ k \leq -10}} 
t^{-\frac{1}{4}-\frac{1}{16}+O(\delta)} 2^{j+\frac{k}{2}} \langle 2^j \rangle^{-\frac{3}{8}} \Vert u \Vert_X \\
&\quad \quad + \sum_{\substack{j \in \mathbb{Z}, \\ 2^j \gtrsim t^{-\frac{1}{36}}, \\ 2^j \lesssim t^{\frac{1}{36}}}} \sum_{\substack{k \in \mathbb{Z}, \\ 2^k \lesssim t^{-\frac{1}{2}+\frac{1}{16}} 2^{-\frac{3j}{2}} \langle 2^j \rangle^{\frac{1}{2}}, \\ k \leq -10}} 
t^{-\frac{1}{4}+O(\delta)} 2^{\frac{5j}{8}+\frac{5k}{4}-\delta k} \langle 2^j \rangle^{-\frac{3}{8}} \Vert u \Vert_X \\
&\quad \lesssim \sum_{\substack{j \in \mathbb{Z}, \\ 2^j \gtrsim t^{-\frac{1}{36}}, \\ 2^j \lesssim t^{\frac{1}{36}}}} 
t^{-\frac{1}{2}-\frac{1}{32}+O(\delta)} 2^{\frac{j}{4}} \langle 2^j \rangle^{-\frac{1}{2}} \Vert u \Vert_X 
\quad + \sum_{\substack{j \in \mathbb{Z}, \\ 2^j \gtrsim t^{-\frac{1}{36}}, \\ 2^j \lesssim t^{\frac{1}{36}}}} 
t^{-\frac{1}{2}-\frac{1}{32}+O(\delta)} 2^{\frac{j}{4}} \langle 2^j \rangle^{-\frac{1}{8}} \Vert u \Vert_X
\quad \\
&\quad \quad + \sum_{\substack{j \in \mathbb{Z}, \\ 2^j \gtrsim t^{-\frac{1}{36}}, \\ 2^j \lesssim t^{\frac{1}{36}}}} 
t^{-\frac{51}{64}+O(\delta)} 2^{-\frac{5j}{4}} \langle 2^j \rangle^{\frac{1}{4}} \Vert u \Vert_X \\
&\quad \lesssim t^{-\frac{1}{2}-\frac{1}{32}+O(\delta)} \Vert u \Vert_X + t^{-\frac{1}{2}-\frac{1}{32}+\frac{1}{36*8}+O(\delta)} \Vert u \Vert_X + t^{-\frac{51}{64}+\frac{5}{36*4}+O(\delta)} \Vert u \Vert_X \\
&\quad \lesssim t^{-\frac{1}{2}-\frac{1}{36}+O(\delta)} \Vert u \Vert_X 
\end{align*}
as wanted. 
\end{Dem}

\subsection{Neighborhood of the plane}

\begin{Lem} Let $t > 0$. Then
\begin{align*}
\Vert e^{-i t \omega(D)} m_{\widehat{\mathcal{P}}}(D) \nabla (h_a(t), h_c(t)) \Vert_{L^4} &\lesssim t^{-\frac{63}{128}} \langle t \rangle^{-\frac{1}{64}} \Vert u \Vert_X \\
\Vert e^{-i t \omega(D)} m_{\widehat{\mathcal{P}}}(D) \partial_x h_b(t) \Vert_{L^4} &\lesssim t^{-\frac{63}{128}} \langle t \rangle^{-\frac{1}{64}} \Vert u \Vert_X
\end{align*}
\end{Lem}

\begin{Dem}
On the one hand, we can estimate using Sobolev's embedding: 
\begin{align*}
&\Vert e^{-i t \omega(D)} \psi_{j, k}^{\widehat{\mathcal{P}}}(D) \nabla h_a(t) \Vert_{L^4} \lesssim 2^{\frac{5j}{4}+\frac{3k}{4}-\delta k - \delta j} \Vert \psi_{j, k}^{\widehat{\mathcal{P}}}(D) \nabla h_a(t) \Vert_{L^{\frac{1}{1-\delta}}_x L^2_y} \\
&\lesssim 2^{\frac{5j}{4}+\frac{3k}{4}-\delta k - \delta j} \Vert \psi_j(D) m_{\widehat{\mathcal{P}}}(D) \nabla h_a(t) \Vert_{L^{\frac{1}{1-\delta}}_x L^2_y} \\
&\lesssim 2^{\frac{5j}{4}+\frac{3k}{4}-\delta k - \delta j} \Vert \langle x \rangle^{\frac{1}{2}} \psi_j(D) m_{\widehat{\mathcal{P}}}(D) \nabla h_a(t) \Vert_{L^2} \\
&\lesssim 2^{\frac{5j}{4}+\frac{3k}{4}-\delta k - \delta j} \Vert \psi_j(D) \nabla h_a(t) \Vert_{L^2}^{\frac{1}{2}} \left( \Vert \psi_j(D) h_a(t) \Vert_{H^1} + \Vert m_{\widehat{\mathcal{P}}}(D) \nabla x h_a(t) \Vert_{L^2} \right)^{\frac{1}{2}} \\
&\lesssim 2^{\frac{7j}{4}+\frac{3k}{4}-\delta k - \delta j} \langle 2^j \rangle^{-\frac{1}{2}} \Vert u \Vert_X \\
\end{align*}

On the other hand, we already saw that if we set 
\begin{align*}
\widehat{K}_{t, j, k}(\overline{\xi}) &= e^{-i t \omega(\overline{\xi})} \psi_{j, k}^{\widehat{\mathcal{P}}}(\overline{\xi}) 
\end{align*}
then $\Vert K_{t, j, k} \Vert_{L^{\infty}} \lesssim t^{-\frac{3}{2}} 2^{-\frac{3j}{2}-\frac{k}{2}}$ if $t 2^{k+3j} \gtrsim 1$, so we can apply Young's convolution inequality and Riesz-Thorin interpolation theorem to get 
\begin{align*}
&\Vert e^{-i t \omega(D)} \psi_{j, k}^{\widehat{\mathcal{P}}}(D) \nabla h_a(t) \Vert_{L^4} \lesssim \left( t^{-\frac{3}{2}} 2^{-\frac{3j}{2}-\frac{k}{2}} \right)^{\frac{1}{2}} \Vert \psi_{j,k}^{\widehat{\mathcal{P}}}(D) \nabla h_a(t) \Vert_{L^{\frac{4}{3}}} \\
&\lesssim t^{-\frac{3}{4}} 2^{-\frac{3j}{4}-\frac{k}{4}} \Vert \psi_j(D) m_{\widehat{\mathcal{P}}}(D) \nabla h_a(t) \Vert_{L^{\frac{4}{3}}} \\
&\lesssim t^{-\frac{3}{4}} 2^{-\frac{3j}{4}-\frac{k}{4}} \Vert \langle x, y \rangle^{\frac{3}{4}+\delta} \psi_j(D) m_{\widehat{\mathcal{P}}}(D) \nabla h_a(t) \Vert_{L^2} \\
&\lesssim t^{-\frac{3}{4}} 2^{-\frac{3j}{4}-\frac{k}{4}} \Vert \psi_j(D) \nabla h_a(t) \Vert_{L^2}^{\frac{1}{4}-\delta} \left( \Vert \psi_j(D) h_a(t) \Vert_{H^1} + \Vert \psi_j(D) m_{\widehat{\mathcal{P}}}(D) \nabla (x, y) h_a(t) \Vert_{L^2} \right)^{\frac{3}{4}+\delta} \\
&\lesssim t^{-\frac{3}{4}} 2^{-\frac{j}{2}-\frac{k}{4}-\delta j} \langle 2^j \rangle^{-\frac{1}{4}+\delta} \Vert u \Vert_X
\end{align*}

First, if $t \geq 1$, we decompose: 
\begin{align*}
&\Vert e^{-i t \omega(D)} m_{\widehat{\mathcal{P}}}(D) \nabla h_a(t) \Vert_{L^4} \\
&\lesssim \sum_{\substack{j \in \mathbb{Z}, \\ 2^j \lesssim t^{-\frac{1}{3}}}} \Vert e^{-i t \omega(D)} \psi_j(D) \nabla h_a(t) \Vert_{L^4} \quad + 
\sum_{\substack{j \in \mathbb{Z}, \\ 2^j \gtrsim t^{-\frac{1}{3}}}} \sum_{\substack{k \in \mathbb{Z}, \\ 2^k \lesssim t^{-\frac{3}{4}} 2^{-\frac{9j}{4}}, \\ k \leq -10}} \Vert e^{-i t \omega(D)} \psi_{j, k}^{\widehat{\mathcal{P}}}(D) \nabla h_a(t) \Vert_{L^4} \\
&\quad \quad + \sum_{\substack{j \in \mathbb{Z}, \\ 2^j \gtrsim t^{-\frac{1}{3}}}} \sum_{\substack{k \in \mathbb{Z}, \\ 2^k \gtrsim t^{-\frac{3}{4}} 2^{-\frac{9j}{4}}, \\ k \leq -10}} \Vert e^{-i t \omega(D)} \psi_{j, k}^{\widehat{\mathcal{P}}}(D) \nabla h_a(t) \Vert_{L^4} \\
&\lesssim \sum_{\substack{j \in \mathbb{Z}, \\ 2^j \lesssim t^{-\frac{1}{3}}}} 2^{\frac{7j}{4}} \Vert u \Vert_X \quad + \sum_{\substack{j \in \mathbb{Z}, \\ 2^j \gtrsim t^{-\frac{1}{3}}}} \sum_{\substack{k \in \mathbb{Z}, \\ 2^k \lesssim t^{-\frac{3}{4}} 2^{-\frac{9j}{4}}, \\ k \leq -10}} 2^{\frac{7j}{4}+\frac{3k}{4}-\delta k - \delta j} \langle 2^j \rangle^{-\frac{1}{2}} \Vert u \Vert_X \\
&\quad \quad + \sum_{\substack{j \in \mathbb{Z}, \\ 2^j \gtrsim t^{-\frac{1}{3}}}} \sum_{\substack{k \in \mathbb{Z}, \\ 2^k \gtrsim t^{-\frac{3}{4}} 2^{-\frac{9j}{4}}, \\ k \leq -10}} t^{-\frac{3}{4}} 2^{-\frac{j}{2}-\frac{k}{4}-\delta j} \langle 2^j \rangle^{-\frac{1}{4}+\delta} \Vert u \Vert_X \\
&\lesssim t^{-\frac{7}{12}} \Vert u \Vert_X + \sum_{\substack{j \in \mathbb{Z}, \\ 2^j \gtrsim t^{-\frac{1}{3}}}} t^{-\frac{9}{16}+\frac{3\delta}{4}} 2^{\frac{j}{16}+\frac{5\delta j}{4}} \langle 2^j \rangle^{-\frac{1}{2}} \Vert u \Vert_X 
\quad + \sum_{\substack{j \in \mathbb{Z}, \\ 2^j \gtrsim t^{-\frac{1}{3}}}} t^{-\frac{9}{16}} 2^{\frac{j}{16}-\delta j} \langle 2^j \rangle^{-\frac{1}{4}+\delta} \Vert u \Vert_X \\
&\lesssim t^{-\frac{9}{16}+\frac{3\delta}{4}} \Vert u \Vert_X 
\end{align*}
as wanted. Here above, we used that if $2^j \gtrsim t^{-\frac{1}{3}}$ and $2^k \gtrsim t^{-\frac{3}{4}} 2^{-\frac{9j}{4}}$, then 
\begin{align*}
t 2^{k+3j} \gtrsim t^{\frac{1}{4}} 2^{\frac{3j}{4}} \gtrsim 1
\end{align*}
which allows to apply the previous estimates. 

For $t \leq 1$, we estimate as
\begin{align*}
&\Vert e^{-i t \omega(D)} m_{\widehat{\mathcal{P}}}(D) \nabla h_a(t) \Vert_{L^4} \\
&\lesssim \sum_{\substack{j \in \mathbb{Z}, \\ 2^j \lesssim t^{-\frac{9}{32}}}} \Vert e^{-i t \omega(D)} \psi_j(D) \nabla h_a(t) \Vert_{L^4} \quad + 
\sum_{\substack{j \in \mathbb{Z}, \\ 2^j \gtrsim t^{-\frac{9}{32}}}} \sum_{\substack{k \in \mathbb{Z}, \\ 2^k \lesssim t^{-\frac{3}{4}} 2^{-2j}, \\ k \leq -10}} \Vert e^{-i t \omega(D)} \psi_{j, k}^{\widehat{\mathcal{P}}}(D) \nabla h_a(t) \Vert_{L^4} \\
&\quad \quad + \sum_{\substack{j \in \mathbb{Z}, \\ 2^j \gtrsim t^{-\frac{9}{32}}}} \sum_{\substack{k \in \mathbb{Z}, \\ 2^k \gtrsim t^{-\frac{3}{4}} 2^{-2j}, \\ k \leq -10}} \Vert e^{-i t \omega(D)} \psi_{j, k}^{\widehat{\mathcal{P}}}(D) \nabla h_a(t) \Vert_{L^4} \\
&\lesssim \sum_{\substack{j \in \mathbb{Z}, \\ 2^j \lesssim t^{-\frac{9}{32}}}} 2^{\frac{7j}{4}} \langle 2^j \rangle^{-1} \Vert u \Vert_X \quad + \sum_{\substack{j \in \mathbb{Z}, \\ 2^j \gtrsim t^{-\frac{9}{32}}}} \sum_{\substack{k \in \mathbb{Z}, \\ 2^k \lesssim t^{-\frac{3}{4}} 2^{-2j}, \\ k \leq -10}} 2^{\frac{5j}{4}+\frac{3k}{4}-\delta k - \delta j} \Vert u \Vert_X \\
&\quad \quad + \sum_{\substack{j \in \mathbb{Z}, \\ 2^j \gtrsim t^{-\frac{9}{32}}}} \sum_{\substack{k \in \mathbb{Z}, \\ 2^k \gtrsim t^{-\frac{3}{4}} 2^{-2j}, \\ k \leq -10}} t^{-\frac{3}{4}} 2^{-\frac{3j}{4}-\frac{k}{4}} \Vert u \Vert_X \\
&\lesssim t^{-\frac{63}{128}} \Vert u \Vert_X + \sum_{\substack{j \in \mathbb{Z}, \\ 2^j \gtrsim t^{-\frac{9}{32}}}} t^{-\frac{9}{16}+\frac{3\delta}{4}} 2^{-\frac{j}{4}+\delta j} \Vert u \Vert_X 
\quad + \sum_{\substack{j \in \mathbb{Z}, \\ 2^j \gtrsim t^{-\frac{9}{32}}}} t^{-\frac{9}{16}} 2^{-\frac{j}{4}} \Vert u \Vert_X \\
&\lesssim t^{-\frac{63}{128}} \Vert u \Vert_X
\end{align*}
which concludes the estimate for $h_a$. 

The estimate for $h_c$ is exactly the same. Finally, for $h_b$, we replace the full derivative $\nabla$ by $\partial_x$: in particular, we can apply the same estimates since $\partial_x h_b$ satisfies the same bounds as $\nabla h_a$. 
\end{Dem}

\subsection{Neighborhood of the cone}

In the neighborhood of the cone, in order to separate properly the $X_b$ weight from the $X_a, X_c$ weights, it is needed to do by hand the stationary phase estimate and use the pseudo-differential localisation lemma \ref{lemlocpseudonontrivial}. 

\begin{Lem} Let $t > 0, j, k \in \mathbb{Z}$ be such that $2^j \gg t^{-\frac{1}{3}}$ and $k \leq -10$. Let $\alpha \in \{ a, b, c \}$. Then 
\begin{align*}
\Vert e^{-i t \omega(D)} \psi_{j, k}^{\widehat{\mathcal{C}}}(D) h_{\alpha}(t) \Vert_{L^4} &\lesssim t^{-\frac{15}{32}} \langle t \rangle^{\delta} 2^{-\frac{21j}{32}+\delta k+10\delta |j|} \Vert u \Vert_X 
\end{align*}
Moreover, 
\begin{align*}
\Vert e^{-i t \omega(D)} \psi_{j, k}^{\widehat{\mathcal{C}}}(D) h_b(t) \Vert_{L^4} &\lesssim t^{-\frac{3}{4}} 2^{-\frac{3j}{2}+\delta k+10\delta |j|} \Vert u \Vert_X + t^{-\frac{1}{2}-\frac{1}{1024}} 2^{-\frac{7j}{8}+\delta k} \langle 2^j \rangle^{-\frac{1}{4}} \Vert u \Vert_X \\
&\quad + \langle t \rangle^{\frac{1}{16}+103\delta} t^{-\frac{45}{64}} 2^{-j-\frac{23j}{64}+\delta k+10\delta |j|} \Vert u \Vert_X 
\end{align*} \label{lemestimeeL4displocalezoneC} 
\end{Lem}

Admitting this lemma, we can already show: 

\begin{Cor} Let $t > 0$. Then
\begin{align*}
\Vert e^{-i t \omega(D)} m_{\widehat{\mathcal{C}}}(D) \nabla h_b(t) \Vert_{L^4} &\lesssim t^{-\frac{31}{64}} \langle t \rangle^{-\frac{1}{64}-\frac{1}{1024}} \Vert u \Vert_X 
\end{align*}
\end{Cor}

\begin{Dem}
We decompose and apply Lemmas \ref{lemestdispL4BF} and \ref{lemestimeeL4displocalezoneC}: if $t \geq 1$, 
\begin{align*}
&\Vert e^{-i t \omega(D)} m_{\widehat{\mathcal{C}}}(D) \nabla h_b(t) \Vert_{L^4} \\
&\lesssim 
\sum_{\substack{j \in \mathbb{Z}, \\ 2^j \lesssim t^{-\frac{1}{3}}}} \Vert e^{-i t \omega(D)} \psi_j(D) \nabla h_b(t) \Vert_{L^4} \quad + \sum_{\substack{j \in \mathbb{Z}, \\ 2^j \gtrsim t^{-\frac{1}{3}}}} \sum_{\substack{k \in \mathbb{Z}, \\ k \leq -10}} \Vert e^{-i t \omega(D)} \psi_{j, k}^{\widehat{\mathcal{C}}}(D) \nabla h_b(t) \Vert_{L^4} \\
&\lesssim \sum_{\substack{j \in \mathbb{Z}, \\ 2^j \lesssim t^{-\frac{1}{3}}}} 2^{\frac{7j}{4}} \langle 2^j \rangle^{-1} \Vert u \Vert_X \\
&\quad \quad 
+ \sum_{\substack{j \in \mathbb{Z}, \\ 2^j \gtrsim t^{-\frac{1}{3}}}} \sum_{\substack{k \in \mathbb{Z}, \\ k \leq -10}} \left( t^{-\frac{3}{4}} 2^{-\frac{j}{2}+10\delta |j|} + t^{-\frac{1}{2}-\frac{1}{1024}} 2^{\frac{j}{8}} \langle 2^j \rangle^{-\frac{1}{4}} + \langle t \rangle^{\frac{1}{16}+103\delta} t^{-\frac{45}{64}} 2^{-\frac{23j}{64}+10\delta |j|} \right) 2^{\delta k} \Vert u \Vert_X \\
&\lesssim t^{-\frac{7}{12}} \Vert u \Vert_X 
+ \left( t^{-\frac{7}{12}+\frac{10\delta}{3}} + t^{-\frac{1}{2}-\frac{1}{1024}} + t^{-\frac{41}{64}+\frac{23}{192}+107\delta} \right) \Vert u \Vert_X \\
&\lesssim \left( t^{-\frac{1}{2}-\frac{1}{1024}} + t^{-\frac{25}{48}+107\delta} \right) \Vert u \Vert_X \\
&\lesssim t^{-\frac{1}{2}-\frac{1}{1024}} \Vert u \Vert_X
\end{align*}
as wanted. On the other hand, if $t \leq 1$, 
\begin{align*}
&\Vert e^{-i t \omega(D)} m_{\widehat{\mathcal{C}}}(D) \nabla h_b(t) \Vert_{L^4} \\
&\lesssim 
\sum_{\substack{j \in \mathbb{Z}, \\ 2^j \lesssim t^{-\frac{5}{8}}}} \Vert e^{-i t \omega(D)} \psi_j(D) \nabla h_b(t) \Vert_{L^4} \quad + \sum_{\substack{j \in \mathbb{Z}, \\ 2^j \gtrsim t^{-\frac{5}{8}}}} \sum_{\substack{k \in \mathbb{Z}, \\ k \leq -10}} \Vert e^{-i t \omega(D)} \psi_{j, k}^{\widehat{\mathcal{C}}}(D) \nabla h_b(t) \Vert_{L^4} \\
&\lesssim \sum_{\substack{j \in \mathbb{Z}, \\ 2^j \lesssim t^{-\frac{5}{8}}}} 2^{\frac{3j}{4}} \Vert u \Vert_X \quad 
+ \sum_{\substack{j \in \mathbb{Z}, \\ 2^j \gtrsim t^{-\frac{5}{8}}}} \sum_{\substack{k \in \mathbb{Z}, \\ k \leq -10}} \left( t^{-\frac{3}{4}} 2^{-\frac{j}{2}+10\delta |j|} + t^{-\frac{1}{2}-\frac{1}{1024}} 2^{-\frac{j}{8}} + t^{-\frac{45}{64}} 2^{-\frac{23j}{64}+10\delta |j|} \right) 2^{\delta k} \Vert u \Vert_X \\
&\lesssim t^{-\frac{15}{32}} \Vert u \Vert_X 
+ \left( t^{-\frac{3}{4}+\frac{5}{16}-\frac{50\delta}{8}} + t^{-\frac{1}{2}-\frac{1}{1024}+\frac{5}{64}} + t^{-\frac{45}{64}+\frac{23*5}{64*8}-\frac{50\delta}{8}} \right) \Vert u \Vert_X \\
&\lesssim t^{-\frac{15}{32}} \Vert u \Vert_X 
+ \left( t^{-\frac{7}{16}-\frac{50\delta}{8}} + t^{-\frac{7}{16}} + t^{-\frac{245}{512}-\frac{50\delta}{8}} \right) \Vert u \Vert_X \\
&\lesssim t^{-\frac{31}{64}} \Vert u \Vert_X 
\end{align*}
as wanted. 
\end{Dem}

We now prove Lemma \ref{lemestimeeL4displocalezoneC}. 

\begin{Dem}
By Fourier's inversion formula, for $(x, y) \in \mathbb{R}^3$: 
\begin{align}
e^{-i t \omega(D)} \psi_{j, k}^{\widehat{\mathcal{C}}}(D) h_{\alpha}(t, x, y) &= \int e^{-i t \Phi\left( \frac{x}{t}, \frac{y}{t}, \overline{\xi} \right)} \psi_{j, k}^{\widehat{\mathcal{C}}}(\overline{\xi}) \widehat{h}_{\alpha}(t, \overline{\xi}) ~ d\overline{\xi} \label{preuveestimeeL4CformuleFourier} 
\end{align}
where $\Phi$ is the usual phase
\begin{align*}
\Phi\left( \frac{x}{t}, \frac{y}{t}, \overline{\xi} \right) &= \xi_0^3 + \xi_0 |\xi|^2 - \frac{\xi_0 x}{t} - \frac{\xi \cdot y}{t} 
\end{align*}

To estimate \eqref{preuveestimeeL4CformuleFourier} in $L^4_{x, y}$, we will separate the $L^4$ estimate depending on the support in $(x, y)$, and notably the direction of $(x, y)$. Furthermore, at fixed $(x, y)$, we can apply integrations by parts on \eqref{preuveestimeeL4CformuleFourier} as for the $L^{\infty}$ dispersive estimate, and then obtain an estimate applying interpolation and Lemma \ref{lemlocpseudonontrivial}. 

Recall the notation $\mathfrak{t} = t 2^{3j}$ and $\overline{\xi_a} = \sqrt{\frac{|x|+|y|}{t}}$. 

For $n \in \mathbb{Z}$ such that $2^n \gg 2^j$, we set
\begin{align*}
A_{HF} &= \{ (x, y), ~ \overline{\xi_a} \ll 2^j \} \\
A_n^{\mathcal{R}} &= \{ (x, y), ~ \overline{\xi_a} \simeq 2^n, ~ |y| \gtrsim |x|+|y|, ~ \left| |x| - \sqrt{3} |y| \right| \gtrsim |x|+|y| \} \\
A_n^{\mathcal{C}} &= \{ (x, y), ~ \overline{\xi_a} \simeq 2^n, ~ \left| |x| - \sqrt{3} |y| \right| \ll |x|+|y| \} \\
A_{res}^{\mathcal{R}} &= \{ (x, y), ~ \overline{\xi_a} \simeq 2^j, ~ |y| \gtrsim |x|+|y|, ~ \left| |x| - \sqrt{3} |y| \right| \gtrsim |x|+|y| \} \\
A_{res}^{\mathcal{C}} &= \{ (x, y), ~ \overline{\xi_a}| \simeq 2^j, ~ \left| |x| - \sqrt{3} |y| \right| \ll |x|+|y| \} \\
A_{res}^{\mathcal{L}} &= \{ (x, y), ~ \overline{\xi_a} \gtrsim 2^j, ~ |y| \ll |x| \}
\end{align*}

Then we can decompose: 
\begin{align*}
\Vert \eqref{preuveestimeeL4CformuleFourier} \Vert_{L^4(\mathbb{R}^3)} &\lesssim \Vert \eqref{preuveestimeeL4CformuleFourier} \Vert_{L^4(A_{HF})} + \Vert \eqref{preuveestimeeL4CformuleFourier} \Vert_{L^4(A_{res}^{\mathcal{R}})} + \Vert \eqref{preuveestimeeL4CformuleFourier} \Vert_{L^4(A_{res}^{\mathcal{C}})} + \Vert \eqref{preuveestimeeL4CformuleFourier} \Vert_{L^4(A_{res}^{\mathcal{L}})} \\
&\quad \quad + \sum_{\substack{n \in \mathbb{Z}, \\ n \geq j+10}} \left( \Vert \eqref{preuveestimeeL4CformuleFourier} \Vert_{L^4(A_n^{\mathcal{R}})} + \Vert \eqref{preuveestimeeL4CformuleFourier} \Vert_{L^4(A_n^{\mathcal{C}})} \right) 
\end{align*}

\paragraph{High frequencies} Let us start by estimating $L^4(A_{HF})$. 

If $(x, y) \in A_1$, $|\overline{\xi}| \simeq 2^j$, we have $\widehat{X}_a(\overline{\xi}) \cdot \nabla_{\overline{\xi}} \Phi \simeq 2^{2j}$, so we can apply an integration by parts: 
\begin{align*}
\eqref{preuveestimeeL4CformuleFourier} &= t^{-1} \int e^{-i t \Phi\left( \frac{x}{t}, \frac{y}{t}, \overline{\xi} \right)} \nabla_{\overline{\xi}} \cdot \left( \frac{\widehat{X}_a(\overline{\xi})}{\widehat{X}_a(\overline{\xi}) \cdot \nabla_{\overline{\xi}} \Phi} \psi_{j, k}^{\widehat{\mathcal{C}}}(\overline{\xi}) \right) \widehat{h}_{\alpha}(t, \overline{\xi}) ~ d\overline{\xi} \\
&\quad + t^{-1} \int e^{-i t \Phi\left( \frac{x}{t}, \frac{y}{t}, \overline{\xi} \right)} \frac{1}{\widehat{X}_a(\overline{\xi}) \cdot \nabla_{\overline{\xi}} \Phi} \psi_{j, k}^{\widehat{\mathcal{C}}}(\overline{\xi}) \widehat{X}_a(\overline{\xi}) \cdot \nabla_{\overline{\xi}} \widehat{h}_{\alpha}(t, \overline{\xi}) ~ d\overline{\xi}
\end{align*}
On the support of $\psi_{j, k}^{\widehat{\mathcal{C}}}$ and for this choice of $(x, y)$, we have $\widehat{X}_{\beta}(\overline{\xi}) \cdot \nabla_{\overline{\xi}} \Phi \lesssim 2^{2j}$ for every $\beta \in \{ a, b, c \}$: in particular, we may add everywhere a symbol
\begin{align*}
\widetilde{\psi}\left( 2^{-2j} \widehat{X}_a(\overline{\xi}) \cdot \nabla_{\overline{\xi}} \Phi \right) \widetilde{\chi}\left( 2^{-2j} \widehat{X}_b(\overline{\xi}) \cdot \nabla_{\overline{\xi}} \Phi \right) \widetilde{\chi}\left( 2^{-2j} \widehat{X}_b(\overline{\xi}) \cdot \nabla_{\overline{\xi}} \Phi \right) 
\end{align*}
without changing the value, where $\widetilde{\psi}$, $\widetilde{\chi}$ have similar properties as $\psi, \chi$, but potentially larger supports. In particular, we may now absorb other symbols into this one and apply Lemma \ref{lemlocpseudonontrivial}. We can therefore estimate: 
\begin{align*}
\Vert \eqref{preuveestimeeL4CformuleFourier} \Vert_{L^4(A_1)} &\lesssim \Vert \eqref{preuveestimeeL4CformuleFourier} \Vert_{L^2}^{\frac{1}{2}} \Vert \eqref{preuveestimeeL4CformuleFourier} \Vert_{L^{\infty}}^{\frac{1}{2}} \\
&\lesssim t^{-1} 2^{-3j} \left( \Vert h_{\alpha}(t) \Vert_{L^2} + \Vert \nabla X_a h_{\alpha}(t) \Vert_{L^2} \right)^{\frac{1}{2}} \left( \Vert \psi_{j, k}^{\widehat{\mathcal{C}}} \Vert_{L^2} \left( \Vert h_{\alpha}(t) \Vert_{L^2} + \Vert \nabla X_a h_{\alpha}(t) \Vert_{L^2} \right) \right)^{\frac{1}{2}} \\
&\lesssim t^{-1} 2^{-\frac{9j}{4}+\frac{k}{4}} \Vert u \Vert_X 
\end{align*}
which is enough for any $\alpha$. 

\paragraph{Resonant frequencies, away from the line and the cone} Let us consider the norm $L^4(A_{res}^{\mathcal{R}})$. 

We now introduce, as for the $L^{\infty}$ dispersive estimate, the localisation symbols 
\begin{align*}
\Psi_{(l_a, l_b, l_c)}^a\left( \overline{\xi}, \frac{x}{t}, \frac{y}{t} \right) &= \psi\left( 2^{-2j-l_a} \widehat{X}_a(\overline{\xi}) \cdot \nabla_{\overline{\xi}} \Phi \right) \chi\left( 2^{-2j-l_b} \widehat{X}_b(\overline{\xi}) \cdot \nabla_{\overline{\xi}} \Phi \right) \chi\left( 2^{-2j-l_c} \widehat{X}_c(\overline{\xi}) \cdot \nabla_{\overline{\xi}} \Phi \right) \\
\Psi_{(l_a, l_b, l_c)}^b\left( \overline{\xi}, \frac{x}{t}, \frac{y}{t} \right) &= \chi\left( 2^{-2j-l_a} \widehat{X}_a(\overline{\xi}) \cdot \nabla_{\overline{\xi}} \Phi \right) \psi\left( 2^{-2j-l_b} \widehat{X}_b(\overline{\xi}) \cdot \nabla_{\overline{\xi}} \Phi \right) \chi\left( 2^{-2j-l_c} \widehat{X}_c(\overline{\xi}) \cdot \nabla_{\overline{\xi}} \Phi \right) \\
\Psi_{(l_a, l_b, l_c)}^c\left( \overline{\xi}, \frac{x}{t}, \frac{y}{t} \right) &= \chi\left( 2^{-2j-l_a} \widehat{X}_a(\overline{\xi}) \cdot \nabla_{\overline{\xi}} \Phi \right) \chi\left( 2^{-2j-l_b} \widehat{X}_b(\overline{\xi}) \cdot \nabla_{\overline{\xi}} \Phi \right) \psi\left( 2^{-2j-l_c} \widehat{X}_c(\overline{\xi}) \cdot \nabla_{\overline{\xi}} \Phi \right) \\
\Psi_{(l_a, l_b, l_c)}^{int}\left( \overline{\xi}, \frac{x}{t}, \frac{y}{t} \right) &= \chi\left( 2^{-2j-l_a} \widehat{X}_a(\overline{\xi}) \cdot \nabla_{\overline{\xi}} \Phi \right) \chi\left( 2^{-2j-l_b} \widehat{X}_b(\overline{\xi}) \cdot \nabla_{\overline{\xi}} \Phi \right) \chi\left( 2^{-2j-l_c} \widehat{X}_c(\overline{\xi}) \cdot \nabla_{\overline{\xi}} \Phi \right) 
\end{align*}

Then, we decompose
\begin{subequations}
\begin{align}
\eqref{preuveestimeeL4CformuleFourier} &= \sum_{l = l_0}^0 \int e^{-i t \Phi\left( \frac{x}{t}, \frac{y}{t}, \overline{\xi} \right)} \Psi_{(l, 0, l)}^a \psi_{j, k}^{\widehat{\mathcal{C}}}(\overline{\xi}) \widehat{h}_{\alpha}(t, \overline{\xi}) ~ d\overline{\xi} \label{estdispL4CresR-a} \\
&\quad + \sum_{l = l_0}^0 \int e^{-i t \Phi\left( \frac{x}{t}, \frac{y}{t}, \overline{\xi} \right)} \Psi_{(l, 0, l)}^c \psi_{j, k}^{\widehat{\mathcal{C}}}(\overline{\xi}) \widehat{h}_{\alpha}(t, \overline{\xi}) ~ d\overline{\xi} \label{estdispL4CresR-c} \\
&\quad + \int e^{-i t \Phi\left( \frac{x}{t}, \frac{y}{t}, \overline{\xi} \right)} \Psi_{(l_0, 0, l_0)}^{int} \psi_{j, k}^{\widehat{\mathcal{C}}}(\overline{\xi}) \widehat{h}_{\alpha}(t, \overline{\xi}) ~ d\overline{\xi} \label{estdispL4CresR-int}
\end{align}
\end{subequations}
for $l_0$ such that $2^{l_0} \simeq \mathfrak{t}^{-\frac{1}{2}+\varrho}$ for some $\varrho > 0$ small enough.  

Note that, if $\left( \widehat{X}_a, \widehat{X}_c \right) \cdot \nabla_{\overline{\xi}} \Phi \ll 2^{2j}$, then necessarily $\widehat{X}_b \cdot \nabla_{\overline{\xi}} \Phi \simeq 2^{2j}$. 

As in the proof of the $L^{\infty}$ dispersive estimate, we can estimate the anisotropic volumes: 
\begin{align*}
\Vert \Psi_{(l, 0, l)} \psi_{j, k} \Vert_{L^{\infty}(A_{res}^{\mathcal{R}}, L^{p_a}_a L^{p_b}_b L^{p_c}_c)} &\lesssim 2^{\frac{l}{p_a}+\frac{l}{p_c}+\frac{k}{p_b}} 2^{\frac{j}{p_a}+\frac{j}{p_b}+\frac{j}{p_c}}
\end{align*}
see \eqref{estimeevolumiquePsilacasCloinLloinC}. 

On \eqref{estdispL4CresR-a}, we apply at most $n$, large with respect to $\varrho^{-1}$ (fixed), integrations by parts along $\widehat{X}_a$: 
\begin{subequations}
\begin{align}
\eqref{estdispL4CresR-a} &= \sum_{l = l_0}^0 t^{-n} 2^{-2nl-3nj} \int e^{-i t \Phi\left( \frac{x}{t}, \frac{y}{t}, \overline{\xi} \right)} \Psi_{(l, 0, l)}^a \psi_{j, k}^{\widehat{\mathcal{C}}}(\overline{\xi}) \widehat{h}_{\alpha}(t, \overline{\xi}) ~ d\overline{\xi} \label{estdispL4CresR-a-1} \\
&\quad + \sum_{i = 1}^n \sum_{l = l_0}^0 t^{-i} 2^{-2il+l-3ij+j} \int e^{-i t \Phi\left( \frac{x}{t}, \frac{y}{t}, \overline{\xi} \right)} \Psi_{(l, 0, l)}^a \psi_{j, k}^{\widehat{\mathcal{C}}}(\overline{\xi}) \widehat{X}_a(\overline{\xi}) \cdot \nabla_{\overline{\xi}} \widehat{h}_{\alpha}(t, \overline{\xi}) ~ d\overline{\xi} \label{estdispL4CresR-a-2} 
\end{align}
\end{subequations}
where the symbols may change as long as they keep similar properties. We can then estimate by interpolation and applying Lemma \ref{lemlocpseudonontrivial}, and denoting by $\eqref{estdispL4CresR-a-1}_l$ the $l$-th term of the sum: 
\begin{align*}
\Vert \eqref{estdispL4CresR-a-1} \Vert_{L^4(A_{res}^{\mathcal{R}})} &\lesssim \sum_{l = l_0}^0 \Vert \eqref{estdispL4CresR-a-1}_l \Vert_{L^2}^{\frac{1}{2}} \Vert \eqref{estdispL4CresR-a-1}_l \Vert_{L^{\infty}(A_{res}^{\mathcal{R}})}^{\frac{1}{2}} \\
&\lesssim \sum_{l = l_0}^0 t^{-n} 2^{-2nl-3nj} \Vert h_{\alpha}(t) \Vert_{L^2}^{\frac{1}{2}} \left( \Vert \Psi_{(l, 0, l)}^a \psi_{j, k} \Vert_{L^{\infty}(A_{res}^{\mathcal{R}}, L^{\frac{1}{1-\kappa}}_{a, c} L^2_b)} \Vert \psi_j \widehat{h}_{\alpha}(t) \Vert_{L^{\frac{1}{\kappa}}_{a, c} L^2_b} \right)^{\frac{1}{2}} \\
&\lesssim \sum_{l = l_0}^0 \mathfrak{t}^{-n} 2^{-2nl+l-\kappa l+\frac{k}{4}+\frac{5j}{4}-\kappa j} \Vert h_{\alpha}(t) \Vert_{L^2}^{\frac{1}{2}} \Vert \psi_j \widehat{h}_{\alpha}(t) \Vert_{H^1_{a, c} L^2_b}^{\frac{1}{2}} \\
&\lesssim \mathfrak{t}^{-2n \varrho - \frac{1}{2}+\varrho+\frac{\kappa}{2}} 2^{\frac{k}{4}+\frac{3j}{4}-\kappa j} \Vert u \Vert_X \\
&\lesssim t^{-1} 2^{-\frac{9j}{4}+\frac{k}{4}-\kappa j} \Vert u \Vert_X \\
\Vert \eqref{estdispL4CresR-a-2} \Vert_{L^4(A_{res}^{\mathcal{R}})} &\lesssim \sum_{i = 1}^n \sum_{l = l_0}^0 \Vert \eqref{estdispL4CresR-a-2}_l \Vert_{L^2}^{\frac{1}{2}} \Vert \eqref{estdispL4CresR-a-2}_l \Vert_{L^{\infty}(A_{res}^{\mathcal{R}})}^{\frac{1}{2}} \\
&\lesssim \sum_{i = 1}^n \sum_{l = l_0}^0 \mathfrak{t}^{-i} 2^{-2il+l} \Vert \nabla X_a h_{\alpha}(t) \Vert_{L^2}^{\frac{1}{2}} \left( \Vert \Psi_{(l, 0, l)} \psi_{j, k} \Vert_{L^{\infty}(A_{res}^{\mathcal{R}}, L^2)} \Vert \nabla X_a h_{\alpha}(t) \Vert_{L^2}^{\frac{1}{2}} \right)^{\frac{1}{2}} \\
&\lesssim \sum_{i = 1}^n \sum_{l = l_0}^0 \mathfrak{t}^{-i} 2^{-2il+\frac{3l}{2}+\frac{k}{4}+\frac{3j}{4}} \Vert u \Vert_X \\
&\lesssim \mathfrak{t}^{-\frac{3}{4}} 2^{\frac{k}{4}+\frac{3j}{4}} \Vert u \Vert_X \\
&\lesssim t^{-\frac{3}{4}} 2^{-\frac{3j}{2}+\frac{k}{4}} \Vert u \Vert_X 
\end{align*}
for $n \varrho \gg 1$, where $\kappa > 0$ is a small enough parameter and where we used $\mathfrak{t} \gg 1$. 

We proceed the same way for \eqref{estdispL4CresR-c} replacing $\widehat{X}_a$ by $\widehat{X}_c$. 

Then, for \eqref{estdispL4CresR-int}, we consider two cases depending on $k$. First, if $2^k \lesssim \mathfrak{t}^{-\frac{1}{64}}$, we can estimate directly:  
\begin{align*}
\Vert \eqref{estdispL4CresR-int} \Vert_{L^4(A_{res}^{\mathcal{R}})} &\lesssim \Vert \eqref{estdispL4CresR-int} \Vert_{L^2}^{\frac{1}{2}} \Vert \eqref{estdispL4CresR-int} \Vert_{L^{\infty}(A_{res}^{\mathcal{R}})}^{\frac{1}{2}} \\
&\lesssim \left( \langle 2^j \rangle^{-1} \Vert \langle \nabla \rangle h_{\alpha}(t) \Vert_{L^2} \right)^{\frac{1}{2}} \left( \Vert \Psi_{(l, 0, l)}^{int} \psi_{j, k} \Vert_{L^{\infty}(A_{res}^{\mathcal{R}}, L^{\frac{1}{1-\kappa}}_{a, c} L^2_b)} \Vert \psi_j \widehat{h}_{\alpha}(t) \Vert_{L^{\frac{1}{\kappa}}_{a, c} L^2_b} \right)^{\frac{1}{2}} \\
&\lesssim 2^{l-\kappa l+\frac{k}{4}+\frac{5j}{4}-\kappa j} \langle 2^j \rangle^{-\frac{1}{2}} \Vert u \Vert_X^{\frac{1}{2}} \Vert \psi_j \widehat{h}_{\alpha}(t) \Vert_{H^1_{a, c} L^2_b}^{\frac{1}{2}} \\
&\lesssim \mathfrak{t}^{-\frac{1}{2}+\varrho+\frac{\kappa}{2}-\frac{1}{512}} 2^{\frac{k}{8}+\frac{3j}{4}-\kappa j} \langle 2^j \rangle^{-\frac{1}{2}} \Vert u \Vert_X \\
&\lesssim t^{-\frac{1}{2}-\frac{1}{1024}} 2^{\frac{k}{8}-\frac{7j}{8}} \langle 2^j \rangle^{-\frac{1}{4}} \Vert u \Vert_X
\end{align*}
which is enough. 

On the other hand, if $2^k \gtrsim \mathfrak{t}^{-\frac{1}{64}}$, since we can assume $l_0 \leq -10$ here, we can apply an integration by parts along $\widehat{X}_b(\overline{\xi})$: 
\begin{subequations}
\begin{align}
\eqref{estdispL4CresR-int} &= t^{-1} 2^{-2j} \left( 2^{-l_0-j} + 2^{-k-j} \right) \int e^{-i t \Phi\left( \frac{x}{t}, \frac{y}{t}, \overline{\xi} \right)} \Psi_{(l_0, 0, l_0)}^{int} \psi_{j, k}^{\widehat{\mathcal{C}}}(\overline{\xi}) \widehat{h}_{\alpha}(t, \overline{\xi}) ~ d\overline{\xi} \label{estdispL4CresR-int-1} \\
&\quad + t^{-1} 2^{-2j-k} \int e^{-i t \Phi\left( \frac{x}{t}, \frac{y}{t}, \overline{\xi} \right)} \Psi_{(l_0, 0, l_0)}^{int} \psi_{j, k}^{\widehat{\mathcal{C}}}(\overline{\xi}) m_b(\overline{\xi}) \widehat{X}_b(\overline{\xi}) \cdot \nabla_{\overline{\xi}} \widehat{h}_{\alpha}(t, \overline{\xi}) ~ d\overline{\xi} \label{estdispL4CresR-int-2} 
\end{align}
\end{subequations}
We then estimate: 
\begin{align*}
\Vert \eqref{estdispL4CresR-int-1} \Vert_{L^4(A_{res}^{\mathcal{R}})} &\lesssim \Vert \eqref{estdispL4CresR-int-1} \Vert_{L^2}^{\frac{1}{2}} \Vert \eqref{estdispL4CresR-int-1} \Vert_{L^{\infty}(A_{res}^{\mathcal{R}})}^{\frac{1}{2}} \\
&\lesssim \mathfrak{t}^{-1} \left( 2^{-l_0} + 2^{-k} \right) \Vert h_{\alpha}(t) \Vert_{L^2}^{\frac{1}{2}} \left( \Vert \Psi_{(l_0, 0, l_0)} \psi_{j, k}^{\widehat{\mathcal{C}}} \Vert_{L^{\infty}(A_{res}^{\mathcal{R}}, L^{\frac{1}{1-\kappa}}_{a, c} L^2_b)} \Vert \psi_j \widehat{h}_{\alpha}(t) \Vert_{L^{\frac{1}{\kappa}}_{a, c} L^2_b} \right)^{\frac{1}{2}} \\
&\lesssim \mathfrak{t}^{-1} 2^{-\kappa l_0+\frac{k}{4}+\frac{3j}{4}-\kappa j} \Vert u \Vert_X \\
&\lesssim \mathfrak{t}^{-1+\frac{\kappa}{2}} 2^{\frac{k}{4}+\frac{3j}{4}-\kappa j} \Vert u \Vert_X \\
&\lesssim t^{-\frac{3}{4}} 2^{-\frac{3j}{2}-\kappa j+\frac{k}{4}} \Vert u \Vert_X \\
\Vert \eqref{estdispL4CresR-int-2} \Vert_{L^4(A_{res}^{\mathcal{R}})} &\lesssim \Vert \eqref{estdispL4CresR-int-2} \Vert_{L^2}^{\frac{1}{2}} \Vert \eqref{estdispL4CresR-int-2} \Vert_{L^{\infty}(A_{res}^{\mathcal{R}})}^{\frac{1}{2}} \\
&\lesssim \mathfrak{t}^{-1} 2^{-k} \Vert \nabla X_a h_{\alpha}(t) \Vert_{L^2}^{\frac{1}{2}} \left( \Vert \Psi_{(l_0, 0, l_0)} \psi_{j, k}^{\widehat{\mathcal{C}}} \Vert_{L^{\infty}(A_{res}^{\mathcal{R}}, L^2)} \Vert \partial_x m_b X_b h_{\alpha}(t) \Vert_{L^2} \right)^{\frac{1}{2}} \\
&\lesssim \mathfrak{t}^{-1} 2^{\frac{l_0}{2}-\frac{3k}{4}+\frac{3j}{4}} \Vert u \Vert_X \\
&\lesssim \mathfrak{t}^{-1} 2^{\frac{k}{4}+\frac{3j}{4}} \Vert u \Vert_X \\
&\lesssim t^{-\frac{3}{4}} 2^{-\frac{3j}{2}+\frac{k}{4}} \Vert u \Vert_X
\end{align*}
which is enough. 

\paragraph{Resonant frequencies, near the cone} Let us consider the norm $L^4(A_{res}^{\mathcal{C}})$. 

In this case, we can reuse the same strategy as previously, but this time we may have 
\begin{align*}
\widehat{X}_{\beta}(\overline{\xi}) \cdot \nabla_{\overline{\xi}} \Phi \ll 2^{2j} 
\end{align*}
for every $\beta \in \{ a, b, c \}$. The beginning of the estimate is exactly the same and we can restrict our attention to \eqref{estdispL4CresR-int} only. Furthermore, we can also estimate in the case $2^k \lesssim \mathfrak{t}^{-\frac{1}{64}}$ as before, so we assume now $2^k \gtrsim \mathfrak{t}^{-\frac{1}{64}}$. 

We then separate the $L^4$ norm in three pieces, depending on the size of $\tau := \sqrt{\frac{||x| - \sqrt{3} |y||}{t}}$ with respect to $2^k$. 

We have
\begin{align*}
\widehat{X}_b(\overline{\xi}) \cdot \nabla_{\overline{\xi}} \Phi &= \widehat{X}_b'(\overline{\xi}) \cdot \nabla_{\overline{\xi}} \Phi + O\left( \widehat{X}_a(\overline{\xi}) \cdot \nabla_{\overline{\xi}} \Phi \right) + O\left( \widehat{X}_c(\overline{\xi}) \cdot \nabla_{\overline{\xi}} \Phi \right) \\
&= \frac{\xi_0}{2 |\xi_0|} \left( \sqrt{3} |\xi_0| - |\xi| \right)^2 - \frac{\xi_0}{|\xi_0|} \left( x - \sqrt{3} \frac{\xi_0}{|\xi_0|} |y| \right) + O\left( 2^{2j+l_0} \right) 
\end{align*}
In particular, if $2^k \gg \tau$, then this quantity is of order $2^{2k+2j}$ without needing further localisation. If $2^k \ll \tau$, then this quantity is of order $\tau^2 2^{2j} \gg 2^{2k+2j}$. In these two cases, we can apply an integration by parts just like in the case away from the cone treated before, and estimate the same way up to losing a factor $2^{-2k}$ everywhere: but precisely as we assumed here $2^{2k} \gtrsim \mathfrak{t}^{-\frac{1}{32}}$, it is easy to see that the same estimates provide a sufficient bound. 

Finally, if $\tau \simeq 2^k$, we can localise also using $l_b$ and decompose into: 
\begin{subequations}
\begin{align}
\eqref{estdispL4CresR-int} &= \sum_{l = l_1}^{2k} \int e^{-i t \Phi\left( \frac{x}{t}, \frac{y}{t}, \overline{\xi} \right)} \Psi_{(l_0, l, l_0)}^b \psi_{j, k}^{\widehat{\mathcal{C}}}(\overline{\xi}) \widehat{h}_{\alpha}(t, \overline{\xi}) ~ d\overline{\xi} \label{estdispL4CresC-int-b} \\
&\quad + \int e^{-i t \Phi\left( \frac{x}{t}, \frac{y}{t}, \overline{\xi} \right)} \Psi_{(l_0, l_1, l_0)}^{int} \psi_{j, k}^{\widehat{\mathcal{C}}}(\overline{\xi}) \widehat{h}_{\alpha}(t, \overline{\xi}) ~ d\overline{\xi} \label{estdispL4CresC-int-int} 
\end{align}
\end{subequations}
for $l_1$ such that $2^{l_1} \simeq \mathfrak{t}^{-\frac{1}{16}}$. 
Note that we have the volume estimate
\begin{align*}
\Vert \Psi_{(l_0, l, l_0)}^{*} \psi_{j, k}^{\widehat{\mathcal{C}}} \Vert_{L^{\infty}(A_{res}^{\mathcal{C}}, L^{p_c}_c L^{p_b}_b L^{p_a}_a)} &\lesssim 2^{\frac{l_0}{p_a}+\frac{l_0}{p_c}+\frac{l-k}{p_b}} 2^{\frac{j}{p_a}+\frac{j}{p_b}+\frac{j}{p_c}} 
\end{align*}
obtained in the proof of the dispersive $L^{\infty}$ estimate, see \eqref{estvolumiquePsilalbcasCCtausimk}. (Note that the symbols used in \eqref{estvolumiquePsilalbcasCCtausimk} were slightly different, as $m_b \widehat{X}_b \cdot \nabla \Phi$ was localised instead of $\widehat{X}_b \cdot \nabla \Phi$ here, hence a gain of a factor $2^k$ in the estimate above.) 

On \eqref{estdispL4CresC-int-b}, we apply an integration by parts along $\widehat{X}_b$: 
\begin{subequations}
\begin{align}
\eqref{estdispL4CresC-int-b} &= \sum_{l = l_1}^{2k} t^{-1} 2^{-l-2j} \left( 2^{-l_0-j} + 2^{-2k-j} \right) \int e^{-i t \Phi\left( \frac{x}{t}, \frac{y}{t}, \overline{\xi} \right)} \Psi_{(l_0, l, l_0)}^b \psi_{j, k}^{\widehat{\mathcal{C}}}(\overline{\xi}) \widehat{h}_{\alpha}(t, \overline{\xi}) ~ d\overline{\xi} \label{estdispL4CresC-int-b-1} \\
&\quad + \sum_{l = l_1}^{2k} t^{-1} 2^{-k-l-2j} \int e^{-i t \Phi\left( \frac{x}{t}, \frac{y}{t}, \overline{\xi} \right)} \Psi_{(l_0, l, l_0)}^b \psi_{j, k}^{\widehat{\mathcal{C}}}(\overline{\xi}) m_b(\overline{\xi}) \widehat{X}_b(\overline{\xi}) \cdot \nabla_{\overline{\xi}} \widehat{h}_{\alpha}(t, \overline{\xi}) ~ d\overline{\xi} \label{estdispL4CresC-int-b-2} 
\end{align}
\end{subequations}
and then estimate: 
\begin{align*}
&\Vert \eqref{estdispL4CresC-int-b-1} \Vert_{L^4(A_{res}^{\mathcal{C}})} \lesssim \sum_{l = l_1}^{2k} \Vert \eqref{estdispL4CresC-int-b-1}_l \Vert_{L^2}^{\frac{1}{2}} \Vert \eqref{estdispL4CresC-int-b-1}_l \Vert_{L^{\infty}}^{\frac{1}{2}} \\
&\lesssim \sum_{l = l_1}^{2k} t^{-1} 2^{-l-3j-l_0} \Vert h_{\alpha}(t) \Vert_{L^2}^{\frac{1}{2}} \left( \Vert \Psi_{(l_0, l, l_0)}^b \psi_{j, k}^{\widehat{\mathcal{C}}} \Vert_{L^{\infty}(A_{res}^{\mathcal{C}}, L^{\frac{1}{1-\kappa}}_c L^2_b L^{\frac{1}{1-\kappa}}_a)} \Vert \psi_j \widehat{h}_{\alpha}(t) \Vert_{L^{\frac{1}{\kappa}}_c L^2_b L^{\frac{1}{\kappa}}_a} \right)^{\frac{1}{2}} \\
&\lesssim \sum_{l = l_1}^{2k} \mathfrak{t}^{-1} 2^{-l-\kappa l_0+\frac{l-k}{4}+\frac{3j}{4}} \Vert u \Vert_X \\
&\lesssim \mathfrak{t}^{-1+\frac{1}{16}+\frac{\delta}{32}+\frac{\kappa}{2}} 2^{\frac{k}{4}+\delta k+\frac{3j}{4}} \Vert u \Vert_X \\
&\lesssim t^{-\frac{3}{4}} 2^{-\frac{3j}{2}+\frac{k}{4}} \Vert u \Vert_X \\
&\Vert \eqref{estdispL4CresC-int-b-2} \Vert_{L^4(A_{res}^{\mathcal{C}})} \lesssim \sum_{l = l_1}^{2k} \Vert \eqref{estdispL4CresC-int-b-2}_l \Vert_{L^2}^{\frac{1}{2}} \Vert \eqref{estdispL4CresC-int-b-2}_l \Vert_{L^{\infty}(A_{res}^{\mathcal{C}})}^{\frac{1}{2}} \\
&\lesssim \sum_{l = l_1}^{2k} \mathfrak{t}^{-1} 2^{-k-l} \Vert \partial_x m_b(D) X_b h_{\alpha}(t) \Vert_{L^2}^{\frac{1}{2}} \left( \Vert \Psi_{(l_0, l, l_0)}^b \psi_{j, k}^{\widehat{\mathcal{C}}} \Vert_{L^{\infty}(A_{res}^{\mathcal{C}}, L^2)} \Vert \partial_x m_b(D) X_b h_{\alpha}(t) \Vert_{L^2} \right)^{\frac{1}{2}} \\
&\lesssim \sum_{l = l_1}^{2k} \mathfrak{t}^{-1} 2^{-k-l+\frac{l_0}{2}+\frac{l-k}{4}+\frac{3j}{4}} \Vert u \Vert_X \\
&\lesssim \mathfrak{t}^{-1+\frac{\delta}{32}} 2^{\frac{k}{4}+\delta k+\frac{3j}{4}} \Vert u \Vert_X \\
&\lesssim t^{-\frac{3}{4}} 2^{-\frac{3j}{2}+\frac{k}{4}} \Vert u \Vert_X \\
&\Vert \eqref{estdispL4CresC-int-int} \Vert_{L^4(A_{res}^{\mathcal{C}})} \lesssim \Vert \eqref{estdispL4CresC-int-int} \Vert_{L^2}^{\frac{1}{2}} \Vert \eqref{estdispL4CresC-int-int} \Vert_{L^{\infty}(A_{res}^{\mathcal{C}})}^{\frac{1}{2}} \\
&\lesssim \left( \langle 2^j \rangle^{-1} \Vert \langle \nabla \rangle h_{\alpha}(t) \Vert_{L^2} \right)^{\frac{1}{2}} \left( \Vert \Psi_{(l_0, l_1, l_0)}^{int} \psi_{j, k}^{\widehat{\mathcal{C}}} \Vert_{L^{\infty}(A_{res}^{\mathcal{C}}, L^{\frac{1}{1-\kappa}}_c L^2_b L^{\frac{1}{1-\kappa}}_a)} \Vert \psi_j \widehat{h}_{\alpha}(t) \Vert_{L^{\frac{1}{\kappa}}_c L^2_b L^{\frac{1}{\kappa}}_a} \right)^{\frac{1}{2}} \\
&\lesssim 2^{l_0-\kappa l_0+\frac{l_1-k}{4}+\frac{3j}{4}-\kappa j} \langle 2^j \rangle^{-\frac{1}{2}} \Vert u \Vert_X \\
&\lesssim \mathfrak{t}^{-\frac{1}{2}+\varrho+\frac{\kappa}{2}-\frac{1}{32}} 2^{\frac{k}{8}+\frac{3j}{4}-\kappa j} \langle 2^j \rangle^{-\frac{1}{2}} \Vert u \Vert_X \\
&\lesssim t^{-\frac{1}{2}-\frac{1}{64}} 2^{-\frac{7j}{8}+\frac{k}{8}} \langle 2^j \rangle^{-\frac{1}{4}} \Vert u \Vert_X 
\end{align*}
which is enough. 

\paragraph{Resonant frequencies, near the line} Let us consider the norm $L^4(A_{res}^{\mathcal{L}})$. 

As in the case $A_{res}^{\mathcal{R}}$, as soon as $(\widehat{X}_a, \widehat{X}_c) \cdot \nabla \Phi \ll 2^{2j}$, then $\widehat{X}_b \cdot \nabla \Phi \simeq 2^{2j}$. 

We therefore decompose as before using the localisation symbols $\Psi_{(l_a, l_b, l_a)}$, and we obtain the terms \eqref{estdispL4CresR-a}, \eqref{estdispL4CresR-c}, \eqref{estdispL4CresR-int}, choosing the same $l_0$. Again, we decompose \eqref{estdispL4CresR-a} into $\eqref{estdispL4CresR-a-1}+\eqref{estdispL4CresR-a-2}$, but this time the volume estimate is
\begin{align*}
\Vert \Psi_{(l, 0, l)}^{*} \psi_{j, k}^{\widehat{\mathcal{C}}} \Vert_{L^{\infty}(A_{res}^{\mathcal{L}}, L^{p_c}_c L^{p_b}_b L^{p_a}_a)} &\lesssim 2^{\frac{l}{p_a}+\frac{k}{p_b}} 2^{\frac{j}{p_a}+\frac{j}{p_b}+\frac{j}{p_c}}
\end{align*}
see \eqref{estimeevolumiquePsilacasCvoisL}. 

We may estimate \eqref{estdispL4CresR-a-1} in a very similar way as before (using $n \varrho \gg 1$), then 
\begin{align*}
&\Vert \eqref{estdispL4CresR-a-2} \Vert_{L^4(A_{res}^{\mathcal{L}})} \lesssim \sum_{i = 1}^n \sum_{l = l_0}^0 \Vert \eqref{estdispL4CresR-a-2}_l \Vert_{L^2}^{\frac{1}{2}} \Vert \eqref{estdispL4CresR-a-2}_l \Vert_{L^{\infty}(A_{res}^{\mathcal{L}})}^{\frac{1}{2}} \\
&\lesssim \sum_{i = 1}^n \sum_{l = l_0}^0 \mathfrak{t}^{-i} 2^{-2il+l} \Vert \partial_x m_b(D) X_b h_{\alpha}(t) \Vert_{L^2}^{\frac{1}{2}} \left( \Vert \Psi_{(l, 0, l)}^a \psi_{j, k}^{\widehat{\mathcal{C}}} \Vert_{L^{\infty}(A_{res}^{\mathcal{L}}, L^2_{c, b} L^{\frac{1}{1-\kappa}}_a)} \Vert \psi_j \widehat{h}_{\alpha}(t) \Vert_{L^2_{c, b} L^{\frac{1}{\kappa}}_a} \right)^{\frac{1}{2}} \\
&\lesssim \sum_{i = 1}^n \sum_{l = l_0}^0 \mathfrak{t}^{-i} 2^{-2il+\frac{3l}{2}-\frac{\kappa l}{2}+\frac{k}{4}+j-\frac{\kappa j}{2}} \Vert u \Vert_X^{\frac{1}{2}} \Vert \psi_j \widehat{h}_{\alpha}(t) \Vert_{L^2_{c, b} H^1_a}^{\frac{1}{4}} \Vert \widehat{h}_{\alpha}(t) \Vert_{L^2}^{\frac{1}{4}} \\
&\lesssim \mathfrak{t}^{-\frac{3}{4}+\frac{\kappa}{4}} 2^{\frac{k}{4}+\frac{3j}{4}-\frac{\kappa j}{2}} \Vert u \Vert_X \\
&\lesssim t^{-\frac{11}{16}} 2^{-\frac{5j}{16}-\frac{\kappa j}{2}+\frac{k}{4}} \Vert u \Vert_X
\end{align*}
which is enough. We estimate \eqref{estdispL4CresR-c} the same way. 

Then, for the internal term, we exhaust the cases depending on $\alpha$ and $k$. First, let us consider $\alpha \in \{ a, c \}$. 

If $2^k \lesssim \mathfrak{t}^{-\frac{7}{8}}$, we can estimate directly: 
\begin{align*}
\Vert \eqref{estdispL4CresR-int} \Vert_{L^4(A_{res}^{\mathcal{L}})} &\lesssim \Vert \eqref{estdispL4CresR-int} \Vert_{L^2}^{\frac{1}{2}} \Vert \eqref{estdispL4CresR-int} \Vert_{L^{\infty}(A_{res}^{\mathcal{L}})}^{\frac{1}{2}} \\
&\lesssim \left( \langle 2^j \rangle^{-1} \Vert \langle \nabla \rangle h_{\alpha}(t) \Vert_{L^2} \right)^{\frac{1}{2}} 
\left( \Vert \Psi_{(l_0, 0, l_0)}^{int} \psi_{j, k}^{\widehat{\mathcal{C}}} \Vert_{L^{\infty}(A_{res}^{\mathcal{C}}, L^2_{c, b} L^{\frac{1}{1-\kappa}}_a)} \Vert \psi_j \widehat{h}_{\alpha}(t) \Vert_{L^2_{c, b} L^{\frac{1}{\kappa}}_a} \right)^{\frac{1}{2}} \\
&\lesssim 2^{\frac{l_0}{2}-\frac{\kappa l_0}{2}+\frac{k}{4}+\frac{3j}{4}-\frac{\kappa j}{2}} \langle 2^j \rangle^{-\frac{1}{2}} \Vert u \Vert_X \\
&\lesssim \mathfrak{t}^{-\frac{15}{32}+\frac{\varrho}{2}+\frac{\kappa}{4}+\frac{7 \delta}{8}} 2^{\frac{3j}{4}-\frac{\kappa j}{2}+\delta k} \langle 2^j \rangle^{-\frac{1}{2}} \Vert u \Vert_X 
\end{align*}
If however $2^k \gtrsim \mathfrak{t}^{-\frac{7}{8}}$, we apply an integration by parts along $\widehat{X}_{b-corr}$, defined so to have 
\begin{align*}
\widehat{X}_{b-corr}(\overline{\xi}) \cdot \nabla_{\overline{\xi}} \left[ \Psi_{(l_0, 0, l_0)}^{int} \right] &\lesssim 2^{-j} 
\end{align*}
as in the proof of the $L^{\infty}$ dispersive estimate. Then, 
\begin{align*}
\eqref{estdispL4CresR-int} &= t^{-1} 2^{-k-3j} \int e^{i t \Phi} \Psi_{(l_0, 0, l_0)}^{int} \psi_{j, k}^{\widehat{\mathcal{C}}} \widehat{h}_{\alpha}(t) ~ d\overline{\xi} \\
&\quad + t^{-1} 2^{-k-3j} \int e^{i t \Phi} \Psi_{(l_0, 0, l_0)}^{int} \psi_{j, k}^{\widehat{\mathcal{C}}} \overline{\xi} m_b(\overline{\xi}) \widehat{X}_{b-corr}(\overline{\xi}) \cdot \nabla \widehat{h}_{\alpha}(t) ~ d\overline{\xi} 
\end{align*}
We can estimate: 
\begin{align*}
\Vert \eqref{estdispL4CresR-int} \Vert_{L^4(A_{res}^{\mathcal{L}})} &\lesssim \Vert \eqref{estdispL4CresR-int} \Vert_{L^2}^{\frac{1}{2}} \Vert \eqref{estdispL4CresR-int} \Vert_{L^{\infty}(A_{res}^{\mathcal{L}})}^{\frac{1}{2}} \\
&\lesssim \mathfrak{t}^{-1} 2^{-k} \left( \Vert h_{\alpha}(t) \Vert_{L^2} + \Vert \partial_x m_b(D) X_b h_{\alpha}(t) \Vert_{L^2} \right) \Vert \Psi_{(l_0, 0, l_0)}^{int} \psi_{j, k}^{\widehat{\mathcal{C}}} \Vert_{L^{\infty}(A_{res}^{\mathcal{L}}, L^2)}^{\frac{1}{2}} \\
&\lesssim \mathfrak{t}^{-1} 2^{\frac{l_0}{4}-\frac{3k}{4}+\frac{3j}{4}} \Vert u \Vert_X \\
&\lesssim \mathfrak{t}^{-\frac{15}{32}+\frac{\varrho}{4}+\frac{7\delta}{8}} 2^{\frac{3j}{4}+\delta k} \Vert u \Vert_X 
\end{align*}
as wanted. 

For $\alpha = b$, we decompose differently. First, if $2^k \gtrsim \langle t \rangle^{-\frac{1}{12}} \mathfrak{t}^{-\frac{9}{16}}$, we reuse the integration by parts here above and get 
\begin{align*}
\Vert \eqref{estdispL4CresR-int} \Vert_{L^4(A_{res}^{\mathcal{L}})} &\lesssim \mathfrak{t}^{-1} 2^{\frac{l_0}{4}-\frac{3k}{4}+\frac{3j}{4}} \Vert u \Vert_X \\
&\lesssim \langle t \rangle^{\frac{1}{16}+\frac{\delta}{12}} \mathfrak{t}^{-1-\frac{1}{8}+\frac{\varrho}{4}+\frac{27}{64}+\frac{9\delta}{16}} 2^{\frac{3j}{4}+\delta k} \Vert u \Vert_X \\
&\lesssim \langle t \rangle^{\frac{1}{16}+\frac{\delta}{12}} t^{-\frac{45}{64}+\frac{\varrho}{4}+\frac{9\delta}{16}} 2^{-j-\frac{23j}{64}+\frac{3\varrho j}{4}+\frac{27\delta j}{16}+\delta k} \Vert u \Vert_X
\end{align*}
However, if $2^k \lesssim \langle t \rangle^{-\frac{1}{12}} \mathfrak{t}^{-\frac{9}{16}}$, we apply no integration by parts and estimate directly as: 
\begin{align*}
&\Vert e^{-i t \omega(D)} \psi_{j, k}^{\widehat{\mathcal{C}}} h_b(t) \Vert_{L^4(A_{res}^{\mathcal{C}})} \lesssim \Vert e^{-i t \omega(D)} \psi_{j, k}^{\widehat{\mathcal{C}}} m_b(D) X_b f(t) \Vert_{L^4(A_{res}^{\mathcal{C}})} + \Vert e^{-i t \omega(D)} \psi_{j, k}^{\widehat{\mathcal{C}}} g_b(t) \Vert_{L^4(A_{res}^{\mathcal{C}})} \\
&\lesssim 2^{\frac{3j}{4}+\frac{5k}{4}} \Vert \psi_j(D) X_b f(t) \Vert_{L^2} + \left( 2^{\frac{k}{4}} \Vert |\nabla|^{\frac{3}{4}} g_b(t) \Vert_{L^2} \right)^{4\delta} \left( \langle 2^j \rangle^{-\frac{3}{2}} \Vert e^{-i t \omega(D)} \langle \nabla \rangle^{\frac{3}{2}} g_b(t) \Vert_{L^4} \right)^{1-4\delta} \\
&\lesssim 2^{\frac{3j}{4}+\frac{5k}{4}} \langle t \rangle^{\frac{1}{6}+101\delta} \Vert u \Vert_X + 2^{\delta k} \langle 2^j \rangle^{-\frac{5}{4}} t^{-\frac{5}{12}} \langle t \rangle^{-\frac{1}{12}-\frac{1}{128}} \Vert u \Vert_X \\
&\lesssim \langle t \rangle^{\frac{1}{16}+102\delta} \mathfrak{t}^{-\frac{45}{64}+\frac{9\delta}{16}} 2^{\frac{3j}{4}+\delta k} \Vert u \Vert_X + 2^{\delta k} \langle 2^j \rangle^{-\frac{5}{4}} t^{-\frac{5}{12}} \langle t \rangle^{-\frac{1}{12}-\frac{1}{128}} \Vert u \Vert_X \\
&\lesssim \langle t \rangle^{\frac{1}{16}+102\delta} t^{-\frac{45}{64}+\frac{9\delta}{16}} 2^{-j-\frac{23j}{64}+\frac{27\delta j}{16}+\delta k} \Vert u \Vert_X + 2^{\delta k} \langle 2^j \rangle^{-\frac{5}{4}} t^{-\frac{5}{12}} \langle t \rangle^{-\frac{1}{12}-\frac{1}{128}} \Vert u \Vert_X
\end{align*}
which is enough. 

\paragraph{Low frequencies, away from the cone} Let $n \geq j+10$ and consider the norm $L^4(A_n^{\mathcal{R}})$. 

Note that, for fixed $(x, y)$, we have
\begin{align*}
\begin{pmatrix} \widehat{X}_a(\overline{\xi}) & \widehat{X}_c(\overline{\xi}) \end{pmatrix}^T \begin{pmatrix} x \\ y \end{pmatrix} &= \begin{pmatrix} \frac{x \xi_0 + y \cdot \xi}{|\overline{\xi}|} \\ \frac{y \cdot J \xi}{|\xi|} \end{pmatrix} 
\end{align*}
We aim to show that at least one of them needs to be of order $2^{2n} t$. Then, since we are localised on low frequencies $n \geq j+10$, this will imply that 
\begin{align*}
\left| \left( \widehat{X}_a(\overline{\xi}) \cdot \nabla_{\overline{\xi}} \Phi, \widehat{X}_c(\overline{\xi}) \cdot \nabla_{\overline{\xi}} \Phi \right) \right| \gtrsim 2^{2n} 
\end{align*}

Assume by contradiction that both of these quantities are small with respect to $2^{2n} t$. Then $y$ is close to alignment with $\frac{\xi}{|\xi|}$, and since $\sqrt{3} |\xi_0| - |\xi| \ll |\overline{\xi}|$, 
\begin{align*}
\left| \frac{x \xi_0 + y \cdot \xi}{|\overline{\xi}|} \right| = \frac{1}{2} \left| |x| \pm \sqrt{3} |y| \right| + o\left( 2^{2n} t \right) 
\end{align*}
But $\left| |x| \pm \sqrt{3} |y| \right| \gtrsim |x|+|y|$ whatever the choice of the sign, since $(x, y) \in A_n^{\widehat{\mathcal{R}}}$. This is a contradiction, and therefore for any $(x, y) \in A_n^{\mathcal{R}}$ we can decompose
\begin{align*}
1 &= \widetilde{\psi}\left( 2^{-2n} \widehat{X}_a(\overline{\xi}) \cdot \nabla_{\overline{\xi}} \Phi \right) \widetilde{\chi}\left( 2^{-2n} \widehat{X}_b(\overline{\xi}) \cdot \nabla_{\overline{\xi}} \Phi \right) \widetilde{\chi}\left( 2^{-2n} \widehat{X}_c(\overline{\xi}) \cdot \nabla_{\overline{\xi}} \Phi \right) \\
&\quad + \widetilde{\chi}\left( 2^{-2n} \widehat{X}_a(\overline{\xi}) \cdot \nabla_{\overline{\xi}} \Phi \right) \widetilde{\chi}\left( 2^{-2n} \widehat{X}_b(\overline{\xi}) \cdot \nabla_{\overline{\xi}} \Phi \right) \widetilde{\psi}\left( 2^{-2n} \widehat{X}_c(\overline{\xi}) \cdot \nabla_{\overline{\xi}} \Phi \right) 
\end{align*}
for some choice of $\widetilde{\psi}, \widetilde{\chi}$. In particular, we can apply integrations by parts as in the high frequencies case, replacing $j$ by $n$ (and possibly $X_a$ by $X_c$). We then get: 
\begin{align*}
\Vert \eqref{preuveestimeeL4CformuleFourier} \Vert_{L^4(A_n^{\mathcal{R}})} &\lesssim t^{-1} 2^{-\frac{9n}{4}+\frac{k}{4}} \Vert u \Vert_X
\end{align*}
Summing over $n \geq j+10$, we recover: 
\begin{align*}
\sum_{n \geq j+10} \Vert \eqref{preuveestimeeL4CformuleFourier} \Vert_{L^4(A_n^{\mathcal{R}})} &\lesssim t^{-1} 2^{-\frac{9j}{4}+\frac{k}{4}} \Vert u \Vert_X
\end{align*}
which is enough. 

\paragraph{Low frequencies, near the cone} Let us finally consider the norm $L^4(A_n^{\mathcal{C}})$. 

This time, reusing the computation above, we can have 
\begin{align*}
\widehat{X}_a(\overline{\xi}) \cdot \nabla_{\overline{\xi}} \Phi, ~ \widehat{X}_c(\overline{\xi}) \cdot \nabla_{\overline{\xi}} \Phi ~ \ll 2^{2n} 
\end{align*}
but this forces $\widehat{X}_b(\overline{\xi}) \cdot \nabla_{\overline{\xi}} \Phi \simeq 2^{2n}$. 

In particular, we can again decompose using the symbols $\Psi_{(l_a, l_b, l_c)}$ and the estimates are simpler than in the resonant frequency case. Just as in the previous low frequency case, we automatically gain a geometric factor in $n$, so that summability in $n$ is not an issue. We skip the details. 
\end{Dem}

\subsection{Integrated \texorpdfstring{$L^4$}{L4} estimate} 

In previous sections, we showed: 

\begin{Prop} Let $t > 0$. Then
\begin{align*}
\Vert e^{-i t \omega(D)} \partial_x h_b(t) \Vert_{L^4} &\lesssim t^{-\frac{511}{1024}} \langle t \rangle^{-\frac{1}{512}} \Vert u \Vert_X \\
\Vert e^{-i t \omega(D)} m_{\widehat{\mathcal{L}}}(D) \nabla (h_a(t), h_c(t)) \Vert_{L^4} &\lesssim t^{-\frac{63}{128}} \Vert u \Vert_X \\
\Vert e^{-i t \omega(D)} \left( 1 - m_{\widehat{\mathcal{C}}}(D) - m_{\widehat{\mathcal{L}}}(D) \right) \nabla (h_a(t), h_c(t)) \Vert_{L^4} &\lesssim t^{-\frac{63}{128}} \langle t \rangle^{-\frac{1}{64}} \Vert u \Vert_X 
\end{align*} \label{propositionestimeedispersiveL4totaleponctuelle} 
\end{Prop}

We will now improve the $L^4$ estimates on $h_a, h_c$ by integrating them in time and using the fact that $X_a, X_c$ correspond to invariances of the equation. 

\begin{Prop} Let $t > 0$ and $\alpha \in \{ a, c \}$. Then
\begin{align*}
\int_0^t s^{\frac{63}{64}} \langle s \rangle^{\frac{1}{50}} \Vert e^{-i s \omega(D)} \nabla h_{\alpha}(s) \Vert_{L^4}^4 ~ ds ~ \lesssim \Vert u \Vert_X^4 
\end{align*} \label{propositionestimeedispersiveL4totaleintegrale} 
\end{Prop}

\begin{Dem}
It is enough to show the property replacing $h_{\alpha}$ by $m_A(D) h_{\alpha}$, for every $A \in \{ \widehat{\mathcal{R}}, \widehat{\mathcal{C}}, \widehat{\mathcal{L}}, \widehat{\mathcal{P}} \}$. If $A = \widehat{\mathcal{R}}$ or $A = \widehat{\mathcal{P}}$, the result is a consequence of Proposition \ref{propositionestimeedispersiveL4totaleponctuelle}. 

Then, for $A = \widehat{\mathcal{L}}$, we develop: 
\begin{align*}
&\int_0^t s^{\frac{63}{64}} \langle s \rangle^{\frac{1}{50}} \Vert e^{-i s \omega(D)} \nabla m_{\widehat{\mathcal{L}}}(D) h_{\alpha}(s) \Vert_{L^4}^4 ~ ds \\
&\quad = \int_0^t \int \int \int e^{i s \varphi_3(\overline{\xi}, \overline{\eta}, \overline{\sigma})} s^{\frac{63}{64}} \langle s \rangle^{\frac{1}{50}} |\overline{\xi}| |\overline{\eta}| |\overline{\sigma}| |\overline{\rho}| m_{\widehat{\mathcal{L}}}^{\otimes 4}(\overline{\xi}, \overline{\eta}, \overline{\sigma}, \overline{\rho}) \\
&\pushright{\widehat{h}_{\alpha}(s, \overline{\xi}) \widehat{h}_{\alpha}(s, \overline{\eta}) \widehat{h}_{\alpha}(s, \overline{\sigma}) \widehat{h}_{\alpha}(s, \overline{\rho}) ~ d\overline{\xi} d\overline{\eta} d\overline{\sigma} ds}
\end{align*}
where $\varphi_3$ is the cubic interaction phase and where we recall the convention $\overline{\rho} = \overline{\xi}-\overline{\eta}-\overline{\sigma}$. Let $\mu_{\leq}$ be a Coifman-Meyer symbol localising on $\{ |\overline{\rho}| \gtrsim |\overline{\xi}| + |\overline{\eta}|+|\overline{\sigma}| \}$: up to interchanging the variables, it is enough to estimate only 
\begin{align*}
\int_0^t \int \int \int e^{i s \varphi_3(\overline{\xi}, \overline{\eta}, \overline{\sigma})} s^{\frac{63}{64}} \langle s \rangle^{\frac{1}{50}} \mu_{\leq}(\overline{\xi}, \overline{\eta}, \overline{\sigma}) |\overline{\xi}| |\overline{\eta}| |\overline{\sigma}| |\overline{\rho}| m_{\widehat{\mathcal{L}}}^{\otimes 4}(\overline{\xi}, \overline{\eta}, \overline{\sigma}, \overline{\rho}) \\
\widehat{h}_{\alpha}(s, \overline{\xi}) \widehat{h}_{\alpha}(s, \overline{\eta}) \widehat{h}_{\alpha}(s, \overline{\sigma}) \widehat{h}_{\alpha}(s, \overline{\rho}) ~ d\overline{\xi} d\overline{\eta} d\overline{\sigma} ds
\end{align*}

We now decompose
\begin{align*}
h_{\alpha}(t) &= m_{\alpha}(D) X_{\alpha} f(t) - g_{\alpha}(t) 
\end{align*}
for the term in $\overline{\rho}$. On the one hand, we estimate: 
\begin{align*}
&\int_0^t \int \int \int e^{i s \varphi_3(\overline{\xi}, \overline{\eta}, \overline{\sigma})} s^{\frac{63}{64}} \langle s \rangle^{\frac{1}{50}} \mu_{\leq}(\overline{\xi}, \overline{\eta}, \overline{\sigma}) |\overline{\xi}| |\overline{\eta}| |\overline{\sigma}| |\overline{\rho}| m_{\widehat{\mathcal{L}}}^{\otimes 4}(\overline{\xi}, \overline{\eta}, \overline{\sigma}, \overline{\rho}) \\
&\pushright{\widehat{h}_{\alpha}(s, \overline{\xi}) \widehat{h}_{\alpha}(s, \overline{\eta}) \widehat{h}_{\alpha}(s, \overline{\sigma}) \widehat{g}_{\alpha}(s, \overline{\rho}) ~ d\overline{\xi} d\overline{\eta} d\overline{\sigma} ds} \\
&\quad \lesssim \int_0^t s^{\frac{63}{64}} \langle s \rangle^{\frac{1}{50}} \Vert e^{-i s \omega(D)} \nabla g_{\alpha}(s) \Vert_{L^4} \Vert e^{-i s \omega(D)} \nabla m_{\widehat{\mathcal{L}}}(D) h_{\alpha}(s) \Vert_{L^4}^3 ~ ds \\
&\quad \lesssim \left( \int_0^t s^{\frac{63}{64}} \langle s \rangle^{\frac{1}{50}} \Vert e^{-i s \omega(D)} \nabla m_{\widehat{\mathcal{L}}}(D) h_{\alpha}(s) \Vert_{L^4}^4 ~ ds \right)^{\frac{3}{4}} \left( \int_0^t s^{\frac{63}{64}} \langle s \rangle^{\frac{1}{50}} \Vert e^{-i s \omega(D)} \nabla g_{\alpha}(s) \Vert_{L^4}^4 ~ ds \right)^{\frac{1}{4}} \\
&\quad \lesssim \left( \int_0^t s^{\frac{63}{64}} \langle s \rangle^{\frac{1}{50}} \Vert e^{-i s \omega(D)} \nabla m_{\widehat{\mathcal{L}}}(D) h_{\alpha}(s) \Vert_{L^4}^4 ~ ds \right)^{\frac{3}{4}} \left( \int_0^t s^{-\frac{2}{3}-\frac{1}{64}} \langle s \rangle^{-\frac{1}{2}+\frac{1}{50}} \Vert u \Vert_X^4 ~ ds \right)^{\frac{1}{4}} \\
&\quad \lesssim \Vert u \Vert_X \left( \int_0^t s^{\frac{63}{64}} \langle s \rangle^{\frac{1}{50}} \Vert e^{-i s \omega(D)} \nabla m_{\widehat{\mathcal{L}}}(D) h_{\alpha}(s) \Vert_{L^4}^4 ~ ds \right)^{\frac{3}{4}}
\end{align*}

Then, where we created the $m_{\alpha}(D) X_{\alpha}$ operator, we can decompose
\begin{align*}
|\overline{\rho}| \widehat{X}_a(\overline{\rho}) \cdot \nabla_{\overline{\xi}} \widehat{f}(s, \overline{\rho}) 
&= \overline{\rho} \cdot \nabla_{\overline{\xi}} \widehat{f}(s, \overline{\rho}) \\
&= \left( \overline{\xi} \cdot \nabla_{\overline{\xi}} + \overline{\eta} \cdot \nabla_{\overline{\eta}} + \overline{\sigma} \cdot \nabla_{\overline{\sigma}} \right) \widehat{f}(s, \overline{\rho}) \\
&= \left( |\overline{\xi}| \widehat{X}_a(\overline{\xi}) \cdot \nabla_{\overline{\xi}} + |\overline{\eta}| \widehat{X}_a(\overline{\eta}) \cdot \nabla_{\overline{\eta}} + |\overline{\sigma}| \widehat{X}_a(\overline{\sigma}) \cdot \nabla_{\overline{\sigma}} \right) \widehat{f}(s, \overline{\rho}) 
\end{align*}
and similarly: 
\begin{align*}
&|\overline{\rho}| m_c(\overline{\rho}) \widehat{X}_c(\overline{\rho}) \cdot \nabla_{\overline{\xi}} \widehat{f}(s, \overline{\rho}) \\
&\quad = \left( |\overline{\xi}| m_c(\overline{\xi}) \widehat{X}_c(\overline{\xi}) \cdot \nabla_{\overline{\xi}} + |\overline{\eta}| m_c(\overline{\eta}) \widehat{X}_c(\overline{\eta}) \cdot \nabla_{\overline{\eta}} + |\overline{\sigma}| m_c(\overline{\sigma}) \widehat{X}_c(\overline{\sigma}) \cdot \nabla_{\overline{\sigma}} \right) \widehat{f}(s, \overline{\rho})
\end{align*}
Moreover, 
\begin{align*}
&\left( |\overline{\xi}| \widehat{X}_a(\overline{\xi}) \cdot \nabla_{\overline{\xi}} + |\overline{\eta}| \widehat{X}_a(\overline{\eta}) \cdot \nabla_{\overline{\eta}} + |\overline{\sigma}| \widehat{X}_a(\overline{\sigma}) \cdot \nabla_{\overline{\sigma}} \right) \varphi_3 = 3 \varphi_3 \\
&\left( |\overline{\xi}| m_c(\overline{\xi}) \widehat{X}_c(\overline{\xi}) \cdot \nabla_{\overline{\xi}} + |\overline{\eta}| m_c(\overline{\eta}) \widehat{X}_c(\overline{\eta}) \cdot \nabla_{\overline{\eta}} + |\overline{\sigma}| m_c(\overline{\sigma}) \widehat{X}_c(\overline{\sigma}) \cdot \nabla_{\overline{\sigma}} \right) \varphi_3 = 0 
\end{align*}
by direct computation. Let us denote by $c_a = 3$, $c_c = 0$. 

We now apply an integration by parts in Fourier: 
\begin{subequations}
\begin{align}
&\begin{aligned}
\int_0^t \int \int \int e^{i s \varphi_3(\overline{\xi}, \overline{\eta}, \overline{\sigma})} s^{\frac{63}{64}} \langle s \rangle^{\frac{1}{50}} \mu_{\leq}(\overline{\xi}, \overline{\eta}, \overline{\sigma}) |\overline{\xi}| |\overline{\eta}| |\overline{\sigma}| m_{\widehat{\mathcal{L}}}^{\otimes 4}(\overline{\xi}, \overline{\eta}, \overline{\sigma}, \overline{\rho}) \\
\widehat{h}_{\alpha}(s, \overline{\xi}) \widehat{h}_{\alpha}(s, \overline{\eta}) \widehat{h}_{\alpha}(s, \overline{\sigma}) |\overline{\rho}| m_{\alpha}(\overline{\rho}) \widehat{X}_{\alpha}(\overline{\rho}) \cdot \nabla_{\overline{\xi}} \widehat{f}(s, \overline{\rho}) ~ d\overline{\xi} d\overline{\eta} d\overline{\sigma} ds
\end{aligned} \notag \\
&\quad \begin{aligned}
= - \int_0^t \int \int \int e^{i s \varphi_3(\overline{\xi}, \overline{\eta}, \overline{\sigma})} i s c_{\alpha} \varphi_3 ~ s^{\frac{63}{64}} \langle s \rangle^{\frac{1}{50}} \mu_{\leq}(\overline{\xi}, \overline{\eta}, \overline{\sigma}) |\overline{\xi}| |\overline{\eta}| |\overline{\sigma}| m_{\widehat{\mathcal{L}}}^{\otimes 4}(\overline{\xi}, \overline{\eta}, \overline{\sigma}, \overline{\rho}) \\
\widehat{h}_{\alpha}(s, \overline{\xi}) \widehat{h}_{\alpha}(s, \overline{\eta}) \widehat{h}_{\alpha}(s, \overline{\sigma}) \widehat{f}(s, \overline{\rho}) ~ d\overline{\xi} d\overline{\eta} d\overline{\sigma} ds
\end{aligned} \label{estdispL4integralepreuve-1} \\
&\quad \quad \begin{aligned}
&- 3 \int_0^t \int \int \int e^{i s \varphi_3(\overline{\xi}, \overline{\eta}, \overline{\sigma})} s^{\frac{63}{64}} \langle s \rangle^{\frac{1}{50}} \mu_{\leq}(\overline{\xi}, \overline{\eta}, \overline{\sigma}) |\overline{\xi}| |\overline{\eta}| |\overline{\sigma}| m_{\widehat{\mathcal{L}}}^{\otimes 4}(\overline{\xi}, \overline{\eta}, \overline{\sigma}, \overline{\rho}) \\
&\quad \quad \quad \quad |\overline{\xi}| m_{\alpha}(\overline{\xi}) \widehat{X}_{\alpha}(\overline{\xi}) \cdot \nabla_{\overline{\xi}} \widehat{h}_{\alpha}(s, \overline{\xi}) \widehat{h}_{\alpha}(s, \overline{\eta}) \widehat{h}_{\alpha}(s, \overline{\sigma}) \widehat{f}(s, \overline{\rho}) ~ d\overline{\xi} d\overline{\eta} d\overline{\sigma} ds 
\end{aligned} \label{estdispL4integralepreuve-2} \\
&\quad \quad \begin{aligned}
&- 3 \int_0^t \int \int \int e^{i s \varphi_3(\overline{\xi}, \overline{\eta}, \overline{\sigma})} s^{\frac{63}{64}} \langle s \rangle^{\frac{1}{50}} \nabla_{\overline{\xi}} \cdot \left( |\overline{\xi}| m_{\alpha}(\overline{\xi}) \widehat{X}_{\alpha}(\overline{\xi}) \mu_{\leq}(\overline{\xi}, \overline{\eta}, \overline{\sigma}) |\overline{\xi}| |\overline{\eta}| |\overline{\sigma}| m_{\widehat{\mathcal{L}}}^{\otimes 4}(\overline{\xi}, \overline{\eta}, \overline{\sigma}, \overline{\rho}) \right) \\
&\pushright{ \widehat{h}_{\alpha}(s, \overline{\xi}) \widehat{h}_{\alpha}(s, \overline{\eta}) \widehat{h}_{\alpha}(s, \overline{\sigma}) |\overline{\rho}| m_{\alpha}(\overline{\rho}) \widehat{X}_{\alpha}(\overline{\rho}) \cdot \nabla_{\overline{\xi}} \widehat{f}(s, \overline{\rho}) ~ d\overline{\xi} d\overline{\eta} d\overline{\sigma} ds }
\end{aligned} \label{estdispL4integralepreuve-3} 
\end{align}
\end{subequations} 
where we symmetrised the Fourier variables. 

Using $\mu_{\leq}$, we can estimate
\begin{align*}
\eqref{estdispL4integralepreuve-2} &\lesssim \int_0^t s^{\frac{63}{64}} \langle s \rangle^{\frac{1}{50}} \Vert \nabla m_{\alpha}(D) X_{\alpha} h_{\alpha}(s) \Vert_{L^2} \Vert e^{-i s \omega(D)} m_{\widehat{\mathcal{L}}}(D) \nabla h_{\alpha}(s) \Vert_{L^4}^2 \Vert \nabla m_{\widehat{\mathcal{L}}}(D) u(s) \Vert_{L^{\infty}} ~ ds \\
&\lesssim \left( \int_0^t s^{\frac{63}{64}} \langle s \rangle^{\frac{1}{50}} \Vert e^{-i s \omega(D)} m_{\widehat{\mathcal{L}}}(D) \nabla h_{\alpha}(s) \Vert_{L^4}^4 ~ ds \right)^{\frac{1}{2}} \\
&\quad \quad \left( \int_0^t s^{\frac{63}{64}} \langle s \rangle^{\frac{1}{50}} \Vert \nabla m_{\alpha}(D) X_{\alpha} h_{\alpha}(s) \Vert_{L^2}^2 \Vert \nabla m_{\widehat{\mathcal{L}}}(D) u(s) \Vert_{L^{\infty}}^2 ~ ds \right)^{\frac{1}{2}} \\
&\lesssim \left( \int_0^t s^{\frac{63}{64}} \langle s \rangle^{\frac{1}{50}} \Vert e^{-i s \omega(D)} m_{\widehat{\mathcal{L}}}(D) \nabla h_{\alpha}(s) \Vert_{L^4}^4 ~ ds \right)^{\frac{1}{2}} \left( \int_0^t s^{-\frac{2}{3}-\frac{1}{64}} \langle s \rangle^{-\frac{1}{2}+\frac{1}{50}} \Vert u \Vert_X^4 ~ ds \right)^{\frac{1}{2}} \\
&\lesssim \Vert u \Vert_X^2 \left( \int_0^t s^{\frac{63}{64}} \langle s \rangle^{\frac{1}{50}} \Vert e^{-i s \omega(D)} m_{\widehat{\mathcal{L}}}(D) \nabla h_{\alpha}(s) \Vert_{L^4}^4 ~ ds \right)^{\frac{1}{2}} \\
\eqref{estdispL4integralepreuve-3} &\lesssim \int_0^t s^{\frac{63}{64}} \langle s \rangle^{\frac{1}{50}} \Vert h_{\alpha}(s) \Vert_{L^2} \Vert e^{-i s \omega(D)} m_{\widehat{\mathcal{L}}}(D) \nabla h_{\alpha}(s) \Vert_{L^4}^2 \Vert \nabla m_{\widehat{\mathcal{L}}}(D) u(s) \Vert_{L^{\infty}} ~ ds \\
&\lesssim \Vert u \Vert_X^2 \left( \int_0^t s^{\frac{63}{64}} \langle s \rangle^{\frac{1}{50}} \Vert e^{-i s \omega(D)} m_{\widehat{\mathcal{L}}}(D) \nabla h_{\alpha}(s) \Vert_{L^4}^4 ~ ds \right)^{\frac{1}{2}}
\end{align*}

Then, on \eqref{estdispL4integralepreuve-1}, we apply an integration by parts in time: 
\begin{subequations}
\begin{align}
&\eqref{estdispL4integralepreuve-1} \notag \\
&\begin{aligned}
= \int_0^t \int \int \int e^{i s \varphi_3(\overline{\xi}, \overline{\eta}, \overline{\sigma})} \partial_s \left( s c_{\alpha} ~ s^{\frac{63}{64}} \langle s \rangle^{\frac{1}{50}} \right) \mu_{\leq}(\overline{\xi}, \overline{\eta}, \overline{\sigma}) |\overline{\xi}| |\overline{\eta}| |\overline{\sigma}| m_{\widehat{\mathcal{L}}}^{\otimes 4}(\overline{\xi}, \overline{\eta}, \overline{\sigma}, \overline{\rho}) \\
\widehat{h}_{\alpha}(s, \overline{\xi}) \widehat{h}_{\alpha}(s, \overline{\eta}) \widehat{h}_{\alpha}(s, \overline{\sigma}) \widehat{f}(s, \overline{\rho}) ~ d\overline{\xi} d\overline{\eta} d\overline{\sigma} ds 
\end{aligned} \label{estdispL4integralepreuve-1-1} \\
&\quad \begin{aligned}
+ 3 \int_0^t \int \int \int e^{i s \varphi_3(\overline{\xi}, \overline{\eta}, \overline{\sigma})} s c_{\alpha} ~ s^{\frac{63}{64}} \langle s \rangle^{\frac{1}{50}} \mu_{\leq}(\overline{\xi}, \overline{\eta}, \overline{\sigma}) |\overline{\xi}| |\overline{\eta}| |\overline{\sigma}| m_{\widehat{\mathcal{L}}}^{\otimes 4}(\overline{\xi}, \overline{\eta}, \overline{\sigma}, \overline{\rho}) \\
\partial_s \widehat{h}_{\alpha}(s, \overline{\xi}) \widehat{h}_{\alpha}(s, \overline{\eta}) \widehat{h}_{\alpha}(s, \overline{\sigma}) \widehat{f}(s, \overline{\rho}) ~ d\overline{\xi} d\overline{\eta} d\overline{\sigma} ds 
\end{aligned} \label{estdispL4integralepreuve-1-2} \\
&\quad \begin{aligned}
+ \int_0^t \int \int \int e^{i s \varphi_3(\overline{\xi}, \overline{\eta}, \overline{\sigma})} s c_{\alpha} ~ s^{\frac{63}{64}} \langle s \rangle^{\frac{1}{50}} \mu_{\leq}(\overline{\xi}, \overline{\eta}, \overline{\sigma}) |\overline{\xi}| |\overline{\eta}| |\overline{\sigma}| m_{\widehat{\mathcal{L}}}^{\otimes 4}(\overline{\xi}, \overline{\eta}, \overline{\sigma}, \overline{\rho}) \\
\widehat{h}_{\alpha}(s, \overline{\xi}) \widehat{h}_{\alpha}(s, \overline{\eta}) \widehat{h}_{\alpha}(s, \overline{\sigma}) \partial_s \widehat{f}(s, \overline{\rho}) ~ d\overline{\xi} d\overline{\eta} d\overline{\sigma} ds 
\end{aligned} \label{estdispL4integralepreuve-1-3} \\
&\quad \begin{aligned}
- \int \int \int e^{i t \varphi_3(\overline{\xi}, \overline{\eta}, \overline{\sigma})} t c_{\alpha} ~ t^{\frac{63}{64}} \langle t \rangle^{\frac{1}{50}} \mu_{\leq}(\overline{\xi}, \overline{\eta}, \overline{\sigma}) |\overline{\xi}| |\overline{\eta}| |\overline{\sigma}| m_{\widehat{\mathcal{L}}}^{\otimes 4}(\overline{\xi}, \overline{\eta}, \overline{\sigma}, \overline{\rho}) \\
\widehat{h}_{\alpha}(t, \overline{\xi}) \widehat{h}_{\alpha}(t, \overline{\eta}) \widehat{h}_{\alpha}(t, \overline{\sigma}) \widehat{f}(t, \overline{\rho}) ~ d\overline{\xi} d\overline{\eta} d\overline{\sigma} 
\end{aligned} \label{estdispL4integralepreuve-1-4} 
\end{align}
\end{subequations}
We can now estimate using Proposition \ref{propositionestimeedispersiveL4totaleponctuelle} and Propositions \ref{prop-estimeeaprioriH1-h} and \ref{propositionestimeesstandardscattering}: 
\begin{align*}
\eqref{estdispL4integralepreuve-1-1} &\lesssim \int_0^t s^{\frac{63}{64}} \langle s \rangle^{\frac{1}{50}} \Vert h_{\alpha}(s) \Vert_{L^2} \Vert e^{-i s \omega(D)} \nabla m_{\widehat{\mathcal{L}}}(D) \nabla h_{\alpha}(s) \Vert_{L^4}^2 \Vert \nabla m_{\widehat{\mathcal{L}}}(D) u(s) \Vert_{L^{\infty}} ~ ds \\
&\lesssim \Vert u \Vert_X^2 \left( \int_0^t s^{\frac{63}{64}} \langle s \rangle^{\frac{1}{50}} \Vert e^{-i s \omega(D)} m_{\widehat{\mathcal{L}}}(D) \nabla h_{\alpha}(s) \Vert_{L^4}^4 ~ ds \right)^{\frac{1}{2}} \\
\eqref{estdispL4integralepreuve-1-2} &\lesssim \int_0^t s^{1+\frac{63}{64}} \langle s \rangle^{\frac{1}{50}} \Vert \partial_s h_{\alpha}(s) \Vert_{L^2} \Vert e^{-i s \omega(D)} \nabla m_{\widehat{\mathcal{L}}}(D) \nabla h_{\alpha}(s) \Vert_{L^4}^2 \Vert \nabla m_{\widehat{\mathcal{L}}}(D) u(s) \Vert_{L^{\infty}} ~ ds \\ 
&\lesssim \left( \int_0^t s^{1-\delta} \langle s \rangle^{2\delta} \Vert \partial_s h_{\alpha}(s) \Vert_{L^2}^2 ~ ds \right)^{\frac{1}{2}} \left( \int_0^t s^{\frac{63}{64}} \langle s \rangle^{\frac{1}{50}} \Vert e^{-i s \omega(D)} \nabla m_{\widehat{\mathcal{L}}}(D) \nabla h_{\alpha}(s) \Vert_{L^4}^4 ~ ds \right)^{\frac{1}{2}} \\
&\quad \quad \quad \left( \sup_{s \in (0, t]} s^{1+\frac{\delta}{2}-\frac{1}{128}} \langle s \rangle^{\frac{1}{100}-\delta} \Vert \nabla m_{\widehat{\mathcal{L}}}(D) u(s) \Vert_{L^{\infty}} \right) \\
&\lesssim \Vert u \Vert_X^2 \left( \Vert u \Vert_X + \sum_{\beta = a, c} \left( \int_0^t s^{\frac{63}{64}} \langle s \rangle^{\frac{1}{50}} \Vert e^{-i s \omega(D)} \nabla h_{\beta}(s) \Vert_{L^4}^4 ~ ds \right)^{\frac{1}{4}} \right) \\
&\quad \quad \quad \left( \int_0^t s^{\frac{63}{64}} \langle s \rangle^{\frac{1}{50}} \Vert e^{-i s \omega(D)} \nabla m_{\widehat{\mathcal{L}}}(D) \nabla h_{\alpha}(s) \Vert_{L^4}^4 ~ ds \right)^{\frac{1}{2}} \\
\eqref{estdispL4integralepreuve-1-3} &\lesssim \int_0^t s^{\frac{63}{64}} \langle s \rangle^{\frac{1}{50}} \Vert e^{-i s \omega(D)} \nabla m_{\widehat{\mathcal{L}}}(D) \nabla h_{\alpha}(s) \Vert_{L^4}^3 \Vert e^{-i s \omega(D)} \partial_s f(s) \Vert_{L^4} ~ ds \\
&\lesssim \left( \int_0^t s^{\frac{63}{64}} \langle s \rangle^{\frac{1}{50}} \Vert e^{-i s \omega(D)} \nabla m_{\widehat{\mathcal{L}}}(D) \nabla h_{\alpha}(s) \Vert_{L^4}^4 ~ ds \right)^{\frac{3}{4}} \left( \int_0^t s^{\frac{63}{64}} \langle s \rangle^{\frac{1}{50}} \Vert e^{-i s \omega(D)} \partial_s f(s) \Vert_{L^4}^4 ~ ds \right)^{\frac{1}{4}} \\
&\lesssim \Vert u \Vert_X \left( \int_0^t s^{\frac{63}{64}} \langle s \rangle^{\frac{1}{50}} \Vert e^{-i s \omega(D)} \nabla m_{\widehat{\mathcal{L}}}(D) \nabla h_{\alpha}(s) \Vert_{L^4}^4 ~ ds \right)^{\frac{3}{4}} \\
\eqref{estdispL4integralepreuve-1-4} &\lesssim t^{\frac{127}{64}} \langle t \rangle^{\frac{1}{50}} \Vert h_{\alpha}(t) \Vert_{L^2} \Vert e^{-i t \omega(D)} \nabla h_{\alpha}(t) \Vert_{L^4}^2 \Vert m_{\widehat{\mathcal{L}}}(D) \nabla u(t) \Vert_{L^{\infty}} \\
&\lesssim t^{\frac{1}{6}} \langle t \rangle^{-\frac{1}{6}-\frac{379}{4800}+O(\delta)} \Vert u \Vert_X^4 \\
&\lesssim \Vert u \Vert_X^4
\end{align*}

Hence, we showed that
\begin{align*}
&\int_0^t s^{\frac{63}{64}} \langle s \rangle^{\frac{1}{50}} \Vert e^{-i s \omega(D)} m_{\widehat{\mathcal{L}}}(D) \nabla h_{\alpha}(s) \Vert_{L^4}^4 ~ ds \\
&\quad \lesssim \Vert u \Vert_X^4 + \Vert u \Vert_X^2 \left( \int_0^t s^{\frac{63}{64}} \langle s \rangle^{\frac{1}{50}} \Vert e^{-i s \omega(D)} m_{\widehat{\mathcal{L}}}(D) \nabla h_{\alpha}(s) \Vert_{L^4}^4 ~ ds \right)^{\frac{1}{2}} \\
&\quad \quad \quad + \Vert u \Vert_X \left( \int_0^t s^{\frac{63}{64}} \langle s \rangle^{\frac{1}{50}} \Vert e^{-i s \omega(D)} m_{\widehat{\mathcal{L}}}(D) \nabla h_{\alpha}(s) \Vert_{L^4}^4 ~ ds \right)^{\frac{3}{4}} \\
&\quad \quad \quad + \Vert u \Vert_X^2 \left( \int_0^t s^{\frac{63}{64}} \langle s \rangle^{\frac{1}{50}} \Vert e^{-i s \omega(D)} m_{\widehat{\mathcal{L}}}(D) \nabla h_{\alpha}(s) \Vert_{L^4}^4 ~ ds \right)^{\frac{1}{2}} \\
&\quad \quad \quad \quad \quad \sum_{\beta = a, c} \left( \int_0^t s^{\frac{63}{64}} \langle s \rangle^{\frac{1}{50}} \Vert e^{-i s \omega(D)} \nabla h_{\beta}(s) \Vert_{L^4}^4 ~ ds \right)^{\frac{1}{2}}
\end{align*}
In particular, up to showing a similar estimate for $A = \widehat{\mathcal{C}}$, we will deduce the wanted inequality by absorbing the last term by bootstrap and noting that the quantity on the left cannot (for fixed $t \in (0, \infty)$) be infinite thanks to Proposition \ref{propositionestimeedispersiveL4totaleponctuelle}. 

Finally, if $A = \widehat{\mathcal{C}}$, we decompose using Minkowski's inequality: 
\begin{align*}
&\left( \int_0^t s^{\frac{63}{64}} \langle s \rangle^{\frac{1}{50}} \Vert e^{-i s \omega(D)} \nabla m_A(D) h_{\alpha}(s) \Vert_{L^4}^4 ~ ds \right)^{\frac{1}{4}} \\
&\lesssim \sum_{j \in \mathbb{Z}} \left( \int_0^t s^{\frac{63}{64}} \langle s \rangle^{\frac{1}{50}} \Vert e^{-i s \omega(D)} \psi_j(D) m_{\widehat{\mathcal{C}}}(D) \nabla h_{\alpha}(s) \Vert_{L^4}^4 ~ ds \right)^{\frac{1}{4}} 
\end{align*}
Then, for fixed $j$, the integral is finite by Lemma \ref{lemestimeeL4displocalezoneC} and it is enough to estimate each term in a way that will be summable in $j$. 

However, for fixed $j$, we can reuse exactly the same decomposition and integrations by parts as in the case $A = \widehat{\mathcal{L}}$. All the estimates are identical up to gaining a factor $2^{\delta j} \langle 2^j \rangle^{-2\delta}$ by applying the $L^{\infty}$ estimate at localised frequency, Lemma \ref{lemestdispvoiscone}. The only different term is the boundary term \eqref{estdispL4integralepreuve-1-4}, for which we apply either Lemma \ref{lemestimeeL4displocalezoneC} or Lemma \ref{lemestdispL4BF}: if $2^j \lesssim t^{-\frac{1}{3}}$, then 
\begin{align*}
\eqref{estdispL4integralepreuve-1-4}_j &\lesssim t^{\frac{127}{64}} \langle t \rangle^{\frac{1}{50}} 2^{2j} \Vert e^{-i t \omega(D)} \psi_j(D) h_{\alpha}(t) \Vert_{L^4}^2 \Vert h_{\alpha}(t) \Vert_{L^2} \Vert \partial_x u(t) \Vert_{L^{\infty}} \\
&\lesssim t^{\frac{127}{64}-\frac{5}{6}} \langle t \rangle^{\frac{1}{50}-\frac{1}{4}+100\delta} 2^{\frac{7j}{2}} \langle 2^j \rangle^{-2} \Vert u \Vert_X^4
\end{align*}
This is summable in $j$, and the sum is
\begin{align*}
&t^{\frac{31}{64}+\frac{1}{6}} \langle t \rangle^{\frac{1}{50}-\frac{11}{12}+100\delta} \Vert u \Vert_X^4 \\
&\quad \lesssim \Vert u \Vert_X^4 
\end{align*}
On the other hand, if $2^j \gtrsim t^{-\frac{1}{3}}$, 
\begin{align*}
\eqref{estdispL4integralepreuve-1-4}_j &\lesssim t^{\frac{127}{64}} \langle t \rangle^{\frac{1}{50}} 2^{\frac{3j}{2}} \langle 2^j \rangle^{-1} \Vert e^{-i t \omega(D)} \psi_j(D) h_{\alpha}(t) \Vert_{L^4}^2 \Vert \langle \nabla \rangle h_{\alpha}(t) \Vert_{L^2} \Vert m_{\widehat{\mathcal{C}}}(D) |\nabla|^{\frac{3}{2}} u(t) \Vert_{L^{\infty}} \\
&\lesssim t^{\frac{127}{64}-\frac{5}{6}} \langle t \rangle^{\frac{1}{50}-\frac{1}{3}+100\delta} 2^{\frac{3j}{2}} \langle 2^j \rangle^{-1} \left( t^{-\frac{15}{32}} \langle t \rangle^{\delta} 2^{-\frac{21j}{32}+10\delta |j|} \right)^2 \Vert u \Vert_X^4 \\
&\lesssim t^{\frac{3}{64}+\frac{1}{6}} \langle t \rangle^{\frac{1}{50}-\frac{1}{3}+102\delta} 2^{\frac{3j}{16}+20\delta |j|} \langle 2^j \rangle^{-1} \Vert u \Vert_X^4 
\end{align*}
This is summable in $j$, and the sum is 
\begin{align*}
&t^{\frac{3}{64}+\frac{1}{6}} \langle t \rangle^{\frac{1}{50}-\frac{1}{3}+102\delta} \Vert u \Vert_X^4 \\
&\quad \lesssim \Vert u \Vert_X^4
\end{align*}
The conclusion is the same. 
\end{Dem}

\section{Weighted decomposition} \label{section-decomplemma} 

The goal of this section is to prove a key decomposition lemma, which will allow to define the decomposition
\[ m_{\alpha}(D) X_{\alpha} f(t) = h_{\alpha}(t) + g_{\alpha}(t) \]
introduced in the main Proposition \ref{prop-estimeeapriori-complete}. 

For $\mu$ a symbol, $F_1, F_2$ two functions, denote by 
\begin{align*}
T_{\mu}[F_1, F_2](t) &= \mathcal{F}^{-1} \left( \int e^{i t \varphi} \mu(\overline{\xi}, \overline{\eta}) \widehat{F}_1(t, \overline{\eta}) \widehat{F}_2(t, \overline{\sigma}) ~ d\overline{\eta} \right) 
\end{align*}
We will denote by $T[F_1, F_2](t) = T_1[F_1, F_2]$ where $1$ means the symbol identically equal to $1$. We then define the interaction operator: 
\[ I_{\mu}[F_1, F_2](t) = \int_0^t T_{\mu}[F_1, F_2](s) ~ ds \]
In a similar way, we denote by $I[F_1, F_2](t) = I_1[F_1, F_2](t)$. Note that we have by Duhamel's formula 
\begin{align*}
f(t) = f(0) + \partial_x I[f, f](t) 
\end{align*}
by definition of $I$. In particular, applying the weights and computing commutators, we see that 
\begin{align*}
m_a(D) X_a f(t) &= m_a(D) X_a f(0) + \partial_x \left( m_a(D) X_a I[f, f](t) \right) + \frac{\partial_x}{|\nabla|} I[f, f](t) \\
m_b(D) X_b f(t) &= m_b(D) X_b f(0) + \partial_x \left( m_b(D) X_b I[f, f](t) \right) + \frac{i |\nabla_y|}{|\nabla|} I[f, f](t) \\
m_c(D) X_c f(t) &= m_c(D) X_c f(0) + \partial_x \left( m_c(D) X_c I[f, f](t) \right) 
\end{align*}

The following bilinear decomposition lemma is the goal of this section: 

\begin{Lem} Let $F_1, F_2$ be two sufficiently regular functions and $\mu$ a symbol of order $0$. Then: 
\begin{enumerate}
\item We can decompose: 
\begin{align*}
&|\nabla| m_a(D) X_a I_{\mu}[F_1, F_2](t) \\
&= I_{\mu}\left[ |\nabla| m_a(D) X_a F_1, ~ F_2 \right](t) + I_{\mu}\left[ F_1, ~ |\nabla| m_a(D) X_a F_2 \right](t) 
+ 3 t T_{\mu}[F_1, F_2](t) \\
&\quad - 3 I_{\mu}\left[ t \partial_t F_1, ~ F_2 \right](t) - 3 I_{\mu}\left[ F_1, ~ t \partial_t F_2 \right](t) 
+ I_{\widetilde{\mu}}[F_1, F_2](t) 
\end{align*}
for some $0$ order symbol $\widetilde{\mu}$. If $\mu \equiv 1$, $\widetilde{\mu} \equiv 0$. 
\item We can decompose: 
\begin{align*}
&|\nabla| m_c(D) X_c I_{\mu}[F_1, F_2](t) \\
&= I_{\mu}\left[ |\nabla| m_c(D) X_c F_1, ~ F_2 \right](t) + I_{\mu}\left[ F_1, ~ |\nabla| m_c(D) X_c F_2 \right](t) 
+ I_{\widetilde{\mu}}[F_1, F_2](t) 
\end{align*}
for some $0$ order symbol $\widetilde{\mu}$. If $\mu \equiv 1$, $\widetilde{\mu} \equiv 0$. 
\item There exist symbols $m, \widetilde{\mu}$ of order $0$ (that can below vary from line to line) and a symbol $m_g$ (the same at each line where it appears) such that 
\begin{subequations}
\label{lemdecequHGtot} 
\begin{align}
&i \xi_0 m_b(\overline{\xi}) \widehat{X}_b(\overline{\xi}) \cdot \nabla_{\overline{\xi}} \widehat{I}_{\mu}[F_1, F_2](t, \overline{\xi}) \notag \\
&= \int_0^t \int e^{i s \varphi} i \mu(\overline{\xi}, \overline{\eta}) \eta_0 m_b(\overline{\eta}) \widehat{X}_b(\overline{\eta}) \cdot \nabla_{\overline{\eta}} \widehat{F}_1(s, \overline{\eta}) \widehat{F}_2(s, \overline{\sigma}) ~ d\overline{\eta} ds \label{lemdecequH-01-symeta} \\
&+ \int_0^t \int e^{i s \varphi} i \mu(\overline{\xi}, \overline{\eta}) \widehat{F}_1(s, \overline{\eta}) \sigma_0 m_b(\overline{\sigma}) \widehat{X}_b(\overline{\sigma}) \cdot \nabla_{\overline{\xi}} \widehat{F}_2(s, \overline{\sigma}) ~ d\overline{\eta} ds \label{lemdecequH-02-symsigma} \\
&+ \sum_{\alpha = a, c} \int_0^t \int e^{i s \varphi} \mu(\overline{\xi}, \overline{\eta}) m(\overline{\xi}, \overline{\eta}) \frac{|\overline{\eta}| \sigma_0}{|\overline{\eta}| + |\overline{\sigma}|} m_{\alpha}(\overline{\eta}) \widehat{X}_{\alpha}(\overline{\eta}) \cdot \nabla_{\overline{\eta}} \widehat{F}_1(s, \overline{\eta}) \widehat{F}_2(s, \overline{\sigma}) ~ d\overline{\eta} ds \label{lemdecequH-03-resxetaac} \\
&+ \int_0^t \int e^{i s \varphi} \mu(\overline{\xi}, \overline{\eta}) m(\overline{\xi}, \overline{\eta}) \frac{\eta_0 \sigma_0}{|\overline{\eta}| + |\overline{\sigma}|} m_b(\overline{\eta}) \widehat{X}_b(\overline{\eta}) \cdot \nabla_{\overline{\eta}} \widehat{F}_1(s, \overline{\eta}) \widehat{F}_2(s, \overline{\sigma}) ~ d\overline{\eta} ds \label{lemdecequH-04-resxetabbon} \\
&+ \int_0^t \int e^{i s \varphi} \mu(\overline{\xi}, \overline{\eta}) m(\overline{\xi}, \overline{\eta}) \frac{\eta_0 \sigma_0 |\overline{\eta}|}{(|\eta_0|+|\sigma_0|)(|\overline{\eta}| + |\overline{\sigma}|)} m_b(\overline{\eta}) \widehat{X}_b(\overline{\eta}) \cdot \nabla_{\overline{\eta}} \widehat{F}_1(s, \overline{\eta}) \widehat{F}_2(s, \overline{\sigma}) ~ d\overline{\eta} ds \label{lemdecequH-05-resxetabbonbis} \\
&+ \int_0^t \int e^{i s \varphi} \mu(\overline{\xi}, \overline{\eta}) m(\overline{\xi}, \overline{\eta}) \frac{\varphi m_{\widehat{\mathcal{P}}}(\overline{\eta})}{(|\overline{\eta}| + |\overline{\sigma}|)^2} m_b(\overline{\eta}) \widehat{X}_b(\overline{\eta}) \cdot \nabla_{\overline{\eta}} \widehat{F}_1(s, \overline{\eta}) \widehat{F}_2(s, \overline{\sigma}) ~ d\overline{\eta} ds \label{lemdecequH-06-resxetabphi} \\
&+ \sum_{\alpha = a, c} \int_0^t \int e^{i s \varphi} \mu(\overline{\xi}, \overline{\eta}) m(\overline{\xi}, \overline{\eta}) \frac{|\overline{\sigma}| \eta_0}{|\overline{\eta}| + |\overline{\sigma}|} \widehat{F}_1(s, \overline{\eta}) m_{\alpha}(\overline{\sigma}) \widehat{X}_{\alpha}(\overline{\sigma}) \cdot \nabla_{\overline{\eta}} \widehat{F}_2(s, \overline{\sigma}) ~ d\overline{\eta} ds \label{lemdecequH-07-resxsigmaac} \\
&+ \int_0^t \int e^{i s \varphi} \mu(\overline{\xi}, \overline{\eta}) m(\overline{\xi}, \overline{\eta}) \frac{\sigma_0 \eta_0}{|\overline{\eta}| + |\overline{\sigma}|} \widehat{F}_1(s, \overline{\eta}) m_b(\overline{\sigma}) \widehat{X}_b(\overline{\sigma}) \cdot \nabla_{\overline{\eta}} \widehat{F}_2(s, \overline{\sigma}) ~ d\overline{\eta} ds \label{lemdecequH-08-resxsigmabbon} \\
&+ \int_0^t \int e^{i s \varphi} \mu(\overline{\xi}, \overline{\eta}) m(\overline{\xi}, \overline{\eta}) \frac{\sigma_0 \eta_0 |\overline{\sigma}|}{(|\eta_0| + |\sigma_0|) (|\overline{\eta}| + |\overline{\sigma}|)} \widehat{F}_1(s, \overline{\eta}) m_b(\overline{\sigma}) \widehat{X}_b(\overline{\sigma}) \cdot \nabla_{\overline{\eta}} \widehat{F}_2(s, \overline{\sigma}) ~ d\overline{\eta} ds \label{lemdecequH-09-resxsigmabbonbis} \\
&+ \int_0^t \int e^{i s \varphi} \mu(\overline{\xi}, \overline{\eta}) m(\overline{\xi}, \overline{\eta}) \frac{\varphi m_{\widehat{\mathcal{P}}}(\overline{\sigma})}{(|\overline{\eta}| + |\overline{\sigma}|)^2} \widehat{F}_1(s, \overline{\eta}) m_b(\overline{\sigma}) \widehat{X}_b(\overline{\sigma}) \cdot \nabla_{\overline{\eta}} \widehat{F}_2(s, \overline{\sigma}) ~ d\overline{\eta} ds \label{lemdecequH-10-resxsigmabphi} \\
&+ \int_0^t \int e^{i s \varphi} \widetilde{\mu}(\overline{\xi}, \overline{\eta}) \widehat{F}_1(s, \overline{\eta}) \widehat{F}_2(s, \overline{\sigma}) ~ d\overline{\eta} ds \label{lemdecequH-11-dersymb} \\
&+ \int_0^t \int s e^{i s \varphi} \mu(\overline{\xi}, \overline{\eta}) m_g(\overline{\xi}, \overline{\eta}) \partial_s \widehat{F}_1(s, \overline{\eta}) \widehat{F}_2(s, \overline{\sigma}) ~ d\overline{\eta} ds \label{lemdecequH-12-resteta} \\
&+ \int_0^t \int s e^{i s \varphi} \mu(\overline{\xi}, \overline{\eta}) m_g(\overline{\xi}, \overline{\eta}) \widehat{F}_1(s, \overline{\eta}) \partial_s \widehat{F}_2(s, \overline{\sigma}) ~ d\overline{\eta} ds \label{lemdecequH-13-restsigma} \\
&- \int t e^{i t \varphi} \mu(\overline{\xi}, \overline{\eta}) m_g(\overline{\xi}, \overline{\eta}) \widehat{F}_1(t, \overline{\eta}) \widehat{F}_2(t, \overline{\sigma}) ~ d\overline{\eta} \label{lemdecequG} 
\end{align}
\end{subequations}
Furthermore, there exists a decomposition 
\begin{align}
m_g(\overline{\xi}, \overline{\eta}) = m_b(\overline{\xi}) m_1(\overline{\xi}, \overline{\eta}) + \nabla_{\overline{\eta}} \varphi(\overline{\xi}, \overline{\eta}) ~ m_2(\overline{\xi}, \overline{\eta}) \label{decompositionfinemg-gainkout} 
\end{align}
for some symbols $m_1, m_2$. 
\end{enumerate} \label{lemdecompositionFGH} 
\end{Lem}

\begin{Rem} The fact that the symbol $m_g$ is exactly the same in \eqref{lemdecequH-12-resteta}, \eqref{lemdecequH-13-restsigma} and \eqref{lemdecequG} will be of particular importance in the estimates. Note that this property is however automatic, since these terms will come from an integration by parts in time, for which \eqref{lemdecequG} is the boundary term. 
\end{Rem}

\subsection{Decomposition for \texorpdfstring{$X_a$}{Xa}} 

Starting from $\widehat{I}$, we apply $|\overline{\xi}| \widehat{X}_a(\overline{\xi}) \cdot \nabla_{\overline{\xi}} = \overline{\xi} \cdot \nabla_{\overline{\xi}}$ (recall $m_a \equiv 1$): 
\begin{subequations}
\begin{align}
|\overline{\xi}| \widehat{X}_a(\overline{\xi}) \cdot \nabla_{\overline{\xi}} \widehat{I}_{\mu}[F_1, F_2](t, \overline{\xi}) 
&= \int_0^t \int i s \overline{\xi} \cdot \nabla_{\overline{\xi}} \varphi e^{i s \varphi} \mu(\overline{\xi}, \overline{\eta}) \widehat{F}_1(s, \overline{\eta}) \widehat{F}_2(s, \overline{\sigma}) ~ d\overline{\eta} ds \label{equdecXaf-1} \\
&+ \int_0^t \int e^{i s \varphi} \mu(\overline{\xi}, \overline{\eta}) \widehat{F}_1(s, \overline{\eta}) \overline{\xi} \cdot \nabla_{\overline{\xi}} \widehat{F}_2(s, \overline{\sigma}) ~ d\overline{\eta} ds \label{equdecXaf-2} \\
&+ \int_0^t \int e^{i s \varphi} \overline{\xi} \cdot \nabla_{\overline{\xi}} \mu(\overline{\xi}, \overline{\eta}) \widehat{F}_1(s, \overline{\eta}) \widehat{F}_2(s, \overline{\sigma}) ~ d\overline{\eta} ds \label{equdecXaf-3} 
\end{align}
\end{subequations}
We then write that: 
\begin{align}
\overline{\xi} \cdot \nabla_{\overline{\xi}} \widehat{F}_2(s, \overline{\sigma}) &= - \overline{\eta} \cdot \nabla_{\overline{\eta}} \widehat{F}_2(s, \overline{\sigma}) + \overline{\sigma} \cdot \nabla_{\overline{\xi}} \widehat{F}_2(s, \overline{\sigma}) \label{decompositionderiveeaxitosigma} 
\end{align}
so that we can decompose: 
\begin{subequations}
\begin{align}
\eqref{equdecXaf-2} &= \int_0^t \int e^{i s \varphi} \mu(\overline{\xi}, \overline{\eta}) \widehat{F}_1(s, \overline{\eta}) |\overline{\sigma}| \widehat{X}_a(\overline{\sigma}) \cdot \nabla_{\overline{\xi}} \widehat{F}_2(s, \overline{\sigma}) ~ d\overline{\eta} ds \label{equdecXaf-2-1} \\
&- \int_0^t \int e^{i s \varphi} \mu(\overline{\xi}, \overline{\eta}) \widehat{F}_1(s, \overline{\eta}) \overline{\eta} \cdot \nabla_{\overline{\eta}} \widehat{F}_2(s, \overline{\sigma}) ~ d\overline{\eta} ds \label{equdecXaf-2-2} 
\end{align}
\end{subequations}
We then apply an integration by parts in $\overline{\eta}$ on \eqref{equdecXaf-2-2}: 
\begin{subequations}
\begin{align}
\eqref{equdecXaf-2-2} &= \int_0^t \int i s \overline{\eta} \cdot \nabla_{\overline{\eta}} \varphi e^{i s \varphi} \mu(\overline{\xi}, \overline{\eta}) \widehat{F}_1(s, \overline{\eta}) \widehat{F}_2(s, \overline{\sigma}) ~ d\overline{\eta} ds \label{equdecXaf-2-2-1} \\
&+ \int_0^t \int e^{i s \varphi} |\overline{\eta}| \widehat{X}_a(\overline{\eta}) \cdot \nabla_{\overline{\eta}} \widehat{F}_1(s, \overline{\eta}) \widehat{F}_2(s, \overline{\sigma}) ~ d\overline{\eta} ds \label{equdecXaf-2-2-2} \\
&+ \int_0^t \int e^{i s \varphi} \nabla_{\overline{\eta}} \cdot (\overline{\eta} \mu(\overline{\xi}, \overline{\eta})) \widehat{F}_1(s, \overline{\eta}) \widehat{F}_2(s, \overline{\sigma}) ~ d\overline{\eta} ds \label{equdecXaf-2-2-3} 
\end{align}
\end{subequations}
Now we group \eqref{equdecXaf-1} and \eqref{equdecXaf-2-2-1}: 
\begin{align*}
\eqref{equdecXaf-1} + \eqref{equdecXaf-2-2-1} &= \int_0^t \int i s \left( \overline{\xi} \cdot \nabla_{\overline{\xi}} \varphi + \overline{\eta} \cdot \nabla_{\overline{\eta}} \varphi \right) e^{i s \varphi} \mu(\overline{\xi}, \overline{\eta}) \widehat{F}_1(s, \overline{\eta}) \widehat{F}_2(s, \overline{\sigma}) ~ d\overline{\eta} ds
\end{align*}
But
\begin{align}
\overline{\xi} \cdot \nabla_{\overline{\xi}} \varphi + \overline{\eta} \cdot \nabla_{\overline{\eta}} \varphi = 3 \varphi \label{relationdecompositionchampa} 
\end{align}
by homogeneity of $\varphi$ (or equivalently by scaling invariance of the equation). We can therefore apply an integration by parts in time: 
\begin{subequations}
\begin{align}
\eqref{equdecXaf-1} + \eqref{equdecXaf-2-2-1} &= 3 \int_0^t \int i s \varphi e^{i s \varphi} \mu(\overline{\xi}, \overline{\eta}) \widehat{F}_1(s, \overline{\eta}) \widehat{F}_2(s, \overline{\sigma}) ~ d\overline{\eta} ds \notag \\
&= - 3 \int_0^t \int e^{i s \varphi} \mu(\overline{\xi}, \overline{\eta}) \widehat{F}_1(s, \overline{\eta}) \widehat{F}_2(s, \overline{\sigma}) ~ d\overline{\eta} ds \label{equdecXaf-1-1} \\
&- 3 \int_0^t \int s e^{i s \varphi} \mu(\overline{\xi}, \overline{\eta}) \partial_s \left( \widehat{F}_1(s, \overline{\eta}) \widehat{F}_2(s, \overline{\sigma}) \right) ~ d\overline{\eta} ds \label{equdecXaf-1-2} \\
&+ 3 \int t e^{i t \varphi} \mu(\overline{\xi}, \overline{\eta}) \widehat{F}_1(t, \overline{\eta}) \widehat{F}_2(t, \overline{\sigma}) ~ d\overline{\eta} \label{equdecXaf-1-3} 
\end{align}
\end{subequations}
Grouping now all the terms and coming back to the physical space, we have indeed: 
\begin{align*}
&|\nabla| X_a I_{\mu}[F_1, F_2](t) \\
&= \mathcal{F}^{-1} \left( \eqref{equdecXaf-2-2-2} + \eqref{equdecXaf-2-1} + \eqref{equdecXaf-1-3} + \eqref{equdecXaf-1-2} + \left( \eqref{equdecXaf-3} + \eqref{equdecXaf-2-2-3} + \eqref{equdecXaf-1-1} \right) \right) \\
&= I_{\mu}\left[ |\nabla| m_a(D) X_a F_1, ~ F_2 \right](t) + I_{\mu}\left[ F_1, ~ |\nabla| m_a(D) X_a F_2 \right](t) 
+ 3 t T_{\mu}[F_1, F_2](t) \\
&\quad - 3 I_{\mu}\left[ \partial_t F_1, ~ F_2 \right](t) - 3 I_{\mu}\left[ F_1, ~ \partial_t F_2 \right](t) 
+ I_{\widetilde{\mu}}[F_1, F_2](t) 
\end{align*}
Furthermore, we have the expression: 
\begin{align*}
\widetilde{\mu}(\overline{\xi}, \overline{\eta}) &= \overline{\xi} \cdot \nabla_{\overline{\xi}} \mu(\overline{\xi}, \overline{\eta}) + \nabla_{\overline{\eta}} \cdot (\overline{\eta} \mu(\overline{\xi}, \overline{\eta})) - 3 \mu(\overline{\xi}, \overline{\eta}) \\
&= \overline{\xi} \cdot \nabla_{\overline{\xi}} \mu(\overline{\xi}, \overline{\eta}) + \overline{\eta} \cdot \nabla_{\overline{\eta}} \mu(\overline{\xi}, \overline{\eta})
\end{align*}
In particular, if $\mu \equiv 1$, then $\widetilde{\mu} \equiv 0$. 

This concludes the proof of Lemma \ref{lemdecompositionFGH} for $\alpha = a$. 

\subsection{Decomposition for \texorpdfstring{$X_c$}{Xc}} 

We proceed the same way as for $X_a$. All steps are the same, rewriting 
\[ J \xi \cdot \nabla_{\xi} \widehat{f}(s, \overline{\sigma}) = - J \eta \cdot \nabla_{\eta} \widehat{f}(s, \overline{\sigma}) + J \sigma \cdot \nabla_{\xi} \widehat{f}(s, \overline{\sigma}) \]
instead of \eqref{decompositionderiveeaxitosigma}, and using the identity
\begin{align*}
J \xi \cdot \nabla_{\xi} \varphi + J \eta \cdot \nabla_{\eta} \varphi = 0 
\end{align*}
instead of \eqref{relationdecompositionchampa}. We also use the fact that
\[ |\overline{\xi}| m_c(\overline{\xi}) \widehat{X}_c(\overline{\xi}) = \begin{pmatrix} 0 \\ J \xi \end{pmatrix} \]
to absorb the singularity near $\xi = 0$. Finally, for the expression of $\widetilde{\mu}$, we find: 
\begin{align*}
\widetilde{\mu}(\overline{\xi}, \overline{\eta}) &= J\xi \cdot \nabla_{\xi} \mu(\overline{\xi}, \overline{\eta}) + \nabla_{\eta} \cdot (J \eta \mu(\overline{\xi}, \overline{\eta})) \\
&= J \xi \cdot \nabla_{\xi} \mu(\overline{\xi}, \overline{\eta}) + J \eta \cdot \nabla_{\eta} \mu(\overline{\xi}, \overline{\eta}) 
\end{align*}
Again, if $\mu \equiv 1$, then $\widetilde{\mu} \equiv 0$. 

This concludes the proof of Lemma \ref{lemdecompositionFGH} for $\alpha = c$. 

\subsection{Decomposition for \texorpdfstring{$X_b$}{Xb}} 

When we apply $m_b \widehat{X}_b$, we cannot find an analogous relationship because $X_b$ is not related to an invariance of the equation. We will therefore try to make use of the cancellations of $m_b$ in order to get nonetheless the announced decomposition. 

To that end, we decompose the frequency space according to the following localisation symbols: 
\begin{align*}
\mu_{BBB}(\overline{\xi}, \overline{\eta}), \quad \mu_{HBH}(\overline{\xi}, \overline{\eta}), \quad \mu_{HHB}(\overline{\xi}, \overline{\eta}), \quad \mu_{BHH}(\overline{\xi}, \overline{\eta}) 
\end{align*}
which localise respectively on $\{ |\overline{\xi}| \simeq |\overline{\eta}| \simeq |\overline{\sigma}| \}$, $\{ |\overline{\xi}| \simeq |\overline{\sigma}| \gg |\overline{\eta}| \}$, $\{ |\overline{\xi}| \simeq |\overline{\eta}| \gg |\overline{\sigma}| \}$, and $\{ |\overline{\xi}| \ll |\overline{\eta}| \simeq |\overline{\sigma}| \}$. 

Furthermore, we will also decompose the interactions depending on the direction of $\overline{\eta}$: 
\[ 1 = m_{\widehat{\mathcal{R}}}(\overline{\eta}) + m_{\widehat{\mathcal{C}}}(\overline{\eta}) + m_{\widehat{\mathcal{L}}}(\overline{\eta}) + m_{\widehat{\mathcal{P}}}(\overline{\eta}) \]
and likewise for $\overline{\sigma}$. 

Let us denote, for $A_1, A_2 \in \{ \widehat{\mathcal{R}}, \widehat{\mathcal{C}}, \widehat{\mathcal{L}}, \widehat{\mathcal{P}} \}$, 
\begin{align*}
\widehat{I}_{\mu}^{(BBB, A_1A_2)}[F_1, F_2](t, \overline{\xi}):= \int_0^t \int e^{i s \varphi} \mu(\overline{\xi}, \overline{\eta}) \mu_{BBB}(\overline{\xi}, \overline{\eta}) m_{A_1}(\overline{\eta}) m_{A_2}(\overline{\sigma}) \widehat{F}_1(s, \overline{\eta}) \widehat{F}_2(s, \overline{\sigma}) ~ d\overline{\eta} ds
\end{align*}
and 
\begin{align*}
\widehat{I}_{\mu}^{(HBH, A_1)}[F_1, F_2](t, \overline{\xi}) := \int_0^t \int e^{i s \varphi} \mu(\overline{\xi}, \overline{\eta}) \mu_{HBH}(\overline{\xi}, \overline{\eta}) m_{A_1}(\overline{\sigma}) \widehat{F}_1(s, \overline{\eta}) \widehat{F}_2(s, \overline{\sigma}) ~ d\overline{\eta} ds
\end{align*}
and in a symmetric way the $I_{\mu}^{(HHB, A_1)}$. If $A_1 \neq \widehat{\mathcal{R}}$, we set: 
\begin{align*}
\widehat{I}_{\mu}^{(BHH, A_1)}[F_1, F_2](t, \overline{\xi}) := \int_0^t \int e^{i s \varphi} \mu(\overline{\xi}, \overline{\eta}) \mu_{BHH}(\overline{\xi}, \overline{\eta}) m_{A_1}(\overline{\sigma}) m_{A_1}(\overline{\eta}) \widehat{F}_1(s, \overline{\eta}) \widehat{F}_2(s, \overline{\sigma}) ~ d\overline{\eta} ds 
\end{align*}
while: 
\begin{align*}
\widehat{I}_{\mu}^{(BHH, \widehat{\mathcal{R}})}[F_1, F_2](t, \overline{\xi}) := \int_0^t \int e^{i s \varphi} \mu(\overline{\xi}, \overline{\eta}) \mu_{BHH}(\overline{\xi}, \overline{\eta}) \mu_{\widehat{\mathcal{R}}}(\overline{\xi}, \overline{\eta}) \widehat{F}_1(s, \overline{\eta}) \widehat{F}_2(s, \overline{\sigma}) ~ d\overline{\eta} ds 
\end{align*}
where
\begin{align*}
\mu_{\widehat{\mathcal{R}}}(\overline{\xi}, \overline{\eta}) := 1 - \sum_{A_1 = \widehat{\mathcal{C}}, \widehat{\mathcal{L}}, \widehat{\mathcal{P}}} m_{A_1}(\overline{\sigma}) m_{A_1}(\overline{\eta}) 
\end{align*}
On the support of $\mu_{BHH} \mu_{\widehat{\mathcal{R}}}$, $\overline{\eta}$ and $\overline{\sigma}$ are localised away from $\widehat{\mathcal{C}} \cup \widehat{\mathcal{L}} \cup \widehat{\mathcal{P}}$. 

In the $BBB$ case, we have therefore 16 interactions to consider: using the fact that the expression in Lemma \ref{lemdecompositionFGH} is symmetric in $F_1, F_2$ (up to the $\mu$ symbol), we can symmetrize the variables $\overline{\eta}$ and $\overline{\sigma}$ and consider only 10 interactions. In the cases $HBH$ and $BHH$, there are each time 4 interactions to consider, while the case $HHB$ is symmetric to the case $HBH$. 

Note that \eqref{lemdecequH-01-symeta} (respectively \eqref{lemdecequH-02-symsigma}), which is the only term for which the symbol is imposed, is of the form \eqref{lemdecequH-04-resxetabbon} (respectively \eqref{lemdecequH-08-resxsigmabbon}) as soon as $|\sigma_0| \gtrsim |\overline{\eta}| + |\overline{\sigma}|$ (respectively $|\eta_0| \gtrsim |\overline{\eta}| + |\overline{\sigma}|$). It is therefore needed to have it appear explicitely only if $\overline{\sigma}$ is in the neighborhood of the plane or if the interaction is of type HHB (respectively if $\overline{\eta}$ is in the neighborhood of the plane or if the interaction is of type $HBH$).  

\paragraph{Notation} We will use the following important notation. We will write 
\begin{align*}
\mu_1 = O(\mu_2)
\end{align*}
if there exist a symbol $\mu_3$ such that $\mu_1 = \mu_2 \mu_3$. We will write 
\begin{align*}
\mu_1 = o(\mu_2)
\end{align*}
if there exist a symbol $\mu_3$ such that $\mu_1 = \mu_2 \mu_3$ and $\Vert \mu_3 \Vert_{L^{\infty}} \ll 1$. 

These notations will be used extensively as they will allow the computations to be (a bit) more compact. 

\subsubsection{Interaction \texorpdfstring{$(BBB, \widehat{\mathcal{R}}\widehat{\mathcal{R}})$}{(BBB, RR)}} 

Let us consider
\begin{subequations}
\begin{align}
&\xi_0 m_b(\overline{\xi}) \widehat{X}_b(\overline{\xi}) \cdot \nabla_{\overline{\xi}} \widehat{I}_{\mu}^{(BBB, \widehat{\mathcal{R}})\widehat{\mathcal{R}}}[F_1, F_2](t, \overline{\xi}) \notag \\
&\quad = \int_0^t \int i s \xi_0 m_b(\overline{\xi}) \widehat{X}_b(\overline{\xi}) \cdot \nabla_{\overline{\xi}} \varphi e^{i s \varphi} \mu(\overline{\xi}, \overline{\eta}) \mu_{BBB}(\overline{\xi}, \overline{\eta}) m_{\widehat{\mathcal{R}}}(\overline{\eta}) m_{\widehat{\mathcal{R}}}(\overline{\sigma}) \widehat{F}_1(s, \overline{\eta}) \widehat{F}_2(s, \overline{\sigma}) ~ d\overline{\eta} ds \label{equdecchampbRR-1} \\
&\quad \quad + \int_0^t \int e^{i s \varphi} \xi_0 \mu(\overline{\xi}, \overline{\eta}) \mu_{BBB}(\overline{\xi}, \overline{\eta}) m_{\widehat{\mathcal{R}}}(\overline{\eta}) m_{\widehat{\mathcal{R}}}(\overline{\sigma}) \widehat{F}_1(s, \overline{\eta}) m_b(\overline{\xi}) \widehat{X}_b(\overline{\xi}) \cdot \nabla_{\overline{\xi}} \widehat{F}_2(s, \overline{\sigma}) ~ d\overline{\eta} ds \label{equdecchampbRR-2} \\
&\quad \quad + \int_0^t \int e^{i s \varphi} \xi_0 m_b(\overline{\xi}) \widehat{X}_b(\overline{\xi}) \cdot \nabla_{\overline{\xi}} \left( \mu(\overline{\xi}, \overline{\eta}) \mu_{BBB}(\overline{\xi}, \overline{\eta}) m_{\widehat{\mathcal{R}}}(\overline{\eta}) m_{\widehat{\mathcal{R}}}(\overline{\sigma}) \right) \widehat{F}_1(s, \overline{\eta}) \widehat{F}_2(s, \overline{\sigma}) ~ d\overline{\eta} ds \label{equdecchampbRR-3} 
\end{align}
\end{subequations}
\eqref{equdecchampbRR-2} and \eqref{equdecchampbRR-3} are of the form $\eqref{lemdecequH-07-resxsigmaac} + \eqref{lemdecequH-08-resxsigmabbon}$ and \eqref{lemdecequH-11-dersymb} respectively. 

Then, for \eqref{equdecchampbRR-1}, we apply Lemma \ref{lem-non-res-loin0Cone} to get 
\begin{align*}
\xi_0 &= m_{-2} \varphi + m_{-1} \nabla_{\overline{\eta}} \varphi 
\end{align*}
for locally smooth symbols $m_{-2}, m_{-1}$, homogeneous and of respective degrees $-2, -1$. It is then possible to apply integrations by parts and recover easily terms of the form \eqref{lemdecequH-03-resxetaac} to \eqref{lemdecequG}: 
\begin{align*}
&\eqref{equdecchampbRR-1} \\
&= \int_0^t \int i s m_{-1} \nabla_{\overline{\eta}} \varphi m_b(\overline{\xi}) \widehat{X}_b(\overline{\xi}) \cdot \nabla_{\overline{\xi}} \varphi e^{i s \varphi} \mu(\overline{\xi}, \overline{\eta}) \mu_{BBB}(\overline{\xi}, \overline{\eta}) m_{\widehat{\mathcal{R}}}(\overline{\eta}) m_{\widehat{\mathcal{R}}}(\overline{\sigma}) \widehat{F}_1(s, \overline{\eta}) \widehat{F}_2(s, \overline{\sigma}) ~ d\overline{\eta} ds \\
&\quad + \int_0^t \int i s m_{-2} \varphi m_b(\overline{\xi}) \widehat{X}_b(\overline{\xi}) \cdot \nabla_{\overline{\xi}} \varphi e^{i s \varphi} \mu(\overline{\xi}, \overline{\eta}) \mu_{BBB}(\overline{\xi}, \overline{\eta}) m_{\widehat{\mathcal{R}}}(\overline{\eta}) m_{\widehat{\mathcal{R}}}(\overline{\sigma}) \widehat{F}_1(s, \overline{\eta}) \widehat{F}_2(s, \overline{\sigma}) ~ d\overline{\eta} ds \\
&= - \int_0^t \int e^{i s \varphi} \nabla_{\overline{\eta}} \cdot \left( m_{-1} m_b(\overline{\xi}) \widehat{X}_b(\overline{\xi}) \cdot \nabla_{\overline{\xi}} \varphi \mu(\overline{\xi}, \overline{\eta}) \mu_{BBB}(\overline{\xi}, \overline{\eta}) m_{\widehat{\mathcal{R}}}(\overline{\eta}) m_{\widehat{\mathcal{R}}}(\overline{\sigma}) \right) \widehat{F}_1(s, \overline{\eta}) \widehat{F}_2(s, \overline{\sigma}) ~ d\overline{\eta} ds \\
&\quad - \int_0^t \int e^{i s \varphi} m_{-1} m_b(\overline{\xi}) \widehat{X}_b(\overline{\xi}) \cdot \nabla_{\overline{\xi}} \varphi \mu(\overline{\xi}, \overline{\eta}) \mu_{BBB}(\overline{\xi}, \overline{\eta}) m_{\widehat{\mathcal{R}}}(\overline{\eta}) m_{\widehat{\mathcal{R}}}(\overline{\sigma}) \nabla_{\overline{\eta}} \left( \widehat{F}_1(s, \overline{\eta}) \widehat{F}_2(s, \overline{\sigma}) \right) ~ d\overline{\eta} ds \\
&\quad - \int_0^t \int m_{-2} m_b(\overline{\xi}) \widehat{X}_b(\overline{\xi}) \cdot \nabla_{\overline{\xi}} \varphi e^{i s \varphi} \mu(\overline{\xi}, \overline{\eta}) \mu_{BBB}(\overline{\xi}, \overline{\eta}) m_{\widehat{\mathcal{R}}}(\overline{\eta}) m_{\widehat{\mathcal{R}}}(\overline{\sigma}) \widehat{F}_1(s, \overline{\eta}) \widehat{F}_2(s, \overline{\sigma}) ~ d\overline{\eta} ds \\
&\quad - \int_0^t \int s m_{-2} m_b(\overline{\xi}) \widehat{X}_b(\overline{\xi}) \cdot \nabla_{\overline{\xi}} \varphi e^{i s \varphi} \mu(\overline{\xi}, \overline{\eta}) \mu_{BBB}(\overline{\xi}, \overline{\eta}) m_{\widehat{\mathcal{R}}}(\overline{\eta}) m_{\widehat{\mathcal{R}}}(\overline{\sigma}) \partial_s \left( \widehat{F}_1(s, \overline{\eta}) \widehat{F}_2(s, \overline{\sigma}) \right) ~ d\overline{\eta} ds \\
&\quad + \int t m_{-2} m_b(\overline{\xi}) \widehat{X}_b(\overline{\xi}) \cdot \nabla_{\overline{\xi}} \varphi e^{i t \varphi} \mu(\overline{\xi}, \overline{\eta}) \mu_{BBB}(\overline{\xi}, \overline{\eta}) m_{\widehat{\mathcal{R}}}(\overline{\eta}) m_{\widehat{\mathcal{R}}}(\overline{\sigma}) \widehat{F}_1(t, \overline{\eta}) \widehat{F}_2(t, \overline{\sigma}) ~ d\overline{\eta} ds
\end{align*}
as wanted. We didn't use the presence of $m_b(\overline{\xi})$ so \eqref{decompositionfinemg-gainkout} is automatic. 

\subsubsection{Interaction \texorpdfstring{$(BBB, \widehat{\mathcal{C}}\widehat{\mathcal{R}})$}{(BBB, CR)}} 

Let us consider 
\begin{subequations}
\begin{align}
&\xi_0 m_b(\overline{\xi}) \widehat{X}_b(\overline{\xi}) \cdot \nabla_{\overline{\xi}} \widehat{I}_{\mu}^{(BBB, \widehat{\mathcal{C}}\widehat{\mathcal{R}})}[F_1, F_2](t, \overline{\xi}) \notag \\
&\quad = \int_0^t \int i s \xi_0 m_b(\overline{\xi}) \widehat{X}_b(\overline{\xi}) \cdot \nabla_{\overline{\xi}} \varphi e^{i s \varphi} \mu(\overline{\xi}, \overline{\eta}) \mu_{BBB}(\overline{\xi}, \overline{\eta}) m_{\widehat{\mathcal{C}}}(\overline{\eta}) m_{\widehat{\mathcal{R}}}(\overline{\sigma}) \widehat{F}_1(s, \overline{\eta}) \widehat{F}_2(s, \overline{\sigma}) ~ d\overline{\eta} ds \label{equdecchampbCR-1} \\
&\quad \quad + \int_0^t \int e^{i s \varphi} \xi_0 \mu(\overline{\xi}, \overline{\eta}) \mu_{BBB}(\overline{\xi}, \overline{\eta}) m_{\widehat{\mathcal{C}}}(\overline{\eta}) m_{\widehat{\mathcal{R}}}(\overline{\sigma}) \widehat{F}_1(s, \overline{\eta}) m_b(\overline{\xi}) \widehat{X}_b(\overline{\xi}) \cdot \nabla_{\overline{\xi}} \widehat{F}_2(s, \overline{\sigma}) ~ d\overline{\eta} ds \label{equdecchampbCR-2} \\
&\quad \quad + \int_0^t \int e^{i s \varphi} \xi_0 m_b(\overline{\xi}) \widehat{X}_b(\overline{\xi}) \cdot \nabla_{\overline{\xi}} \left( \mu(\overline{\xi}, \overline{\eta}) \mu_{BBB}(\overline{\xi}, \overline{\eta}) m_{\widehat{\mathcal{C}}}(\overline{\eta}) m_{\widehat{\mathcal{R}}}(\overline{\sigma}) \right) \widehat{F}_1(s, \overline{\eta}) \widehat{F}_2(s, \overline{\sigma}) ~ d\overline{\eta} ds \label{equdecchampbCR-3} 
\end{align}
\end{subequations}
Again, \eqref{equdecchampbCR-2} and \eqref{equdecchampbCR-3} are already in the form $\eqref{lemdecequH-07-resxsigmaac} + \eqref{lemdecequH-08-resxsigmabbon}$ and \eqref{lemdecequH-11-dersymb} respectively.  

Then, for \eqref{equdecchampbCR-1}, we can apply integrations by parts as soon as
\begin{align*}
1 = O(\varphi) + O(\widehat{X}_c(\overline{\eta}) \cdot \nabla_{\overline{\eta}} \varphi) + O(\widehat{X}_a(\overline{\eta}) \cdot \nabla_{\overline{\eta}} \varphi) + O\left( m_b(\overline{\eta}) \nabla_{\overline{\eta}} \varphi \right) 
\end{align*}
Let us therefore assume that 
\begin{align*}
\varphi = o(1), \quad \widehat{X}_c(\overline{\eta}) \cdot \nabla_{\overline{\eta}} \varphi = o(1), \quad \widehat{X}_a(\overline{\eta}) \cdot \nabla_{\overline{\eta}} \varphi = o(1)
\end{align*}

Moreover, we may localise $\overline{\eta}$ by some finer $\widetilde{m}_{\widehat{\mathcal{C}}}$ symbol whose symbol has support disjoint from the one of $m_{\widehat{\mathcal{R}}}$ and such that $\frac{\overline{\eta}_b^{\overline{\eta}}}{|\overline{\eta}|}$ is very small with respect to every other localisation constant: otherwise, we can apply the same decomposition as in the case $(BBB, \widehat{\mathcal{R}}\widehat{\mathcal{R}})$ already treated. 

Since $\overline{\sigma}$ is localised by $m_{\widehat{\mathcal{R}}}$ and $\overline{\eta}$ by $\widetilde{m}_{\widehat{\mathcal{C}}}$, we have that
\begin{align*}
1 &= O\left( \overline{\sigma}_a^{\overline{\sigma}} \right), \quad 1 = O\left( \overline{\sigma}_b^{\overline{\sigma}} \right) 
\end{align*}
and therefore 
\begin{align*}
1 &= O\left( \overline{\sigma}_b^{\overline{\eta}} \right) + O\left( \xi_t^{\overline{\eta} \overline{\sigma}} \right) 
\end{align*}
Hence, if $\widetilde{m}_{\widehat{\mathcal{C}}}$ localises finely enough, 
\begin{align*}
1 = O\left( \nabla_{\overline{\eta}} \varphi \right) 
\end{align*}
by derivating in the direction $\left( \frac{\eta_0}{|\eta_0|} \partial_{\eta_0} - \sqrt{3} \frac{\eta}{|\eta|} \cdot \nabla_{\eta} \right)$. Finally, by Lemma \ref{lemcalculsconecoordonneesconiquesvarphi}, 
\begin{align*}
\xi_t^{\overline{\eta} \overline{\sigma}} &= O\left( \widehat{X}_c(\overline{\eta}) \cdot \nabla_{\overline{\eta}} \varphi \right) \\
\overline{\eta}_b^{\overline{\eta}} &= O\left( m_b(\overline{\eta}) \nabla_{\overline{\eta}} \varphi \right) \\
\left( \overline{\eta}_a^{\overline{\eta}} \right)^2 - \left( \overline{\sigma}_a^{\overline{\eta}} \right)^2 &= O\left( \widehat{X}_a(\overline{\eta}) \cdot \nabla_{\overline{\eta}} \varphi \right) + O\left( \xi_t^{\overline{\eta} \overline{\sigma}} \right) + O\left( \overline{\eta}_b^{\overline{\eta}} \right) 
\end{align*}
There exist then a sign $\epsilon \in \{ -1, 1 \}$ such that, locally, 
\begin{align*}
\overline{\sigma}_a^{\overline{\eta}} &= \epsilon \overline{\eta}_a^{\overline{\eta}} + O\left( \left( \overline{\eta}_a^{\overline{\eta}} \right)^2 - \left( \overline{\sigma}_a^{\overline{\eta}} \right)^2 \right) 
\end{align*}

Note that here, 
\begin{align*}
|\xi| \xi_t^{\overline{\xi} \overline{\eta}} &= \frac{J \xi \cdot \eta}{|\eta|} = \frac{J \sigma \cdot \eta}{|\eta|} = O\left( \xi_t^{\overline{\eta} \overline{\sigma}} \right) \\
|\xi| \xi_t^{\overline{\xi} \overline{\sigma}} &= O\left( \xi_t^{\overline{\eta} \overline{\sigma}} \right) 
\end{align*}
Therefore, by Lemma \ref{lemcalculsconecoordonneesconiquesvarphi} again, and using the linearity $\overline{\xi}_{\alpha}^{\overline{\eta}} = \overline{\eta}_{\alpha}^{\overline{\eta}} + \overline{\sigma}_{\alpha}^{\overline{\eta}}$, we have 
\begin{align*}
6 \sqrt{3} \frac{\eta_0}{|\eta_0|} \varphi &= \left( \overline{\xi}_a^{\overline{\eta}} \right)^3 + \left( \overline{\xi}_b^{\overline{\xi}} \right)^3 - \left( \overline{\eta}_a^{\overline{\eta}} \right)^3 - \left( \overline{\sigma}_a^{\overline{\eta}} \right)^3 - \left( \overline{\sigma}_b^{\overline{\eta}} \right)^3 + O\left( \xi_t^{\overline{\eta} \overline{\sigma}} \right) + O\left( \overline{\eta}_b^{\overline{\eta}} \right) \\
&= \left( 1 + \epsilon \right) \left( 1 + 2 \epsilon \right) \left( \overline{\eta}_a^{\overline{\eta}} \right)^3 + O\left( \xi_t^{\overline{\eta} \overline{\sigma}} \right) + O\left( \overline{\eta}_b^{\overline{\eta}} \right) + O\left( \left( \overline{\eta}_a^{\overline{\eta}} \right)^2 - \left( \overline{\sigma}_a^{\overline{\eta}} \right)^2 \right)
\end{align*}

If $\epsilon = 1$, then $(1 + \epsilon) (1 + 2 \epsilon) = 6$ so locally
\begin{align*}
1 = O\left( \varphi \right) + O\left( \xi_t^{\overline{\eta} \overline{\sigma}} \right) + O\left( \overline{\eta}_b^{\overline{\eta}} \right) + O\left( \left( \overline{\eta}_a^{\overline{\eta}} \right)^2 - \left( \overline{\sigma}_a^{\overline{\eta}} \right)^2 \right)
\end{align*}
which allows to apply integrations by parts in a similar way as what was done in the case $(BBB, \widehat{\mathcal{R}}\widehat{\mathcal{R}})$, up to only differentiating in directions compatible with $\widehat{F}_1(s, \overline{\eta})$ or with a factor $m_b(\overline{\eta})$. 

If $\epsilon = -1$, then 
\begin{align*}
\overline{\xi}_a^{\overline{\eta}} &= \overline{\eta}_a^{\overline{\eta}} + \overline{\sigma}_a^{\overline{\eta}} \\
&= O\left( \left( \overline{\eta}_a^{\overline{\eta}} \right)^2 - \left( \overline{\sigma}_a^{\overline{\eta}} \right)^2 \right) 
\end{align*}
and therefore 
\begin{align*}
3 \xi_0^2 - |\xi|^2 &= 3 \xi_0^2 - \left( \frac{\eta \cdot \xi}{|\eta|} \right)^2 + O\left( \xi_t^{\overline{\eta} \overline{\sigma}} \right) \\
&= \overline{\xi}_a^{\overline{\eta}} \overline{\xi}_b^{\overline{\eta}} + O\left( \xi_t^{\overline{\eta} \overline{\sigma}} \right) \\
&= O\left( \left( \overline{\eta}_a^{\overline{\eta}} \right)^2 - \left( \overline{\sigma}_a^{\overline{\eta}} \right)^2 \right) + O\left( \xi_t^{\overline{\eta} \overline{\sigma}} \right)
\end{align*}
But precisely here 
\begin{align*}
m_b(\overline{\xi}) &= O\left( 3 \xi_0^2 - |\xi|^2 \right) 
\end{align*}
which justifies that the presence of $m_b(\overline{\xi})$ allows the integrations by parts nedded to bring \eqref{equdecchampbCR-1} back to terms of \eqref{lemdecequHGtot}. We saw above that $1 = O(\nabla_{\overline{\eta}} \varphi)$ here, so \eqref{decompositionfinemg-gainkout} is also true. 

\subsubsection{Interaction \texorpdfstring{$(BBB, \widehat{\mathcal{L}}\widehat{\mathcal{R}})$}{(BBB, LR)}}

Let us consider
\begin{subequations}
\begin{align}
&\xi_0 m_b(\overline{\xi}) \widehat{X}_b(\overline{\xi}) \cdot \nabla_{\overline{\xi}} \widehat{I}_{\mu}^{(BBB, \widehat{\mathcal{L}}\widehat{\mathcal{R}})}[F_1, F_2](t, \overline{\xi}) \notag \\
&\quad = \int_0^t \int i s \xi_0 m_b(\overline{\xi}) \widehat{X}_b(\overline{\xi}) \cdot \nabla_{\overline{\xi}} \varphi e^{i s \varphi} \mu(\overline{\xi}, \overline{\eta}) \mu_{BBB}(\overline{\xi}, \overline{\eta}) m_{\widehat{\mathcal{L}}}(\overline{\eta}) m_{\widehat{\mathcal{R}}}(\overline{\sigma}) \widehat{F}_1(s, \overline{\eta}) \widehat{F}_2(s, \overline{\sigma}) ~ d\overline{\eta} ds \label{equdecchampbLR-1} \\
&\quad \quad + \int_0^t \int e^{i s \varphi} \xi_0 \mu(\overline{\xi}, \overline{\eta}) \mu_{BBB}(\overline{\xi}, \overline{\eta}) m_{\widehat{\mathcal{L}}}(\overline{\eta}) m_{\widehat{\mathcal{R}}}(\overline{\sigma}) \widehat{F}_1(s, \overline{\eta}) m_b(\overline{\xi}) \widehat{X}_b(\overline{\xi}) \cdot \nabla_{\overline{\xi}} \widehat{F}_2(s, \overline{\sigma}) ~ d\overline{\eta} ds \label{equdecchampbLR-2} \\
&\quad \quad + \int_0^t \int e^{i s \varphi} \xi_0 m_b(\overline{\xi}) \widehat{X}_b(\overline{\xi}) \cdot \nabla_{\overline{\xi}} \left( \mu(\overline{\xi}, \overline{\eta}) \mu_{BBB}(\overline{\xi}, \overline{\eta}) m_{\widehat{\mathcal{L}}}(\overline{\eta}) m_{\widehat{\mathcal{R}}}(\overline{\sigma}) \right) \widehat{F}_1(s, \overline{\eta}) \widehat{F}_2(s, \overline{\sigma}) ~ d\overline{\eta} ds \label{equdecchampbLR-3} 
\end{align}
\end{subequations}
Again, \eqref{equdecchampbLR-2} and \eqref{equdecchampbLR-3} are already of the form $\eqref{lemdecequH-07-resxsigmaac} + \eqref{lemdecequH-08-resxsigmabbon}$ and \eqref{lemdecequH-11-dersymb} respectively. 

As before, up to reusing the case $(BBB, \widehat{\mathcal{R}}\widehat{\mathcal{R}})$, we can localise $\overline{\eta}$ arbitrarily close to $\widehat{\mathcal{L}}$. 

Then, for \eqref{equdecchampbLR-1}, we have that
\begin{align*}
\nabla_{\eta} \varphi &= 2 \sigma_0 \sigma - 2 \eta_0 \eta = 2 \sigma_0 \sigma + o(1) 
\end{align*}
hence $1 = O(\nabla_{\eta} \varphi)$ (and \eqref{decompositionfinemg-gainkout} is automatic), and 
\[ |\eta| = O(m_b(\overline{\eta}) \widehat{X}_b(\overline{\eta}) \cdot \nabla_{\overline{\eta}} \varphi) + O(m_c(\overline{\eta}) \widehat{X}_c(\overline{\eta}) \cdot \nabla_{\overline{\eta}} \varphi) \]
Then, 
\begin{align*}
\partial_{\eta_0} \varphi &= 3 \eta_0^2 + |\eta|^2 - 3 \sigma_0^2 - |\sigma|^2 \\
&= O(|\eta|) + 3 \eta_0^2 - 3 \sigma_0^2 - |\sigma|^2
\end{align*}
therefore, 
\begin{align*}
3 \eta_0^2 - 3 \sigma_0^2 - |\sigma|^2 = O(\widehat{X}_a(\overline{\eta}) \cdot \nabla_{\overline{\eta}} \varphi) + O(m_b(\overline{\eta}) \widehat{X}_b(\overline{\eta}) \cdot \nabla_{\overline{\eta}} \varphi) + O(m_c(\overline{\eta}) \widehat{X}_c(\overline{\eta}) \cdot \nabla_{\overline{\eta}} \varphi)
\end{align*}
Finally, 
\begin{align*}
\varphi &= \xi_0^3 + \xi_0 |\xi|^2 - \eta_0^3 - \eta_0 |\eta|^2 - \sigma_0^3 - \sigma_0 |\sigma|^2 \\
&= 3 \eta_0^2 \xi_0 + O(\eta) + O(\partial_{\eta_0} \varphi)
\end{align*}
We thus have 
\begin{align*}
\xi_0 = O(\varphi) + O(\widehat{X}_a(\overline{\eta}) \cdot \nabla_{\overline{\eta}} \varphi) + O(m_b(\overline{\eta}) \widehat{X}_b(\overline{\eta}) \cdot \nabla_{\overline{\eta}} \varphi) + O(m_c(\overline{\eta}) \widehat{X}_c(\overline{\eta}) \cdot \nabla_{\overline{\eta}} \varphi)
\end{align*}
and this relation is enough to apply the integrations by parts needed to bring \eqref{equdecchampbLR-1} back to terms of \eqref{lemdecequHGtot}. 

\subsubsection{Interaction \texorpdfstring{$(BBB, \widehat{\mathcal{P}}\widehat{\mathcal{R}})$}{(BBB, PR)}}

Let us consider 
\begin{subequations}
\begin{align}
&\xi_0 m_b(\overline{\xi}) \widehat{X}_b(\overline{\xi}) \cdot \nabla_{\overline{\xi}} \widehat{I}_{\mu}^{(BBB, \widehat{\mathcal{P}}\widehat{\mathcal{R}})}[F_1, F_2](t, \overline{\xi}) \notag \\
&\quad = \int_0^t \int i s \xi_0 m_b(\overline{\xi}) \widehat{X}_b(\overline{\xi}) \cdot \nabla_{\overline{\xi}} \varphi e^{i s \varphi} \mu(\overline{\xi}, \overline{\eta}) \mu_{BBB}(\overline{\xi}, \overline{\eta}) m_{\widehat{\mathcal{P}}}(\overline{\eta}) m_{\widehat{\mathcal{R}}}(\overline{\sigma}) \widehat{F}_1(s, \overline{\eta}) \widehat{F}_2(s, \overline{\sigma}) ~ d\overline{\eta} ds \label{equdecchampbPR-1} \\
&\quad \quad + \int_0^t \int e^{i s \varphi} \xi_0 \mu(\overline{\xi}, \overline{\eta}) \mu_{BBB}(\overline{\xi}, \overline{\eta}) m_{\widehat{\mathcal{P}}}(\overline{\eta}) m_{\widehat{\mathcal{R}}}(\overline{\sigma}) \widehat{F}_1(s, \overline{\eta}) m_b(\overline{\xi}) \widehat{X}_b(\overline{\xi}) \cdot \nabla_{\overline{\xi}} \widehat{F}_2(s, \overline{\sigma}) ~ d\overline{\eta} ds \label{equdecchampbPR-2} \\
&\quad \quad + \int_0^t \int e^{i s \varphi} \xi_0 m_b(\overline{\xi}) \widehat{X}_b(\overline{\xi}) \cdot \nabla_{\overline{\xi}} \left( \mu(\overline{\xi}, \overline{\eta}) \mu_{BBB}(\overline{\xi}, \overline{\eta}) m_{\widehat{\mathcal{P}}}(\overline{\eta}) m_{\widehat{\mathcal{R}}}(\overline{\sigma}) \right) \widehat{F}_1(s, \overline{\eta}) \widehat{F}_2(s, \overline{\sigma}) ~ d\overline{\eta} ds \label{equdecchampbPR-3} 
\end{align}
\end{subequations}
Again, \eqref{equdecchampbPR-3} is already of the form \eqref{lemdecequH-11-dersymb}. Then, we need to get the term \eqref{lemdecequH-02-symsigma} from the decomposition, so we remove it from \eqref{equdecchampbPR-2} (with the same localisation): 
\begin{align*}
\eqref{equdecchampbPR-2} - \eqref{lemdecequH-02-symsigma} = - \int_0^t \int e^{i s \varphi} &\mu(\overline{\xi}, \overline{\eta}) \mu_{BBB}(\overline{\xi}, \overline{\eta}) m_{\widehat{\mathcal{P}}}(\overline{\eta}) m_{\widehat{\mathcal{R}}}(\overline{\sigma}) \widehat{F}_1(s, \overline{\eta}) \\
&\left( \xi_0 m_b(\overline{\xi}) \widehat{X}_b(\overline{\xi}) - \sigma_0 m_b(\overline{\sigma}) \widehat{X}_b(\overline{\sigma}) \right) \cdot \nabla_{\overline{\eta}} \widehat{F}_2(s, \overline{\sigma}) ~ d\overline{\eta} ds
\end{align*}
We then apply an integration by parts in $\overline{\eta}$ on this term: 
\begin{subequations}
\begin{align}
&\eqref{equdecchampbPR-2} - \eqref{lemdecequH-02-symsigma} \notag \\
&\begin{aligned}
= \int_0^t \int i s \left( \xi_0 m_b(\overline{\xi}) \widehat{X}_b(\overline{\xi}) - \sigma_0 m_b(\overline{\sigma}) \widehat{X}_b(\overline{\sigma}) \right) \cdot \nabla_{\overline{\eta}} \varphi e^{i s \varphi} \mu(\overline{\xi}, \overline{\eta}) \mu_{BBB}(\overline{\xi}, \overline{\eta}) m_{\widehat{\mathcal{P}}}(\overline{\eta}) \\
m_{\widehat{\mathcal{R}}}(\overline{\sigma}) \widehat{F}_1(s, \overline{\eta}) \widehat{F}_2(s, \overline{\sigma}) ~ d\overline{\eta} ds 
\end{aligned} \label{equdecchampbPR-2-1} \\
&\begin{aligned}
+ \int_0^t \int e^{i s \varphi} \mu(\overline{\xi}, \overline{\eta}) \mu_{BBB}(\overline{\xi}, \overline{\eta}) \left( \xi_0 m_b(\overline{\xi}) \widehat{X}_b(\overline{\xi}) - \sigma_0 m_b(\overline{\sigma}) \widehat{X}_b(\overline{\sigma}) \right) \cdot \nabla_{\overline{\eta}} \widehat{F}_1(s, \overline{\eta}) \\
m_{\widehat{\mathcal{P}}}(\overline{\eta}) m_{\widehat{\mathcal{R}}}(\overline{\sigma}) \widehat{F}_2(s, \overline{\sigma}) ~ d\overline{\eta} ds 
\end{aligned} \label{equdecchampbPR-2-2} \\
&\begin{aligned}
+ \int_0^t \int \nabla_{\overline{\eta}} \cdot \left( \left( \xi_0 m_b(\overline{\xi}) \widehat{X}_b(\overline{\xi}) - \sigma_0 m_b(\overline{\sigma}) \widehat{X}_b(\overline{\sigma}) \right) \mu(\overline{\xi}, \overline{\eta}) \mu_{BBB}(\overline{\xi}, \overline{\eta}) m_{\widehat{\mathcal{P}}}(\overline{\eta}) m_{\widehat{\mathcal{R}}}(\overline{\sigma}) \right) \\
e^{i s \varphi} \widehat{F}_1(s, \overline{\eta}) \widehat{F}_2(s, \overline{\sigma}) ~ d\overline{\eta} ds 
\end{aligned} \label{equdecchampbPR-2-3} 
\end{align}
\end{subequations}
\eqref{equdecchampbPR-2-3} is of the form \eqref{lemdecequH-11-dersymb}. 

In \eqref{equdecchampbPR-2-2}, we decompose by Lemma \ref{lem-projection-chpsvect}: 
\begin{align*}
&\xi_0 m_b(\overline{\xi}) \widehat{X}_b(\overline{\xi}) - \sigma_0 m_b(\overline{\sigma}) \widehat{X}_b(\overline{\sigma}) \\
&= O(\eta_0) + O(\widehat{X}_a(\overline{\eta})) + O(\widehat{X}_c(\overline{\eta})) + \xi_0 \left( m_b(\overline{\xi}) P_b^b(\overline{\xi}, \overline{\eta}) - m_b(\overline{\sigma}) P_b^b(\overline{\sigma}, \overline{\eta}) \right) \widehat{X}_b(\overline{\eta}) \\
&= O(\eta_0) + O(\widehat{X}_a(\overline{\eta})) + O(\widehat{X}_c(\overline{\eta})) + \xi_0 \left( m_b(\overline{\xi}) \frac{|\xi|}{|\overline{\xi}|} - m_b(\overline{\sigma}) \frac{|\sigma|}{|\overline{\sigma}|} \right) \frac{|\eta|}{|\overline{\eta}|} \widehat{X}_b(\overline{\eta}) \\
&= O(\eta_0) + O(\widehat{X}_a(\overline{\eta})) + O(\widehat{X}_c(\overline{\eta})) + O(|\xi| - |\sigma|) \widehat{X}_b(\overline{\eta})
\end{align*}
The term in $O(\eta_0)$ contributes as $\eqref{lemdecequH-03-resxetaac} + \eqref{lemdecequH-04-resxetabbon}$, the terms in $O(\widehat{X}_a(\overline{\eta}))$ or $O(\widehat{X}_c(\overline{\eta}))$ contribute as \eqref{lemdecequH-03-resxetaac}. Finally, for the last one, we use that 
\begin{align*}
\varphi &= O(\eta_0) + \sigma_0 (|\xi|^2 - |\sigma|^2) 
\end{align*}
hence 
\begin{align*}
|\xi| - |\sigma| = O(\eta_0) + O(\varphi)
\end{align*}
But here the term in $O(\varphi)$ contributes as \eqref{lemdecequH-06-resxetabphi}. In particular, we can rewrite \eqref{equdecchampbPR-2-2} entirely as terms from \eqref{lemdecequHGtot}. 

Finally, we group \eqref{equdecchampbPR-1} and \eqref{equdecchampbPR-2-1}: 
\begin{align*}
\eqref{equdecchampbPR-1} + \eqref{equdecchampbPR-2-1} = \int_0^t \int &i s \left( \xi_0 m_b(\overline{\xi}) \widehat{X}_b(\overline{\xi}) \cdot (\nabla_{\overline{\xi}} + \nabla_{\overline{\eta}}) - \sigma_0 m_b(\overline{\sigma}) \widehat{X}_b(\overline{\sigma}) \cdot \nabla_{\overline{\eta}} \right) \varphi e^{i s \varphi} \\
&\mu(\overline{\xi}, \overline{\eta}) \mu_{BBB}(\overline{\xi}, \overline{\eta}) m_{\widehat{\mathcal{P}}}(\overline{\eta}) m_{\widehat{\mathcal{R}}}(\overline{\sigma}) \widehat{F}_1(s, \overline{\eta}) \widehat{F}_2(s, \overline{\sigma}) ~ d\overline{\eta} ds
\end{align*}
We then compute: 
\begin{align*}
&\left( \xi_0 m_b(\overline{\xi}) \widehat{X}_b(\overline{\xi}) \cdot (\nabla_{\overline{\xi}} + \nabla_{\overline{\eta}}) - \sigma_0 m_b(\overline{\sigma}) \widehat{X}_b(\overline{\sigma}) \cdot \nabla_{\overline{\eta}} \right) \varphi \\
&\quad = O(\eta_0) + O(|\xi| - |\sigma|) \\
&\quad \quad \quad + m_b(\overline{\xi}) \left( \frac{\xi_0|\xi|}{|\overline{\xi}|} (3 \xi_0^2 + |\xi|^2 - |\eta|^2) - 2 \frac{\xi_0^3 |\xi|}{|\overline{\xi}|} - \frac{\sigma_0 |\sigma|}{|\overline{\sigma}|} (3 \sigma_0^2 + |\sigma|^2 - |\eta|^2) + 2 \frac{\sigma_0^3 |\sigma|}{|\overline{\sigma}|} \right) \\
&\quad = O(\eta_0) + O(|\xi| - |\sigma|) 
\end{align*}
On the other hand, 
\begin{align*}
\nabla_{\eta} \varphi &= 2 \sigma_0 \sigma - 2 \eta_0 \eta = 2 \sigma_0 \sigma + o(1) 
\end{align*}
so that $1 = O(\nabla_{\eta} \varphi)$: not only is \eqref{decompositionfinemg-gainkout} automatic, but also the factor $O(\eta_0)$ allows to apply an integration by parts in $\eta$ and only obtain terms of the form \eqref{lemdecequH-07-resxsigmaac}, \eqref{lemdecequH-08-resxsigmabbon}, \eqref{lemdecequH-03-resxetaac}, \eqref{lemdecequH-04-resxetabbon} or \eqref{lemdecequH-11-dersymb} (using that $\nabla_{\eta}$ decomposes into a combination of $\widehat{X}_a(\overline{\eta}), \eta_0 \widehat{X}_b(\overline{\eta}), \widehat{X}_c(\overline{\eta})$). 

We already saw that $\varphi = O(\eta_0) + \sigma_0 (|\xi|^2 - |\sigma|^2)$, so we can also apply integrations by parts on the term in $O(|\xi| - |\sigma|)$. 

\subsubsection{Interaction \texorpdfstring{$(BBB, \widehat{\mathcal{C}}\widehat{\mathcal{C}})$}{(BBB, CC)}} 

Let us consider
\begin{subequations}
\begin{align}
&\xi_0 m_b(\overline{\xi}) \widehat{X}_b(\overline{\xi}) \cdot \nabla_{\overline{\xi}} \widehat{I}_{\mu}^{(BBB, \widehat{\mathcal{C}}\widehat{\mathcal{C}})}[F_1, F_2](t, \overline{\xi}) \notag \\
&\quad = \int_0^t \int i s \xi_0 m_b(\overline{\xi}) \widehat{X}_b(\overline{\xi}) \cdot \nabla_{\overline{\xi}} \varphi e^{i s \varphi} \mu(\overline{\xi}, \overline{\eta}) \mu_{BBB}(\overline{\xi}, \overline{\eta}) m_{\widehat{\mathcal{C}}}(\overline{\eta}) m_{\widehat{\mathcal{C}}}(\overline{\sigma}) \widehat{F}_1(s, \overline{\eta}) \widehat{F}_2(s, \overline{\sigma}) ~ d\overline{\eta} ds \label{equdecchampbCC-1} \\
&\quad \quad + \int_0^t \int e^{i s \varphi} \xi_0 \mu(\overline{\xi}, \overline{\eta}) \mu_{BBB}(\overline{\xi}, \overline{\eta}) m_{\widehat{\mathcal{C}}}(\overline{\eta}) m_{\widehat{\mathcal{C}}}(\overline{\sigma}) \widehat{F}_1(s, \overline{\eta}) m_b(\overline{\xi}) \widehat{X}_b(\overline{\xi}) \cdot \nabla_{\overline{\xi}} \widehat{F}_2(s, \overline{\sigma}) ~ d\overline{\eta} ds \label{equdecchampbCC-2} \\
&\quad \quad + \int_0^t \int e^{i s \varphi} \xi_0 m_b(\overline{\xi}) \widehat{X}_b(\overline{\xi}) \cdot \nabla_{\overline{\xi}} \left( \mu(\overline{\xi}, \overline{\eta}) \mu_{BBB}(\overline{\xi}, \overline{\eta}) m_{\widehat{\mathcal{C}}}(\overline{\eta}) m_{\widehat{\mathcal{C}}}(\overline{\sigma}) \right) \widehat{F}_1(s, \overline{\eta}) \widehat{F}_2(s, \overline{\sigma}) ~ d\overline{\eta} ds \label{equdecchampbCC-3} 
\end{align}
\end{subequations}
\eqref{equdecchampbCC-3} is under the form \eqref{lemdecequH-11-dersymb}. 

We consider two subcases: either $\epsilon^{\overline{\eta} \overline{\sigma}} \theta^{\overline{\eta} \overline{\sigma}}$ is away enough from $1$ (which means that $\overline{\eta}$ and $\overline{\sigma}$ are not aligned, and in particular $\overline{\xi}$ cannot be too close to the cone); or $\epsilon^{\overline{\eta} \overline{\sigma}} \theta^{\overline{\eta} \overline{\sigma}}$ is close to $1$. We can obtain such a decomposition by adding Coifman-Meyer type symbols which, for simplicity, we will omit in what follows. 

\paragraph{1.} Localise first to have $\epsilon^{\overline{\eta} \overline{\sigma}} \theta^{\overline{\eta} \overline{\sigma}}$ away enough from $1$. 

To deal with \eqref{equdecchampbCC-2}, we decompose 
\[ \widehat{X}_b(\overline{\xi}) = \sum_{\alpha = a, b, c} P_b^{\alpha}(\overline{\xi}, \overline{\sigma}) \widehat{X}_{\alpha}(\overline{\sigma}) \]
by Lemma \ref{lem-projection-chpsvect}. If $\alpha = a$ or $c$, the contribution is of the form \eqref{lemdecequH-07-resxsigmaac} and therefore it only remains $\alpha = b$. By Lemma \ref{lemlienentrebaseancbasenouv}, we can also replace $\widehat{X}_b(\overline{\sigma})$ by $\frac{\sigma_0}{|\sigma_0|} \widehat{X}_b'(\overline{\sigma})$, and by Lemma \ref{lemchampmodifiebCC-propfond} we can then replace it by $\widehat{X}_{b-\widehat{\mathcal{C}}}(\overline{\eta}, \overline{\sigma})$, and finally apply an integration by parts in $\overline{\eta}$. Indeed, the modified vector field $\widehat{X}_{b-\widehat{\mathcal{C}}}(\overline{\eta}, \overline{\sigma})$ is not singular in this situation as $\epsilon^{\overline{\eta} \overline{\sigma}} \theta^{\overline{\eta} \overline{\sigma}}$ is not close to $1$. We then get (with a change of sign induced by $\nabla_{\overline{\xi}} \widehat{F}_2(s, \overline{\sigma}) = - \nabla_{\overline{\eta}} \widehat{F}_2(s, \overline{\sigma})$): 
\begin{subequations}
\begin{align}
&\begin{aligned}
-\int_0^t \int e^{i s \varphi} \xi_0 \mu(\overline{\xi}, \overline{\eta}) \mu_{BBB}(\overline{\xi}, \overline{\eta}) m_{\widehat{\mathcal{C}}}(\overline{\eta}) m_{\widehat{\mathcal{C}}}(\overline{\sigma}) \widehat{F}_1(s, \overline{\eta}) m_b(\overline{\xi}) P_b^b(\overline{\xi}, \overline{\sigma}) \frac{\sigma_0}{|\sigma_0|} \\
\widehat{X}_{b-\widehat{\mathcal{C}}}(\overline{\eta}, \overline{\sigma}) \cdot \nabla_{\overline{\eta}} \widehat{F}_2(s, \overline{\sigma}) ~ d\overline{\eta} ds 
\end{aligned} \notag \\
&\quad \begin{aligned}
= \int_0^t \int i s \xi_0 m_b(\overline{\xi}) P_b^b(\overline{\xi}, \overline{\sigma}) \frac{\sigma_0}{|\sigma_0|} \widehat{X}_{b-\widehat{\mathcal{C}}}(\overline{\eta}, \overline{\sigma}) \cdot \nabla_{\overline{\eta}} \varphi e^{i s \varphi} \mu(\overline{\xi}, \overline{\eta}) \mu_{BBB}(\overline{\xi}, \overline{\eta}) m_{\widehat{\mathcal{C}}}(\overline{\eta}) \\
m_{\widehat{\mathcal{C}}}(\overline{\sigma}) \widehat{F}_1(s, \overline{\eta}) \widehat{F}_2(s, \overline{\sigma}) ~ d\overline{\eta} ds
\end{aligned} \label{equdecchampbCC-2-1} \\
&\quad \quad \begin{aligned}
+ \int_0^t \int e^{i s \varphi} \xi_0 \mu(\overline{\xi}, \overline{\eta}) \mu_{BBB}(\overline{\xi}, \overline{\eta}) m_{\widehat{\mathcal{C}}}(\overline{\eta}) m_b(\overline{\xi}) P_b^b(\overline{\xi}, \overline{\sigma}) \widehat{X}_{b-\widehat{\mathcal{C}}}(\overline{\eta}, \overline{\sigma}) \cdot \nabla_{\overline{\eta}} \widehat{F}_1(s, \overline{\eta}) \\
m_{\widehat{\mathcal{C}}}(\overline{\sigma}) \frac{\sigma_0}{|\sigma_0|} \widehat{F}_2(s, \overline{\sigma}) ~ d\overline{\eta} ds 
\end{aligned} \label{equdecchampbCC-2-2} \\
&\quad \quad \begin{aligned}
+ \int_0^t \int \nabla_{\overline{\eta}} \cdot \left( m_b(\overline{\xi}) P_b^b(\overline{\xi}, \overline{\sigma}) \frac{\sigma_0}{|\sigma_0|} \widehat{X}_{b-\widehat{\mathcal{C}}}(\overline{\eta}, \overline{\sigma}) \xi_0 \mu(\overline{\xi}, \overline{\eta}) \mu_{BBB}(\overline{\xi}, \overline{\eta}) m_{\widehat{\mathcal{C}}}(\overline{\eta}) m_{\widehat{\mathcal{C}}}(\overline{\sigma}) \right) \\
e^{i s \varphi} \widehat{F}_1(s, \overline{\eta}) \widehat{F}_2(s, \overline{\sigma}) ~ d\overline{\eta} ds 
\end{aligned} \label{equdecchampbCC-2-3} 
\end{align}
\end{subequations}
\eqref{equdecchampbCC-2-3} is of the form \eqref{lemdecequH-11-dersymb}, and \eqref{equdecchampbCC-2-2} is of the form $\eqref{lemdecequH-03-resxetaac} + \eqref{lemdecequH-04-resxetabbon}$ by Lemma \ref{lemchampmodifiebCC-propfond}. 

Finally, we group \eqref{equdecchampbCC-1} and \eqref{equdecchampbCC-2-1}: 
\begin{align}
\eqref{equdecchampbCC-1} + \eqref{equdecchampbCC-2-1} &= \int_0^t \int i s \xi_0 m_b(\overline{\xi}) \left( \widehat{X}_b(\overline{\xi}) \cdot \nabla_{\overline{\xi}} + P_b^b(\overline{\xi}, \overline{\sigma}) \frac{\sigma_0}{|\sigma_0|} \widehat{X}_{b-\widehat{\mathcal{C}}}(\overline{\eta}, \overline{\sigma}) \cdot \nabla_{\overline{\eta}} \right) \varphi \notag \\
&e^{i s \varphi} \mu(\overline{\xi}, \overline{\eta}) \mu_{BBB}(\overline{\xi}, \overline{\eta}) m_{\widehat{\mathcal{C}}}(\overline{\eta}) m_{\widehat{\mathcal{C}}}(\overline{\sigma}) \widehat{F}_1(s, \overline{\eta}) \widehat{F}_2(s, \overline{\sigma}) ~ d\overline{\eta} ds \label{equdecchampbCC-1-1} 
\end{align}

By Lemma \ref{lemcalculsconecoordonneesconiquesvarphi}, 
\begin{align*}
\xi_t^{\overline{\eta} \overline{\sigma}} &= O\left( \widehat{X}_c(\overline{\eta}) \cdot \nabla_{\overline{\eta}} \varphi \right) \\
\widehat{X}_a(\overline{\eta}) \cdot \nabla_{\overline{\eta}} \varphi &= \frac{\eta_0}{|\overline{\eta}|} \left( \left(  \overline{\sigma}_a^{\overline{\eta}} \right)^2 
- \left( \overline{\eta}_a^{\overline{\eta}} \right)^2 \right) 
+ O\left( \xi_t^{\overline{\eta} \overline{\sigma}} \right) + O\left( \overline{\eta}_b^{\overline{\eta}} \right) 
\end{align*}
Since we localised to have $\epsilon^{\overline{\eta} \overline{\sigma}} \theta^{\overline{\eta} \overline{\sigma}}$ away enough from $1$, this forces 
\begin{align*}
1 + \epsilon^{\overline{\eta} \overline{\sigma}} \theta^{\overline{\eta} \overline{\sigma}} &= \frac{1 - \left( \theta^{\overline{\eta} \overline{\sigma}} \right)^2}{1 - \epsilon^{\overline{\eta} \overline{\sigma}} \theta^{\overline{\eta} \overline{\sigma}}} = O\left( \xi_t^{\overline{\eta} \overline{\sigma}} \right) 
\end{align*}
and therefore 
\begin{align*}
\overline{\sigma}_a^{\overline{\eta}} &= \sqrt{3} \epsilon^{\overline{\eta} \overline{\sigma}} |\sigma_0| + \theta^{\overline{\eta} \overline{\sigma}} |\sigma| = O\left( \xi_t^{\overline{\eta} \overline{\sigma}} \right) + O\left( \overline{\sigma}_b^{\overline{\sigma}} \right) 
\end{align*}
We finally get that 
\begin{align*}
\widehat{X}_a(\overline{\eta}) \cdot \nabla_{\overline{\eta}} \varphi &= - \frac{\eta_0}{|\overline{\eta}|} 
- \left( \overline{\eta}_a^{\overline{\eta}} \right)^2 
+ O\left( \xi_t^{\overline{\eta} \overline{\sigma}} \right) + O\left( \overline{\eta}_b^{\overline{\eta}} \right) + O\left( \overline{\sigma}_b^{\overline{\sigma}} \right) 
\end{align*}
In particular, assuming that $m_{\widehat{\mathcal{C}}}$ localises finely enough, 
\begin{align*}
1 &= O\left( \widehat{X}_a(\overline{\eta}) \cdot \nabla_{\overline{\eta}} \varphi \right) + O\left( \widehat{X}_c(\overline{\eta}) \cdot \nabla_{\overline{\eta}} \varphi \right) 
\end{align*}
This computations justifies that, in \eqref{equdecchampbCC-1-1}, we can apply integrations by parts and decompose as soon as we have a factor $O\left( \overline{\sigma}_b^{\overline{\sigma}} \right)$: the presence of $\overline{\sigma}_b^{\overline{\sigma}}$ allows to gain a symbol $m_b(\overline{\sigma})$ and hence to differentiate $\widehat{F}_2$ in any direction, while above the directions are compatible with $\widehat{F}_1$. Note also that \eqref{decompositionfinemg-gainkout} is automatic. 

In a symmetric way, we can control any contribution of the form $O\left( \overline{\eta}_b^{\overline{\eta}} \right)$. 

Finally, 
\begin{align*}
\varphi &= \xi_0^3 + \xi_0 |\xi|^2 - \eta_0^3 - \eta_0 |\eta|^2 - \sigma_0^3 - \sigma_0 |\sigma|^2 \\
&= 6 \xi_0 \eta_0 \sigma_0 \left( 1 + \epsilon^{\overline{\eta} \overline{\sigma}} \theta^{\overline{\eta} \overline{\sigma}} \right) + O\left( \overline{\eta}_b^{\overline{\eta}} \right) + O\left( \overline{\sigma}_b^{\overline{\sigma}} \right)
\end{align*}

Now we compute 
\begin{align*}
&m_b(\overline{\xi}) \left( \widehat{X}_b(\overline{\xi}) \cdot \nabla_{\overline{\xi}} + P_b^b(\overline{\xi}, \overline{\sigma}) \frac{\sigma_0}{|\sigma_0|} \widehat{X}_{b-\widehat{\mathcal{C}}}(\overline{\eta}, \overline{\sigma}) \cdot \nabla_{\overline{\eta}} \right) \varphi \xi_0 \\
&\quad = O\left( \overline{\eta}_b^{\overline{\eta}} \right) + O\left( \overline{\sigma}_b^{\overline{\sigma}} \right) + O\left( \xi_0 \epsilon^{\overline{\eta} \overline{\sigma}} \theta^{\overline{\eta} \overline{\sigma}} \right) \\
&\quad = O\left( m_b(\overline{\sigma}) \widehat{X}_a(\overline{\eta}) \cdot \nabla_{\overline{\eta}} \varphi \right) 
+ O\left( m_b(\overline{\sigma}) \widehat{X}_c(\overline{\eta}) \cdot \nabla_{\overline{\eta}} \varphi \right) 
+ O\left( m_b(\overline{\eta}) \widehat{X}_a(\overline{\sigma}) \cdot \nabla_{\overline{\eta}} \varphi \right) 
+ O\left( m_b(\overline{\eta}) \widehat{X}_c(\overline{\sigma}) \cdot \nabla_{\overline{\eta}} \varphi \right) 
+ O(\varphi) 
\end{align*}
This justifies that we can apply all necessary integrations by parts to bring \eqref{equdecchampbCC-1-1} back to terms of \eqref{lemdecequHGtot}. 

\paragraph{2.} Let us now localise $\epsilon^{\overline{\eta} \overline{\sigma}} \theta^{\overline{\eta} \overline{\sigma}}$ close to $1$. This means that $\overline{\eta}$ and $\overline{\sigma}$ are close to alignment: in particular, $\overline{\xi}$ is also localised near $\widehat{\mathcal{C}}$. 

We have here that 
\begin{align*}
\varphi &= \xi_0^3 + \xi_0 |\xi|^2 - \eta_0^3 - \eta_0 |\eta|^2 - \sigma_0^3 - \sigma_0 |\sigma|^2 \\
&= o(1) + 12 \xi_0 \eta_0 \sigma_0
\end{align*}
hence $1 = O(\varphi)$. We can directly apply an integration by parts in time to \eqref{equdecchampbCC-1} and get back to terms of \eqref{lemdecequHGtot}, with \eqref{decompositionfinemg-gainkout} being satisfied since we did not use the presence of $m_b(\overline{\xi})$. The only difficulty that remains now is to rewrite \eqref{equdecchampbCC-2}. To that end, we rewrite 
\begin{align*}
3 \xi_0^2 - |\xi|^2 &= 3 \eta_0^2 - |\eta|^2 + 3 \sigma_0^2 - |\sigma|^2 + 6 \eta_0 \sigma_0 - 2 \eta \cdot \sigma \\
&= O\left( \overline{\eta}_b^{\overline{\eta}} \right) + O\left( \overline{\sigma}_b^{\overline{\sigma}} \right) + 6 \eta_0 \sigma_0 \left( 1 - \epsilon^{\overline{\eta} \overline{\sigma}} \theta^{\overline{\eta} \overline{\sigma}} \right) 
\end{align*}
hence we deduce 
\begin{align*}
m_b(\overline{\xi}) = O(m_b(\overline{\sigma})) + O(m_b(\overline{\eta})) + O\left( 1 - \epsilon^{\overline{\eta} \overline{\sigma}} \theta^{\overline{\eta} \overline{\sigma}} \right) 
\end{align*}
Therefore, 
\begin{subequations}
\begin{align}
&\eqref{equdecchampbCC-2} \notag \\
&= \int_0^t \int e^{i s \varphi} \xi_0 \mu(\overline{\xi}, \overline{\eta}) \mu_{BBB}(\overline{\xi}, \overline{\eta}) m_{\widehat{\mathcal{C}}}(\overline{\eta}) m_{\widehat{\mathcal{C}}}(\overline{\sigma}) \widehat{F}_1(s, \overline{\eta}) O(m_b(\overline{\sigma})) \widehat{X}_b(\overline{\xi}) \cdot \nabla_{\overline{\xi}} \widehat{F}_2(s, \overline{\sigma}) ~ d\overline{\eta} ds \label{equdecchampbCC-2-bis-1} \\
&+ \int_0^t \int e^{i s \varphi} \xi_0 \mu(\overline{\xi}, \overline{\eta}) \mu_{BBB}(\overline{\xi}, \overline{\eta}) m_{\widehat{\mathcal{C}}}(\overline{\eta}) m_{\widehat{\mathcal{C}}}(\overline{\sigma}) \widehat{F}_1(s, \overline{\eta}) O(m_b(\overline{\eta})) \widehat{X}_b(\overline{\xi}) \cdot \nabla_{\overline{\xi}} \widehat{F}_2(s, \overline{\sigma}) ~ d\overline{\eta} ds \label{equdecchampbCC-2-bis-2} \\
&+ \int_0^t \int e^{i s \varphi} \xi_0 \mu(\overline{\xi}, \overline{\eta}) \mu_{BBB}(\overline{\xi}, \overline{\eta}) m_{\widehat{\mathcal{C}}}(\overline{\eta}) m_{\widehat{\mathcal{C}}}(\overline{\sigma}) \widehat{F}_1(s, \overline{\eta}) O\left( 1 - \epsilon^{\overline{\eta} \overline{\sigma}} \theta^{\overline{\eta} \overline{\sigma}} \right) \widehat{X}_b(\overline{\xi}) \cdot \nabla_{\overline{\xi}} \widehat{F}_2(s, \overline{\sigma}) ~ d\overline{\eta} ds \label{equdecchampbCC-2-bis-3} 
\end{align}
\end{subequations}
But \eqref{equdecchampbCC-2-bis-1} is already of the form $\eqref{lemdecequH-07-resxsigmaac} + \eqref{lemdecequH-08-resxsigmabbon}$, while we can apply an integration by parts on \eqref{equdecchampbCC-2-bis-2} to bring it back to \eqref{lemdecequH-03-resxetaac}, \eqref{lemdecequH-04-resxetabbon}, \eqref{lemdecequH-11-dersymb} or to a term similar to \eqref{equdecchampbCC-1} and treated the same way using $1 = O(\varphi)$, and we automatically have a factor $\nabla_{\overline{\eta}} \varphi$ so that \eqref{decompositionfinemg-gainkout} holds. 

Finally, for \eqref{equdecchampbCC-2-bis-3}, we use the fact that, even if $\widehat{X}_{b-\widehat{\mathcal{C}}}(\overline{\eta}, \overline{\sigma})$ is singular, 
\begin{align*}
\left( 1 - \epsilon^{\overline{\eta} \overline{\sigma}} \theta^{\overline{\eta} \overline{\sigma}} \right) \widehat{X}_{b-\widehat{\mathcal{C}}}(\overline{\eta}, \overline{\sigma})
\end{align*}
is not, so we can apply an integration by parts as in the previous case, and then use $1 = O(\varphi)$ and the presence of $\nabla_{\overline{\eta}} \varphi$ to deal with the term where the exponential is differentiated. Indeed: 
\begin{align*}
&\left( 1 - \epsilon^{\overline{\eta} \overline{\sigma}} \theta^{\overline{\eta} \overline{\sigma}} \right) \widehat{X}_{b-\widehat{\mathcal{C}}}(\overline{\eta}, \overline{\sigma}) \\
&= \left( 1 - \epsilon^{\overline{\eta} \overline{\sigma}} \theta^{\overline{\eta} \overline{\sigma}} \right) \widehat{X}_b'(\overline{\sigma}) - \left( 1 - \epsilon^{\overline{\eta} \overline{\sigma}} \theta^{\overline{\eta} \overline{\sigma}} \right) \frac{\widetilde{P}_b^b(\overline{\sigma}, \overline{\eta})}{\widetilde{P}_a^b(\overline{\eta}, \overline{\sigma})^2 + \widetilde{P}_c^b(\overline{\eta}, \overline{\sigma})^2} \left( \widetilde{P}_a^b(\overline{\sigma}, \overline{\eta}) \widehat{X}_a'(\overline{\sigma}) + \widetilde{P}_c^b(\overline{\sigma}, \overline{\eta}) \widehat{X}_c'(\overline{\eta}^n) \right)
\end{align*}
and 
\begin{align*}
\frac{1 - \epsilon^{\overline{\eta} \overline{\sigma}} \theta^{\overline{\eta} \overline{\sigma}}}{\widetilde{P}_a^b(\overline{\eta}, \overline{\sigma})^2 + \widetilde{P}_c^b(\overline{\eta}, \overline{\sigma})^2} &= \frac{1 - \epsilon^{\overline{\eta} \overline{\sigma}} \theta^{\overline{\eta} \overline{\sigma}}}{\left( \frac{\sqrt{3}}{4} \epsilon^{\overline{\eta} \overline{\sigma}} - \frac{\sqrt{3}}{4} \theta^{\overline{\eta} \overline{\sigma}} \right)^2 + \left( \frac{\sqrt{3}}{2} \xi_t^{\overline{\eta} \overline{\sigma}} \right)^2} 
\end{align*}
We then use that $1 - \epsilon^{\overline{\eta} \overline{\sigma}} \theta^{\overline{\eta} \overline{\sigma}} = O\left( \left( \xi_t^{\overline{\eta} \overline{\sigma}} \right)^2 \right)$ to conclude that this quantity is not singular. 

\subsubsection{Interaction \texorpdfstring{$(BBB, \widehat{\mathcal{L}}\widehat{\mathcal{C}})$}{(BBB, LC)}}

Let us consider 
\begin{subequations}
\begin{align}
&\xi_0 m_b(\overline{\xi}) \widehat{X}_b(\overline{\xi}) \cdot \nabla_{\overline{\xi}} \widehat{I}_{\mu}^{(BBB, \widehat{\mathcal{L}}\widehat{\mathcal{C}})}[F_1, F_2](t, \overline{\xi}) \notag \\
&\quad = \int_0^t \int i s \xi_0 m_b(\overline{\xi}) \widehat{X}_b(\overline{\xi}) \cdot \nabla_{\overline{\xi}} \varphi e^{i s \varphi} \mu(\overline{\xi}, \overline{\eta}) \mu_{BBB}(\overline{\xi}, \overline{\eta}) m_{\widehat{\mathcal{L}}}(\overline{\eta}) m_{\widehat{\mathcal{C}}}(\overline{\sigma}) \widehat{F}_1(s, \overline{\eta}) \widehat{F}_2(s, \overline{\sigma}) ~ d\overline{\eta} ds \label{equdecchampbLC-1} \\
&\quad \quad + \int_0^t \int e^{i s \varphi} \xi_0 \mu(\overline{\xi}, \overline{\eta}) \mu_{BBB}(\overline{\xi}, \overline{\eta}) m_{\widehat{\mathcal{L}}}(\overline{\eta}) m_{\widehat{\mathcal{C}}}(\overline{\sigma}) \widehat{F}_1(s, \overline{\eta}) m_b(\overline{\xi}) \widehat{X}_b(\overline{\xi}) \cdot \nabla_{\overline{\xi}} \widehat{F}_2(s, \overline{\sigma}) ~ d\overline{\eta} ds \label{equdecchampbLC-2} \\
&\quad \quad + \int_0^t \int e^{i s \varphi} \xi_0 m_b(\overline{\xi}) \widehat{X}_b(\overline{\xi}) \cdot \nabla_{\overline{\xi}} \left( \mu(\overline{\xi}, \overline{\eta}) \mu_{BBB}(\overline{\xi}, \overline{\eta}) m_{\widehat{\mathcal{L}}}(\overline{\eta}) m_{\widehat{\mathcal{C}}}(\overline{\sigma}) \right) \widehat{F}_1(s, \overline{\eta}) \widehat{F}_2(s, \overline{\sigma}) ~ d\overline{\eta} ds \label{equdecchampbLC-3} 
\end{align}
\end{subequations}
\eqref{equdecchampbLC-3} is of the form \eqref{lemdecequH-11-dersymb}. 

Then, we can project
\begin{align*}
\widehat{X}_b(\overline{\xi}) &= \sum_{\alpha = a, b, c} P_b^{\alpha}(\overline{\xi}, \overline{\sigma}) \widehat{X}_{\alpha}(\overline{\sigma}) 
\end{align*}
The contribution in \eqref{equdecchampbLC-2} of the terms $\alpha = a$ or $\alpha = c$ is of the form \eqref{lemdecequH-07-resxsigmaac}. Up to such terms, we may as well replace $\widehat{X}_b(\overline{\sigma})$ by $\frac{\sigma_0}{|\sigma_0|} \widehat{X}_b'(\overline{\sigma})$, and then by $\frac{\sigma_0}{|\sigma_0|} \widehat{X}_{b-\widehat{\mathcal{L}}}(\overline{\eta}, \overline{\sigma})$, before applying an integration by parts. 
\begin{subequations}
\begin{align}
&\begin{aligned}
\int_0^t \int e^{i s \varphi} \xi_0 \mu(\overline{\xi}, \overline{\eta}) \mu_{BBB}(\overline{\xi}, \overline{\eta}) m_{\widehat{\mathcal{L}}}(\overline{\eta}) m_{\widehat{\mathcal{C}}}(\overline{\sigma}) \widehat{F}_1(s, \overline{\eta}) m_b(\overline{\xi}) P_b^b(\overline{\xi}, \overline{\sigma}) \frac{\sigma_0}{|\sigma_0|} \\
\widehat{X}_{b-\widehat{\mathcal{L}}}(\overline{\eta}, \overline{\sigma}) \cdot \nabla_{\overline{\xi}} \widehat{F}_2(s, \overline{\sigma}) ~ d\overline{\eta} ds 
\end{aligned} \notag \\
&\quad \begin{aligned}
= \int_0^t \int i s \xi_0 m_b(\overline{\xi}) P_b^b(\overline{\xi}, \overline{\sigma}) \frac{\sigma_0}{|\sigma_0|} \widehat{X}_{b-\widehat{\mathcal{L}}}(\overline{\eta}, \overline{\sigma}) \cdot \nabla_{\overline{\eta}} \varphi e^{i s \varphi} \mu(\overline{\xi}, \overline{\eta}) \mu_{BBB}(\overline{\xi}, \overline{\eta}) m_{\widehat{\mathcal{L}}}(\overline{\eta}) \\
m_{\widehat{\mathcal{C}}}(\overline{\sigma}) \widehat{F}_1(s, \overline{\eta}) \widehat{F}_2(s, \overline{\sigma}) ~ d\overline{\eta} ds 
\end{aligned} \label{equdecchampbLC-2-1} \\
&\quad \begin{aligned}
+ \int_0^t \int e^{i s \varphi} \xi_0 \mu(\overline{\xi}, \overline{\eta}) \mu_{BBB}(\overline{\xi}, \overline{\eta}) m_{\widehat{\mathcal{L}}}(\overline{\eta}) m_{\widehat{\mathcal{C}}}(\overline{\sigma}) \widehat{X}_{b-\widehat{\mathcal{L}}}(\overline{\eta}, \overline{\sigma}) \cdot \nabla_{\overline{\eta}} \widehat{F}_1(s, \overline{\eta}) m_b(\overline{\xi}) \\
P_b^b(\overline{\xi}, \overline{\sigma}) \frac{\sigma_0}{|\sigma_0|} \widehat{F}_2(s, \overline{\sigma}) ~ d\overline{\eta} ds 
\end{aligned} \label{equdecchampbLC-2-2} \\
&\quad \begin{aligned}
+ \int_0^t \int \nabla_{\overline{\eta}} \cdot \left( m_b(\overline{\xi}) P_b^b(\overline{\xi}, \overline{\sigma}) \frac{\sigma_0}{|\sigma_0|} \widehat{X}_{b-\widehat{\mathcal{L}}}(\overline{\eta}, \overline{\sigma}) \xi_0 \mu(\overline{\xi}, \overline{\eta}) \mu_{BBB}(\overline{\xi}, \overline{\eta}) m_{\widehat{\mathcal{L}}}(\overline{\eta}) m_{\widehat{\mathcal{C}}}(\overline{\sigma}) \right) \\
e^{i s \varphi} \widehat{F}_1(s, \overline{\eta}) \widehat{F}_2(s, \overline{\sigma}) ~ d\overline{\eta} ds 
\end{aligned} \label{equdecchampbLC-2-3} 
\end{align}
\end{subequations}
Lemma \ref{lemchampmodifiebLCpropfond} ensures that \eqref{equdecchampbLC-2-2} is of the form $\eqref{lemdecequH-03-resxetaac} + \eqref{lemdecequH-04-resxetabbon}$, and \eqref{equdecchampbLC-2-3} is of the form \eqref{lemdecequH-11-dersymb}. Then, we group \eqref{equdecchampbLC-2-1} with \eqref{equdecchampbLC-1}: 
\begin{align}
\eqref{equdecchampbLC-2-1} + \eqref{equdecchampbLC-1} = \int_0^t \int &i s \xi_0 m_b(\overline{\xi}) \left( \widehat{X}_b(\overline{\xi}) \cdot \nabla_{\overline{\xi}} + P_b^b(\overline{\xi}, \overline{\sigma}) \frac{\sigma_0}{|\sigma_0|} \widehat{X}_{b-\widehat{\mathcal{L}}}(\overline{\eta}, \overline{\sigma}) \cdot \nabla_{\overline{\eta}} \right) \varphi e^{i s \varphi} \notag \\
&\mu(\overline{\xi}, \overline{\eta}) \mu_{BBB}(\overline{\xi}, \overline{\eta}) m_{\widehat{\mathcal{L}}}(\overline{\eta}) m_{\widehat{\mathcal{C}}}(\overline{\sigma}) \widehat{F}_1(s, \overline{\eta}) \widehat{F}_2(s, \overline{\sigma}) ~ d\overline{\eta} ds \label{equdecchampbLC-1-1} 
\end{align}

Note that, here, $\nabla_{\eta} \varphi = o(1) + \sigma_0 \sigma$ so that $1 = O(\nabla_{\overline{\eta}} \varphi)$ and \eqref{decompositionfinemg-gainkout} is automatic. 

\begin{figure}
\centering
\begin{subfigure}{0.4\textwidth}
\centering
\begin{tikzpicture}
        \draw[->] (-0.5, 0) -- (2.3, 0) node[above] {$\xi$};
        \draw[->] (0, -1.5) -- (0, 1.5) node[right] {$\xi_0$};
        \draw[->] (0, 0) -- ({0.97*sqrt(3)}, 0.97) node[above] {$\overline{\sigma}$};
        \draw[->] ({sqrt(3)}, 1) -- ({sqrt(3)}, -0.95) node[right] {$\overline{\eta}$};
        \draw[->] (0, 0) -- ({sqrt(3)}, -1) node[below] {$\overline{\xi}$};
        \draw[dotted] ({-0.25*sqrt(3)}, -0.25) -- ({1.25*sqrt(3)}, 1.25) node[right] {$\widehat{\mathcal{C}}$};
        \draw[dotted] ({-0.25*sqrt(3)}, 0.25) -- ({1.25*sqrt(3)}, -1.25);
    \end{tikzpicture}
    \subcaption{Remaining case for the $(BBB, \widehat{\mathcal{L}}\widehat{\mathcal{C}})$ interaction} \label{figurecasrestantBBBLC}
\end{subfigure} \hfill
\begin{subfigure}{0.4\textwidth}
\centering
\begin{tikzpicture}
		\draw[->] (-2.0, 0) -- (2., 0) node[above] {$\xi$};
		\draw[->] (0, -0.5) -- (0, 1.5) node[right] {$\xi_0$};
        \draw[->] (0, 0) -- ({0.97*sqrt(3)}, 0.97) node[below] {$\overline{\sigma}$};
        \draw[->] ({sqrt(3)}, 1) -- (-{0.95*sqrt(3)}, 1) node[above] {$\overline{\eta}$};
        \draw[->] (0, 0) -- (-{sqrt(3)}, 1) node[below] {$\overline{\xi}$};
        \draw[dotted] ({-0.25*sqrt(3)}, -0.25) -- ({1.25*sqrt(3)}, 1.25) node[right] {$\widehat{\mathcal{C}}$};
        \draw[dotted] ({0.25*sqrt(3)}, -0.25) -- ({-1.25*sqrt(3)}, 1.25);
\end{tikzpicture}
\subcaption{Remaining case for the $(BBB, \widehat{\mathcal{P}}\widehat{\mathcal{C}})$ interaction} \label{figurecasrestantBBBPC} 
\end{subfigure} 
\caption{Remaining cases of $BBB$ interactions} \label{figurecasrestantBBB} 
\end{figure}
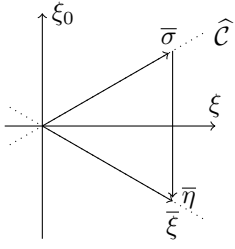
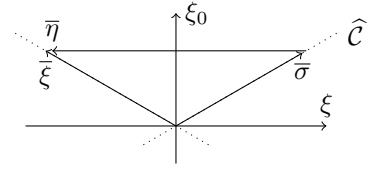

We now compute that 
\begin{align*}
\varphi &= \xi_0^3 + \xi_0 |\xi|^2 - \eta_0^3 - \eta_0 |\eta|^2 - \sigma_0^3 - \sigma_0 |\sigma|^2 \\
&= o(1) + 3 \sigma_0 \eta_0 (\xi_0 + \sigma_0) 
\end{align*}
In particular, as soon as $|\xi_0 + \sigma_0| \gtrsim |\overline{\xi}|$, we have $1 = O(\varphi)$ and we may thus control \eqref{equdecchampbLC-1-1} easily by integration by parts in time. We thus restrict our attention to the case $\xi_0 + \sigma_0 = o(1)$. In particular, $\overline{\xi}$ is close to $\widehat{\mathcal{C}}$, and $\eta_0 + 2 \sigma_0 = o(1)$, see Figure \ref{figurecasrestantBBBLC}. 

We now compute that
\begin{align*}
\partial_{\eta_0} \varphi &= 3 \sigma_0^2 + |\sigma|^2 - 3 \eta_0^2 - |\eta|^2 \\
&= o(1) - 6 \sigma_0^2 
\end{align*}
so $1 = O(\partial_{\eta_0} \varphi)$, hence
\[ \overline{\sigma}_b^{\overline{\sigma}} = O\left( m_b(\overline{\sigma}) m_a(\overline{\eta}) \widehat{X}_a(\overline{\eta}) \cdot \nabla_{\overline{\eta}} \varphi \right) \]
On the other hand, by Lemma \ref{lemcalculsconecoordonneesconiquesvarphi}, 
\begin{align*}
\xi_t^{\overline{\xi} \overline{\sigma}} = \frac{|\eta|}{|\xi|} \xi_t^{\overline{\eta} \overline{\sigma}} &= O\left( \widehat{X}_c(\overline{\sigma}) \cdot \nabla_{\overline{\eta}} \varphi \right) 
\end{align*}
In particular, since $\xi_t^{\overline{\xi} \overline{\sigma}} = O\left( |\eta| \right)$ by the previous computation, we deduce that 
\begin{align*}
\left( \xi_t^{\overline{\xi} \overline{\sigma}} \right)^2 &= O\left( |\eta| \widehat{X}_c(\overline{\sigma}) \cdot \nabla_{\overline{\eta}} \varphi \right)
\end{align*}
We may express the same way $|\eta|^2 \left( \xi_t^{\overline{\eta} \overline{\xi}} \right)^2$. Moreover, here, $\epsilon^{\overline{\xi} \overline{\sigma}} = -1$ and $\theta^{\overline{\xi} \overline{\sigma}}$ is close to $1$, hence 
\begin{align*}
1 + \epsilon^{\overline{\xi} \overline{\sigma}} \theta^{\overline{\xi} \overline{\sigma}} &= O\left( \left( \xi_t^{\overline{\xi} \overline{\sigma}} \right)^2 \right) 
\end{align*}

Then, 
\begin{align*}
\overline{\sigma}_a^{\overline{\xi}} &= \sqrt{3} \epsilon^{\overline{\xi} \overline{\sigma}} |\sigma_0| + \theta^{\overline{\xi} \overline{\sigma}} |\sigma| \\
&= O\left( \left( \xi_t^{\overline{\xi} \overline{\sigma}} \right)^2 \right) + O\left( \overline{\sigma}_b^{\overline{\sigma}} \right) 
\end{align*}

Finally, by Lemma \ref{lemcalculsconecoordonneesconiquesvarphi}, and using also Lemma \ref{lemcalculvarphiKdV1D}, we may simplify 
\begin{align*}
6 \sqrt{3} \frac{\xi_0}{|\xi_0|} \varphi &= \left( \overline{\xi}_a^{\overline{\xi}} \right)^3 + \left( \overline{\xi}_b^{\overline{\xi}} \right)^3 - \left( \overline{\eta}_a^{\overline{\xi}} \right)^3 - \left( \overline{\eta}_b^{\overline{\xi}} \right)^3 - \left( \overline{\sigma}_a^{\overline{\xi}} \right)^3 - \left( \overline{\sigma}_b^{\overline{\xi}} \right)^3 + O\left( \left( \xi_t^{\overline{\xi} \overline{\sigma}} \right)^2 \right) \\
&= 3 \overline{\xi}_b^{\overline{\xi}} \overline{\eta}_b^{\overline{\xi}} \overline{\sigma}_b^{\overline{\xi}} + O\left( \left( \xi_t^{\overline{\xi} \overline{\sigma}} \right)^2 \right) + O\left( \overline{\sigma}_b^{\overline{\sigma}} \right) 
\end{align*}
But $\overline{\eta}_b^{\overline{\xi}} \simeq |\eta_0|$, $\overline{\sigma}_b^{\overline{\xi}} \simeq |\overline{\sigma}|$, so
\begin{align*}
\overline{\xi}_b^{\overline{\xi}} &= O(\varphi) + O\left( \left( \xi_t^{\overline{\xi} \overline{\sigma}} \right)^2 \right) + O\left( \overline{\sigma}_b^{\overline{\sigma}} \right)
\end{align*}
which is enough to apply integrations by parts and decompose \eqref{equdecchampbLC-1-1}. 

\subsubsection{Interaction \texorpdfstring{$(BBB, \widehat{\mathcal{P}}\widehat{\mathcal{C}})$}{(BBB, PC)}}

Let us consider 
\begin{subequations}
\begin{align}
&\xi_0 m_b(\overline{\xi}) \widehat{X}_b(\overline{\xi}) \cdot \nabla_{\overline{\xi}} \widehat{I}_{\mu}^{(BBB, \widehat{\mathcal{P}}\widehat{\mathcal{C}})}[F_1, F_2](t, \overline{\xi}) \notag \\
&\quad = \int_0^t \int i s \xi_0 m_b(\overline{\xi}) \widehat{X}_b(\overline{\xi}) \cdot \nabla_{\overline{\xi}} \varphi e^{i s \varphi} \mu(\overline{\xi}, \overline{\eta}) \mu_{BBB}(\overline{\xi}, \overline{\eta}) m_{\widehat{\mathcal{P}}}(\overline{\eta}) m_{\widehat{\mathcal{C}}}(\overline{\sigma}) \widehat{F}_1(s, \overline{\eta}) \widehat{F}_2(s, \overline{\sigma}) ~ d\overline{\eta} ds \label{equdecchampbPC-1} \\
&\quad \quad + \int_0^t \int e^{i s \varphi} \xi_0 \mu(\overline{\xi}, \overline{\eta}) \mu_{BBB}(\overline{\xi}, \overline{\eta}) m_{\widehat{\mathcal{P}}}(\overline{\eta}) m_{\widehat{\mathcal{C}}}(\overline{\sigma}) \widehat{F}_1(s, \overline{\eta}) m_b(\overline{\xi}) \widehat{X}_b(\overline{\xi}) \cdot \nabla_{\overline{\xi}} \widehat{F}_2(s, \overline{\sigma}) ~ d\overline{\eta} ds \label{equdecchampbPC-2} \\
&\quad \quad + \int_0^t \int e^{i s \varphi} \xi_0 m_b(\overline{\xi}) \widehat{X}_b(\overline{\xi}) \cdot \nabla_{\overline{\xi}} \left( \mu(\overline{\xi}, \overline{\eta}) \mu_{BBB}(\overline{\xi}, \overline{\eta}) m_{\widehat{\mathcal{P}}}(\overline{\eta}) m_{\widehat{\mathcal{C}}}(\overline{\sigma}) \right) \widehat{F}_1(s, \overline{\eta}) \widehat{F}_2(s, \overline{\sigma}) ~ d\overline{\eta} ds \label{equdecchampbPC-3} 
\end{align}
\end{subequations}
\eqref{equdecchampbPC-3} is of the form \eqref{lemdecequH-11-dersymb}. 

We remove from \eqref{equdecchampbPC-2} the term \eqref{lemdecequH-02-symsigma}, and then apply an integration by parts: 
\begin{subequations}
\begin{align}
&\eqref{equdecchampbPC-2} - \eqref{lemdecequH-02-symsigma} \notag \\
&\begin{aligned}
= 
\int_0^t \int e^{i s \varphi} \mu(\overline{\xi}, \overline{\eta}) \mu_{BBB}(\overline{\xi}, \overline{\eta}) m_{\widehat{\mathcal{C}}}(\overline{\sigma}) \left( \xi_0 m_b(\overline{\xi}) \widehat{X}_b(\overline{\xi}) - \sigma_0 m_b(\overline{\sigma}) \widehat{X}_b(\overline{\sigma}) \right) \cdot \nabla_{\overline{\xi}} \widehat{F}_2(s, \overline{\sigma}) \\
m_{\widehat{\mathcal{P}}}(\overline{\eta}) \widehat{F}_1(s, \overline{\eta}) ~ d\overline{\eta} ds 
\end{aligned} \notag \\
&\begin{aligned}
= \int_0^t \int i s \left( \xi_0 m_b(\overline{\xi}) \widehat{X}_b(\overline{\xi}) - \sigma_0 m_b(\overline{\sigma}) \widehat{X}_b(\overline{\sigma}) \right) \cdot \nabla_{\overline{\eta}} \varphi e^{i s \varphi} \mu(\overline{\xi}, \overline{\eta}) \mu_{BBB}(\overline{\xi}, \overline{\eta}) m_{\widehat{\mathcal{P}}}(\overline{\eta}) \\
m_{\widehat{\mathcal{C}}}(\overline{\sigma}) \widehat{F}_1(s, \overline{\eta}) \widehat{F}_2(s, \overline{\sigma}) ~ d\overline{\eta} ds 
\end{aligned} \label{equdecchampbPC-2-1} \\
&\begin{aligned}
+ \int_0^t \int e^{i s \varphi} \mu(\overline{\xi}, \overline{\eta}) \mu_{BBB}(\overline{\xi}, \overline{\eta})  \left( \xi_0 m_b(\overline{\xi}) \widehat{X}_b(\overline{\xi}) - \sigma_0 m_b(\overline{\sigma}) \widehat{X}_b(\overline{\sigma}) \right) \cdot \nabla_{\overline{\eta}} \widehat{F}_1(s, \overline{\eta}) \\
m_{\widehat{\mathcal{P}}}(\overline{\eta})m_{\widehat{\mathcal{C}}}(\overline{\sigma}) \widehat{F}_2(s, \overline{\sigma}) ~ d\overline{\eta} ds 
\end{aligned} \label{equdecchampbPC-2-2} \\
&\begin{aligned}
+ \int_0^t \int \nabla_{\overline{\eta}} \cdot \left( \left( \xi_0 m_b(\overline{\xi}) \widehat{X}_b(\overline{\xi}) - \sigma_0 m_b(\overline{\sigma}) \widehat{X}_b(\overline{\sigma}) \right) \mu(\overline{\xi}, \overline{\eta}) \mu_{BBB}(\overline{\xi}, \overline{\eta}) m_{\widehat{\mathcal{P}}}(\overline{\eta}) m_{\widehat{\mathcal{C}}}(\overline{\sigma}) \right) \\
e^{i s \varphi} \widehat{F}_1(s, \overline{\eta}) \widehat{F}_2(s, \overline{\sigma}) ~ d\overline{\eta} ds 
\end{aligned} \label{equdecchampbPC-2-3} 
\end{align}
\end{subequations}
\eqref{equdecchampbPC-2-3} is of the form \eqref{lemdecequH-11-dersymb}. 

For \eqref{equdecchampbPC-2-2}, we project: 
\begin{align*}
\xi_0 m_b(\overline{\xi}) \widehat{X}_b(\overline{\xi}) - \sigma_0 m_b(\overline{\sigma}) \widehat{X}_b(\overline{\sigma}) 
&= O(\eta_0) + O(\widehat{X}_a(\overline{\eta})) + O(\widehat{X}_c(\overline{\eta})) + O(|\xi| - |\sigma|) \widehat{X}_b(\overline{\eta}) 
\end{align*}
by the same computation as in the $(BBB, \widehat{\mathcal{P}}\widehat{\mathcal{R}})$ interaction. But we have again that 
\begin{align*}
\varphi &= O(\eta_0) + \sigma_0 (|\xi|^2 - |\sigma|^2) 
\end{align*}
so that, here, $|\xi| - |\sigma| = O(\eta_0) + O(\varphi)$. That way, we can bring \eqref{equdecchampbPC-2-2} back to terms of \eqref{lemdecequHGtot}. 

Then, we group \eqref{equdecchampbPC-2-1} and \eqref{equdecchampbPC-1}: 
\begin{align}
\eqref{equdecchampbPC-1} + \eqref{equdecchampbPC-2-1} = 
&\int_0^t \int i s \left( \xi_0 m_b(\overline{\xi}) \widehat{X}_b(\overline{\xi}) \cdot \nabla_{\overline{\xi}} + \left( \xi_0 m_b(\overline{\xi}) \widehat{X}_b(\overline{\xi}) - \sigma_0 m_b(\overline{\sigma}) \widehat{X}_b(\overline{\sigma}) \right) \cdot \nabla_{\overline{\eta}} \right) \varphi \notag \\
&e^{i s \varphi} \mu(\overline{\xi}, \overline{\eta}) \mu_{BBB}(\overline{\xi}, \overline{\eta}) m_{\widehat{\mathcal{P}}}(\overline{\eta}) m_{\widehat{\mathcal{C}}}(\overline{\sigma}) \widehat{F}_1(s, \overline{\eta}) \cdot \nabla_{\overline{\xi}} \widehat{F}_2(s, \overline{\sigma}) ~ d\overline{\eta} ds \label{equdecchampbPC-1-1} 
\end{align}

We have that 
\begin{align*}
\widehat{X}_a(\overline{\sigma}) \cdot \nabla_{\overline{\eta}} \varphi &= \frac{\sigma_0}{|\overline{\sigma}|} \left( 3 \sigma_0^2 + |\sigma|^2 - 3 \eta_0^2 - |\eta|^2 \right) + \frac{\sigma}{|\overline{\sigma}|} \cdot (2 \sigma_0 \sigma - 2 \eta_0 \eta) \\
&= o(1) + \frac{\sigma_0}{|\overline{\sigma}|} \left( 4 |\sigma|^2 - |\eta|^2 \right) \\
\varphi &= \xi_0^3 + \xi_0 |\xi|^2 - \eta_0^3 - \eta_0 |\eta|^2 - \sigma_0^3 - \sigma_0 |\sigma|^2 \\
&= o(1) + \xi_0 (|\xi|^2 - |\sigma|^2) 
\end{align*}
In particular, as soon as $|\xi| - |\sigma| \gtrsim |\overline{\xi}|$, we have $1 = O(\varphi)$ and there is no difficulty to decompose, as again $\nabla_{\eta} \varphi \simeq \sigma_0 \sigma$ so $1 = O(\nabla_{\overline{\eta}} \varphi)$ and \eqref{decompositionfinemg-gainkout} is automatic. 

If $|\xi| - |\sigma| \ll |\overline{\xi}|$, then $\overline{\xi}$ is close to the cone. We separate into two cases. 

\paragraph{1.} If $1 = O\left( 4 |\sigma|^2 - |\eta|^2 \right) = O\left( \widehat{X}_a(\overline{\sigma}) \cdot \nabla_{\overline{\eta}} \varphi \right)$, then we can apply an integration by parts on any term having a factor $O(\eta_0)$ using $\widehat{X}_a(\overline{\sigma})$ and get terms from \eqref{lemdecequHGtot}. 

We then have that 
\begin{align*}
\varphi &= O(\eta_0) + \xi_0 (|\xi|^2 - |\sigma|^2) 
\end{align*}
On the other hand, 
\begin{align*}
|\overline{\xi}| - |\overline{\sigma}| &= \frac{|\overline{\xi}|^2 - |\overline{\sigma}|^2}{|\overline{\xi}| + |\overline{\sigma}|} \\
&= O(\eta_0) + O(|\xi|^2 - |\sigma|^2) \\
(3 \xi_0^2 - |\xi|^2) - (3 \sigma_0^2 - |\sigma|^2) &= O(\eta_0) + O(|\xi|^2 - |\sigma|^2) 
\end{align*}
We deduce that 
\begin{align*}
m_b(\overline{\xi}) = m_b(\overline{\sigma}) + O(\eta_0) + O(\varphi) 
\end{align*}
Then, we can also compute that 
\begin{align*}
&\left( \xi_0 m_b(\overline{\xi}) \widehat{X}_b(\overline{\xi}) \cdot \nabla_{\overline{\xi}} +  \left( \xi_0 m_b(\overline{\xi}) \widehat{X}_b(\overline{\xi}) - \sigma_0 m_b(\overline{\sigma}) \widehat{X}_b(\overline{\sigma}) \right) \cdot \nabla_{\overline{\eta}} \right) \varphi \\
&\quad = O(\eta_0) + O(\varphi)
\end{align*}
this justifies that we can decompose \eqref{equdecchampbPC-1-1}. 

\paragraph{2.} Let us now consider the case $4 |\sigma|^2 - |\eta|^2 = o(1)$. By the triangular inequality, since $|\xi| - |\sigma| = o(1)$, we must have $2 \sigma + \eta = \xi + \sigma = o(1)$, $\theta^{\overline{\eta} \overline{\sigma}}$ close to $-1$, see Figure \ref{figurecasrestantBBBPC}. 

We have that
\begin{align*}
\nabla_{\eta} \varphi &= 2 \sigma_0 \sigma - 2 \eta_0 \eta = o(1) + 2 \sigma_0 \sigma 
\end{align*}
so $1 = O(\nabla_{\eta} \varphi)$, and we may apply integrations by parts in $\eta$ whenever we have a factor $O\left( \eta_0 \overline{\sigma}_b^{\overline{\sigma}} \right)$ to get terms of \eqref{lemdecequHGtot}. 

By Lemma \ref{lemcalculsconecoordonneesconiquesvarphi}, 
\begin{align*}
\eta_0^2 \xi_t^{\overline{\eta} \overline{\sigma}} &= O\left( \eta_0 \widehat{X}_c(\overline{\sigma}) \cdot \nabla_{\overline{\eta}} \varphi \right) 
\end{align*}
and we can compute that
\begin{align*}
\varphi &= 3 \xi_0 \eta_0 \sigma_0 + \xi_0 (|\xi|^2 - |\sigma|^2) - \eta_0 |\eta|^2 + \eta_0 |\sigma|^2 \\
&= O\left( \eta_0 \overline{\sigma}_b^{\overline{\sigma}} \right) + \xi_0 (|\xi|^2 - |\sigma|^2) + \eta_0 \left( 3 (\xi_0 + \sigma_0) \sigma_0 - |\eta|^2 \right) 
\end{align*}
which implies 
\begin{align*}
\xi_0 (|\xi|^2 - |\sigma|^2) &= O(\varphi) + O\left( \eta_0 \overline{\sigma}_b^{\overline{\sigma}} \right) + \eta_0 \left( |\eta|^2 - 3 (\xi_0 + \sigma_0) \sigma_0 \right) 
\end{align*}
In particular, 
\begin{align*}
\overline{\sigma}_b^{\overline{\sigma}} (|\xi|^2 - |\sigma|^2) &= O(\varphi) + O\left( \eta_0 \overline{\sigma}_b^{\overline{\sigma}} \right)
\end{align*}
Finally, we can also compute that 
\begin{align*}
&\eta_0 \widehat{X}_a(\overline{\sigma}) \cdot \nabla_{\overline{\eta}} \varphi \\
&= O\left( \eta_0 \overline{\sigma}_b^{\overline{\sigma}} \right) + O\left( \eta_0^2 \xi_t^{\overline{\eta} \overline{\sigma}} \right) + \eta_0 \frac{\sigma_0}{|\overline{\sigma}|} \left( 2 \sqrt{3} |\sigma_0| + |\eta| - \sqrt{3} \eta_0 \frac{\sigma_0}{|\sigma_0|} \right) \left( \sqrt{3} |\sigma_0| + \sqrt{3} |\xi_0| - |\eta| \right) \\
\end{align*}
We deduce that we can also decompose in the presence of a factor $\eta_0 \left( \sqrt{3} |\sigma_0| + \sqrt{3} |\xi_0| - |\eta| \right)$. 

Note now that 
\begin{align*}
3 \xi_0^2 - |\xi|^2 - 3 \sigma_0^2 + |\sigma|^2 &= 3 \eta_0 (\xi_0 + \sigma_0) + O(\varphi) + O\left( \eta_0 \overline{\sigma}_b^{\overline{\sigma}} \right) + \frac{\eta_0}{\xi_0} \left( 3 (\xi_0 + \sigma_0) \sigma_0 - |\eta|^2 \right) \\
&= O(\varphi) + O\left( \eta_0 \overline{\sigma}_b^{\overline{\sigma}} \right) + \frac{\eta_0}{\xi_0} \left( 3 (\xi_0 + \sigma_0)^2 - |\eta|^2 \right) \\
&= O(\varphi) + O\left( \eta_0 \overline{\sigma}_b^{\overline{\sigma}} \right) + O\left(\eta_0 (\sqrt{3} |\sigma_0| + \sqrt{3} |\xi_0| - |\eta|) \right) 
\end{align*}

It only remains to compute
\begin{align*}
&m_b(\overline{\xi}) \xi_0 \widehat{X}_b(\overline{\xi}) \cdot (\nabla_{\overline{\xi}} + \nabla_{\overline{\eta}}) \varphi - m_b(\overline{\sigma}) \sigma_0 \widehat{X}_b(\overline{\sigma}) \cdot \nabla_{\overline{\eta}} \varphi \\
&\quad = O\left( \eta_0 \overline{\sigma}_b^{\overline{\sigma}} \right) + O\left( \overline{\sigma}_b^{\overline{\sigma}} (|\xi| - |\sigma|) \right) + O(3 \xi_0^2 - |\xi|^2 - 3 \sigma_0^2 + |\sigma|^2) 
\end{align*} 
which justifies that we can apply integrations by parts and get the desired decomposition. 

\subsubsection{Interaction \texorpdfstring{$(BBB, \widehat{\mathcal{L}}\widehat{\mathcal{L}})$}{(BBB, LL)}}

Let us consider
\begin{subequations}
\begin{align}
&\xi_0 m_b(\overline{\xi}) \widehat{X}_b(\overline{\xi}) \cdot \nabla_{\overline{\xi}} \widehat{I}_{\mu}^{(BBB, \widehat{\mathcal{L}}\widehat{\mathcal{L}})}[F_1, F_2](t, \overline{\xi}) \notag \\
&\quad = \int_0^t \int i s \xi_0 m_b(\overline{\xi}) \widehat{X}_b(\overline{\xi}) \cdot \nabla_{\overline{\xi}} \varphi e^{i s \varphi} \mu(\overline{\xi}, \overline{\eta}) \mu_{BBB}(\overline{\xi}, \overline{\eta}) m_{\widehat{\mathcal{L}}}(\overline{\eta}) m_{\widehat{\mathcal{L}}}(\overline{\sigma}) \widehat{F}_1(s, \overline{\eta}) \widehat{F}_2(s, \overline{\sigma}) ~ d\overline{\eta} ds \label{equdecchampbLL-1} \\
&\quad \quad + \int_0^t \int e^{i s \varphi} \xi_0 \mu(\overline{\xi}, \overline{\eta}) \mu_{BBB}(\overline{\xi}, \overline{\eta}) m_{\widehat{\mathcal{L}}}(\overline{\eta}) m_{\widehat{\mathcal{L}}}(\overline{\sigma}) \widehat{F}_1(s, \overline{\eta}) m_b(\overline{\xi}) \widehat{X}_b(\overline{\xi}) \cdot \nabla_{\overline{\xi}} \widehat{F}_2(s, \overline{\sigma}) ~ d\overline{\eta} ds \label{equdecchampbLL-2} \\
&\quad \quad + \int_0^t \int e^{i s \varphi} \xi_0 m_b(\overline{\xi}) \widehat{X}_b(\overline{\xi}) \cdot \nabla_{\overline{\xi}} \left( \mu(\overline{\xi}, \overline{\eta}) \mu_{BBB}(\overline{\xi}, \overline{\eta}) m_{\widehat{\mathcal{L}}}(\overline{\eta}) m_{\widehat{\mathcal{L}}}(\overline{\sigma}) \right) \widehat{F}_1(s, \overline{\eta}) \widehat{F}_2(s, \overline{\sigma}) ~ d\overline{\eta} ds \label{equdecchampbLL-3} 
\end{align}
\end{subequations}
\eqref{equdecchampbLL-3} is of the form \eqref{lemdecequH-11-dersymb}. 

Note that 
\begin{align*}
\varphi &= \xi_0^3 + \xi_0 |\xi|^2 - \eta_0^3 - \eta_0 |\eta|^2 - \sigma_0^3 - \sigma_0 |\sigma|^2 \\
&= o(1) + 3 \xi_0 \eta_0 \sigma_0 
\end{align*}
But $|\overline{\xi}| \simeq |\overline{\eta}| \simeq |\overline{\sigma}| \gg |\eta| + |\sigma| \gtrsim |\xi|$, so $|\xi_0| \simeq |\eta_0| \simeq |\sigma_0| \simeq |\overline{\xi}|$. Hence, $1 = O(\varphi)$. There is no difficulty to decompose \eqref{equdecchampbLL-1}, where we keep the factor $m_b(\overline{\xi})$ to have \eqref{decompositionfinemg-gainkout}. 

Then, for \eqref{equdecchampbLL-2}, we decompose 
\begin{align*}
m_b(\overline{\xi}) \widehat{X}_b(\overline{\xi}) &= O(e_0) + O(\xi) = O(e_0) + O(\eta) + O(\sigma) 
\end{align*}
where $e_0 = \begin{pmatrix} 1 \\ 0 \\ 0 \end{pmatrix}$ is the first vector of the canonical basis. Therefore, 
\begin{subequations}
\begin{align}
\eqref{equdecchampbLL-2} &= \int_0^t \int e^{i s \varphi} \xi_0 \mu(\overline{\xi}, \overline{\eta}) \mu_{BBB}(\overline{\xi}, \overline{\eta}) m_{\widehat{\mathcal{L}}}(\overline{\eta}) m_{\widehat{\mathcal{L}}}(\overline{\sigma}) \widehat{F}_1(s, \overline{\eta}) O(1) \partial_{\xi_0} \widehat{F}_2(s, \overline{\sigma}) ~ d\overline{\eta} ds \label{equdecchampbLL-2-1} \\
&+ \int_0^t \int e^{i s \varphi} \xi_0 \mu(\overline{\xi}, \overline{\eta}) \mu_{BBB}(\overline{\xi}, \overline{\eta}) m_{\widehat{\mathcal{L}}}(\overline{\eta}) m_{\widehat{\mathcal{L}}}(\overline{\sigma}) \widehat{F}_1(s, \overline{\eta}) O(\sigma) \nabla_{\xi} \widehat{F}_2(s, \overline{\sigma}) ~ d\overline{\eta} ds \label{equdecchampbLL-2-2} \\
&+ \int_0^t \int e^{i s \varphi} \xi_0 \mu(\overline{\xi}, \overline{\eta}) \mu_{BBB}(\overline{\xi}, \overline{\eta}) m_{\widehat{\mathcal{L}}}(\overline{\eta}) m_{\widehat{\mathcal{L}}}(\overline{\sigma}) \widehat{F}_1(s, \overline{\eta}) O(\eta) \nabla_{\xi} \widehat{F}_2(s, \overline{\sigma}) ~ d\overline{\eta} ds \label{equdecchampbLL-2-3} 
\end{align}
\end{subequations}
\eqref{equdecchampbLL-2-1} and \eqref{equdecchampbLL-2-2} are of the form $\eqref{lemdecequH-07-resxsigmaac} + \eqref{lemdecequH-08-resxsigmabbon}$. For \eqref{equdecchampbLL-2-3}, we apply an integration by parts in $\eta$ and get, thanks to the $O(\eta)$ factor, uniquely terms of the form $\eqref{lemdecequH-03-resxetaac} + \eqref{lemdecequH-04-resxetabbon}$, \eqref{lemdecequH-11-dersymb} or a term similar to \eqref{equdecchampbLL-1} that can also be treated using $1 = O(\varphi)$, and that the derivation of the exponential gives automatically a factor $\nabla_{\overline{\eta}} \varphi$ so that \eqref{decompositionfinemg-gainkout} holds. 

\subsubsection{Interaction \texorpdfstring{$(BBB, \widehat{\mathcal{P}}\widehat{\mathcal{L}})$}{(BBB, PL)}}

Let us consider 
\begin{subequations}
\begin{align}
&\xi_0 m_b(\overline{\xi}) \widehat{X}_b(\overline{\xi}) \cdot \nabla_{\overline{\xi}} \widehat{I}_{\mu}^{(BBB, \widehat{\mathcal{P}}\widehat{\mathcal{L}})}[F_1, F_2](t, \overline{\xi}) \notag \\
&\quad = \int_0^t \int i s \xi_0 m_b(\overline{\xi}) \widehat{X}_b(\overline{\xi}) \cdot \nabla_{\overline{\xi}} \varphi e^{i s \varphi} \mu(\overline{\xi}, \overline{\eta}) \mu_{BBB}(\overline{\xi}, \overline{\eta}) m_{\widehat{\mathcal{P}}}(\overline{\eta}) m_{\widehat{\mathcal{L}}}(\overline{\sigma}) \widehat{F}_1(s, \overline{\eta}) \widehat{F}_2(s, \overline{\sigma}) ~ d\overline{\eta} ds \label{equdecchampbPL-1} \\
&\quad \quad + \int_0^t \int e^{i s \varphi} \xi_0 \mu(\overline{\xi}, \overline{\eta}) \mu_{BBB}(\overline{\xi}, \overline{\eta}) m_{\widehat{\mathcal{P}}}(\overline{\eta}) m_{\widehat{\mathcal{L}}}(\overline{\sigma}) \widehat{F}_1(s, \overline{\eta}) m_b(\overline{\xi}) \widehat{X}_b(\overline{\xi}) \cdot \nabla_{\overline{\xi}} \widehat{F}_2(s, \overline{\sigma}) ~ d\overline{\eta} ds \label{equdecchampbPL-2} \\
&\quad \quad + \int_0^t \int e^{i s \varphi} \xi_0 m_b(\overline{\xi}) \widehat{X}_b(\overline{\xi}) \cdot \nabla_{\overline{\xi}} \left( \mu(\overline{\xi}, \overline{\eta}) \mu_{BBB}(\overline{\xi}, \overline{\eta}) m_{\widehat{\mathcal{P}}}(\overline{\eta}) m_{\widehat{\mathcal{L}}}(\overline{\sigma}) \right) \widehat{F}_1(s, \overline{\eta}) \widehat{F}_2(s, \overline{\sigma}) ~ d\overline{\eta} ds \label{equdecchampbPL-3} 
\end{align}
\end{subequations}
\eqref{equdecchampbPL-3} is of the form \eqref{lemdecequH-11-dersymb}. 

We have that
\begin{align*}
\varphi &= \xi_0^3 + \xi_0 |\xi|^2 - \eta_0^3 - \eta_0 |\eta|^2 - \sigma_0^3 - \sigma_0 |\sigma|^2 \\
&= o(1) + \xi_0 |\xi|^2 
\end{align*}
But since $\xi_0 \simeq \sigma_0 \simeq |\overline{\sigma}| \simeq |\overline{\eta}| \simeq |\eta| \simeq |\xi|$, we deduce that $1 = O(\varphi)$. Again, there is no difficulty bringing back \eqref{equdecchampbPL-1} to terms of \eqref{lemdecequHGtot}, and \eqref{decompositionfinemg-gainkout} holds by the presence of $m_b(\overline{\xi})$. 

On the other hand, every term where the frequential derivative hits $F_1$ is automatically of the form \eqref{lemdecequH-06-resxetabphi}. 

The rest of the decomposition is similar to the cases $(BBB, \widehat{\mathcal{P}}\widehat{\mathcal{R}})$ or $(BBB, \widehat{\mathcal{P}}\widehat{\mathcal{C}})$: we remove from \eqref{equdecchampbPL-2} the term \eqref{lemdecequH-02-symsigma}, apply an integration by parts in $\overline{\eta}$, and get either a term similar to \eqref{equdecchampbPL-1} that we treat the same way (this time using the factor $\nabla_{\overline{\eta}} \varphi$ coming from the differentiation of the exponential to check \eqref{decompositionfinemg-gainkout}), a term of the form \eqref{lemdecequH-11-dersymb}, or a term where the derivative hits $\widehat{F}_1$, and therefore of the form \eqref{lemdecequH-06-resxetabphi}. 

\subsubsection{Interaction \texorpdfstring{$(BBB, \widehat{\mathcal{P}}\widehat{\mathcal{P}})$}{(BBB, PP)}}

Let us consider
\begin{align*}
&\xi_0 m_b(\overline{\xi}) \widehat{X}_b(\overline{\xi}) \cdot \nabla_{\overline{\xi}} \widehat{I}_{\mu}^{(BBB, \widehat{\mathcal{P}}\widehat{\mathcal{P}})}[F_1, F_2](t, \overline{\xi}) \\
&= \xi_0 m_b(\overline{\xi}) \widehat{X}_b(\overline{\xi}) \cdot \nabla_{\overline{\xi}} \left( \int_0^t \int e^{i s \varphi} \mu(\overline{\xi}, \overline{\eta}) \mu_{BBB}(\overline{\xi}, \overline{\eta}) m_{\widehat{\mathcal{P}}}(\overline{\eta}) m_{\widehat{\mathcal{P}}}(\overline{\sigma}) \widehat{F}_1(s, \overline{\eta}) \widehat{F}_2(s, \overline{\sigma}) ~ d\overline{\eta} ds \right)
\end{align*}
By symmetry of $\overline{\eta}$ and $\overline{\sigma}$ in this case, we can decompose again the integral: either we are in $\{ |\eta_0| \simeq |\sigma_0| \}$ or in $\{ |\eta_0| \ll |\sigma_0| \}$ (symmetrizing with $(\overline{\xi}, \overline{\eta}) \mapsto (\overline{\xi}, \overline{\sigma})$ the area $\{ |\sigma_0| \ll |\eta_0| \}$). 

Let us denote by $\mu_{BB}(\xi_0, \eta_0)$ the symbol localising on $\{ |\eta_0| \simeq |\sigma_0| \}$ and $\mu_{BH}(\xi_0, \eta_0)$ localising on $\{ |\eta_0| \ll |\sigma_0| \}$. These symbols depend only on $(\xi_0, \eta_0)$ and are of Coifman-Meyer type in these coordinates, hence they keep the Hölder property. 

Note that here, $|\overline{\xi}| \gtrsim |\overline{\eta}|+|\overline{\sigma}|$ and $\xi_0 = \eta_0 + \sigma_0$ forces $\overline{\xi}$ to be localised by $m_{\widehat{\mathcal{P}}}+m_{\widehat{\mathcal{R}}}$. In particular, $1 = O(m_b(\overline{\xi}))$ and \eqref{decompositionfinemg-gainkout} holds automatically. 

\paragraph{1.} Consider first the integral over $\{ |\eta_0| \simeq |\sigma_0| \}$. We get
\begin{subequations}
\begin{align}
&\begin{aligned}
\int_0^t \int i s \xi_0 m_b(\overline{\xi}) \widehat{X}_b(\overline{\xi}) \cdot \nabla_{\overline{\xi}} \varphi e^{i s \varphi} \mu(\overline{\xi}, \overline{\eta}) \mu_{BBB}(\overline{\xi}, \overline{\eta}) \mu_{BB}(\xi_0, \eta_0) m_{\widehat{\mathcal{P}}}(\overline{\eta}) m_{\widehat{\mathcal{P}}}(\overline{\sigma}) \widehat{F}_1(s, \overline{\eta}) \\
\widehat{F}_2(s, \overline{\sigma}) ~ d\overline{\eta} ds 
\end{aligned} \label{equdecchampbPP-BB-1} \\
&\quad \begin{aligned}
+ \int_0^t \int e^{i s \varphi} \xi_0 \mu(\overline{\xi}, \overline{\eta}) \mu_{BBB}(\overline{\xi}, \overline{\eta}) \mu_{BB}(\xi_0, \eta_0) m_{\widehat{\mathcal{P}}}(\overline{\eta}) m_{\widehat{\mathcal{P}}}(\overline{\sigma}) \widehat{F}_1(s, \overline{\eta}) m_b(\overline{\xi}) \\
\widehat{X}_b(\overline{\xi}) \cdot \nabla_{\overline{\xi}} \widehat{F}_2(s, \overline{\sigma}) ~ d\overline{\eta} ds 
\end{aligned} \label{equdecchampbPP-BB-2} \\
&\quad \begin{aligned}
+ \int_0^t \int e^{i s \varphi} \xi_0 m_b(\overline{\xi}) \widehat{X}_b(\overline{\xi}) \cdot \nabla_{\overline{\xi}} \left( \mu(\overline{\xi}, \overline{\eta}) \mu_{BBB}(\overline{\xi}, \overline{\eta}) \mu_{BB}(\xi_0, \eta_0) m_{\widehat{\mathcal{P}}}(\overline{\eta}) m_{\widehat{\mathcal{P}}}(\overline{\sigma}) \right) \\
\widehat{F}_1(s, \overline{\eta}) \widehat{F}_2(s, \overline{\sigma}) ~ d\overline{\eta} ds 
\end{aligned} \label{equdecchampbPP-BB-3} 
\end{align}
\end{subequations}
\eqref{equdecchampbPP-BB-3} is of the form \eqref{lemdecequH-11-dersymb}. (Indeed, even if $\mu_{BB}$ is differentiated, the singularity behaves like $\frac{1}{|\eta_0| + |\sigma_0|}$ which is compensated by $\xi_0$.) 

On the other hand, since $\xi_0 = \eta_0 + \sigma_0$, and that $\frac{\sigma_0}{\eta_0}, \frac{\eta_0}{\sigma_0}, \frac{|\eta_0|+|\sigma_0|}{|\eta_0|}, \frac{|\eta_0|+|\sigma_0|}{|\sigma_0|}$ are not singular here, \eqref{equdecchampbPP-BB-2} is of the form \eqref{lemdecequH-09-resxsigmabbonbis}. Likewise, \eqref{lemdecequH-01-symeta} and \eqref{lemdecequH-02-symsigma} are already of the form \eqref{lemdecequH-05-resxetabbonbis} and \eqref{lemdecequH-09-resxsigmabbonbis} and we don't need to make them appear. 

Finally, for \eqref{equdecchampbPP-BB-1}, we compute that 
\begin{align*}
\partial_{\eta_0} \varphi &= o(1) + |\sigma|^2 - |\eta|^2
\end{align*}
In particular, as soon as $1 = O(|\sigma| - |\eta|)$, we can apply an integration by parts along $\partial_{\eta_0}$ and use that $\xi_0 = O(\eta_0) = O(\sigma_0)$ to obtain terms from \eqref{lemdecequHGtot}. Let us therefore only consider the case localised on $|\sigma| - |\eta| = o(1)$. We consider two subcases. 

\paragraph{1.1} If $1 = O(J \eta \cdot \sigma)$ (i.e. $\eta$ and $\sigma$ are not aligned), then 
\begin{align*}
J \eta \cdot \nabla_{\eta} \varphi = 2 \sigma_0 J \eta \cdot \sigma
\end{align*}
and we can decompose in the presence of a factor $O(\eta_0^2)$. Then, 
\begin{align*}
\partial_{\eta_0} \varphi &= O(\eta_0^2) + |\sigma|^2 - |\eta|^2 \\
\varphi &= O(\eta_0^2) + \xi_0 |\xi|^2 - \eta_0 |\eta|^2 - \sigma_0 |\sigma|^2 \\
&= O(\eta_0^2) + O(\eta_0 \partial_{\eta_0} \varphi) + \xi_0 (|\xi|^2 - |\sigma|^2) 
\end{align*}
so we can decompose also in the presence of a factor $O(\xi_0 (|\xi|^2 - |\sigma|^2))$. Finally, 
\begin{align*}
\xi_0 |\overline{\xi}| \widehat{X}_b(\overline{\xi}) \cdot \nabla_{\overline{\xi}} \varphi 
&= O(\eta_0^2) + O(\varphi) + O(\eta_0 \partial_{\eta_0} \varphi) 
\end{align*}
which concludes. 

\paragraph{1.2} If $J \eta \cdot \sigma = o(1)$ (i.e. $\eta$ and $\sigma$ are close to alignment), then since $|\sigma| - |\eta| = o(1)$ and $1 = O(\xi)$, we must have $\sigma = \eta + o(1)$, so $1 = O(|\xi|^2 - |\sigma|^2)$. Therefore, from the computation of $\varphi$: 
\begin{align*}
\varphi + \sigma_0 \partial_{\eta_0} \varphi &= \xi_0^3 + \xi_0 |\xi|^2 - \eta_0^3 - \eta_0 |\eta|^2 - \sigma_0^3 - \sigma_0 |\sigma|^2 + \sigma_0 (3 \sigma_0^2 + |\sigma|^2 - 3 \eta_0^2 - |\eta|^2) \\
&= \xi_0 (|\xi|^2 - |\eta|^2 + o(1)) 
\end{align*}
This means that we can decompose whenever we have a factor $O(\xi_0)$, which is enough. 

\paragraph{2.} Let us consider the case $\{ |\eta_0| \ll |\sigma_0| \}$. We then have
\begin{subequations}
\begin{align}
&\begin{aligned}
\int_0^t \int i s \xi_0 m_b(\overline{\xi}) \widehat{X}_b(\overline{\xi}) \cdot \nabla_{\overline{\xi}} \varphi e^{i s \varphi} \mu(\overline{\xi}, \overline{\eta}) \mu_{BBB}(\overline{\xi}, \overline{\eta}) \mu_{BH}(\xi_0, \eta_0) m_{\widehat{\mathcal{P}}}(\overline{\eta}) m_{\widehat{\mathcal{P}}}(\overline{\sigma}) \widehat{F}_1(s, \overline{\eta}) \\
\widehat{F}_2(s, \overline{\sigma}) ~ d\overline{\eta} ds 
\end{aligned} \label{equdecchampbPP-BH-1} \\
&\quad \begin{aligned}
+ \int_0^t \int e^{i s \varphi} \xi_0 \mu(\overline{\xi}, \overline{\eta}) \mu_{BBB}(\overline{\xi}, \overline{\eta}) \mu_{BH}(\xi_0, \eta_0) m_{\widehat{\mathcal{P}}}(\overline{\eta}) m_{\widehat{\mathcal{P}}}(\overline{\sigma}) \widehat{F}_1(s, \overline{\eta}) m_b(\overline{\xi}) \\
\widehat{X}_b(\overline{\xi}) \cdot \nabla_{\overline{\xi}} \widehat{F}_2(s, \overline{\sigma}) ~ d\overline{\eta} ds 
\end{aligned} \label{equdecchampbPP-BH-2} \\
&\quad \begin{aligned}
+ \int_0^t \int e^{i s \varphi} \xi_0 m_b(\overline{\xi}) \widehat{X}_b(\overline{\xi}) \cdot \nabla_{\overline{\xi}} \left( \mu(\overline{\xi}, \overline{\eta}) \mu_{BBB}(\overline{\xi}, \overline{\eta}) \mu_{BH}(\xi_0, \eta_0) m_{\widehat{\mathcal{P}}}(\overline{\eta}) m_{\widehat{\mathcal{P}}}(\overline{\sigma}) \right) \\
\widehat{F}_1(s, \overline{\eta}) \widehat{F}_2(s, \overline{\sigma}) ~ d\overline{\eta} ds 
\end{aligned} \label{equdecchampbPP-BH-3} 
\end{align}
\end{subequations}
Again, \eqref{equdecchampbPP-BH-3} is of the form \eqref{lemdecequH-11-dersymb}. 

In this situation, \eqref{lemdecequH-01-symeta} is of the form \eqref{lemdecequH-05-resxetabbonbis} and we don't need to make it appear explicitely. 

We need however \eqref{lemdecequH-02-symsigma} that we remove from \eqref{equdecchampbPP-BH-2} before applying an integration by parts: 
\begin{subequations}
\begin{align}
&\eqref{equdecchampbPP-BH-2} - \eqref{lemdecequH-02-symsigma} \notag \\
&\quad \begin{aligned}
= \int_0^t \int e^{i s \varphi} \mu(\overline{\xi}, \overline{\eta}) \mu_{BBB}(\overline{\xi}, \overline{\eta}) \mu_{BH}(\xi_0, \eta_0) m_{\widehat{\mathcal{P}}}(\overline{\eta}) m_{\widehat{\mathcal{P}}}(\overline{\sigma}) \widehat{F}_1(s, \overline{\eta}) \\
\left( \xi_0 m_b(\overline{\xi}) \widehat{X}_b(\overline{\xi}) - \sigma_0 m_b(\overline{\sigma}) \widehat{X}_b(\overline{\sigma}) \right) \cdot \nabla_{\overline{\xi}} \widehat{F}_2(s, \overline{\sigma}) ~ d\overline{\eta} ds 
\end{aligned} \notag \\
&\quad \begin{aligned}
= \int_0^t \int i s \left( \xi_0 m_b(\overline{\xi}) \widehat{X}_b(\overline{\xi}) - \sigma_0 m_b(\overline{\sigma}) \widehat{X}_b(\overline{\sigma}) \right) \cdot \nabla_{\overline{\eta}} \varphi e^{i s \varphi} \mu(\overline{\xi}, \overline{\eta}) \mu_{BBB}(\overline{\xi}, \overline{\eta}) \\
\mu_{BH}(\xi_0, \eta_0) m_{\widehat{\mathcal{P}}}(\overline{\eta}) m_{\widehat{\mathcal{P}}}(\overline{\sigma}) \widehat{F}_1(s, \overline{\eta}) \widehat{F}_2(s, \overline{\sigma}) ~ d\overline{\eta} ds 
\end{aligned} \label{equdecchampbPP-BH-2-1} \\
&\quad \quad \begin{aligned}
+ \int_0^t \int e^{i s \varphi} \mu(\overline{\xi}, \overline{\eta}) \mu_{BBB}(\overline{\xi}, \overline{\eta}) \left( \xi_0 m_b(\overline{\xi}) \widehat{X}_b(\overline{\xi}) - \sigma_0 m_b(\overline{\sigma}) \widehat{X}_b(\overline{\sigma}) \right) \cdot \nabla_{\overline{\eta}} \widehat{F}_1(s, \overline{\eta}) \\
\mu_{BH}(\xi_0, \eta_0) m_{\widehat{\mathcal{P}}}(\overline{\eta}) m_{\widehat{\mathcal{P}}}(\overline{\sigma}) \widehat{F}_2(s, \overline{\sigma}) ~ d\overline{\eta} ds 
\end{aligned} \label{equdecchampbPP-BH-2-2} \\
&
+ \int_0^t \int \nabla_{\overline{\eta}} \cdot \left( \left( \xi_0 m_b(\overline{\xi}) \widehat{X}_b(\overline{\xi}) - \sigma_0 m_b(\overline{\sigma}) \widehat{X}_b(\overline{\sigma}) \right) \mu(\overline{\xi}, \overline{\eta}) \mu_{BBB}(\overline{\xi}, \overline{\eta}) \mu_{BH}(\xi_0, \eta_0) m_{\widehat{\mathcal{P}}}(\overline{\eta}) m_{\widehat{\mathcal{P}}}(\overline{\sigma}) \right) \notag \\
&\quad \quad \quad \quad \quad \quad \quad \quad \quad \quad \quad \quad \quad \quad \quad \quad \quad \quad \quad \quad \quad \quad \quad \quad e^{i s \varphi} \widehat{F}_1(s, \overline{\eta}) \widehat{F}_2(s, \overline{\sigma}) ~ d\overline{\eta} ds \label{equdecchampbPP-BH-2-3} 
\end{align}
\end{subequations}
\eqref{equdecchampbPP-BH-2-3} is of the form \eqref{lemdecequH-11-dersymb}. 

For \eqref{equdecchampbPP-BH-2-2}, just as in the case $(BBB, \widehat{\mathcal{P}}\widehat{\mathcal{R}})$, we project: 
\begin{align*}
&\xi_0 m_b(\overline{\xi}) \widehat{X}_b(\overline{\xi}) - \sigma_0 m_b(\overline{\sigma}) \widehat{X}_b(\overline{\sigma}) 
\\
&= O(\eta_0) + O(\sigma_0 \widehat{X}_a(\overline{\eta})) + O(\sigma_0 \widehat{X}_c(\overline{\eta})) + \sigma_0 \left( m_b(\overline{\xi}) P_b^b(\overline{\xi}, \overline{\eta}) - m_b(\overline{\sigma}) P_b^b(\overline{\sigma}, \overline{\eta}) \right) \widehat{X}_b(\overline{\eta}) \\
&= O(\eta_0) + O(\sigma_0 \widehat{X}_a(\overline{\eta})) + O(\sigma_0 \widehat{X}_c(\overline{\eta})) + O(\sigma_0 (|\xi| - |\sigma|)) \widehat{X}_b(\overline{\eta})
\end{align*}
But again, $\varphi = O(\eta_0) + \sigma_0 (|\xi|^2 - |\sigma|^2)$, so finally
\begin{align*}
\xi_0 m_b(\overline{\xi}) \widehat{X}_b(\overline{\xi}) - \sigma_0 m_b(\overline{\sigma}) \widehat{X}_b(\overline{\sigma}) 
&= O(\eta_0) + O(\sigma_0 \widehat{X}_a(\overline{\eta})) + O(\sigma_0 \widehat{X}_c(\overline{\eta})) + O(\varphi) \widehat{X}_b(\overline{\eta})
\end{align*}
which justifies that we can decompose 
\eqref{equdecchampbPP-BH-2-2} into \eqref{lemdecequH-03-resxetaac}, \eqref{lemdecequH-05-resxetabbonbis} and \eqref{lemdecequH-06-resxetabphi}. 

Finally, we group \eqref{equdecchampbPP-BH-2-1} and \eqref{equdecchampbPP-BH-1}: 
\begin{align*}
&\eqref{equdecchampbPP-BH-1} + \eqref{equdecchampbPP-BH-2-1} \\
&\quad = 
\int_0^t \int i s \left( \xi_0 m_b(\overline{\xi}) \widehat{X}_b(\overline{\xi}) \cdot (\nabla_{\overline{\xi}} + \nabla_{\overline{\eta}}) - \sigma_0 m_b(\overline{\sigma}) \widehat{X}_b(\overline{\sigma}) \cdot \nabla_{\overline{\eta}} \right) \varphi e^{i s \varphi} \\
&\quad \quad \quad \quad \quad \quad \mu(\overline{\xi}, \overline{\eta}) \mu_{BBB}(\overline{\xi}, \overline{\eta}) \mu_{BH}(\xi_0, \eta_0) m_{\widehat{\mathcal{P}}}(\overline{\eta}) m_{\widehat{\mathcal{P}}}(\overline{\sigma}) \widehat{F}_1(s, \overline{\eta}) \widehat{F}_2(s, \overline{\sigma}) ~ d\overline{\eta} ds
\end{align*}

We have that 
\begin{align*}
\nabla_{\eta} \varphi &= 2 \sigma_0 \sigma - 2 \eta_0 \eta \\
&= 2 \sigma_0 \sigma (1 + o(1)) 
\end{align*}
which means we can decompose whenever we have a factor $O(\sigma_0 \eta_0)$. Then, 
\begin{align*}
\varphi - \eta_0 \partial_{\eta_0} \varphi &= \xi_0^3 + \xi_0 |\xi|^2 - \eta_0^3 - \eta_0 |\eta|^2 - \sigma_0^3 - \sigma_0 |\sigma|^2 - \eta_0 (3 \sigma_0^2 + |\sigma|^2 - 3 \eta_0^2 - |\eta|^2) \\
&= 3 \xi_0 \eta_0 \sigma_0 + \xi_0 |\xi|^2 - \xi_0 |\sigma|^2 - 3 \eta_0 (\sigma_0 - \eta_0) \xi_0 \\
&= O(\eta_0 \sigma_0) + \xi_0 (|\xi|^2 - |\sigma|^2) 
\end{align*}
so we can also decompose if we have a factor $O(\xi_0 (|\xi|^2 - |\sigma|^2))$. Finally, we compute that 
\begin{align*}
m_b(\overline{\xi}) - m_b(\overline{\sigma}) &= O(\eta_0) + O(|\xi|^2 - |\sigma|^2) 
\end{align*}
hence 
\begin{align*}
&\left( \xi_0 m_b(\overline{\xi}) \widehat{X}_b(\overline{\xi}) \cdot (\nabla_{\overline{\xi}} + \nabla_{\overline{\eta}}) - \sigma_0 m_b(\overline{\sigma}) \widehat{X}_b(\overline{\sigma}) \cdot \nabla_{\overline{\eta}} \right) \varphi \\
&\quad = O(\sigma_0 \eta_0) + O(\xi_0 (|\xi|^2 - |\sigma|^2)) + m_b(\overline{\xi}) \eta_0 (|\sigma|^2 - |\eta|^2) \\
&\quad = O(\sigma_0 \eta_0) + O(\xi_0 (|\xi|^2 - |\sigma|^2)) + O(\eta_0 \partial_{\eta_0} \varphi) 
\end{align*} 
which justifies that we can apply all the desired integrations by parts to decompose $\eqref{equdecchampbPP-BH-1} + \eqref{equdecchampbPP-BH-2-1}$. 

\subsubsection{Interaction \texorpdfstring{$(HBH, \widehat{\mathcal{R}})$}{(HBH, R)}} 

Let us consider
\begin{subequations}
\begin{align}
&\xi_0 m_b(\overline{\xi}) \widehat{X}_b(\overline{\xi}) \cdot \nabla_{\overline{\xi}} \widehat{I}_{\mu}^{(HBH, \widehat{\mathcal{R}})}[F_1, F_2](t, \overline{\xi}) \notag \\
&\quad = \int_0^t \int i s \xi_0 m_b(\overline{\xi}) \widehat{X}_b(\overline{\xi}) \cdot \nabla_{\overline{\xi}} \varphi e^{i s \varphi} \mu(\overline{\xi}, \overline{\eta}) \mu_{HBH}(\overline{\xi}, \overline{\eta}) m_{\widehat{\mathcal{R}}}(\overline{\sigma}) \widehat{F}_1(s, \overline{\eta}) \widehat{F}_2(s, \overline{\sigma}) ~ d\overline{\eta} ds \label{equdecchampbHBHR-1} \\
&\quad \quad + \int_0^t \int e^{i s \varphi} \xi_0 \mu(\overline{\xi}, \overline{\eta}) \mu_{HBH}(\overline{\xi}, \overline{\eta}) m_{\widehat{\mathcal{R}}}(\overline{\sigma}) \widehat{F}_1(s, \overline{\eta}) m_b(\overline{\xi}) \widehat{X}_b(\overline{\xi}) \cdot \nabla_{\overline{\xi}} \widehat{F}_2(s, \overline{\sigma}) ~ d\overline{\eta} ds \label{equdecchampbHBHR-2} \\
&\quad \quad + \int_0^t \int e^{i s \varphi} \xi_0 m_b(\overline{\xi}) \widehat{X}_b(\overline{\xi}) \cdot \nabla_{\overline{\xi}} \left( \mu(\overline{\xi}, \overline{\eta}) \mu_{HBH}(\overline{\xi}, \overline{\eta}) m_{\widehat{\mathcal{R}}}(\overline{\sigma}) \right) \widehat{F}_1(s, \overline{\eta}) \widehat{F}_2(s, \overline{\sigma}) ~ d\overline{\eta} ds \label{equdecchampbHBHR-3} 
\end{align}
\end{subequations}
\eqref{equdecchampbHBHR-3} is of the form \eqref{lemdecequH-11-dersymb}. 

Note that in any case, here, as well as for all the $HBH$ interactions, $1 = O(\nabla_{\overline{\eta}} \varphi)$ so that \eqref{decompositionfinemg-gainkout} is automatic. 

\paragraph{1.} Let us first localise $\overline{\eta}$ away enough from $\widehat{\mathcal{P}}$. Then 
\begin{align}
m_b(\overline{\xi}) \xi_0 \widehat{X}_b(\overline{\xi}) = m_b(\overline{\sigma}) \sigma_0 \widehat{X}_b(\overline{\sigma}) + O(\overline{\eta}) \label{equationdecompositionHBHsymsigma} 
\end{align}
so that \eqref{equdecchampbHBHR-2} can be decomposed into \eqref{lemdecequH-02-symsigma} and \eqref{lemdecequH-08-resxsigmabbon}. 

Then, 
\begin{align*}
\widehat{X}_b(\overline{\xi}) \cdot \nabla_{\overline{\xi}} \varphi &= O\left( \nabla_{\overline{\xi}} \varphi \right) \\
&= O(\overline{\eta}) 
\end{align*}
If we furthermore localise to have 
\[ 1 = O\left( \widehat{X}_a(\overline{\eta}) \cdot \nabla_{\overline{\eta}} \varphi \right) + O\left( m_b(\overline{\eta}) \widehat{X}_b(\overline{\eta}) \cdot \nabla_{\overline{\eta}} \varphi \right) + O\left( m_c(\overline{\eta}) \widehat{X}_c(\overline{\eta}) \cdot \nabla_{\overline{\eta}} \varphi \right) \]
then it is easy to apply an integration by parts in $\overline{\eta}$ on \eqref{equdecchampbHBHR-1} and bring it back to terms of \eqref{lemdecequHGtot} thanks to the presence of the factor $O(\overline{\eta})$. 

Localise then to have
\begin{align}
\widehat{X}_a(\overline{\eta}) \cdot \nabla_{\overline{\eta}} \varphi = o(1), \quad m_b(\overline{\eta}) \widehat{X}_b(\overline{\eta}) \cdot \nabla_{\overline{\eta}} \varphi = o(1), \quad m_c(\overline{\eta}) \widehat{X}_c(\overline{\eta}) \cdot \nabla_{\overline{\eta}} \varphi = o(1) \label{hypothesetravailHBHR1} 
\end{align}
Since $\partial_{\eta_0} \varphi = 3 \sigma_0^2 + |\sigma|^2 - 3 \eta_0^2 - |\eta|^2$ is not $o(1)$, this forces $\overline{\eta}$ to be localised near $\widehat{\mathcal{C}}$. But then, by Lemma \ref{lemcalculsconecoordonneesconiquesvarphi}, 
\begin{align*}
\xi_t^{\overline{\eta} \overline{\sigma}} &= O\left( \widehat{X}_c(\overline{\eta}) \cdot \nabla_{\overline{\eta}} \varphi \right) \\
\widehat{X}_a(\overline{\eta}) \cdot \nabla_{\overline{\eta}} \varphi &= \frac{\eta_0}{|\overline{\eta}|} \left( \overline{\sigma}_a^{\overline{\eta}} \right)^2 + O\left( \overline{\eta} \right) + O\left( \xi_t^{\overline{\eta} \overline{\sigma}} \right) 
\end{align*}
But since $\overline{\sigma}$ is localised by $m_{\widehat{\mathcal{R}}}$, 
\begin{align*}
\overline{\sigma}_a^{\overline{\eta}} &= \sqrt{3} \epsilon^{\overline{\eta} \overline{\sigma}} |\sigma_0| + \theta^{\overline{\eta} \overline{\sigma}} |\sigma| \\
&= \epsilon^{\overline{\eta} \overline{\sigma}} \left( \sqrt{3} |\sigma_0| \pm |\sigma| \right) + O\left( \xi_t^{\overline{\eta} \overline{\sigma}} \right) \\
&\simeq |\overline{\sigma}|^2 
\end{align*}
We deduce that hypothesis \eqref{hypothesetravailHBHR1} is never satisfied. 

\paragraph{2.} Assume now that $\overline{\eta}$ is localised near $\widehat{\mathcal{P}}$. 

We then have
\begin{align*}
\nabla_{\eta} \varphi &= 2 \sigma_0 \sigma + o(1) 
\end{align*}
so that $1 = O(\nabla_{\eta} \varphi)$. We may therefore easily decompose any term having a factor $O(\eta_0)$ in \eqref{equdecchampbHBHR-1}. Then, 
\begin{align*}
\varphi &= O(\eta_0) + \xi_0 (|\xi|^2 - |\sigma|^2) 
\end{align*}
so we can also decompose having a factor $O(|\xi| - |\sigma|)$. 

For \eqref{equdecchampbHBHR-2}, we need to remove \eqref{lemdecequH-02-symsigma} and apply an integration by parts: 
\begin{subequations}
\begin{align}
&\eqref{equdecchampbHBHR-2} - \eqref{lemdecequH-02-symsigma} \notag \\
&\begin{aligned}
= \int_0^t \int e^{i s \varphi} \mu(\overline{\xi}, \overline{\eta}) \mu_{HBH}(\overline{\xi}, \overline{\eta}) m_{\widehat{\mathcal{R}}}(\overline{\sigma}) \left( \xi_0 m_b(\overline{\xi}) \widehat{X}_b(\overline{\xi}) - \sigma_0 m_b(\overline{\sigma}) \widehat{X}_b(\overline{\sigma}) \right) \cdot \nabla_{\overline{\xi}} \widehat{F}_2(s, \overline{\sigma}) \\
m_{\widehat{\mathcal{P}}}(\overline{\eta}) \widehat{F}_1(s, \overline{\eta}) ~ d\overline{\eta} ds 
\end{aligned} \notag \\
&\begin{aligned}
= \int_0^t \int i s \left( \xi_0 m_b(\overline{\xi}) \widehat{X}_b(\overline{\xi}) - \sigma_0 m_b(\overline{\sigma}) \widehat{X}_b(\overline{\sigma}) \right) \cdot \nabla_{\overline{\eta}} \varphi e^{i s \varphi} \mu(\overline{\xi}, \overline{\eta}) \mu_{HBH}(\overline{\xi}, \overline{\eta}) m_{\widehat{\mathcal{P}}}(\overline{\eta}) \\
m_{\widehat{\mathcal{R}}}(\overline{\sigma}) \widehat{F}_1(s, \overline{\eta}) \widehat{F}_2(s, \overline{\sigma}) ~ d\overline{\eta} ds 
\end{aligned} \label{equdecchampbHBHR-2-1} \\
&\quad \begin{aligned}
+ \int_0^t \int e^{i s \varphi} \mu(\overline{\xi}, \overline{\eta}) \mu_{HBH}(\overline{\xi}, \overline{\eta}) \left( \xi_0 m_b(\overline{\xi}) \widehat{X}_b(\overline{\xi}) - \sigma_0 m_b(\overline{\sigma}) \widehat{X}_b(\overline{\sigma}) \right) \cdot \nabla_{\overline{\eta}} \widehat{F}_1(s, \overline{\eta}) \\
m_{\widehat{\mathcal{P}}}(\overline{\eta}) m_{\widehat{\mathcal{R}}}(\overline{\sigma}) \widehat{F}_2(s, \overline{\sigma}) ~ d\overline{\eta} ds 
\end{aligned} \label{equdecchampbHBHR-2-2} \\
&\begin{aligned}
+ \int_0^t \int \nabla_{\overline{\xi}} \cdot \left( \left( \xi_0 m_b(\overline{\xi}) \widehat{X}_b(\overline{\xi}) - \sigma_0 m_b(\overline{\sigma}) \widehat{X}_b(\overline{\sigma}) \right) \mu(\overline{\xi}, \overline{\eta}) \mu_{HBH}(\overline{\xi}, \overline{\eta}) m_{\widehat{\mathcal{P}}}(\overline{\eta}) m_{\widehat{\mathcal{R}}}(\overline{\sigma}) \right) \\
e^{i s \varphi} \widehat{F}_1(s, \overline{\eta}) \widehat{F}_2(s, \overline{\sigma}) ~ d\overline{\eta} ds 
\end{aligned} \label{equdecchampbHBHR-2-3} 
\end{align}
\end{subequations}
\eqref{equdecchampbHBHR-2-3} is of the form \eqref{lemdecequH-11-dersymb}.

For \eqref{equdecchampbHBHR-2-2}, we project: 
\begin{align*}
\xi_0 m_b(\overline{\xi}) \widehat{X}_b(\overline{\xi}) - \sigma_0 m_b(\overline{\sigma}) \widehat{X}_b(\overline{\sigma}) 
&= O(\eta_0) + O(\widehat{X}_a(\overline{\eta})) + O(\widehat{X}_c(\overline{\eta})) + O(|\xi| - |\sigma|) \widehat{X}_b(\overline{\eta}) 
\end{align*}
As before, $\varphi = O(\eta_0) + \sigma_0 (|\xi|^2 - |\sigma|^2)$, so we may rewrite \eqref{equdecchampbHBHR-2-2} in the form of terms of \eqref{lemdecequHGtot}. 

Finally, we group \eqref{equdecchampbHBHR-1} and \eqref{equdecchampbHBHR-2-1} and get:  
\begin{align*}
\int_0^t \int i s \left( \xi_0 m_b(\overline{\xi}) \widehat{X}_b(\overline{\xi}) \cdot \left( \nabla_{\overline{\xi}} + \nabla_{\overline{\eta}} \right) - \sigma_0 m_b(\overline{\sigma}) \widehat{X}_b(\overline{\sigma}) \cdot \nabla_{\overline{\eta}} \right) \varphi e^{i s \varphi} \mu(\overline{\xi}, \overline{\eta}) \mu_{HBH}(\overline{\xi}, \overline{\eta}) m_{\widehat{\mathcal{P}}}(\overline{\eta}) \\
m_{\widehat{\mathcal{R}}}(\overline{\sigma}) \widehat{F}_1(s, \overline{\eta}) \widehat{F}_2(s, \overline{\sigma}) ~ d\overline{\eta} ds
\end{align*}
We then compute
\begin{align*}
&\left( \xi_0 m_b(\overline{\xi}) \widehat{X}_b(\overline{\xi}) \cdot \left( \nabla_{\overline{\xi}} + \nabla_{\overline{\eta}} \right) - \sigma_0 m_b(\overline{\sigma}) \widehat{X}_b(\overline{\sigma}) \cdot \nabla_{\overline{\eta}} \right) \varphi \\
&\quad = O(\eta_0) + O(|\xi| - |\sigma|)
\end{align*}
which concludes. 

\subsubsection{Interaction \texorpdfstring{$(HBH, \widehat{\mathcal{C}})$}{(HBH, C)}}

Let us consider
\begin{subequations}
\begin{align}
&\xi_0 m_b(\overline{\xi}) \widehat{X}_b(\overline{\xi}) \cdot \nabla_{\overline{\xi}} \widehat{I}_{\mu}^{(HBH, \widehat{\mathcal{C}})}[F_1, F_2](t, \overline{\xi}) \notag \\
&\quad = \int_0^t i s \xi_0 m_b(\overline{\xi}) \widehat{X}_b(\overline{\xi}) \cdot \nabla_{\overline{\xi}} \varphi \int e^{i s \varphi} \mu(\overline{\xi}, \overline{\eta}) \mu_{HBH}(\overline{\xi}, \overline{\eta}) m_{\widehat{\mathcal{C}}}(\overline{\sigma}) \widehat{F}_1(s, \overline{\eta}) \widehat{F}_2(s, \overline{\sigma}) ~ d\overline{\eta} ds \label{equdecchampbHBHC-1} \\
&\quad \quad + \int_0^t \int e^{i s \varphi} \xi_0 \mu(\overline{\xi}, \overline{\eta}) \mu_{HBH}(\overline{\xi}, \overline{\eta}) m_{\widehat{\mathcal{C}}}(\overline{\sigma}) \widehat{F}_1(s, \overline{\eta}) m_b(\overline{\xi}) \widehat{X}_b(\overline{\xi}) \cdot \nabla_{\overline{\xi}} \widehat{F}_2(s, \overline{\sigma}) ~ d\overline{\eta} ds \label{equdecchampbHBHC-2} \\
&\quad \quad + \int_0^t \int e^{i s \varphi} \xi_0 m_b(\overline{\xi}) \widehat{X}_b(\overline{\xi}) \cdot \nabla_{\overline{\xi}} \left( \mu(\overline{\xi}, \overline{\eta}) \mu_{HBH}(\overline{\xi}, \overline{\eta}) m_{\widehat{\mathcal{C}}}(\overline{\sigma}) \right) \widehat{F}_1(s, \overline{\eta}) \widehat{F}_2(s, \overline{\sigma}) ~ d\overline{\eta} ds \label{equdecchampbHBHC-3} 
\end{align}
\end{subequations}
\eqref{equdecchampbHBHC-3} is of the form \eqref{lemdecequH-11-dersymb}. 

Here, we need to further localise
\[ 1 = m_{\widehat{\mathcal{R}}}(\overline{\eta}) + m_{\widehat{\mathcal{C}}}(\overline{\eta}) + m_{\widehat{\mathcal{L}}}(\overline{\eta}) + m_{\widehat{\mathcal{P}}}(\overline{\eta}) \]

\paragraph{1.} Let us localise $\overline{\eta}$ by $m_{\widehat{\mathcal{R}}}(\overline{\eta})$. Then, 
\begin{align*}
\widehat{X}_a(\overline{\sigma}) \cdot \nabla_{\overline{\eta}} \varphi &= \frac{\sigma_0}{|\overline{\sigma}|} (3 \sigma_0^2 + |\sigma|^2) + o(1) 
\end{align*}
Therefore, using $\nabla_{\overline{\xi}} \varphi = O(\overline{\eta})$, we may apply an integration by parts along $\widehat{X}_a(\overline{\sigma})$ on \eqref{equdecchampbHBHC-1} and recover terms from \eqref{lemdecequHGtot}. 

On the other hand, we decompose \eqref{equdecchampbHBHC-2} using
\begin{subequations} \label{decompositionequdecchampbHBHC-2} 
\begin{align}
m_b(\overline{\xi}) \xi_0 \widehat{X}_b(\overline{\xi}) &= m_b(\overline{\sigma}) \sigma_0 \widehat{X}_b(\overline{\sigma}) \label{decompositionequdecchampbHBHC-2-1} \\
&+ O(\overline{\eta}) \widehat{X}_a(\overline{\sigma}) + O(\overline{\eta}) \widehat{X}_c(\overline{\sigma}) \label{decompositionequdecchampbHBHC-2-2} \\
&+ \left( m_b(\overline{\xi}) \xi_0 P_b^b(\overline{\xi}, \overline{\sigma}) - m_b(\overline{\sigma}) \sigma_0 \right) \widehat{X}_b(\overline{\sigma}) \label{decompositionequdecchampbHBHC-2-3} 
\end{align}
\end{subequations}
\eqref{decompositionequdecchampbHBHC-2-1} contributes as \eqref{lemdecequH-02-symsigma}, \eqref{decompositionequdecchampbHBHC-2-2} as \eqref{lemdecequH-07-resxsigmaac}. 

We then further rewrite: 
\begin{align*}
\eqref{decompositionequdecchampbHBHC-2-3} &= O\left( \overline{\eta} \overline{\sigma}_b^{\overline{\sigma}} \right) \widehat{X}_b(\overline{\sigma}) + O(1) \left( \overline{\xi}_b^{\overline{\xi}} - \overline{\sigma}_b^{\overline{\sigma}} \right) \widehat{X}_b(\overline{\sigma}) 
\end{align*}
The first term contributes like \eqref{lemdecequH-08-resxsigmabbon}. For the second, we apply an integration by parts in $\overline{\eta}$ and get: 
\begin{align*}
&\int_0^t \int i s O\left( \overline{\xi}_b^{\overline{\xi}} - \overline{\sigma}_b^{\overline{\sigma}} \right) \nabla_{\overline{\eta}} \varphi e^{i s \varphi} \xi_0 \mu(\overline{\xi}, \overline{\eta}) \mu_{HBH}(\overline{\xi}, \overline{\eta}) m_{\widehat{\mathcal{R}}}(\overline{\eta}) m_{\widehat{\mathcal{C}}}(\overline{\sigma}) \widehat{F}_1(s, \overline{\eta}) \widehat{F}_2(s, \overline{\sigma}) ~ d\overline{\eta} ds
\end{align*}
up to terms from \eqref{lemdecequHGtot} since $\overline{\xi}_b^{\overline{\xi}} - \overline{\sigma}_b^{\overline{\sigma}} = O(\overline{\eta})$. We can then decompose the same way as \eqref{equdecchampbHBHC-1}. 

\paragraph{2.} Localise now $\overline{\eta}$ by $m_{\widehat{\mathcal{C}}}(\overline{\eta})$. 

We start with \eqref{equdecchampbHBHC-2}. Again we apply \eqref{decompositionequdecchampbHBHC-2} to get back to the contribution of \eqref{decompositionequdecchampbHBHC-2-3}; then we can rewrite: 
\begin{align*}
\eqref{decompositionequdecchampbHBHC-2-3} &= O\left( \overline{\eta} \overline{\sigma}_b^{\overline{\sigma}} \right) \widehat{X}_b(\overline{\sigma}) + P_b^b(\overline{\xi}, \overline{\sigma}) \frac{(3 \xi_0^2 - |\xi|^2 - 3 \sigma_0^2 + |\sigma|^2) |\xi|}{|\overline{\xi}|^3} \xi_0 \widehat{X}_b(\overline{\sigma}) \\
&= O\left( \overline{\eta} \overline{\sigma}_b^{\overline{\sigma}} \right) \widehat{X}_b(\overline{\sigma}) + O(\overline{\eta}) \widehat{X}_a(\overline{\sigma}) + O(\overline{\eta}) \widehat{X}_c(\overline{\sigma}) \\
&\quad + P_b^b(\overline{\xi}, \overline{\sigma}) \frac{(3 \xi_0^2 - |\xi|^2 - 3 \sigma_0^2 + |\sigma|^2) |\xi|}{|\overline{\xi}|^3} |\xi_0| \widehat{X}_{b-\widehat{\mathcal{C}}}(\overline{\eta}, \overline{\sigma}) 
\end{align*}
and again it is enough to restrict our attention to the last term, for which we apply an integration by parts in $\overline{\eta}$. Since $3 \xi_0^2 - |\xi|^2 - 3 \sigma_0^2 + |\sigma|^2 = O(\overline{\eta})$, we get terms from the decomposition plus
\begin{align*}
\int_0^t \int i s P_b^b(\overline{\xi}, \overline{\sigma}) \frac{(3 \xi_0^2 - |\xi|^2 - 3 \sigma_0^2 + |\sigma|^2)|\xi|}{|\overline{\xi}|^3} |\xi_0| \widehat{X}_{b-\widehat{\mathcal{C}}}(\overline{\eta}, \overline{\sigma}) \cdot \nabla_{\overline{\eta}} \varphi e^{i s \varphi} \mu(\overline{\xi}, \overline{\eta}) \mu_{HBH}(\overline{\xi}, \overline{\eta}) m_{\widehat{\mathcal{C}}}(\overline{\eta}) \\
m_{\widehat{\mathcal{C}}}(\overline{\sigma}) \widehat{F}_1(s, \overline{\eta}) \widehat{F}_2(s, \overline{\sigma}) ~ d\overline{\eta} ds
\end{align*}
that we group with \eqref{equdecchampbHBHC-1} and get: 
\begin{align} 
\int_0^t \int i s &\left( m_b(\overline{\xi}) \xi_0 \widehat{X}_b(\overline{\xi}) \cdot \nabla_{\overline{\xi}} + P_b^b(\overline{\xi}, \overline{\sigma}) \frac{(3 \xi_0^2 - |\xi|^2 - 3 \sigma_0^2 + |\sigma|^2)|\xi|}{|\overline{\xi}|^3} |\xi_0| \widehat{X}_{b-\widehat{\mathcal{C}}}(\overline{\eta}, \overline{\sigma}) \cdot \nabla_{\overline{\eta}} \right) \varphi e^{i s \varphi} \notag \\
&\mu(\overline{\xi}, \overline{\eta}) \mu_{HBH}(\overline{\xi}, \overline{\eta}) m_{\widehat{\mathcal{C}}}(\overline{\eta}) m_{\widehat{\mathcal{C}}}(\overline{\sigma}) \widehat{F}_1(s, \overline{\eta}) \widehat{F}_2(s, \overline{\sigma}) ~ d\overline{\eta} ds \label{equdecchampbHBHCC-1} 
\end{align}

We now apply a second localisation depending on $\epsilon^{\overline{\eta} \overline{\sigma}} \theta^{\overline{\eta} \overline{\sigma}}$, to know whether this quantity is close or not to $-1$. 

\paragraph{2.1} If $\epsilon^{\overline{\eta} \overline{\sigma}} \theta^{\overline{\eta} \overline{\sigma}}$ is away enough from $-1$, then either $\xi_t^{\overline{\eta} \overline{\sigma}}$ is away enough from $0$, or $\epsilon^{\overline{\eta} \overline{\sigma}} \theta^{\overline{\eta} \overline{\sigma}}$ is close to $1$. In either case, 
\begin{align*}
1 = O\left( \widehat{Y}(\overline{\eta}, \overline{\sigma}) \cdot \nabla_{\overline{\eta}} \varphi \right) 
\end{align*}
by Lemma \ref{lemcalculchampmodifie}, which means we can decompose whenever we have a factor $O(\overline{\eta})$ in \eqref{equdecchampbHBHCC-1}. But precisely $\nabla_{\overline{\xi}} \varphi = O(\overline{\eta})$ and $3 \xi_0^2 - |\xi|^2 - 3 \sigma_0^2 + |\sigma|^2 = O(\overline{\eta})$. 

\paragraph{2.2} If $\epsilon^{\overline{\eta} \overline{\sigma}} \theta^{\overline{\eta} \overline{\sigma}}$ is close to $-1$, then Lemma \ref{lemcalculchampmodifie} only provides a way to decompose when we have a factor $O\left( \xi_t^{\overline{\eta} \overline{\sigma}} |\overline{\eta}| \right)$, that is $O\left( \xi_t^{\overline{\xi} \overline{\sigma}} \right)$. Then, we have
\begin{align*}
\widehat{X}_a(\overline{\sigma}) \cdot \nabla_{\overline{\eta}} \varphi &= \frac{\sigma_0}{|\overline{\sigma}|} (3 \sigma_0^2 + 3 |\sigma|^2) + o(1) 
\end{align*}
so we also know how to decompose if a factor $O\left( \overline{\eta}_b^{\overline{\eta}} \right)$ is present. Note now that
\begin{align*}
\overline{\eta}_a^{\overline{\sigma}} &= \sqrt{3} \epsilon^{\overline{\eta} \overline{\sigma}} |\eta_0| + \theta^{\overline{\eta} \overline{\sigma}} |\eta| \\
&= \epsilon^{\overline{\eta} \overline{\sigma}} \left( \sqrt{3} |\eta_0| - |\eta| \right) + O\left( |\eta| \xi_t^{\overline{\eta} \overline{\sigma}} \right) \\
&= O\left( \overline{\eta}_b^{\overline{\eta}} \right) + O\left( \xi_t^{\overline{\xi} \overline{\sigma}} \right) \\
\overline{\eta}_b^{\overline{\sigma}} &= \epsilon^{\overline{\eta} \overline{\sigma}} \overline{\eta}_a^{\overline{\eta}} + O\left( \xi_t^{\overline{\xi} \overline{\sigma}} \right) \\
\overline{\xi}_a^{\overline{\sigma}} &= \overline{\xi}_a^{\overline{\xi}} + O\left( \xi_t^{\overline{\xi} \overline{\sigma}} \right) \\
\overline{\xi}_b^{\overline{\sigma}} &= \overline{\xi}_b^{\overline{\xi}} + O\left( \xi_t^{\overline{\xi} \overline{\sigma}} \right)
\end{align*}
Finally, by Lemmas \ref{lemcalculsconecoordonneesconiquesvarphi} and \ref{lemcalculvarphiKdV1D}, 
\begin{align*}
6 \sqrt{3} \frac{\sigma_0}{|\sigma_0|} \varphi &= \left( \overline{\xi}_a^{\overline{\sigma}} \right)^3 + \left( \overline{\xi}_b^{\overline{\sigma}} \right)^3 - \left( \overline{\eta}_a^{\overline{\sigma}} \right)^3 - \left( \overline{\eta}_b^{\overline{\sigma}} \right)^3 - \left( \overline{\sigma}_a^{\overline{\sigma}} \right)^3 - \left( \overline{\sigma}_b^{\overline{\sigma}} \right)^3 + O\left( \xi_t^{\overline{\xi} \overline{\sigma}} \right) + O\left( |\overline{\eta}| \xi_t^{\overline{\eta} \overline{\sigma}} \right) \\
&= 3 \epsilon^{\overline{\eta} \overline{\sigma}} \overline{\xi}_b^{\overline{\xi}} \overline{\eta}_a^{\overline{\eta}} \overline{\sigma}_b^{\overline{\sigma}} + O\left( \overline{\eta}_b^{\overline{\eta}} \right) + O\left( \xi_t^{\overline{\xi} \overline{\sigma}} \right) 
\end{align*}
hence
\begin{align*}
|\overline{\eta}| \overline{\xi}_b^{\overline{\xi}} \overline{\sigma}_b^{\overline{\sigma}} &= O(\varphi) + O\left( \overline{\eta}_b^{\overline{\eta}} \right) + O\left( \xi_t^{\overline{\xi} \overline{\sigma}} \right)
\end{align*}

It only remains to compute 
\begin{align*}
&\left( m_b(\overline{\xi}) \xi_0 \widehat{X}_b(\overline{\xi}) \cdot \nabla_{\overline{\xi}} + P_b^b(\overline{\xi}, \overline{\sigma}) \frac{(3 \xi_0^2 - |\xi|^2 - 3 \sigma_0^2 + |\sigma|^2)|\xi|}{|\overline{\xi}|^3} |\xi_0| \widehat{X}_{b-\widehat{\mathcal{C}}}(\overline{\eta}, \overline{\sigma}) \cdot \nabla_{\overline{\eta}} \right) \varphi \\
&\quad = O\left( \overline{\xi}_b^{\overline{\sigma}} \overline{\sigma}_b^{\overline{\sigma}} \overline{\eta}_b^{\overline{\sigma}} \right) + O\left( \xi_t^{\overline{\xi} \overline{\sigma}} \right) + O\left( \overline{\eta}_b^{\overline{\eta}} \right)
\end{align*} 
which concludes. 

\paragraph{3.} Localise now $\overline{\eta}$ by $m_{\widehat{\mathcal{L}}}(\overline{\eta})$. 

We have
\begin{align*}
\varphi &= \xi_0^3 + \xi_0 |\xi|^2 - \eta_0^3 - \eta_0 |\eta|^2 - \sigma_0^3 - \sigma_0 |\sigma|^2 \\
&= o(\eta_0) + \eta_0 (3 \xi_0^2 + |\xi|^2) 
\end{align*}
so that $\overline{\eta} = O(\eta_0) = O(\varphi)$, which allows to directly decompose \eqref{equdecchampbHBHC-1}. 

Then, for \eqref{equdecchampbHBHC-2}, we apply the same strategy as in the case 2. above where $\overline{\eta}$ was close to $\widehat{\mathcal{C}}$: we decompose according to \eqref{decompositionequdecchampbHBHC-2} then project on $\widehat{X}_{b-\widehat{\mathcal{L}}}(\overline{\eta}, \overline{\sigma})$, keeping a factor $O\left( \overline{\xi}_b^{\overline{\xi}} - \overline{\sigma}_b^{\overline{\sigma}} \right) = O(\overline{\eta})$. We can then apply an integration by parts in $\overline{\eta}$ and get terms from \eqref{lemdecequHGtot} or of the form \eqref{equdecchampbHBHC-1} that can be decomposed the same way. 

\paragraph{4.} Localise finally $\overline{\eta}$ by $m_{\widehat{\mathcal{P}}}(\overline{\eta})$. 

We start by removing from \eqref{equdecchampbHBHC-2} the desired term \eqref{lemdecequH-02-symsigma}, then we apply an integration by parts in $\overline{\eta}$: 
\begin{subequations}
\begin{align}
&\begin{aligned}
\int_0^t \int e^{i s \varphi} \mu(\overline{\xi}, \overline{\eta}) \mu_{HBH}(\overline{\xi}, \overline{\eta}) m_{\widehat{\mathcal{C}}}(\overline{\sigma}) \left( m_b(\overline{\xi}) \xi_0 \widehat{X}_b(\overline{\xi}) - m_b(\overline{\sigma}) \sigma_0 \widehat{X}_b(\overline{\sigma}) \right) \cdot \nabla_{\overline{\xi}} \widehat{F}_2(s, \overline{\sigma}) \\
m_{\widehat{\mathcal{P}}}(\overline{\eta}) \widehat{F}_1(s, \overline{\eta}) ~ d\overline{\eta} ds 
\end{aligned} \notag \\
&\quad \begin{aligned}
= \int_0^t \int i s \left( m_b(\overline{\xi}) \xi_0 \widehat{X}_b(\overline{\xi}) - m_b(\overline{\sigma}) \sigma_0 \widehat{X}_b(\overline{\sigma}) \right) \cdot \nabla_{\overline{\eta}} \varphi e^{i s \varphi} \mu(\overline{\xi}, \overline{\eta}) \mu_{HBH}(\overline{\xi}, \overline{\eta}) m_{\widehat{\mathcal{P}}}(\overline{\eta}) \\
m_{\widehat{\mathcal{C}}}(\overline{\sigma}) \widehat{F}_1(s, \overline{\eta}) \widehat{F}_2(s, \overline{\sigma}) ~ d\overline{\eta} ds 
\end{aligned} \label{equdecchampbHBHPC-2-1} \\
&\quad \begin{aligned}
+ \int_0^t \int e^{i s \varphi} \mu(\overline{\xi}, \overline{\eta}) \mu_{HBH}(\overline{\xi}, \overline{\eta}) \left( m_b(\overline{\xi}) \xi_0 \widehat{X}_b(\overline{\xi}) - m_b(\overline{\sigma}) \sigma_0 \widehat{X}_b(\overline{\sigma}) \right) \cdot \nabla_{\overline{\eta}} \widehat{F}_1(s, \overline{\eta}) \\
m_{\widehat{\mathcal{C}}}(\overline{\sigma}) \widehat{F}_2(s, \overline{\sigma}) ~ d\overline{\eta} ds 
\end{aligned} \label{equdecchampbHBHPC-2-2} \\
&\quad \begin{aligned}
+ \int_0^t \int \nabla_{\overline{\eta}} \cdot \left( \left( m_b(\overline{\xi}) \xi_0 \widehat{X}_b(\overline{\xi}) - m_b(\overline{\sigma}) \sigma_0 \widehat{X}_b(\overline{\sigma}) \right) \mu(\overline{\xi}, \overline{\eta}) \mu_{HBH}(\overline{\xi}, \overline{\eta}) m_{\widehat{\mathcal{P}}}(\overline{\eta}) m_{\widehat{\mathcal{C}}}(\overline{\sigma}) \right) \\
e^{i s \varphi} \widehat{F}_1(s, \overline{\eta}) \widehat{F}_2(s, \overline{\sigma}) ~ d\overline{\eta} ds 
\end{aligned} \label{equdecchampbHBHPC-2-3} 
\end{align}
\end{subequations}
\eqref{equdecchampbHBHPC-2-3} is of the form \eqref{lemdecequH-11-dersymb}.

\eqref{equdecchampbHBHPC-2-2} is of the form $\eqref{lemdecequH-03-resxetaac} + \eqref{lemdecequH-04-resxetabbon} + \eqref{lemdecequH-05-resxetabbonbis} + \eqref{lemdecequH-06-resxetabphi}$ by projecting: 
\begin{align*} m_b(\overline{\xi}) \xi_0 \widehat{X}_b(\overline{\xi}) - m_b(\overline{\sigma}) \sigma_0 \widehat{X}_b(\overline{\sigma}) &= O(\eta_0) + O(|\overline{\eta}|) \widehat{X}_a(\overline{\eta}) + O(|\overline{\eta}|) \widehat{X}_c(\overline{\eta}) + O(|\xi| - |\sigma|) \widehat{X}_b(\overline{\eta}) 
\end{align*}
and using $\varphi = O(\eta_0) + \xi_0 (|\xi| - |\sigma|)$. 

Finally, we group \eqref{equdecchampbHBHPC-2-1} with \eqref{equdecchampbHBHC-1}. 

Note now that 
\begin{align*}
\widehat{X}_a(\overline{\sigma}) \cdot \nabla_{\overline{\eta}} \varphi &= \frac{\sigma_0}{|\overline{\sigma}|} (3 \sigma_0^2 + |\sigma|^2 - 3 \eta_0^2 - |\eta|^2) + \frac{\sigma}{|\overline{\sigma}|} \cdot (2 \sigma_0 \sigma - 2 \eta_0 \eta) \\
&= o(1) + \frac{\sigma_0}{|\overline{\sigma}|} (3 \sigma_0^2 + 3 |\sigma|^2) 
\end{align*}
In particular, in \eqref{equdecchampbHBHC-1}, we can decompose whenever we have a factor $O(\eta_0)$ applying an integration by parts along $\widehat{X}_a(\overline{\sigma})$. Then, 
\begin{align*}
\varphi &= \xi_0^3 + \xi_0 |\xi|^2 - \eta_0^3 - \eta_0 |\eta|^2 - \sigma_0^3 - \sigma_0 |\sigma|^2 \\
&= O(\eta_0) + \xi_0 (|\xi|^2 - |\sigma|^2) 
\end{align*}
so we can also decompose terms with a factor $O(|\xi| - |\sigma|)$. In particular, 
\begin{align*}
m_b(\overline{\xi}) - m_b(\overline{\sigma}) &= O(\eta_0) + O(|\xi| - |\sigma|) \\
|\overline{\xi}| - |\overline{\sigma}| &= O(\eta_0) + O(|\xi| - |\sigma|) 
\end{align*}

We then have that 
\begin{align*}
&m_b(\overline{\xi}) \xi_0 \widehat{X}_b(\overline{\xi}) \cdot (\nabla_{\overline{\xi}} + \nabla_{\overline{\eta}}) \varphi - m_b(\overline{\sigma}) \sigma_0 \widehat{X}_b(\overline{\sigma}) \cdot \nabla_{\overline{\eta}} \varphi \\
&\quad = O(\eta_0) + O(|\xi| - |\sigma|) + \frac{m_b(\overline{\xi})}{|\overline{\xi}|} \left( \xi_0 |\xi| (\xi_0^2 + |\xi|^2 - |\eta|^2) - \sigma_0 |\sigma| (\sigma_0^2 + |\sigma|^2 - |\eta|^2) \right) \\
&\quad = O(\eta_0) + O(|\xi| - |\sigma|)
\end{align*}
which concludes. 

\subsubsection{Interaction \texorpdfstring{$(HBH, \widehat{\mathcal{L}})$}{(HBH, L)}}

Let us consider
\begin{subequations}
\begin{align}
&\xi_0 m_b(\overline{\xi}) \widehat{X}_b(\overline{\xi}) \cdot \nabla_{\overline{\xi}} \widehat{I}_{\mu}^{(HBH, \widehat{\mathcal{L}})}[F_1, F_2](t, \overline{\xi}) \notag \\
&\quad = \int_0^t \int i s \xi_0 m_b(\overline{\xi}) \widehat{X}_b(\overline{\xi}) \cdot \nabla_{\overline{\xi}} \varphi e^{i s \varphi} \mu(\overline{\xi}, \overline{\eta}) \mu_{HBH}(\overline{\xi}, \overline{\eta}) m_{\widehat{\mathcal{L}}}(\overline{\sigma}) \widehat{F}_1(s, \overline{\eta}) \widehat{F}_2(s, \overline{\sigma}) ~ d\overline{\eta} ds \label{equdecchampbHBHL-1} \\
&\quad \quad + \int_0^t \int e^{i s \varphi} \xi_0 \mu(\overline{\xi}, \overline{\eta}) \mu_{HBH}(\overline{\xi}, \overline{\eta}) m_{\widehat{\mathcal{L}}}(\overline{\sigma}) \widehat{F}_1(s, \overline{\eta}) m_b(\overline{\xi}) \widehat{X}_b(\overline{\xi}) \cdot \nabla_{\overline{\xi}} \widehat{F}_2(s, \overline{\sigma}) ~ d\overline{\eta} ds \label{equdecchampbHBHL-2} \\
&\quad \quad + \int_0^t \int e^{i s \varphi} \xi_0 m_b(\overline{\xi}) \widehat{X}_b(\overline{\xi}) \cdot \nabla_{\overline{\xi}} \left( \mu(\overline{\xi}, \overline{\eta}) \mu_{HBH}(\overline{\xi}, \overline{\eta}) m_{\widehat{\mathcal{L}}}(\overline{\sigma}) \right) \widehat{F}_1(s, \overline{\eta}) \widehat{F}_2(s, \overline{\sigma}) ~ d\overline{\eta} ds \label{equdecchampbHBHL-3} 
\end{align}
\end{subequations}
\eqref{equdecchampbHBHL-3} is of the form \eqref{lemdecequH-11-dersymb}. 

We then localise $\overline{\eta}$ more finely: 
\[ 1 = m_{\widehat{\mathcal{R}}}(\overline{\eta}) + m_{\widehat{\mathcal{C}}}(\overline{\eta}) + m_{\widehat{\mathcal{L}}}(\overline{\eta}) + m_{\widehat{\mathcal{P}}}(\overline{\eta}) \]

\paragraph{1.} Localise first $\overline{\eta}$ by $m_{\widehat{\mathcal{R}}}$. 

On \eqref{equdecchampbHBHL-2}, we can apply the same decompositions as in the case $(HBH, \widehat{\mathcal{C}})$, removing \eqref{lemdecequH-02-symsigma} and then applying integrations by parts, to get down to the analysis of 
\begin{align*}
\int_0^t \int i s \left( m_b(\overline{\xi}) \xi_0 \widehat{X}_b(\overline{\xi}) - m_b(\overline{\sigma}) \sigma_0 \widehat{X}_b(\overline{\sigma}) \right) \cdot \nabla_{\overline{\eta}} \varphi e^{i s \varphi} \mu(\overline{\xi}, \overline{\eta}) \mu_{HBH}(\overline{\xi}, \overline{\eta}) m_{\widehat{\mathcal{L}}}(\overline{\sigma}) \widehat{F}_1(s, \overline{\eta}) \\
m_b(\overline{\xi}) \widehat{X}_b(\overline{\xi}) \cdot \nabla_{\overline{\xi}} \widehat{F}_2(s, \overline{\sigma}) ~ d\overline{\eta} ds
\end{align*} 
which we group with \eqref{equdecchampbHBHL-1}. 

We now compute that 
\begin{align*}
\partial_{\eta_0} \varphi &= 3 \sigma_0^2 + |\sigma|^2 - 3 \eta_0^2 - |\eta|^2 \\
&= 3 \sigma_0^2 + |\sigma|^2 + o(1)
\end{align*}
so any factor $O(\overline{\eta})$ allows a decomposition by integration by parts along $\widehat{X}_a(\overline{\sigma})$. 

But
\begin{align*}
&m_b(\overline{\xi}) \xi_0 \widehat{X}_b(\overline{\xi}) \cdot (\nabla_{\overline{\xi}} + \nabla_{\overline{\eta}}) \varphi - m_b(\overline{\sigma}) \sigma_0 \widehat{X}_b(\overline{\sigma}) \cdot \nabla_{\overline{\eta}} \varphi \\
&\quad = O(\overline{\eta}) 
\end{align*}
which concludes. 

\paragraph{2.} Localise now $\overline{\eta}$ by $m_{\widehat{\mathcal{C}}}$. 

For \eqref{equdecchampbHBHL-2}, we decompose
\begin{subequations}
\begin{align}
m_b(\overline{\xi}) \xi_0 \widehat{X}_b(\overline{\xi}) &= m_b(\overline{\sigma}) \sigma_0 \widehat{X}_b(\overline{\sigma}) \label{decompositionchampbHBHCL-1} \\
&+ O(\overline{\eta}) e_0 \label{decompositionchampbHBHCL-2} \\
&+ \frac{3 \xi_0^2 - |\xi|^2}{|\overline{\xi}|^4} \begin{pmatrix} 0 \\ - 2 \xi_0^2 \xi \end{pmatrix} - \frac{3 \sigma_0^2 - |\sigma|^2}{|\overline{\sigma}|^4} \begin{pmatrix} 0 \\ - 2 \sigma_0^2 \sigma \end{pmatrix} \label{decompositionchampbHBHCL-3} 
\end{align}
\end{subequations}
\eqref{decompositionchampbHBHCL-1} contributes like \eqref{lemdecequH-02-symsigma}, and \eqref{decompositionchampbHBHCL-2} like $\eqref{lemdecequH-07-resxsigmaac} + \eqref{lemdecequH-08-resxsigmabbon}$. For \eqref{decompositionchampbHBHCL-3}, we rewrite that 
\begin{align*}
\begin{pmatrix} 0 \\ \xi_0^2 \xi \end{pmatrix} &= \begin{pmatrix} 0 \\ \xi_0^2 \frac{J \eta \cdot \xi}{|\eta|} \frac{J \eta}{|\eta|} \end{pmatrix} + \begin{pmatrix} \xi_0^2 \frac{\eta \cdot \xi}{|\eta|} \frac{|\overline{\eta}|}{|\eta|} \frac{\eta_0}{|\overline{\eta}|} \\ \xi_0^2 \frac{\eta \cdot \xi}{|\eta|} \frac{|\overline{\eta}|}{|\eta|} \frac{\eta}{|\overline{\eta}|} \end{pmatrix} - \xi_0^2 \frac{\eta \cdot \xi}{|\eta|} \frac{|\overline{\eta}|}{|\eta|} \frac{\eta_0}{|\overline{\eta}|} e_0 \\
&= \xi_0^2 \frac{J \eta \cdot \xi}{|\eta|} \widehat{X}_c(\overline{\eta}) + \xi_0^2 \frac{\eta \cdot \xi}{|\eta|} \frac{|\overline{\eta}|}{|\eta|} \widehat{X}_a(\overline{\eta}) - \xi_0^2 \frac{\eta_0 \eta \cdot \xi}{|\eta|^2} e_0 
\end{align*}
Similarly, 
\begin{align*}
\begin{pmatrix} 0 \\ \sigma_0^2 \sigma \end{pmatrix} &= \sigma_0^2 \frac{J \eta \cdot \sigma}{|\eta|} \widehat{X}_c(\overline{\eta}) + \sigma_0^2 \frac{\eta \cdot \sigma}{|\eta|} \frac{|\overline{\eta}|}{|\eta|} \widehat{X}_a(\overline{\eta}) - \sigma_0^2 \frac{\eta_0 \eta \cdot \sigma}{|\eta|^2} e_0 
\end{align*}
Therefore, 
\begin{subequations}
\begin{align}
\eqref{decompositionchampbHBHCL-3} &= \left( \frac{3 \xi_0^2 - |\xi|^2}{|\overline{\xi}|^4} \xi_0^2 \frac{J \eta \cdot \xi}{|\eta|} - \frac{3 \sigma_0^2 - |\sigma|^2}{|\overline{\sigma}|^4} \sigma_0^2 \frac{J \eta \cdot \sigma}{|\eta|} \right) \widehat{X}_c(\overline{\eta}) \label{decompositionchampbHBHCL-3-1} \\
&+ \left( \frac{3 \xi_0^2 - |\xi|^2}{|\overline{\xi}|^4} \xi_0^2 \frac{\eta \cdot \xi}{|\eta|} \frac{|\overline{\eta}|}{|\eta|} - \frac{3 \sigma_0^2 - |\sigma|^2}{|\overline{\sigma}|^4} \sigma_0^2 \frac{\eta \cdot \sigma}{|\eta|} \frac{|\overline{\eta}|}{|\eta|} \right) \widehat{X}_a(\overline{\eta}) \label{decompositionchampbHBHCL-3-2} \\
&- \left( \frac{3 \xi_0^2 - |\xi|^2}{|\overline{\xi}|^4} \xi_0^2 \frac{\eta_0 \eta \cdot \xi}{|\eta|^2} - \frac{3 \sigma_0^2 - |\sigma|^2}{|\overline{\sigma}|^4} \sigma_0^2 \frac{\eta_0 \eta \cdot \sigma}{|\eta|^2} \right) e_0 \label{decompositionchampbHBHCL-3-3} 
\end{align}
\end{subequations}
\eqref{decompositionchampbHBHCL-3-3} contributes like $\eqref{lemdecequH-07-resxsigmaac} + \eqref{lemdecequH-08-resxsigmabbon}$: indeed, 
\begin{align*}
&\frac{3 \xi_0^2 - |\xi|^2}{|\overline{\xi}|^4} \xi_0^2 \frac{\eta_0 \eta \cdot \xi}{|\eta|^2} - \frac{3 \sigma_0^2 - |\sigma|^2}{|\overline{\sigma}|^4} \sigma_0^2 \frac{\eta_0 \eta \cdot \sigma}{|\eta|^2} = O(\overline{\eta})
\end{align*}
Then, on \eqref{decompositionchampbHBHCL-3-1} and \eqref{decompositionchampbHBHCL-3-2}, we apply an integration by parts in $\overline{\eta}$. Since
\begin{align*}
&\frac{3 \xi_0^2 - |\xi|^2}{|\overline{\xi}|^4} \xi_0^2 \frac{J \eta \cdot \xi}{|\eta|} - \frac{3 \sigma_0^2 - |\sigma|^2}{|\overline{\sigma}|^4} \sigma_0^2 \frac{J \eta \cdot \sigma}{|\eta|} = O(\overline{\eta}) \\
&\frac{3 \xi_0^2 - |\xi|^2}{|\overline{\xi}|^4} \xi_0^2 \frac{\eta \cdot \xi}{|\eta|} \frac{|\overline{\eta}|}{|\eta|} - \frac{3 \sigma_0^2 - |\sigma|^2}{|\overline{\sigma}|^4} \sigma_0^2 \frac{\eta \cdot \sigma}{|\eta|} \frac{|\overline{\eta}|}{|\eta|} = O(\overline{\eta})
\end{align*}
we deduce that we only obtain terms of the form \eqref{lemdecequH-03-resxetaac}, \eqref{lemdecequH-11-dersymb} and 
\begin{align*}
&\int_0^t \int i s O(\overline{\eta}) \nabla_{\overline{\eta}} \varphi e^{i s \varphi} \mu(\overline{\xi}, \overline{\eta}) \mu_{HBH}(\overline{\xi}, \overline{\eta}) m_{\widehat{\mathcal{L}}}(\overline{\sigma}) \widehat{F}_1(s, \overline{\eta}) \widehat{F}_2(s, \overline{\sigma}) ~ d\overline{\eta} ds
\end{align*}
that we group with \eqref{equdecchampbHBHL-1}. 

Note now that
\begin{align*}
\varphi &= \xi_0^3 + \xi_0 |\xi|^2 - \eta_0^3 - \eta_0 |\eta|^2 - \sigma_0^3 - \sigma_0 |\sigma|^2 \\
&= o(\overline{\eta}) + 3 \xi_0 \eta_0 \sigma_0
\end{align*}
In particular, $\overline{\eta} = O(\eta_0) = O(\varphi)$. 

But we already saw that, in the expressions of \eqref{decompositionchampbHBHCL-3-1} and \eqref{decompositionchampbHBHCL-3-2}, all remaining terms have a factor $O(\overline{\eta})$; while we always have $\nabla_{\overline{\xi}} \varphi = O(\overline{\eta})$. This concludes the decomposition in this case. 

\paragraph{3.} Localise now $\overline{\eta}$ by $m_{\widehat{\mathcal{L}}}$. 

For \eqref{equdecchampbHBHL-2}, we remove \eqref{lemdecequH-02-symsigma} and compute
\begin{align*}
&m_b(\overline{\xi}) \xi_0 \widehat{X}_b(\overline{\xi}) - m_b(\overline{\sigma}) \sigma_0 \widehat{X}_b(\overline{\sigma}) \\
&\quad = O(\overline{\eta}) e_0 + \begin{pmatrix} 0 \\ 
- \frac{3 \xi_0^2 - |\xi|^2}{|\overline{\xi}|^4} \xi_0^2 \xi + \frac{3 \sigma_0^2 - |\sigma|^2}{|\overline{\sigma}|^4} \sigma_0^2 \sigma \end{pmatrix} \\
&\quad = O(\overline{\eta}) e_0 + O(\eta) + \begin{pmatrix} 0 \\ 
- \frac{3 \xi_0^2 - |\sigma|^2}{|\overline{\xi}|^4} \xi_0^2 \xi + \frac{3 \sigma_0^2 - |\sigma|^2}{|\overline{\sigma}|^4} \sigma_0^2 \sigma \end{pmatrix} \\
&\quad = O(\overline{\eta}) e_0 + O(\eta) + \begin{pmatrix} 0 \\
\frac{1}{|\overline{\xi}|^4 |\overline{\sigma}|^4} \left( (3 \sigma_0^2 - |\sigma|^2) \sigma_0^2 |\overline{\xi}|^4 - (3 \xi_0^2 - |\sigma|^2) \xi_0^2 |\overline{\sigma}|^4 \right) \sigma \end{pmatrix}
\end{align*}
We may then simplify 
\begin{align*}
|\sigma|^2 \sigma (\sigma_0^2 |\overline{\xi}|^4 - \xi_0^2 |\overline{\sigma}|^4) &= O(|\sigma| \overline{\eta}) \\
\sigma_0^4 |\overline{\xi}|^4 - \xi_0^4 |\overline{\sigma}|^4 &= \left( \sigma_0^2 |\overline{\xi}|^2 - \xi_0^2 |\overline{\sigma}|^2 \right) \left( \sigma_0^2 |\overline{\xi}|^2 + \xi_0^2 |\overline{\sigma}|^2 \right) \\
&= \left( \sigma_0^2 |\xi|^2 - \xi_0^2 |\sigma|^2 \right) O(1) \\
&= O(\eta) + O(|\sigma| \overline{\eta})
\end{align*}
This means that 
\begin{align*}
&m_b(\overline{\xi}) \xi_0 \widehat{X}_b(\overline{\xi}) - m_b(\overline{\sigma}) \sigma_0 \widehat{X}_b(\overline{\sigma}) \\
&\quad = O(\overline{\eta}) e_0 + O(\eta) + O(|\sigma| |\overline{\eta}|) 
\end{align*}
The contributions of $O(\overline{\eta}) e_0$ is of the form $\eqref{lemdecequH-07-resxsigmaac} + \eqref{lemdecequH-08-resxsigmabbon}$. For the contribution of $O(\eta)$, we apply an integration by parts in $\overline{\eta}$ and get terms of the form \eqref{lemdecequH-03-resxetaac}, \eqref{lemdecequH-04-resxetabbon}, \eqref{lemdecequH-11-dersymb} and 
\begin{align}
\int_0^t \int i s O(\eta) e^{i s \varphi} \mu(\overline{\xi}, \overline{\eta}) m_{HBH}(\overline{\xi}, \overline{\eta}) m_{\widehat{\mathcal{L}}}(\overline{\eta}) m_{\widehat{\mathcal{L}}}(\overline{\sigma}) \widehat{F}_1(s, \overline{\eta}) \widehat{F}_2(s, \overline{\sigma}) ~ d\overline{\eta} ds \label{equdecchampbHBHLL-reste}
\end{align}

Note that 
\begin{align*}
\partial_{\eta_0} \varphi &= o(1) + 3 \sigma_0^2 + |\sigma|^2
\end{align*}
so that it is easy to decompose any term with a factor $O(\overline{\eta})$, hence also $O(\nabla_{\overline{\xi}} \varphi)$ which is the case of \eqref{equdecchampbHBHL-1}, and also \eqref{equdecchampbHBHLL-reste}. 

Finally, the contribution of $O(|\sigma| \overline{\eta})$ is of the form \eqref{lemdecequH-08-resxsigmabbon}. 

\paragraph{4.} Localise finally $\overline{\eta}$ by $m_{\widehat{\mathcal{P}}}$. 

We have that 
\begin{align*}
\partial_{\eta_0} \varphi &= o(1) + 3 \sigma_0^2 
\end{align*}
which allows to decompose any term having a factor $O(\eta_0)$; then, 
\begin{align*}
\varphi &= O(\eta_0) + \xi_0 (|\xi|^2 - |\sigma|^2) 
\end{align*}
which allows to decompose as well terms having a factor $O(|\xi| - |\sigma|)$. Note also that 
\begin{align*}
|\overline{\xi}| - |\overline{\sigma}| &= \frac{|\overline{\xi}|^2 - |\overline{\sigma}|^2}{|\overline{\xi}| + |\overline{\sigma}|} \\
&= O(\eta_0) + O(|\xi| - |\sigma|) 
\end{align*}

Then, we remove from \eqref{equdecchampbHBHL-2} the desired term \eqref{lemdecequH-02-symsigma}, and then apply an integration by parts and group with \eqref{equdecchampbHBHL-1}: 
\begin{subequations}
\begin{align}
&\eqref{equdecchampbHBHL-1} + \eqref{equdecchampbHBHL-2} - \eqref{lemdecequH-02-symsigma} \notag \\
&\begin{aligned} 
= \int_0^t \int i s \left( \xi_0 m_b(\overline{\xi}) \widehat{X}_b(\overline{\xi}) \cdot \left( \nabla_{\overline{\xi}} + \nabla_{\overline{\eta}} \right) - \sigma_0 m_b(\overline{\sigma}) \widehat{X}_b(\overline{\sigma}) \cdot \nabla_{\overline{\eta}} \right) \varphi e^{i s \varphi} \mu(\overline{\xi}, \overline{\eta}) \\
\mu_{HBH}(\overline{\xi}, \overline{\eta}) m_{\widehat{\mathcal{L}}}(\overline{\sigma}) m_{\widehat{\mathcal{P}}}(\overline{\eta}) \widehat{F}_1(s, \overline{\eta}) \widehat{F}_2(s, \overline{\sigma}) ~ d\overline{\eta} ds 
\end{aligned} \label{equdecchampbHBHPL-1} \\
&\begin{aligned}
+ \int_0^t \int e^{i s \varphi} \mu(\overline{\xi}, \overline{\eta}) \mu_{HBH}(\overline{\xi}, \overline{\eta}) \left( \xi_0 m_b(\overline{\xi}) \widehat{X}_b(\overline{\xi}) - \sigma_0 m_b(\overline{\sigma}) \widehat{X}_b(\overline{\sigma}) \right) \cdot \nabla_{\overline{\eta}} \widehat{F}_1(s, \overline{\eta}) \\
m_{\widehat{\mathcal{P}}}(\overline{\eta}) m_{\widehat{\mathcal{L}}}(\overline{\sigma}) \widehat{F}_2(s, \overline{\sigma}) ~ d\overline{\eta} ds 
\end{aligned} \label{equdecchampbHBHPL-2} \\
&\begin{aligned}
+ \int_0^t \int \nabla_{\overline{\eta}} \cdot \left( \left( \xi_0 m_b(\overline{\xi}) \widehat{X}_b(\overline{\xi}) - \sigma_0 m_b(\overline{\sigma}) \widehat{X}_b(\overline{\sigma}) \right) \mu(\overline{\xi}, \overline{\eta}) \mu_{HBH}(\overline{\xi}, \overline{\eta}) m_{\widehat{\mathcal{L}}}(\overline{\sigma}) m_{\widehat{\mathcal{P}}}(\overline{\eta}) \right) \\
e^{i s \varphi} \widehat{F}_1(s, \overline{\eta}) \widehat{F}_2(s, \overline{\sigma}) ~ d\overline{\eta} ds 
\end{aligned} \label{equdecchampbHBHPL-3} 
\end{align}
\end{subequations}
\eqref{equdecchampbHBHPL-3} is of the form \eqref{lemdecequH-11-dersymb}. 

We have that
\begin{align*}
&m_b(\overline{\xi}) \xi_0 \widehat{X}_b(\overline{\xi}) \cdot (\nabla_{\overline{\xi}} + \nabla_{\overline{\eta}}) \varphi - m_b(\overline{\sigma}) \sigma_0 \widehat{X}_b(\overline{\sigma}) \cdot \nabla_{\overline{\eta}} \varphi \\
&\quad = O(\eta_0) + \frac{m_b(\overline{\xi})}{|\overline{\xi}|} \xi_0 |\xi| (\xi_0^2 + |\xi|^2 - |\eta|^2) - \frac{m_b(\overline{\sigma})}{|\overline{\sigma}|} \sigma_0 |\sigma| (\sigma_0^2 + |\sigma|^2 - |\eta|^2) \\
&\quad = O(\eta_0) + O(|\xi| - |\sigma|) 
\end{align*}
so that we can decompose \eqref{equdecchampbHBHPL-1}. 

Finally, for \eqref{equdecchampbHBHPL-2}, we project: 
\begin{align*}
&\xi_0 m_b(\overline{\xi}) \widehat{X}_b(\overline{\xi}) - \sigma_0 m_b(\overline{\sigma}) \widehat{X}_b(\overline{\sigma}) \\
&= O(\eta_0) + O(|\overline{\eta}| \widehat{X}_a(\overline{\eta})) + O(|\overline{\eta}| \widehat{X}_c(\overline{\eta})) 
+ \xi_0 \left( m_b(\overline{\xi}) P_b^b(\overline{\xi}, \overline{\eta}) - m_b(\overline{\sigma}) P_b^b(\overline{\sigma}) \right) \widehat{X}_b(\overline{\eta}) \\
&= O(\eta_0) + O(|\overline{\eta}| \widehat{X}_a(\overline{\eta})) + O(|\overline{\eta}| \widehat{X}_c(\overline{\eta})) 
+ O(|\xi| - |\sigma|) \widehat{X}_b(\overline{\eta}) \\
&= O(\eta_0) + O((|\overline{\eta}| \widehat{X}_a(\overline{\eta})) + O(|\overline{\eta}| \widehat{X}_c(\overline{\eta})) 
+ O(\varphi \widehat{X}_b(\overline{\eta}))
\end{align*}
This is enough to rewrite \eqref{equdecchampbHBHPL-2} as a combination of \eqref{lemdecequH-03-resxetaac}, \eqref{lemdecequH-04-resxetabbon} and \eqref{lemdecequH-06-resxetabphi}. 

\subsubsection{Interaction \texorpdfstring{$(HBH, \widehat{\mathcal{P}})$}{(HBH, P)}}

Let us consider 
\begin{subequations}
\begin{align}
&\xi_0 m_b(\overline{\xi}) \widehat{X}_b(\overline{\xi}) \cdot \nabla_{\overline{\xi}} \widehat{I}_{\mu}^{(HBH, \widehat{\mathcal{P}})}[F_1, F_2](t, \overline{\xi}) \notag \\
&\quad = \int_0^t \int i s \xi_0 m_b(\overline{\xi}) \widehat{X}_b(\overline{\xi}) \cdot \nabla_{\overline{\xi}} \varphi e^{i s \varphi} \mu(\overline{\xi}, \overline{\eta}) \mu_{HBH}(\overline{\xi}, \overline{\eta}) m_{\widehat{\mathcal{P}}}(\overline{\sigma}) \widehat{F}_1(s, \overline{\eta}) \widehat{F}_2(s, \overline{\sigma}) ~ d\overline{\eta} ds \label{equdecchampbHBHP-1} \\
&\quad \quad + \int_0^t \int e^{i s \varphi} \xi_0 \mu(\overline{\xi}, \overline{\eta}) \mu_{HBH}(\overline{\xi}, \overline{\eta}) m_{\widehat{\mathcal{P}}}(\overline{\sigma}) \widehat{F}_1(s, \overline{\eta}) m_b(\overline{\xi}) \widehat{X}_b(\overline{\xi}) \cdot \nabla_{\overline{\xi}} \widehat{F}_2(s, \overline{\sigma}) ~ d\overline{\eta} ds \label{equdecchampbHBHP-2} \\
&\quad \quad + \int_0^t \int e^{i s \varphi} \xi_0 m_b(\overline{\xi}) \widehat{X}_b(\overline{\xi}) \cdot \nabla_{\overline{\xi}} \left( \mu(\overline{\xi}, \overline{\eta}) \mu_{HBH}(\overline{\xi}, \overline{\eta}) m_{\widehat{\mathcal{P}}}(\overline{\sigma}) \right) \widehat{F}_1(s, \overline{\eta}) \widehat{F}_2(s, \overline{\sigma}) ~ d\overline{\eta} ds \label{equdecchampbHBHP-3} 
\end{align}
\end{subequations}
\eqref{equdecchampbHBHP-3} is already of the form \eqref{lemdecequH-11-dersymb}. 

We now introduce symbols $\mu_{BH}, \mu_{HH}, \mu_{HB}$ as in the previous case $(BBB, \widehat{\mathcal{P}}\widehat{\mathcal{P}})$, that localise on $\{ |\eta_0| \ll |\sigma_0| \}$, $\{ |\eta_0| \simeq |\sigma_0| \}$ et $\{ |\eta_0| \gg |\sigma_0| \}$ respectively. 

\paragraph{1.} Localise first on $\{ |\eta_0| \simeq |\sigma_0| \}$. Then \eqref{lemdecequH-02-symsigma} is of the form \eqref{lemdecequH-09-resxsigmabbonbis} and we don't need to remove it explicitely. Furthermore, \eqref{equdecchampbHBHP-2} is also already of the form \eqref{lemdecequH-09-resxsigmabbonbis}. 

\paragraph{1.1} Let us first localise $\overline{\eta}$ away enough from $\widehat{\mathcal{P}}$. Then 
\begin{align*}
\eta_0 \widehat{X}_a(\overline{\eta}) \cdot \nabla_{\overline{\eta}} \varphi &= \frac{\eta_0^2}{|\overline{\eta}|} \left( 3 \sigma_0^2 + |\sigma|^2 - 3 \eta_0^2 - |\eta|^2 \right) + \frac{\eta_0 \eta}{|\overline{\eta}|} \cdot (2 \sigma_0 \sigma - 2 \eta_0 \eta) \\
&= \frac{\eta_0}{|\overline{\eta}|} \left( \eta_0 |\sigma|^2 + o(\eta_0) \right) 
\end{align*}
so that we can decompose if we have a factor $O(\eta_0)$. But
\begin{align*}
\nabla_{\overline{\xi}} \varphi &= O(\overline{\eta}) = O(\eta_0)
\end{align*}
and we can apply the necessary integrations by parts for \eqref{equdecchampbHBHP-1}. 

On the other hand, to get \eqref{lemdecequH-01-symeta}, we can apply an integration by parts on this term and get terms of the form \eqref{lemdecequH-07-resxsigmaac}, \eqref{lemdecequH-09-resxsigmabbonbis}, \eqref{lemdecequH-11-dersymb} or a term similar to \eqref{equdecchampbHBHP-1}, that we can treat the same way. 

\paragraph{1.2} Localise now $\overline{\eta}$ in the neighborhood of $\widehat{\mathcal{P}}$. Then
\begin{align*}
\partial_{\eta_0} \varphi &= 3 \sigma_0^2 + |\sigma|^2 - 3 \eta_0^2 - |\eta|^2 \\
&\simeq |\sigma|^2
\end{align*}
and so we can again decompose if we have a factor $O(\eta_0)$. The rest of the argument is identical. 

\paragraph{2.} Localise now on $\{ |\eta_0| \gg |\sigma_0| \}$. Then \eqref{lemdecequH-02-symsigma} is of the form \eqref{lemdecequH-09-resxsigmabbonbis}; moreover, \eqref{equdecchampbHBHP-2} is of the form \eqref{lemdecequH-09-resxsigmabbonbis}. We then remove \eqref{lemdecequH-01-symeta} and apply an integration by parts: 
\begin{subequations}
\begin{align}
&- \int_0^t \int e^{i s \varphi} \mu(\overline{\xi}, \overline{\eta}) \mu_{HBH}(\overline{\xi}, \overline{\eta}) \mu_{HB}(\xi_0, \eta_0) m_{\widehat{\mathcal{P}}}(\overline{\sigma}) m_b(\overline{\eta}) \eta_0 \widehat{X}_b(\overline{\eta}) \cdot \nabla_{\overline{\eta}} \widehat{F}_1(s, \overline{\eta}) \widehat{F}_2(s, \overline{\sigma}) ~ d\overline{\eta} ds \notag \\
&\quad = \int_0^t \int i s m_b(\overline{\eta}) \eta_0 \widehat{X}_b(\overline{\eta}) \cdot \nabla_{\overline{\eta}} \varphi e^{i s \varphi} \mu(\overline{\xi}, \overline{\eta}) \mu_{HBH}(\overline{\xi}, \overline{\eta}) \mu_{HB}(\xi_0, \eta_0) m_{\widehat{\mathcal{P}}}(\overline{\sigma}) \widehat{F}_1(s, \overline{\eta}) \widehat{F}_2(s, \overline{\sigma}) ~ d\overline{\eta} ds \label{equdecchampbHBHP-HB-2-1} \\
&\quad \quad \int_0^t \int e^{i s \varphi} \mu(\overline{\xi}, \overline{\eta}) \mu_{HBH}(\overline{\xi}, \overline{\eta}) \mu_{HB}(\xi_0, \eta_0) m_{\widehat{\mathcal{P}}}(\overline{\sigma}) \widehat{F}_1(s, \overline{\eta}) m_b(\overline{\eta}) \eta_0 \widehat{X}_b(\overline{\eta}) \cdot \nabla_{\overline{\eta}} \widehat{F}_2(s, \overline{\sigma}) ~ d\overline{\eta} ds \label{equdecchampbHBHP-HB-2-2} \\
&\quad \quad \int_0^t \int e^{i s \varphi} \nabla_{\overline{\eta}} \cdot \left( m_b(\overline{\eta}) \eta_0 \widehat{X}_b(\overline{\eta}) \mu(\overline{\xi}, \overline{\eta}) \mu_{HBH}(\overline{\xi}, \overline{\eta}) \mu_{BH}(\xi_0, \eta_0) m_{\widehat{\mathcal{P}}}(\overline{\sigma}) \right) \widehat{F}_1(s, \overline{\eta}) \widehat{F}_2(s, \overline{\sigma}) ~ d\overline{\eta} ds \label{equdecchampbHBHP-HB-2-3} 
\end{align}
\end{subequations}
\eqref{equdecchampbHBHP-HB-2-2} is of the form \eqref{lemdecequH-09-resxsigmabbonbis}, and \eqref{equdecchampbHBHP-HB-2-3} of the form \eqref{lemdecequH-11-dersymb}. Then, we group \eqref{equdecchampbHBHP-1} and \eqref{equdecchampbHBHP-HB-2-1}: 
\begin{align*}
\eqref{equdecchampbHBHP-1} + \eqref{equdecchampbHBHP-HB-2-1} = \int_0^t \int i s \left( m_b(\overline{\xi}) \xi_0 \widehat{X}_b(\overline{\xi}) \cdot \nabla_{\overline{\xi}} + m_b(\overline{\eta}) \eta_0 \widehat{X}_b(\overline{\eta}) \cdot \nabla_{\overline{\eta}} \right) \varphi \\e^{i s \varphi} \mu(\overline{\xi}, \overline{\eta}) \mu_{HBH}(\overline{\xi}, \overline{\eta}) \mu_{HB}(\xi_0, \eta_0) m_{\widehat{\mathcal{P}}}(\overline{\sigma}) \widehat{F}_1(s, \overline{\eta}) \widehat{F}_2(s, \overline{\sigma}) ~ d\overline{\eta} ds
\end{align*}

\paragraph{2.1} Localise first $\overline{\eta}$ away enough from $\widehat{\mathcal{P}}$. We then compute that 
\begin{align*}
\widehat{X}_a(\overline{\eta}) \cdot \nabla_{\overline{\eta}} \varphi &= \frac{\eta_0}{|\overline{\eta}|} (3 \sigma_0^2 + |\sigma|^2 - 3 \eta_0^2 - |\eta|^2) + \frac{\eta}{|\overline{\eta}|} \cdot (2 \sigma_0 \sigma - 2 \eta_0 \eta) \\
&= \frac{\eta_0}{|\overline{\eta}|} |\sigma|^2 + o(1) 
\end{align*}
which allows to decompose whenever we have a factor $O(\sigma_0)$. On the other hand, 
\begin{align*}
\varphi &= \xi_0^3 + \xi_0 |\xi|^2 - \eta_0^3 - \eta_0 |\eta|^2 - \sigma_0^3 - \sigma_0 |\sigma|^2 \\
&= O(\sigma_0) + \xi_0 (|\xi|^2 - |\eta|^2) 
\end{align*}
and $|\xi|^2 - |\eta|^2 \simeq |\xi|^2$, so we can also decompose if we have a factor $O(\xi_0)$. Finally, it is clear that 
\begin{align*}
&\left( m_b(\overline{\xi}) \xi_0 \widehat{X}_b(\overline{\xi}) \cdot \nabla_{\overline{\xi}} + m_b(\overline{\eta}) \eta_0 \widehat{X}_b(\overline{\eta}) \cdot \nabla_{\overline{\eta}} \right) \varphi \\
&\quad = O(\xi_0) + O(\eta_0) = O(\xi_0) + O(\sigma_0)
\end{align*}
which concludes. 

\paragraph{2.2} Localise now $\overline{\eta}$ near $\widehat{\mathcal{P}}$. We compute that 
\begin{align*}
\partial_{\eta_0} \varphi &= 3 \sigma_0^2 + |\sigma|^2 - 3 \eta_0^2 - |\eta|^2 \\
&\simeq |\sigma|^2
\end{align*}
so again we can decompose if we have a factor $O(\sigma_0)$. The rest of the argument is identical. 

\paragraph{3.} Localise finally on $\{ |\eta_0| \ll |\sigma_0| \}$. 

In this case, we need to remove \eqref{lemdecequH-02-symsigma} from \eqref{equdecchampbHBHP-2}: 
\begin{align*}
\eqref{equdecchampbHBHP-2} - \eqref{lemdecequH-02-symsigma} = \int_0^t \int e^{i s \varphi} \mu(\overline{\xi}, \overline{\eta}) \mu_{HBH}(\overline{\xi}, \overline{\eta}) \mu_{BH}(\xi_0, \eta_0) m_{\widehat{\mathcal{P}}}(\overline{\sigma}) \widehat{F}_1(s, \overline{\eta}) \\
\left( m_b(\overline{\xi}) \xi_0 \widehat{X}_b(\overline{\xi}) - m_b(\overline{\sigma}) \sigma_0 \widehat{X}_b(\overline{\sigma}) \right) \cdot \nabla_{\overline{\xi}} \widehat{F}_2(s, \overline{\sigma}) ~ d\overline{\eta} ds
\end{align*}
Here above, we have clearly that 
\begin{align*}
m_b(\overline{\xi}) \xi_0 \widehat{X}_b(\overline{\xi}) - m_b(\overline{\sigma}) \sigma_0 \widehat{X}_b(\overline{\sigma}) &= O(\eta_0) + \xi_0 \left( m_b(\overline{\xi}) \widehat{X}_b(\overline{\xi}) - m_b(\overline{\sigma}) \widehat{X}_b(\overline{\sigma}) \right) 
\end{align*}
The term with a factor $O(\eta_0)$ in $\eqref{equdecchampbHBHP-2} - \eqref{lemdecequH-02-symsigma}$ is of the form \eqref{lemdecequH-09-resxsigmabbonbis}, so we only keep the second one. 

We also need to get \eqref{lemdecequH-01-symeta}. But, up to terms of the form \eqref{lemdecequH-09-resxsigmabbonbis} or \eqref{lemdecequH-11-dersymb}, we can rewrite after an integration by parts \eqref{lemdecequH-01-symeta} as 
\begin{align}
- \int_0^t \int i s \eta_0 m_b(\overline{\eta}) \widehat{X}_b(\overline{\eta}) \cdot \nabla_{\overline{\eta}} \varphi e^{i s \varphi} \mu(\overline{\xi}, \overline{\eta}) \mu_{HBH}(\overline{\xi}, \overline{\eta}) \mu_{HBH}(\xi_0, \eta_0) m_{\widehat{\mathcal{P}}}(\overline{\sigma}) \widehat{F}_1(s, \overline{\eta}) \widehat{F}_2(s, \overline{\sigma}) ~ d\overline{\eta} ds \label{equdecchampbHBHP-1-bissym-cas3} 
\end{align}

\paragraph{3.1} Localise first $\overline{\eta}$ away enough from $\widehat{\mathcal{P}}$. Then 
\begin{align*}
m_b(\overline{\xi}) \widehat{X}_b(\overline{\xi}) - m_b(\overline{\sigma}) \widehat{X}_b(\overline{\sigma}) &= O(\overline{\eta}) = O(\eta_0)
\end{align*}
Therefore, we can rewrite the remaining term in $\eqref{equdecchampbHBHP-2} - \eqref{lemdecequH-02-symsigma}$ in the form \eqref{lemdecequH-09-resxsigmabbonbis}. On the other hand, 
\begin{align*}
\widehat{X}_a(\overline{\eta}) \cdot \nabla_{\overline{\eta}} \varphi &= \frac{\eta_0}{|\overline{\eta}|} |\sigma|^2 + o(1) 
\end{align*}
so we know how to decompose having a factor $O(\eta_0)$, hence $O(\overline{\eta})$. But 
\begin{align*}
\nabla_{\overline{\xi}} \varphi &= O(\overline{\eta})
\end{align*}
which is enough to decompose \eqref{equdecchampbHBHP-1}. 

Finally, for \eqref{equdecchampbHBHP-1-bissym-cas3}, it is of the same form as \eqref{equdecchampbHBHP-1} and can be decomposed the same way. 

\paragraph{3.2} Localise now $\overline{\eta}$ near $\widehat{\mathcal{P}}$. 

This time, we apply an integration by parts in $\overline{\eta}$ on $\eqref{equdecchampbHBHP-2} - \eqref{lemdecequH-02-symsigma}$: 
\begin{subequations}
\begin{align}
&\eqref{equdecchampbHBHP-2} - \eqref{lemdecequH-02-symsigma} \notag \\
&\begin{aligned}
= \int_0^t \int i s \xi_0 \left( m_b(\overline{\xi}) \widehat{X}_b(\overline{\xi}) - m_b(\overline{\sigma}) \widehat{X}_b(\overline{\sigma}) \right) \cdot \nabla_{\overline{\eta}} \varphi e^{i s \varphi} \mu(\overline{\xi}, \overline{\eta}) \mu_{HBH}(\overline{\xi}, \overline{\eta}) \mu_{BH}(\xi_0, \eta_0) \\
m_{\widehat{\mathcal{P}}}(\overline{\sigma}) m_{\widehat{\mathcal{P}}}(\overline{\eta}) \widehat{F}_1(s, \overline{\eta}) \widehat{F}_2(s, \overline{\sigma}) ~ d\overline{\eta} ds 
\end{aligned} \label{equdecchampbHBHP-BH-2-1} \\
&\begin{aligned}
+ \int_0^t \int e^{i s \varphi} \mu(\overline{\xi}, \overline{\eta}) \mu_{HBH}(\overline{\xi}, \overline{\eta}) \xi_0 \left( m_b(\overline{\xi}) \widehat{X}_b(\overline{\xi}) - m_b(\overline{\sigma}) \widehat{X}_b(\overline{\sigma}) \right) \cdot \nabla_{\overline{\eta}} \widehat{F}_1(s, \overline{\eta}) \\
\mu_{BH}(\xi_0, \eta_0) m_{\widehat{\mathcal{P}}}(\overline{\sigma}) m_{\widehat{\mathcal{P}}}(\overline{\eta}) \widehat{F}_2(s, \overline{\sigma}) ~ d\overline{\eta} ds 
\end{aligned} \label{equdecchampbHBHP-BH-2-2} \\
&
+ \int_0^t \int \nabla_{\overline{\eta}} \cdot \left( \xi_0 \left( m_b(\overline{\xi}) \widehat{X}_b(\overline{\xi}) - m_b(\overline{\sigma}) \widehat{X}_b(\overline{\sigma}) \right) \mu(\overline{\xi}, \overline{\eta}) \mu_{HBH}(\overline{\xi}, \overline{\eta}) \mu_{BH}(\xi_0, \eta_0) m_{\widehat{\mathcal{P}}}(\overline{\sigma}) m_{\widehat{\mathcal{P}}}(\overline{\eta}) \right) \notag \\
&\quad \quad \quad \quad \quad \quad \quad \quad \quad \quad \quad \quad \quad \quad \quad \quad \quad \quad \quad \quad \quad \quad \quad \quad e^{i s \varphi} \widehat{F}_1(s, \overline{\eta}) \widehat{F}_2(s, \overline{\sigma}) ~ d\overline{\eta} ds \label{equdecchampbHBHP-BH-2-3} 
\end{align}
\end{subequations}
\eqref{equdecchampbHBHP-BH-2-3} is of the form \eqref{lemdecequH-11-dersymb}. 

For \eqref{equdecchampbHBHP-BH-2-2}, we project: 
\begin{align*}
\xi_0 \left( m_b(\overline{\xi}) \widehat{X}_b(\overline{\xi}) - m_b(\overline{\sigma}) \widehat{X}_b(\overline{\sigma}) \right) 
&= O(\eta_0) + O(|\overline{\eta}| \widehat{X}_a(\overline{\eta})) + O(|\overline{\eta}| \widehat{X}_c(\overline{\eta})) 
+ O(\xi_0(|\xi| - |\sigma|) \widehat{X}_b(\overline{\eta})) 
\end{align*}
as we already did previously. But, in the same way, $\varphi = \xi_0 (|\xi|^2 - |\sigma|^2) + O(\eta_0)$, and so we can rewrite
\begin{align*}
\xi_0 \left( m_b(\overline{\xi}) \widehat{X}_b(\overline{\xi}) - m_b(\overline{\sigma}) \widehat{X}_b(\overline{\sigma}) \right) 
&= O(\eta_0 \widehat{X}_b(\overline{\eta})) + O(|\overline{\eta}| \widehat{X}_a(\overline{\eta})) + O(|\overline{\eta}| \widehat{X}_c(\overline{\eta})) 
+ O(\varphi \widehat{X}_b(\overline{\eta})) 
\end{align*}
The factor $O(\eta_0)$ contributes like \eqref{lemdecequH-05-resxetabbonbis}; those in $O(|\overline{\eta}| \widehat{X}_a(\overline{\eta}))$ and $O(|\overline{\eta}| \widehat{X}_c(\overline{\eta}))$ like \eqref{lemdecequH-03-resxetaac}; and the one in $O(\varphi \widehat{X}_b(\overline{\eta}))$ like \eqref{lemdecequH-06-resxetabphi}. 

Finally, we group \eqref{equdecchampbHBHP-BH-2-1} with \eqref{equdecchampbHBHP-1} and \eqref{equdecchampbHBHP-1-bissym-cas3}: 
\begin{align*}
&\eqref{equdecchampbHBHP-1} + \eqref{equdecchampbHBHP-BH-2-1} + \eqref{equdecchampbHBHP-1-bissym-cas3} \\
&\quad = \int_0^t \int i s \left( m_b(\overline{\xi}) \xi_0 \widehat{X}_b(\overline{\xi}) \cdot (\nabla_{\overline{\xi}} + \nabla_{\overline{\eta}}) - m_b(\overline{\sigma}) \xi_0 \widehat{X}_b(\overline{\sigma}) \cdot \nabla_{\overline{\eta}} - m_b(\overline{\eta}) \eta_0 \widehat{X}_b(\overline{\eta}) \cdot \nabla_{\overline{\eta}} \right) \varphi e^{i s \varphi} \\
&\quad \quad \quad \quad \quad \quad \frac{\xi_0}{|\overline{\xi}|} \mu(\overline{\xi}, \overline{\eta}) \mu_{HBH}(\overline{\xi}, \overline{\eta}) \mu_{BH}(\xi_0, \eta_0) m_{\widehat{\mathcal{P}}}(\overline{\sigma}) m_{\widehat{\mathcal{P}}}(\overline{\eta}) \widehat{F}_1(s, \overline{\eta}) \widehat{F}_2(s, \overline{\sigma}) ~ d\overline{\eta} ds
\end{align*}

We now compute that 
\begin{align*}
\partial_{\eta_0} \varphi &= |\sigma|^2 + o(1) 
\end{align*}
so we know how to decompose if we have a factor $O(\eta_0)$, and then we already saw that $\varphi = \xi_0 (|\xi|^2 - |\sigma|^2) + O(\eta_0)$, so we can also decompose if we have $O(\xi_0 (|\xi| - |\sigma|))$. It is then easy to check that 
\begin{align*}
\xi_0 (m_b(\overline{\xi}) - m_b(\overline{\sigma})) &= O(\eta_0) + O(\xi_0 (|\xi| - |\sigma|)) 
\end{align*}
therefore
\begin{align*}
&\left( m_b(\overline{\xi}) \xi_0 \widehat{X}_b(\overline{\xi}) \cdot (\nabla_{\overline{\xi}} + \nabla_{\overline{\eta}}) - m_b(\overline{\sigma}) \xi_0 \widehat{X}_b(\overline{\sigma}) \cdot \nabla_{\overline{\eta}} - \eta_0 m_b(\overline{\eta}) \widehat{X}_b(\overline{\eta}) \cdot \nabla_{\overline{\eta}} \right) \varphi \\
&\quad = O(\eta_0) + O(\xi_0 (|\xi| - |\sigma|)) + \xi_0 m_b(\overline{\xi}) \Biggl( |\xi| (3 \xi_0^2 + |\xi|^2 - 3 \eta_0^2 - |\eta|^2) - \xi_0 \frac{\xi}{|\xi|} \cdot (2 \xi_0 \xi - 2 \eta_0 \eta) \\
&\quad \quad \quad - |\sigma| (3 \sigma_0^2 + |\sigma|^2 - 3 \eta_0^2 - |\eta|^2) + \sigma_0 \frac{\sigma}{|\sigma|} \cdot (2 \sigma_0 \sigma - 2 \eta_0 \eta) \Biggl) \\
&\quad = O(\eta_0) + O(\xi_0 (|\xi| - |\sigma|))
\end{align*}
which is enough. 

\subsubsection{Interaction \texorpdfstring{$(BHH, \widehat{\mathcal{R}})$}{(BHH, R)}}

Let us consider
\begin{subequations}
\begin{align}
&\xi_0 m_b(\overline{\xi}) \widehat{X}_b(\overline{\xi}) \cdot \nabla_{\overline{\xi}} \widehat{I}_{\mu}^{(BHH, \widehat{\mathcal{R}})}[F_1, F_2](t, \overline{\xi}) \notag \\
&\quad = \int_0^t \int i s \xi_0 m_b(\overline{\xi}) \widehat{X}_b(\overline{\xi}) \cdot \nabla_{\overline{\xi}} \varphi e^{i s \varphi} \mu(\overline{\xi}, \overline{\eta}) \mu_{BHH}(\overline{\xi}, \overline{\eta}) \mu_{\widehat{\mathcal{R}}}(\overline{\xi}, \overline{\eta}) \widehat{F}_1(s, \overline{\eta}) \widehat{F}_2(s, \overline{\sigma}) ~ d\overline{\eta} ds \label{equdecchampbBHHR-1} \\
&\quad \quad + \int_0^t \int e^{i s \varphi} \xi_0 \mu(\overline{\xi}, \overline{\eta}) \mu_{BHH}(\overline{\xi}, \overline{\eta}) \mu_{\widehat{\mathcal{R}}}(\overline{\xi}, \overline{\eta}) \widehat{F}_1(s, \overline{\eta}) m_b(\overline{\xi}) \widehat{X}_b(\overline{\xi}) \cdot \nabla_{\overline{\xi}} \widehat{F}_2(s, \overline{\sigma}) ~ d\overline{\eta} ds \label{equdecchampbBHHR-2} \\
&\quad \quad + \int_0^t \int e^{i s \varphi} \xi_0 m_b(\overline{\xi}) \widehat{X}_b(\overline{\xi}) \cdot \nabla_{\overline{\xi}} \left( \mu(\overline{\xi}, \overline{\eta}) \mu_{BHH}(\overline{\xi}, \overline{\eta}) \mu_{\widehat{\mathcal{R}}}(\overline{\xi}, \overline{\eta}) \right) \widehat{F}_1(s, \overline{\eta}) \widehat{F}_2(s, \overline{\sigma}) ~ d\overline{\eta} ds \label{equdecchampbBHHR-3} 
\end{align}
\end{subequations}
\eqref{equdecchampbBHHR-3} is of the form \eqref{lemdecequH-11-dersymb}, and \eqref{equdecchampbBHHR-2} of the form $\eqref{lemdecequH-07-resxsigmaac} + \eqref{lemdecequH-08-resxsigmabbon}$. Finally, for \eqref{equdecchampbBHHR-1}, we apply Lemma \ref{lem-non-res-loin0Cone} and get that
\[ \xi_0 = O(\nabla_{\overline{\eta}} \varphi) + O(\varphi) \]
which is enough to apply integrations by parts on \eqref{equdecchampbBHHR-1} and get terms from \eqref{lemdecequHGtot}. We didn't use the presence of $m_b(\overline{\xi})$, so \eqref{decompositionfinemg-gainkout} holds. 

\subsubsection{Interaction \texorpdfstring{$(BHH, \widehat{\mathcal{C}})$}{(BHH, C)}} 

Let us consider
\begin{subequations}
\begin{align}
&\xi_0 m_b(\overline{\xi}) \widehat{X}_b(\overline{\xi}) \cdot \nabla_{\overline{\xi}} \widehat{I}_{\mu}^{(BHH, \widehat{\mathcal{C}})}[F_1, F_2](t, \overline{\xi}) \notag \\
&\quad = \int_0^t \int i s \xi_0 m_b(\overline{\xi}) \widehat{X}_b(\overline{\xi}) \cdot \nabla_{\overline{\xi}} \varphi e^{i s \varphi} \mu(\overline{\xi}, \overline{\eta}) \mu_{BHH}(\overline{\xi}, \overline{\eta}) m_{\widehat{\mathcal{C}}}(\overline{\eta}) m_{\widehat{\mathcal{C}}}(\overline{\sigma}) \widehat{F}_1(s, \overline{\eta}) \widehat{F}_2(s, \overline{\sigma}) ~ d\overline{\eta} ds \label{equdecchampbBHHC-1} \\
&\quad \quad + \int_0^t \int e^{i s \varphi} \xi_0 \mu(\overline{\xi}, \overline{\eta}) \mu_{BHH}(\overline{\xi}, \overline{\eta}) m_{\widehat{\mathcal{C}}}(\overline{\eta}) m_{\widehat{\mathcal{C}}}(\overline{\sigma}) \widehat{F}_1(s, \overline{\eta}) m_b(\overline{\xi}) \widehat{X}_b(\overline{\xi}) \cdot \nabla_{\overline{\xi}} \widehat{F}_2(s, \overline{\sigma}) ~ d\overline{\eta} ds \label{equdecchampbBHHC-2} \\
&\quad \quad + \int_0^t \int e^{i s \varphi} \xi_0 m_b(\overline{\xi}) \widehat{X}_b(\overline{\xi}) \cdot \nabla_{\overline{\xi}} \left( \mu(\overline{\xi}, \overline{\eta}) \mu_{BHH}(\overline{\xi}, \overline{\eta}) m_{\widehat{\mathcal{C}}}(\overline{\eta}) m_{\widehat{\mathcal{C}}}(\overline{\sigma}) \right) \widehat{F}_1(s, \overline{\eta}) \widehat{F}_2(s, \overline{\sigma}) ~ d\overline{\eta} ds \label{equdecchampbBHHC-3} 
\end{align}
\end{subequations}
\eqref{equdecchampbBHHC-3} is of the form \eqref{lemdecequH-11-dersymb}. 

For \eqref{equdecchampbBHHC-2}, we can project $\widehat{X}_b(\overline{\xi})$ on the basis $(\widehat{X}_a(\overline{\sigma}), \widehat{X}_b(\overline{\sigma}), \widehat{X}_c(\overline{\sigma}))$ and all contributions except the one of $\widehat{X}_b(\overline{\sigma})$ are of the form \eqref{lemdecequH-07-resxsigmaac}. Likewise, up to terms of the form \eqref{lemdecequH-07-resxsigmaac}, we can replace $\widehat{X}_b$ by $\frac{\sigma_0}{|\sigma_0|} \widehat{X}_b'$. We therefore replace \eqref{equdecchampbBHHC-2} by
\begin{align}
\int_0^t \int e^{i s \varphi} \xi_0 \mu(\overline{\xi}, \overline{\eta}) \mu_{BHH}(\overline{\xi}, \overline{\eta}) m_{\widehat{\mathcal{C}}}(\overline{\eta}) m_{\widehat{\mathcal{C}}}(\overline{\sigma}) \widehat{F}_1(s, \overline{\eta}) m_b(\overline{\xi}) \frac{\sigma_0}{|\sigma_0|} P_b^b(\overline{\xi}, \overline{\sigma}) \widehat{X}_b'(\overline{\sigma}) \cdot \nabla_{\overline{\xi}} \widehat{F}_2(s, \overline{\sigma}) ~ d\overline{\eta} ds \label{equdecchampbBHHC-2-bis} 
\end{align}

The aim is now to apply an integration by parts along $\widehat{X}_{b-\widehat{\mathcal{C}}}(\overline{\eta}, \overline{\sigma})$, but unfortunately this vector field is singular in this situation. In order to compensate the singularity, we write that 
\begin{align*}
\overline{\xi}_b^{\overline{\eta}} &= \overline{\eta}_b^{\overline{\eta}} + \overline{\sigma}_b^{\overline{\eta}} \\
&= \overline{\eta}_b^{\overline{\eta}} - \overline{\sigma}_b^{\overline{\sigma}} - |\sigma| - \frac{\eta \cdot \sigma}{|\eta|} \\
&= \overline{\eta}_b^{\overline{\eta}} - \overline{\sigma}_b^{\overline{\sigma}} - |\sigma| \left( 1 + \theta^{\overline{\eta} \overline{\sigma}} \right) 
\end{align*}
Yet
\begin{align*}
3 \xi_0^2 - |\xi|^2 &= 3 \xi_0^2 - \left( \frac{\eta \cdot \xi}{|\eta|} \right)^2 - \left( \frac{J \eta \cdot \xi}{|\eta|} \right)^2 \\
&= \overline{\xi}_a^{\overline{\eta}} \overline{\xi}_b^{\overline{\eta}} - |\sigma|^2 \left( \xi_t^{\overline{\sigma} \overline{\eta}} \right)^2 \\
&= \overline{\xi}_a^{\overline{\eta}} \left( \overline{\eta}_b^{\overline{\eta}} - \overline{\sigma}_b^{\overline{\sigma}} - |\sigma| \left( \theta^{\overline{\eta} \overline{\sigma}} + 1 \right) \right) - |\sigma|^2 \left( \xi_t^{\overline{\sigma} \overline{\eta}} \right)^2
\end{align*}
We can then separate into: 
\begin{subequations}
\begin{align}
&\eqref{equdecchampbBHHC-2-bis} \notag \\
&\begin{aligned}
= \int_0^t \int e^{i s \varphi} \mu(\overline{\xi}, \overline{\eta}) \mu_{BHH}(\overline{\xi}, \overline{\eta}) m_{\widehat{\mathcal{C}}}(\overline{\eta}) m_{\widehat{\mathcal{C}}}(\overline{\sigma}) \widehat{F}_1(s, \overline{\eta}) \frac{\xi_0 |\xi| \overline{\xi}_a^{\overline{\eta}}}{|\overline{\xi}|^3} \overline{\eta}_b^{\overline{\eta}} \frac{\sigma_0}{|\sigma_0|} P_b^b(\overline{\xi}, \overline{\sigma}) \\
\widehat{X}_b'(\overline{\sigma}) \cdot \nabla_{\overline{\xi}} \widehat{F}_2(s, \overline{\sigma}) ~ d\overline{\eta} ds 
\end{aligned} \label{equdecchampbBHHC-2-1} \\
&\begin{aligned}
- \int_0^t \int e^{i s \varphi} \mu(\overline{\xi}, \overline{\eta}) \mu_{BHH}(\overline{\xi}, \overline{\eta}) m_{\widehat{\mathcal{C}}}(\overline{\eta}) m_{\widehat{\mathcal{C}}}(\overline{\sigma}) \widehat{F}_1(s, \overline{\eta}) \frac{\xi_0 |\xi| \overline{\xi}_a^{\overline{\eta}}}{|\overline{\xi}|^3} \overline{\sigma}_b^{\overline{\sigma}} \frac{\sigma_0}{|\sigma_0|} P_b^b(\overline{\xi}, \overline{\sigma}) \\
\widehat{X}_b'(\overline{\sigma}) \cdot \nabla_{\overline{\xi}} \widehat{F}_2(s, \overline{\sigma}) ~ d\overline{\eta} ds 
\end{aligned} \label{equdecchampbBHHC-2-2} \\
&\begin{aligned}
- \int_0^t \int e^{i s \varphi} \mu(\overline{\xi}, \overline{\eta}) \mu_{BHH}(\overline{\xi}, \overline{\eta}) m_{\widehat{\mathcal{C}}}(\overline{\eta}) m_{\widehat{\mathcal{C}}}(\overline{\sigma}) \left( \frac{\xi_0 |\xi| \overline{\xi}_a^{\overline{\eta}}}{|\overline{\xi}|^3} \frac{|\sigma| \left( 1 + \theta^{\overline{\eta} \overline{\sigma}} \right)}{\xi_t^{\overline{\eta} \overline{\sigma}}} + \frac{\xi_0 |\sigma|^2 \xi_t^{\overline{\eta} \overline{\sigma}} |\xi|}{|\overline{\xi}|^3} \right) \\
\widehat{F}_1(s, \overline{\eta}) \xi_t^{\overline{\eta} \overline{\sigma}} \frac{\sigma_0}{|\sigma_0|} P_b^b(\overline{\xi}, \overline{\sigma}) \widehat{X}_b'(\overline{\sigma}) \cdot \nabla_{\overline{\xi}} \widehat{F}_2(s, \overline{\sigma}) ~ d\overline{\eta} ds
\end{aligned} \label{equdecchampbBHHC-2-3} 
\end{align}
\end{subequations}
\eqref{equdecchampbBHHC-2-2} is of the form \eqref{lemdecequH-08-resxsigmabbon}. 

On \eqref{equdecchampbBHHC-2-1}, we directly apply an integration by parts in $\overline{\eta}$: 
\begin{subequations}
\begin{align}
&\eqref{equdecchampbBHHC-2-1} \notag \\
&\begin{aligned}
= \int_0^t \int i s \frac{\xi_0 |\xi| \overline{\xi}_a^{\overline{\eta}}}{|\overline{\xi}|^3} \overline{\eta}_b^{\overline{\eta}} \frac{\sigma_0}{|\sigma_0|} P_b^b(\overline{\xi}, \overline{\sigma}) \widehat{X}_b'(\overline{\sigma}) \cdot \nabla_{\overline{\eta}} \varphi e^{i s \varphi} \mu(\overline{\xi}, \overline{\eta}) \mu_{BHH}(\overline{\xi}, \overline{\eta}) m_{\widehat{\mathcal{C}}}(\overline{\eta}) m_{\widehat{\mathcal{C}}}(\overline{\sigma}) \\
\widehat{F}_1(s, \overline{\eta}) \widehat{F}_2(s, \overline{\sigma}) ~ d\overline{\eta} ds 
\end{aligned} \label{equdecchampbBHHC-2-1-1} \\
&\begin{aligned}
+ \int_0^t \int e^{i s \varphi} \mu(\overline{\xi}, \overline{\eta}) \mu_{BHH}(\overline{\xi}, \overline{\eta}) m_{\widehat{\mathcal{C}}}(\overline{\eta}) \frac{\xi_0 |\xi| \overline{\xi}_a^{\overline{\eta}}}{|\overline{\xi}|^3} \overline{\eta}_b^{\overline{\eta}} \frac{\sigma_0}{|\sigma_0|} P_b^b(\overline{\xi}, \overline{\sigma}) \widehat{X}_b'(\overline{\sigma}) \cdot \nabla_{\overline{\eta}} \widehat{F}_1(s, \overline{\eta}) \\
m_{\widehat{\mathcal{C}}}(\overline{\sigma}) \widehat{F}_2(s, \overline{\sigma}) ~ d\overline{\eta} ds 
\end{aligned} \label{equdecchampbBHHC-2-1-2} \\
&\begin{aligned}
+ \int_0^t \int e^{i s \varphi} \nabla_{\overline{\eta}} \cdot \left( \frac{\xi_0 |\xi| \overline{\xi}_a^{\overline{\eta}}}{|\overline{\xi}|^3} \overline{\eta}_b^{\overline{\eta}} \frac{\sigma_0}{|\sigma_0|} P_b^b(\overline{\xi}, \overline{\sigma}) \widehat{X}_b'(\overline{\sigma}) \mu(\overline{\xi}, \overline{\eta}) \mu_{BHH}(\overline{\xi}, \overline{\eta}) m_{\widehat{\mathcal{C}}}(\overline{\eta}) m_{\widehat{\mathcal{C}}}(\overline{\sigma}) \right) \\
\widehat{F}_1(s, \overline{\eta}) \widehat{F}_2(s, \overline{\sigma}) ~ d\overline{\eta} ds 
\end{aligned} \label{equdecchampbBHHC-2-1-3} 
\end{align}
\end{subequations}
The presence of $\overline{\eta}_b^{\overline{\eta}}$ ensures that \eqref{equdecchampbBHHC-2-1-2} is of the form $\eqref{lemdecequH-03-resxetaac} + \eqref{lemdecequH-04-resxetabbon}$, and \eqref{equdecchampbBHHC-2-1-3} is of the form \eqref{lemdecequH-11-dersymb}. 

On the other hand, on \eqref{equdecchampbBHHC-2-3}, the presence of a factor $\xi_t^{\overline{\eta} \overline{\sigma}}$ allows to consider the field $\widehat{X}_{b-\widehat{\mathcal{C}}}(\overline{\eta}, \overline{\sigma})$: indeed, in the bases $(\widehat{X}_a', \widehat{X}_b', \widehat{X}_c')(\overline{\sigma})$, it has coordinates
\begin{align*}
- \frac{\widetilde{P}_b^b(\overline{\eta}, \overline{\sigma}) \widetilde{P}_a^b(\overline{\eta}, \overline{\sigma})}{\widetilde{P}_a^b(\overline{\eta}, \overline{\sigma})^2 + \widetilde{P}_c^b(\overline{\eta}, \overline{\sigma})^2}, \quad 1, \quad \frac{\widetilde{P}_b^b(\overline{\eta}, \overline{\sigma}) \widetilde{P}_c^b(\overline{\eta}, \overline{\sigma})}{\widetilde{P}_a^b(\overline{\eta}, \overline{\sigma})^2 + \widetilde{P}_c^b(\overline{\eta}, \overline{\sigma})^2}
\end{align*}
But $\widetilde{P}_b^b(\overline{\eta}, \overline{\sigma})$ is close to $-1$, while 
\begin{align*}
\widetilde{P}_c^b(\overline{\eta}, \overline{\sigma}) 
&= - \frac{\xi_t^{\overline{\eta} \overline{\sigma}}}{2} \\
\widetilde{P}_a^b(\overline{\eta}, \overline{\sigma}) 
&= - \frac{\sqrt{3}}{4} \left( 1 + \theta^{\overline{\eta} \overline{\sigma}} \right) = O\left( \left( \xi_t^{\overline{\eta} \overline{\sigma}} \right)^2 \right) 
\end{align*}
We deduce that 
\begin{align*}
- \frac{\widetilde{P}_b^b(\overline{\eta}, \overline{\sigma}) \widetilde{P}_a^b(\overline{\eta}, \overline{\sigma})}{\widetilde{P}_a^b(\overline{\eta}, \overline{\sigma})^2 + \widetilde{P}_c^b(\overline{\eta}, \overline{\sigma})^2} 
&= O(1) \\
\frac{\widetilde{P}_b^b(\overline{\eta}, \overline{\sigma}) \widetilde{P}_c^b(\overline{\eta}, \overline{\sigma})}{\widetilde{P}_a^b(\overline{\eta}, \overline{\sigma})^2 + \widetilde{P}_c^b(\overline{\eta}, \overline{\sigma})^2} 
&= O\left( \left( \xi_t^{\overline{\eta} \overline{\sigma}} \right)^{-1} \right) 
\end{align*}
In particular, the vector field $\xi_t^{\overline{\eta} \overline{\sigma}} \widehat{X}_{b-\widehat{\mathcal{C}}}(\overline{\eta}, \overline{\sigma})$ is not singular anymore, and we can therefore replace \eqref{equdecchampbBHHC-2-3}, up to terms of the form \eqref{lemdecequH-07-resxsigmaac}, by 
\begin{subequations}
\begin{align}
&\begin{aligned}
&- \int_0^t \int e^{i s \varphi} \mu(\overline{\xi}, \overline{\eta}) \mu_{BHH}(\overline{\xi}, \overline{\eta}) m_{\widehat{\mathcal{C}}}(\overline{\eta}) m_{\widehat{\mathcal{C}}}(\overline{\sigma}) \widehat{F}_1(s, \overline{\eta}) \left( \frac{\xi_0 |\xi| \overline{\xi}_a^{\overline{\eta}}}{|\overline{\xi}|^3} \frac{|\sigma| \left( 1 + \theta^{\overline{\eta} \overline{\sigma}} \right)}{\xi_t^{\overline{\eta} \overline{\sigma}}} + \frac{\xi_0 |\sigma|^2 \xi_t^{\overline{\eta} \overline{\sigma}} |\xi|}{|\overline{\xi}|^3} \right) \\
&\pushright{\xi_t^{\overline{\eta} \overline{\sigma}} \frac{\sigma_0}{|\sigma_0|} P_b^b(\overline{\xi}, \overline{\sigma}) \widehat{X}_{b-\widehat{\mathcal{C}}}(\overline{\sigma}) \cdot \nabla_{\overline{\xi}} \widehat{F}_2(s, \overline{\sigma}) ~ d\overline{\eta} ds}
\end{aligned} \notag \\
&\begin{aligned}
= - \int_0^t \int i s \left( \frac{\xi_0 |\xi| \overline{\xi}_a^{\overline{\eta}}}{|\overline{\xi}|^3} \frac{|\sigma| \left( 1 + \theta^{\overline{\eta} \overline{\sigma}} \right)}{\xi_t^{\overline{\eta} \overline{\sigma}}} + \frac{\xi_0 |\sigma|^2 \xi_t^{\overline{\eta} \overline{\sigma}} |\xi|}{|\overline{\xi}|^3} \right) \xi_t^{\overline{\eta} \overline{\sigma}} \frac{\sigma_0}{|\sigma_0|} P_b^b(\overline{\xi}, \overline{\sigma}) \widehat{X}_{b-\widehat{\mathcal{C}}}(\overline{\sigma}) \cdot \nabla_{\overline{\eta}} \varphi \\
e^{i s \varphi} \mu(\overline{\xi}, \overline{\eta}) \mu_{BHH}(\overline{\xi}, \overline{\eta}) m_{\widehat{\mathcal{C}}}(\overline{\eta}) m_{\widehat{\mathcal{C}}}(\overline{\sigma}) \widehat{F}_1(s, \overline{\eta}) \widehat{F}_2(s, \overline{\sigma}) ~ d\overline{\eta} ds 
\end{aligned} \label{equdecchampbBHHC-2-3-1} \\
&\begin{aligned} 
- \int_0^t \int e^{i s \varphi} \mu(\overline{\xi}, \overline{\eta}) \mu_{BHH}(\overline{\xi}, \overline{\eta}) m_{\widehat{\mathcal{C}}}(\overline{\eta}) m_{\widehat{\mathcal{C}}}(\overline{\sigma}) \left( \frac{\xi_0 |\xi| \overline{\xi}_a^{\overline{\eta}}}{|\overline{\xi}|^3} \frac{|\sigma| \left( 1 + \theta^{\overline{\eta} \overline{\sigma}} \right)}{\xi_t^{\overline{\eta} \overline{\sigma}}} + \frac{\xi_0 |\sigma|^2 \xi_t^{\overline{\eta} \overline{\sigma}} |\xi|}{|\overline{\xi}|^3} \right) \\
\xi_t^{\overline{\eta} \overline{\sigma}} \frac{\sigma_0}{|\sigma_0|} P_b^b(\overline{\xi}, \overline{\sigma}) \widehat{X}_{b-\widehat{\mathcal{C}}}(\overline{\sigma}) \cdot \nabla_{\overline{\eta}} \widehat{F}_1(s, \overline{\eta}) \widehat{F}_2(s, \overline{\sigma}) ~ d\overline{\eta} ds 
\end{aligned} \label{equdecchampbBHHC-2-3-2} \\
&\begin{aligned} 
- \int_0^t \int e^{i s \varphi} \nabla_{\overline{\eta}} \cdot \Biggl( \left( \frac{\xi_0 |\xi| \overline{\xi}_a^{\overline{\eta}}}{|\overline{\xi}|^3} \frac{|\sigma| \left( 1 + \theta^{\overline{\eta} \overline{\sigma}} \right)}{\xi_t^{\overline{\eta} \overline{\sigma}}} + \frac{\xi_0 |\sigma|^2 \xi_t^{\overline{\eta} \overline{\sigma}} |\xi|}{|\overline{\xi}|^3} \right) \xi_t^{\overline{\eta} \overline{\sigma}} \frac{\sigma_0}{|\sigma_0|} P_b^b(\overline{\xi}, \overline{\sigma}) \\
\widehat{X}_{b-\widehat{\mathcal{C}}}(\overline{\sigma}) \mu(\overline{\xi}, \overline{\eta}) \mu_{BHH}(\overline{\xi}, \overline{\eta}) m_{\widehat{\mathcal{C}}}(\overline{\eta}) m_{\widehat{\mathcal{C}}}(\overline{\sigma}) \Biggl) \widehat{F}_1(s, \overline{\eta}) \widehat{F}_2(s, \overline{\sigma}) ~ d\overline{\eta} ds 
\end{aligned} \label{equdecchampbBHHC-2-3-3} 
\end{align}
\end{subequations}
\eqref{equdecchampbBHHC-2-3-2} is of the form \eqref{lemdecequH-03-resxetaac} by Lemma \ref{lemchampmodifiebCC-propfond}, and \eqref{equdecchampbBHHC-2-3-3} is of the form \eqref{lemdecequH-11-dersymb}. 

We now try to find good symbols that would allow to apply integrations by parts on the remaining terms and prove the decomposition. First, we could try to apply the field $\widehat{Y}(\overline{\eta}, \overline{\sigma})$: by Lemma \ref{lemcalculchampmodifie}, since here $\epsilon^{\overline{\eta} \overline{\sigma}} = -1$ and $\theta^{\overline{\eta} \overline{\sigma}}$ is close to $-1$, 
\begin{align*}
\widehat{Y}(\overline{\eta}, \overline{\sigma}) \cdot \nabla_{\overline{\eta}} \varphi = 2 \frac{\eta_0}{|\eta_0|} \frac{1 + \epsilon^{\overline{\eta} \overline{\sigma}} \theta^{\overline{\eta} \overline{\sigma}}}{7 + \epsilon^{\overline{\eta} \overline{\sigma}} \theta^{\overline{\eta} \overline{\sigma}}} \left( \left( \overline{\sigma}_a^{\overline{\sigma}} \right)^2 - \left( \overline{\eta}_a^{\overline{\eta}} \right)^2 \right) 
\end{align*}
so we can control $\overline{\sigma}_a^{\overline{\sigma}} - \overline{\eta}_a^{\overline{\eta}}$. 

Moreover, 
\begin{align*}
\overline{\sigma}_b^{\overline{\sigma}} \xi_t^{\overline{\eta} \overline{\sigma}} &= O\left( m_b(\overline{\sigma}) \widehat{X}_c(\overline{\eta}) \cdot \nabla_{\overline{\eta}} \varphi \right) 
\end{align*}
and in a symmetric way we also have $\overline{\eta}_b^{\overline{\eta}} \xi_t^{\overline{\eta} \overline{\sigma}}$. 

We will also use that 
\begin{align*}
|\eta| - |\sigma| &= \frac{1}{2} \left( \overline{\sigma}_a^{\overline{\sigma}} - \overline{\sigma}_b^{\overline{\sigma}} - \overline{\eta}_a^{\overline{\eta}} + \overline{\eta}_b^{\overline{\eta}} \right) \\
&= \frac{1}{2} \overline{\xi}_b^{\overline{\eta}} + O\left( \xi_t^{\overline{\eta} \overline{\sigma}} \right) + O\left( \overline{\eta}_a^{\overline{\eta}} - \overline{\sigma}_a^{\overline{\sigma}} \right) \\
|\xi|^2 &= |\eta|^2 + |\sigma|^2 + 2 \eta \cdot \sigma \\
&= \left( |\eta| - |\sigma| \right)^2 + O\left( \xi_t^{\overline{\eta} \overline{\sigma}} \right) \\
&= \frac{1}{4} \left( \overline{\xi}_b^{\overline{\eta}} \right)^2 + O\left( \xi_t^{\overline{\eta} \overline{\sigma}} \right) + O\left( \overline{\eta}_a^{\overline{\eta}} - \overline{\sigma}_a^{\overline{\sigma}} \right) \\
\xi \cdot \sigma &= \eta \cdot \sigma + |\sigma|^2 \\
&= |\sigma| \left( |\sigma| - |\eta| \right) + O\left( \xi_t^{\overline{\eta} \overline{\sigma}} \right) \\
&= - \frac{1}{2} |\sigma| \overline{\xi}_b^{\overline{\eta}} + O\left( \xi_t^{\overline{\eta} \overline{\sigma}} \right) + O\left( \overline{\eta}_a^{\overline{\eta}} - \overline{\sigma}_a^{\overline{\sigma}} \right) \\
\overline{\xi}_a^{\overline{\eta}} &= \overline{\eta}_a^{\overline{\eta}} + \overline{\sigma}_a^{\overline{\eta}} \\
&= O\left( \xi_t^{\overline{\eta} \overline{\sigma}} \right) + O\left( \overline{\eta}_a^{\overline{\eta}} - \overline{\sigma}_a^{\overline{\sigma}} \right) \\
\overline{\eta}_b^{\overline{\eta}} - \overline{\sigma}_b^{\overline{\sigma}} &= \sqrt{3} |\eta_0| - \sqrt{3} |\sigma_0| - |\eta| + |\sigma| = O\left( \overline{\xi} \right) 
\end{align*}

For \eqref{equdecchampbBHHC-1}, we write that 
\begin{align*}
\xi_0 m_b(\overline{\xi}) &= O\left( \frac{\xi_0 (3 \xi_0^2 - |\xi|^2)}{|\overline{\xi}|^2} \right) \\
&= O\left( \frac{\xi_0}{|\overline{\xi}|} \xi_t^{\overline{\eta} \overline{\sigma}} \right) + O\left( \frac{\xi_0 \overline{\xi}_a^{\overline{\eta}}}{|\overline{\xi}|^2} \left( \overline{\eta}_b^{\overline{\eta}} - \overline{\sigma}_b^{\overline{\sigma}} \right) \right) \\
&= O\left( \frac{\xi_0}{|\overline{\xi}|} \xi_t^{\overline{\eta} \overline{\sigma}} \right) + O\left( \overline{\eta}_a^{\overline{\eta}} - \overline{\sigma}_a^{\overline{\sigma}} \right) 
\end{align*}

For \eqref{equdecchampbBHHC-2-1-1}, we write that 
\begin{align*}
\frac{\overline{\xi}_a^{\overline{\xi}}}{|\overline{\xi}|} \overline{\eta}_b^{\overline{\eta}} \widehat{X}_b'(\overline{\sigma}) \cdot \nabla_{\overline{\eta}} \varphi 
&= \frac{\overline{\xi}_a^{\overline{\xi}}}{|\overline{\xi}|} \overline{\eta}_b^{\overline{\eta}} O\left( \overline{\xi} \right) \\
&= O\left( \overline{\xi}_a^{\overline{\xi}} \overline{\eta}_b^{\overline{\eta}} \right) \\
&= O\left( \overline{\eta}_b^{\overline{\eta}} \xi_t^{\overline{\eta} \overline{\sigma}} \right) + O\left( \overline{\eta}_a^{\overline{\eta}} - \overline{\sigma}_a^{\overline{\sigma}} \right) 
\end{align*}

Finally, for \eqref{equdecchampbBHHC-2-3-1}, the symbol appearing is already of the form $O\left( \frac{\xi_0}{|\overline{\xi}|} \xi_t^{\overline{\eta} \overline{\sigma}} \right)$. 

To conclude, we only need to explain how to decompose if we have a factor $O\left( \frac{\xi_0}{|\overline{\xi}|} \xi_t^{\overline{\eta} \overline{\sigma}} \right)$. 

To that end, we exhaust the cases depending on the angular localisations of $\overline{\xi}, \overline{\eta}, \overline{\sigma}$. Even if $\overline{\xi}$ is very small with respect to the two other Fourier variables, we can obtain such localisations decomposing the unit sphere $S^2$ into a sufficiently fine (finite) partition of unity, and then using Hörmander-Mikhlin symbols in all three Fourier variables to localise according to this partition of $S^2$, to the third power. 

\paragraph{1.} If $\overline{\xi}_a^{\overline{\eta}} \simeq |\overline{\xi}|$, then we may compute more finely 
\begin{align*}
\overline{\xi}_a^{\overline{\eta}} &= \overline{\eta}_a^{\overline{\eta}} + \overline{\sigma}_a^{\overline{\eta}} \\
&= \overline{\eta}_a^{\overline{\eta}} - \overline{\sigma}_a^{\overline{\sigma}} + O\left( 1 + \theta^{\overline{\eta} \overline{\sigma}} \right) 
\end{align*}
But $\xi_t^{\overline{\eta} \overline{\sigma}} = O\left( |\xi| \right) = O\left( |\overline{\xi}| \right) = o(1)$, so that 
\begin{align*}
1 + \theta^{\overline{\eta} \overline{\sigma}} &= O\left( \left( \xi_t^{\overline{\eta} \overline{\sigma}} \right)^2 \right) \\
&= o\left( |\overline{\xi}| \right) 
\end{align*}
We deduce that 
\begin{align*}
O\left( \frac{\xi_0}{|\overline{\xi}|} \xi_t^{\overline{\eta} \overline{\sigma}} \right) &= O\left( \overline{\xi} \right) = O\left( \overline{\xi}_a^{\overline{\eta}} \right) = O\left( \overline{\eta}_a^{\eta} - \overline{\sigma}_a^{\overline{\sigma}} \right) 
\end{align*}
which is enough. 

\paragraph{2.} If $\overline{\xi}_a^{\overline{\eta}} \ll |\overline{\xi}|$, we can write that 
\begin{align*}
\xi_0 &= \frac{1}{2 \sqrt{3}} \frac{\eta_0}{|\eta_0|} \left( 2 \sqrt{3} |\eta_0| - 2 \sqrt{3} |\sigma_0| \right) \\
&= \frac{1}{2 \sqrt{3}} \frac{\eta_0}{|\eta_0|} \left( \overline{\eta}_a^{\overline{\eta}} + \overline{\eta}_b^{\overline{\eta}} - \overline{\sigma}_a^{\overline{\sigma}} - \overline{\sigma}_b^{\overline{\sigma}} \right) \\
&= O\left( \overline{\eta}_b^{\overline{\eta}} - \overline{\sigma}_b^{\overline{\sigma}} \right) + O\left( \overline{\eta}_a^{\overline{\eta}} - \overline{\sigma}_a^{\overline{\sigma}} \right) 
\end{align*}
Therefore, 
\begin{align*}
\frac{\xi_0}{|\overline{\xi}|} \xi_t^{\overline{\eta} \overline{\sigma}} &= \frac{O\left( \overline{\eta}_b^{\overline{\eta}} - \overline{\sigma}_b^{\overline{\sigma}} \right)}{|\overline{\xi}|} \xi_t^{\overline{\eta} \overline{\sigma}} + O\left( \overline{\eta}_a^{\overline{\eta}} - \overline{\sigma}_a^{\overline{\sigma}} \right)
\end{align*}
and we already know how to decompose the second term. We therefore only have to treat 
\begin{align}
\int_0^t \int e^{i s \varphi} s \frac{\overline{\eta}_b^{\overline{\eta}} - \overline{\sigma}_b^{\overline{\sigma}}}{|\overline{\xi}|} \xi_t^{\overline{\eta} \overline{\sigma}} |\overline{\eta}|^3 \mu \widehat{F}_1(s, \overline{\eta}) \widehat{F}_2(s, \overline{\sigma}) ~ d\overline{\eta} ds \label{equdecchampbBHHCC-termefinalrestant} 
\end{align}
where $\mu$ now contains all previous localisation symbols and $O(1)$ remainders whose form is irrelevant. $\mu$ is still of order $0$. 

We know that 
\begin{align*}
\xi_t^{\overline{\eta} \overline{\sigma}} &= O\left( \widehat{X}_c(\overline{\eta}) \cdot \nabla_{\overline{\eta}} \varphi \right) 
\end{align*}
so that we may apply an integration by parts in this direction. Thus, allowing the symbol $\mu$ to vary from line to line, 
\begin{subequations}
\begin{align}
\eqref{equdecchampbBHHCC-termefinalrestant} &= \int_0^t \int e^{i s \varphi} \frac{\overline{\eta}_b^{\overline{\eta}} - \overline{\sigma}_b^{\overline{\sigma}}}{|\overline{\xi}|} |\overline{\eta}| \mu \widehat{X}_c(\overline{\eta}) \cdot \nabla_{\overline{\eta}} \widehat{F}_1(s, \overline{\eta}) \widehat{F}_2(s, \overline{\sigma}) ~ d\overline{\eta} ds \label{equdecchampbBHHCC-termefinalrestant-1} \\
&+ \int_0^t \int e^{i s \varphi} \frac{\overline{\eta}_b^{\overline{\eta}} - \overline{\sigma}_b^{\overline{\sigma}}}{|\overline{\xi}|} |\overline{\eta}| \mu \widehat{F}_1(s, \overline{\eta}) \left( P_c^a(\overline{\eta}, \overline{\sigma}) \widehat{X}_a(\overline{\sigma}) + P_c^c(\overline{\eta}, \overline{\sigma}) \widehat{X}_c(\overline{\sigma}) \right) \cdot \nabla_{\overline{\eta}} \widehat{F}_2(s, \overline{\sigma}) ~ d\overline{\eta} ds \label{equdecchampbBHHCC-termefinalrestant-2} \\
&+ \int_0^t \int e^{i s \varphi} \frac{\overline{\eta}_b^{\overline{\eta}} - \overline{\sigma}_b^{\overline{\sigma}}}{|\overline{\xi}|} |\overline{\eta}| \mu \widehat{F}_1(s, \overline{\eta}) P_c^b(\overline{\eta}, \overline{\sigma}) \widehat{X}_b(\overline{\sigma}) \cdot \nabla_{\overline{\eta}} \widehat{F}_2(s, \overline{\sigma}) ~ d\overline{\eta} ds \label{equdecchampbBHHCC-termefinalrestant-3} \\
&+ \int_0^t \int e^{i s \varphi} \nabla_{\overline{\eta}} \cdot \left( \widehat{X}_c(\overline{\eta}) \frac{\overline{\eta}_b^{\overline{\eta}} - \overline{\sigma}_b^{\overline{\sigma}}}{|\overline{\xi}|} |\overline{\eta}| \mu \right) \widehat{F}_1(s, \overline{\eta}) \widehat{F}_2(s, \overline{\sigma}) ~ d\overline{\eta} ds \label{equdecchampbBHHCC-termefinalrestant-4}
\end{align}
\end{subequations}
\eqref{equdecchampbBHHCC-termefinalrestant-1} is of the form \eqref{lemdecequH-03-resxetaac}, \eqref{equdecchampbBHHCC-termefinalrestant-2} of the form \eqref{lemdecequH-07-resxsigmaac}, \eqref{equdecchampbBHHCC-termefinalrestant-4} of the form \eqref{lemdecequH-11-dersymb}. 

For \eqref{equdecchampbBHHCC-termefinalrestant-3}, we have that
\begin{align*}
P_c^b(\overline{\eta}, \overline{\sigma}) &= O\left( \xi_t^{\overline{\eta} \overline{\sigma}} \right) 
\end{align*}
Therefore, we can rewrite the symbol as 
\begin{align*}
\frac{\overline{\eta}_b^{\overline{\eta}} - \overline{\sigma}_b^{\overline{\sigma}}}{|\overline{\xi}|} P_b^b(\overline{\eta}, \overline{\sigma}) &= \overline{\eta}_b^{\overline{\eta}} \frac{\xi_t^{\overline{\eta} \overline{\sigma}}}{|\overline{\xi}|} + O\left( \overline{\sigma}_b^{\overline{\sigma}} \right) 
\end{align*}
The second term contributes like \eqref{lemdecequH-08-resxsigmabbon}. On the first, we recover a structure analogous to \eqref{equdecchampbBHHC-2-1} and we can apply another integration by parts in $\overline{\eta}$: the presence of the factor $\overline{\eta}_b^{\overline{\eta}}$ ensures that we only get terms of the form \eqref{lemdecequH-03-resxetaac}, \eqref{lemdecequH-04-resxetabbon} or \eqref{lemdecequH-11-dersymb}, plus 
\begin{align*}
\int_0^t \int e^{i s \varphi} s \overline{\eta}_b^{\overline{\eta}} \frac{\xi_t^{\overline{\eta} \overline{\sigma}}}{|\overline{\xi}|} \nabla_{\overline{\eta}} \varphi |\overline{\eta}|^2 \mu \widehat{F}_1(s, \overline{\eta}) \widehat{F}_2(s, \overline{\sigma}) ~ d\overline{\eta} ds
\end{align*}
But precisely $\nabla_{\overline{\eta}} \varphi = O\left( \overline{\xi} \right)$, so we only get now a factor $O\left( \overline{\eta}_b^{\overline{\eta}} \xi_t^{\overline{\eta} \overline{\sigma}} \right)$, and we know how to decompose it. 

Note that we applied no integrations by parts in time for this interaction, so that locally $m_g \equiv 0$ and thus \eqref{decompositionfinemg-gainkout} holds. 

\subsubsection{Interaction \texorpdfstring{$(BHH, \widehat{\mathcal{L}})$}{(BHH, L)}}

Let us consider 
\begin{subequations}
\begin{align}
&\xi_0 m_b(\overline{\xi}) \widehat{X}_b(\overline{\xi}) \cdot \nabla_{\overline{\xi}} \widehat{I}_{\mu}^{(BHH, \widehat{\mathcal{L}})}[F_1, F_2](t, \overline{\xi}) \notag \\
&\quad = \int_0^t \int i s \xi_0 m_b(\overline{\xi}) \widehat{X}_b(\overline{\xi}) \cdot \nabla_{\overline{\xi}} \varphi e^{i s \varphi} \mu(\overline{\xi}, \overline{\eta}) \mu_{BHH}(\overline{\xi}, \overline{\eta}) m_{\widehat{\mathcal{L}}}(\overline{\eta}) m_{\widehat{\mathcal{L}}}(\overline{\sigma}) \widehat{F}_1(s, \overline{\eta}) \widehat{F}_2(s, \overline{\sigma}) ~ d\overline{\eta} ds \label{equdecchampbBHHL-1} \\
&\quad \quad + \int_0^t \int e^{i s \varphi} \xi_0 \mu(\overline{\xi}, \overline{\eta}) \mu_{BHH}(\overline{\xi}, \overline{\eta}) m_{\widehat{\mathcal{L}}}(\overline{\eta}) m_{\widehat{\mathcal{L}}}(\overline{\sigma}) \widehat{F}_1(s, \overline{\eta}) m_b(\overline{\xi}) \widehat{X}_b(\overline{\xi}) \cdot \nabla_{\overline{\xi}} \widehat{F}_2(s, \overline{\sigma}) ~ d\overline{\eta} ds \label{equdecchampbBHHL-2} \\
&\quad \quad + \int_0^t \int e^{i s \varphi} \xi_0 m_b(\overline{\xi}) \widehat{X}_b(\overline{\xi}) \cdot \nabla_{\overline{\xi}} \left( \mu(\overline{\xi}, \overline{\eta}) \mu_{BHH}(\overline{\xi}, \overline{\eta}) m_{\widehat{\mathcal{L}}}(\overline{\eta}) m_{\widehat{\mathcal{L}}}(\overline{\sigma}) \right) \widehat{F}_1(s, \overline{\eta}) \widehat{F}_2(s, \overline{\sigma}) ~ d\overline{\eta} ds \label{equdecchampbBHHL-3} 
\end{align}
\end{subequations}
\eqref{equdecchampbBHHL-3} is of the form \eqref{lemdecequH-11-dersymb}. 

For \eqref{equdecchampbBHHL-2}, we have
\begin{align*}
m_b(\overline{\xi}) |\overline{\xi}| \widehat{X}_b(\overline{\xi}) &= O(1) e_0 - \frac{3 \xi_0^2 - |\xi|^2}{|\overline{\xi}|^3} \xi_0 \xi \\
&= O(1) e_0 + O(\sigma) - \frac{3 \xi_0^2 - |\xi|^2}{|\overline{\xi}|^3} \xi_0 \begin{pmatrix} 0 \\ \eta \end{pmatrix}
\end{align*}
The two first terms contribute like \eqref{lemdecequH-07-resxsigmaac} or \eqref{lemdecequH-08-resxsigmabbon}; on the third, we apply an integration by parts in $\eta$ and get terms of the form \eqref{lemdecequH-03-resxetaac}, \eqref{lemdecequH-04-resxetabbon}, \eqref{lemdecequH-11-dersymb} plus
\begin{align*}
- \int_0^t \int i s \frac{3 \xi_0^2 - |\xi|^2}{|\overline{\xi}|^3} \xi_0 \eta \cdot \nabla_{\eta} \varphi e^{i s \varphi} \frac{\xi_0}{|\overline{\xi}|} \mu(\overline{\xi}, \overline{\eta}) \mu_{BHH}(\overline{\xi}, \overline{\eta}) m_{\widehat{\mathcal{L}}}(\overline{\eta}) m_{\widehat{\mathcal{L}}}(\overline{\sigma}) \widehat{F}_1(s, \overline{\eta}) \widehat{F}_2(s, \overline{\sigma}) ~ d\overline{\eta} ds
\end{align*}
that we group with \eqref{equdecchampbBHHL-1} to get: 
\begin{align*}
\int_0^t \int i s \left( m_b(\overline{\xi}) |\overline{\xi}| \widehat{X}_b(\overline{\xi}) \cdot \nabla_{\overline{\xi}} - \frac{3 \xi_0^2 - |\xi|^2}{|\overline{\xi}|^3} \xi_0 \eta \cdot \nabla_{\eta} \right) \varphi e^{i s \varphi} \frac{\xi_0}{|\overline{\xi}|} \mu(\overline{\xi}, \overline{\eta}) \mu_{BHH}(\overline{\xi}, \overline{\eta}) m_{\widehat{\mathcal{L}}}(\overline{\eta}) \\
m_{\widehat{\mathcal{L}}}(\overline{\sigma}) \widehat{F}_1(s, \overline{\eta}) \widehat{F}_2(s, \overline{\sigma}) ~ d\overline{\eta} ds
\end{align*}

Note that 
\begin{align*}
\varphi + \sigma_0 \partial_{\eta_0} \varphi &= \xi_0^3 + \xi_0 |\xi|^2 - \eta_0^3 - \eta_0 |\eta|^2 - \sigma_0^3 - \sigma_0 |\sigma|^2 + \sigma_0 (3 \sigma_0^2 + |\sigma|^2 - 3 \eta_0^2 - |\eta|^2) \\
&= \xi_0 (3 \eta_0 \sigma_0 + |\xi|^2 - |\eta|^2 + 3 \sigma_0 (\sigma_0 - \eta_0)) \\
&= \xi_0 (3 \sigma_0^2 + o(1)) 
\end{align*}
so we know how to decompose if we have a factor $O(\xi_0)$. But
\begin{align*}
m_b(\overline{\xi}) |\overline{\xi}| \widehat{X}_b(\overline{\xi}) \cdot \nabla_{\overline{\xi}} \varphi \frac{\xi_0}{|\overline{\xi}|} = O(\xi_0 m_b(\overline{\xi}))
\end{align*}
and 
\begin{align*}
\frac{3 \xi_0^2 - |\xi|^2}{|\overline{\xi}|^3} \xi_0 \eta \cdot \nabla_{\eta} \varphi = O(\xi_0 m_b(\overline{\xi})) 
\end{align*}
which concludes, the presence of $m_b(\overline{\xi})$ ensuring that \eqref{decompositionfinemg-gainkout} holds. 

\subsubsection{Interaction \texorpdfstring{$(BHH, \widehat{\mathcal{P}})$}{(BHH, P)}}

Here, we need to separate into three subcases depending on the relative sizes of $\eta_0$ and $\sigma_0$. Due to the symmetry between $\overline{\eta}$ and $\overline{\sigma}$ here, it is enough to treat only two: either $|\sigma_0| \simeq |\eta_0|$, or $|\eta_0| \ll |\sigma_0|$. 

\paragraph{1.} Let us first localise on $\{ |\eta_0| \simeq |\sigma_0| \}$, and consider: 
\begin{subequations}
\begin{align}
&\xi_0 m_b(\overline{\xi}) \widehat{X}_b(\overline{\xi}) \cdot \nabla_{\overline{\xi}} \left( \int_0^t \int e^{i s \varphi} \mu(\overline{\xi}, \overline{\eta}) \mu_{BHH}(\overline{\xi}, \overline{\eta}) \mu_{BB}(\xi_0, \eta_0) m_{\widehat{\mathcal{P}}}(\overline{\eta}) m_{\widehat{\mathcal{P}}}(\overline{\sigma}) \widehat{F}_1(s, \overline{\eta}) \widehat{F}_2(s, \overline{\sigma}) ~ d\overline{\eta} ds \right) \notag \\
&= \int_0^t \int i s \xi_0 m_b(\overline{\xi}) \widehat{X}_b(\overline{\xi}) \cdot \nabla_{\overline{\xi}} \varphi e^{i s \varphi} \mu(\overline{\xi}, \overline{\eta}) \mu_{BHH}(\overline{\xi}, \overline{\eta}) \mu_{BB}(\xi_0, \eta_0) m_{\widehat{\mathcal{P}}}(\overline{\eta}) m_{\widehat{\mathcal{P}}}(\overline{\sigma}) \widehat{F}_1(s, \overline{\eta}) \widehat{F}_2(s, \overline{\sigma}) ~ d\overline{\eta} ds \label{equdecchampbBHHPBB-1} \\
&+ \int_0^t \int e^{i s \varphi} \xi_0 \mu(\overline{\xi}, \overline{\eta}) \mu_{BHH}(\overline{\xi}, \overline{\eta}) \mu_{BB}(\xi_0, \eta_0) m_{\widehat{\mathcal{P}}}(\overline{\eta}) m_{\widehat{\mathcal{P}}}(\overline{\sigma}) \widehat{F}_1(s, \overline{\eta}) m_b(\overline{\xi}) \widehat{X}_b(\overline{\xi}) \cdot \nabla_{\overline{\xi}} \widehat{F}_2(s, \overline{\sigma}) ~ d\overline{\eta} ds \label{equdecchampbBHHPBB-2} \\
&+ \int_0^t \int e^{i s \varphi} \xi_0 m_b(\overline{\xi}) \widehat{X}_b(\overline{\xi}) \cdot \nabla_{\overline{\xi}} \left( \mu(\overline{\xi}, \overline{\eta}) \mu_{BHH}(\overline{\xi}, \overline{\eta}) \mu_{BB}(\xi_0, \eta_0) m_{\widehat{\mathcal{P}}}(\overline{\eta}) m_{\widehat{\mathcal{P}}}(\overline{\sigma}) \right) \widehat{F}_1(s, \overline{\eta}) \widehat{F}_2(s, \overline{\sigma}) ~ d\overline{\eta} ds \label{equdecchampbBHHPBB-3} 
\end{align}
\end{subequations}
\eqref{equdecchampbBHHPBB-3} is of the form \eqref{lemdecequH-11-dersymb} (even if $\mu_{BB}$ is differentiated, the presence of $\xi_0$ compensates the singularity). \eqref{equdecchampbBHHPBB-2} is of the form \eqref{lemdecequH-09-resxsigmabbonbis} since $\xi_0 = O(\eta_0)$. 

Then, \eqref{lemdecequH-01-symeta} and \eqref{lemdecequH-02-symsigma} are both already of the form \eqref{lemdecequH-05-resxetabbonbis} and \eqref{lemdecequH-09-resxsigmabbonbis} so we don't need to remove them. 

Finally, for \eqref{equdecchampbBHHPBB-1}, we compute that
\begin{align*}
\varphi - \eta_0 \partial_{\eta_0} \varphi 
&= \xi_0 \left( 3 \eta_0^2 + |\xi|^2 - |\sigma|^2 \right) \\
&= \xi_0 (|\sigma|^2 + o(1)) 
\end{align*}
so we know how to decompose in the presence of $O(\xi_0)$, and hence \eqref{equdecchampbBHHPBB-1} can be decomposed into terms of \eqref{lemdecequHGtot}. We did not use $m_b(\overline{\xi})$ so \eqref{decompositionfinemg-gainkout} holds. 

\paragraph{2.} Localise now on $\{ |\sigma_0| \gg |\eta_0| \}$. Then we consider: 
\begin{subequations}
\begin{align}
&\xi_0 m_b(\overline{\xi}) \widehat{X}_b(\overline{\xi}) \cdot \nabla_{\overline{\xi}} \left( \int_0^t \int e^{i s \varphi} \mu(\overline{\xi}, \overline{\eta}) \mu_{BHH}(\overline{\xi}, \overline{\eta}) \mu_{BH}(\xi_0, \eta_0) m_{\widehat{\mathcal{P}}}(\overline{\eta}) m_{\widehat{\mathcal{P}}}(\overline{\sigma}) \widehat{F}_1(s, \overline{\eta}) \widehat{F}_2(s, \overline{\sigma}) ~ d\overline{\eta} ds \right) \notag \\
&= \int_0^t \int i s \xi_0 m_b(\overline{\xi}) \widehat{X}_b(\overline{\xi}) \cdot \nabla_{\overline{\xi}} \varphi e^{i s \varphi} \mu(\overline{\xi}, \overline{\eta}) \mu_{BHH}(\overline{\xi}, \overline{\eta}) \mu_{BH}(\xi_0, \eta_0) m_{\widehat{\mathcal{P}}}(\overline{\eta}) m_{\widehat{\mathcal{P}}}(\overline{\sigma}) \widehat{F}_1(s, \overline{\eta}) \widehat{F}_2(s, \overline{\sigma}) ~ d\overline{\eta} ds \label{equdecchampbBHHPBH-1} \\
&+ \int_0^t \int e^{i s \varphi} \xi_0 \mu(\overline{\xi}, \overline{\eta}) \mu_{BHH}(\overline{\xi}, \overline{\eta}) \mu_{BH}(\xi_0, \eta_0) m_{\widehat{\mathcal{P}}}(\overline{\eta}) m_{\widehat{\mathcal{P}}}(\overline{\sigma}) \widehat{F}_1(s, \overline{\eta}) m_b(\overline{\xi}) \widehat{X}_b(\overline{\xi}) \cdot \nabla_{\overline{\xi}} \widehat{F}_2(s, \overline{\sigma}) ~ d\overline{\eta} ds \label{equdecchampbBHHPBH-2} \\
&+ \int_0^t \int e^{i s \varphi} \xi_0 m_b(\overline{\xi}) \widehat{X}_b(\overline{\xi}) \cdot \nabla_{\overline{\xi}} \left( \mu(\overline{\xi}, \overline{\eta}) \mu_{BHH}(\overline{\xi}, \overline{\eta}) \mu_{BH}(\xi_0, \eta_0) m_{\widehat{\mathcal{P}}}(\overline{\eta}) m_{\widehat{\mathcal{P}}}(\overline{\sigma}) \right) \widehat{F}_1(s, \overline{\eta}) \widehat{F}_2(s, \overline{\sigma}) ~ d\overline{\eta} ds \label{equdecchampbBHHPBH-3} 
\end{align}
\end{subequations}
\eqref{equdecchampbBHHPBH-3} is of the form \eqref{lemdecequH-11-dersymb}. 

Then, \eqref{lemdecequH-01-symeta} is already of the form \eqref{lemdecequH-05-resxetabbonbis}. However, we remove from \eqref{equdecchampbBHHPBH-2} the desired term \eqref{lemdecequH-02-symsigma} and then apply an integration by parts and group with \eqref{equdecchampbBHHPBH-1}: 
\begin{subequations}
\begin{align}
&\eqref{equdecchampbBHHPBH-1} + \eqref{equdecchampbBHHPBH-2} - \eqref{lemdecequH-02-symsigma} \notag \\
&\begin{aligned}
= \int_0^t \int i s \left( \xi_0 m_b(\overline{\xi}) \widehat{X}_b(\overline{\xi}) \cdot \left( \nabla_{\overline{\xi}} + \nabla_{\overline{\eta}} \right) - \sigma_0 m_b(\overline{\sigma}) \widehat{X}_b(\overline{\sigma}) \cdot \nabla_{\overline{\eta}} \right) \varphi e^{i s \varphi} \mu(\overline{\xi}, \overline{\eta}) \\
\mu_{BHH}(\overline{\xi}, \overline{\eta}) \mu_{BH}(\xi_0, \eta_0) m_{\widehat{\mathcal{P}}}(\overline{\eta}) m_{\widehat{\mathcal{P}}}(\overline{\sigma}) \widehat{F}_1(s, \overline{\eta}) \widehat{F}_2(s, \overline{\sigma}) ~ d\overline{\eta} ds \end{aligned} \label{equdecchampbBHHPBH-1-1} \\
&\begin{aligned}
+ \int_0^t \int e^{i s \varphi} \mu(\overline{\xi}, \overline{\eta}) \mu_{BHH}(\overline{\xi}, \overline{\eta}) \left( \xi_0 m_b(\overline{\xi}) \widehat{X}_b(\overline{\xi}) - \sigma_0 m_b(\overline{\sigma}) \widehat{X}_b(\overline{\sigma}) \right) \cdot \nabla_{\overline{\eta}} \widehat{F}_1(s, \overline{\eta}) \\
\mu_{BH}(\xi_0, \eta_0) m_{\widehat{\mathcal{P}}}(\overline{\eta}) m_{\widehat{\mathcal{P}}}(\overline{\sigma}) \widehat{F}_2(s, \overline{\sigma}) ~ d\overline{\eta} ds 
\end{aligned} \label{equdecchampbBHHPBH-1-2} \\
&\begin{aligned} 
+ \int_0^t \int \nabla_{\overline{\eta}} \cdot \Biggl( \left( \xi_0 m_b(\overline{\xi}) \widehat{X}_b(\overline{\xi}) - \sigma_0 m_b(\overline{\sigma}) \widehat{X}_b(\overline{\sigma}) \right) \mu(\overline{\xi}, \overline{\eta}) \mu_{BHH}(\overline{\xi}, \overline{\eta}) \mu_{BH}(\xi_0, \eta_0) \\
m_{\widehat{\mathcal{P}}}(\overline{\eta}) m_{\widehat{\mathcal{P}}}(\overline{\sigma}) \Biggl) e^{i s \varphi} \widehat{F}_1(s, \overline{\eta}) \widehat{F}_2(s, \overline{\sigma}) ~ d\overline{\eta} ds 
\end{aligned} \label{equdecchampbBHHPBH-1-3} 
\end{align}
\end{subequations}
\eqref{equdecchampbBHHPBH-1-3} is of the form \eqref{lemdecequH-11-dersymb}. 

For \eqref{equdecchampbBHHPBH-1-2}, as before, we project: 
\begin{align*}
\xi_0 m_b(\overline{\xi}) \widehat{X}_b(\overline{\xi}) - \sigma_0 m_b(\overline{\sigma}) \widehat{X}_b(\overline{\sigma})
&= O(\eta_0) + O(\sigma_0 \widehat{X}_a(\overline{\eta})) + O(\sigma_0 \widehat{X}_c(\overline{\eta})) + O(\varphi \widehat{X}_b(\overline{\eta})) 
\end{align*}
which allows to decompose into \eqref{lemdecequH-03-resxetaac}, \eqref{lemdecequH-05-resxetabbonbis} and \eqref{lemdecequH-06-resxetabphi}. 

Finally, for \eqref{equdecchampbBHHPBH-1-3}, we first compute that
\begin{align*}
\varphi - \eta_0 \partial_{\eta_0} \varphi &= \xi_0 \left( 3 \eta_0^2 + |\xi|^2 - |\sigma|^2 \right) \\
&= - \xi_0 \left( |\sigma|^2 + o(1) \right) 
\end{align*}
so we know how to decompose $O(\xi_0)$. But here, $\sigma_0 = O(\xi_0)$, so that 
\begin{align*}
&\left( \xi_0 m_b(\overline{\xi}) \widehat{X}_b(\overline{\xi}) \cdot \left( \nabla_{\overline{\xi}} + \nabla_{\overline{\eta}} \right) - \sigma_0 m_b(\overline{\sigma}) \widehat{X}_b(\overline{\sigma}) \cdot \nabla_{\overline{\eta}} \right) \varphi \\
&\quad = O(\xi_0m_b(\overline{\xi})) + O(\sigma_0 m_b(\overline{\xi})) = O(\xi_0 m_b(\overline{\xi}))
\end{align*}
which concludes.

\section{Structure of the weighted profiles} \label{section-structurepoids} 

\subsection{Structure of \texorpdfstring{$h, g$}{h, g}} 

\begin{Lem} Let $\alpha \in \{ a, b, c \}$. Then we can decompose
\begin{align*}
m_{\alpha}(D) X_{\alpha} f(t) &= h_{\alpha}(t) + g_{\alpha}(t) \\
h_{\alpha}(t) &= m_{\alpha}(D) X_{\alpha} f(0) + h_{\alpha, 2}(t) + h_{\alpha, 3}(t) 
\end{align*}
and we have the expressions: 
\begin{align*}
h_{a, 2}(t) &= 2 \frac{\partial_x}{|\nabla|} I\left[ |\nabla| h_a, ~ f \right](t) + \frac{\partial_x}{|\nabla|} I[f, f](t) \\
g_a(t) &= 3 t \frac{\partial_x}{|\nabla|} T[f, f](t) \\
h_{c, 2}(t) &= 2 \frac{\partial_x}{|\nabla|} I\left[ |\nabla| h_c, ~ f \right](t) \\
g_c(t) &= h_{a, 3}(t) = h_{c, 3}(t) = 0
\end{align*}
as well as for $\alpha = b$, in Fourier space: 
\begin{subequations}
\begin{align}
\widehat{h}_{b, 2}(t, \overline{\xi}) 
&= 2 \int_0^t \int e^{i s \varphi} i \eta_0 \widehat{h}_b(s, \overline{\eta}) \widehat{f}(s, \overline{\sigma}) ~ d\overline{\eta} ds \label{lemstructurehb-02-symeta} \\
&+ \sum_{\alpha = a, c} \int_0^t \int e^{i s \varphi} \mu(\overline{\xi}, \overline{\eta}) \frac{|\overline{\eta}| \sigma_0}{|\overline{\eta}| + |\overline{\sigma}|} \widehat{h}_{\alpha}(s, \overline{\eta}) \widehat{f}(s, \overline{\sigma}) ~ d\overline{\eta} ds \label{lemstructurehb-03-resxetaac} \\
&+ \int_0^t \int e^{i s \varphi} \mu(\overline{\xi}, \overline{\eta}) \frac{\eta_0 \sigma_0}{|\overline{\eta}| + |\overline{\sigma}|} \widehat{h}_b(s, \overline{\eta}) \widehat{f}(s, \overline{\sigma}) ~ d\overline{\eta} ds \label{lemstructurehb-04-resxetabbon} \\
&+ \int_0^t \int e^{i s \varphi} \mu(\overline{\xi}, \overline{\eta}) \frac{\eta_0 \sigma_0 |\overline{\eta}|}{(|\eta_0|+|\sigma_0|)(|\overline{\eta}| + |\overline{\sigma}|)} \widehat{h}_b(s, \overline{\eta}) \widehat{f}(s, \overline{\sigma}) ~ d\overline{\eta} ds \label{lemstructurehb-05-resxetabbonbis} \\
&+ \int_0^t \int e^{i s \varphi} \mu(\overline{\xi}, \overline{\eta}) \frac{\varphi m_{\widehat{\mathcal{P}}}(\overline{\eta})}{(|\overline{\eta}| + |\overline{\sigma}|)^2} \widehat{h}_b(s, \overline{\eta}) \widehat{f}(s, \overline{\sigma}) ~ d\overline{\eta} ds \label{lemstructurehb-06-resxetabphi} \\
&+ \int_0^t \int e^{i s \varphi} \mu(\overline{\xi}, \overline{\eta}) \widehat{f}(s, \overline{\eta}) \widehat{f}(s, \overline{\sigma}) ~ d\overline{\eta} \label{lemstructurehb-11-dersymb} \\
\widehat{h}_{b, 3}(t, \overline{\xi}) &= \int_0^t \int \int e^{i s \varphi_3} \mu(\overline{\xi}, \overline{\eta}, \overline{\sigma}) i s \eta_0 \frac{|\overline{\sigma}|+|\overline{\rho}|}{|\overline{\eta}|+|\overline{\sigma}|+|\overline{\rho}|} \widehat{f}(s, \overline{\eta}) \widehat{f}(s, \overline{\sigma}) \widehat{f}(s, \overline{\rho}) ~ d\overline{\eta} d\overline{\sigma} ds \label{lemstructurehb-13-termecubique} \\
\widehat{g}_b(t, \overline{\xi}) &= t \int e^{i t \varphi} \mu(\overline{\xi}, \overline{\eta}) \widehat{f}(t, \overline{\eta}) \widehat{f}(t, \overline{\sigma}) ~ d\overline{\eta} \label{lemstructurehb-14-g} 
\end{align}
\end{subequations} 
where the symbol $\mu$ may change from line to line, and where we recall that $\varphi_3$ stands for the cubic interaction phase with the cubic convention $\overline{\rho} = \overline{\xi} - \overline{\eta}-\overline{\sigma}$. 
\label{lemstructurehalphagalpha} 
\end{Lem}

\begin{Dem}
By Duhamel's formula and the definition of $I$, we have 
\begin{align*}
f(t) &= f(0) + \partial_x I[f, f](t) 
\end{align*}

Let us apply $m_a(D) X_a$ and get: 
\begin{align*}
m_a(D) X_a f(t) &= m_a(D) X_a f(0) + \frac{\partial_x}{|\nabla|} \left( |\nabla| m_a(D) X_a I[f, f](t) \right) + \left[ m_a(D) X_a, \partial_x \right] I[f, f](t) \\
&= m_a(D) X_a f(0) + 
2 \frac{\partial_x}{|\nabla|} I\left[ |\nabla| m_a(D) X_a f, ~ f \right](t) 
+ 3 t \frac{\partial_x}{|\nabla|} T[f, f](t) \\
&\quad - 6 \frac{\partial_x}{|\nabla|} I\left[ t \partial_t f, ~ f \right](t) 
+ \frac{\partial_x}{|\nabla|} I[f, f](t) 
\end{align*}
by the decomposition lemma \ref{lemdecompositionFGH}. We thus already have 
\begin{align*}
g_a(t) &= 3 t \frac{\partial_x}{|\nabla|} T[f, f](t) = 3 t |\nabla|^{-1} \partial_t f(t) 
\end{align*}
We put the rest into $h_a(t)$, and we may decompose further
\begin{align*}
&m_a(D) X_a f(0) + 2 \frac{\partial_x}{|\nabla|} I\left[ |\nabla| m_a(D) X_a f, ~ f \right](t) 
- 6 \frac{\partial_x}{|\nabla|} I\left[ t \partial_t f, ~ f \right](t) 
+ \frac{\partial_x}{|\nabla|} I[f, f](t) \\
&= m_a(D) X_a f(0) + 2 \frac{\partial_x}{|\nabla|} I\left[ |\nabla| h_a, ~ f \right](t) + 2 \frac{\partial_x}{|\nabla|} I\left[ |\nabla| g_a, ~ f \right](t) 
- 6 \frac{\partial_x}{|\nabla|} I\left[ t \partial_t f, ~ f \right](t) 
+ \frac{\partial_x}{|\nabla|} I[f, f](t) \\
&= m_a(D) X_a f(0) + 2 \frac{\partial_x}{|\nabla|} I\left[ |\nabla| h_a, ~ f \right](t) 
+ \frac{\partial_x}{|\nabla|} I[f, f](t)
\end{align*}
as wanted. 

If we apply $m_c(D) X_c$, the reasoning is the same and even simpler, noting that $X_c$ commutes with $\partial_x$: 
\begin{align*}
m_c(D) X_c f(t) &= m_c(D) X_c f(0) + \frac{\partial_x}{|\nabla|} \left( |\nabla| m_c(D) X_c I[f, f](t) \right) \\
&= m_c(D) X_c f(0) + 2 \frac{\partial_x}{|\nabla|} I\left[ |\nabla| m_c(D) X_c f, ~ f \right](t) \\
&= m_c(D) X_c f(0) + 2 \frac{\partial_x}{|\nabla|} I\left[ |\nabla| h_c, ~ f \right](t) 
\end{align*}
and we set $g_c(t) = 0$. 

Finally, for $m_b(D) X_b$, we start by commuting: 
\begin{align*}
m_b(D) X_b f(t) &= m_b(D) X_b f(0) + \partial_x m_b(D) X_b I[f, f](t) + \frac{i |\nabla_y|}{|\nabla|} m_b(D) I[f, f](t) 
\end{align*}
The term $\frac{i |\nabla_y|}{|\nabla|} m_b(D) I[f, f](t)$ is of the form \eqref{lemstructurehb-11-dersymb}. For the other, we can apply the decomposition lemma \ref{lemdecompositionFGH} and get: 
\begin{subequations}
\begin{align}
&i \xi_0 m_b(\overline{\xi}) \widehat{X}_b(\overline{\xi}) \cdot \nabla_{\overline{\xi}} \widehat{I}[f, f](t, \overline{\xi}) \notag \\ 
&= 2 \int_0^t \int e^{i s \varphi} i \eta_0 m_b(\overline{\eta}) \widehat{X}_b(\overline{\eta}) \cdot \nabla_{\overline{\eta}} \widehat{f}(s, \overline{\eta}) \widehat{f}(s, \overline{\sigma}) ~ d\overline{\eta} ds \label{preuvelemstructurehb-01-symeta} \\
&+ \sum_{\alpha = a, c} \int_0^t \int e^{i s \varphi} \mu(\overline{\xi}, \overline{\eta}) \frac{|\overline{\eta}| \sigma_0}{|\overline{\eta}| + |\overline{\sigma}|} m_{\alpha}(\overline{\eta}) \widehat{X}_{\alpha}(\overline{\eta}) \cdot \nabla_{\overline{\eta}} \widehat{f}(s, \overline{\eta}) \widehat{f}(s, \overline{\sigma}) ~ d\overline{\eta} ds \label{preuvelemstructurehb-03-resxetaac} \\
&+ \int_0^t \int e^{i s \varphi} \mu(\overline{\xi}, \overline{\eta}) \frac{\eta_0 \sigma_0}{|\overline{\eta}| + |\overline{\sigma}|} m_b(\overline{\eta}) \widehat{X}_b(\overline{\eta}) \cdot \nabla_{\overline{\eta}} \widehat{f}(s, \overline{\eta}) \widehat{f}(s, \overline{\sigma}) ~ d\overline{\eta} ds \label{preuvelemstructurehb-04-resxetabbon} \\
&+ \int_0^t \int e^{i s \varphi} \mu(\overline{\xi}, \overline{\eta}) \frac{\eta_0 \sigma_0 |\overline{\eta}|}{(|\eta_0|+|\sigma_0|)(|\overline{\eta}| + |\overline{\sigma}|)} m_b(\overline{\eta}) \widehat{X}_b(\overline{\eta}) \cdot \nabla_{\overline{\eta}} \widehat{f}(s, \overline{\eta}) \widehat{f}(s, \overline{\sigma}) ~ d\overline{\eta} ds \label{preuvelemstructurehb-05-resxetabbonbis} \\
&+ \int_0^t \int e^{i s \varphi} \mu(\overline{\xi}, \overline{\eta}) \frac{\varphi m_{\widehat{\mathcal{P}}}(\overline{\eta})}{(|\overline{\eta}| + |\overline{\sigma}|)^2} m_b(\overline{\eta}) \widehat{X}_b(\overline{\eta}) \cdot \nabla_{\overline{\eta}} \widehat{f}(s, \overline{\eta}) \widehat{f}(s, \overline{\sigma}) ~ d\overline{\eta} ds \label{preuvelemstructurehb-06-resxetabphi} \\
&+ \int_0^t \int e^{i s \varphi} \mu(\overline{\xi}, \overline{\eta}) \widehat{f}(s, \overline{\eta}) \widehat{f}(s, \overline{\sigma}) ~ d\overline{\eta} ds \label{preuvelemstructurehb-11-dersymb} \\
&+ \int_0^t \int s e^{i s \varphi} \left( m_g(\overline{\xi}, \overline{\eta}) + m_g(\overline{\xi}, \overline{\sigma}) \right) \partial_s \widehat{f}(s, \overline{\eta}) \widehat{f}(s, \overline{\sigma}) ~ d\overline{\eta} ds \label{preuvelemstructurehb-12-resteta} \\
&- \int t e^{i t \varphi} m_g(\overline{\xi}, \overline{\eta}) \widehat{f}(t, \overline{\eta}) \widehat{f}(t, \overline{\sigma}) ~ d\overline{\eta} \label{preuvelemstructurehb-g} 
\end{align}
\end{subequations}
by the symmetry of $\overline{\eta}$ and $\overline{\sigma}$ in this situation. 
Above, we can symmetrize: 
\begin{align}
\eqref{preuvelemstructurehb-g} &= - \frac{1}{2} \int t e^{i t \varphi} \left( m_g(\overline{\xi}, \overline{\eta}) + m_g(\overline{\xi}, \overline{\sigma}) \right) \widehat{f}(t, \overline{\eta}) \widehat{f}(t, \overline{\sigma}) ~ d\overline{\eta} \label{preuvelemstructurehb-gsym} 
\end{align}
We then define $\widehat{g}_b$ as only \eqref{preuvelemstructurehb-g}, which has indeed the desired structure \eqref{lemstructurehb-14-g}, and the rest as $h_b(t)$. To obtain the formula for $h_b(t)$, we develop everywhere $m_{\alpha} X_{\alpha} f(t) = h_{\alpha}(t) + g_{\alpha}(t)$: 
\begin{subequations}
\begin{align}
\eqref{preuvelemstructurehb-03-resxetaac} &= \sum_{\alpha = a, c} \int_0^t \int e^{i s \varphi} \mu(\overline{\xi}, \overline{\eta}) \frac{|\overline{\eta}| \sigma_0}{|\overline{\eta}| + |\overline{\sigma}|} \widehat{h}_{\alpha}(s, \overline{\eta}) \widehat{f}(s, \overline{\sigma}) ~ d\overline{\eta} ds \label{preuvelemstructurehb-03-1} \\
&+ \sum_{\alpha = a, c} \int_0^t \int e^{i s \varphi} \mu(\overline{\xi}, \overline{\eta}) \frac{|\overline{\eta}| \sigma_0}{|\overline{\eta}| + |\overline{\sigma}|} \widehat{g}_{\alpha}(s, \overline{\eta}) \widehat{f}(s, \overline{\sigma}) ~ d\overline{\eta} ds \label{preuvelemstructurehb-03-2} 
\end{align}
\end{subequations} 
\eqref{preuvelemstructurehb-03-1} is of the form \eqref{lemstructurehb-03-resxetaac}, \eqref{preuvelemstructurehb-03-2} is of the form \eqref{lemstructurehb-13-termecubique} by developing $g_{\alpha}(t)$ as a quadratic term in $f$ and changing the variables. The same way, developing \eqref{preuvelemstructurehb-04-resxetabbon}, \eqref{preuvelemstructurehb-05-resxetabbonbis}, \eqref{preuvelemstructurehb-06-resxetabphi} respectively, we get in a term of the form \eqref{lemstructurehb-04-resxetabbon}, \eqref{lemstructurehb-05-resxetabbonbis}, \eqref{lemstructurehb-06-resxetabphi} respectively plus in any case a term of the form \eqref{lemstructurehb-13-termecubique}, using in the last one that 
\begin{align*}
\varphi &= O(\eta_0 \overline{\sigma}) + O(\sigma_0 \overline{\eta}) 
\end{align*}
by Lemma \ref{lemcalculvarphiKdV1D}. \eqref{preuvelemstructurehb-11-dersymb} is already of the form \eqref{lemstructurehb-11-dersymb}. 

It thus only remains
\begin{subequations}
\begin{align}
&\eqref{preuvelemstructurehb-01-symeta} + \eqref{preuvelemstructurehb-12-resteta} = 2 \int_0^t \int e^{i s \varphi} i \eta_0 \widehat{h}_b(s, \overline{\eta}) \widehat{f}(s, \overline{\sigma}) ~ d\overline{\eta} ds + 2 \int_0^t \int e^{i s \varphi} i \eta_0 \widehat{g}_b(s, \overline{\eta}) \widehat{f}(s, \overline{\sigma}) ~ d\overline{\eta} ds \notag \\
&\pushright{+ \int_0^t \int s e^{i s \varphi} \left( m_g(\overline{\xi}, \overline{\eta}) + m_g(\overline{\xi}, \overline{\sigma}) \right) \partial_s \widehat{f}(s, \overline{\eta}) \widehat{f}(s, \overline{\sigma}) ~ d\overline{\eta} ds} \notag \\
&= 2 \int_0^t \int e^{i s \varphi} i \eta_0 \widehat{h}_b(s, \overline{\eta}) \widehat{f}(s, \overline{\sigma}) ~ d\overline{\eta} ds \label{preuvelemstructurehb-01-1} \\
&\begin{aligned}
+ \int_0^t \int \int e^{i s \varphi_3} i s \left( m_g(\overline{\xi}, \overline{\eta}+\overline{\sigma}) + m_g(\overline{\xi}, \overline{\rho}) - m_g(\overline{\eta}+\overline{\sigma}, \overline{\eta}) - m_g(\overline{\eta}+\overline{\sigma}, \overline{\sigma}) \right) \\
(\eta_0+\sigma_0) \widehat{f}(s, \overline{\eta}) \widehat{f}(s, \overline{\sigma}) \widehat{f}(s, \overline{\rho}) ~ d\overline{\eta} ds 
\end{aligned} \label{preuvelemstructurehb-01-2} 
\end{align}
\end{subequations}
where we developed using \eqref{preuvelemstructurehb-gsym}. \eqref{preuvelemstructurehb-01-1} corresponds exactly to \eqref{lemstructurehb-02-symeta}. Finally, for \eqref{preuvelemstructurehb-01-2}, we can rewrite
\begin{align*}
\eqref{preuvelemstructurehb-01-2}= 2 \int_0^t \int \int e^{i s \varphi_3} i s \eta_0 \left( m_g(\overline{\xi}, \overline{\eta}+\overline{\sigma}) + m_g(\overline{\xi}, \overline{\rho}) - m_g(\overline{\eta}+\overline{\sigma}, \overline{\eta}) - m_g(\overline{\eta}+\overline{\sigma}, \overline{\sigma}) \right) \\
\widehat{f}(s, \overline{\eta}) \widehat{f}(s, \overline{\sigma}) \widehat{f}(s, \overline{\rho}) ~ d\overline{\eta} ds
\end{align*}
by permutation of the variables. If we localise in order to have $|\overline{\sigma}|+|\overline{\rho}| \gtrsim |\overline{\eta}|$, we have automatically the desired form; else, $|\overline{\sigma}|+|\overline{\rho}| \ll |\overline{\eta}|$ and we can develop by Taylor's formula: 
\begin{align*}
&m_g(\overline{\xi}, \overline{\eta}+\overline{\sigma}) + m_g(\overline{\xi}, \overline{\rho}) - m_g(\overline{\eta}+\overline{\sigma}, \overline{\eta}) - m_g(\overline{\eta}+\overline{\sigma}, \overline{\sigma}) \\
&= m_g(\overline{\xi}, \overline{\xi}) + m_g(\overline{\xi}, 0) - m_g(\overline{\xi}, \overline{\xi}) - m_g(\overline{\xi}, 0) + O(\overline{\sigma}) + O(\overline{\rho}) \\
&= O(\overline{\sigma}) + O(\overline{\rho}) 
\end{align*}
which is what we wanted. 
\end{Dem}

\subsection{Quadratic structure} 

\begin{Lem} Let $\alpha \in \{ a, b, c \}$, $\beta \in \{ a, c \}$. Then we can decompose
\begin{align*}
|\nabla| m_{\alpha}(D) X_{\alpha} h_{\beta}(t) &= h_{\alpha \beta}(t) + g_{\alpha \beta}(t) \\
h_{\alpha \beta}(t) &= |\nabla| m_{\alpha}(D) X_{\alpha} m_{\beta}(D) X_{\beta} f(0) + h_{\alpha \beta, 2}(t) + h_{\alpha \beta, 3}(t)  + h_{\alpha \beta, 4}(t) 
\end{align*}
with the following expressions: 
\begin{subequations}
\begin{align}
\widehat{h}_{\alpha \beta, 2}(t, \overline{\xi}) &= 2 m_{0, \alpha}(\overline{\xi}) \int_0^t \int e^{i s \varphi} m_{0, \overline{\alpha}}(\overline{\eta}) |\overline{\eta}| \widehat{h}_{\alpha \beta}(s, \overline{\eta}) \widehat{f}(s, \overline{\sigma}) ~ d\overline{\sigma} ds \label{lemstructhquadac-2-sym} \\
&\quad + \int_0^t \int e^{i s \varphi} \mu(\overline{\xi}, \overline{\eta}) \widehat{h}_{b \beta}(s, \overline{\eta}) \sigma_0 \widehat{f}(s, \overline{\sigma}) ~ d\overline{\eta} ds \label{lemstructhquadac-2-doubleresx} \\
&\quad + \int_0^t \int e^{i s \varphi} \mu(\overline{\xi}, \overline{\eta}) \frac{\varphi}{(|\overline{\eta}| + |\overline{\sigma}|)^2} \widehat{h}_{b \beta}(s, \overline{\eta}) \widehat{f}(s, \overline{\sigma}) ~ d\overline{\eta} ds \label{lemstructhquadac-2-doubleresxvarphi} \\
&+ \sum_{\gamma = a, c} \int_0^t \int e^{i s \varphi} \mu(\overline{\xi}, \overline{\eta}) |\overline{\eta}| \widehat{h}_{\beta}(s, \overline{\eta}) |\overline{\sigma}| \widehat{h}_{\gamma}(s, \overline{\sigma}) ~ d\overline{\eta} ds \label{lemstructhquadac-2-resxresx} \\
&+ \int_0^t \int e^{i s \varphi} \mu(\overline{\xi}, \overline{\eta}) |\overline{\eta}| \widehat{h}_{\beta}(s, \overline{\eta}) \sigma_0 \widehat{h}_b(s, \overline{\sigma}) ~ d\overline{\eta} ds \label{lemstructhquadac-2-resxresxb} \\
&+ \int_0^t \widetilde{\chi}(s) \int e^{i s \varphi} \mu(\overline{\xi}, \overline{\eta}) \langle (|\overline{\eta}|+|\overline{\sigma}|)^{-1} \rangle \langle \overline{\eta} \rangle \widehat{h}_{\beta}(s, \overline{\eta}) \langle \overline{\sigma} \rangle \widehat{h}_b(s, \overline{\sigma}) ~ d\overline{\eta} ds \label{lemstructhquadac-2-resxresxvarphi} \\
&+ \sum_{\gamma = a, c} \int_0^t \int e^{i s \varphi} \mu(\overline{\xi}, \overline{\eta}) |\overline{\eta}| \widehat{h}_{\gamma}(s, \overline{\eta}) \widehat{f}(s, \overline{\sigma}) ~ d\overline{\eta} ds \label{lemstructhquadac-2-dersymb-resx} \\
&+ \int_0^t \int e^{i s \varphi} \mu(\overline{\xi}, \overline{\eta}) \eta_0 \widehat{h}_b(s, \overline{\eta}) \widehat{f}(s, \overline{\sigma}) ~ d\overline{\eta} ds \label{lemstructhquadac-2-dersymb-resxb} \\
&+ \sum_{\gamma = a, b, c} \int_0^t \int e^{i s \varphi} \mu(\overline{\xi}, \overline{\eta}) \widehat{h}_{\gamma}(s, \overline{\eta}) \sigma_0 \widehat{f}(s, \overline{\sigma}) ~ d\overline{\eta} ds \label{lemstructhquadac-2-dersymb-resxbis} \\
&+ \int_0^t \int e^{i s \varphi} \mu(\overline{\xi}, \overline{\eta}) \widehat{f}(s, \overline{\eta}) \widehat{f}(s, \overline{\sigma}) ~ d\overline{\eta} ds \label{lemstructhquadac-2-doubledersymb} 
 \\
\widehat{h}_{\alpha \beta, 3}(t, \overline{\xi}) &= \int_0^t \int \int e^{i s \varphi_3(\overline{\xi}, \overline{\eta}, \overline{\sigma})} s \mu(\overline{\xi}, \overline{\eta}, \overline{\sigma}) |\overline{\eta}| \widehat{h}_{\beta}(s, \overline{\eta}) |\overline{\sigma}| \widehat{f}(s, \overline{\sigma}) \widehat{f}(s, \overline{\rho}) ~ d\overline{\eta} d\overline{\sigma} ds \label{lemstructhquadac-3-resx} \\
&+ \int_0^t \int \int e^{i s \varphi_3(\overline{\xi}, \overline{\eta}, \overline{\sigma})} s \mu(\overline{\xi}, \overline{\eta}, \overline{\sigma}) |\overline{\eta}| \widehat{f}(s, \overline{\eta}) \widehat{f}(s, \overline{\sigma}) \widehat{f}(s, \overline{\rho}) ~ d\overline{\eta} d\overline{\sigma} ds \label{lemstructhquadac-3-dersymb} \\
&\quad + \sum_{\gamma_1, \gamma_2 = a, b, c} \int_0^t \int \int e^{i s \varphi_3(\overline{\xi}, \overline{\eta}, \overline{\sigma})} \mu(\overline{\xi}, \overline{\eta}, \overline{\sigma}) |\overline{\eta}|^{-1} \widehat{h}_{\gamma_1}(s, \overline{\eta}) \widehat{h}_{\gamma_2}(s, \overline{\sigma}) |\overline{\rho}| \widehat{f}(s, \overline{\rho}) ~ d\overline{\eta} d\overline{\sigma} ds \label{lemstructhquadac-3-varphi} \\
&\quad \begin{aligned}
+ \sum_{\gamma_1, \gamma_2 = a, b, c} \int_0^t \int \int e^{i s \varphi_3(\overline{\xi}, \overline{\eta}, \overline{\sigma})} \mu(\overline{\xi}, \overline{\eta}, \overline{\sigma}) |\overline{\eta}|^{-1} \widehat{h}_{\gamma_1}(s, \overline{\eta}) m_{0, \overline{\gamma_2}}(\overline{\sigma}) |\overline{\sigma}| \\
\widehat{h}_{\gamma_2}(s, \overline{\sigma}) \widehat{f}(s, \overline{\rho}) ~ d\overline{\eta} d\overline{\sigma} ds 
\end{aligned} \label{lemstructhquadac-3-varphibis} \\
&\quad + \sum_{\gamma = a, b, c} \int_0^t \int \int e^{i s \varphi_3(\overline{\xi}, \overline{\eta}, \overline{\sigma})} \mu(\overline{\xi}, \overline{\eta}, \overline{\sigma}) |\overline{\eta}|^{-1} \widehat{h}_{\gamma}(s, \overline{\eta}) \widehat{f}(s, \overline{\sigma}) \widehat{f}(s, \overline{\rho}) ~ d\overline{\eta} ds \label{lemstructhquadac-3-varphidersymb} \\
\widehat{h}_{\alpha \beta, 4}(t, \overline{\xi}) &= \sum_{\gamma = a, b, c} \int_0^t \int \int \int e^{i s \varphi_4(\overline{\xi}, \overline{\eta}^1, \overline{\eta}^2, \overline{\eta}^3)} \mu(\overline{\xi}, \overline{\eta}^1, \overline{\eta}^2, \overline{\eta}^3) s |\overline{\eta}^1|^{-1} \widehat{h}_{\gamma}(s, \overline{\eta}^1) |\overline{\eta}^2| \widehat{f}(s, \overline{\eta}^2) \notag \\
&\quad \quad \quad \quad \quad \quad \quad \quad \quad \quad \quad \quad \quad \quad \quad \quad \quad \quad \widehat{f}(s, \overline{\eta}^3) \widehat{f}(s, \overline{\eta}^4) ~ d\overline{\eta}^1 d\overline{\eta}^2 d\overline{\eta}^3 ds \label{lemstructhquadac-4-dersymb} \\
\widehat{g}_{\alpha \beta}(t, \overline{\xi}) &= t \int e^{i t \varphi} \mu(\overline{\xi}, \overline{\eta}) m_{0, \overline{\beta}}(\overline{\eta}) |\overline{\eta}| \widehat{h}_{\beta}(t, \overline{\eta}) \widehat{f}(t, \overline{\sigma}) ~ d\overline{\eta} \label{lemstructhquadac-gresx} \\
&\quad + t \int e^{i t \varphi} \mu(\overline{\xi}, \overline{\eta}) \widehat{f}(t, \overline{\eta}) \widehat{f}(t, \overline{\sigma}) ~ d\overline{\eta} \label{lemstructhquadac-gdersymb} \\
&\quad + (1 - \chi(t)) \int e^{i t \varphi} \mu(\overline{\xi}, \overline{\eta}) \frac{1}{|\overline{\eta}| + |\overline{\sigma}|} \widehat{h}_{\beta}(t, \overline{\eta}) \widehat{h}_b(t, \overline{\sigma}) ~ d\overline{\eta} \label{lemstructhquadac-gvarphi} 
\end{align}
\end{subequations} 
where the symbol $\mu$ can change from line to line, and where we denoted by $m_{0, \alpha}(\overline{\xi}) = m_{0, \overline{b}}(\overline{\xi}) = \frac{i \xi_0}{|\overline{\xi}|}$, $m_{0, b}(\overline{\xi}) = m_{0, \overline{\alpha}}(\overline{\xi}) = 1$ for every $\alpha \neq b$, and where $\widetilde{\chi}$ is some compactly supported function and $\chi$ the usual Littlewood-Paley localisation function. \label{lemmestructurehquadabcsurac}
\end{Lem}

\begin{Dem}
By Lemma \ref{lemstructurehalphagalpha}, we have for $\beta \in \{ a, c \}$: 
\begin{align*}
h_{\beta}(t) &= m_{\beta}(D) X_{\beta} f(0) + 2 \frac{\partial_x}{|\nabla|} I\left[ |\nabla| h_{\beta}, ~ f \right](t) + c_{\beta} \frac{\partial_x}{|\nabla|} I[f, f](t)
\end{align*}
for some constant $c_{\beta}$ ($c_a = 1, c_c = 0$). 

In particular, we can develop
\begin{align*}
&|\nabla| m_{\alpha}(D) X_{\alpha} h_{\beta}(t) = |\nabla| m_{\alpha}(D) X_{\alpha} m_{\beta}(D) X_{\beta} f(0) 
+ 2 \partial_x m_{\alpha}(D) X_{\alpha} I\left[ |\nabla| h_{\beta}, ~ f \right](t) \\
&+ c_{\beta} \partial_x m_{\alpha}(D) X_{\alpha} I[f, f](t) + \left[ |\nabla| m_{\alpha}(D) X_{\alpha}, ~ \frac{\partial_x}{|\nabla|} \right] \left( 2 I\left[ |\nabla| h_{\beta}, ~ f \right](t) + c_{\beta} I[f, f](t) \right) 
\end{align*}
We can then apply the decomposition Lemma to develop all terms. 

Note first that $\left[ |\nabla| m_{\alpha}(D) X_{\alpha}, ~ \frac{\partial_x}{|\nabla|} \right]$ is a Hörmander-Mikhlin type symbol, so the second line above corresponds to \eqref{lemstructhquadac-2-dersymb-resx} and \eqref{lemstructhquadac-2-doubledersymb}. 

If $\alpha = a$, we have {\scriptsize 
\begin{subequations}
\begin{align}
&2 \partial_x m_a(D) X_a I\left[ |\nabla| h_{\beta}, ~ f \right](t) + c_{\beta} \partial_x m_a(D) X_a I[f, f](t) \notag \\
&\quad = 2 \frac{\partial_x}{|\nabla|} I\left[ |\nabla|^2 X_a h_{\beta}, ~ f \right](t) 
+ 2 \frac{\partial_x}{|\nabla|} I\left[ |\nabla| h_{\beta}, ~ f \right](t) 
+ 2 \frac{\partial_x}{|\nabla|} I\left[ |\nabla| h_{\beta}, ~ |\nabla| X_a f \right](t) 
+ 3 t \frac{\partial_x}{|\nabla|} T\left[ |\nabla| h_{\beta}, ~ f \right](t) \notag \\
&\quad \quad - 3 \frac{\partial_x}{|\nabla|} I\left[ t |\nabla| \partial_t h_{\beta}, ~ f \right](t) 
- 3 \frac{\partial_x}{|\nabla|} I\left[ |\nabla| h_{\beta}, ~ t \partial_t f \right](t) 
+ 2 c_{\beta} \frac{\partial_x}{|\nabla|} I\left[ |\nabla| X_a f, ~ f \right](t) 
+ 3 t \frac{\partial_x}{|\nabla|} T[f, f](t) 
- 6 \frac{\partial_x}{|\nabla|} I\left[ t \partial_t f, ~ f \right](t) \notag \\
&\quad = 2 \frac{\partial_x}{|\nabla|} I\left[ |\nabla| h_{\alpha \beta}, ~ f \right](t) 
+ 2 \frac{\partial_x}{|\nabla|} I\left[ |\nabla| h_{\beta}, ~ f \right](t)
+ 2 \frac{\partial_x}{|\nabla|} I\left[ |\nabla| h_{\beta}, ~ |\nabla| h_a \right](t) 
+ 2 c_{\beta} \frac{\partial_x}{|\nabla|} I\left[ |\nabla| h_a, ~ f \right](t) \label{developpementlemstructhquad-a-2ok} \\
&\quad \quad + 2 \frac{\partial_x}{|\nabla|} I\left[ |\nabla| g_{\alpha \beta}, ~ f \right](t)
+ 2 \frac{\partial_x}{|\nabla|} I\left[ |\nabla| h_{\beta}, ~ |\nabla| g_a \right](t) 
- 6 \frac{\partial_x}{|\nabla|} I\left[ t \partial_x T\left[ |\nabla| h_{\beta}, ~ f \right], ~ f \right](t) \label{developpementlemstructhquad-a-3comp} \\
&\quad \quad - 3 c_{\beta} \frac{\partial_x}{|\nabla|} I\left[ t \partial_x T[f, f], ~ f \right](t) 
- 3 \frac{\partial_x}{|\nabla|} I\left[ |\nabla| h_{\beta}, ~ t \partial_x T[f, f] \right](t) 
+ 2 c_{\beta} \frac{\partial_x}{|\nabla|} I\left[ |\nabla| g_a, ~ f \right](t) 
- 6 \frac{\partial_x}{|\nabla|} I\left[ t \partial_x T[f, f], ~ f \right](t) \label{developpementlemstructhquad-a-3ok} \\
&\quad \quad
+ 3 t \frac{\partial_x}{|\nabla|} T\left[ |\nabla| h_{\beta}, ~ f \right](t) 
+ 3 t \frac{\partial_x}{|\nabla|} T[f, f](t) \label{developpementlemstructhquad-a-g}
\end{align}
\end{subequations} }
We can readily identify \eqref{developpementlemstructhquad-a-g} as $\widehat{g}_{a \beta}(t)$, both terms being of the form \eqref{lemstructhquadac-gresx} et \eqref{lemstructhquadac-gdersymb}. Furthermore, the terms of \eqref{developpementlemstructhquad-a-2ok} correspond already to \eqref{lemstructhquadac-2-sym} and terms from \eqref{lemstructhquadac-2-resxresx}, \eqref{lemstructhquadac-2-dersymb-resx}, \eqref{lemstructhquadac-2-doubledersymb}. Finally, \eqref{developpementlemstructhquad-a-3ok} and the central term from line \eqref{developpementlemstructhquad-a-3comp} are of the form \eqref{lemstructhquadac-3-resx} and \eqref{lemstructhquadac-3-dersymb} (up to developing $g_a$ using Lemma \ref{lemstructurehalphagalpha}). It thus only remains: 
\begin{align*}
&2 \frac{\partial_x}{|\nabla|} I\left[ |\nabla| g_{\alpha \beta}, ~ f \right](t)
- 6 \frac{\partial_x}{|\nabla|} I\left[ t \partial_x T\left[ |\nabla| h_{\beta}, ~ f \right], ~ f \right](t) \\
&= 6 \frac{\partial_x}{|\nabla|} I\left[ t \partial_x T\left[ |\nabla| h_{\beta}, ~ f \right], ~ f \right](t)
+ 6 \frac{\partial_x}{|\nabla|} I\left[ t \partial_x T\left[ f, ~ f \right], ~ f \right](t)
- 6 \frac{\partial_x}{|\nabla|} I\left[ t \partial_x T\left[ |\nabla| h_{\beta}, ~ f \right], ~ f \right](t) \\
&= 6 \frac{\partial_x}{|\nabla|} I\left[ t \partial_x T\left[ f, ~ f \right], ~ f \right](t)
\end{align*}
injecting the exact formula found above for $g_{\alpha \beta}$, that is \eqref{developpementlemstructhquad-a-g}. We get a term of the form \eqref{lemstructhquadac-3-dersymb}. 

If $\alpha = c$, we can proceed to the same development which is simpler: 
\begin{align*}
&2 \partial_x m_c(D) X_c I\left[ |\nabla| h_{\beta}, ~ f \right](t) + c_{\beta} \partial_x m_c(D) X_c I[f, f](t) \\
&= 2 \frac{\partial_x}{|\nabla|} I\left[ |\nabla|^2 m_c(D) X_c h_{\beta}, ~ f \right](t)
+ 2 \frac{\partial_x}{|\nabla|} I\left[ |\nabla| h_{\beta}, ~ |\nabla| m_c(D) X_c f \right](t)
+ 2 c_{\beta} \frac{\partial_x}{|\nabla|} I\left[ |\nabla| m_c(D) X_c f, ~ f \right](t) 
\end{align*}
We may then identify $g_{c \beta}(t) = 0$, and thus 
\begin{align*}
&2 \partial_x m_c(D) X_c I\left[ |\nabla| h_{\beta}, ~ f \right](t) + c_{\beta} \partial_x m_c(D) X_c I[f, f](t) \\
&\quad = 2 \frac{\partial_x}{|\nabla|} I\left[ |\nabla| h_{c \beta}, ~ f \right](t)
+ 2 \frac{\partial_x}{|\nabla|} I\left[ |\nabla| h_{\beta}, ~ |\nabla| h_c \right](t)
+ 2 c_{\beta} \frac{\partial_x}{|\nabla|} I\left[ |\nabla| h_c, ~ f \right](t) 
\end{align*}
which corresponds exactly to \eqref{lemstructhquadac-2-sym} and terms of the form \eqref{lemstructhquadac-2-resxresx} and \eqref{lemstructhquadac-2-dersymb-resx}. 

Finally, if $\alpha = b$, we develop in Fourier: 
\begin{subequations}
\begin{align}
&2 i \xi_0 m_b(\overline{\xi}) \widehat{X}_b(\overline{\xi}) \cdot \nabla_{\overline{\xi}} \widehat{I}\left[ |\nabla| h_{\beta}, ~ f \right](t, \overline{\xi}) + c_{\beta} i \xi_0 m_b(\overline{\xi}) \widehat{X}_b(\overline{\xi}) \cdot \nabla_{\overline{\xi}} \widehat{I}[f, f](t, \overline{\xi}) \notag \\
&\quad = 2 \int_0^t \int e^{i s \varphi} i \eta_0 |\overline{\eta}| m_b(\overline{\eta}) \widehat{X}_b(\overline{\eta}) \cdot \nabla_{\overline{\eta}} \widehat{h}_{\beta}(s, \overline{\eta}) \widehat{f}(s, \overline{\sigma}) ~ d\overline{\eta} ds \label{developpementlemstructquadhac-1} \\
&\quad \quad + \sum_{\gamma = a, b, c} \int_0^t \int e^{i s \varphi} \mu(\overline{\xi}, \overline{\eta}) |\overline{\eta}| m_{\gamma}(\overline{\eta}) \widehat{X}_{\gamma}(\overline{\eta}) \cdot \nabla_{\overline{\eta}} \widehat{h}_{\beta}(s, \overline{\eta}) \sigma_0 \widehat{f}(s, \overline{\sigma}) ~ d\overline{\eta} ds \label{developpementlemstructquadhac-2} \\
&\quad \quad + \int_0^t \int e^{i s \varphi} \mu(\overline{\xi}, \overline{\eta}) \frac{\varphi |\overline{\eta}| m_{\widehat{\mathcal{P}}}(\overline{\eta})}{(|\overline{\eta}| + |\overline{\sigma}|)^2} m_b(\overline{\eta}) \widehat{X}_b(\overline{\eta}) \cdot \nabla_{\overline{\eta}} \widehat{h}_{\beta}(s, \overline{\eta}) \widehat{f}(s, \overline{\sigma}) ~ d\overline{\eta} ds \label{developpementlemstructquadhac-3} \\
&\quad \quad + \int_0^t \int e^{i s \varphi} \mu(\overline{\xi}, \overline{\eta}) |\overline{\eta}| \widehat{h}_{\beta}(s, \overline{\eta}) \sigma_0 m_b(\overline{\sigma}) \widehat{X}_b(\overline{\sigma}) \cdot \nabla_{\overline{\xi}} \widehat{f}(s, \overline{\sigma}) ~ d\overline{\eta} ds \label{developpementlemstructquadhac-4} \\
&\quad \quad + \sum_{\gamma = a, c} \int_0^t \int e^{i s \varphi} \mu(\overline{\xi}, \overline{\eta}) |\overline{\eta}| \widehat{h}_{\beta}(s, \overline{\eta}) |\overline{\sigma}| m_{\gamma}(\overline{\sigma}) \widehat{X}_{\gamma}(\overline{\sigma}) \cdot \nabla_{\overline{\eta}} \widehat{f}(s, \overline{\sigma}) ~ d\overline{\eta} ds \label{developpementlemstructquadhac-5} \\
&\quad \quad + \int_0^t \int e^{i s \varphi} \mu(\overline{\xi}, \overline{\eta}) \frac{\varphi |\overline{\eta}| m_{\widehat{\mathcal{P}}}(\overline{\sigma})}{(|\overline{\eta}| + |\overline{\sigma}|)^2} \widehat{h}_{\beta}(s, \overline{\eta}) m_b(\overline{\sigma}) \widehat{X}_b(\overline{\sigma}) \cdot \nabla_{\overline{\eta}} \widehat{f}(s, \overline{\sigma}) ~ d\overline{\eta} ds \label{developpementlemstructquadhac-6} \\
&\quad \quad + \int_0^t \int e^{i s \varphi} \mu(\overline{\xi}, \overline{\eta}) |\overline{\eta}| \widehat{h}_{\beta}(s, \overline{\eta}) \widehat{f}(s, \overline{\sigma}) ~ d\overline{\eta} ds \label{developpementlemstructquadhac-7} \\
&\quad \quad + \int_0^t \int s e^{i s \varphi} m_g(\overline{\xi}, \overline{\eta}) |\overline{\eta}| \partial_s \widehat{h}_{\beta}(s, \overline{\eta}) \widehat{f}(s, \overline{\sigma}) ~ d\overline{\eta} ds \label{developpementlemstructquadhac-8} \\
&\quad \quad + \int_0^t \int s e^{i s \varphi} m_g(\overline{\xi}, \overline{\eta}) |\overline{\eta}| \widehat{h}_{\beta}(s, \overline{\eta}) \partial_s \widehat{f}(s, \overline{\sigma}) ~ d\overline{\eta} ds \label{developpementlemstructquadhac-9} \\
&\quad \quad - \int t e^{i t \varphi} m_g(\overline{\xi}, \overline{\eta}) |\overline{\eta}| \widehat{h}_{\beta}(t, \overline{\eta}) \widehat{f}(t, \overline{\sigma}) ~ d\overline{\eta} \label{developpementlemstructquadhac-10} \\
&\quad \quad + \int_0^t \int e^{i s \varphi} \mu(\overline{\xi}, \overline{\eta}) \eta_0 m_b(\overline{\eta}) \widehat{X}_b(\overline{\eta}) \cdot \nabla_{\overline{\eta}} \widehat{f}(s, \overline{\eta}) \widehat{f}(s, \overline{\sigma}) ~ d\overline{\eta} ds \label{developpementlemstructquadhac-11} \\
&\quad \quad + \sum_{\gamma = a, b, c} \int_0^t \int e^{i s \varphi} \mu(\overline{\xi}, \overline{\eta}) m_{\gamma}(\overline{\eta}) \widehat{X}_{\gamma}(\overline{\eta}) \cdot \nabla_{\overline{\eta}} \widehat{f}(s, \overline{\eta}) \sigma_0 \widehat{f}(s, \overline{\sigma}) ~ d\overline{\eta} ds \label{developpementlemstructquadhac-12} \\
&\quad \quad + \int_0^t \int e^{i s \varphi} \mu(\overline{\xi}, \overline{\eta}) \widehat{f}(s, \overline{\eta}) \widehat{f}(s, \overline{\sigma}) ~ d\overline{\eta} ds \label{developpementlemstructquadhac-13} \\
&\quad \quad + \int_0^t \int s e^{i s \varphi} \mu(\overline{\xi}, \overline{\eta}) \partial_s \widehat{f}(s, \overline{\eta}) \widehat{f}(s, \overline{\sigma}) ~ d\overline{\eta} ds \label{developpementlemstructquadhac-14} \\
&\quad \quad + \frac{\xi_0}{|\overline{\xi}|} \int t e^{i t \varphi} \mu(\overline{\xi}, \overline{\eta}) \widehat{f}(t, \overline{\eta}) \widehat{f}(t, \overline{\sigma}) ~ d\overline{\eta} \label{developpementlemstructquadhac-15} 
\end{align}
\end{subequations} 
In order to start with defining $g_{b \beta}$, we can separate in \eqref{developpementlemstructquadhac-6}: 
\begin{subequations}
\begin{align}
\eqref{developpementlemstructquadhac-6} &= \int_0^t (1 - \chi(s)) \int e^{i s \varphi} \mu(\overline{\xi}, \overline{\eta}) \frac{\varphi |\overline{\eta}| m_{\widehat{\mathcal{P}}}(\overline{\sigma})}{(|\overline{\eta}| + |\overline{\sigma}|)^2} \widehat{h}_{\beta}(s, \overline{\eta}) \widehat{h}_b(s, \overline{\sigma}) ~ d\overline{\eta} ds \label{developpementlemstructquadhac-6-1} \\
&\quad + \int_0^t \chi(s) \int e^{i s \varphi} \mu(\overline{\xi}, \overline{\eta}) \frac{\varphi |\overline{\eta}| m_{\widehat{\mathcal{P}}}(\overline{\sigma})}{(|\overline{\eta}| + |\overline{\sigma}|)^2} \widehat{h}_{\beta}(s, \overline{\eta}) \widehat{h}_b(s, \overline{\sigma}) ~ d\overline{\eta} ds \label{developpementlemstructquadhac-6-2} \\
&\quad + \int_0^t \int e^{i s \varphi} \mu(\overline{\xi}, \overline{\eta}) \frac{\varphi |\overline{\eta}| m_{\widehat{\mathcal{P}}}(\overline{\sigma})}{(|\overline{\eta}| + |\overline{\sigma}|)^2} \widehat{h}_{\beta}(s, \overline{\eta}) \widehat{g}_b(s, \overline{\sigma}) ~ d\overline{\eta} ds \label{developpementlemstructquadhac-6-3}
\end{align}
\end{subequations} 
for $\chi$ the usual Littlewood-Paley localisation function. We then apply an integration by parts in time on the first term: 
\begin{subequations}
\begin{align}
\eqref{developpementlemstructquadhac-6-1} &= \int_0^t \chi'(s) \int e^{i s \varphi} \mu(\overline{\xi}, \overline{\eta}) \frac{|\overline{\eta}| m_{\widehat{\mathcal{P}}}(\overline{\sigma})}{(|\overline{\eta}| + |\overline{\sigma}|)^2} \widehat{h}_{\beta}(s, \overline{\eta}) \widehat{h}_b(s, \overline{\sigma}) ~ d\overline{\eta} ds \label{developpementlemstructquadhac-6-1-1} \\
&\quad + \int_0^t (1 - \chi(s)) \int e^{i s \varphi} \mu(\overline{\xi}, \overline{\eta}) \frac{|\overline{\eta}| m_{\widehat{\mathcal{P}}}(\overline{\sigma})}{(|\overline{\eta}| + |\overline{\sigma}|)^2} \partial_s \widehat{h}_{\beta}(s, \overline{\eta}) \widehat{h}_b(s, \overline{\sigma}) ~ d\overline{\eta} ds \label{developpementlemstructquadhac-6-1-2} \\
&\quad + \int_0^t (1 - \chi(s)) \int e^{i s \varphi} \mu(\overline{\xi}, \overline{\eta}) \frac{|\overline{\eta}| m_{\widehat{\mathcal{P}}}(\overline{\sigma})}{(|\overline{\eta}| + |\overline{\sigma}|)^2} \widehat{h}_{\beta}(s, \overline{\eta}) \partial_s \widehat{h}_b(s, \overline{\sigma}) ~ d\overline{\eta} ds \label{developpementlemstructquadhac-6-1-3} \\
&\quad + (1 - \chi(t)) \int e^{i t \varphi} \mu(\overline{\xi}, \overline{\eta}) \frac{|\overline{\eta}| m_{\widehat{\mathcal{P}}}(\overline{\sigma})}{(|\overline{\eta}| + |\overline{\sigma}|)^2} \widehat{h}_{\beta}(t, \overline{\eta}) \widehat{h}_b(t, \overline{\sigma}) ~ d\overline{\eta} \label{developpementlemstructquadhac-6-1-4}
\end{align}
\end{subequations} 
where the symbol $\mu$ can change from line to line. 

We then define
\begin{align*}
g_{b \beta}(t) &:= \eqref{developpementlemstructquadhac-10} + \eqref{developpementlemstructquadhac-15} + \eqref{developpementlemstructquadhac-6-1-4} 
\end{align*}
and $h_{b \beta}(t)$ corresponds to all other terms. 

First, $g_{b \beta}(t)$ has indeed the desired form, given that \eqref{developpementlemstructquadhac-10} is of the form \eqref{lemstructhquadac-gresx}, \eqref{developpementlemstructquadhac-15} of the form \eqref{lemstructhquadac-gdersymb} and \eqref{developpementlemstructquadhac-6-1-4} of the form \eqref{lemstructhquadac-gvarphi}. 

In every other term, we may now systematically develop
\begin{align*}
|\nabla| m_{\gamma}(D) X_{\gamma} h_{\beta} = h_{\gamma \beta}(t) + g_{\gamma \beta}(t), \quad m_{\gamma}(D) X_{\gamma} f = h_{\gamma} + g_{\gamma}
\end{align*}

Let us first consider the quadratic terms obtained that way, that is the ones without any temporal derivative nor any $g$. We can already see that the quadratic term from \eqref{developpementlemstructquadhac-1} is exactly \eqref{lemstructhquadac-2-sym}, that the one from \eqref{developpementlemstructquadhac-2} is of the form \eqref{lemstructhquadac-2-doubleresx}, the one of \eqref{developpementlemstructquadhac-3} of the form \eqref{lemstructhquadac-2-doubleresxvarphi}, the one of \eqref{developpementlemstructquadhac-4} of the form \eqref{lemstructhquadac-2-resxresxb}, the one of \eqref{developpementlemstructquadhac-5} of the form \eqref{lemstructhquadac-2-resxresx}, the one of \eqref{developpementlemstructquadhac-11} of the form \eqref{lemstructhquadac-2-dersymb-resxb}, the one of \eqref{developpementlemstructquadhac-12} of the form \eqref{lemstructhquadac-2-dersymb-resxbis}. \eqref{developpementlemstructquadhac-7} is of the form \eqref{lemstructhquadac-2-dersymb-resx}, \eqref{developpementlemstructquadhac-13} of the form \eqref{lemstructhquadac-2-doubledersymb}, while \eqref{developpementlemstructquadhac-8}, \eqref{developpementlemstructquadhac-9}, \eqref{developpementlemstructquadhac-14} are purely cubic. Finally, in the decomposition of \eqref{developpementlemstructquadhac-6}, \eqref{developpementlemstructquadhac-6-3}, \eqref{developpementlemstructquadhac-6-1-2} and \eqref{developpementlemstructquadhac-6-1-3} are purely cubic, while \eqref{developpementlemstructquadhac-6-2} and \eqref{developpementlemstructquadhac-6-1-1} are of the form \eqref{lemstructhquadac-2-resxresxvarphi}. 

If we now consider the cubic terms, using the form of $g_{\alpha \beta}$ obtained previously, we already have that the cubic terms from \eqref{developpementlemstructquadhac-2} and \eqref{developpementlemstructquadhac-3} are of the form $\eqref{lemstructhquadac-3-resx}+\eqref{lemstructhquadac-3-resx} + \eqref{lemstructhquadac-3-varphi}$, and then using the form of $g_{\alpha}$ and $\partial_t f$, that the cubic terms from \eqref{developpementlemstructquadhac-4}, \eqref{developpementlemstructquadhac-5}, \eqref{developpementlemstructquadhac-6-3}, \eqref{developpementlemstructquadhac-9}, \eqref{developpementlemstructquadhac-11}, \eqref{developpementlemstructquadhac-12} and \eqref{developpementlemstructquadhac-14} are also of the desired form. For \eqref{developpementlemstructquadhac-6-1-2} and \eqref{developpementlemstructquadhac-6-1-3}, we can use that, for any $\alpha$, 
\begin{align*}
\partial_t h_{\alpha}(t) &= \sum_{\gamma = a, b, c} T_{\mu}\left[ m_{0, \overline{\gamma}}(D) |\nabla| h_{\gamma}, ~ f \right](t) 
+ \sum_{\gamma = a, b, c} T_{\mu}\left[ h_{\gamma}, ~ |\nabla| f \right](t) 
+ T_{\mu}[f, f](t) \\
&\quad + T_{\mu}\left[ f, ~ |\nabla| T_{\mu'}[f, f] \right](t) 
\end{align*}
for some $\mu$ that can change from term to term. Therefore, we deduce that the cubic terms obtained that way are also of the form \eqref{lemstructhquadac-3-varphi}, \eqref{lemstructhquadac-3-varphibis} or \eqref{lemstructhquadac-3-varphidersymb}, and that only a quadratic term of the form \eqref{lemstructhquadac-4-dersymb} remains. 

It then only remains to consider the cubic contribution of \eqref{developpementlemstructquadhac-1} and \eqref{developpementlemstructquadhac-8}, that is 
\begin{align*}
&2 \int_0^t \int e^{i s \varphi} i \eta_0 \widehat{g}_{b \beta}(s, \overline{\eta}) \widehat{f}(s, \overline{\sigma}) ~ d\overline{\eta} ds \\
&\quad + \int_0^t \int s e^{i s \varphi} m_g(\overline{\xi}, \overline{\eta}) |\overline{\eta}| \partial_s \widehat{h}_{\beta}(s, \overline{\eta}) \widehat{f}(s, \overline{\sigma}) ~ d\overline{\eta} ds 
\end{align*}
Replacing $g_{b \beta}(t)$ by its explicit expression $\eqref{developpementlemstructquadhac-10} + \eqref{developpementlemstructquadhac-15} + \eqref{developpementlemstructquadhac-6-1-4}$, the contribution of $\eqref{developpementlemstructquadhac-15} + \eqref{developpementlemstructquadhac-6-1-4}$ is of the form $\eqref{lemstructhquadac-3-dersymb}+\eqref{lemstructhquadac-3-varphi}+\eqref{lemstructhquadac-3-varphibis}$. Likewise, since 
\begin{align*}
\partial_t h_{\beta}(t) &= 2 \frac{\partial_x}{|\nabla|} T\left[ |\nabla| h_{\beta}, ~ f \right](t) + c_{\beta} \frac{\partial_x}{|\nabla|} T[f, f](t)
\end{align*}
the contribution of the $T[f, f]$ in \eqref{developpementlemstructquadhac-8} is of the form \eqref{lemstructhquadac-3-dersymb}. We therefore have only remaining
\begin{align*}
&2 \int_0^t \int e^{i s \varphi} i \eta_0 \eqref{developpementlemstructquadhac-10}(s, \overline{\eta}) \widehat{f}(s, \overline{\sigma}) ~ d\overline{\eta} ds + 2 \int_0^t \int s e^{i s \varphi} m_g(\overline{\xi}, \overline{\eta}) i \eta_0 \widehat{T}[|\nabla| h_{\beta}, ~ f](s, \overline{\eta}) \widehat{f}(s, \overline{\sigma}) ~ d\overline{\eta} ds \\
&= - 2 \int_0^t \int \int e^{i s \varphi_3(\overline{\xi}, \overline{\eta}, \overline{\sigma})} i s (\eta_0+\sigma_0) m_g(\overline{\eta}+\overline{\sigma}, \overline{\eta}) |\overline{\eta}| \widehat{h}_{\beta}(s, \overline{\eta}) \widehat{f}(s, \overline{\sigma}) \widehat{f}(s, \overline{\rho}) ~ d\overline{\eta} d\overline{\sigma} ds \\
&\quad + 2 \int_0^t \int \int e^{i s \varphi_3(\overline{\xi}, \overline{\eta}, \overline{\sigma})} i s (\eta_0+\sigma_0) m_g(\overline{\xi}, \overline{\eta}+\overline{\sigma}) |\overline{\eta}| \widehat{h}_{\beta}(s, \overline{\eta}) \widehat{f}(s, \overline{\sigma}) \widehat{f}(s, \overline{\rho}) ~ d\overline{\eta} d\overline{\sigma} ds \\
&= 2 \int_0^t \int \int e^{i s \varphi_3(\overline{\xi}, \overline{\eta}, \overline{\sigma})} i s (\eta_0+\sigma_0) \left( m_g(\overline{\xi}, \overline{\eta}+\overline{\sigma}) - m_g(\overline{\eta}+\overline{\sigma}, \overline{\eta}) \right) |\overline{\eta}| \widehat{h}_{\beta}(s, \overline{\eta}) \widehat{f}(s, \overline{\sigma}) \widehat{f}(s, \overline{\rho}) ~ d\overline{\eta} d\overline{\sigma} ds
\end{align*}
We then separate the integral into two pieces by localisation: either $|\overline{\sigma}|+|\overline{\rho}| \gtrsim |\overline{\eta}|$, and we automatically have the form \eqref{lemstructhquadac-3-resx}; or $|\overline{\sigma}|+|\overline{\rho}| \ll |\overline{\eta}|$ and we have by Taylor's formula
\begin{align*}
m_g(\overline{\xi}, \overline{\eta}+\overline{\sigma}) - m_g(\overline{\eta}+\overline{\sigma}, \overline{\eta})
&= O(\overline{\sigma}) + O(\overline{\rho})
\end{align*}
which gives exactly the structure of \eqref{lemstructhquadac-3-resx}. This concludes. 
\end{Dem}

\begin{Lem} Let $\alpha \in \{ a, c \}$. Then we can decompose 
\begin{align*}
|\nabla| m_{\alpha}(D) X_{\alpha} h_b(t) &= h_{\alpha b}(t) + g_{\alpha b}(t) \\
h{\alpha b}(t) &= |\nabla| m_{\alpha}(D) X_{\alpha} m_b(D) X_b f(0) + h_{\alpha b, 2}(t) + h_{\alpha b, 3}(t) + h_{\alpha b, 4}(t) 
\end{align*}
and we have the expressions
\begin{subequations}
\begin{align}
\widehat{h}_{\alpha b, 2}(t, \overline{\xi}) &= 2 \int_0^t \int e^{i s \varphi} i \eta_0 \widehat{h}_{\alpha b}(s, \overline{\eta}) \widehat{f}(s, \overline{\sigma}) ~ d\overline{\eta} ds \label{lemstructhquadacsurb-2-sym} \\
&\quad + \eqref{lemstructhquadac-2-doubleresx} 
+ \eqref{lemstructhquadac-2-doubleresxvarphi} 
+ \eqref{lemstructhquadac-2-resxresx}
+ \eqref{lemstructhquadac-2-resxresxb}
+ \eqref{lemstructhquadac-2-resxresxvarphi}
+ \eqref{lemstructhquadac-2-dersymb-resx}
+ \eqref{lemstructhquadac-2-dersymb-resxb} \notag \\
&\quad + \eqref{lemstructhquadac-2-dersymb-resxbis} 
+ \eqref{lemstructhquadac-2-doubledersymb} \notag \\
\widehat{h}_{\alpha b, 3}(t, \overline{\xi}) &= 
\eqref{lemstructhquadac-3-resx} + \eqref{lemstructhquadac-3-dersymb} + \eqref{lemstructhquadac-3-varphi} + \eqref{lemstructhquadac-3-varphibis} + \eqref{lemstructhquadac-3-varphidersymb} \notag \\
\widehat{h}_{\alpha b, 4}(t, \overline{\xi}) &= \int_0^t \int \int \int e^{i s \varphi_4(\overline{\xi}, \overline{\eta}^1, \overline{\eta}^2, \overline{\eta}^3)} \mu(\overline{\xi}, \overline{\eta}^1, \overline{\eta}^2, \overline{\eta}^3) s^2 |\overline{\eta}^1| \widehat{f}(s, \overline{\eta}^1) |\overline{\eta}^2|\widehat{f}(s, \overline{\eta}^2) \widehat{f}(s, \overline{\eta}^3) \notag \\
&\quad \quad \quad \quad \quad \quad \quad \quad \quad \quad \quad \quad \quad \quad \quad \quad \quad \quad \widehat{f}(s, \overline{\eta}^4) ~ d\overline{\eta}^1 d\overline{\eta}^2 d\overline{\eta}^3 ds \label{lemstructhquadacsurb-4-s2} \\
&\quad + \eqref{lemstructhquadac-4-dersymb} \notag \\
\widehat{g}_{\alpha b}(t, \overline{\xi}) &= t^2 \int \int e^{i t \varphi_3(\overline{\xi}, \overline{\eta}, \overline{\sigma})} \mu(\overline{\xi}, \overline{\eta}, \overline{\sigma}) |\overline{\eta}| \widehat{f}(t, \overline{\eta}) \widehat{f}(t, \overline{\sigma}) \widehat{f}(t, \overline{\rho}) ~ d\overline{\eta} d\overline{\sigma} \label{lemstructhquadacsurb-g-t2} \\
&\quad + 
\eqref{lemstructhquadac-gresx} + \eqref{lemstructhquadac-gdersymb} + \eqref{lemstructhquadac-gvarphi} \notag 
\end{align}
\end{subequations} 
for $\beta$ replaced by $\alpha$. 
\label{lemmestructurehquadacsurb} 
\end{Lem}

\begin{Dem}
We start from the formula of $h_b(t)$ given by Lemma \ref{lemstructurehalphagalpha}, decomposing $h_b(t) = h_{b, 2}(t) + h_{b, 3}(t)$. 

All terms in the decomposition of $h_{b, 2}(t)$ are time-integrated quadratic interactions, and we can therefore apply again the decomposition lemma \ref{lemdecompositionFGH}. The decompositions are very close to those already applied in the proof of Lemma \ref{lemmestructurehquadabcsurac} (in the opposite order) and we skip the details. The only different step is to check that the new term in $g_{\alpha b}$, that is \eqref{lemstructhquadacsurb-g-t2}, only creates well-controlated cubic terms when developing $m_{\gamma}(D) X_{\gamma} h_b = h_{\gamma b} + g_{\gamma b}$, but it is easy to check that all theses contributions are of the form \eqref{lemstructhquadacsurb-4-s2}. 

It remains to consider the effect of $|\nabla| m_{\alpha}(D) X_{\alpha}$ on the cubic term $h_{b, 3}(t)$. We can then develop in Fourier, first for $\alpha = a$: 
\begin{subequations}
\begin{align}
&|\overline{\xi}| \widehat{X}_a(\overline{\xi}) \cdot \nabla_{\overline{\xi}} \widehat{h}_{b, 3}(t, \overline{\xi}) \notag \\
&= \int_0^t \int \int \overline{\xi} \cdot \nabla_{\overline{\xi}} \left( e^{i s \varphi_3(\overline{\xi}, \overline{\eta}, \overline{\sigma})} \mu(\overline{\xi}, \overline{\eta}, \overline{\sigma}) s \eta_0 \frac{|\overline{\sigma}|+|\overline{\rho}|}{|\overline{\eta}|+|\overline{\sigma}|+|\overline{\rho}|} \widehat{f}(s, \overline{\eta}) \widehat{f}(s, \overline{\sigma}) \widehat{f}(s, \overline{\rho}) \right) ~ d\overline{\eta} d\overline{\sigma} ds \notag \\
&= \int_0^t \int \int e^{i s \varphi_3(\overline{\xi}, \overline{\eta}, \overline{\sigma})} i s \overline{\xi} \cdot \nabla_{\overline{\xi}} \varphi_3(\overline{\xi}, \overline{\eta}, \overline{\sigma}) \mu(\overline{\xi}, \overline{\eta}, \overline{\sigma}) s \eta_0 \frac{|\overline{\sigma}|+|\overline{\rho}|}{|\overline{\eta}|+|\overline{\sigma}|+|\overline{\rho}|} \widehat{f}(s, \overline{\eta}) \widehat{f}(s, \overline{\sigma}) \widehat{f}(s, \overline{\rho}) ~ d\overline{\eta} d\overline{\sigma} ds \label{preuvelemmestructhquadacsurb-cubique-1} \\
&\quad + \int_0^t \int \int e^{i s \varphi_3(\overline{\xi}, \overline{\eta}, \overline{\sigma})} \mu(\overline{\xi}, \overline{\eta}, \overline{\sigma}) s \eta_0 \frac{|\overline{\sigma}|+|\overline{\rho}|}{|\overline{\eta}|+|\overline{\sigma}|+|\overline{\rho}|} \widehat{f}(s, \overline{\eta}) \widehat{f}(s, \overline{\sigma}) \overline{\xi} \cdot \nabla_{\overline{\xi}} \widehat{f}(s, \overline{\rho}) ~ d\overline{\eta} d\overline{\sigma} ds \label{preuvelemmestructhquadacsurb-cubique-2} \\
&\quad + \int_0^t \int \int e^{i s \varphi_3(\overline{\xi}, \overline{\eta}, \overline{\sigma})} \overline{\xi} \cdot \nabla_{\overline{\xi}} \left( \mu(\overline{\xi}, \overline{\eta}, \overline{\sigma}) s \eta_0 \frac{|\overline{\sigma}|+|\overline{\rho}|}{|\overline{\eta}|+|\overline{\sigma}|+|\overline{\rho}|} \right) \widehat{f}(s, \overline{\eta}) \widehat{f}(s, \overline{\sigma}) \widehat{f}(s, \overline{\rho}) ~ d\overline{\eta} d\overline{\sigma} ds \label{preuvelemmestructhquadacsurb-cubique-3}
\end{align}
\end{subequations}
Note that, here above, the symbol $\mu$ has no precised form but is the same at each line, and will remain so in what follows. \eqref{preuvelemmestructhquadacsurb-cubique-3} is of the form \eqref{lemstructhquadac-3-dersymb}. In \eqref{preuvelemmestructhquadacsurb-cubique-2}, we can decompose: 
\begin{align*}
\overline{\xi} \cdot \nabla_{\overline{\xi}} \widehat{f}(s, \overline{\rho}) 
&= \overline{\rho} \cdot \nabla_{\overline{\xi}} \widehat{f}(s, \overline{\rho}) - \overline{\eta} \cdot \nabla_{\overline{\eta}} \widehat{f}(s, \overline{\rho}) - \overline{\sigma} \cdot \nabla_{\overline{\sigma}} \widehat{f}(s, \overline{\rho}) 
\end{align*}
and then apply integrations by parts on the two last terms. We then get:  
\begin{subequations}
\begin{align}
&\eqref{preuvelemmestructhquadacsurb-cubique-1}+\eqref{preuvelemmestructhquadacsurb-cubique-2} \notag \\
&\begin{aligned}
= \int_0^t \int \int e^{i s \varphi_3(\overline{\xi}, \overline{\eta}, \overline{\sigma})} i s \left( \overline{\xi} \cdot \nabla_{\overline{\xi}} + \overline{\eta} \cdot \nabla_{\overline{\eta}} + \overline{\sigma} \cdot \nabla_{\overline{\sigma}} \right) \varphi_3(\overline{\xi}, \overline{\eta}, \overline{\sigma}) \mu(\overline{\xi}, \overline{\eta}, \overline{\sigma}) s \eta_0 \\
\frac{|\overline{\sigma}|+|\overline{\rho}|}{|\overline{\eta}|+|\overline{\sigma}|+|\overline{\rho}|} \widehat{f}(s, \overline{\eta}) \widehat{f}(s, \overline{\sigma}) \widehat{f}(s, \overline{\rho}) ~ d\overline{\eta} d\overline{\sigma} ds 
\end{aligned} \label{preuvelemmestructhquadacsurb-cubique-1-1} \\
&+ \int_0^t \int \int e^{i s \varphi_3(\overline{\xi}, \overline{\eta}, \overline{\sigma})} \widetilde{\mu}(\overline{\xi}, \overline{\eta}, \overline{\sigma}) s |\overline{\eta}| \widehat{f}(s, \overline{\eta}) \widehat{f}(s, \overline{\sigma}) |\overline{\rho}| \widehat{X}_a(\overline{\rho}) \cdot \nabla_{\overline{\xi}} \widehat{f}(s, \overline{\rho}) ~ d\overline{\eta} d\overline{\sigma} ds \label{preuvelemmestructhquadacsurb-cubique-1-2} \\
&+ \int_0^t \int \int e^{i s \varphi_3(\overline{\xi}, \overline{\eta}, \overline{\sigma})} \widetilde{\mu}(\overline{\xi}, \overline{\eta}, \overline{\sigma}) s |\overline{\eta}| \widehat{f}(s, \overline{\eta}) \widehat{f}(s, \overline{\sigma}) \widehat{f}(s, \overline{\rho}) ~ d\overline{\eta} d\overline{\sigma} ds \label{preuvelemmestructhquadacsurb-cubique-1-3}
\end{align}
\end{subequations} 
For some symbol $\widetilde{\mu}$ that can change from line to line. \eqref{preuvelemmestructhquadacsurb-cubique-1-2} is of the form \eqref{lemstructhquadac-3-resx}, \eqref{preuvelemmestructhquadacsurb-cubique-1-3} of the form \eqref{lemstructhquadac-3-dersymb}. Finally, for \eqref{preuvelemmestructhquadacsurb-cubique-1-1}, we use the relation
\begin{align*}
\left( \overline{\xi} \cdot \nabla_{\overline{\xi}} + \overline{\eta} \cdot \nabla_{\overline{\eta}} + \overline{\sigma} \cdot \nabla_{\overline{\sigma}} \right) \varphi_3(\overline{\xi}, \overline{\eta}, \overline{\sigma}) 
&= 3 \varphi_3(\overline{\xi}, \overline{\eta}, \overline{\sigma}) 
\end{align*}
which allows to apply an integration by parts in time: 
\begin{subequations}
\begin{align}
\eqref{preuvelemmestructhquadacsurb-cubique-1-1} &= - 3 \int_0^t \int \int e^{i s \varphi_3(\overline{\xi}, \overline{\eta}, \overline{\sigma})} i s^2 \mu(\overline{\xi}, \overline{\eta}, \overline{\sigma}) \eta_0 \frac{|\overline{\sigma}|+|\overline{\rho}|}{|\overline{\eta}|+|\overline{\sigma}|+|\overline{\rho}|} \partial_s \left( \widehat{f}(s, \overline{\eta}) \widehat{f}(s, \overline{\sigma}) \widehat{f}(s, \overline{\rho}) \right) ~ d\overline{\eta} d\overline{\sigma} ds \label{preuvelemmestructhquadacsurb-cubique-1-1-1} \\
&\quad - 6 \int_0^t \int \int e^{i s \varphi_3(\overline{\xi}, \overline{\eta}, \overline{\sigma})} i s \mu(\overline{\xi}, \overline{\eta}, \overline{\sigma}) \eta_0 \frac{|\overline{\sigma}|+|\overline{\rho}|}{|\overline{\eta}|+|\overline{\sigma}|+|\overline{\rho}|} \widehat{f}(s, \overline{\eta}) \widehat{f}(s, \overline{\sigma}) \widehat{f}(s, \overline{\rho}) ~ d\overline{\eta} d\overline{\sigma} ds \label{preuvelemmestructhquadacsurb-cubique-1-1-2} \\
&\quad + 3 \int \int e^{i t \varphi_3(\overline{\xi}, \overline{\eta}, \overline{\sigma})} i t^2 \mu(\overline{\xi}, \overline{\eta}, \overline{\sigma}) \eta_0 \frac{|\overline{\sigma}|+|\overline{\rho}|}{|\overline{\eta}|+|\overline{\sigma}|+|\overline{\rho}|} \widehat{f}(t, \overline{\eta}) \widehat{f}(t, \overline{\sigma}) \widehat{f}(t, \overline{\rho}) ~ d\overline{\eta} d\overline{\sigma} ds \label{preuvelemmestructhquadacsurb-cubique-1-1-3}
\end{align}
\end{subequations}
\eqref{preuvelemmestructhquadacsurb-cubique-1-1-2} is of the form \eqref{lemstructhquadac-3-dersymb}, and by developing $\partial_s f$, \eqref{preuvelemmestructhquadacsurb-cubique-1-1-1} is of the form \eqref{lemstructhquadacsurb-4-s2}. 

Finally, we can place \eqref{preuvelemmestructhquadacsurb-cubique-1-1-3} in $g_{\alpha b}(t)$ since it is of the form \eqref{lemstructhquadacsurb-g-t2}. 

If $\alpha = c$, the decomposition is the same and is simpler, using again 
\begin{align*}
J \xi \cdot \nabla_{\xi} \widehat{f}(s, \overline{\rho}) &= J \rho \cdot \nabla_{\xi} \widehat{f}(s, \overline{\rho}) - J \eta \cdot \nabla_{\eta} \widehat{f}(s, \overline{\rho}) - J \sigma \cdot \nabla_{\sigma} \widehat{f}(s, \overline{\rho}) 
\end{align*}
to apply integrations by parts and the relation 
\begin{align*}
\left( J \xi \cdot \nabla_{\xi} + J \eta \cdot \nabla_{\eta} + J \sigma \cdot \nabla_{\sigma} \right) \varphi_3 = 0 
\end{align*}
to simplify the remainder. 
\end{Dem}

\section{Estimates on \texorpdfstring{$f, u, g$}{f, u, g}} \label{section-estimees-simples-fug}

\begin{Prop} \textbf{Linear estimates.} Let $\kappa > 0$. There exist $C > 0$ such that, for any $t > 0$, 
\begin{align*}
\Vert \langle \nabla \rangle u(t) \Vert_{L^4} &\leq C t^{-\frac{5}{12}} \langle t \rangle^{-\frac{1}{8}+50\delta} \Vert u \Vert_X \\
\Vert \langle \nabla \rangle^{\frac{1}{4}} u(t) \Vert_{L^4} &\leq C \langle t \rangle^{-\frac{13}{24}+50\delta} \Vert u \Vert_X \\
\Vert \langle \nabla \rangle^{\frac{3}{2}+\kappa} |\partial_x|^{-\frac{1}{2}+\kappa} |\nabla|^{-\frac{1}{2}-\kappa} u(t) \Vert_{L^2} &\leq C \Vert u \Vert_X 
\end{align*}
\textbf{A priori estimate and scattering.} There exist $C' > 0$ such that, for every $t \geq 0$, for every $\alpha \in \{ a, b, c \}$, 
\begin{align*}
\Vert \partial_t f(t) \Vert_{H^1} &\leq C' t^{-\frac{5}{6}} \langle t \rangle^{-\frac{1}{4}+100\delta} \Vert u \Vert_X^2 \\
\Vert \partial_t f(t) \Vert_{L^4} &\leq C' t^{-\frac{5}{6}} \langle t \rangle^{-\frac{19}{24}+150\delta} \Vert u \Vert_X^2 \\
\Vert \nabla \partial_t f(t) \Vert_{L^4} &\leq C' t^{-\frac{41}{36}} \langle t \rangle^{-\frac{35}{72}+150\delta} \Vert u \Vert_X^2 
\end{align*}
In particular, $\partial_t f \in L^1_t H^1_{x, y}$ and
\begin{align*}
\Vert u(t) \Vert_{H^1} &\leq \Vert u_0 \Vert_{H^1} + C' \Vert u \Vert_X^2 
\end{align*}
 \label{propositionestimeesstandardscattering} 
\end{Prop}

\begin{Dem}
We start with the linear estimates. 

For the $L^2$ estimate, we apply Parseval's identity and we can separate: 
\begin{align*}
&\Vert \langle \nabla \rangle^{\frac{3}{2}+\kappa} |\partial_x|^{-\frac{1}{2}+\kappa} |\nabla|^{-\frac{1}{2}-\kappa} u(t) \Vert_{L^2} \\ 
&\lesssim \Vert \chi(D) m_{\widehat{\mathcal{P}}}(D) |\partial_x|^{-\frac{1}{2}+\kappa} |\nabla_y|^{-\frac{1}{2}-\kappa} f(t) \Vert_{L^2} + \Vert (1 - \chi(D)) m_{\widehat{\mathcal{P}}}(D) \nabla |\partial_x|^{-\frac{1}{2}+\kappa} f(t) \Vert_{L^2} \\
&\quad + \Vert \left( 1 - m_{\widehat{\mathcal{P}}}(D) \right) \langle \nabla \rangle^{\frac{3}{2}+\kappa} |\partial_x|^{-\frac{1}{2}+\kappa} |\nabla|^{-\frac{1}{2}-\kappa} f(t) \Vert_{L^2} 
\end{align*}
The part localised away from $\widehat{\mathcal{P}}$ can be estimated as $\Vert |\nabla|^{-1} f(t) \Vert_{L^2}$ by the Hardy inequality. For the part localised at low frequencies, we can first apply the fractional Hardy inequality in dimension $1$ and in dimension $2$, noting that $|x|^r$ and $|\nabla_y|^p$ commute: 
\begin{align*}
\Vert \chi(D) m_{\widehat{\mathcal{P}}}(D) |\partial_x|^{-\frac{1}{2}+\kappa} |\nabla_y|^{-\frac{1}{2}-\kappa} f(t) \Vert_{L^2} 
&\lesssim \Vert |x|^{\frac{1}{2}-\kappa} \chi(D) m_{\widehat{\mathcal{P}}}(D) |\nabla|^{-\frac{1}{2}-\kappa} f(t) \Vert_{L^2} \\
&\lesssim \Vert |y|^{\frac{1}{2}+\kappa} |x|^{\frac{1}{2}-\kappa} \chi(D) m_{\widehat{\mathcal{P}}}(D) f(t) \Vert_{L^2} \\
&\lesssim \Vert (x, y) \chi(D) m_{\widehat{\mathcal{P}}}(D) f(t) \Vert_{L^2} \\
&\lesssim \Vert |\nabla|^{-1} f(t) \Vert_{L^2} + \Vert m_{\widehat{\mathcal{P}}}(D) (x, y) f(t) \Vert_{L^2} \\
&\lesssim \Vert u \Vert_X
\end{align*}
Finally, for the high frequency part, we proceed the same way: 
\begin{align*}
\Vert (1 - \chi(D)) m_{\widehat{\mathcal{P}}}(D) \nabla |\partial_x|^{-\frac{1}{2}+\kappa} f(t) \Vert_{L^2} 
&\lesssim \Vert |x|^{\frac{1}{2}-\kappa} (1 - \chi(D)) m_{\widehat{\mathcal{P}}}(D) \nabla f(t) \Vert_{L^2} \\
&\lesssim \Vert f(t) \Vert_{H^1} + \Vert m_{\widehat{\mathcal{P}}}(D) x f(t) \Vert_{L^2} \\
&\lesssim \Vert u \Vert_X
\end{align*}

Then, for the $L^4$ estimate, we interpolate and use the above linear estimates we just showed, with $\kappa = \frac{\delta}{2}$: 
\begin{align*}
\Vert \langle \nabla \rangle u(t) \Vert_{L^4} &\lesssim \Vert m_{\widehat{\mathcal{P}}}(D) \langle \nabla \rangle u(t) \Vert_{L^4} + \Vert (1 - m_{\widehat{\mathcal{P}}}(D)) \langle \nabla \rangle u(t) \Vert_{L^4} \\
&\lesssim \sum_{\substack{j, k \in \mathbb{Z}, \\ k \leq -10}} \langle 2^j \rangle \Vert \psi_{j, k}^{\widehat{\mathcal{P}}}(D) u(t) \Vert_{L^4} \quad + \sum_{j \in \mathbb{Z}} \langle 2^j \rangle \Vert \psi_j(D) (1 - m_{\widehat{\mathcal{P}}}(D)) u(t) \Vert_{L^4} \\
&\lesssim \sum_{\substack{j, k \in \mathbb{Z}, \\ k \leq -10}} \langle 2^j \rangle \Vert \psi_{j, k}^{\widehat{\mathcal{P}}}(D) u(t) \Vert_{L^2}^{\frac{1}{2}} \Vert \psi_{j, k}^{\widehat{\mathcal{P}}}(D) u(t) \Vert_{L^{\infty}}^{\frac{1}{2}} \\
&\quad \quad + \sum_{j \in \mathbb{Z}} \langle 2^j \rangle \Vert \psi_j(D) (1 - m_{\widehat{\mathcal{P}}}(D)) u(t) \Vert_{L^2}^{\frac{1}{2}} \Vert \psi_j(D) (1 - m_{\widehat{\mathcal{P}}}(D)) u(t) \Vert_{L^{\infty}}^{\frac{1}{2}} \\
&\lesssim \sum_{\substack{j, k \in \mathbb{Z}, \\ k \leq -10}} t^{-\frac{5}{12}} \langle t \rangle^{-\frac{1}{8}+50\delta} 2^{\frac{\delta j+\delta k}{2}} 2^{-\frac{j}{2}-\frac{k}{4}} \langle 2^j \rangle^{\frac{3}{4}-\frac{1}{100}} \Vert \psi_{j, k}^{\widehat{\mathcal{P}}}(D) u(t) \Vert_{L^2}^{\frac{1}{2}} \Vert u \Vert_X^{\frac{1}{2}} \\
&\quad \quad + \sum_{j \in \mathbb{Z}} t^{-\frac{5}{12}} \langle t \rangle^{-\frac{1}{8}+50\delta} 2^{-\frac{j}{2}+\frac{\delta j}{2}} \langle 2^j \rangle^{1-\frac{1}{100}} \Vert \psi_j(D) (1 - m_{\widehat{\mathcal{P}}}(D)) u(t) \Vert_{L^2}^{\frac{1}{2}} \Vert u \Vert_X^{\frac{1}{2}} \\
&\lesssim \sum_{\substack{j, k \in \mathbb{Z}, \\ k \leq -10}} t^{-\frac{5}{12}} \langle t \rangle^{-\frac{1}{8}+50\delta} 2^{\frac{2\delta j+\delta k}{4}} \langle 2^j \rangle^{-\frac{1}{100}} \Vert \langle \nabla \rangle^{\frac{3}{2}} |\nabla|^{-\frac{1}{2}-\frac{\delta}{2}} |\partial_x|^{-\frac{1}{2}+\frac{\delta}{2}} u(t) \Vert_{L^2}^{\frac{1}{2}} \Vert u \Vert_X^{\frac{1}{2}} \\
&\quad \quad + \sum_{j \in \mathbb{Z}} t^{-\frac{5}{12}} \langle t \rangle^{-\frac{1}{8}+50\delta} 2^{\frac{\delta j}{2}} \langle 2^j \rangle^{-\frac{1}{100}} \Vert \langle \nabla \rangle^2 |\nabla|^{-1} u(t) \Vert_{L^2}^{\frac{1}{2}} \Vert u \Vert_X^{\frac{1}{2}} \\
&\lesssim t^{-\frac{5}{12}} \langle t \rangle^{-\frac{1}{8}+50\delta} \Vert u \Vert_X
\end{align*}
as wanted. The estimate of $\Vert \langle \nabla \rangle^{\frac{1}{4}} u(t) \Vert_{L^4}$ is identical for $t \geq 1$, and for $t \leq 1$ it follows from Sobolev's embedding: 
\begin{align*}
\Vert \langle \nabla \rangle^{\frac{1}{4}} u(t) \Vert_{L^4} &\lesssim \Vert u(t) \Vert_{H^1} \lesssim \Vert u \Vert_X
\end{align*}

We now prove the a priori estimates. 

For the $L^2$ estimate, we can compute directly 
\begin{align*}
\Vert \partial_t f(t) \Vert_{L^2} &= \Vert u(t) \partial_x u(t) \Vert_{L^2} \\
&\leq \Vert u(t) \Vert_{L^2} \Vert \partial_x u(t) \Vert_{L^{\infty}} \\
&\lesssim t^{-\frac{5}{6}} \langle t \rangle^{-\frac{1}{4}+100\delta} \Vert u \Vert_X^2 
\end{align*}
In fact it is not needed to use the presence of the $\partial_x$ factor here: 
\begin{align*}
\Vert \frac{\langle \nabla \rangle}{|\partial_x|} \partial_t f(t) \Vert_{L^2} &\lesssim \Vert \langle \nabla \rangle u(t) \Vert_{L^4} \Vert u(t) \Vert_{L^4} \\
&\lesssim t^{-\frac{5}{6}} \langle t \rangle^{-\frac{1}{4}+100\delta} \Vert u \Vert_X^2
\end{align*}
as well. 

Then, for the $\dot{H}^1$ estimate, we separate: 
\begin{align*}
\Vert \nabla \partial_t f(t) \Vert_{L^2} &\lesssim \Vert \chi\left( 2^{-j_0} \nabla \right) \nabla \partial_t f(t) \Vert_{L^2} + \Vert \left( 1 - \chi\left( 2^{-j_0} \nabla \right) \right) \nabla \partial_t f(t) \Vert_{L^2} 
\end{align*}
for some $j_0 \geq 0$ such that $2^{j_0} \simeq t^{-\frac{1}{3}}$ if $t \leq 1$, $j_0 = 0$ if $t \geq 1$. 

On the one hand, for the low frequency part, if $t \geq 1$ we localised on $|\overline{\xi}| \leq 1$ so we can simply bound by the $L^2$ estimate. If $t \leq 1$, we separate the $BHH$ interaction from the rest: 
\begin{align*}
\Vert \chi\left( 2^{-j_0} \nabla \right) \nabla \partial_t f(t) \Vert_{L^2} &\lesssim 2^{2j_0} \Vert \chi\left( 2^{-j_0-10} \nabla \right) u(t) \Vert_{L^4}^2 + \Vert \left( 1 - \chi\left( 2^{-j_0-10} \nabla \right) \right) \nabla u(t) \Vert_{L^4}^2 \\
&\lesssim 2^{\frac{5j_0}{2}} \Vert \nabla u(t) \Vert_{L^2}^2 + t^{-\frac{5}{6}} \Vert u \Vert_X^2 \\
&\lesssim t^{-\frac{5}{6}} \Vert u \Vert_X^2 
\end{align*}

For the high frequency part, we decompose the interaction as: 
\begin{subequations}
\begin{align}
\left( 1 - \chi\left( 2^{-j_0} \overline{\xi} \right) \right) \overline{\xi} \partial_t \widehat{f}(t, \overline{\xi}) &= \left( 1 - \chi\left( 2^{-j_0} \overline{\xi} \right) \right) \overline{\xi} \int e^{i t \varphi} \mu_{BBB}(\overline{\xi}, \overline{\eta}) i \xi_0 \widehat{f}(t, \overline{\eta}) \widehat{f}(t, \overline{\sigma}) ~ d\overline{\eta} \label{estapriori-dtfBBBH1} \\
&\quad + 2 \left( 1 - \chi\left( 2^{-j_0} \overline{\xi} \right) \right) \overline{\xi} \int e^{i t \varphi} \mu_{HBH}(\overline{\xi}, \overline{\eta}) i \xi_0 \widehat{f}(t, \overline{\eta}) \widehat{f}(t, \overline{\sigma}) ~ d\overline{\eta} \label{estapriori-dtfBHBH1} 
\end{align}
\end{subequations} 
where we symmetrized $\overline{\eta}$ and $\overline{\sigma}$. We can easily estimate
\begin{align*}
\Vert \eqref{estapriori-dtfBBBH1} \Vert_{L^2} &\lesssim \Vert \nabla u(t) \Vert_{L^2} \Vert \partial_x u(t) \Vert_{L^{\infty}} \\
&\lesssim t^{-\frac{5}{6}} \langle t \rangle^{-\frac{1}{4}+100\delta} \Vert u \Vert_X^2
\end{align*}
Then, if $|\overline{\eta}| \ll |\overline{\xi}| \simeq |\overline{\sigma}|$, we have
\begin{align*}
\partial_{\eta_0} \varphi &= 3 \sigma_0^2 + |\sigma|^2 - 3 \eta_0^2 - |\eta|^2 \simeq |\overline{\sigma}|^2 
\end{align*}
and therefore we can apply an integration by parts: 
\begin{subequations}
\begin{align}
\eqref{estapriori-dtfBHBH1} &= \left( 1 - \chi\left( 2^{-j_0} \overline{\xi} \right) \right) t^{-1} \int e^{i t \varphi} \mu(\overline{\xi}, \overline{\eta}) \mu_{HBH}(\overline{\xi}, \overline{\eta}) \partial_{\eta_0} \widehat{f}(t, \overline{\eta}) \widehat{f}(t, \overline{\sigma}) ~ d\overline{\eta} \label{estapriori-dtfBHBH1-1} \\
&\quad + \left( 1 - \chi\left( 2^{-j_0} \overline{\xi} \right) \right) t^{-1} \int e^{i t \varphi} \mu(\overline{\xi}, \overline{\eta}) \mu_{HBH}(\overline{\xi}, \overline{\eta}) \widehat{f}(t, \overline{\eta}) \partial_{\eta_0} \widehat{f}(t, \overline{\sigma}) ~ d\overline{\eta} \label{estapriori-dtfBHBH1-2} \\
&\quad + \left( 1 - \chi\left( 2^{-j_0} \overline{\xi} \right) \right) t^{-1} \int e^{i t \varphi} \mu(\overline{\xi}, \overline{\eta}) \mu_{HBH}(\overline{\xi}, \overline{\eta}) |\overline{\xi}|^{-1} \widehat{f}(t, \overline{\eta}) \widehat{f}(t, \overline{\sigma}) ~ d\overline{\eta} \label{estapriori-dtfBHBH1-3}
\end{align}
\end{subequations}
for some symbol $\mu$ of order $0$, that can change from line to line. Then, if $t \geq 1$: 
\begin{align*}
\Vert \eqref{estapriori-dtfBHBH1-1} \Vert_{L^2} &\lesssim t^{-1} 2^{-j_0} \Vert (x, y) f(t) \Vert_{L^2} \Vert \nabla u(t) \Vert_{L^{\infty}} \\
&\lesssim t^{-\frac{11}{6}} \langle t \rangle^{\frac{1}{3}+201\delta} 2^{-j_0} \Vert u \Vert_X^2 \\
&\lesssim t^{-\frac{3}{2}+201\delta} \Vert u \Vert_X^2 \\
\Vert \eqref{estapriori-dtfBHBH1-2} \Vert_{L^2} &\lesssim t^{-1} \Vert u(t) \Vert_{L^{\infty}} \Vert \nabla (x, y) f(t) \Vert_{L^2} \\
&\lesssim t^{-\frac{11}{6}} \langle t \rangle^{\frac{1}{4}+201\delta} \Vert u \Vert_X^2 \\
&\lesssim t^{-\frac{3}{2}+201\delta} \Vert u \Vert_X^2 \\
\Vert \eqref{estapriori-dtfBHBH1-3} \Vert_{L^2} &\lesssim t^{-1} \Vert u(t) \Vert_{L^4}^2 \\
&\lesssim t^{-\frac{13}{6}+100\delta} \Vert u \Vert_X^2 
\end{align*}
If $t \leq 1$: 
\begin{align*}
\Vert \eqref{estapriori-dtfBHBH1-1} \Vert_{L^2} &\lesssim t^{-1} 2^{-\frac{j_0}{2}} \Vert e^{-i t \omega(D)} (x, y) f(t) \Vert_{L^6} \Vert |\nabla|^{\frac{1}{2}} u(t) \Vert_{L^3} \\
&\lesssim t^{-1} 2^{-\frac{j_0}{2}} \Vert \nabla (x, y) f(t) \Vert_{L^2} \Vert \nabla u(t) \Vert_{L^2} \\
&\lesssim t^{-\frac{5}{6}} \Vert u \Vert_X^2 \\
\Vert \eqref{estapriori-dtfBHBH1-2} \Vert_{L^2} &\lesssim t^{-1} 2^{-\frac{j_0}{2}} \Vert u(t) \Vert_{L^6} \Vert e^{-i t \omega(D)} |\nabla|^{\frac{1}{2}} (x, y) f(t) \Vert_{L^3} \\
&\lesssim t^{-\frac{5}{6}} \Vert u \Vert_X^2 \\
\Vert \eqref{estapriori-dtfBHBH1-3} \Vert_{L^2} &\lesssim t^{-1} 2^{-\frac{3j_0}{2}} \Vert u(t) \Vert_{L^6} \Vert |\nabla|^{\frac{1}{2}} u(t) \Vert_{L^3} \\
&\lesssim t^{-\frac{5}{6}} \Vert u \Vert_X^2 
\end{align*}

We therefore get globally: 
\begin{align*}
\Vert \partial_t f(t) \Vert_{\dot{H}^1} &\lesssim t^{-\frac{5}{6}} \langle t \rangle^{-\frac{1}{4}+100\delta} \Vert u \Vert_X^2
\end{align*}
as wanted. 

For the $L^4$ estimate, we have for $t \geq 1$:
\begin{align*}
\Vert e^{-i t \omega(D)} \partial_t f(t) \Vert_{L^4} &\lesssim \Vert \partial_x u(t) \Vert_{L^{\infty}} \Vert u(t) \Vert_{L^4} \\
&\lesssim t^{-\frac{13}{8}+150\delta} \Vert u \Vert_X^2 
\end{align*}
For $t \leq 1$, we can absorb the singularity by Sobolev's embedding and interpolation: 
\begin{align*}
\Vert e^{-it\omega(D)} \partial_t f(t) \Vert_{L^4} &\lesssim \Vert |\nabla|^{\frac{3}{4}} \partial_t f(t) \Vert_{L^2} \\
&\lesssim t^{-\frac{5}{6}} \Vert u \Vert_X^2 
\end{align*}
We get then for any $t > 0$
\begin{align*}
\Vert e^{-i t \omega(D)} \partial_t f(t) \Vert_{L^4} &\lesssim t^{-\frac{5}{6}} \langle t \rangle^{-\frac{19}{24}+150\delta} \Vert u \Vert_X^2 
\end{align*}

If we add one derivative, we can reuse the high and low frequencies decomposition used above for the $\dot{H}^1$ estimate, dealing differently with the cases $t \geq 1$ and $t \leq 1$. We choose $j_0 = 0$ for any $t$ this time, so that the low frequencies part satisfies the same estimates as without the additional derivative. 

If $t \geq 1$, we can apply very similar transformations and estimates to the high frequencies part as for the $\dot{H}^1$ estimate. We skip the details. 

If $t \leq 1$, using the same decomposition and the same integral by parts as in the case $t \geq 1$ (needed to absorb the additional derivative), 
\begin{align*}
\Vert e^{-it \omega(D)} \mathcal{F}^{-1} \eqref{estapriori-dtfBHBH1-1} \Vert_{L^4} &\lesssim t^{-1} \Vert e^{-i t \omega(D)} (x, y) f(t) \Vert_{L^6} \Vert u(t) \Vert_{L^{12}} \\
&\lesssim t^{-1} \Vert \nabla (x, y) f(t) \Vert_{L^2} \Vert \nabla u(t) \Vert_{L^2}^{\frac{5}{6}} \Vert \nabla u(t) \Vert_{L^{\infty}}^{\frac{1}{6}} \\
&\lesssim t^{-\frac{41}{36}} \Vert u \Vert_X^2 \\
\Vert e^{-i t \omega(D)} \mathcal{F}^{-1} \eqref{estapriori-dtfBHBH1-2} \Vert_{L^4} &\lesssim t^{-1} \Vert u(t) \Vert_{L^{12}} \Vert e^{-i t \omega(D)} (x, y) f(t) \Vert_{L^6} \\
&\lesssim t^{-\frac{41}{36}} \Vert u \Vert_X^2 \\
\Vert e^{-i t \omega(D)} \mathcal{F}^{-1} \eqref{estapriori-dtfBHBH1-3} \Vert_{L^4} &\lesssim t^{-1} \Vert u(t) \Vert_{L^4}^2 \\
&\lesssim t^{-1} \Vert u \Vert_X^2
\end{align*}
We therefore get for any $t > 0$ 
\begin{align*}
\Vert e^{-i t \omega(D)} \nabla \partial_t f(t) \Vert_{L^4} &\lesssim t^{-\frac{41}{36}} \langle t \rangle^{-\frac{35}{72}+150\delta} \Vert u \Vert_X^2
\end{align*}
as wanted. 
\end{Dem}

We note that the scattering results propagate to higher regularities, as shown in the following lemma: 

\begin{Lem} Let $N \geq 2$. Assume that $u_0 \in H^N$ and that $T$ is chosen small enough so that $u \in C^0([0, T], H^N)$. Then 
\begin{align*}
\Vert u(t) \Vert_{H^N}^2 &\leq \Vert u_0 \Vert_{H^N}^2 + \Vert u \Vert_X \Vert u \Vert_{L^{\infty}_t H^N}^2 \left( 1 + \Vert u \Vert_X \right)^2 
\end{align*}
\end{Lem}

In particular, if the initial data is more regular, the a priori estimate can be extended changing $H^1$ into $H^N$. The scattering then holds in $H^s$ for any $s < N$: since $f(t)$ is bounded in $H^N$, up to extracting a subsequence it converges weakly in $H^N$ to a limit, which has to be $f_{\infty}$, and therefore $f_{\infty} \in H^N$; then by interpolation one has 
\begin{align*}
\Vert f(t) - f_{\infty} \Vert_{H^s} &\lesssim \Vert f(t) - f_{\infty} \Vert_{L^2}^{\frac{N-s}{N}} \Vert f(t) - f_{\infty} \Vert_{H^N}^{\frac{s}{N}}
\end{align*}
The first factor converges to $0$ as $t \to \infty$, while the second is uniformly bounded, hence the convergence in $H^s$ for any $s < N$. 

\begin{Dem}
It is well-known that the $L^2$ is preserved by the equation. 

Then, 
\begin{align*}
\Vert \partial_x^N u(t) \Vert_{L^2}^2 &= \Vert \partial_x^N u_0 \Vert_{L^2}^2 - 2 \int_0^t \int \partial_x^N u(s, x, y) \partial_x^{N+1} \Delta u(s, x, y) ~ dx dy ds \\
&\quad - 2 \int_0^t \int \partial_x^N u(s, x, y) \partial_x^N \left( u(s, x, y) \partial_x u(s, x, y) \right) ~ dx dy ds \\
&= \Vert \partial_x^N u_0 \Vert_{L^2}^2 + 2 \int_0^t \int \nabla \partial_x^N u(s, x, y) \cdot \partial_x \nabla \partial_x^N u(s, x, y) ~ dx dy ds \\
&\quad - 2 \int_0^t \int \partial_x^N u(s, x, y) \partial_x^{N+1} u(s, x, y) ~ u(s, x, y) ~ dx dy ds 
\\
&\quad - 2 \int_0^t \int \partial_x^N u(s, x, y) \left[ \partial_x^N, u(s, x, y) \right] \partial_x u(s, x, y) ~ dx dy ds \\
&= \Vert \partial_x^N u_0 \Vert_{L^2}^2 + \int_0^t \int \left( \partial_x^N u(s, x, y) \right)^2 \partial_x u(s, x, y) ~ dx dy ds 
\\
&\quad - 2 \int_0^t \int \partial_x^N u(s, x, y) \left[ \partial_x^N, u(s, x, y) \right] \partial_x u(s, x, y) ~ dx dy ds 
\end{align*}
by integration by parts. But 
\begin{align*}
\left| \int_0^t \int \left( \partial_x^N u(s, x, y) \right)^2 \partial_x u(s, x, y) ~ dx dy ds \right| &\leq \Vert u \Vert_{L^{\infty} H^N}^2 \int_0^t \Vert \partial_x u(s) \Vert_{L^{\infty}} ~ ds \\
&\lesssim \Vert u \Vert_{L^{\infty} H^N}^2 \Vert u \Vert_X \\
\left| \int_0^t \int \partial_x^N u(s, x, y) \left[ \partial_x^N, u(s, x, y) \right] \partial_x u(s, x, y) ~ dx dy ds \right| &\lesssim \int_0^t \Vert \partial_x^N u(s) \Vert_{L^2} \Vert \partial_x^N u(s) \Vert_{L^2} \Vert \partial_x u(s) \Vert_{L^{\infty}} ~ ds \\
&\lesssim \Vert u \Vert_{L^{\infty} H^N}^2 \Vert u \Vert_X
\end{align*}
by Kato-Ponce's estimate. 

We then only need to estimate $\Vert \nabla_y^N u(t) \Vert_{L^2}$. If we localise on $|\xi_0| \gtrsim |\xi|$, the above estimate is enough, so we can simply consider 
\begin{align*}
&\Vert \nabla_y^N m_{\widehat{\mathcal{P}}}(D) u(t) \Vert_{L^2}^2 = \Vert \nabla_y^N m_{\widehat{\mathcal{P}}}(D) f(t) \Vert_{L^2}^2 \\
&= \Vert \nabla_y^N m_{\widehat{\mathcal{P}}}(D) u_0 \Vert_{L^2}^2 + 2 \int_0^t \int \int |\xi|^{2N} m_{\widehat{\mathcal{P}}}(\overline{\xi}) \widehat{f}(s, -\overline{\xi}) e^{i s \varphi} m_{\widehat{\mathcal{P}}}(\overline{\xi}) i \xi_0 \widehat{f}(s, \overline{\eta}) \widehat{f}(s, \overline{\sigma}) ~ d\overline{\eta} d\overline{\xi} ds 
\end{align*}
by Parseval's identity. We now localise on different areas. 

If we localise on $\{ |\overline{\eta}| \simeq |\overline{\sigma}| \}$ by some symbol $\mu = \mu_{BHH} + \mu_{BBB}$, then we can directly estimate
\begin{align*}
&\left| \int_0^t \int \int |\xi|^{2N} m_{\widehat{\mathcal{P}}}(\overline{\xi}) \widehat{f}(s, -\overline{\xi}) e^{i s \varphi} \mu(\overline{\xi}, \overline{\eta}) i m_{\widehat{\mathcal{P}}}(\overline{\xi}) \xi_0 \widehat{f}(s, \overline{\eta}) \widehat{f}(s, \overline{\sigma}) ~ d\overline{\eta} d\overline{\xi} ds \right| \\
&\quad \lesssim \int_0^t \Vert \nabla^N f(s) \Vert_{L^2}^2 \Vert \partial_x u(s) \Vert_{L^{\infty}} ~ ds \\
&\quad \lesssim \Vert u \Vert_{L^{\infty} H^N}^2 \Vert u \Vert_X
\end{align*}

By symmetry, it suffices to consider the case $\{ |\overline{\eta}| \ll |\overline{\sigma}| \simeq |\overline{\xi}| \}$. 

We now symmetrise 
\begin{align*}
&2 \int_0^t \int \int |\xi|^{2N} m_{\widehat{\mathcal{P}}}(\overline{\xi}) \widehat{f}(s, -\overline{\xi}) e^{i s \varphi} \mu_{HBH}(\overline{\xi}, \overline{\sigma}) m_{\widehat{\mathcal{P}}}(\overline{\xi}) i \xi_0 \widehat{f}(s, \overline{\eta}) \widehat{f}(s, \overline{\sigma}) ~ d\overline{\eta} d\overline{\xi} ds \\
&\quad = - 2 \int_0^t \int \int |\sigma|^{2N} m_{\widehat{\mathcal{P}}}(\overline{\sigma}) \widehat{f}(s, \overline{\sigma}) e^{i s \varphi} \mu_{HBH}(\overline{\xi}, \overline{\sigma}) m_{\widehat{\mathcal{P}}}(\overline{\sigma}) i \sigma_0 \widehat{f}(s, \overline{\eta}) \widehat{f}(s, -\overline{\xi}) ~ d\overline{\eta} d\overline{\xi} ds \\
&\quad = \int_0^t \int \int \left( \xi_0 |\xi|^{2N} m_{\widehat{\mathcal{P}}}(\overline{\xi})^2 - \sigma_0 |\sigma|^{2N} m_{\widehat{\mathcal{P}}}(\overline{\sigma})^2 \right) \widehat{f}(s, -\overline{\xi}) e^{i s \varphi} \mu_{HBH}(\overline{\xi}, \overline{\sigma}) i \widehat{f}(s, \overline{\eta}) \widehat{f}(s, \overline{\sigma}) ~ d\overline{\eta} d\overline{\xi} ds
\end{align*}

Note that 
\begin{align*}
\varphi &= O(\eta_0) + \xi_0 (|\xi|^2 - |\sigma|^2) 
\end{align*}
hence
\begin{align*}
\xi_0 |\xi|^{2N} m_{\widehat{\mathcal{P}}}(\overline{\xi})^2 - \sigma_0 |\sigma|^{2N} m_{\widehat{\mathcal{P}}}(\overline{\sigma})^2 &= O(\eta_0) + O\left( \xi_0 \left( m_{\widehat{\mathcal{P}}}(\overline{\xi}) - m_{\widehat{\mathcal{P}}}(\overline{\sigma}) \right) \right) + O\left( \xi_0 (|\xi|^{2N} - |\sigma|^{2N}) \right) \\
&= O(\eta_0) + O\left( \xi_0 (|\xi| - |\sigma|) \right) \\
&= O(\eta_0) + O(\varphi)
\end{align*}

On the one hand, 
\begin{align*}
&\left| \int_0^t \int \int \eta_0 |\overline{\xi}|^{2N} \widehat{f}(s, -\overline{\xi}) e^{i s \varphi} \mu_{HBH}(\overline{\xi}, \overline{\sigma}) i \widehat{f}(s, \overline{\eta}) \widehat{f}(s, \overline{\sigma}) ~ d\overline{\eta} d\overline{\xi} ds \right| \\
&\quad \lesssim \int_0^t \Vert \nabla^N f(s) \Vert_{L^2}^2 \Vert \partial_x u(s) \Vert_{L^{\infty}} ~ ds \\
&\quad \lesssim \Vert u \Vert_{L^{\infty} H^N}^2 \Vert u \Vert_X
\end{align*}
so we know how to control the contribution of $O(\eta_0)$. 

On the other hand, if we have a factor $\varphi$, we can use $\varphi = O(\overline{\eta})$ to estimate the short time contribution: 
\begin{align*}
&\left| \int_0^{\min(1, t)} \int \int \varphi |\overline{\xi}|^{2N-2} \widehat{f}(s, -\overline{\xi}) e^{i s \varphi} \mu_{HBH}(\overline{\xi}, \overline{\sigma}) i \widehat{f}(s, \overline{\eta}) \widehat{f}(s, \overline{\sigma}) ~ d\overline{\eta} d\overline{\xi} ds \right| \\
&\quad \lesssim \int_0^{\min(1, t)} \Vert \nabla^N f(s) \Vert_{L^2}^2 \Vert \nabla u(s) \Vert_{L^{\infty}} ~ ds \\
&\quad \lesssim \Vert u \Vert_{L^{\infty} H^N}^2 \Vert u \Vert_X
\end{align*}
For $t \geq 1$ we apply an integration by parts in time: 
\begin{subequations}
\begin{align}
&\int_1^t \int \int \varphi |\overline{\xi}|^{2N-2} \widehat{f}(s, -\overline{\xi}) e^{i s \varphi} \mu_{HBH}(\overline{\xi}, \overline{\sigma}) i \widehat{f}(s, \overline{\eta}) \widehat{f}(s, \overline{\sigma}) ~ d\overline{\eta} d\overline{\xi} ds \notag \\
&\quad = - \int_1^t \int \int e^{i s \varphi} \mu_{HBH}(\overline{\xi}, \overline{\sigma}) |\overline{\xi}|^{2N-2} \partial_s \left( \widehat{f}(s, -\overline{\xi}) \widehat{f}(s, \overline{\sigma}) \right) \widehat{f}(s, \overline{\eta}) ~ d\overline{\eta} d\overline{\xi} ds \label{esthauteregHNvarphi-1} \\
&\quad \quad - \int_1^t \int \int e^{i s \varphi} \mu_{HBH}(\overline{\xi}, \overline{\sigma}) |\overline{\xi}|^{2N-2} \widehat{f}(s, -\overline{\xi}) \widehat{f}(s, \overline{\sigma}) \partial_s \widehat{f}(s, \overline{\eta}) ~ d\overline{\eta} d\overline{\xi} ds \label{esthauteregHNvarphi-2} \\
&\quad \quad + \int \int e^{i t \varphi} \mu_{HBH}(\overline{\xi}, \overline{\sigma}) |\overline{\xi}|^{2N-2} \widehat{f}(t, -\overline{\xi}) \widehat{f}(t, \overline{\sigma}) \widehat{f}(t, \overline{\eta}) ~ d\overline{\eta} d\overline{\xi} \label{esthauteregHNvarphi-3} \\
&\quad \quad - \int \int e^{i \varphi} \mu_{HBH}(\overline{\xi}, \overline{\sigma}) |\overline{\xi}|^{2N-2} \widehat{f}(1, -\overline{\xi}) \widehat{f}(1, \overline{\sigma}) \widehat{f}(1, \overline{\eta}) ~ d\overline{\eta} d\overline{\xi} \label{esthauteregHNvarphi-4}
\end{align}
\end{subequations}
Finally, we can estimate
\begin{align*}
|\eqref{esthauteregHNvarphi-1}| &\lesssim \int_1^t \Vert e^{-i s \omega(D)} \partial_s f(s) \Vert_{W^{N-2, 4}} \Vert u(s) \Vert_{H^N} \Vert u(s) \Vert_{L^4} ~ ds \\
&\lesssim \int_1^t \Vert u(s) \Vert_{L^{\infty}} \Vert u(s) \Vert_{W^{N-1, 4}} \Vert u(s) \Vert_{H^N} \Vert u(s) \Vert_{L^4} ~ ds \\
&\lesssim \int_1^t s^{-\frac{31}{24}+200\delta} \Vert u \Vert_X^2 \Vert u(s) \Vert_{H^N}^2 ~ ds \\
&\lesssim \Vert u \Vert_X^2 \Vert u \Vert_{L^{\infty} H^N}^2 \\
|\eqref{esthauteregHNvarphi-2}| &\lesssim \int_1^t \Vert u(s) \Vert_{H^N}^2 \Vert e^{-i s \omega(D)} |\nabla|^{-1} \partial_s f(s) \Vert_{L^{\infty}} ~ ds \\
&\lesssim \int_1^t \Vert u \Vert_{L^{\infty} H^N}^2 \Vert \partial_s f(s) \Vert_{H^1} ~ ds \\
&\lesssim \Vert u \Vert_{L^{\infty} H^N}^2 \Vert u \Vert_X^2 \\
|\eqref{esthauteregHNvarphi-3}| &\lesssim \Vert u(t) \Vert_{H^N}^2 \Vert u(t) \Vert_{L^{\infty}} \\
&\lesssim \Vert u \Vert_{L^{\infty} H^N}^2 \Vert u \Vert_X \\
|\eqref{esthauteregHNvarphi-4}| &\lesssim \Vert u(1) \Vert_{H^N}^2 \Vert u(1) \Vert_{L^{\infty}} \\
&\lesssim \Vert u \Vert_{L^{\infty} H^N}^2 \Vert u \Vert_X
\end{align*}
by applying Sobolev's embedding and estimates from Proposition \ref{propositionestimeesstandardscattering}. 
\end{Dem}

For $\alpha = a, b, c$, by Lemma \ref{lemstructurehalphagalpha}, $g_{\alpha}$ has the form
\begin{align*}
\widehat{g}_{\alpha}(t, \overline{\xi}) &= t \int e^{i t \varphi(\overline{\xi}, \overline{\eta})} \mu(\overline{\xi}, \overline{\eta}) \widehat{f}(t, \overline{\eta}) \widehat{f}(t, \overline{\sigma}) ~ d\overline{\eta} 
\end{align*}
for some symbol $\mu$ (depending on $\alpha$) whose form is not important for the proof of the following lemma. 

\begin{Lem} For any $\alpha \in \{ a, b, c \}$, any $t > 0$, 
\begin{subequations}
\begin{align}
\Vert g_{\alpha}(t) \Vert_{H^1} &\lesssim \langle t \rangle^{-\frac{1}{12}+100\delta} \Vert u \Vert_X^2 \label{estimeegL2-simple} \\
\Vert |\nabla|^{\frac{1}{2}} g_{\alpha}(t) \Vert_{H^1} &\lesssim \langle t \rangle^{-\frac{1}{6}+100\delta} \Vert u \Vert_X^2 \label{estimeegL2-dec1.6} \\
\Vert e^{-i t \omega(D)} \nabla g_{\alpha}(t) \Vert_{L^4} &\lesssim t^{\frac{1}{6}} \langle t \rangle^{-\frac{17}{24}+150\delta} \Vert u \Vert_X^2 \label{estimeegL4}
\end{align}
\end{subequations} 
\end{Lem}

\begin{Dem}
For \eqref{estimeegL2-simple}, we can simply estimate: 
\begin{align*}
\Vert g_{\alpha}(t) \Vert_{H^1} &\lesssim t \Vert u(t) \Vert_{L^4} \Vert \langle \nabla \rangle u(t) \Vert_{L^4} \\
&\lesssim \langle t \rangle^{-\frac{1}{12}+100\delta} \Vert u \Vert_X^2 
\end{align*}

For \eqref{estimeegL2-dec1.6}, we may separate the interaction intro three pieces: by symmetry of the expression, we can always assume that we localise to have $|\overline{\eta}| \lesssim |\overline{\sigma}|$, and 
\begin{enumerate}
\item either from the two interacting frequencies $\overline{\eta}$ and $\overline{\sigma}$, one is away enough from $\widehat{\mathcal{P}}$ and is not small with respect to the other; 
\item or the interaction is of type $BBB$ or $BHH$ with $\overline{\eta}$ and $\overline{\sigma}$ near $\widehat{\mathcal{P}}$; 
\item or the interaction is of type $HBH$ with $\overline{\sigma}$ near the plane. 
\end{enumerate}
In the first case, denoting by $\mu$ again a symbol localising on such an area, we may estimate 
\begin{align*}
\Vert t \langle \nabla \rangle |\nabla|^{\frac{1}{2}} T_{\mu}[f, f](t) \Vert_{L^2} &\lesssim t \Vert \langle \nabla \rangle |\nabla|^{-1} u(t) \Vert_{L^2} \Vert |\nabla|^{\frac{3}{2}} (1 - m_{\widehat{\mathcal{P}}}(D)) u(t) \Vert_{L^{\infty}} \\
&\lesssim \langle t \rangle^{-\frac{1}{6}+100\delta} \Vert u \Vert_X^2 
\end{align*}
In the second case, denoting again by $\mu$ a symbol having the correct localisation, we can further localise to have $|\eta_0| \lesssim |\sigma_0|$ up to symmetry, and we estimate 
\begin{align*}
\Vert t \langle \nabla \rangle |\nabla|^{\frac{1}{2}} T_{\mu}[f, f](t) \Vert_{L^2} &\lesssim t \Vert |\partial_x|^{-\frac{1}{3}} \langle \nabla \rangle |\nabla|^{-\frac{2}{3}} u(t) \Vert_{L^{\infty}} \Vert |\partial_x|^{\frac{1}{3}} |\nabla|^{\frac{7}{6}} u(t) \Vert_{L^{\infty}} \\
&\lesssim \langle t \rangle^{-\frac{1}{6}+100\delta} \Vert u \Vert_X^2 
\end{align*}

Finally, for the last case, we need to decompose a bit further. Note that we now reduced to study the following term in Fourier
\begin{align}
t \int e^{i t \varphi} \mu(\overline{\xi}, \overline{\eta}) \mu_{HBH}(\overline{\xi}, \overline{\eta}) m_{\widehat{\mathcal{P}}}(\overline{\sigma}) \widehat{f}(t, \overline{\eta}) \widehat{f}(t, \overline{\sigma}) ~ d\overline{\eta} \label{preuveestimeegtermeHBHPreste} 
\end{align}
If $\overline{\eta}$ is not near the plane, we have that 
\begin{align*}
\widehat{X}_a(\overline{\eta}) \cdot \nabla_{\overline{\eta}} \varphi &\simeq |\overline{\sigma}|^2 
\end{align*}
while if $\overline{\eta}$ is near the plane, 
\begin{align*}
\partial_{\eta_0} \varphi &\simeq |\overline{\sigma}|^2 
\end{align*}
In any case, we can apply an integration by parts in $\overline{\eta}$, and we get: 
\begin{subequations}
\begin{align}
\eqref{preuveestimeegtermeHBHPreste} &= \int e^{i t \varphi} |\overline{\sigma}|^{-2} |\overline{\eta}|^{-1} \mu(\overline{\xi}, \overline{\eta}) \mu_{HBH}(\overline{\xi}, \overline{\eta}) m_{\widehat{\mathcal{P}}}(\overline{\sigma}) \widehat{f}(t, \overline{\eta}) \widehat{f}(t, \overline{\sigma}) ~ d\overline{\eta} \label{preuveestimeegtermeHBHPreste1} \\
&\quad + \int e^{i t \varphi} |\overline{\sigma}|^{-2} \mu(\overline{\xi}, \overline{\eta}) \mu_{HBH}(\overline{\xi}, \overline{\eta}) m_{\widehat{\mathcal{P}}}(\overline{\sigma}) \widehat{Y}(\overline{\eta}) \cdot \nabla_{\overline{\eta}} \widehat{f}(t, \overline{\eta}) \widehat{f}(t, \overline{\sigma}) ~ d\overline{\eta} \label{preuveestimeegtermeHBHPreste2} \\
&\quad + \int e^{i t \varphi} |\overline{\sigma}|^{-2} \mu(\overline{\xi}, \overline{\eta}) \mu_{HBH}(\overline{\xi}, \overline{\eta}) m_{\widehat{\mathcal{P}}}(\overline{\sigma}) \widehat{f}(t, \overline{\eta}) \nabla_{\overline{\eta}} \widehat{f}(t, \overline{\sigma}) ~ d\overline{\eta} \label{preuveestimeegtermeHBHPreste3}
\end{align}
\end{subequations} 
where the symbol $\mu$ may change from line to line, and where $\widehat{Y}$ is a vector field equal either to $\widehat{X}_a$ or to $(1, 0, 0)$, such that $\Vert Y f \Vert_{L^2} \lesssim \Vert u \Vert_X$. We then estimate for $t \geq 1$: 
\begin{align*}
\Vert \langle \overline{\xi} \rangle |\overline{\xi}|^{\frac{1}{2}} \eqref{preuveestimeegtermeHBHPreste1} \Vert_{L^2} &\lesssim \Vert |\nabla|^{-\frac{3}{2}} u(t) \Vert_{L^3} \Vert \langle \nabla \rangle |\nabla|^{-1} u(t) \Vert_{L^6} \\
&\lesssim \Vert |\nabla|^{-1} u(t) \Vert_{L^2} \left( \Vert u(t) \Vert_{L^6} + \Vert m_{\widehat{\mathcal{P}}}(D) |\nabla|^{-1} u(t) \Vert_{L^6} \right) \\
&\lesssim \Vert u \Vert_X \left( \langle t \rangle^{-\frac{1}{2}} \Vert u \Vert_X + \Vert |\nabla|^{-1} u(t) \Vert_{L^2}^{\frac{2}{3}} \Vert m_{\widehat{\mathcal{P}}}(D) |\nabla|^{\frac{1}{2}} u(t) \Vert_{L^{\infty}}^{\frac{1}{3}} \right) \\
&\lesssim t^{-\frac{1}{3}+100\delta} \Vert u \Vert_X^2 \\
\Vert \langle \overline{\xi} \rangle |\overline{\xi}|^{\frac{1}{2}} \eqref{preuveestimeegtermeHBHPreste2} \Vert_{L^2} &\lesssim \Vert |\nabla|^{-\frac{1}{2}} e^{-i t \omega(D)} Yf(t) \Vert_{L^3} \Vert \langle \nabla \rangle |\nabla|^{-1} u(t) \Vert_{L^6} \\
&\lesssim \Vert t^{-\frac{1}{3}+100\delta} \Vert u \Vert_X^2 \\
\Vert \langle \overline{\xi} \rangle |\overline{\xi}|^{\frac{1}{2}} \eqref{preuveestimeegtermeHBHPreste3} \Vert_{L^2} &\lesssim \Vert |\nabla|^{-\frac{1}{2}} u(t) \Vert_{L^3} \Vert \langle \nabla \rangle |\nabla|^{-1} m_{\widehat{\mathcal{P}}}(D) e^{-i t \omega(D)} (x, y) f(t) \Vert_{L^6} \\
&\lesssim \Vert |\nabla|^{-1} u(t) \Vert_{L^2}^{\frac{3}{7}} \Vert u(t) \Vert_{L^4}^{\frac{4}{7}} \Vert m_{\widehat{\mathcal{P}}}(D) (x, y) f(t) \Vert_{H^1} \\
&\lesssim t^{-\frac{2}{7}+100\delta} \Vert u \Vert_X^2 
\end{align*}
by Gagliardo-Nirenberg's inequality. If $t \leq 1$, we have directly that 
\begin{align*}
\Vert \langle \overline{\xi} \rangle |\overline{\xi}|^{\frac{1}{2}} \eqref{preuveestimeegtermeHBHPreste1} \Vert_{L^2} &\lesssim \Vert |\nabla|^{-1} u(t) \Vert_{L^6} \Vert \langle \nabla \rangle |\nabla|^{-\frac{3}{2}} u(t) \Vert_{L^3} \\
&\lesssim \Vert u(t) \Vert_{L^2} \Vert \langle \nabla \rangle |\nabla|^{-1} u(t) \Vert_{L^2} \\
&\lesssim \Vert u \Vert_X^2 \\
\Vert \langle \overline{\xi} \rangle |\overline{\xi}|^{\frac{1}{2}} \eqref{preuveestimeegtermeHBHPreste2} \Vert_{L^2} &\lesssim \Vert |\nabla|^{-\frac{1}{2}} e^{-i t \omega(D)} Yf(t) \Vert_{L^3} \Vert \langle \nabla \rangle |\nabla|^{-1} u(t) \Vert_{L^6} \\
&\lesssim \Vert u \Vert_X^2 \\
\Vert \langle \overline{\xi} \rangle |\overline{\xi}|^{\frac{1}{2}} \eqref{preuveestimeegtermeHBHPreste3} \Vert_{L^2} &\lesssim \Vert |\nabla|^{-\frac{1}{2}} u(t) \Vert_{L^3} \Vert \langle \nabla \rangle |\nabla|^{-1} m_{\widehat{\mathcal{P}}}(D) e^{-i t \omega(D)} (x, y) f(t) \Vert_{L^6} \\
&\lesssim \Vert u(t) \Vert_{L^2} \Vert m_{\widehat{\mathcal{P}}}(D) (x, y) f(t) \Vert_{H^1} \\
&\lesssim \Vert u \Vert_X^2
\end{align*}
which concludes. 

Finally, for \eqref{estimeegL4}, we notice that $\nabla g_{\alpha}(t)$ has a similar structure to $t \partial_t f(t)$, up to the $x$-derivative. We can therefore apply similar estimate, replacing only $\Vert \partial_x u(t) \Vert_{L^{\infty}}$ by $\Vert \nabla u(t) \Vert_{L^{\infty}}$ and thus losing $\langle t \rangle^{\frac{1}{12}}$, which corresponds to the announced estimate. 
\end{Dem} 

\begin{Lem} Let $j, k \in \mathbb{Z}$ such that $k \leq -10$, for any $t \geq 1$, we have 
\begin{align*}
\Vert \psi_{j, k}^{\widehat{\mathcal{C}}}(D) g_b(t) \Vert_{L^2} &\lesssim 2^{k-\frac{j}{2}} \langle 2^j \rangle^{-1} t^{-\frac{1}{6}+100\delta} \Vert u \Vert_X^2 + t^{-\frac{1}{2}+201\delta} 2^{-3j} \Vert u \Vert_X^2 
\end{align*}
\end{Lem}

\begin{Dem}
We apply Equation \ref{decompositionfinemg-gainkout} from Lemma \ref{lemdecompositionFGH} and decompose 
\begin{align*}
m_g(\overline{\xi}, \overline{\eta}) &= m_b(\overline{\xi}) m_1(\overline{\xi}, \overline{\eta}) + \nabla_{\overline{\eta}} \varphi(\overline{\xi}, \overline{\eta}) ~ m_2(\overline{\xi}, \overline{\eta})
\end{align*}
where $m_1$ is of degree $0$ and $m_2$ of degree $-2$. On the support of $\psi_{j, k}^{\widehat{\mathcal{C}}}(D)$, $m_b \lesssim 2^k$ and so the contribution of the first part follows immediately from the previous estimates. 

Now, for the part gaining a factor $\nabla_{\overline{\eta}} \varphi$, we can apply an integration by parts: 
\begin{subequations}
\begin{align}
&t \int e^{i t \varphi} \nabla_{\overline{\eta}} \varphi m_2 \widehat{f}(t, \overline{\eta}) \widehat{f}(t, \overline{\sigma}) ~ d\overline{\eta} \notag \\
&\quad = \int e^{i t \varphi} \nabla_{\overline{\eta}} m_2 \widehat{f}(t, \overline{\eta}) \widehat{f}(t, \overline{\sigma}) ~ d\overline{\eta} \label{preuveestimeefineg1} \\
&\quad \quad + \int e^{i t \varphi} m_2 \nabla_{\overline{\eta}} \left( \widehat{f}(t, \overline{\eta}) \widehat{f}(t, \overline{\sigma}) \right) ~ d\overline{\eta} \label{preuveestimeefineg2}
\end{align}
\end{subequations}
Then, we can estimate: 
\begin{align*}
\Vert \psi_{j, k}^{\widehat{\mathcal{C}}} \eqref{preuveestimeefineg1} \Vert_{L^2} 
&\lesssim 2^{-3j} \Vert \xi_0^3 \psi_{j, k}^{\widehat{\mathcal{C}}} \eqref{preuveestimeefineg1} \Vert_{L^2} \\
&\lesssim 2^{-3j} \Vert |\nabla|^{-1} u(t) \Vert_{L^2} \Vert \partial_x u(t) \Vert_{L^{\infty}} \\
&\lesssim 2^{-3j} t^{-\frac{13}{12}+100\delta} \Vert u \Vert_X^2 \\
\Vert \psi_{j, k}^{\widehat{\mathcal{C}}} \eqref{preuveestimeefineg2} \Vert_{L^2} &\lesssim 2^{-3j} \left( \Vert \nabla (x, y) f(t) \Vert_{L^2} \Vert u(t) \Vert_{L^{\infty}} + \Vert (x, y) f(t) \Vert_{L^2} \Vert \nabla u(t) \Vert_{L^{\infty}} \right) \\
&\lesssim 2^{-3j} t^{-\frac{1}{2}+201\delta} \Vert u \Vert_X^2 
\end{align*}
as wanted. 
\end{Dem}

\section{\texorpdfstring{$H^1$}{H1} estimate of \texorpdfstring{$h$}{h}} \label{section-estimeepoidssimple}

The goal of this section is to show the following proposition: 

\begin{Prop} There exist $C > 0$ such that, for every $\alpha \in \{ a, b, c \}$ and $t \geq 0$, 
\begin{align*}
\Vert h_{\alpha}(t) \Vert_{H^1}^2 &\leq \Vert (x, y) u_0 \Vert_{H^1}^2 + C \Vert u \Vert_X^3 \left( 1 + \Vert u \Vert_X \right) 
\end{align*} 
Furthermore, 
\begin{align}
\int_0^t s^{1-\delta} \langle s \rangle^{2\delta} \Vert \partial_s h_{\alpha}(s) \Vert_{L^2}^2 ~ ds &\lesssim \Vert u \Vert_X^2 \left( \Vert u \Vert_X + \sum_{\beta = a, c} \left( \int_0^t s^{\frac{63}{64}} \langle s \rangle^{\frac{1}{50}} \Vert e^{-i s \omega(D)} \nabla h_{\beta}(s) \Vert_{L^4}^4 ~ ds \right)^{\frac{1}{4}} \right)^2 \label{estimeeaprioriscatteringL2halpha} 
\end{align}
\label{prop-estimeeaprioriH1-h} 
\end{Prop}

\begin{Dem}
We apply Lemma \ref{lemstructurehalphagalpha} which gives an expression fo $h_{\alpha}$, $\alpha = a, b, c$, and we differentiate it in time to get $\partial_t h_{\alpha}$. 

We start with $h_b$. 

We estimate: 
\begin{align*}
\Vert \partial_t \eqref{lemstructurehb-02-symeta} \Vert_{L^2} &\lesssim \Vert e^{-i t \omega(D)} \partial_x h_b(t) \Vert_{L^4} \Vert u(t) \Vert_{L^4} \\
&\lesssim t^{-\frac{511}{512}} \langle t \rangle^{-\frac{1}{256}} \Vert u \Vert_X^2 \\
\Vert \partial_t \mathcal{F}^{-1} \eqref{lemstructurehb-03-resxetaac} \Vert_{H^1} &\lesssim \sum_{\beta = a, c} \Vert h_{\beta}(t) \Vert_{H^1} \Vert \partial_x u(t) \Vert_{L^{\infty}} \\
&\lesssim t^{-\frac{5}{6}} \langle t \rangle^{-\frac{1}{4}+100\delta} \Vert u \Vert_X^2 \\
\Vert \partial_t \mathcal{F}^{-1} \eqref{lemstructurehb-04-resxetabbon} \Vert_{H^1} + \Vert \partial_t \mathcal{F}^{-1} \eqref{lemstructurehb-05-resxetabbonbis} \Vert_{H^1} &\lesssim \Vert h_b(t) \Vert_{H^1} \Vert \partial_x u(t) \Vert_{L^{\infty}} \\
&\lesssim t^{-\frac{5}{6}} \langle t \rangle^{-\frac{1}{4}+100\delta} \Vert u \Vert_X^2 \\
\Vert \partial_t \mathcal{F}^{-1} \eqref{lemstructurehb-06-resxetabphi} \Vert_{H^1} &\lesssim \Vert h_b(t) \Vert_{H^1} \Vert \partial_x u(t) \Vert_{L^{\infty}} + \Vert e^{-i t \omega(D)} \partial_x h_b(t) \Vert_{L^4} \Vert \nabla u(t) \Vert_{L^4} \\
&\lesssim t^{-\frac{5}{6}} \langle t \rangle^{-\frac{5}{24}+100\delta} \Vert u \Vert_X^2 \\
\Vert \partial_t \mathcal{F}^{-1} \eqref{lemstructurehb-11-dersymb} \Vert_{H^1} &\lesssim \Vert \langle \nabla \rangle u(t) \Vert_{L^4} \Vert u(t) \Vert_{L^4} \\
&\lesssim t^{-\frac{5}{6}} \langle t \rangle^{-\frac{1}{4}+100\delta} \Vert u \Vert_X^2 \\
\Vert \partial_t \mathcal{F}^{-1} \eqref{lemstructurehb-13-termecubique} \Vert_{H^1} &\lesssim t \Vert \partial_x u(t) \Vert_{L^{\infty}} \Vert \langle \nabla \rangle u(t) \Vert_{L^4} \Vert u(t) \Vert_{L^4} \\
&\lesssim t^{-\frac{2}{3}} \langle t \rangle^{-\frac{1}{2}+200\delta} \Vert u \Vert_X^3 
\end{align*}
Here above, for \eqref{lemstructurehb-06-resxetabphi}, we used that 
\begin{align*}
\varphi &= O(\eta_0 \overline{\sigma}) + O(\sigma_0 \overline{\eta}) 
\end{align*}
which follows from Lemma \ref{lemcalculvarphiKdV1D}. 

This proves \eqref{estimeeaprioriscatteringL2halpha} for $\alpha = b$. 

If $\alpha = a$ or $c$, we can proceed the same way, noting that the expressions given by Lemma \ref{lemstructurehalphagalpha} are all of the form of one already estimated in the case $\alpha = b$ (or simpler), except for the symmetric term \eqref{lemstructurehb-02-symeta}, for which it is needed to use the integration in time: 
\begin{align*}
&\int_0^t s^{1-\delta} \langle s \rangle^{2\delta} \Vert \partial_t \eqref{lemstructurehb-02-symeta} \Vert_{L^2}^2 ~ ds \lesssim \int_0^t s^{1-\delta} \langle s \rangle^{2\delta} \Vert e^{-i s \omega(D)} \nabla h_{\alpha}(s) \Vert_{L^4}^2 \Vert u(s) \Vert_{L^4}^2 ~ ds \\
&\lesssim \left( \int_0^t s^{1-\delta} \langle s \rangle^{2\delta} \Vert e^{-i s \omega(D)} \nabla h_{\alpha}(s) \Vert_{L^4}^4 ~ ds \right)^{\frac{1}{2}} \left( \int_0^t s^{1-\delta} \langle s \rangle^{2\delta} \Vert u(s) \Vert_{L^4}^4 ~ ds \right)^{\frac{1}{2}} \\
&\lesssim \left( \int_0^t s^{\frac{63}{64}} \langle s \rangle^{\frac{1}{50}} \Vert e^{-i s \omega(D)} \nabla h_{\alpha}(s) \Vert_{L^4}^4 ~ ds \right)^{\frac{1}{2}} \Vert u \Vert_X^2 
\end{align*}
This proves \eqref{estimeeaprioriscatteringL2halpha}. 

To get the $H^1$ estimate, we apply Proposition \ref{propositionestimeedispersiveL4totaleintegrale} to improve the previous estimate into 
\begin{align*}
\int_0^t \Vert \partial_s h_{\alpha}(s) \Vert_{L^2} ~ ds &\lesssim \Vert u \Vert_X 
\end{align*}
for any $\alpha$. Then, we apply Duhamel's formula directly on the norm: for instance, for $\alpha = b$, 
\begin{align*}
&\Vert h_b(t) \Vert_{H^1}^2 = \Vert h_b(0) \Vert_{H^1}^2 + 2 \int_0^t \int \langle \nabla \rangle h_b(s, x, y) ~ \langle \nabla \rangle \partial_s h_b(s, x, y) ~ dx dy ds \\
&\leq \Vert m_b(D) X_b f(0) \Vert_{H^1}^2 + 2 \int_0^t \int \langle \overline{\xi} \rangle \widehat{h}_b(s, -\overline{\xi}) ~ \langle \overline{\xi} \rangle \partial_s \eqref{lemstructurehb-02-symeta}(s, \overline{\xi}) ~ d\overline{\xi} ds \\
&\quad + 2 \int_0^t \Vert h_b(s) \Vert_{H^1} \Vert \partial_s \left( h_b - \mathcal{F}^{-1} \eqref{lemstructurehb-02-symeta} \right) \Vert_{H^1} ~ ds \\
&\leq \Vert (x, y) f(0) \Vert_{H^1}^2 + 4 \int_0^t \int \int \langle \overline{\xi} \rangle^2 \widehat{h}_b(s, -\overline{\xi}) ~ e^{i s \varphi} i \eta_0 \widehat{h}_b(s, \overline{\eta}) \widehat{f}(s, \overline{\sigma}) ~ d\overline{\eta} d\overline{\xi} ds 
+ C \Vert u \Vert_X^3 \left( 1 + \Vert u \Vert_X \right) 
\end{align*}
by using that every term except \eqref{lemstructurehb-02-symeta} satisfy a $L^1_t H^1_{x, y}$ estimate for their time derivative. 

For the term coming from \eqref{lemstructurehb-02-symeta}, we apply the symmetry $\mathfrak{s}: (\overline{\xi}, \overline{\eta}) \mapsto (-\overline{\eta}, -\overline{\xi})$: 
\begin{align*}
&4 \int_0^t \int \int \langle \overline{\xi} \rangle^2 \widehat{h}_b(s, -\overline{\xi}) ~ e^{i s \varphi} i \eta_0 \widehat{h}_b(s, \overline{\eta}) \widehat{f}(s, \overline{\sigma}) ~ d\overline{\eta} d\overline{\xi} ds \\
&\quad = - 4 \int_0^t \int \int \langle \overline{\eta} \rangle^2 \widehat{h}_b(s, \overline{\eta}) ~ e^{i s \varphi} i \xi_0 \widehat{h}_b(s, -\overline{\xi}) \widehat{f}(s, \overline{\sigma}) ~ d\overline{\eta} d\overline{\xi} ds \\
&\quad = 2 i \int_0^t \int \left( \eta_0 \int \langle \overline{\xi} \rangle^2 - \xi_0 \langle \overline{\eta} \rangle^2 \right) \widehat{h}_b(s, -\overline{\xi}) ~ e^{i s \varphi} \widehat{h}_b(s, \overline{\eta}) \widehat{f}(s, \overline{\sigma}) ~ d\overline{\eta} d\overline{\xi} ds
\end{align*}
But
\begin{align*}
\eta_0 \langle \overline{\xi} \rangle^2 - \xi_0 \langle \overline{\eta} \rangle^2 &= - \sigma_0 + \eta_0 \overline{\xi} \cdot \overline{\sigma} 
+ \eta_0 \overline{\xi} \cdot \overline{\eta} 
- \xi_0 \overline{\eta} \cdot \overline{\xi} + \xi_0 \overline{\eta} \cdot \overline{\sigma} \\
&= - \sigma_0 + \eta_0 \overline{\xi} \cdot \overline{\sigma} - \sigma_0 \overline{\xi} \cdot \overline{\eta} + \xi_0 \overline{\eta} \cdot \overline{\sigma}
\end{align*}
We can therefore estimate
\begin{align*}
&4 \int_0^t \int \int \langle \overline{\xi} \rangle^2 \widehat{h}_b(s, -\overline{\xi}) ~ e^{i s \varphi} i \eta_0 \widehat{h}_b(s, \overline{\eta}) \widehat{f}(s, \overline{\sigma}) ~ d\overline{\eta} d\overline{\xi} ds \\
&\quad \lesssim \int_0^t \left( \Vert h_b(s) \Vert_{H^1}^2 \Vert \partial_x u(s) \Vert_{L^{\infty}} + \Vert e^{-i s \omega(D)} \partial_x h_b(s) \Vert_{L^4} \Vert h_b(s) \Vert_{H^1} \Vert \nabla u(s) \Vert_{L^4} \right) ~ ds \\
&\quad \lesssim \int_0^t \Vert u \Vert_X^3 s^{-\frac{5}{6}} \langle s \rangle^{-\frac{5}{24}+100\delta} ~ ds \\
&\quad \lesssim \Vert u \Vert_X^3 
\end{align*}
which is enough. In the case $\alpha = a$ or $c$, we can apply exactly the same steps for the symmetric term, and use the integrated $L^4$ estimate at the end instead. 
\end{Dem}

\section{Quadratic weight estimates except (b, b)} \label{section-estimees-quadratiques-horsbb} 

In this section, we show: 

\begin{Prop} Assume $\Vert u \Vert_X$ is small enough and let $(\alpha, \beta) \in \{ a, b, c \}$ be such that $(\alpha, \beta) \neq (b, b)$. Then
\begin{align*}
\Vert \nabla m_{\alpha}(D) X_{\alpha} h_{\beta}(t) \Vert_{L^2}^2 &\lesssim \Vert (1 + x^2 + |y|^2) u_0 \Vert_{H^1}^2 
+ \Vert u \Vert_X^3 
\end{align*}
\end{Prop}

\begin{Dem}
We apply Lemmas \ref{lemmestructurehquadabcsurac} and \ref{lemmestructurehquadacsurb} and we estimate every term. Since we can decompose as 
\begin{align*}
|\nabla| m_{\alpha}(D) X_{\alpha} h_{\beta}(t) &= h_{\alpha \beta}(t) + g_{\alpha \beta}(t) 
\end{align*}
we can estimate separately each contribution. 

\paragraph{1. Estimate on $g$} We first show: 
\begin{align*}
\Vert g_{\alpha \beta}(t) \Vert_{L^2} &\lesssim \Vert u \Vert_X^2 
\end{align*}

Indeed: 
\begin{align*}
\Vert \eqref{lemstructhquadac-gdersymb} \Vert_{L^2} &\lesssim t \Vert u(t) \Vert_{L^4}^2 \\
&\lesssim t^{\frac{1}{6}} \langle t \rangle^{-\frac{1}{4}+100\delta} \Vert u \Vert_X^2 \\
\Vert \eqref{lemstructhquadac-gvarphi} \Vert_{L^2} &\lesssim |1 - \chi(t)| \sum_{\gamma_1, \gamma_2 = a, b, c} \Vert e^{-i t \omega(D)} |\nabla|^{-1} h_{\gamma_1}(t) \Vert_{L^6} \Vert e^{-i t \omega(D)} h_{\gamma_2}(t) \Vert_{L^3} \\
&\lesssim |1 - \chi(t)| \sum_{\gamma = a, b, c} \Vert h_{\gamma}(t) \Vert_{H^1}^2 \\
&\lesssim |1 - \chi(t)| \Vert u \Vert_X^2 \\
\Vert \eqref{lemstructhquadacsurb-g-t2} \Vert_{L^2} &\lesssim t^2 \Vert \nabla u(t) \Vert_{L^{\infty}} \Vert u(t) \Vert_{L^4}^2 \\
&\lesssim t^{\frac{1}{3}} \langle t \rangle^{-\frac{5}{12}+200\delta} \Vert u \Vert_X^3 
\end{align*}
which is enough assuming $\Vert u \Vert_X \lesssim 1$. 

It only remains \eqref{lemstructhquadac-gresx}. If $\beta = b$, we can estimate as
\begin{align*}
\Vert \eqref{lemstructhquadac-gresx} \Vert_{L^2} &\lesssim t \Vert e^{-i t \omega(D)} \partial_x h_b(t) \Vert_{L^4} \Vert u(t) \Vert_{L^4} \\
&\lesssim t^{\frac{1}{12}+\frac{1}{1024}} \langle t \rangle^{-\frac{1}{8}-\frac{1}{512}+50\delta} \Vert u \Vert_X^2
\end{align*}
which is enough. If $\beta = a, c$, we separate: 
\begin{subequations}
\begin{align}
\eqref{lemstructhquadac-gresx} &= t (1 - \chi(\overline{\xi})) \int e^{i t \varphi} \mu(\overline{\xi}, \overline{\eta}) \mu_{HHB}(\overline{\xi}, \overline{\eta}) m_{\widehat{\mathcal{C}}}(\overline{\eta}) |\overline{\eta}| \widehat{h}_{\beta}(t, \overline{\eta}) \widehat{f}(t, \overline{\sigma}) ~ d\overline{\eta} \label{preuveestimeequadratiquegresx1} \\
&\quad + t \chi(\overline{\xi}) \int e^{i t \varphi} \mu(\overline{\xi}, \overline{\eta}) \mu_{HHB}(\overline{\xi}, \overline{\eta}) m_{\widehat{\mathcal{C}}}(\overline{\eta}) |\overline{\eta}| \widehat{h}_{\beta}(t, \overline{\eta}) \widehat{f}(t, \overline{\sigma}) ~ d\overline{\eta} \label{preuveestimeequadratiquegresx1bis} \\
&\quad + t \int e^{i t \varphi} \mu(\overline{\xi}, \overline{\eta}) (1 - \mu_{HHB}(\overline{\xi}, \overline{\eta})) m_{\widehat{\mathcal{C}}}(\overline{\eta}) |\overline{\eta}| \widehat{h}_{\beta}(t, \overline{\eta}) \widehat{f}(t, \overline{\sigma}) ~ d\overline{\eta} \label{preuveestimeequadratiquegresx2} \\
&\quad + t \int e^{i t \varphi} \mu(\overline{\xi}, \overline{\eta}) \left( 1 - m_{\widehat{\mathcal{C}}}(\overline{\eta}) \right) |\overline{\eta}| \widehat{h}_{\beta}(t, \overline{\eta}) \widehat{f}(t, \overline{\sigma}) ~ d\overline{\eta} \label{preuveestimeequadratiquegresx3}
\end{align}
\end{subequations} 
First, we can estimate: 
\begin{align*}
\Vert \eqref{preuveestimeequadratiquegresx1bis} \Vert_{L^2} &\lesssim t \Vert e^{-i t \omega(D)} \chi(D) \nabla h_{\beta}(t) \Vert_{L^4} \Vert u(t) \Vert_{L^4} \\
&\lesssim t^{\frac{11}{96}} \langle t \rangle^{-\frac{1}{8}+51\delta} \Vert u \Vert_X^2 \\
\Vert \eqref{preuveestimeequadratiquegresx2} \Vert_{L^2} &\lesssim t \Vert e^{-i t \omega(D)} |\nabla|^{\frac{11}{16}} \langle \nabla \rangle^{-\frac{1}{16}} h_{\beta}(t) \Vert_{L^4} \Vert |\nabla|^{\frac{5}{16}} \langle \nabla \rangle^{\frac{1}{16}} u(t) \Vert_{L^4} \\
&\lesssim t^{\frac{11}{96}} \langle t \rangle^{-\frac{1}{8}+51\delta} \Vert u \Vert_X^2 \\
\Vert \eqref{preuveestimeequadratiquegresx3} \Vert_{L^2} &\lesssim t \Vert e^{-i t \omega(D)} \nabla (1 - m_{\widehat{\mathcal{C}}}(D)) h_{\beta}(t) \Vert_{L^4} \Vert u(t) \Vert_{L^4} \\
&\lesssim t^{\frac{65}{128}} \langle t \rangle^{-\frac{13}{24}+50\delta} \Vert u \Vert_X^2 
\end{align*}
where we used Lemma \ref{lemestimeeL4displocalezoneC}. 

For \eqref{preuveestimeequadratiquegresx1}, we apply an integration by parts in the direction $\widehat{X}_a(\overline{\eta})$, noting that 
\begin{align*}
\widehat{X}_a(\overline{\eta}) \cdot \nabla_{\overline{\eta}} \varphi \simeq |\overline{\eta}|^2 
\end{align*}
when $|\overline{\sigma}| \ll |\overline{\eta}|$. This leads to: 
\begin{subequations}
\begin{align}
\eqref{preuveestimeequadratiquegresx1} &= (1 - \chi(\overline{\xi})) \int e^{i t \varphi} \mu(\overline{\xi}, \overline{\eta}) \mu_{HHB}(\overline{\xi}, \overline{\eta}) m_{\widehat{\mathcal{C}}}(\overline{\eta}) |\overline{\eta}|^{-1} |\overline{\sigma}|^{-1} \widehat{h}_{\beta}(t, \overline{\eta}) \widehat{f}(t, \overline{\sigma}) ~ d\overline{\eta} \label{preuveestimeequadratiquegresx1-1} \\
&\quad + (1 - \chi(\overline{\xi})) \int e^{i t \varphi} \mu(\overline{\xi}, \overline{\eta}) \mu_{HHB}(\overline{\xi}, \overline{\eta}) m_{\widehat{\mathcal{C}}}(\overline{\eta}) |\overline{\eta}|^{-1} \widehat{X}_a(\overline{\eta}) \cdot \nabla_{\overline{\eta}} \widehat{h}_{\beta}(t, \overline{\eta}) \widehat{f}(t, \overline{\sigma}) ~ d\overline{\eta} \label{preuveestimeequadratiquegresx1-2} \\
&\quad + (1 - \chi(\overline{\xi})) \int e^{i t \varphi} \mu(\overline{\xi}, \overline{\eta}) \mu_{HHB}(\overline{\xi}, \overline{\eta}) m_{\widehat{\mathcal{C}}}(\overline{\eta}) |\overline{\eta}|^{-1} \widehat{h}_{\beta}(t, \overline{\eta}) \nabla_{\overline{\eta}} \widehat{f}(t, \overline{\sigma}) ~ d\overline{\eta} \label{preuveestimeequadratiquegresx1-3} 
\end{align}
\end{subequations} 
where the symbol can change from line to line. We then estimate: 
\begin{align*}
\Vert \eqref{preuveestimeequadratiquegresx1-1} \Vert_{L^2} &\lesssim \Vert e^{-i t \omega(D)} (1 - \chi(D)) h_{\beta}(t) \Vert_{L^4} \Vert |\nabla|^{-1} u(t) \Vert_{L^4} \\
&\lesssim \langle t \rangle^{-\frac{15}{32}} \Vert u \Vert_X^2 \\
\Vert \eqref{preuveestimeequadratiquegresx1-2} \Vert_{L^2} &\lesssim \Vert |\overline{\xi}|^{\frac{3}{4}} \eqref{preuveestimeequadratiquegresx1-2} \Vert_{L^{\frac{4}{3}}} \\
&\lesssim \Vert \nabla X_a h_{\beta}(t) \Vert_{L^2} \Vert u(t) \Vert_{L^4} \\
&\lesssim \langle t \rangle^{-\frac{1}{2}} \Vert u \Vert_X^2 \\
\Vert \eqref{preuveestimeequadratiquegresx1-3} \Vert_{L^2} &\lesssim \Vert e^{-i t \omega(D)} (1 - \chi(D)) h_{\beta}(t) \Vert_{L^3} \Vert e^{-i t \omega(D)} (x, y) f(t) \Vert_{L^6} \\
&\lesssim \Vert h_{\beta}(t) \Vert_{L^2}^{\frac{1}{3}} \Vert e^{-i t \omega(D)} (1 - \chi(D)) h_{\beta}(t) \Vert_{L^4}^{\frac{2}{3}} \Vert \nabla (x, y) f(t) \Vert_{L^2} \\
&\lesssim \langle t \rangle^{-\frac{7}{48}+101\delta} \Vert u \Vert_X^2 
\end{align*}
which is enough. Above, we used again Lemma \ref{lemestimeeL4displocalezoneC}. 

\paragraph{2. Estimate on $h$} To estimate $h_{\alpha \beta}$, we apply Duhamel's formula: 
\begin{align*}
\Vert h_{\alpha \beta}(t) \Vert_{L^2}^2 &= \Vert h_{\alpha \beta}(0) \Vert_{L^2}^2 + 2 \int_0^t \int h_{\alpha \beta}(s, x, y) ~ \partial_s h_{\alpha \beta}(s, x, y) ~ dx dy ds 
\end{align*}
We have indeed: 
\begin{align*}
\Vert h_{\alpha \beta}(0) \Vert_{L^2} &= \Vert \nabla m_{\alpha}(D) X_{\alpha} m_{\beta}(D) X_{\beta} u_0 \Vert_{L^2} \lesssim \Vert (1 + x^2 + |y|^2) u_0 \Vert_{H^1}
\end{align*}
as wanted. 

Then, we develop $\partial_s h_{\alpha \beta}$ using the expression from Lemmas \ref{lemmestructurehquadabcsurac} and \ref{lemmestructurehquadacsurb} and estimate each contribution. 

\paragraph{2.1. Direct estimates} Some of the term can be estimated directly: 
\begin{align*}
&\begin{aligned}
\int_0^t s^{1-\delta} \langle s \rangle^{2\delta} \Vert \partial_s \eqref{lemstructhquadac-2-doubleresx} \Vert_{L^2}^2 ~ ds
&\lesssim \sum_{\gamma = a, c} \int_0^t s^{1-\delta} \langle s \rangle^{2\delta} \Vert h_{b \gamma}(s) \Vert_{L^2}^2 \Vert \partial_x u(s) \Vert_{L^{\infty}}^2 ~ ds \\
&\lesssim \int_0^t s^{-\frac{5}{6}} \langle s \rangle^{-\frac{1}{3}+201\delta} \Vert u \Vert_X^4 ~ ds ~ \lesssim \Vert u \Vert_X^4
\end{aligned} \\
&\begin{aligned}
&\int_0^t s^{1-\delta} \langle s \rangle^{2\delta} \Vert \partial_s (\eqref{lemstructhquadac-2-resxresx} + \eqref{lemstructhquadac-2-resxresxb}) \Vert_{L^2}^2 ~ ds \\
&\quad \lesssim \sum_{\gamma_1 = a, b, c} \sum_{\gamma_2 = a, c} \int_0^t s^{1-\delta} \langle s \rangle^{2\delta} \Vert e^{-i s \omega(D)} m_{0, \overline{\gamma_1}}(D) \nabla h_{\gamma_1}(s) \Vert_{L^4}^2 \Vert e^{-i s \omega(D)} \nabla h_{\gamma_2}(s) \Vert_{L^4}^2 ~ ds \\
&\quad \lesssim \sum_{\gamma_1 = a, b, c} \left( \int_0^t s^{1-\delta} \langle s \rangle^{2\delta} \Vert e^{-i s \omega(D)} m_{0, \overline{\gamma_1}}(D) \nabla h_{\gamma_1}(s) \Vert_{L^4}^4 ~ ds \right)^{\frac{1}{2}} \\
&\quad \quad \quad \times \sum_{\gamma_2 = a, c} \left( \int_0^t s^{1-\delta} \langle s \rangle^{2\delta} \Vert e^{-i s \omega(D)} \nabla h_{\gamma_2}(s) \Vert_{L^4}^4 ~ ds \right)^{\frac{1}{2}} \\
&\quad \lesssim \Vert u \Vert_X^4 
\end{aligned} \\
&\begin{aligned} 
&\int_0^t s^{1-\delta} \langle s \rangle^{2\delta} \Vert \partial_s \eqref{lemstructhquadac-2-resxresxvarphi} \Vert_{L^2}^2 ~ ds \\
&\quad \lesssim \sum_{\gamma_1, \gamma_2 = a, b, c} \int_0^t s^{1-\delta} \langle s \rangle^{2\delta} |\widetilde{\chi}(s)| \left( \Vert e^{-i s \omega(D)} |\nabla|^{-1} h_{\gamma_1}(s) \Vert_{L^4} + \Vert e^{-i s \omega(D)} \nabla h_{\gamma_1}(s) \Vert_{L^4} \right)^2 \\
&\pushright{\Vert e^{-i s \omega(D)} \nabla h_{\gamma_2}(s) \Vert_{L^4}^2 ~ ds} \\
&\quad \lesssim \int_0^t s^{-\frac{5}{6}} |\widetilde{\chi}(s)| \Vert u \Vert_X^4 ~ ds ~ \lesssim \Vert u \Vert_X^4 
\end{aligned} \\
&\begin{aligned} 
&\int_0^t s^{1-\delta} \langle s \rangle^{2\delta} \Vert \partial_s (\eqref{lemstructhquadac-2-dersymb-resx} + \eqref{lemstructhquadac-2-dersymb-resxb}) \Vert_{L^2}^2 ~ ds \\
&\quad \lesssim \sum_{\gamma = a, b, c} \int_0^t s^{1-\delta} \langle s \rangle^{2\delta} \Vert e^{-i s \omega(D)} m_{0, \overline{\gamma_1}}(D) \nabla h_{\gamma_1}(s) \Vert_{L^4}^2 \Vert u(s) \Vert_{L^4}^2 ~ ds \\
&\quad \lesssim \int_0^t s^{-\frac{5}{6}} \langle s \rangle^{-\frac{1}{4}+201\delta} \Vert u \Vert_X^4 ~ ds ~ \lesssim \Vert u \Vert_X^4 
\end{aligned} \\
&\begin{aligned} 
\int_0^t s^{1-\delta} \langle s \rangle^{2\delta} \Vert \partial_s \eqref{lemstructhquadac-2-dersymb-resxbis} \Vert_{L^2}^2 ~ ds &\lesssim \sum_{\gamma = a, b, c} \int_0^t s^{1-\delta} \langle s \rangle^{2\delta} \Vert h_{\gamma}(s) \Vert_{L^2}^2 \Vert \partial_x u(s) \Vert_{L^{\infty}}^2 ~ ds \\
&\lesssim \int_0^t s^{-\frac{5}{6}} \langle s \rangle^{-\frac{1}{3}+201\delta} \Vert u \Vert_X^4 ~ ds ~ \lesssim \Vert u \Vert_X^4 
\end{aligned} \\
&\begin{aligned} 
\int_0^t s^{1-\delta} \langle s \rangle^{2\delta} \Vert \partial_s \eqref{lemstructhquadac-2-doubledersymb} \Vert_{L^2}^2 ~ ds &\lesssim \int_0^t s^{1-\delta} \langle s \rangle^{2\delta} \Vert u(s) \Vert_{L^4}^4 ~ ds \\
&\lesssim \int_0^t s^{-\frac{5}{6}} \langle s \rangle^{-\frac{1}{3}+201\delta} \Vert u \Vert_X^4 ~ ds ~ \lesssim \Vert u \Vert_X^4 
\end{aligned} \\
&\begin{aligned} 
&\int_0^t s^{1-\delta} \langle s \rangle^{2\delta} \Vert \partial_s \eqref{lemstructhquadac-3-resx} \Vert_{L^2}^2 ~ ds
\lesssim \sum_{\gamma = a, c} \int_0^t s^{3-\delta} \langle s \rangle^{2\delta} \Vert e^{-i s \omega(D)} \nabla h_{\gamma}(s) \Vert_{L^4}^2 \Vert \nabla u(s) \Vert_{L^{\infty}}^2 \Vert u(s) \Vert_{L^4}^2 ~ ds \\
&\quad \lesssim \sum_{\gamma = a, c} \left( \int_0^t s^{1-\delta} \langle s \rangle^{2\delta} \Vert e^{-i s \omega(D)} \nabla h_{\gamma}(s) \Vert_{L^4}^4 ~ ds \right)^{\frac{1}{2}} \left( \int_0^t s^{5-\delta} \langle s \rangle^{2\delta} \Vert \nabla u(s) \Vert_{L^{\infty}}^4 \Vert u(s) \Vert_{L^4}^4 ~ ds \right)^{\frac{1}{2}} \\
&\quad \lesssim \Vert u \Vert_X^2 \left( \int_0^t s^{-\frac{5}{6}} \langle s \rangle^{-\frac{1}{4}+601\delta} \Vert u \Vert_X^8 ~ ds \right)^{\frac{1}{2}} \lesssim \Vert u \Vert_X^6 
\end{aligned} \\
&\begin{aligned} 
\int_0^t s^{1-\delta} \langle s \rangle^{2\delta} \Vert \partial_s \eqref{lemstructhquadac-3-dersymb} \Vert_{L^2}^2 ~ ds
&\lesssim \int_0^t s^{3-\delta} \langle s \rangle^{2\delta} \Vert \nabla u(s) \Vert_{L^{\infty}}^2 \Vert u(s) \Vert_{L^4}^4 ~ ds \\
&\lesssim \int_0^t s^{-\frac{5}{6}} \langle s \rangle^{-\frac{1}{4}+401\delta} \Vert u \Vert_X^6 ~ ds ~ \lesssim \Vert u \Vert_X^6 
\end{aligned} \\
&\begin{aligned} 
\int_0^t s^{\frac{5}{6}} \langle s \rangle^{\frac{1}{6}+\delta} \Vert \partial_s \eqref{lemstructhquadac-3-varphibis} \Vert_{L^2}^2 ~ ds
&\lesssim \sum_{\gamma_1, \gamma_2 = a, b, c} \int_0^t s^{\frac{5}{6}} \langle s \rangle^{\frac{1}{6}+\delta} \Vert e^{-i s \omega(D)} |\nabla|^{-1} h_{\gamma_1}(s) \Vert_{L^{\infty}}^2 \\
&\pushright{\Vert e^{-i s \omega(D)} m_{0, \overline{\gamma_2}}(D) \nabla h_{\gamma_2}(s) \Vert_{L^4}^2 \Vert u(s) \Vert_{L^4}^2 ~ ds} \\
&\lesssim \sum_{\gamma_1 = a, b, c} \int_0^t s^{-\frac{5}{6}} \langle s \rangle^{-\frac{1}{4}+201\delta} \Vert h_{\gamma_1}(s) \Vert_{H^1}^2 \Vert u \Vert_X^4 ~ ds ~ \lesssim \Vert u \Vert_X^6 
\end{aligned} \\
&\begin{aligned} 
\int_0^t s^{1-\delta} \langle s \rangle^{2\delta} \Vert \partial_s \eqref{lemstructhquadac-3-varphidersymb} \Vert_{L^2}^2 ~ ds
&\lesssim \sum_{\gamma = a, b, c} \int_0^t s^{1-\delta} \langle s \rangle^{2\delta} \Vert e^{-i s \omega(D)} |\nabla|^{-1} h_{\gamma}(s) \Vert_{L^{\infty}}^2 \Vert u(s) \Vert_{L^4}^4 ~ ds \\
&\lesssim \sum_{\gamma = a, b, c} \int_0^t s^{-\frac{5}{6}} \langle s \rangle^{-\frac{1}{3}+201\delta} \Vert h_{\gamma}(s) \Vert_{H^1}^2 \Vert u \Vert_X^4 ~ ds ~ \lesssim \Vert u \Vert_X^6 
\end{aligned} \\
&\begin{aligned} 
\int_0^t s^{1-\delta} \langle s \rangle^{2\delta} \Vert \partial_s \eqref{lemstructhquadac-4-dersymb} \Vert_{L^2}^2 ~ ds
&\lesssim \sum_{\gamma = a, b, c} \int_0^t s^{3-\delta} \langle s \rangle^{2\delta} \Vert e^{- i s \omega(D)} |\nabla|^{-1} h_{\gamma}(s) \Vert_{L^{\infty}}^2 \Vert \nabla u(s) \Vert_{L^{\infty}}^2 \Vert u(s) \Vert_{L^4}^4 ~ ds \\
&\lesssim \int_0^t s^{-\frac{5}{6}} \langle s \rangle^{-\frac{1}{4}+401\delta} \Vert u \Vert_X^8 ~ ds ~ \lesssim \Vert u \Vert_X^8 
\end{aligned} \\
&\begin{aligned} 
\int_0^t s^{1-\delta} \langle s \rangle^{2\delta} \Vert \partial_s \eqref{lemstructhquadacsurb-4-s2} \Vert_{L^2}^2 ~ ds 
&\lesssim \int_0^t s^{5-\delta} \langle s \rangle^{2\delta} \Vert \nabla u(s) \Vert_{L^{\infty}}^4 \Vert u(s) \Vert_{L^4}^4 ~ ds \\
&\lesssim \int_0^t s^{-\delta} \langle s \rangle^{-\frac{7}{6}+601\delta} \Vert u \Vert_X^8 ~ ds ~ \lesssim \Vert u \Vert_X^8 
\end{aligned} 
\end{align*}
For \eqref{lemstructhquadac-3-varphi}, we will use that: 
\begin{align*}
\Vert e^{-i t \omega(D)} h_{\gamma_2}(t) \Vert_{L^3} &\lesssim \Vert e^{-i t \omega(D)} |\nabla|^{\frac{1}{6}} |\partial_x|^{\frac{1}{12}} h_{\gamma_2}(t) \Vert_{L^{\frac{12}{5}}} \\
&\lesssim \Vert e^{-i t \omega(D)} \partial_x h_{\gamma_2}(t) \Vert_{L^{\frac{12}{5}}}^{\frac{1}{12}} \Vert e^{-i t \omega(D)} |\nabla|^{\frac{2}{11}} h_{\gamma_2}(t) \Vert_{L^{\frac{12}{5}}}^{\frac{11}{12}} \\
&\lesssim \Vert e^{-i t \omega(D)} \partial_x h_{\gamma_2}(t) \Vert_{L^4}^{\frac{1}{36}} \Vert h_{\gamma_2}(t) \Vert_{H^1}^{\frac{35}{36}} 
\end{align*}
Therefore: 
\begin{align*}
&\int_0^t s^{1-\delta} \langle s \rangle^{2\delta} \Vert \partial_s \eqref{lemstructhquadac-3-varphi} \Vert_{L^2}^2 ~ ds 
\\
&\lesssim \sum_{\gamma_1, \gamma_2 = a, b, c} \int_0^t s^{1-\delta} \langle s \rangle^{2\delta} \Vert e^{-i s \omega(D)} |\nabla|^{-1} h_{\gamma_1}(s) \Vert_{L^6}^2 \Vert e^{-i s \omega(D)} h_{\gamma_2}(s) \Vert_{L^3}^2 \Vert \nabla u(s) \Vert_{L^{\infty}}^2 ~ ds \\
&\lesssim \sum_{\gamma_1, \gamma_2 = a, b, c} \int_0^t s^{1-\delta} \langle s \rangle^{2\delta} \Vert h_{\gamma_1}(s) \Vert_{L^2}^2 \Vert h_{\gamma_2}(s) \Vert_{H^1}^{\frac{35}{18}} \Vert e^{-i s \omega(D)} \partial_x h_{\gamma_2}(s) \Vert_{L^4}^{\frac{1}{18}} \Vert \nabla u(s) \Vert_{L^{\infty}}^2 ~ ds \\
&\lesssim \sum_{\gamma_2 = a, b, c} \int_0^t s^{-\frac{5}{6}} \langle s \rangle^{-\frac{1}{6}+101\delta} \Vert u \Vert_X^{4+\frac{35}{18}} \Vert e^{-i s \omega(D)} \partial_x h_{\gamma_2}(s) \Vert_{L^4}^{\frac{1}{18}} ~ ds \\
&\lesssim \sum_{\gamma_2 = a, b, c} \left( \int_0^t s^{\frac{511}{512}} \langle s \rangle^{\frac{1}{256}} \Vert e^{-i s \omega(D)} \partial_x h_{\gamma_2}(s) \Vert_{L^4}^4 ~ ds \right)^{\frac{1}{72}} \\
&\quad \quad \quad \left( \int_0^t s^{-\frac{5}{6}*\frac{72}{71}-\frac{511}{512}*\frac{1}{72}} \langle s \rangle^{\left(-\frac{1}{6}+101\delta\right)*\frac{72}{71}-\frac{1}{256}*\frac{1}{72}} ~ ds \right)^{\frac{71}{72}} \Vert u \Vert_X^{6-\frac{1}{18}} \\
&\lesssim \Vert u \Vert_X^6 
\end{align*}
where $\delta$ is small enough. 

These estimates are enough in order to control the corresponding contributions: indeed, if $A(s)$ satisfies 
\begin{align}
\int_0^t s^{1-\delta} \langle s \rangle^{2\delta} \Vert A(s) \Vert_{L^2}^2 ~ ds \lesssim \Vert u \Vert_X^4 \label{estimeegeneriquederiveehquadactermessimples} 
\end{align}
which is the case for every temporal derivative above if $\Vert u \Vert_X \leq 1$, then 
\begin{align*}
&\int_0^t \int h_{\alpha \beta}(s, x, y) ~ A(s, x, y) ~ dx dy ds \\
&\lesssim \left( \int_0^t s^{-1+\delta} \langle s \rangle^{-2\delta} \Vert h_{\alpha \beta}(s) \Vert_{L^2}^2 ~ ds \right)^{\frac{1}{2}} \left( \int_0^t s^{1-\delta} \langle s \rangle^{2\delta} \Vert A(s) \Vert_{L^2}^2 ~ ds \right)^{\frac{1}{2}} \\
&\lesssim \left( \int_0^t s^{-1+\delta} \langle s \rangle^{-2\delta} \Vert u \Vert_X^2 ~ ds \right)^{\frac{1}{2}} \Vert u \Vert_X^2 \\
&\lesssim \Vert u \Vert_X^3 
\end{align*}
as wanted. 

\paragraph{2.2. Symmetric term} For \eqref{lemstructhquadac-2-sym} or \eqref{lemstructhquadacsurb-2-sym},  we separate the cases depending on $\alpha, \beta$. 

If $\alpha = b$, $\beta \in \{ a, c \}$, we have 
\begin{align*}
\eqref{lemstructhquadac-2-sym} &= 2 \int_0^t \int e^{i s \varphi} i \eta_0 \widehat{h}_{b \beta}(s, \overline{\eta}) \widehat{f}(s, \overline{\sigma}) ~ d\overline{\eta} ds 
\end{align*}
Therefore, we need to estimate: 
\begin{align*}
&\int_0^t \int h_{b \beta}(s, x, y) ~ \partial_s \mathcal{F}^{-1} \eqref{lemstructhquadac-2-sym}(s, x, y) ~ dx dy ds \\
&\quad = 2 \int_0^t \int \int \widehat{h}_{b \beta}(s, -\overline{\xi}) e^{i s \varphi} i \eta_0 \widehat{h}_{b \beta}(s, \overline{\eta}) \widehat{f}(s, \overline{\sigma}) d\overline{\eta} d\overline{\xi} ds \\
&\quad = - 2 \int_0^t \int \int \widehat{h}_{b \beta}(s, \overline{\eta}) e^{i s \varphi} i \xi_0 \widehat{h}_{b \beta}(s, -\overline{\xi}) \widehat{f}(s, \overline{\sigma}) d\overline{\eta} d\overline{\xi} ds \\
&\quad = - \int_0^t \int \int \widehat{h}_{b \beta}(s, \overline{\eta}) e^{i s \varphi} i \sigma_0 \widehat{h}_{b \beta}(s, -\overline{\xi}) \widehat{f}(s, \overline{\sigma}) d\overline{\eta} d\overline{\xi} ds \\
&\quad \lesssim \int_0^t \Vert h_{b \beta}(s) \Vert_{L^2}^2 \Vert \partial_x u(s) \Vert_{L^{\infty}} ~ ds \\
&\quad \lesssim \int_0^t s^{-\frac{5}{6}} \langle s \rangle^{-\frac{1}{4}+100\delta} \Vert u \Vert_X^3 ~ ds ~ \lesssim \Vert u \Vert_X^3
\end{align*}
using the symmetry $\mathfrak{s}: (\overline{\xi}, \overline{\eta}) \mapsto (-\overline{\eta}, -\overline{\xi})$. 

Likewise, if $\alpha \in \{ a, c \}$ and $\beta = b$, we need to consider: 
\begin{align*}
\eqref{lemstructhquadacsurb-2-sym} &= 2 \int_0^t \int e^{i s \varphi} i \eta_0 \widehat{h}_{\alpha b}(s, \overline{\eta}) \widehat{f}(s, \overline{\sigma}) ~ d\overline{\eta} ds
\end{align*}
that has exactly the same structure. 

Finally, if $\alpha, \beta \in \{ a, c \}$, we separate depending on the size of the frequencies: 
\begin{subequations}
\begin{align}
\eqref{lemstructhquadac-2-sym} &= 2 \frac{i \xi_0}{|\overline{\xi}|} \int_0^t \int e^{i s \varphi} |\overline{\eta}| \widehat{h}_{\alpha \beta}(s, \overline{\eta}) \widehat{f}(s, \overline{\sigma}) ~ d\overline{\eta} ds \notag \\
&= 2 \frac{i \xi_0}{|\overline{\xi}|} \int_0^t \int e^{i s \varphi} \left( \mu_{BBB}(\overline{\xi}, \overline{\eta}) + \mu_{HHB}(\overline{\xi}, \overline{\eta}) \right) |\overline{\eta}| \widehat{h}_{\alpha \beta}(s, \overline{\eta}) \widehat{f}(s, \overline{\sigma}) ~ d\overline{\eta} ds \label{estimeequadhac-symHHB} \\
&\quad + 2 \frac{i \xi_0}{|\overline{\xi}|} \int_0^t \int e^{i s \varphi} \left( \mu_{BHH}(\overline{\xi}, \overline{\eta}) + \mu_{HBH}(\overline{\xi}, \overline{\eta}) \right) |\overline{\eta}| \widehat{h}_{\alpha \beta}(s, \overline{\eta}) \widehat{f}(s, \overline{\sigma}) ~ d\overline{\eta} ds \label{estimeequadhac-symBHH} 
\end{align}
\end{subequations}

On \eqref{estimeequadhac-symHHB}, we apply the same symmetry $\mathfrak{s}$ as before: 
\begin{align*}
&\int_0^t \int h_{\alpha \beta}(s, x, y) ~ \partial_s \mathcal{F}^{-1} \eqref{estimeequadhac-symHHB}(s, x, y) ~ dx dy ds \\
&\quad = 2 \int_0^t \int \int \widehat{h}_{\alpha \beta}(s, -\overline{\xi}) e^{i s \varphi} \frac{i \xi_0 |\overline{\eta}|}{|\overline{\xi}|} \left( \mu_{BBB}(\overline{\xi}, \overline{\eta}) + \mu_{HHB}(\overline{\xi}, \overline{\eta}) \right) \widehat{h}_{\alpha \beta}(s, \overline{\eta}) \widehat{f}(s, \overline{\sigma}) ~ d\overline{\eta} d\overline{\xi} ds \\
&\quad = \int_0^t \int \int \widehat{h}_{\alpha \beta}(s, -\overline{\xi}) e^{i s \varphi} \left( \frac{i \xi_0 |\overline{\eta}|}{|\overline{\xi}|} - \frac{i \eta_0 |\overline{\xi}|}{|\overline{\eta}|} \right) \left( \mu_{BBB}(\overline{\xi}, \overline{\eta}) + \mu_{HHB}(\overline{\xi}, \overline{\eta}) \right) \widehat{h}_{\alpha \beta}(s, \overline{\eta}) \widehat{f}(s, \overline{\sigma}) ~ d\overline{\eta} d\overline{\xi} ds
\end{align*}
up to choosing $\mu_{BBB}, \mu_{HHB}$ invariant by $\mathfrak{s}$. But
\begin{align*}
\varphi &= O(\sigma_0) + \xi_0 (|\xi|^2 - |\eta|^2) = O(\sigma_0) + \eta_0 |\overline{\xi}|^2 - \xi_0 |\overline{\eta}|^2 
\end{align*}
so that
\begin{align*}
\frac{i \xi_0 |\overline{\eta}|}{|\overline{\xi}|} - \frac{i \eta_0 |\overline{\xi}|}{|\overline{\eta}|} &= O(\varphi) + O(\sigma_0) 
\end{align*}
The term in $O(\sigma_0)$ can be estimated as before, while the term in $O(\varphi)$ contributes as \eqref{lemstructhquadac-2-doubleresxvarphi} and its estimate will be done below. 

On \eqref{estimeequadhac-symBHH}, we develop more precisely: 
\begin{align*}
\varphi &= \xi_0^3 + \xi_0 |\xi|^2 - \eta_0^3 - \eta_0 |\eta|^2 - \sigma_0^3 - \sigma_0 |\sigma|^2 \\
&= 3 \xi_0 \eta_0 \sigma_0 + \xi_0 |\xi|^2 - \xi_0 |\eta|^2 + \sigma_0 (|\eta|^2 - |\sigma|^2) \\
&= O(\sigma_0 \overline{\xi}) + \xi_0 (|\xi|^2 - |\eta|^2) 
\end{align*}
In particular, on the support of $\mu_{BHH} + \mu_{HBH}$, $1 = O(|\xi|^2 - |\eta|^2)$, so that 
\begin{align*}
\frac{\xi_0 |\overline{\eta}|}{|\overline{\xi}|} &= O(\sigma_0) + O\left( \frac{\varphi}{|\overline{\xi}|} \right) 
\end{align*}
Again, the term in $O(\sigma_0)$ can be estimated as before, while the term in $O\left( \frac{\varphi}{|\overline{\xi}|} \right)$ contributes as 
\begin{align}
\int_0^t \int e^{i s \varphi} |\overline{\xi}|^{-1} \mu(\overline{\xi}, \overline{\eta}) \frac{\varphi}{|\overline{\eta}|+|\overline{\sigma}|} \widehat{h}_{\alpha \beta}(s, \overline{\eta}) \widehat{f}(s, \overline{\sigma}) ~ d\overline{\eta} ds \label{estimeequadhac-varphibis} 
\end{align}
for some symbol $\mu$, and we estimate it below. 

\paragraph{2.3. Term with $\varphi$} It only remains \eqref{lemstructhquadac-2-doubleresxvarphi} and \eqref{estimeequadhac-varphibis}: but since the former has a simpler form than the latter, we only treat the second one. In Duhamel's formula, for $(\gamma_1, \gamma_2) \in \{ a, b, c \}^2 \setminus \{ (b, b) \}$: 
\begin{subequations}
\begin{align}
&\int_0^t \int \int |\overline{\xi}|^{-1} \widehat{h}_{\alpha \beta}(s, -\overline{\xi}) e^{i s \varphi} \mu(\overline{\xi}, \overline{\eta}) \frac{\varphi}{|\overline{\eta}|+|\overline{\sigma}|} \widehat{h}_{\gamma_1, \gamma_2}(s, \overline{\eta}) \widehat{f}(s, \overline{\sigma}) ~ d\overline{\eta} ds \notag \\
&\quad = \int_0^t \int \int |\overline{\xi}|^{-1} \partial_s \widehat{h}_{\alpha \beta}(s, -\overline{\xi}) e^{i s \varphi} \mu(\overline{\xi}, \overline{\eta}) \frac{1}{|\overline{\eta}|+|\overline{\sigma}|} \widehat{h}_{\gamma_1 \gamma_2}(s, \overline{\eta}) \widehat{f}(s, \overline{\sigma}) ~ d\overline{\eta} ds \label{estimeequadhac-varphi-1} \\
&\quad \quad + \int_0^t \int \int |\overline{\xi}|^{-1} \widehat{h}_{\alpha \beta}(s, -\overline{\xi}) e^{i s \varphi} \mu(\overline{\xi}, \overline{\eta}) \frac{1}{|\overline{\eta}|+|\overline{\sigma}|} \partial_s \widehat{h}_{\gamma_1 \gamma_2}(s, \overline{\eta}) \widehat{f}(s, \overline{\sigma}) ~ d\overline{\eta} ds \label{estimeequadhac-varphi-2} \\
&\quad \quad + \int_0^t \int \int |\overline{\xi}|^{-1} \widehat{h}_{\alpha \beta}(s, -\overline{\xi}) e^{i s \varphi} \mu(\overline{\xi}, \overline{\eta}) \frac{1}{|\overline{\eta}|+|\overline{\sigma}|} \widehat{h}_{\gamma_1 \gamma_2}(s, \overline{\eta}) \partial_s \widehat{f}(s, \overline{\sigma}) ~ d\overline{\eta} ds \label{estimeequadhac-varphi-3} \\
&\quad \quad + \int \int |\overline{\xi}|^{-1} \widehat{h}_{\alpha \beta}(t, -\overline{\xi}) e^{i t \varphi} \mu(\overline{\xi}, \overline{\eta}) \frac{1}{|\overline{\eta}|+|\overline{\sigma}|} \widehat{h}_{\gamma_1 \gamma_2}(t, \overline{\eta}) \widehat{f}(t, \overline{\sigma}) ~ d\overline{\eta} ds \label{estimeequadhac-varphi-4} \\
&\quad \quad + \int \int |\overline{\xi}|^{-1} \widehat{h}_{\alpha \beta}(0, -\overline{\xi}) e^{i s \varphi} \mu(\overline{\xi}, \overline{\eta}) \frac{\varphi}{|\overline{\eta}|+|\overline{\sigma}|} \widehat{h}_{\gamma_1 \gamma_2}(0, \overline{\eta}) \widehat{f}(0, \overline{\sigma}) ~ d\overline{\eta} ds \label{estimeequadhac-varphi-5}
\end{align}
\end{subequations} 
by integration by parts in time (the symbol $\mu$ may change from line to line). Already, we can estimate: 
\begin{align*}
\eqref{estimeequadhac-varphi-3} &\lesssim \int_0^t \Vert e^{-i s \omega(D)} |\nabla|^{-1} h_{\alpha \beta}(s) \Vert_{L^6} \Vert h_{\gamma_1 \gamma_2}(s) \Vert_{L^2} \Vert e^{-i s \omega(D)} |\nabla|^{-1} \partial_s f(s) \Vert_{L^3} ~ ds \\
&\lesssim \int_0^t \Vert u \Vert_X^2 \Vert u(s) \Vert_{L^6}^2 ~ ds \\
&\lesssim \int_0^t \langle s \rangle^{-\frac{13}{12}+100\delta} \Vert u \Vert_X^4 ~ ds ~ \lesssim \Vert u \Vert_X^4 \\
\eqref{estimeequadhac-varphi-4} &\lesssim \Vert e^{-i t \omega(D)} |\nabla|^{-1} h_{\alpha \beta}(t) \Vert_{L^6} \Vert h_{\gamma_1 \gamma_2}(t) \Vert_{L^2} \Vert |\nabla|^{-1} u(t) \Vert_{L^3} \\
&\lesssim \Vert u \Vert_X^3 
\end{align*}
and \eqref{estimeequadhac-varphi-5} can be estimated like \eqref{estimeequadhac-varphi-4}. 

For the terms containing a derivative $\partial_s h_{\alpha \beta}$ or $\partial_s h_{\gamma_1 \gamma_2}$, we can reuse the decomposition we just used: we saw that 
\begin{align*}
\partial_s h_{\alpha \beta}(s) &= A(s) + \partial_s \mathcal{F}^{-1} \eqref{lemstructhquadac-2-sym} + \partial_s \mathcal{F}^{-1} \eqref{lemstructhquadac-2-doubleresxvarphi} 
\end{align*}
for $A(s)$ satisfying \eqref{estimeegeneriquederiveehquadactermessimples}. We will use the weaker bound
\begin{align}
\int_0^t s^{1-\delta} \langle s \rangle^{-1000\delta} \Vert A(s) \Vert_{L^2}^2 ~ ds ~ \lesssim \Vert u \Vert_X^4 \label{estimeegeneriquederiveehquadactermessimples-faible} 
\end{align}
Then by developing, we can already estimate the terms containing $A(s)$: 
\begin{align*}
&\int_0^t \int \int |\overline{\xi}|^{-1} \widehat{A}(s, -\overline{\xi}) e^{i s \varphi} \mu(\overline{\xi}, \overline{\eta}) \frac{1}{|\overline{\eta}|+|\overline{\sigma}|} \widehat{h}_{\gamma_1 \gamma_2}(s, \overline{\eta}) \widehat{f}(s, \overline{\sigma}) ~ d\overline{\eta} ds \\
&\quad \lesssim \int_0^t \Vert e^{- i s \omega(D)} |\nabla|^{-1} A(s) \Vert_{L^6} \Vert h_{\gamma_1 \gamma_2}(s) \Vert_{L^2} \Vert |\nabla|^{-1} u(s) \Vert_{L^3} ~ ds \\
&\quad \lesssim \left( \int_0^t s^{1-\delta} \langle s \rangle^{-1000\delta} \Vert A(s) \Vert_{L^2}^2 ~ ds \right)^{\frac{1}{2}} \left( \int_0^t s^{-1+\delta} \langle s \rangle^{1000\delta} \Vert u \Vert_X^2 \Vert |\nabla|^{-1} u(s) \Vert_{L^2}^{\frac{10}{7}} \Vert u(s) \Vert_{L^4}^{\frac{4}{7}} ~ ds \right)^{\frac{1}{2}} \\
&\quad \lesssim \Vert u \Vert_X^2 \left( \int_0^t s^{-1+\delta} \langle s \rangle^{-\frac{2}{7}+1000\delta} \Vert u \Vert_X^4 ~ ds \right)^{\frac{1}{2}} \\
&\quad \lesssim \Vert u \Vert_X^4 \\
&\int_0^t \int \int |\overline{\xi}|^{-1} \widehat{h}_{\alpha \beta}(s, -\overline{\xi}) e^{i s \varphi} \mu(\overline{\xi}, \overline{\eta}) \frac{1}{|\overline{\eta}|+|\overline{\sigma}|} \widehat{A}(s, \overline{\eta}) \widehat{f}(s, \overline{\sigma}) ~ d\overline{\eta} ds \\
&\quad \lesssim \int_0^t \Vert e^{-i s \omega(D)} |\nabla|^{-1} h_{\alpha \beta}(s) \Vert_{L^6} \Vert A(s) \Vert_{L^2} \Vert |\nabla|^{-1} u(s) \Vert_{L^3} ~ ds \\
&\quad \lesssim \Vert u \Vert_X^4 
\end{align*}

For the term containing \eqref{lemstructhquadac-2-doubleresxvarphi}, we estimate: 
\begin{align*}
\Vert \partial_s \eqref{lemstructhquadac-2-doubleresxvarphi} \Vert_{L^2} &\lesssim \Vert h_{b \beta}(s) \Vert_{L^2} \Vert \nabla u(s) \Vert_{L^{\infty}} \\
&\lesssim s^{-\frac{5}{6}} \langle s \rangle^{-\frac{1}{6}+100\delta} \Vert u \Vert_X^2 
\end{align*}
Hence, $\partial_s \eqref{lemstructhquadac-2-doubleresxvarphi}$ satisfies \eqref{estimeegeneriquederiveehquadactermessimples-faible} and we can estimate its contribution like the one from $A(s)$. 

Finally, for the term \eqref{lemstructhquadac-2-sym}, we decompose: 
\begin{subequations}
\begin{align}
\eqref{lemstructhquadac-2-sym} &= \int_0^t \int e^{i s \varphi} \mu(\overline{\xi}, \overline{\eta}) \mu_{HHB}(\overline{\eta}, \overline{\sigma}) |\overline{\eta}| \widehat{h}_{\alpha \beta}(s, \overline{\eta}) \widehat{f}(s, \overline{\sigma}) ~ d\overline{\sigma} ds \label{estimeehquadhac-symdec2HHB} \\
&\quad + \int_0^t \int e^{i s \varphi} \mu(\overline{\xi}, \overline{\eta}) \left( 1 - \mu_{HHB}(\overline{\eta}, \overline{\sigma}) \right) |\overline{\eta}| \widehat{h}_{\alpha \beta}(s, \overline{\eta}) \widehat{f}(s, \overline{\sigma}) ~ d\overline{\sigma} ds \label{estimeehquadhac-symdec2reste}
\end{align}
\end{subequations}
for some symbol $\mu$. We then have
\begin{align*}
\Vert \partial_s \eqref{estimeehquadhac-symdec2reste} \Vert_{L^2} &\lesssim \Vert h_{\alpha \beta}(s) \Vert_{L^2} \Vert \nabla u(s) \Vert_{L^{\infty}} \\
&\lesssim s^{-\frac{5}{6}} \langle s \rangle^{-\frac{1}{6}+100\delta} \Vert u \Vert_X^2 
\end{align*}
so we can estimate the contribution of $\partial_s \eqref{estimeehquadhac-symdec2reste}$ like the one from $A(s)$. Finally, 
\begin{align*}
\Vert e^{-i s \omega(D)} |\nabla|^{-1} \partial_s \mathcal{F}^{-1} \eqref{estimeehquadhac-symdec2HHB} \Vert_{L^{\frac{4}{3}}} &\lesssim \Vert h_{\alpha \beta}(s) \Vert_{L^2} \Vert u(s) \Vert_{L^4} \\
&\lesssim \langle s \rangle^{-\frac{13}{24}+50\delta} \Vert u \Vert_X^2 
\end{align*}
and so we can estimate in a symmetric way the two terms: 
\begin{align*}
&\int_0^t \int \int |\overline{\xi}|^{-1} \partial_s \eqref{estimeehquadhac-symdec2HHB}(s, -\overline{\xi}) e^{i s \varphi} \mu(\overline{\xi}, \overline{\eta}) \frac{1}{|\overline{\eta}|+|\overline{\sigma}|} \widehat{h}_{\gamma_1 \gamma_2}(s, \overline{\eta}) \widehat{f}(s, \overline{\sigma}) ~ d\overline{\eta} ds \\
&\quad + \int_0^t \int \int |\overline{\xi}|^{-1} \widehat{h}_{\alpha \beta}(s, -\overline{\xi}) e^{i s \varphi} \mu(\overline{\xi}, \overline{\eta}) \frac{1}{|\overline{\eta}|+|\overline{\sigma}|} \partial_s \eqref{estimeehquadhac-symdec2HHB}(s, \overline{\eta}) \widehat{f}(s, \overline{\sigma}) ~ d\overline{\eta} ds \\
&\lesssim \int_0^t \Vert e^{-i s \omega(D)} |\nabla|^{-1} \partial_s \mathcal{F}^{-1} \eqref{estimeehquadhac-symdec2HHB} \Vert_{L^{\frac{4}{3}}} \Vert e^{-i s \omega(D)} |\nabla|^{-1} h_{\gamma_1 \gamma_2}(s) \Vert_{L^6} \Vert u(s) \Vert_{L^{12}} ~ ds \\
&\lesssim \int_0^t \langle s \rangle^{-\frac{13}{12}+100\delta} \Vert u \Vert_X^4 ~ ds ~ \lesssim \Vert u \Vert_X^4 
\end{align*}
as wanted. This concludes. 
\end{Dem}

\begin{Lem} Let $t > 0$, $\alpha \in \{ a, b, c \}$. Then
\begin{align*}
\Vert \nabla X_c g_{\alpha}(t) \Vert_{L^2} &\lesssim \Vert u \Vert_X^2  
\end{align*}
Furthermore, for $t \geq 1$, 
\begin{align*}
\Vert \psi_j(D) m_{\widehat{\mathcal{C}}}(D) X_c g_{\alpha}(t) \Vert_{L^2} &\lesssim \left( t^{-\frac{1}{12}+\frac{1}{32}+50\delta} + t^{-\frac{1}{8}} 2^{-\frac{7j}{6}+1000\delta j} \right) \Vert u \Vert_X^2
\end{align*}
\end{Lem}

\begin{Dem}
We apply Lemma \ref{lemstructurehalphagalpha} to get that, for any $\alpha$, $g_{\alpha}$ has the structure \ref{lemstructurehb-14-g} for some symbol depending on $\alpha$. 

Now we may run the proof of the decomposition lemma for $X_c$ again, noting that for $X_c$ we do not use the fact that there is a time integral, and so we deduce that 
\begin{subequations}
\begin{align}
|\overline{\xi}| \widehat{X}_c(\overline{\xi}) \cdot \nabla_{\overline{\xi}} \widehat{g}_{\alpha}(t, \overline{\xi}) 
&= t \int e^{i t \varphi} \mu(\overline{\xi}, \overline{\eta}) \widehat{f}(t, \overline{\eta}) \widehat{f}(t, \overline{\sigma}) ~ d\overline{\eta} \label{estimeeaprioriXcg-1} \\
&\quad + t \int e^{i t \varphi} \mu(\overline{\xi}, \overline{\eta}) |\overline{\eta}| \widehat{X}_c(\overline{\eta}) \cdot \nabla_{\overline{\eta}} \widehat{f}(t, \overline{\eta}) \widehat{f}(t, \overline{\sigma}) ~ d\overline{\eta} \label{estimeeaprioriXcg-2} 
\end{align}
\end{subequations} 
for some symbol $\mu$ of order $0$ (that can change from line to line). But \eqref{estimeeaprioriXcg-1} is similar to \eqref{lemstructhquadac-gdersymb} and \eqref{estimeeaprioriXcg-2} is similar to \eqref{lemstructhquadac-gresx}, so the same estimates apply. This proves the first estimate of the lemma. 

In order to improve the estimate in the presence of an additional derivative, we reuse the previous disjonction. First, we have by the same computation as before that 
\begin{align*}
\Vert \psi_j(\overline{\xi}) m_{\widehat{\mathcal{C}}}(\overline{\xi}) \eqref{estimeeaprioriXcg-1} \Vert_{L^2} &\lesssim t 2^{-\frac{j}{2}} \Vert u(t) \Vert_{L^2} \Vert |\nabla|^{\frac{1}{2}} \partial_x u(t) \Vert_{L^{\infty}} \\
&\lesssim t^{-\frac{1}{6}+100\delta} 2^{-\frac{j}{2}} \Vert u \Vert_X^2 \\
&\lesssim t^{-\frac{1}{12}+\frac{1}{32}+50\delta} \Vert u \Vert_X^2 + t^{-\frac{1}{8}} 2^{-\frac{7j}{6}+1000\delta j} \Vert u \Vert_X^2 
\end{align*}
Then, for \eqref{estimeeaprioriXcg-2}, we separate again: 
\begin{subequations}
\begin{align}
&\psi_j(\overline{\xi}) m_{\widehat{\mathcal{C}}}(\overline{\xi}) \eqref{estimeeaprioriXcg-2} \notag \\
&= \psi_j(\overline{\xi}) m_{\widehat{\mathcal{C}}}(\overline{\xi}) t \int e^{i t \varphi} \mu(\overline{\xi}, \overline{\eta}) \mu_{HHB}(\overline{\xi}, \overline{\eta}) |\overline{\eta}| \widehat{X}_c(\overline{\eta}) \cdot \nabla_{\overline{\eta}} \widehat{f}(t, \overline{\eta}) \widehat{f}(t, \overline{\sigma}) ~ d\overline{\eta} \label{estimeeaprioriXcg-2-1} \\
&\quad + \psi_j(\overline{\xi}) m_{\widehat{\mathcal{C}}}(\overline{\xi}) t \int e^{i t \varphi} \mu(\overline{\xi}, \overline{\eta}) \left( 1 - \mu_{HHB}(\overline{\xi}, \overline{\eta}) \right) \left( 1 - m_{\widehat{\mathcal{P}}}(\overline{\sigma}) \right) |\overline{\eta}| \widehat{X}_c(\overline{\eta}) \cdot \nabla_{\overline{\eta}} \widehat{f}(t, \overline{\eta}) \widehat{f}(t, \overline{\sigma}) ~ d\overline{\eta} \label{estimeeaprioriXcg-2-2} \\
&\quad + \psi_j(\overline{\xi}) m_{\widehat{\mathcal{C}}}(\overline{\xi}) t \int e^{i t \varphi} \mu(\overline{\xi}, \overline{\eta}) \left( 1 - \mu_{HHB}(\overline{\xi}, \overline{\eta}) \right) \left( 1 - m_{\widehat{\mathcal{C}}}(\overline{\eta}) \right) m_{\widehat{\mathcal{P}}}(\overline{\sigma}) |\overline{\eta}| \widehat{X}_c(\overline{\eta}) \cdot \nabla_{\overline{\eta}} \widehat{f}(t, \overline{\eta}) \widehat{f}(t, \overline{\sigma}) ~ d\overline{\eta} \label{estimeeaprioriXcg-2-3} \\
&\quad + \psi_j(\overline{\xi}) m_{\widehat{\mathcal{C}}}(\overline{\xi}) t \int e^{i t \varphi} \mu(\overline{\xi}, \overline{\eta}) \left( 1 - \mu_{HHB}(\overline{\xi}, \overline{\eta}) \right) m_{\widehat{\mathcal{C}}}(\overline{\eta}) m_{\widehat{\mathcal{P}}}(\overline{\sigma}) |\overline{\eta}| \widehat{X}_c(\overline{\eta}) \cdot \nabla_{\overline{\eta}} \widehat{f}(t, \overline{\eta}) \widehat{f}(t, \overline{\sigma}) ~ d\overline{\eta} \label{estimeeaprioriXcg-2-4}
\end{align}
\end{subequations}
We then estimate: 
\begin{align*}
\Vert \eqref{estimeeaprioriXcg-2-2} \Vert_{L^2} &\lesssim t 2^{-\frac{j}{2}} \Vert h_c(t) \Vert_{L^2} \Vert |\nabla|^{\frac{1}{2}} \partial_x u(t) \Vert_{L^{\infty}} \\
&\lesssim t^{-\frac{1}{6}+100\delta} 2^{-\frac{j}{2}} \Vert u \Vert_X^2 \\
\Vert \eqref{estimeeaprioriXcg-2-3} \Vert_{L^2} &\lesssim t \Vert e^{-i t \omega(D)} (1 - m_{\widehat{\mathcal{C}}}(D)) \nabla h_c(t) \Vert_{L^4} \Vert |\nabla|^{\frac{1}{2}} u(t) \Vert_{L^4} \\
&\lesssim t^{-\frac{1}{12}+\frac{1}{128}+50\delta} \Vert u \Vert_X^2 \\
\Vert \eqref{estimeeaprioriXcg-2-4} \Vert_{L^2} &\lesssim t \Vert e^{-i t \omega(D)} m_{\widehat{\mathcal{C}}}(D) |\nabla|^{\frac{11}{16}} \langle \nabla \rangle^{-\frac{1}{16}} h_c(t) \Vert_{L^4} \Vert m_{\widehat{\mathcal{P}}}(D) |\nabla|^{\frac{5}{16}} \langle \nabla \rangle^{\frac{1}{16}} u(t) \Vert_{L^4} \\
&\lesssim t^{-\frac{1}{12}+\frac{1}{32}+50\delta} \Vert u \Vert_X^2 
\end{align*}
Then, for \eqref{estimeeaprioriXcg-2-1}, we estimate directly if $2^j \lesssim t^{-\frac{13}{171}}$: 
\begin{align*}
\Vert \eqref{estimeeaprioriXcg-2-1} \Vert_{L^2} &\lesssim t 2^{\frac{11j}{32}-1000\delta} \Vert e^{-i t \omega(D)} m_{\widehat{\mathcal{C}}}(D) |\nabla|^{\frac{21}{32}+1000\delta} h_c(t) \Vert_{L^4} \Vert u(t) \Vert_{L^4} \\
&\lesssim t^{-\frac{1}{96}+50\delta} 2^{\frac{11j}{32}-1000\delta} \Vert u \Vert_X^2 \\
&\lesssim t^{-\frac{257}{2052}+2000\delta} 2^{-\frac{7j}{6}+1000\delta} \Vert u \Vert_X^2 \\
&\lesssim t^{-\frac{1}{8}} 2^{-\frac{7j}{16}+1000\delta} \Vert u \Vert_X^2 
\end{align*}
On the other hand, if $2^j \gtrsim t^{-\frac{13}{171}}$, we apply an integration by parts as before: 
\begin{subequations}
\begin{align}
\eqref{estimeeaprioriXcg-2-1} &= \psi_j(\overline{\xi}) m_{\widehat{\mathcal{C}}}(\overline{\xi}) \int e^{i t \varphi} \mu(\overline{\xi}, \overline{\eta}) \mu_{HHB}(\overline{\xi}, \overline{\eta}) |\overline{\eta}|^{-1} |\overline{\sigma}|^{-1} \widehat{X}_c(\overline{\eta}) \cdot \nabla_{\overline{\eta}} \widehat{f}(t, \overline{\eta}) \widehat{f}(t, \overline{\sigma}) ~ d\overline{\eta} \label{estimeeaprioriXcg-2-1-1} \\
&\quad + \psi_j(\overline{\xi}) m_{\widehat{\mathcal{C}}}(\overline{\xi}) \int e^{i t \varphi} \mu(\overline{\xi}, \overline{\eta}) \mu_{HHB}(\overline{\xi}, \overline{\eta}) |\overline{\eta}|^{-1} \widehat{X}_a(\overline{\eta}) \cdot \nabla_{\overline{\eta}} \widehat{h}_c(t, \overline{\eta}) \widehat{f}(t, \overline{\sigma}) ~ d\overline{\eta} \label{estimeeaprioriXcg-2-1-2} \\
&\quad + \psi_j(\overline{\xi}) m_{\widehat{\mathcal{C}}}(\overline{\xi}) \int e^{i t \varphi} \mu(\overline{\xi}, \overline{\eta}) \mu_{HHB}(\overline{\xi}, \overline{\eta}) \left( 1 - m_{\widehat{\mathcal{C}}}(\overline{\eta}) \right) |\overline{\eta}|^{-1} \widehat{X}_c(\overline{\eta}) \cdot \nabla_{\overline{\eta}} \widehat{f}(t, \overline{\eta}) \nabla_{\overline{\eta}} \widehat{f}(t, \overline{\sigma}) ~ d\overline{\eta} \label{estimeeaprioriXcg-2-1-3} \\
&\quad + \psi_j(\overline{\xi}) m_{\widehat{\mathcal{C}}}(\overline{\xi}) \int e^{i t \varphi} \mu(\overline{\xi}, \overline{\eta}) \mu_{HHB}(\overline{\xi}, \overline{\eta}) m_{\widehat{\mathcal{C}}}(\overline{\eta}) |\overline{\eta}|^{-1} \widehat{X}_c(\overline{\eta}) \cdot \nabla_{\overline{\eta}} \widehat{f}(t, \overline{\eta}) \nabla_{\overline{\eta}} \widehat{f}(t, \overline{\sigma}) ~ d\overline{\eta} \label{estimeeaprioriXcg-2-1-4}
\end{align}
\end{subequations} 
We then estimate: 
\begin{align*}
\Vert \eqref{estimeeaprioriXcg-2-1-2} \Vert_{L^2} &\lesssim 2^{-\frac{5j}{4}} \Vert e^{-i t \omega(D)} |\nabla|^{\frac{1}{4}} X_a h_c(t) \Vert_{L^4} \Vert u(t) \Vert_{L^4} \\
&\lesssim t^{-\frac{1}{2}} 2^{-\frac{5j}{4}} \Vert u \Vert_X^2 \\
&\lesssim t^{-\frac{1}{8}} 2^{-\frac{7j}{6}+1000\delta j} \Vert u \Vert_X^2 \\
\Vert \eqref{estimeeaprioriXcg-2-1-3} \Vert_{L^2} &\lesssim 2^{-j} \Vert e^{-i t \omega(D)} (1 - m_{\widehat{\mathcal{C}}}(D)) \psi_j(D) h_c(t) \Vert_{L^3} \Vert e^{-i t \omega(D)} (x, y) f(t) \Vert_{L^6} \\
&\lesssim 2^{-\frac{5j}{3}} \Vert h_c(t) \Vert_{L^2}^{\frac{1}{3}} \Vert e^{-i t \omega(D)} (1 - m_{\widehat{\mathcal{C}}}(D)) \nabla h_c(t) \Vert_{L^4}^{\frac{2}{3}} \Vert \nabla (x, y) f(t) \Vert_{L^2} \\
&\lesssim t^{-\frac{47}{96}+101\delta} 2^{-\frac{5j}{3}} \Vert u \Vert_X^2 \\
&\lesssim t^{-\frac{1}{8}} 2^{-\frac{7j}{6}+1000\delta j} \Vert u \Vert_X^2 \\
\Vert \eqref{estimeeaprioriXcg-2-1-4} \Vert_{L^2} &\lesssim 2^{-j} \Vert e^{-i t \omega(D)} m_{\widehat{\mathcal{C}}}(D) \psi_j(D) h_c(t) \Vert_{L^3} \Vert e^{-i t \omega(D)} (x, y) f(t) \Vert_{L^6} \\
&\lesssim 2^{-j} \Vert \psi_j(D) h_c(t) \Vert_{L^2}^{\frac{1}{3}} \Vert e^{-i t \omega(D)} m_{\widehat{\mathcal{C}}}(D) \psi_j(D) h_c(t) \Vert_{L^4}^{\frac{2}{3}} \Vert \nabla (x, y) f(t) \Vert_{L^2} \\
&\lesssim t^{-\frac{7}{48}+101\delta} 2^{-\frac{23j}{16}-1000\delta j} \langle 2^j \rangle^{-\frac{13}{48}} \Vert u \Vert_X^2 \\
&\lesssim t^{-\frac{257}{2052}+2000\delta} 2^{-\frac{7j}{6}+1000 \delta j} \Vert u \Vert_X^2 \lesssim t^{-\frac{1}{8}} 2^{-\frac{7j}{6}+1000\delta j} \Vert u \Vert_X^2 
\end{align*}
\eqref{estimeeaprioriXcg-2-1-1} is simpler than \eqref{estimeeaprioriXcg-2-1-3} and \eqref{estimeeaprioriXcg-2-1-4}. 
\end{Dem}

\section{Quadratic weighted estimate, case (b, b)} \label{section-estimee-quad-bb} 

In this section, we prove: 

\begin{Prop} For any $t > 0$, 
\begin{align*}
\Vert \partial_x m_b(D) X_b h_b(t) \Vert_{L^2}^2 &\lesssim \Vert (1 + x^2 + |y|^2) u_0 \Vert_{H^1}^2 + \langle t \rangle^{\frac{1}{12}+301\delta} \Vert u \Vert_X^3 \\
\Vert (1 - m_{\widehat{\mathcal{R}}}(D)) \partial_x m_b(D) X_b h_b(t) \Vert_{L^2}^2 &\lesssim \Vert (1 + x^2 + |y|^2) u_0 \Vert_{H^1}^2 + \Vert u \Vert_X^3 
\end{align*} \label{propestimeeaprioriquadhbb} 
\end{Prop}

Again, we will separate into
\begin{align*}
\partial_x m_b(D) X_b h_b(t) &= h_{bb}(t) + g_{bb}(t) 
\end{align*}
where $h_{bb}(t)$ contains the initial data as well as the time-integrated terms, while $g_{bb}(t)$ contains boundary terms. We will then estimate 
\begin{align*}
\Vert g_{bb}(t) \Vert_{L^2} &\lesssim \Vert u \Vert_X^2 \\
\Vert h_{bb}(t) \Vert_{L^2}^2 &\lesssim \Vert h_{bb}(0) \Vert_{L^2}^2 + \Vert u \Vert_X^3 
\end{align*}
In order to get this decomposition, we start from Lemma \ref{lemstructurehalphagalpha} which gives the decomposition $h_b(t) = h_b(0) + h_{b, 2}(t) + h_{b, 3}(t)$. We may then compute directly $\partial_x m_b(D) X_b h_{b, 2}(t)$ applying the decomposition lemma \ref{lemdecompositionFGH}, but we need to prove a new decomposition for the cubic term $\partial_x m_b(D) X_b h_{b, 3}(t)$. The main difference is that, now, such a decomposition cannot be as smooth as before, as there are new resonances in the cubic term, and these resonances explain the time growth of this weighted estimate. It turns out these resonances are well localised, in the sense that they can be restricted to the case where the output frequency is localised sufficiently away from $\widehat{\mathcal{C}}, \widehat{\mathcal{L}}$ and $\widehat{\mathcal{P}}$. Also, even in the cubic estimate without time growth, some decompositions require more singular symbols than before. 

For the cubic term, we will prove the following lemma: 

\begin{Lem} It is possible to decompose
\begin{align*}
\partial_x m_b(D) X_b h_{b, 3}(t) &= h_{bb, 3}(t) + g_{bb, 3}(t) 
\end{align*}
such that 
\begin{align*}
\int_0^t \langle s \rangle^{-\frac{1}{24}-301\delta} \Vert \partial_s h_{bb, 3}(s) ~ ds \Vert_{L^2} &\lesssim \Vert u \Vert_X^3 \\
\int_0^t \Vert (1 - m_{\widehat{\mathcal{R}}}(D)) \partial_s h_{bb, 3}(s) ~ ds \Vert_{L^2} &\lesssim \Vert u \Vert_X^3 \\
\Vert g_{bb, 3}(t) \Vert_{L^2} &\lesssim \Vert u \Vert_X^3 \\
\Vert e^{-i t \omega(D)} \nabla g_{bb, 3}(t) \Vert_{L^4} &\lesssim t^{-\frac{5}{12}} \langle t \rangle^{-\frac{1}{8}+O(\delta)} \Vert u \Vert_X^3 
\end{align*} \label{lemmedecompositioncubiquehquadbb} 
\end{Lem} 

For now, we admit this lemma and prove Proposition \ref{propestimeeaprioriquadhbb}. 

\subsection{Proof of Proposition \ref{propestimeeaprioriquadhbb}}

By Lemma \ref{lemmedecompositioncubiquehquadbb}, we can already decompose $\partial_x m_b(D) X_b h_{b, 3}(t) = h_{bb, 3}(t) + g_{bb, 3}(t)$. It thus only remains to decompose also $\partial_x m_b(D) X_b h_{b, 2}(t)$ and to estimate every term. 

To that end, we apply Lemma \ref{lemstructurehalphagalpha} that gives an expression for $h_{b, 2}(t)$, and then the decomposition lemma \ref{lemdecompositionFGH} to compute an expression for $\partial_x m_b(D) X_b h_{b, 2}(t)$. We then get: 
\begin{subequations}
\begin{align}
&\xi_0 m_b(\overline{\xi}) \widehat{X}_b(\overline{\xi}) \cdot \nabla_{\overline{\xi}} \widehat{h}_{b, 2}(t, \overline{\xi}) \notag \\
&= - 2 \int_0^t \int e^{i s \varphi} \eta_0^2 m_b(\overline{\eta}) \widehat{X}_b(\overline{\eta}) \cdot \nabla_{\overline{\eta}} \widehat{h}_b(s, \overline{\eta}) \widehat{f}(s, \overline{\sigma}) ~ d\overline{\eta} ds \label{decompositionhb2-de-sym} \\ 
&+ \int_0^t \int e^{i s \varphi} \mu(\overline{\xi}, \overline{\eta}) \eta_0 m_b(\overline{\eta}) \widehat{X}_b(\overline{\eta}) \cdot \nabla_{\overline{\eta}} \widehat{h}_b(s, \overline{\eta}) \sigma_0 \widehat{f}(s, \overline{\sigma}) ~ d\overline{\eta} ds \label{decompositionhb2-de-b2bon} \\
&+ \int_0^t \int e^{i s \varphi} \mu(\overline{\xi}, \overline{\eta}) \frac{\eta_0 \varphi m_{\widehat{\mathcal{P}}}(\overline{\eta})}{(|\overline{\eta}| + |\overline{\sigma}|)^2} m_b(\overline{\eta}) \widehat{X}_b(\overline{\eta}) \cdot \nabla_{\overline{\eta}} \widehat{h}_b(s, \overline{\eta}) \widehat{f}(s, \overline{\sigma}) ~ d\overline{\eta} ds \label{decompositionhb2-de-b2eta0varphi} \\
&+ \int_0^t \int e^{i s \varphi} \mu(\overline{\xi}, \overline{\eta}) \frac{\varphi^2 m_{\widehat{\mathcal{P}}}(\overline{\eta})}{(|\overline{\eta}| + |\overline{\sigma}|)^4} m_b(\overline{\eta}) \widehat{X}_b(\overline{\eta}) \cdot \nabla_{\overline{\eta}} \widehat{h}_b(s, \overline{\eta}) \widehat{f}(s, \overline{\sigma}) ~ d\overline{\eta} ds \label{decompositionhb2-de-b2sigma0varphi} \\
&+ \sum_{\alpha = a, c} \int_0^t \int e^{i s \varphi} \mu(\overline{\xi}, \overline{\eta}) |\overline{\eta}| m_{\alpha}(\overline{\eta}) \widehat{X}_{\alpha}(\overline{\eta}) \cdot \nabla_{\overline{\eta}} \widehat{h}_b(s, \overline{\eta}) \sigma_0 \widehat{f}(s, \overline{\sigma}) ~ d\overline{\eta} ds \label{decompositionhb2-de-acb} \\ 
&+ \sum_{\alpha = a, b, c} \sum_{\beta = a, c} \int_0^t \int e^{i s \varphi} \mu(\overline{\xi}, \overline{\eta}) |\overline{\eta}| m_{\alpha}(\overline{\eta}) \widehat{X}_{\alpha}(\overline{\eta}) \cdot \nabla_{\overline{\eta}} \widehat{h}_{\beta}(s, \overline{\eta}) \sigma_0 \widehat{f}(s, \overline{\sigma}) ~ d\overline{\eta} ds \label{decompositionhb2-de-acac} 
\\
&+ \int_0^t \int e^{i s \varphi} \mu(\overline{\xi}, \overline{\eta}) \eta_0 \widehat{h}_b(s, \overline{\eta}) \sigma_0 m_b(\overline{\sigma}) \widehat{X}_b(\overline{\sigma}) \cdot \nabla_{\overline{\xi}} \widehat{f}(s, \overline{\sigma}) ~ d\overline{\eta} ds \label{decompositionhb2-ee-bbbon} \\ 
&+ \int_0^t \int e^{i s \varphi} \mu(\overline{\xi}, \overline{\eta}) \frac{\eta_0 \varphi}{(|\overline{\eta}| + |\overline{\sigma}|)^2} \widehat{h}_b(s, \overline{\eta}) m_b(\overline{\sigma}) \widehat{X}_b(\overline{\sigma}) \cdot \nabla_{\overline{\eta}} \widehat{f}(s, \overline{\sigma}) ~ d\overline{\eta} ds \label{decompositionhb2-ee-bbeta0varphi} \\
&+ \int_0^t \int e^{i s \varphi} \mu(\overline{\xi}, \overline{\eta}) \frac{\sigma_0 \varphi}{(|\overline{\eta}| + |\overline{\sigma}|)^2} \widehat{h}_b(s, \overline{\eta}) m_b(\overline{\sigma}) \widehat{X}_b(\overline{\sigma}) \cdot \nabla_{\overline{\eta}} \widehat{f}(s, \overline{\sigma}) ~ d\overline{\eta} ds \label{decompositionhb2-ee-bbsigma0varphi} \\
&+ \sum_{\alpha = a, c} \int_0^t \int e^{i s \varphi} \mu(\overline{\xi}, \overline{\eta}) \eta_0 \widehat{h}_b(s, \overline{\eta}) |\overline{\sigma}| m_{\alpha}(\overline{\sigma}) \widehat{X}_{\alpha}(\overline{\sigma}) \cdot \nabla_{\overline{\eta}} \widehat{f}(s, \overline{\sigma}) ~ d\overline{\eta} ds \label{decompositionhb2-ee-bac} \\
&+ \sum_{\alpha = a, c} \int_0^t \int e^{i s \varphi} \mu(\overline{\xi}, \overline{\eta}) |\overline{\eta}| \widehat{h}_{\alpha}(s, \overline{\eta}) \sigma_0 m_b(\overline{\sigma}) \widehat{X}_b(\overline{\sigma}) \cdot \nabla_{\overline{\xi}} \widehat{f}(s, \overline{\sigma}) ~ d\overline{\eta} ds \label{decompositionhb2-ee-acb} \\
&+ \sum_{\alpha, \beta = a, c} \int_0^t \int e^{i s \varphi} \mu(\overline{\xi}, \overline{\eta}) |\overline{\eta}| \widehat{h}_{\beta}(s, \overline{\eta}) |\overline{\sigma}| m_{\alpha}(\overline{\sigma}) \widehat{X}_{\alpha}(\overline{\sigma}) \cdot \nabla_{\overline{\eta}} \widehat{f}(s, \overline{\sigma}) ~ d\overline{\eta} ds \label{decompositionhb2-ee-acac} 
\\
&+ \int_0^t \int e^{i s \varphi} \mu(\overline{\xi}, \overline{\eta}) \eta_0 \widehat{h}_b(s, \overline{\eta}) \widehat{f}(s, \overline{\sigma}) ~ d\overline{\eta} ds \label{decompositionhb2-eders-b} \\
&+ \sum_{\alpha = a, b, c} \int_0^t \int e^{i s \varphi} \mu(\overline{\xi}, \overline{\eta}) \widehat{h}_{\alpha}(s, \overline{\eta}) \sigma_0 \widehat{f}(s, \overline{\sigma}) ~ d\overline{\eta} ds \label{decompositionhb2-eders-ac} \\
&+ \int_0^t \int e^{i s \varphi} \mu(\overline{\xi}, \overline{\eta}) \eta_0 m_b(\overline{\eta}) \widehat{X}_b(\overline{\eta}) \cdot \nabla_{\overline{\eta}} \widehat{f}(s, \overline{\eta}) \widehat{f}(s, \overline{\sigma}) ~ d\overline{\eta} ds \label{decompositionhb2-eders-bf} \\
&+ \sum_{\alpha = a, b, c} \int_0^t \int e^{i s \varphi} \mu(\overline{\xi}, \overline{\eta}) m_{\alpha}(\overline{\eta}) \widehat{X}_{\alpha}(\overline{\eta}) \cdot \nabla_{\overline{\eta}} \widehat{f}(s, \overline{\eta}) \sigma_0 \widehat{f}(s, \overline{\sigma}) ~ d\overline{\eta} ds \label{decompositionhb2-eders-acf} 
\\
&+ \int_0^t \int e^{i s \varphi} \mu(\overline{\xi}, \overline{\eta}) \widehat{f}(s, \overline{\eta}) \widehat{f}(s, \overline{\sigma}) ~ d\overline{\eta} ds \label{decompositionhb2-ders2} 
\\
&+ \int_0^t \int e^{i s \varphi} s m_g(\overline{\xi}, \overline{\eta}) \eta_0 \partial_s \widehat{h}_b(s, \overline{\eta}) \widehat{f}(s, \overline{\sigma}) ~ d\overline{\eta} ds \label{decompositionhb2-etr-b} \\
&+ \sum_{\alpha = a, b, c} \int_0^t \int e^{i s \varphi} s \mu(\overline{\xi}, \overline{\eta}) \partial_s \widehat{h}_{\alpha}(s, \overline{\eta}) |\overline{\sigma}| \widehat{f}(s, \overline{\sigma}) ~ d\overline{\eta} ds \label{decompositionhb2-etr-ac} 
\\
&+ \int_0^t \int e^{i s \varphi} s \mu(\overline{\xi}, \overline{\eta}) \eta_0 \widehat{h}_b(s, \overline{\eta}) \partial_s \widehat{f}(s, \overline{\sigma}) ~ d\overline{\eta} ds \label{decompositionhb2-et-b} \\
&+ \sum_{\alpha = a, c} \int_0^t \int e^{i s \varphi} s \mu(\overline{\xi}, \overline{\eta}) |\overline{\eta}| \widehat{h}_{\alpha}(s, \overline{\eta}) \partial_s \widehat{f}(s, \overline{\sigma}) ~ d\overline{\eta} ds \label{decompositionhb2-et-ac} \\
&+ \int_0^t \int e^{i s \varphi} s \mu(\overline{\xi}, \overline{\eta}) \frac{\varphi}{(|\overline{\eta}| + |\overline{\sigma}|)^2} \widehat{h}_b(s, \overline{\eta}) \partial_s \widehat{f}(s, \overline{\sigma}) ~ d\overline{\eta} ds \label{decompositionhb2-et-varphi} 
\\
&+ \int_0^t \int e^{i s \varphi} s \mu(\overline{\xi}, \overline{\eta}) \partial_s \widehat{f}(s, \overline{\eta}) \widehat{f}(s, \overline{\sigma}) ~ d\overline{\eta} ds \label{decompositionhb2-tders} 
\\
&- t \int e^{i t \varphi} m_g(\overline{\xi}, \overline{\eta}) \eta_0 \widehat{h}_b(t, \overline{\eta}) \widehat{f}(t, \overline{\sigma}) ~ d\overline{\eta} \label{decompositionhb2-bord-b} \\
&+ \sum_{\alpha = a, c} t \int e^{i t \varphi} \mu(\overline{\xi}, \overline{\eta}) \frac{|\overline{\eta}|}{|\overline{\eta}|+|\overline{\sigma}|} \widehat{h}_{\alpha}(t, \overline{\eta}) |\overline{\sigma}| \widehat{f}(t, \overline{\sigma}) ~ d\overline{\eta} \label{decompositionhb2-bord-ac} \\
&+ t \int e^{i t \varphi} \mu(\overline{\xi}, \overline{\eta}) \frac{\eta_0 \sigma_0}{|\overline{\eta}|+|\overline{\sigma}|} \widehat{h}_b(t, \overline{\eta}) \widehat{f}(t, \overline{\sigma}) ~ d\overline{\eta} \label{decompositionhb2-bord-bbon} \\
&+ t \int e^{i t \varphi} \mu(\overline{\xi}, \overline{\eta}) \frac{\eta_0 \sigma_0 |\overline{\eta}|}{(|\eta_0|+|\sigma_0|)(|\overline{\eta}|+|\overline{\sigma}|)} \widehat{h}_b(t, \overline{\eta}) \widehat{f}(t, \overline{\sigma}) ~ d\overline{\eta} \label{decompositionhb2-bord-bbonbis} \\
&+ t \int e^{i t \varphi} \mu(\overline{\xi}, \overline{\eta}) \frac{\varphi}{|\overline{\eta}|^2+|\overline{\sigma}|^2} \widehat{h}_b(t, \overline{\eta}) \widehat{f}(t, \overline{\sigma}) ~ d\overline{\eta} \label{decompositionhb2-bord-bvarphi} \\
&+ t \int e^{i t \varphi} \mu(\overline{\xi}, \overline{\eta}) \widehat{f}(t, \overline{\eta}) \widehat{f}(t, \overline{\sigma}) ~ d\overline{\eta} \label{decompositionhb2-bord-ders} 
\end{align}
\end{subequations}

We organised the terms by the following ordering: 
\begin{itemize}
\item from \eqref{decompositionhb2-de-sym} to \eqref{decompositionhb2-de-acac}, these are interactions terms between a quadratic weight term and an unweighted term; 
\item from \eqref{decompositionhb2-ee-bbbon} to \eqref{decompositionhb2-ee-acac}, these are interactions terms between two simply weighted terms; 
\item from \eqref{decompositionhb2-eders-b} to \eqref{decompositionhb2-eders-acf}, these are interactions terms between a simply weighted term and an unweighted term; 
\item \eqref{decompositionhb2-ders2} is the interaction term of two unweighted terms; 
\item \eqref{decompositionhb2-etr-b} and \eqref{decompositionhb2-etr-ac} are the interaction terms between the time derivative of a simply weighted term and an unweighted term; 
\item from \eqref{decompositionhb2-et-b} to \eqref{decompositionhb2-et-varphi}, these are interaction terms between a simply weighted term and the time derivative of an unweighted term; 
\item \eqref{decompositionhb2-tders} is an interaction term between an unweighted time derivative and an unweighted term; 
\item from \eqref{decompositionhb2-bord-b} to \eqref{decompositionhb2-bord-ders}, these are boundary terms in time. 
\end{itemize}

To understand where each term is coming from, we provide the following table: the vertical axis corresponds to a line in the abstract expression given by the decomposition lemma, applied on a term given by Lemma \ref{lemstructurehalphagalpha} on the horizontal axis, and we read an output which is one of the above terms (or two sometimes). 

\begin{tabular}{|c||c|c|c|c|c|c|}
\hline 
& \eqref{lemstructurehb-02-symeta} & \eqref{lemstructurehb-03-resxetaac} & \eqref{lemstructurehb-04-resxetabbon} & \eqref{lemstructurehb-05-resxetabbonbis} & \eqref{lemstructurehb-06-resxetabphi} & \eqref{lemstructurehb-11-dersymb} 
\\ \hline \hline 
\eqref{lemdecequH-01-symeta} & \eqref{decompositionhb2-de-sym} & \eqref{decompositionhb2-de-acac} & \eqref{decompositionhb2-de-b2bon} & \eqref{decompositionhb2-de-b2bon} & \eqref{decompositionhb2-de-b2eta0varphi} & \eqref{decompositionhb2-eders-bf} 
\\ \hline
\eqref{lemdecequH-02-symsigma} & \eqref{decompositionhb2-ee-bbbon} & \eqref{decompositionhb2-ee-acb} & \eqref{decompositionhb2-ee-bbbon} & \eqref{decompositionhb2-ee-bbbon} & \eqref{decompositionhb2-ee-bbeta0varphi} & \eqref{decompositionhb2-eders-bf} 
\\ \hline
\eqref{lemdecequH-03-resxetaac} & \eqref{decompositionhb2-de-acb} & \eqref{decompositionhb2-de-acac} & \eqref{decompositionhb2-de-acac} & \eqref{decompositionhb2-de-acac} & \eqref{decompositionhb2-de-acac} & \eqref{decompositionhb2-eders-acf} 
\\ \hline
\eqref{lemdecequH-04-resxetabbon} & \eqref{decompositionhb2-de-b2bon} & \eqref{decompositionhb2-de-acac} & \eqref{decompositionhb2-de-b2bon} & \eqref{decompositionhb2-de-b2bon} & \eqref{decompositionhb2-de-b2bon} & \eqref{decompositionhb2-eders-bf} 
\\ \hline
\eqref{lemdecequH-05-resxetabbonbis} & \eqref{decompositionhb2-de-b2bon} & \eqref{decompositionhb2-de-acac} & \eqref{decompositionhb2-de-b2bon} & \eqref{decompositionhb2-de-b2bon} & \eqref{decompositionhb2-de-b2bon} & \eqref{decompositionhb2-eders-bf} 
\\ \hline
\eqref{lemdecequH-06-resxetabphi} & \eqref{decompositionhb2-de-b2eta0varphi} & \eqref{decompositionhb2-de-acac} & \eqref{decompositionhb2-de-b2bon} & \eqref{decompositionhb2-de-b2bon} & \eqref{decompositionhb2-de-b2sigma0varphi} & $\eqref{decompositionhb2-eders-bf} + \eqref{decompositionhb2-eders-acf}$ 
\\ \hline
\eqref{lemdecequH-07-resxsigmaac} & \eqref{decompositionhb2-ee-bac} & \eqref{decompositionhb2-ee-acac} & \eqref{decompositionhb2-ee-bac} & \eqref{decompositionhb2-ee-bac} & \eqref{decompositionhb2-ee-bac} & \eqref{decompositionhb2-eders-acf} 
\\ \hline
\eqref{lemdecequH-08-resxsigmabbon} & \eqref{decompositionhb2-ee-bbbon} & \eqref{decompositionhb2-ee-acb} & \eqref{decompositionhb2-ee-bbbon} & \eqref{decompositionhb2-ee-bbbon} & \eqref{decompositionhb2-ee-bbbon} & \eqref{decompositionhb2-eders-bf} 
\\ \hline
\eqref{lemdecequH-09-resxsigmabbonbis} & \eqref{decompositionhb2-ee-bbbon} & \eqref{decompositionhb2-ee-acb} & \eqref{decompositionhb2-ee-bbbon} & \eqref{decompositionhb2-ee-bbbon} & \eqref{decompositionhb2-ee-bbbon} & \eqref{decompositionhb2-eders-bf} 
\\ \hline
\eqref{lemdecequH-10-resxsigmabphi} & \eqref{decompositionhb2-ee-bbeta0varphi} & \eqref{decompositionhb2-ee-acb} & \eqref{decompositionhb2-ee-bbbon} & \eqref{decompositionhb2-ee-bbbon} & $\eqref{decompositionhb2-ee-bbeta0varphi}+\eqref{decompositionhb2-ee-bbsigma0varphi}$ & $\eqref{decompositionhb2-eders-bf} + \eqref{decompositionhb2-eders-acf}$
\\ \hline
\eqref{lemdecequH-11-dersymb} & \eqref{decompositionhb2-eders-b} & \eqref{decompositionhb2-eders-ac} & \eqref{decompositionhb2-eders-b} & \eqref{decompositionhb2-eders-b} & $\eqref{decompositionhb2-eders-b}+\eqref{decompositionhb2-eders-bf}$ & \eqref{decompositionhb2-ders2} 
\\ \hline
\eqref{lemdecequH-12-resteta} & \eqref{decompositionhb2-etr-b} & \eqref{decompositionhb2-etr-ac} & \eqref{decompositionhb2-etr-ac} & \eqref{decompositionhb2-etr-ac} & \eqref{decompositionhb2-etr-ac} & \eqref{decompositionhb2-tders} 
\\ \hline
\eqref{lemdecequH-13-restsigma} & \eqref{decompositionhb2-et-b} & \eqref{decompositionhb2-et-ac} & \eqref{decompositionhb2-et-b} & \eqref{decompositionhb2-et-b} & \eqref{decompositionhb2-et-varphi} & \eqref{decompositionhb2-tders} 
\\ \hline
\eqref{lemdecequG} & \eqref{decompositionhb2-bord-b} & \eqref{decompositionhb2-bord-ac} & \eqref{decompositionhb2-bord-bbon} & \eqref{decompositionhb2-bord-bbonbis} & \eqref{decompositionhb2-bord-bvarphi} & \eqref{decompositionhb2-bord-ders} 
\\ \hline 
\end{tabular}

\paragraph{Direct estimates} We can directly estimate: 
\begin{align*}
&\begin{aligned}
\Vert \partial_t \eqref{decompositionhb2-de-b2bon} \Vert_{L^2} &\lesssim \Vert \partial_x X_b h_b(t) \Vert_{L^2} \Vert \partial_x u(t) \Vert_{L^{\infty}} \lesssim \langle t \rangle^{\frac{1}{24}} \Vert u \Vert_X t^{-\frac{5}{6}} \langle t \rangle^{-\frac{1}{4}+100\delta} \Vert u \Vert_X \\
&\lesssim t^{-\frac{5}{6}} \langle t \rangle^{-\frac{5}{24}+100\delta} \Vert u \Vert_X^2 
\end{aligned} \\
&\begin{aligned} 
\Vert \partial_t \eqref{decompositionhb2-de-acb} \Vert_{L^2} &\lesssim \sum_{\alpha = a, c} \Vert \nabla X_{\alpha} h_b(t) \Vert_{L^2} \Vert \partial_x u(t) \Vert_{L^{\infty}} \lesssim t^{-\frac{5}{6}} \langle t \rangle^{-\frac{1}{4}+100\delta} \Vert u \Vert_X^2 
\end{aligned} \\
&\begin{aligned} 
\Vert \partial_t \eqref{decompositionhb2-de-acac} \Vert_{L^2} &\lesssim \sum_{\alpha = a, b, c} \sum_{\beta = a, c} \Vert \nabla X_{\alpha} h_{\beta}(t) \Vert_{L^2} \Vert \partial_x u(t) \Vert_{L^{\infty}} \lesssim t^{-\frac{5}{6}} \langle t \rangle^{-\frac{1}{4}+100\delta} \Vert u \Vert_X^2 
\end{aligned} \\
&\begin{aligned} 
\Vert \partial_t \eqref{decompositionhb2-ee-bbbon} \Vert_{L^2} &\lesssim \Vert e^{-i s \omega(D)} \partial_x h_b(t) \Vert_{L^4} \Vert e^{-i s \omega(D)} \partial_x X_b f(t) \Vert_{L^4} \lesssim t^{-\frac{511}{512}} \langle t \rangle^{-\frac{1}{256}} \Vert u \Vert_X^2 
\end{aligned} \\
&\begin{aligned} 
&\int_0^t s^{1-\delta} \langle s \rangle^{2\delta} \Vert \partial_s \eqref{decompositionhb2-ee-bac} \Vert_{L^2}^2 ~ ds \lesssim \sum_{\alpha = a, c} \int_0^t s^{1-\delta} \langle s \rangle^{2\delta} \Vert e^{-i s \omega(D)} \partial_x h_b(s) \Vert_{L^4}^2 \Vert e^{-i s \omega(D)} \nabla h_{\alpha}(s) \Vert_{L^4}^2 ~ ds \\
&\quad \lesssim \sum_{\alpha = a, c} \left( \int_0^t s^{-1+\frac{1}{256}-\delta} \langle s \rangle^{2\delta-\frac{1}{128}} \Vert u \Vert_X^4 ~ ds \right)^{\frac{1}{2}} \left( \int_0^t s^{1-\delta} \langle s \rangle^{2\delta} \Vert e^{-i s \omega(D)} \nabla h_{\alpha}(s) \Vert_{L^4}^4 ~ ds \right)^{\frac{1}{2}} \\
&\quad \lesssim \Vert u \Vert_X^4 
\end{aligned} \\
&\begin{aligned} 
\int_0^t s^{1-\delta} \langle s \rangle^{2\delta} \Vert \partial_s \eqref{decompositionhb2-ee-acb} \Vert_{L^2}^2 ~ ds &\lesssim \sum_{\alpha = a, c} \int_0^t s^{1-\delta} \langle s \rangle^{2\delta} \Vert e^{-i s \omega(D)} \nabla h_{\alpha}(s) \Vert_{L^4}^2 \Vert e^{-i s \omega(D)} \partial_x X_b f(s) \Vert_{L^4}^2 ~ ds \\
&\lesssim \Vert u \Vert_X^4 
\end{aligned} \\
&\begin{aligned} 
\int_0^t s^{1-\delta} \langle s \rangle^{2\delta} \Vert \partial_s \eqref{decompositionhb2-ee-acac} \Vert_{L^2}^2 ~ ds &\lesssim \int_0^t s^{1-\delta} \langle s \rangle^{2\delta} \left( \sum_{\alpha = a, c} \Vert e^{-i s \omega(D)} \nabla h_{\alpha}(s) \Vert_{L^4} \right)^4 ~ ds \lesssim \Vert u \Vert_X^4 
\end{aligned} \\
&\begin{aligned} 
\Vert \partial_t \eqref{decompositionhb2-eders-b} \Vert_{L^2} &\lesssim \Vert e^{-i t \omega(D)} \partial_x h_b(t) \Vert_{L^4} \Vert u(t) \Vert_{L^4} \lesssim t^{-\frac{5}{6}} \langle t \rangle^{-\frac{5}{24}+100\delta} \Vert u \Vert_X^2 
\end{aligned} \\
&\begin{aligned} 
\Vert \partial_t \eqref{decompositionhb2-eders-ac} \Vert_{L^2} &\lesssim \sum_{\alpha = a, b, c} \Vert h_{\alpha}(t) \Vert_{L^2} \Vert \partial_x u(t) \Vert_{L^{\infty}} \lesssim t^{-\frac{5}{6}} \langle t \rangle^{-\frac{1}{4}+100\delta} \Vert u \Vert_X^2 
\end{aligned} \\
&\begin{aligned} 
\Vert \partial_t \eqref{decompositionhb2-eders-bf} \Vert_{L^2} &\lesssim \Vert e^{-i t \omega(D)} \partial_x X_b f(t) \Vert_{L^4} \Vert u(t) \Vert_{L^4} \lesssim t^{-\frac{5}{6}} \langle t \rangle^{-\frac{5}{24}+100\delta} \Vert u \Vert_X^2 
\end{aligned} \\
&\begin{aligned} 
\Vert \partial_t \eqref{decompositionhb2-eders-acf} \Vert_{L^2} &\lesssim \sum_{\alpha = a, b, c} \Vert X_{\alpha} f(t) \Vert_{L^2} \Vert \partial_x u(t) \Vert_{L^{\infty}} \lesssim t^{-\frac{5}{6}} \langle t \rangle^{-\frac{1}{4}+100\delta} \Vert u \Vert_X^2 
\end{aligned} \\
&\begin{aligned} 
\Vert \partial_t \eqref{decompositionhb2-ders2} \Vert_{L^2} &\lesssim \Vert u(t) \Vert_{L^4}^2 \lesssim t^{-\frac{5}{6}} \langle t \rangle^{-\frac{1}{4}+100\delta} \Vert u \Vert_X^2 
\end{aligned} \\
&\begin{aligned} 
\int_0^t \Vert \partial_s \eqref{decompositionhb2-etr-ac} \Vert_{L^2} ~ ds &\lesssim \int_0^t s \sum_{\alpha = a, b, c} \Vert \partial_s h_{\alpha}(s) \Vert_{L^2} \Vert \partial_x u(s) \Vert_{L^{\infty}} ~ ds \\
&\lesssim \int_0^t s^{\frac{1}{6}-\delta} \langle s \rangle^{-\frac{1}{4}+102\delta} \Vert \partial_s h_{\alpha}(s) \Vert_{L^2} \Vert u \Vert_X ~ ds \lesssim \Vert u \Vert_X^3 
\end{aligned} \\
&\begin{aligned} 
\Vert \partial_t \eqref{decompositionhb2-et-b} \Vert_{L^2} &\lesssim t \Vert e^{-i t \omega(D)} \partial_x h_b(t) \Vert_{L^4} \Vert e^{-i t \omega(D)} \partial_t f(t) \Vert_{L^4} \lesssim t^{-\frac{1}{4}} \langle t \rangle^{-\frac{7}{8}+250\delta} \Vert u \Vert_X^3 
\end{aligned} \\
&\begin{aligned} 
\Vert \partial_t \eqref{decompositionhb2-et-ac} \Vert_{L^2} &\lesssim t \sum_{\alpha = a, c} \Vert e^{-i t \omega(D)} \nabla h_{\alpha}(t) \Vert_{L^4} \Vert e^{-i t \omega(D)} \partial_t f(t) \Vert_{L^4} \lesssim t^{-\frac{2}{3}} \langle t \rangle^{-\frac{19}{24}+250\delta} \Vert u \Vert_X^3 
\end{aligned} \\
&\begin{aligned} 
\Vert \partial_t \eqref{decompositionhb2-tders} \Vert_{L^2} &\lesssim t \Vert e^{-i t \omega(D)} \partial_t f(t) \Vert_{L^4} \Vert u(t) \Vert_{L^4} \lesssim t^{-\frac{1}{4}} \langle t \rangle^{-\frac{11}{12}+200\delta} \Vert u \Vert_X^3 
\end{aligned} 
\end{align*}

We can also estimate boundary terms: 
\begin{align*}
\Vert \eqref{decompositionhb2-bord-b} \Vert_{L^2} &\lesssim t \Vert e^{-i t \omega(D)} \partial_x h_b(t) \Vert_{L^4} \Vert u(t) \Vert_{L^4} \lesssim t^{\frac{1}{6}} \langle t \rangle^{-\frac{5}{24}+200\delta} \Vert u \Vert_X^2 \\
\Vert \eqref{decompositionhb2-bord-ac} \Vert_{L^2} &\lesssim t \sum_{\alpha = a, b, c} \Vert h_{\alpha}(t) \Vert_{L^2} \Vert \nabla u(t) \Vert_{L^{\infty}} \lesssim t^{\frac{1}{6}} \langle t \rangle^{-\frac{1}{6}+100\delta} \Vert u \Vert_X^2 \\
\Vert \eqref{decompositionhb2-bord-ders} \Vert_{L^2} &\lesssim t \Vert u(t) \Vert_{L^4}^2 \lesssim t^{\frac{1}{6}} \langle t \rangle^{-\frac{1}{4}+100\delta} \Vert u \Vert_X^2 
\end{align*}
It is possible to estimate \eqref{decompositionhb2-bord-bbon}, \eqref{decompositionhb2-bord-bbonbis} and \eqref{decompositionhb2-bord-bvarphi} like \eqref{decompositionhb2-bord-ac}. 

\paragraph{Estimate with $\varphi$} Let us now consider the terms having a factor $\varphi$ in their expression, i.e. \eqref{decompositionhb2-de-b2eta0varphi}, \eqref{decompositionhb2-de-b2sigma0varphi}, \eqref{decompositionhb2-ee-bbeta0varphi}, \eqref{decompositionhb2-ee-bbsigma0varphi} and \eqref{decompositionhb2-et-varphi}. 

For \eqref{decompositionhb2-et-varphi}, we develop it into a cubic term and decompose: 
\begin{subequations}
\begin{align}
\eqref{decompositionhb2-et-varphi} &= \int_0^t \int \int e^{i s \varphi_3(\overline{\xi}, \overline{\eta}, \overline{\sigma})} s \mu(\overline{\xi}, \overline{\eta}, \overline{\sigma}) \chi\left( \frac{|\overline{\rho}|}{|\overline{\sigma}|} \right) \frac{\sigma_0 \varphi(\overline{\xi}, \overline{\eta})}{(|\overline{\eta}|+|\overline{\sigma}+\overline{\rho}|)^2} \widehat{h}_b(s, \overline{\eta}) \widehat{f}(s, \overline{\sigma}) \widehat{f}(s, \overline{\rho}) ~ d\overline{\eta} d\overline{\sigma} ds \label{decompositionhb2-et-varphi-1} \\
&\quad + \int_0^t \int \int e^{i s \varphi_3(\overline{\xi}, \overline{\eta}, \overline{\sigma})} s \mu(\overline{\xi}, \overline{\eta}, \overline{\sigma}) \chi\left( \frac{|\overline{\sigma}|}{|\overline{\rho}|} \right)\frac{\sigma_0 \varphi(\overline{\xi}, \overline{\eta})}{(|\overline{\eta}|+|\overline{\sigma}+\overline{\rho}|)^2} \widehat{h}_b(s, \overline{\eta}) \widehat{f}(s, \overline{\sigma}) \widehat{f}(s, \overline{\rho}) ~ d\overline{\eta} d\overline{\sigma} ds \label{decompositionhb2-et-varphi-2}
\end{align}
\end{subequations} 
where the function $\chi$ is used here to localise on $\{ |\overline{\rho}| \lesssim |\overline{\sigma}| \}$ or on $\{ |\overline{\sigma}| \gtrsim |\overline{\rho}| \}$. We then note that 
\begin{align*}
\frac{\varphi}{(|\overline{\eta}|+|\overline{\sigma}+\overline{\rho}|)^2} \chi\left( \frac{|\overline{\sigma}|}{|\overline{\rho}|} \right) &= O(\overline{\rho}) 
\end{align*}
so that we can directly estimate: 
\begin{align*}
\Vert \partial_t \eqref{decompositionhb2-et-varphi-2} \Vert_{L^2} &\lesssim t \Vert h_b(t) \Vert_{L^2} \Vert \partial_x u(t) \Vert_{L^{\infty}} \Vert \nabla u(t) \Vert_{L^{\infty}} \\
&\lesssim t^{-\frac{2}{3}} \langle t \rangle^{-\frac{5}{12}+200\delta} \Vert u \Vert_X^3 
\end{align*}
On the other hand, for \eqref{decompositionhb2-et-varphi-1}, 
\begin{align*}
\varphi(\overline{\xi}, \overline{\eta}) &= \varphi_3(\overline{\xi}, \overline{\eta}, \overline{\sigma}) + O(\overline{\rho}) 
\end{align*}
The term having a factor $O(\overline{\rho})$ can be estimated like \eqref{decompositionhb2-et-varphi-2}. On the other, we can apply an integration by parts in time: 
\begin{subequations}
\begin{align}
&\int_0^t \int \int e^{i s \varphi_3(\overline{\xi}, \overline{\eta}, \overline{\sigma})} s \mu(\overline{\xi}, \overline{\eta}, \overline{\sigma}) \chi\left( \frac{|\overline{\rho}|}{|\overline{\sigma}|} \right) \frac{\sigma_0 \varphi_3(\overline{\xi}, \overline{\eta}, \overline{\sigma})}{(|\overline{\eta}|+|\overline{\sigma}+\overline{\rho}|)^2} \widehat{h}_b(s, \overline{\eta}) \widehat{f}(s, \overline{\sigma}) \widehat{f}(s, \overline{\rho}) ~ d\overline{\eta} d\overline{\sigma} ds \notag \\
&\quad = - \int_0^t \int \int e^{i s \varphi_3(\overline{\xi}, \overline{\eta}, \overline{\sigma})} \mu(\overline{\xi}, \overline{\eta}, \overline{\sigma}) \chi\left( \frac{|\overline{\rho}|}{|\overline{\sigma}|} \right) \frac{\sigma_0}{(|\overline{\eta}|+|\overline{\sigma}+\overline{\rho}|)^2} \widehat{h}_b(s, \overline{\eta}) \widehat{f}(s, \overline{\sigma}) \widehat{f}(s, \overline{\rho}) ~ d\overline{\eta} d\overline{\sigma} ds \label{decompositionhb2-et-varphi-1-2} \\
&\quad \quad - \int_0^t \int \int e^{i s \varphi_3(\overline{\xi}, \overline{\eta}, \overline{\sigma})} s \mu(\overline{\xi}, \overline{\eta}, \overline{\sigma}) \chi\left( \frac{|\overline{\rho}|}{|\overline{\sigma}|} \right) \frac{\sigma_0}{(|\overline{\eta}|+|\overline{\sigma}+\overline{\rho}|)^2} \partial_s \widehat{h}_b(s, \overline{\eta}) \widehat{f}(s, \overline{\sigma}) \widehat{f}(s, \overline{\rho}) ~ d\overline{\eta} d\overline{\sigma} ds \label{decompositionhb2-et-varphi-1-3} \\
&\quad \quad - \int_0^t \int \int e^{i s \varphi_3(\overline{\xi}, \overline{\eta}, \overline{\sigma})} s \mu(\overline{\xi}, \overline{\eta}, \overline{\sigma}) \chi\left( \frac{|\overline{\rho}|}{|\overline{\sigma}|} \right) \frac{\sigma_0}{(|\overline{\eta}|+|\overline{\sigma}+\overline{\rho}|)^2} \widehat{h}_b(s, \overline{\eta}) \partial_s \left( \widehat{f}(s, \overline{\sigma}) \widehat{f}(s, \overline{\rho}) \right) ~ d\overline{\eta} d\overline{\sigma} ds \label{decompositionhb2-et-varphi-1-4} \\
&\quad \quad + \int \int e^{i t \varphi_3(\overline{\xi}, \overline{\eta}, \overline{\sigma})} t \mu(\overline{\xi}, \overline{\eta}, \overline{\sigma}) \chi\left( \frac{|\overline{\rho}|}{|\overline{\sigma}|} \right) \frac{\sigma_0}{(|\overline{\eta}|+|\overline{\sigma}+\overline{\rho}|)^2} \widehat{h}_b(t, \overline{\eta}) \widehat{f}(t, \overline{\sigma}) \widehat{f}(t, \overline{\rho}) ~ d\overline{\eta} d\overline{\sigma} \label{decompositionhb2-et-varphi-1-5}
\end{align}
\end{subequations}
We then estimate:
\begin{align*}
\Vert \partial_t \eqref{decompositionhb2-et-varphi-1-2} \Vert_{L^2} &\lesssim \Vert e^{-i t \omega(D)} |\nabla|^{-1} h_b(t) \Vert_{L^6} \Vert \partial_x u(t) \Vert_{L^{\infty}} \Vert |\nabla|^{-1} u(t) \Vert_{L^3} \\
&\lesssim t^{-\frac{5}{6}} \langle t \rangle^{-\frac{1}{4}+100\delta} \Vert u \Vert_X^3 \\
\int_0^t \Vert \partial_s \eqref{decompositionhb2-et-varphi-1-3} \Vert_{L^2} ~ ds &\lesssim \int_0^t s \Vert e^{-i s \omega(D)} |\nabla|^{-1} \partial_s h_b(s) \Vert_{L^6} \Vert \partial_x u(s) \Vert_{L^{\infty}} \Vert |\nabla|^{-1} u(s) \Vert_{L^3} ~ ds \\
&\lesssim \int_0^t s^{\frac{1}{6}} \langle s \rangle^{-\frac{1}{4}+100\delta} \Vert \partial_t h_b(t) \Vert_{L^2} \Vert u \Vert_X^2 ~ ds \lesssim \Vert u \Vert_X^4 \\
\Vert \partial_t \eqref{decompositionhb2-et-varphi-1-4} \Vert_{L^2} &\lesssim t \Vert e^{-i t \omega(D)} |\nabla|^{-1} h_b(t) \Vert_{L^{12}} \Vert \partial_t f(t) \Vert_{L^4} \Vert u(t) \Vert_{L^6} \\
&\lesssim \Vert h_b(t) \Vert_{H^1} t^{\frac{1}{6}} \langle t \rangle^{-\frac{19}{24}+150\delta} \Vert u \Vert_X^2 \Vert \langle \nabla \rangle u(t) \Vert_{L^4} \\
&\lesssim t^{-\frac{1}{4}} \langle t \rangle^{-\frac{11}{12}+150\delta} \Vert u \Vert_X^4 \\
\Vert \eqref{decompositionhb2-et-varphi-1-5} \Vert_{L^2} &\lesssim t \Vert e^{- i t \omega(D)} |\nabla|^{-1} h_b(t) \Vert_{L^6} \Vert \partial_x u(t) \Vert_{L^{\infty}} \Vert |\nabla|^{-1} u(t) \Vert_{L^3} \\
&\lesssim t^{\frac{1}{6}} \langle t \rangle^{-\frac{1}{4}+100\delta} \Vert u \Vert_X^3 
\end{align*}

For \eqref{decompositionhb2-ee-bbsigma0varphi}, we develop $m_b(D) X_b f(t) = h_b(t) + g_b(t)$: 
\begin{subequations}
\begin{align}
\eqref{decompositionhb2-ee-bbsigma0varphi} &= \int_0^t \int e^{i s \varphi} \mu(\overline{\xi}, \overline{\eta}) \frac{\sigma_0 \varphi}{(|\overline{\eta}| + |\overline{\sigma}|)^2} \widehat{h}_b(s, \overline{\eta}) \widehat{h}_b(s, \overline{\sigma}) ~ d\overline{\eta} ds \label{decompositionhb2-ee-bbsigma0varphi-1} \\
&\quad + \int_0^t \int e^{i s \varphi} \mu(\overline{\xi}, \overline{\eta}) \frac{\sigma_0 \varphi}{(|\overline{\eta}| + |\overline{\sigma}|)^2} \widehat{h}_b(s, \overline{\eta}) \widehat{g}_b(s, \overline{\sigma}) ~ d\overline{\eta} ds\label{decompositionhb2-ee-bbsigma0varphi-2}
\end{align}
\end{subequations}
\eqref{decompositionhb2-ee-bbsigma0varphi-2} is similar to \eqref{decompositionhb2-et-varphi} (up to the symbol) and can be estimated the same way. On \eqref{decompositionhb2-ee-bbsigma0varphi-1}, we apply an integration by parts in time: 
\begin{subequations}
\begin{align}
\eqref{decompositionhb2-ee-bbsigma0varphi-1} &= - \int_0^t \int e^{i s \varphi} \mu(\overline{\xi}, \overline{\eta}) \frac{\sigma_0}{(|\overline{\eta}| + |\overline{\sigma}|)^2} \partial_s \left( \widehat{h}_b(s, \overline{\eta}) \widehat{h}_b(s, \overline{\sigma}) \right) ~ d\overline{\eta} ds \label{decompositionhb2-ee-bbsigma0varphi-1-3} \\
&\quad + \int e^{i t \varphi} \mu(\overline{\xi}, \overline{\eta}) \frac{\sigma_0}{(|\overline{\eta}| + |\overline{\sigma}|)^2} \widehat{h}_b(t, \overline{\eta}) \widehat{h}_b(t, \overline{\sigma}) ~ d\overline{\eta} \label{decompositionhb2-ee-bbsigma0varphi-1-4} \\
&\quad + \int \mu(\overline{\xi}, \overline{\eta}) \frac{\sigma_0}{(|\overline{\eta}| + |\overline{\sigma}|)^2} \widehat{h}_b(0, \overline{\eta}) \widehat{h}_b(0, \overline{\sigma}) ~ d\overline{\eta} \label{decompositionhb2-ee-bbsigma0varphi-1-5}
\end{align}
\end{subequations} 
We then estimate: 
\begin{align*}
\int_0^t \Vert \partial_s \eqref{decompositionhb2-ee-bbsigma0varphi-1-3} \Vert_{L^2} ~ ds &\lesssim \int_0^t \Vert e^{-i s \omega(D)} |\nabla|^{-1} \partial_s h_b(s) \Vert_{L^6} \Vert e^{-i s \omega(D)} h_b(s) \Vert_{L^3} ~ ds \\
&\lesssim \int_0^t \Vert \partial_s h_b(s) \Vert_{L^2} \Vert h_b(s) \Vert_{H^1} ~ ds \lesssim \Vert u \Vert_X^2 \\
\Vert \eqref{decompositionhb2-ee-bbsigma0varphi-1-4} \Vert_{L^2} &\lesssim \Vert e^{-i t \omega(D)} |\nabla|^{-1} h_b(t) \Vert_{L^6} \Vert e^{-i t \omega(D)} h_b(t) \Vert_{L^3} \\
&\lesssim \Vert h_b(t) \Vert_{L^2} \Vert h_b(t) \Vert_{H^1} ~ \lesssim \Vert u \Vert_X^2 
\end{align*}
We estimate \eqref{decompositionhb2-ee-bbsigma0varphi-1-5} like \eqref{decompositionhb2-ee-bbsigma0varphi-1-4}. 

For \eqref{decompositionhb2-ee-bbeta0varphi}, we develop as well: 
\begin{subequations}
\begin{align}
\eqref{decompositionhb2-ee-bbeta0varphi} &= \int_0^t \int e^{i s \varphi} \mu(\overline{\xi}, \overline{\eta}) \frac{\eta_0 \varphi}{(|\overline{\eta}| + |\overline{\sigma}|)^2} \widehat{h}_b(s, \overline{\eta}) \widehat{h}_b(s, \overline{\sigma}) ~ d\overline{\eta} ds \label{decompositionhb2-ee-bbeta0varphi-1} \\
&\quad + \int_0^t \int e^{i s \varphi} \mu(\overline{\xi}, \overline{\eta}) \frac{\eta_0 \varphi}{(|\overline{\eta}| + |\overline{\sigma}|)^2} \widehat{h}_b(s, \overline{\eta}) \widehat{g}_b(s, \overline{\sigma}) ~ d\overline{\eta} ds \label{decompositionhb2-ee-bbeta0varphi-2} 
\end{align}
\end{subequations}
\eqref{decompositionhb2-ee-bbeta0varphi-1} is similar to \eqref{decompositionhb2-ee-bbsigma0varphi-1} and can be estimated the same way. For \eqref{decompositionhb2-ee-bbeta0varphi-2}, we estimate directly: 
\begin{align*}
\Vert \partial_t \eqref{decompositionhb2-ee-bbeta0varphi-2} \Vert_{L^2} &\lesssim \Vert e^{-i t \omega(D)} \partial_x h_b(t) \Vert_{L^4} \Vert e^{-i t \omega(D)} \nabla g_b(t) \Vert_{L^4} \\
&\lesssim t^{-\frac{2}{3}} \langle t \rangle^{-\frac{3}{8}+200\delta} \Vert u \Vert_X^3 
\end{align*}

For \eqref{decompositionhb2-de-b2sigma0varphi}, we separate into three pieces depending on the localisations of $\overline{\eta}, \overline{\sigma}$. First, if $\overline{\eta}$ and $\overline{\sigma}$ are localised by $1 - m_{\widehat{\mathcal{P}}}$, we can use that
\begin{align*}
\varphi &= O\left( \overline{\sigma} \overline{\eta} \right) = O(\sigma_0 \eta_0) 
\end{align*}
and we recover a term of the form \eqref{decompositionhb2-de-b2bon}. Then, if $\overline{\sigma}$ is localised by $m_{\widehat{\mathcal{R}}}+m_{\widehat{\mathcal{P}}}$, we use that
\begin{align*}
\widehat{X}_b(\overline{\eta}) &= \sum_{\alpha = a, b, c} \mu_{\alpha}(\overline{\xi}, \overline{\eta}) m_{\alpha}(\overline{\sigma}) \widehat{X}_{\alpha}(\overline{\sigma}) 
\end{align*}
since in this case, $m_{\alpha}(\overline{\sigma}) \simeq 1$ for every $\alpha$. If finally $\overline{\sigma}$ is localised by $m_{\widehat{\mathcal{C}}}+m_{\widehat{\mathcal{L}}}$ and $\overline{\eta}$ by $m_{\widehat{\mathcal{P}}}$, this means that $\left( X_a(\overline{\eta}), X_c(\overline{\eta}), X_a(\overline{\sigma}) \right)$ is a local basis, and we may decompose 
\begin{align*}
\widehat{X}_b(\overline{\eta}) &= \mu_a(\overline{\xi}, \overline{\eta}) \widehat{X}_a(\overline{\sigma}) + \mu_a'(\overline{\xi}, \overline{\eta}) \widehat{X}_a(\overline{\eta}) + \mu_c'(\overline{\xi}, \overline{\eta}) \widehat{X}_c(\overline{\eta}) 
\end{align*}
Thus, in the two last cases, we can decompose as
\begin{align*}
\widehat{X}_b(\overline{\eta}) &= \mu_b(\overline{\xi}, \overline{\eta}) m_b(\overline{\sigma}) \widehat{X}_b(\overline{\sigma}) + \sum_{\alpha = a, c} \left( \mu_{\alpha}(\overline{\xi}, \overline{\eta}) m_{\alpha}(\overline{\sigma}) \widehat{X}_{\alpha}(\overline{\sigma}) + \mu_{\alpha}'(\overline{\xi}, \overline{\eta}) m_{\alpha}(\overline{\eta}) \widehat{X}_{\alpha}(\overline{\eta}) \right) 
\end{align*}
Since $\varphi = O(\overline{\eta})$, the contribution in \eqref{decompositionhb2-de-b2sigma0varphi} of terms having $m_{\alpha} X_{\alpha}(\overline{\eta})$, $\alpha = a, c$, is of the form \eqref{decompositionhb2-de-acb}. On every other contributions, we apply an integration by parts: 
\begin{subequations}
\begin{align}
&\int_0^t \int e^{i s \varphi} \mu(\overline{\xi}, \overline{\eta}) \frac{\varphi^2}{(|\overline{\eta}| + |\overline{\sigma}|)^4} \mu_{\alpha}(\overline{\xi}, \overline{\eta}) m_{\alpha}(\overline{\sigma}) \widehat{X}_{\alpha}(\overline{\sigma}) \cdot \nabla_{\overline{\eta}} \widehat{h}_b(s, \overline{\eta}) \widehat{f}(s, \overline{\sigma}) ~ d\overline{\eta} ds \notag \\
&\quad = - \int_0^t \int e^{i s \varphi} \mu(\overline{\xi}, \overline{\eta}) \frac{\varphi^2}{(|\overline{\eta}| + |\overline{\sigma}|)^4} \mu_{\alpha}(\overline{\xi}, \overline{\eta}) \widehat{h}_b(s, \overline{\eta}) m_{\alpha}(\overline{\sigma}) \widehat{X}_{\alpha}(\overline{\sigma}) \cdot \nabla_{\overline{\eta}} \widehat{f}(s, \overline{\sigma}) ~ d\overline{\eta} ds \label{decompositionhb2-de-b2sigma0varphi-1} \\
&\quad \quad + \int_0^t \int e^{i s \varphi} \mu(\overline{\xi}, \overline{\eta}) s \frac{\varphi^2}{(|\overline{\eta}|+|\overline{\sigma}|)^2} \mu_{\alpha}(\overline{\xi}, \overline{\eta}) \widehat{h}_b(s, \overline{\eta}) \widehat{f}(s, \overline{\sigma}) ~ d\overline{\eta} ds \label{decompositionhb2-de-b2sigma0varphi-2} \\
&\quad \quad - \int_0^t \int e^{i s \varphi} \nabla_{\overline{\eta}} \cdot \left( \mu_{\alpha}(\overline{\xi}, \overline{\eta}) m_{\alpha}(\overline{\sigma}) \widehat{X}_{\alpha}(\overline{\sigma}) \mu(\overline{\xi}, \overline{\eta}) \frac{\varphi^2}{(|\overline{\eta}| + |\overline{\sigma}|)^4} \right) \widehat{h}_b(s, \overline{\eta}) \widehat{f}(s, \overline{\sigma}) ~ d\overline{\eta} ds \label{decompositionhb2-de-b2sigma0varphi-3}
\end{align}
\end{subequations} 
where the symbol $\mu$ may change from line to line. \eqref{decompositionhb2-de-b2sigma0varphi-1} can be estimated like $\eqref{decompositionhb2-ee-bbeta0varphi}+\eqref{decompositionhb2-ee-bbsigma0varphi}$, up to replacing $m_b X_b$ by $m_{\alpha} X_{\alpha}$, for which every estimate can be done the same way (up to integrating in time and using the integrated $L^4$ estimate). \eqref{decompositionhb2-de-b2sigma0varphi-3} is of the form $\eqref{decompositionhb2-eders-b}+\eqref{decompositionhb2-eders-ac}$. Finally, on \eqref{decompositionhb2-de-b2sigma0varphi-2}, we can apply an integration by parts in time: 
\begin{subequations}
\begin{align}
\eqref{decompositionhb2-de-b2sigma0varphi-2} &= - \int_0^t \int e^{i s \varphi} \mu(\overline{\xi}, \overline{\eta}) \frac{\varphi}{(|\overline{\eta}|+|\overline{\sigma}|)^2} \mu_{\alpha}(\overline{\xi}, \overline{\eta}) \widehat{h}_b(s, \overline{\eta}) \widehat{f}(s, \overline{\sigma}) ~ d\overline{\eta} ds \label{decompositionhb2-de-b2sigma0varphi-2-1} \\
&\quad - \int_0^t \int e^{i s \varphi} \mu(\overline{\xi}, \overline{\eta}) s \frac{\varphi}{(|\overline{\eta}|+|\overline{\sigma}|)^2} \mu_{\alpha}(\overline{\xi}, \overline{\eta}) \partial_s \left( \widehat{h}_b(s, \overline{\eta}) \widehat{f}(s, \overline{\sigma}) \right) ~ d\overline{\eta} ds \label{decompositionhb2-de-b2sigma0varphi-2-2} \\
&\quad + \int e^{i t \varphi} \mu(\overline{\xi}, \overline{\eta}) t \frac{\varphi}{(|\overline{\eta}|+|\overline{\sigma}|)^2} \mu_{\alpha}(\overline{\xi}, \overline{\eta}) \widehat{h}_b(t, \overline{\eta}) \widehat{f}(t, \overline{\sigma}) ~ d\overline{\eta} \label{decompositionhb2-de-b2sigma0varphi-2-3}
\end{align}
\end{subequations} 
Using that $\varphi = O(\eta_0 \overline{\sigma}) + O(\sigma_0 \overline{\eta})$, \eqref{decompositionhb2-de-b2sigma0varphi-2-1} is of the form $\eqref{decompositionhb2-eders-b}+\eqref{decompositionhb2-eders-ac}$, \eqref{decompositionhb2-de-b2sigma0varphi-2-2} is of the form $\eqref{decompositionhb2-etr-ac} + \eqref{decompositionhb2-et-varphi}$ and \eqref{decompositionhb2-de-b2sigma0varphi-2-3} is of the form \eqref{decompositionhb2-bord-bvarphi}. 

Finally, for \eqref{decompositionhb2-de-b2eta0varphi}, if $\overline{\sigma}$ is localised by $1 - m_{\widehat{\mathcal{P}}}$, we can write $\varphi = O(\overline{\sigma}) = O(\sigma_0)$ and we get back to \eqref{decompositionhb2-de-b2bon}. If however $\overline{\sigma}$ is localised by $m_{\widehat{\mathcal{P}}}$, we can apply an integration by parts in frequency: 
\begin{subequations}
\begin{align}
&\int_0^t \int e^{i s \varphi} \mu(\overline{\xi}, \overline{\eta}) \frac{\eta_0 \varphi m_{\widehat{\mathcal{P}}}(\overline{\sigma})}{(|\overline{\eta}| + |\overline{\sigma}|)^2} m_b(\overline{\eta}) \widehat{X}_b(\overline{\eta}) \cdot \nabla_{\overline{\eta}} \widehat{h}_b(s, \overline{\eta}) \widehat{f}(s, \overline{\sigma}) ~ d\overline{\eta} ds \notag \\
&\quad = - \int_0^t \int e^{i s \varphi} \mu(\overline{\xi}, \overline{\eta}) \frac{\eta_0 \varphi m_{\widehat{\mathcal{P}}}(\overline{\sigma})}{(|\overline{\eta}| + |\overline{\sigma}|)^2} m_b(\overline{\eta}) \widehat{h}_b(s, \overline{\eta}) \widehat{X}_b(\overline{\eta}) \cdot \nabla_{\overline{\eta}} \widehat{f}(s, \overline{\sigma}) ~ d\overline{\eta} ds \label{decompositionhb2-de-b2eta0varphi-1} \\
&\quad \quad - \int_0^t \int e^{i s \varphi} \mu(\overline{\xi}, \overline{\eta}) s \eta_0 \varphi \widehat{h}_b(s, \overline{\eta}) \widehat{f}(s, \overline{\sigma}) ~ d\overline{\eta} ds \label{decompositionhb2-de-b2eta0varphi-2} \\
&\quad \quad - \int_0^t \int e^{i s \varphi} \nabla_{\overline{\eta}} \cdot \left( \widehat{X}_b(\overline{\eta}) \mu(\overline{\xi}, \overline{\eta}) \frac{\eta_0 \varphi m_{\widehat{\mathcal{P}}}(\overline{\sigma})}{(|\overline{\eta}| + |\overline{\sigma}|)^2} m_b(\overline{\eta}) \right) \widehat{h}_b(s, \overline{\eta}) \widehat{f}(s, \overline{\sigma}) ~ d\overline{\eta} ds \label{decompositionhb2-de-b2eta0varphi-3}
\end{align}
\end{subequations} 
\eqref{decompositionhb2-de-b2eta0varphi-1} is of the form \eqref{decompositionhb2-ee-bbeta0varphi} (up to replacing $m_b X_b$ by $m_{\alpha} X_{\alpha}$, as for \eqref{decompositionhb2-de-b2sigma0varphi-1}), and \eqref{decompositionhb2-de-b2eta0varphi-3} is of the form $\eqref{decompositionhb2-eders-b}+\eqref{decompositionhb2-eders-ac}$. Finally, for \eqref{decompositionhb2-de-b2eta0varphi-2}, we apply an integration by parts in time, and in the computations below we use the fact that the symbol $\mu$ is the same at each line: 
\begin{subequations}
\begin{align}
\eqref{decompositionhb2-de-b2eta0varphi-2} &= - \int_0^t \int e^{i s \varphi} \mu(\overline{\xi}, \overline{\eta}) \eta_0 \widehat{h}_b(s, \overline{\eta}) \widehat{f}(s, \overline{\sigma}) ~ d\overline{\eta} ds \label{decompositionhb2-de-b2eta0varphi-2-1} \\
&\quad - \int_0^t \int e^{i s \varphi} \mu(\overline{\xi}, \overline{\eta}) s \eta_0 \widehat{h}_b(s, \overline{\eta}) \partial_s \widehat{f}(s, \overline{\sigma}) ~ d\overline{\eta} ds \label{decompositionhb2-de-b2eta0varphi-2-2} \\
&\quad - \int_0^t \int e^{i s \varphi} \mu(\overline{\xi}, \overline{\eta}) s \eta_0 \partial_s \widehat{h}_b(s, \overline{\eta}) \widehat{f}(s, \overline{\sigma}) ~ d\overline{\eta} ds \label{decompositionhb2-de-b2eta0varphi-2-3} \\
&\quad + \int e^{i t \varphi} \mu(\overline{\xi}, \overline{\eta}) t \eta_0 \widehat{h}_b(t, \overline{\eta}) \widehat{f}(t, \overline{\sigma}) ~ d\overline{\eta} \label{decompositionhb2-de-b2eta0varphi-2-4} 
\end{align}
\end{subequations}
\eqref{decompositionhb2-de-b2eta0varphi-2-1} is of the form \eqref{decompositionhb2-eders-b}, \eqref{decompositionhb2-de-b2eta0varphi-2-2} of the form \eqref{decompositionhb2-et-b}. Finally, for \eqref{decompositionhb2-de-b2eta0varphi-2-3} and \eqref{decompositionhb2-de-b2eta0varphi-2-4}, we can add them respectively to \eqref{decompositionhb2-etr-b} and \eqref{decompositionhb2-bord-b}, noting that this changes the symbol $m_g$ given by the decomposition lemma, but that this symbol remains the same in \eqref{decompositionhb2-etr-b} and \eqref{decompositionhb2-bord-b}. 

\paragraph{Symmetric estimate} It only remains \eqref{decompositionhb2-de-sym} and \eqref{decompositionhb2-etr-b}. 

Let us first define 
\begin{align*}
g_{bb}(t) &= \mathcal{F}^{-1} \Biggl( \eqref{decompositionhb2-bord-b}+\eqref{decompositionhb2-bord-ac}+\eqref{decompositionhb2-bord-bbon}+\eqref{decompositionhb2-bord-bbonbis}+\eqref{decompositionhb2-bord-bvarphi}+
\eqref{decompositionhb2-bord-ders}+\eqref{decompositionhb2-et-varphi-1-5}\\
&\quad \quad +\eqref{decompositionhb2-ee-bbsigma0varphi-1-4} \Biggl) 
\quad + g_{bb, 3}(t) 
\end{align*}
that is every boundary term obtained in the argument (not the initial data at $t = 0$), and we define $h_{bb}(t)$ as the rest of the terms. 

On the one hand, all previous estimates and Lemma \ref{lemmedecompositioncubiquehquadbb} ensure that 
\begin{align*}
\Vert g_{bb}(t) \Vert_{L^2} &\lesssim \langle t \rangle^{\frac{1}{24}+301\delta} \Vert u \Vert_X^2 \\
\Vert (1 - m_{\widehat{\mathcal{R}}}(D)) g_{bb}(t) \Vert_{L^2} &\lesssim \Vert u \Vert_X^2 
\end{align*}
assuming $\Vert u \Vert_X \leq 1$. 

On the other hand, to estimate $h_{bb}(t)$ in $L^2$, we write by Duhamel's formula 
\begin{align*}
\Vert h_{bb}(t) \Vert_{L^2}^2 &= \Vert h_{bb}(0) \Vert_{L^2}^2 + 2 \int_0^t \int h_{bb}(s, x, y) ~ \partial_s h_{bb}(s, x, y) ~ dx dy ds 
\end{align*}
and then develop $\partial_s h_{bb}(s)$. We then have 
\begin{align*}
&\int_0^t \int h_{bb}(s, x, y) ~ \partial_s h_{bb}(s, x, y) ~ dx dy ds \\
&\leq \int_0^t \int h_{bb}(s, x, y) ~ \partial_s \mathcal{F}^{-1} \left( \eqref{decompositionhb2-de-sym} + \eqref{decompositionhb2-etr-b} \right)(s, x, y) ~ dx dy ds \\
&\quad + \int_0^t \Vert h_{bb}(s) \Vert_{L^2} \Vert \partial_s h_{bb, 3}(s) \Vert_{L^2} ~ ds + \int_0^t \Vert h_{bb}(s) \Vert_{L^2} \Vert \partial_s \left( \widehat{h}_{bb}(s) - \eqref{decompositionhb2-de-sym} - \eqref{decompositionhb2-etr-b} - \widehat{h}_{bb, 3}(s) \right) \Vert_{L^2} ~ ds \\
&\lesssim \int_0^t \int h_{bb}(s, x, y) ~ \partial_s \mathcal{F}^{-1} \left( \eqref{decompositionhb2-de-sym} + \eqref{decompositionhb2-etr-b} \right)(s, x, y) ~ dx dy ds \\
&\quad + \langle t \rangle^{\frac{1}{12}+301\delta} \Vert u \Vert_X \int_0^t \Vert \partial_s h_{bb, 3}(s) \Vert_{L^2} ~ ds + \Vert u \Vert_X^3 \\
&\lesssim \int_0^t \int h_{bb}(s, x, y) ~ \partial_s \mathcal{F}^{-1} \left( \eqref{decompositionhb2-de-sym} + \eqref{decompositionhb2-etr-b} \right)(s, x, y) ~ dx dy ds ~ + \langle t \rangle^{\frac{1}{12}+602\delta} \Vert u \Vert_X^3
\end{align*}
using that, for all the terms estimated before, they satisfy some $L^1_t L^2_x$ estimate for their time derivative. The same way, if we localise by $1-m_{\widehat{\mathcal{R}}}$, 
\begin{align*}
&\Vert (1 - m_{\widehat{\mathcal{R}}}(D)) h_{bb}(t) \Vert_{L^2}^2 = \Vert (1 - m_{\widehat{\mathcal{R}}}(D)) h_{bb}(0) \Vert_{L^2}^2 \\
&\quad \quad + 2 \int_0^t \int (1 - m_{\widehat{\mathcal{R}}}(D)) h_{bb}(s, x, y) ~ \partial_s (1 - m_{\widehat{\mathcal{R}}}(D)) h_{bb}(s, x, y) ~ dx dy ds \\
&\lesssim \Vert (1 + x^2 + |y|^2) u_0 \Vert_{H^1}^2 \\
&\quad + \int_0^t \int (1 - m_{\widehat{\mathcal{R}}}(D)) h_{bb}(s, x, y) ~ \partial_s (1 - m_{\widehat{\mathcal{R}}}(D)) \mathcal{F}^{-1} \left( \eqref{decompositionhb2-de-sym} + \eqref{decompositionhb2-etr-b} \right)(s, x, y) ~ dx dy ds ~ + \Vert u \Vert_X^3 
\end{align*}
by a similar computation. 

We can then develop \eqref{decompositionhb2-de-sym}, for $m$ a (real-valued, symmetric with respect to $0$) Hörmander-Mikhlin symbol: 
\begin{subequations}
\begin{align}
\eqref{decompositionhb2-de-sym} &= 2 \int_0^t \int e^{i s \varphi} i \eta_0 \widehat{h}_{bb}(s, \overline{\eta}) \widehat{f}(s, \overline{\sigma}) ~ d\overline{\eta} ds \label{decompositionhb2-de-sym-1} \\
&\quad + 2 \int_0^t \int e^{i s \varphi} i \eta_0 \widehat{g}_{bb}(s, \overline{\eta}) \widehat{f}(s, \overline{\sigma}) ~ d\overline{\eta} ds \label{decompositionhb2-de-sym-2} 
\end{align}  
\end{subequations}
Therefore, for $m \equiv 1$ or $m = 1 - m_{\widehat{\mathcal{R}}}$, 
\begin{align*}
&\int_0^t \int m(D) h_{bb}(s, x, y) ~ \partial_s m(D) \mathcal{F}^{-1} \eqref{decompositionhb2-de-sym-1} (s, x, y) ~ dx dy ds \\
&\quad = 2 \int_0^t \int \int m(\overline{\xi})^2 \widehat{h}_{bb}(s, -\overline{\xi}) ~ e^{i s \varphi} i \eta_0 \widehat{h}_{bb}(s, \overline{\eta}) \widehat{f}(s, \overline{\sigma}) ~ d\overline{\eta} d\overline{\xi} ds \\
&\quad = - i \int_0^t \int \int \left( m(\overline{\eta})^2 \xi_0 - m(\overline{\xi})^2 \eta_0 \right) \widehat{h}_{bb}(s, -\overline{\xi}) ~ e^{i s \varphi} \widehat{h}_{bb}(s, \overline{\eta}) \widehat{f}(s, \overline{\sigma}) ~ d\overline{\eta} d\overline{\xi} ds
\end{align*}
using the symmetry $\mathfrak{s} : (\overline{\xi}, \overline{\eta}) \mapsto (-\overline{\eta}, -\overline{\xi})$. 

If $m \equiv 1$, $m(\overline{\eta})^2 \xi_0 - m(\overline{\xi})^2 \eta_0 = \sigma_0$ and we recover an estimate similar to the one of \eqref{decompositionhb2-de-b2bon}. If $m = 1 - m_{\widehat{\mathcal{R}}}$, then since $m$ only depends on $\frac{|\xi|^2}{\xi_0^2 + |\xi|^2}$, we can write 
\begin{align*}
m(\overline{\xi}) &= M\left( \frac{|\xi|^2}{\xi_0^2 + |\xi|^2} \right)  
\end{align*}
for some smooth function $M$. Then, 
\begin{align*}
m(\overline{\xi}) - m(\overline{\eta}) &= M\left( \frac{|\xi|^2}{\xi_0^2 + |\xi|^2} \right) - M\left( \frac{|\eta|^2}{\eta_0^2 + |\eta|^2} \right) \\
&= \int_0^1 \left( \frac{|\xi|^2}{\xi_0^2 + |\xi|^2} - \frac{|\eta|^2}{\eta_0^2 + |\eta|^2} \right) M'\left( (1-\tau) \frac{|\eta|^2}{\eta_0^2 + |\eta|^2} + \tau \frac{|\xi|^2}{\xi_0^2 + |\xi|^2} \right) ~ d\tau \\
&= O(\sigma_0) + O(|\xi|^2 - |\eta|^2) 
\end{align*}
Since on the other hand $\varphi = O(\sigma_0) + \eta_0 (|\xi|^2 - |\eta|^2)$, we deduce that 
\begin{align*}
m(\overline{\eta})^2 \xi_0 - m(\overline{\xi})^2 \eta_0 &= O\left( m(\overline{\xi}) \sigma_0 \right) + O\left( m(\overline{\eta}) \sigma_0 \right) + O\left( m(\overline{\xi}) \varphi \right) + O\left( m(\overline{\eta}) \varphi \right) 
\end{align*}
The term in $O(\sigma_0 m(\overline{\xi}))$ can be estimated like \eqref{decompositionhb2-de-b2bon}, the term in $O(\sigma_0 m(\overline{\eta}))$ is symmetric. 

It then only remains for $\mu$ an order $0$ symbol: 
\begin{subequations}
\begin{align}
&\int_0^t \int \int \mu(\overline{\xi}, \overline{\eta}) m(\overline{\xi}) \frac{\varphi}{|\overline{\eta}|^2 + |\overline{\sigma}|^2} e^{i s \varphi} \widehat{h}_{bb}(s, -\overline{\xi}) \widehat{h}_{bb}(s, \overline{\eta}) \widehat{f}(s, \overline{\sigma}) ~ d\overline{\eta} d\overline{\xi} ds \notag \\
&\quad = \int_0^t \int \int \mu(\overline{\xi}, \overline{\eta}) m(\overline{\xi}) \frac{1}{|\overline{\eta}|^2 + |\overline{\sigma}|^2} e^{i s \varphi} \partial_s \widehat{h}_{bb}(s, -\overline{\xi}) \widehat{h}_{bb}(s, \overline{\eta}) \widehat{f}(s, \overline{\sigma}) ~ d\overline{\eta} d\overline{\xi} ds \label{decompositionhb2-de-sym-1-1} \\
&\quad \quad + \int_0^t \int \int \mu(\overline{\xi}, \overline{\eta}) m(\overline{\xi}) \frac{1}{|\overline{\eta}|^2 + |\overline{\sigma}|^2} e^{i s \varphi} \widehat{h}_{bb}(s, -\overline{\xi}) \partial_s \widehat{h}_{bb}(s, \overline{\eta}) \widehat{f}(s, \overline{\sigma}) ~ d\overline{\eta} d\overline{\xi} ds \label{decompositionhb2-de-sym-1-1bis} \\
&\quad \quad + \int_0^t \int \int \mu(\overline{\xi}, \overline{\eta}) m(\overline{\xi}) \frac{1}{|\overline{\eta}|^2 + |\overline{\sigma}|^2} e^{i s \varphi} \widehat{h}_{bb}(s, -\overline{\xi}) \widehat{h}_{bb}(s, \overline{\eta}) \partial_s \widehat{f}(s, \overline{\sigma}) ~ d\overline{\eta} d\overline{\xi} ds \label{decompositionhb2-de-sym-1-2} \\
&\quad \quad + \int \int \mu(\overline{\xi}, \overline{\eta}) m(\overline{\xi}) \frac{1}{|\overline{\eta}|^2 + |\overline{\sigma}|^2} e^{i t \varphi} \widehat{h}_{bb}(t, -\overline{\xi}) \widehat{h}_{bb}(t, \overline{\eta}) \widehat{f}(t, \overline{\sigma}) ~ d\overline{\eta} d\overline{\xi} \label{decompositionhb2-de-sym-1-3} \\
&\quad \quad + \int \int \mu(\overline{\xi}, \overline{\eta}) m(\overline{\xi}) \frac{1}{|\overline{\eta}|^2 + |\overline{\sigma}|^2} \widehat{h}_{bb}(0, -\overline{\xi}) \widehat{h}_{bb}(0, \overline{\eta}) \widehat{f}(0, \overline{\sigma}) ~ d\overline{\eta} d\overline{\xi} \label{decompositionhb2-de-sym-1-4} 
\end{align}
\end{subequations} 
(where $\mu$ can change from line to line). We then estimate: 
\begin{align*}
\eqref{decompositionhb2-de-sym-1-2} &\lesssim \int_0^t \Vert (1 - m_{\widehat{\mathcal{R}}}(D)) h_{bb}(s) \Vert_{L^2} \Vert e^{-i s \omega(D)} |\nabla|^{-1} h_{bb}(s) \Vert_{L^6} \Vert e^{-i s \omega(D)} |\nabla|^{-1} \partial_s f(s) \Vert_{L^3} ~ ds \\
&\lesssim \int_0^t \Vert u \Vert_X^2 \Vert h_{bb}(s) \Vert_{L^2} \Vert u(s) \Vert_{L^6}^2 ~ ds \lesssim \int_0^t t^{-\frac{5}{6}} \langle t \rangle^{-\frac{1}{4}+\frac{1}{24}+501\delta} \Vert u \Vert_X^4 ~ ds \\
&\lesssim \Vert u \Vert_X^4 \\
\eqref{decompositionhb2-de-sym-1-3} &\lesssim \Vert (1 - m_{\widehat{\mathcal{R}}}(D)) h_{bb}(t) \Vert_{L^2} \Vert e^{-i t \omega(D)} |\nabla|^{-1} h_{bb}(t) \Vert_{L^6} \Vert u(t) \Vert_{L^3} \lesssim \Vert u \Vert_X^3 
\end{align*}
\eqref{decompositionhb2-de-sym-1-4} can be estimated like \eqref{decompositionhb2-de-sym-1-3}. Finally, for \eqref{decompositionhb2-de-sym-1-1} and \eqref{decompositionhb2-de-sym-1-1bis}, we develop again $\partial_s h_{bb}(s)$. On every term except \eqref{decompositionhb2-de-sym} and \eqref{decompositionhb2-etr-b}, we have an estimate of the form 
\begin{align*}
\int_0^t \langle s \rangle^{-\frac{1}{12}-500\delta} \Vert \partial_s A(s) \Vert_{L^2} ~ ds \lesssim \Vert u \Vert_X^2
\end{align*}
so we can estimate: 
\begin{align*}
&\int_0^t \int \int \mu(\overline{\xi}, \overline{\eta}) \frac{1}{|\overline{\eta}|^2 + |\overline{\sigma}|^2} e^{i s \varphi} \partial_s \widehat{A}(s, -\overline{\xi}) \widehat{h}_{bb}(s, \overline{\eta}) \widehat{f}(s, \overline{\sigma}) ~ d\overline{\eta} d\overline{\xi} ds \\
&\quad \lesssim \int_0^t \Vert \partial_s A(s) \Vert_{L^2} \Vert e^{-i s \omega(D)} |\nabla|^{-1} h_{bb}(s) \Vert_{L^6} \Vert u(s) \Vert_{L^3} ~ ds \\
&\quad \lesssim \int_0^t \langle s \rangle^{\frac{1}{24}} \Vert \partial_s A(s) \Vert_{L^2} \Vert u \Vert_X \Vert u(s) \Vert_{L^2}^{\frac{1}{3}} \Vert u(s) \Vert_{L^4}^{\frac{2}{3}} ~ ds \\
&\quad \lesssim \int_0^t \langle s \rangle^{-\frac{23}{72}+50\delta} \Vert \partial_s A(s) \Vert_{L^2} \Vert u \Vert_X^2 ~ ds \lesssim \Vert u \Vert_X^4
\end{align*}
Finally, for \eqref{decompositionhb2-de-sym} and \eqref{decompositionhb2-etr-b}, we can extend
\begin{subequations}
\begin{align}
\eqref{decompositionhb2-de-sym} &= - 2 \int_0^t \int e^{i s \varphi} \mu_{HHB}(\overline{\xi}, \overline{\eta}) \eta_0^2 m_b(\overline{\eta}) \widehat{X}_b(\overline{\eta}) \cdot \nabla_{\overline{\eta}} \widehat{h}_b(s, \overline{\eta}) \widehat{f}(s, \overline{\sigma}) ~ d\overline{\eta} ds \label{decompositionhb2-de-sym-termeHHB} \\
&\quad - 2 \int_0^t \int e^{i s \varphi} \left( 1 - \mu_{HHB}(\overline{\xi}, \overline{\eta}) \right) \eta_0^2 m_b(\overline{\eta}) \widehat{X}_b(\overline{\eta}) \cdot \nabla_{\overline{\eta}} \widehat{h}_b(s, \overline{\eta}) \widehat{f}(s, \overline{\sigma}) ~ d\overline{\eta} ds \label{decompositionhb2-de-sym-termeok} 
\end{align}
\end{subequations}
and we have 
\begin{align*}
\Vert |\overline{\xi}|^{-1} \partial_t \eqref{decompositionhb2-de-sym-termeHHB} \Vert_{L^2} &\lesssim \Vert \partial_x m_b(D) X_b h_b(t) \Vert_{L^2} \Vert u(t) \Vert_{L^{\infty}} \lesssim t^{-\frac{5}{6}} \langle t \rangle^{\frac{1}{8}+401\delta} \Vert u \Vert_X^2 \\
\Vert \partial_t \eqref{decompositionhb2-de-sym-termeok} \Vert_{L^2} &\lesssim \Vert \partial_x m_b(D) X_b h_b(t) \Vert_{L^2} \Vert \nabla u(t) \Vert_{L^{\infty}} \lesssim t^{-\frac{5}{6}} \langle t \rangle^{-\frac{1}{8}+401\delta} \Vert u \Vert_X^2
\end{align*}
Therefore, we can estimate
\begin{align*}
&\int_0^t \int \int \mu(\overline{\xi}, \overline{\eta}) \frac{1}{|\overline{\eta}|^2 + |\overline{\sigma}|^2} e^{i s \varphi} \partial_s \eqref{decompositionhb2-de-sym-termeHHB}(s, -\overline{\xi}) \widehat{h}_{bb}(s, \overline{\eta}) \widehat{f}(s, \overline{\sigma}) ~ d\overline{\eta} d\overline{\xi} ds \\
&\quad \lesssim \int_0^t \Vert |\overline{\xi}|^{-1} \partial_s \eqref{decompositionhb2-de-sym-termeHHB} \Vert_{L^2} \Vert e^{-i s \omega(D)} |\nabla|^{-1} h_{bb}(s) \Vert_{L^2} \Vert u(t) \Vert_{L^3} ~ ds \\
&\quad \lesssim \int_0^t t^{-\frac{5}{6}} \langle t \rangle^{\frac{1}{8}+\frac{1}{24}-\frac{13}{36}+401\delta} \Vert u \Vert_X^4 ~ ds \\
&\quad \lesssim \int_0^t t^{-\frac{5}{6}} \langle t \rangle^{-\frac{1}{6}-\frac{1}{36}+401\delta} \Vert u \Vert_X^4 ~ ds ~ \lesssim \Vert u \Vert_X^4 \\
&\int_0^t \int \int \mu(\overline{\xi}, \overline{\eta}) \frac{1}{|\overline{\eta}|^2 + |\overline{\sigma}|^2} e^{i s \varphi} \partial_s \eqref{decompositionhb2-de-sym-termeok}(s, -\overline{\xi}) \widehat{h}_{bb}(s, \overline{\eta}) \widehat{f}(s, \overline{\sigma}) ~ d\overline{\eta} d\overline{\xi} ds \\
&\quad \lesssim \int_0^t \Vert \partial_s \eqref{decompositionhb2-de-sym-termeok} \Vert_{L^2} \Vert e^{-i s \omega(D)} |\nabla|^{-1} h_{bb}(s) \Vert_{L^6} \Vert |\nabla|^{-1} u(s) \Vert_{L^3} ~ ds \\
&\quad \lesssim \int_0^t s^{-\frac{5}{6}} \langle s \rangle^{-\frac{1}{8}+\frac{1}{24}+401\delta} \Vert u \Vert_X^3 \Vert |\nabla|^{-1} u(s) \Vert_{L^2}^{\frac{5}{7}} \Vert u(s) \Vert_{L^4}^{\frac{2}{7}} ~ ds \\
&\quad \lesssim \int_0^t s^{-\frac{5}{6}} \langle s \rangle^{-\frac{5}{21}+401\delta} \Vert u \Vert_X^4 ~ ds ~ \lesssim \Vert u \Vert_X^4 
\end{align*}
We proceed the same way for the contribution of \eqref{decompositionhb2-etr-b}, which is simpler. 

It then remains \eqref{decompositionhb2-de-sym-2} and \eqref{decompositionhb2-etr-b}. In \eqref{decompositionhb2-de-sym-2}, we develop $g_{bb}$ and we can estimate every term except \eqref{decompositionhb2-bord-b} and \eqref{decompositionhb2-bord-bvarphi}: 
\begin{align*}
&\begin{aligned}
&\int_0^t \langle s \rangle^{-\frac{1}{2}-\frac{1}{48}} \Vert e^{-i s \omega(D)} \nabla \mathcal{F}^{-1} \eqref{decompositionhb2-bord-ac} \Vert_{L^4} ~ ds \lesssim \int_0^t s \langle s \rangle^{-\frac{1}{2}-\frac{1}{48}} \sum_{\alpha = a, c} \Vert e^{-i s \omega(D)} \nabla h_{\alpha}(s) \Vert_{L^4} \Vert \nabla u(s) \Vert_{L^{\infty}} ~ ds \\
&\quad \lesssim \sum_{\alpha = a, c} \left( \int_0^t s^{\frac{63}{64}} \langle s \rangle^{\frac{1}{50}} \Vert e^{-i s \omega(D)} \nabla h_{\alpha}(s) \Vert_{L^4}^4 ~ ds \right)^{\frac{1}{4}} \left( \int_0^t s^{1+\frac{1}{192}} \langle s \rangle^{-\frac{1}{150}-\frac{2}{3}-\frac{1}{36}} \Vert \nabla u(s) \Vert_{L^{\infty}}^{\frac{4}{3}} ~ ds \right)^{\frac{3}{4}} \\
&\quad \lesssim \Vert u \Vert_X \left( \int_0^t s^{-\frac{1}{9}+\frac{1}{192}} \langle s \rangle^{-\frac{1}{150}-\frac{2}{3}-\frac{1}{36}-\frac{2}{9}+\frac{400\delta}{3}} \Vert u \Vert_X^{\frac{4}{3}} ~ ds \right)^{\frac{3}{4}} \lesssim \Vert u \Vert_X^2 
\end{aligned} \\
&\begin{aligned} 
\Vert e^{-i t \omega(D)} \nabla \mathcal{F}^{-1} \eqref{decompositionhb2-bord-bbon} \Vert_{L^4} &\lesssim t \Vert e^{-i t \omega(D)} \partial_x h_b(t) \Vert_{L^4} \Vert \nabla u(t) \Vert_{L^{\infty}} \lesssim t^{-\frac{1}{4}} \langle t \rangle^{-\frac{1}{4}+100\delta} \Vert u \Vert_X^2 
\end{aligned} \\
&\begin{aligned} 
\Vert e^{-i t \omega(D)} \partial_x \mathcal{F}^{-1} \eqref{decompositionhb2-bord-bbonbis} \Vert_{L^4} &\lesssim t \Vert e^{-i t \omega(D)} \partial_x h_b(t) \Vert_{L^4} \Vert \nabla u(t) \Vert_{L^{\infty}} \lesssim t^{-\frac{1}{4}} \langle t \rangle^{-\frac{1}{4}+100\delta} \Vert u \Vert_X^2 
\end{aligned} \\
&\begin{aligned} 
\Vert e^{-i t \omega(D)} \nabla \mathcal{F}^{-1} \eqref{decompositionhb2-bord-ders} \Vert_{L^4} &\lesssim t \Vert u(t) \Vert_{L^4} \Vert \nabla u(t) \Vert_{L^{\infty}} \lesssim t^{-\frac{1}{4}} \langle t \rangle^{-\frac{3}{8}+150\delta} \Vert u \Vert_X^2 
\end{aligned} \\
&\begin{aligned} 
\Vert e^{-i t \omega(D)} \nabla \mathcal{F}^{-1} \eqref{decompositionhb2-et-varphi-1-5} \Vert_{L^4} &\lesssim t \Vert e^{-i t \omega(D)} |\nabla|^{-1} h_b(t) \Vert_{L^6} \Vert \nabla u(t) \Vert_{L^{\infty}} \Vert u(t) \Vert_{L^{12}} \\
&\lesssim t^{\frac{1}{6}} \langle t \rangle^{-\frac{1}{6}+100\delta} \Vert h_b(t) \Vert_{L^2} \Vert u \Vert_X \Vert |\nabla|^{\frac{1}{2}} u(t) \Vert_{L^4} \\
&\lesssim t^{-\frac{1}{4}} \langle t \rangle^{-\frac{7}{24}+150\delta} \Vert u \Vert_X^3 
\end{aligned} \\
&\begin{aligned} 
\Vert e^{-i t \omega(D)} \nabla \mathcal{F}^{-1} \eqref{decompositionhb2-ee-bbsigma0varphi-1-4} \Vert_{L^4} &\lesssim \Vert e^{-i t \omega(D)} \partial_x h_b(t) \Vert_{L^4} \Vert e^{-i t \omega(D)} |\nabla|^{-1} h_b(t) \Vert_{L^{\infty}} \\
&\lesssim t^{-\frac{5}{12}} \langle t \rangle^{-\frac{1}{12}+O(\delta)} \Vert u \Vert_X^2 
\end{aligned} \\
&\begin{aligned} 
\Vert e^{-i t \omega(D)} \nabla g_{bb, 3}(t) \Vert_{L^4} &\lesssim t^2 \Vert \partial_x u(t) \Vert_{L^{\infty}} \Vert \nabla u(t) \Vert_{L^{\infty}} \Vert u(t) \Vert_{L^4} \lesssim t^{-\frac{1}{12}} \langle t \rangle^{-\frac{13}{24}+250\delta} \Vert u \Vert_X^3
\end{aligned} 
\end{align*}
Therefore, for any of these terms, denoting them by $A(s)$, we can estimate 
\begin{align*}
&\int_0^t \Vert \partial_s \left[ - 2 \int_0^s \int e^{i s' \varphi} \eta_0 \widehat{A}(s', \overline{\eta}) \widehat{f}(s', \overline{\sigma}) ~ d\overline{\eta} ds' \right] \Vert_{L^2} ~ ds \\
&\quad \lesssim \int_0^t \Vert e^{-i s \omega(D)} \partial_x A(s) \Vert_{L^4} \Vert u(s) \Vert_{L^4} ~ ds \lesssim \Vert u \Vert_X^3 
\end{align*}

For \eqref{decompositionhb2-bord-bvarphi}, we develop
\begin{subequations}
\begin{align}
&2 i \int_0^t \int \int \widehat{h}_{bb}(s, -\overline{\xi}) ~ e^{i s \varphi} \eta_0 \eqref{decompositionhb2-bord-bvarphi}(s, \overline{\eta}) \widehat{f}(s, \overline{\sigma}) ~ d\overline{\eta} ds \notag \\
&\quad = \int_0^t \int \int \int \widehat{h}_{bb}(s, -\overline{\xi}) ~ e^{i s \varphi_3(\overline{\xi}, \overline{\eta}, \overline{\sigma})} s (\eta_0+\sigma_0) \mu(\overline{\xi}, \overline{\eta}, \overline{\sigma}) \frac{\varphi(\overline{\eta}+\overline{\sigma}, \overline{\eta})}{|\overline{\eta}|^2+|\overline{\sigma}|^2} \widehat{h}_b(s, \overline{\eta}) \widehat{f}(s, \overline{\sigma}) \widehat{f}(s, \overline{\rho}) ~ d\overline{\eta} ds \notag \\
&\quad = \int_0^t \int \int \int \widehat{h}_{bb}(s, -\overline{\xi}) ~ e^{i s \varphi_3(\overline{\xi}, \overline{\eta}, \overline{\sigma})} s \mu(\overline{\xi}, \overline{\eta}, \overline{\sigma}) \eta_0 \widehat{h}_b(s, \overline{\eta}) \widehat{f}(s, \overline{\sigma}) |\overline{\rho}| \widehat{f}(s, \overline{\rho}) ~ d\overline{\eta} ds \label{decompositionhb2-symb-bordbvarphi-1} \\
&\quad \quad + \int_0^t \int \int \int \widehat{h}_{bb}(s, -\overline{\xi}) ~ e^{i s \varphi_3(\overline{\xi}, \overline{\eta}, \overline{\sigma})} s \mu(\overline{\xi}, \overline{\eta}, \overline{\sigma}) \widehat{h}_b(s, \overline{\eta}) \sigma_0 \widehat{f}(s, \overline{\sigma}) |\overline{\rho}| \widehat{f}(s, \overline{\rho}) ~ d\overline{\eta} ds \label{decompositionhb2-symb-bordbvarphi-2} \\
&\quad \quad \begin{aligned}
+ \int_0^t \int \int \int \widehat{h}_{bb}(s, -\overline{\xi}) ~ e^{i s \varphi_3(\overline{\xi}, \overline{\eta}, \overline{\sigma})} s (\eta_0+\sigma_0) \mu(\overline{\xi}, \overline{\eta}, \overline{\sigma}) \frac{\varphi_3(\overline{\xi}, \overline{\eta}, \overline{\sigma})}{|\overline{\eta}|^2+|\overline{\sigma}|^2+|\overline{\rho}|^2} \\
\widehat{h}_b(s, \overline{\eta}) \widehat{f}(s, \overline{\sigma}) \widehat{f}(s, \overline{\rho}) ~ d\overline{\eta} ds 
\end{aligned} \label{decompositionhb2-symb-bordbvarphi-3}
\end{align}
\end{subequations}
for some symbol $\mu$ that can change from line to line. Above, we used that 
\begin{align*}
\frac{\varphi(\overline{\eta}+\overline{\sigma}, \overline{\eta})}{|\overline{\eta}|^2 + |\overline{\sigma}|^2} &= O\left( \frac{\varphi_3(\overline{\xi}, \overline{\eta}, \overline{\sigma})}{|\overline{\eta}|^2+|\overline{\sigma}|^2+|\overline{\rho}|^2} \right) + O(\overline{\rho}) 
\end{align*}
decomposing depending on whether $|\overline{\rho}| \gtrsim |\overline{\eta}|+|\overline{\sigma}|$ or $|\overline{\rho}| \ll |\overline{\eta}|+|\overline{\sigma}|$. We can already estimate: 
\begin{align*}
\eqref{decompositionhb2-symb-bordbvarphi-1} &\lesssim \int_0^t \Vert h_{bb}(s) \Vert_{L^2} s \Vert e^{-i s \omega(D)} \partial_x h_b(s) \Vert_{L^4} \Vert u(s) \Vert_{L^4} \Vert \nabla u(s) \Vert_{L^{\infty}} ~ ds \\
&\lesssim \int_0^t s^{-\frac{2}{3}} \langle s \rangle^{-\frac{1}{3}-\frac{1}{1024}+451\delta} \Vert u \Vert_X^4 ~ ds \lesssim \Vert u \Vert_X^4 \\
\eqref{decompositionhb2-symb-bordbvarphi-2} &\lesssim \int_0^t \Vert h_{bb}(s) \Vert_{L^2} s \Vert h_b(s) \Vert_{L^2} \Vert \partial_x u(s) \Vert_{L^{\infty}} \Vert \nabla u(s) \Vert_{L^{\infty}} ~ ds \\
&\lesssim \int_0^t s^{-\frac{2}{3}} \langle s \rangle^{-\frac{3}{8}+501\delta} \Vert u \Vert_X^4 ~ ds ~ \lesssim \Vert u \Vert_X^4 
\end{align*}
On \eqref{decompositionhb2-symb-bordbvarphi-3}, we can apply an integration by parts in time: {\footnotesize 
\begin{subequations}
\begin{align}
\eqref{decompositionhb2-symb-bordbvarphi-3} &= \int_0^t \int \int \int \widehat{h}_{bb}(s, -\overline{\xi}) ~ e^{i s \varphi_3(\overline{\xi}, \overline{\eta}, \overline{\sigma})} \mu(\overline{\xi}, \overline{\eta}, \overline{\sigma}) \frac{1}{|\overline{\eta}|+|\overline{\sigma}|+|\overline{\rho}|} \widehat{h}_b(s, \overline{\eta}) \widehat{f}(s, \overline{\sigma}) \widehat{f}(s, \overline{\rho}) ~ d\overline{\eta} d\overline{\sigma} ds \label{decompositionhb2-symb-bordbvarphi-3-1} \\
&\quad + \int_0^t \int \int \int \partial_s \widehat{h}_{bb}(s, -\overline{\xi}) ~ e^{i s \varphi_3(\overline{\xi}, \overline{\eta}, \overline{\sigma})} s \mu(\overline{\xi}, \overline{\eta}, \overline{\sigma}) \frac{1}{|\overline{\eta}|+|\overline{\sigma}|+|\overline{\rho}|} \widehat{h}_b(s, \overline{\eta}) \widehat{f}(s, \overline{\sigma}) \widehat{f}(s, \overline{\rho}) ~ d\overline{\eta} d\overline{\sigma} ds \label{decompositionhb2-symb-bordbvarphi-3-2} \\
&\quad + \int_0^t \int \int \int \widehat{h}_{bb}(s, -\overline{\xi}) ~ e^{i s \varphi_3(\overline{\xi}, \overline{\eta}, \overline{\sigma})} s \mu(\overline{\xi}, \overline{\eta}, \overline{\sigma}) \frac{1}{|\overline{\eta}|+|\overline{\sigma}|+|\overline{\rho}|} \partial_s \widehat{h}_b(s, \overline{\eta}) \widehat{f}(s, \overline{\sigma}) \widehat{f}(s, \overline{\rho}) ~ d\overline{\eta} d\overline{\sigma} ds \label{decompositionhb2-symb-bordbvarphi-3-3} \\
&\quad + \int_0^t \int \int \int \widehat{h}_{bb}(s, -\overline{\xi}) ~ e^{i s \varphi_3(\overline{\xi}, \overline{\eta}, \overline{\sigma})} s \mu(\overline{\xi}, \overline{\eta}, \overline{\sigma}) \frac{1}{|\overline{\eta}|+|\overline{\sigma}|+|\overline{\rho}|} \widehat{h}_b(s, \overline{\eta}) \partial_s \left( \widehat{f}(s, \overline{\sigma}) \widehat{f}(s, \overline{\rho}) \right) ~ d\overline{\eta} d\overline{\sigma} ds \label{decompositionhb2-symb-bordbvarphi-3-4} \\
&\quad + \int \int \int \widehat{h}_{bb}(t, -\overline{\xi}) ~ e^{i t \varphi_3(\overline{\xi}, \overline{\eta}, \overline{\sigma})} t \mu(\overline{\xi}, \overline{\eta}, \overline{\sigma}) \frac{1}{|\overline{\eta}|+|\overline{\sigma}|+|\overline{\rho}|} \widehat{h}_b(t, \overline{\eta}) \widehat{f}(t, \overline{\sigma}) \widehat{f}(t, \overline{\rho}) ~ d\overline{\eta} d\overline{\sigma} \label{decompositionhb2-symb-bordbvarphi-3-5}
\end{align}
\end{subequations} }
where $\mu$ can change from line to line. Then: 
\begin{align*}
\eqref{decompositionhb2-symb-bordbvarphi-3-1} &\lesssim \int_0^t \Vert h_{bb}(s) \Vert_{L^2} \Vert e^{-i s \omega(D)} |\nabla|^{-1} h_b(s) \Vert_{L^6} \Vert u(s) \Vert_{L^6}^2 ~ ds \\
&\lesssim \int_0^t s^{-\frac{5}{6}} \langle s \rangle^{-\frac{5}{24}+501\delta} \Vert u \Vert_X^4 ~ ds ~ \lesssim \Vert u \Vert_X^4 \\
\eqref{decompositionhb2-symb-bordbvarphi-3-3} &\lesssim \int_0^t \Vert h_{bb}(s) \Vert_{L^2} s \Vert e^{-i s \omega(D)} |\nabla|^{-1} \partial_s h_b(s) \Vert_{L^6} \Vert u(s) \Vert_{L^6}^2 ~ ds \\
&\lesssim \int_0^t \langle s \rangle^{-\frac{1}{24}+501\delta} \Vert u \Vert_X^3 \Vert \partial_s h_b(s) \Vert_{L^2} ~ ds \lesssim \Vert u \Vert_X^5 \\
\eqref{decompositionhb2-symb-bordbvarphi-3-4} &\lesssim \int_0^t \Vert h_{bb}(s) \Vert_{L^2} s \Vert e^{-i s \omega(D)} |\nabla|^{-1} h_b(s) \Vert_{L^6} \Vert e^{-i s \omega(D)} \partial_s f(s) \Vert_{L^6} \Vert u(s) \Vert_{L^6} ~ ds \\
&\lesssim \int_0^t s \langle s \rangle^{\frac{1}{24}+301\delta} \Vert u \Vert_X \Vert h_b(s) \Vert_{L^2} \Vert \partial_x u(s) \Vert_{L^{\infty}} \Vert u(s) \Vert_{L^6}^2 ~ ds \\
&\lesssim \int_0^t s^{-\frac{2}{3}} \langle s \rangle^{-\frac{11}{24}+601\delta} \Vert u \Vert_X^5 ~ ds ~ \lesssim \Vert u \Vert_X^5 
\end{align*}
Finally, for \eqref{decompositionhb2-symb-bordbvarphi-3-2}, we can reuse the computation done for \eqref{decompositionhb2-de-sym-1-1} and \eqref{decompositionhb2-de-sym-1-1bis}, noting that 
\begin{align*}
&\Vert \langle \overline{\xi} \rangle \int \int e^{i s \varphi_3(\overline{\xi}, \overline{\eta}, \overline{\sigma})} s \mu(\overline{\xi}, \overline{\eta}, \overline{\sigma}) \frac{1}{|\overline{\eta}|+|\overline{\sigma}|+|\overline{\rho}|} \widehat{h}_b(s, \overline{\eta}) \widehat{f}(s, \overline{\sigma}) \widehat{f}(s, \overline{\rho}) ~ d\overline{\eta} d\overline{\sigma} \Vert_{L^2_{\overline{\xi}}} \\
&\quad \lesssim s \Vert e^{-i s \omega(D)} \langle \nabla \rangle |\nabla|^{-1} h_b(s) \Vert_{L^6} \Vert u(s) \Vert_{L^6}^2 \lesssim s^{\frac{1}{6}} \langle s \rangle^{-\frac{1}{4}} \Vert u \Vert_X^3 
\end{align*}
which allows to apply the same integrations by parts (and even simpler). 

Before treating the contribution of \eqref{decompositionhb2-bord-b} in \eqref{decompositionhb2-de-sym-2}, we estimate a part of \eqref{decompositionhb2-etr-b}. By developing $\partial_s h_b(s)$ using Lemma \ref{lemstructurehalphagalpha}, we can first estimate: 
\begin{align*}
\Vert e^{-i t \omega(D)} \partial_x \mathcal{F}^{-1} \partial_t \eqref{lemstructurehb-03-resxetaac} \Vert_{L^4} 
&\lesssim \sum_{\alpha = a, c} \Vert e^{-i t \omega(D)} \nabla h_{\alpha}(t) \Vert_{L^4} \Vert \partial_x u(t) \Vert_{L^{\infty}} \\
&\lesssim t^{-\frac{5}{4}} \langle t \rangle^{-\frac{1}{3}+O(\delta)} \Vert u \Vert_X^2 \\
\Vert e^{-i t \omega(D)} \partial_x \mathcal{F}^{-1} \partial_t \left( \eqref{lemstructurehb-04-resxetabbon} + \eqref{lemstructurehb-05-resxetabbonbis} \right) \Vert_{L^4}
&\lesssim \Vert e^{-i t \omega(D)} \partial_x h_b(t) \Vert_{L^4} \Vert \partial_x u(t) \Vert_{L^{\infty}} \\
&\lesssim t^{-\frac{5}{4}} \langle t \rangle^{-\frac{1}{3}+O(\delta)} \Vert u \Vert_X^2 \\
\Vert e^{-i t \omega(D)} \partial_x \mathcal{F}^{-1} \partial_t \eqref{lemstructurehb-11-dersymb} \Vert_{L^4} 
&\lesssim \Vert u(t) \Vert_{L^4} \Vert \partial_x u(t) \Vert_{L^{\infty}} \\
&\lesssim t^{-\frac{5}{4}} \langle t \rangle^{-\frac{3}{8}+150\delta} \Vert u \Vert_X^2 \\
\Vert e^{-i t \omega(D)} \partial_x \mathcal{F}^{-1} \partial_t \eqref{lemstructurehb-13-termecubique} \Vert_{L^4}
&\lesssim t \Vert \partial_x u(t) \Vert_{L^{\infty}} \Vert u(t) \Vert_{L^{\infty}} \Vert u(t) \Vert_{L^4} \\
&\lesssim t^{-\frac{13}{12}} \langle t \rangle^{-\frac{13}{24}+250\delta} \Vert u \Vert_X^3
\end{align*}
Therefore, denoting by $A(t)$ one of those terms, we can estimate: 
\begin{align*}
&\int_0^t \int \int \widehat{h}_{bb}(s, -\overline{\xi}) ~ e^{i s \varphi} s m_g(\overline{\xi}, \overline{\eta}) \eta_0 \partial_s \widehat{A}(s, \overline{\eta}) \widehat{f}(s, \overline{\sigma}) ~ d\overline{\eta} ds \\
&\quad \lesssim \int_0^t \Vert h_{bb}(s) \Vert_{L^2} s \Vert e^{-i s \omega(D)} \partial_x \partial_s A(s) \Vert_{L^4} \Vert u(s) \Vert_{L^4} ~ ds \\
&\quad \lesssim \int_0^t s^{-\frac{2}{3}} \langle s \rangle^{-\frac{5}{12}+1000\delta} \Vert u \Vert_X^4 ~ ds ~ \lesssim \Vert u \Vert_X^4
\end{align*}
For the contribution of \eqref{lemstructurehb-06-resxetabphi}, we develop and get the term: 
\begin{align*}
&\int_0^t \int \int \widehat{h}_{bb}(s, -\overline{\xi}) ~ e^{i s \varphi_3(\overline{\xi}, \overline{\eta}, \overline{\sigma})} s \mu(\overline{\xi}, \overline{\eta}, \overline{\sigma}) (\eta_0+\sigma_0) \frac{\varphi(\overline{\eta}+\overline{\sigma}, \overline{\eta})}{(|\overline{\eta}| + |\overline{\sigma}|)^2} \widehat{h}_b(s, \overline{\eta}) \widehat{f}(s, \overline{\sigma}) \widehat{f}(s, \overline{\rho}) ~ d\overline{\eta} d\overline{\sigma} d\overline{\xi} ds 
\end{align*}
But this term is identical to the contribution of \eqref{decompositionhb2-bord-bvarphi} in \eqref{decompositionhb2-de-sym-2}, and can therefore be estimated the same way. 

Finally, we group the contribution of \eqref{decompositionhb2-bord-b} in \eqref{decompositionhb2-de-sym-2}, and the one of \eqref{lemstructurehb-02-symeta} in \eqref{decompositionhb2-etr-b}, using the fact that the symbol $m_g$ is the same: 
\begin{align}
\int_0^t \int \int e^{i s \varphi_3(\overline{\xi}, \overline{\eta}, \overline{\sigma})} 2 i s \eta_0 (\eta_0+\sigma_0) \left( m_g(\overline{\xi}, \overline{\eta}+\overline{\sigma}) - m_g(\overline{\eta}+\overline{\sigma}, \overline{\eta}) \right) \widehat{h}_b(s, \overline{\eta}) \widehat{f}(s, \overline{\sigma}) \widehat{f}(s, \overline{\rho}) ~ d\overline{\eta} d\overline{\sigma} ds \label{decompositionhb2-de-sym-2-bordbcompensation} 
\end{align}
But by Taylor's formula, we have that 
\begin{align*}
&\eta_0 (\eta_0+\sigma_0) \left( m_g(\overline{\xi}, \overline{\eta}+\overline{\sigma}) - m_g(\overline{\eta}+\overline{\sigma}, \overline{\eta}) \right) \\
&= O(\eta_0 \overline{\sigma}) + O(\eta_0 \overline{\rho}) 
+ \eta_0^2 \chi\left( 2^{10} \frac{|\overline{\sigma}|+|\overline{\rho}|}{|\overline{\eta}|} \right) \left( m_g(\overline{\xi}, \overline{\eta}+\overline{\sigma}) - m_g(\overline{\eta}+\overline{\sigma}, \overline{\eta}) \right) \\
&= O(\eta_0 \overline{\sigma}) + O(\eta_0 \overline{\rho}) 
\end{align*}
Therefore, we can estimate: 
\begin{align*}
\Vert \partial_t \eqref{decompositionhb2-de-sym-2-bordbcompensation} \Vert_{L^2} &\lesssim t \Vert e^{-i t \omega(D)} \partial_x h_b(t) \Vert_{L^4} \Vert \nabla u(t) \Vert_{L^{\infty}} \Vert u(t) \Vert_{L^4} \\
&\lesssim t^{-\frac{2}{3}} \langle t \rangle^{-\frac{3}{8}+150\delta} \Vert u \Vert_X^3 
\end{align*}
which is integrable in time. 

This concludes the proof of Proposition \ref{propestimeeaprioriquadhbb}. 

\subsection{Proof of Lemma \ref{lemmedecompositioncubiquehquadbb}}

It only remains to prove Lemma \ref{lemmedecompositioncubiquehquadbb} to conclude for the quadratic weighted estimate. As in the proof of the decomposition lemma \ref{lemdecompositionFGH}, we decompose $h_{b, 3}(t)$ depending on the directions of the interacting frequencies. However, the main contribution is the direction, since the form of $h_{b, 3}(t)$ ensures that derivatives are already easier to distribute. Thus, we only separate the cases depending on the directions (and only in certain specific cases depending on their relative sizes). 

We will denote by $I_{\mu}^{A_1 A_2 A_3}[F_1, F_2, F_3]$ the interaction 
\begin{align*}
&\widehat{I}_{\mu}^{A_1 A_2 A_3}[F_1, F_2, F_3](t, \overline{\xi}) \\
&\quad = \int_0^t \int \int e^{i s \varphi_3(\overline{\xi}, \overline{\eta}, \overline{\sigma})} \mu(\overline{\xi}, \overline{\eta}, \overline{\sigma}) m_{A_1}(\overline{\eta}) \widehat{F}_1(s, \overline{\eta}) m_{A_2}(\overline{\sigma}) \widehat{F}_2(s, \overline{\sigma}) m_{A_3}(\overline{\rho}) \widehat{F}_3(s, \overline{\rho}) ~ d\overline{\eta} d\overline{\sigma} ds
\end{align*}
for $A_1, A_2, A_3 \in \{ \widehat{\mathcal{R}}, \widehat{\mathcal{C}}, \widehat{\mathcal{L}}, \widehat{\mathcal{P}} \}$. In each of the following subsections, we will consider objects of the form 
\begin{align*}
\xi_0 m_b(\overline{\xi}) \widehat{X}_b(\overline{\xi}) \cdot \nabla_{\overline{\xi}} \widehat{I}_{s \mu}^{A_1 A_2 A_3}[F_1, F_2, F_3](t, \overline{\xi})
\end{align*}
where the $F_i$ are all equal to $f$, except one of them which is $\partial_x f$, and $\mu$ is a symbol of order $0$ that allows to distribute derivatives like in the form from Lemma \ref{lemstructurehalphagalpha}, that is we can avoid that if, say, $F_1 = \partial_x f$, we can avoid using the structure of $\mu$ that a second derivative falls on $F_1$. 

We will treat in this symmetric way $F_1, F_2, F_3$ in order to symmetrize the various geometric areas and the interactions, which reduces the number of cases to consider to only 20. 

We will also denote by $\mu_{BBBB}, \mu_{HHHB}, \mu_{HHBH}, ...$ symbols localising the frequencies depending on their relative sizes, that will be useful in certain specific cases. 

\subsubsection{Generic estimates}

The following lemma gives estimates we will use repeatedly: 

\begin{Lem} Let $\mu$ be a symbol allowing for Hölder-type estimates. Consider the following quantities: 
\begin{subequations} \label{qteslemestimeesgeneriquesdeccubiquehb2} 
\begin{align}
&\int_0^t \int \int e^{i s \varphi_3(\overline{\xi}, \overline{\eta}, \overline{\sigma})} s \mu(\overline{\xi}, \overline{\eta}, \overline{\sigma}) |\overline{\eta}| \widehat{f}(s, \overline{\eta}) \widehat{f}(s, \overline{\sigma}) |\overline{\rho}| m_{\alpha}(\overline{\rho}) \widehat{X}_{\alpha}(\overline{\rho}) \cdot \nabla_{\overline{\xi}} \widehat{f}(s, \overline{\rho}) ~ d\overline{\eta} d\overline{\sigma} ds \label{lemestimeesgeneriquesdeccubiquehb2-eta0rhoresxrho} \\
&\int_0^t \int \int e^{i s \varphi_3(\overline{\xi}, \overline{\eta}, \overline{\sigma})} s \mu(\overline{\xi}, \overline{\eta}, \overline{\sigma}) \widehat{f}(s, \overline{\eta}) |\overline{\sigma}| \widehat{f}(s, \overline{\sigma}) \rho_0 m_{\beta}(\overline{\rho}) \widehat{X}_{\beta}(\overline{\rho}) \cdot \nabla_{\overline{\xi}} \widehat{f}(s, \overline{\rho}) ~ d\overline{\eta} d\overline{\sigma} ds \label{lemestimeesgeneriquesdeccubiquehb2-rho0sigmaresxrho} \\
&\int_0^t \int \int e^{i s \varphi_3(\overline{\xi}, \overline{\eta}, \overline{\sigma})} s \mu(\overline{\xi}, \overline{\eta}, \overline{\sigma}) \eta_0 \widehat{f}(s, \overline{\eta}) |\overline{\sigma}| \widehat{f}(s, \overline{\sigma}) m_{\beta}(\overline{\rho}) \widehat{X}_{\beta}(\overline{\rho}) \cdot \nabla_{\overline{\xi}} \widehat{f}(s, \overline{\rho}) ~ d\overline{\eta} d\overline{\sigma} ds \label{lemestimeesgeneriquesdeccubiquehb2-eta0sigmaresxrho} \\
&\int_0^t \int \int e^{i s \varphi_3(\overline{\xi}, \overline{\eta}, \overline{\sigma})} s \mu(\overline{\xi}, \overline{\eta}, \overline{\sigma}) \eta_0 \widehat{f}(s, \overline{\eta}) \sigma_0 \widehat{f}(s, \overline{\sigma}) \left( 1 - m_{\widehat{\mathcal{C}}}(\overline{\rho}) \right) \nabla_{\overline{\eta}} \widehat{f}(s, \overline{\rho}) ~ d\overline{\eta} d\overline{\sigma} ds \label{lemestimeesgeneriquesdeccubiquehb2-eta0sigma0resxdegrho} \\
&\int_0^t \int \int e^{i s \varphi_3(\overline{\xi}, \overline{\eta}, \overline{\sigma})} s^2 \mu(\overline{\xi}, \overline{\eta}, \overline{\sigma}) \partial_s \widehat{f}(s, \overline{\eta}) |\overline{\sigma}| \widehat{f}(s, \overline{\sigma}) \widehat{f}(s, \overline{\rho}) ~ d\overline{\eta} d\overline{\sigma} ds \label{lemestimeesgeneriquesdeccubiquehb2-sigmadseta} \\
&\int_0^t \int \int e^{i s \varphi_3(\overline{\xi}, \overline{\eta}, \overline{\sigma})} s \mu(\overline{\xi}, \overline{\eta}, \overline{\sigma}) |\overline{\eta}| \widehat{f}(s, \overline{\eta}) \widehat{f}(s, \overline{\sigma}) \widehat{f}(s, \overline{\rho}) ~ d\overline{\eta} d\overline{\sigma} ds \label{lemestimeesgeneriquesdeccubiquehb2-dersymb} \\
&\int_0^t \int \int e^{i s \varphi_3(\overline{\xi}, \overline{\eta}, \overline{\sigma})} s \mu(\overline{\xi}, \overline{\eta}, \overline{\sigma}) |\overline{\eta}| \widehat{f}(s, \overline{\eta}) |\overline{\sigma}|^{-1} \widehat{f}(s, \overline{\sigma}) |\overline{\rho}| \widehat{f}(s, \overline{\rho}) ~ d\overline{\eta} d\overline{\sigma} ds \label{lemestimeesgeneriquesdeccubiquehb2-dersymbbis} 
\end{align}
\end{subequations}
where $\alpha \in \{ a, c \}$ and $\beta \in \{ a, b, c \}$ are arbitrary. For every $\gamma$, we have 
\begin{align*}
\int_0^t \Vert \partial_s (\ref{qteslemestimeesgeneriquesdeccubiquehb2}.\gamma) \Vert_{L^2} ~ ds &\lesssim \Vert u \Vert_X^3 \left( 1 + \Vert u \Vert_X \right) 
\end{align*}
Furthermore, for every term of the form 
\begin{align}
t^2 \int \int e^{i s \varphi_3(\overline{\xi}, \overline{\eta}, \overline{\sigma})} \mu(\overline{\xi}, \overline{\eta}, \overline{\sigma}) \frac{\eta_0 |\overline{\sigma}|}{|\overline{\eta}|+|\overline{\sigma}|+|\overline{\rho}|} \widehat{f}(s, \overline{\eta}) \widehat{f}(s, \overline{\sigma}) \widehat{f}(s, \overline{\rho}) ~ d\overline{\eta} d\overline{\sigma} ds \label{lemestimeesgeneriquesdeccubiquehb2-g}
\end{align}
we have 
\begin{align*}
\Vert \eqref{lemestimeesgeneriquesdeccubiquehb2-g} \Vert_{L^2} &\lesssim \Vert u \Vert_X^3 \\
\Vert e^{-i t \omega(D)} \nabla \eqref{lemestimeesgeneriquesdeccubiquehb2-g} \Vert_{L^4} &\lesssim t^{-\frac{5}{12}} \langle t \rangle^{-\frac{1}{8}+O(\delta)} \Vert u \Vert_X^3 
\end{align*}
\end{Lem}

\begin{Dem}
We can estimate directly: 
\begin{align*}
&\begin{aligned}
&\int_0^t \Vert \partial_s \eqref{lemestimeesgeneriquesdeccubiquehb2-eta0rhoresxrho} \Vert_{L^2} ~ ds \lesssim \int_0^t s \Vert \nabla u(s) \Vert_{L^{\infty}} \Vert u(s) \Vert_{L^4} \sum_{\alpha = a, c} \Vert e^{-is \omega(D)} \nabla m_{\alpha}(D) X_{\alpha} f(s) \Vert_{L^4} ~ ds \\
&\quad \lesssim \sum_{\alpha = a, c} \left( \int_0^t s \Vert \nabla u(s) \Vert_{L^{\infty}}^{\frac{4}{3}} \Vert u(s) \Vert_{L^4}^{\frac{4}{3}} ~ ds \right)^{\frac{3}{4}} \left( \int_0^t s \Vert e^{-i s \omega(D)} \nabla m_{\alpha}(D) X_{\alpha} f(s) \Vert_{L^4}^4 ~ ds \right)^{\frac{1}{4}} \\
&\quad \lesssim \left( \int_0^t s^{-\frac{2}{3}} \langle s \rangle^{-\frac{2}{9}-\frac{1}{6}+200\delta} \Vert u \Vert_X^{\frac{8}{3}} ~ ds \right)^{\frac{3}{4}} \Vert u \Vert_X \lesssim \Vert u \Vert_X^3 
\end{aligned} \\
&\begin{aligned} 
\int_0^t \Vert \partial_s \eqref{lemestimeesgeneriquesdeccubiquehb2-rho0sigmaresxrho} \Vert_{L^2} ~ ds &\lesssim \int_0^t s \Vert u(s) \Vert_{L^4} \Vert \nabla u(s) \Vert_{L^{\infty}} \sum_{\beta = a, b, c} \Vert e^{-is \omega(D)} \partial_x m_{\beta}(D) X_{\beta} f(s) \Vert_{L^4} ~ ds \\
&\lesssim \Vert u \Vert_X^3 
\end{aligned} \\
&\begin{aligned} 
\Vert \partial_t \eqref{lemestimeesgeneriquesdeccubiquehb2-eta0sigmaresxrho} \Vert_{L^2} &\lesssim t \Vert \partial_x u(t) \Vert_{L^{\infty}} \Vert \nabla u(t) \Vert_{L^{\infty}} \sum_{\beta = a, b, c} \Vert m_{\beta}(D) X_{\beta} f(t) \Vert_{L^2} \\
&\lesssim t^{-\frac{2}{3}} \langle t \rangle^{-\frac{5}{12}+200\delta} \Vert u \Vert_X^3 
\end{aligned} \\
&\begin{aligned} 
\Vert \partial_t \eqref{lemestimeesgeneriquesdeccubiquehb2-eta0sigma0resxdegrho} \Vert_{L^2} &\lesssim t \Vert \partial_x u(t) \Vert_{L^{\infty}}^2 \Vert \left( 1 - m_{\widehat{\mathcal{C}}}(D) \right) (x, y) f(t) \Vert_{L^2} \lesssim t^{-\frac{2}{3}} \langle t \rangle^{-\frac{19}{48}+402\delta} \Vert u \Vert_X^3 
\end{aligned} \\
&\begin{aligned} 
\Vert \partial_t \eqref{lemestimeesgeneriquesdeccubiquehb2-sigmadseta} \Vert_{L^2} &\lesssim t^2 \Vert e^{-i t \omega(D)} \partial_t f(t) \Vert_{L^4} \Vert \nabla u(t) \Vert_{L^{\infty}} \Vert u(t) \Vert_{L^4} \lesssim t^{-\frac{1}{2}} \langle t \rangle^{-\frac{2}{3}+300\delta} \Vert u \Vert_X^4 
\end{aligned} \\
&\begin{aligned} 
\Vert \partial_t \eqref{lemestimeesgeneriquesdeccubiquehb2-dersymb} \Vert_{L^2} &\lesssim t \Vert \nabla u(t) \Vert_{L^{\infty}} \Vert u(t) \Vert_{L^4}^2 \lesssim t^{-\frac{2}{3}} \langle t \rangle^{-\frac{5}{12}+200\delta} \Vert u \Vert_X^3 
\end{aligned} \\
&\begin{aligned} 
\Vert \partial_t \eqref{lemestimeesgeneriquesdeccubiquehb2-dersymbbis} \Vert_{L^2} &\lesssim t \Vert \nabla u(t) \Vert_{L^{\infty}}^2 \Vert |\nabla|^{-1} u(t) \Vert_{L^2} \lesssim t^{-\frac{2}{3}} \langle t \rangle^{-\frac{1}{3}+200\delta} \Vert u \Vert_X^3 
\end{aligned} \\
&\begin{aligned} 
\Vert \eqref{lemestimeesgeneriquesdeccubiquehb2-g} \Vert_{L^2} &\lesssim t^2 \Vert \nabla u(t) \Vert_{L^{\infty}} \Vert u(t) \Vert_{L^4}^2 \lesssim t^{\frac{1}{3}} \langle t \rangle^{-\frac{5}{12}+200\delta} \Vert u \Vert_X^3 
\end{aligned} \\
&\begin{aligned} 
\Vert e^{-i t \omega(D)} \nabla \eqref{lemestimeesgeneriquesdeccubiquehb2-g} \Vert_{L^2} &\lesssim t^2 \Vert \nabla u(t) \Vert_{L^{\infty}}^2 \Vert u(t) \Vert_{L^4} \lesssim t^{-\frac{1}{12}} \langle t \rangle^{-\frac{11}{24}+250\delta} \Vert u \Vert_X^3 
\end{aligned} 
\end{align*}
In all the cases where we proved pointwise in $t$ estimates, these are integrable in $t$, so it concludes. 
\end{Dem}

\subsubsection{Interaction \texorpdfstring{$\widehat{\mathcal{R}}\widehat{\mathcal{R}}\widehat{\mathcal{R}}$}{RRR}}

Let us consider: {\footnotesize 
\begin{subequations}
\begin{align}
&\xi_0 m_b(\overline{\xi}) \widehat{X}_b(\overline{\xi}) \cdot \nabla_{\overline{\xi}} \widehat{I}_{s \mu}^{\widehat{\mathcal{R}}\widehat{\mathcal{R}}\widehat{\mathcal{R}}}[F_1, F_2, F_3](t, \overline{\xi}) \notag \\
&\quad
= \int_0^t \int \int i s^2 \xi_0 m_b(\overline{\xi}) \widehat{X}_b(\overline{\xi}) \cdot \nabla_{\overline{\xi}} \varphi_3 e^{i s \varphi_3} \mu(\overline{\xi}, \overline{\eta}, \overline{\sigma}) m_{\widehat{\mathcal{R}}}(\overline{\eta}) m_{\widehat{\mathcal{R}}}(\overline{\sigma}) m_{\widehat{\mathcal{R}}}(\overline{\rho}) \widehat{F}_1(s, \overline{\eta}) \widehat{F}_2(s, \overline{\sigma}) \widehat{F}_3(s, \overline{\rho}) ~ d\overline{\eta} d\overline{\sigma} ds \label{equdecchampbRRR-1} \\
&\quad \quad + \int_0^t \int \int e^{i s \varphi_3} s \xi_0 \mu(\overline{\xi}, \overline{\eta}, \overline{\sigma}) m_{\widehat{\mathcal{R}}}(\overline{\eta}) m_{\widehat{\mathcal{R}}}(\overline{\sigma}) m_{\widehat{\mathcal{R}}}(\overline{\rho}) \widehat{F}_1(s, \overline{\eta}) \widehat{F}_2(s, \overline{\sigma}) m_b(\overline{\xi}) \widehat{X}_b(\overline{\xi}) \cdot \nabla_{\overline{\xi}} \widehat{F}_3(s, \overline{\rho}) ~ d\overline{\eta} d\overline{\sigma} ds \label{equdecchampbRRR-2} \\
&\quad \quad + \int_0^t \int \int e^{i s \varphi_3} s \xi_0 m_b(\overline{\xi}) \widehat{X}_b(\overline{\xi}) \cdot \nabla_{\overline{\xi}} \left( \mu(\overline{\xi}, \overline{\eta}, \overline{\sigma}) m_{\widehat{\mathcal{R}}}(\overline{\eta}) m_{\widehat{\mathcal{R}}}(\overline{\sigma}) m_{\widehat{\mathcal{R}}}(\overline{\rho}) \right) \widehat{F}_1(s, \overline{\eta}) \widehat{F}_2(s, \overline{\sigma}) \widehat{F}_3(s, \overline{\rho}) ~ d\overline{\eta} d\overline{\sigma} ds \label{equdecchampbRRR-3} 
\end{align}
\end{subequations} }
\eqref{equdecchampbRRR-3} is of the form $\eqref{lemestimeesgeneriquesdeccubiquehb2-dersymb}+\eqref{lemestimeesgeneriquesdeccubiquehb2-dersymbbis}$, and \eqref{equdecchampbRRR-2} is of the form $\eqref{lemestimeesgeneriquesdeccubiquehb2-eta0sigmaresxrho}+\eqref{lemestimeesgeneriquesdeccubiquehb2-rho0sigmaresxrho}$. 

Finally, for \eqref{equdecchampbRRR-1}, we can apply a partition of unity on the sphere $S^8$ to decompose it into a finite number of neighborhoods such that, either $(\varphi_3, \nabla_{\overline{\eta}, \overline{\sigma}} \varphi_3) \neq 0$ on the considered neighborhood, either this neighborhood is centered around a point such that $(\varphi_3, \nabla_{\overline{\eta}, \overline{\sigma}} \varphi_3) = 0$. 

In the first case, we have locally
\begin{align*}
1 = O(\varphi_3) + O(\nabla_{\overline{\eta}, \overline{\sigma}} \varphi_3) 
\end{align*}

In the second case, by Lemma \ref{lemcalculresonancescubiquesgen}, up to interchanging $\overline{\eta}, \overline{\sigma}, \overline{\rho}$, this means that at the center of the neighborhood we have $\overline{\eta} + \overline{\sigma} = 0$, $\overline{\rho} = \overline{\xi}$. 

Note now that $\nabla_{\overline{\eta}} \varphi_3 - \nabla_{\overline{\sigma}} \varphi_3$ vanishes on $\{ \overline{\eta} + \overline{\sigma} = 0 \}$ and satisfies
\begin{align*}
\nabla_{\overline{\eta}} \left( \nabla_{\overline{\eta}} - \nabla_{\overline{\sigma}} \right)^T \varphi_3 &= \begin{pmatrix}
6 \eta_0 & 2 \eta^T \\ 2 \eta & 2 \eta_0 I_2 
\end{pmatrix}
\end{align*}
whose determinant is $8 \eta_0 (3 \eta_0^2 - |\eta|^2)$ which does not vanish since we assumed $\overline{\eta}$ localised by $m_{\widehat{\mathcal{R}}}$. In particular, this means that $\{ \nabla_{\overline{\eta}} \varphi_3 = \nabla_{\overline{\sigma}} \varphi_3 \}$ defines locally a $5$-dimensional submanifold of $S^8$, that coincides with $\{ \overline{\eta} + \overline{\sigma} = 0 \} \cap S^8$. Furthermore, every smooth function vanishing on $\{ \overline{\eta} + \overline{\sigma} = 0 \} \cap S^8$ belongs to the ideal of smooth functions generated by the components of $\nabla_{\overline{\eta}} \varphi_3 - \nabla_{\overline{\sigma}} \varphi_3$. 

But precisely $\nabla_{\overline{\xi}} \varphi_3$ vanishes on $\{ \overline{\eta} + \overline{\sigma} = 0 \} = \{ \overline{\xi} = \overline{\rho} \}$, so we can factorize it by $\nabla_{\overline{\eta}} \varphi_3 - \nabla_{\overline{\sigma}} \varphi_3$. 

Therefore, in any case, we have that 
\begin{align*}
\nabla_{\overline{\xi}} \varphi &= O(\varphi_3) + O(\nabla_{\overline{\eta}} \varphi_3) + O(\nabla_{\overline{\sigma}} \varphi_3)
\end{align*}
We can then apply integrations by parts on \eqref{equdecchampbRRR-1}: {\footnotesize 
\begin{align*}
\eqref{equdecchampbRRR-1} &= \int_0^t \int \int \xi_0 e^{i s \varphi_3} s \mu(\overline{\xi}, \overline{\eta}, \overline{\sigma}) m_{\widehat{\mathcal{R}}}(\overline{\eta}) m_{\widehat{\mathcal{R}}}(\overline{\sigma}) m_{\widehat{\mathcal{R}}}(\overline{\rho}) \nabla_{\overline{\eta}, \overline{\sigma}} \left( \widehat{F}_1(s, \overline{\eta}) \widehat{F}_2(s, \overline{\sigma}) \widehat{F}_3(s, \overline{\rho}) \right) ~ d\overline{\eta} d\overline{\sigma} ds \\
&\quad + \int_0^t \int \int e^{i s \varphi_3} s (|\overline{\eta}|+|\overline{\sigma}|+|\overline{\rho}|) \mu(\overline{\xi}, \overline{\eta}, \overline{\sigma}) m_{\widehat{\mathcal{R}}}(\overline{\eta}) m_{\widehat{\mathcal{R}}}(\overline{\sigma}) m_{\widehat{\mathcal{R}}}(\overline{\rho}) \widehat{F}_1(s, \overline{\eta}) \widehat{F}_2(s, \overline{\sigma}) \widehat{F}_3(s, \overline{\rho}) ~ d\overline{\eta} d\overline{\sigma} ds \\
&\quad + \int_0^t \int \int s^2 e^{i s \varphi_3} \mu(\overline{\xi}, \overline{\eta}, \overline{\sigma}) m_{\widehat{\mathcal{R}}}(\overline{\eta}) m_{\widehat{\mathcal{R}}}(\overline{\sigma}) m_{\widehat{\mathcal{R}}}(\overline{\rho}) \partial_s \left( \widehat{F}_1(s, \overline{\eta}) \widehat{F}_2(s, \overline{\sigma}) \widehat{F}_3(s, \overline{\rho}) \right) ~ d\overline{\eta} d\overline{\sigma} ds \\
&\quad + \int \int t^2 e^{i t \varphi_3} \mu(\overline{\xi}, \overline{\eta}, \overline{\sigma}) m_{\widehat{\mathcal{R}}}(\overline{\eta}) m_{\widehat{\mathcal{R}}}(\overline{\sigma}) m_{\widehat{\mathcal{R}}}(\overline{\rho}) \widehat{F}_1(t, \overline{\eta}) \widehat{F}_2(t, \overline{\sigma}) \widehat{F}_3(t, \overline{\rho}) ~ d\overline{\eta} d\overline{\sigma} 
\end{align*} }
where the symbol $\mu$ can change from line to line as long as it keeps similar properties. We recognize then terms of the form \eqref{qteslemestimeesgeneriquesdeccubiquehb2} or \eqref{lemestimeesgeneriquesdeccubiquehb2-g}.  

\subsubsection{Interaction \texorpdfstring{$\widehat{\mathcal{C}}\widehat{\mathcal{R}}\widehat{\mathcal{R}}$}{CRR}}

Let us consider {\footnotesize 
\begin{subequations}
\begin{align}
&\xi_0 m_b(\overline{\xi}) \widehat{X}_b(\overline{\xi}) \cdot \nabla_{\overline{\xi}} \widehat{I}_{s \mu}^{\widehat{\mathcal{C}}\widehat{\mathcal{R}}\widehat{\mathcal{R}}}[F_1, F_2, F_3](t, \overline{\xi}) \notag \\
&\quad = \int_0^t \int \int i s^2 \xi_0 m_b(\overline{\xi}) \widehat{X}_b(\overline{\xi}) \cdot \nabla_{\overline{\xi}} \varphi_3 e^{i s \varphi_3} \mu(\overline{\xi}, \overline{\eta}, \overline{\sigma}) m_{\widehat{\mathcal{C}}}(\overline{\eta}) m_{\widehat{\mathcal{R}}}(\overline{\sigma}) m_{\widehat{\mathcal{R}}}(\overline{\rho}) \widehat{F}_1(s, \overline{\eta}) \widehat{F}_2(s, \overline{\sigma}) \widehat{F}_3(s, \overline{\rho}) ~ d\overline{\eta} d\overline{\sigma} ds \label{equdecchampbCRR-1} \\
&\quad \quad + \int_0^t \int \int e^{i s \varphi_3} s \xi_0 \mu(\overline{\xi}, \overline{\eta}, \overline{\sigma}) m_{\widehat{\mathcal{C}}}(\overline{\eta}) m_{\widehat{\mathcal{R}}}(\overline{\sigma}) m_{\widehat{\mathcal{R}}}(\overline{\rho}) \widehat{F}_1(s, \overline{\eta}) \widehat{F}_2(s, \overline{\sigma}) m_b(\overline{\xi}) \widehat{X}_b(\overline{\xi}) \cdot \nabla_{\overline{\xi}} \widehat{F}_3(s, \overline{\rho}) ~ d\overline{\eta} d\overline{\sigma} ds \label{equdecchampbCRR-2} \\
&\quad \quad + \int_0^t \int \int e^{i s \varphi_3} s \xi_0 m_b(\overline{\xi}) \widehat{X}_b(\overline{\xi}) \cdot \nabla_{\overline{\xi}} \left( \mu(\overline{\xi}, \overline{\eta}, \overline{\sigma}) m_{\widehat{\mathcal{C}}}(\overline{\eta}) m_{\widehat{\mathcal{R}}}(\overline{\sigma}) m_{\widehat{\mathcal{R}}}(\overline{\rho}) \right) \widehat{F}_1(s, \overline{\eta}) \widehat{F}_2(s, \overline{\sigma}) \widehat{F}_3(s, \overline{\rho}) ~ d\overline{\eta} d\overline{\sigma} ds \label{equdecchampbCRR-3} 
\end{align}
\end{subequations} }
\eqref{equdecchampbCRR-3} is of the form $\eqref{lemestimeesgeneriquesdeccubiquehb2-dersymb}+\eqref{lemestimeesgeneriquesdeccubiquehb2-dersymbbis}$, and \eqref{equdecchampbCRR-2} is of the form $\eqref{lemestimeesgeneriquesdeccubiquehb2-eta0sigmaresxrho}+\eqref{lemestimeesgeneriquesdeccubiquehb2-rho0sigmaresxrho}$. 

For \eqref{equdecchampbCRR-1}, as before, we can restrict our attention to the angular neighborhood of a point 
\begin{align*}
Z = (\xi_0^Z, \xi^Z, \eta_0^Z, \eta^Z, \sigma_0^Z, \sigma^Z) \in S^8
\end{align*} 
such that 
\begin{align*}
\varphi_3 = 0, \quad \widehat{X}_a(\overline{\eta}) \cdot \nabla_{\overline{\eta}} \varphi_3 = 0, \quad \widehat{X}_c(\overline{\eta}) \cdot \nabla_{\overline{\eta}} \varphi_3 = 0, \quad \nabla_{\overline{\sigma}} \varphi_3 = 0 
\end{align*}
and $\overline{\eta}^Z \in \widehat{\mathcal{C}}$. (If we are on the support of $m_{\widehat{\mathcal{C}}}$ in $\overline{\eta}$ but in the neighborhood of a point such that $\overline{\eta}^Z \notin \widehat{\mathcal{C}}$, up to choosing the neighborhood small enough we can reuse the decomposition of the case $\widehat{\mathcal{R}} \widehat{\mathcal{R}} \widehat{\mathcal{R}}$.) 

In particular, $\nabla_{\sigma_0} \varphi_3(Z) = 0$ implies that, locally $|\overline{\sigma}| \simeq |\overline{\rho}|$, while $\widehat{X}_c(\overline{\eta}) \cdot \nabla_{\overline{\eta}} \varphi_3(Z) = 0$ forces $\eta^Z$ and $\rho^Z$ to be aligned, then  
\begin{align*}
0 = \widehat{X}_a(\overline{\eta}) \cdot \nabla_{\overline{\eta}} \varphi_3(Z) &= \frac{\eta_0^Z}{|\overline{\eta}^Z|} \left( 3 (\rho_0^Z)^2 + |\rho^Z|^2 - 3 (\eta_0^Z)^2 - |\eta^Z|^2 \right) + \frac{\eta^Z}{|\overline{\eta}^Z|} \cdot (2 \rho_0^Z \rho^Z - 2 \eta_0^Z \eta^Z ) \\
&= \frac{\eta_0^Z}{|\overline{\eta}^Z|} \left( 3 (\rho_0^Z)^2 + |\rho^Z|^2 \pm 2 \sqrt{3} |\rho_0^Z| |\rho^Z| - 12 (\eta_0^Z)^2 \right) 
\end{align*}
for some choice of a sign, which corresponds to the sign of $\rho_0^Z \eta_0^Z \rho^Z \cdot \eta^Z$. Since $\overline{\rho}$ is localised by $m_{\widehat{\mathcal{R}}}$, we deduce that, locally, $|\overline{\eta}| \simeq |\overline{\rho}|$. 

Then, by Lemma \ref{lemcalculsconecoordonneesconiquesvarphi}, we have that 
\begin{align*}
\xi_t^{\overline{\eta} \overline{\rho}} &= O\left( \widehat{X}_c(\overline{\eta}) \cdot \nabla_{\overline{\eta}} \varphi_3 \right) \\
\xi_t^{\overline{\eta} \overline{\sigma}} &= O\left( \widehat{X}_c(\overline{\eta}) \cdot (\nabla_{\overline{\eta}}-\nabla_{\overline{\sigma}}) \varphi_3 \right) \\
1 &= O\left( \nabla_{\overline{\eta}} \varphi_3 \right) \\
\left( \overline{\rho}_a^{\overline{\eta}} \right)^2 - \left( \overline{\eta}_a^{\overline{\eta}} \right)^2 &= O\left( \widehat{X}_a(\overline{\eta}) \cdot \nabla_{\overline{\eta}} \varphi_3 \right) + O\left( \xi_t^{\overline{\eta} \overline{\rho}} \right) + O\left( \overline{\eta}_b^{\overline{\eta}} \right) \\
\left( \overline{\sigma}_a^{\overline{\eta}} \right)^2 - \left( \overline{\eta}_a^{\overline{\eta}} \right)^2 &= O\left( \widehat{X}_a(\overline{\eta}) \cdot (\nabla_{\overline{\eta}}-\nabla_{\overline{\sigma}}) \varphi_3 \right) + O\left( \xi_t^{\overline{\eta} \overline{\sigma}} \right) + O\left( \overline{\eta}_b^{\overline{\eta}} \right) 
\end{align*}
We also have 
\begin{align*}
\widehat{X}_b'(\overline{\eta}) \cdot \nabla_{\overline{\sigma}} \varphi_3 &= \frac{\eta_0}{|\eta_0|} \left( 3 \rho_0^2 + |\rho|^2 - 3 \sigma_0^2 - |\sigma|^2 \right) - \sqrt{3} \frac{\eta}{|\eta|} \cdot \left( 2 \rho_0 \rho - 2 \sigma_0 \sigma \right) \\
&= \frac{\eta_0}{|\eta_0|} \left( 3 \rho_0^2 + \left( \frac{\eta \cdot \rho}{|\eta|} \right)^2 - 2 \sqrt{3} \frac{\eta_0 \rho_0}{|\eta_0|} \frac{\eta \cdot \rho}{|\eta|} - 3 \sigma_0^2 - \left( \frac{\eta \cdot \sigma}{|\eta|} \right)^2 + 2 \sqrt{3} \frac{\eta_0 \sigma_0}{|\eta_0|} \frac{\eta \cdot \sigma}{|\eta|} \right) \\
&\quad + \frac{\eta_0}{|\eta_0|} \left( \left( \frac{J \eta \cdot \rho}{|\eta|} \right)^2 - \left( \frac{J \eta \cdot \sigma}{|\eta|} \right)^2 \right) \\
&= \frac{\eta_0}{|\eta_0|} \left( \left( \overline{\rho}_b^{\overline{\eta}} \right)^2 - \left( \overline{\sigma}_b^{\overline{\eta}} \right)^2 \right) + O\left( \xi_t^{\overline{\eta} \overline{\sigma}} \right) + O\left( \xi_t^{\overline{\eta} \overline{\rho}} \right) 
\end{align*}
In particular, near the considered point $Z$, there exists signs $\epsilon_1, \epsilon_2, \epsilon_3 \in \{ -1, 1 \}$ such that, locally, 
\begin{align*}
\overline{\sigma}_a^{\overline{\eta}} &= \epsilon_1 \overline{\eta}_a^{\overline{\eta}} + O\left( \widehat{X}_a(\overline{\eta}) \cdot \nabla_{\overline{\eta}} \varphi_3 \right) + O\left( \widehat{X}_c(\overline{\eta}) \cdot \nabla_{\overline{\eta}} \varphi_3 \right) + O\left( m_b(\overline{\eta}) \widehat{X}_b(\overline{\eta}) \cdot \nabla_{\overline{\eta}} \varphi_3 \right) + O\left( \nabla_{\overline{\sigma}} \varphi_3 \right) \\
\overline{\rho}_a^{\overline{\eta}} &= \epsilon_2 \overline{\eta}_a^{\overline{\eta}} + O\left( \widehat{X}_a(\overline{\eta}) \cdot \nabla_{\overline{\eta}} \varphi_3 \right) + O\left( \widehat{X}_c(\overline{\eta}) \cdot \nabla_{\overline{\eta}} \varphi_3 \right) + O\left( m_b(\overline{\eta}) \widehat{X}_b(\overline{\eta}) \cdot \nabla_{\overline{\eta}} \varphi_3 \right) + O\left( \nabla_{\overline{\sigma}} \varphi_3 \right) \\
\overline{\rho}_b^{\overline{\eta}} &= \epsilon_3 \overline{\sigma}_b^{\overline{\eta}} + O\left( \widehat{X}_a(\overline{\eta}) \cdot \nabla_{\overline{\eta}} \varphi_3 \right) + O\left( \widehat{X}_c(\overline{\eta}) \cdot \nabla_{\overline{\eta}} \varphi_3 \right) + O\left( m_b(\overline{\eta}) \widehat{X}_b(\overline{\eta}) \cdot \nabla_{\overline{\eta}} \varphi_3 \right) + O\left( \nabla_{\overline{\sigma}} \varphi_3 \right)
\end{align*}
In particular, 
\begin{align*}
&\xi_t^{\overline{\eta} \overline{\xi}} = O\left( \widehat{X}_c(\overline{\eta}) \cdot \nabla_{\overline{\eta}} \varphi_3 \right) + O\left( \nabla_{\overline{\sigma}} \varphi_3 \right) \\
&\overline{\xi}_a^{\overline{\eta}} = (1 + \epsilon_1 + \epsilon_2) \overline{\eta}_a^{\overline{\eta}} + O\left( \widehat{X}_a(\overline{\eta}) \cdot \nabla_{\overline{\eta}} \varphi_3 \right) + O\left( \widehat{X}_c(\overline{\eta}) \cdot \nabla_{\overline{\eta}} \varphi_3 \right) + O\left( m_b(\overline{\eta}) \widehat{X}_b(\overline{\eta}) \cdot \nabla_{\overline{\eta}} \varphi_3 \right) + O\left( \nabla_{\overline{\sigma}} \varphi_3 \right) \\
&\overline{\xi}_b^{\overline{\eta}} = (1 + \epsilon_3) \overline{\sigma}_b^{\overline{\eta}} + O\left( \widehat{X}_a(\overline{\eta}) \cdot \nabla_{\overline{\eta}} \varphi_3 \right) + O\left( \widehat{X}_c(\overline{\eta}) \cdot \nabla_{\overline{\eta}} \varphi_3 \right) + O\left( m_b(\overline{\eta}) \widehat{X}_b(\overline{\eta}) \cdot \nabla_{\overline{\eta}} \varphi_3 \right) + O\left( \nabla_{\overline{\sigma}} \varphi_3 \right)
\end{align*}
This implies that $|\overline{\xi}| \simeq |\overline{\eta}|$, so we can localise by $\mu_{BBBB}$. 

We can then compute by Lemma \ref{lemcalculsconecoordonneesconiquesvarphi} that
\begin{align*}
6 \sqrt{3} \frac{\eta_0}{|\eta_0|} \varphi_3 &= \left( \overline{\xi}_a^{\overline{\eta}} \right)^3 + \left( \overline{\xi}_b^{\overline{\eta}} \right)^3 - \left( \overline{\eta}_a^{\overline{\eta}} \right)^3 - \left( \overline{\eta}_b^{\overline{\eta}} \right)^3 - \left( \overline{\sigma}_a^{\overline{\eta}} \right)^3 - \left( \overline{\sigma}_b^{\overline{\eta}} \right)^3 - \left( \overline{\rho}_a^{\overline{\eta}} \right)^3 - \left( \overline{\rho}_b^{\overline{\eta}} \right)^3 \\
&\quad + O\left( \widehat{X}_c(\overline{\eta}) \cdot \nabla_{\overline{\eta}} \varphi_3 \right) + O\left( \nabla_{\overline{\sigma}} \varphi_3 \right) \\
&= 6 \left( \overline{\eta}_a^{\overline{\eta}} \right)^3 (1 + \epsilon_1) (1 + \epsilon_2) 
+ 3 \left( \overline{\sigma}_b^{\overline{\eta}} \right)^3 (1 + \epsilon_3) \\
&\quad + O\left( \widehat{X}_a(\overline{\eta}) \cdot \nabla_{\overline{\eta}} \varphi_3 \right) + O\left( \widehat{X}_c(\overline{\eta}) \cdot \nabla_{\overline{\eta}} \varphi_3 \right) + O\left( m_b(\overline{\eta}) \widehat{X}_b(\overline{\eta}) \cdot \nabla_{\overline{\eta}} \varphi_3 \right) + O\left( \nabla_{\overline{\sigma}} \varphi_3 \right)
\end{align*}
For $\varphi_3$ to vanish at $Z$, there are only to possibilities: 
\begin{itemize}
\item either $\epsilon_3 = -1$ and $(\epsilon_1, \epsilon_2) \neq (1, 1)$; 
\item or $\epsilon_3 = 1 = \epsilon_1 = \epsilon_2$ and 
\begin{align*}
\overline{\sigma}_b^{\overline{\eta}} &= - 2^{\frac{2}{3}} \overline{\eta}_a^{\overline{\eta}} 
\end{align*}
at $Z$. 
\end{itemize}

If $\epsilon_3 = -1$, then 
\begin{align*}
m_b(\overline{\xi}) &= O\left( \overline{\xi}_b^{\overline{\xi}} \right) = O\left( \overline{\xi}_b^{\overline{\eta}} \right) + O\left( \xi_t^{\overline{\eta} \overline{\xi}} \right) \\
&= O\left( \widehat{X}_a(\overline{\eta}) \cdot \nabla_{\overline{\eta}} \varphi_3 \right) + O\left( \widehat{X}_c(\overline{\eta}) \cdot \nabla_{\overline{\eta}} \varphi_3 \right) + O\left( m_b(\overline{\eta}) \widehat{X}_b(\overline{\eta}) \cdot \nabla_{\overline{\eta}} \varphi_3 \right) + O\left( \nabla_{\overline{\sigma}} \varphi_3 \right)
\end{align*}
and we can apply integrations by parts as in the case $\widehat{\mathcal{R}}\widehat{\mathcal{R}}\widehat{\mathcal{R}}$ to get only terms of \eqref{qteslemestimeesgeneriquesdeccubiquehb2} or of the form $g_{bb, 3}(t)$. 

If $\epsilon_3 = \epsilon_1 = \epsilon_2 = 1$, then at $Z$ we can compute that 
\begin{align*}
\overline{\xi}_b^{\overline{\eta}} &= - 2^{\frac{5}{3}} \overline{\eta}_a^{\overline{\eta}}, \quad \quad \quad \overline{\xi}_a^{\overline{\eta}} &= 3 \overline{\eta}_a^{\overline{\eta}} 
\end{align*}
In particular, $\overline{\xi}$ is localised by $m_{\widehat{\mathcal{R}}}$, and even outside of the support of $m_{\widehat{\mathcal{C}}}+m_{\widehat{\mathcal{L}}}+m_{\widehat{\mathcal{P}}}$ (assuming these symbols are chosen fine enough). We then always have that 
\begin{align*}
1 &= O\left( \nabla_{\overline{\eta}} \varphi_3 \right) 
\end{align*}
so we can apply an integration by parts along $\overline{\eta}$. We then obtain terms from \eqref{qteslemestimeesgeneriquesdeccubiquehb2}, plus a degenerate term: 
\begin{align}
&\int_0^t \int \int e^{i s \varphi_3} s \mu_{BBBB} \mu |\overline{\xi}|^2 m_{\widehat{\mathcal{C}}}(\overline{\eta}) m_{\widehat{\mathcal{R}}}(\overline{\sigma}) m_{\widehat{\mathcal{R}}}(\overline{\rho}) \nabla_{\overline{\eta}} \widehat{f}(s, \overline{\eta}) \widehat{f}(s, \overline{\sigma}) \widehat{f}(s, \overline{\rho}) ~ d\overline{\eta} d\overline{\sigma} ds \label{termedegenerecasCRRcubiquehb2} 
\end{align}

We then consider two subcases depending on the size of $\overline{\xi}$. First, if we localise to have $|\overline{\xi}| \lesssim t^{-\frac{5}{24}+10\delta}$, then we do not apply the above integration by parts and estimate directly: 
\begin{align*}
\Vert \chi\left( t^{\frac{5}{24}-10\delta} \overline{\xi} \right) \partial_t \eqref{equdecchampbCRR-1} \Vert_{L^2} 
&\lesssim t^{2-\frac{5}{8}+30\delta} \Vert |\nabla|^{-1} m_{\widehat{\mathcal{C}}}(D) f(t) \Vert_{L^2} \Vert \nabla m_{\widehat{\mathcal{R}}}(D) u(t) \Vert_{L^{\infty}}^2 \\
&\lesssim t^{-\frac{7}{24}} \langle t \rangle^{-\frac{2}{3}+230\delta} \Vert u \Vert_X^3 
\end{align*}
which is enough. On the other hand, if we localise $|\overline{\xi}| \gtrsim t^{-\frac{5}{24}+1O\delta}$, then we rather estimate \eqref{termedegenerecasCRRcubiquehb2} and get: 
\begin{align*}
\Vert \left( 1 - \chi\left( t^{\frac{5}{24}-10\delta} \overline{\xi} \right) \right) \partial_t \eqref{termedegenerecasCRRcubiquehb2} \Vert_{L^2}
&\lesssim t^{1+\frac{5}{24}-10\delta} \Vert \nabla (x, y) f(t) \Vert_{L^2} \Vert \nabla m_{\widehat{\mathcal{R}}}(D) u(t) \Vert_{L^{\infty}}^2 \\
&\lesssim t^{-\frac{11}{24}} \langle t \rangle^{-\frac{1}{2}+291\delta} \Vert u \Vert_X^3 
\end{align*}
which is also enough. 

\subsubsection{Interaction \texorpdfstring{$\widehat{\mathcal{L}}\widehat{\mathcal{R}}\widehat{\mathcal{R}}$}{LRR}} 
 
Let us consider: {\footnotesize 
\begin{subequations}
\begin{align}
&\xi_0 m_b(\overline{\xi}) \widehat{X}_b(\overline{\xi}) \cdot \nabla_{\overline{\xi}} \widehat{I}_{s\mu}^{\widehat{\mathcal{L}}\widehat{\mathcal{R}}\widehat{\mathcal{R}}}[F_1, F_2, F_3](t, \overline{\xi}) \notag \\
&\quad = \int_0^t \int \int i s^2 \xi_0 m_b(\overline{\xi}) \widehat{X}_b(\overline{\xi}) \cdot \nabla_{\overline{\xi}} \varphi_3 e^{i s \varphi_3} \mu(\overline{\xi}, \overline{\eta}, \overline{\sigma}) m_{\widehat{\mathcal{L}}}(\overline{\eta}) m_{\widehat{\mathcal{R}}}(\overline{\sigma}) m_{\widehat{\mathcal{R}}}(\overline{\rho}) \widehat{F}_1(s, \overline{\eta}) \widehat{F}_2(s, \overline{\sigma}) \widehat{F}_3(s, \overline{\rho}) ~ d\overline{\eta} d\overline{\sigma} ds \label{equdecchampbLRR-1} \\
&\quad \quad + \int_0^t \int \int e^{i s \varphi_3} \xi_0 s \mu(\overline{\xi}, \overline{\eta}, \overline{\sigma}) m_{\widehat{\mathcal{L}}}(\overline{\eta}) m_{\widehat{\mathcal{R}}}(\overline{\sigma}) m_{\widehat{\mathcal{R}}}(\overline{\rho}) \widehat{F}_1(s, \overline{\eta}) \widehat{F}_2(s, \overline{\sigma}) m_b(\overline{\xi}) \widehat{X}_b(\overline{\xi}) \cdot \nabla_{\overline{\xi}} \widehat{F}_3(s, \overline{\rho}) ~ d\overline{\eta} d\overline{\sigma} ds \label{equdecchampbLRR-2} \\
&\quad \quad + \int_0^t \int \int e^{i s \varphi_3} s \xi_0 m_b(\overline{\xi}) \widehat{X}_b(\overline{\xi}) \cdot \nabla_{\overline{\xi}} \left( \mu(\overline{\xi}, \overline{\eta}, \overline{\sigma}) m_{\widehat{\mathcal{L}}}(\overline{\eta}) m_{\widehat{\mathcal{R}}}(\overline{\sigma}) m_{\widehat{\mathcal{R}}}(\overline{\rho}) \right) \widehat{F}_1(s, \overline{\eta}) \widehat{F}_2(s, \overline{\sigma}) \widehat{F}_3(s, \overline{\rho}) ~ d\overline{\eta} d\overline{\sigma} ds \label{equdecchampbLRR-3} 
\end{align}
\end{subequations} }
\eqref{equdecchampbLRR-3} is of the form $\eqref{lemestimeesgeneriquesdeccubiquehb2-dersymb}+\eqref{lemestimeesgeneriquesdeccubiquehb2-dersymbbis}$, and \eqref{equdecchampbLRR-2} is of the form $\eqref{lemestimeesgeneriquesdeccubiquehb2-eta0sigmaresxrho}+\eqref{lemestimeesgeneriquesdeccubiquehb2-rho0sigmaresxrho}$. 

Finally, for \eqref{equdecchampbLRR-1}, as before, we can restrict our attention to the angular neighborhood of a point $Z \in S^8$ such that 
\begin{align*}
\varphi_3 = 0, \quad \nabla_{\overline{\sigma}} \varphi_3 = 0, \quad \partial_{\eta_0} \varphi_3 = 0, \quad \eta^Z = 0 
\end{align*}
and $\overline{\sigma}^Z, \overline{\rho}^Z$ are in the support of $m_{\widehat{\mathcal{R}}}$, in particular $\sigma^Z, \rho^Z \neq 0$. We then have $|\overline{\eta}| \simeq |\overline{\sigma}| \simeq |\overline{\rho}|$. But at the point $Z$, 
\begin{align*}
\nabla_{\eta} \varphi_3 &= 2 \rho_0^Z \rho^Z - 2 \eta_0^Z \eta^Z = 2 \rho_0^Z \rho^Z \neq 0
\end{align*}
so that, locally, $1 = O(\nabla_{\eta} \varphi_3)$. We can then apply an integration by parts along $\overline{\eta}$ and obtain terms from \eqref{qteslemestimeesgeneriquesdeccubiquehb2}. 

\subsubsection{Interactions \texorpdfstring{$\widehat{\mathcal{P}}\widehat{\mathcal{R}}\widehat{\mathcal{R}}$}{PRR}, \texorpdfstring{$\widehat{\mathcal{P}}\widehat{\mathcal{C}}\widehat{\mathcal{R}}$}{PCR}, \texorpdfstring{$\widehat{\mathcal{P}}\widehat{\mathcal{L}}\widehat{\mathcal{R}}$}{PLR}}

Let us consider {\footnotesize 
\begin{subequations}
\begin{align}
&\xi_0 m_b(\overline{\xi}) \widehat{X}_b(\overline{\xi}) \cdot \nabla_{\overline{\xi}} \widehat{I}_{s \mu}^{\widehat{\mathcal{P}}A_2\widehat{\mathcal{R}}}[F_1, F_2, F_3](t, \overline{\xi}) \notag \\
&\quad = \int_0^t \int \int i s^2 \xi_0 m_b(\overline{\xi}) \widehat{X}_b(\overline{\xi}) \cdot \nabla_{\overline{\xi}} \varphi_3 e^{i s \varphi_3} \mu(\overline{\xi}, \overline{\eta}, \overline{\sigma}) m_{\widehat{\mathcal{P}}}(\overline{\eta}) m_{A_2}(\overline{\sigma}) m_{\widehat{\mathcal{R}}}(\overline{\rho}) \widehat{F}_1(s, \overline{\eta}) \widehat{F}_2(s, \overline{\sigma}) \widehat{F}_3(s, \overline{\rho}) ~ d\overline{\eta} d\overline{\sigma} ds \label{equdecchampbPRR-1} \\
&\quad \quad + \int_0^t \int \int e^{i s \varphi_3} s \xi_0 \mu(\overline{\xi}, \overline{\eta}, \overline{\sigma}) m_{\widehat{\mathcal{P}}}(\overline{\eta}) m_{A_2}(\overline{\sigma}) m_{\widehat{\mathcal{R}}}(\overline{\rho}) \widehat{F}_1(s, \overline{\eta}) \widehat{F}_2(s, \overline{\sigma}) m_b(\overline{\xi}) \widehat{X}_b(\overline{\xi}) \cdot \nabla_{\overline{\xi}} \widehat{F}_3(s, \overline{\rho}) ~ d\overline{\eta} d\overline{\sigma} ds \label{equdecchampbPRR-2} \\
&\quad \quad + \int_0^t \int \int e^{i s \varphi_3} s \xi_0 m_b(\overline{\xi}) \widehat{X}_b(\overline{\xi}) \cdot \nabla_{\overline{\xi}} \left( \mu(\overline{\xi}, \overline{\eta}, \overline{\sigma}) m_{\widehat{\mathcal{P}}}(\overline{\eta}) m_{A_2}(\overline{\sigma}) m_{\widehat{\mathcal{R}}}(\overline{\rho}) \right) \widehat{F}_1(s, \overline{\eta}) \widehat{F}_2(s, \overline{\sigma}) \widehat{F}_3(s, \overline{\rho}) ~ d\overline{\eta} d\overline{\sigma} ds \label{equdecchampbPRR-3} 
\end{align}
\end{subequations} }
for $A_2 \in \{ \widehat{\mathcal{R}}, \widehat{\mathcal{C}}, \widehat{\mathcal{L}} \}$. 
\eqref{equdecchampbPRR-3} is of the form $\eqref{lemestimeesgeneriquesdeccubiquehb2-dersymb}+\eqref{lemestimeesgeneriquesdeccubiquehb2-dersymbbis}$, and \eqref{equdecchampbPRR-2} is of the form $\eqref{lemestimeesgeneriquesdeccubiquehb2-eta0sigmaresxrho}+\eqref{lemestimeesgeneriquesdeccubiquehb2-rho0sigmaresxrho}$ (up to symmetrizing the variables). 

As before, for \eqref{equdecchampbPRR-1}, we can consider only the angular neighborhood of a point $Z \in S^8$ such that, at $Z$, 
\begin{align*}
\varphi_3 = 0, \quad \nabla_{\eta} \varphi_3 = 0, \quad \eta_0^Z = 0, \quad \widehat{X}_a(\overline{\sigma}^Z) \cdot \nabla_{\overline{\sigma}} \varphi_3 = 0 
\end{align*}
and $\overline{\rho}^Z$ is in the support of $m_{\widehat{\mathcal{R}}}$. In particular, $\rho_0^Z \rho^Z = 0$, so $\overline{\rho}^Z = 0$, and $\overline{\sigma}^Z = 0$ as well. Therefore, in the neighborhood of $Z$, $|\overline{\sigma}|+|\overline{\rho}| \ll |\overline{\eta}|$. In particular, $1 = O(\partial_{\eta_0} \varphi_3)$. 

We have that 
\begin{align*}
\varphi_3 &= O(\eta_0) + \xi_0^3 + \xi_0 |\xi|^2 - \sigma_0^3 - \sigma_0 |\sigma|^2 - \rho_0^3 - \rho_0 |\rho|^2 \\
&= O(\eta_0) + O\left( \overline{\rho} \overline{\sigma} \right) + \xi_0 (|\xi|^2 - |\sigma|^2- |\rho|^2) 
\end{align*}
Since here $|\xi| \simeq |\eta| \simeq |\overline{\eta}| \gg |\overline{\sigma}|+|\overline{\rho}|$, we deduce that 
\begin{align*}
\xi_0 &= O(\varphi_3) + O(\eta_0) + O\left( \overline{\rho} \overline{\sigma} \right) 
\end{align*}
Recall that $\mu$ allows to distribute one derivative. We can then apply an integration by parts in time on the term having a factor $O(\varphi_3)$ and recover terms from \eqref{qteslemestimeesgeneriquesdeccubiquehb2} or \eqref{lemestimeesgeneriquesdeccubiquehb2-g}; for the term having a factor $O(\eta_0)$, we can apply an integration by parts along $\eta_0$ and recover also terms from \eqref{qteslemestimeesgeneriquesdeccubiquehb2}; finally, on the term having a factor $O(\overline{\rho} \overline{\sigma})$, we can also apply an integration by parts in $\eta_0$ and we recognize terms from \eqref{qteslemestimeesgeneriquesdeccubiquehb2}. 

\subsubsection{Interaction \texorpdfstring{$\widehat{\mathcal{C}}\widehat{\mathcal{C}}\widehat{\mathcal{R}}$}{CCR}}

Let us consider: {\footnotesize 
\begin{subequations}
\begin{align}
&\xi_0 m_b(\overline{\xi}) \widehat{X}_b(\overline{\xi}) \cdot \nabla_{\overline{\xi}} \widehat{I}_{s \mu}^{\widehat{\mathcal{C}}\widehat{\mathcal{C}}\widehat{\mathcal{R}}}[F_1, F_2, F_3](t, \overline{\xi}) \notag \\
&\quad = \int_0^t \int \int i s^2 \xi_0 m_b(\overline{\xi}) \widehat{X}_b(\overline{\xi}) \cdot \nabla_{\overline{\xi}} \varphi_3 e^{i s \varphi_3} \mu(\overline{\xi}, \overline{\eta}, \overline{\sigma}) m_{\widehat{\mathcal{C}}}(\overline{\eta}) m_{\widehat{\mathcal{C}}}(\overline{\sigma}) m_{\widehat{\mathcal{R}}}(\overline{\rho}) \widehat{F}_1(s, \overline{\eta}) \widehat{F}_2(s, \overline{\sigma}) \widehat{F}_3(s, \overline{\rho}) ~ d\overline{\eta} d\overline{\sigma} ds \label{equdecchampbCCR-1} \\
&\quad \quad + \int_0^t \int \int e^{i s \varphi_3} s \xi_0 \mu(\overline{\xi}, \overline{\eta}, \overline{\sigma}) m_{\widehat{\mathcal{C}}}(\overline{\eta}) m_{\widehat{\mathcal{C}}}(\overline{\sigma}) m_{\widehat{\mathcal{R}}}(\overline{\rho}) \widehat{F}_1(s, \overline{\eta}) \widehat{F}_2(s, \overline{\sigma}) m_b(\overline{\xi}) \widehat{X}_b(\overline{\xi}) \cdot \nabla_{\overline{\xi}} \widehat{F}_3(s, \overline{\rho}) ~ d\overline{\eta} d\overline{\sigma} ds \label{equdecchampbCCR-2} \\
&\quad \quad + \int_0^t \int \int e^{i s \varphi_3} s \xi_0 m_b(\overline{\xi}) \widehat{X}_b(\overline{\xi}) \cdot \nabla_{\overline{\xi}} \left( \mu(\overline{\xi}, \overline{\eta}, \overline{\sigma}) m_{\widehat{\mathcal{C}}}(\overline{\eta}) m_{\widehat{\mathcal{C}}}(\overline{\sigma}) m_{\widehat{\mathcal{R}}}(\overline{\rho}) \right) \widehat{F}_1(s, \overline{\eta}) \widehat{F}_2(s, \overline{\sigma}) \widehat{F}_3(s, \overline{\rho}) ~ d\overline{\eta} d\overline{\sigma} ds \label{equdecchampbCCR-3} 
\end{align}
\end{subequations} }
\eqref{equdecchampbCCR-3} is of the form $\eqref{lemestimeesgeneriquesdeccubiquehb2-dersymb}+\eqref{lemestimeesgeneriquesdeccubiquehb2-dersymbbis}$, and \eqref{equdecchampbCCR-2} is of the form $\eqref{lemestimeesgeneriquesdeccubiquehb2-eta0sigmaresxrho}+\eqref{lemestimeesgeneriquesdeccubiquehb2-rho0sigmaresxrho}$. 

For \eqref{equdecchampbCCR-1}, we need only to consider the angular neighborhood of a point $Z \in S^8$ such that, at $Z$, 
\begin{align*}
\varphi_3 = \widehat{X}_a(\overline{\eta}) \cdot \nabla_{\overline{\eta}} \varphi_3 = \widehat{X}_c(\overline{\eta}) \cdot \nabla_{\overline{\eta}} \varphi_3 = \widehat{X}_a(\overline{\sigma}) \cdot \nabla_{\overline{\sigma}} \varphi_3 = \widehat{X}_c(\overline{\sigma}) \cdot \nabla_{\overline{\sigma}} \varphi_3 = m_b(\overline{\eta}^Z) = m_b(\overline{\sigma}^Z) = 0 
\end{align*}
and $\overline{\rho}^Z$ is in the support of $m_{\widehat{\mathcal{R}}}$. In particular, this implies that, near $Z$, $|\overline{\eta}| \simeq |\overline{\sigma}| \simeq |\overline{\rho}|$. Furthermore, by Lemma \ref{lemcalculsconecoordonneesconiquesvarphi}, 
\begin{align*}
&\left( \xi_t^{\overline{\eta} \overline{\rho}}, \xi_t^{\overline{\sigma} \overline{\rho}}, \overline{\eta}_b^{\overline{\eta}}, \overline{\sigma}_b^{\overline{\sigma}}, \left( \overline{\eta}_a^{\overline{\eta}} \right)^2 - \left( \overline{\rho}_a^{\overline{\eta}} \right)^2, \left( \overline{\sigma}_a^{\overline{\sigma}} \right)^2 - \left( \overline{\rho}_a^{\overline{\eta}} \right)^2 \right) \\
&= O\left( \widehat{X}_a(\overline{\eta}) \cdot \nabla_{\overline{\eta}} \varphi_3 \right) + O\left( \widehat{X}_c(\overline{\eta}) \cdot \nabla_{\overline{\eta}} \varphi_3 \right) + O\left( \widehat{X}_a(\overline{\sigma}) \cdot \nabla_{\overline{\sigma}} \varphi_3 \right) \\
&\quad + O\left( \widehat{X}_c(\overline{\sigma}) \cdot \nabla_{\overline{\sigma}} \varphi_3 \right) + O\left( m_b(\overline{\eta}) \widehat{X}_b(\overline{\eta}) \cdot \nabla_{\overline{\eta}} \varphi_3 \right) + O\left( m_b(\overline{\sigma}) \widehat{X}_b(\overline{\sigma}) \cdot \nabla_{\overline{\sigma}} \varphi_3 \right)
\end{align*}
Locally, there are signs $\epsilon_1, \epsilon_2 \in \{ 1, -1 \}$ such that 
\begin{align*}
\overline{\rho}_a^{\overline{\eta}} = \epsilon_1 \overline{\eta}_a^{\overline{\eta}} + O\left( \left( \overline{\eta}_a^{\overline{\eta}} \right)^2 - \left( \overline{\rho}_a^{\overline{\eta}} \right)^2 \right), \quad \quad \overline{\rho}_a^{\overline{\sigma}} = \epsilon_2 \overline{\sigma}_a^{\overline{\sigma}} + O\left( \left( \overline{\sigma}_a^{\overline{\sigma}} \right)^2 - \left( \overline{\rho}_a^{\overline{\sigma}} \right)^2 \right) 
\end{align*}

\paragraph{1.} Assume first that, near $Z$, $\epsilon^{\overline{\eta} \overline{\sigma}} \theta^{\overline{\eta} \overline{\sigma}}$ is close to $1$, i.e. $\overline{\eta}$ and $\overline{\sigma}$ are close to alignment. We then have 
\begin{align*}
\overline{\rho}_a^{\overline{\sigma}} &= \epsilon^{\overline{\eta} \overline{\sigma}} \overline{\rho}_a^{\overline{\eta}} + O\left( \xi_t^{\overline{\rho} \overline{\sigma}} \right) + O\left( \xi_t^{\overline{\rho} \overline{\eta}} \right) 
\end{align*}
and likewise replacing $\overline{\rho}$ by one of the other variables, or $a$ by $b$. We then have $\epsilon_2 = \epsilon_1 \epsilon^{\overline{\eta} \overline{\sigma}}$, and 
\begin{align*}
\overline{\xi}_a^{\overline{\eta}} &= \left( 1 + \epsilon^{\overline{\eta} \overline{\sigma}} + \epsilon_1 \right)^3 \overline{\eta}_a^{\overline{\eta}} + O\left( \widehat{X}_a(\overline{\eta}) \cdot \nabla_{\overline{\eta}} \varphi_3 \right) + O\left( \widehat{X}_c(\overline{\eta}) \cdot \nabla_{\overline{\eta}} \varphi_3 \right) + O\left( \widehat{X}_a(\overline{\sigma}) \cdot \nabla_{\overline{\sigma}} \varphi_3 \right) \\
&\quad + O\left( \widehat{X}_c(\overline{\sigma}) \cdot \nabla_{\overline{\sigma}} \varphi_3 \right) + O\left( m_b(\overline{\eta}) \widehat{X}_b(\overline{\eta}) \cdot \nabla_{\overline{\eta}} \varphi_3 \right) + O\left( m_b(\overline{\sigma}) \widehat{X}_b(\overline{\sigma}) \cdot \nabla_{\overline{\sigma}} \varphi_3 \right) \\
\overline{\xi}_b^{\overline{\eta}} &= \overline{\rho}_b^{\overline{\eta}} + O\left( \widehat{X}_a(\overline{\eta}) \cdot \nabla_{\overline{\eta}} \varphi_3 \right) + O\left( \widehat{X}_c(\overline{\eta}) \cdot \nabla_{\overline{\eta}} \varphi_3 \right) + O\left( \widehat{X}_a(\overline{\sigma}) \cdot \nabla_{\overline{\sigma}} \varphi_3 \right) \\
&\quad + O\left( \widehat{X}_c(\overline{\sigma}) \cdot \nabla_{\overline{\sigma}} \varphi_3 \right) + O\left( m_b(\overline{\eta}) \widehat{X}_b(\overline{\eta}) \cdot \nabla_{\overline{\eta}} \varphi_3 \right) + O\left( m_b(\overline{\sigma}) \widehat{X}_b(\overline{\sigma}) \cdot \nabla_{\overline{\sigma}} \varphi_3 \right)
\end{align*}
Therefore, by Lemma \ref{lemcalculsconecoordonneesconiquesvarphi}, {\footnotesize 
\begin{align*}
6 \sqrt{3} \frac{\eta_0}{|\eta_0|} \varphi_3 &= \left( \overline{\eta}_a^{\overline{\eta}} \right)^3 \left( \left( 1 + \epsilon^{\overline{\eta} \overline{\sigma}} + \epsilon_1 \right)^3 - 1 - \epsilon^{\overline{\eta} \overline{\sigma}} - \epsilon_1 \right) + O\left( \widehat{X}_a(\overline{\eta}) \cdot \nabla_{\overline{\eta}} \varphi_3 \right) + O\left( \widehat{X}_c(\overline{\eta}) \cdot \nabla_{\overline{\eta}} \varphi_3 \right) \\
&+ O\left( \widehat{X}_a(\overline{\sigma}) \cdot \nabla_{\overline{\sigma}} \varphi_3 \right) + O\left( \widehat{X}_c(\overline{\sigma}) \cdot \nabla_{\overline{\sigma}} \varphi_3 \right) + O\left( m_b(\overline{\eta}) \widehat{X}_b(\overline{\eta}) \cdot \nabla_{\overline{\eta}} \varphi_3 \right) + O\left( m_b(\overline{\sigma}) \widehat{X}_b(\overline{\sigma}) \cdot \nabla_{\overline{\sigma}} \varphi_3 \right) \\
&= 6 \left( \overline{\eta}_a^{\overline{\eta}} \right)^3 (1 + \epsilon_1) \left( 1 + \epsilon^{\overline{\eta} \overline{\sigma}} \right) + O\left( \widehat{X}_a(\overline{\eta}) \cdot \nabla_{\overline{\eta}} \varphi_3 \right) + O\left( \widehat{X}_c(\overline{\eta}) \cdot \nabla_{\overline{\eta}} \varphi_3 \right) + O\left( \widehat{X}_a(\overline{\sigma}) \cdot \nabla_{\overline{\sigma}} \varphi_3 \right) \\
&\quad + O\left( \widehat{X}_c(\overline{\sigma}) \cdot \nabla_{\overline{\sigma}} \varphi_3 \right) + O\left( m_b(\overline{\eta}) \widehat{X}_b(\overline{\eta}) \cdot \nabla_{\overline{\eta}} \varphi_3 \right) + O\left( m_b(\overline{\sigma}) \widehat{X}_b(\overline{\sigma}) \cdot \nabla_{\overline{\sigma}} \varphi_3 \right)
\end{align*} }
In particular, for $\varphi_3$ to vanish at $Z$, we need $\left( \epsilon_1, \epsilon^{\overline{\eta} \overline{\sigma}} \right) \neq (1, 1)$. In any case, this implies that 
\begin{align*}
\left( \overline{\xi}_a^{\overline{\eta}} \right)^2 &= \left( \overline{\rho}_a^{\overline{\eta}} \right)^2 + O\left( \widehat{X}_a(\overline{\eta}) \cdot \nabla_{\overline{\eta}} \varphi_3 \right) + O\left( \widehat{X}_c(\overline{\eta}) \cdot \nabla_{\overline{\eta}} \varphi_3 \right) + O\left( \widehat{X}_a(\overline{\sigma}) \cdot \nabla_{\overline{\sigma}} \varphi_3 \right) \\
&\quad + O\left( \widehat{X}_c(\overline{\sigma}) \cdot \nabla_{\overline{\sigma}} \varphi_3 \right) + O\left( m_b(\overline{\eta}) \widehat{X}_b(\overline{\eta}) \cdot \nabla_{\overline{\eta}} \varphi_3 \right) + O\left( m_b(\overline{\sigma}) \widehat{X}_b(\overline{\sigma}) \cdot \nabla_{\overline{\sigma}} \varphi_3 \right)
\end{align*}
We deduce that
\begin{align*}
\partial_{\xi_0} \varphi_3 &= 3 \xi_0^2 + |\xi|^2 - 3 \rho_0^2 - |\rho|^2 \\
&= \frac{1}{2} \left( \left( \overline{\xi}_a^{\overline{\eta}} \right)^2 + \left( \overline{\xi}_b^{\overline{\eta}} \right)^2 - \left( \overline{\rho}_a^{\overline{\eta}} \right)^2 - \left( \overline{\rho}_b^{\overline{\eta}} \right)^2 \right) + O\left( \xi_t^{\overline{\xi} \overline{\eta}} \right) + O\left( \xi_t^{\overline{\rho} \overline{\eta}} \right) \\
&= O\left( \widehat{X}_a(\overline{\eta}) \cdot \nabla_{\overline{\eta}} \varphi_3 \right) + O\left( \widehat{X}_c(\overline{\eta}) \cdot \nabla_{\overline{\eta}} \varphi_3 \right) + O\left( \widehat{X}_a(\overline{\sigma}) \cdot \nabla_{\overline{\sigma}} \varphi_3 \right) \\
&\quad + O\left( \widehat{X}_c(\overline{\sigma}) \cdot \nabla_{\overline{\sigma}} \varphi_3 \right) + O\left( m_b(\overline{\eta}) \widehat{X}_b(\overline{\eta}) \cdot \nabla_{\overline{\eta}} \varphi_3 \right) + O\left( m_b(\overline{\sigma}) \widehat{X}_b(\overline{\sigma}) \cdot \nabla_{\overline{\sigma}} \varphi_3 \right) \\
\frac{J \eta}{|\eta|} \cdot \nabla_{\xi} \varphi_3 &= O\left( \xi_t^{\overline{\xi} \overline{\eta}} \right) + O\left( \xi_t^{\overline{\rho} \overline{\eta}} \right) \\
&= O\left( \widehat{X}_c(\overline{\eta}) \cdot \nabla_{\overline{\eta}} \varphi_3 \right) + O\left( \widehat{X}_c(\overline{\sigma}) \cdot \nabla_{\overline{\sigma}} \varphi_3 \right) \\
\frac{\eta_0 \eta}{|\eta_0| |\eta|} \cdot \nabla_{\xi} \varphi_3 &= 2 \frac{\eta_0 \xi_0}{|\eta_0|} \frac{\eta \cdot \xi}{|\eta|} - 2 \frac{\eta_0 \rho_0}{|\eta_0|} \frac{\eta \cdot \rho}{|\eta|} \\
&= \frac{1}{2 \sqrt{3}} \left( \left( \overline{\xi}_a^{\overline{\eta}} \right)^2 - \left( \overline{\xi}_b^{\overline{\eta}} \right)^2 - \left( \overline{\rho}_a^{\overline{\eta}} \right)^2 + \left( \overline{\rho}_b^{\overline{\eta}} \right)^2 \right) \\
&= O\left( \widehat{X}_a(\overline{\eta}) \cdot \nabla_{\overline{\eta}} \varphi_3 \right) + O\left( \widehat{X}_c(\overline{\eta}) \cdot \nabla_{\overline{\eta}} \varphi_3 \right) + O\left( \widehat{X}_a(\overline{\sigma}) \cdot \nabla_{\overline{\sigma}} \varphi_3 \right) \\
&\quad + O\left( \widehat{X}_c(\overline{\sigma}) \cdot \nabla_{\overline{\sigma}} \varphi_3 \right) + O\left( m_b(\overline{\eta}) \widehat{X}_b(\overline{\eta}) \cdot \nabla_{\overline{\eta}} \varphi_3 \right) + O\left( m_b(\overline{\sigma}) \widehat{X}_b(\overline{\sigma}) \cdot \nabla_{\overline{\sigma}} \varphi_3 \right)
\end{align*}
which allows to apply sufficient integrations by parts from $\nabla_{\overline{\xi}} \varphi_3$ in \eqref{equdecchampbCCR-1}. 

\paragraph{2.} Assume now that, near $Z$, $\epsilon^{\overline{\eta} \overline{\sigma}} \theta^{\overline{\eta} \overline{\sigma}}$ is close to $-1$. Then 
\begin{align*}
\overline{\rho}_a^{\overline{\sigma}} &= \epsilon^{\overline{\eta} \overline{\sigma}} \overline{\rho}_b^{\overline{\eta}} + O\left( \xi_t^{\overline{\rho} \overline{\eta}} \right) + O\left( \xi_t^{\overline{\rho} \overline{\sigma}} \right) 
\end{align*}
and likewise replacing $\overline{\rho}$ by another variable or exchanging the role of $a, b$. Therefore, 
\begin{align*}
\overline{\xi}_a^{\overline{\eta}} &= (1 + \epsilon_1) \overline{\eta}_a^{\overline{\eta}} + O\left( \widehat{X}_a(\overline{\eta}) \cdot \nabla_{\overline{\eta}} \varphi_3 \right) + O\left( \widehat{X}_c(\overline{\eta}) \cdot \nabla_{\overline{\eta}} \varphi_3 \right) + O\left( \widehat{X}_a(\overline{\sigma}) \cdot \nabla_{\overline{\sigma}} \varphi_3 \right) \\
&\quad + O\left( \widehat{X}_c(\overline{\sigma}) \cdot \nabla_{\overline{\sigma}} \varphi_3 \right) + O\left( m_b(\overline{\eta}) \widehat{X}_b(\overline{\eta}) \cdot \nabla_{\overline{\eta}} \varphi_3 \right) + O\left( m_b(\overline{\sigma}) \widehat{X}_b(\overline{\sigma}) \cdot \nabla_{\overline{\sigma}} \varphi_3 \right) \\
\overline{\xi}_b^{\overline{\eta}} &= \epsilon^{\overline{\eta} \overline{\sigma}} (1 + \epsilon_2) \overline{\sigma}_a^{\overline{\sigma}} + O\left( \widehat{X}_a(\overline{\eta}) \cdot \nabla_{\overline{\eta}} \varphi_3 \right) + O\left( \widehat{X}_c(\overline{\eta}) \cdot \nabla_{\overline{\eta}} \varphi_3 \right) + O\left( \widehat{X}_a(\overline{\sigma}) \cdot \nabla_{\overline{\sigma}} \varphi_3 \right) \\
&\quad + O\left( \widehat{X}_c(\overline{\sigma}) \cdot \nabla_{\overline{\sigma}} \varphi_3 \right) + O\left( m_b(\overline{\eta}) \widehat{X}_b(\overline{\eta}) \cdot \nabla_{\overline{\eta}} \varphi_3 \right) + O\left( m_b(\overline{\sigma}) \widehat{X}_b(\overline{\sigma}) \cdot \nabla_{\overline{\sigma}} \varphi_3 \right)
\end{align*}
We then have by Lemma \ref{lemcalculsconecoordonneesconiquesvarphi}: {\footnotesize 
\begin{align*}
6 \sqrt{3} \frac{\eta_0}{|\eta_0|} \varphi_3 &= \left( \overline{\eta}_a^{\overline{\eta}} \right)^3 \left( \left( 1 + \epsilon_1 \right)^3 - 1 - \epsilon_1 \right) + \epsilon^{\overline{\eta} \overline{\sigma}} \left( \overline{\sigma}_a^{\overline{\sigma}} \right)^3 \left( \left( 1 + \epsilon_2 \right)^2 - 1 - \epsilon_2 \right) + O\left( \widehat{X}_a(\overline{\eta}) \cdot \nabla_{\overline{\eta}} \varphi_3 \right) + O\left( \widehat{X}_c(\overline{\eta}) \cdot \nabla_{\overline{\eta}} \varphi_3 \right) \\
&\quad + O\left( \widehat{X}_a(\overline{\sigma}) \cdot \nabla_{\overline{\sigma}} \varphi_3 \right) + O\left( \widehat{X}_c(\overline{\sigma}) \cdot \nabla_{\overline{\sigma}} \varphi_3 \right) + O\left( m_b(\overline{\eta}) \widehat{X}_b(\overline{\eta}) \cdot \nabla_{\overline{\eta}} \varphi_3 \right) + O\left( m_b(\overline{\sigma}) \widehat{X}_b(\overline{\sigma}) \cdot \nabla_{\overline{\sigma}} \varphi_3 \right) \\
&= 3 \left( \overline{\eta}_a^{\overline{\eta}} \right)^3 \left( 1 + \epsilon_1 \right) + 3 \epsilon^{\overline{\eta} \overline{\sigma}} \left( \overline{\sigma}_a^{\overline{\sigma}} \right)^3 \left( 1 + \epsilon_2 \right) + O\left( \widehat{X}_a(\overline{\eta}) \cdot \nabla_{\overline{\eta}} \varphi_3 \right) + O\left( \widehat{X}_c(\overline{\eta}) \cdot \nabla_{\overline{\eta}} \varphi_3 \right) + O\left( \widehat{X}_a(\overline{\sigma}) \cdot \nabla_{\overline{\sigma}} \varphi_3 \right) \\
&\quad + O\left( \widehat{X}_c(\overline{\sigma}) \cdot \nabla_{\overline{\sigma}} \varphi_3 \right) + O\left( m_b(\overline{\eta}) \widehat{X}_b(\overline{\eta}) \cdot \nabla_{\overline{\eta}} \varphi_3 \right) + O\left( m_b(\overline{\sigma}) \widehat{X}_b(\overline{\sigma}) \cdot \nabla_{\overline{\sigma}} \varphi_3 \right)
\end{align*} }
For $\varphi_3$ to vanish at $Z$, it is needed that $\epsilon_1 = \epsilon_2 = -1$, or that $\epsilon_1 = \epsilon_2 = 1 = - \epsilon^{\overline{\eta} \overline{\sigma}}$ and $\overline{\eta}_a^{\overline{\eta}} = \overline{\sigma}_a^{\overline{\sigma}}$ at $Z$. In the first case, 
\begin{align*}
\overline{\xi} &= O\left( \overline{\xi}_a^{\overline{\eta}} \right) + O\left( \overline{\xi}_b^{\overline{\eta}} \right) + O\left( \xi_t^{\overline{\eta} \overline{\xi}} \right) \\
&= O\left( \widehat{X}_a(\overline{\eta}) \cdot \nabla_{\overline{\eta}} \varphi_3 \right) + O\left( \widehat{X}_c(\overline{\eta}) \cdot \nabla_{\overline{\eta}} \varphi_3 \right) + O\left( \widehat{X}_a(\overline{\sigma}) \cdot \nabla_{\overline{\sigma}} \varphi_3 \right) \\
&\quad + O\left( \widehat{X}_c(\overline{\sigma}) \cdot \nabla_{\overline{\sigma}} \varphi_3 \right) + O\left( m_b(\overline{\eta}) \widehat{X}_b(\overline{\eta}) \cdot \nabla_{\overline{\eta}} \varphi_3 \right) + O\left( m_b(\overline{\sigma}) \widehat{X}_b(\overline{\sigma}) \cdot \nabla_{\overline{\sigma}} \varphi_3 \right)
\end{align*}
which is enough. In the second case, this implies that 
\begin{align*}
\rho_0 &= \frac{\eta_0}{2 \sqrt{3} |\eta_0|} \left( \overline{\rho}_a^{\overline{\eta}} + \overline{\rho}_b^{\overline{\eta}} \right) \\
&= \frac{\eta_0}{2 \sqrt{3} |\eta_0|} \left( \overline{\rho}_a^{\overline{\eta}} - \overline{\rho}_b^{\overline{\sigma}} \right) + O\left( \widehat{X}_c(\overline{\eta}) \cdot \nabla_{\overline{\eta}} \varphi_3 \right) + O\left( \widehat{X}_c(\overline{\sigma}) \cdot \nabla_{\overline{\sigma}} \varphi_3 \right) \\
&= O\left( \widehat{X}_a(\overline{\eta}) \cdot \nabla_{\overline{\eta}} \varphi_3 \right) + O\left( \widehat{X}_c(\overline{\eta}) \cdot \nabla_{\overline{\eta}} \varphi_3 \right) + O\left( \widehat{X}_a(\overline{\sigma}) \cdot \nabla_{\overline{\sigma}} \varphi_3 \right) \\
&\quad + O\left( \widehat{X}_c(\overline{\sigma}) \cdot \nabla_{\overline{\sigma}} \varphi_3 \right) + O\left( m_b(\overline{\eta}) \widehat{X}_b(\overline{\eta}) \cdot \nabla_{\overline{\eta}} \varphi_3 \right) + O\left( m_b(\overline{\sigma}) \widehat{X}_b(\overline{\sigma}) \cdot \nabla_{\overline{\sigma}} \varphi_3 \right)
\end{align*}
In particular, ar $Z$, $\rho_0^Z = 0$, which is a contradiction given the localisation. 

\subsubsection{Interaction \texorpdfstring{$\widehat{\mathcal{C}}\widehat{\mathcal{L}}\widehat{\mathcal{R}}$}{CLR}}

Let us consider: {\footnotesize 
\begin{subequations}
\begin{align}
&\xi_0 m_b(\overline{\xi}) \widehat{X}_b(\overline{\xi}) \cdot \nabla_{\overline{\xi}} \widehat{I}_{s \mu}^{\widehat{\mathcal{C}}\widehat{\mathcal{L}}\widehat{\mathcal{R}}}[F_1, F_2, F_3](t, \overline{\xi}) \notag \\
&\quad = \int_0^t \int \int i s^2 \xi_0 m_b(\overline{\xi}) \widehat{X}_b(\overline{\xi}) \cdot \nabla_{\overline{\xi}} \varphi_3 e^{i s \varphi_3} \mu(\overline{\xi}, \overline{\eta}, \overline{\sigma}) m_{\widehat{\mathcal{C}}}(\overline{\eta}) m_{\widehat{\mathcal{L}}}(\overline{\sigma}) m_{\widehat{\mathcal{R}}}(\overline{\rho}) \widehat{F}_1(s, \overline{\eta}) \widehat{F}_2(s, \overline{\sigma}) \widehat{F}_3(s, \overline{\rho}) ~ d\overline{\eta} d\overline{\sigma} ds \label{equdecchampbLCR-1} \\
&\quad \quad + \int_0^t \int \int e^{i s \varphi_3} s \xi_0 \mu(\overline{\xi}, \overline{\eta}, \overline{\sigma}) m_{\widehat{\mathcal{C}}}(\overline{\eta}) m_{\widehat{\mathcal{L}}}(\overline{\sigma}) m_{\widehat{\mathcal{R}}}(\overline{\rho}) \widehat{F}_1(s, \overline{\eta}) \widehat{F}_2(s, \overline{\sigma}) m_b(\overline{\xi}) \widehat{X}_b(\overline{\xi}) \cdot \nabla_{\overline{\xi}} \widehat{F}_3(s, \overline{\rho}) ~ d\overline{\eta} d\overline{\sigma} ds \label{equdecchampbLCR-2} \\
&\quad \quad + \int_0^t \int \int e^{i s \varphi_3} s \xi_0 m_b(\overline{\xi}) \widehat{X}_b(\overline{\xi}) \cdot \nabla_{\overline{\xi}} \left( \mu(\overline{\xi}, \overline{\eta}, \overline{\sigma}) m_{\widehat{\mathcal{C}}}(\overline{\eta}) m_{\widehat{\mathcal{L}}}(\overline{\sigma}) m_{\widehat{\mathcal{R}}}(\overline{\rho}) \right) \widehat{F}_1(s, \overline{\eta}) \widehat{F}_2(s, \overline{\sigma}) \widehat{F}_3(s, \overline{\rho}) ~ d\overline{\eta} d\overline{\sigma} ds \label{equdecchampbLCR-3} 
\end{align}
\end{subequations} }
\eqref{equdecchampbLCR-3} is of the form $\eqref{lemestimeesgeneriquesdeccubiquehb2-dersymb}+\eqref{lemestimeesgeneriquesdeccubiquehb2-dersymbbis}$, and \eqref{equdecchampbLCR-2} is of the form $\eqref{lemestimeesgeneriquesdeccubiquehb2-eta0sigmaresxrho}+\eqref{lemestimeesgeneriquesdeccubiquehb2-rho0sigmaresxrho}$. 

For \eqref{equdecchampbLCR-1}, we need only to consider the angular neighborhood of a point $Z$ such that, at $Z$, 
\begin{align*}
\varphi_3 = 0, \quad \partial_{\sigma_0} \varphi_3 = 0, \quad \sigma^Z = 0, \quad \widehat{X}_a(\overline{\eta}) \cdot \nabla_{\overline{\eta}} \varphi_3 = \widehat{X}_c(\overline{\eta}) \cdot \nabla_{\overline{\eta}} \varphi_3 = m_b(\overline{\eta}^Z) = 0 
\end{align*}
and $\overline{\rho}^Z$ is in the support of $m_{\widehat{\mathcal{R}}}$. In particular, this implies $|\overline{\eta}| \simeq |\overline{\sigma}| \simeq |\overline{\rho}|$ near $Z$. Also, 
\begin{align*}
1 = O\left( \nabla_{\sigma} \varphi_3 \right) 
\end{align*}
therefore we can apply integrations by parts in the direction $\sigma$. We then get terms from  \eqref{qteslemestimeesgeneriquesdeccubiquehb2}. 

\subsubsection{Interaction \texorpdfstring{$\widehat{\mathcal{L}}\widehat{\mathcal{L}}\widehat{\mathcal{R}}$}{LLR}}

Let us consider: {\footnotesize 
\begin{subequations}
\begin{align}
&\xi_0 m_b(\overline{\xi}) \widehat{X}_b(\overline{\xi}) \cdot \nabla_{\overline{\xi}} \widehat{I}_{s \mu}^{\widehat{\mathcal{L}}\widehat{\mathcal{L}}\widehat{\mathcal{R}}}[F_1, F_2, F_3](t, \overline{\xi}) \notag \\
&\quad = \int_0^t \int \int i s^2 \xi_0 m_b(\overline{\xi}) \widehat{X}_b(\overline{\xi}) \cdot \nabla_{\overline{\xi}} \varphi_3 e^{i s \varphi_3} \mu(\overline{\xi}, \overline{\eta}, \overline{\sigma}) m_{\widehat{\mathcal{L}}}(\overline{\eta}) m_{\widehat{\mathcal{L}}}(\overline{\sigma}) m_{\widehat{\mathcal{R}}}(\overline{\rho}) \widehat{F}_1(s, \overline{\eta}) \widehat{F}_2(s, \overline{\sigma}) \widehat{F}_3(s, \overline{\rho}) ~ d\overline{\eta} d\overline{\sigma} ds \label{equdecchampbLLR-1} \\
&\quad \quad + \int_0^t \int \int e^{i s \varphi_3} s \xi_0 \mu(\overline{\xi}, \overline{\eta}, \overline{\sigma}) m_{\widehat{\mathcal{L}}}(\overline{\eta}) m_{\widehat{\mathcal{L}}}(\overline{\sigma}) m_{\widehat{\mathcal{R}}}(\overline{\rho}) \widehat{F}_1(s, \overline{\eta}) \widehat{F}_2(s, \overline{\sigma}) m_b(\overline{\xi}) \widehat{X}_b(\overline{\xi}) \cdot \nabla_{\overline{\xi}} \widehat{F}_3(s, \overline{\rho}) ~ d\overline{\eta} d\overline{\sigma} ds \label{equdecchampbLLR-2} \\
&\quad \quad + \int_0^t \int \int e^{i s \varphi_3} s \xi_0 m_b(\overline{\xi}) \widehat{X}_b(\overline{\xi}) \cdot \nabla_{\overline{\xi}} \left( \mu(\overline{\xi}, \overline{\eta}, \overline{\sigma}) m_{\widehat{\mathcal{L}}}(\overline{\eta}) m_{\widehat{\mathcal{L}}}(\overline{\sigma}) m_{\widehat{\mathcal{R}}}(\overline{\rho}) \right) \widehat{F}_1(s, \overline{\eta}) \widehat{F}_2(s, \overline{\sigma}) \widehat{F}_3(s, \overline{\rho}) ~ d\overline{\eta} d\overline{\sigma} ds \label{equdecchampbLLR-3} 
\end{align}
\end{subequations} }
\eqref{equdecchampbLLR-3} is of the form $\eqref{lemestimeesgeneriquesdeccubiquehb2-dersymb}+\eqref{lemestimeesgeneriquesdeccubiquehb2-dersymbbis}$, and \eqref{equdecchampbLLR-2} is of the form $\eqref{lemestimeesgeneriquesdeccubiquehb2-eta0sigmaresxrho}+\eqref{lemestimeesgeneriquesdeccubiquehb2-rho0sigmaresxrho}$. 

For \eqref{equdecchampbLLR-1}, we need only to consider the neighborhood of a point $Z \in S^8$ such that, at $Z$, 
\begin{align*}
\varphi_3 = 0, \quad \partial_{\eta_0} \varphi_3 = \partial_{\sigma_0} \varphi_3 = 0, \quad \eta^Z = \sigma^Z = 0 
\end{align*}
and $\overline{\rho}^Z$ is in the support of $m_{\widehat{\mathcal{R}}}$. In particular, $|\overline{\eta}| \simeq |\overline{\sigma}| \simeq |\overline{\rho}|$. This implies that, at $Z$, $\nabla_{\eta} \varphi_3 = 2 \rho_0^Z \rho^Z \neq 0$, so locally $1 = O(\nabla_{\eta} \varphi)$. We can thus apply an integration by parts in $\eta$ and obtain terms from \eqref{qteslemestimeesgeneriquesdeccubiquehb2}. 

\subsubsection{Interaction \texorpdfstring{$\widehat{\mathcal{P}}\widehat{\mathcal{P}}\widehat{\mathcal{R}}$}{PPR}}

Let us consider: {\footnotesize 
\begin{subequations}
\begin{align}
&\xi_0 m_b(\overline{\xi}) \widehat{X}_b(\overline{\xi}) \cdot \nabla_{\overline{\xi}} \widehat{I}_{s \mu}^{\widehat{\mathcal{P}}\widehat{\mathcal{P}}\widehat{\mathcal{R}}}[F_1, F_2, F_3](t, \overline{\xi}) \notag \\
&\quad = \int_0^t \int \int i s^2 \xi_0 m_b(\overline{\xi}) \widehat{X}_b(\overline{\xi}) \cdot \nabla_{\overline{\xi}} \varphi_3 e^{i s \varphi_3} \mu(\overline{\xi}, \overline{\eta}, \overline{\sigma}) m_{\widehat{\mathcal{P}}}(\overline{\eta}) m_{\widehat{\mathcal{P}}}(\overline{\sigma}) m_{\widehat{\mathcal{R}}}(\overline{\rho}) \widehat{F}_1(s, \overline{\eta}) \widehat{F}_2(s, \overline{\sigma}) \widehat{F}_3(s, \overline{\rho}) ~ d\overline{\eta} d\overline{\sigma} ds \label{equdecchampbPPR-1} \\
&\quad \quad + \int_0^t \int \int e^{i s \varphi_3} s \xi_0 \mu(\overline{\xi}, \overline{\eta}, \overline{\sigma}) m_{\widehat{\mathcal{P}}}(\overline{\eta}) m_{\widehat{\mathcal{P}}}(\overline{\sigma}) m_{\widehat{\mathcal{R}}}(\overline{\rho}) \widehat{F}_1(s, \overline{\eta}) \widehat{F}_2(s, \overline{\sigma}) m_b(\overline{\xi}) \widehat{X}_b(\overline{\xi}) \cdot \nabla_{\overline{\xi}} \widehat{F}_3(s, \overline{\rho}) ~ d\overline{\eta} d\overline{\sigma} ds \label{equdecchampbPPR-2} \\
&\quad \quad + \int_0^t \int \int e^{i s \varphi_3} s \xi_0 m_b(\overline{\xi}) \widehat{X}_b(\overline{\xi}) \cdot \nabla_{\overline{\xi}} \left( \mu(\overline{\xi}, \overline{\eta}, \overline{\sigma}) m_{\widehat{\mathcal{P}}}(\overline{\eta}) m_{\widehat{\mathcal{P}}}(\overline{\sigma}) m_{\widehat{\mathcal{R}}}(\overline{\rho}) \right) \widehat{F}_1(s, \overline{\eta}) \widehat{F}_2(s, \overline{\sigma}) \widehat{F}_3(s, \overline{\rho}) ~ d\overline{\eta} d\overline{\sigma} ds \label{equdecchampbPPR-3} 
\end{align}
\end{subequations} }
\eqref{equdecchampbPPR-3} is of the form $\eqref{lemestimeesgeneriquesdeccubiquehb2-dersymb}+\eqref{lemestimeesgeneriquesdeccubiquehb2-dersymbbis}$, and \eqref{equdecchampbPPR-2} is of the form $\eqref{lemestimeesgeneriquesdeccubiquehb2-eta0sigmaresxrho}+\eqref{lemestimeesgeneriquesdeccubiquehb2-rho0sigmaresxrho}$. 

For \eqref{equdecchampbPPR-1}, we only need to consider the neighborhood of a point $Z \in S^8$ such that, at $Z$, 
\begin{align*}
\varphi_3 = 0, \quad \nabla_{\eta} \varphi = \nabla_{\sigma} \varphi = 0, \quad \eta_0^Z = \sigma_0^Z = 0
\end{align*}
In particular, $\overline{\rho}^Z = 0$, so $|\overline{\rho}| \ll |\overline{\eta}|+|\overline{\sigma}|$ near $Z$. 

We separate into two subcases depending on the relative sizes of $\overline{\eta}, \overline{\sigma}$ (which are symmetric here). 

\paragraph{1.} Assume first that, near $Z$, $|\overline{\eta}| \simeq |\overline{\sigma}|$. Then we have
\begin{align*}
\xi_0 &= O(\eta_0) + O(\sigma_0) + O(\overline{\rho})
\end{align*}
On the term with a factor $\eta_0$, we can apply an integration by parts in $\eta_0$ and, using $|\overline{\eta}| \simeq |\overline{\sigma}|$ to distribute the derivatives, we get terms from \eqref{qteslemestimeesgeneriquesdeccubiquehb2}. We can treat the term having a factor $\sigma_0$ in a symmetric way. Finally, on the term having a factor $O(\overline{\rho})$, we can also apply an integration by parts in $\eta_0$ and, using $|\overline{\sigma}| \simeq |\overline{\eta}|$ again to distribute the derivatives, we only get terms from \eqref{qteslemestimeesgeneriquesdeccubiquehb2}. 

\paragraph{2.} Assume that, near $Z$, $|\overline{\eta}| \gg |\overline{\sigma}|$. Then 
\begin{align*}
\varphi_3 &= O(\eta_0) + \xi_0^3 + \xi_0 |\xi|^2 - \sigma_0^3 - \sigma_0 |\sigma|^2 - \rho_0^3 - \rho_0 |\rho|^2 \\
&= O(\eta_0) + O(\overline{\sigma} \overline{\rho}) + \xi_0 \left( |\xi|^2 - |\sigma|^2 - |\rho|^2 \right) 
\end{align*}
Since $|\xi| \simeq |\eta| \simeq |\overline{\eta}| \gg |\overline{\sigma}|+|\overline{\rho}|$ here, we have that 
\begin{align*}
\xi_0 &= O(\varphi_3) + O(\eta_0) + O(\overline{\sigma} \overline{\rho}) 
\end{align*}
Moreover, $\mu$ allows to distribute one derivative. On the term with a factor $O(\varphi_3)$, we apply an integration by parts in time; on the ones having a factor $O(\eta_0)$ or $O(\overline{\sigma} \overline{\rho})$, we apply an integration by parts in $\eta_0$ and get in both cases terms from \eqref{qteslemestimeesgeneriquesdeccubiquehb2} ot \eqref{lemestimeesgeneriquesdeccubiquehb2-g}. 

\subsubsection{Interaction \texorpdfstring{$\widehat{\mathcal{C}}\widehat{\mathcal{C}}\widehat{\mathcal{C}}$}{CCC}}

Let us consider {\footnotesize 
\begin{subequations}
\begin{align}
&\xi_0 m_b(\overline{\xi}) \widehat{X}_b(\overline{\xi}) \cdot \nabla_{\overline{\xi}} \widehat{I}_{s \mu}^{\widehat{\mathcal{C}}\widehat{\mathcal{C}}\widehat{\mathcal{C}}}[F_1, F_2, F_3](t, \overline{\xi}) \notag \\
&\quad = \int_0^t \int \int i s^2 \xi_0 m_b(\overline{\xi}) \widehat{X}_b(\overline{\xi}) \cdot \nabla_{\overline{\xi}} \varphi_3 e^{i s \varphi_3} \mu(\overline{\xi}, \overline{\eta}, \overline{\sigma}) m_{\widehat{\mathcal{C}}}(\overline{\eta}) m_{\widehat{\mathcal{C}}}(\overline{\sigma}) m_{\widehat{\mathcal{C}}}(\overline{\rho}) \widehat{F}_1(s, \overline{\eta}) \widehat{F}_2(s, \overline{\sigma}) \widehat{F}_3(s, \overline{\rho}) ~ d\overline{\eta} d\overline{\sigma} ds \label{equdecchampbCCC-1} \\
&\quad \quad + \int_0^t \int \int e^{i s \varphi_3} s \xi_0 \mu(\overline{\xi}, \overline{\eta}, \overline{\sigma}) m_{\widehat{\mathcal{C}}}(\overline{\eta}) m_{\widehat{\mathcal{C}}}(\overline{\sigma}) m_{\widehat{\mathcal{C}}}(\overline{\rho}) \widehat{F}_1(s, \overline{\eta}) \widehat{F}_2(s, \overline{\sigma}) m_b(\overline{\xi}) \widehat{X}_b(\overline{\xi}) \cdot \nabla_{\overline{\xi}} \widehat{F}_3(s, \overline{\rho}) ~ d\overline{\eta} d\overline{\sigma} ds \label{equdecchampbCCC-2} \\
&\quad \quad + \int_0^t \int \int e^{i s \varphi_3} s \xi_0 m_b(\overline{\xi}) \widehat{X}_b(\overline{\xi}) \cdot \nabla_{\overline{\xi}} \left( \mu(\overline{\xi}, \overline{\eta}, \overline{\sigma}) m_{\widehat{\mathcal{C}}}(\overline{\eta}) m_{\widehat{\mathcal{C}}}(\overline{\sigma}) m_{\widehat{\mathcal{C}}}(\overline{\rho}) \right) \widehat{F}_1(s, \overline{\eta}) \widehat{F}_2(s, \overline{\sigma}) \widehat{F}_3(s, \overline{\rho}) ~ d\overline{\eta} d\overline{\sigma} ds \label{equdecchampbCCC-3} 
\end{align}
\end{subequations} }
\eqref{equdecchampbCCC-3} is of the form $\eqref{lemestimeesgeneriquesdeccubiquehb2-dersymb}+\eqref{lemestimeesgeneriquesdeccubiquehb2-dersymbbis}$. 

We separate three subcases, depending on the values of $\epsilon^{\overline{\eta} \overline{\sigma}} \theta^{\overline{\eta} \overline{\sigma}}, \epsilon^{\overline{\eta} \overline{\rho}} \theta^{\overline{\eta} \overline{\rho}}, \epsilon^{\overline{\sigma} \overline{\rho}} \theta^{\overline{\sigma} \overline{\rho}}$: 
\begin{enumerate}
\item either all of these values are close to $1$ (alignment); 
\item either two at least are sufficiently away from $1, -1$; 
\item or two at least are close to $-1$. 
\end{enumerate}
By symmetry of $\overline{\eta}, \overline{\sigma}, \overline{\rho}$, in the last case, we can assume that $\epsilon^{\overline{\eta} \overline{\sigma}} \theta^{\overline{\eta} \overline{\sigma}}, \epsilon^{\overline{\eta} \overline{\rho}} \theta^{\overline{\eta} \overline{\rho}}$ are close $-1$. In particular, this means that $\eta, \sigma, \rho$ are close to colinearity, so $\theta^{\overline{\sigma}, \overline{\rho}}$ is close to $\{ -1, 1 \}$: but the sign of 
\begin{align*}
\epsilon^{\overline{\eta} \overline{\sigma}} \theta^{\overline{\eta} \overline{\sigma}} \epsilon^{\overline{\eta} \overline{\rho}} \theta^{\overline{\eta} \overline{\rho}} \epsilon^{\overline{\sigma} \overline{\rho}} \theta^{\overline{\sigma} \overline{\rho}}
\end{align*}
needs to be $+1$, so this forces $\epsilon^{\overline{\sigma} \overline{\rho}} \theta^{\overline{\sigma} \overline{\rho}}$ to be close to $1$. 

\paragraph{1.} Localise first such that $\epsilon^{\overline{\eta} \overline{\sigma}} \theta^{\overline{\eta} \overline{\sigma}}, \epsilon^{\overline{\eta} \overline{\rho}} \theta^{\overline{\eta} \overline{\rho}}, \epsilon^{\overline{\sigma} \overline{\rho}} \theta^{\overline{\sigma} \overline{\rho}}$ are all close to $1$, i.e. $\overline{\eta}, \overline{\sigma}, \overline{\rho}$ are close to alignment. 

\paragraph{1.1.} Let us first also localise to have $|\overline{\xi}| \gtrsim |\overline{\eta}| + |\overline{\sigma}| + |\overline{\rho}|$. Then {\footnotesize 
\begin{align*}
&3 \xi_0^2 - |\xi|^2 = 3 \left( \eta_0^2 + \sigma_0^2 + \rho_0^2 + 2 \eta_0 \rho_0 + 2 \eta_0 \sigma_0 + 2 \sigma_0 \rho_0 \right) - \left( |\eta|^2 + |\sigma|^2 + |\rho|^2 + 2 \eta \cdot \sigma + 2 \eta \cdot \rho + 2 \sigma \cdot \rho \right) \\
&\quad = \left( 3 \eta_0^2 - |\eta|^2 \right) + \left( 3 \sigma_0^2 - |\sigma|^2 \right) + \left( 3 \rho_0^2 - |\rho|^2 \right) + \epsilon^{\overline{\eta} \overline{\rho}} \left( \sqrt{3} |\eta_0| + |\eta| \right) \left( \sqrt{3} |\rho_0| - |\rho| \right) + \epsilon^{\overline{\eta} \overline{\rho}} \left( \sqrt{3} |\eta_0| - |\eta| \right) \left( \sqrt{3} |\rho_0| + |\rho| \right) \\
&\quad + \epsilon^{\overline{\sigma} \overline{\rho}} \left( \sqrt{3} |\sigma_0| + |\sigma| \right) \left( \sqrt{3} |\rho_0| - |\rho| \right)
+ \epsilon^{\overline{\sigma} \overline{\rho}} \left( \sqrt{3} |\sigma_0| - |\sigma| \right) \left( \sqrt{3} |\rho_0| + |\rho| \right)
+ \epsilon^{\overline{\eta} \overline{\sigma}} \left( \sqrt{3} |\eta_0| + |\eta| \right) \left( \sqrt{3} |\sigma_0| - |\sigma| \right) \\
&\quad + \epsilon^{\overline{\eta} \overline{\sigma}} \left( \sqrt{3} |\eta_0| - |\eta| \right) \left( \sqrt{3} |\sigma_0| + |\sigma| \right) 
+ 2 \left( \epsilon^{\overline{\eta} \overline{\rho}} - \theta^{\overline{\eta} \overline{\rho}} \right) |\eta| |\rho|
+ 2 \left( \epsilon^{\overline{\sigma} \overline{\rho}} - \theta^{\overline{\sigma} \overline{\rho}} \right) |\sigma| |\rho|  
+ 2 \left( \epsilon^{\overline{\eta} \overline{\sigma}} - \theta^{\overline{\eta} \overline{\sigma}} \right) |\eta| |\sigma| 
\end{align*} }
Therefore, we can rewrite \eqref{equdecchampbCCC-2} as {\footnotesize 
\begin{subequations}
\begin{align}
&\eqref{equdecchampbCCC-2} = \int_0^t \int \int e^{i s \varphi_3} s \mu(\overline{\xi}, \overline{\eta}, \overline{\sigma}) m_{\widehat{\mathcal{C}}}(\overline{\eta}) m_{\widehat{\mathcal{C}}}(\overline{\sigma}) m_{\widehat{\mathcal{C}}}(\overline{\rho}) \widehat{F}_1(s, \overline{\eta}) \widehat{F}_2(s, \overline{\sigma}) O\left( \overline{\rho}_b^{\overline{\rho}} \right) \widehat{X}_b(\overline{\rho}) \cdot \nabla_{\overline{\xi}} \widehat{F}_3(s, \overline{\rho}) ~ d\overline{\eta} d\overline{\sigma} ds \label{equdecchampbCCC-2-cas1-1} \\
&+ \int_0^t \int \int e^{i s \varphi_3} s \mu(\overline{\xi}, \overline{\eta}, \overline{\sigma}) m_{\widehat{\mathcal{C}}}(\overline{\eta}) m_{\widehat{\mathcal{C}}}(\overline{\sigma}) m_{\widehat{\mathcal{C}}}(\overline{\rho}) \widehat{F}_1(s, \overline{\eta}) \widehat{F}_2(s, \overline{\sigma}) O(1) \widehat{X}_a(\overline{\rho}) \cdot \nabla_{\overline{\xi}} \widehat{F}_3(s, \overline{\rho}) ~ d\overline{\eta} d\overline{\sigma} ds \label{equdecchampbCCC-2-cas1-2} \\
&+ \int_0^t \int \int e^{i s \varphi_3} s \mu(\overline{\xi}, \overline{\eta}, \overline{\sigma}) m_{\widehat{\mathcal{C}}}(\overline{\eta}) m_{\widehat{\mathcal{C}}}(\overline{\sigma}) m_{\widehat{\mathcal{C}}}(\overline{\rho}) \widehat{F}_1(s, \overline{\eta}) \widehat{F}_2(s, \overline{\sigma}) O(1) \widehat{X}_c(\overline{\rho}) \cdot \nabla_{\overline{\xi}} \widehat{F}_3(s, \overline{\rho}) ~ d\overline{\eta} d\overline{\sigma} ds \label{equdecchampbCCC-2-cas1-3} \\
&\begin{aligned}
&+ \int_0^t \int \int e^{i s \varphi_3} s \mu(\overline{\xi}, \overline{\eta}, \overline{\sigma}) m_{\widehat{\mathcal{C}}}(\overline{\eta}) m_{\widehat{\mathcal{C}}}(\overline{\sigma}) m_{\widehat{\mathcal{C}}}(\overline{\rho}) \widehat{F}_1(s, \overline{\eta}) \widehat{F}_2(s, \overline{\sigma}) \\
&\quad \quad \quad \frac{\xi_0 |\xi|}{|\overline{\xi}|^3} \overline{\eta}_b^{\overline{\eta}} \left( \overline{\eta}_a^{\overline{\eta}} + \epsilon^{\overline{\eta} \overline{\sigma}} \overline{\sigma}_a^{\overline{\sigma}} + \epsilon^{\overline{\eta} \overline{\rho}} \overline{\rho}_a^{\overline{\rho}} \right) \widehat{X}_b(\overline{\xi}) \cdot \nabla_{\overline{\xi}} \widehat{F}_3(s, \overline{\rho}) ~ d\overline{\eta} d\overline{\sigma} ds \end{aligned} \label{equdecchampbCCC-2-cas1-4} \\
&\begin{aligned}
&+ \int_0^t \int \int e^{i s \varphi_3} s \mu(\overline{\xi}, \overline{\eta}, \overline{\sigma}) m_{\widehat{\mathcal{C}}}(\overline{\eta}) m_{\widehat{\mathcal{C}}}(\overline{\sigma}) m_{\widehat{\mathcal{C}}}(\overline{\rho}) \widehat{F}_1(s, \overline{\eta}) \widehat{F}_2(s, \overline{\sigma}) \\
&\quad \quad \quad \frac{\xi_0 |\xi|}{|\overline{\xi}|^3} \overline{\sigma}_b^{\overline{\sigma}} \left( \overline{\sigma}_a^{\overline{\sigma}} + \epsilon^{\overline{\sigma} \overline{\eta}} \overline{\eta}_a^{\overline{\eta}} + \epsilon^{\overline{\sigma} \overline{\rho}} \overline{\rho}_a^{\overline{\rho}} \right) \widehat{X}_b(\overline{\xi}) \cdot \nabla_{\overline{\xi}} \widehat{F}_3(s, \overline{\rho}) ~ d\overline{\eta} d\overline{\sigma} ds \end{aligned} \label{equdecchampbCCC-2-cas1-5} \\
&\begin{aligned}
&+ \int_0^t \int \int e^{i s \varphi_3} s \mu(\overline{\xi}, \overline{\eta}, \overline{\sigma}) m_{\widehat{\mathcal{C}}}(\overline{\eta}) m_{\widehat{\mathcal{C}}}(\overline{\sigma}) m_{\widehat{\mathcal{C}}}(\overline{\rho}) \widehat{F}_1(s, \overline{\eta}) \widehat{F}_2(s, \overline{\sigma}) \\
&\quad \quad 2 \frac{\xi_0 |\xi|}{|\overline{\xi}|^3} \left( \epsilon^{\overline{\eta} \overline{\rho}} - \theta^{\overline{\eta} \overline{\rho}} \right) |\eta| |\rho| \frac{\rho_0}{|\rho_0|} P_b^b(\overline{\xi}, \overline{\rho}) \widehat{X}_{b-\widehat{\mathcal{C}}}(\overline{\eta}, \overline{\rho}) \cdot \nabla_{\overline{\xi}} \widehat{F}_3(s, \overline{\rho}) ~ d\overline{\eta} d\overline{\sigma} ds 
\end{aligned} \label{equdecchampbCCC-2-cas1-6} \\
&\begin{aligned}
&+ \int_0^t \int \int e^{i s \varphi_3} s \mu(\overline{\xi}, \overline{\eta}, \overline{\sigma}) m_{\widehat{\mathcal{C}}}(\overline{\eta}) m_{\widehat{\mathcal{C}}}(\overline{\sigma}) m_{\widehat{\mathcal{C}}}(\overline{\rho}) \widehat{F}_1(s, \overline{\eta}) \widehat{F}_2(s, \overline{\sigma}) \\
&\quad \quad 2 \frac{\xi_0 |\xi|}{|\overline{\xi}|^3} \left( \epsilon^{\overline{\sigma} \overline{\rho}} - \theta^{\overline{\sigma} \overline{\rho}} \right) |\sigma| |\rho| \frac{\rho_0}{|\rho_0|} P_b^b(\overline{\xi}, \overline{\rho}) \widehat{X}_{b-\widehat{\mathcal{C}}}(\overline{\sigma}, \overline{\rho}) \cdot \nabla_{\overline{\xi}} \widehat{F}_3(s, \overline{\rho}) ~ d\overline{\eta} d\overline{\sigma} ds 
\end{aligned} \label{equdecchampbCCC-2-cas1-7} \\
&\begin{aligned}
&+ \int_0^t \int \int e^{i s \varphi_3} s \mu(\overline{\xi}, \overline{\eta}, \overline{\sigma}) m_{\widehat{\mathcal{C}}}(\overline{\eta}) m_{\widehat{\mathcal{C}}}(\overline{\sigma}) m_{\widehat{\mathcal{C}}}(\overline{\rho}) \widehat{F}_1(s, \overline{\eta}) \widehat{F}_2(s, \overline{\sigma}) 2 \frac{\xi_0 |\xi|}{|\overline{\xi}|^3} \left( \epsilon^{\overline{\eta} \overline{\sigma}} - \theta^{\overline{\eta} \overline{\sigma}} \right) |\eta| |\sigma| \widehat{X}_{b}(\overline{\xi}) \cdot \nabla_{\overline{\xi}} \widehat{F}_3(s, \overline{\rho}) ~ d\overline{\eta} d\overline{\sigma} ds 
\end{aligned} \label{equdecchampbCCC-2-cas1-8} 
\end{align}
\end{subequations} }
\eqref{equdecchampbCCC-2-cas1-1}, \eqref{equdecchampbCCC-2-cas1-2} and \eqref{equdecchampbCCC-2-cas1-3} are already of the form \eqref{qteslemestimeesgeneriquesdeccubiquehb2}. On \eqref{equdecchampbCCC-2-cas1-4}, we can apply an integration by parts along $\widehat{X}_b(\overline{\xi}) \cdot \nabla_{\overline{\eta}}$, on \eqref{equdecchampbCCC-2-cas1-5} along $\widehat{X}_b(\overline{\xi}) \cdot \nabla_{\overline{\sigma}}$, on \eqref{equdecchampbCCC-2-cas1-6} along $\widehat{X}_{b-\widehat{\mathcal{C}}}(\overline{\rho}, \overline{\eta}) \cdot \nabla_{\overline{\eta}}$, on \eqref{equdecchampbCCC-2-cas1-7} along $\widehat{X}_{b-\widehat{\mathcal{C}}}(\overline{\rho}, \overline{\sigma}) \cdot \nabla_{\overline{\xi}}$. Here we use that $(1 - \epsilon^{\overline{\eta} \overline{\rho}} \theta^{\overline{\eta} \overline{\rho}}) \widehat{X}_{b-\widehat{\mathcal{C}}}(\overline{\rho}, \overline{\eta})$ is not singular. We then get terms from \eqref{qteslemestimeesgeneriquesdeccubiquehb2}, plus {\footnotesize 
\begin{align}
\begin{aligned}
&\int_0^t \int \int i s^2 \frac{\xi_0 |\xi|}{|\overline{\xi}|^3} \Biggl( \overline{\eta}_b^{\overline{\eta}} \left( \overline{\eta}_a^{\overline{\eta}} + \epsilon^{\overline{\eta} \overline{\sigma}} \overline{\sigma}_a^{\overline{\sigma}} + \epsilon^{\overline{\eta} \overline{\rho}} \overline{\rho}_a^{\overline{\rho}} \right) \widehat{X}_b(\overline{\xi}) \cdot \nabla_{\overline{\eta}} + \overline{\sigma}_b^{\overline{\sigma}} \left( \overline{\sigma}_a^{\overline{\sigma}} + \epsilon^{\overline{\sigma} \overline{\eta}} \overline{\eta}_a^{\overline{\eta}} + \epsilon^{\overline{\sigma} \overline{\rho}} \overline{\rho}_a^{\overline{\rho}} \right) \widehat{X}_b(\overline{\xi}) \cdot \nabla_{\overline{\sigma}} \\
&\quad \quad \quad + 2 \left( \epsilon^{\overline{\eta} \overline{\rho}} - \theta^{\overline{\eta} \overline{\rho}} \right) |\eta| |\rho| \frac{\rho_0}{|\rho_0|} P_b^b(\overline{\xi}, \overline{\rho}) \widehat{X}_{b-\widehat{\mathcal{C}}}(\overline{\eta}, \overline{\rho}) \cdot \nabla_{\overline{\eta}} + 2 \left( \epsilon^{\overline{\sigma} \overline{\rho}} - \theta^{\overline{\sigma} \overline{\rho}} \right) |\sigma| |\rho| \frac{\rho_0}{|\rho_0|} P_b^b(\overline{\xi}, \overline{\rho}) \widehat{X}_{b-\widehat{\mathcal{C}}}(\overline{\sigma}, \overline{\rho}) \cdot \nabla_{\overline{\sigma}} \Biggl) \varphi_3 \\
&\quad \quad \quad \quad e^{i s \varphi_3} \mu(\overline{\xi}, \overline{\eta}, \overline{\sigma}) m_{\widehat{\mathcal{C}}}(\overline{\eta}) m_{\widehat{\mathcal{C}}}(\overline{\sigma}) m_{\widehat{\mathcal{C}}}(\overline{\rho}) \widehat{F}_1(s, \overline{\eta}) \widehat{F}_2(s, \overline{\sigma}) \widehat{F}_3(s, \overline{\rho}) ~ d\overline{\eta} d\overline{\sigma} ds
\end{aligned}
\label{equdecchampbCCC-2-cas1-termereste} 
\end{align} }
Finally, for \eqref{equdecchampbCCC-2-cas1-8}, we start by an integration by parts along $\widehat{X}_b(\overline{\xi}) \cdot \nabla_{\overline{\eta}}$, to get terms from \eqref{qteslemestimeesgeneriquesdeccubiquehb2} plus {\footnotesize 
\begin{align*}
&\int_0^t \int \int is^2 2 \frac{\xi_0 |\xi|}{|\overline{\xi}|^3} \left( \epsilon^{\overline{\eta} \overline{\sigma}} - \theta^{\overline{\eta} \overline{\sigma}} \right) |\eta| |\sigma| \widehat{X}_{b}(\overline{\xi}) \cdot \nabla_{\overline{\eta}} \varphi_3 e^{i s \varphi_3} \mu(\overline{\xi}, \overline{\eta}, \overline{\sigma}) m_{\widehat{\mathcal{C}}}(\overline{\eta}) m_{\widehat{\mathcal{C}}}(\overline{\sigma}) m_{\widehat{\mathcal{C}}}(\overline{\rho}) \widehat{F}_1(s, \overline{\eta}) \widehat{F}_2(s, \overline{\sigma}) \widehat{F}_3(s, \overline{\rho}) ~ d\overline{\eta} d\overline{\sigma} ds \\
&\begin{aligned}
&+ \int_0^t \int \int e^{i s \varphi_3} s \mu(\overline{\xi}, \overline{\eta}, \overline{\sigma}) m_{\widehat{\mathcal{C}}}(\overline{\eta}) m_{\widehat{\mathcal{C}}}(\overline{\sigma}) m_{\widehat{\mathcal{C}}}(\overline{\rho}) \widehat{F}_2(s, \overline{\sigma}) \widehat{F}_3(s, \overline{\rho}) \\
&\quad \quad 2 \frac{\xi_0 |\xi|}{|\overline{\xi}|^3} \left( \epsilon^{\overline{\eta} \overline{\sigma}} - \theta^{\overline{\eta} \overline{\sigma}} \right) |\eta| |\sigma| \frac{\eta_0}{|\eta_0|} P_b^b(\overline{\xi}, \overline{\eta}) \widehat{X}_{b-\widehat{\mathcal{C}}}(\overline{\sigma}, \overline{\eta}) \cdot \left( \nabla_{\overline{\eta}} - \nabla_{\overline{\sigma}} \right) \widehat{F}_1(s, \overline{\eta}) ~ d\overline{\eta} d\overline{\sigma} ds
\end{aligned} 
\end{align*} }
On the second term, we can apply an integration by parts along $\widehat{X}_{b-\widehat{\mathcal{C}}}(\overline{\eta}, \overline{\sigma}) \cdot \left( \nabla_{\overline{\eta}} - \nabla_{\overline{\sigma}} \right)$ (using that $\left( \nabla_{\overline{\eta}} - \nabla_{\overline{\sigma}} \right) \overline{\rho} = 0$) and get back to {\footnotesize 
\begin{align}
\begin{aligned}
&\int_0^t \int \int is^2 2 \frac{\xi_0 |\xi|}{|\overline{\xi}|^3} \left( \epsilon^{\overline{\eta} \overline{\sigma}} - \theta^{\overline{\eta} \overline{\sigma}} \right) |\eta| |\sigma| \left( \widehat{X}_{b}(\overline{\xi}) \cdot \nabla_{\overline{\eta}} - \frac{\eta_0}{|\eta_0|} P_b^b(\overline{\xi}, \overline{\eta}) \widehat{X}_{b-\widehat{\mathcal{C}}}(\overline{\sigma}, \overline{\eta}) \cdot \left( \nabla_{\overline{\eta}} - \nabla_{\overline{\sigma}} \right) \right) \varphi_3 \\
&\quad \quad  e^{i s \varphi_3} \mu(\overline{\xi}, \overline{\eta}, \overline{\sigma}) m_{\widehat{\mathcal{C}}}(\overline{\eta}) m_{\widehat{\mathcal{C}}}(\overline{\sigma}) m_{\widehat{\mathcal{C}}}(\overline{\rho}) \widehat{F}_1(s, \overline{\eta}) \widehat{F}_2(s, \overline{\sigma}) \widehat{F}_3(s, \overline{\rho}) ~ d\overline{\eta} d\overline{\sigma} ds 
\end{aligned}
\label{equdecchampbCCC-2-cas1-termereste2} 
\end{align} }
up to well-controlled terms. 

We now group \eqref{equdecchampbCCC-1}, \eqref{equdecchampbCCC-2-cas1-termereste} and \eqref{equdecchampbCCC-2-cas1-termereste2}: {\footnotesize 
\begin{align}
\begin{aligned}
&\int_0^t \int \int i s^2 \frac{\xi_0 |\xi|}{|\overline{\xi}|^3} \Biggl( \left( 3 \xi_0^2 - |\xi|^2 \right) \widehat{X}_b(\overline{\xi}) \cdot \nabla_{\overline{\xi}} + \overline{\eta}_b^{\overline{\eta}} \left( \overline{\eta}_a^{\overline{\eta}} + \epsilon^{\overline{\eta} \overline{\sigma}} \overline{\sigma}_a^{\overline{\sigma}} + \epsilon^{\overline{\eta} \overline{\rho}} \overline{\rho}_a^{\overline{\rho}} \right) \widehat{X}_b(\overline{\xi}) \cdot \nabla_{\overline{\eta}} + \overline{\sigma}_b^{\overline{\sigma}} \left( \overline{\sigma}_a^{\overline{\sigma}} + \epsilon^{\overline{\sigma} \overline{\eta}} \overline{\eta}_a^{\overline{\eta}} + \epsilon^{\overline{\sigma} \overline{\rho}} \overline{\rho}_a^{\overline{\rho}} \right) \widehat{X}_b(\overline{\xi}) \cdot \nabla_{\overline{\sigma}} \\
&\quad \quad + 2 \left( \epsilon^{\overline{\eta} \overline{\rho}} - \theta^{\overline{\eta} \overline{\rho}} \right) |\eta| |\rho| \frac{\rho_0}{|\rho_0|} P_b^b(\overline{\xi}, \overline{\rho}) \widehat{X}_{b-\widehat{\mathcal{C}}}(\overline{\eta}, \overline{\rho}) \cdot \nabla_{\overline{\eta}} + 2 \left( \epsilon^{\overline{\sigma} \overline{\rho}} - \theta^{\overline{\sigma} \overline{\rho}} \right) |\sigma| |\rho| \frac{\rho_0}{|\rho_0|} P_b^b(\overline{\xi}, \overline{\rho}) \widehat{X}_{b-\widehat{\mathcal{C}}}(\overline{\sigma}, \overline{\rho}) \cdot \nabla_{\overline{\sigma}} \\
&\quad \quad + 2 \left( \epsilon^{\overline{\eta} \overline{\sigma}} - \theta^{\overline{\eta} \overline{\sigma}} \right) |\eta| |\sigma| \left( \widehat{X}_{b}(\overline{\xi}) \cdot \nabla_{\overline{\eta}} - \frac{\eta_0}{|\eta_0|} P_b^b(\overline{\xi}, \overline{\eta}) \widehat{X}_{b-\widehat{\mathcal{C}}}(\overline{\sigma}, \overline{\eta}) \cdot \left( \nabla_{\overline{\eta}} - \nabla_{\overline{\sigma}} \right) \right)
\Biggl) \varphi_3 \\
&\quad \quad e^{i s \varphi_3} \mu(\overline{\xi}, \overline{\eta}, \overline{\sigma}) m_{\widehat{\mathcal{C}}}(\overline{\eta}) m_{\widehat{\mathcal{C}}}(\overline{\sigma}) m_{\widehat{\mathcal{C}}}(\overline{\rho}) \widehat{F}_1(s, \overline{\eta}) \widehat{F}_2(s, \overline{\sigma}) \widehat{F}_3(s, \overline{\rho}) ~ d\overline{\eta} d\overline{\sigma} ds
\end{aligned}
\end{align} }

Let us now compute symbols allowing for good integrations by parts. By Lemma \ref{lemcalculchampmodifie}, we have 
\begin{align*}
\overline{\eta}_a^{\overline{\eta}} - \overline{\rho}_a^{\overline{\rho}} &= O\left( \widehat{Y}(\overline{\eta}, \overline{\rho}) \cdot \nabla_{\overline{\eta}} \varphi_3 \right) \\
\overline{\sigma}_a^{\overline{\sigma}} - \overline{\rho}_a^{\overline{\rho}} &= O\left( \widehat{Y}(\overline{\sigma}, \overline{\rho}) \cdot \nabla_{\overline{\eta}} \varphi_3 \right) 
\end{align*}
In particular, we easily localise on $|\overline{\eta}| \simeq |\overline{\rho}| \simeq |\overline{\sigma}|$, and even on a neighborhood where
\begin{align*}
\widehat{Y}(\overline{\eta}, \overline{\rho}) \cdot \nabla_{\overline{\eta}} \varphi_3 &= o(1) \\
\widehat{Y}(\overline{\sigma}, \overline{\rho}) \cdot \nabla_{\overline{\sigma}} \varphi_3 &= o(1)
\end{align*}
Since $\overline{\eta}, \overline{\sigma}, \overline{\rho}$ are close to alignment, we know that $\overline{\xi}_a^{\overline{\xi}} = \sqrt{3} |\xi_0| + |\xi|$ needs to be close to $\overline{\eta}_a^{\overline{\eta}}$ or to $3 \overline{\eta}_a^{\overline{\eta}}$. In particular, in this second case, we have by Lemma \ref{lemcalculssimplescoordonneesconiquevarphi} that
\begin{align*}
6 \sqrt{3} \varphi_3 &= O\left( \overline{\xi}_b^{\overline{\xi}} \right) + O\left( \overline{\eta}_b^{\overline{\eta}} \right) + O\left( \overline{\sigma}_b^{\overline{\sigma}} \right) + O\left( \overline{\rho}_b^{\overline{\rho}} \right) + O\left( \overline{\eta}_a^{\overline{\eta}} - \overline{\rho}_a^{\overline{\rho}} \right) + O\left( \overline{\sigma}_a^{\overline{\sigma}} - \overline{\rho}_a^{\overline{\rho}} \right) 
+ \left( \overline{\xi}_a^{\overline{\xi}} \right)^3 - 3 \left( \overline{\eta}_a^{\overline{\eta}} \right)^3 \\
&= o(1) + 24 \left( \overline{\eta}_a^{\overline{\eta}} \right)^3 
\end{align*}
and so we can apply an integration by parts immediately. Then, up to exchanging $\overline{\eta}, \overline{\sigma}, \overline{\rho}$, we need only to consider the case 
\begin{align*}
\epsilon^{\overline{\eta} \overline{\sigma}} = -1 = \epsilon^{\overline{\sigma} \overline{\rho}}, \quad \epsilon^{\overline{\eta} \overline{\rho}} = 1 = \epsilon^{\overline{\xi} \overline{\rho}}
\end{align*}

Then, by Lemma \ref{lemcalculsconecoordonneesconiquesvarphi}, 
\begin{align*}
\xi_t^{\overline{\eta} \overline{\rho}} &= O\left( \widehat{X}_c(\overline{\eta}) \cdot \nabla_{\overline{\eta}} \varphi_3 \right) 
\end{align*}
We can therefore bound any term having a factor $O\left( \overline{\rho}_b^{\overline{\rho}} \xi_t^{\overline{\eta} \overline{\rho}} \right)$, and symmetric-wise
\begin{align*}
O\left( \overline{\eta}_b^{\overline{\eta}} \xi_t^{\overline{\eta} \overline{\rho}} \right), \quad O\left( \overline{\rho}_b^{\overline{\rho}} \xi_t^{\overline{\sigma} \overline{\rho}} \right), \quad O\left( \overline{\sigma}_b^{\overline{\sigma}} \xi_t^{\overline{\eta} \overline{\sigma}} \right), \quad O\left( \overline{\eta}_b^{\overline{\eta}} \xi_t^{\overline{\eta} \overline{\sigma}} \right), \quad O\left( \overline{\sigma}_b^{\overline{\sigma}} \xi_t^{\overline{\sigma} \overline{\rho}} \right)
\end{align*}
Note now that $\xi_t^{\overline{\eta} \overline{\sigma}} = O\left( \xi_t^{\overline{\eta} \overline{\rho}} \right) + O\left( \xi_t^{\overline{\sigma} \overline{\rho}} \right)$ (since $|\eta|, |\sigma|, |\rho|$ are all of the same order here), so we have good control over any factor being product of some $\xi_t$ with some $\overline{\eta}_b^{\overline{\eta}}, \overline{\sigma}_b^{\overline{\sigma}}, \overline{\rho}_b^{\overline{\rho}}$. In what follows, to simplify the notations, we will denote such symbols by $O(\xi_t m_b)$. 

We can then also compute that {\footnotesize 
\begin{align*}
&\left( \overline{\xi}_a^{\overline{\xi}} \right)^2 + \left( \overline{\xi}_b^{\overline{\xi}} \right)^2 = 6 \xi_0^2 + 2 |\xi|^2 \\
&= 6 \eta_0^2 + 2 |\eta|^2 + 6 \sigma_0^2 + 2 |\sigma|^2 + 6 \rho_0^2 + 2 |\rho|^2 + 12 \eta_0 \sigma_0 + 4 \eta \cdot \sigma + 12 \eta_0 \rho_0 + 4 \eta \cdot \rho + 12 \sigma_0 \rho_0 + 4 \sigma \cdot \rho \\
&= \left( \overline{\eta}_a^{\overline{\eta}} \right)^2 + \left( \overline{\eta}_b^{\overline{\eta}} \right)^2 + \left( \overline{\sigma}_a^{\overline{\sigma}} \right)^2 + \left( \overline{\sigma}_b^{\overline{\sigma}} \right)^2 + \left( \overline{\rho}_a^{\overline{\rho}} \right)^2 + \left( \overline{\rho}_b^{\overline{\rho}} \right)^2 - 2 \overline{\eta}_a^{\overline{\eta}} \overline{\sigma}_a^{\overline{\sigma}} - 2 \overline{\eta}_b^{\overline{\eta}} \overline{\sigma}_b^{\overline{\sigma}} + 4 |\eta| |\sigma| \left( \theta^{\overline{\eta} \overline{\sigma}} + 1 \right) \\
&\quad + 2 \overline{\eta}_a^{\overline{\eta}} \overline{\rho}_a^{\overline{\rho}} + 2 \overline{\eta}_b^{\overline{\eta}} \overline{\rho}_b^{\overline{\rho}} + 4 |\eta| |\rho| \left( \theta^{\overline{\eta} \overline{\rho}} - 1 \right) 
- 2 \overline{\sigma}_a^{\overline{\sigma}} \overline{\rho}_a^{\overline{\rho}} - 2 \overline{\sigma}_b^{\overline{\sigma}} \overline{\rho}_b^{\overline{\rho}} + 4 |\sigma| |\rho| \left( \theta^{\overline{\eta} \overline{\sigma}} + 1 \right) \\
&= O\left( \overline{\eta}_a^{\overline{\eta}} - \overline{\sigma}_a^{\overline{\sigma}} \right) + O\left( \overline{\eta}_a^{\overline{\eta}} - \overline{\rho}_a^{\overline{\rho}} \right) + O(m_b \xi_t) + \left( \overline{\eta}_a^{\overline{\eta}} \right)^2 + \left( \overline{\eta}_b^{\overline{\eta}} - \overline{\sigma}_b^{\overline{\sigma}} + \overline{\rho}_b^{\overline{\rho}} \right)^2
+ \left( \overline{\eta}_a^{\overline{\eta}} \right)^2 \left( 1 + \theta^{\overline{\eta} \overline{\sigma}} + \theta^{\overline{\eta} \overline{\rho}} + \theta^{\overline{\eta} \overline{\sigma}} \right) \\
&2 \overline{\xi}_a^{\overline{\xi}} \overline{\xi}_b^{\overline{\xi}} 
= 6 \xi_0^2 - 2 |\xi|^2 \\
&= 6 \eta_0^2 - 2 |\eta|^2 + 6 \sigma_0^2 - 2 |\sigma|^2 + 6 \rho_0^2 - 2 |\rho|^2 + 12 \eta_0 \sigma_0 - 4 \eta \cdot \sigma + 12 \eta_0 \rho_0 - 4 \eta \cdot \rho + 12 \sigma_0 \rho_0 - 4 \sigma \cdot \rho \\
&= 2 \overline{\eta}_a^{\overline{\eta}} \overline{\eta}_b^{\overline{\eta}} + 2 \overline{\sigma}_a^{\overline{\sigma}} \overline{\sigma}_b^{\overline{\sigma}} + 2 \overline{\rho}_a^{\overline{\rho}} \overline{\rho}_b^{\overline{\rho}} - 2 \overline{\eta}_a^{\overline{\eta}} \overline{\sigma}_b^{\overline{\sigma}} - 2 \overline{\eta}_b^{\overline{\eta}} \overline{\sigma}_a^{\overline{\sigma}} - 4 |\eta| |\sigma| \left( \theta^{\overline{\eta} \overline{\sigma}} + 1 \right) \\
&\quad + 2 \overline{\eta}_a^{\overline{\eta}} \overline{\rho}_b^{\overline{\rho}} + 2 \overline{\eta}_b^{\overline{\eta}} \overline{\rho}_a^{\overline{\rho}} - 4 |\eta| |\rho| \left( \theta^{\overline{\eta} \overline{\rho}} - 1 \right) 
- 2 \overline{\sigma}_a^{\overline{\sigma}} \overline{\rho}_b^{\overline{\rho}} - 2 \overline{\sigma}_b^{\overline{\sigma}} \overline{\rho}_a^{\overline{\rho}} - 4 |\sigma| |\rho| \left( \theta^{\overline{\sigma} \overline{\rho}} + 1 \right) \\
&= O\left( \overline{\eta}_a^{\overline{\eta}} - \overline{\sigma}_a^{\overline{\sigma}} \right) + O\left( \overline{\eta}_a^{\overline{\eta}} - \overline{\rho}_a^{\overline{\rho}} \right) + O(m_b \xi_t) \\
&\quad + 2 \overline{\eta}_a^{\overline{\eta}} \left( \overline{\eta}_b^{\overline{\eta}} - \overline{\sigma}_b^{\overline{\sigma}} + \overline{\rho}_b^{\overline{\rho}} \right) 
- \left( \overline{\eta}_a^{\overline{\eta}} \right)^2 \left( 1 + \theta^{\overline{\eta} \overline{\sigma}} + \theta^{\overline{\eta} \overline{\rho}} + \theta^{\overline{\eta} \overline{\sigma}} \right)
\end{align*} }
Therefore, 
\begin{align*}
\overline{\xi}_a^{\overline{\xi}} &= O\left( \overline{\eta}_a^{\overline{\eta}} - \overline{\sigma}_a^{\overline{\sigma}} \right) + O\left( \overline{\eta}_a^{\overline{\eta}} - \overline{\rho}_a^{\overline{\rho}} \right) + O(m_b \xi_t) \\
&\quad + \frac{1}{2} \left( \overline{\eta}_a^{\overline{\eta}} + \overline{\eta}_b^{\overline{\eta}} - \overline{\sigma}_b^{\overline{\sigma}} + \overline{\rho}_b^{\overline{\rho}} \right) 
+ \frac{1}{2} \left( \left( \overline{\eta}_a^{\overline{\eta}} - \overline{\eta}_b^{\overline{\eta}} + \overline{\sigma}_b^{\overline{\sigma}} - \overline{\rho}_b^{\overline{\rho}} \right)^2
+ 2 \left( \overline{\eta}_a^{\overline{\eta}} \right)^2 \left( 1 + \theta^{\overline{\eta} \overline{\sigma}} + \theta^{\overline{\eta} \overline{\rho}} + \theta^{\overline{\eta} \overline{\sigma}} \right) \right)^{\frac{1}{2}} \\
&= O\left( \overline{\eta}_a^{\overline{\eta}} - \overline{\sigma}_a^{\overline{\sigma}} \right) + O\left( \overline{\eta}_a^{\overline{\eta}} - \overline{\rho}_a^{\overline{\rho}} \right) + O(m_b \xi_t) + O\left( \left( 1 + \theta^{\overline{\eta} \overline{\sigma}} + \theta^{\overline{\eta} \overline{\rho}} + \theta^{\overline{\eta} \overline{\sigma}} \right)^2 \right) \\
&\quad + \frac{1}{2} \left( \overline{\eta}_a^{\overline{\eta}} + \overline{\eta}_b^{\overline{\eta}} - \overline{\sigma}_b^{\overline{\sigma}} + \overline{\rho}_b^{\overline{\rho}} \right) 
+ \frac{1}{2} \left( \overline{\eta}_a^{\overline{\eta}} - \overline{\eta}_b^{\overline{\eta}} + \overline{\sigma}_b^{\overline{\sigma}} - \overline{\rho}_b^{\overline{\rho}} \right)
+ \frac{1}{2} \overline{\eta}_a^{\overline{\eta}} \left( 1 + \theta^{\overline{\eta} \overline{\sigma}} + \theta^{\overline{\eta} \overline{\rho}} + \theta^{\overline{\eta} \overline{\sigma}} \right) \\
&= O\left( \overline{\eta}_a^{\overline{\eta}} - \overline{\sigma}_a^{\overline{\sigma}} \right) + O\left( \overline{\eta}_a^{\overline{\eta}} - \overline{\rho}_a^{\overline{\rho}} \right) + O(m_b \xi_t) + O\left( \left( 1 + \theta^{\overline{\eta} \overline{\sigma}} + \theta^{\overline{\eta} \overline{\rho}} + \theta^{\overline{\eta} \overline{\sigma}} \right)^2 \right) \\
&\quad + \overline{\eta}_a^{\overline{\eta}} 
+ \frac{1}{2} \overline{\eta}_a^{\overline{\eta}} \left( 1 + \theta^{\overline{\eta} \overline{\sigma}} + \theta^{\overline{\eta} \overline{\rho}} + \theta^{\overline{\eta} \overline{\sigma}} \right) \\
\overline{\xi}_b^{\overline{\xi}} &= O\left( \overline{\eta}_a^{\overline{\eta}} - \overline{\sigma}_a^{\overline{\sigma}} \right) + O\left( \overline{\eta}_a^{\overline{\eta}} - \overline{\rho}_a^{\overline{\rho}} \right) + O(m_b \xi_t) + O\left( \left( 1 + \theta^{\overline{\eta} \overline{\sigma}} + \theta^{\overline{\eta} \overline{\rho}} + \theta^{\overline{\eta} \overline{\sigma}} \right) \right) \\
&\quad + \overline{\eta}_b^{\overline{\eta}} - \overline{\sigma}_b^{\overline{\sigma}} + \overline{\rho}_b^{\overline{\rho}} 
\end{align*}
In particular, 
\begin{align*}
\left( \overline{\xi}_b^{\overline{\xi}} \right)^3 &= \left( \overline{\eta}_b^{\overline{\eta}} - \overline{\sigma}_b^{\overline{\sigma}} + \overline{\rho}_b^{\overline{\rho}} \right)^3 + O\left( \left( 1 + \theta^{\overline{\eta} \overline{\sigma}} + \theta^{\overline{\eta} \overline{\rho}} + \theta^{\overline{\eta} \overline{\sigma}} \right)^2 \right) \\
&\quad + O\left( \overline{\eta}_a^{\overline{\eta}} - \overline{\sigma}_a^{\overline{\sigma}} \right) + O\left( \overline{\eta}_a^{\overline{\eta}} - \overline{\rho}_a^{\overline{\rho}} \right) + O(m_b \xi_t)
\end{align*}

Now, we compute 
\begin{align*}
6 \sqrt{3} \frac{\eta_0}{|\eta_0|} \varphi_3
&= \left( \overline{\xi}_a^{\overline{\xi}} \right)^3 + \left( \overline{\xi}_b^{\overline{\xi}} \right)^3 - \left( \overline{\eta}_a^{\overline{\eta}} \right)^3 - \left( \overline{\eta}_b^{\overline{\eta}} \right)^3 + \left( \overline{\sigma}_a^{\overline{\sigma}} \right)^3 + \left( \overline{\sigma}_b^{\overline{\sigma}} \right)^3 - \left( \overline{\rho}_a^{\overline{\rho}} \right)^3 - \left( \overline{\rho}_b^{\overline{\rho}} \right)^3 \\
&= O\left( \overline{\eta}_a^{\overline{\eta}} - \overline{\sigma}_a^{\overline{\sigma}} \right) + O\left( \overline{\eta}_a^{\overline{\eta}} - \overline{\rho}_a^{\overline{\rho}} \right) + O(m_b \xi_t) \\
&\quad + \frac{3}{2} \left( \overline{\eta}_a^{\overline{\eta}} \right)^3 \left( 1 + \theta^{\overline{\eta} \overline{\sigma}} + \theta^{\overline{\eta} \overline{\rho}} + \theta^{\overline{\eta} \overline{\sigma}} \right) \left( 1 + O\left( 1 + \theta^{\overline{\eta} \overline{\sigma}} + \theta^{\overline{\eta} \overline{\rho}} + \theta^{\overline{\eta} \overline{\sigma}} \right) \right) \\
&\quad + \left( \overline{\xi}_b^{\overline{\xi}} \right)^3 - \left( \overline{\eta}_b^{\overline{\eta}} \right)^3 + \left( \overline{\sigma}_b^{\overline{\sigma}} \right)^3 - \left( \overline{\rho}_b^{\overline{\rho}} \right)^3
\end{align*}
In particular, we deduce that 
\begin{align*}
&\frac{3}{2} \left( \overline{\eta}_a^{\overline{\eta}} \right)^3 \left( 1 + \theta^{\overline{\eta} \overline{\sigma}} + \theta^{\overline{\eta} \overline{\rho}} + \theta^{\overline{\eta} \overline{\sigma}} \right) \left( 1 + O\left( 1 + \theta^{\overline{\eta} \overline{\sigma}} + \theta^{\overline{\eta} \overline{\rho}} + \theta^{\overline{\eta} \overline{\sigma}} \right) \right) \\
&= O\left( \overline{\eta}_a^{\overline{\eta}} - \overline{\sigma}_a^{\overline{\sigma}} \right) + O\left( \overline{\eta}_a^{\overline{\eta}} - \overline{\rho}_a^{\overline{\rho}} \right) + O(m_b \xi_t) + O(\varphi_3) \\
&\quad - \left( \left( \overline{\eta}_b^{\overline{\eta}} - \overline{\sigma}_b^{\overline{\sigma}} + \overline{\rho}_b^{\overline{\rho}} \right)^3 - \left( \overline{\eta}_b^{\overline{\eta}} \right)^3 + \left( \overline{\sigma}_b^{\overline{\sigma}} \right)^3 - \left( \overline{\rho}_b^{\overline{\rho}} \right)^3 \right)
\end{align*}
In particular, 
\begin{align*}
\left( 1 + \theta^{\overline{\eta} \overline{\sigma}} + \theta^{\overline{\eta} \overline{\rho}} + \theta^{\overline{\eta} \overline{\sigma}} \right)^2 &= O\left( \overline{\eta}_a^{\overline{\eta}} - \overline{\sigma}_a^{\overline{\sigma}} \right) + O\left( \overline{\eta}_a^{\overline{\eta}} - \overline{\rho}_a^{\overline{\rho}} \right) + O(m_b \xi_t) + O(\varphi_3)
\end{align*}
and therefore 
\begin{align*}
\left( \overline{\xi}_b^{\overline{\xi}} \right)^3 &= \left( \overline{\eta}_b^{\overline{\eta}} - \overline{\sigma}_b^{\overline{\sigma}} + \overline{\rho}_b^{\overline{\rho}} \right)^3 + O\left( \overline{\eta}_a^{\overline{\eta}} - \overline{\sigma}_a^{\overline{\sigma}} \right) + O\left( \overline{\eta}_a^{\overline{\eta}} - \overline{\rho}_a^{\overline{\rho}} \right) + O(m_b \xi_t) + O(\varphi_3) 
\end{align*}
Moreover, any quantity having the form $\overline{\xi}_b^{\overline{\xi}} \xi_t^{\alpha \beta}$ for any $\alpha, \beta$ can be rewritten as $O(m_b \xi_t)$, $O(\varphi_3)$ or $O\left( \overline{\eta}_a^{\overline{\eta}} - \overline{\sigma}_a^{\overline{\sigma}} \right) + O\left( \overline{\eta}_a^{\overline{\eta}} - \overline{\rho}_a^{\overline{\rho}} \right)$. 

Finally, we compute {\footnotesize 
\begin{align*}
&\Biggl( \left( 3 \xi_0^2 - |\xi|^2 \right) \widehat{X}_b(\overline{\xi}) \cdot \nabla_{\overline{\xi}} + \overline{\eta}_b^{\overline{\eta}} \left( \overline{\eta}_a^{\overline{\eta}} + \epsilon^{\overline{\eta} \overline{\sigma}} \overline{\sigma}_a^{\overline{\sigma}} + \epsilon^{\overline{\eta} \overline{\rho}} \overline{\rho}_a^{\overline{\rho}} \right) \widehat{X}_b(\overline{\xi}) \cdot \nabla_{\overline{\eta}} + \overline{\sigma}_b^{\overline{\sigma}} \left( \overline{\sigma}_a^{\overline{\sigma}} + \epsilon^{\overline{\sigma} \overline{\eta}} \overline{\eta}_a^{\overline{\eta}} + \epsilon^{\overline{\sigma} \overline{\rho}} \overline{\rho}_a^{\overline{\rho}} \right) \widehat{X}_b(\overline{\xi}) \cdot \nabla_{\overline{\sigma}} \\
&\quad \quad + 2 \left( \epsilon^{\overline{\eta} \overline{\rho}} - \theta^{\overline{\eta} \overline{\rho}} \right) |\eta| |\rho| \frac{\rho_0}{|\rho_0|} P_b^b(\overline{\xi}, \overline{\rho}) \widehat{X}_{b-\widehat{\mathcal{C}}}(\overline{\eta}, \overline{\rho}) \cdot \nabla_{\overline{\eta}} + 2 \left( \epsilon^{\overline{\sigma} \overline{\rho}} - \theta^{\overline{\sigma} \overline{\rho}} \right) |\sigma| |\rho| \frac{\rho_0}{|\rho_0|} P_b^b(\overline{\xi}, \overline{\rho}) \widehat{X}_{b-\widehat{\mathcal{C}}}(\overline{\sigma}, \overline{\rho}) \cdot \nabla_{\overline{\sigma}} \\
&\quad \quad + 2 \left( \epsilon^{\overline{\eta} \overline{\sigma}} - \theta^{\overline{\eta} \overline{\sigma}} \right) |\eta| |\sigma| \left( \widehat{X}_{b}(\overline{\xi}) \cdot \nabla_{\overline{\eta}} - \frac{\eta_0}{|\eta_0|} P_b^b(\overline{\xi}, \overline{\eta}) \widehat{X}_{b-\widehat{\mathcal{C}}}(\overline{\sigma}, \overline{\eta}) \cdot \left( \nabla_{\overline{\eta}} - \nabla_{\overline{\sigma}} \right) \right)
\Biggl) \varphi_3 \\
&= \frac{1}{|\overline{\xi}| |\xi|} \Biggl( \overline{\xi}_a^{\overline{\xi}} \overline{\xi}_b^{\overline{\xi}} \left( |\xi|^2 \left( 3 \xi_0^2 + |\xi|^2 - 3 \rho_0^2 - |\rho|^2 \right) - \xi_0 \xi \cdot \left( 2 \xi_0 \xi - 2 \rho_0 \rho \right) \right) \\
&\quad + \overline{\eta}_b^{\overline{\eta}} \left( \overline{\eta}_a^{\overline{\eta}} - \overline{\sigma}_a^{\overline{\sigma}} + \overline{\rho}_a^{\overline{\rho}} \right) 
\left( |\xi|^2 \left( 3 \rho_0^2 + |\rho|^2 - 3 \eta_0^2 - |\eta|^2 \right) - \xi_0 \xi \cdot \left( 2 \rho_0 \rho - 2 \eta_0 \eta \right) \right)  \\
&\quad + \overline{\sigma}_b^{\overline{\sigma}} \left( \overline{\sigma}_a^{\overline{\sigma}} - \overline{\eta}_a^{\overline{\eta}} - \overline{\rho}_a^{\overline{\rho}} \right) 
\left( |\xi|^2 \left( 3 \rho_0^2 + |\rho|^2 - 3 \sigma_0^2 - |\sigma|^2 \right) - \xi_0 \xi \cdot \left( 2 \rho_0 \rho - 2 \sigma_0 \sigma \right) \right) \\
&\quad \quad + 2 \left( 1 - \theta^{\overline{\eta} \overline{\rho}} \right) |\eta| |\rho| \frac{|\xi|^2 |\rho|^2 + \xi_0 \rho_0 \xi \cdot \rho}{|\overline{\rho}| |\rho|} \frac{2 \sqrt{3}}{7 + \theta^{\overline{\eta} \overline{\rho}}} \left( \left( \overline{\rho}_b^{\overline{\rho}} \right)^2 - \left( \overline{\eta}_a^{\overline{\eta}} \right)^2 - \frac{1 + \theta^{\overline{\eta} \overline{\rho}}}{6} \left( \left( \overline{\rho}_a^{\overline{\rho}} \right)^2 - \left( \overline{\rho}_b^{\overline{\rho}} \right)^2 \right) \right) \\
&\quad \quad - 2 \left( 1 + \theta^{\overline{\sigma} \overline{\rho}} \right) |\sigma| |\rho| \frac{|\xi|^2 |\rho|^2 + \xi_0 \rho_0 \xi \cdot \rho}{|\overline{\rho}| |\rho|} 
\frac{2 \sqrt{3}}{7 - \theta^{\overline{\sigma} \overline{\rho}}} \left( \left( \overline{\rho}_b^{\overline{\rho}} \right)^2 - \left( \overline{\sigma}_a^{\overline{\sigma}} \right)^2 - \frac{1 - \theta^{\overline{\sigma} \overline{\rho}}}{6} \left( \left( \overline{\rho}_a^{\overline{\rho}} \right)^2 - \left( \overline{\rho}_b^{\overline{\rho}} \right)^2 \right) \right) \\
&\quad \quad - 2 \left( 1 + \theta^{\overline{\eta} \overline{\sigma}} \right) |\eta| |\sigma| \left( 
|\xi|^2 \left( 3 \rho_0^2 + |\rho|^2 - 3 \eta_0^2 - |\eta|^2 \right) - \xi_0 \xi \cdot \left( 2 \rho_0 \rho - 2 \eta_0 \eta \right) \right) \\
&\quad \quad - 2 \left( 1 + \theta^{\overline{\eta} \overline{\sigma}} \right) |\eta| |\sigma| \frac{|\xi|^2 |\eta|^2 + \xi_0 \eta_0 \xi \cdot \eta}{|\overline{\eta}| |\eta|} 
\frac{2 \sqrt{3}}{7 - \theta^{\overline{\eta} \overline{\sigma}}} \left( \left( \overline{\eta}_b^{\overline{\eta}} \right)^2 - \left( \overline{\sigma}_a^{\overline{\sigma}} \right)^2 - \frac{1 - \theta^{\overline{\eta} \overline{\sigma}}}{6} \left( \left( \overline{\eta}_a^{\overline{\eta}} \right)^2 - \left( \overline{\eta}_b^{\overline{\eta}} \right)^2 \right) \right) 
\Biggl) \\
&= O\left( \overline{\eta}_a^{\overline{\eta}} - \overline{\sigma}_a^{\overline{\sigma}} \right) + O\left( \overline{\eta}_a^{\overline{\eta}} - \overline{\rho}_a^{\overline{\rho}} \right) + O(m_b \xi_t) + O(\varphi_3) \\
&\quad + \frac{1}{|\overline{\xi}| |\xi|} \Biggl( \overline{\eta}_a^{\overline{\eta}} \overline{\xi}_b^{\overline{\xi}} \left( |\xi|^2 \frac{1}{2} \left( \left( \overline{\xi}_a^{\overline{\xi}} \right)^2 + \left( \overline{\xi}_b^{\overline{\xi}} \right)^2 - \left( \overline{\rho}_a^{\overline{\rho}} \right)^2 - \left( \overline{\rho}_b^{\overline{\rho}} \right)^2 \right) - \frac{1}{2 \sqrt{3}} |\xi_0| |\xi| \left( \left( \overline{\xi}_a^{\overline{\xi}} \right)^2 - \left( \overline{\xi}_b^{\overline{\xi}} \right)^2 - \left( \overline{\rho}_a^{\overline{\rho}} \right)^2 + \left( \overline{\rho}_b^{\overline{\rho}} \right)^2 \right) \right) \\
&\quad + \overline{\eta}_b^{\overline{\eta}} \overline{\eta}_a^{\overline{\eta}}
\left( |\xi|^2 \frac{1}{2} \left( \left( \overline{\rho}_a^{\overline{\rho}} \right)^2 + \left( \overline{\rho}_b^{\overline{\rho}} \right)^2 - \left( \overline{\eta}_a^{\overline{\eta}} \right)^2 - \left( \overline{\eta}_b^{\overline{\eta}} \right)^2 \right) - \frac{1}{2 \sqrt{3}} |\xi_0| |\xi| \left( \left( \overline{\rho}_a^{\overline{\rho}} \right)^2 - \left( \overline{\rho}_b^{\overline{\rho}} \right)^2 - \left( \overline{\eta}_a^{\overline{\eta}} \right)^2 + \left( \overline{\eta}_b^{\overline{\eta}} \right)^2 \right) \right)  \\
&\quad - \overline{\sigma}_b^{\overline{\sigma}} \overline{\eta}_a^{\overline{\eta}} 
\left( |\xi|^2 \frac{1}{2} \left( \left( \overline{\rho}_a^{\overline{\rho}} \right)^2 + \left( \overline{\rho}_b^{\overline{\rho}} \right)^2 - \left( \overline{\sigma}_a^{\overline{\sigma}} \right)^2 - \left( \overline{\sigma}_b^{\overline{\sigma}} \right)^2 \right) - \frac{1}{2 \sqrt{3}} |\xi_0| |\xi| \left( \left( \overline{\rho}_a^{\overline{\rho}} \right)^2 - \left( \overline{\rho}_b^{\overline{\rho}} \right)^2 - \left( \overline{\sigma}_a^{\overline{\sigma}} \right)^2 + \left( \overline{\sigma}_b^{\overline{\sigma}} \right)^2 \right) \right) \\
&\quad \quad - \frac{1}{2} \left( 1 - \theta^{\overline{\eta} \overline{\rho}} \right) \left( \overline{\eta}_a^{\overline{\eta}} \right)^2 \frac{\frac{1}{16} \left( \overline{\eta}_a^{\overline{\eta}} \right)^4 + \frac{1}{48} \left( \overline{\eta}_a^{\overline{\eta}} \right)^4 \left( \theta^{\overline{\eta} \overline{\rho}} + \theta^{\overline{\sigma} \overline{\rho}} + 1 \right)}{\frac{1}{2 \sqrt{3}} \left( \overline{\eta}_a^{\overline{\eta}} \right)^2} 2 \sqrt{3} \frac{1}{6} \left( \overline{\eta}_a^{\overline{\eta}} \right)^2 \\
&\quad \quad + \frac{1}{2} \left( 1 + \theta^{\overline{\sigma} \overline{\rho}} \right) \left( \overline{\eta}_a^{\overline{\eta}} \right)^2 \frac{\frac{1}{16} \left( \overline{\eta}_a^{\overline{\eta}} \right)^4 + \frac{1}{48} \left( \overline{\eta}_a^{\overline{\eta}} \right)^4 \left( \theta^{\overline{\eta} \overline{\rho}} + \theta^{\overline{\sigma} \overline{\rho}} + 1 \right)}{\frac{1}{2 \sqrt{3}} \left( \overline{\eta}_a^{\overline{\eta}} \right)^2} 
2 \sqrt{3} \frac{1}{6} \left( \overline{\eta}_a^{\overline{\eta}} \right)^2 \\
&\quad \quad - \frac{1}{2} \left( 1 + \theta^{\overline{\eta} \overline{\sigma}} \right) \left( \overline{\eta}_a^{\overline{\eta}} \right)^2 \left(
|\xi|^2 \frac{1}{2} \left( \left( \overline{\rho}_a^{\overline{\rho}} \right)^2 + \left( \overline{\rho}_b^{\overline{\rho}} \right)^2 - \left( \overline{\eta}_a^{\overline{\eta}} \right)^2 - \left( \overline{\eta}_b^{\overline{\eta}} \right)^2 \right) - \frac{1}{24} \left( \overline{\eta}_a^{\overline{\eta}} \right)^4 \left( \theta^{\overline{\sigma} \overline{\rho}} - \theta^{\overline{\sigma} \overline{\eta}} \right) \right) \\
&\quad \quad + \frac{1}{2} \left( 1 + \theta^{\overline{\eta} \overline{\sigma}} \right) \left( \overline{\eta}_a^{\overline{\eta}} \right)^2 \frac{\frac{1}{16} \left( \overline{\eta}_a^{\overline{\eta}} \right)^4 + \frac{1}{48} \left( \overline{\eta}_a^{\overline{\eta}} \right)^4 \left( 1 + \theta^{\overline{\eta} \overline{\sigma}} + \theta^{\overline{\eta} \overline{\rho}} \right)}{\frac{1}{2 \sqrt{3}} \left( \overline{\eta}_a^{\overline{\eta}} \right)^2} 
2 \sqrt{3} \frac{1}{6} \left( \overline{\eta}_a^{\overline{\eta}} \right)^2 
\Biggl) \\
&= O\left( \overline{\eta}_a^{\overline{\eta}} - \overline{\sigma}_a^{\overline{\sigma}} \right) + O\left( \overline{\eta}_a^{\overline{\eta}} - \overline{\rho}_a^{\overline{\rho}} \right) + O(m_b \xi_t) + O(\varphi_3) \\
&\quad + \frac{1}{|\overline{\xi}| |\xi|} \Biggl( \frac{1}{6} \left( \overline{\eta}_a^{\overline{\eta}} \right)^3 \overline{\xi}_b^{\overline{\xi}} \left( \left( \overline{\xi}_b^{\overline{\xi}} \right)^2 - \left( \overline{\rho}_b^{\overline{\rho}} \right)^2 \right) + \frac{1}{6} \overline{\eta}_b^{\overline{\eta}} \left( \overline{\eta}_a^{\overline{\eta}} \right)^3
\left( \left( \overline{\rho}_b^{\overline{\rho}} \right)^2 - \left( \overline{\eta}_b^{\overline{\eta}} \right)^2 \right) - \frac{1}{6} \overline{\sigma}_b^{\overline{\sigma}} \left( \overline{\eta}_a^{\overline{\eta}} \right)^3 
\left( \left( \overline{\rho}_b^{\overline{\rho}} \right)^2 - \left( \overline{\sigma}_b^{\overline{\sigma}} \right)^2 \right) \\
&\quad \quad - \frac{1}{48} \left( 1 - \theta^{\overline{\eta} \overline{\rho}} \right) \left( \overline{\eta}_a^{\overline{\eta}} \right)^6 \left( 3 + \left( \theta^{\overline{\eta} \overline{\rho}} + \theta^{\overline{\sigma} \overline{\rho}} + 1 \right)\right)  + \frac{1}{48} \left( 1 + \theta^{\overline{\sigma} \overline{\rho}} \right) \left( \overline{\eta}_a^{\overline{\eta}} \right)^6 \left( 3 + \left( \theta^{\overline{\eta} \overline{\rho}} + \theta^{\overline{\sigma} \overline{\rho}} + 1 \right) \right) \\
&\quad \quad + \frac{1}{48} \left( 1 + \theta^{\overline{\eta} \overline{\sigma}} \right) \left( \overline{\eta}_a^{\overline{\eta}} \right)^6 \left( \theta^{\overline{\sigma} \overline{\rho}} - \theta^{\overline{\sigma} \overline{\eta}} \right) + \frac{1}{48} \left( 1 + \theta^{\overline{\eta} \overline{\sigma}} \right) \left( \overline{\eta}_a^{\overline{\eta}} \right)^6 \left( 3 + \left( 1 + \theta^{\overline{\eta} \overline{\sigma}} + \theta^{\overline{\eta} \overline{\rho}} \right) \right) 
\Biggl) \\
&= O\left( \overline{\eta}_a^{\overline{\eta}} - \overline{\sigma}_a^{\overline{\sigma}} \right) + O\left( \overline{\eta}_a^{\overline{\eta}} - \overline{\rho}_a^{\overline{\rho}} \right) + O(m_b \xi_t) + O(\varphi_3) \\
&\quad + \frac{1}{48 |\overline{\xi}| |\xi|} \left( \overline{\eta}_a^{\overline{\eta}} \right)^3 \Biggl( 
8 \left( \left( \overline{\xi}_b^{\overline{\xi}} \right)^3 - \left( \overline{\eta}_b^{\overline{\eta}} \right)^3 + \left( \overline{\sigma}_b^{\overline{\sigma}} \right)^3 - \left( \overline{\rho}_b^{\overline{\rho}} \right)^2 \left( \overline{\xi}_b^{\overline{\xi}} - \overline{\eta}_b^{\overline{\eta}} + \overline{\sigma}_b^{\overline{\sigma}} \right) \right) \\
&\quad \quad + \left( \overline{\eta}_a^{\overline{\eta}} \right)^3 \left( 3 + \left( \theta^{\overline{\eta} \overline{\rho}} + \theta^{\overline{\sigma} \overline{\rho}} + 1 \right) \right) \left( 1 + \theta^{\overline{\eta} \overline{\sigma}} + \theta^{\overline{\eta} \overline{\rho}} + \theta^{\overline{\sigma} \overline{\rho}} \right) 
\Biggl) \\
&= O\left( \overline{\eta}_a^{\overline{\eta}} - \overline{\sigma}_a^{\overline{\sigma}} \right) + O\left( \overline{\eta}_a^{\overline{\eta}} - \overline{\rho}_a^{\overline{\rho}} \right) + O(m_b \xi_t) + O(\varphi_3) \\
&\quad + \frac{1}{9 |\overline{\xi}| |\xi|} \left( \overline{\eta}_a^{\overline{\eta}} \right)^3 \left( \left( \overline{\xi}_b^{\overline{\xi}} \right)^3 - \left( \overline{\eta}_b^{\overline{\eta}} \right)^3 + \left( \overline{\sigma}_b^{\overline{\sigma}} \right)^3 - \left( \overline{\rho}_b^{\overline{\rho}} \right)^3 \right) 
\end{align*} }

In particular, up to terms from \eqref{qteslemestimeesgeneriquesdeccubiquehb2}, we reduced to
\begin{align*}
&\int_0^t \int \int e^{i s \varphi_3} s^2 m_{\widehat{\mathcal{C}}}(\overline{\xi}) \mu_{BBBB}(\overline{\xi}, \overline{\eta}, \overline{\sigma}) \mu(\overline{\xi}, \overline{\eta}, \overline{\sigma}) \left( \left( \overline{\xi}_b^{\overline{\xi}} \right)^3 - \left( \overline{\eta}_b^{\overline{\eta}} \right)^3 + \left( \overline{\sigma}_b^{\overline{\sigma}} \right)^3 - \left( \overline{\rho}_b^{\overline{\rho}} \right)^3 \right) \\
&\quad \quad \quad m_{\widehat{\mathcal{C}}}(\overline{\eta}) m_{\widehat{\mathcal{C}}}(\overline{\sigma}) m_{\widehat{\mathcal{C}}}(\overline{\rho}) \widehat{F}_1(s, \overline{\eta}) \widehat{F}_2(s, \overline{\sigma}) \widehat{F}_3(s, \overline{\rho}) ~ d\overline{\eta} d\overline{\sigma} ds
\end{align*}

We now explain how to get a decomposition for the single remaining term, having the above form. To that end, we add a localisation $\mu_{sing}$ which is singular at $\Gamma = \{ (\overline{\xi}, \overline{\eta}, \overline{\sigma}), \overline{\xi} = \overline{\eta} = - \overline{\sigma} \}$, bounded but with derivatives growing like $\frac{1}{d(\cdot, \Gamma)}$. Note now that 
\begin{align*}
d\left( (\overline{\xi}, \overline{\eta}, \overline{\sigma}), \Gamma \right) &\simeq \left| \overline{\eta}_a^{\overline{\eta}} - \overline{\sigma}_a^{\overline{\sigma}} \right| + \left| \overline{\eta}_a^{\overline{\eta}} - \overline{\rho}_a^{\overline{\rho}} \right| + \overline{\eta}_a^{\overline{\eta}} \left| \xi_t^{\overline{\eta} \overline{\sigma}} \right| + \overline{\eta}_a^{\overline{\eta}} \left| \xi_t^{\overline{\eta} \overline{\rho}} \right| + \left| \overline{\eta}_b^{\overline{\eta}} - \overline{\rho}_b^{\overline{\rho}} \right| + \left| \overline{\eta}_b^{\overline{\eta}} - \overline{\sigma}_b^{\overline{\sigma}} \right| 
\end{align*}
In particular, we can choose $\mu_{sing}$ localising depending on which of these scalare quantities is the largest. The boundedness of the associated multilinear operator will follow from Theorem \ref{theoMTTmultilinsingGamma}. 

We only rewrite {\footnotesize 
\begin{align*}
&\left( \overline{\xi}_b^{\overline{\xi}} \right)^3 - \left( \overline{\eta}_b^{\overline{\eta}} \right)^3 + \left( \overline{\sigma}_b^{\overline{\sigma}} \right)^3 - \left( \overline{\rho}_b^{\overline{\rho}} \right)^3 \\
&\quad = O\left( \overline{\eta}_a^{\overline{\eta}} - \overline{\sigma}_a^{\overline{\sigma}} \right) + O\left( \overline{\eta}_a^{\overline{\eta}} - \overline{\rho}_a^{\overline{\rho}} \right) + O(m_b \xi_t) + O(\varphi_3) + \left( \overline{\eta}_b^{\overline{\eta}} - \overline{\sigma}_b^{\overline{\sigma}} + \overline{\rho}_b^{\overline{\rho}} \right)^3 - \left( \overline{\eta}_b^{\overline{\eta}} \right)^3 + \left( \overline{\sigma}_b^{\overline{\sigma}} \right)^3 - \left( \overline{\rho}_b^{\overline{\rho}} \right)^3 \\
&\quad = O\left( \overline{\eta}_a^{\overline{\eta}} - \overline{\sigma}_a^{\overline{\sigma}} \right) + O\left( \overline{\eta}_a^{\overline{\eta}} - \overline{\rho}_a^{\overline{\rho}} \right) + O(m_b \xi_t) + O(\varphi_3) \\
&\quad \quad + \left( \overline{\eta}_b^{\overline{\eta}} - \overline{\sigma}_b^{\overline{\sigma}} \right) 
\left( \left( \overline{\eta}_b^{\overline{\eta}} - \overline{\sigma}_b^{\overline{\sigma}} + \overline{\rho}_b^{\overline{\rho}} \right)^2 + \left( \overline{\eta}_b^{\overline{\eta}} - \overline{\sigma}_b^{\overline{\sigma}} + \overline{\rho}_b^{\overline{\rho}} \right) \overline{\rho}_b^{\overline{\rho}} + \left( \overline{\rho}_b^{\overline{\rho}} \right)^2 \right) 
- \left( \overline{\eta}_b^{\overline{\eta}} - \overline{\sigma}_b^{\overline{\sigma}} \right) \left( \left( \overline{\eta}_b^{\overline{\eta}} \right)^2 + \overline{\eta}_b^{\overline{\eta}} \overline{\sigma}_b^{\overline{\sigma}} + \left( \overline{\sigma}_b^{\overline{\sigma}} \right)^2 \right) \\
&\quad = O\left( \overline{\eta}_a^{\overline{\eta}} - \overline{\sigma}_a^{\overline{\sigma}} \right) + O\left( \overline{\eta}_a^{\overline{\eta}} - \overline{\rho}_a^{\overline{\rho}} \right) + O(m_b \xi_t) + O(\varphi_3) \\
&\quad \quad + 3 \left( \overline{\eta}_b^{\overline{\eta}} - \overline{\sigma}_b^{\overline{\sigma}} \right) 
\left( \left( \overline{\rho}_b^{\overline{\rho}} \right)^2 
- \overline{\eta}_b^{\overline{\eta}} \overline{\sigma}_b^{\overline{\sigma}}
+ \overline{\eta}_b^{\overline{\eta}} \overline{\rho}_b^{\overline{\rho}}
- \overline{\sigma}_b^{\overline{\sigma}} \overline{\rho}_b^{\overline{\rho}}   
\right) \\
&\quad = O\left( \overline{\eta}_a^{\overline{\eta}} - \overline{\sigma}_a^{\overline{\sigma}} \right) + O\left( \overline{\eta}_a^{\overline{\eta}} - \overline{\rho}_a^{\overline{\rho}} \right) + O(m_b \xi_t) + O(\varphi_3) + 3 \left( \overline{\eta}_b^{\overline{\eta}} - \overline{\sigma}_b^{\overline{\sigma}} \right) \left( \overline{\rho}_b^{\overline{\rho}} - \overline{\sigma}_b^{\overline{\sigma}} \right) 
\left( \overline{\rho}_b^{\overline{\rho}} 
+ \overline{\eta}_b^{\overline{\eta}}  
\right)
\end{align*} }
This allows to reduce to quantities dominated by the distance to $\Gamma$. Let us now consider 
\begin{align}
\begin{aligned} 
&\int_0^t \int \int e^{i s \varphi_3} s^2 m_{\widehat{\mathcal{C}}}(\overline{\xi}) \mu_{BBBB}(\overline{\xi}, \overline{\eta}, \overline{\sigma}) \mu(\overline{\xi}, \overline{\eta}, \overline{\sigma}) \left( \left( \overline{\xi}_b^{\overline{\xi}} \right)^3 - \left( \overline{\eta}_b^{\overline{\eta}} \right)^3 + \left( \overline{\sigma}_b^{\overline{\sigma}} \right)^3 - \left( \overline{\rho}_b^{\overline{\rho}} \right)^3 \right) \\
&\quad \quad \quad m_{\widehat{\mathcal{C}}}(\overline{\eta}) m_{\widehat{\mathcal{C}}}(\overline{\sigma}) m_{\widehat{\mathcal{C}}}(\overline{\rho}) |\overline{\eta}| \widehat{f}(s, \overline{\eta}) \widehat{f}(s, \overline{\sigma}) \widehat{f}(s, \overline{\rho}) ~ d\overline{\eta} d\overline{\sigma} ds 
\end{aligned} \label{termedegenereCCCcubiquehb2-11} 
\end{align}

\paragraph{1.1.1.} First, we localise to have 
\begin{align*}
d\left( (\overline{\xi}, \overline{\eta}, \overline{\sigma}), \Gamma \right) &\simeq \left| \overline{\eta}_a^{\overline{\eta}} - \overline{\sigma}_a^{\overline{\sigma}} \right| + \left| \overline{\eta}_a^{\overline{\eta}} - \overline{\rho}_a^{\overline{\rho}} \right|
\end{align*}
Then, we can write on the support of $\mu_{sing}$ that
\begin{align*} 
\left( \overline{\eta}_b^{\overline{\eta}} - \overline{\sigma}_b^{\overline{\sigma}} \right) \left( \overline{\rho}_b^{\overline{\rho}} - \overline{\sigma}_b^{\overline{\sigma}} \right) 
\left( \overline{\rho}_b^{\overline{\rho}} 
+ \overline{\eta}_b^{\overline{\eta}}  
\right) &= O\left( \left( \overline{\eta}_b^{\overline{\eta}} - \overline{\sigma}_b^{\overline{\sigma}} \right) \left( \overline{\eta}_a^{\overline{\eta}} - \overline{\sigma}_a^{\overline{\sigma}} \right) \right) 
+ O\left( \left( \overline{\eta}_b^{\overline{\eta}} - \overline{\sigma}_b^{\overline{\sigma}} \right) \left( \overline{\eta}_a^{\overline{\eta}} - \overline{\rho}_a^{\overline{\rho}} \right) \right) 
\end{align*}
where the $O$ allows for symbols having the same type of singularity than $\mu_{sing}$. The two terms are symmetric and we only need to consider $O\left( \left( \overline{\eta}_b^{\overline{\eta}} - \overline{\sigma}_b^{\overline{\sigma}} \right) \left( \overline{\eta}_a^{\overline{\eta}} - \overline{\rho}_a^{\overline{\rho}} \right) \right)$. 
The factor $\overline{\eta}_a^{\overline{\eta}} - \overline{\rho}_a^{\overline{\rho}}$ allows to apply an integration by parts along $\widehat{Y}(\overline{\eta}, \overline{\rho}) \cdot \nabla_{\overline{\eta}}$. We get: 
\begin{subequations}
\begin{align}
&\eqref{termedegenereCCCcubiquehb2-11} \notag \\
&\begin{aligned}
= 
\int_0^t \int \int e^{i s \varphi_3} s m_{\widehat{\mathcal{C}}}(\overline{\xi}) \mu_{BBBB}(\overline{\xi}, \overline{\eta}, \overline{\sigma}) \mu_{sing}(\overline{\xi}, \overline{\eta}, \overline{\sigma}) \left( \overline{\eta}_b^{\overline{\eta}} - \overline{\sigma}_b^{\overline{\sigma}} \right) m_{\widehat{\mathcal{C}}}(\overline{\eta}) m_{\widehat{\mathcal{C}}}(\overline{\sigma}) m_{\widehat{\mathcal{C}}}(\overline{\rho}) \\
|\overline{\eta}| \widehat{Y}(\overline{\eta}, \overline{\rho}) \cdot \nabla_{\overline{\eta}} \left( \widehat{f}(s, \overline{\eta}) \widehat{f}(s, \overline{\sigma}) \widehat{f}(s, \overline{\rho}) \right) ~ d\overline{\eta} d\overline{\sigma} ds 
\end{aligned} \label{termedegenereCCCcubiquehb2-111-1} \\
&
+ \int_0^t \int \int \nabla_{\overline{\eta}} \cdot \left( \widehat{Y}(\overline{\eta}, \overline{\rho}) m_{\widehat{\mathcal{C}}}(\overline{\xi}) \mu_{BBBB}(\overline{\xi}, \overline{\eta}, \overline{\sigma}) \mu_{sing}(\overline{\xi}, \overline{\eta}, \overline{\sigma}) \left( \overline{\eta}_b^{\overline{\eta}} - \overline{\sigma}_b^{\overline{\sigma}} \right) m_{\widehat{\mathcal{C}}}(\overline{\eta}) m_{\widehat{\mathcal{C}}}(\overline{\sigma}) m_{\widehat{\mathcal{C}}}(\overline{\rho}) |\overline{\eta}| \right) \notag \\
&\quad \quad \quad \quad \quad \quad \quad \quad \quad \quad \quad \quad \quad \quad \quad \quad \quad \quad e^{i s \varphi_3} s \widehat{f}(s, \overline{\eta}) \widehat{f}(s, \overline{\sigma}) \widehat{f}(s, \overline{\rho}) ~ d\overline{\eta} d\overline{\sigma} ds \label{termedegenereCCCcubiquehb2-111-2}
\end{align}
\end{subequations}
On the one hand, we can estimate: 
\begin{align*}
\Vert \partial_t \eqref{termedegenereCCCcubiquehb2-111-1} \Vert_{L^2} &\lesssim t \left( \sum_{\alpha = a, b, c} \Vert m_{\alpha}(D) X_{\alpha} f(t) \Vert_{L^2} \right) \Vert \partial_x u(t) \Vert_{L^{\infty}} \Vert \partial_x u(t) \Vert_{L^{\infty}} \\
&\lesssim t^{-\frac{2}{3}} \langle t \rangle^{-\frac{1}{2}+200\delta} \Vert u \Vert_X^3
\end{align*}
On the other hand, on \eqref{termedegenereCCCcubiquehb2-111-2}, either the derivative $\nabla_{\overline{\eta}}$ does not touch $\mu_{sing}$, or it hits $\mu_{sing}$ but the factor $\left( \overline{\eta}_b^{\overline{\eta}} - \overline{\sigma}_b^{\overline{\sigma}} \right)$ allows to absorb the singularity, and we can therefore estimate in any case: 
\begin{align*}
\Vert \partial_t \eqref{termedegenereCCCcubiquehb2-111-2} \Vert_{L^2} &\lesssim t \Vert u(t) \Vert_{L^4}^2 \Vert \partial_x u(t) \Vert_{L^{\infty}} \\
&\lesssim t^{-\frac{2}{3}} \langle t \rangle^{-\frac{1}{2}+200\delta} \Vert u \Vert_X^3
\end{align*}
which is enough. 

\paragraph{1.1.2.} We now localise on 
\begin{align*}
d\left( (\overline{\xi}, \overline{\eta}, \overline{\sigma}), \Gamma \right) &\simeq \overline{\eta}_a^{\overline{\eta}} \left| \xi_t^{\overline{\eta} \overline{\sigma}} \right| + \overline{\eta}_a^{\overline{\eta}} \left| \xi_t^{\overline{\eta} \overline{\rho}} \right| \gg \left| \overline{\eta}_a^{\overline{\eta}} - \overline{\sigma}_a^{\overline{\sigma}} \right| + \left| \overline{\eta}_a^{\overline{\eta}} - \overline{\rho}_a^{\overline{\rho}} \right|
\end{align*}
Then we can decompose
\begin{align*}
\left( \overline{\eta}_b^{\overline{\eta}} - \overline{\sigma}_b^{\overline{\sigma}} \right) \left( \overline{\rho}_b^{\overline{\rho}} - \overline{\sigma}_b^{\overline{\sigma}} \right) 
\left( \overline{\rho}_b^{\overline{\rho}} 
+ \overline{\eta}_b^{\overline{\eta}}  
\right) &= O\left( \left( \overline{\eta}_b^{\overline{\eta}} - \overline{\sigma}_b^{\overline{\sigma}} \right) m_b \xi_t \right) 
\end{align*}
and again we can apply an integration by parts, this time in a direction $\widehat{X}_c$, absorbing the potential singularity if the derivative hits $\mu_{sing}$ thanks to the presence of an additional factor $\left( \overline{\eta}_b^{\overline{\eta}} - \overline{\sigma}_b^{\overline{\sigma}} \right)$. The estimates are very similar to the above case 1.1.1. and we skip the details. 

\paragraph{1.1.3.} Finally, we localise on 
\begin{align*}
d\left( (\overline{\xi}, \overline{\eta}, \overline{\sigma}), \Gamma \right) &\simeq \left| \overline{\eta}_b^{\overline{\eta}} - \overline{\rho}_b^{\overline{\eta}} \right| + \left| \overline{\eta}_b^{\overline{\eta}} - \overline{\sigma}_b^{\overline{\sigma}} \right| \gg \left| \overline{\eta}_a^{\overline{\eta}} - \overline{\sigma}_a^{\overline{\sigma}} \right| + \left| \overline{\eta}_a^{\overline{\eta}} - \overline{\rho}_a^{\overline{\rho}} \right| + \overline{\eta}_a^{\overline{\eta}} \left| \xi_t^{\overline{\eta} \overline{\sigma}} \right| + \overline{\eta}_a^{\overline{\eta}} \left| \xi_t^{\overline{\eta} \overline{\rho}} \right|
\end{align*}

We can decompose the frequency space using 
\begin{align*}
1 &= \sum_{j \in \mathbb{Z}} \sum_{k \in \mathbb{Z}, k \leq -10} \psi_{j, k}^{\widehat{\mathcal{C}}}(\overline{\xi}) 
\end{align*}
(almost everywhere), and likewise for $\overline{\eta}, \overline{\sigma}, \overline{\rho}$. The presence of $\mu_{BBBB}$ allows to reduce to only one sum in $j$. Then, on the support of $\mu_{sing}$, we have 
\begin{align*}
\left| \overline{\eta}_a^{\overline{\eta}} - \overline{\sigma}_a^{\overline{\sigma}} \right| + \left| \overline{\eta}_a^{\overline{\eta}} - \overline{\rho}_a^{\overline{\rho}} \right| + \overline{\eta}_a^{\overline{\eta}} \left| \xi_t^{\overline{\eta} \overline{\sigma}} \right| + \overline{\eta}_a^{\overline{\eta}} \left| \xi_t^{\overline{\eta} \overline{\rho}} \right|
&\ll \left| \overline{\eta}_b^{\overline{\eta}} - \overline{\rho}_b^{\overline{\eta}} \right| + \left| \overline{\eta}_b^{\overline{\eta}} - \overline{\sigma}_b^{\overline{\sigma}} \right| \lesssim |\overline{\eta}_b^{\overline{\eta}}| + |\overline{\rho}_b^{\overline{\rho}}| + |\overline{\sigma}_b^{\overline{\sigma}}| 
\end{align*}
In particular, since we saw that 
\begin{align*}
\overline{\xi}_b^{\overline{\xi}} &= \overline{\eta}_b^{\overline{\eta}} - \overline{\sigma}_b^{\overline{\sigma}} + \overline{\rho}_b^{\overline{\rho}} + O\left( \overline{\eta}_a^{\overline{\eta}} - \overline{\sigma}_a^{\overline{\sigma}} \right) + O\left( \overline{\eta}_a^{\overline{\eta}} - \overline{\rho}_a^{\overline{\rho}} \right) + O\left( \xi_t^{\overline{\eta} \overline{\sigma}} \right) + O\left( \xi_t^{\overline{\eta} \overline{\rho}} \right) 
\end{align*}
we deduce that the sums in $k$ can be simplified following a Bony-type decomposition in the pseudo-coordinate $b$ (which is not linear, but here behaves linearly up to much smaller terms). We can then reduce to only two cases, the other being symmetric: either $\overline{\eta}_b^{\overline{\eta}} \simeq \overline{\sigma}_b^{\overline{\sigma}} \gtrsim \overline{\rho}_b^{\overline{\rho}}, \overline{\xi}_b^{\overline{\xi}}$; or $\overline{\eta}_b^{\overline{\eta}} \simeq \overline{\xi}_b^{\overline{\xi}} \gg |\overline{\sigma}_b^{\overline{\sigma}}| + |\overline{\rho}_b^{\overline{\rho}}|$. We will use 
\begin{align}
\chi_{j, k}^{\widehat{\mathcal{C}}}(\overline{\xi}) := \psi_j(\overline{\xi}) \chi\left( 2^{-k} \frac{\overline{\xi}_b^{\overline{\xi}}}{\overline{\xi}_a^{\overline{\xi}}} \right) \label{localisationCtaillejtotaleexactetaillekfinemaj} 
\end{align}
that localises the direction $b$ smaller than $2^{j+k}$, and the frequency of size $2^j$. 

\paragraph{1.1.3.1.} Let us first consider 
\begin{align}
\begin{aligned}
&\sum_{j \in \mathbb{Z}} \sum_{\substack{k \in \mathbb{Z}, \\ k \leq -10}} &\int_0^t \int \int e^{i s \varphi_3} s^2 \chi_{j, k}^{\widehat{\mathcal{C}}}(\overline{\xi}) \mu_{sing}(\overline{\xi}, \overline{\eta}, \overline{\sigma}) \left( \overline{\eta}_b^{\overline{\eta}} - \overline{\sigma}_b^{\overline{\sigma}} \right) \left( \overline{\rho}_b^{\overline{\rho}} - \overline{\sigma}_b^{\overline{\sigma}} \right) 
\left( \overline{\rho}_b^{\overline{\rho}} 
+ \overline{\eta}_b^{\overline{\eta}}  
\right) \\
&&\quad \quad \quad \psi_{j, k}^{\widehat{\mathcal{C}}}(\overline{\eta}) \psi_{j, k}^{\widehat{\mathcal{C}}}(\overline{\sigma}) \chi_{j, k}^{\widehat{\mathcal{C}}}(\overline{\rho}) |\overline{\eta}| \widehat{f}(s, \overline{\eta}) \widehat{f}(s, \overline{\sigma}) \widehat{f}(s, \overline{\rho}) ~ d\overline{\eta} d\overline{\sigma} ds 
\end{aligned} \label{termedegenereCCCcubiquehb2-1131}
\end{align}

We can then use that 
\begin{align*}
\overline{\eta}_a^{\overline{\xi}} &= \overline{\eta}_a^{\overline{\eta}} + \frac{\xi \cdot \eta}{|\xi|} - |\eta| \\
&= \overline{\eta}_a^{\overline{\eta}} + O\left( \left( \xi_t^{\overline{\xi} \overline{\eta}} \right)^2 \right) 
\end{align*}
and likewise (up to changing the signs) by replacing $\overline{\eta}$ by another variable, or $a$ by $b$. Thus, 
\begin{align*}
\overline{\xi}_{\alpha}^{\overline{\xi}} &= \overline{\eta}_{\alpha}^{\overline{\eta}} - \overline{\sigma}_{\alpha}^{\overline{\eta}} + \overline{\rho}_{\alpha}^{\overline{\eta}} + O(\xi_t^2) 
\end{align*}
for $\alpha = a, b$, where $O(\xi_t^2)$ is any quadratic expression in the $\xi_t$. By Lemma \ref{lemcalculssimplescoordonneesconiquevarphi}, we therefore have 
\begin{align*}
6 \sqrt{3} \frac{\xi_0}{|\xi_0|} \varphi_3 &= \sum_{\alpha = a, b} \left[ \left( \overline{\eta}_{\alpha}^{\overline{\eta}} - \overline{\sigma}_{\alpha}^{\overline{\eta}} + \overline{\rho}_{\alpha}^{\overline{\eta}} \right)^3 - \left( \overline{\eta}_{\alpha}^{\overline{\eta}} \right)^3 + \left( \overline{\sigma}_{\alpha}^{\overline{\eta}} \right)^3 - \left( \overline{\rho}_{\alpha}^{\overline{\eta}} \right)^3 \right] \quad + O(\xi_t^2) \\
&= 3 \left( \overline{\eta}_a^{\overline{\eta}} - \overline{\sigma}_a^{\overline{\sigma}} \right) \left( \overline{\rho}_a^{\overline{\rho}} - \overline{\sigma}_a^{\overline{\sigma}} \right) \left( \overline{\eta}_a^{\overline{\eta}} + \overline{\rho}_a^{\overline{\rho}} \right) 
+ 3 \left( \overline{\eta}_b^{\overline{\eta}} - \overline{\sigma}_b^{\overline{\sigma}} \right) \left( \overline{\rho}_b^{\overline{\rho}} - \overline{\sigma}_b^{\overline{\sigma}} \right) \left( \overline{\eta}_b^{\overline{\eta}} + \overline{\rho}_b^{\overline{\rho}} \right) + O(\xi_t^2) 
\end{align*}
using also Lemma \ref{lemcalculvarphiSchrodcub1D}. 

In particular, we deduce that 
\begin{align*}
\left( \overline{\eta}_b^{\overline{\eta}} - \overline{\sigma}_b^{\overline{\sigma}} \right) \left( \overline{\rho}_b^{\overline{\rho}} - \overline{\sigma}_b^{\overline{\sigma}} \right) \left( \overline{\eta}_b^{\overline{\eta}} + \overline{\rho}_b^{\overline{\rho}} \right) 
&= O(\varphi_3) + O(\xi_t^2) + O\left( \left( \overline{\eta}_b^{\overline{\eta}} - \overline{\sigma}_b^{\overline{\sigma}} \right) \left( \overline{\rho}_b^{\overline{\rho}} - \overline{\sigma}_b^{\overline{\sigma}} \right) \right) 
\end{align*}
When we have the factor $O(\varphi_3)$, we apply an integration by parts in time; when we have the factor $O(\xi_t^2)$, we can use one of the $\xi_t$ factor to apply an integration by parts along $\widehat{X}_c(\overline{\rho}) \cdot \nabla_{\overline{\eta}}$ or $\widehat{X}_c(\overline{\rho}) \cdot \nabla_{\overline{\sigma}}$; finally, when we have a factor $O\left( \left( \overline{\eta}_b^{\overline{\eta}} - \overline{\sigma}_b^{\overline{\sigma}} \right) \left( \overline{\rho}_b^{\overline{\rho}} - \overline{\sigma}_b^{\overline{\sigma}} \right) \right)$, we use $\left( \overline{\rho}_b^{\overline{\rho}} - \overline{\sigma}_b^{\overline{\sigma}} \right)$ to apply an integration by parts along $\widehat{Y}(\overline{\sigma}, \overline{\rho}) \cdot \nabla_{\overline{\sigma}}$. We then get: {\footnotesize 
\begin{subequations}
\begin{align}
&\eqref{termedegenereCCCcubiquehb2-1131} 
= \sum_{j \in \mathbb{Z}} \sum_{\substack{k \in \mathbb{Z}, \\ k \leq -10}} \int_0^t \int \int e^{i s \varphi_3} s \chi_{j, k}^{\widehat{\mathcal{C}}}(\overline{\xi}) \mu_{sing}(\overline{\xi}, \overline{\eta}, \overline{\sigma}) \psi_{j, k}^{\widehat{\mathcal{C}}}(\overline{\eta}) \psi_{j, k}^{\widehat{\mathcal{C}}}(\overline{\sigma}) \chi_{j, k}^{\widehat{\mathcal{C}}}(\overline{\rho}) |\overline{\eta}| \widehat{f}(s, \overline{\eta}) \widehat{f}(s, \overline{\sigma}) \widehat{f}(s, \overline{\rho}) ~ d\overline{\eta} d\overline{\sigma} ds \label{termedegenereCCCcubiquehb2-1131-1} \\
&\quad + \sum_{j \in \mathbb{Z}} \sum_{\substack{k \in \mathbb{Z}, \\ k \leq -10}} \int_0^t \int \int e^{i s \varphi_3} s^2 \chi_{j, k}^{\widehat{\mathcal{C}}}(\overline{\xi}) \mu_{sing}(\overline{\xi}, \overline{\eta}, \overline{\sigma}) \psi_{j, k}^{\widehat{\mathcal{C}}}(\overline{\eta}) \psi_{j, k}^{\widehat{\mathcal{C}}}(\overline{\sigma}) \chi_{j, k}^{\widehat{\mathcal{C}}}(\overline{\rho}) |\overline{\eta}| \partial_s \left( \widehat{f}(s, \overline{\eta}) \widehat{f}(s, \overline{\sigma}) \widehat{f}(s, \overline{\rho}) \right) ~ d\overline{\eta} d\overline{\sigma} ds \label{termedegenereCCCcubiquehb2-1131-2} \\
&\quad + \sum_{j \in \mathbb{Z}} \sum_{\substack{k \in \mathbb{Z}, \\ k \leq -10}} \int \int e^{i t \varphi_3} t^2 \chi_{j, k}^{\widehat{\mathcal{C}}}(\overline{\xi}) \mu_{sing}(\overline{\xi}, \overline{\eta}, \overline{\sigma}) \psi_{j, k}^{\widehat{\mathcal{C}}}(\overline{\eta}) \psi_{j, k}^{\widehat{\mathcal{C}}}(\overline{\sigma}) \chi_{j, k}^{\widehat{\mathcal{C}}}(\overline{\rho}) |\overline{\eta}| \widehat{f}(t, \overline{\eta}) \widehat{f}(t, \overline{\sigma}) \widehat{f}(t, \overline{\rho}) ~ d\overline{\eta} d\overline{\sigma} \label{termedegenereCCCcubiquehb2-1131-3} \\
&\quad + \sum_{j \in \mathbb{Z}} \sum_{\substack{k \in \mathbb{Z}, \\ k \leq -10}} \int_0^t \int \int e^{i s \varphi_3} s \chi_{j, k}^{\widehat{\mathcal{C}}}(\overline{\xi}) \mu_{sing}(\overline{\xi}, \overline{\eta}, \overline{\sigma}) \psi_{j, k}^{\widehat{\mathcal{C}}}(\overline{\eta}) \psi_{j, k}^{\widehat{\mathcal{C}}}(\overline{\sigma}) \chi_{j, k}^{\widehat{\mathcal{C}}}(\overline{\rho}) |\overline{\eta}|^2 \widehat{f}(s, \overline{\eta}) \widehat{f}(s, \overline{\sigma}) \widehat{X}_c(\overline{\rho}) \cdot \nabla_{\overline{\xi}} \widehat{f}(s, \overline{\rho}) ~ d\overline{\eta} d\overline{\sigma} ds \label{termedegenereCCCcubiquehb2-1131-4} \\
&\quad + \sum_{j \in \mathbb{Z}} \sum_{\substack{k \in \mathbb{Z}, \\ k \leq -10}} \int_0^t \int \int e^{i s \varphi_3} s \chi_{j, k}^{\widehat{\mathcal{C}}}(\overline{\xi}) \mu_{sing}(\overline{\xi}, \overline{\eta}, \overline{\sigma}) \psi_{j, k}^{\widehat{\mathcal{C}}}(\overline{\eta}) \psi_{j, k}^{\widehat{\mathcal{C}}}(\overline{\sigma}) \chi_{j, k}^{\widehat{\mathcal{C}}}(\overline{\rho}) |\overline{\eta}|^2 O(\xi_t) \nabla_{\overline{\eta}} \widehat{f}(s, \overline{\eta}) \widehat{f}(s, \overline{\sigma}) \widehat{f}(s, \overline{\rho}) ~ d\overline{\eta} d\overline{\sigma} ds \label{termedegenereCCCcubiquehb2-1131-5} \\
&\quad + \sum_{j \in \mathbb{Z}} \sum_{\substack{k \in \mathbb{Z}, \\ k \leq -10}} \int_0^t \int \int e^{i s \varphi_3} s \chi_{j, k}^{\widehat{\mathcal{C}}}(\overline{\xi}) O(\xi_t) \nabla \mu_{sing}(\overline{\xi}, \overline{\eta}, \overline{\sigma}) \psi_{j, k}^{\widehat{\mathcal{C}}}(\overline{\eta}) \psi_{j, k}^{\widehat{\mathcal{C}}}(\overline{\sigma}) \chi_{j, k}^{\widehat{\mathcal{C}}}(\overline{\rho}) |\overline{\eta}|^2 \widehat{f}(s, \overline{\eta}) \widehat{f}(s, \overline{\sigma}) \widehat{f}(s, \overline{\rho}) ~ d\overline{\eta} d\overline{\sigma} ds \label{termedegenereCCCcubiquehb2-1131-6} \\
&\quad + \sum_{j \in \mathbb{Z}} \sum_{\substack{k \in \mathbb{Z}, \\ k \leq -10}} \int_0^t \int \int e^{i s \varphi_3} s \chi_{j, k}^{\widehat{\mathcal{C}}}(\overline{\xi}) \mu_{sing}(\overline{\xi}, \overline{\eta}, \overline{\sigma}) \widehat{X}_c(\overline{\rho}) \cdot \nabla_{\overline{\eta}} \left( \psi_{j, k}^{\widehat{\mathcal{C}}}(\overline{\eta}) \right) \psi_{j, k}^{\widehat{\mathcal{C}}}(\overline{\sigma}) \chi_{j, k}^{\widehat{\mathcal{C}}}(\overline{\rho}) |\overline{\eta}|^2 \widehat{f}(s, \overline{\eta}) \widehat{f}(s, \overline{\sigma}) \widehat{f}(s, \overline{\rho}) ~ d\overline{\eta} d\overline{\sigma} ds \label{termedegenereCCCcubiquehb2-1131-7} \\
&\quad + \sum_{j \in \mathbb{Z}} \sum_{\substack{k \in \mathbb{Z}, \\ k \leq -10}} \int_0^t \int \int e^{i s \varphi_3} s \chi_{j, k}^{\widehat{\mathcal{C}}}(\overline{\xi}) \mu_{sing}(\overline{\xi}, \overline{\eta}, \overline{\sigma}) \psi_{j, k}^{\widehat{\mathcal{C}}}(\overline{\eta}) \psi_{j, k}^{\widehat{\mathcal{C}}}(\overline{\sigma}) \chi_{j, k}^{\widehat{\mathcal{C}}}(\overline{\rho}) |\overline{\eta}|^2 \widehat{f}(s, \overline{\eta}) \widehat{Y}(\overline{\sigma}, \overline{\rho}) \cdot \nabla_{\overline{\sigma}} \left( \widehat{f}(s, \overline{\sigma}) \widehat{f}(s, \overline{\rho}) \right) ~ d\overline{\eta} d\overline{\sigma} ds \label{termedegenereCCCcubiquehb2-1131-8} \\
&\quad + \sum_{j \in \mathbb{Z}} \sum_{\substack{k \in \mathbb{Z}, \\ k \leq -10}} \int_0^t \int \int e^{i s \varphi_3} s \chi_{j, k}^{\widehat{\mathcal{C}}}(\overline{\xi}) \left( \overline{\eta}_a^{\overline{\eta}} - \overline{\sigma}_a^{\overline{\sigma}} \right) \nabla \mu_{sing}(\overline{\xi}, \overline{\eta}, \overline{\sigma}) \psi_{j, k}^{\widehat{\mathcal{C}}}(\overline{\eta}) \psi_{j, k}^{\widehat{\mathcal{C}}}(\overline{\sigma}) \chi_{j, k}^{\widehat{\mathcal{C}}}(\overline{\rho}) |\overline{\eta}| \widehat{f}(s, \overline{\eta}) \widehat{f}(s, \overline{\sigma}) \widehat{f}(s, \overline{\rho}) ~ d\overline{\eta} d\overline{\sigma} ds \label{termedegenereCCCcubiquehb2-1131-9} \\
&\quad + \sum_{j \in \mathbb{Z}} \sum_{\substack{k \in \mathbb{Z}, \\ k \leq -10}} \int_0^t \int \int e^{i s \varphi_3} s \chi_{j, k}^{\widehat{\mathcal{C}}}(\overline{\xi}) \mu_{sing}(\overline{\xi}, \overline{\eta}, \overline{\sigma}) \psi_{j, k}^{\widehat{\mathcal{C}}}(\overline{\eta}) \widehat{Y}_(\overline{\sigma}, \overline{\rho}) \cdot \nabla_{\overline{\sigma}} \left( \psi_{j, k}^{\widehat{\mathcal{C}}}(\overline{\sigma}) \chi_{j, k}^{\widehat{\mathcal{C}}}(\overline{\rho}) \right) |\overline{\eta}|^2 \widehat{f}(s, \overline{\eta}) \widehat{f}(s, \overline{\sigma}) \widehat{f}(s, \overline{\rho}) ~ d\overline{\eta} d\overline{\sigma} ds \label{termedegenereCCCcubiquehb2-1131-10} 
\end{align}
\end{subequations} }
where the symbols can change from line to line as long as they keep similar properties, and where we symmetrized between $\overline{\eta}$ and $\overline{\sigma}$. We also used that $\widehat{X}_c(\overline{\rho}) \cdot \nabla \chi_{j, k}^{\widehat{\mathcal{C}}}(\overline{\rho}) = 0$. 

First we can estimate
\begin{align*}
\Vert \partial_t \eqref{termedegenereCCCcubiquehb2-1131-1} \Vert_{L^2} &\lesssim \sum_{j \in \mathbb{Z}} \sum_{\substack{k \in \mathbb{Z}, \\ k \leq -10}} t \Vert \nabla \psi_{j, k}^{\widehat{\mathcal{C}}}(D) u(t) \Vert_{L^{\infty}}^2 \Vert |\nabla|^{-1} u(t) \Vert_{L^2} \\
&\lesssim \sum_{j \in \mathbb{Z}} \sum_{\substack{k \in \mathbb{Z}, \\ k \leq -10}} t^{-\frac{2}{3}} \langle t \rangle^{-\frac{1}{2}+200\delta} 2^{\delta j+\delta k} \langle 2^j \rangle^{-\frac{1}{2}} \Vert u \Vert_X^3 \\
&\lesssim t^{-\frac{2}{3}} \langle t \rangle^{-\frac{1}{2}+200\delta} \Vert u \Vert_X^3 \\
\Vert \partial_t \eqref{termedegenereCCCcubiquehb2-1131-2} \Vert_{L^2} &\lesssim \sum_{j \in \mathbb{Z}} \sum_{\substack{k \in \mathbb{Z}, \\ k \leq -10}} t^2 \Vert \partial_t f(t) \Vert_{L^4} \Vert \nabla \psi_{j, k}^{\widehat{\mathcal{C}}}(D) u(t) \Vert_{L^{\infty}} \Vert u(t) \Vert_{L^4} \\
&\lesssim \sum_{j \in \mathbb{Z}} \sum_{\substack{k \in \mathbb{Z}, \\ k \leq -10}} t^{-\frac{1}{2}} \langle t \rangle^{-\frac{3}{4}+300\delta} 2^{\delta j+\delta k} \langle 2^j \rangle^{-\frac{1}{2}} \Vert u \Vert_X^4 \\
&\lesssim t^{-\frac{1}{2}} \langle t \rangle^{-\frac{3}{4}+300\delta} \Vert u \Vert_X^4 \\
\Vert \partial_t (\eqref{termedegenereCCCcubiquehb2-1131-4}+\eqref{termedegenereCCCcubiquehb2-1131-8}) \Vert_{L^2} &\lesssim \sum_{j \in \mathbb{Z}} \sum_{\substack{k \in \mathbb{Z}, \\ k \leq -10}} t \Vert \nabla \psi_{j, k}^{\widehat{\mathcal{C}}}(D) u(t) \Vert_{L^{\infty}} \Vert \partial_x u(t) \Vert_{L^{\infty}} \sum_{\alpha = a, b, c} \Vert m_{\alpha}(D) X_{\alpha} f(t) \Vert_{L^2} \\
&\lesssim \sum_{j \in \mathbb{Z}} \sum_{\substack{k \in \mathbb{Z}, \\ k \leq -10}} t^{-\frac{2}{3}} \langle t \rangle^{-\frac{1}{2}+200\delta} 2^{\delta j+\delta k} \langle 2^j \rangle^{-\frac{1}{2}} \Vert u \Vert_X^3 \\
&\lesssim t^{-\frac{2}{3}} \langle t \rangle^{-\frac{1}{2}+200\delta} \Vert u \Vert_X^3 
\end{align*}
Then, for \eqref{termedegenereCCCcubiquehb2-1131-5}, by the localisation $\mu_{sing}$ and the $(j, k)$-localisation, we have
\begin{align*}
O(\xi_t) &= O\left( \overline{\eta}_b^{\overline{\eta}} - \overline{\sigma}_b^{\overline{\sigma}} \right) + O\left( \overline{\eta}_b^{\overline{\eta}} - \overline{\rho}_b^{\overline{\rho}} \right) \\
&= O\left( \overline{\eta}_b^{\overline{\eta}} \right) 
\end{align*}
which allows to estimate \eqref{termedegenereCCCcubiquehb2-1131-5} like \eqref{termedegenereCCCcubiquehb2-1131-8}. The same way, $O(\xi_t)$ can absorb the singularity of $\nabla \mu_{sing}$ and we can estimate \eqref{termedegenereCCCcubiquehb2-1131-6} like \eqref{termedegenereCCCcubiquehb2-1131-1}, then \eqref{termedegenereCCCcubiquehb2-1131-9} as well. Finally, we compute that
\begin{align*}
\widehat{X}_c(\overline{\rho}) \cdot \nabla_{\overline{\eta}} \psi_{j, k}^{\widehat{\mathcal{C}}}(\overline{\eta}) 
&= O(1) + O\left( \frac{J \rho \cdot \eta}{|\rho|} 2^{-k} \right) \\
&= O(1) 
\end{align*}
by the same argument, so we can estimate \eqref{termedegenereCCCcubiquehb2-1131-7} like \eqref{termedegenereCCCcubiquehb2-1131-1}. Finally, by Lemma \ref{lemchampmodifieYpropfond}, we have
\begin{align*}
&\widehat{Y}(\overline{\sigma}, \overline{\rho}) \cdot \nabla_{\overline{\sigma}} \left[ \psi_{j, k}^{\widehat{\mathcal{C}}}(\overline{\sigma}) \chi_{j, k}^{\widehat{\mathcal{C}}}(\overline{\rho}) \right] \\
&\quad = \sum_{\alpha = a, b, c} \left[ O\left( m_{\alpha}(\overline{\sigma}) \widehat{X}_{\alpha}(\overline{\sigma}) \cdot \nabla_{\overline{\sigma}} \psi_{j, k}^{\widehat{\mathcal{C}}}(\overline{\sigma}) \right) 
+ O\left( m_{\alpha}(\overline{\rho}) \widehat{X}_{\alpha}(\overline{\rho}) \cdot \nabla_{\overline{\sigma}} \chi_{j, k}^{\widehat{\mathcal{C}}}(\overline{\rho}) \right) \right] \\
&\quad = O(1) 
\end{align*}
because $\widehat{X}_a(\overline{\sigma}) \cdot \nabla \overline{\sigma}_b^{\overline{\sigma}} = \widehat{X}_c(\overline{\sigma}) \cdot \nabla \overline{\sigma}_b^{\overline{\sigma}} = 0$. We can thus also estimate \eqref{termedegenereCCCcubiquehb2-1131-10} like \eqref{termedegenereCCCcubiquehb2-1131-1}. 

It only remains \eqref{termedegenereCCCcubiquehb2-1131-3}, but it is a boundary term that can be placed in $g_{bb, 3}(t)$, since we can estimate it as: 
\begin{align*}
\Vert \eqref{termedegenereCCCcubiquehb2-1131-3} \Vert_{L^2} &\lesssim \sum_{j \in \mathbb{Z}} \sum_{\substack{k \in \mathbb{Z}, \\ k \leq -10}} t^2 \Vert \psi_{j, k}^{\widehat{\mathcal{C}}}(D) \nabla u(t) \Vert_{L^{\infty}} \Vert u(t) \Vert_{L^4}^2 \\
&\lesssim \sum_{j \in \mathbb{Z}} \sum_{\substack{k \in \mathbb{Z}, \\ k \leq -10}} t^{\frac{1}{3}} \langle t \rangle^{-\frac{1}{2}+200\delta} 2^{\delta j+\delta k} \langle 2^j \rangle^{-\frac{1}{2}} \Vert u \Vert_X^3 \\
&\lesssim \Vert u \Vert_X^3 \\
\Vert e^{-i t \omega(D)} \nabla \mathcal{F}^{-1} \eqref{termedegenereCCCcubiquehb2-1131-3} \Vert_{L^4} &\lesssim \sum_{j \in \mathbb{Z}} \sum_{\substack{k \in \mathbb{Z}, \\ k \leq -10}} t^2 \Vert \psi_{j, k}^{\widehat{\mathcal{C}}}(D) \nabla u(t) \Vert_{L^{\infty}} \Vert \partial_x u(t) \Vert_{L^{\infty}} \Vert u(t) \Vert_{L^4} \\
&\lesssim \sum_{j \in \mathbb{Z}} \sum_{\substack{k \in \mathbb{Z}, \\ k \leq -10}} t^{-\frac{1}{12}} \langle t \rangle^{-\frac{5}{8}+250\delta} 2^{\delta j+\delta k} \langle 2^j \rangle^{-\frac{1}{2}} \Vert u \Vert_X^3 \\
&\lesssim t^{-\frac{1}{12}} \langle t \rangle^{-\frac{5}{8}+250\delta} \Vert u \Vert_X^3 
\end{align*}

\paragraph{1.1.3.2} Let us now consider: 
\begin{align}
\begin{aligned}
&\sum_{j \in \mathbb{Z}} \sum_{\substack{k \in \mathbb{Z}, \\ k \leq -10}} &\int_0^t \int \int e^{i s \varphi_3} s^2 \psi_{j, k}^{\widehat{\mathcal{C}}}(\overline{\xi}) \mu_{sing}(\overline{\xi}, \overline{\eta}, \overline{\sigma}) \left( \overline{\eta}_b^{\overline{\eta}} - \overline{\sigma}_b^{\overline{\sigma}} \right) \left( \overline{\rho}_b^{\overline{\rho}} - \overline{\sigma}_b^{\overline{\sigma}} \right) 
\left( \overline{\rho}_b^{\overline{\rho}} 
+ \overline{\eta}_b^{\overline{\eta}}  
\right) \\
&&\quad \quad \quad \psi_{j, k}^{\widehat{\mathcal{C}}}(\overline{\eta}) \chi_{j, k-10}^{\widehat{\mathcal{C}}}(\overline{\sigma}) \chi_{j, k-10}^{\widehat{\mathcal{C}}}(\overline{\rho}) |\overline{\eta}| \widehat{f}(s, \overline{\eta}) \widehat{f}(s, \overline{\sigma}) \widehat{f}(s, \overline{\rho}) ~ d\overline{\eta} d\overline{\sigma} ds 
\end{aligned} \label{termedegenereCCCcubiquehb2-1132}
\end{align}

In particular, we have that 
\begin{align*}
d\left( (\overline{\xi}, \overline{\eta}, \overline{\sigma}), \Gamma \right) &\simeq \overline{\eta}_b^{\overline{\eta}} \gg \left( \overline{\eta}_a^{\overline{\eta}} - \overline{\rho}_a^{\overline{\rho}}, \overline{\rho}_b^{\overline{\rho}} \right) 
\end{align*}
By Lemma \ref{lemcalculssimplescoordonneesconiquevarphi}, 
\begin{align*}
2 \partial_{\eta_0} \varphi_3 &= \left( \overline{\rho}_a^{\overline{\rho}} \right)^2 + \left( \overline{\rho}_b^{\overline{\rho}} \right)^2 - \left( \overline{\eta}_a^{\overline{\eta}} \right)^2 - \left( \overline{\eta}_b^{\overline{\eta}} \right)^2 \\
&\simeq - \left( \overline{\eta}_b^{\overline{\eta}} \right)^2 
\end{align*}
and likewise for $\left( \partial_{\eta_0} - \partial_{\sigma_0} \right) \varphi_3$. We can then decompose: 
\begin{align*}
\left( \overline{\eta}_b^{\overline{\eta}} - \overline{\sigma}_b^{\overline{\sigma}} \right) \left( \overline{\rho}_b^{\overline{\rho}} - \overline{\sigma}_b^{\overline{\sigma}} \right) 
\left( \overline{\rho}_b^{\overline{\rho}} 
+ \overline{\eta}_b^{\overline{\eta}}  
\right) &= O\left( \left( \overline{\eta}_b^{\overline{\eta}} \right)^2 \overline{\rho}_b^{\overline{\rho}} \right) + O\left( \left( \overline{\eta}_b^{\overline{\eta}} \right)^2 \overline{\sigma}_b^{\overline{\sigma}} \right) 
\end{align*}
The two terms are symmetric, so we only need to consider the first one. We can use the factor $\left( \overline{\eta}_b^{\overline{\eta}} \right)^2$ to apply an integration by parts in $\eta_0$, and we get: {\scriptsize 
\begin{subequations}
\begin{align}
\eqref{termedegenereCCCcubiquehb2-1132} &= \sum_{j \in \mathbb{Z}} \sum_{\substack{k \in \mathbb{Z}, \\ k \leq -10}} \int_0^t \int \int e^{i s \varphi_3} s \psi_{j, k}^{\widehat{\mathcal{C}}}(\overline{\xi}) \mu_{sing}(\overline{\xi}, \overline{\eta}, \overline{\sigma}) \psi_{j, k}^{\widehat{\mathcal{C}}}(\overline{\eta}) \chi_{j, k-10}^{\widehat{\mathcal{C}}}(\overline{\sigma}) \chi_{j, k-10}^{\widehat{\mathcal{C}}}(\overline{\rho}) |\overline{\eta}| \widehat{f}(s, \overline{\eta}) \widehat{f}(s, \overline{\sigma}) \widehat{f}(s, \overline{\rho}) ~ d\overline{\eta} d\overline{\sigma} ds \label{termedegenereCCCcubiquehb2-1132-1} \\
&\quad + \sum_{j \in \mathbb{Z}} \sum_{\substack{k \in \mathbb{Z}, \\ k \leq -10}} \int_0^t \int \int e^{i s \varphi_3} s \psi_{j, k}^{\widehat{\mathcal{C}}}(\overline{\xi}) \mu_{sing}(\overline{\xi}, \overline{\eta}, \overline{\sigma}) \psi_{j, k}^{\widehat{\mathcal{C}}}(\overline{\eta}) \chi_{j, k-10}^{\widehat{\mathcal{C}}}(\overline{\sigma}) \chi_{j, k-10}^{\widehat{\mathcal{C}}}(\overline{\rho}) |\overline{\eta}| \overline{\rho}_b^{\overline{\rho}} \partial_{\eta_0} \widehat{f}(s, \overline{\eta}) \widehat{f}(s, \overline{\sigma}) \widehat{f}(s, \overline{\rho}) ~ d\overline{\eta} d\overline{\sigma} ds \label{termedegenereCCCcubiquehb2-1132-2} \\
&\quad + \sum_{j \in \mathbb{Z}} \sum_{\substack{k \in \mathbb{Z}, \\ k \leq -10}} \int_0^t \int \int e^{i s \varphi_3} s \psi_{j, k}^{\widehat{\mathcal{C}}}(\overline{\xi}) \mu_{sing}(\overline{\xi}, \overline{\eta}, \overline{\sigma}) \psi_{j, k}^{\widehat{\mathcal{C}}}(\overline{\eta}) \chi_{j, k-10}^{\widehat{\mathcal{C}}}(\overline{\sigma}) \chi_{j, k-10}^{\widehat{\mathcal{C}}}(\overline{\rho}) |\overline{\eta}| \widehat{f}(s, \overline{\eta}) \widehat{f}(s, \overline{\sigma}) \overline{\rho}_b^{\overline{\rho}} \partial_{\eta_0} \widehat{f}(s, \overline{\rho}) ~ d\overline{\eta} d\overline{\sigma} ds \label{termedegenereCCCcubiquehb2-1132-3} \\
&\quad + \sum_{j \in \mathbb{Z}} \sum_{\substack{k \in \mathbb{Z}, \\ k \leq -10}} \int_0^t \int \int e^{i s \varphi_3} s \psi_{j, k}^{\widehat{\mathcal{C}}}(\overline{\xi}) \overline{\rho}_b^{\overline{\rho}} \nabla \mu_{sing}(\overline{\xi}, \overline{\eta}, \overline{\sigma}) \psi_{j, k}^{\widehat{\mathcal{C}}}(\overline{\eta}) \chi_{j, k-10}^{\widehat{\mathcal{C}}}(\overline{\sigma}) \chi_{j, k-10}^{\widehat{\mathcal{C}}}(\overline{\rho}) |\overline{\eta}| \widehat{f}(s, \overline{\eta}) \widehat{f}(s, \overline{\sigma}) \widehat{f}(s, \overline{\rho}) ~ d\overline{\eta} d\overline{\sigma} ds \label{termedegenereCCCcubiquehb2-1132-4} \\
&\quad + \sum_{j \in \mathbb{Z}} \sum_{\substack{k \in \mathbb{Z}, \\ k \leq -10}} \int_0^t \int \int e^{i s \varphi_3} s \psi_{j, k}^{\widehat{\mathcal{C}}}(\overline{\xi}) \mu_{sing}(\overline{\xi}, \overline{\eta}, \overline{\sigma}) \overline{\rho}_b^{\overline{\rho}} \partial_{\eta_0} \left( \psi_{j, k}^{\widehat{\mathcal{C}}}(\overline{\eta}) \chi_{j, k-10}^{\widehat{\mathcal{C}}}(\overline{\sigma}) \chi_{j, k-10}^{\widehat{\mathcal{C}}}(\overline{\rho}) \right) |\overline{\eta}| \widehat{f}(s, \overline{\eta}) \widehat{f}(s, \overline{\sigma}) \widehat{f}(s, \overline{\rho}) ~ d\overline{\eta} d\overline{\sigma} ds \label{termedegenereCCCcubiquehb2-1132-5} 
\end{align}
\end{subequations} }
\eqref{termedegenereCCCcubiquehb2-1132-1} can be estimated like \eqref{termedegenereCCCcubiquehb2-1131-1}, \eqref{termedegenereCCCcubiquehb2-1132-3} like \eqref{termedegenereCCCcubiquehb2-1131-8}, and using that $\overline{\rho}_b^{\overline{\rho}} \lesssim \overline{\eta}_b^{\overline{\eta}} \simeq d\left( (\overline{\xi}, \overline{\sigma}, \overline{\rho}), \Gamma \right)$, \eqref{termedegenereCCCcubiquehb2-1132-4} and \eqref{termedegenereCCCcubiquehb2-1132-5} can also be estimated like \eqref{termedegenereCCCcubiquehb2-1131-1}. Finally, for \eqref{termedegenereCCCcubiquehb2-1132-2}, we use that the supports of the $\psi_{j, k}^{\widehat{\mathcal{C}}}$ form a locally finite covering of a neighborhood of $\widehat{\mathcal{C}}$ in order to write: 
\begin{align*}
\Vert \partial_t \eqref{termedegenereCCCcubiquehb2-1132-2} \Vert_{L^2} &\lesssim \sum_{\alpha = a, b, c} \left( \sum_{j \in \mathbb{Z}} \sum_{k \in \mathbb{Z}, k \leq -10} \left[ t \Vert \psi_{j, k}^{\widehat{\mathcal{C}}}(D) m_{\alpha}(D) X_{\alpha} f(t) \Vert_{L^2} \Vert \partial_x u(t) \Vert_{L^{\infty}}^2 \right]^2 \right)^{\frac{1}{2}} \\
&\lesssim \sum_{\alpha = a, b, c} t^{-\frac{2}{3}} \langle t \rangle^{-\frac{1}{2}+200\delta} \Vert u \Vert_X^2 \left( \sum_{j \in \mathbb{Z}} \sum_{k \in \mathbb{Z}, k \leq -10} \Vert \psi_{j, k}^{\widehat{\mathcal{C}}}(D) m_{\alpha}(D) X_{\alpha} f(t) \Vert_{L^2}^2 \right)^{\frac{1}{2}} \\
&\lesssim \sum_{\alpha = a, b, c} t^{-\frac{2}{3}} \langle t \rangle^{-\frac{1}{2}+200\delta} \Vert u \Vert_X^2 \Vert m_{\alpha}(D) X_{\alpha} f(t) \Vert_{L^2} \\
&\lesssim t^{-\frac{2}{3}} \langle t \rangle^{-\frac{1}{2}+200\delta} \Vert u \Vert_X^3
\end{align*}
as wanted. 

\paragraph{1.2.} Let us now consider the case where $|\overline{\xi}| \ll |\overline{\eta}| + |\overline{\sigma}| + |\overline{\rho}|$. 

As before, we can apply integrations by parts as soon as we have a factor $\overline{\eta}_a^{\overline{\eta}} - \overline{\sigma}_a^{\overline{\sigma}}$ or $\overline{\eta}_a^{\overline{\eta}} - \overline{\rho}_a^{\overline{\rho}}$: but here all these quantities cannot be $o(1)$ simultaneously. In particular, there is no particular difficulty to rewrite \eqref{equdecchampbCCC-1} as terms of \eqref{qteslemestimeesgeneriquesdeccubiquehb2}. It thus suffices to rewrite \eqref{equdecchampbCCC-2} to \eqref{equdecchampbCCC-1}. 

To that end, we write that 
\begin{align*}
\xi_0 m_b(\overline{\xi}) &= \frac{\xi_0 |\xi|}{|\overline{\xi}|^2} \frac{3 \xi_0^2 - |\xi|^2}{|\overline{\xi}|} \\
&= \frac{\xi_0 |\xi|}{|\overline{\xi}|^2} \left( \frac{\overline{\xi}_b^{\overline{\rho}} \overline{\xi}_a^{\overline{\rho}}}{|\overline{\xi}|} + O\left( \xi_t^{\overline{\xi} \overline{\rho}} \right) \right) \\
&= O\left( \overline{\rho}_b^{\overline{\rho}} \right) + O\left( \overline{\sigma}_b^{\overline{\sigma}} \right) + O\left( \overline{\eta}_b^{\overline{\eta}} \right) + O\left( \xi_t^{\overline{\eta} \overline{\rho}} \right) + O\left( \xi_t^{\overline{\sigma} \overline{\rho}} \right) 
\end{align*}
The term with factor $O\left( \overline{\rho}_b^{\overline{\rho}} \right)$ already has the desired form; on the one in $O\left( \overline{\eta}_b^{\overline{\eta}} \right)$, we can apply an integration by parts $\nabla_{\overline{\eta}}$, and symmetrically on the one in $O\left( \overline{\sigma}_b^{\overline{\sigma}} \right)$; finally, on the one in $O\left( \xi_t^{\overline{\eta} \overline{\rho}} \right)$, we can rewrite it like $\widehat{X}_{b-\widehat{\mathcal{C}}}$ without singularity and then apply an integration by parts, and symmetrically for $O\left( \xi_t^{\overline{\sigma} \overline{\rho}} \right)$. (All these transformations are similar to the ones already applied in case 1.1.) 

\paragraph{2.} Let us now consider the case where $\eta, \sigma, \rho$ are not all close to alignment. 

In particular, up to exchanging them, we can assume that $\rho$ is aligned neither with $\eta$ nor with $\sigma$. In particular, $\widehat{X}_{b-\widehat{\mathcal{C}}}(\overline{\eta}, \overline{\rho})$ is not singular here and we can decompose \eqref{equdecchampbCCC-2} by projecting $\widehat{X}_b(\overline{\xi})$ on the basis $\left( \widehat{X}_{\alpha}(\overline{\rho}) \right)_{\alpha = a, b, c}$ to get terms from \eqref{qteslemestimeesgeneriquesdeccubiquehb2} plus one on which we apply an integration by parts: {\scriptsize 
\begin{subequations}
\begin{align}
&\int_0^t \int \int e^{i s \varphi_3} s \xi_0 \mu(\overline{\xi}, \overline{\eta}, \overline{\sigma}) m_{\widehat{\mathcal{C}}}(\overline{\eta}) m_{\widehat{\mathcal{C}}}(\overline{\sigma}) m_{\widehat{\mathcal{C}}}(\overline{\rho}) \widehat{F}_1(s, \overline{\eta}) \widehat{F}_2(s, \overline{\sigma}) m_b(\overline{\xi}) \frac{\rho_0}{|\rho_0|} P_b^b(\overline{\xi}, \overline{\rho}) \widehat{X}_{b-\widehat{\mathcal{C}}}(\overline{\eta}, \overline{\rho}) \cdot \nabla_{\overline{\xi}} \widehat{F}_3(s, \overline{\rho}) ~ d\overline{\eta} d\overline{\sigma} ds \notag \\
&= \int_0^t \int \int i s^2 m_b(\overline{\xi}) \frac{\rho_0}{|\rho_0|} P_b^b(\overline{\xi}, \overline{\rho}) \widehat{X}_{b-\widehat{\mathcal{C}}}(\overline{\eta}, \overline{\rho}) \cdot \nabla_{\overline{\eta}} \varphi_3 e^{i s \varphi_3} \xi_0 \mu(\overline{\xi}, \overline{\eta}, \overline{\sigma}) m_{\widehat{\mathcal{C}}}(\overline{\eta}) m_{\widehat{\mathcal{C}}}(\overline{\sigma}) m_{\widehat{\mathcal{C}}}(\overline{\rho}) \widehat{F}_1(s, \overline{\eta}) \widehat{F}_2(s, \overline{\sigma}) \widehat{F}_3(s, \overline{\rho}) ~ d\overline{\eta} d\overline{\sigma} ds \label{equdecchampbCCC-cas2-2-1} \\
&+ \int_0^t \int \int e^{i s \varphi_3} s \xi_0 \mu(\overline{\xi}, \overline{\eta}, \overline{\sigma}) m_{\widehat{\mathcal{C}}}(\overline{\eta}) m_{\widehat{\mathcal{C}}}(\overline{\sigma}) m_{\widehat{\mathcal{C}}}(\overline{\rho}) m_b(\overline{\xi}) \frac{\rho_0}{|\rho_0|} P_b^b(\overline{\xi}, \overline{\rho}) \widehat{X}_{b-\widehat{\mathcal{C}}}(\overline{\eta}, \overline{\rho}) \cdot \nabla_{\overline{\eta}} \widehat{F}_1(s, \overline{\eta}) \widehat{F}_2(s, \overline{\sigma}) \widehat{F}_3(s, \overline{\rho}) ~ d\overline{\eta} d\overline{\sigma} ds \label{equdecchampbCCC-cas2-2-2} \\
&+ \int_0^t \int \int e^{i s \varphi_3} s \nabla_{\overline{\eta}} \cdot \left( m_b(\overline{\xi}) \frac{\rho_0}{|\rho_0|} P_b^b(\overline{\xi}, \overline{\rho}) \widehat{X}_{b-\widehat{\mathcal{C}}}(\overline{\eta}, \overline{\rho}) \xi_0 \mu(\overline{\xi}, \overline{\eta}, \overline{\sigma}) m_{\widehat{\mathcal{C}}}(\overline{\eta}) m_{\widehat{\mathcal{C}}}(\overline{\sigma}) m_{\widehat{\mathcal{C}}}(\overline{\rho}) \right) \widehat{F}_1(s, \overline{\eta}) \widehat{F}_2(s, \overline{\sigma}) \widehat{F}_3(s, \overline{\rho}) ~ d\overline{\eta} d\overline{\sigma} ds \label{equdecchampbCCC-cas2-2-3} 
\end{align}
\end{subequations} }
\eqref{equdecchampbCCC-cas2-2-2} and \eqref{equdecchampbCCC-cas2-2-3} are terms from \eqref{qteslemestimeesgeneriquesdeccubiquehb2}. 

Finally, we group \eqref{equdecchampbCCC-cas2-2-1} with \eqref{equdecchampbCCC-1} and we get 
\begin{align}
\begin{aligned}
&\int_0^t \int \int i s^2 m_b(\overline{\xi}) \left( \widehat{X}_b(\overline{\xi}) \cdot \nabla_{\overline{\xi}} + \frac{\rho_0}{|\rho_0|} P_b^b(\overline{\xi}, \overline{\rho}) \widehat{X}_{b-\widehat{\mathcal{C}}}(\overline{\eta}, \overline{\rho}) \cdot \nabla_{\overline{\eta}} \right) \varphi_3 e^{i s \varphi_3} \xi_0 \mu(\overline{\xi}, \overline{\eta}, \overline{\sigma}) \\
&\quad \quad m_{\widehat{\mathcal{C}}}(\overline{\eta}) m_{\widehat{\mathcal{C}}}(\overline{\sigma}) m_{\widehat{\mathcal{C}}}(\overline{\rho}) \widehat{F}_1(s, \overline{\eta}) \widehat{F}_2(s, \overline{\sigma}) \widehat{F}_3(s, \overline{\rho}) ~ d\overline{\eta} d\overline{\sigma} ds 
\end{aligned} \label{equdecchampbCCC-cas2-final} 
\end{align} 

Since $\rho$ is not aligned with $\eta$ or $\sigma$, $\xi_t^{\overline{\eta} \overline{\rho}}$ and $\xi_t^{\overline{\sigma} \overline{\rho}}$ cannot be close to $0$, and $\theta^{\overline{\eta} \overline{\rho}}$ or $\theta^{\overline{\sigma} \overline{\rho}}$ cannot be close to $\pm 1$. In particular, we deduce that 
\begin{align*}
\left( \overline{\eta}_a^{\overline{\eta}} \right)^2 - \left( \overline{\rho}_a^{\overline{\eta}} \right)^2 &= O\left( \widehat{Y}(\overline{\eta}, \overline{\rho}) \cdot \nabla_{\overline{\eta}} \varphi_3 \right) \\
\left( \overline{\sigma}_a^{\overline{\sigma}} \right)^2 - \left( \overline{\rho}_a^{\overline{\eta}} \right)^2 &= O\left( \widehat{Y}(\overline{\sigma}, \overline{\rho}) \cdot \nabla_{\overline{\sigma}} \varphi_3 \right) 
\end{align*}
We can therefore only consider the angular neighborhood of a point where all these quantities are $o(1)$, and thus $|\overline{\rho}| \simeq |\overline{\eta}| \simeq |\overline{\sigma}|$. Furthermore, since $\xi_0 = \eta_0 + \sigma_0 + \rho_0$, we deduce that $|\overline{\xi}| \simeq |\overline{\eta}|$ as well. By symmetry of $\overline{\eta}$ and $\overline{\sigma}$, we can assume that $\xi_0$ has the same sign as $\eta_0$, i.e. $\epsilon^{\overline{\xi} \overline{\eta}} = 1$. We then have 
\begin{align*}
\overline{\eta}_a^{\overline{\eta}} - \overline{\rho}_a^{\overline{\eta}} &= O\left( \widehat{Y}(\overline{\eta}, \overline{\rho}) \cdot \nabla_{\overline{\eta}} \varphi_3 \right) \\
\overline{\sigma}_a^{\overline{\sigma}} - \overline{\rho}_a^{\overline{\eta}} &= O\left( \widehat{Y}(\overline{\sigma}, \overline{\rho}) \cdot \nabla_{\overline{\sigma}} \varphi_3 \right) \\
1 &= O\left( \widehat{X}_c(\overline{\rho}) \cdot \nabla_{\overline{\eta}} \varphi_3 \right) \\
1 &= O\left( \widehat{X}_c(\overline{\eta}) \cdot \nabla_{\overline{\eta}} \varphi_3 \right) \\
1 &= O\left( \widehat{X}_c(\overline{\rho}) \cdot \nabla_{\overline{\sigma}} \varphi_3 \right) 
\end{align*}
We can apply integrations by parts whenever we have a factor $O(\overline{\eta}_b^{\overline{\eta}}), O(\overline{\sigma}_b^{\overline{\sigma}}), O(\overline{\rho}_b^{\overline{\rho}})$. We then compute 
\begin{align*}
\frac{\eta_0}{|\eta_0|} \varphi_3 &= |\xi_0|^3 + |\xi_0| |\xi|^2 - |\eta_0|^3 - |\eta_0| - \epsilon^{\overline{\sigma} \overline{\eta}} |\sigma_0|^3 - \epsilon^{\overline{\sigma} \overline{\eta}} |\sigma_0| |\sigma|^2 - \epsilon^{\overline{\rho} \overline{\eta}} |\rho_0|^3
- \epsilon^{\overline{\rho} \overline{\eta}} |\rho_0| |\rho|^2 \\
&= |\xi_0| \left( \xi_0^2 + |\xi|^2 \right) - \frac{1}{6 \sqrt{3}} \left( \left( \overline{\eta}_a^{\overline{\eta}} \right)^3 + \left( \overline{\eta}_b^{\overline{\eta}} \right)^3 + \epsilon^{\overline{\sigma} \overline{\eta}} \left( \overline{\sigma}_a^{\overline{\sigma}} \right)^3 + \epsilon^{\overline{\sigma} \overline{\eta}} \left( \overline{\sigma}_b^{\overline{\sigma}} \right)^3 + \epsilon^{\overline{\rho} \overline{\eta}} \left( \overline{\rho}_a^{\overline{\rho}} \right)^3 + \epsilon^{\overline{\rho} \overline{\eta}} \left( \overline{\rho}_b^{\overline{\rho}} \right)^3 \right) \\
&= O(\overline{\eta}_b^{\overline{\eta}}) + O(\overline{\sigma}_b^{\overline{\sigma}}) + O(\overline{\rho}_b^{\overline{\rho}}) + O(\overline{\eta}_a^{\overline{\eta}} - \overline{\rho}_a^{\overline{\eta}}) + O(\overline{\sigma}_a^{\overline{\sigma}} - \overline{\rho}_a^{\overline{\eta}}) \\
&\quad + \left( |\eta_0| + \epsilon^{\overline{\eta} \overline{\sigma}} |\sigma_0| + \epsilon^{\overline{\eta} \overline{\rho}} |\rho_0| \right) \left( \xi_0^2 + |\xi|^2 \right) - \frac{1}{6 \sqrt{3}} \left( 1 + \epsilon^{\overline{\sigma} \overline{\eta}} + \epsilon^{\overline{\rho} \overline{\eta}} \right) \left( \overline{\eta}_a^{\overline{\eta}} \right)^3 \\
&= O(\overline{\eta}_b^{\overline{\eta}}) + O(\overline{\sigma}_b^{\overline{\sigma}}) + O(\overline{\rho}_b^{\overline{\rho}}) + O(\overline{\eta}_a^{\overline{\eta}} - \overline{\rho}_a^{\overline{\eta}}) + O(\overline{\sigma}_a^{\overline{\sigma}} - \overline{\rho}_a^{\overline{\eta}}) \\
&\quad + \frac{1}{2 \sqrt{3}} \overline{\eta}_a^{\overline{\eta}} \left( 1 + \epsilon^{\overline{\eta} \overline{\sigma}} + \epsilon^{\overline{\eta} \overline{\rho}} \right) \left( \left( \xi_0^2 + |\xi|^2 \right) - \frac{1}{3} \left( \overline{\eta}_a^{\overline{\eta}} \right)^2 \right) 
\end{align*}
Since $\overline{\eta}_a^{\overline{\eta}} \simeq |\overline{\eta}|$ and $1 + \epsilon^{\overline{\eta} \overline{\sigma}} + \epsilon^{\overline{\eta} \overline{\rho}} \in \{ -1, 1, 3 \}$, we deduce that 
\begin{align*}
|\xi|^2 = \frac{1}{3} \left( \overline{\eta}_a^{\overline{\eta}} \right)^2 - \xi_0^2 + O(\overline{\eta}_b^{\overline{\eta}}) + O(\overline{\sigma}_b^{\overline{\sigma}}) + O(\overline{\rho}_b^{\overline{\rho}}) + O(\overline{\eta}_a^{\overline{\eta}} - \overline{\rho}_a^{\overline{\eta}}) + O(\overline{\sigma}_a^{\overline{\sigma}} - \overline{\rho}_a^{\overline{\eta}}) + O(\varphi_3) 
\end{align*}
As before, we only need to consider the neighborhood of a point where $\varphi_3 = o(1)$. 

Then, 
\begin{align*}
3 \xi_0^2 - |\xi|^2 
&= O(\overline{\eta}_b^{\overline{\eta}}) + O(\overline{\sigma}_b^{\overline{\sigma}}) + O(\overline{\rho}_b^{\overline{\rho}}) + O(\overline{\eta}_a^{\overline{\eta}} - \overline{\rho}_a^{\overline{\eta}}) + O(\overline{\sigma}_a^{\overline{\sigma}} - \overline{\rho}_a^{\overline{\eta}}) + O(\varphi_3) \\
&\quad + 4 \left( |\eta_0| + \epsilon^{\overline{\eta} \overline{\sigma}} |\sigma_0| + \epsilon^{\overline{\eta} \overline{\rho}} |\rho_0| \right)^2 - \frac{1}{3} \left( \overline{\eta}_a^{\overline{\eta}} \right)^2 \\
&= O(\overline{\eta}_b^{\overline{\eta}}) + O(\overline{\sigma}_b^{\overline{\sigma}}) + O(\overline{\rho}_b^{\overline{\rho}}) + O(\overline{\eta}_a^{\overline{\eta}} - \overline{\rho}_a^{\overline{\eta}}) + O(\overline{\sigma}_a^{\overline{\sigma}} - \overline{\rho}_a^{\overline{\eta}}) + O(\varphi_3) \\
&\quad + \frac{1}{3} \left( \left( 1 + \epsilon^{\overline{\eta} \overline{\sigma}} + \epsilon^{\overline{\eta} \overline{\rho}} \right)^2 - 1 \right) \left( \overline{\eta}_a^{\overline{\eta}} \right)^2
\end{align*}
In particular, if $\epsilon^{\overline{\eta} \overline{\sigma}} + \epsilon^{\overline{\eta} \overline{\rho}} = 0$, the presence of $m_b$ allows to conclude. 

On the other hand, if $\epsilon^{\overline{\eta} \overline{\sigma}} + \epsilon^{\overline{\eta} \overline{\rho}} = 2$, then $|\xi_0| \simeq |\overline{\xi}| \simeq |\overline{\eta}|$ and the previous computations ensures that $3 \xi_0^2 - |\xi|^2 \simeq |\overline{\xi}|^2$; finally, 
\begin{align*}
\xi_0^2 &= \left( |\eta_0| + |\sigma_0| + |\rho_0| \right)^2 \\
&= \frac{3}{4} \left( \overline{\eta}_a^{\overline{\eta}} \right)^2 + o(1) \\
\varphi_3 &\simeq \overline{\eta}_a^{\overline{\eta}} \left( \xi_0^2 + |\xi|^2 - \frac{1}{3} \left( \overline{\eta}_a^{\overline{\eta}} \right)^2 \right) \\
&\simeq \overline{\eta}_a^{\overline{\eta}} \left( |\xi|^2 + \left( \frac{3}{4} - \frac{1}{3} \right) \left( \overline{\eta}_a^{\overline{\eta}} \right)^2 \right) \\
&\simeq \left( \overline{\eta}_a^{\overline{\eta}} \right)^3 
\end{align*}
which also concludes. 

\paragraph{3.} Finally, let us consider the case where, up to exchanging the variables, we have $\epsilon^{\overline{\eta} \overline{\sigma}} \theta^{\overline{\eta} \overline{\sigma}}$ and $\epsilon^{\overline{\eta} \overline{\rho}} \theta^{\overline{\eta} \overline{\rho}}$ close to $-1$. As we saw, this forces $\epsilon^{\overline{\sigma} \overline{\rho}} \theta^{\overline{\sigma} \overline{\rho}}$ to be close to $1$. 

We can then apply the same integrations by parts as in the case 2. above to recover from \eqref{equdecchampbCCC-1} and \eqref{equdecchampbCCC-2} the term \eqref{equdecchampbCCC-cas2-final}. 

This time, we can compute that 
\begin{align*}
\widehat{Y}(\overline{\sigma}, \overline{\rho}) \cdot \nabla_{\overline{\sigma}} \varphi_3 &\simeq \left( \overline{\sigma}_a^{\overline{\sigma}} \right)^2 - \left( \overline{\rho}_a^{\overline{\rho}} \right)^2 \\
\widehat{Y}(\overline{\eta}, \overline{\rho}) \cdot \nabla_{\overline{\eta}} \varphi_3 &\simeq \xi_t^{\overline{\eta} \overline{\rho}} \left( \left( \overline{\eta}_a^{\overline{\eta}} \right)^2 - \left( \overline{\rho}_a^{\overline{\rho}} \right)^2 \right) 
\end{align*}
We only need to consider a neighborhood where these quantities are $o(1)$. 

We separate into subcases depending on the relative sizes. 

\paragraph{3.1.} Let us first localise on $|\overline{\eta}| \simeq |\overline{\sigma}| \simeq |\overline{\rho}|$. Then
\begin{align*}
\widehat{X}_a(\overline{\eta}) \cdot \nabla_{\overline{\eta}} \varphi_3 &= \frac{\eta_0}{|\overline{\eta}|} \left( 3 \eta_0^2 + |\eta|^2 - 3 \rho_0^2 - |\rho|^2 \right) + \frac{\eta}{|\overline{\eta}|} \cdot \left( 2 \eta_0 \eta - 2 \rho_0 \rho \right) \\
&= o(1) + \frac{\eta_0}{|\overline{\eta}|} \left( 3 \eta_0^2 + 3 |\eta|^2 - 3 \rho_0^2 - |\rho|^2 + 2 \sqrt{3} |\rho_0| |\rho| \right) \\
&= o(1) + 3 \eta_0 |\overline{\eta}| 
\end{align*}
so we can apply an integration by parts and get good terms if we have a factor $O(\overline{\rho}_b^{\overline{\rho}})$. Exchanging the role of $\overline{\eta}$ and $\overline{\rho}$, we can treat the same way $O(\overline{\eta}_b^{\overline{\eta}})$, and by symmetry $O(\overline{\sigma}_b^{\overline{\sigma}})$. 

\paragraph{3.1.1.} Assume that $\epsilon^{\overline{\sigma} \overline{\rho}} = -1$. Then $\epsilon^{\overline{\xi} \overline{\eta}} = 1$. 

We then have 
\begin{align*}
&\overline{\xi}_a^{\overline{\xi}} \overline{\xi}_b^{\overline{\xi}} = 3 \xi_0^2 - |\xi|^2 \\
&= 3 \eta_0^2 - |\eta|^2 + 3 \sigma_0^2 - |\sigma|^2 + 3 \rho_0^2 - |\rho|^2 + 6 \eta_0 \sigma_0 - 2 \eta \cdot \sigma + 6 \eta_0 \rho_0 - 2 \eta \cdot \rho + 6 \sigma_0 \rho_0 - 2 \sigma \cdot \rho \\
&= O(\overline{\eta}_b^{\overline{\eta}}) + O(\overline{\sigma}_b^{\overline{\sigma}}) + O(\overline{\rho}_b^{\overline{\rho}}) 
+ O((\overline{\eta}_a^{\overline{\eta}} - \overline{\rho}_a^{\overline{\rho}}) \xi_t^{\overline{\eta} \overline{\rho}}) 
+ O((\overline{\eta}_a^{\overline{\eta}} - \overline{\sigma}_a^{\overline{\sigma}}) \xi_t^{\overline{\eta} \overline{\sigma}}) 
+ O(\overline{\rho}_a^{\overline{\rho}} - \overline{\sigma}_a^{\overline{\sigma}}) \\
&\quad - \frac{1}{2} \left( \overline{\eta}_a^{\overline{\eta}} \right)^2 \left( 1 + \epsilon^{\overline{\eta} \overline{\sigma}} \theta^{\overline{\eta} \overline{\sigma}} + \epsilon^{\overline{\eta} \overline{\rho}} \theta^{\overline{\eta} \overline{\rho}} + \epsilon^{\overline{\sigma} \overline{\rho}} \theta^{\overline{\sigma} \overline{\rho}} \right) \\
&\left( \overline{\xi}_a^{\overline{\xi}} \right)^2 + \left( \overline{\xi}_b^{\overline{\xi}} \right)^2 = 6 \xi_0^2 + 2 |\xi|^2 \\
&= 6 \eta_0^2 + 2 |\eta|^2 + 6 \sigma_0^2 + 2 |\sigma|^2 + 6 \rho_0^2 + 2 |\rho|^2 + 12 \eta_0 \sigma_0 + 4 \eta \cdot \sigma + 12 \eta_0 \rho_0 + 4 \eta \cdot \rho + 12 \sigma_0 \rho_0 + 4 \sigma \cdot \rho \\
&= O(\overline{\eta}_b^{\overline{\eta}}) + O(\overline{\sigma}_b^{\overline{\sigma}}) + O(\overline{\rho}_b^{\overline{\rho}}) + O((\overline{\eta}_a^{\overline{\eta}} - \overline{\rho}_a^{\overline{\rho}}) \xi_t^{\overline{\eta} \overline{\rho}}) 
+ O((\overline{\eta}_a^{\overline{\eta}} - \overline{\sigma}_a^{\overline{\sigma}}) \xi_t^{\overline{\eta} \overline{\sigma}}) 
+ O(\overline{\rho}_a^{\overline{\rho}} - \overline{\sigma}_a^{\overline{\sigma}}) \\
&\quad + \left( \overline{\eta}_a^{\overline{\eta}} \right)^2 + \left( \overline{\eta}_a^{\overline{\eta}} \right)^2 \left( 1 + \epsilon^{\overline{\eta} \overline{\sigma}} \theta^{\overline{\eta} \overline{\sigma}} + \epsilon^{\overline{\eta} \overline{\rho}} \theta^{\overline{\eta} \overline{\rho}} + \epsilon^{\overline{\sigma} \overline{\rho}} \theta^{\overline{\sigma} \overline{\rho}} \right) 
\end{align*}
hence
\begin{align*}
\overline{\xi}_a^{\overline{\xi}} + \overline{\xi}_b^{\overline{\xi}} &= O(\overline{\eta}_b^{\overline{\eta}}) + O(\overline{\sigma}_b^{\overline{\sigma}}) + O(\overline{\rho}_b^{\overline{\rho}}) + O((\overline{\eta}_a^{\overline{\eta}} - \overline{\rho}_a^{\overline{\rho}}) \xi_t^{\overline{\eta} \overline{\rho}}) 
+ O((\overline{\eta}_a^{\overline{\eta}} - \overline{\sigma}_a^{\overline{\sigma}}) \xi_t^{\overline{\eta} \overline{\sigma}}) 
+ O(\overline{\rho}_a^{\overline{\rho}} - \overline{\sigma}_a^{\overline{\sigma}}) \\
&\quad + \overline{\eta}_a^{\overline{\eta}} \\
\overline{\xi}_a^{\overline{\xi}} - \overline{\xi}_b^{\overline{\xi}} &= O(\overline{\eta}_b^{\overline{\eta}}) + O(\overline{\sigma}_b^{\overline{\sigma}}) + O(\overline{\rho}_b^{\overline{\rho}}) + O((\overline{\eta}_a^{\overline{\eta}} - \overline{\rho}_a^{\overline{\rho}}) \xi_t^{\overline{\eta} \overline{\rho}}) 
+ O((\overline{\eta}_a^{\overline{\eta}} - \overline{\sigma}_a^{\overline{\sigma}}) \xi_t^{\overline{\eta} \overline{\sigma}}) 
+ O(\overline{\rho}_a^{\overline{\rho}} - \overline{\sigma}_a^{\overline{\sigma}}) \\
&\quad + \overline{\eta}_a^{\overline{\eta}} + \overline{\eta}_a^{\overline{\eta}} \left( 1 + \epsilon^{\overline{\eta} \overline{\sigma}} \theta^{\overline{\eta} \overline{\sigma}} + \epsilon^{\overline{\eta} \overline{\rho}} \theta^{\overline{\eta} \overline{\rho}} + \epsilon^{\overline{\sigma} \overline{\rho}} \theta^{\overline{\sigma} \overline{\rho}} \right) 
+ O\left( \left( 1 + \epsilon^{\overline{\eta} \overline{\sigma}} \theta^{\overline{\eta} \overline{\sigma}} + \epsilon^{\overline{\eta} \overline{\rho}} \theta^{\overline{\eta} \overline{\rho}} + \epsilon^{\overline{\sigma} \overline{\rho}} \theta^{\overline{\sigma} \overline{\rho}} \right)^2 \right) \\
\overline{\xi}_a^{\overline{\xi}} &= \overline{\eta}_a^{\overline{\eta}} + \frac{1}{2} \overline{\eta}_a^{\overline{\eta}} \left( 1 + \epsilon^{\overline{\eta} \overline{\sigma}} \theta^{\overline{\eta} \overline{\sigma}} + \epsilon^{\overline{\eta} \overline{\rho}} \theta^{\overline{\eta} \overline{\rho}} + \epsilon^{\overline{\sigma} \overline{\rho}} \theta^{\overline{\sigma} \overline{\rho}} \right) 
+ O\left( \left( 1 + \epsilon^{\overline{\eta} \overline{\sigma}} \theta^{\overline{\eta} \overline{\sigma}} + \epsilon^{\overline{\eta} \overline{\rho}} \theta^{\overline{\eta} \overline{\rho}} + \epsilon^{\overline{\sigma} \overline{\rho}} \theta^{\overline{\sigma} \overline{\rho}} \right)^2 \right) \\
&\quad + O(\overline{\eta}_b^{\overline{\eta}}) + O(\overline{\sigma}_b^{\overline{\sigma}}) + O(\overline{\rho}_b^{\overline{\rho}}) + O((\overline{\eta}_a^{\overline{\eta}} - \overline{\rho}_a^{\overline{\rho}}) \xi_t^{\overline{\eta} \overline{\rho}}) 
+ O((\overline{\eta}_a^{\overline{\eta}} - \overline{\sigma}_a^{\overline{\sigma}}) \xi_t^{\overline{\eta} \overline{\sigma}}) 
+ O(\overline{\rho}_a^{\overline{\rho}} - \overline{\sigma}_a^{\overline{\sigma}}) \\
\overline{\xi}_b^{\overline{\xi}} &= - \frac{1}{2} \overline{\eta}_a^{\overline{\eta}} \left( 1 + \epsilon^{\overline{\eta} \overline{\sigma}} \theta^{\overline{\eta} \overline{\sigma}} + \epsilon^{\overline{\eta} \overline{\rho}} \theta^{\overline{\eta} \overline{\rho}} + \epsilon^{\overline{\sigma} \overline{\rho}} \theta^{\overline{\sigma} \overline{\rho}} \right) 
+ O\left( \left( 1 + \epsilon^{\overline{\eta} \overline{\sigma}} \theta^{\overline{\eta} \overline{\sigma}} + \epsilon^{\overline{\eta} \overline{\rho}} \theta^{\overline{\eta} \overline{\rho}} + \epsilon^{\overline{\sigma} \overline{\rho}} \theta^{\overline{\sigma} \overline{\rho}} \right)^2 \right) \\
&\quad + O(\overline{\eta}_b^{\overline{\eta}}) + O(\overline{\sigma}_b^{\overline{\sigma}}) + O(\overline{\rho}_b^{\overline{\rho}}) + O((\overline{\eta}_a^{\overline{\eta}} - \overline{\rho}_a^{\overline{\rho}}) \xi_t^{\overline{\eta} \overline{\rho}}) 
+ O((\overline{\eta}_a^{\overline{\eta}} - \overline{\sigma}_a^{\overline{\sigma}}) \xi_t^{\overline{\eta} \overline{\sigma}}) 
+ O(\overline{\rho}_a^{\overline{\rho}} - \overline{\sigma}_a^{\overline{\sigma}})
\end{align*}
Therefore, 
\begin{align*}
6 \sqrt{3} \frac{\eta_0}{|\eta_0|} \varphi_3 &= \left( \overline{\xi}_a^{\overline{\xi}} \right)^3 + \left( \overline{\xi}_b^{\overline{\xi}} \right)^3 - \left( \overline{\eta}_a^{\overline{\eta}} \right)^3 - \left( \overline{\eta}_b^{\overline{\eta}} \right)^3 - \epsilon^{\overline{\eta} \overline{\rho}} \left( \overline{\rho}_a^{\overline{\rho}} \right)^3 - \epsilon^{\overline{\eta} \overline{\rho}} \left( \overline{\rho}_b^{\overline{\rho}} \right)^3 - \epsilon^{\overline{\eta} \overline{\sigma}} \left( \overline{\sigma}_a^{\overline{\sigma}} \right)^3 - \epsilon^{\overline{\eta} \overline{\sigma}} \left( \overline{\sigma}_b^{\overline{\sigma}} \right)^3 \\
&= O(\overline{\eta}_b^{\overline{\eta}}) + O(\overline{\sigma}_b^{\overline{\sigma}}) + O(\overline{\rho}_b^{\overline{\rho}}) + O((\overline{\eta}_a^{\overline{\eta}} - \overline{\rho}_a^{\overline{\rho}}) \xi_t^{\overline{\eta} \overline{\rho}}) 
+ O((\overline{\eta}_a^{\overline{\eta}} - \overline{\sigma}_a^{\overline{\sigma}}) \xi_t^{\overline{\eta} \overline{\sigma}}) 
+ O(\overline{\rho}_a^{\overline{\rho}} - \overline{\sigma}_a^{\overline{\sigma}}) \\
&\quad + \frac{3}{2} \left( \overline{\eta}_a^{\overline{\eta}} \right)^3 \left( 1 + \epsilon^{\overline{\eta} \overline{\sigma}} \theta^{\overline{\eta} \overline{\sigma}} + \epsilon^{\overline{\eta} \overline{\rho}} \theta^{\overline{\eta} \overline{\rho}} + \epsilon^{\overline{\sigma} \overline{\rho}} \theta^{\overline{\sigma} \overline{\rho}} \right) \left( 1 + o(1) \right) 
\end{align*}
We deduce that 
\begin{align*}
\overline{\xi}_b^{\overline{\xi}} &= O(\varphi_3) + O(\overline{\eta}_b^{\overline{\eta}}) + O(\overline{\sigma}_b^{\overline{\sigma}}) + O(\overline{\rho}_b^{\overline{\rho}}) + O((\overline{\eta}_a^{\overline{\eta}} - \overline{\rho}_a^{\overline{\rho}}) \xi_t^{\overline{\eta} \overline{\rho}}) 
+ O((\overline{\eta}_a^{\overline{\eta}} - \overline{\sigma}_a^{\overline{\sigma}}) \xi_t^{\overline{\eta} \overline{\sigma}}) 
+ O(\overline{\rho}_a^{\overline{\rho}} - \overline{\sigma}_a^{\overline{\sigma}})
\end{align*}
which concludes since, here, $|\overline{\xi}| \simeq |\overline{\eta}|$. 

\paragraph{3.1.2.} Assume now that $\epsilon^{\overline{\sigma} \overline{\rho}} = 1$. 

We then have
\begin{align*}
\rho_0 &= \sigma_0 + o(1) \\
\rho &= \sigma + o(1) \\
|\sigma|^2 &= 3 \sigma_0^2 + o(1) \\
|\eta|^2 &= 3 \eta_0^2 + o(1) \\
\sigma \cdot \eta &= |\sigma| |\eta| \theta^{\overline{\sigma} \overline{\eta}} = - 3 \epsilon^{\overline{\sigma} \overline{\eta}} |\sigma_0| |\eta_0| + o(1) = - 3 \eta_0 \sigma_0 + o(1) 
\end{align*}
Therefore, 
\begin{align*}
\varphi_3 &= \xi_0^3 + \xi_0 |\xi|^2 - \eta_0^3 - \eta_0 |\eta|^2 - \sigma_0^3 - \sigma_0 |\sigma|^2 - \rho_0^3 - \rho_0 |\rho|^2 \\
&= \left( \eta_0 + 2 \sigma_0 \right)^3 + \left( \eta_0 + 2 \sigma_0 \right) \left( |\eta|^2 + 4 \eta \cdot \sigma + 4 |\sigma|^2 \right) - 4 \eta_0^3 - 8 \sigma_0^3 + o(1) \\
&= \left( \eta_0 + 2 \sigma_0 \right) \left( 4 \eta_0^2 - 8 \eta_0 \sigma_0 + 16 \sigma_0^2 \right) - 4 \eta_0^3 - 8 \sigma_0^3 + o(1) \\
&= 24 \sigma_0^3 + o(1) 
\end{align*}
In particular, if we localised to have the $o(1)$ small enough (which is always possible up to reusing already treated cases), we have locally $1 = O(\varphi_3)$ and we get terms from \eqref{qteslemestimeesgeneriquesdeccubiquehb2}. 

\paragraph{3.2.} Let us now localise on $|\overline{\eta}| \gg |\overline{\sigma}| + |\overline{\rho}|$. Then using $\widehat{Y}(\overline{\eta}, \overline{\rho}) \cdot \nabla_{\overline{\eta}}$ we have a control over terms with a factor $\xi_t^{\overline{\eta} \overline{\rho}}$, and likewise $\xi_t^{\overline{\eta} \overline{\sigma}}$, so all $\xi_t$. 

Moreover, computing
\begin{align*}
\widehat{X}_a(\overline{\eta}) \cdot \nabla_{\overline{\eta}} \varphi_3 &= o(1) - \frac{\eta_0}{|\overline{\eta}|} \left( 3 \eta_0^2 + 3 |\eta|^2 \right) 
\end{align*}
we deduce that we have a good control if we have a factor $O(\overline{\rho}_b^{\overline{\rho}})$, and likewise $O(\overline{\sigma}_b^{\overline{\sigma}})$. Finally, 
\begin{align*}
\widehat{X}_a(\overline{\rho}) \cdot \nabla_{\overline{\eta}} \varphi_3 &= \frac{\rho_0}{|\overline{\rho}|} \left( 3 \rho_0^2 + |\rho|^2 - 3 \eta_0^2 - |\eta|^2 \right) + \frac{\rho}{|\overline{\rho}|} \cdot \left( 2 \rho_0 \rho - 2 \eta_0 \eta \right) \\
&= O\left( \overline{\rho}_b^{\overline{\rho}} \right) + O\left( \xi_t^{\overline{\eta} \overline{\rho}} \right)
+ \frac{\rho_0}{|\overline{\rho}|} \left( 3 \rho_0^2 + 3 |\rho|^2 - 3 \eta_0^2 - |\eta|^2 + 2 \sqrt{3} |\eta_0| |\eta| \right) \\
&= O\left( \overline{\rho}_b^{\overline{\rho}} \right) + O\left( \xi_t^{\overline{\eta} \overline{\rho}} \right)
+ \frac{\rho_0}{|\overline{\rho}|} \left( \left( \overline{\rho}_a^{\overline{\rho}} \right)^2 - \left( \overline{\eta}_b^{\overline{\eta}} \right)^2 \right) 
\end{align*}
which allows to control $O\left( \overline{\eta}_b^{\overline{\eta}} \left( \left( \overline{\rho}_a^{\overline{\rho}} \right)^2 - \left( \overline{\eta}_b^{\overline{\eta}} \right)^2 \right) \right)$. 

On the other hand, since $\overline{\xi}$ is close to $\overline{\eta}$, we can compute that 
\begin{align*}
\overline{\xi}_a^{\overline{\xi}} &= \overline{\eta}_a^{\overline{\xi}} + \overline{\sigma}_a^{\overline{\xi}} + \overline{\rho}_a^{\overline{\xi}} \\
&= \overline{\eta}_a^{\overline{\eta}} + O(\overline{\sigma}_b^{\overline{\sigma}}) + O(\overline{\rho}_b^{\overline{\rho}}) + O(\xi_t^{\overline{\eta} \overline{\sigma}}) + O(\xi_t^{\overline{\eta} \overline{\rho}}) \\
\overline{\xi}_b^{\overline{\xi}} &= \overline{\eta}_b^{\overline{\eta}} + \epsilon^{\overline{\eta} \overline{\sigma}} \overline{\sigma}_a^{\overline{\sigma}} + \epsilon^{\overline{\eta} \overline{\rho}} \overline{\rho}_a^{\overline{\rho}} + O(\xi_t^{\overline{\eta} \overline{\sigma}}) + O(\xi_t^{\overline{\eta} \overline{\rho}})
\end{align*}

Now, 
\begin{align*}
&6 \sqrt{3} \frac{\eta_0}{|\eta_0|} \varphi_3 = \left( \overline{\xi}_a^{\overline{\xi}} \right)^3 + \left( \overline{\xi}_b^{\overline{\xi}} \right)^3 - \left( \overline{\eta}_a^{\overline{\eta}} \right)^3 - \left( \overline{\eta}_b^{\overline{\eta}} \right)^3 - \epsilon^{\overline{\eta} \overline{\sigma}} \left( \overline{\sigma}_a^{\overline{\sigma}} \right)^3 - \epsilon^{\overline{\eta} \overline{\sigma}} \left( \overline{\sigma}_b^{\overline{\sigma}} \right)^3 - \epsilon^{\overline{\eta} \overline{\rho}} \left( \overline{\rho}_a^{\overline{\rho}} \right)^3 - \epsilon^{\overline{\eta} \overline{\rho}} \left( \overline{\rho}_b^{\overline{\rho}} \right)^3 \\
&= O(\overline{\sigma}_b^{\overline{\sigma}}) + O(\overline{\rho}_b^{\overline{\rho}}) + O(\xi_t^{\overline{\eta} \overline{\rho}}) + O(\xi_t^{\overline{\eta} \overline{\sigma}}) + \left( \overline{\xi}_b^{\overline{\xi}} \right)^3 - \left( \overline{\eta}_b^{\overline{\eta}} \right)^3 - \epsilon^{\overline{\eta} \overline{\sigma}} \left( \overline{\sigma}_a^{\overline{\sigma}} \right)^3 - \epsilon^{\overline{\eta} \overline{\rho}} \left( \overline{\rho}_a^{\overline{\rho}} \right)^3 \\
&= O(\overline{\sigma}_b^{\overline{\sigma}}) + O(\overline{\rho}_b^{\overline{\rho}}) + O(\xi_t^{\overline{\eta} \overline{\rho}}) + O(\xi_t^{\overline{\eta} \overline{\sigma}}) + O\left( \overline{\eta}_b^{\overline{\eta}} \left( \left( \overline{\eta}_b^{\overline{\eta}} \right)^2 - \left( \overline{\rho}_a^{\overline{\rho}} \right)^2 \right) \right) 
+ O\left( \left( \overline{\sigma}_a^{\overline{\sigma}} \right)^2 - \left( \overline{\rho}_a^{\overline{\rho}} \right)^2 \right) \\
&\quad + \left( \overline{\xi}_b^{\overline{\xi}} \right)^3 - \left( \overline{\sigma}_a^{\overline{\sigma}} \right)^2 \left( \overline{\eta}_b^{\overline{\eta}} + \epsilon^{\overline{\eta} \overline{\sigma}} \overline{\sigma}_a^{\overline{\sigma}} + \epsilon^{\overline{\eta} \overline{\rho}} \overline{\rho}_a^{\overline{\rho}} \right) \\
&= O(\overline{\sigma}_b^{\overline{\sigma}}) + O(\overline{\rho}_b^{\overline{\rho}}) + O(\xi_t^{\overline{\eta} \overline{\rho}}) + O(\xi_t^{\overline{\eta} \overline{\sigma}}) + O\left( \overline{\eta}_b^{\overline{\eta}} \left( \left( \overline{\eta}_b^{\overline{\eta}} \right)^2 - \left( \overline{\rho}_a^{\overline{\rho}} \right)^2 \right) \right) 
+ O\left( \left( \overline{\sigma}_a^{\overline{\sigma}} \right)^2 - \left( \overline{\rho}_a^{\overline{\rho}} \right)^2 \right) \\
&\quad + \overline{\xi}_b^{\overline{\xi}} \left( \left( \overline{\xi}_b^{\overline{\xi}} \right)^2 - \left( \overline{\sigma}_a^{\overline{\sigma}} \right)^2 \right) 
\end{align*} 

But we also have {\footnotesize 
\begin{align*}
&m_b(\overline{\xi}) \left( \widehat{X}_b(\overline{\xi}) \cdot \nabla_{\overline{\xi}} + \frac{\rho_0}{|\rho_0|} P_b^b(\overline{\xi}, \overline{\rho}) \widehat{X}_{b-\widehat{\mathcal{C}}}(\overline{\eta}, \overline{\rho}) \cdot \nabla_{\overline{\eta}} \right) \varphi_3 \\
&\quad = \frac{1}{|\overline{\xi}|^4} \overline{\xi}_a^{\overline{\xi}} \overline{\xi}_b^{\overline{\xi}} \Biggl( |\xi|^2 \left( 3 \xi_0^2 + |\xi|^2 - 3 \rho_0^2 - |\rho|^2 \right) - \xi_0 \xi \cdot \left( 2 \xi_0 \xi - 2 \rho_0 \rho \right) \\
&\quad \quad + \frac{|\xi|^2 |\rho|^2 + \xi_0 \rho_0 \xi \cdot \rho}{|\overline{\rho}| |\rho|} \frac{2 \sqrt{3}}{7 + \epsilon^{\overline{\eta} \overline{\rho}} \theta^{\overline{\eta} \overline{\rho}}} \left( (\sqrt{3} |\rho_0| - |\rho|)^2 - (\sqrt{3} |\eta_0| + |\eta|)^2 - 2 \sqrt{3} \frac{1 + \epsilon^{\overline{\eta} \overline{\rho}} \theta^{\overline{\eta} \overline{\rho}}}{3} |\rho_0| |\rho| \right) \Biggl) \\
&\quad = O(\overline{\rho}_b^{\overline{\rho}}) + O(\xi_t^{\overline{\eta} \overline{\rho}}) + O(\xi_t^{\overline{\eta} \overline{\sigma}}) + \frac{1}{|\overline{\xi}|^4} \overline{\xi}_a^{\overline{\xi}} \overline{\xi}_b^{\overline{\xi}} \Biggl( |\xi|^2 \left( 3 \xi_0^2 + |\xi|^2 - 3 \rho_0^2 - |\rho|^2 \right) - 2 \xi_0^2 |\xi|^2 
+ 2 |\xi_0| |\rho_0| |\rho| \left( - |\eta| + \epsilon^{\overline{\eta} \overline{\sigma}} |\sigma| + \epsilon^{\overline{\eta} \overline{\rho}} |\rho| \right)  \\
&\quad \quad - \frac{1}{6} \left( 3 |\xi|^2 + \sqrt{3} |\xi_0| \left( - |\eta| + \epsilon^{\overline{\eta} \overline{\sigma}} |\sigma| + \epsilon^{\overline{\eta} \overline{\rho}}  |\rho| \right) \right) \left( \overline{\eta}_a^{\overline{\eta}} \right)^2 \Biggl) \\
&\quad = O(\overline{\sigma}_b^{\overline{\sigma}}) + O(\overline{\rho}_b^{\overline{\rho}}) + O(\xi_t^{\overline{\eta} \overline{\rho}}) + O(\xi_t^{\overline{\eta} \overline{\sigma}}) + O\left( \left( \overline{\sigma}_a^{\overline{\sigma}} \right)^2 - \left( \overline{\rho}_a^{\overline{\rho}} \right)^2 \right) + \frac{1}{|\overline{\xi}|^4} \overline{\xi}_a^{\overline{\xi}} \overline{\xi}_b^{\overline{\xi}} \Biggl( 
\frac{1}{8} \left( \overline{\xi}_a^{\overline{\xi}} - \overline{\xi}_b^{\overline{\xi}} \right)^2 \left( \left( \overline{\xi}_a^{\overline{\xi}} \right)^2 + \left( \overline{\xi}_b^{\overline{\xi}} \right)^2 - \left( \overline{\rho}_a^{\overline{\rho}} \right)^2 \right) \\
&\quad \quad - \frac{1}{24} \left( \left( \overline{\xi}_a^{\overline{\xi}} \right)^2 - \left( \overline{\xi}_b^{\overline{\xi}} \right)^2 \right)^2 
+ \frac{1}{24} \left( \overline{\xi}_a^{\overline{\xi}} + \overline{\xi}_b^{\overline{\xi}} \right) \left( \overline{\sigma}_a^{\overline{\sigma}} \right)^2 \left( - \overline{\eta}_a^{\overline{\eta}} + \overline{\eta}_b^{\overline{\eta}} + \epsilon^{\overline{\eta} \overline{\sigma}} \overline{\sigma}_a^{\overline{\sigma}} + \epsilon^{\overline{\eta} \overline{\rho}} \overline{\rho}_a^{\overline{\rho}} \right)  \\
&\quad \quad - \frac{1}{24} \left( 3 \left( \overline{\xi}_a^{\overline{\xi}} - \overline{\xi}_b^{\overline{\xi}} \right)^2 + \left( \overline{\xi}_a^{\overline{\xi}} + \overline{\xi}_b^{\overline{\xi}} \right) \left( - \overline{\eta}_a^{\overline{\eta}} + \overline{\eta}_b^{\overline{\eta}} + \epsilon^{\overline{\eta} \overline{\sigma}} \overline{\sigma}_a^{\overline{\sigma}} + \epsilon^{\overline{\eta} \overline{\rho}}  \overline{\rho}_a^{\overline{\rho}} \right) \right) \left( \overline{\eta}_a^{\overline{\eta}} \right)^2 \Biggl) \\
&\quad = O(\overline{\sigma}_b^{\overline{\sigma}}) + O(\overline{\rho}_b^{\overline{\rho}}) + O(\xi_t^{\overline{\eta} \overline{\rho}}) + O(\xi_t^{\overline{\eta} \overline{\sigma}}) + O\left( \left( \overline{\sigma}_a^{\overline{\sigma}} \right)^2 - \left( \overline{\rho}_a^{\overline{\rho}} \right)^2 \right) + \frac{1}{|\overline{\xi}|^4} \overline{\xi}_a^{\overline{\xi}} \overline{\xi}_b^{\overline{\xi}} \Biggl( 
\frac{1}{8} \left( \overline{\eta}_a^{\overline{\eta}} - \overline{\xi}_b^{\overline{\xi}} \right)^2 \left( \left( \overline{\xi}_b^{\overline{\xi}} \right)^2 - \left( \overline{\sigma}_a^{\overline{\sigma}} \right)^2 \right) 
\\
&\quad \quad - \frac{1}{24} \left( \left( \overline{\eta}_a^{\overline{\eta}} \right)^2 - \left( \overline{\xi}_b^{\overline{\xi}} \right)^2 \right)^2 
+ \frac{1}{24} \left( \overline{\eta}_a^{\overline{\eta}} + \overline{\xi}_b^{\overline{\xi}} \right) \left( \overline{\sigma}_a^{\overline{\sigma}} \right)^2 \left( - \overline{\eta}_a^{\overline{\eta}} + \overline{\xi}_b^{\overline{\xi}} \right) - \frac{1}{24} \left( \overline{\eta}_a^{\overline{\eta}} + \overline{\xi}_b^{\overline{\xi}} \right) \left( - \overline{\eta}_a^{\overline{\eta}} + \overline{\xi}_b^{\overline{\xi}} \right) \left( \overline{\eta}_a^{\overline{\eta}} \right)^2 \Biggl) \\
&\quad = O(\overline{\sigma}_b^{\overline{\sigma}}) + O(\overline{\rho}_b^{\overline{\rho}}) + O(\xi_t^{\overline{\eta} \overline{\rho}}) + O(\xi_t^{\overline{\eta} \overline{\sigma}}) + O\left( \left( \overline{\sigma}_a^{\overline{\sigma}} \right)^2 - \left( \overline{\rho}_a^{\overline{\rho}} \right)^2 \right) \\
&\quad \quad + \frac{1}{24 |\overline{\xi}|^4} \overline{\xi}_a^{\overline{\xi}} \overline{\xi}_b^{\overline{\xi}} \left( \overline{\eta}_a^{\overline{\eta}} - \overline{\xi}_b^{\overline{\xi}} \right) \Biggl( 
3 \left( \overline{\eta}_a^{\overline{\eta}} - \overline{\xi}_b^{\overline{\xi}} \right) \left( \left( \overline{\xi}_b^{\overline{\xi}} \right)^2 - \left( \overline{\sigma}_a^{\overline{\sigma}} \right)^2 \right) - \left( \overline{\eta}_a^{\overline{\eta}} + \overline{\xi}_b^{\overline{\xi}} \right) \left( \overline{\sigma}_a^{\overline{\sigma}} \right)^2 + \left( \overline{\eta}_a^{\overline{\eta}} + \overline{\xi}_b^{\overline{\xi}} \right) \left( \overline{\xi}_b^{\overline{\xi}} \right)^2 \Biggl) \\
&\quad = O(\overline{\sigma}_b^{\overline{\sigma}}) + O(\overline{\rho}_b^{\overline{\rho}}) + O(\xi_t^{\overline{\eta} \overline{\rho}}) + O(\xi_t^{\overline{\eta} \overline{\sigma}}) + O\left( \left( \overline{\sigma}_a^{\overline{\sigma}} \right)^2 - \left( \overline{\rho}_a^{\overline{\rho}} \right)^2 \right) \\
&\quad \quad + \frac{1}{12 |\overline{\xi}|^4} \overline{\xi}_a^{\overline{\xi}} \overline{\xi}_b^{\overline{\xi}} \left( \overline{\eta}_a^{\overline{\eta}} - \overline{\xi}_b^{\overline{\xi}} \right) \left( 2 \overline{\eta}_a^{\overline{\eta}} - \overline{\xi}_b^{\overline{\xi}} \right) \Biggl( 
\left( \overline{\xi}_b^{\overline{\xi}} \right)^2 
- \left( \overline{\sigma}_a^{\overline{\sigma}} \right)^2 \Biggl) \\
&\quad = O(\overline{\sigma}_b^{\overline{\sigma}}) + O(\overline{\rho}_b^{\overline{\rho}}) + O(\xi_t^{\overline{\eta} \overline{\rho}}) + O(\xi_t^{\overline{\eta} \overline{\sigma}}) + O\left( \left( \overline{\sigma}_a^{\overline{\sigma}} \right)^2 - \left( \overline{\rho}_a^{\overline{\rho}} \right)^2 \right) + O(\varphi_3) 
\end{align*} }
which concludes, using that $\mu$ allows to distribute derivatives.  

\paragraph{3.3.} Finally, let us localise on $|\overline{\eta}| \ll |\overline{\sigma}| \simeq |\overline{\rho}|$. Then using $\widehat{Y}$ we have good control over the factor $O(\overline{\sigma}_a^{\overline{\sigma}} - \overline{\rho}_a^{\overline{\rho}})$, $O(\xi_t^{\overline{\eta} \overline{\rho}})$, $O(\xi_t^{\overline{\eta} \overline{\sigma}})$, and so every $\xi_t$. 

On the other hand, as before, using $\widehat{X}_a(\overline{\rho}) \cdot \nabla_{\overline{\eta}} \varphi_3$, we also have a control over terms with a factor $\overline{\eta}_b^{\overline{\eta}}$. 

Then, if $\epsilon^{\overline{\sigma} \overline{\rho}} = 1$, 
\begin{align*}
\varphi_3 &= \xi_0^3 + \xi_0 |\xi|^2 - \eta_0^3 - \eta_0 |\eta|^2 - \sigma_0^3 - \sigma_0 |\sigma|^2 - \rho_0^3 - \rho_0 |\rho|^2 \\
&= o(1) + 4 \xi_0^3 - 4 \sigma_0^3 - 4 \rho_0^3 \\
&= o(1) + 24 \rho_0^3 
\end{align*}
so $1 = O(\varphi_3)$, which concludes. 

Consider now the case $\epsilon^{\overline{\sigma} \overline{\rho}} = -1$, which forces $|\overline{\xi}| \ll |\overline{\rho}| \simeq |\overline{\sigma}|$. Then 
\begin{align*}
\xi_0 m_b(\overline{\xi}) &= \frac{\xi_0 |\xi|}{|\overline{\xi}|^2} \frac{3 \xi_0^2 - |\xi|^2}{|\overline{\xi}|} \\
&= \frac{\xi_0 |\xi|}{|\overline{\xi}|^2} \left( \frac{\overline{\xi}_a^{\overline{\rho}} \overline{\xi}_b^{\overline{\rho}}}{|\overline{\xi}|} + \frac{1}{|\overline{\xi}|} \left( \frac{J \rho \cdot \xi}{|\rho|} \right)^2 \right) \\
&= O(\xi_t^{\overline{\eta} \overline{\rho}}) + O(\xi_t^{\overline{\eta} \overline{\sigma}}) + O(\overline{\xi}_a^{\overline{\rho}}) 
\end{align*}
But 
\begin{align*}
\overline{\xi}_a^{\overline{\rho}} &= \overline{\eta}_a^{\overline{\rho}} + \overline{\sigma}_a^{\overline{\rho}} + \overline{\rho}_a^{\overline{\rho}} \\
&= O(\xi_t^{\overline{\eta} \overline{\sigma}}) + O(\xi_t^{\overline{\eta} \overline{\rho}}) + O(\overline{\eta}_b^{\overline{\eta}}) + O(\overline{\sigma}_a^{\overline{\sigma}} - \overline{\rho}_a^{\overline{\rho}}) 
\end{align*}
which concludes. 

\subsubsection{Interaction \texorpdfstring{$\widehat{\mathcal{L}}\widehat{\mathcal{C}}\widehat{\mathcal{C}}$}{LCC}}

Let us consider: {\footnotesize 
\begin{subequations}
\begin{align}
&\xi_0 m_b(\overline{\xi}) \widehat{X}_b(\overline{\xi}) \cdot \nabla_{\overline{\xi}} \widehat{I}_{s \mu}^{\widehat{\mathcal{L}}\widehat{\mathcal{C}}\widehat{\mathcal{C}}}[F_1, F_2, F_3](t, \overline{\xi}) \notag \\
&\quad = \int_0^t \int \int i s^2 \xi_0 m_b(\overline{\xi}) \widehat{X}_b(\overline{\xi}) \cdot \nabla_{\overline{\xi}} \varphi_3 e^{i s \varphi_3} \mu(\overline{\xi}, \overline{\eta}, \overline{\sigma}) m_{\widehat{\mathcal{L}}}(\overline{\eta}) m_{\widehat{\mathcal{C}}}(\overline{\sigma}) m_{\widehat{\mathcal{C}}}(\overline{\rho}) \widehat{F}_1(s, \overline{\eta}) \widehat{F}_2(s, \overline{\sigma}) \widehat{F}_3(s, \overline{\rho}) ~ d\overline{\eta} d\overline{\sigma} ds \label{equdecchampbLCC-1} \\
&\quad \quad + \int_0^t \int \int e^{i s \varphi_3} s \xi_0 \mu(\overline{\xi}, \overline{\eta}, \overline{\sigma}) m_{\widehat{\mathcal{L}}}(\overline{\eta}) m_{\widehat{\mathcal{C}}}(\overline{\sigma}) m_{\widehat{\mathcal{C}}}(\overline{\rho}) \widehat{F}_1(s, \overline{\eta}) \widehat{F}_2(s, \overline{\sigma}) m_b(\overline{\xi}) \widehat{X}_b(\overline{\xi}) \cdot \nabla_{\overline{\xi}} \widehat{F}_3(s, \overline{\rho}) ~ d\overline{\eta} d\overline{\sigma} ds \label{equdecchampbLCC-2} \\
&\quad \quad + \int_0^t \int \int e^{i s \varphi_3} s \xi_0 m_b(\overline{\xi}) \widehat{X}_b(\overline{\xi}) \cdot \nabla_{\overline{\xi}} \left( \mu(\overline{\xi}, \overline{\eta}, \overline{\sigma}) m_{\widehat{\mathcal{L}}}(\overline{\eta}) m_{\widehat{\mathcal{C}}}(\overline{\sigma}) m_{\widehat{\mathcal{C}}}(\overline{\rho}) \right) \widehat{F}_1(s, \overline{\eta}) \widehat{F}_2(s, \overline{\sigma}) \widehat{F}_3(s, \overline{\rho}) ~ d\overline{\eta} d\overline{\sigma} ds \label{equdecchampbLCC-3} 
\end{align}
\end{subequations} }
\eqref{equdecchampbLCC-3} is of the form $\eqref{lemestimeesgeneriquesdeccubiquehb2-dersymb}+\eqref{lemestimeesgeneriquesdeccubiquehb2-dersymbbis}$. 

Furthermore, we can use the field $\widehat{X}_{b-\widehat{\mathcal{L}}}(\overline{\eta}, \overline{\rho})$ and apply an integration by parts to bring back \eqref{equdecchampbLCC-2} to terms of \eqref{qteslemestimeesgeneriquesdeccubiquehb2} plus {\footnotesize 
\begin{align*}
\int_0^t \int \int i s^2 m_b(\overline{\xi}) P_b^b(\overline{\xi}, \overline{\rho}) \widehat{X}_{b-\widehat{\mathcal{L}}}(\overline{\eta}, \overline{\rho}) \cdot \nabla_{\overline{\eta}} \varphi_3 e^{i s \varphi_3} \xi_0 \mu(\overline{\xi}, \overline{\eta}, \overline{\sigma}) m_{\widehat{\mathcal{L}}}(\overline{\eta}) m_{\widehat{\mathcal{C}}}(\overline{\sigma}) m_{\widehat{\mathcal{C}}}(\overline{\rho}) \widehat{F}_1(s, \overline{\eta}) \widehat{F}_2(s, \overline{\sigma}) \widehat{F}_3(s, \overline{\rho}) ~ d\overline{\eta} d\overline{\sigma} ds
\end{align*} }
that we group with \eqref{equdecchampbLCC-1} to get {\footnotesize 
\begin{align*}
\int_0^t \int \int i s^2 m_b(\overline{\xi}) \left( \widehat{X}_b(\overline{\xi}) \cdot \nabla_{\overline{\xi}} + P_b^b(\overline{\xi}, \overline{\rho}) \widehat{X}_{b-\widehat{\mathcal{L}}}(\overline{\eta}, \overline{\rho}) \cdot \nabla_{\overline{\eta}} \right) \varphi_3 e^{i s \varphi_3} \xi_0 \mu(\overline{\xi}, \overline{\eta}, \overline{\sigma}) m_{\widehat{\mathcal{L}}}(\overline{\eta}) m_{\widehat{\mathcal{C}}}(\overline{\sigma}) m_{\widehat{\mathcal{C}}}(\overline{\rho}) \\
\widehat{F}_1(s, \overline{\eta}) \widehat{F}_2(s, \overline{\sigma}) \widehat{F}_3(s, \overline{\rho}) ~ d\overline{\eta} d\overline{\sigma} ds
\end{align*} }

\paragraph{1.} Let us first localise on $|\overline{\eta}| \simeq |\overline{\sigma}| \simeq |\overline{\rho}|$. 

Note that 
\begin{align*}
\nabla_{\eta} \varphi_3 &= o(1) + 2 \rho_0 \rho \simeq |\overline{\rho}|^2 
\end{align*}
In particular, if we have a factor $O(\overline{\rho}_b^{\overline{\rho}})$ we can apply an integration by parts along $\eta$ and get terms from \eqref{qteslemestimeesgeneriquesdeccubiquehb2}. Likewise with $O(\overline{\sigma}_b^{\overline{\sigma}})$. Then by Lemma \ref{lemcalculsconecoordonneesconiquesvarphi}, 
\begin{align*}
|\eta| \xi_t^{\overline{\eta} \overline{\rho}} &= O\left( \widehat{X}_c(\overline{\rho}) \cdot \nabla_{\overline{\eta}} \varphi_3 \right) \\
\left( \overline{\rho}_a^{\overline{\rho}} \right)^2 - \left( \overline{\eta}_a^{\overline{\rho}} \right)^2 &= O\left( \widehat{X}_a(\overline{\rho}) \cdot \nabla_{\overline{\eta}} \varphi_3 \right) + O\left( |\eta| \xi_t^{\overline{\eta} \overline{\rho}} \right) + O\left( \overline{\rho}_b^{\overline{\rho}} \right)
\end{align*}
So, up to symmetrizing, we have a good control on terms with a factor
\begin{align*}
O(|\eta| \xi_t^{\overline{\eta} \overline{\rho}}), \quad O\left( \left( \overline{\rho}_a^{\overline{\rho}} \right)^2 - \left( \overline{\eta}_a^{\overline{\rho}} \right)^2 \right), \quad O(|\eta| \xi_t^{\overline{\eta} \overline{\sigma}}), \quad O\left( \left( \overline{\sigma}_a^{\overline{\sigma}} \right)^2 - \left( \overline{\eta}_a^{\overline{\sigma}} \right)^2 \right)
\end{align*} 

Without loss of generality, we only consider a neighborhood where all these quantities are $o(1)$. 

Note that, locally, 
\begin{align*}
\overline{\eta}_a^{\overline{\sigma}} &= \epsilon^{\overline{\eta} \overline{\sigma}} \overline{\sigma}_a^{\overline{\sigma}} + O\left( \overline{\sigma}_b^{\overline{\sigma}} \right) + O\left( |\eta| \xi_t^{\overline{\eta} \overline{\sigma}} \right) + O\left( \left( \overline{\sigma}_a^{\overline{\sigma}} \right)^2 - \left( \overline{\eta}_a^{\overline{\sigma}} \right)^2 \right) \\
\overline{\rho}_a^{\overline{\rho}} &= \overline{\sigma}_a^{\overline{\sigma}} + O\left( \overline{\sigma}_b^{\overline{\sigma}} \right) + O\left( |\eta| \xi_t^{\overline{\eta} \overline{\sigma}} \right) + O\left( \left( \overline{\sigma}_a^{\overline{\sigma}} \right)^2 - \left( \overline{\eta}_a^{\overline{\sigma}} \right)^2 \right) + O\left( \overline{\rho}_b^{\overline{\rho}} \right) + O\left( |\eta| \xi_t^{\overline{\eta} \overline{\rho}} \right) + O\left( \left( \overline{\rho}_a^{\overline{\rho}} \right)^2 - \left( \overline{\eta}_a^{\overline{\rho}} \right)^2 \right)
\end{align*}

\paragraph{1.1.} Firts, we localise to have $\rho, \sigma$ not aligned, that is $\theta^{\overline{\rho} \overline{\sigma}}$ away enough from $\pm 1$. 

In this case, $\eta = O(|\eta| \xi_t^{\overline{\eta} \overline{\sigma}}) + O(|\eta| \xi_t^{\overline{\eta} \overline{\rho}})$. We can then simplify
\begin{align*}
\eta_0 &= 2 \epsilon^{\overline{\eta} \overline{\sigma}} \sigma_0 + O\left( \overline{\sigma}_b^{\overline{\sigma}} \right) + O\left( \left( \overline{\sigma}_a^{\overline{\sigma}} \right)^2 - \left( \overline{\eta}_a^{\overline{\sigma}} \right)^2 \right) + O(\eta) \\
\rho_0 &= \epsilon^{\overline{\rho} \overline{\sigma}} \sigma_0 + O\left( \overline{\sigma}_b^{\overline{\sigma}} \right) + O(\eta) + O\left( \left( \overline{\sigma}_a^{\overline{\sigma}} \right)^2 - \left( \overline{\eta}_a^{\overline{\sigma}} \right)^2 \right) + O\left( \overline{\rho}_b^{\overline{\rho}} \right) + O\left( \left( \overline{\rho}_a^{\overline{\rho}} \right)^2 - \left( \overline{\eta}_a^{\overline{\rho}} \right)^2 \right)
\end{align*}

In particular, $\xi_0 = \alpha \sigma_0 + R$ for some well-controlled remainder $R$, and $\alpha \in \{ -2, 0, 2, 4 \}$ entirely determined by the signs $\epsilon^{\overline{\eta} \overline{\sigma}}, \epsilon^{\overline{\rho} \overline{\sigma}}$. If $\alpha = 0$, the factor $\xi_0$ allows to conclude immediately. Then, if $\alpha \neq 0$, we can compute that {\footnotesize 
\begin{align*}
\varphi_3 &= \xi_0^3 + \xi_0 |\xi|^2 - \eta_0^3 - \eta_0 |\eta|^2 - \sigma_0^3 - \sigma_0 |\sigma|^2 - \rho_0^3 - \rho_0 |\rho|^2 \\
&= \xi_0 \left( \xi_0^2 + |\xi|^2 \right) - 4 \eta_0 \sigma_0^2 - 4 \sigma_0^3 - 4 \rho_0 \sigma_0^2 + O\left( \overline{\sigma}_b^{\overline{\sigma}} \right) + O(\eta) + O\left( \left( \overline{\sigma}_a^{\overline{\sigma}} \right)^2 - \left( \overline{\eta}_a^{\overline{\sigma}} \right)^2 \right) + O\left( \overline{\rho}_b^{\overline{\rho}} \right) + O\left( \left( \overline{\rho}_a^{\overline{\rho}} \right)^2 - \left( \overline{\eta}_a^{\overline{\rho}} \right)^2 \right) \\
&= \xi_0 \left( \alpha^2 \sigma_0^2 + |\xi|^2 - 4 \sigma_0^2 \right) + O\left( \overline{\sigma}_b^{\overline{\sigma}} \right) + O(\eta) + O\left( \left( \overline{\sigma}_a^{\overline{\sigma}} \right)^2 - \left( \overline{\eta}_a^{\overline{\sigma}} \right)^2 \right) + O\left( \overline{\rho}_b^{\overline{\rho}} \right) + O\left( \left( \overline{\rho}_a^{\overline{\rho}} \right)^2 - \left( \overline{\eta}_a^{\overline{\rho}} \right)^2 \right)
\end{align*} }
If $\alpha = 4$, $1 = O(\varphi)$ which concludes. If $\alpha = \pm 2$, then $|\xi|^2$ allows for good integrations by parts. But precisely, {\footnotesize 
\begin{align*}
&m_b(\overline{\xi}) \left( \widehat{X}_b(\overline{\xi}) \cdot \nabla_{\overline{\xi}} + P_b^b(\overline{\xi}, \overline{\rho}) \frac{\rho_0}{|\rho_0|} \widehat{X}_{b-\widehat{\mathcal{L}}}(\overline{\eta}, \overline{\rho}) \cdot \nabla_{\overline{\eta}} \right) \varphi_3 \\
&\quad = O(|\xi|) \left( \frac{|\xi|}{|\overline{\xi}|} \Biggl( 3 \xi_0^2 + |\xi|^2 - 3 \rho_0^2 - |\rho|^2 \right) - \frac{\xi_0 \xi}{|\overline{\xi}| |\xi|} \cdot \left( 2 \xi_0 \xi - 2 \rho_0 \rho \right) \\
&\quad \quad \quad + \frac{|\xi|^2 |\rho|^2 + \xi_0 \rho_0 \xi \cdot \rho}{|\overline{\xi}| |\overline{\rho}| |\xi| |\rho|} \frac{2 \sqrt{3}}{3} \left( 3 \rho_0^2 + |\rho|^2 - 3 \eta_0^2 - |\eta|^2 \right) \Biggl) \\
&\quad = O(|\xi|^2) + O(|\xi|) \left( 2 \rho_0 \frac{\xi_0 \xi \cdot \rho}{|\overline{\xi}| |\xi|} - 4 \sqrt{3} \sigma_0^2 \frac{\xi_0 \rho_0 \xi \cdot \rho}{|\overline{\xi}| |\overline{\rho}| |\xi| |\rho|} \right) + O\left( \overline{\sigma}_b^{\overline{\sigma}} \right) + O(\eta) + O\left( \left( \overline{\sigma}_a^{\overline{\sigma}} \right)^2 - \left( \overline{\eta}_a^{\overline{\sigma}} \right)^2 \right) \\
&\quad \quad \quad + O\left( \overline{\rho}_b^{\overline{\rho}} \right) + O\left( \left( \overline{\rho}_a^{\overline{\rho}} \right)^2 - \left( \overline{\eta}_a^{\overline{\rho}} \right)^2 \right) \\
&\quad = O(|\xi|^2) + O\left( \overline{\sigma}_b^{\overline{\sigma}} \right) + O(\eta) + O\left( \left( \overline{\sigma}_a^{\overline{\sigma}} \right)^2 - \left( \overline{\eta}_a^{\overline{\sigma}} \right)^2 \right) + O\left( \overline{\rho}_b^{\overline{\rho}} \right) + O\left( \left( \overline{\rho}_a^{\overline{\rho}} \right)^2 - \left( \overline{\eta}_a^{\overline{\rho}} \right)^2 \right)
\end{align*} }
which concludes. 

\paragraph{1.2.} Let us now localise to have $\epsilon^{\overline{\rho} \overline{\sigma}} \theta^{\overline{\rho} \overline{\sigma}}$ close to $1$. 

We then have that 
\begin{align*}
\frac{\sigma \cdot \eta}{|\sigma|} &= \overline{\eta}_a^{\overline{\sigma}} - \sqrt{3} \frac{\sigma_0 \eta_0}{|\sigma_0|} \\
&= \epsilon^{\overline{\sigma} \overline{\eta}} \overline{\sigma}_a^{\overline{\sigma}} - \sqrt{3} \frac{\sigma_0 \eta_0}{|\sigma_0|} + O\left( \left( \overline{\sigma}_a^{\overline{\sigma}} \right)^2 - \left( \overline{\eta}_a^{\overline{\sigma}} \right)^2 \right) \\
&= 2 \sqrt{3} \frac{\eta_0 \sigma_0}{|\eta_0|} - \sqrt{3} \frac{\sigma_0 \eta_0}{|\sigma_0|} + O\left( \left( \overline{\sigma}_a^{\overline{\sigma}} \right)^2 - \left( \overline{\eta}_a^{\overline{\sigma}} \right)^2 \right) + O\left( \overline{\sigma}_b^{\overline{\sigma}} \right) \\
\xi_0 &= \eta_0 + \sigma_0 + \rho_0 \\
&= \eta_0 + \sigma_0 + \frac{\rho_0}{2 \sqrt{3} |\rho_0|} \overline{\rho}_a^{\overline{\rho}} + O\left( \overline{\rho}_b^{\overline{\rho}} \right) \\
&= \eta_0 + \sigma_0 + \frac{\rho_0}{2 \sqrt{3} |\rho_0|} \overline{\sigma}_a^{\overline{\sigma}} + O\left( \overline{\sigma}_a^{\overline{\sigma}} - \overline{\rho}_a^{\overline{\rho}} \right) + O\left( \overline{\rho}_b^{\overline{\rho}} \right) \\
&= \eta_0 + \sigma_0 \left( 1 + \epsilon^{\overline{\rho} \overline{\sigma}} \right) + O\left( \overline{\sigma}_a^{\overline{\sigma}} - \overline{\rho}_a^{\overline{\rho}} \right) + O\left( \overline{\rho}_b^{\overline{\rho}} \right) + O\left( \overline{\sigma}_b^{\overline{\sigma}} \right) \\
\rho_0 &= \epsilon^{\overline{\rho} \overline{\sigma}} \sigma_0 + O\left( \overline{\sigma}_a^{\overline{\sigma}} - \overline{\rho}_a^{\overline{\rho}} \right) + O\left( \overline{\rho}_b^{\overline{\rho}} \right) + O\left( \overline{\sigma}_b^{\overline{\sigma}} \right)
\end{align*}
Therefore, 
\begin{align*}
\varphi_3 &= \xi_0^3 + \xi_0 |\xi|^2 - \eta_0^3 - \eta_0 |\eta|^2 - \sigma_0^3 - \sigma_0 |\sigma|^2 - \rho_0^3 - \rho_0 |\rho|^2 \\
&= \xi_0 \left( \xi_0^2 + |\xi|^2 \right) - \eta_0^3 - 3 \eta_0 \left( 2 \frac{\eta_0 \sigma_0}{|\eta_0|} - \frac{\sigma_0 \eta_0}{|\sigma_0|} \right)^2 - 4 \sigma_0^3 - 4 \rho_0^3 \\
&\quad + O\left( \overline{\sigma}_a^{\overline{\sigma}} - \overline{\rho}_a^{\overline{\rho}} \right) + O\left( \overline{\rho}_b^{\overline{\rho}} \right) + O\left( \overline{\sigma}_b^{\overline{\sigma}} \right) + O\left( |\eta| \xi_t^{\overline{\eta} \overline{\sigma}} \right) \\
&= \xi_0 \left( \xi_0^2 + |\xi|^2 \right) - 4 \left( \eta_0^3 + 3 \eta_0 \sigma_0^2 - 3 \epsilon^{\overline{\eta} \overline{\sigma}} \eta_0^2 \sigma_0 + \left( 1 + \epsilon^{\overline{\sigma} \overline{\rho}} \right) \sigma_0^3 \right) \\
&\quad + O\left( \overline{\sigma}_a^{\overline{\sigma}} - \overline{\rho}_a^{\overline{\rho}} \right) + O\left( \overline{\rho}_b^{\overline{\rho}} \right) + O\left( \overline{\sigma}_b^{\overline{\sigma}} \right) + O\left( |\eta| \xi_t^{\overline{\eta} \overline{\sigma}} \right) 
\end{align*}

\paragraph{1.2.1.} If $\epsilon^{\overline{\rho} \overline{\sigma}} = \epsilon^{\overline{\eta} \overline{\sigma}} = 1$, then we can directly write that 
\begin{align*}
\eta_0 &= 2 \sigma_0 + o(1) = 2 \rho_0 + o(1) \\
\xi_0 &= 4 \sigma_0 + o(1) \\
|\xi|^2 &= 4 |\sigma|^2 + o(1) = 12 \sigma_0^2 + o(1) \\
\varphi_3 &= o(1) + 4 \sigma_0 \left( 16 \sigma_0^2 + 12 \sigma_0^2 \right) - 4 \left( 8 \sigma_0^3 + 6 \sigma_0^3 - 12 \sigma_0^3 + 2 \sigma_0^3 \right) \\
&= o(1) + 24 \sigma_0^3 
\end{align*}
so $1 = O(\varphi_3)$. 

\paragraph{1.2.2.} If $\epsilon^{\overline{\rho} \overline{\sigma}} = 1$, $\epsilon^{\overline{\eta} \overline{\sigma}} = -1$, then we can simplify the expression of $\varphi_3$: 
\begin{align*}
\varphi_3 &= \xi_0 \left( \xi_0^2 + |\xi|^2 \right) - 4 \left( \eta_0^3 + 3 \eta_0 \sigma_0^2 + 3 \eta_0^2 \sigma_0 + 2 \sigma_0^3 \right) \\
&\quad + O\left( \overline{\sigma}_a^{\overline{\sigma}} - \overline{\rho}_a^{\overline{\rho}} \right) + O\left( \overline{\rho}_b^{\overline{\rho}} \right) + O\left( \overline{\sigma}_b^{\overline{\sigma}} \right) + O\left( |\eta| \xi_t^{\overline{\eta} \overline{\sigma}} \right) \\
&= \xi_0 \left( \xi_0^2 + |\xi|^2 \right) - 4 (\eta_0 + 2 \sigma_0) \left( \eta_0^2 + \sigma_0^2 + \eta_0 \sigma_0 \right) \\
&\quad + O\left( \overline{\sigma}_a^{\overline{\sigma}} - \overline{\rho}_a^{\overline{\rho}} \right) + O\left( \overline{\rho}_b^{\overline{\rho}} \right) + O\left( \overline{\sigma}_b^{\overline{\sigma}} \right) + O\left( |\eta| \xi_t^{\overline{\eta} \overline{\sigma}} \right) \\
&= \xi_0 \left( \xi_0^2 + |\xi|^2 - 4 \left( \eta_0^2 + \sigma_0^2 + \eta_0 \sigma_0 \right) \right) \\
&\quad + O\left( \overline{\sigma}_a^{\overline{\sigma}} - \overline{\rho}_a^{\overline{\rho}} \right) + O\left( \overline{\rho}_b^{\overline{\rho}} \right) + O\left( \overline{\sigma}_b^{\overline{\sigma}} \right) + O\left( |\eta| \xi_t^{\overline{\eta} \overline{\sigma}} \right)
\end{align*}
Note then that $\xi_0 = o(1)$ here, and that 
\begin{align*}
\frac{\sigma \cdot \eta}{|\sigma|} &= 2 \sqrt{3} \frac{\eta_0 \sigma_0}{|\eta_0|} - \sqrt{3} \frac{\sigma_0 \eta_0}{|\sigma_0|} + O\left( \left( \overline{\sigma}_a^{\overline{\sigma}} \right)^2 - \left( \overline{\eta}_a^{\overline{\sigma}} \right)^2 \right) + O\left( \overline{\sigma}_b^{\overline{\sigma}} \right) \\
&= \sqrt{3} \frac{\eta_0}{|\eta_0|} (2 \sigma_0 + \eta_0) + O\left( \left( \overline{\sigma}_a^{\overline{\sigma}} \right)^2 - \left( \overline{\eta}_a^{\overline{\sigma}} \right)^2 \right) + O\left( \overline{\sigma}_b^{\overline{\sigma}} \right) \\
&= \sqrt{3} \frac{\eta_0}{|\eta_0|} \xi_0 + O\left( \left( \overline{\sigma}_a^{\overline{\sigma}} \right)^2 - \left( \overline{\eta}_a^{\overline{\sigma}} \right)^2 \right) + O\left( \overline{\sigma}_b^{\overline{\sigma}} \right) + O\left( \overline{\sigma}_a^{\overline{\sigma}} - \overline{\rho}_a^{\overline{\rho}} \right) + O\left( \overline{\rho}_b^{\overline{\rho}} \right) 
\end{align*}
We deduce that 
\begin{align*}
|\eta| &= \sqrt{3} \frac{\eta_0}{|\eta_0|} \xi_0 + O\left( \left( \overline{\sigma}_a^{\overline{\sigma}} \right)^2 - \left( \overline{\eta}_a^{\overline{\sigma}} \right)^2 \right) + O\left( \overline{\sigma}_b^{\overline{\sigma}} \right) + O\left( \overline{\sigma}_a^{\overline{\sigma}} - \overline{\rho}_a^{\overline{\rho}} \right) + O\left( \overline{\rho}_b^{\overline{\rho}} \right) + O(|\eta| \xi_t^{\overline{\eta} \overline{\sigma}}) 
\end{align*}
Moreover, 
\begin{align*}
\overline{\eta}_a^{\overline{\sigma}} - \overline{\eta}_a^{\overline{\rho}} &= \overline{\rho}_a^{\overline{\rho}} - \overline{\sigma}_a^{\overline{\sigma}} + O\left( \left( \overline{\eta}_a^{\overline{\sigma}} \right)^2 - \left( \overline{\sigma}_a^{\overline{\sigma}} \right)^2 \right) + O\left( \left( \overline{\eta}_a^{\overline{\rho}} \right)^2 - \left( \overline{\rho}_a^{\overline{\rho}} \right)^2 \right) \\
&= O\left( \overline{\rho}_a^{\overline{\rho}} - \overline{\sigma}_a^{\overline{\sigma}} \right) + O\left( \left( \overline{\eta}_a^{\overline{\sigma}} \right)^2 - \left( \overline{\sigma}_a^{\overline{\sigma}} \right)^2 \right) + O\left( \left( \overline{\eta}_a^{\overline{\rho}} \right)^2 - \left( \overline{\rho}_a^{\overline{\rho}} \right)^2 \right)
\end{align*}
But on the other hand, 
\begin{align*}
\overline{\eta}_a^{\overline{\sigma}} - \overline{\eta}_a^{\overline{\rho}} &= \frac{\sigma \cdot \eta}{|\sigma|} - \frac{\rho \cdot \eta}{|\rho|} \\
&= \frac{\sigma}{|\sigma|} \cdot \left( \frac{\rho}{|\rho|} \frac{\eta \cdot \rho}{|\rho|} + \frac{J \rho}{|\rho|} \frac{\eta \cdot J \rho}{|\rho|} \right) - \frac{\rho \cdot \eta}{|\rho|} \\
&= O(|\eta| \xi_t^{\overline{\eta} \overline{\rho}}) + \frac{\eta \cdot \rho}{|\rho|} \left( \theta^{\overline{\sigma} \overline{\rho}} - 1 \right) 
\end{align*}
We can now simplify
\begin{align*}
\varphi_3 &= \xi_0 \left( \xi_0^2 + |\sigma|^2 + |\rho|^2 + |\eta|^2 + 2 \sigma \cdot \rho + 2 \sigma \cdot \eta + 2 \rho \cdot \eta - 4 \left( \eta_0^2 + \sigma_0^2 + \eta_0 \sigma_0 \right) \right) \\
&\quad + O\left( \overline{\sigma}_a^{\overline{\sigma}} - \overline{\rho}_a^{\overline{\rho}} \right) + O\left( \overline{\rho}_b^{\overline{\rho}} \right) + O\left( \overline{\sigma}_b^{\overline{\sigma}} \right) + O\left( |\eta| \xi_t^{\overline{\eta} \overline{\sigma}} \right) \\
&= \xi_0 \left( \xi_0^2 + |\eta|^2 + 8 \sigma_0^2 - 4 \eta_0^2 - 4 \eta_0 \sigma_0 + 12 \frac{\eta_0}{|\eta_0|} |\sigma_0| \xi_0 \right) \\
&\quad + O\left( \overline{\sigma}_a^{\overline{\sigma}} - \overline{\rho}_a^{\overline{\rho}} \right) + O\left( \overline{\rho}_b^{\overline{\rho}} \right) + O\left( \overline{\sigma}_b^{\overline{\sigma}} \right) + O\left( |\eta| \xi_t^{\overline{\eta} \overline{\sigma}} \right) \\
&= 4 \xi_0 \left( \xi_0^2 + |\eta|^2 + \xi_0 (\sigma_0 - \eta_0) + 3 \frac{\eta_0}{|\eta_0|} |\sigma_0| \xi_0 \right) \\
&\quad + O\left( \overline{\sigma}_a^{\overline{\sigma}} - \overline{\rho}_a^{\overline{\rho}} \right) + O\left( \overline{\rho}_b^{\overline{\rho}} \right) + O\left( \overline{\sigma}_b^{\overline{\sigma}} \right) + O\left( |\eta| \xi_t^{\overline{\eta} \overline{\sigma}} \right) \\
&= 4 \xi_0 ( \xi_0^2 + |\eta|^2 ) + O\left( \overline{\sigma}_a^{\overline{\sigma}} - \overline{\rho}_a^{\overline{\rho}} \right) + O\left( \overline{\rho}_b^{\overline{\rho}} \right) + O\left( \overline{\sigma}_b^{\overline{\sigma}} \right) + O\left( |\eta| \xi_t^{\overline{\eta} \overline{\sigma}} \right)
\end{align*}
so we can apply integrations by parts if we have a factor $O(\xi_0^3)$. 

But we can compute that {\footnotesize 
\begin{align*}
&\xi_0 \left( \widehat{X}_a(\overline{\xi}) \cdot \nabla_{\overline{\xi}} + \frac{\rho_0}{|\rho_0|} P_b^b(\overline{\xi}, \overline{\rho}) \widehat{X}_{b-\widehat{\mathcal{L}}}(\overline{\eta}, \overline{\rho}) \cdot \nabla_{\overline{\eta}} \right) \varphi_3 \\
&\quad = \xi_0 \frac{|\xi|}{|\overline{\xi}|} \left( 3 \xi_0^2 + |\xi|^2 - 3 \rho_0^2 - |\rho|^2 \right) - \frac{\xi_0^2 \xi}{|\overline{\xi}| |\xi|} \cdot \left( 2 \xi_0 \xi - 2 \rho_0 \rho \right) 
+ \frac{2 \sqrt{3}}{3} \xi_0 \frac{|\xi|^2 |\rho|^2 + \xi_0 \rho_0 \xi \cdot \rho}{|\overline{\xi}| |\overline{\rho}| |\xi| |\rho|} \left( 3 \rho_0^2 + |\rho|^2 - 3 \eta_0^2 - |\eta|^2 \right) \\
&\quad = \frac{|\xi|}{|\overline{\xi}|} \xi_0 \left( |\xi|^2 - 3 \eta_0^2 \right) 
+ \frac{2 \sqrt{3}}{3} \frac{\xi_0 \rho_0 \xi \cdot \rho}{|\overline{\xi}| |\overline{\rho}| |\xi| |\rho|} \xi_0 \left( 
\left( \overline{\rho}_a^{\overline{\rho}} \right)^2 - 3 \eta_0^2 \right) 
\\
&\quad + O(\xi_0^3) + O\left( |\eta| \xi_t^{\overline{\eta} \overline{\sigma}} \right) + O\left( \overline{\rho}_b^{\overline{\rho}} \right) + O\left( \overline{\sigma}_b^{\overline{\sigma}} \right) + O\left( \left( \overline{\sigma}_a^{\overline{\sigma}} \right)^2 - \left( \overline{\eta}_a^{\overline{\sigma}} \right)^2 \right) + O\left( \left( \overline{\rho}_a^{\overline{\rho}} \right)^2 - \left( \overline{\eta}_a^{\overline{\rho}} \right)^2 \right) + O\left( |\eta| \xi_t^{\overline{\eta} \overline{\rho}} \right) \\
&\quad = \frac{|\xi|}{|\overline{\xi}|} \xi_0 \left( |\sigma|^2 + |\rho|^2 + 2 |\sigma| |\rho| + 2 \sigma \cdot \eta + 2 \rho \cdot \eta - 3 \eta_0^2 \right)  
\\
&\quad + O(\xi_0^3) + O\left( |\eta| \xi_t^{\overline{\eta} \overline{\sigma}} \right) + O\left( \overline{\rho}_b^{\overline{\rho}} \right) + O\left( \overline{\sigma}_b^{\overline{\sigma}} \right) + O\left( \left( \overline{\sigma}_a^{\overline{\sigma}} \right)^2 - \left( \overline{\eta}_a^{\overline{\sigma}} \right)^2 \right) + O\left( \left( \overline{\rho}_a^{\overline{\rho}} \right)^2 - \left( \overline{\eta}_a^{\overline{\rho}} \right)^2 \right) + O\left( |\eta| \xi_t^{\overline{\eta} \overline{\rho}} \right) \\
&\quad = \frac{|\xi|}{|\overline{\xi}|} \xi_0 \left( 12 \sigma_0^2 + 12 \frac{\eta_0}{|\eta_0|} |\sigma_0| \xi_0 - 3 \eta_0^2 \right)  
\\
&\quad + O(\xi_0^3) + O\left( |\eta| \xi_t^{\overline{\eta} \overline{\sigma}} \right) + O\left( \overline{\rho}_b^{\overline{\rho}} \right) + O\left( \overline{\sigma}_b^{\overline{\sigma}} \right) + O\left( \left( \overline{\sigma}_a^{\overline{\sigma}} \right)^2 - \left( \overline{\eta}_a^{\overline{\sigma}} \right)^2 \right) + O\left( \left( \overline{\rho}_a^{\overline{\rho}} \right)^2 - \left( \overline{\eta}_a^{\overline{\rho}} \right)^2 \right) + O\left( |\eta| \xi_t^{\overline{\eta} \overline{\rho}} \right) \\
&\quad = 3 \frac{|\xi|}{|\overline{\xi}|} \xi_0^2 \left( (2 \sigma_0 - \eta_0) - 4 \sigma_0 \right)  
\\
&\quad + O(\xi_0^3) + O\left( |\eta| \xi_t^{\overline{\eta} \overline{\sigma}} \right) + O\left( \overline{\rho}_b^{\overline{\rho}} \right) + O\left( \overline{\sigma}_b^{\overline{\sigma}} \right) + O\left( \left( \overline{\sigma}_a^{\overline{\sigma}} \right)^2 - \left( \overline{\eta}_a^{\overline{\sigma}} \right)^2 \right) + O\left( \left( \overline{\rho}_a^{\overline{\rho}} \right)^2 - \left( \overline{\eta}_a^{\overline{\rho}} \right)^2 \right) + O\left( |\eta| \xi_t^{\overline{\eta} \overline{\rho}} \right) \\
&\quad = O(\xi_0^3) + O\left( |\eta| \xi_t^{\overline{\eta} \overline{\sigma}} \right) + O\left( \overline{\rho}_b^{\overline{\rho}} \right) + O\left( \overline{\sigma}_b^{\overline{\sigma}} \right) + O\left( \left( \overline{\sigma}_a^{\overline{\sigma}} \right)^2 - \left( \overline{\eta}_a^{\overline{\sigma}} \right)^2 \right) + O\left( \left( \overline{\rho}_a^{\overline{\rho}} \right)^2 - \left( \overline{\eta}_a^{\overline{\rho}} \right)^2 \right) + O\left( |\eta| \xi_t^{\overline{\eta} \overline{\rho}} \right) 
\end{align*} }
which concludes. 

\paragraph{1.2.3.} If $\epsilon^{\overline{\eta} \overline{\sigma}} = 1$ and $\epsilon^{\overline{\sigma} \overline{\rho}} = -1$, then $\xi = o(1)$, $\xi_0 = \eta_0 + o(1)$, and more precisely, 
\begin{align*}
\xi_0 - \eta_0 &= \sigma_0 + \rho_0 \\
&= \frac{\sigma_0}{2 \sqrt{3} |\sigma_0|} \left( \overline{\sigma}_a^{\overline{\sigma}} - \overline{\rho}_a^{\overline{\rho}} \right) + O\left( \overline{\sigma}_b^{\overline{\sigma}} \right) + O\left( \overline{\rho}_b^{\overline{\rho}} \right) \\
&= O\left( \overline{\sigma}_a^{\overline{\sigma}} - \overline{\rho}_a^{\overline{\rho}} \right) + O\left( \overline{\sigma}_b^{\overline{\sigma}} \right) + O\left( \overline{\rho}_b^{\overline{\rho}} \right)
\end{align*}
so we can simplify 
\begin{align*}
\varphi_3 &= \xi_0 \left( \xi_0^2 + |\xi|^2 \right) - 4 \eta_0 \left( \eta_0^2 + 3 \sigma_0^2 - 3 \eta_0 \sigma_0 \right) \\
&\quad + O\left( \overline{\sigma}_a^{\overline{\sigma}} - \overline{\rho}_a^{\overline{\rho}} \right) + O\left( \overline{\rho}_b^{\overline{\rho}} \right) + O\left( \overline{\sigma}_b^{\overline{\sigma}} \right) + O\left( |\eta| \xi_t^{\overline{\eta} \overline{\sigma}} \right) \\
&= \xi_0 \left( |\xi|^2 - 3 \xi_0^2 - 12 \sigma_0^2 + 12 \xi_0 \sigma_0 \right) \\
&\quad + O\left( \overline{\sigma}_a^{\overline{\sigma}} - \overline{\rho}_a^{\overline{\rho}} \right) + O\left( \overline{\rho}_b^{\overline{\rho}} \right) + O\left( \overline{\sigma}_b^{\overline{\sigma}} \right) + O\left( |\eta| \xi_t^{\overline{\eta} \overline{\sigma}} \right) \\
&= \xi_0 \left( |\xi|^2 - 3 (\xi_0 - 2 \sigma_0)^2 \right) \\
&\quad + O\left( \overline{\sigma}_a^{\overline{\sigma}} - \overline{\rho}_a^{\overline{\rho}} \right) + O\left( \overline{\rho}_b^{\overline{\rho}} \right) + O\left( \overline{\sigma}_b^{\overline{\sigma}} \right) + O\left( |\eta| \xi_t^{\overline{\eta} \overline{\sigma}} \right)
\end{align*}
As before, up to the sign, 
\begin{align*}
\overline{\eta}_a^{\overline{\sigma}} + \overline{\eta}_a^{\overline{\rho}} 
&= O\left( \overline{\rho}_a^{\overline{\rho}} - \overline{\sigma}_a^{\overline{\sigma}} \right) + O\left( \left( \overline{\eta}_a^{\overline{\sigma}} \right)^2 - \left( \overline{\sigma}_a^{\overline{\sigma}} \right)^2 \right) + O\left( \left( \overline{\eta}_a^{\overline{\rho}} \right)^2 - \left( \overline{\rho}_a^{\overline{\rho}} \right)^2 \right)
\end{align*}
and also 
\begin{align*}
\overline{\eta}_a^{\overline{\sigma}} + \overline{\eta}_a^{\overline{\rho}} 
&= O(|\eta| \xi_t^{\overline{\eta} \overline{\rho}}) + \frac{\eta \cdot \rho}{|\rho|} \left( \theta^{\overline{\sigma} \overline{\rho}} + 1 \right) 
\end{align*}
Therefore, it is now sufficient to compute that {\footnotesize 
\begin{align*}
&|\xi| \left( \widehat{X}_a(\overline{\xi}) \cdot \nabla_{\overline{\xi}} + \frac{\rho_0}{|\rho_0|} P_b^b(\overline{\xi}, \overline{\rho}) \widehat{X}_{b-\widehat{\mathcal{L}}}(\overline{\eta}, \overline{\rho}) \cdot \nabla_{\overline{\eta}} \right) \varphi_3 \\
&\quad = \frac{|\xi|^2}{|\overline{\xi}|} \left( 3 \xi_0^2 + |\xi|^2 - 3 \rho_0^2 - |\rho|^2 \right) - \frac{\xi_0^2 \xi}{|\overline{\xi}|} \cdot \left( 2 \xi_0 \xi - 2 \rho_0 \rho \right) 
+ \frac{2 \sqrt{3}}{3} \frac{|\xi|^2 |\rho|^2 + \xi_0 \rho_0 \xi \cdot \rho}{|\overline{\xi}| |\overline{\rho}| |\rho|} \left( 3 \rho_0^2 + |\rho|^2 - 3 \eta_0^2 - |\eta|^2 \right) \\
&\quad = \frac{|\xi|^2}{|\overline{\xi}|} \left( \xi_0^2 - 3 \eta_0^2 \right) 
+ \frac{2 \sqrt{3}}{3} \frac{\xi_0 \rho_0 \xi \cdot \rho}{|\overline{\xi}| |\overline{\rho}| |\rho|} 2 \sqrt{3} \frac{\sigma_0 \eta_0}{|\sigma_0|} \frac{\sigma \cdot \eta}{|\sigma|} 
\\
&\quad + O(\varphi_3) + O\left( |\eta| \xi_t^{\overline{\eta} \overline{\sigma}} \right) + O\left( \overline{\rho}_b^{\overline{\rho}} \right) + O\left( \overline{\sigma}_b^{\overline{\sigma}} \right) + O\left( \left( \overline{\sigma}_a^{\overline{\sigma}} \right)^2 - \left( \overline{\eta}_a^{\overline{\sigma}} \right)^2 \right) + O\left( \left( \overline{\rho}_a^{\overline{\rho}} \right)^2 - \left( \overline{\eta}_a^{\overline{\rho}} \right)^2 \right) + O\left( |\eta| \xi_t^{\overline{\eta} \overline{\rho}} \right) \\
&\quad = - 2 \xi_0^2 \frac{|\xi|^2}{|\overline{\xi}|} 
- 2 \xi_0^2 \frac{\xi \cdot \rho}{|\overline{\xi}| |\rho|} \frac{\sigma \cdot \eta}{|\sigma|} 
\\
&\quad + O(\varphi_3) + O\left( |\eta| \xi_t^{\overline{\eta} \overline{\sigma}} \right) + O\left( \overline{\rho}_b^{\overline{\rho}} \right) + O\left( \overline{\sigma}_b^{\overline{\sigma}} \right) + O\left( \left( \overline{\sigma}_a^{\overline{\sigma}} \right)^2 - \left( \overline{\eta}_a^{\overline{\sigma}} \right)^2 \right) + O\left( \left( \overline{\rho}_a^{\overline{\rho}} \right)^2 - \left( \overline{\eta}_a^{\overline{\rho}} \right)^2 \right) + O\left( |\eta| \xi_t^{\overline{\eta} \overline{\rho}} \right) \\
&\quad = - 2 \frac{\xi_0^2}{|\overline{\xi}|} \left( |\xi|^2 + \frac{\xi \cdot \rho}{|\rho|} \frac{\sigma \cdot \eta}{|\sigma|} \right)
\\
&\quad + O(\varphi_3) + O\left( |\eta| \xi_t^{\overline{\eta} \overline{\sigma}} \right) + O\left( \overline{\rho}_b^{\overline{\rho}} \right) + O\left( \overline{\sigma}_b^{\overline{\sigma}} \right) + O\left( \left( \overline{\sigma}_a^{\overline{\sigma}} \right)^2 - \left( \overline{\eta}_a^{\overline{\sigma}} \right)^2 \right) + O\left( \left( \overline{\rho}_a^{\overline{\rho}} \right)^2 - \left( \overline{\eta}_a^{\overline{\rho}} \right)^2 \right) + O\left( |\eta| \xi_t^{\overline{\eta} \overline{\rho}} \right)
\end{align*} }
But
\begin{align*}
|\xi|^2 + \frac{\xi \cdot \rho}{|\rho|} \frac{\sigma \cdot \eta}{|\sigma|} &= |\xi|^2 + \frac{\eta \cdot \rho}{|\rho|} \frac{\sigma \cdot \eta}{|\sigma|} + \left( |\rho| + \frac{\sigma \cdot \rho}{|\rho|} \right) \frac{\sigma \cdot \eta}{|\sigma|} \\
&= |\xi|^2 - |\eta|^2 + O(|\eta| \xi_t^{\overline{\eta} \overline{\sigma}}) + O(|\eta| \xi_t^{\overline{\eta} \overline{\rho}}) + O(|\eta| (1 + \theta)) \\
&= O(\varphi_3) + O(|\eta| \xi_t^{\overline{\eta} \overline{\sigma}}) + O(|\eta| \xi_t^{\overline{\eta} \overline{\rho}}) + O(|\eta| (1 + \theta))
\end{align*}
which concludes. 

\paragraph{1.2.4.} If $\epsilon^{\overline{\eta} \overline{\sigma}} = -1 = \epsilon^{\overline{\sigma} \overline{\rho}}$, then $\epsilon^{\overline{\eta} \overline{\rho}} = 1$, so up to exchanging $\overline{\sigma}$ and $\overline{\rho}$ we recover case 1.2.3. above.  

\paragraph{1.3.} Localise now to have $\epsilon^{\overline{\sigma} \overline{\rho}} \theta^{\overline{\sigma} \overline{\rho}}$ close to $-1$. 

\paragraph{1.3.1.} If $\epsilon^{\overline{\eta} \overline{\sigma}} = \epsilon^{\overline{\rho} \overline{\sigma}} = 1$, then $\xi = o(1)$, $\xi_0 = 4 \sigma_0 + o(1)$ and $\eta_0 = 2 \sigma_0 + o(1) = 2 \rho_0 + o(1)$, thus 
\begin{align*}
\varphi_3 &= o(1) + \xi_0^3 - \eta_0^3 - 4 \sigma_0^3 - 4 \rho_0^3 \\
&= o(1) + \sigma_0^3 \left( 4^3 - 2^3 - 4 - 4 \right) \\
&= o(1) + 48 \sigma_0^3
\end{align*}
so $1 = O(\varphi_3)$. 

\paragraph{1.3.2.} If $\epsilon^{\overline{\eta} \overline{\sigma}} = -1$ and $\epsilon^{\overline{\rho} \overline{\sigma}} = 1$, then $\overline{\xi} = o(1)$. We can then compute that {\footnotesize 
\begin{align*}
\varphi_3 &= \xi_0^3 + \xi_0 |\xi|^2 - \eta_0^3 - \eta_0 |\eta|^2 - \sigma_0^3 - \sigma_0 |\sigma|^2 - \rho_0^3 - \rho_0 |\rho|^2 \\
&= o(\xi_0) - \eta_0^3 - \frac{1}{2} \eta_0 \left( \frac{\eta \cdot \sigma}{|\sigma|} \right)^2 - \frac{1}{2} \eta_0 \left( \frac{\eta \cdot \rho}{|\rho|} \right)^2 - 4 \sigma_0^3 - 4 \rho_0^3 \\
&\quad + O\left( |\eta| \xi_t^{\overline{\eta} \overline{\sigma}} \right) + O\left( |\eta| \xi_t^{\overline{\eta} \overline{\rho}} \right) + O\left( \overline{\rho}_b^{\overline{\rho}} \right) + O\left( \overline{\sigma}_b^{\overline{\sigma}} \right) \\
&= o(\xi_0) - \eta_0^3 - \frac{1}{2} \eta_0 \left( - 2 \sqrt{3} |\sigma_0| + \sqrt{3} |\eta_0| \right)^2 - \frac{1}{2} \eta_0 \left( - 2 \sqrt{3} |\rho_0| + \sqrt{3} |\eta_0| \right)^2 - 4 \sigma_0^3 - 4 \rho_0^3 \\
&\quad + O\left( |\eta| \xi_t^{\overline{\eta} \overline{\sigma}} \right) + O\left( |\eta| \xi_t^{\overline{\eta} \overline{\rho}} \right) + O\left( \overline{\rho}_b^{\overline{\rho}} \right) + O\left( \overline{\sigma}_b^{\overline{\sigma}} \right) + O\left( \left( \overline{\eta}_a^{\overline{\sigma}} \right)^2 - \left( \overline{\sigma}_a^{\overline{\sigma}} \right)^2 \right) + O\left( \left( \overline{\eta}_a^{\overline{\rho}} \right)^2 - \left( \overline{\rho}_a^{\overline{\rho}} \right)^2 \right) \\
&= o(\xi_0) - 4 \eta_0^3 - 6 \eta_0 \sigma_0^2 - 6 \eta_0 \rho_0^2 - 6 \eta_0^2 \sigma_0 - 6 \eta_0^2 \rho_0 - 4 \sigma_0^3 - 4 \rho_0^3 \\
&\quad + O\left( |\eta| \xi_t^{\overline{\eta} \overline{\sigma}} \right) + O\left( |\eta| \xi_t^{\overline{\eta} \overline{\rho}} \right) + O\left( \overline{\rho}_b^{\overline{\rho}} \right) + O\left( \overline{\sigma}_b^{\overline{\sigma}} \right) + O\left( \left( \overline{\eta}_a^{\overline{\sigma}} \right)^2 - \left( \overline{\sigma}_a^{\overline{\sigma}} \right)^2 \right) + O\left( \left( \overline{\eta}_a^{\overline{\rho}} \right)^2 - \left( \overline{\rho}_a^{\overline{\rho}} \right)^2 \right) \\
&= o(\xi_0) - 4 \eta_0^2 \xi_0 - 2 \eta_0 \sigma_0 \xi_0 - 2 \eta_0 \rho_0 \xi_0 - 4 \sigma_0^2 \xi_0 - 4 \rho_0^2 \xi_0 + 4 \rho_0 \sigma_0 \xi_0 \\
&\quad + O\left( |\eta| \xi_t^{\overline{\eta} \overline{\sigma}} \right) + O\left( |\eta| \xi_t^{\overline{\eta} \overline{\rho}} \right) + O\left( \overline{\rho}_b^{\overline{\rho}} \right) + O\left( \overline{\sigma}_b^{\overline{\sigma}} \right) + O\left( \left( \overline{\eta}_a^{\overline{\sigma}} \right)^2 - \left( \overline{\sigma}_a^{\overline{\sigma}} \right)^2 \right) + O\left( \left( \overline{\eta}_a^{\overline{\rho}} \right)^2 - \left( \overline{\rho}_a^{\overline{\rho}} \right)^2 \right) \\
&= o(\xi_0) + \xi_0 \left( - 4 \eta_0^2 - 2 \eta_0 \sigma_0 - 2 \eta_0 \rho_0 - 4 \sigma_0^2 - 4 \rho_0^2 + 4 \rho_0 \sigma_0 \right) \\
&\quad + O\left( |\eta| \xi_t^{\overline{\eta} \overline{\sigma}} \right) + O\left( |\eta| \xi_t^{\overline{\eta} \overline{\rho}} \right) + O\left( \overline{\rho}_b^{\overline{\rho}} \right) + O\left( \overline{\sigma}_b^{\overline{\sigma}} \right) + O\left( \left( \overline{\eta}_a^{\overline{\sigma}} \right)^2 - \left( \overline{\sigma}_a^{\overline{\sigma}} \right)^2 \right) + O\left( \left( \overline{\eta}_a^{\overline{\rho}} \right)^2 - \left( \overline{\rho}_a^{\overline{\rho}} \right)^2 \right) \\
\end{align*} }
But here $\eta_0 = - 2 \sigma_0 + o(1) = - 2 \rho_0 + o(1)$, so 
\begin{align*}
- 4 \eta_0^2 - 2 \eta_0 \sigma_0 - 2 \eta_0 \rho_0 - 4 \sigma_0^2 - 4 \rho_0^2 + 4 \rho_0 \sigma_0 
&= - 12 \sigma_0^2
\end{align*}
This justifies that {\footnotesize 
\begin{align*}
\xi_0 &= O(\varphi_3) + O\left( |\eta| \xi_t^{\overline{\eta} \overline{\sigma}} \right) + O\left( |\eta| \xi_t^{\overline{\eta} \overline{\rho}} \right) + O\left( \overline{\rho}_b^{\overline{\rho}} \right) + O\left( \overline{\sigma}_b^{\overline{\sigma}} \right) + O\left( \left( \overline{\eta}_a^{\overline{\sigma}} \right)^2 - \left( \overline{\sigma}_a^{\overline{\sigma}} \right)^2 \right) + O\left( \left( \overline{\eta}_a^{\overline{\rho}} \right)^2 - \left( \overline{\rho}_a^{\overline{\rho}} \right)^2 \right)
\end{align*} }
and it concludes. 

\paragraph{1.3.3.} If $\epsilon^{\overline{\eta} \overline{\sigma}} = 1$ and $\epsilon^{\overline{\rho} \overline{\sigma}} = -1$, then $\eta_0 = 2 \sigma_0 + o(1) = - 2 \rho_0 + o(1)$, $\xi_0 = 2 \sigma_0 + o(1)$, $\xi = 2 \sigma + o(1)$. 

We can directly compute that
\begin{align*}
\varphi_3 &= \xi_0^3 + \xi_0 |\xi|^2 - \eta_0^3 - \eta_0 |\eta|^2 - \sigma_0^3 - \sigma_0 |\sigma|^2 - \rho_0^3 - \rho_0 |\rho|^2 \\
&= o(1) + 4 \xi_0^3 - \eta_0^3 - 4 \sigma_0^3 - 4 \rho_0^3 \\
&= o(1) + 16 \sigma_0^3 
\end{align*}
so $1 = O(\varphi_3)$. 

\paragraph{1.3.4.} If $\epsilon^{\overline{\eta} \overline{\sigma}} = -1$ and $\epsilon^{\overline{\rho} \overline{\sigma}} = -1$, then up to exchanging $\overline{\sigma}$ and $\overline{\rho}$ we recover case 1.3.3. above. 

\paragraph{2.} Let us now consider the case $|\overline{\eta}| \ll |\overline{\sigma}| \simeq |\overline{\rho}|$. 

We have directly 
\begin{align*}
1 = O\left( \widehat{X}_a(\overline{\rho}) \cdot \nabla_{\overline{\eta}} \varphi_3 \right) 
\end{align*}
so we can automatically apply an integration by parts and get terms from \eqref{qteslemestimeesgeneriquesdeccubiquehb2}. 

\paragraph{3.} Let us now consider the case $|\overline{\sigma}| \ll |\overline{\eta}| \simeq |\overline{\rho}|$. 

Then a factor $O(\overline{\sigma})$ is enough: indeed, we always have 
\begin{align*}
1 = O\left( \widehat{X}_a(\overline{\sigma}) \cdot \left( \nabla_{\overline{\eta}} - \nabla_{\overline{\sigma}} \right) \varphi_3 \right) 
\end{align*}
so we can apply an integration by parts in this direction and we recover then terms from \eqref{qteslemestimeesgeneriquesdeccubiquehb2}. 

But now, up to terms $O(\overline{\sigma})$, we recover a quadratic interaction and we can apply the same decomposition as in the proof of the decomposition lemma. This easily brings terms from \eqref{qteslemestimeesgeneriquesdeccubiquehb2}. 

\paragraph{4.} Let us now consider the case $|\overline{\eta}| + |\overline{\sigma}| \ll |\overline{\rho}|$. 

We then control any term containing a factor $O(\overline{\sigma})$. Indeed, $1 = O\left( \widehat{X}_a(\overline{\rho}) \cdot \nabla_{\overline{\eta}} \varphi_3 \right)$ so we can apply an integration by parts in this direction and the presence of $O(\overline{\sigma})$ ensures that we have terms from \eqref{qteslemestimeesgeneriquesdeccubiquehb2}. 

As in case 3. above, we recover a quadratic case and can proceed the same way, using $\mu$ in order to distribute derivatives in the case of the symmetric term of the decomposition lemma. 

\paragraph{5.} Consider finally the term $|\overline{\sigma}| + |\overline{\rho}| \ll |\overline{\eta}|$. 

In this case, we control terms with a factor $O(\overline{\rho}_b^{\overline{\rho}})$: we have indeed
\begin{align*}
\partial_{\eta_0} \varphi_3 &= o(1) - 3 \eta_0^2 - |\eta|^2
\end{align*}
so $1 = O(\partial_{\eta_0} \varphi_3)$ and if we have an additional factor $\overline{\rho}_b^{\overline{\rho}}$ we only get terms from \eqref{qteslemestimeesgeneriquesdeccubiquehb2}. 

We can proceed similarly with a factor $O(\overline{\sigma}_b^{\overline{\sigma}})$. 

Furthermore, we can also control terms with a factor $O(\eta \overline{\sigma}), O(\eta \overline{\rho})$: since $1 = O\left( \widehat{X}_a(\overline{\rho}) \cdot \nabla_{\overline{\eta}} \varphi_3 \right)$, we can apply an integration by parts in this direction and indeed get terms from \eqref{qteslemestimeesgeneriquesdeccubiquehb2}. 

Finally, we also control terms with a factor $O\left( \overline{\sigma} \overline{\rho} \right)$ applying an integration by parts along $\widehat{X}_a(\overline{\rho}) \cdot \nabla_{\overline{\eta}}$ as well. 

We now write {\footnotesize 
\begin{align*}
\varphi_3 &= \xi_0^3 + \xi_0 |\xi|^2 - \eta_0^3 - \eta_0 |\eta|^2 - \sigma_0^3 - \sigma_0 |\sigma|^2 - \rho_0^3 - \rho_0 |\rho|^2 \\
&= (\xi_0 - \eta_0) \left( \xi_0^2 + \xi_0 \eta_0 + \eta_0^2 - |\eta|^2 \right) - (\sigma_0 + \rho_0) \left( \sigma_0^2 + \rho_0^2 \right) + \xi_0 \left( |\xi|^2 - |\eta|^2 \right) - (\sigma_0 + \rho_0) \left( |\sigma|^2 + |\rho|^2 \right) + O\left( \overline{\sigma} \overline{\rho} \right) \\
&= (\sigma_0+\rho_0) \left( \xi_0^2 + \xi_0 \eta_0 + \eta_0^2 - \sigma_0^2 - \rho_0^2 - |\sigma|^2 - |\rho|^2 \right) + \xi_0 \left( \xi - \eta \right) \cdot \left( \sigma + \rho \right) + O\left( \overline{\sigma} \overline{\rho} \right) + O\left( \eta \overline{\sigma} \right) + O\left( \eta \overline{\rho} \right) \\
&= (\sigma_0+\rho_0) \left( \xi_0^2 + \xi_0 \eta_0 + \eta_0^2 - \sigma_0^2 - \rho_0^2 - |\sigma|^2 - |\rho|^2 \right) + \xi_0 \left( |\sigma|^2 + |\rho|^2 \right) + O\left( \overline{\sigma} \overline{\rho} \right) + O\left( \eta \overline{\sigma} \right) + O\left( \eta \overline{\rho} \right) \\
&= (\sigma_0+\rho_0) \left( 3 \xi_0^2 + (\sigma_0 + \rho_0)^2 - 3 (\sigma_0 + \rho_0) \xi_0 - 4 \sigma_0^2 - 4 \rho_0^2 \right) + 3 \xi_0 \left( \sigma_0^2 + \rho_0^2 \right) + O\left( \overline{\sigma} \overline{\rho} \right) + O\left( \eta \overline{\sigma} \right) + O\left( \eta \overline{\rho} \right) \\
&\quad \quad + O\left( \overline{\rho}_b^{\overline{\rho}} \right) + O\left( \overline{\sigma}_b^{\overline{\sigma}} \right) \\
&= 3 \xi_0^2 (\sigma_0+\rho_0) 
+ (\sigma_0 + \rho_0)^3 - 4 (\sigma_0 + \rho_0) \left( \sigma_0^2 + \rho_0^2 \right)  + O\left( \overline{\sigma} \overline{\rho} \right) + O\left( \eta \overline{\sigma} \right) + O\left( \eta \overline{\rho} \right) + O\left( \overline{\rho}_b^{\overline{\rho}} \right) + O\left( \overline{\sigma}_b^{\overline{\sigma}} \right) \\
&= (\sigma_0 + \rho_0) \left( 3 \xi_0^2 + o(1) \right) + O\left( \overline{\sigma} \overline{\rho} \right) + O\left( \eta \overline{\sigma} \right) + O\left( \eta \overline{\rho} \right) + O\left( \overline{\rho}_b^{\overline{\rho}} \right) + O\left( \overline{\sigma}_b^{\overline{\sigma}} \right)
\end{align*} }
Therefore, 
\begin{align*}
\sigma_0 + \rho_0 &= O(\varphi_3) + O\left( \overline{\sigma} \overline{\rho} \right) + O\left( \eta \overline{\sigma} \right) + O\left( \eta \overline{\rho} \right) + O\left( \overline{\rho}_b^{\overline{\rho}} \right) + O\left( \overline{\sigma}_b^{\overline{\sigma}} \right)
\end{align*}

Finally, using that $\mu$ allows to distribute derivatives, and that $m_b(\overline{\xi})$ contrains a factor $|\xi|$, we get (symmetrizing $\overline{\sigma}$ and $\overline{\rho}$) 
\begin{align*}
|\xi| |\sigma| &= O\left( \eta \overline{\sigma} \right) + O\left( \overline{\rho} \overline{\sigma} \right) + O\left( \overline{\sigma}_b^{\overline{\sigma}} \right) 
+ 3 \sigma_0^2 \\
&= O\left( \eta \overline{\sigma} \right) + O\left( \overline{\rho} \overline{\sigma} \right) 
+ O\left( \eta \overline{\rho} \right) + O\left( \overline{\rho}_b^{\overline{\rho}} \right) + O\left( \overline{\sigma}_b^{\overline{\sigma}} \right) + O(\sigma_0 \varphi_3) 
\end{align*}
which concludes. 

\subsubsection{Interaction \texorpdfstring{$\widehat{\mathcal{C}}\widehat{\mathcal{C}}\widehat{\mathcal{P}}$}{CCP}}

Let us consider: {\footnotesize 
\begin{subequations}
\begin{align}
&\xi_0 m_b(\overline{\xi}) \widehat{X}_b(\overline{\xi}) \cdot \nabla_{\overline{\xi}} \widehat{I}_{s \mu}^{\widehat{\mathcal{C}}\widehat{\mathcal{C}}\widehat{\mathcal{P}}}[F_1, F_2, F_3](t, \overline{\xi}) \notag \\
&\quad = \int_0^t \int \int i s^2 \xi_0 m_b(\overline{\xi}) \widehat{X}_b(\overline{\xi}) \cdot \nabla_{\overline{\xi}} \varphi_3 e^{i s \varphi_3} \mu(\overline{\xi}, \overline{\eta}, \overline{\sigma}) m_{\widehat{\mathcal{C}}}(\overline{\eta}) m_{\widehat{\mathcal{C}}}(\overline{\sigma}) m_{\widehat{\mathcal{P}}}(\overline{\rho}) \widehat{F}_1(s, \overline{\eta}) \widehat{F}_2(s, \overline{\sigma}) \widehat{F}_3(s, \overline{\rho}) ~ d\overline{\eta} d\overline{\sigma} ds \label{equdecchampbCCP-1} \\
&\quad \quad + \int_0^t \int \int e^{i s \varphi_3} s \xi_0 \mu(\overline{\xi}, \overline{\eta}, \overline{\sigma}) m_{\widehat{\mathcal{C}}}(\overline{\eta}) m_{\widehat{\mathcal{C}}}(\overline{\sigma}) m_{\widehat{\mathcal{P}}}(\overline{\rho}) \widehat{F}_1(s, \overline{\eta}) \widehat{F}_2(s, \overline{\sigma}) m_b(\overline{\xi}) \widehat{X}_b(\overline{\xi}) \cdot \nabla_{\overline{\xi}} \widehat{F}_3(s, \overline{\rho}) ~ d\overline{\eta} d\overline{\sigma} ds \label{equdecchampbCCP-2} \\
&\quad \quad + \int_0^t \int \int e^{i s \varphi_3} s \xi_0 m_b(\overline{\xi}) \widehat{X}_b(\overline{\xi}) \cdot \nabla_{\overline{\xi}} \left( \mu(\overline{\xi}, \overline{\eta}, \overline{\sigma}) m_{\widehat{\mathcal{C}}}(\overline{\eta}) m_{\widehat{\mathcal{C}}}(\overline{\sigma}) m_{\widehat{\mathcal{P}}}(\overline{\rho}) \right) \widehat{F}_1(s, \overline{\eta}) \widehat{F}_2(s, \overline{\sigma}) \widehat{F}_3(s, \overline{\rho}) ~ d\overline{\eta} d\overline{\sigma} ds \label{equdecchampbCCP-3} 
\end{align}
\end{subequations} }
\eqref{equdecchampbCCP-3} is of the form $\eqref{lemestimeesgeneriquesdeccubiquehb2-dersymb}+\eqref{lemestimeesgeneriquesdeccubiquehb2-dersymbbis}$. 

Then, for \eqref{equdecchampbCCP-1} and \eqref{equdecchampbCCP-2}, we separate into several subcases depending on the relative sizes of frequencies. 

\paragraph{1.} Let us first localise to have $|\overline{\eta}| \simeq |\overline{\sigma}| \simeq |\overline{\rho}|$. 

Then \eqref{equdecchampbCCP-2} is already of the form \eqref{qteslemestimeesgeneriquesdeccubiquehb2}. 

Note that $\nabla_{\eta} \varphi_3 = 2 \rho_0 \rho - 2 \eta_0 \eta \simeq |\overline{\eta}|^2$, so we have a control over terms with a factor $O\left( \overline{\eta}_b^{\overline{\eta}} \right)$, or likewise $O\left( \overline{\sigma}_b^{\overline{\sigma}} \right)$. 

Then, by Lemma \ref{lemcalculsconecoordonneesconiquesvarphi}, 
\begin{align*}
\rho_0 \xi_t^{\overline{\eta} \overline{\rho}} &= O\left( \widehat{X}_c(\overline{\eta}) \cdot \nabla_{\overline{\eta}} \varphi_3 \right) 
\end{align*}
and we may compute that 
\begin{align*}
\widehat{X}_a(\overline{\eta}) \cdot \nabla_{\overline{\eta}} \varphi_3 &= \frac{\eta_0}{|\overline{\eta}|} \left( 3 \rho_0^2 + |\rho|^2 - 3 \eta_0^2 - |\eta|^2 \right) + \frac{\eta}{|\overline{\eta}|} \cdot \left( 2 \rho_0 \rho - 2 \eta_0 \eta \right) \\
&= O\left( \overline{\eta}_b^{\overline{\eta}} \right) + O\left( \rho_0 \xi_t^{\overline{\eta} \overline{\rho}} \right) + \frac{\eta_0}{|\overline{\eta}|} \left( 3 \rho_0^2 + |\rho|^2 + 2 \sqrt{3} \frac{\eta_0}{|\eta_0|} \mbox{sgn}(\theta^{\overline{\eta} \overline{\rho}}) \rho_0 |\rho| - \left( \overline{\eta}_a^{\overline{\eta}} \right)^2 \right) 
\end{align*}
Without loss of generality, we may thus consider a neighborhood of $|\rho| = 2 |\eta| = 2 |\sigma|$. 

\paragraph{1.1.} If $\xi_t^{\overline{\eta} \overline{\rho}}$ or $\xi_t^{\overline{\sigma} \overline{\rho}}$ is away enough from $0$, we have a good control over terms with a factor $O(\rho_0)$, and then as well $|\rho| - \overline{\eta}_a^{\overline{\eta}}$ and $|\rho| - \overline{\sigma}_a^{\overline{\sigma}}$. If $\epsilon^{\overline{\sigma} \overline{\eta}} = -1$, then 
\begin{align*}
\xi_0 &= \rho_0 + \eta_0 + \sigma_0 \\
&= O(\rho_0) + O\left( \overline{\eta}_a^{\overline{\eta}} - \overline{\sigma}_a^{\overline{\sigma}} \right) + O\left( \overline{\eta}_b^{\overline{\eta}} \right) + O\left( \overline{\sigma}_b^{\overline{\sigma}} \right) 
\end{align*}
which concludes. 

If $\epsilon^{\overline{\sigma} \overline{\eta}} = 1$, then 
\begin{align*}
\varphi_3 &= O(\rho_0) + \xi_0^3 + \xi_0 |\xi|^2 - \eta_0^3 - \eta_0 |\eta|^2 - \sigma_0^3 - \sigma_0 |\sigma|^2 \\
&= O(\rho_0) + O\left( \overline{\eta}_b^{\overline{\eta}} \right) + O\left( \overline{\sigma}_b^{\overline{\sigma}} \right) 
+ \xi_0 \left( \xi_0^2 + |\xi|^2 - 4 \eta_0^2 - 4 \sigma_0^2 + 4 \eta_0 \sigma_0 \right) 
\end{align*}
But $\xi_0 = 2 \eta_0 + o(1)$, $\eta_0 = \sigma_0 + o(1)$, so 
\begin{align*}
\varphi_3 &= o(1) + 2 \eta_0 |\xi|^2 
\end{align*}
However, $\xi = \rho + \sigma + \eta$ and here $|\rho| = 2 |\eta| + o(1)$, $|\eta| = |\sigma| + o(1)$, with $\rho$ and $\eta$, or $\rho$ and $\sigma$ not close to alignement by hypothesis, so that $1 = O(\xi)$, hence $1 = O(\varphi_3)$. 

\paragraph{1.2.} Let us now localise to have $\rho, \eta, \sigma$ close to alignment. 

Here, we have $\xi_0 = o(1) + \eta_0 + \sigma_0$ and $|\eta_0| = |\sigma_0| + o(1)$, so $\xi_0 = o(1)$ or $\xi_0 = 2 \eta_0 + o(1)$. Likewise, $|\xi| = o(1)$ or $|\xi| = 2 |\eta| + o(1)$ or $|\xi| = 4 |\eta| + o(1)$. In fact, 
\begin{align*}
|\xi| = |\eta| \left( 1 + \theta^{\overline{\eta} \overline{\sigma}} + 2 \theta^{\overline{\eta} \overline{\rho}} \right) + o(1) 
\end{align*}
so the only case where $\xi = o(1)$ is $\theta^{\overline{\eta} \overline{\sigma}} = 1 + o(1)$, $\theta^{\overline{\eta} \overline{\rho}} = -1 + o(1)$. 

Note also that 
\begin{align*}
\varphi_3 &= \xi_0^3 + \xi_0 |\xi|^2 - \eta_0^3 - \eta_0 |\eta|^2 - \sigma_0^3 - \sigma_0 |\sigma|^2 - \rho_0^3 - \rho_0 |\rho|^2 \\
&= o(1) + (\eta_0 + \sigma_0)^3 - \eta_0^3 - \sigma_0^3 - \eta_0 |\eta|^2 - \sigma_0 |\sigma|^2 + \xi_0 |\xi|^2 \\
&= o(1) + 3 \eta_0 \sigma_0^2 - \eta_0 |\eta|^2 + 3 \sigma_0 \eta_0^2 - \sigma_0 |\sigma|^2 + \xi_0 |\xi|^2 \\
&= o(1) + \xi_0 |\xi|^2
\end{align*}
Therefore, if $1 = O(\xi)$ and $1 = O(\xi_0)$, we have $1 = O(\varphi_3)$ and we can immediately conclude. 

There are now two subcases to consider: either $\epsilon^{\overline{\eta} \overline{\sigma}} = -1$ (and $\xi_0 = o(1)$), or $\epsilon^{\overline{\eta} \overline{\sigma}} = 1$, $\theta^{\overline{\eta} \overline{\sigma}}$ close to $1$, $\theta^{\overline{\eta} \overline{\rho}}$ close to $-1$ (and $\xi = o(1)$). 

\paragraph{1.2.1.} Consider first the case $\epsilon^{\overline{\eta} \overline{\sigma}} = -1$. 

\paragraph{1.2.1.1.} Assume first that $\theta^{\overline{\eta} \overline{\sigma}}$ is close to $-1$. 

Then using $\widehat{Y}(\overline{\eta}, \overline{\sigma})$ we have a good control over terms with a factor $\overline{\eta}_a^{\overline{\eta}} - \overline{\sigma}_a^{\overline{\sigma}}$. In particular, 
\begin{align*}
\xi_0 &= \eta_0 + \sigma_0 + \rho_0 \\
&= \rho_0 + O\left( \overline{\eta}_a^{\overline{\eta}} - \overline{\sigma}_a^{\overline{\sigma}} \right) + O\left( \overline{\eta}_b^{\overline{\eta}} \right) + O\left( \overline{\sigma}_b^{\overline{\sigma}} \right) 
\end{align*}

We can then compute more precisely: 
\begin{align*}
\varphi_3 &= \xi_0^3 + \xi_0 |\xi|^2 - \eta_0^3 - \eta_0 |\eta|^2 - \sigma_0^3 - \sigma_0 |\sigma|^2 - \rho_0^3 - \rho_0 |\rho|^2 \\
&= O\left( \overline{\eta}_a^{\overline{\eta}} - \overline{\sigma}_a^{\overline{\sigma}} \right) + O\left( \overline{\eta}_b^{\overline{\eta}} \right) + O\left( \overline{\sigma}_b^{\overline{\sigma}} \right)
+ \xi_0 \left( |\xi|^2 - |\rho|^2 \right) 
\end{align*}
Finally, we also have 
\begin{align*}
\xi \cdot \rho &= \eta \cdot \rho + \sigma \cdot \rho + |\rho|^2 \\
&= O(\xi_t^{\overline{\eta} \overline{\rho}}) + O(\xi_t^{\overline{\sigma} \overline{\rho}}) + |\rho| \theta^{\overline{\rho} \overline{\eta}} \left( |\eta| + \theta^{\overline{\eta} \overline{\sigma}} |\sigma| \right) + |\rho|^2 \\
&= O(\xi_t^{\overline{\eta} \overline{\rho}}) + O(\xi_t^{\overline{\sigma} \overline{\rho}}) + O\left( \overline{\eta}_a^{\overline{\eta}} - \overline{\sigma}_a^{\overline{\sigma}} \right) + O\left( \overline{\eta}_b^{\overline{\eta}} \right) + O\left( \overline{\sigma}_b^{\overline{\sigma}} \right) + |\rho|^2 
\end{align*}

But 
\begin{align*}
\xi_0 \widehat{X}_b(\overline{\xi}) \cdot \nabla_{\overline{\xi}} \varphi_3 
&= \frac{|\xi|}{|\overline{\xi}|} \left( 3 \xi_0^2 + |\xi|^2 - 3 \rho_0^2 - |\rho|^2 \right) - \frac{\xi_0^2 \xi}{|\overline{\xi}| |\xi|} \cdot \left( 2 \xi_0 \xi - 2 \rho_0 \rho \right) \\
&= O\left( \overline{\eta}_a^{\overline{\eta}} - \overline{\sigma}_a^{\overline{\sigma}} \right) + O\left( \overline{\eta}_b^{\overline{\eta}} \right) + O\left( \overline{\sigma}_b^{\overline{\sigma}} \right) + O(\varphi_3) - \frac{2 \xi_0^3}{|\overline{\xi}| |\xi|} \left( |\xi|^2 - \xi \cdot \rho \right) \\
&= O\left( \overline{\eta}_a^{\overline{\eta}} - \overline{\sigma}_a^{\overline{\sigma}} \right) + O\left( \overline{\eta}_b^{\overline{\eta}} \right) + O\left( \overline{\sigma}_b^{\overline{\sigma}} \right) + O(\varphi_3) + O(\rho_0 \xi_t^{\overline{\eta} \overline{\rho}}) + O(\rho_0 \xi_t^{\overline{\sigma} \overline{\rho}})
\end{align*}
which concludes. 

\paragraph{1.2.1.2.} Assume now that $\theta^{\overline{\eta} \overline{\sigma}}$ is close to $1$ and that $\theta^{\overline{\eta} \overline{\rho}}$ is close to $-1$. 

Then 
\begin{align*}
|\rho| - \sqrt{3} \frac{\eta_0}{|\eta_0|} \rho_0 &= 2 \sqrt{3} |\eta_0| + O\left( \widehat{X}_a(\overline{\eta}) \cdot \nabla_{\overline{\eta}} \varphi_3 \right) + O\left( \overline{\eta}_b^{\overline{\eta}} \right) \\
|\rho| - \sqrt{3} \frac{\sigma_0}{|\sigma_0|} \rho_0 &= 2 \sqrt{3} |\sigma_0| + O\left( \widehat{X}_a(\overline{\sigma}) \cdot \nabla_{\overline{\sigma}} \varphi_3 \right) + O\left( \overline{\sigma}_b^{\overline{\sigma}} \right)
\end{align*}
so that 
\begin{align*}
|\eta_0| - |\sigma_0| &= - \frac{\eta_0}{|\eta_0|} \rho_0 + O\left( \widehat{X}_a(\overline{\eta}) \cdot \nabla_{\overline{\eta}} \varphi_3 \right) + O\left( \overline{\eta}_b^{\overline{\eta}} \right) + O\left( \widehat{X}_a(\overline{\sigma}) \cdot \nabla_{\overline{\sigma}} \varphi_3 \right) + O\left( \overline{\sigma}_b^{\overline{\sigma}} \right)
\end{align*}
Finally, 
\begin{align*}
\xi_0 &= \frac{\eta_0}{|\eta_0|} \left( |\eta_0| - |\sigma_0| + \frac{\eta_0}{|\eta_0|} \rho_0 \right) \\
&= O\left( \widehat{X}_a(\overline{\eta}) \cdot \nabla_{\overline{\eta}} \varphi_3 \right) + O\left( \overline{\eta}_b^{\overline{\eta}} \right) + O\left( \widehat{X}_a(\overline{\sigma}) \cdot \nabla_{\overline{\sigma}} \varphi_3 \right) + O\left( \overline{\sigma}_b^{\overline{\sigma}} \right)
\end{align*}
which concludes. 

\paragraph{1.2.1.3.} Assume now that $\theta^{\overline{\eta} \overline{\sigma}}$ and $\theta^{\overline{\eta} \overline{\rho}}$ are close to $1$. 

This time, reusing computations from the case 1.2.1.2., we get 
\begin{align*}
|\rho| + \sqrt{3} \frac{\eta_0}{|\eta_0|} \rho_0 &= 2 \sqrt{3} |\eta_0| + O\left( \widehat{X}_a(\overline{\eta}) \cdot \nabla_{\overline{\eta}} \varphi_3 \right) + O\left( \overline{\eta}_b^{\overline{\eta}} \right) \\
|\rho| + \sqrt{3} \frac{\sigma_0}{|\sigma_0|} \rho_0 &= 2 \sqrt{3} |\sigma_0| + O\left( \widehat{X}_a(\overline{\sigma}) \cdot \nabla_{\overline{\sigma}} \varphi_3 \right) + O\left( \overline{\sigma}_b^{\overline{\sigma}} \right) \\
|\eta_0| - |\sigma_0| &= \frac{\eta_0}{|\eta_0|} \rho_0 + O\left( \widehat{X}_a(\overline{\eta}) \cdot \nabla_{\overline{\eta}} \varphi_3 \right) + O\left( \overline{\eta}_b^{\overline{\eta}} \right) + O\left( \widehat{X}_a(\overline{\sigma}) \cdot \nabla_{\overline{\sigma}} \varphi_3 \right) + O\left( \overline{\sigma}_b^{\overline{\sigma}} \right) \\
\xi_0 &= \frac{\eta_0}{|\eta_0|} \left( |\eta_0| - |\sigma_0| + \frac{\eta_0}{|\eta_0|} \rho_0 \right) \\
&= 2 \rho_0 + O\left( \widehat{X}_a(\overline{\eta}) \cdot \nabla_{\overline{\eta}} \varphi_3 \right) + O\left( \overline{\eta}_b^{\overline{\eta}} \right) + O\left( \widehat{X}_a(\overline{\sigma}) \cdot \nabla_{\overline{\sigma}} \varphi_3 \right) + O\left( \overline{\sigma}_b^{\overline{\sigma}} \right)
\end{align*}

We can then compute that {\footnotesize 
\begin{align*}
\varphi_3 &= \xi_0^3 + \xi_0 |\xi|^2 - \eta_0^3 - \eta_0 |\eta|^2 - \sigma_0^3 - \sigma_0 |\sigma|^2 - \rho_0^3 - \rho_0 |\rho|^2 \\
&= \rho_0 \left( 2 |\xi|^2 - |\rho|^2 \right) - 4 \eta_0^3 - 4 \sigma_0^3 
+ o(\rho_0) + O\left( \widehat{X}_a(\overline{\eta}) \cdot \nabla_{\overline{\eta}} \varphi_3 \right) + O\left( \overline{\eta}_b^{\overline{\eta}} \right) + O\left( \widehat{X}_a(\overline{\sigma}) \cdot \nabla_{\overline{\sigma}} \varphi_3 \right) + O\left( \overline{\sigma}_b^{\overline{\sigma}} \right) \\
&= \rho_0 \left( 2 |\xi|^2 - |\rho|^2 - 4 \left( \eta_0^2 + \sigma_0^2 - \eta_0 \sigma_0 \right) \right) 
+ o(\rho_0) + O\left( \widehat{X}_a(\overline{\eta}) \cdot \nabla_{\overline{\eta}} \varphi_3 \right) + O\left( \overline{\eta}_b^{\overline{\eta}} \right) + O\left( \widehat{X}_a(\overline{\sigma}) \cdot \nabla_{\overline{\sigma}} \varphi_3 \right) + O\left( \overline{\sigma}_b^{\overline{\sigma}} \right) \\
&= 72 \rho_0 \eta_0^2
+ o(\rho_0) + O\left( \widehat{X}_a(\overline{\eta}) \cdot \nabla_{\overline{\eta}} \varphi_3 \right) + O\left( \overline{\eta}_b^{\overline{\eta}} \right) + O\left( \widehat{X}_a(\overline{\sigma}) \cdot \nabla_{\overline{\sigma}} \varphi_3 \right) + O\left( \overline{\sigma}_b^{\overline{\sigma}} \right)
\end{align*} }
since here $|\xi| = |\rho| + |\eta| + |\sigma| + o(1) = 4 |\eta| + o(1) = 4 \sqrt{3} |\eta_0| + o(1)$, $|\rho| = 2 \sqrt{3} |\eta_0| + o(1)$, $\sigma_0 = - \eta_0 + o(1)$. Therefore, 
\begin{align*}
\xi_0 &= O(\varphi_3) + O\left( \widehat{X}_a(\overline{\eta}) \cdot \nabla_{\overline{\eta}} \varphi_3 \right) + O\left( \overline{\eta}_b^{\overline{\eta}} \right) + O\left( \widehat{X}_a(\overline{\sigma}) \cdot \nabla_{\overline{\sigma}} \varphi_3 \right) + O\left( \overline{\sigma}_b^{\overline{\sigma}} \right)
\end{align*}
which concludes. 

\paragraph{1.2.2.} Consider now the case $\epsilon^{\overline{\eta} \overline{\sigma}} = 1$, $\theta^{\overline{\eta} \overline{\sigma}}$ close to $1$, $\theta^{\overline{\eta} \overline{\rho}}$ close to $-1$. 

This time, we have
\begin{align*}
|\rho| - \sqrt{3} \frac{\eta_0}{|\eta_0|} \rho_0 &= 2 \sqrt{3} |\eta_0| + O\left( \widehat{X}_a(\overline{\eta}) \cdot \nabla_{\overline{\eta}} \varphi_3 \right) + O\left( \overline{\eta}_b^{\overline{\eta}} \right) \\
|\rho| - \sqrt{3} \frac{\sigma_0}{|\sigma_0|} \rho_0 &= 2 \sqrt{3} |\sigma_0| + O\left( \widehat{X}_a(\overline{\sigma}) \cdot \nabla_{\overline{\sigma}} \varphi_3 \right) + O\left( \overline{\sigma}_b^{\overline{\sigma}} \right) \\
|\eta_0| - |\sigma_0| &= O\left( \widehat{X}_a(\overline{\eta}) \cdot \nabla_{\overline{\eta}} \varphi_3 \right) + O\left( \overline{\eta}_b^{\overline{\eta}} \right) + O\left( \widehat{X}_a(\overline{\sigma}) \cdot \nabla_{\overline{\sigma}} \varphi_3 \right) + O\left( \overline{\sigma}_b^{\overline{\sigma}} \right) \\
\xi_0 &= 2 \eta_0 + \rho_0 + O\left( \widehat{X}_a(\overline{\eta}) \cdot \nabla_{\overline{\eta}} \varphi_3 \right) + O\left( \overline{\eta}_b^{\overline{\eta}} \right) + O\left( \widehat{X}_a(\overline{\sigma}) \cdot \nabla_{\overline{\sigma}} \varphi_3 \right) + O\left( \overline{\sigma}_b^{\overline{\sigma}} \right)
\end{align*}
We may therefore compute that {\footnotesize 
\begin{align*}
\varphi_3 &= \xi_0^3 + \xi_0 |\xi|^2 - \eta_0^3 - \eta_0 |\eta|^2 - \sigma_0^3 - \sigma_0 |\sigma|^2 - \rho_0^3 - \rho_0 |\rho|^2 \\
&= \left( 2 \eta_0 + \rho_0 \right)^3 + \xi_0 |\xi|^2 - 8 \eta_0^3 - \rho_0^3 - \rho_0 \left( 2 \sqrt{3} \eta_0 + \sqrt{3} \rho_0 \right)^2 + O\left( \widehat{X}_a(\overline{\eta}) \cdot \nabla_{\overline{\eta}} \varphi_3 \right) + O\left( \overline{\eta}_b^{\overline{\eta}} \right) \\
&\quad \quad + O\left( \widehat{X}_a(\overline{\sigma}) \cdot \nabla_{\overline{\sigma}} \varphi_3 \right) + O\left( \overline{\sigma}_b^{\overline{\sigma}} \right) \\
&= \xi_0 |\xi|^2 - 6 \eta_0 \rho_0^2 - 3 \rho_0^3 + O\left( \widehat{X}_a(\overline{\eta}) \cdot \nabla_{\overline{\eta}} \varphi_3 \right) + O\left( \overline{\eta}_b^{\overline{\eta}} \right) + O\left( \widehat{X}_a(\overline{\sigma}) \cdot \nabla_{\overline{\sigma}} \varphi_3 \right) + O\left( \overline{\sigma}_b^{\overline{\sigma}} \right) \\
&= \xi_0 \left( |\xi|^2 - 3 \rho_0^2 \right) + O\left( \widehat{X}_a(\overline{\eta}) \cdot \nabla_{\overline{\eta}} \varphi_3 \right) + O\left( \overline{\eta}_b^{\overline{\eta}} \right) + O\left( \widehat{X}_a(\overline{\sigma}) \cdot \nabla_{\overline{\sigma}} \varphi_3 \right) + O\left( \overline{\sigma}_b^{\overline{\sigma}} \right)
\end{align*} }
Since $1 = O(\xi_0)$, we have 
\begin{align*}
|\xi|^2 - 3 \rho_0^2 &= O(\varphi_3) + O\left( \widehat{X}_a(\overline{\eta}) \cdot \nabla_{\overline{\eta}} \varphi_3 \right) + O\left( \overline{\eta}_b^{\overline{\eta}} \right) + O\left( \widehat{X}_a(\overline{\sigma}) \cdot \nabla_{\overline{\sigma}} \varphi_3 \right) + O\left( \overline{\sigma}_b^{\overline{\sigma}} \right)
\end{align*}

We will also need {\footnotesize 
\begin{align*}
\rho_0 \xi \cdot \rho &= \rho_0 |\rho|^2 + \rho_0 \eta \cdot \rho + \rho_0 \sigma \cdot \rho \\
&= \rho_0 \left( 2 \sqrt{3} \eta_0 + \sqrt{3} \rho_0 \right)^2 - \rho_0 |\eta| |\rho| - \rho_0 |\sigma| |\rho| \\
&\quad + O\left( \rho_0 \xi_t^{\overline{\eta} \overline{\rho}} \right) + O\left( \rho_0 \xi_t^{\overline{\sigma} \overline{\rho}} \right) + O\left( \widehat{X}_a(\overline{\eta}) \cdot \nabla_{\overline{\eta}} \varphi_3 \right) + O\left( \overline{\eta}_b^{\overline{\eta}} \right) + O\left( \widehat{X}_a(\overline{\sigma}) \cdot \nabla_{\overline{\sigma}} \varphi_3 \right) + O\left( \overline{\sigma}_b^{\overline{\sigma}} \right) \\
&= 3 \rho_0^2 \left( 2 \eta_0 + \rho_0 \right) + O\left( \rho_0 \xi_t^{\overline{\eta} \overline{\rho}} \right) + O\left( \rho_0 \xi_t^{\overline{\sigma} \overline{\rho}} \right) + O\left( \widehat{X}_a(\overline{\eta}) \cdot \nabla_{\overline{\eta}} \varphi_3 \right) + O\left( \overline{\eta}_b^{\overline{\eta}} \right) + O\left( \widehat{X}_a(\overline{\sigma}) \cdot \nabla_{\overline{\sigma}} \varphi_3 \right) + O\left( \overline{\sigma}_b^{\overline{\sigma}} \right)
\end{align*} }

Finally, we compute that {\footnotesize 
\begin{align*}
|\xi| \widehat{X}_b(\overline{\xi}) \cdot \nabla_{\overline{\xi}} \varphi_3 
&= \frac{|\xi|^2}{|\overline{\xi}|} \left( 3 \xi_0^2 + |\xi|^2 - 3 \rho_0^2 - |\rho|^2 \right) - \frac{\xi_0 \xi}{|\overline{\xi}|} \cdot \left( 2 \xi_0 \xi - 2 \rho_0 \rho \right) \\
&= - 2 \frac{\xi_0}{|\overline{\xi}|} \left( \xi_0 |\xi|^2 - 3 \rho_0^2 (2 \eta_0 + \rho_0) \right) + O\left( \rho_0 \xi_t^{\overline{\eta} \overline{\rho}} \right) + O\left( \rho_0 \xi_t^{\overline{\sigma} \overline{\rho}} \right) + O\left( \widehat{X}_a(\overline{\eta}) \cdot \nabla_{\overline{\eta}} \varphi_3 \right) + O\left( \overline{\eta}_b^{\overline{\eta}} \right) \\
&\quad \quad + O\left( \widehat{X}_a(\overline{\sigma}) \cdot \nabla_{\overline{\sigma}} \varphi_3 \right) + O\left( \overline{\sigma}_b^{\overline{\sigma}} \right) + O(\varphi_3) \\
&= O\left( \rho_0 \xi_t^{\overline{\eta} \overline{\rho}} \right) + O\left( \rho_0 \xi_t^{\overline{\sigma} \overline{\rho}} \right) + O\left( \widehat{X}_a(\overline{\eta}) \cdot \nabla_{\overline{\eta}} \varphi_3 \right) + O\left( \overline{\eta}_b^{\overline{\eta}} \right) + O\left( \widehat{X}_a(\overline{\sigma}) \cdot \nabla_{\overline{\sigma}} \varphi_3 \right) + O\left( \overline{\sigma}_b^{\overline{\sigma}} \right) + O(\varphi_3)
\end{align*} }
and this concludes. 

\paragraph{2.} Let us now localise to have $|\overline{\eta}| \ll |\overline{\sigma}| \simeq |\overline{\rho}|$. 

In this case, we decompose \eqref{equdecchampbCCP-2} as: 
\begin{align*}
\widehat{X}_b(\overline{\xi}) &= \frac{|\xi| |\overline{\sigma}|}{|\overline{\xi}| \sigma_0} \widehat{X}_a(\overline{\sigma}) + O\left( \widehat{X}_a(\overline{\rho}) \right) + O\left( \widehat{X}_c(\overline{\rho}) \right) + O\left( \rho_0 \widehat{X}_b(\overline{\rho}) \right) 
\end{align*}
The contribution of the three last terms is of the form of terms from \eqref{qteslemestimeesgeneriquesdeccubiquehb2}. On the first, we can apply an integration by parts and get terms from \eqref{qteslemestimeesgeneriquesdeccubiquehb2}, plus a term having the form
\begin{align*}
\int_0^t \int \int e^{i s \varphi_3} s^2 m_b(\overline{\xi}) |\overline{\sigma}|^4 \mu(\overline{\xi}, \overline{\eta}, \overline{\sigma}) \widehat{f}(s, \overline{\eta}) \widehat{f}(s, \overline{\sigma}) \widehat{f}(s, \overline{\rho}) ~ d\overline{\eta} d\overline{\sigma} ds 
\end{align*}
\eqref{equdecchampbCCP-1} also has this form. 

Note that $\xi_0 = \sigma_0 + o(1)$ so $|\overline{\xi}| \simeq |\overline{\sigma}| \simeq |\xi_0|$. 

But we have that
\begin{align*}
\varphi_3 &= \xi_0^3 + \xi_0 |\xi|^2 - \eta_0^3 - \eta_0 |\eta|^2 - \sigma_0^3 - \sigma_0 |\sigma|^2 - \rho_0^3 - \rho_0 |\rho|^2 \\
&= O(\rho_0) + O(\overline{\eta}) + \xi_0 \left( |\xi|^2 - |\sigma|^2 \right) 
\end{align*}
and therefore
\begin{align*}
m_b(\overline{\xi}) &= O\left( 3 \xi_0^2 - |\xi|^2 \right) \\
&= O(\rho_0) + O\left( \overline{\eta} \right) + O(\varphi_3) + O\left( \overline{\sigma}_b^{\overline{\sigma}} \right) 
\end{align*}
Besides, $\nabla_{\sigma} \varphi_3 = o(1) - 2 \sigma_0 \sigma$ so $1 = O(\nabla_{\sigma} \varphi_3)$, and since $\overline{\eta}$ is small with respect to the other variables, we easily get $1 = O\left( \widehat{X}_a(\overline{\eta}) \cdot \nabla_{\overline{\eta}} \varphi_3 \right)$. If we have a factor $O(\rho_0)$, we can apply an integration by parts along $\widehat{X}_a(\overline{\eta})$, as well as on the term with factor $O(\overline{\eta})$, and along $\nabla_{\sigma}$ if we have a factor $O\left( \overline{\sigma}_b^{\overline{\sigma}} \right)$, while we can apply an integration by parts in time if we have a factor $O(\varphi_3)$. In any case, we get terms from \eqref{qteslemestimeesgeneriquesdeccubiquehb2}. 

\paragraph{3.} Let us now localise to have $|\overline{\rho}| \ll |\overline{\eta}| \simeq |\overline{\sigma}|$. 

Then \eqref{equdecchampbCCP-2} is already of the form \eqref{lemestimeesgeneriquesdeccubiquehb2-eta0sigmaresxrho}. 

For \eqref{equdecchampbCCP-1}, we have $1 = O\left( \widehat{X}_a(\overline{\eta}) \cdot \nabla_{\overline{\eta}} \varphi_3 \right)$ and it is easy to apply an integration by parts to get terms from \eqref{qteslemestimeesgeneriquesdeccubiquehb2}. 

\paragraph{4.} Let us now localise to have $|\overline{\eta}| + |\overline{\rho}| \ll |\overline{\sigma}|$. 

As in the case 2. above, we can rewrite \eqref{equdecchampbCCP-2} as terms from \eqref{qteslemestimeesgeneriquesdeccubiquehb2} plus 
\begin{align*}
\int_0^t \int \int e^{i s \varphi_3} s^2 m_b(\overline{\xi}) |\overline{\sigma}|^4 \mu \widehat{f}(s, \overline{\eta}) \widehat{f}(s, \overline{\sigma}) \widehat{f}(s, \overline{\rho}) ~ d\overline{\eta} d\overline{\sigma} ds 
\end{align*}
where $\mu$ allows to distribute a derivative. \eqref{equdecchampbCCP-1} can also be written under this form. $\mu$ then allows to apply the same computations as in the case 2. above, we skip the details. 

\paragraph{5.} Let us finally localise to have $|\overline{\eta}| + |\overline{\sigma}| \ll |\overline{\rho}|$. 

Again, by integration by parts we can rewrite \eqref{equdecchampbCCP-2} as terms of \eqref{qteslemestimeesgeneriquesdeccubiquehb2} (using that $\mu$ allows to distribute a derivative) plus a term having the form 
\begin{align*}
\int_0^t \int \int e^{i s \varphi_3} s^2 \xi_0 |\overline{\rho}|^2 |\overline{\eta}| \mu m_{\widehat{\mathcal{C}}}(\overline{\eta}) m_{\widehat{\mathcal{C}}}(\overline{\sigma}) m_{\widehat{\mathcal{P}}}(\overline{\rho}) \widehat{f}(s, \overline{\eta}) \widehat{f}(s, \overline{\sigma}) \widehat{f}(s, \overline{\rho}) ~ d\overline{\eta} d\overline{\sigma} ds 
\end{align*}
up to symmetrising $\overline{\eta}, \overline{\sigma}$. \eqref{equdecchampbCCP-1} also has this form. Here above, $\mu$ does not allow to distribute any more derivative (as we already used it to get the above form). 

We can now compute that {\footnotesize 
\begin{align*}
\varphi_3 &= \xi_0^3 + \xi_0 |\xi|^2 - \eta_0^3 - \eta_0 |\eta|^2 - \sigma_0^3 - \sigma_0 |\sigma|^2 - \rho_0^3 - \rho_0 |\rho|^2 \\
&= \left( \eta_0 + \sigma_0 + \rho_0 \right)^3 - \eta_0^3 - \sigma_0^3 - \rho_0^3 
+ \xi_0 \left( |\eta|^2 + |\sigma|^2 + |\rho|^2 + 2 \eta \cdot \sigma + 2 \eta \cdot \rho + 2 \sigma \cdot \rho \right) 
- \eta_0 |\eta|^2 - \sigma_0 |\sigma|^2 - \rho_0 |\rho|^2 \\
&= O(\rho_0 \overline{\eta}) + O(\rho_0 \overline{\sigma}) + O(\overline{\eta} \overline{\sigma}) 
+ (\xi_0 - \eta_0) |\eta|^2 + (\xi_0 - \sigma_0) |\sigma|^2 + (\xi_0 - \rho_0) |\rho|^2 + 2 \xi_0 \eta \cdot \rho + 2 \xi_0 \sigma \cdot \rho \\
&= O(\rho_0 \overline{\eta}) + O(\rho_0 \overline{\sigma}) + O(\overline{\eta} \overline{\sigma}) 
+ (\xi_0 - \rho_0) |\rho|^2 + 2 (\xi_0 - \rho_0) (\eta + \sigma) \cdot \rho \\
&= O(\rho_0 \overline{\eta}) + O(\rho_0 \overline{\sigma}) + O(\overline{\eta} \overline{\sigma}) 
+ (\xi_0 - \rho_0) \left( |\rho|^2 + o(1) \right) 
\end{align*} }
Therefore, 
\begin{align*}
\eta_0 + \sigma_0 = O(\varphi_3) + O(\rho_0 \overline{\eta}) + O(\rho_0 \overline{\sigma}) + O(\overline{\eta} \overline{\sigma})
\end{align*}
We can then develop
\begin{align*}
\xi_0 |\overline{\eta}| &= O((\eta_0 + \sigma_0)|\overline{\eta}|) + \rho_0 |\overline{\eta}| \\
&= O(|\overline{\eta}| \varphi_3) + O(\rho_0 \overline{\eta}) + O(\rho_0 \overline{\sigma}) + O(\overline{\eta} \overline{\sigma})
\end{align*}
Now we can apply an integration by parts in time on the term having a factor $O(|\overline{\eta}|\varphi_3)$, along $\widehat{X}_a(\overline{\eta}) \cdot \nabla_{\overline{\eta}}$ for the factor $O(\rho_0 \overline{\eta})$, along $\widehat{X}_a(\overline{\sigma}) \cdot \nabla_{\overline{\sigma}}$ for the factor $O(\rho_0 \overline{\sigma})$ and the factor $O(\overline{\eta} \overline{\sigma})$. In all cases, we get terms from \eqref{qteslemestimeesgeneriquesdeccubiquehb2}. 

\subsubsection{Interaction \texorpdfstring{$\widehat{\mathcal{L}}\widehat{\mathcal{L}}\widehat{\mathcal{C}}$}{LLC}}

Let us consider: {\footnotesize 
\begin{subequations}
\begin{align}
&\xi_0 m_b(\overline{\xi}) \widehat{X}_b(\overline{\xi}) \cdot \nabla_{\overline{\xi}} \widehat{I}_{s\mu}^{\widehat{\mathcal{L}}\widehat{\mathcal{L}}\widehat{\mathcal{C}}}[F_1, F_2, F_3](t, \overline{\xi}) \notag \\
&\quad = \int_0^t \int \int i s^2 \xi_0 m_b(\overline{\xi}) \widehat{X}_b(\overline{\xi}) \cdot \nabla_{\overline{\xi}} \varphi_3 e^{i s \varphi_3} \mu(\overline{\xi}, \overline{\eta}, \overline{\sigma}) m_{\widehat{\mathcal{L}}}(\overline{\eta}) m_{\widehat{\mathcal{L}}}(\overline{\sigma}) m_{\widehat{\mathcal{C}}}(\overline{\rho}) \widehat{F}_1(s, \overline{\eta}) \widehat{F}_2(s, \overline{\sigma}) \widehat{F}_3(s, \overline{\rho}) ~ d\overline{\eta} d\overline{\sigma} ds \label{equdecchampbLLC-1} \\
&\quad \quad + \int_0^t \int \int e^{i s \varphi_3} s \xi_0 \mu(\overline{\xi}, \overline{\eta}, \overline{\sigma}) m_{\widehat{\mathcal{L}}}(\overline{\eta}) m_{\widehat{\mathcal{L}}}(\overline{\sigma}) m_{\widehat{\mathcal{C}}}(\overline{\rho}) \widehat{F}_1(s, \overline{\eta}) \widehat{F}_2(s, \overline{\sigma}) m_b(\overline{\xi}) \widehat{X}_b(\overline{\xi}) \cdot \nabla_{\overline{\xi}} \widehat{F}_3(s, \overline{\rho}) ~ d\overline{\eta} d\overline{\sigma} ds \label{equdecchampbLLC-2} \\
&\quad \quad + \int_0^t \int \int e^{i s \varphi_3} s \xi_0 m_b(\overline{\xi}) \widehat{X}_b(\overline{\xi}) \cdot \nabla_{\overline{\xi}} \left( \mu(\overline{\xi}, \overline{\eta}, \overline{\sigma}) m_{\widehat{\mathcal{L}}}(\overline{\eta}) m_{\widehat{\mathcal{L}}}(\overline{\sigma}) m_{\widehat{\mathcal{C}}}(\overline{\rho}) \right) \widehat{F}_1(s, \overline{\eta}) \widehat{F}_2(s, \overline{\sigma}) \widehat{F}_3(s, \overline{\rho}) ~ d\overline{\eta} d\overline{\sigma} ds \label{equdecchampbLLC-3} 
\end{align}
\end{subequations} }
\eqref{equdecchampbLLC-3} is of the form $\eqref{lemestimeesgeneriquesdeccubiquehb2-dersymb}+\eqref{lemestimeesgeneriquesdeccubiquehb2-dersymbbis}$. 

On \eqref{equdecchampbLLC-2}, we can project $\widehat{X}_b(\overline{\xi})$ on the basis $\left( \widehat{X}_a(\overline{\rho}), \widehat{X}_c(\overline{\rho}), \widehat{X}_{b-\widehat{\mathcal{L}}}(\overline{\eta}, \overline{\rho}) \right)$ to get terms from \eqref{qteslemestimeesgeneriquesdeccubiquehb2} and one with frequency derivative along $\widehat{X}_{b-\widehat{\mathcal{L}}}$, for which we apply an integration by parts in $\overline{\eta}$ to get, up to terms from \eqref{qteslemestimeesgeneriquesdeccubiquehb2}, {\footnotesize 
\begin{align*}
\int_0^t \int \int i s^2 m_b(\overline{\xi}) \xi_0 \frac{\rho_0}{|\rho_0|} P_b^b(\overline{\xi}, \overline{\rho}) \widehat{X}_{b-\widehat{\mathcal{L}}}(\overline{\eta}, \overline{\rho}) \cdot \nabla_{\overline{\eta}} \varphi_3 ~ e^{i s \varphi_3} \mu(\overline{\xi}, \overline{\eta}, \overline{\sigma}) m_{\widehat{\mathcal{L}}}(\overline{\eta}) m_{\widehat{\mathcal{L}}}(\overline{\sigma}) m_{\widehat{\mathcal{C}}}(\overline{\rho}) \widehat{F}_1(s, \overline{\eta}) \widehat{F}_2(s, \overline{\sigma})  \widehat{F}_3(s, \overline{\rho}) ~ d\overline{\eta} d\overline{\sigma} ds
\end{align*} }
that we can group with \eqref{equdecchampbLLC-1} to reduce to 
\begin{align*}
&\int_0^t \int \int i s^2 \xi_0 m_b(\overline{\xi}) \left( \widehat{X}_b(\overline{\xi}) \cdot \nabla_{\overline{\xi}} + \frac{\rho_0}{|\rho_0|} P_b^b(\overline{\xi}, \overline{\rho}) \widehat{X}_{b-\widehat{\mathcal{L}}}(\overline{\eta}, \overline{\rho}) \cdot \nabla_{\overline{\eta}} \right) \varphi_3 \\
&\quad \quad \quad e^{i s \varphi_3} \mu(\overline{\xi}, \overline{\eta}, \overline{\sigma}) m_{\widehat{\mathcal{L}}}(\overline{\eta}) m_{\widehat{\mathcal{L}}}(\overline{\sigma}) m_{\widehat{\mathcal{C}}}(\overline{\rho}) \widehat{F}_1(s, \overline{\eta}) \widehat{F}_2(s, \overline{\sigma}) \widehat{F}_3(s, \overline{\rho}) ~ d\overline{\eta} d\overline{\sigma} ds
\end{align*}

We now separate into subcases depending on the relative size of frequencies. 

Note that we always have 
\begin{align*}
\left( \partial_{\eta_0} - \partial_{\sigma_0} \right) \varphi_3 &= 3 \sigma_0^2 + |\sigma|^2 - 3 \eta_0^2 - |\eta|^2
\end{align*}
so we only need to consider the cases $|\overline{\eta}| \simeq |\overline{\sigma}| \gtrsim |\overline{\rho}|$ and $|\overline{\eta}|+|\overline{\sigma}| \ll |\overline{\rho}|$. 

\paragraph{1.} Let us first localise to have $|\overline{\eta}| \simeq |\overline{\sigma}| \simeq |\overline{\rho}|$. 

Then since
\begin{align*}
\nabla_{\eta} \varphi_3 &= o(1) + 2 \rho_0 \rho \simeq |\overline{\rho}|^2 
\end{align*}
we deduce that we have a good control if we have a factor $O(\overline{\rho}_b^{\overline{\rho}})$. Then, computing $\widehat{X}_c(\overline{\rho}) \cdot \nabla_{\overline{\eta}} \varphi_3$, we can proceed the same way with $|\eta| \xi_t^{\overline{\eta} \overline{\rho}}$, and then 
\begin{align*}
\widehat{X}_a(\overline{\rho}) \cdot \nabla_{\overline{\eta}} \varphi_3 &= \frac{\rho_0}{|\overline{\rho}|} \left( \left( \overline{\rho}_a^{\overline{\rho}} \right)^2 - \left( \overline{\eta}_a^{\overline{\rho}} \right)^2 \right) + O\left( \overline{\rho}_b^{\overline{\rho}} \right) + O\left( |\eta| \xi_t^{\overline{\eta} \overline{\rho}} \right) 
\end{align*}
so we also have a good control if we have a factor $O\left( \overline{\rho}_a^{\overline{\rho}} - \epsilon^{\overline{\eta} \overline{\rho}} \overline{\eta}_a^{\overline{\rho}} \right)$, or likewise $O\left( |\sigma| \xi_t^{\overline{\sigma} \overline{\rho}} \right)$ or $O\left( \overline{\rho}_a^{\overline{\rho}} - \epsilon^{\overline{\sigma} \overline{\rho}} \overline{\sigma}_a^{\overline{\rho}} \right)$. 

Finally, 
\begin{align*}
\left( \partial_{\eta_0} - \partial_{\sigma_0} \right) \varphi_3 &= 3 \sigma_0^2 + |\sigma|^2 - 3 \eta_0^2 - |\eta|^2 \\
&= \frac{1}{2} \left( \left( \overline{\sigma}_a^{\overline{\rho}} \right)^2 + \left( \overline{\sigma}_b^{\overline{\rho}} \right)^2 - \left( \overline{\eta}_a^{\overline{\rho}} \right)^2 - \left( \overline{\eta}_b^{\overline{\rho}} \right)^2 \right) + O\left( |\eta| \xi_t^{\overline{\eta} \overline{\rho}} \right) + O\left( |\sigma| \xi_t^{\overline{\sigma} \overline{\rho}} \right) \\
&= \frac{1}{2} \left( \left( \overline{\sigma}_b^{\overline{\rho}} \right)^2 - \left( \overline{\eta}_b^{\overline{\rho}} \right)^2 \right) + O\left( \overline{\rho}_a^{\overline{\rho}} - \epsilon^{\overline{\eta} \overline{\rho}} \overline{\eta}_a^{\overline{\rho}} \right) + O\left( \overline{\rho}_a^{\overline{\rho}} - \epsilon^{\overline{\sigma} \overline{\rho}} \overline{\sigma}_a^{\overline{\rho}} \right) + O\left( |\eta| \xi_t^{\overline{\eta} \overline{\rho}} \right) + O\left( |\sigma| \xi_t^{\overline{\sigma} \overline{\rho}} \right)
\end{align*}
so we have a control over terms with a factor $O\left( \overline{\eta}_b^{\overline{\rho}} - \epsilon^{\overline{\eta} \overline{\sigma}} \overline{\sigma}_b^{\overline{\rho}} \right)$. 

We will also need that 
\begin{align*}
|\xi|^2 &= \left( \frac{\rho \cdot \xi}{|\rho|} \right)^2 + \left( \frac{J \rho \cdot \xi}{|\rho|} \right)^2 \\
&= \left( \frac{\rho \cdot \xi}{|\rho|} \right)^2 + O\left( |\eta| \xi_t^{\overline{\eta} \overline{\rho}} \right) + O\left( |\sigma| \xi_t^{\overline{\sigma} \overline{\rho}} \right)
\end{align*}

Now, we can compute that {\footnotesize 
\begin{align*}
6 \sqrt{3} \frac{\rho_0}{|\rho_0|} \varphi_3 &= \left( \overline{\xi}_a^{\overline{\rho}} \right)^3 + \left( \overline{\xi}_b^{\overline{\rho}} \right)^3 - \left( \overline{\eta}_a^{\overline{\rho}} \right)^3 - \left( \overline{\eta}_b^{\overline{\rho}} \right)^3 - \left( \overline{\rho}_a^{\overline{\rho}} \right)^3 - \left( \overline{\rho}_b^{\overline{\rho}} \right)^3 - \left( \overline{\sigma}_a^{\overline{\rho}} \right)^3 - \left( \overline{\sigma}_b^{\overline{\rho}} \right)^3 + O\left( |\eta| \xi_t^{\overline{\eta} \overline{\rho}} \right) + O\left( |\sigma| \xi_t^{\overline{\sigma} \overline{\rho}} \right) \\
&= \left( \overline{\rho}_a^{\overline{\rho}} \right)^3 \left( \left( 1 + \epsilon^{\overline{\rho} \overline{\eta}} + \epsilon^{\overline{\rho} \overline{\sigma}} \right)^3 - \epsilon^{\overline{\rho} \overline{\eta}} - \epsilon^{\overline{\rho} \overline{\sigma}} - 1 \right) 
+ \left( \overline{\eta}_b^{\overline{\rho}} \right)^3 \left( \left( 1 + \epsilon^{\overline{\eta} \overline{\sigma}} \right)^3 - 1 - \epsilon^{\overline{\eta} \overline{\sigma}} \right) \\
&\quad + O\left( |\eta| \xi_t^{\overline{\eta} \overline{\rho}} \right) + O\left( |\sigma| \xi_t^{\overline{\sigma} \overline{\rho}} \right) + O\left( \overline{\eta}_b^{\overline{\rho}} - \epsilon^{\overline{\eta} \overline{\sigma}} \overline{\sigma}_b^{\overline{\rho}} \right) + O\left( \overline{\rho}_a^{\overline{\rho}} - \epsilon^{\overline{\eta} \overline{\rho}} \overline{\eta}_a^{\overline{\rho}} \right) + O\left( \overline{\rho}_a^{\overline{\rho}} - \epsilon^{\overline{\sigma} \overline{\rho}} \overline{\sigma}_a^{\overline{\rho}} \right) + O\left( \overline{\rho}_b^{\overline{\rho}} \right) \\
&= 6 \left( \overline{\rho}_a^{\overline{\rho}} \right)^3 \left( 
1 + \epsilon^{\overline{\eta} \overline{\rho}} + \epsilon^{\overline{\sigma} \overline{\rho}} + \epsilon^{\overline{\eta} \overline{\sigma}} \right) 
+ 3 \left( \overline{\eta}_b^{\overline{\rho}} \right)^3 \left( 
1 + \epsilon^{\overline{\eta} \overline{\sigma}} \right) \\
&\quad + O\left( |\eta| \xi_t^{\overline{\eta} \overline{\rho}} \right) + O\left( |\sigma| \xi_t^{\overline{\sigma} \overline{\rho}} \right) + O\left( \overline{\eta}_b^{\overline{\rho}} - \epsilon^{\overline{\eta} \overline{\sigma}} \overline{\sigma}_b^{\overline{\rho}} \right) + O\left( \overline{\rho}_a^{\overline{\rho}} - \epsilon^{\overline{\eta} \overline{\rho}} \overline{\eta}_a^{\overline{\rho}} \right) + O\left( \overline{\rho}_a^{\overline{\rho}} - \epsilon^{\overline{\sigma} \overline{\rho}} \overline{\sigma}_a^{\overline{\rho}} \right) + O\left( \overline{\rho}_b^{\overline{\rho}} \right)
\end{align*} }

We now separate into subcases depending on the signs. 

\paragraph{1.1.} If $\epsilon^{\overline{\eta} \overline{\rho}} = \epsilon^{\overline{\sigma} \overline{\rho}} = 1 = \epsilon^{\overline{\eta} \overline{\sigma}}$, then 
\begin{align*}
\varphi_3 &= 24 \left( \overline{\rho}_a^{\overline{\rho}} \right)^3 
+ 6 \left( \overline{\eta}_b^{\overline{\rho}} \right)^3 + o(1) 
\end{align*}
But $\overline{\eta}_b^{\overline{\rho}} = \overline{\eta}_a^{\overline{\rho}} + o(1) = \overline{\rho}_a^{\overline{\rho}} + o(1)$ so $1 = O(\varphi_3)$. 

\paragraph{1.2.} If $\epsilon^{\overline{\eta} \overline{\rho}} = 1$, $\epsilon^{\overline{\sigma} \overline{\rho}} = -1 = \epsilon^{\overline{\eta} \overline{\sigma}}$, then we have $\overline{\rho}_a^{\overline{\rho}} = \overline{\rho}_a^{\overline{\rho}} + o(1)$, so $|\overline{\xi}| \simeq |\overline{\rho}|$, and therefore 
\begin{align*}
3 \xi_0^2 - |\xi|^2 &= \overline{\xi}_a^{\overline{\rho}} \overline{\xi}_b^{\overline{\rho}} + O\left( |\eta| \xi_t^{\overline{\eta} \overline{\rho}} \right) + O\left( |\sigma| \xi_t^{\overline{\sigma} \overline{\rho}} \right) \\
&= \overline{\xi}_a^{\overline{\rho}} \left( \overline{\rho}_b^{\overline{\rho}} + \overline{\eta}_b^{\overline{\rho}} + \overline{\sigma}_b^{\overline{\rho}} \right) + O\left( |\eta| \xi_t^{\overline{\eta} \overline{\rho}} \right) + O\left( |\sigma| \xi_t^{\overline{\sigma} \overline{\rho}} \right) \\
&= O\left( \overline{\rho}_b^{\overline{\rho}} \right) + O\left( \overline{\eta}_b^{\overline{\rho}} + \overline{\sigma}_b^{\overline{\rho}} \right) + O\left( |\eta| \xi_t^{\overline{\eta} \overline{\rho}} \right) + O\left( |\sigma| \xi_t^{\overline{\sigma} \overline{\rho}} \right)
\end{align*}
which concludes. 

\paragraph{1.3.} If $\epsilon^{\overline{\sigma} \overline{\rho}} = 1$, $\epsilon^{\overline{\eta} \overline{\rho}} = -1 = \epsilon^{\overline{\eta} \overline{\sigma}}$, then up to exchanging $\overline{\eta}$ and $\overline{\sigma}$ we recover case 1.2. above. 

\paragraph{1.4.} If $\epsilon^{\overline{\eta} \overline{\rho}} = \epsilon^{\overline{\sigma} \overline{\rho}} = -1$, $\epsilon^{\overline{\eta} \overline{\sigma}} = 1$, then 
\begin{align*}
\varphi_3 &= 6 \left( \overline{\eta}_b^{\overline{\rho}} \right)^3 + o(1) 
\end{align*}
so that $1 = O(\varphi_3)$. 

\paragraph{2.} Let us now localise to have $|\overline{\rho}| \ll |\overline{\eta}| \simeq |\overline{\sigma}|$. 

We then have that $\xi_0 m_b(\overline{\xi}) = O(\xi) = O(\eta) + O(\sigma) + O(\rho)$. Moreover, $1 = O\left( \widehat{X}_a(\overline{\rho}) \cdot \nabla_{\overline{\eta}} \varphi_3 \right)$ and $1 = O\left( \widehat{X}_a(\overline{\rho}) \cdot \nabla_{\overline{\sigma}} \varphi_3 \right)$, so we can apply an integration by parts along $\widehat{X}_a(\overline{\rho}) \cdot \nabla_{\overline{\eta}}$ if we have a factor $O(\eta)$ or $O(\rho)$, and along $\widehat{X}_a(\overline{\rho}) \cdot \nabla_{\overline{\sigma}}$ if we have a factor $O(\sigma)$, so we only recover terms from \eqref{qteslemestimeesgeneriquesdeccubiquehb2}. 

\paragraph{3.} Let us finally localise to have $|\overline{\eta}| + |\overline{\sigma}| \ll |\overline{\rho}|$. 

We can then distribute a derivative using $\mu$, and by symmetry me may assume we obtain an additional factor $|\overline{\eta}|$. Since $1 = O\left( \widehat{X}_a(\overline{\rho}) \cdot \nabla_{\overline{\sigma}} \varphi_3 \right)$ again, we can apply an integration by parts in this direction and get terms from \eqref{qteslemestimeesgeneriquesdeccubiquehb2}. 

\subsubsection{Interaction \texorpdfstring{$\widehat{\mathcal{L}}\widehat{\mathcal{C}}\widehat{\mathcal{P}}$}{LCP}}

Let us now consider: {\footnotesize 
\begin{subequations}
\begin{align}
&\xi_0 m_b(\overline{\xi}) \widehat{X}_b(\overline{\xi}) \cdot \nabla_{\overline{\xi}} \widehat{I}_{s \mu}^{\widehat{\mathcal{L}}\widehat{\mathcal{C}}\widehat{\mathcal{P}}}[F_1, F_2, F_3](t, \overline{\xi}) \notag \\
&\quad = \int_0^t \int \int i s^2 \xi_0 m_b(\overline{\xi}) \widehat{X}_b(\overline{\xi}) \cdot \nabla_{\overline{\xi}} \varphi_3 e^{i s \varphi_3} \mu(\overline{\xi}, \overline{\eta}, \overline{\sigma}) m_{\widehat{\mathcal{L}}}(\overline{\eta}) m_{\widehat{\mathcal{C}}}(\overline{\sigma}) m_{\widehat{\mathcal{P}}}(\overline{\rho}) \widehat{F}_1(s, \overline{\eta}) \widehat{F}_2(s, \overline{\sigma}) \widehat{F}_3(s, \overline{\rho}) ~ d\overline{\eta} d\overline{\sigma} ds \label{equdecchampbLCP-1} \\
&\quad \quad + \int_0^t \int \int e^{i s \varphi_3} s \xi_0 \mu(\overline{\xi}, \overline{\eta}, \overline{\sigma}) m_{\widehat{\mathcal{L}}}(\overline{\eta}) m_{\widehat{\mathcal{C}}}(\overline{\sigma}) m_{\widehat{\mathcal{P}}}(\overline{\rho}) \widehat{F}_1(s, \overline{\eta}) \widehat{F}_2(s, \overline{\sigma}) m_b(\overline{\xi}) \widehat{X}_b(\overline{\xi}) \cdot \nabla_{\overline{\xi}} \widehat{F}_3(s, \overline{\rho}) ~ d\overline{\eta} d\overline{\sigma} ds \label{equdecchampbLCP-2} \\
&\quad \quad + \int_0^t \int \int e^{i s \varphi_3} s \xi_0 m_b(\overline{\xi}) \widehat{X}_b(\overline{\xi}) \cdot \nabla_{\overline{\xi}} \left( \mu(\overline{\xi}, \overline{\eta}, \overline{\sigma}) m_{\widehat{\mathcal{L}}}(\overline{\eta}) m_{\widehat{\mathcal{C}}}(\overline{\sigma}) m_{\widehat{\mathcal{P}}}(\overline{\rho}) \right) \widehat{F}_1(s, \overline{\eta}) \widehat{F}_2(s, \overline{\sigma}) \widehat{F}_3(s, \overline{\rho}) ~ d\overline{\eta} d\overline{\sigma} ds \label{equdecchampbLCP-3} 
\end{align}
\end{subequations} }
\eqref{equdecchampbLCP-3} is of the form $\eqref{lemestimeesgeneriquesdeccubiquehb2-dersymb}+\eqref{lemestimeesgeneriquesdeccubiquehb2-dersymbbis}$. 

We now separate into subcases depending on the relative sizes of the frequencies. 

\paragraph{1.} Let us first localise to have $|\overline{\eta}| \simeq |\overline{\sigma}| \simeq |\overline{\rho}|$. 

Then \eqref{equdecchampbLCP-2} is of the form \eqref{lemestimeesgeneriquesdeccubiquehb2-eta0sigmaresxrho}. 

We now compute that 
\begin{align*}
\partial_{\eta_0} \varphi_3 &= 3 \eta_0^2 + |\eta|^2 - 3 \rho_0^2 - |\rho|^2 \\
\widehat{X}_a(\overline{\sigma}) \cdot \nabla_{\overline{\sigma}} \varphi_3 &= \frac{\sigma_0}{|\overline{\sigma}|} \left( 3 \rho_0^2 + |\rho|^2 - 3 |\overline{\sigma}|^2 \right) + 2 \rho_0 \frac{\sigma \cdot \rho}{|\overline{\sigma}|} 
\end{align*}
Therefore, we only need to consider the neighborhood of $3 \eta_0^2 = |\rho|^2 = 12 \sigma_0^2$. 

If $\epsilon^{\overline{\eta} \overline{\sigma}} = 1$, $\eta_0 = 2 \sigma_0 + o(1)$ and therefore $\xi_0 = 3 \sigma_0 + o(1)$, then 
\begin{align*}
\varphi_3 &= \xi_0^3 + \xi_0 |\xi|^2 - \eta_0^3 - 4 \sigma_0^3 + o(1) \\
&= 27 \sigma_0^3 + 3 \sigma_0 |\xi|^2 - 8 \sigma_0^3 - 4 \sigma_0^3 + o(1) \\
&= 3 \sigma_0 \left( 5 \sigma_0^2 + |\xi|^2 \right) + o(1) 
\end{align*}
so $1 = O(\varphi_3)$. 

Assume therefore $\epsilon^{\overline{\eta} \overline{\sigma}} = -1$, so $\eta_0 = - 2 \sigma_0 + o(1)$. This time, $\xi_0 = - \sigma_0 + o(1)$, and then 
\begin{align*}
\varphi_3 &= - \sigma_0^3 - \sigma_0 |\xi|^2 + 8 \sigma_0^3 - 4 \sigma_0^3 + o(1) \\
&= \sigma_0 \left( |\sigma|^2 - |\xi|^2 \right) + o(1) 
\end{align*}
We may therefore also restrict to a neighborhood of $|\sigma| = |\xi|$. But $\xi = o(1) + \sigma + \rho$ and $|\rho| = 2 |\sigma| + o(1) = |\sigma| + |\xi| + o(1)$, so this forces $\theta^{\overline{\sigma} \overline{\rho}}$ close to $-1$, $\xi = - \sigma + o(1)$. 

In this situation, we note that, if we apply an integration by parts in the direction $m_{\alpha}(\overline{\sigma}) \widehat{X}_{\alpha}(\overline{\sigma}) \cdot \left( \nabla_{\overline{\eta}} - \nabla_{\overline{\sigma}} \right)$, we obtain terms from \eqref{qteslemestimeesgeneriquesdeccubiquehb2}, and a term of the form 
\begin{align}
\int_0^t \int \int e^{i s \varphi_3} s \mu m_{\widehat{\mathcal{L}}}(\overline{\eta}) m_{\widehat{\mathcal{C}}}(\overline{\sigma}) m_{\widehat{\mathcal{P}}}(\overline{\rho}) \nabla_{\overline{\eta}} \widehat{f}(s, \overline{\eta}) |\overline{\sigma}|^{\frac{3}{2}} \widehat{f}(s, \overline{\sigma}) |\overline{\rho}|^{\frac{1}{2}} \widehat{f}(s, \overline{\rho}) ~ d\overline{\eta} d\overline{\sigma} ds \label{termedegenereokcubicbbLCPBBB} 
\end{align}
that satisfies: 
\begin{align*}
\Vert \partial_t \eqref{termedegenereokcubicbbLCPBBB} \Vert_{L^2} &\lesssim t \Vert m_{\widehat{\mathcal{L}}}(D) (x, y) f(t) \Vert_{L^2} \Vert m_{\widehat{\mathcal{C}}}(D) |\nabla|^{\frac{3}{2}} u(t) \Vert_{L^{\infty}} \Vert m_{\widehat{\mathcal{P}}}(D) |\nabla|^{\frac{1}{2}} u(t) \Vert_{L^{\infty}} \\
&\lesssim t^{-\frac{2}{3}} \langle t \rangle^{\frac{5}{48}-\frac{1}{3}-\frac{1}{6}+402\delta} \Vert u \Vert_X^3 
\end{align*}
and so such a term also satisfies a sufficient estimate for the proof of the Lemma. We therefore have a good control over terms with a factor
\begin{align*}
\widehat{X}_c(\overline{\sigma}) \cdot \left( \nabla_{\overline{\eta}} - \nabla_{\overline{\sigma}} \right) \varphi_3 &= - 2 \eta_0 \frac{J \sigma \cdot \eta}{|\sigma|} \simeq |\overline{\eta}| |\eta| \xi_t^{\overline{\eta} \overline{\sigma}} \\
\widehat{X}_a(\overline{\sigma}) \cdot \left( \nabla_{\overline{\eta}} - \nabla_{\overline{\sigma}} \right) \varphi_3 &= \frac{\sigma_0}{|\overline{\sigma}|} \left( \left( \overline{\sigma}_a^{\overline{\sigma}} \right)^2 - \left( \overline{\eta}_a^{\overline{\sigma}} \right)^2 \right) + O\left( |\eta| \xi_t^{\overline{\eta} \overline{\sigma}} \right) + O\left( \overline{\sigma}_b^{\overline{\sigma}} \right) 
\end{align*}
i.e. $O\left( \overline{\sigma}_a^{\overline{\sigma}} + \overline{\eta}_a^{\overline{\sigma}} \right)$. 

Then, we can rewrite 
\begin{align*}
\widehat{X}_a(\overline{\sigma}) \cdot \nabla_{\overline{\sigma}} \varphi_3 &= \frac{\sigma_0}{|\overline{\sigma}|} \left( 3 \rho_0^2 + |\rho|^2 - 2 \sqrt{3} \frac{\sigma_0}{|\sigma_0|} \rho_0 |\rho| - 12 \sigma_0^2 \right) + O\left( \rho_0 \xi_t^{\overline{\rho} \overline{\sigma}} \right) + O\left( \overline{\sigma}_b^{\overline{\sigma}} \right) 
\end{align*}
so we can also control terms with a factor $O\left( |\rho| - 2 \sqrt{3} |\sigma_0| - \sqrt{3} \frac{\sigma_0}{|\sigma_0|} \rho_0 \right)$. 

Finally, we can combine with {\footnotesize 
\begin{align*}
\partial_{\eta_0} \varphi_3 &= 3 \rho_0^2 + |\rho|^2 - 3 \eta_0^2 - |\eta|^2 \\
&= 3 \rho_0^2 + \left( 2 \sqrt{3} \sigma_0 + \sqrt{3} \rho_0 \right)^2 - \frac{1}{2} \left( \left( \overline{\eta}_a^{\overline{\sigma}} \right)^2 + \left( \overline{\eta}_b^{\overline{\sigma}} \right)^2 \right) + O\left( |\eta| \xi_t^{\overline{\eta} \overline{\sigma}} \right) + O\left( |\rho| - 2 \sqrt{3} |\sigma_0| - \sqrt{3} \frac{\sigma_0}{|\sigma_0|} \rho_0 \right) \\
&= 6 \left( \sigma_0 + \rho_0 \right)^2 - \frac{1}{2} \left( \overline{\eta}_b^{\overline{\sigma}} \right)^2 + O\left( |\eta| \xi_t^{\overline{\eta} \overline{\sigma}} \right) + O\left( |\rho| - 2 \sqrt{3} |\sigma_0| - \sqrt{3} \frac{\sigma_0}{|\sigma_0|} \rho_0 \right) + O\left( \overline{\sigma}_a^{\overline{\sigma}} + \overline{\eta}_a^{\overline{\sigma}} \right) + O\left( \overline{\sigma}_b^{\overline{\sigma}} \right)
\end{align*} }
which allows to control terms with a factor $O\left( \overline{\eta}_b^{\overline{\sigma}} + 2 \sqrt{3} \frac{\sigma_0}{|\sigma_0|} (\sigma_0 + \rho_0) \right)$. We now write {\footnotesize 
\begin{align*}
\eta_0 &= \frac{1}{2 \sqrt{3}} \frac{\sigma_0}{|\sigma_0|} \left( \overline{\eta}_a^{\overline{\sigma}} + \overline{\eta}_b^{\overline{\sigma}} \right) \\
&= \frac{1}{2 \sqrt{3}} \frac{\sigma_0}{|\sigma_0|} \left( - 2 \sqrt{3} |\sigma_0| - 2 \sqrt{3} \frac{\sigma_0}{|\sigma_0|} (\sigma_0 + \rho_0) \right) + O\left( \overline{\eta}_a^{\overline{\sigma}} + \overline{\sigma}_a^{\overline{\sigma}} \right) + O\left( \overline{\sigma}_b^{\overline{\sigma}} \right) + O\left( \overline{\eta}_b^{\overline{\sigma}} + 2 \sqrt{3} \frac{\sigma_0}{|\sigma_0|} (\sigma_0 + \rho_0) \right) \\
&= - 2 \sigma_0 - \rho_0 + O\left( \overline{\eta}_a^{\overline{\sigma}} + \overline{\sigma}_a^{\overline{\sigma}} \right) + O\left( \overline{\sigma}_b^{\overline{\sigma}} \right) + O\left( \overline{\eta}_b^{\overline{\sigma}} + 2 \sqrt{3} \frac{\sigma_0}{|\sigma_0|} (\sigma_0 + \rho_0) \right) \\
\xi_0 &= \sigma_0 + \eta_0 + \rho_0 \\
&= - \sigma_0 + O\left( \overline{\eta}_a^{\overline{\sigma}} + \overline{\sigma}_a^{\overline{\sigma}} \right) + O\left( \overline{\sigma}_b^{\overline{\sigma}} \right) + O\left( \overline{\eta}_b^{\overline{\sigma}} + 2 \sqrt{3} \frac{\sigma_0}{|\sigma_0|} (\sigma_0 + \rho_0) \right) \\
\frac{\sigma \cdot \eta}{|\sigma|} &= \frac{1}{2} \left( \overline{\eta}_a^{\overline{\sigma}} - \overline{\eta}_b^{\overline{\sigma}} \right) \\
&= \sqrt{3} \frac{\sigma_0}{|\sigma_0|} \rho_0 + O\left( \overline{\eta}_a^{\overline{\sigma}} + \overline{\sigma}_a^{\overline{\eta}} \right) + O\left( \overline{\sigma}_b^{\overline{\sigma}} \right) + O\left( \overline{\eta}_b^{\overline{\sigma}} + 2 \sqrt{3} \frac{\sigma_0}{|\sigma_0|} (\sigma_0 + \rho_0) \right) \\
|\xi|^2 &= \left( \frac{\sigma \cdot \xi}{|\sigma|} \right)^2 + \left( \frac{J \sigma \cdot \xi}{|\sigma|} \right)^2 = \left( |\sigma| + \frac{\sigma \cdot \eta}{|\sigma|} + \frac{\sigma \cdot \rho}{|\sigma|} \right)^2 + \left( \frac{J \sigma \cdot \rho}{|\sigma|} \right)^2 + O\left( |\eta| \xi_t^{\overline{\eta} \overline{\sigma}} \right) \\
&= \left( \sqrt{3} |\sigma_0| + \sqrt{3} \frac{\sigma_0}{|\sigma_0|} \rho_0 + \frac{\sigma \cdot \rho}{|\sigma|} \right)^2 + \left( \frac{J \sigma \cdot \rho}{|\sigma|} \right)^2 + O\left( |\eta| \xi_t^{\overline{\eta} \overline{\sigma}} \right) + O\left( \overline{\eta}_b^{\overline{\sigma}} + 2 \sqrt{3} \frac{\sigma_0}{|\sigma_0|} (\sigma_0 + \rho_0) \right) \\
&\quad + O\left( \overline{\eta}_a^{\overline{\sigma}} + \overline{\sigma}_a^{\overline{\eta}} \right) + O\left( \overline{\sigma}_b^{\overline{\sigma}} \right) \\
&= 3 \sigma_0^2 + 3 \rho_0^2 + 6 \sigma_0 \rho_0 + |\rho|^2 + 2 \sqrt{3} \frac{\sigma_0}{|\sigma_0|} \left( \sigma_0 + \rho_0 \right) \frac{\sigma \cdot \rho}{|\sigma|} + O\left( |\eta| \xi_t^{\overline{\eta} \overline{\sigma}} \right) + O\left( \overline{\eta}_b^{\overline{\sigma}} + 2 \sqrt{3} \frac{\sigma_0}{|\sigma_0|} (\sigma_0 + \rho_0) \right) \\
&\quad + O\left( \overline{\eta}_a^{\overline{\sigma}} + \overline{\sigma}_a^{\overline{\eta}} \right) + O\left( \overline{\sigma}_b^{\overline{\sigma}} \right) \\
&= 3 \sigma_0^2 + 3 \rho_0^2 + 6 \sigma_0 \rho_0 + 3 \left( 2 \sigma_0 + \rho_0 \right)^2 
- 2 \sqrt{3} \frac{\sigma_0}{|\sigma_0|} \left( \sigma_0 + \rho_0 \right) |\rho| 
+ 2 \sqrt{3} \frac{\sigma_0}{|\sigma_0|} \left( \sigma_0 + \rho_0 \right) |\rho| \left( 1 + \theta^{\overline{\sigma} \overline{\rho}} \right) \\
&\quad + O\left( |\eta| \xi_t^{\overline{\eta} \overline{\sigma}} \right) + O\left( \overline{\eta}_a^{\overline{\sigma}} + \overline{\sigma}_a^{\overline{\eta}} \right) + O\left( \overline{\sigma}_b^{\overline{\sigma}} \right) + O\left( \overline{\eta}_b^{\overline{\sigma}} + 2 \sqrt{3} \frac{\sigma_0}{|\sigma_0|} (\sigma_0 + \rho_0) \right) + O\left( |\rho| - 2 \sqrt{3} |\sigma_0| - \sqrt{3} \frac{\sigma_0}{|\sigma_0|} \rho_0 \right) \\
&= 3 \sigma_0^2 + 12 \sigma_0^2 \left( 1 + \theta^{\overline{\sigma} \overline{\rho}} \right) + O\left( |\eta| \xi_t^{\overline{\eta} \overline{\sigma}} \right) + O\left( \overline{\eta}_b^{\overline{\sigma}} + 2 \sqrt{3} \frac{\sigma_0}{|\sigma_0|} (\sigma_0 + \rho_0) \right) + O\left( |\rho| - 2 \sqrt{3} |\sigma_0| - \sqrt{3} \frac{\sigma_0}{|\sigma_0|} \rho_0 \right) \\
&\quad + O\left( \rho_0 \xi_t^{\overline{\sigma} \overline{\rho}} \right) + O\left( \overline{\eta}_a^{\overline{\sigma}} + \overline{\sigma}_a^{\overline{\eta}} \right) + O\left( \overline{\sigma}_b^{\overline{\sigma}} \right)
\end{align*} }

It thus only remains to simplify {\footnotesize 
\begin{align*}
\varphi_3 &= \xi_0^3 + \xi_0 |\xi|^2 - \eta_0^3 - \eta_0 |\eta|^2 - \sigma_0^3 - \sigma_0 |\sigma|^2 - \rho_0^3 - \rho_0 |\rho|^2 \\
&= - \sigma_0^3 - \sigma_0 \left( 3 \sigma_0^2 + 12 \sigma_0^2 \left( 1 + \theta^{\overline{\sigma} \overline{\rho}} \right) \right) 
+ (2 \sigma_0 + \rho_0)^3 + 3 (2 \sigma_0 + \rho_0) \rho_0^2 - 4 \sigma_0^3 - \rho_0^3 - 3 \rho_0 \left( 2 \sigma_0 + \rho_0 \right)^2 + O\left( |\eta| \xi_t^{\overline{\eta} \overline{\sigma}} \right) \\
&\quad + O\left( \overline{\eta}_a^{\overline{\sigma}} + \overline{\sigma}_a^{\overline{\eta}} \right) + O\left( \overline{\sigma}_b^{\overline{\sigma}} \right) + O\left( \overline{\eta}_b^{\overline{\sigma}} + 2 \sqrt{3} \frac{\sigma_0}{|\sigma_0|} (\sigma_0 + \rho_0) \right) + O\left( |\rho| - 2 \sqrt{3} |\sigma_0| - \sqrt{3} \frac{\sigma_0}{|\sigma_0|} \rho_0 \right) + O\left( \rho_0 \xi_t^{\overline{\sigma} \overline{\rho}} \right) \\
&= - 12 \sigma_0^3 \left( 1 + \theta^{\overline{\sigma} \overline{\rho}} \right) + O\left( \overline{\eta}_b^{\overline{\sigma}} + 2 \sqrt{3} \frac{\sigma_0}{|\sigma_0|} (\sigma_0 + \rho_0) \right) + O\left( |\rho| - 2 \sqrt{3} |\sigma_0| - \sqrt{3} \frac{\sigma_0}{|\sigma_0|} \rho_0 \right) + O\left( \rho_0 \xi_t^{\overline{\sigma} \overline{\rho}} \right) \\
&\quad + O\left( |\eta| \xi_t^{\overline{\eta} \overline{\sigma}} \right) + O\left( \overline{\eta}_a^{\overline{\sigma}} + \overline{\sigma}_a^{\overline{\eta}} \right) + O\left( \overline{\sigma}_b^{\overline{\sigma}} \right)
\end{align*} }
so we have a control over terms with a factor $O\left( 1 + \theta^{\overline{\sigma} \overline{\rho}} \right)$. 

But {\footnotesize 
\begin{align*}
3 \xi_0^2 - |\xi|^2 &= 3 \sigma_0^2 - \left( 3 \sigma_0^2 + 12 \sigma_0^2 \left( 1 + \theta^{\overline{\sigma} \overline{\rho}} \right) \right) + O\left( |\eta| \xi_t^{\overline{\eta} \overline{\sigma}} \right) + O\left( \overline{\eta}_a^{\overline{\sigma}} + \overline{\sigma}_a^{\overline{\eta}} \right) + O\left( \overline{\sigma}_b^{\overline{\sigma}} \right) \\
&\quad + O\left( \overline{\eta}_b^{\overline{\sigma}} + 2 \sqrt{3} \frac{\sigma_0}{|\sigma_0|} (\sigma_0 + \rho_0) \right) + O\left( |\rho| - 2 \sqrt{3} |\sigma_0| - \sqrt{3} \frac{\sigma_0}{|\sigma_0|} \rho_0 \right) + O\left( \rho_0 \xi_t^{\overline{\sigma} \overline{\rho}} \right) \\
&= O(\varphi_3) + O\left( |\eta| \xi_t^{\overline{\eta} \overline{\sigma}} \right) + O\left( \overline{\eta}_a^{\overline{\sigma}} + \overline{\sigma}_a^{\overline{\eta}} \right) \\
&\quad + O\left( \overline{\sigma}_b^{\overline{\sigma}} \right) + O\left( \overline{\eta}_b^{\overline{\sigma}} + 2 \sqrt{3} \frac{\sigma_0}{|\sigma_0|} (\sigma_0 + \rho_0) \right) + O\left( |\rho| - \sqrt{3} \frac{\sigma_0}{|\sigma_0|} \left( 2 \sigma_0 + \rho_0 \right) \right) + O\left( \rho_0 \xi_t^{\overline{\sigma} \overline{\rho}} \right)
\end{align*} }
which concludes. 

\paragraph{2.} Let us now localise to have $|\overline{\eta}| \simeq |\overline{\sigma}| \gg |\overline{\rho}|$. 

Again, \eqref{equdecchampbLCP-2} is of the form \eqref{lemestimeesgeneriquesdeccubiquehb2-eta0sigmaresxrho}, and since 
\begin{align*}
1 = O(\partial_{\eta_0} \varphi_3) 
\end{align*}
we can easily rewrite \eqref{equdecchampbLCP-1} as terms from \eqref{qteslemestimeesgeneriquesdeccubiquehb2}. 

\paragraph{3.} Let us now localise to have $|\overline{\eta}| \simeq |\overline{\rho}| \gg |\overline{\sigma}|$. 

Then we can develop in \eqref{equdecchampbLCP-2}
\begin{align*}
\widehat{X}_b(\overline{\xi}) \cdot \nabla_{\overline{\xi}} &= \frac{|\xi|}{|\overline{\xi}|} \partial_{\xi_0} - \frac{\xi_0 \xi}{|\overline{\xi}| |\xi|} \cdot \nabla_{\xi} 
\end{align*}
The second term contributes like terms from \eqref{qteslemestimeesgeneriquesdeccubiquehb2}. For the first, we can apply an integration by parts in $\eta_0$ and recover terms from \eqref{qteslemestimeesgeneriquesdeccubiquehb2} plus a term with the form
\begin{align*}
\int_0^t \int \int e^{i s \varphi_3} s^2 |\overline{\eta}|^4 \mu \widehat{f}(s, \overline{\eta}) \widehat{f}(s, \overline{\sigma}) \widehat{f}(s, \overline{\rho}) ~ d\overline{\eta} d\overline{\sigma} ds 
\end{align*}
\eqref{equdecchampbLCP-1} also has the same form. But on the one hand, 
\begin{align*}
1 = O\left( \widehat{X}_a(\overline{\sigma}) \cdot \nabla_{\overline{\sigma}} \varphi_3 \right), \quad 1 = O\left( \widehat{X}_a(\overline{\sigma}) \cdot \left( \nabla_{\overline{\eta}} - \nabla_{\overline{\sigma}} \right) \varphi_3 \right) 
\end{align*}
and on the other
\begin{align*}
\varphi_3 &= \xi_0^3 + \xi_0 |\xi|^2 - \eta_0^3 - \eta_0 |\eta|^2 + O\left( \overline{\sigma} \right) + O(\rho_0) \\
&= \eta_0 |\rho|^2 + O\left( \overline{\sigma} \right) + O(\rho_0) + O(\eta)
\end{align*}
so that 
\begin{align*}
1 = O(\varphi_3) + O\left( \overline{\sigma} \right) + O(\rho_0) + O(\eta)
\end{align*}
We may then apply an integration by parts in time when we have a factor $O(\varphi_3)$, or an integration by parts along $\widehat{X}_a(\overline{\sigma}) \cdot \nabla_{\overline{\eta}}$ if we have a factor $O(\overline{\sigma}) + O(\rho_0)$, or along $\widehat{X}_a(\overline{\sigma}) \cdot \left( \nabla_{\overline{\eta}} - \nabla_{\overline{\sigma}} \right)$ for the factor $O(\eta)$, and we recover indeed terms from \eqref{qteslemestimeesgeneriquesdeccubiquehb2}. 

\paragraph{4.} Let us now localise to have $|\overline{\sigma}| \simeq |\overline{\rho}| \gg |\overline{\eta}|$. 

In this case, we can again apply an integration by parts on \eqref{equdecchampbLCP-2} to rewrite it as terms from \eqref{qteslemestimeesgeneriquesdeccubiquehb2} plus 
\begin{align*}
\int_0^t \int \int e^{i s \varphi_3} s^2 m_b(\overline{\xi}) |\overline{\sigma}|^4 \mu \widehat{f}(s, \overline{\eta}) \widehat{f}(s, \overline{\sigma}) \widehat{f}(s, \overline{\rho}) ~ d\overline{\eta} d\overline{\sigma} ds 
\end{align*}
\eqref{equdecchampbLCP-1} also has this form

Here, $\xi_0 = \sigma_0 + o(1)$ so $|\overline{\xi}| \simeq |\overline{\sigma}| \simeq |\xi_0|$. 

But 
\begin{align*}
\varphi_3 &= O(\rho_0) + O(\overline{\eta}) + \xi_0^3 + \xi_0 |\xi|^2 - \sigma_0^3 - \sigma_0 |\sigma|^2 \\
&= O(\rho_0) + O(\overline{\eta}) + \xi_0 \left( |\xi|^2 - |\sigma|^2 \right) 
\end{align*}
We deduce that 
\begin{align*}
m_b(\overline{\xi}) &= O\left( 3 \xi_0^2 - |\xi|^2 \right) = O(\varphi_3) + O(\rho_0) + O(\overline{\eta}) + O(m_b(\overline{\sigma}))
\end{align*}
We can now apply an integration by parts in time for the factor $O(\varphi_3)$, along $\partial_{\eta_0}$ for the factor $O(\rho_0)$ and the factor $O(\overline{\eta})$, and along $\partial_{\eta_0}-\partial_{\sigma_0}$ for the factor $O(m_b(\overline{\sigma}))$. In all cases, we get terms from \eqref{qteslemestimeesgeneriquesdeccubiquehb2}. 

\paragraph{5.} Let us now localise to have $|\overline{\eta}| \gg |\overline{\sigma}|+|\overline{\rho}|$. 

In this case, we can decompose $\widehat{X}_b(\overline{\xi})$ in \eqref{equdecchampbLCP-2} as terms from \eqref{qteslemestimeesgeneriquesdeccubiquehb2} (using that $\mu$ allows to distribute a derivative), and a term on which we apply an integration by parts in $\eta_0$ as before, to get terms from \eqref{qteslemestimeesgeneriquesdeccubiquehb2} plus a term explicitely equal to: 
\begin{align*}
\int_0^t \int \int e^{i s \varphi_3} i s^2 \xi_0 m_b(\overline{\xi}) \frac{|\xi|}{|\overline{\xi}|} \partial_{\eta_0} \varphi_3 \mu(\overline{\xi}, \overline{\eta}, \overline{\sigma}) \widehat{F}_1(\overline{\eta}) \widehat{F}_2(\overline{\sigma}) \widehat{F}_3(\overline{\rho}) ~ d\overline{\eta} d\overline{\sigma} ds 
\end{align*}
and that we group with \eqref{equdecchampbLCP-1} to get: 
\begin{align}
\int_0^t \int \int e^{i s \varphi_3} i s^2 \xi_0 m_b(\overline{\xi}) \left( \frac{|\xi|}{|\overline{\xi}|} \left( \partial_{\xi_0} + \partial_{\eta_0} \right) - \frac{\xi_0 \xi}{|\overline{\xi}| |\xi|} \cdot \nabla_{\xi} \right) \varphi_3 \mu(\overline{\xi}, \overline{\eta}, \overline{\sigma}) \widehat{F}_1(\overline{\eta}) \widehat{F}_2(\overline{\sigma}) \widehat{F}_3(\overline{\rho}) ~ d\overline{\eta} d\overline{\sigma} ds \label{termefinalcubiquehb2-LCP-HHBB} 
\end{align}

Here, we have that 
\begin{align*}
1 &= O(\partial_{\eta_0} \varphi_3) \\
1 &= O\left( \widehat{X}_a(\overline{\sigma}) \cdot \left( \nabla_{\overline{\eta}} - \nabla_{\overline{\sigma}} \right) \varphi_3 \right) 
\end{align*}
Furthermore, $\xi_0 = \eta_0 + o(1)$ so $|\xi_0| \simeq |\overline{\xi}| \simeq |\overline{\eta}|$, and 
\begin{align*}
\varphi_3 &= O(\rho_0) + O(\overline{\sigma}) + \xi_0^3 + \xi_0 |\xi|^2 - \eta_0^3 - \eta_0 |\eta|^2 \\
&= O(\rho_0) + O(\overline{\sigma}) + \xi_0 \left( |\xi|^2 - |\eta|^2 \right) 
\end{align*}
hence $|\xi|^2 - |\eta|^2 = O(\varphi_3) + O(\rho_0) + O(\overline{\sigma})$. Note then that, in \eqref{termefinalcubiquehb2-LCP-HHBB}, if we obtain a factor $O(|\overline{\rho}|\varphi_3)$, $O(\rho_0)$ or $O(\overline{\sigma})$, we can apply an integration by parts in time, or in $\eta_0$ in the two last cases, and we recover in any case terms from \eqref{qteslemestimeesgeneriquesdeccubiquehb2}. In particular, we can assume that $F_1 = \partial_x f$, $F_2 = F_3 = f$, and that $\mu$ gives an additional factor $|\overline{\rho}|$. We can then compute that, since $m_b(\overline{\xi}) = O(\xi)$, 
\begin{align*}
&|\xi| |\overline{\rho}| \left( \frac{|\xi|}{|\overline{\xi}|} \left( \partial_{\xi_0} + \partial_{\eta_0} \right) - \frac{\xi_0 \xi}{|\overline{\xi}| |\xi|} \cdot \nabla_{\xi} \right) \varphi_3 \\
&\quad = |\overline{\rho}| \frac{|\xi|^2}{|\overline{\xi}|} \left( 3 \xi_0^2 + |\xi|^2 - 3 \eta_0^2 - |\eta|^2 \right) - |\overline{\rho}| \frac{\xi_0}{|\overline{\xi}|} \left( 2 \xi_0 |\xi|^2 - 2 \rho_0 |\rho|^2 \right) \\
&\quad = O\left( |\overline{\rho}| \varphi_3 \right) + O(\rho_0) + O(\overline{\sigma}) + O\left( |\eta| |\overline{\rho}| \right) 
\end{align*}
But we can also apply an integration by parts along $\widehat{X}_a(\overline{\sigma}) \cdot \left( \nabla_{\overline{\eta}} - \nabla_{\overline{\sigma}} \right)$ when we have the factor $O\left( |\eta| |\overline{\rho}| \right)$ and we recover terms from \eqref{qteslemestimeesgeneriquesdeccubiquehb2}. 

\paragraph{6.} Localise now to have $|\overline{\sigma}| \gg |\overline{\eta}|+|\overline{\rho}|$. 

In this case, as before, we can rewrite \eqref{equdecchampbLCP-2} as terms from \eqref{qteslemestimeesgeneriquesdeccubiquehb2} (using that $\mu$ allows to distribute a derivative), plus a term of the form 
\begin{align*}
\int_0^t \int \int e^{i s \varphi_3} s^2 m_b(\overline{\xi}) |\overline{\sigma}|^4 \mu \widehat{f}(s, \overline{\eta}) \widehat{f}(s, \overline{\sigma}) \widehat{f}(s, \overline{\rho}) ~ d\overline{\eta} d\overline{\sigma} ds 
\end{align*}
\eqref{equdecchampbLCP-1} also has this form. But 
\begin{align*}
\varphi_3 &= O(\rho_0) + O(\overline{\eta}) + \xi_0 (|\xi|^2 - |\sigma|^2) 
\end{align*}
so 
\begin{align*}
m_b(\overline{\xi}) &= O\left( 3 \xi_0^2 - |\xi|^2 \right) \\
&= O\left( m_b(\overline{\sigma}) \right) + O(\varphi_3) + O(\rho_0) + O(\overline{\eta}) 
\end{align*}
On the other hand we have that 
\begin{align*}
1 = O\left( \widehat{X}_a(\overline{\sigma}) \cdot \nabla_{\overline{\sigma}} \varphi_3 \right), \quad 1 = O\left( \widehat{X}_a(\overline{\sigma}) \cdot \left( \nabla_{\overline{\eta}} - \nabla_{\overline{\sigma}} \right) \varphi_3 \right), \quad 1 = O\left( \nabla_{\sigma} \varphi_3 \right) 
\end{align*}
We can use $\mu$ to distribute a derivative, and apply an integration by parts in time for the factor $O(\varphi_3)$; along $\nabla_{\sigma}$ for the factor $O\left( m_b(\overline{\sigma}) \right)$; along $\widehat{X}_a(\overline{\sigma}) \cdot \nabla_{\overline{\sigma}} \varphi_3$ for the factors $O(\rho_0)$ and $O(\overline{\eta})$. In all cases, we recover terms from \eqref{qteslemestimeesgeneriquesdeccubiquehb2}. 

\paragraph{7.} Let us finally localise to have $|\overline{\rho}| \gg |\overline{\eta}|+|\overline{\sigma}|$. 

Again, we rewrite \eqref{equdecchampbLCP-2}, up to terms from \eqref{qteslemestimeesgeneriquesdeccubiquehb2}, as 
\begin{align*}
\int_0^t \int \int e^{i s \varphi_3} s^2 \xi_0 |\overline{\rho}|^3 \mu \widehat{f}(s, \overline{\eta}) \widehat{f}(s, \overline{\sigma}) \widehat{f}(s, \overline{\rho}) ~ d\overline{\eta} d\overline{\sigma} ds 
\end{align*}
and \eqref{equdecchampbLCP-1} also has this form, where $\mu$ allows to distribute a derivative. 

We then have, by the same computation as for the interaction $\widehat{\mathcal{C}}\widehat{\mathcal{C}}\widehat{\mathcal{P}}$, case 5., 
\begin{align*}
\varphi_3 &= (\xi_0 - \rho_0) \left( |\rho|^2 + o(1) \right) + O\left( \rho_0 \overline{\eta} \right) + O\left( \rho_0 \overline{\sigma} \right) + O\left( \overline{\eta} \overline{\sigma} \right) 
\end{align*}
and then 
\begin{align*}
\xi_0 \mu &= (\xi_0 - \rho_0) \mu + \rho_0 \mu \\
&= O(\mu \varphi_3) + O\left( \rho_0 \overline{\eta} \right) + O\left( \rho_0 \overline{\sigma} \right) + O\left( \overline{\eta} \overline{\sigma} \right)
\end{align*}
We may then apply an integration by parts in time for the factor $O(\mu \varphi_3)$; along $\partial_{\eta_0}$ for the factor $O\left( \rho_0 \overline{\eta} \right)$; along $\widehat{X}_a(\overline{\sigma}) \cdot \nabla_{\overline{\sigma}}$ for the factors $O\left( \rho_0 \overline{\sigma} \right)$ and $O\left( \overline{\eta} \overline{\sigma} \right)$, and we recover terms from \eqref{qteslemestimeesgeneriquesdeccubiquehb2}. 

\subsubsection{Interaction \texorpdfstring{$\widehat{\mathcal{C}}\widehat{\mathcal{P}}\widehat{\mathcal{P}}$}{CPP}}

Let us now consider: {\footnotesize 
\begin{subequations}
\begin{align}
&\xi_0 m_b(\overline{\xi}) \widehat{X}_b(\overline{\xi}) \cdot \nabla_{\overline{\xi}} \widehat{I}_{s\mu}^{\widehat{\mathcal{C}}\widehat{\mathcal{P}}\widehat{\mathcal{P}}}[F_1, F_2, F_3](t, \overline{\xi}) \notag \\
&\quad = \int_0^t \int \int i s^2 \xi_0 m_b(\overline{\xi}) \widehat{X}_b(\overline{\xi}) \cdot \nabla_{\overline{\xi}} \varphi_3 e^{i s \varphi_3} \mu(\overline{\xi}, \overline{\eta}, \overline{\sigma}) m_{\widehat{\mathcal{C}}}(\overline{\eta}) m_{\widehat{\mathcal{P}}}(\overline{\sigma}) m_{\widehat{\mathcal{P}}}(\overline{\rho}) \widehat{F}_1(s, \overline{\eta}) \widehat{F}_2(s, \overline{\sigma}) \widehat{F}_3(s, \overline{\rho}) ~ d\overline{\eta} d\overline{\sigma} ds \label{equdecchampbCPP-1} \\
&\quad \quad + \int_0^t \int \int e^{i s \varphi_3} s \xi_0 \mu(\overline{\xi}, \overline{\eta}, \overline{\sigma}) m_{\widehat{\mathcal{C}}}(\overline{\eta}) m_{\widehat{\mathcal{P}}}(\overline{\sigma}) m_{\widehat{\mathcal{P}}}(\overline{\rho}) \widehat{F}_1(s, \overline{\eta}) \widehat{F}_2(s, \overline{\sigma}) m_b(\overline{\xi}) \widehat{X}_b(\overline{\xi}) \cdot \nabla_{\overline{\xi}} \widehat{F}_3(s, \overline{\rho}) ~ d\overline{\eta} d\overline{\sigma} ds \label{equdecchampbCPP-2} \\
&\quad \quad + \int_0^t \int \int e^{i s \varphi_3} s \xi_0 m_b(\overline{\xi}) \widehat{X}_b(\overline{\xi}) \cdot \nabla_{\overline{\xi}} \left( \mu(\overline{\xi}, \overline{\eta}, \overline{\sigma}) m_{\widehat{\mathcal{C}}}(\overline{\eta}) m_{\widehat{\mathcal{P}}}(\overline{\sigma}) m_{\widehat{\mathcal{P}}}(\overline{\rho}) \right) \widehat{F}_1(s, \overline{\eta}) \widehat{F}_2(s, \overline{\sigma}) \widehat{F}_3(s, \overline{\rho}) ~ d\overline{\eta} d\overline{\sigma} ds \label{equdecchampbCPP-3} 
\end{align}
\end{subequations} }
\eqref{equdecchampbCPP-3} is of the form $\eqref{lemestimeesgeneriquesdeccubiquehb2-dersymb}+\eqref{lemestimeesgeneriquesdeccubiquehb2-dersymbbis}$. 

We then decompose depending on the relative sizes. Without loss of generality, up to applying this decomposition before applying $\xi_0 m_b(\overline{\xi}) \widehat{X}_b(\overline{\xi}) \cdot \nabla_{\overline{\xi}}$, we can always assume that $|\overline{\sigma}| \gtrsim |\overline{\rho}|$. 

\paragraph{1.} Let us first localise to have $|\overline{\eta}| \simeq |\overline{\sigma}| \simeq |\overline{\rho}|$. \eqref{equdecchampbCPP-2} is then already of the form \eqref{lemestimeesgeneriquesdeccubiquehb2-eta0sigmaresxrho}. 

Then, 
\begin{align*}
\partial_{\sigma_0} \varphi_3 &= o(1) + |\rho|^2 - |\sigma|^2 \\
\widehat{X}_a(\overline{\eta}) \cdot \nabla_{\overline{\eta}} \varphi_3 &= o(1) + \frac{\eta_0}{|\overline{\eta}|} \left( |\rho|^2 - 12 \eta_0^2 \right) 
\end{align*}
so we only need to consider the neighborhood of $|\rho| = |\sigma| = 2 \sqrt{3} |\eta_0|$. In particular, $\xi_0 = \eta_0 + o(1)$ so $|\overline{\xi}| \simeq |\overline{\eta}|$ as well. But then 
\begin{align*}
\varphi_3 &= \xi_0^3 + \xi_0 |\xi|^2 - \eta_0^3 - \eta_0 |\eta|^2 + o(1) \\
&= \eta_0 (|\xi|^2 - |\eta|^2) + o(1) 
\end{align*}
So we can restrict also to the neighborhood of $|\xi| = |\eta|$. Therefore $\overline{\xi}$ is in the neighborhood of $\widehat{\mathcal{C}}$. 

In all cases, $\nabla_{\eta} \varphi_3 \simeq |\overline{\eta}|^2$, so we have a good control over terms with a factor $O\left( \overline{\eta}_b^{\overline{\eta}} \right)$. 

\paragraph{1.1.} Let us first localise to have $\xi_t^{\overline{\eta} \overline{\rho}}$ away enough from $0$, that is $\eta$ and $\rho$ not close to alignment. Then $\xi_t^{\overline{\eta} \overline{\sigma}}$ cannot be close to $0$: if it were the case, we would have 
\begin{align*}
|\xi|^2 &= |\eta|^2 + o(1) \\
&= |\eta|^2 + |\sigma|^2 + |\rho|^2 + 2 \eta \cdot \sigma + 2 \eta \cdot \rho + 2 \sigma \cdot \rho \\
&= 9 |\eta|^2 + 4 |\eta|^2 \theta^{\overline{\eta} \overline{\sigma}} + 4 |\eta|^2 \theta^{\overline{\eta} \overline{\rho}} + 8 |\eta|^2 \theta^{\overline{\eta} \overline{\sigma}} \theta^{\overline{\sigma} \overline{\rho}} + o(1) \\
&= |\eta|^2 \left( 9 + 4 \theta^{\overline{\eta} \overline{\sigma}} + 4 \theta^{\overline{\eta} \overline{\rho}} + 8 \theta^{\overline{\eta} \overline{\sigma}} \theta^{\overline{\sigma} \overline{\rho}} \right) + o(1) 
\end{align*}
But $\theta^{\overline{\eta} \overline{\sigma}}$ needs to be $1 + o(1)$ or $-1 + o(1)$, while $\theta^{\overline{\eta} \overline{\rho}}$ is away enough from $\pm 1$. If $\theta^{\overline{\eta} \overline{\sigma}} = 1 + o(1)$, we can then simplify 
\begin{align*}
|\xi|^2 = |\eta|^2 \left( 13 + 12 \theta^{\overline{\eta} \overline{\rho}} \right) + o(1)
\end{align*}
which forces $|\xi|^2 - |\eta|^2 \gtrsim |\overline{\eta}|^2$, and this is a contradiction. If $\theta^{\overline{\eta} \overline{\sigma}} = -1 + o(1)$, we have instead
\begin{align*}
|\xi|^2 
&= |\eta|^2 \left( 5 - 4 \theta^{\overline{\eta} \overline{\rho}} \right) + o(1) 
\end{align*}
and again $|\xi|^2 - |\eta|^2 \gtrsim |\overline{\eta}|^2$, a contradiction. We deduce that $\xi_t^{\overline{\eta} \overline{\sigma}} \gtrsim 1$. 

Finally, we write that 
\begin{align*}
\widehat{X}_c(\overline{\eta}) \cdot \nabla_{\overline{\eta}} \varphi_3 &= 2 \rho_0 \frac{J \eta \cdot \rho}{|\eta|} \simeq |\overline{\eta}| \rho_0 
\end{align*}
In particular, having a factor $O(\rho_0)$ allows a good control by integration by parts along $\widehat{X}_c(\overline{\eta}) \cdot \nabla_{\overline{\eta}}$ to recover terms from \eqref{qteslemestimeesgeneriquesdeccubiquehb2}. Symmetric-wise, we can treat the factor $O(\sigma_0)$. But now 
\begin{align*}
\varphi_3 &= \xi_0^3 + \xi_0 |\xi|^2 - \eta_0^3 - \eta_0 |\eta|^2 + O(\rho_0) + O(\sigma_0) \\
&= \xi_0 (|\xi|^2 - |\eta|^2) + O(\rho_0) + O(\sigma_0)
\end{align*}
so we have a good control over terms with a factor $O(|\xi| - |\eta|)$. Finally, 
\begin{align*}
3 \xi_0^2 - |\xi|^2 &= 3 \eta_0^2 - |\eta|^2 + O(\sigma_0) + O(\rho_0) + O(|\xi| - |\eta|) \\
&= O(\sigma_0) + O(\rho_0) + O(|\xi| - |\eta|) + O\left( \overline{\eta}_b^{\overline{\eta}} \right) 
\end{align*}
and we precisely have such a factor in $m_b(\overline{\xi})$. 

\paragraph{1.2.} Let us now localise to have $\xi_t^{\overline{\eta} \overline{\sigma}}$ and $\xi_t^{\overline{\eta} \overline{\rho}}$ close to $0$. In particular, up to exchanging $\overline{\sigma}$ and $\overline{\rho}$, we can assume that $\theta^{\overline{\eta} \overline{\sigma}}$ is close to $1$ and $\theta^{\overline{\eta} \overline{\rho}}$ is close to $-1$. 

Let us compute explicitely some factors allowing for good control. $\nabla_{\eta} \varphi_3 \simeq |\overline{\eta}|^2$ so the factor $O\left( \overline{\eta}_b^{\overline{\eta}} \right)$ allows for a good control. Then, 
\begin{align*}
\widehat{X}_c(\overline{\eta}) \cdot \nabla_{\overline{\eta}} \varphi_3 &= 2 \rho_0 \frac{J \eta \cdot \rho}{|\eta|} \simeq \rho_0 |\rho| \xi_t^{\overline{\eta} \overline{\rho}} \\
\widehat{X}_a(\overline{\eta}) \cdot \nabla_{\overline{\eta}} \varphi_3 &= \frac{\eta_0}{|\overline{\eta}|} \left( 3 \rho_0^2 + |\rho|^2 - 2 \sqrt{3} \frac{\eta_0}{|\eta_0|} \rho_0 |\rho| - 12 \eta_0^2 \right) + O\left( \overline{\eta}_b^{\overline{\eta}} \right) + O\left( \rho_0 \xi_t^{\overline{\eta} \overline{\rho}} \right) 
\end{align*}
We may thus estimate terms with a factor $O\left( |\rho| - 2 \sqrt{3} |\eta_0| - \sqrt{3} \frac{\eta_0}{|\eta_0|} \rho_0 \right)$ or $O\left( \rho_0 \xi_t^{\overline{\eta} \overline{\rho}} \right)$; likewise, replacing $\overline{\rho}$ by $\overline{\sigma}$, we can also estimate terms with a factor $O\left( |\sigma| - 2 \sqrt{3} |\eta_0| + \sqrt{3} \frac{\eta_0}{|\eta_0|} \sigma_0 \right)$ or $O\left( \sigma_0 \xi_t^{\overline{\eta} \overline{\sigma}} \right)$. Then, we compute that {\footnotesize 
\begin{align*}
&\partial_{\sigma_0} \varphi_3 = 3 \rho_0^2 + |\rho|^2 - 3 \sigma_0^2 - |\sigma|^2 \\
&\quad = 3 \rho_0^2 + 3 \left( 2 \eta_0 + \rho_0 \right)^2 - 3 \sigma_0^2 - 3 \left( 2 \eta_0 - \sigma_0 \right)^2 + O\left( |\rho| - 2 \sqrt{3} |\eta_0| - \sqrt{3} \frac{\eta_0}{|\eta_0|} \rho_0 \right) + O\left( |\sigma| - 2 \sqrt{3} |\eta_0| + \sqrt{3} \frac{\eta_0}{|\eta_0|} \sigma_0 \right) \\
&\quad = 6 (\rho_0 + \sigma_0) \left( 2 \eta_0 + \rho_0 - \sigma_0 \right) 
+ O\left( |\rho| - 2 \sqrt{3} |\eta_0| - \sqrt{3} \frac{\eta_0}{|\eta_0|} \rho_0 \right) + O\left( |\sigma| - 2 \sqrt{3} |\eta_0| + \sqrt{3} \frac{\eta_0}{|\eta_0|} \sigma_0 \right)
\end{align*} }
In particular, we can estimate terms with a factor $O(\rho_0 + \sigma_0) = O(\xi_0 - \eta_0)$. We can thus rewrite {\footnotesize 
\begin{align*}
\varphi_3 &= \xi_0^3 + \xi_0 |\xi|^2 - \eta_0^3 - \eta_0 |\eta|^2 - \rho_0^3 - \rho_0 |\rho|^2 - \sigma_0^3 - \sigma_0 |\sigma|^2 \\
&= \xi_0 (|\xi|^2 - |\eta|^2) - \rho_0 \left( |\rho|^2 - |\sigma|^2 \right) + O(\rho_0 + \sigma_0) \\
&= \xi_0 (|\xi|^2 - |\eta|^2) - \rho_0 (|\rho| + |\sigma|) \left( 2 \sqrt{3} |\eta_0| + \sqrt{3} \frac{\eta_0}{|\eta_0|} \rho_0 - 2 \sqrt{3} |\eta_0| + \sqrt{3} \frac{\eta_0}{|\eta_0|} \sigma_0 \right) \\
&\quad + O(\rho_0 + \sigma_0) + O\left( |\rho| - 2 \sqrt{3} |\eta_0| - \sqrt{3} \frac{\eta_0}{|\eta_0|} \rho_0 \right) + O\left( |\sigma| - 2 \sqrt{3} |\eta_0| + \sqrt{3} \frac{\eta_0}{|\eta_0|} \sigma_0 \right) \\
&= \xi_0 (|\xi|^2 - |\eta|^2) + O(\rho_0 + \sigma_0) + O\left( |\rho| - 2 \sqrt{3} |\eta_0| - \sqrt{3} \frac{\eta_0}{|\eta_0|} \rho_0 \right) + O\left( |\sigma| - 2 \sqrt{3} |\eta_0| + \sqrt{3} \frac{\eta_0}{|\eta_0|} \sigma_0 \right)
\end{align*} }
and so we have a control over terms with a factor $O(|\xi| - |\eta|)$. Finally, we have as before that 
\begin{align*}
3 \xi_0^2 - |\xi|^2 &= O\left( \overline{\eta}_b^{\overline{\eta}} \right) + O\left( \rho_0 + \sigma_0 \right) + O\left( |\xi| - |\eta| \right) 
\end{align*}
which concludes. 

\paragraph{2.} Let us now localise to have $|\overline{\eta}| \simeq |\overline{\eta}| \gg |\overline{\rho}|$. Then \eqref{qteslemestimeesgeneriquesdeccubiquehb2} is already of the form \eqref{lemestimeesgeneriquesdeccubiquehb2-eta0sigmaresxrho}. 

On the other hand, $1 = O\left( \widehat{X}_a(\overline{\eta}) \cdot \nabla_{\overline{\eta}} \varphi_3 \right)$ so we can easily rewrite \eqref{equdecchampbCPP-1} as terms from \eqref{qteslemestimeesgeneriquesdeccubiquehb2}. 

\paragraph{3.} Let us now localise to have $|\overline{\sigma}| \simeq |\overline{\rho}| \gg |\overline{\eta}|$. 

Then we can decompote
\begin{align*}
\widehat{X}_b(\overline{\xi}) \cdot \nabla_{\overline{\eta}} &= \frac{|\xi| |\overline{\eta}|}{|\overline{\xi}| \eta_0} \widehat{X}_a(\overline{\eta}) \cdot \nabla_{\overline{\eta}} + O(\nabla_{\eta}) 
\end{align*}
to rewrite \eqref{equdecchampbCPP-2} as terms from \eqref{qteslemestimeesgeneriquesdeccubiquehb2} plus a term of the form 
\begin{align*}
\int_0^t \int \int e^{i s \varphi_3} s^2 \xi_0 |\overline{\sigma}|^3 \mu \widehat{f}(s, \overline{\eta}) \widehat{f}(s, \overline{\sigma}) \widehat{f}(s, \overline{\rho}) ~ d\overline{\eta} d\overline{\sigma} ds 
\end{align*}
\eqref{equdecchampbCPP-1} also has this form. But 
\begin{align*}
\xi_0 &= O(\rho_0) + O(\sigma_0) + O(\overline{\eta}) 
\end{align*}
We can then apply an integration by parts along $\widehat{X}_a(\overline{\eta}) \cdot \nabla_{\overline{\eta}}$ for the factors $O(\rho_0)$ and $O(\overline{\eta})$, and along $\widehat{X}_a(\overline{\eta}) \cdot \left( \nabla_{\overline{\eta}} - \nabla_{\overline{\sigma}} \right)$ for the factor $O(\sigma_0)$, and we only recover terms from \eqref{qteslemestimeesgeneriquesdeccubiquehb2}. 

\paragraph{4.} Let us now localise to have $|\overline{\eta}| \gg |\overline{\sigma}|+|\overline{\rho}|$. Again, we can rewrite \eqref{equdecchampbCPP-2} as terms of \eqref{qteslemestimeesgeneriquesdeccubiquehb2} plus a term of the form 
\begin{align*}
\int_0^t \int \int e^{i s \varphi_3} s^2 |\overline{\eta}|^4 \mu \widehat{f}(s, \overline{\eta}) \widehat{f}(s, \overline{\sigma}) \widehat{f}(s, \overline{\rho}) ~ d\overline{\eta} d\overline{\sigma} ds 
\end{align*}
where $\mu$ allows to distribute a derivative. \eqref{equdecchampbCPP-1} also has this form. 

We then have that 
\begin{align*}
1 = O\left( \widehat{X}_a(\overline{\eta}) \cdot \nabla_{\overline{\eta}} \varphi_3 \right), \quad 1 = O\left( \widehat{X}_a(\overline{\eta}) \cdot \left( \nabla_{\overline{\eta}} - \nabla_{\overline{\sigma}} \right) \varphi_3 \right)
\end{align*}
If we have a factor $O(\overline{\sigma})$, we can apply an integration by parts along $\widehat{X}_a(\overline{\eta}) \cdot \nabla_{\overline{\eta}}$, and along $\widehat{X}_a(\overline{\eta}) \cdot \left( \nabla_{\overline{\eta}} - \nabla_{\overline{\sigma}} \right)$ if we have a factor $O(\overline{\rho})$, to recover terms from \eqref{qteslemestimeesgeneriquesdeccubiquehb2}. But $\mu$ exactly gives such factors. 

\paragraph{5.} Let us finally localise to have $|\overline{\sigma}| \gg |\overline{\eta}|+|\overline{\rho}|$. 

Again, we rewrite \eqref{equdecchampbCPP-2} as terms of \eqref{qteslemestimeesgeneriquesdeccubiquehb2} plus a term of the form
\begin{align*}
\int_0^t \int \int e^{i s \varphi_3} s^2 \xi_0 |\overline{\sigma}|^2 \mu \widehat{F}_1(\overline{\eta}) \widehat{F}_2(\overline{\sigma}) \widehat{F}_3(\overline{\rho}) ~ d\overline{\eta} d\overline{\sigma} ds 
\end{align*}
\eqref{equdecchampbCPP-1} also has this form. 

If $F_2 = \partial_x f$, $\mu$ allows to distribute a derivative and we get a symbol of the form $O\left( \xi_0 \sigma_0 (|\overline{\eta}|+|\overline{\rho}|) \right) = O(\sigma_0 \overline{\eta}) + O(\sigma_0 \overline{\rho})$. 

If $F_1 = \partial_x f$ or $F_3 = \partial_x f$, we get a symbol of the form $O(\xi_0 \overline{\eta}) + O(\xi_0 \rho_0) = O((\xi_0 - \sigma_0) |\overline{\eta}|) + O((\xi_0 - \sigma_0) |\overline{\rho}|) + O(\sigma_0 \overline{\eta}) + O(\sigma_0 \overline{\rho})$. 

In any case, it is enough to treat the factors $O((\xi_0-\sigma_0)|\overline{\eta}|)$, $O((\xi_0-\sigma_0)|\overline{\rho}|)$, $O(\sigma_0 \overline{\eta})$ or $O(\sigma_0 \overline{\rho})$. For the factors $O(\sigma_0 \overline{\eta})$ and $O(\sigma_0 \overline{\rho})$, we can apply an integration by parts along $\widehat{X}_a(\overline{\eta}) \cdot \left( \nabla_{\overline{\eta}} - \nabla_{\overline{\sigma}} \right)$ using that 
\begin{align*}
1 = O\left( \widehat{X}_a(\overline{\eta}) \cdot \left( \nabla_{\overline{\eta}} - \nabla_{\overline{\sigma}} \right) \varphi_3 \right) 
\end{align*}
and we recover indeed terms of \eqref{qteslemestimeesgeneriquesdeccubiquehb2}. Moreover, a factor $O(\overline{\eta} \overline{\rho})$ also allows for an integration by parts in this same direction and brings only terms of \eqref{qteslemestimeesgeneriquesdeccubiquehb2}. 

On the other hand, 
\begin{align*}
\varphi_3 &= \xi_0^3 - \eta_0^3 - \sigma_0^3 - \rho_0^3 + \xi_0 |\xi|^2 - \eta_0 |\eta|^2 - \sigma_0 |\sigma|^2 - \rho_0 |\rho|^2 \\
&= O(\sigma_0 \overline{\eta}) + O(\sigma_0 \overline{\rho}) + O(\overline{\eta} \overline{\rho}) + \xi_0 \left( |\eta|^2 + |\sigma|^2 + |\rho|^2 + 2 \sigma \cdot (\eta + \rho) \right) - \eta_0 |\eta|^2 - \sigma_0 |\sigma|^2 - \rho_0 |\rho|^2 \\
&= O(\sigma_0 \overline{\eta}) + O(\sigma_0 \overline{\rho}) + O(\overline{\eta} \overline{\rho}) + (\xi_0 - \sigma_0) \left( |\sigma|^2 + 2 \sigma \cdot (\eta + \rho) \right) 
\end{align*}
so that 
\begin{align*}
(\xi_0 - \sigma_0) |\overline{\eta}| &= O(|\overline{\eta}| \varphi_3) + O(\sigma_0 \overline{\eta}) + O(\sigma_0 \overline{\rho}) + O(\overline{\eta} \overline{\rho}) \\
(\xi_0 - \sigma_0) |\overline{\rho}| &= O(|\overline{\rho}| \varphi_3) + O(\sigma_0 \overline{\eta}) + O(\sigma_0 \overline{\rho}) + O(\overline{\eta} \overline{\rho})
\end{align*}
which concludes. 

\subsubsection{Interaction \texorpdfstring{$\widehat{\mathcal{L}}\widehat{\mathcal{L}}\widehat{\mathcal{L}}$}{LLL}}

Let us consider: {\footnotesize 
\begin{subequations}
\begin{align}
&\xi_0 m_b(\overline{\xi}) \widehat{X}_b(\overline{\xi}) \cdot \nabla_{\overline{\xi}} \widehat{I}_{s\mu}^{\widehat{\mathcal{L}}\widehat{\mathcal{L}}\widehat{\mathcal{L}}}[F_1, F_2, F_3](t, \overline{\xi}) \notag \\
&\quad = \int_0^t \int \int i s^2 \xi_0 m_b(\overline{\xi}) \widehat{X}_b(\overline{\xi}) \cdot \nabla_{\overline{\xi}} \varphi_3 e^{i s \varphi_3} \mu(\overline{\xi}, \overline{\eta}, \overline{\sigma}) m_{\widehat{\mathcal{L}}}(\overline{\eta}) m_{\widehat{\mathcal{L}}}(\overline{\sigma}) m_{\widehat{\mathcal{L}}}(\overline{\rho}) \widehat{F}_1(s, \overline{\eta}) \widehat{F}_2(s, \overline{\sigma}) \widehat{F}_3(s, \overline{\rho}) ~ d\overline{\eta} d\overline{\sigma} ds \label{equdecchampbLLL-1} \\
&\quad \quad + \int_0^t \int \int e^{i s \varphi_3} s \xi_0 \mu(\overline{\xi}, \overline{\eta}, \overline{\sigma}) m_{\widehat{\mathcal{L}}}(\overline{\eta}) m_{\widehat{\mathcal{L}}}(\overline{\sigma}) m_{\widehat{\mathcal{L}}}(\overline{\rho}) \widehat{F}_1(s, \overline{\eta}) \widehat{F}_2(s, \overline{\sigma}) m_b(\overline{\xi}) \widehat{X}_b(\overline{\xi}) \cdot \nabla_{\overline{\xi}} \widehat{F}_3(s, \overline{\rho}) ~ d\overline{\eta} d\overline{\sigma} ds \label{equdecchampbLLL-2} \\
&\quad \quad + \int_0^t \int \int e^{i s \varphi_3} s \xi_0 m_b(\overline{\xi}) \widehat{X}_b(\overline{\xi}) \cdot \nabla_{\overline{\xi}} \left( \chi(\overline{\xi}, \overline{\eta}, \overline{\sigma}) m_{\widehat{\mathcal{L}}}(\overline{\eta}) m_{\widehat{\mathcal{L}}}(\overline{\sigma}) m_{\widehat{\mathcal{L}}}(\overline{\rho}) \right) \widehat{F}_1(s, \overline{\eta}) \widehat{F}_2(s, \overline{\sigma}) \widehat{F}_3(s, \overline{\rho}) ~ d\overline{\eta} d\overline{\sigma} ds \label{equdecchampbLLL-3} 
\end{align}
\end{subequations} }
\eqref{equdecchampbLLL-3} is of the form $\eqref{lemestimeesgeneriquesdeccubiquehb2-dersymb}+\eqref{lemestimeesgeneriquesdeccubiquehb2-dersymbbis}$. 

Without loss of generality, we need only to consider three cases for the relative sizes of frequencies: $|\overline{\eta}| \simeq |\overline{\sigma}| \simeq |\overline{\rho}|$, $|\overline{\eta}| \simeq |\overline{\sigma}| \gg |\overline{\rho}|$ or $|\overline{\eta}| \gg |\overline{\sigma}|+|\overline{\rho}|$ (by symmetry of the frequencies in this case). In the first two cases, \eqref{equdecchampbLLL-2} is of the form \eqref{lemestimeesgeneriquesdeccubiquehb2-eta0sigmaresxrho}. In the last case, we can distribute a derivative using $\mu$, and we get terms of \eqref{qteslemestimeesgeneriquesdeccubiquehb2} except if we obtain an additional factor $|\overline{\rho}|$, but then we can apply an integration by parts in $\sigma$ to recover terms from \eqref{qteslemestimeesgeneriquesdeccubiquehb2} plus a term having the form
\begin{align*}
\int_0^t \int \int e^{i s \varphi_3} s^2 |\overline{\eta}|^3 |\overline{\rho}| \mu \widehat{f}(s, \overline{\eta}) \widehat{f}(s, \overline{\sigma}) \widehat{f}(s, \overline{\rho}) ~ d\overline{\eta} d\overline{\sigma} ds 
\end{align*}
However $\partial_{\eta_0} \varphi_3 \simeq |\overline{\eta}|^2$ in this case so we easily recover terms from \eqref{qteslemestimeesgeneriquesdeccubiquehb2}. 

It thus only remains \eqref{equdecchampbLLL-1}. We have that 
\begin{align*}
\partial_{\eta_0} \varphi_3 &= 3 \rho_0^2 + |\rho|^2 - 3 \eta_0^2 - |\eta|^2 
\end{align*}
so we may restrict to the neighborhood of $|\rho_0| = |\eta_0|$, and symmetrically of $|\rho_0| = |\sigma_0|$. We thus have $|\overline{\eta}| \simeq |\overline{\sigma}| \simeq |\overline{\rho}|$. 

If $\epsilon^{\overline{\eta} \overline{\sigma}} = \epsilon^{\overline{\eta} \overline{\rho}} = 1$, then $\xi_0 = 3 \eta_0 + o(1)$, $\eta_0 = \sigma_0 + o(1) = \rho_0 + o(1)$, so 
\begin{align*}
\varphi_3 &= \xi_0^3 - \eta_0^3 - \sigma_0^3 - \rho_0^3 + o(1) \\
&= 24 \eta_0^3 + o(1) 
\end{align*}
and $1 = O(\varphi_3)$. Otherwise, up to exchanging $\overline{\eta}, \overline{\sigma}, \overline{\rho}$, we may assume that $\epsilon^{\overline{\eta} \overline{\sigma}} = 1$, $\epsilon^{\overline{\sigma} \overline{\rho}} = \epsilon^{\overline{\eta} \overline{\rho}} = -1$. In particular, $|\overline{\eta}| \simeq |\overline{\sigma}| \simeq |\overline{\rho}| \simeq |\overline{\xi}|$. We then have 
\begin{align*}
\eta &= \frac{\rho_0}{\eta_0} \rho + O(\nabla_{\eta} \varphi_3) \\
\partial_{\eta_0} \varphi_3 &= 3 \rho_0^2 + |\rho|^2 - 3 \eta_0^3 - |\eta|^2 \\
&= 3 \left( \rho_0^2 - \eta_0^2 \right) + |\rho|^2 \left( 1 - \frac{\rho_0^2}{\eta_0^2} \right) + O(\nabla_{\eta} \varphi_3) \\
&= \left( \rho_0^2 - \eta_0^2 \right) \left( 3 - \frac{|\rho|^2}{\eta_0^2} \right) + O(\nabla_{\eta} \varphi_3) 
\end{align*}
Therefore $\rho_0 + \eta_0 = O(\partial_{\eta_0} \varphi_3) + O(\nabla_{\eta} \varphi_3)$. In particular, 
\begin{align*}
\eta &= \frac{\rho_0}{\eta_0} \rho + O(\nabla_{\eta} \varphi_3) \\
&= - \rho + O(\nabla_{\eta} \varphi_3) + O(\partial_{\eta_0} \varphi_3) 
\end{align*}
Likewise, we show 
\begin{align*}
\rho_0 + \sigma_0 &= O(\partial_{\sigma_0} \varphi_3) + O(\nabla_{\sigma} \varphi_3) \\
\rho + \sigma &= O(\partial_{\sigma_0} \varphi_3) + O(\nabla_{\sigma} \varphi_3) 
\end{align*}
We deduce that 
\begin{align*}
\xi_0 + \rho_0 &= O(\partial_{\eta_0} \varphi_3) + O(\nabla_{\eta} \varphi_3) + O(\partial_{\sigma_0} \varphi_3) + O(\nabla_{\sigma} \varphi_3) \\
\xi + \rho &= O(\partial_{\eta_0} \varphi_3) + O(\nabla_{\eta} \varphi_3) + O(\partial_{\sigma_0} \varphi_3) + O(\nabla_{\sigma} \varphi_3)
\end{align*}
and thus 
\begin{align*}
\nabla_{\overline{\xi}} \varphi_3 &= O(\partial_{\eta_0} \varphi_3) + O(\nabla_{\eta} \varphi_3) + O(\partial_{\sigma_0} \varphi_3) + O(\nabla_{\sigma} \varphi_3)
\end{align*}
We can then apply integrations by parts and get terms from \eqref{qteslemestimeesgeneriquesdeccubiquehb2}. 

\subsubsection{Interaction \texorpdfstring{$\widehat{\mathcal{L}}\widehat{\mathcal{L}}\widehat{\mathcal{P}}$}{LLP}}

Let us consider: {\footnotesize 
\begin{subequations}
\begin{align}
&\xi_0 m_b(\overline{\xi}) \widehat{X}_b(\overline{\xi}) \cdot \nabla_{\overline{\xi}} \widehat{I}_{s \mu}^{\widehat{\mathcal{L}}\widehat{\mathcal{L}}\widehat{\mathcal{P}}}[F_1, F_2, F_3](t, \overline{\xi}) \notag \\
&\quad = \int_0^t \int \int i s^2 \xi_0 m_b(\overline{\xi}) \widehat{X}_b(\overline{\xi}) \cdot \nabla_{\overline{\xi}} \varphi_3 e^{i s \varphi_3} \mu(\overline{\xi}, \overline{\eta}, \overline{\sigma}) m_{\widehat{\mathcal{L}}}(\overline{\eta}) m_{\widehat{\mathcal{L}}}(\overline{\sigma}) m_{\widehat{\mathcal{P}}}(\overline{\rho}) \widehat{F}_1(s, \overline{\eta}) \widehat{F}_2(s, \overline{\sigma}) \widehat{F}_3(s, \overline{\rho}) ~ d\overline{\eta} d\overline{\sigma} ds \label{equdecchampbLLP-1} \\
&\quad \quad + \int_0^t \int \int e^{i s \varphi_3} s \xi_0 \mu(\overline{\xi}, \overline{\eta}, \overline{\sigma}) m_{\widehat{\mathcal{L}}}(\overline{\eta}) m_{\widehat{\mathcal{L}}}(\overline{\sigma}) m_{\widehat{\mathcal{P}}}(\overline{\rho}) \widehat{F}_1(s, \overline{\eta}) \widehat{F}_2(s, \overline{\sigma}) m_b(\overline{\xi}) \widehat{X}_b(\overline{\xi}) \cdot \nabla_{\overline{\xi}} \widehat{F}_3(s, \overline{\rho}) ~ d\overline{\eta} d\overline{\sigma} ds \label{equdecchampbLLP-2} \\
&\quad \quad + \int_0^t \int \int e^{i s \varphi_3} s \xi_0 m_b(\overline{\xi}) \widehat{X}_b(\overline{\xi}) \cdot \nabla_{\overline{\xi}} \left( \mu(\overline{\xi}, \overline{\eta}, \overline{\sigma}) m_{\widehat{\mathcal{L}}}(\overline{\eta}) m_{\widehat{\mathcal{L}}}(\overline{\sigma}) m_{\widehat{\mathcal{P}}}(\overline{\rho}) \right) \widehat{F}_1(s, \overline{\eta}) \widehat{F}_2(s, \overline{\sigma}) \widehat{F}_3(s, \overline{\rho}) ~ d\overline{\eta} d\overline{\sigma} ds \label{equdecchampbLLP-3} 
\end{align}
\end{subequations} }
\eqref{equdecchampbLLP-3} is of the form $\eqref{lemestimeesgeneriquesdeccubiquehb2-dersymb}+\eqref{lemestimeesgeneriquesdeccubiquehb2-dersymbbis}$. 

We now consider subcases depending on the relative sizes of frequencies, noting that we may symmetrize $\overline{\eta}$ and $\overline{\sigma}$ here. 

\paragraph{1.} Let us first localise to have $|\overline{\eta}| \simeq |\overline{\sigma}| \gtrsim |\overline{\rho}|$. Then \eqref{equdecchampbLLP-2} is of the form \eqref{lemestimeesgeneriquesdeccubiquehb2-eta0sigmaresxrho}. 

For \eqref{equdecchampbLLP-1}, we have 
\begin{align*}
\partial_{\eta_0} \varphi_3 &= o(1) + |\rho|^2 - 3 \eta_0^2 \\
\partial_{\sigma_0} \varphi_3 &= o(1) + |\rho|^2 - 3 \sigma_0^2 
\end{align*}
so we can restrict to the neighborhood of $|\rho| = \sqrt{3} |\eta_0| = \sqrt{3} |\sigma_0|$. 

If $\epsilon^{\overline{\eta} \overline{\sigma}} = 1$, then $\xi_0 = 2 \eta_0 + o(1)$ and 
\begin{align*}
\varphi_3 &= \xi_0^3 + \xi_0 |\xi|^2 - \eta_0^3 - \eta_0 |\eta|^2 - \sigma_0^3 - \sigma_0 |\sigma|^2 - \rho_0^3 - \rho_0 |\rho|^2 \\
&= 8 \eta_0^3 + 6 \eta_0^3 - \eta_0^3 - \eta_0^3 + o(1) \\
&= 12 \eta_0^3 + o(1) 
\end{align*}
so $1 = O(\varphi_3)$. 

Assume now $\epsilon^{\overline{\eta} \overline{\sigma}} = -1$. In this case, by similar computations as for the interaction $\widehat{\mathcal{L}}\widehat{\mathcal{L}}\widehat{\mathcal{L}}$ above, we have that 
\begin{align*}
\sigma_0 + \eta_0 &= O\left( \left( \nabla_{\overline{\eta}} - \partial_{\overline{\sigma}} \right) \varphi_3 \right) \\
\sigma + \eta &= O\left( \left( \nabla_{\overline{\eta}} - \partial_{\overline{\sigma}} \right) \varphi_3 \right)
\end{align*}
In particular, we can apply integrations by parts on terms having a factor $O(\rho_0 (\eta + \sigma))$ or $O(\rho_0 (\eta_0 + \sigma_0))$ and get terms from \eqref{qteslemestimeesgeneriquesdeccubiquehb2}. On the other hand, 
\begin{align*}
|\eta|^2 - |\rho|^2 &= 3 \rho_0^2 - 3 \eta_0^2 + O(\partial_{\eta_0} \varphi_3) \\
|\sigma|^2 - |\rho|^2 &= 3 \rho_0^2 - 3 \sigma_0^2 + O(\partial_{\sigma_0} \varphi_3) \\
\rho_0 |\xi|^2 &= \rho_0 |\eta + \sigma + \rho|^2 \\
&= \rho_0 |\rho|^2 + O(\rho_0 (\eta + \sigma)) 
\end{align*}
Therefore, {\footnotesize 
\begin{align*}
\varphi_3 &= \xi_0^3 + \xi_0 |\xi|^2 - \eta_0^3 - \eta_0 |\eta|^2 - \sigma_0^3 - \sigma_0 |\sigma|^2 - \rho_0^3 - \rho_0 |\rho|^2 \\
&= (\eta_0 + \sigma_0)^3 + \rho_0^3 + (\eta_0 + \sigma_0) |\xi|^2 + \rho_0 |\xi|^2 - (\eta_0 + \sigma_0) (\eta_0^2 - \eta_0 \sigma_0 + \sigma_0^2) - (\eta_0 + \sigma_0) |\eta|^2 \\
&\quad + \sigma_0 (|\eta|^2 - |\sigma|^2) - \rho_0^3 - \rho_0 |\rho|^2 + O(\rho_0 (\eta_0 + \sigma_0)) \\
&= (\eta_0 + \sigma_0) \left( (\eta_0 + \sigma_0)^2 + |\xi|^2 - \eta_0^2 + \eta_0 \sigma_0 - \sigma_0^2 - |\eta|^2 \right) + \sigma_0 (3 \sigma_0^2 - 3 \eta_0^2) \\
&\quad + O(\rho_0 (\eta_0 + \sigma_0)) + O(\rho_0 (\eta + \sigma)) + O(\partial_{\eta_0} \varphi_3) + O(\partial_{\sigma_0} \varphi_3) \\
&= (\eta_0 + \sigma_0) \left( (\eta_0 + \sigma_0)^2 + |\xi|^2 - \eta_0^2 + \eta_0 \sigma_0 - \sigma_0^2 - |\eta|^2 + 3 \sigma_0 (\sigma_0 - \eta_0) \right) \\
&\quad + O(\rho_0 (\eta_0 + \sigma_0)) + O(\rho_0 (\eta + \sigma)) + O(\partial_{\eta_0} \varphi_3) + O(\partial_{\sigma_0} \varphi_3)
\end{align*} }
But 
\begin{align*}
(\eta_0 + \sigma_0)^2 + |\xi|^2 - \eta_0^2 + \eta_0 \sigma_0 - \sigma_0^2 - |\eta|^2 + 3 \sigma_0 (\sigma_0 - \eta_0) 
&= o(1) + 3 \eta_0^2 
\end{align*}
We deduce that 
\begin{align*}
\eta_0 + \sigma_0 &= O(\varphi_3) + O(\rho_0 (\eta_0 + \sigma_0)) + O(\rho_0 (\eta + \sigma)) + O(\partial_{\eta_0} \varphi_3) + O(\partial_{\sigma_0} \varphi_3)
\end{align*}
Note that $\xi_0 - \rho_0 = \eta_0 + \sigma_0$. 

Now we only need to compute 
\begin{align*}
\xi_0 \widehat{X}_b(\overline{\xi}) \cdot \nabla_{\overline{\xi}} \varphi_3 &= \xi_0 \frac{|\xi|}{|\overline{\xi}|} \left( 3 \xi_0^2 + |\xi|^2 - 3 \rho_0^2 - |\rho|^2 \right) - \xi_0 \frac{\xi_0 \xi}{|\overline{\xi}| |\xi|} \cdot \left( 2 \xi_0 \xi - 2 \rho_0 \rho \right) \\
&= O(\xi_0 - \rho_0) + O(\rho_0 (\xi - \rho)) \\
&= O(\xi_0 - \rho_0) + O(\rho_0 (\eta + \sigma)) 
\end{align*}
which concludes. 

\paragraph{2.} Let us now localise to have $|\overline{\eta}| \simeq |\overline{\rho}| \gg |\overline{\sigma}|$, so we can as for previous interactions decompose and apply an integration by parts on \eqref{equdecchampbLLP-2} to rewrite it as terms from \eqref{qteslemestimeesgeneriquesdeccubiquehb2} plus a term having the form
\begin{align*}
\int_0^t \int \int e^{i s \varphi_3} s^2 |\overline{\eta}|^4 \mu \widehat{f}(s, \overline{\eta}) \widehat{f}(s, \overline{\sigma}) \widehat{f}(s, \overline{\rho}) ~ d\overline{\eta} d\overline{\sigma} ds 
\end{align*}
and \eqref{equdecchampbLLP-1} is also of this form. But here
\begin{align*}
1 = O\left( (\partial_{\eta_0}-\partial_{\sigma_0}) \varphi_3 \right) 
\end{align*}
so we can easily recover only terms of \eqref{qteslemestimeesgeneriquesdeccubiquehb2}. 

\paragraph{3.} Let us now localise to have $|\overline{\eta}| \gg |\overline{\sigma}| + |\overline{\rho}|$, then as before we can rewrite \eqref{equdecchampbLLP-2} as terms of \eqref{qteslemestimeesgeneriquesdeccubiquehb2} plus a term of the form 
\begin{align*}
\int_0^t \int \int e^{i s \varphi_3} s^2 |\overline{\eta}|^4 \mu \widehat{f}(s, \overline{\eta}) \widehat{f}(s, \overline{\sigma}) \widehat{f}(s, \overline{\rho}) ~ d\overline{\eta} d\overline{\sigma} ds 
\end{align*}
where $\mu$ allows to distribute a derivative. \eqref{equdecchampbLLP-1} also has this form. 

In particular, the symbol $\mu$ gives an additional factor $O\left( \overline{\sigma} \right)$ or $O\left( \overline{\rho} \right)$. In both cases, we can apply an integration by parts along $\left( \partial_{\eta_0} - \partial_{\sigma_0} \right)$, since 
\begin{align*}
1 = O\left( \left( \partial_{\eta_0} - \partial_{\sigma_0} \right) \varphi_3 \right) 
\end{align*}
and we get terms of \eqref{qteslemestimeesgeneriquesdeccubiquehb2}. 

\paragraph{4.} Finally, let us localise to have $|\overline{\rho}| \gg |\overline{\eta}|+|\overline{\sigma}|$. Again, we rewrite \eqref{equdecchampbLLP-2} as terms of \eqref{qteslemestimeesgeneriquesdeccubiquehb2} plus a term of the form 
\begin{align*}
\int_0^t \int \int e^{i s \varphi_3} s^2 \xi_0 |\overline{\rho}|^3 \mu \widehat{f}(s, \overline{\eta}) \widehat{f}(s, \overline{\sigma}) \widehat{f}(s, \overline{\rho}) ~ d\overline{\eta} d\overline{\sigma} ds
\end{align*}
where $\mu$ allows to distribute a derivative. \eqref{equdecchampbLLP-1} has the same form. 

By the same compotuation as for the interaction $\widehat{\mathcal{C}}\widehat{\mathcal{C}}\widehat{\mathcal{P}}$, we have
\begin{align*}
\xi_0 - \rho_0 &= O(\varphi_3) + O\left( \rho_0 \overline{\eta} \right) + O\left( \rho_0 \overline{\sigma} \right) + O\left( \overline{\eta} \overline{\sigma} \right) 
\end{align*}
and therefore 
\begin{align*}
\xi_0 \mu &= ( \xi_0 - \rho_0 ) \mu + \rho_0 \mu \\
&= O(\mu \varphi_3) + O\left( \rho_0 \overline{\eta} \right) + O\left( \rho_0 \overline{\sigma} \right) + O\left( \overline{\eta} \overline{\sigma} \right)
\end{align*}
We also have automatically $1 = O(\partial_{\eta_0} \varphi_3)$. 
We can then apply an integration by parts in time for the factor $O(\mu \varphi_3)$ and distribute derivatives, along $\partial_{\eta_0}$ for the factors $O\left( \rho_0 \overline{\eta} \right)$, $O\left( \rho_0 \overline{\sigma} \right)$, $O\left( \overline{\eta} \overline{\sigma} \right)$, and we recover terms from \eqref{qteslemestimeesgeneriquesdeccubiquehb2}. 

\subsubsection{Interaction \texorpdfstring{$\widehat{\mathcal{L}}\widehat{\mathcal{P}}\widehat{\mathcal{P}}$}{LPP}}

Let us consider: {\footnotesize 
\begin{subequations}
\begin{align}
&\xi_0 m_b(\overline{\xi}) \widehat{X}_b(\overline{\xi}) \cdot \nabla_{\overline{\xi}} \widehat{I}_{s \mu}^{\widehat{\mathcal{L}}\widehat{\mathcal{P}}\widehat{\mathcal{P}}}[F_1, F_2, F_3](t, \overline{\xi}) \notag \\
&\quad = \int_0^t \int \int i s^2 \xi_0 m_b(\overline{\xi}) \widehat{X}_b(\overline{\xi}) \cdot \nabla_{\overline{\xi}} \varphi_3 e^{i s \varphi_3} \mu(\overline{\xi}, \overline{\eta}, \overline{\sigma}) m_{\widehat{\mathcal{L}}}(\overline{\eta}) m_{\widehat{\mathcal{P}}}(\overline{\sigma}) m_{\widehat{\mathcal{P}}}(\overline{\rho}) \widehat{F}_1(s, \overline{\eta}) \widehat{F}_2(s, \overline{\sigma}) \widehat{F}_3(s, \overline{\rho}) ~ d\overline{\eta} d\overline{\sigma} ds \label{equdecchampbLPP-1} \\
&\quad \quad + \int_0^t \int \int e^{i s \varphi_3} s \xi_0 \mu(\overline{\xi}, \overline{\eta}, \overline{\sigma}) m_{\widehat{\mathcal{L}}}(\overline{\eta}) m_{\widehat{\mathcal{P}}}(\overline{\sigma}) m_{\widehat{\mathcal{P}}}(\overline{\rho}) \widehat{F}_1(s, \overline{\eta}) \widehat{F}_2(s, \overline{\sigma}) m_b(\overline{\xi}) \widehat{X}_b(\overline{\xi}) \cdot \nabla_{\overline{\xi}} \widehat{F}_3(s, \overline{\rho}) ~ d\overline{\eta} d\overline{\sigma} ds \label{equdecchampbLPP-2} \\
&\quad \quad + \int_0^t \int \int e^{i s \varphi_3} s \xi_0 m_b(\overline{\xi}) \widehat{X}_b(\overline{\xi}) \cdot \nabla_{\overline{\xi}} \left( \mu(\overline{\xi}, \overline{\eta}, \overline{\sigma}) m_{\widehat{\mathcal{L}}}(\overline{\eta}) m_{\widehat{\mathcal{P}}}(\overline{\sigma}) m_{\widehat{\mathcal{P}}}(\overline{\rho}) \right) \widehat{F}_1(s, \overline{\eta}) \widehat{F}_2(s, \overline{\sigma}) \widehat{F}_3(s, \overline{\rho}) ~ d\overline{\eta} d\overline{\sigma} ds \label{equdecchampbLPP-3} 
\end{align}
\end{subequations} }
\eqref{equdecchampbLPP-3} has the form $\eqref{lemestimeesgeneriquesdeccubiquehb2-dersymb}+\eqref{lemestimeesgeneriquesdeccubiquehb2-dersymbbis}$. 

We separate into subcases depending on the relative sizes of frequencies. Without loss of generality, we can assume that $|\overline{\rho}| \lesssim |\overline{\sigma}|$ ou $|\overline{\eta}| \gg |\overline{\sigma}|+|\overline{\rho}|$. 

\paragraph{1.} Let us first localise to have $|\overline{\eta}| \simeq |\overline{\sigma}| \simeq |\overline{\rho}|$. Then \eqref{equdecchampbLPP-2} already has the form \eqref{lemestimeesgeneriquesdeccubiquehb2-eta0sigmaresxrho}. 

Then, for \eqref{equdecchampbLPP-1}, 
\begin{align*}
\partial_{\eta_0} \varphi_3 &= o(1) + |\rho|^2 - 3 \eta_0^2 \\
\partial_{\sigma_0} \varphi_3 &= o(1) + |\rho|^2 - |\sigma|^2 
\end{align*}
and so we only need to consider the neighborhood of $|\rho| = |\sigma| = \sqrt{3} |\eta_0|$. In particular, $|\overline{\eta}| \simeq |\overline{\sigma}| \simeq |\overline{\rho}| \simeq |\overline{\xi}|$. Then, 
\begin{align*}
\varphi_3 &= \xi_0^3 + \xi_0 |\xi|^2 - \eta_0^3 + o(1) \\
&= \xi_0 |\xi|^2 + o(1) 
\end{align*}
where $\xi_0 = \eta_0 + o(1)$, so $|\xi_0| \simeq |\overline{\xi}|$. We can therefore reduce to the neighborhood of $\xi = 0$. But since $\xi = \eta + \sigma + \rho$ with $\eta = o(1)$, $|\sigma| = |\rho| + o(1)$, this forces $\theta^{\overline{\sigma} \overline{\rho}} = -1 + o(1)$. 

Let us now compute factors that allow for good estimates. 
\begin{align*}
\frac{J \sigma}{|\sigma|} \cdot \nabla_{\sigma} \varphi_3 &= 2 \rho_0 \frac{J \sigma \cdot \rho}{|\sigma|} 
\end{align*}
so we can estimate any term having a factor $O\left( \rho_0 \xi_t^{\overline{\sigma} \overline{\rho}} \right)$, and symmetric-wise $O\left( \sigma_0 \xi_t^{\overline{\sigma} \overline{\rho}} \right)$. Then, 
\begin{align*}
\frac{\sigma}{|\sigma|} \cdot \nabla_{\sigma} \varphi_3 &= 2 \rho_0 \rho \cdot \sigma - 2 \sigma_0 |\sigma|^2 \\
&= - 2 \rho_0 |\rho| |\sigma| - 2 \sigma_0 |\sigma|^2 + O\left( \rho_0 \xi_t^{\overline{\sigma} \overline{\rho}} \right) 
\end{align*}
so we can estimate terms having a factor $O\left( \rho_0 |\rho| + \sigma_0 |\sigma| \right)$. Finally, 
\begin{align*}
\partial_{\sigma_0} \varphi_3 &= 3 \rho_0^2 + |\rho|^2 - 3 \sigma_0^2 - |\sigma|^2 \\
&= 3 \rho_0^2 + 2 \sqrt{3} \rho_0 |\rho| + |\rho|^2 - 3 \sigma_0^2 + 2 \sqrt{3} \sigma_0 |\sigma| - |\sigma|^2 + O\left( \rho_0 |\rho| + \sigma_0 |\sigma| \right) \\
&= \left( |\rho| + \sqrt{3} \rho_0 \right)^2 - \left( |\sigma| - \sqrt{3} \sigma_0 \right)^2 + O\left( \rho_0 |\rho| + \sigma_0 |\sigma| \right)
\end{align*}
so we can estimate any term having a factor $O\left( |\rho| - |\sigma| + \sqrt{3} (\rho_0 + \sigma_0) \right)$. Doing the same computation but changing the sign of $\rho_0 |\rho| - \sigma_0 |\sigma|$, we also get 
\begin{align*}
O\left( |\rho| - |\sigma| - \sqrt{3} (\rho_0 + \sigma_0) \right)
\end{align*} 
and therefore $O(|\rho| - |\sigma|)$ and $O(\rho_0 + \sigma_0) = O(\xi_0 - \eta_0)$. 

Note also that we can apply an integration by parts along $\eta \cdot \nabla_{\eta}$ and recover only terms from \eqref{qteslemestimeesgeneriquesdeccubiquehb2} since $\eta$ compensate the potential singular weight on $\widehat{f}(s, \overline{\eta})$, and 
\begin{align*}
2 \rho_0 \rho \cdot \eta &= 2 \eta_0 |\eta|^2 + O\left( \eta \cdot \nabla_{\eta} \varphi_3 \right) 
\end{align*}
We also have 
\begin{align*}
3 \rho_0^2 + |\rho|^2 &= 3 \eta_0^2 + |\eta|^2 + O\left( \partial_{\eta_0} \varphi_3 \right) 
\end{align*}
Finally, we compute that 
\begin{align*}
\varphi_3 &= \xi_0^3 + \xi_0 |\xi|^2 - \eta_0^3 - \eta_0 |\eta|^2 - \sigma_0^3 - \sigma_0 |\sigma|^2 - \rho_0^3 - \rho_0 |\rho|^2 \\
&= \xi_0 (|\xi|^2 - |\eta|^2) + O(\xi_0 - \eta_0) + O(\sigma_0 + \rho_0) + O(|\sigma| - |\rho|) 
\end{align*}
so we can estimate terms having a factor $O(|\xi|^2 - |\eta|^2)$. 

Now we only need to compute {\footnotesize 
\begin{align*}
&|\xi| \widehat{X}_b(\overline{\xi}) \cdot \nabla_{\overline{\xi}} \varphi_3 = \frac{|\xi|^2}{|\overline{\xi}|} \left( 3 \xi_0^2 + |\xi|^2 - 3 \rho_0^2 - |\rho|^2 \right) - \frac{\xi_0 \xi}{|\overline{\xi}|} \cdot \left( 2 \xi_0 \xi - 2 \rho_0 \rho \right) \\
&\quad = \frac{|\xi|^2}{|\overline{\xi}|} \left( 3 \xi_0^2 + |\xi|^2 - 3 \eta_0^2 - |\eta|^2 \right) - \frac{\xi_0}{|\overline{\xi}|} \left( 2 \xi_0 |\xi|^2 - 2 \rho_0 \rho \cdot \xi \right) + O\left( \partial_{\eta_0} \varphi_3 \right) \\
&\quad = - 2 \frac{\xi_0}{|\overline{\xi}|} \left( \xi_0 |\xi|^2 - \rho_0 \rho \cdot \eta - \rho_0 \rho \cdot \sigma - \rho_0 |\rho|^2 \right) + O\left( \partial_{\eta_0} \varphi_3 \right) + O(\xi_0 - \eta_0) + O(|\xi|^2 - |\eta|^2) \\
&\quad = - 2 \frac{\xi_0}{|\overline{\xi}|} \left( \xi_0 |\xi|^2 - \eta_0 |\eta|^2 + \rho_0 |\rho| |\sigma| - \rho_0 |\rho|^2 \right) + O\left( \partial_{\eta_0} \varphi_3 \right) + O(\xi_0 - \eta_0) + O(|\xi|^2 - |\eta|^2) \\
&\quad \quad + O\left( \eta \cdot \nabla_{\eta} \varphi_3 \right) + O\left( \rho_0 \xi_t^{\overline{\rho} \overline{\sigma}} \right) \\
&\quad = O\left( \partial_{\eta_0} \varphi_3 \right) + O(\xi_0 - \eta_0) + O(|\xi|^2 - |\eta|^2) + O\left( \eta \cdot \nabla_{\eta} \varphi_3 \right) + O\left( \rho_0 \xi_t^{\overline{\rho} \overline{\sigma}} \right) + O(|\sigma| - |\rho|) 
\end{align*} }
which concludes. 

\paragraph{2.} Let us localise to have $|\overline{\eta}| \simeq |\overline{\sigma}| \gg |\overline{\rho}|$. Then \eqref{equdecchampbLPP-2} is already of the form \eqref{lemestimeesgeneriquesdeccubiquehb2-eta0sigmaresxrho}. 

On the other hand, we have $1 = O(\partial_{\eta_0} \varphi_3)$, so we can easily rewrite \eqref{equdecchampbLPP-1} as terms of \eqref{qteslemestimeesgeneriquesdeccubiquehb2}. 

\paragraph{3.} Let us localise to have $|\overline{\sigma}| \simeq |\overline{\rho}| \gg |\overline{\eta}|$. 

If $F_1 = \partial_x f$, then \eqref{equdecchampbLPP-2} is already of the form \eqref{lemestimeesgeneriquesdeccubiquehb2-eta0sigmaresxrho}, and since $1 = O(\partial_{\eta_0} \varphi_3)$ we can easily rewrite \eqref{equdecchampbLPP-1} as terms of \eqref{qteslemestimeesgeneriquesdeccubiquehb2}. 

Else, by symmetry between $\overline{\sigma}$ and $\overline{\rho}$ we may assume that $F_3 = \partial_x f$, so \eqref{equdecchampbLPP-2} is of the form \eqref{lemestimeesgeneriquesdeccubiquehb2-rho0sigmaresxrho}, and again we can apply an integration by parts along $\partial_{\eta_0}$ on \eqref{equdecchampbLPP-1}. 

\paragraph{4.} Let us localise to have $|\overline{\eta}| \gg |\overline{\sigma}|+|\overline{\rho}|$. 

Then as for the interaction $\widehat{\mathcal{L}}\widehat{\mathcal{L}}\widehat{\mathcal{P}}$ above we can rewrite \eqref{equdecchampbLPP-2} as terms of \eqref{qteslemestimeesgeneriquesdeccubiquehb2} plus a term of the form 
\begin{align*}
\int_0^t \int \int e^{i s \varphi_3} s^2 |\overline{\eta}|^4 \mu \widehat{f}(s, \overline{\eta}) \widehat{f}(s, \overline{\sigma}) \widehat{f}(s, \overline{\rho}) ~ d\overline{\eta} d\overline{\sigma} ds 
\end{align*}
where $\mu$ allows to distribute a derivative. \eqref{equdecchampbLPP-1} also has this form. 

But 
\begin{align*}
1 = O\left( \partial_{\eta_0} \varphi_3 \right), \quad 1 = O\left( \left( \partial_{\eta_0} - \partial_{\sigma_0} \right) \varphi_3 \right) 
\end{align*}
so we can apply an integration by parts along $\partial_{\eta_0}$ if we have a factor $O(\overline{\sigma})$, and along $\left( \partial_{\eta_0} - \partial_{\sigma_0} \right)$ if we have a factor $O(\overline{\rho})$, which is always the case using $\mu$, so we recover in all cases terms of \eqref{qteslemestimeesgeneriquesdeccubiquehb2}. 

\paragraph{5.} Finally, let us localise to have $|\overline{\sigma}| \gg |\overline{\eta}|+|\overline{\rho}|$. We can rewrite \eqref{equdecchampbLPP-2} and \eqref{equdecchampbLPP-1}, up to terms of \eqref{qteslemestimeesgeneriquesdeccubiquehb2}, as 
\begin{align*}
\int_0^t \int \int e^{i s \varphi_3} s^2 \xi_0 |\overline{\sigma}|^3 \mu \widehat{f}(s, \overline{\eta}) \widehat{f}(s, \overline{\sigma}) \widehat{f}(s, \overline{\rho}) ~ d\overline{\eta} d\overline{\sigma} ds 
\end{align*}
where $\mu$ allows to distribute a derivative. But here $\xi_0 = \eta_0 + \sigma_0 + \rho_0$ and we have 
\begin{align*}
1 = O\left( \partial_{\sigma_0} \varphi_3 \right), \quad 1 = O\left( \left( \partial_{\eta_0} - \partial_{\sigma_0} \right) \varphi_3 \right) 
\end{align*}
Therefore, using $\mu$, we obtain either an additional factor $|\overline{\eta}|$ or $|\overline{\rho}|$. We can then apply an integration by parts along $\partial_{\sigma_0}$ if we have a factor $O(\eta_0)$ or $O(\sigma_0 |\overline{\eta}|)$ or $O(\rho_0)$, and along $\left( \partial_{\eta_0} - \partial_{\sigma_0} \right)$ for the factor $O(\sigma_0 |\overline{\rho}|)$. We then get only terms of \eqref{qteslemestimeesgeneriquesdeccubiquehb2}. 

\subsubsection{Interaction \texorpdfstring{$\widehat{\mathcal{P}}\widehat{\mathcal{P}}\widehat{\mathcal{P}}$}{PPP}}

Let us finally consider: {\footnotesize 
\begin{subequations}
\begin{align}
&\xi_0 m_b(\overline{\xi}) \widehat{X}_b(\overline{\xi}) \cdot \nabla_{\overline{\xi}} \widehat{I}_{s \mu}^{\widehat{\mathcal{P}}\widehat{\mathcal{P}}\widehat{\mathcal{P}}}[F_1, F_2, F_3](t, \overline{\xi}) \notag \\
&\quad = \int_0^t \int \int i s^2 \xi_0 m_b(\overline{\xi}) \widehat{X}_b(\overline{\xi}) \cdot \nabla_{\overline{\xi}} \varphi_3 e^{i s \varphi_3} \mu(\overline{\xi}, \overline{\eta}, \overline{\sigma}) m_{\widehat{\mathcal{P}}}(\overline{\eta}) m_{\widehat{\mathcal{P}}}(\overline{\sigma}) m_{\widehat{\mathcal{P}}}(\overline{\rho}) \widehat{F}_1(s, \overline{\eta}) \widehat{F}_2(s, \overline{\sigma}) \widehat{F}_3(s, \overline{\rho}) ~ d\overline{\eta} d\overline{\sigma} ds \label{equdecchampbPPP-1} \\
&\quad \quad + \int_0^t \int \int e^{i s \varphi_3} s \xi_0 \mu(\overline{\xi}, \overline{\eta}, \overline{\sigma}) m_{\widehat{\mathcal{P}}}(\overline{\eta}) m_{\widehat{\mathcal{P}}}(\overline{\sigma}) m_{\widehat{\mathcal{P}}}(\overline{\rho}) \widehat{F}_1(s, \overline{\eta}) \widehat{F}_2(s, \overline{\sigma}) m_b(\overline{\xi}) \widehat{X}_b(\overline{\xi}) \cdot \nabla_{\overline{\xi}} \widehat{F}_3(s, \overline{\rho}) ~ d\overline{\eta} d\overline{\sigma} ds \label{equdecchampbPPP-2} \\
&\quad \quad + \int_0^t \int \int e^{i s \varphi_3} s \xi_0 m_b(\overline{\xi}) \widehat{X}_b(\overline{\xi}) \cdot \nabla_{\overline{\xi}} \left( \mu(\overline{\xi}, \overline{\eta}, \overline{\sigma}) m_{\widehat{\mathcal{P}}}(\overline{\eta}) m_{\widehat{\mathcal{P}}}(\overline{\sigma}) m_{\widehat{\mathcal{P}}}(\overline{\rho}) \right) \widehat{F}_1(s, \overline{\eta}) \widehat{F}_2(s, \overline{\sigma}) \widehat{F}_3(s, \overline{\rho}) ~ d\overline{\eta} d\overline{\sigma} ds \label{equdecchampbPPP-3} 
\end{align}
\end{subequations} }
\eqref{equdecchampbPPP-3} is of the form $\eqref{lemestimeesgeneriquesdeccubiquehb2-dersymb}+\eqref{lemestimeesgeneriquesdeccubiquehb2-dersymbbis}$. 

In this case, by symmetry, we can always assume that $F_3 = \partial_x f$, $F_1 = F_2 = f$. In particular, \eqref{equdecchampbPPP-2} is of the form \eqref{lemestimeesgeneriquesdeccubiquehb2-rho0sigmaresxrho}. 

For \eqref{equdecchampbPPP-1}, we have 
\begin{align*}
\partial_{\eta_0} \varphi_3 &= o(1) + |\rho|^2 - |\eta|^2 \\
\partial_{\sigma_0} \varphi_3 &= o(1) + |\rho|^2 - |\sigma|^2 
\end{align*}
In particular, 
\begin{itemize}
\item if $|\overline{\rho}| \simeq |\overline{\eta}| \gg |\overline{\sigma}|$, we can apply an integration by parts along $\partial_{\sigma_0}$ and (since $F_3 = \partial_x f$) we recover terms from \eqref{qteslemestimeesgeneriquesdeccubiquehb2}; symmetric-wise if $|\overline{\rho}| \simeq |\overline{\sigma}| \gg |\overline{\eta}|$; 
\item if $|\overline{\eta}| \simeq |\overline{\sigma}| \gg |\overline{\rho}|$, we can apply an integration by parts along $\partial_{\eta_0}$ and (since $F_3 = \partial_x f$) we recover terms from \eqref{qteslemestimeesgeneriquesdeccubiquehb2}; 
\item if $|\overline{\rho}| \gg |\overline{\sigma}|+|\overline{\eta}|$, $\mu$ allows to gain a factor $O(\overline{\eta})$ or $O(\overline{\sigma})$, and we can apply an integration by parts along $\partial_{\sigma_0}$, respectively $\partial_{\eta_0}$, to recover terms from \eqref{qteslemestimeesgeneriquesdeccubiquehb2}. 
\end{itemize}
Therefore, only two cases remain: $|\overline{\eta}| \simeq |\overline{\sigma}| \simeq |\overline{\rho}|$, or $|\overline{\eta}| \gg |\overline{\sigma}|+|\overline{\rho}|$. 

\paragraph{1.} If we localise to have $|\overline{\eta}| \simeq |\overline{\sigma}| \simeq |\overline{\rho}|$, the previous computation ensures that we only need to consider the neighborhood of $|\rho| = |\eta| = |\sigma|$. 

We can then compute that 
\begin{align*}
\rho_0 \xi_t^{\overline{\eta} \overline{\rho}} &= O\left( \widehat{X}_c(\overline{\eta}) \cdot \nabla_{\overline{\eta}} \varphi_3 \right) 
\end{align*}
In particular, we can estimate any term having a factor $O\left( \rho_0 \xi_t^{\overline{\eta} \overline{\rho}} \right)$ by integration by parts along $\widehat{X}_c(\overline{\eta}) \cdot \nabla_{\overline{\eta}}$ and get only terms of \eqref{qteslemestimeesgeneriquesdeccubiquehb2}. 

Likewise, we can estimate any terms having a factor 
\begin{align*}
O\left( \eta_0 \xi_t^{\overline{\eta} \overline{\rho}} \right), \quad O\left( \eta_0 \xi_t^{\overline{\eta} \overline{\sigma}} \right), \quad O\left( \sigma_0 \xi_t^{\overline{\eta} \overline{\sigma}} \right), \quad O\left( \rho_0 \xi_t^{\overline{\eta} \overline{\sigma}} \right)  
\end{align*}
In particular, due to the automatic presence of the factor $\rho_0$ ($F_3 = \partial_x f$), we can always reduce to a neighborhood where $\xi_t^{\overline{\eta} \overline{\sigma}}, \xi_t^{\overline{\eta} \overline{\rho}}$ are close to $0$, so $\theta^{\overline{\eta} \overline{\sigma}}, \theta^{\overline{\eta} \overline{\rho}}$ are close to $\pm 1$. 

Then, 
\begin{align*}
\eta \cdot \nabla_{\eta} \varphi_3 &= 2 \rho_0 \eta \cdot \rho - 2 \eta_0 |\eta|^2 \\
&= O\left( \rho_0 \xi_t^{\overline{\eta} \overline{\rho}} \right) + 2 |\eta| \left( \rho_0 \mbox{sgn}\left( \theta^{\overline{\eta} \overline{\rho}} \right) |\rho| - \eta_0 |\eta| \right) 
\end{align*}
so 
\begin{align*}
\rho_0 &= \frac{|\eta|}{|\rho|} \mbox{sgn}\left( \theta^{\overline{\eta} \overline{\rho}} \right) \eta_0 + O\left( \nabla_{\eta} \varphi_3 \right) 
\end{align*}
Note that we can estimate terms having a factor $O\left( \nabla_{\eta} \varphi_3 \right)$. Therefore, 
\begin{align*}
\eta_0 \partial_{\eta_0} \varphi_3 &= 3 \eta_0 \left( \rho_0^2 - \eta_0^2 \right) + \eta_0 \left( |\rho|^2 - |\eta|^2 \right) \\
&= \frac{3 \eta_0^3}{|\rho|^2} \left( |\eta|^2 - |\rho|^2 \right) + \eta_0 \left( |\rho|^2 - |\eta|^2 \right) + O\left( \nabla_{\eta} \varphi_3 \right) \\
&= \eta_0 \left( |\rho|^2 - |\eta|^2 \right) \left( 1 - \frac{3 \eta_0^2}{|\rho|^2} \right) + O\left( \nabla_{\eta} \varphi_3 \right) 
\end{align*}
We know how to estimate terms having a factor $O\left( \eta_0 \partial_{\eta_0} \varphi_3 \right)$, so we also get $O\left( \eta_0 \left( |\rho| - |\eta| \right) \right)$, and going back to $\nabla_{\eta} \varphi_3$ we also get $O\left( \rho_0 - \mbox{sgn}\left( \theta^{\overline{\eta} \overline{\rho}} \right) \eta_0 \right)$. Symmetric-wise, we can estimate if we have a factor 
\begin{align*}
O\left( \sigma_0 - \mbox{sgn}\left( \theta^{\overline{\eta} \overline{\sigma}} \right) \eta_0 \right), \quad O\left( \eta_0 \left( |\sigma| - |\eta| \right) \right) 
\end{align*}

Now, we can compute that {\footnotesize 
\begin{align*}
&\varphi_3 = \xi_0^3 - \eta_0^3 - \sigma_0^3 - \rho_0^3 + \xi_0 |\xi|^2 - \eta_0 |\eta|^2 - \sigma_0 |\sigma|^2 - \rho_0 |\rho|^2 \\
&\quad = O\left( \rho_0 - \mbox{sgn}\left( \theta^{\overline{\eta} \overline{\rho}} \right) \eta_0 \right) + O\left( \sigma_0 - \mbox{sgn}\left( \theta^{\overline{\eta} \overline{\sigma}} \right) \eta_0 \right) + O\left( \eta_0 \left( |\sigma| - |\eta| \right) \right) + O\left( \eta_0 \left( |\rho| - |\eta| \right) \right) \\
&\quad + \eta_0^3 \left[ \left( 1 + \mbox{sgn}\left( \theta^{\overline{\eta} \overline{\sigma}} \right) + \mbox{sgn}\left( \theta^{\overline{\eta} \overline{\rho}} \right) \right)^3 - 1 - \mbox{sgn}\left( \theta^{\overline{\eta} \overline{\sigma}} \right) - \mbox{sgn}\left( \theta^{\overline{\eta} \overline{\rho}} \right) \right] \\
&\quad + \eta_0 \left( 1 + \mbox{sgn}\left( \theta^{\overline{\eta} \overline{\sigma}} \right) + \mbox{sgn}\left( \theta^{\overline{\eta} \overline{\rho}} \right) \right) \left( |\eta|^2 + |\sigma|^2 + |\rho|^2 + 2 \eta \cdot \sigma + 2 \eta \cdot \rho + 2 \sigma \cdot \rho \right) - \eta_0 \left( 1 + \mbox{sgn}\left( \theta^{\overline{\eta} \overline{\sigma}} \right) + \mbox{sgn}\left( \theta^{\overline{\eta} \overline{\rho}} \right) \right) |\eta|^2 \\
&\quad = O\left( \rho_0 - \mbox{sgn}\left( \theta^{\overline{\eta} \overline{\rho}} \right) \eta_0 \right) + O\left( \sigma_0 - \mbox{sgn}\left( \theta^{\overline{\eta} \overline{\sigma}} \right) \eta_0 \right) + O\left( \eta_0 \left( |\sigma| - |\eta| \right) \right) + O\left( \eta_0 \left( |\rho| - |\eta| \right) \right) + O\left( \eta_0 \xi_t^{\overline{\eta} \overline{\rho}} \right) + O\left( \eta_0 \xi_t^{\overline{\eta} \overline{\sigma}} \right) \\
&\quad +
\eta_0 \left( 1 + \mbox{sgn}\left( \theta^{\overline{\eta} \overline{\sigma}} \right) + \mbox{sgn}\left( \theta^{\overline{\eta} \overline{\rho}} \right) \right) \Bigl[ \eta_0^2 \left( \left( 1 + \mbox{sgn}\left( \theta^{\overline{\eta} \overline{\sigma}} \right) + \mbox{sgn}\left( \theta^{\overline{\eta} \overline{\rho}} \right) \right)^2 - 1 \right) \\
&\quad \quad \quad + |\eta|^2 \left( 2 + 2 \mbox{sgn}\left( \theta^{\overline{\eta} \overline{\sigma}} \right) + 2 \mbox{sgn}\left( \theta^{\overline{\eta} \overline{\rho}} \right) + 2 \mbox{sgn}\left( \theta^{\overline{\eta} \overline{\sigma}} \right) \mbox{sgn}\left( \theta^{\overline{\eta} \overline{\rho}} \right) \right) \Bigl] \\
&\quad = O\left( \rho_0 - \mbox{sgn}\left( \theta^{\overline{\eta} \overline{\rho}} \right) \eta_0 \right) + O\left( \sigma_0 - \mbox{sgn}\left( \theta^{\overline{\eta} \overline{\sigma}} \right) \eta_0 \right) + O\left( \eta_0 \left( |\sigma| - |\eta| \right) \right) + O\left( \eta_0 \left( |\rho| - |\eta| \right) \right) + O\left( \eta_0 \xi_t^{\overline{\eta} \overline{\rho}} \right) + O\left( \eta_0 \xi_t^{\overline{\eta} \overline{\sigma}} \right) \\
&\quad +
\xi_0 \left[ \eta_0^2 \left( \left( 1 + \mbox{sgn}\left( \theta^{\overline{\eta} \overline{\sigma}} \right) + \mbox{sgn}\left( \theta^{\overline{\eta} \overline{\rho}} \right) \right)^2 - 1 \right) 
+ 2 |\eta|^2 \left( 1 + \mbox{sgn}\left( \theta^{\overline{\eta} \overline{\sigma}} \right) \right) \left( 1 + \mbox{sgn}\left( \theta^{\overline{\eta} \overline{\rho}} \right) \right) \right]
\end{align*} }
In particular, if $\mbox{sgn}\left( \theta^{\overline{\eta} \overline{\sigma}} \right) = \mbox{sgn}\left( \theta^{\overline{\eta} \overline{\rho}} \right) = 1$, then {\footnotesize 
\begin{align*}
\xi_0 &= O\left( \varphi_3 \right) + O\left( \rho_0 - \mbox{sgn}\left( \theta^{\overline{\eta} \overline{\rho}} \right) \eta_0 \right) + O\left( \sigma_0 - \mbox{sgn}\left( \theta^{\overline{\eta} \overline{\sigma}} \right) \eta_0 \right) + O\left( \eta_0 \left( |\sigma| - |\eta| \right) \right) + O\left( \eta_0 \left( |\rho| - |\eta| \right) \right) \\
&\quad + O\left( \eta_0 \xi_t^{\overline{\eta} \overline{\rho}} \right) + O\left( \eta_0 \xi_t^{\overline{\eta} \overline{\sigma}} \right)
\end{align*} }
which concludes. 

Now, we can assume that at least one of the signs of the $\theta$'s is $-1$. In particular, by an analogous computation, since $1 + \mbox{sgn}\left( \theta^{\overline{\eta} \overline{\sigma}} \right) + \mbox{sgn}\left( \theta^{\overline{\eta} \overline{\rho}} \right) \in \{ -1, 1 \}$, {\footnotesize 
\begin{align*}
|\xi_0| &= O\left( \rho_0 - \mbox{sgn}\left( \theta^{\overline{\eta} \overline{\rho}} \right) \eta_0 \right) + O\left( \sigma_0 - \mbox{sgn}\left( \theta^{\overline{\eta} \overline{\sigma}} \right) \eta_0 \right) + |\eta_0| \\
&= O\left( \rho_0 - \mbox{sgn}\left( \theta^{\overline{\eta} \overline{\rho}} \right) \eta_0 \right) + O\left( \sigma_0 - \mbox{sgn}\left( \theta^{\overline{\eta} \overline{\sigma}} \right) \eta_0 \right) + |\sigma_0| \\
&= O\left( \rho_0 - \mbox{sgn}\left( \theta^{\overline{\eta} \overline{\rho}} \right) \eta_0 \right) + O\left( \sigma_0 - \mbox{sgn}\left( \theta^{\overline{\eta} \overline{\sigma}} \right) \eta_0 \right) + |\rho_0| \\
\xi_0 |\xi|^2 &= O\left( \rho_0 - \mbox{sgn}\left( \theta^{\overline{\eta} \overline{\rho}} \right) \eta_0 \right) + O\left( \sigma_0 - \mbox{sgn}\left( \theta^{\overline{\eta} \overline{\sigma}} \right) \eta_0 \right) + O\left( \eta_0 \left( |\sigma| - |\eta| \right) \right) + O\left( \eta_0 \left( |\rho| - |\eta| \right) \right) \\
&\quad + O\left( \eta_0 \xi_t^{\overline{\eta} \overline{\rho}} \right) + O\left( \eta_0 \xi_t^{\overline{\eta} \overline{\sigma}} \right) + \xi_0 |\eta|^2 
\end{align*} }
We can now compute that {\footnotesize 
\begin{align*}
&\xi_0 \widehat{X}_b(\overline{\xi}) \cdot \nabla_{\overline{\xi}} \varphi_3 
= \xi_0 \frac{|\xi|}{|\overline{\xi}|} \left( 3 \xi_0^2 + |\xi|^2 - 3 \eta_0^2 - |\eta|^2 \right) - \xi_0 \frac{\xi_0 \xi}{|\overline{\xi}| |\xi|} \cdot \left( 2 \xi_0 \xi - 2 \rho_0 \rho \right) \\
&= O\left( |\xi_0| - |\eta_0| \right) + O\left( \xi_0 \left( |\xi|^2 - |\eta|^2 \right) \right) 
- 2 \xi_0 \frac{\xi_0}{|\overline{\xi}| |\xi|} \left( \xi_0 |\xi|^2 - \rho_0 \eta \cdot \rho - \rho_0 \sigma \cdot \rho - \rho_0 |\rho|^2 \right) \\
&= O\left( \rho_0 - \mbox{sgn}\left( \theta^{\overline{\eta} \overline{\rho}} \right) \eta_0 \right) + O\left( \sigma_0 - \mbox{sgn}\left( \theta^{\overline{\eta} \overline{\sigma}} \right) \eta_0 \right) + O\left( \eta_0 \left( |\sigma| - |\eta| \right) \right) + O\left( \eta_0 \left( |\rho| - |\eta| \right) \right) + O\left( \eta_0 \xi_t^{\overline{\eta} \overline{\rho}} \right) + O\left( \eta_0 \xi_t^{\overline{\eta} \overline{\sigma}} \right) \\
&- 2 \xi_0 \frac{\xi_0}{|\overline{\xi}| |\xi|} |\eta|^2 \left( \left( 1 + \mbox{sgn}\left( \theta^{\overline{\eta} \overline{\sigma}} \right) + \mbox{sgn}\left( \theta^{\overline{\eta} \overline{\rho}} \right) \right) \eta_0 - \rho_0 \left( 1 + \mbox{sgn}\left( \theta^{\overline{\eta} \overline{\rho}} \right) + \mbox{sgn}\left( \theta^{\overline{\eta} \overline{\rho}} \right) \mbox{sgn}\left( \theta^{\overline{\eta} \overline{\sigma}} \right) \right) \right) \\
&= O\left( \rho_0 - \mbox{sgn}\left( \theta^{\overline{\eta} \overline{\rho}} \right) \eta_0 \right) + O\left( \sigma_0 - \mbox{sgn}\left( \theta^{\overline{\eta} \overline{\sigma}} \right) \eta_0 \right) + O\left( \eta_0 \left( |\sigma| - |\eta| \right) \right) + O\left( \eta_0 \left( |\rho| - |\eta| \right) \right) + O\left( \eta_0 \xi_t^{\overline{\eta} \overline{\rho}} \right) + O\left( \eta_0 \xi_t^{\overline{\eta} \overline{\sigma}} \right) \\
&- 2 \xi_0 \frac{\xi_0}{|\overline{\xi}| |\xi|} |\eta|^2 \eta_0 \left( 1 + \mbox{sgn}\left( \theta^{\overline{\eta} \overline{\sigma}} \right) + \mbox{sgn}\left( \theta^{\overline{\eta} \overline{\rho}} \right) - \mbox{sgn}\left( \theta^{\overline{\eta} \overline{\rho}} \right) \left( 1 + \mbox{sgn}\left( \theta^{\overline{\eta} \overline{\rho}} \right) + \mbox{sgn}\left( \theta^{\overline{\eta} \overline{\rho}} \right) \mbox{sgn}\left( \theta^{\overline{\eta} \overline{\sigma}} \right) \right) \right) \\
&= O\left( \rho_0 - \mbox{sgn}\left( \theta^{\overline{\eta} \overline{\rho}} \right) \eta_0 \right) + O\left( \sigma_0 - \mbox{sgn}\left( \theta^{\overline{\eta} \overline{\sigma}} \right) \eta_0 \right) + O\left( \eta_0 \left( |\sigma| - |\eta| \right) \right) + O\left( \eta_0 \left( |\rho| - |\eta| \right) \right) + O\left( \eta_0 \xi_t^{\overline{\eta} \overline{\rho}} \right) + O\left( \eta_0 \xi_t^{\overline{\eta} \overline{\sigma}} \right)
\end{align*} }
which concludes. 

\paragraph{2.} Let us now localise to have $|\overline{\eta}| \gg |\overline{\sigma}|+|\overline{\rho}|$. 

We have that {\footnotesize 
\begin{align*}
&\varphi_3 = \xi_0^3 - \eta_0^3 - \sigma_0^3 - \rho_0^3 + \xi_0 |\xi|^2 - \eta_0 |\eta|^2 - \sigma_0 |\sigma|^2 - \rho_0 |\rho|^2 \\
&= O\left( \rho_0 \overline{\sigma} \right) + O\left( \sigma_0 \overline{\rho} \right) + O\left( \eta_0 \overline{\sigma} \right) + O\left( \eta_0 \overline{\rho} \right) + \xi_0 \left( |\eta|^2 + |\sigma|^2 + |\rho|^2 + 2 \eta \cdot (\rho + \sigma) + 2 \rho \cdot \sigma \right) - \eta_0 |\eta|^2 - \sigma_0 |\sigma|^2 - \rho_0 |\rho|^2 \\
&= O\left( \rho_0 \overline{\sigma} \right) + O\left( \sigma_0 \overline{\rho} \right) + O\left( \eta_0 \overline{\sigma} \right) + O\left( \eta_0 \overline{\rho} \right) + \xi_0 \left( |\eta|^2 + 2 \eta \cdot (\rho + \sigma) \right) + (\xi_0 - \sigma_0) |\sigma|^2 + (\xi_0 - \rho_0) |\rho|^2 - \eta_0 |\eta|^2 \\
&= O\left( \rho_0 \overline{\sigma} \right) + O\left( \sigma_0 \overline{\rho} \right) + O\left( \eta_0 \overline{\sigma} \right) + O\left( \eta_0 \overline{\rho} \right) + (\xi_0-\eta_0) \left( |\eta|^2 + 2 \eta \cdot (\rho + \sigma) \right)
\end{align*} }
so 
\begin{align*}
\xi_0-\eta_0 &= O\left( \varphi_3 \right) + O\left( \rho_0 \overline{\sigma} \right) + O\left( \sigma_0 \overline{\rho} \right) + O\left( \eta_0 \overline{\sigma} \right) + O\left( \eta_0 \overline{\rho} \right)
\end{align*}
Therefore, since we have in \eqref{equdecchampbPPP-1} a factor $\xi_0$ and a factor $\rho_0$, we can write 
\begin{align*}
\xi_0 \rho_0 &= (\xi_0 - \eta_0) \rho_0 + \eta_0 \rho_0 = O\left( \rho_0 \varphi_3 \right) + O\left( \rho_0 \overline{\sigma} \right) + O\left( \sigma_0 \overline{\rho} \right) + O\left( \eta_0 \overline{\sigma} \right) + O\left( \eta_0 \overline{\rho} \right)
\end{align*}
But here, 
\begin{align*}
1 = O(\partial_{\eta_0} \varphi_3), \quad 1 = O\left( (\partial_{\eta_0}-\partial_{\sigma_0}) \varphi_3 \right) 
\end{align*}
We can therefore apply an integration by parts in time on the term having a factor $O\left( \rho_0 \varphi_3 \right)$, along $\partial_{\eta_0}$ for the factors $O\left( \rho_0 \overline{\sigma} \right)$, $O\left( \eta_0 \overline{\sigma} \right)$, and along $(\partial_{\eta_0}-\partial_{\sigma_0})$ for the factors $O\left( \sigma_0 \overline{\rho} \right)$ and $O\left( \eta_0 \overline{\rho} \right)$. In all cases, we only recover terms from \eqref{qteslemestimeesgeneriquesdeccubiquehb2}. 

This concludes the proof of Lemma \ref{lemmedecompositioncubiquehquadbb}.

\section{Estimate for the singular b-weight} \label{section-estimee-mauvaise-b} 

In this section, we prove the following: 

\begin{Prop} We can decompose 
\begin{align*}
X_b f(t) &= h_{b*}(t) + g_{b*}(t) \\
X_c f(t) &= h_{c*}(t) + g_{c*}(t) 
\end{align*}
such that, for every $t > 0$, 
\begin{align*}
\Vert m_{\widehat{\mathcal{C}}}(D) h_{b*}(t) \Vert_{L^2}^2 &\lesssim \Vert (x, y) u_0 \Vert_{L^2}^2 + \langle t \rangle^{1+202\delta} \Vert u \Vert_X^3 \\
\Vert m_{\widehat{\mathcal{C}}}(D) \nabla h_{b*}(t) \Vert_{L^2}^2 &\lesssim \Vert (x, y) u_0 \Vert_{H^1}^2 + \langle t \rangle^{\frac{1}{3}+202\delta} \Vert u \Vert_X^3 \\
\Vert m_{\widehat{\mathcal{C}}}(D) g_{b*}(t) \Vert_{L^2} &\lesssim \langle t \rangle^{\frac{1}{2}+101\delta} \Vert u \Vert_X^2 \\
\Vert e^{-it \omega(D)} m_{\widehat{\mathcal{C}}}(D) \nabla g_{b*}(t) \Vert_{L^4} &\lesssim t^{-\frac{1}{12}} \langle t \rangle^{101\delta} \Vert u \Vert_X^2 \\
\Vert m_{\widehat{\mathcal{C}}}(D) \nabla g_{b*}(t) \Vert_{L^2} &\lesssim \langle t \rangle^{\frac{1}{6}+101\delta} \Vert u \Vert_X^2 \\
\Vert e^{-it \omega(D)} m_{\widehat{\mathcal{C}}}(D) \langle \nabla \rangle^{-1} \nabla^2 g_{b*}(t) \Vert_{L^4} &\lesssim t^{-\frac{5}{12}} \langle t \rangle^{101\delta} \Vert u \Vert_X^2 \\
\Vert m_{\widehat{\mathcal{L}}}(D) h_{b*}(t) \Vert_{L^2}^2 &\lesssim \Vert (x, y) u_0 \Vert_{H^1}^2 + \langle t \rangle^{\frac{5}{24}+404\delta} \Vert u \Vert_X^3 \left( 1 + \Vert u \Vert_X \right) \\
\Vert m_{\widehat{\mathcal{L}}}(D) h_{c*}(t) \Vert_{L^2}^2 &\lesssim \Vert (x, y) u_0 \Vert_{H^1}^2 + \langle t \rangle^{\frac{5}{24}+404\delta} \Vert u \Vert_X^3 \left( 1 + \Vert u \Vert_X \right) \\
\Vert m_{\widehat{\mathcal{L}}}(D) |\nabla|^{\frac{1}{2}} \langle \nabla \rangle^{\frac{1}{2}} h_{b*}(t) \Vert_{L^2}^2 &\lesssim \Vert \langle x, y \rangle^2 u_0 \Vert_{L^2}^2 + \langle t \rangle^{404\delta} \Vert u \Vert_X^3 \left( 1 + \Vert u \Vert_X \right) \\
\Vert m_{\widehat{\mathcal{L}}}(D) |\nabla|^{\frac{1}{2}} \langle \nabla \rangle^{\frac{1}{2}} h_{c*}(t) \Vert_{L^2}^2 &\lesssim \Vert \langle x, y \rangle^2 u_0 \Vert_{L^2}^2 + \langle t \rangle^{404\delta} \Vert u \Vert_X^3 \left( 1 + \Vert u \Vert_X \right) \\
\Vert m_{\widehat{\mathcal{L}}}(D) g_{b*}(t) \Vert_{L^2} &\lesssim \langle t \rangle^{\frac{5}{48}+202\delta} \Vert u \Vert_X^2 \\
\Vert m_{\widehat{\mathcal{L}}}(D) g_{c*}(t) \Vert_{L^2} &\lesssim \langle t \rangle^{\frac{5}{48}+202\delta} \Vert u \Vert_X^2 \\
\Vert e^{-it\omega(D)} m_{\widehat{\mathcal{L}}}(D) \nabla g_{b*}(t) \Vert_{L^4} &\lesssim t^{-\frac{23}{48}} \langle t \rangle^{202\delta} \Vert u \Vert_X^2 \\
\Vert e^{-it\omega(D)} m_{\widehat{\mathcal{L}}}(D) \nabla g_{c*}(t) \Vert_{L^4} &\lesssim t^{-\frac{23}{48}} \langle t \rangle^{202\delta} \Vert u \Vert_X^2 \\
\Vert m_{\widehat{\mathcal{L}}}(D) |\nabla|^{\frac{1}{2}} \langle \nabla \rangle^{\frac{1}{2}} g_{b*}(t) \Vert_{L^2} &\lesssim \langle t \rangle^{202\delta} \Vert u \Vert_X^2 \\
\Vert m_{\widehat{\mathcal{L}}}(D) |\nabla|^{\frac{1}{2}} \langle \nabla \rangle^{\frac{1}{2}} g_{c*}(t) \Vert_{L^2} &\lesssim \langle t \rangle^{202\delta} \Vert u \Vert_X^2 \\
\Vert e^{-it \omega(D)} m_{\widehat{\mathcal{L}}}(D) \langle \nabla \rangle^{-\frac{1}{2}} \nabla^{\frac{3}{2}} g_{b*}(t) \Vert_{L^4} &\lesssim t^{-\frac{5}{12}} \langle t \rangle^{-\frac{1}{8}+202\delta} \Vert u \Vert_X^2 \\
\Vert e^{-it \omega(D)} m_{\widehat{\mathcal{L}}}(D) \langle \nabla \rangle^{-\frac{1}{2}} \nabla^{\frac{3}{2}} g_{c*}(t) \Vert_{L^4} &\lesssim t^{-\frac{5}{12}} \langle t \rangle^{-\frac{1}{8}+202\delta} \Vert u \Vert_X^2 
\end{align*}
\end{Prop}

As before, in order to get the decomposition $h_{b*}+g_{b*}$, we apply Duhamel's formula in Fourier: 
\begin{align*}
\widehat{X}_b(\overline{\xi}) \cdot \nabla_{\overline{\xi}} \widehat{f}(t) &= \widehat{X}_b(\overline{\xi}) \cdot \nabla_{\overline{\xi}} \widehat{f}(0) + \widehat{X}_b(\overline{\xi}) \cdot \nabla_{\overline{\xi}} \left( i \xi_0 \widehat{I}[f, f](t, \overline{\xi}) \right) 
\end{align*}
We always place the initial data in $h_{b*}$, then we separate the interaction term $I[f, f]$ depending on the geometric directions of the interacting frequencies $\overline{\eta}, \overline{\sigma}$, and on the relative sizes through the symbols $\mu_{BBB}, \mu_{BHH}, \mu_{HBH}, \mu_{HHB}$ already introduced. Note that, if we localise $\overline{\xi}$ by $m_{\widehat{\mathcal{R}}}$ or $m_{\widehat{\mathcal{P}}}$, since on these areas $m_b(\overline{\xi}) \simeq 1$, we can simply apply the already obtained decomposition and we already have (better) estimates. 

By symmetry of $\overline{\eta}, \overline{\sigma}$, we only need to consider the interaction cases $BBB, HBH$ and $BHH$. Moreover, in the $BBB$ case, we can symmetrize again the geometric areas. 

Moreover, for $h_{b*}$, we can apply the Duhamel formula on the norm and estimate for $k = 0, 1$, 
\begin{align}
\Vert \nabla^k h_{b*}(t) \Vert_{L^2}^2 &= \Vert \nabla^k h_{b*}(0) \Vert_{L^2}^2 + 2 \mbox{Re} \int_0^t \int |\overline{\xi}|^{2k} \widehat{h}_{b*}(s, -\overline{\xi}) ~ \partial_s \widehat{h}_{b*}(s, \overline{\xi}) ~ d\overline{\xi} ds \label{estapriorimauvaisbDuhamelhcarre}
\end{align}
by Parseval's identity. 

The initial data is exactly
\begin{align*}
h_{b*}(0) = X_b f(0) = X_b u_0
\end{align*}
and its estimate is direct. 

In \eqref{estapriorimauvaisbDuhamelhcarre}, in the case $\widehat{\mathcal{C}}$, when we develop $\partial_s \widehat{h}_{b*}$, any term that can be estimated in $L^2$ (respectively in $\dot{H}^1$) by 
\begin{align*} \langle s \rangle^{\frac{1}{2}+101\delta} r(s) \Vert u \Vert_X^2 \left( 1 + \Vert u \Vert_X \right) \quad \left( \mbox{ respectively } ~ \langle s \rangle^{\frac{1}{6}+101\delta} r(s) \Vert u \Vert_X^2 \left( 1 + \Vert u \Vert_X \right) ~ \right) 
\end{align*}
for $r \in L^1(\mathbb{R}^{+})$ is automatically compatible with the desired a priori estimate since, if we denote by $F(s)$ such a term, 
\begin{align*}
&\int_0^t \int \widehat{h}_{b*}(s, -\overline{\xi}) ~ \widehat{F}(s, \overline{\xi}) ~ d\overline{\xi} ds \lesssim \int_0^t \Vert h_{b*}(s) \Vert_{L^2} ~ \Vert F(s) \Vert_{L^2} ~ ds \\
&\quad \lesssim \int_0^t \langle s \rangle^{1+202\delta} r(s) \Vert u \Vert_X^3 \left( 1 + \Vert u \Vert_X \right) ~ ds \lesssim \langle t \rangle^{1+202\delta} \Vert u \Vert_X^3 \left( 1 + \Vert u \Vert_X \right) 
\end{align*}
(where the implicit constant depends on the function $r$), here above for the case $k = 0$, and likewise for the case $k = 1$. Therefore, we will only come back to the full equation \eqref{estapriorimauvaisbDuhamelhcarre} if we need some symmetry argument to distribute derivatives. 

We can estimate similarly for $\widehat{\mathcal{L}}$, estimating $\partial_s \widehat{h}_{b*}$ in $L^2$ (respectively in $L^2$ with $|\nabla|^{\frac{1}{2}} \langle \nabla \rangle^{\frac{1}{2}}$) by 
\begin{align*}
\langle s \rangle^{\frac{5}{48}+202\delta} r(s) \Vert u \Vert_X^2 \left( 1 + \Vert u \Vert_X \right) \quad \left( \mbox{ respectively } ~ \langle s \rangle^{202\delta} r(s) \Vert u \Vert_X^2 \left( 1 + \Vert u \Vert_X \right) ~ \right) 
\end{align*}

\subsection{Case of the cone}

In this subsection, we localise $\overline{\xi}$ by $m_{\widehat{\mathcal{C}}}$. 

As we will see, all the estimates will be better than the desired one, except in the case $(BHH, \widehat{\mathcal{C}})$ which is treated last and where the growth exponent are saturated. 

\subsubsection{Generic interaction \texorpdfstring{$BBB$}{BBB}}

Let us first localise by 
\begin{align*}
\mu_{BBB}(\overline{\xi}, \overline{\eta}) \left( 1 - m_{\widehat{\mathcal{P}}}(\overline{\eta}) \right) \left( 1 - m_{\widehat{\mathcal{P}}}(\overline{\sigma}) \right) 
\end{align*} 

In this case, we can develop {\footnotesize 
\begin{subequations}
\begin{align}
&m_{\widehat{\mathcal{C}}}(\overline{\xi}) \widehat{X}_b(\overline{\xi}) \cdot \nabla_{\overline{\xi}} \left( i \xi_0 \widehat{I}_{\mu_{BBB} (1 - m_{\widehat{\mathcal{P}}}(D))^{\otimes 2}}[f, f](t, \overline{\xi}) \right) \notag \\
&\quad = \int_0^t \int e^{i s \varphi} i s m_{\widehat{\mathcal{C}}}(\overline{\xi}) \widehat{X}_b(\overline{\xi}) \cdot \nabla_{\overline{\xi}} \varphi ~ i \xi_0 \mu_{BBB}(\overline{\xi}, \overline{\eta}) \left( 1 - m_{\widehat{\mathcal{P}}}(\overline{\eta}) \right) \left( 1 - m_{\widehat{\mathcal{P}}}(\overline{\sigma}) \right) \widehat{f}(s, \overline{\eta}) \widehat{f}(s, \overline{\sigma}) ~ d\overline{\eta} ds \label{estapriorimauvaisbCgenBBB-s} \\
&\quad \quad + \int_0^t \int e^{i s \varphi} i \xi_0 \mu_{BBB}(\overline{\xi}, \overline{\eta}) \left( 1 - m_{\widehat{\mathcal{P}}}(\overline{\eta}) \right) \left( 1 - m_{\widehat{\mathcal{P}}}(\overline{\sigma}) \right) \widehat{f}(s, \overline{\eta}) m_{\widehat{\mathcal{C}}}(\overline{\xi}) \widehat{X}_b(\overline{\xi}) \cdot \nabla_{\overline{\xi}} \widehat{f}(s, \overline{\sigma}) ~ d\overline{\eta} ds \label{estapriorimauvaisbCgenBBB-eta} \\
&\quad \quad + \int_0^t \int e^{i s \varphi} m_{\widehat{\mathcal{C}}}(\overline{\xi}) \widehat{X}_b(\overline{\xi}) \cdot \nabla_{\overline{\xi}} \left( i \xi_0 \mu_{BBB}(\overline{\xi}, \overline{\eta}) \left( 1 - m_{\widehat{\mathcal{P}}}(\overline{\eta}) \right) \left( 1 - m_{\widehat{\mathcal{P}}}(\overline{\sigma}) \right) \right) \widehat{f}(s, \overline{\eta}) \widehat{f}(s, \overline{\sigma}) ~ d\overline{\eta} ds \label{estapriorimauvaisbCgenBBB-dersymb} 
\end{align}
\end{subequations} }

By Lemma \ref{lem-non-res-loin0Cone}, there exists locally $\mu_j$ symbols of order $j$, $j = -2, -3$, such that
\begin{align*}
1 = \mu_{-2} \nabla_{\overline{\eta}} \varphi + \mu_{-3} \varphi 
\end{align*}
Therefore, we can develop \eqref{estapriorimauvaisbCgenBBB-s} and apply integrations by parts: {\footnotesize 
\begin{subequations}
\begin{align}
&\eqref{estapriorimauvaisbCgenBBB-s} = \int_0^t \int e^{i s \varphi} \mu_{-2}(\overline{\xi}, \overline{\eta}) \nabla_{\overline{\eta}} \varphi ~ i s m_{\widehat{\mathcal{C}}}(\overline{\xi}) \widehat{X}_b(\overline{\xi}) \cdot \nabla_{\overline{\xi}} \varphi ~ i \xi_0 \mu_{BBB}(\overline{\xi}, \overline{\eta}) \left( 1 - m_{\widehat{\mathcal{P}}}(\overline{\eta}) \right) \left( 1 - m_{\widehat{\mathcal{P}}}(\overline{\sigma}) \right) \widehat{f}(s, \overline{\eta}) \widehat{f}(s, \overline{\sigma}) ~ d\overline{\eta} ds \notag \\
&\quad \quad + \int_0^t \int e^{i s \varphi} \mu_{-3}(\overline{\xi}, \overline{\eta}) \varphi ~ i s m_{\widehat{\mathcal{C}}}(\overline{\xi}) \widehat{X}_b(\overline{\xi}) \cdot \nabla_{\overline{\xi}} \varphi ~ i \xi_0 \mu_{BBB}(\overline{\xi}, \overline{\eta}) \left( 1 - m_{\widehat{\mathcal{P}}}(\overline{\eta}) \right) \left( 1 - m_{\widehat{\mathcal{P}}}(\overline{\sigma}) \right) \widehat{f}(s, \overline{\eta}) \widehat{f}(s, \overline{\sigma}) ~ d\overline{\eta} ds \notag \\
&= - \int_0^t \int e^{i s \varphi} \mu_{-2}(\overline{\xi}, \overline{\eta}) m_{\widehat{\mathcal{C}}}(\overline{\xi}) \widehat{X}_b(\overline{\xi}) \cdot \nabla_{\overline{\xi}} \varphi ~ i \xi_0 \mu_{BBB}(\overline{\xi}, \overline{\eta}) \left( 1 - m_{\widehat{\mathcal{P}}}(\overline{\eta}) \right) \left( 1 - m_{\widehat{\mathcal{P}}}(\overline{\sigma}) \right) \nabla_{\overline{\eta}} \left( \widehat{f}(s, \overline{\eta}) \widehat{f}(s, \overline{\sigma}) \right) ~ d\overline{\eta} ds \label{estapriorimauvaisbCgenBBB-s-1} \\
&- \int_0^t \int e^{i s \varphi} \nabla_{\overline{\eta}} \cdot \left( \mu_{-2}(\overline{\xi}, \overline{\eta}) m_{\widehat{\mathcal{C}}}(\overline{\xi}) \widehat{X}_b(\overline{\xi}) \cdot \nabla_{\overline{\xi}} \varphi ~ i \xi_0 \mu_{BBB}(\overline{\xi}, \overline{\eta}) \left( 1 - m_{\widehat{\mathcal{P}}}(\overline{\eta}) \right) \left( 1 - m_{\widehat{\mathcal{P}}}(\overline{\sigma}) \right) \right) \widehat{f}(s, \overline{\eta}) \widehat{f}(s, \overline{\sigma}) ~ d\overline{\eta} ds \label{estapriorimauvaisbCgenBBB-s-2} \\
&- \int_0^t \int e^{i s \varphi} \mu_{-3}(\overline{\xi}, \overline{\eta}) s m_{\widehat{\mathcal{C}}}(\overline{\xi}) \widehat{X}_b(\overline{\xi}) \cdot \nabla_{\overline{\xi}} \varphi ~ i \xi_0 \mu_{BBB}(\overline{\xi}, \overline{\eta}) \left( 1 - m_{\widehat{\mathcal{P}}}(\overline{\eta}) \right) \left( 1 - m_{\widehat{\mathcal{P}}}(\overline{\sigma}) \right) \partial_s \left( \widehat{f}(s, \overline{\eta}) \widehat{f}(s, \overline{\sigma}) \right) ~ d\overline{\eta} ds \label{estapriorimauvaisbCgenBBB-s-3} \\
&- \int_0^t \int e^{i s \varphi} \mu_{-3}(\overline{\xi}, \overline{\eta}) m_{\widehat{\mathcal{C}}}(\overline{\xi}) \widehat{X}_b(\overline{\xi}) \cdot \nabla_{\overline{\xi}} \varphi ~ i \xi_0 \mu_{BBB}(\overline{\xi}, \overline{\eta}) \left( 1 - m_{\widehat{\mathcal{P}}}(\overline{\eta}) \right) \left( 1 - m_{\widehat{\mathcal{P}}}(\overline{\sigma}) \right) \widehat{f}(s, \overline{\eta}) \widehat{f}(s, \overline{\sigma}) ~ d\overline{\eta} ds \label{estapriorimauvaisbCgenBBB-s-4} \\
&+ \int e^{i t \varphi} \mu_{-3}(\overline{\xi}, \overline{\eta}) t m_{\widehat{\mathcal{C}}}(\overline{\xi}) \widehat{X}_b(\overline{\xi}) \cdot \nabla_{\overline{\xi}} \varphi ~ i \xi_0 \mu_{BBB}(\overline{\xi}, \overline{\eta}) \left( 1 - m_{\widehat{\mathcal{P}}}(\overline{\eta}) \right) \left( 1 - m_{\widehat{\mathcal{P}}}(\overline{\sigma}) \right) \widehat{f}(s, \overline{\eta}) \widehat{f}(s, \overline{\sigma}) ~ d\overline{\eta} ds \label{estapriorimauvaisbCgenBBB-s-5} 
\end{align}
\end{subequations} }
We only put the boundary term \eqref{estapriorimauvaisbCgenBBB-s-5} in $g_{b*}$ and the rest in $h_{b*}$. 

We can then estimate 
\begin{align*}
\Vert \eqref{estapriorimauvaisbCgenBBB-s-5} \Vert_{L^2} &\lesssim t \Vert u(t) \Vert_{L^4}^2 \lesssim t^{\frac{1}{6}} \langle t \rangle^{-\frac{1}{4}+100\delta} \Vert u \Vert_X^2 \lesssim \Vert u \Vert_X^2 \\
\Vert \overline{\xi} \eqref{estapriorimauvaisbCgenBBB-s-5} \Vert_{L^2} &\lesssim \Vert \xi_0 \eqref{estapriorimauvaisbCgenBBB-s-5} \Vert_{L^2} \lesssim t \Vert u(t) \Vert_{L^2} \Vert \partial_x u(t) \Vert_{L^{\infty}} \lesssim \Vert u \Vert_X^2 \\
\Vert e^{-i t \omega(D)} \nabla \mathcal{F}^{-1} \eqref{estapriorimauvaisbCgenBBB-s-5} \Vert_{L^4} &\lesssim t \Vert u(t) \Vert_{L^4} \Vert \partial_x u(t) \Vert_{L^{\infty}} \lesssim t^{-\frac{1}{4}} \langle t \rangle^{-\frac{3}{8}+150\delta} \Vert u \Vert_X^2 
\end{align*}
where we used that $\overline{\xi}$ was away from $\widehat{\mathcal{P}}$. Note that we did not use the localisation $\mu_{BBB}$ for these estimates. 

We can directly estimate the derivative in $t$ of the terms added to $h_{b*}$ as well: 
\begin{align*}
\Vert \partial_t \eqref{estapriorimauvaisbCgenBBB-eta} \Vert_{L^2} &\lesssim \Vert \partial_x u(t) \Vert_{L^{\infty}} \Vert (x, y) f(t) \Vert_{L^2} \lesssim t^{-\frac{5}{6}} \langle t \rangle^{-\frac{1}{4}+100\delta} ~ \langle t \rangle^{\frac{1}{2}+101\delta} \Vert u \Vert_X^2 \\
\Vert \overline{\xi} \partial_t \eqref{estapriorimauvaisbCgenBBB-eta} \Vert_{L^2} &\lesssim \Vert \partial_x u(t) \Vert_{L^{\infty}} \Vert \nabla (x, y) f(t) \Vert_{L^2} \lesssim t^{-\frac{5}{6}} \langle t \rangle^{-\frac{1}{4}+100\delta} ~ \langle t \rangle^{\frac{1}{6}+101\delta} \Vert u \Vert_X^2 \\
\Vert \langle \overline{\xi} \rangle \partial_t \eqref{estapriorimauvaisbCgenBBB-dersymb} \Vert_{L^2} &\lesssim \Vert \partial_x u(t) \Vert_{L^{\infty}} \Vert \langle \nabla \rangle |\nabla|^{-1} u(t) \Vert_{L^2} \lesssim t^{-\frac{5}{6}} \langle t \rangle^{-\frac{1}{4}+100\delta} \Vert u \Vert_X^2 \\
\Vert \partial_t \eqref{estapriorimauvaisbCgenBBB-s-3} \Vert_{L^2} &\lesssim t \Vert e^{-i t \omega(D)} \partial_t f(t) \Vert_{L^4} \Vert u(t) \Vert_{L^4} \lesssim t \Vert u(t) \Vert_{L^4}^2 \Vert \partial_x u(t) \Vert_{L^{\infty}} \\
&\lesssim t^{-\frac{2}{3}} \langle t \rangle^{-\frac{1}{2}+200\delta} \Vert u \Vert_X^3 \\
\Vert \overline{\xi} \partial_t \eqref{estapriorimauvaisbCgenBBB-s-3} \Vert_{L^2} &\lesssim t \Vert e^{-it \omega(D)} \partial_x \partial_t f(t) \Vert_{L^4} \Vert u(t) \Vert_{L^4} + t \Vert \partial_t f(t) \Vert_{L^2} \Vert \partial_x u(t) \Vert_{L^{\infty}} \\
&\lesssim t \Vert \partial_x u \Vert_{L^{\infty}} \Vert u(t) \Vert_{L^4}^2 + t \Vert u(t) \Vert_{L^2} \Vert \partial_x u(t) \Vert_{L^{\infty}}^2 + t \Vert u(t) \Vert_{L^4}^2 \Vert \partial_x u(t) \Vert_{L^{\infty}} \\
&\quad + \Vert \nabla (x, y) f(t) \Vert_{L^2} \Vert u(t) \Vert_{L^{\infty}} \Vert u(t) \Vert_{L^4} \\
&\lesssim t^{-\frac{2}{3}} \langle t \rangle^{-\frac{1}{2}+200\delta} \Vert u \Vert_X^3 
+ t^{-\frac{11}{12}} \langle t \rangle^{-\frac{5}{24}+301\delta} \Vert u \Vert_X^3 \lesssim t^{-\frac{11}{12}} \langle t \rangle^{-\frac{5}{24}+301\delta} \Vert u \Vert_X^3 
\end{align*}
\eqref{estapriorimauvaisbCgenBBB-s-1} is similar to \eqref{estapriorimauvaisbCgenBBB-eta}, \eqref{estapriorimauvaisbCgenBBB-s-2} and \eqref{estapriorimauvaisbCgenBBB-s-4} are similar to \eqref{estapriorimauvaisbCgenBBB-dersymb}. Note that, for \eqref{estapriorimauvaisbCgenBBB-s-3}, we did not use the localisation $\mu_{BBB}$. 

\subsubsection{Interaction \texorpdfstring{$(BBB, \widehat{\mathcal{R}}\widehat{\mathcal{P}})$}{(BBB, RP)} or \texorpdfstring{$(BBB, \widehat{\mathcal{L}}\widehat{\mathcal{P}})$}{(BBB, LP)}} 

Let us localise by $\mu_{BBB}(\overline{\xi}, \overline{\eta}) m_{A}(\overline{\eta}) m_{\widehat{\mathcal{P}}}(\overline{\sigma})$, $A = \widehat{\mathcal{R}}$ ou $A = \widehat{\mathcal{L}}$. 

We can then develop 
\begin{subequations}
\begin{align}
&m_{\widehat{\mathcal{C}}}(\overline{\xi}) \widehat{X}_b(\overline{\xi}) \cdot \nabla_{\overline{\xi}} \left( i \xi_0 \widehat{I}_{\mu_{BBB} m_A m_{\widehat{\mathcal{P}}}}[f, f](t, \overline{\xi}) \right) \notag \\
&\quad = \int_0^t \int e^{i s \varphi} i s m_{\widehat{\mathcal{C}}}(\overline{\xi}) \widehat{X}_b(\overline{\xi}) \cdot \nabla_{\overline{\xi}} \varphi ~ i \xi_0 \mu_{BBB}(\overline{\xi}, \overline{\eta}) m_{A}(\overline{\eta}) m_{\widehat{\mathcal{P}}}(\overline{\sigma}) \widehat{f}(s, \overline{\eta}) \widehat{f}(s, \overline{\sigma}) ~ d\overline{\eta} ds \label{estapriorimauvaisbCBBBP-s} \\
&\quad \quad + \int_0^t \int e^{i s \varphi} i \xi_0 \mu_{BBB}(\overline{\xi}, \overline{\eta}) m_{A}(\overline{\eta}) m_{\widehat{\mathcal{P}}}(\overline{\sigma}) \widehat{f}(s, \overline{\eta}) m_{\widehat{\mathcal{C}}}(\overline{\xi}) \widehat{X}_b(\overline{\xi}) \cdot \nabla_{\overline{\xi}} \widehat{f}(s, \overline{\sigma}) ~ d\overline{\eta} ds \label{estapriorimauvaisbCBBBP-eta} \\
&\quad \quad + \int_0^t \int e^{i s \varphi} m_{\widehat{\mathcal{C}}}(\overline{\xi}) \widehat{X}_b(\overline{\xi}) \cdot \nabla_{\overline{\xi}} \left( i \xi_0 \mu_{BBB}(\overline{\xi}, \overline{\eta}) m_{A}(\overline{\eta}) m_{\widehat{\mathcal{P}}}(\overline{\sigma}) \right) \widehat{f}(s, \overline{\eta}) \widehat{f}(s, \overline{\sigma}) ~ d\overline{\eta} ds \label{estapriorimauvaisbCBBBP-dersymb} 
\end{align}
\end{subequations}

We can again apply Lemma \ref{lem-non-res-loin0Cone} in order to apply integrations by parts on \eqref{estapriorimauvaisbCBBBP-s}, and we apply a decomposition similar to the generic $BBB$ case between $h_{b*}$ and $g_{b*}$, all the estimates being the same except for the one of the form 
\begin{align}
\int_0^t \int e^{i s \varphi} \mu_1(\overline{\xi}, \overline{\eta}) \nabla_{\overline{\eta}} \widehat{f}(s, \overline{\eta}) ~ \widehat{f}(s, \overline{\sigma}) ~ d\overline{\eta} ds \label{estapriorimauvaisbCBBBP-termesup} 
\end{align}
where $\mu_1$ is a symbol of order $1$ containing all localisation. But we can estimate 
\begin{align*}
\Vert \langle \overline{\xi} \rangle \partial_t \eqref{estapriorimauvaisbCBBBP-termesup} \Vert_{L^2} &\lesssim \Vert m_A(D) \langle \nabla \rangle (x, y) f(t) \Vert_{L^2} \Vert \nabla u(t) \Vert_{L^{\infty}} \lesssim \langle t \rangle^{\frac{5}{48}+202\delta} t^{-\frac{5}{6}} \langle t \rangle^{-\frac{1}{6}+100\delta} \Vert u \Vert_X^2 \\
&\lesssim \langle t \rangle^{\frac{5}{48}+303\delta} \left( t^{-\frac{5}{6}} \langle t \rangle^{-\frac{1}{6}-\delta} \right) \Vert u \Vert_X^2 
\end{align*}
which is enough. 

\subsubsection{Interaction \texorpdfstring{$(BBB, \widehat{\mathcal{P}}\widehat{\mathcal{C}})$}{(BBB, PC)}}

Let us localise by $\mu_{BBB}(\overline{\xi}, \overline{\eta}) m_{\widehat{\mathcal{P}}}(\overline{\eta}) m_{\widehat{\mathcal{C}}}(\overline{\sigma})$. 

We now develop {\footnotesize 
\begin{subequations}
\begin{align}
&m_{\widehat{\mathcal{C}}}(\overline{\xi}) |\overline{\xi}|^k \widehat{X}_b(\overline{\xi}) \cdot \nabla_{\overline{\xi}} \left( i \xi_0 \widehat{I}_{\mu_{BBB} m_{\widehat{\mathcal{P}}} m_{\widehat{\mathcal{C}}}}[f, f](t, \overline{\xi}) \right) \notag \\
&\quad = \int_0^t \int e^{i s \varphi} i s m_{\widehat{\mathcal{C}}}(\overline{\xi}) |\overline{\xi}|^k \widehat{X}_b(\overline{\xi}) \cdot \nabla_{\overline{\xi}} \varphi ~ i \xi_0 \mu_{BBB}(\overline{\xi}, \overline{\eta}) m_{\widehat{\mathcal{P}}}(\overline{\eta}) m_{\widehat{\mathcal{C}}}(\overline{\sigma}) \widehat{f}(s, \overline{\eta}) \widehat{f}(s, \overline{\sigma}) ~ d\overline{\eta} ds \label{estapriorimauvaisbCBBBPC-s} \\
&\quad \quad + \int_0^t \int e^{i s \varphi} i \xi_0 \mu_{BBB}(\overline{\xi}, \overline{\eta}) m_{\widehat{\mathcal{P}}}(\overline{\eta}) m_{\widehat{\mathcal{C}}}(\overline{\sigma}) \widehat{f}(s, \overline{\eta}) m_{\widehat{\mathcal{C}}}(\overline{\xi}) |\overline{\xi}|^k \widehat{X}_b(\overline{\xi}) \cdot \nabla_{\overline{\xi}} \widehat{f}(s, \overline{\sigma}) ~ d\overline{\eta} ds \label{estapriorimauvaisbCBBBPC-eta} \\
&\quad \quad + \int_0^t \int e^{i s \varphi} m_{\widehat{\mathcal{C}}}(\overline{\xi}) |\overline{\xi}|^k \widehat{X}_b(\overline{\xi}) \cdot \nabla_{\overline{\xi}} \left( i \xi_0 \mu_{BBB}(\overline{\xi}, \overline{\eta}) m_{\widehat{\mathcal{P}}}(\overline{\eta}) m_{\widehat{\mathcal{C}}}(\overline{\sigma}) \right) \widehat{f}(s, \overline{\eta}) \widehat{f}(s, \overline{\sigma}) ~ d\overline{\eta} ds \label{estapriorimauvaisbCBBBPC-dersymb} 
\end{align}
\end{subequations} }

We can separate \eqref{estapriorimauvaisbCBBBPC-eta} into two pieces: {\footnotesize 
\begin{subequations}
\begin{align}
\eqref{estapriorimauvaisbCBBBPC-eta} &= \int_0^t \int e^{i s \varphi} i \mu_{BBB}(\overline{\xi}, \overline{\eta}) m_{\widehat{\mathcal{P}}}(\overline{\eta}) m_{\widehat{\mathcal{C}}}(\overline{\sigma}) \widehat{f}(s, \overline{\eta}) \frac{m_{\widehat{\mathcal{C}}}(\overline{\sigma}) + m_{\widehat{\mathcal{C}}}(\overline{\xi})}{2} \frac{\xi_0 + \sigma_0}{2} |\overline{\sigma}|^k \widehat{X}_b(\overline{\sigma}) \cdot \nabla_{\overline{\xi}} \widehat{f}(s, \overline{\sigma}) ~ d\overline{\eta} ds \label{estapriorimauvaisbCBBBPC-eta-sym} \\
&\quad \begin{aligned}
+ \int_0^t \int e^{i s \varphi} i \left( \xi_0 m_{\widehat{\mathcal{C}}}(\overline{\xi}) |\overline{\xi}|^k \widehat{X}_b(\overline{\xi}) - \frac{\xi_0 + \sigma_0}{2} \frac{m_{\widehat{\mathcal{C}}}(\overline{\sigma}) + m_{\widehat{\mathcal{C}}}(\overline{\xi})}{2} |\overline{\sigma}|^k \widehat{X}_b(\overline{\sigma}) \right) \cdot \nabla_{\overline{\xi}} \widehat{f}(s, \overline{\sigma}) \\
\mu_{BBB}(\overline{\xi}, \overline{\eta}) m_{\widehat{\mathcal{P}}}(\overline{\eta}) m_{\widehat{\mathcal{C}}}(\overline{\sigma}) \widehat{f}(s, \overline{\eta}) ~ d\overline{\eta} ds 
\end{aligned} \label{estapriorimauvaisbCBBBPC-eta-2} 
\end{align}
\end{subequations} }
and apply an integration by parts on \eqref{estapriorimauvaisbCBBBPC-eta-2}: {\footnotesize 
\begin{subequations}
\begin{align}
&\eqref{estapriorimauvaisbCBBBPC-eta-2} \notag \\
&\begin{aligned}
&= \int_0^t \int e^{i s \varphi} i s \left( \xi_0 m_{\widehat{\mathcal{C}}}(\overline{\xi}) |\overline{\xi}|^k \widehat{X}_b(\overline{\xi}) - \frac{\xi_0 + \sigma_0}{2} \frac{m_{\widehat{\mathcal{C}}}(\overline{\sigma}) + m_{\widehat{\mathcal{C}}}(\overline{\xi})}{2} |\overline{\sigma}|^k \widehat{X}_b(\overline{\sigma}) \right) \cdot \nabla_{\overline{\eta}} \varphi \\
&\pushright{ ~ i \mu_{BBB}(\overline{\xi}, \overline{\eta}) m_{\widehat{\mathcal{P}}}(\overline{\eta}) m_{\widehat{\mathcal{C}}}(\overline{\sigma}) \widehat{f}(s, \overline{\eta}) \widehat{f}(s, \overline{\sigma}) ~ d\overline{\eta} ds}
\end{aligned} \label{estapriorimauvaisbCBBBPC-eta-2-1} \\
&\begin{aligned} 
&+ \int_0^t \int e^{i s \varphi} i \mu_{BBB}(\overline{\xi}, \overline{\eta}) m_{\widehat{\mathcal{P}}}(\overline{\eta}) m_{\widehat{\mathcal{C}}}(\overline{\sigma}) \\
&\pushright{ \left( \xi_0 m_{\widehat{\mathcal{C}}}(\overline{\xi}) |\overline{\xi}|^k \widehat{X}_b(\overline{\xi}) - \frac{\xi_0 + \sigma_0}{2} \frac{m_{\widehat{\mathcal{C}}}(\overline{\sigma}) + m_{\widehat{\mathcal{C}}}(\overline{\xi})}{2} |\overline{\sigma}|^k \widehat{X}_b(\overline{\sigma}) \right) \cdot \nabla_{\overline{\eta}} \widehat{f}(s, \overline{\eta}) \widehat{f}(s, \overline{\sigma}) ~ d\overline{\eta} ds}
\end{aligned} \label{estapriorimauvaisbCBBBPC-eta-2-2} \\
&
+ \int_0^t \int \nabla_{\overline{\eta}} \cdot \left( \left( \xi_0 m_{\widehat{\mathcal{C}}}(\overline{\xi}) |\overline{\xi}|^k \widehat{X}_b(\overline{\xi}) - \frac{\xi_0 + \sigma_0}{2} \frac{m_{\widehat{\mathcal{C}}}(\overline{\sigma}) + m_{\widehat{\mathcal{C}}}(\overline{\xi})}{2} |\overline{\sigma}|^k \widehat{X}_b(\overline{\sigma}) \right) i \mu_{BBB}(\overline{\xi}, \overline{\eta}) m_{\widehat{\mathcal{P}}}(\overline{\eta}) m_{\widehat{\mathcal{C}}}(\overline{\sigma}) \right) \notag \\
&\quad \quad \quad \quad \quad \quad \quad \quad \quad \quad \quad \quad \quad \quad \quad \quad \quad \quad e^{i s \varphi} \widehat{f}(s, \overline{\eta}) \widehat{f}(s, \overline{\sigma}) ~ d\overline{\eta} ds \label{estapriorimauvaisbCBBBPC-eta-2-3} 
\end{align}
\end{subequations} }

We now group \eqref{estapriorimauvaisbCBBBPC-s} with \eqref{estapriorimauvaisbCBBBPC-eta-2-1}: {\footnotesize 
\begin{align}
&\eqref{estapriorimauvaisbCBBBPC-s} + \eqref{estapriorimauvaisbCBBBPC-eta-2-1} \notag \\
&\quad \begin{aligned}
&= \int_0^t \int e^{i s \varphi} i s \left( \xi_0 m_{\widehat{\mathcal{C}}}(\overline{\xi}) |\overline{\xi}|^k \widehat{X}_b(\overline{\xi}) \cdot \left( \nabla_{\overline{\xi}} + \nabla_{\overline{\eta}} \right) - \frac{\xi_0 + \sigma_0}{2} \frac{m_{\widehat{\mathcal{C}}}(\overline{\sigma}) + m_{\widehat{\mathcal{C}}}(\overline{\xi})}{2} |\overline{\sigma}|^k \widehat{X}_b(\overline{\sigma}) \cdot \nabla_{\overline{\eta}} \right) \varphi \\
&\pushright{~ i \mu_{BBB}(\overline{\xi}, \overline{\eta}) m_{\widehat{\mathcal{P}}}(\overline{\eta}) m_{\widehat{\mathcal{C}}}(\overline{\sigma}) \widehat{f}(s, \overline{\eta}) \widehat{f}(s, \overline{\sigma}) ~ d\overline{\eta} ds} 
\end{aligned} \label{estapriorimauvaisbCBBBPC-s-final} 
\end{align} }
We can now compute that 
\begin{align*}
\nabla_{\eta} \varphi &= 2 \sigma_0 \sigma - 2 \eta_0 \eta 
\end{align*}
But here $|\eta_0| = o(1)$ and therefore $\sqrt{3} |\xi_0| = |\xi| + o(1) = \sqrt{3} |\sigma_0| + o(1) = |\sigma| + o(1)$, so $|\overline{\xi}| = |\overline{\sigma}| + o(1)$, hence $|\overline{\eta}| \leq 3 |\overline{\sigma}|$ (where the constant is explicit). In particular, $|\eta_0 \eta| \ll |\sigma_0 \sigma|$, so $1 = O(\nabla_{\eta} \varphi)$. On the other hand, 
\begin{align*}
\varphi &= \xi_0^3 + \xi_0 |\xi|^2 - \eta_0^3 - \eta_0 |\eta|^2 - \sigma_0^3 - \sigma_0 |\sigma|^2 \\
&= O(\eta_0) + \xi_0 (|\xi|^2 - |\sigma|^2) 
\end{align*}
so $|\xi| - |\sigma| = O(\varphi) + O(\eta_0) = O(\varphi) + O(\eta_0 \nabla_{\eta} \varphi)$. Finally, {\footnotesize 
\begin{align*}
&\left( \xi_0 m_{\widehat{\mathcal{C}}}(\overline{\xi}) |\overline{\xi}|^k \widehat{X}_b(\overline{\xi}) \cdot \left( \nabla_{\overline{\xi}} + \nabla_{\overline{\eta}} \right) - \frac{\xi_0 + \sigma_0}{2} \frac{m_{\widehat{\mathcal{C}}}(\overline{\sigma}) + m_{\widehat{\mathcal{C}}}(\overline{\xi})}{2} |\overline{\sigma}|^k \widehat{X}_b(\overline{\sigma}) \cdot \nabla_{\overline{\eta}} \right) \varphi \\
&\quad = O(\eta_0) + O\left( m_{\widehat{\mathcal{C}}}(\overline{\xi}) - m_{\widehat{\mathcal{C}}}(\overline{\sigma}) \right) + O\left( |\overline{\xi}| - |\overline{\sigma}| \right) 
+ \xi_0 m_{\widehat{\mathcal{C}}}(\overline{\xi}) |\overline{\xi}|^k \left( \widehat{X}_b(\overline{\xi}) \cdot \left( \nabla_{\overline{\xi}} + \nabla_{\overline{\eta}} \right) - \widehat{X}_b(\overline{\sigma}) \cdot \nabla_{\overline{\eta}} \right) \varphi \\
&\quad = O(\eta_0) + O(|\xi| - |\sigma|) + \xi_0 m_{\widehat{\mathcal{C}}}(\overline{\xi}) |\overline{\xi}|^k \Biggl( \frac{|\xi|}{|\overline{\xi}|} \left( 3 \xi_0^2 + |\xi|^2 - 3 \eta_0^2 - |\eta|^2 \right) - \frac{\xi_0 \xi}{|\overline{\xi}| |\xi|} \cdot \left( 2 \xi_0 \xi - 2 \eta_0 \eta \right) \\
&\pushright{- \frac{|\sigma|}{|\overline{\sigma}|} \left( 3 \sigma_0^2 + |\sigma|^2 - 3 \eta_0^2 - |\eta|^2 \right) 
+ \frac{\sigma_0 \sigma}{|\overline{\sigma}| |\sigma|} \cdot \left( 2 \sigma_0 \sigma - 2 \eta_0 \eta \right) \Biggl)} \\
&\quad = O(\eta_0) + O(|\xi| - |\sigma|) \\
&\quad = O(\eta_0 \nabla_{\eta} \varphi) + O(\varphi)
\end{align*} }

We can therefore apply on \eqref{estapriorimauvaisbCBBBPC-s-final} integrations by parts in $\eta$ or in $s$: thanks to the presence of the factor $\eta_0$, in the case of an integration by parts in $\eta$, we only get terms we can estimate like \eqref{estapriorimauvaisbCgenBBB-eta}, and the terms obtained by integration by parts in time are similar to those obtained in the generic $BBB$ case. We use the same decomposition between $h_{b*}$ and $g_{b*}$. 

Likewise, \eqref{estapriorimauvaisbCBBBPC-dersymb}, \eqref{estapriorimauvaisbCBBBPC-eta-2-2} and \eqref{estapriorimauvaisbCBBBPC-eta-2-3} can be estimated in a similar way as \eqref{estapriorimauvaisbCgenBBB-dersymb} or \eqref{estapriorimauvaisbCgenBBB-eta}. 

We now only have to treat the symmetric term \eqref{estapriorimauvaisbCBBBPC-eta-sym}, in which we separate again {\footnotesize 
\begin{subequations}
\begin{align}
\eqref{estapriorimauvaisbCBBBPC-eta-sym} &= \int_0^t \int e^{i s \varphi} i \mu_{BBB}(\overline{\xi}, \overline{\eta}) m_{\widehat{\mathcal{P}}}(\overline{\eta}) m_{\widehat{\mathcal{C}}}(\overline{\sigma}) \widehat{f}(s, \overline{\eta}) \frac{\xi_0 + \sigma_0}{2} \frac{m_{\widehat{\mathcal{C}}}(\overline{\sigma}) + m_{\widehat{\mathcal{C}}}(\overline{\xi})}{2} |\overline{\sigma}|^k \widehat{h}_{b*}(s, \overline{\sigma}) ~ d\overline{\eta} ds \label{estapriorimauvaisbCBBBPC-eta-symh} \\
&\quad + \int_0^t \int e^{i s \varphi} i \mu_{BBB}(\overline{\xi}, \overline{\eta}) m_{\widehat{\mathcal{P}}}(\overline{\eta}) m_{\widehat{\mathcal{C}}}(\overline{\sigma}) \widehat{f}(s, \overline{\eta}) \frac{\xi_0 + \sigma_0}{2} \frac{m_{\widehat{\mathcal{C}}}(\overline{\sigma}) + m_{\widehat{\mathcal{C}}}(\overline{\xi})}{2} |\overline{\sigma}|^k \widehat{g}_{b*}(s, \overline{\sigma}) ~ d\overline{\eta} ds \label{estapriorimauvaisbCBBBPC-eta-symg} 
\end{align}
\end{subequations} }
On the one hand, 
\begin{align*}
\Vert \partial_t \eqref{estapriorimauvaisbCBBBPC-eta-symg}_{k = 0} \Vert_{L^2} &\lesssim \Vert u \Vert_{L^4} \Vert e^{-it \omega(D)} m_{\widehat{\mathcal{C}}}(D) \nabla g_{b*}(t) \Vert_{L^4} \\
&\lesssim t^{-\frac{5}{6}} \langle t \rangle^{-\frac{1}{4}+100\delta} \langle t \rangle^{\frac{1}{2}+101\delta} \Vert u \Vert_X^2 \\
\Vert \partial_t \eqref{estapriorimauvaisbCBBBPC-eta-symg}_{k = 1} \Vert_{L^2} &\lesssim \Vert u \Vert_{L^4} \Vert e^{-it \omega(D)} m_{\widehat{\mathcal{C}}}(D) \nabla^2 g_{b*}(t) \Vert_{L^4} \\
&\lesssim t^{-\frac{5}{6}} \langle t \rangle^{-\frac{1}{4}+100\delta} \langle t \rangle^{\frac{1}{6}+101\delta} \Vert u \Vert_X^2 
\end{align*}

On the other hand, we cannot estimate directly the time derivative of \eqref{estapriorimauvaisbCBBBPC-eta-symh} so we get back to \eqref{estapriorimauvaisbDuhamelhcarre}: the corresponding contribution is {\footnotesize 
\begin{align}
\int_0^t \int \int |\overline{\xi}|^k m_{\widehat{\mathcal{C}}}(\overline{\xi}) \widehat{h}_{b*}(s, -\overline{\xi}) ~ e^{i s \varphi} i \mu_{BBB}(\overline{\xi}, \overline{\eta}) m_{\widehat{\mathcal{P}}}(\overline{\eta}) m_{\widehat{\mathcal{C}}}(\overline{\sigma}) \widehat{f}(s, \overline{\eta}) \frac{\xi_0 + \sigma_0}{2} \frac{m_{\widehat{\mathcal{C}}}(\overline{\sigma}) + m_{\widehat{\mathcal{C}}}(\overline{\xi})}{2} |\overline{\sigma}|^k \widehat{h}_{b*}(s, \overline{\sigma}) ~ d\overline{\eta} d\overline{\xi} ds \label{estapriorimauvaisbCBBBPC-eta-symh-Duhamel} 
\end{align} }
on which we can apply the change of variables $\mathfrak{s} : (\overline{\xi}, \overline{\eta}) \mapsto (\overline{\eta}-\overline{\xi}, \overline{\eta})$, that can be rewritten in the variables $(\overline{\xi}, \overline{\sigma})$ as $(\overline{\xi}, \overline{\sigma}) \mapsto (-\overline{\sigma}, -\overline{\xi})$. We then have: {\footnotesize 
\begin{align*}
&\eqref{estapriorimauvaisbCBBBPC-eta-symh-Duhamel} \\
&= - \int_0^t \int \int |\overline{\sigma}|^k m_{\widehat{\mathcal{C}}}(\overline{\sigma}) \widehat{h}_{b*}(s, \overline{\sigma}) ~ e^{i s \varphi} i \mu_{BBB}(\overline{\xi}, \overline{\eta}) m_{\widehat{\mathcal{P}}}(\overline{\eta}) m_{\widehat{\mathcal{C}}}(\overline{\xi}) \widehat{f}(s, \overline{\eta}) \frac{\xi_0 + \sigma_0}{2} \frac{m_{\widehat{\mathcal{C}}}(\overline{\sigma}) + m_{\widehat{\mathcal{C}}}(\overline{\xi})}{2} |\overline{\xi}|^k \widehat{h}_{b*}(s, - \overline{\xi}) ~ d\overline{\eta} d\overline{\xi} ds \\
&= - \eqref{estapriorimauvaisbCBBBPC-eta-symh-Duhamel} 
\end{align*} }
where we used that $m_{\widehat{\mathcal{C}}}$ is invariant by $\overline{\xi} \mapsto -\overline{\xi}$, and that $\mu_{BBB}$ can be chosen invariant by $\mathfrak{s}$ as well (more generally, $\mu_{BBB}$ can be chosen to be independant of the ordering of the variables and depends only on their norms). In particular, 
\begin{align*}
\eqref{estapriorimauvaisbCBBBPC-eta-symh-Duhamel} = 0 
\end{align*}
which concludes the estimate. 

\subsubsection{Interaction \texorpdfstring{$(BBB, \widehat{\mathcal{P}}\widehat{\mathcal{P}})$}{(BBB, PP)}} 

Let us localise by $\mu_{BBB}(\overline{\xi}, \overline{\eta}) m_{\widehat{\mathcal{P}}}(\overline{\eta}) m_{\widehat{\mathcal{P}}}(\overline{\sigma})$. 

We then have
\begin{align*}
|\overline{\eta}| + |\overline{\sigma}| \simeq |\overline{\xi}| \simeq |\xi_0| \lesssim |\eta_0| + |\sigma_0| 
\end{align*}
so either $|\sigma_0| \gtrsim |\overline{\sigma}|$ or $|\eta_0| \gtrsim |\overline{\eta}|$ (with constants potentially larger than those in the definition of $m_{\widehat{\mathcal{P}}}$). Up to separating into two subcases and use a localisation function $\widetilde{m}_{\widehat{\mathcal{R}}}$ with larger support than $m_{\widehat{\mathcal{R}}}$, we recover the case $(BBB, \widehat{\mathcal{R}}\widehat{\mathcal{P}})$. 

\subsubsection{Generic interaction \texorpdfstring{$HBH$}{HBH}}

Let us localise by $\mu_{HBH}(\overline{\xi}, \overline{\eta}) m_A(\overline{\sigma})$, $A \neq \widehat{\mathcal{C}}$. 

Then, for $I_{\mu_{HBH} m_A}$ not to be identically zero, we must have $A = \widehat{\mathcal{R}}$. Furthermore, 
\begin{align*}
|\overline{\eta}| + \left| \sqrt{3} |\xi_0| - |\xi| \right| \gtrsim \left| \sqrt{3} |\sigma_0| - |\sigma| \right| \gtrsim |\overline{\sigma}| 
\end{align*}
Therefore, either $|\overline{\eta}| \gtrsim |\overline{\sigma}|$ and we recover the case $BBB$ (up to choosing a localisation $\widetilde{\mu}_{BBB}$ with larger support), or $1 = O(m_b(\overline{\xi}))$ and we recover the non-singular weight $b$, for which the estimates are simpler and the decomposition follows from the decomposition lemma. 

\subsubsection{Interaction \texorpdfstring{$(HBH, \widehat{\mathcal{C}})$}{(HBH, C)}}

Let us localise by $\mu_{HBH}(\overline{\xi}, \overline{\eta}) m_{\widehat{\mathcal{C}}}(\overline{\sigma})$. 

We can then develop {\footnotesize 
\begin{subequations}
\begin{align}
&m_{\widehat{\mathcal{C}}}(\overline{\xi}) |\overline{\xi}|^k \widehat{X}_b(\overline{\xi}) \cdot \nabla_{\overline{\xi}} \left( i \xi_0 \widehat{I}_{\mu_{HBH} m_{\widehat{\mathcal{C}}}}[f, f](t, \overline{\xi}) \right) \notag \\
&\quad = \int_0^t \int e^{i s \varphi} i s m_{\widehat{\mathcal{C}}}(\overline{\xi}) |\overline{\xi}|^k \widehat{X}_b(\overline{\xi}) \cdot \nabla_{\overline{\xi}} \varphi ~ i \xi_0 \mu_{HBH}(\overline{\xi}, \overline{\eta}) m_{\widehat{\mathcal{C}}}(\overline{\sigma}) \widehat{f}(s, \overline{\eta}) \widehat{f}(s, \overline{\sigma}) ~ d\overline{\eta} ds \label{estapriorimauvaisbCHBHC-s} \\
&\quad \quad + \int_0^t \int e^{i s \varphi} i \xi_0 \mu_{HBH}(\overline{\xi}, \overline{\eta}) m_{\widehat{\mathcal{C}}}(\overline{\sigma}) \widehat{f}(s, \overline{\eta}) m_{\widehat{\mathcal{C}}}(\overline{\xi}) |\overline{\xi}|^k \widehat{X}_b(\overline{\xi}) \cdot \nabla_{\overline{\xi}} \widehat{f}(s, \overline{\sigma}) ~ d\overline{\eta} ds \label{estapriorimauvaisbCHBHC-eta} \\
&\quad \quad + \int_0^t \int e^{i s \varphi} m_{\widehat{\mathcal{C}}}(\overline{\xi}) |\overline{\xi}|^k \widehat{X}_b(\overline{\xi}) \cdot \nabla_{\overline{\xi}} \left( i \xi_0 \mu_{HBH}(\overline{\xi}, \overline{\eta}) m_{\widehat{\mathcal{C}}}(\overline{\sigma}) \right) \widehat{f}(s, \overline{\eta}) \widehat{f}(s, \overline{\sigma}) ~ d\overline{\eta} ds \label{estapriorimauvaisbCHBHC-dersymb} 
\end{align}
\end{subequations} }

We now separate \eqref{estapriorimauvaisbCHBHC-eta} (in a similar way as in the case $(BBB, \widehat{\mathcal{P}}\widehat{\mathcal{C}})$): {\footnotesize 
\begin{subequations}
\begin{align}
&\eqref{estapriorimauvaisbCHBHC-eta} = \int_0^t \int e^{i s \varphi} i \mu_{HBH}(\overline{\xi}, \overline{\eta}) m_{\widehat{\mathcal{C}}}(\overline{\sigma}) \widehat{f}(s, \overline{\eta}) \frac{m_{\widehat{\mathcal{C}}}(\overline{\xi}) + m_{\widehat{\mathcal{C}}}(\overline{\sigma})}{2} \frac{\xi_0+\sigma_0}{2} |\overline{\sigma}|^k \widehat{h}_{b*}(s, \overline{\sigma}) ~ d\overline{\eta} ds \label{estapriorimauvaisbCHBHC-eta-symh} \\
&\quad + \int_0^t \int e^{i s \varphi} i \mu_{HBH}(\overline{\xi}, \overline{\eta}) m_{\widehat{\mathcal{C}}}(\overline{\sigma}) \widehat{f}(s, \overline{\eta}) \frac{m_{\widehat{\mathcal{C}}}(\overline{\xi}) + m_{\widehat{\mathcal{C}}}(\overline{\sigma})}{2} \frac{\xi_0+\sigma_0}{2} |\overline{\sigma}|^k \widehat{g}_{b*}(s, \overline{\sigma}) ~ d\overline{\eta} ds \label{estapriorimauvaisbCHBHC-eta-symg} \\
&\quad \begin{aligned}
+ \int_0^t \int e^{i s \varphi} i \left( m_{\widehat{\mathcal{C}}}(\overline{\xi}) \xi_0 |\overline{\xi}|^k \widehat{X}_b(\overline{\xi}) - \frac{m_{\widehat{\mathcal{C}}}(\overline{\xi}) + m_{\widehat{\mathcal{C}}}(\overline{\sigma})}{2} \frac{\xi_0+\sigma_0}{2} |\overline{\sigma}|^k \widehat{X}_b(\overline{\sigma}) \right) \cdot \nabla_{\overline{\xi}} \widehat{f}(s, \overline{\sigma}) \\
\mu_{HBH}(\overline{\xi}, \overline{\eta}) m_{\widehat{\mathcal{C}}}(\overline{\sigma}) \widehat{f}(s, \overline{\eta}) ~ d\overline{\eta} ds 
\end{aligned} \label{estapriorimauvaisbCHBHC-eta-2} 
\end{align}
\end{subequations} }
As in the case $(BBB, \widehat{\mathcal{P}}\widehat{\mathcal{C}})$, choosing $\mu_{HBH}$ invariant by $\mathfrak{s}: (\overline{\xi}, \overline{\eta}) \mapsto (\overline{\eta}-\overline{\xi}, \overline{\eta})$, the contribution of \eqref{estapriorimauvaisbCHBHC-eta-symh} in \eqref{estapriorimauvaisbDuhamelhcarre} is identically zero, while we can estimate \eqref{estapriorimauvaisbCHBHC-eta-symg} like \eqref{estapriorimauvaisbCBBBPC-eta-symg}. Finally, on \eqref{estapriorimauvaisbCHBHC-eta-2}, we apply an integration by parts: {\footnotesize 
\begin{subequations}
\begin{align}
&\begin{aligned}
\eqref{estapriorimauvaisbCHBHC-eta-2} = \int_0^t \int e^{i s \varphi} i s \left( m_{\widehat{\mathcal{C}}}(\overline{\xi}) \xi_0 |\overline{\xi}|^k \widehat{X}_b(\overline{\xi}) - \frac{m_{\widehat{\mathcal{C}}}(\overline{\xi}) + m_{\widehat{\mathcal{C}}}(\overline{\sigma})}{2} \frac{\xi_0+\sigma_0}{2} |\overline{\sigma}|^k \widehat{X}_b(\overline{\sigma}) \right) \cdot \nabla_{\overline{\eta}} \varphi \\
~ i \mu_{HBH}(\overline{\xi}, \overline{\eta}) m_{\widehat{\mathcal{C}}}(\overline{\sigma}) \widehat{f}(s, \overline{\eta}) \widehat{f}(s, \overline{\sigma}) ~ d\overline{\eta} ds 
\end{aligned} \label{estapriorimauvaisbCHBHC-eta-2-1} \\
&\quad \begin{aligned}
+ \int_0^t \int e^{i s \varphi} i \mu_{HBH}(\overline{\xi}, \overline{\eta}) \left( m_{\widehat{\mathcal{C}}}(\overline{\xi}) \xi_0 |\overline{\xi}|^k \widehat{X}_b(\overline{\xi}) - \frac{m_{\widehat{\mathcal{C}}}(\overline{\xi}) + m_{\widehat{\mathcal{C}}}(\overline{\sigma})}{2} \frac{\xi_0+\sigma_0}{2} |\overline{\sigma}|^k \widehat{X}_b(\overline{\sigma}) \right) \cdot \nabla_{\overline{\eta}} \widehat{f}(s, \overline{\eta}) \\
m_{\widehat{\mathcal{C}}}(\overline{\sigma}) \widehat{f}(s, \overline{\sigma}) ~ d\overline{\eta} ds 
\end{aligned} \label{estapriorimauvaisbCHBHC-eta-2-2} \\
&\quad \begin{aligned}
+ \int_0^t \int \nabla_{\overline{\eta}} \cdot \left( \left( m_{\widehat{\mathcal{C}}}(\overline{\xi}) \xi_0 |\overline{\xi}|^k \widehat{X}_b(\overline{\xi}) - \frac{m_{\widehat{\mathcal{C}}}(\overline{\xi}) + m_{\widehat{\mathcal{C}}}(\overline{\sigma})}{2} \frac{\xi_0+\sigma_0}{2} |\overline{\sigma}|^k \widehat{X}_b(\overline{\sigma}) \right) i \mu_{HBH}(\overline{\xi}, \overline{\eta}) m_{\widehat{\mathcal{C}}}(\overline{\sigma}) \right) \\
e^{i s \varphi} \widehat{f}(s, \overline{\eta}) \widehat{f}(s, \overline{\sigma}) ~ d\overline{\eta} ds
\end{aligned} \label{estapriorimauvaisbCHBHC-eta-2-3}
\end{align}
\end{subequations} }
We now group \eqref{estapriorimauvaisbCHBHC-s} with \eqref{estapriorimauvaisbCHBHC-eta-2-1}: {\footnotesize 
\begin{align}
&\eqref{estapriorimauvaisbCHBHC-s} + \eqref{estapriorimauvaisbCHBHC-eta-2-1} \notag \\
&\quad \begin{aligned}
&= \int_0^t \int e^{i s \varphi} i s \left( m_{\widehat{\mathcal{C}}}(\overline{\xi}) \xi_0 |\overline{\xi}|^k \widehat{X}_b(\overline{\xi}) \cdot \left( \nabla_{\overline{\xi}} + \nabla_{\overline{\eta}} \right) - \frac{m_{\widehat{\mathcal{C}}}(\overline{\xi}) + m_{\widehat{\mathcal{C}}}(\overline{\sigma})}{2} \frac{\xi_0+\sigma_0}{2} |\overline{\sigma}|^k \widehat{X}_b(\overline{\sigma}) \cdot \nabla_{\overline{\eta}} \right) \varphi \\
&\pushright{ ~ i \mu_{HBH}(\overline{\xi}, \overline{\eta}) m_{\widehat{\mathcal{C}}}(\overline{\sigma}) \widehat{f}(s, \overline{\eta}) \widehat{f}(s, \overline{\sigma}) ~ d\overline{\eta} ds} \end{aligned} \label{estapriorimauvaisbCHBHC-s-final} 
\end{align} }

We can estimate \eqref{estapriorimauvaisbCHBHC-dersymb} and \eqref{estapriorimauvaisbCHBHC-eta-2-3} like \eqref{estapriorimauvaisbCgenBBB-dersymb}. Moreover, in \eqref{estapriorimauvaisbCHBHC-eta-2-2}, 
\begin{align*}
m_{\widehat{\mathcal{C}}}(\overline{\xi}) \xi_0 |\overline{\xi}|^k \widehat{X}_b(\overline{\xi}) - \frac{m_{\widehat{\mathcal{C}}}(\overline{\xi}) + m_{\widehat{\mathcal{C}}}(\overline{\sigma})}{2} \frac{\xi_0+\sigma_0}{2} |\overline{\sigma}|^k \widehat{X}_b(\overline{\sigma}) = O(\overline{\eta}) 
\end{align*}
which justifies that we can also estimate \eqref{estapriorimauvaisbCHBHC-eta-2-2} like \eqref{estapriorimauvaisbCgenBBB-eta} (using the derivative $\overline{\eta}$ in the case $k = 1$). 

Finally, in order to apply integrations by parts on \eqref{estapriorimauvaisbCHBHC-s-final}, we have 
\begin{align*}
\nabla_{\eta} \varphi &= 2 \sigma_0 \sigma - 2 \eta_0 \eta 
\end{align*}
so $1 = O(\nabla_{\eta} \varphi)$, and then 
\begin{align*}
\varphi &= O(\eta_0) + \xi_0 (|\xi|^2 - |\sigma|^2) 
\end{align*}
Therefore, $\eta_0 = O(\eta_0 \nabla_{\eta} \varphi)$ and $|\xi| - |\sigma| = O(\eta_0 \nabla_{\eta} \varphi) + O(\varphi)$. On the other hand, by a computation analogous to one from the case $(BBB, \widehat{\mathcal{P}}\widehat{\mathcal{C}})$, 
\begin{align*}
&\left( m_{\widehat{\mathcal{C}}}(\overline{\xi}) \xi_0 |\overline{\xi}|^k \widehat{X}_b(\overline{\xi}) \cdot \left( \nabla_{\overline{\xi}} + \nabla_{\overline{\eta}} \right) - \frac{m_{\widehat{\mathcal{C}}}(\overline{\xi}) + m_{\widehat{\mathcal{C}}}(\overline{\sigma})}{2} \frac{\xi_0+\sigma_0}{2} |\overline{\sigma}|^k \widehat{X}_b(\overline{\sigma}) \cdot \nabla_{\overline{\eta}} \right) \varphi \\
&\quad = O(\eta_0) + O(|\xi| - |\sigma|) = O(\eta_0 \nabla_{\eta} \varphi) + O(\varphi) 
\end{align*}
We can therefore apply either an integration by parts in time and estimate as in the generic $BBB$ case, or apply an integration by parts in $\eta$ keeping a factor $O(\eta_0)$, which allows again to estimate like in the $BBB$ case. 

\subsubsection{Generic interaction \texorpdfstring{$BHH$}{BHH}}

Let us localise by $\mu_{BHH}(\overline{\xi}, \overline{\eta}) m_{A_1}(\overline{\eta}) m_{A_2}(\overline{\sigma})$, $(A_1, A_2) \in \{ \widehat{\mathcal{R}}, \widehat{\mathcal{C}}, \widehat{\mathcal{L}}, \widehat{\mathcal{P}} \}^2 \setminus \{ (\widehat{\mathcal{C}}, \widehat{\mathcal{C}}) \}$. By symmetry of $\overline{\eta}, \overline{\sigma}$, we can assume $A_1 \neq \widehat{\mathcal{C}}$. Moreover, either $A_1 = A_2$, or one of them is $\widehat{\mathcal{R}}$. 

We can now develop {\footnotesize 
\begin{subequations}
\begin{align}
&m_{\widehat{\mathcal{C}}}(\overline{\xi}) \widehat{X}_b(\overline{\xi}) \cdot \nabla_{\overline{\xi}} \left( i \xi_0 \widehat{I}_{\mu_{BHH} m_{A_1} m_{A_2}}[f, f](t, \overline{\xi}) \right) \notag \\
&\quad = \int_0^t \int e^{i s \varphi} i s m_{\widehat{\mathcal{C}}}(\overline{\xi}) \widehat{X}_b(\overline{\xi}) \cdot \nabla_{\overline{\xi}} \varphi ~ i \xi_0 \mu_{BHH}(\overline{\xi}, \overline{\eta}) m_{A_1}(\overline{\eta}) m_{A_2}(\overline{\sigma}) \widehat{f}(s, \overline{\eta}) \widehat{f}(s, \overline{\sigma}) ~ d\overline{\eta} ds \label{estapriorimauvaisbCBHHgen-s} \\
&\quad \quad + \int_0^t \int e^{i s \varphi} i \xi_0 \mu_{BHH}(\overline{\xi}, \overline{\eta}) m_{A_1}(\overline{\eta}) m_{A_2}(\overline{\sigma}) \widehat{f}(s, \overline{\eta}) m_{\widehat{\mathcal{C}}}(\overline{\xi}) \widehat{X}_b(\overline{\xi}) \cdot \nabla_{\overline{\xi}} \widehat{f}(s, \overline{\sigma}) ~ d\overline{\eta} ds \label{estapriorimauvaisbCBHHgen-eta} \\
&\quad \quad + \int_0^t \int e^{i s \varphi} m_{\widehat{\mathcal{C}}}(\overline{\xi}) \widehat{X}_b(\overline{\xi}) \cdot \nabla_{\overline{\xi}} \left( i \xi_0 \mu_{BHH}(\overline{\xi}, \overline{\eta}) m_{A_1}(\overline{\eta}) m_{A_2}(\overline{\sigma}) \right) \widehat{f}(s, \overline{\eta}) \widehat{f}(s, \overline{\sigma}) ~ d\overline{\eta} ds \label{estapriorimauvaisbCBHHgen-dersymb} 
\end{align}
\end{subequations} }
Here, we can estimate \eqref{estapriorimauvaisbCBHHgen-dersymb} like \eqref{estapriorimauvaisbCgenBBB-dersymb}. Furthermore, we can estimate \eqref{estapriorimauvaisbCBHHgen-eta} in a similar way as \eqref{estapriorimauvaisbCgenBBB-eta}: indeed, either $A_1 \neq \widehat{\mathcal{P}}$ and it is exactly the same estimate, or $A_1 = \widehat{\mathcal{P}}$ and then $A_2 \in \{ \widehat{\mathcal{R}}, \widehat{\mathcal{P}} \}$ so 
\begin{align*}
\Vert \langle \nabla \rangle \partial_t \eqref{estapriorimauvaisbCBHHgen-eta} \Vert_{L^2} &\lesssim \Vert \nabla u(t) \Vert_{L^{\infty}} \Vert \nabla m_{A_2}(D) (x, y) f(t) \Vert_{L^2} \\
&\lesssim t^{-\frac{5}{6}} \langle t \rangle^{-\frac{1}{6}+100\delta} \Vert u \Vert_X^2
\end{align*}
which is enough.  

It thus only remains \eqref{estapriorimauvaisbCBHHgen-s}: but since $A_1 \neq \widehat{\mathcal{C}}$ by hypothesis, we can apply Lemma \ref{lem-non-res-loin0Cone} to write locally 
\begin{align*}
\xi_0 = O(\nabla_{\overline{\eta}} \varphi) + O(\varphi) 
\end{align*}
which is enough to apply sufficient integrations by parts to get terms that can be estimated as above or as in the generic $BBB$ case.  

\subsubsection{Interaction \texorpdfstring{$(BHH, \widehat{\mathcal{C}})$}{(BHH, C)}} 

Let us finally localise by $\mu_{BHH}(\overline{\xi}, \overline{\eta}) m_{\widehat{\mathcal{C}}}(\overline{\eta}) m_{\widehat{\mathcal{C}}}(\overline{\sigma})$. 

We can then develop {\footnotesize 
\begin{subequations}
\begin{align}
&m_{\widehat{\mathcal{C}}}(\overline{\xi}) \widehat{X}_b(\overline{\xi}) \cdot \nabla_{\overline{\xi}} \left( i \xi_0 \widehat{I}_{\mu_{BHH} m_{\widehat{\mathcal{C}}} m_{\widehat{\mathcal{C}}}}[f, f](t, \overline{\xi}) \right) \notag \\
&\quad = \int_0^t \int e^{i s \varphi} i s m_{\widehat{\mathcal{C}}}(\overline{\xi}) \widehat{X}_b(\overline{\xi}) \cdot \nabla_{\overline{\xi}} \varphi ~ i \xi_0 \mu_{BHH}(\overline{\xi}, \overline{\eta}) m_{\widehat{\mathcal{C}}}(\overline{\eta}) m_{\widehat{\mathcal{C}}}(\overline{\sigma}) \widehat{f}(s, \overline{\eta}) \widehat{f}(s, \overline{\sigma}) ~ d\overline{\eta} ds \label{estapriorimauvaisbCBHHC-s} \\
&\quad \quad + \int_0^t \int e^{i s \varphi} i \xi_0 \mu_{BHH}(\overline{\xi}, \overline{\eta}) m_{\widehat{\mathcal{C}}}(\overline{\eta}) m_{\widehat{\mathcal{C}}}(\overline{\sigma}) \widehat{f}(s, \overline{\eta}) m_{\widehat{\mathcal{C}}}(\overline{\xi}) \widehat{X}_b(\overline{\xi}) \cdot \nabla_{\overline{\xi}} \widehat{f}(s, \overline{\sigma}) ~ d\overline{\eta} ds \label{estapriorimauvaisbCBHHC-eta} \\
&\quad \quad + \int_0^t \int e^{i s \varphi} m_{\widehat{\mathcal{C}}}(\overline{\xi}) \widehat{X}_b(\overline{\xi}) \cdot \nabla_{\overline{\xi}} \left( i \xi_0 \mu_{BHH}(\overline{\xi}, \overline{\eta}) m_{\widehat{\mathcal{C}}}(\overline{\eta}) m_{\widehat{\mathcal{C}}}(\overline{\sigma}) \right) \widehat{f}(s, \overline{\eta}) \widehat{f}(s, \overline{\sigma}) ~ d\overline{\eta} ds \label{estapriorimauvaisbCBHHC-dersymb} 
\end{align}
\end{subequations} }
We can estimate \eqref{estapriorimauvaisbCBHHC-eta} and \eqref{estapriorimauvaisbCBHHC-dersymb} like \eqref{estapriorimauvaisbCgenBBB-eta} and \eqref{estapriorimauvaisbCBHHC-dersymb} of the generic $BBB$ case respectively. 

Reusing computations done in the proof of the decomposition lemma, we have 
\begin{align*}
\overline{\eta}_a^{\overline{\eta}} - \overline{\sigma}_a^{\overline{\sigma}} &= O\left( \widehat{Y}(\overline{\eta}, \overline{\sigma}) \cdot \nabla_{\overline{\eta}} \varphi \right) 
\end{align*}
and we know how to estimate any term having a factor $\xi_t^{\overline{\eta}\overline{\sigma}}$ by a succession of integrations by parts (here, $\frac{\xi_0}{|\overline{\xi}|} \simeq 1$). 

In particular, if we localise to have 
\begin{align*}
\xi_0 = O\left( \overline{\eta}_a^{\overline{\eta}} - \overline{\sigma}_a^{\overline{\sigma}} \right) \quad \mbox{ or } \quad \xi_0 = O\left( \xi_t^{\overline{\eta} \overline{\sigma}} \right) 
\end{align*}
then we can apply integrations by parts as in the case of the decomposition lemma and estimate everything in a similar way as in the generic $BBB$ case, we skip the details. 

Up to using a sufficiently fine partition of unity of the sphere $S^2$ to localise $\frac{\overline{\xi}}{|\overline{\xi}|}, \frac{\overline{\eta}}{|\overline{\eta}|}$ (only adding Hörmander-Mikhlin symbols in $\overline{\xi}, \overline{\eta}$, that we will omit in the following to simplify, and a finite sum), we can localise also $\xi_t^{\overline{\xi} \overline{\eta}}$. But if $1 = O(\xi_t^{\overline{\xi} \overline{\eta}})$, then $\xi_0 = O(\overline{\xi}) = O(\overline{\xi} \xi_t^{\overline{\xi} \overline{\eta}}) = O(\xi_t^{\overline{\eta} \overline{\sigma}})$: we may therefore assume $\xi_t^{\overline{\xi} \overline{\eta}} = o(1)$. Then either $\epsilon^{\overline{\xi} \overline{\eta}} \theta^{\overline{\xi} \overline{\eta}}$ is close to $1$, or to $-1$. But 
\begin{align*}
\overline{\eta}_a^{\overline{\eta}} - \overline{\sigma}_a^{\overline{\sigma}} &= \epsilon^{\overline{\xi} \overline{\eta}} \overline{\eta}_a^{\overline{\eta}} + \overline{\sigma}_a^{\overline{\eta}} + O\left( \xi_t^{\overline{\eta} \overline{\sigma}} \right) = \overline{\xi}_a^{\overline{\eta}} + o(\overline{\xi}) 
\end{align*}
If $\epsilon^{\overline{\xi} \overline{\eta}} \theta^{\overline{\xi} \overline{\eta}}$ is close to $1$, then 
\begin{align*}
\overline{\xi}_a^{\overline{\eta}} &= \epsilon^{\overline{\xi} \overline{\eta}} \overline{\xi}_a^{\overline{\xi}} + O\left( |\xi| \xi_t^{\overline{\xi} \overline{\eta}} \right) = \epsilon^{\overline{\xi} \overline{\eta}} \overline{\xi}_a^{\overline{\xi}} + o(\overline{\xi})
\end{align*}
and therefore 
\begin{align*}
\xi_0 = O(\overline{\xi}) = O\left( \overline{\eta}_a^{\overline{\eta}} - \overline{\sigma}_a^{\overline{\sigma}} \right) 
\end{align*}
We can thus reduce to $\epsilon^{\overline{\xi} \overline{\eta}} \theta^{\overline{\xi} \overline{\eta}}$ close to $-1$. 

\begin{figure}
\centering
\begin{subfigure}{0.45\textwidth}
\centering
\begin{tikzpicture}
        \draw[->] (-2.3, 0) -- (2.3, 0) node[above] {$\xi$};
        \draw[->] (0, -0.5) -- (0, 1.5) node[right] {$\xi_0$};
        \draw[->] (0, 0) -- ({0.97*sqrt(3)}, 1.03) node[above] {$\overline{\eta}$};
        \draw[->] ({0.95*sqrt(3)}, 1.05) -- ({-0.13*sqrt(3)}, 0.17) node[above] {$\overline{\sigma}$};
        \draw[->] (0, 0) -- ({-0.15*sqrt(3)}, 0.15) node[left] {$\overline{\xi}$};
        \draw[dotted] ({-0.25*sqrt(3)}, -0.25) -- ({1.25*sqrt(3)}, 1.25) node[right] {$\widehat{\mathcal{C}}$};
        \draw[dotted] ({-0.5*sqrt(3)}, 0.5) -- ({0.25*sqrt(3)}, -0.25);
    \end{tikzpicture}
    \subcaption{Degenerate $(BHH, \widehat{\mathcal{C}}\widehat{\mathcal{C}})$ interaction} \label{figurecasdegeneremaxCbBHHC}
\end{subfigure} \hfill
\begin{subfigure}{0.45\textwidth}
\centering
\begin{tikzpicture}
		\draw[->] (-1., 0) -- (2.8, 0) node[above] {$\xi$};
		\draw[->] (0, -0.5) -- (0, 2.5) node[above] {$\xi_0$};
        \draw[->] (0, 0) -- ({0.98*sqrt(3)}, 0.98) node[below] {$\overline{\eta}$};
        \draw[->] ({sqrt(3)}, 1) -- ({0.1}, 1.98) node[above right] {$\overline{\sigma}$};
        \draw[->] (0, 0) -- (0.08, 2.) node[left] {$\overline{\xi}$};
        \draw[dotted] ({-0.25*sqrt(3)}, -0.25) -- ({1.5*sqrt(3)}, 1.5) node[right] {$\widehat{\mathcal{C}}$};
        \draw[dotted] ({-0.5*sqrt(3)}, 0.5) -- ({0.3*sqrt(3)}, -0.3);
\end{tikzpicture}
\subcaption{Degenerate $(BBB, \widehat{\mathcal{C}}\widehat{\mathcal{C}})$ interaction} \label{figurecasdegeneremaxLbBBBCC} 
\end{subfigure} 
\caption{Degenerate interactions for the singular $b$-weight} \label{figurecasdegenb} 
\end{figure}
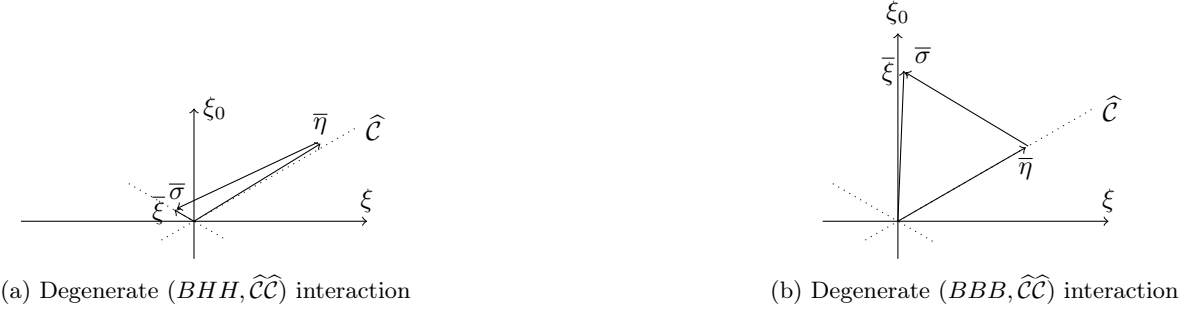

Up to exchanging $\overline{\eta}$ and $\overline{\sigma}$, we will assume $\epsilon^{\overline{\xi} \overline{\eta}} = 1 = -\epsilon^{\overline{\xi} \overline{\sigma}} = - \epsilon^{\overline{\eta} \overline{\sigma}}$. We represent this case on Figure \ref{figurecasdegeneremaxCbBHHC}. 

In this case, 
\begin{align*}
\overline{\xi}_b^{\overline{\xi}} &= \overline{\xi}_a^{\overline{\eta}} + O\left( |\xi| \left( \xi_t^{\overline{\xi} \overline{\eta}} \right)^2 \right) = \overline{\xi}_a^{\overline{\eta}} + O\left( |\xi|^{-1} \left( \xi_t^{\overline{\eta} \overline{\sigma}} \right)^2 \right) \\
\overline{\xi}_a^{\overline{\xi}} &= \overline{\xi}_b^{\overline{\eta}} + O\left( |\xi|^{-1} \left( \xi_t^{\overline{\eta} \overline{\sigma}} \right)^2 \right) \\
\overline{\sigma}_a^{\overline{\sigma}} &= - \overline{\sigma}_a^{\overline{\eta}} + O\left( \left( \xi_t^{\overline{\eta} \overline{\sigma}} \right)^2 \right) \\
\overline{\sigma}_b^{\overline{\sigma}} &= - \overline{\sigma}_b^{\overline{\eta}} + O\left( \left( \xi_t^{\overline{\eta} \overline{\sigma}} \right)^2 \right) 
\end{align*}

Then, by Lemmas \ref{lemcalculssimplescoordonneesconiquevarphi} and \ref{lemcalculvarphiKdV1D}, we can compute
\begin{align*}
6 \sqrt{3} \frac{\eta_0}{|\eta_0|} \varphi &= \left( \overline{\xi}_a^{\overline{\xi}} \right)^3 + \left( \overline{\xi}_b^{\overline{\xi}} \right)^3 - \left( \overline{\eta}_a^{\overline{\eta}} \right)^3 - \left( \overline{\eta}_b^{\overline{\eta}} \right)^3 + \left( \overline{\sigma}_a^{\overline{\sigma}} \right)^3 + \left( \overline{\sigma}_b^{\overline{\sigma}} \right)^3 \\
&= \left( \overline{\xi}_a^{\overline{\eta}} \right)^3 - \left( \overline{\eta}_b^{\overline{\eta}} \right)^3 - \left( \overline{\sigma}_b^{\overline{\eta}} \right)^3
+ O\left( \overline{\eta}_a^{\overline{\eta}} - \overline{\sigma}_a^{\overline{\sigma}} \right) + O\left( \left( \xi_t^{\overline{\eta} \overline{\sigma}} \right)^2 \right) \\
&= 3 \overline{\xi}_a^{\overline{\eta}} \overline{\eta}_b^{\overline{\eta}} \overline{\sigma}_b^{\overline{\eta}} 
+ O\left( \overline{\eta}_a^{\overline{\eta}} - \overline{\sigma}_a^{\overline{\sigma}} \right) + O\left( \left( \xi_t^{\overline{\eta} \overline{\sigma}} \right)^2 \right)
\end{align*}

We now separate the integral \eqref{estapriorimauvaisbCBHHC-s} localising $\overline{\xi}$ and $\overline{\eta}, \overline{\sigma}$ depending on their respective sizes, as well as $\overline{\eta}_b^{\overline{\eta}}$ and $\overline{\sigma}_b^{\overline{\sigma}}$ and $s$: {\footnotesize 
\begin{align*}
\eqref{estapriorimauvaisbCBHHC-s} = \sum_{n \in \mathbb{N}} \sum_{\substack{j, j' \in \mathbb{Z}, \\|j-j'| \leq 2}} \sum_{\substack{l \in \mathbb{Z}, \\ l \leq -10}} \sum_{\substack{k, k' \in \mathbb{Z}, \\k, k' \leq -10}} \int_0^t \psi_n\left( \langle s \rangle \right) \int \psi_{j+l}(\overline{\xi}) \psi_{j, k}^{\widehat{\mathcal{C}}}(\overline{\eta}) \psi_{j', k'}^{\widehat{\mathcal{C}}}(\overline{\sigma}) e^{i s \varphi} i s m_{\widehat{\mathcal{C}}}(\overline{\xi}) \widehat{X}_b(\overline{\xi}) \cdot \nabla_{\overline{\xi}} \varphi ~ i \xi_0 \\
\widehat{f}(s, \overline{\eta}) \widehat{f}(s, \overline{\sigma}) ~ d\overline{\eta} ds
\end{align*} }

We will estimate for fixed $j, j', l, k, k'$ and then sum. To simplify the notations, we only treat the case $j = j'$, the others being straight-forward. We start with the estimates in $L^2$. 

We will denote by 
\begin{align*}
k_0 = k_0(n, j):= -\frac{n}{3}-\frac{29j}{12}+\frac{11j_{+}}{6}
\end{align*}
where $j_{+}$ stands for the positive part of $j$. 

\paragraph{1.} Let us start with the sum restricted to 
\begin{align*}
l - 100 \leq k_0
\end{align*}
In particular, here, there is no dependance in $k, k'$, so we may remove the summation in $k, k'$ and only estimate 
\begin{align}
\int_0^t \psi_n\left( \langle s \rangle \right) \int \psi_{l+j}(\overline{\xi}) \psi_j(\overline{\eta}) \psi_j(\overline{\sigma}) e^{i s \varphi} i s m_{\widehat{\mathcal{C}}}(\overline{\xi}) \widehat{X}_b(\overline{\xi}) \cdot \nabla_{\overline{\xi}} \varphi ~ i \xi_0 \widehat{f}(s, \overline{\eta}) \widehat{f}(s, \overline{\sigma}) ~ d\overline{\eta} ds \label{estapriorimauvaisbCBHHCs1} 
\end{align}
We apply no other transformation here and estimate directly 
\begin{align*}
\Vert \partial_t \eqref{estapriorimauvaisbCBHHCs1} \Vert_{L^2} &\lesssim \psi_n\left( \langle t \rangle \right) t 2^{l+4j-2j_{+}} \Vert |\nabla|^{-1} \langle \nabla \rangle^2 u(t) \Vert_{L^2} \Vert \psi_j(D) u(t) \Vert_{L^{\infty}} \\
&\lesssim t^{\frac{1}{6}} \langle t \rangle^{-\frac{1}{3}+100\delta} \psi_n\left( \langle t \rangle \right) 2^{l+\frac{5j}{2}-2j_{+}} \Vert u \Vert_X^2 
\end{align*}
Therefore, summing: 
\begin{align*}
&\sum_{n \in \mathbb{N}} \sum_{j \in \mathbb{Z}} \sum_{l \in \mathbb{Z}} 1_{l \leq -10} 1_{l-100 \leq k_0} t^{\frac{1}{6}} \langle t \rangle^{-\frac{1}{3}+100\delta} \psi_n\left( \langle t \rangle \right) 2^{l+\frac{5j}{2}-2j_{+}} \Vert u \Vert_X^2 \\
&\quad \lesssim \sum_{n \in \mathbb{N}} \sum_{j \in \mathbb{Z}} t^{\frac{1}{6}} \langle t \rangle^{-\frac{1}{3}+100\delta} \psi_n\left( \langle t \rangle \right) 2^{-\frac{n}{3}+\frac{j}{12}-\frac{j_{+}}{6}} \Vert u \Vert_X^2 \\
&\quad \lesssim t^{\frac{1}{6}} \langle t \rangle^{-\frac{2}{3}+100\delta} \Vert u \Vert_X^2 \lesssim \langle t \rangle^{-1-\delta} \langle t \rangle^{\frac{1}{2}+101\delta} \Vert u \Vert_X^2
\end{align*}
as wanted. 

Note then that, for $k_0 \leq -10$, this forces $2^j \gtrsim 2^{-\frac{n}{3}}$. Therefore, in all subsequent sums, we will have $|\overline{\eta}| \simeq 2^j \gtrsim \langle t \rangle^{-\frac{1}{3}}$. 

\paragraph{2.} Let us then consider the sum restricted to 
\begin{align*}
l -100 \geq k_0, \quad k + 100 \geq l, \quad k'+100 \geq l
\end{align*}
In this case, we can artificially add factors $\overline{\eta}_b^{\overline{\eta}}, \overline{\sigma}_b^{\overline{\sigma}}$: 
\begin{align*}
\xi_0 \widehat{X}_b(\overline{\xi}) \cdot \nabla_{\overline{\xi}} \varphi &= \mu_0 \xi_0 2^{2j} = 2^{-k-k'} \mu_0 \xi_0 \overline{\eta}_b^{\overline{\eta}} \overline{\sigma}_b^{\overline{\sigma}} \\
&= 2^{-k-k'} \mu_0 \varphi + 2^{3j-k-k'} \mu_0' \left( \xi_t^{\overline{\eta} \overline{\sigma}} \right)^2 + 2^{2j-k-k'} \mu_0'' \left( \overline{\eta}_a^{\overline{\eta}} - \overline{\sigma}_a^{\overline{\sigma}} \right) 
\end{align*}
for some symbols $\mu_0, \mu_0', \mu_0''$ of order $0$, that can change from line to line. The symbol $\mu_0$ is not necessarily of Coifman-Meyer type here, but can be decomposed as products of Coifman-Meyer symbols and symbols of the form $\frac{2^{k+j}}{\overline{\eta}_b^{\overline{\eta}}}$ (or $\overline{\sigma}_b^{\overline{\sigma}}$, or $\xi_0$) which, once multiplied by $\psi_{j, k}(\overline{\eta})$, correspond to some $\widetilde{\psi}_{j, k}^{\widehat{\mathcal{C}}}(\overline{\eta})$ for a similar function $\psi$. 

When we have the factor $\varphi$, we apply an integration by parts in $s$: {\footnotesize 
\begin{subequations}
\begin{align}
&\int_0^t \psi_n\left( \langle s \rangle \right) \int \psi_{l+j}(\overline{\xi}) \psi_{j, k}^{\widehat{\mathcal{C}}}(\overline{\eta}) \psi_{j, k'}^{\widehat{\mathcal{C}}}(\overline{\sigma}) e^{i s \varphi} i s 2^{-k-k'} \mu_0 \varphi m_{\widehat{\mathcal{C}}}(\overline{\xi}) ~ i \widehat{f}(s, \overline{\eta}) \widehat{f}(s, \overline{\sigma}) ~ d\overline{\eta} ds \notag \\
&\quad = - \int_0^t \psi_n\left( \langle s \rangle \right) \int \psi_{l+j}(\overline{\xi}) \psi_{j, k}^{\widehat{\mathcal{C}}}(\overline{\eta}) \psi_{j, k'}^{\widehat{\mathcal{C}}}(\overline{\sigma}) e^{i s \varphi} i s 2^{-k-k'} \mu_0 m_{\widehat{\mathcal{C}}}(\overline{\xi}) ~ \partial_s \left( \widehat{f}(s, \overline{\eta}) \widehat{f}(s, \overline{\sigma}) \right) ~ d\overline{\eta} ds \label{estapriorimauvaisbCBHHCs2varphi-1} \\
&\quad \quad - \int_0^t \partial_s \left( s \psi_n\left( \langle s \rangle \right) \right) \int \psi_{l+j}(\overline{\xi}) \psi_{j, k}^{\widehat{\mathcal{C}}}(\overline{\eta}) \psi_{j, k'}^{\widehat{\mathcal{C}}}(\overline{\sigma}) e^{i s \varphi} i 2^{-k-k'} \mu_0 m_{\widehat{\mathcal{C}}}(\overline{\xi}) ~ \widehat{f}(s, \overline{\eta}) \widehat{f}(s, \overline{\sigma}) ~ d\overline{\eta} ds \label{estapriorimauvaisbCBHHCs2varphi-2} \\
&\quad \quad + \psi_n\left( \langle t \rangle \right) \int \psi_{l+j}(\overline{\xi}) \psi_{j, k}^{\widehat{\mathcal{C}}}(\overline{\eta}) \psi_{j, k'}^{\widehat{\mathcal{C}}}(\overline{\sigma}) e^{i t \varphi} i t 2^{-k-k'} \mu_0 m_{\widehat{\mathcal{C}}}(\overline{\xi}) ~ \widehat{f}(t, \overline{\eta}) \widehat{f}(t, \overline{\sigma}) ~ d\overline{\eta} \label{estapriorimauvaisbCBHHCs2varphi-3}
\end{align}
\end{subequations} }
As before, only the boundary term \eqref{estapriorimauvaisbCBHHCs2varphi-3} in $g_{b*}$ and the rest in $h_{b*}$. 

We then estimate: {\footnotesize 
\begin{align*}
&\sum_{n \in \mathbb{N}} \sum_{\substack{j, l, k, k' \in \mathbb{Z}, \\ k, k', l \leq -10}} 1_{l -100 \geq k_0} 1_{\min(k, k') + 100 \geq l} \Vert \partial_t \eqref{estapriorimauvaisbCBHHCs2varphi-1} \Vert_{L^2} \\
&\quad \lesssim \sum_{n \in \mathbb{N}} \sum_{\substack{j, l, k, k' \in \mathbb{Z}, \\ k, k', l \leq -10}} 1_{l -100 \geq k_0} 1_{\min(k, k') + 100 \geq l} \psi_n\left( \langle t \rangle \right) t 2^{-k-k'+j-\frac{3j_{+}}{2}} \Vert |\nabla|^{-1} \langle \nabla \rangle^{\frac{3}{2}} \partial_t f(t) \Vert_{L^2} \Vert \psi_{j, k}^{\widehat{\mathcal{C}}}(D) u(t) \Vert_{L^{\infty}} \\
&\quad \lesssim \sum_{n \in \mathbb{N}} \sum_{\substack{j, l, k, k' \in \mathbb{Z}, \\ k, k', l \leq -10}} 1_{l -100 \geq k_0} 1_{\min(k, k') + 100 \geq l} \psi_n\left( \langle t \rangle \right) t^{-\frac{2}{3}} \langle t \rangle^{-\frac{7}{12}+200\delta} 2^{-k-k'-\frac{j}{2}-\frac{3j_{+}}{2}} \Vert u \Vert_X^2 \\
&\quad \lesssim \sum_{n \in \mathbb{N}} \sum_{j \in \mathbb{Z}} \psi_n\left( \langle t \rangle \right) t^{-\frac{2}{3}} \langle t \rangle^{\frac{1}{12}+200\delta} 2^{\frac{13j}{3}-\frac{29j_{+}}{6}} \Vert u \Vert_X^2 \lesssim t^{-\frac{2}{3}} \langle t \rangle^{\frac{1}{12}+200\delta+\frac{\delta}{3}} \Vert u \Vert_X^2 \\
&\quad \lesssim \left( t^{-\frac{2}{3}} \langle t \rangle^{-\frac{1}{3}-\delta} \right) \langle t \rangle^{\frac{5}{12}+202\delta} \Vert u \Vert_X^2 \\
&\sum_{n \in \mathbb{N}} \sum_{\substack{j, l, k, k' \in \mathbb{Z}, \\ k, k', l \leq -10}} 1_{l-100 \geq k_0} 1_{\min(k, k') + 100 \geq l} \Vert \partial_t \eqref{estapriorimauvaisbCBHHCs2varphi-2} \Vert_{L^2} \\
&\quad \lesssim \sum_{n \in \mathbb{N}} \sum_{\substack{j, l, k, k' \in \mathbb{Z}, \\ k, k', l \leq -10}} 1_{l-100 \geq k_0} 1_{\min(k, k') + 100 \geq l} \left| \partial_t \left( t \psi_n\left( \langle t \rangle \right) \right) \right| 2^{-k-k'+j-2 j_{+}} \Vert \langle \nabla \rangle^2 |\nabla|^{-1} f(t) \Vert_{L^2} \Vert \psi_{j, k}^{\widehat{\mathcal{C}}}(D) u(t) \Vert_{L^{\infty}} \\
&\quad \lesssim \sum_{n \in \mathbb{N}} \sum_{\substack{j, l \in \mathbb{Z}, \\ l \leq -10}} 1_{l-100 \geq k_0} \widetilde{\psi}_n\left( \langle t \rangle \right) t^{-\frac{5}{6}} \langle t \rangle^{-\frac{1}{3}+100\delta} 2^{-2l-\frac{j}{2}-2 j_{+}} \Vert u \Vert_X^2 \\
&\quad \lesssim \sum_{n \in \mathbb{N}} \sum_{\substack{j \in \mathbb{Z}}} \widetilde{\psi}_n\left( \langle t \rangle \right) t^{-\frac{5}{6}} \langle t \rangle^{\frac{1}{3}+100\delta} 2^{\frac{13j}{3}-\frac{17j_{+}}{3}} \Vert u \Vert_X^2 \lesssim \left( t^{-\frac{5}{6}} \langle t \rangle^{-\frac{1}{6}-\delta} \right) \langle t \rangle^{\frac{1}{2}+101\delta} \Vert u \Vert_X^2 \\
&\sum_{n \in \mathbb{N}} \sum_{\substack{j, l, k, k' \in \mathbb{Z}, \\ k, k', l \leq -10}} 1_{l-100 \geq k_0} 1_{\min(k, k') + 100 \geq l} \Vert \eqref{estapriorimauvaisbCBHHCs2varphi-3} \Vert_{L^2} \\
&\quad \lesssim \sum_{n \in \mathbb{N}} \sum_{\substack{j, l, k, k' \in \mathbb{Z}, \\ k, k', l \leq -10}} 1_{l-100 \geq k_0} 1_{\min(k, k') + 100 \geq l} \psi_n\left( \langle t \rangle \right) t 2^{-k-k'+j-2j_{+}} \Vert \langle \nabla \rangle^2 |\nabla|^{-1} f(t) \Vert_{L^2} \Vert \psi_{j, k}^{\widehat{\mathcal{C}}}(D) u(t) \Vert_{L^{\infty}} \\
&\quad \lesssim \sum_{n \in \mathbb{N}} \sum_{\substack{j, l \in \mathbb{Z}, \\ l \leq -10}} 1_{l-100 \geq k_0} \psi_n\left( \langle t \rangle \right) t^{\frac{1}{6}} \langle t \rangle^{-\frac{1}{3}+100\delta} 2^{-2l-\frac{j}{2}-2j_{+}} \Vert u \Vert_X^2 \\
&\quad \lesssim \sum_{n \in \mathbb{N}} \sum_{\substack{j \in \mathbb{Z}}} \psi_n\left( \langle t \rangle \right) t^{\frac{1}{6}} \langle t \rangle^{\frac{1}{3}+100\delta} 2^{\frac{13j}{3}-\frac{17j_{+}}{3}} \Vert u \Vert_X^2 \lesssim \langle t \rangle^{\frac{1}{2}+100\delta} \Vert u \Vert_X^2 \\
&\sum_{n \in \mathbb{N}} \sum_{\substack{j, l, k, k' \in \mathbb{Z}, \\ k, k', l \leq -10}} 1_{l-100 \geq k_0} 1_{\min(k, k') + 100 \geq l} \Vert e^{-it \omega(D)} \nabla \mathcal{F}^{-1} \eqref{estapriorimauvaisbCBHHCs2varphi-3} \Vert_{L^4} \\
&\quad \lesssim \sum_{n \in \mathbb{N}} \sum_{\substack{j, l, k, k' \in \mathbb{Z}, \\ k, k', l \leq -10}} 1_{l-100 \geq k_0} 1_{\min(k, k') + 100 \geq l} \psi_n\left( \langle t \rangle \right) t 2^{-k-k'+l+j} \Vert \psi_{j, k}^{\widehat{\mathcal{C}}}(D) u(t) \Vert_{L^4} \Vert \psi_{j, k}^{\widehat{\mathcal{C}}}(D) u(t) \Vert_{L^{\infty}} \\
&\quad \lesssim \sum_{n \in \mathbb{N}} \sum_{\substack{j, l, k, k' \in \mathbb{Z}, \\ k, k', l \leq -10}} 1_{l-100 \geq k_0} 1_{\min(k, k') + 100 \geq l} \psi_n\left( \langle t \rangle \right) t 2^{-k-k'+l+j} \Vert \psi_j(D) u(t) \Vert_{L^2}^{\frac{1}{2}} \Vert \psi_{j, k}^{\widehat{\mathcal{C}}}(D) u(t) \Vert_{L^{\infty}}^{\frac{3}{2}} \\
&\quad \lesssim \sum_{n \in \mathbb{N}} \sum_{\substack{j, l \in \mathbb{Z}, \\ l \leq -10}} 1_{l-100 \geq k_0} \psi_n\left( \langle t \rangle \right) t^{-\frac{1}{4}} \langle t \rangle^{-\frac{1}{2}+150\delta} 2^{-l-\frac{3j}{4}-j_{+}} \Vert u \Vert_X^2 \\
&\quad \lesssim \sum_{n \in \mathbb{N}} \sum_{\substack{j \in \mathbb{Z}}} \psi_n\left( \langle t \rangle \right) t^{-\frac{1}{4}} \langle t \rangle^{-\frac{1}{6}+150\delta} 2^{\frac{5j}{3}-\frac{17j_{+}}{6}} \Vert u \Vert_X^2 \lesssim t^{-\frac{1}{4}} \langle t \rangle^{-\frac{1}{6}+150\delta} \Vert u \Vert_X^2 
\end{align*} }
Above, $\widetilde{\psi}$ is a function with slightly larger support than $\psi$. 

For the factor $\left( \xi_t^{\overline{\eta} \overline{\sigma}} \right)^2$, recall that
\begin{align*}
\xi_t^{\overline{\eta} \overline{\sigma}} &= O\left( \widehat{X}_c(\overline{\eta}) \cdot \nabla_{\overline{\eta}} \varphi \right) 
\end{align*}
and so we can apply an integration by parts along $\widehat{X}_c(\overline{\eta})$: {\footnotesize 
\begin{subequations}
\begin{align}
&\int_0^t \psi_n\left( \langle s \rangle \right) \int \psi_{l+j}(\overline{\xi}) \psi_{j, k}^{\widehat{\mathcal{C}}}(\overline{\eta}) \psi_{j, k'}^{\widehat{\mathcal{C}}}(\overline{\sigma}) e^{i s \varphi} i s 2^{-k-k'+j} \mu_0 \xi_t^{\overline{\eta} \overline{\sigma}} \widehat{X}_c(\overline{\eta}) \cdot \nabla_{\overline{\eta}} \varphi m_{\widehat{\mathcal{C}}}(\overline{\xi}) ~ i \widehat{f}(s, \overline{\eta}) \widehat{f}(s, \overline{\sigma}) ~ d\overline{\eta} ds \notag \\
&\quad = - \int_0^t \psi_n\left( \langle s \rangle \right) \int \psi_{l+j}(\overline{\xi}) \psi_{j, k}^{\widehat{\mathcal{C}}}(\overline{\eta}) \psi_{j, k'}^{\widehat{\mathcal{C}}}(\overline{\sigma}) e^{i s \varphi} 2^{-k-k'+j} \mu_0 \xi_t^{\overline{\eta} \overline{\sigma}} m_{\widehat{\mathcal{C}}}(\overline{\xi}) ~ i \widehat{f}(s, \overline{\eta}) \widehat{X}_c(\overline{\eta}) \cdot \nabla_{\overline{\eta}} \widehat{f}(s, \overline{\sigma}) ~ d\overline{\eta} ds \label{estapriorimauvaisbCBHHCs2xit-1} \\
&\quad \quad - \int_0^t \psi_n\left( \langle s \rangle \right) \int \psi_{l+j}(\overline{\xi}) \psi_{j, k}^{\widehat{\mathcal{C}}}(\overline{\eta}) \psi_{j, k'}^{\widehat{\mathcal{C}}}(\overline{\sigma}) e^{i s \varphi} 2^{-k-k'+j} \mu_0 \xi_t^{\overline{\eta} \overline{\sigma}} m_{\widehat{\mathcal{C}}}(\overline{\xi}) ~ i \widehat{X}_c(\overline{\eta}) \cdot \nabla_{\overline{\eta}} \widehat{f}(s, \overline{\eta}) \widehat{f}(s, \overline{\sigma}) ~ d\overline{\eta} ds \label{estapriorimauvaisbCBHHCs2xit-2} \\
&\quad \quad - \int_0^t \psi_n\left( \langle s \rangle \right) \int e^{i s \varphi} \nabla_{\overline{\eta}} \cdot \left( \widehat{X}_c(\overline{\eta}) \psi_{l+j}(\overline{\xi}) \psi_{j, k}^{\widehat{\mathcal{C}}}(\overline{\eta}) \psi_{j, k'}^{\widehat{\mathcal{C}}}(\overline{\sigma}) 2^{-k-k'+j} \mu_0 \xi_t^{\overline{\eta} \overline{\sigma}} m_{\widehat{\mathcal{C}}}(\overline{\xi}) \right) ~ i \widehat{f}(s, \overline{\eta}) \widehat{f}(s, \overline{\sigma}) ~ d\overline{\eta} ds \label{estapriorimauvaisbCBHHCs2xit-3} 
\end{align}
\end{subequations} }
(for some symbol $\mu_0$ of order $0$.) 

In \eqref{estapriorimauvaisbCBHHCs2xit-1}, we project: 
\begin{align*}
\widehat{X}_c(\overline{\eta}) &= P_c^a(\overline{\eta}, \overline{\sigma}) \widehat{X}_a(\overline{\sigma}) + P_c^c(\overline{\eta}, \overline{\sigma}) \widehat{X}_c(\overline{\sigma}) + P_c^b(\overline{\eta}, \overline{\sigma}) \widehat{X}_b(\overline{\sigma}) 
\end{align*}
The contributions of $X_a, X_c$ of the first line are similar to \eqref{estapriorimauvaisbCBHHCs2xit-2}. For the contribution of $X_b$, we have that $P_c^b(\overline{\eta}, \overline{\sigma}) = O\left( \xi_t^{\overline{\eta} \overline{\sigma}} \right)$ and thus we get a term of the form {\footnotesize 
\begin{align}
\int_0^t \psi_n\left( \langle s \rangle \right) \int \psi_{l+j}(\overline{\xi}) \psi_{j, k}^{\widehat{\mathcal{C}}}(\overline{\eta}) \psi_{j, k'}^{\widehat{\mathcal{C}}}(\overline{\sigma}) e^{i s \varphi} 2^{-k-k'+j} \mu_0 \left( \xi_t^{\overline{\eta} \overline{\sigma}} \right)^2 m_{\widehat{\mathcal{C}}}(\overline{\xi}) ~ i \widehat{f}(s, \overline{\eta}) \widehat{X}_b(\overline{\sigma}) \cdot \nabla_{\overline{\eta}} \widehat{f}(s, \overline{\sigma}) ~ d\overline{\eta} ds \label{estapriorimauvaisbCBHHCs2xit-1bis} 
\end{align} }
But here, $\xi_t^{\overline{\eta} \overline{\sigma}} \lesssim \frac{\xi_0}{|\overline{\eta}|} \simeq 2^l$. We can therefore estimate, using that $k_0 \leq -10$ implies $2^j \gtrsim 2^{-\frac{n}{3}}$: {\footnotesize 
\begin{align*}
&\sum_{n \in \mathbb{N}} \sum_{\substack{j, l, k, k' \in \mathbb{Z}, \\ l, k, k' \leq -10}} 1_{l-100 \geq k_0} 1_{\min(k, k')+100 \geq l} \Vert \partial_t \eqref{estapriorimauvaisbCBHHCs2xit-1bis} \Vert_{L^2} \\
&\quad \lesssim \sum_{n \in \mathbb{N}} \sum_{\substack{j, l, k, k' \in \mathbb{Z}, \\ l, k, k' \leq -10}} 1_{l-100 \geq k_0} 1_{\min(k, k')+100 \geq l} \psi_n\left( \langle t \rangle \right) 2^{-k-k'+2l} \Vert \psi_{j, k}^{\widehat{\mathcal{C}}}(D) u(t) \Vert_{L^{\infty}} \Vert \nabla X_b f(t) \Vert_{L^2} \\
&\quad \lesssim \sum_{n \in \mathbb{N}} \sum_{\substack{j, l \in \mathbb{Z}, \\ l \leq -10}} 1_{l-100 \geq k_0} \psi_n\left( \langle t \rangle \right) t^{-\frac{5}{6}} \langle t \rangle^{-\frac{1}{12}+201\delta} 2^{-j-\frac{j_{+}}{2}} \Vert u \Vert_X^2 \\
&\quad \lesssim \sum_{n \in \mathbb{N}} \sum_{\substack{j \in \mathbb{Z}}} 1_{-10 \geq k_0} \psi_n\left( \langle t \rangle \right) t^{-\frac{5}{6}} \langle t \rangle^{-\frac{1}{12}+201\delta+\frac{\delta}{3}} 2^{\frac{29\delta j}{12}-\frac{11\delta j_{+}}{6}-j-\frac{j_{+}}{2}} \Vert u \Vert_X^2 \\
&\quad \lesssim t^{-\frac{5}{6}} \langle t \rangle^{\frac{1}{4}+202\delta} \Vert u \Vert_X^2 = \left( t^{-\frac{5}{6}} \langle t \rangle^{-\frac{1}{6}-\delta} \right) \langle t \rangle^{\frac{5}{12}+203\delta} \Vert u \Vert_X^2 \\
&\sum_{n \in \mathbb{N}} \sum_{\substack{j, l, k, k' \in \mathbb{Z}, \\ l, k, k' \leq -10}} 1_{l-100 \geq k_0} 1_{\min(k, k')+100 \geq l} \Vert \partial_t \eqref{estapriorimauvaisbCBHHCs2xit-2} \Vert_{L^2} \\
&\quad \lesssim \sum_{n \in \mathbb{N}} \sum_{\substack{j, l, k, k' \in \mathbb{Z}, \\ l, k, k' \leq -10}} 1_{l-100 \geq k_0} 1_{\min(k, k')+100 \geq l} \psi_n\left( \langle t \rangle \right) 2^{-k-k'+j-j_{+}} \Vert \langle \nabla \rangle X_c f(t) \Vert_{L^2} \Vert \psi_{j, k'}^{\widehat{\mathcal{C}}}(D) u(t) \Vert_{L^{\infty}} \\
&\quad \lesssim \sum_{n \in \mathbb{N}} \sum_{\substack{j, l \in \mathbb{Z}, \\ l \leq -10}} 1_{l-100 \geq k_0} \psi_n\left( \langle t \rangle \right) t^{-\frac{5}{6}} \langle t \rangle^{-\frac{1}{3}+100\delta} 2^{-2l-\frac{j}{2}-j_{+}} \Vert u \Vert_X^2 \\
&\quad \lesssim \sum_{n \in \mathbb{N}} \sum_{j \in \mathbb{Z}} \psi_n\left( \langle t \rangle \right) t^{-\frac{5}{6}} \langle t \rangle^{\frac{1}{3}+100\delta} 2^{\frac{13j}{3}-\frac{14j_{+}}{3}} \Vert u \Vert_X^2 \lesssim \left( t^{-\frac{5}{6}} \langle t \rangle^{-\frac{1}{6}-\delta} \right) \langle t \rangle^{\frac{1}{2}+101\delta} \Vert u \Vert_X^2 
\end{align*} }
Note that, in \eqref{estapriorimauvaisbCBHHCs2xit-2}, we did not use the additional factor $\xi_t^{\overline{\eta} \overline{\sigma}}$ (that was bounded by $1$). 

Finally, for \eqref{estapriorimauvaisbCBHHCs2xit-3}, we have here that
\begin{align*}
&\nabla_{\overline{\eta}} \cdot \left( \widehat{X}_c(\overline{\eta}) \psi_{l+j}(\overline{\xi}) \mu_0 m_{\widehat{\mathcal{C}}}(\overline{\xi}) \right) \lesssim 2^{-j} \\
&\widehat{X}_c(\overline{\eta}) \cdot \nabla_{\overline{\eta}} \psi_{j, k}^{\widehat{\mathcal{C}}}(\overline{\eta}) = 0 \\
&\widehat{X}_c(\overline{\eta}) \cdot \nabla_{\overline{\eta}} \psi_{j, k'}^{\widehat{\mathcal{C}}}(\overline{\sigma}) \lesssim 2^{-j} + 2^{-j-k'} \xi_t^{\overline{\eta} \overline{\sigma}} \lesssim 2^{-j} \\
&\widehat{X}_c(\overline{\eta}) \cdot \nabla_{\overline{\eta}} \xi_t^{\overline{\eta} \overline{\sigma}} \lesssim 2^{-j} 
\end{align*}
where we use that $\xi_t^{\overline{\eta} \overline{\sigma}} \lesssim 2^l$ and that $k' \geq l-100$ here. In particular, we can estimate this term like \eqref{estapriorimauvaisbCBHHCs2varphi-2}. 

It only remains the factor $O\left( \overline{\eta}_a^{\overline{\eta}} - \overline{\sigma}_a^{\overline{\sigma}} \right)$. We can then replace it by $\widehat{Y}(\overline{\eta}, \overline{\sigma}) \cdot \nabla_{\overline{\eta}} \varphi$ and apply an integration by parts. We then only get terms of the form \eqref{estapriorimauvaisbCBHHCs2xit-2} or \eqref{estapriorimauvaisbCBHHCs2xit-3} using the fundamental property of the vector field $\widehat{Y}(\overline{\eta}, \overline{\sigma})$ , recall Lemma \ref{lemchampmodifieYpropfond}, and we estimate the same way. 

\paragraph{3.} Let us now consider the sum restricted to
\begin{align*}
l - 100 \geq k_0, \quad k' + 100 \leq l
\end{align*}
(we can deal in a symmetric way with the case $k+100 \leq l$.) Then $k \geq k_0$ and $k \geq k'+10$. Indeed, 
\begin{align*}
2^l &\simeq \frac{|\overline{\xi}|}{|\overline{\eta}|} \simeq 2^{-j} \overline{\xi}_b^{\overline{\eta}} \simeq 2^{-j} \left( \overline{\eta}_b^{\overline{\eta}} + \overline{\sigma}_b^{\overline{\eta}} \right) 
\end{align*}
But $\overline{\eta}_b^{\overline{\eta}} \simeq 2^{j+k}$ while $\overline{\sigma}_b^{\overline{\eta}} \simeq 2^{j+k'} \ll 2^{j+l}$ here, so that $k \geq k_0$. More precisely, we can even restrict the sum to $|k -l |\leq 10$ (for instance), and to simplify the notations in the following we only consider the diagonal $k = l$ (the other cases being similar). 

Also, here, we group the summation in $k'$: 
\begin{align*}
\sum_{k' + 100 \leq l} \psi_{j, k'}(\overline{\sigma}) &= \psi_j(\overline{\sigma}) \chi\left( 2^{-j-l+100} \overline{\sigma}_b^{\overline{\sigma}} \right) = \chi_{j, l-100}^{\widehat{\mathcal{C}}}(\overline{\xi})
\end{align*}
where we reuse the notation \eqref{localisationCtaillejtotaleexactetaillekfinemaj}. To simplify the notations, we will drop the $-100$ in what follows. 

Let us consider
\begin{align*}
\widehat{X}_b'(\overline{\eta}) := \begin{pmatrix} \dfrac{\eta_0}{|\eta_0|} \\[3mm] - \sqrt{3} \dfrac{\eta}{|\eta|} \end{pmatrix}
\end{align*}
We can compute that 
\begin{align*}
\widehat{X}_b'(\overline{\eta}) \cdot \nabla_{\overline{\eta}} \varphi &= \frac{\eta_0}{|\eta_0|} \left( 3 \sigma_0^2 + |\sigma|^2 - 3 \eta_0^2 - |\eta|^2 \right) - \sqrt{3} \frac{\eta}{|\eta|} \cdot \left( 2 \sigma_0 \sigma - 2 \eta_0 \eta \right) \\
&= \frac{\eta_0}{|\eta_0|} \left( 3 \sigma_0^2 + |\sigma|^2 - 3 \eta_0^2 - |\eta|^2 \right) - \sqrt{3} \frac{\eta_0}{|\eta_0|} \left( 2 |\sigma_0| |\sigma| - 2 |\eta_0| |\eta| \right)
+ O\left( \left( \xi_t^{\overline{\eta} \overline{\sigma}} \right)^2 \right) \\
&= \frac{\eta_0}{|\eta_0|} \left( \left( \overline{\sigma}_b^{\overline{\sigma}} \right)^2 - \left( \overline{\eta}_b^{\overline{\eta}} \right)^2 \right) + O\left( \left( \xi_t^{\overline{\eta} \overline{\sigma}} \right)^2 \right) 
\end{align*}
Note that $\left| \xi_t^{\overline{\eta} \overline{\sigma}} \right| = \left| \xi_t^{\overline{\xi} \overline{\eta}} \right| |\xi| \ll |\xi| \simeq 2^{j+l}$, $\overline{\sigma}_b^{\overline{\sigma}} \ll 2^{j+l}$ but $\overline{\eta}_b^{\overline{\eta}} \simeq 2^{j+l}$. Therefore, 
\begin{align*}
\widehat{X}_b'(\overline{\eta}) \cdot \nabla_{\overline{\eta}} \varphi \simeq 2^{2j+2l} 
\end{align*}

We now create artificially this factor to apply an integration by parts: {\footnotesize 
\begin{subequations}
\begin{align}
&\int_0^t \psi_n\left( \langle s \rangle \right) \int \psi_{j+l}(\overline{\xi}) \psi_{j, l}^{\widehat{\mathcal{C}}}(\overline{\eta}) \chi_{j, l}^{\widehat{\mathcal{C}}}(\overline{\sigma}) e^{i s \varphi} i s \mu_0 \frac{2^{3j+l}}{\widehat{X}_b'(\overline{\eta}) \cdot \nabla_{\overline{\eta}} \varphi} \widehat{X}_b'(\overline{\eta}) \cdot \nabla_{\overline{\eta}} \varphi ~ m_{\widehat{\mathcal{C}}}(\overline{\xi}) i \widehat{f}(s, \overline{\eta}) \widehat{f}(s, \overline{\sigma}) ~ d\overline{\eta} ds \notag \\
&\quad = - \int_0^t \psi_n\left( \langle s \rangle \right) \int \psi_{j+l}(\overline{\xi}) \psi_{j, l}^{\widehat{\mathcal{C}}}(\overline{\eta}) \chi_{j, l}^{\widehat{\mathcal{C}}}(\overline{\sigma}) e^{i s \varphi} ~ \mu_0 \frac{2^{3j+l}}{\widehat{X}_b'(\overline{\eta}) \cdot \nabla_{\overline{\eta}} \varphi} m_{\widehat{\mathcal{C}}}(\overline{\xi}) i \widehat{X}_b'(\overline{\eta}) \cdot \nabla_{\overline{\eta}} \left( \widehat{f}(s, \overline{\eta}) \widehat{f}(s, \overline{\sigma}) \right) ~ d\overline{\eta} ds \label{estapriorimauvaisbCBHHCs3-1} \\
&\quad \quad - \int_0^t \psi_n\left( \langle s \rangle \right) \int e^{i s \varphi} \nabla_{\overline{\eta}} \cdot \left( \widehat{X}_b'(\overline{\eta}) \frac{2^{3j+l}}{\widehat{X}_b'(\overline{\eta}) \cdot \nabla_{\overline{\eta}} \varphi} \mu_0 \psi_{j+l}(\overline{\xi}) \psi_{j, l}^{\widehat{\mathcal{C}}}(\overline{\eta}) \chi_{j, l}^{\widehat{\mathcal{C}}}(\overline{\sigma}) \right) ~ m_{\widehat{\mathcal{C}}}(\overline{\xi}) i \widehat{f}(s, \overline{\eta}) \widehat{f}(s, \overline{\sigma}) ~ d\overline{\eta} ds \label{estapriorimauvaisbCBHHCs3-2} 
\end{align}
\end{subequations} }
for some symbol $\mu_0$ of order $0$, that can be decomposed as a product of Coifman-Meyer and Hörmander-Mikhlin symbols (in $\overline{\xi}$). 

On the one hand, we can estimate: {\footnotesize 
\begin{align*}
&\sum_{n \in \mathbb{N}} \sum_{\substack{j, l \in \mathbb{Z}, \\ l \leq -10}} 1_{l-100 \geq k_0} \Vert \partial_t \eqref{estapriorimauvaisbCBHHCs3-1} \Vert_{L^2} \\
&\quad \lesssim \sum_{n \in \mathbb{N}} \sum_{\substack{j, l \in \mathbb{Z}, \\ l \leq -10}} 1_{l-100 \geq k_0} \psi_n\left( \langle t \rangle \right) 2^{-l} \Vert m_{\widehat{\mathcal{C}}}(D) \nabla (x, y) f(t) \Vert_{L^2} \Vert m_{\widehat{\mathcal{C}}}(D) \psi_j(D) u(t) \Vert_{L^{\infty}} \\
&\quad \lesssim \sum_{n \in \mathbb{N}} \sum_{\substack{j, l \in \mathbb{Z}, \\ l \leq -10}} 1_{l-100 \geq k_0} \psi_n\left( \langle t \rangle \right) t^{-\frac{5}{6}} \langle t \rangle^{-\frac{1}{12}+201\delta} 2^{-l-j-\frac{j_{+}}{2}} \Vert u \Vert_X^2 \\
&\quad \lesssim \sum_{n \in \mathbb{N}} \sum_{\substack{j \in \mathbb{Z}}} \psi_n\left( \langle t \rangle \right) t^{-\frac{5}{6}} \langle t \rangle^{\frac{1}{4}+201\delta} 2^{\frac{17j}{12}-\frac{14j_{+}}{6}} \Vert u \Vert_X^2 \lesssim \left( t^{-\frac{5}{6}} \langle t \rangle^{-\frac{1}{6}-\delta} \right) \langle t \rangle^{\frac{5}{12}+202\delta} \Vert u \Vert_X^2 
\end{align*} }
which is enough. 

For \eqref{estapriorimauvaisbCBHHCs3-2}, we have that {\footnotesize 
\begin{align*}
\nabla_{\overline{\eta}} \cdot \left( \widehat{X}_b'(\overline{\eta}) \mu_0 \psi_{j+l}(\overline{\xi}) \right) &\lesssim 2^{-j} \\
\widehat{X}_b'(\overline{\eta}) \cdot \nabla_{\overline{\eta}} \left( \widehat{X}_b'(\overline{\eta}) \cdot \nabla_{\overline{\eta}} \varphi \right) &= \left( \frac{\eta_0}{|\eta_0|} \partial_{\eta_0} - \sqrt{3} \frac{\eta}{|\eta|} \cdot \nabla_{\eta} \right) \left[ \frac{\eta_0}{|\eta_0|} \left( 3 \sigma_0^2 + |\sigma|^2 - 3 \eta_0^2 - |\eta|^2 \right) - \sqrt{3} \frac{\eta}{|\eta|} \cdot \left( 2 \sigma_0 \sigma - 2 \eta_0 \eta \right) \right] \\
&=  -6 \xi_0 + 2 \sqrt{3} \frac{\eta_0 \eta}{|\eta_0| |\eta|} \cdot \xi - 2 \sqrt{3} \frac{\eta_0 \eta}{|\eta_0| |\eta|} \cdot \xi + 6 \xi_0 = 0 \\
\widehat{X}_b'(\overline{\eta}) \cdot \nabla_{\overline{\eta}} \left( \psi_{j, l}^{\widehat{\mathcal{C}}}(\overline{\eta}) \right) &\lesssim 2^{-l-j} \\
\widehat{X}_b'(\overline{\eta}) \cdot \nabla_{\overline{\eta}} \left( \chi_{j, l}^{\widehat{\mathcal{C}}}(\overline{\sigma}) \right) &\lesssim 2^{-l-j} 
\end{align*} }
Therefore, 
\begin{align*}
\nabla_{\overline{\eta}} \cdot \left( \widehat{X}_b'(\overline{\eta}) \frac{2^{3j+l}}{\widehat{X}_b'(\overline{\eta}) \cdot \nabla_{\overline{\eta}} \varphi} \mu_0 \psi_{j+l}(\overline{\xi}) \psi_{j, l}^{\widehat{\mathcal{C}}}(\overline{\eta}) \chi_{j, l}^{\widehat{\mathcal{C}}}(\overline{\sigma}) \right) &\lesssim 2^{-2l} 
\end{align*}
We can thus estimate \eqref{estapriorimauvaisbCBHHCs3-2} in a similar way as \eqref{estapriorimauvaisbCBHHCs2varphi-2}. 

This concludes the $L^2$ estimate. 

Let us now consider the $\dot{H}^1$ estimate (and the $W^{2, 4}$ estimate for $g$), which means we gain a factor $2^{l+j}$ for every term above. We can reuse essentially the same computations, changing $k_0$ by 
\begin{align*}
k_0 := -\frac{n}{3} - \frac{41j}{24} + \frac{11j_{+}}{12} 
\end{align*}
All the computations are similar and we skip the details. 

\subsection{Case of the line} 

We now localise $\overline{\xi}$ by $m_{\widehat{\mathcal{L}}}$. In this case, we will denote by $\beta \in \{ b, c \}$. 

Again, we will typically get better estimates than the desired one, except for the critial interaction which here comes from the $(BBB, \widehat{\mathcal{C}}\widehat{\mathcal{C}})$ case. 

\subsubsection{Generic interaction \texorpdfstring{$BBB$}{BBB}}

Let us first localise by $\mu_{BBB}(\overline{\xi}, \overline{\eta}) m_{A_1}(\overline{\eta}) m_{A_2}(\overline{\sigma})$, $A_1, A_2 \in \{ \widehat{\mathcal{L}}, \widehat{\mathcal{R}} \}$. 

We can then develop {\footnotesize 
\begin{subequations}
\begin{align}
&m_{\widehat{\mathcal{L}}}(\overline{\xi}) \widehat{X}_{\beta}(\overline{\xi}) \cdot \nabla_{\overline{\xi}} \widehat{I}_{\mu_{BBB} m_{A_1} m_{A_2}}[f, f](t, \overline{\xi}) \notag \\
&= \int_0^t \int e^{i s \varphi} i s \widehat{X}_{\beta}(\overline{\eta}) \cdot \nabla_{\overline{\xi}} \varphi ~ \mu_{BBB}(\overline{\xi}, \overline{\eta}) m_{\widehat{\mathcal{L}}}(\overline{\xi}) m_{A_1}(\overline{\eta}) m_{A_2}(\overline{\sigma}) i \xi_0 \widehat{f}(s, \overline{\eta}) \widehat{f}(s, \overline{\sigma}) ~ d\overline{\eta} ds \label{estapriorimauvaisbLgenBBB-1} \\
&\quad + \int_0^t \int e^{i s \varphi} \mu_{BBB}(\overline{\xi}, \overline{\eta}) m_{\widehat{\mathcal{L}}}(\overline{\xi}) m_{A_1}(\overline{\eta}) m_{A_2}(\overline{\sigma}) i \xi_0 \widehat{f}(s, \overline{\eta}) \widehat{X}_{\beta}(\overline{\eta}) \cdot \nabla_{\overline{\xi}} \widehat{f}(s, \overline{\sigma}) ~ d\overline{\eta} ds \label{estapriorimauvaisbLgenBBB-2} \\
&\quad + \int_0^t \int e^{i s \varphi} m_{\widehat{\mathcal{L}}}(\overline{\xi}) \widehat{X}_{\beta}(\overline{\xi}) \cdot \nabla_{\overline{\xi}} \left( \mu_{BBB}(\overline{\xi}, \overline{\eta}) m_{A_1}(\overline{\eta}) m_{A_2}(\overline{\sigma}) i \xi_0 \right) \widehat{f}(s, \overline{\eta}) \widehat{f}(s, \overline{\sigma}) ~ d\overline{\eta} ds \label{estapriorimauvaisbLgenBBB-3} 
\end{align}
\end{subequations} }

By Lemma \ref{lem-non-res-loin0Cone}, as in the generic $BBB$ case for the cone, we can always rewrite \eqref{estapriorimauvaisbLgenBBB-1} by integrations by parts to \eqref{estapriorimauvaisbLgenBBB-2}, \eqref{estapriorimauvaisbLgenBBB-3} or {\footnotesize 
\begin{subequations}
\begin{align}
&\int_0^t \int e^{i s \varphi} s \mu_0 ~ \mu_{BBB}(\overline{\xi}, \overline{\eta}) m_{\widehat{\mathcal{L}}}(\overline{\xi}) m_{A_1}(\overline{\eta}) m_{A_2}(\overline{\sigma}) i \partial_s \left( \widehat{f}(s, \overline{\eta}) \widehat{f}(s, \overline{\sigma}) \right) ~ d\overline{\eta} ds \label{estapriorimauvaisbLgenBBB-4} \\
&\quad + \int e^{i t \varphi} t \mu_0 ~ \mu_{BBB}(\overline{\xi}, \overline{\eta}) m_{\widehat{\mathcal{L}}}(\overline{\xi}) m_{A_1}(\overline{\eta}) m_{A_2}(\overline{\sigma}) i \widehat{f}(t, \overline{\eta}) \widehat{f}(t, \overline{\sigma}) ~ d\overline{\eta} \label{estapriorimauvaisbLgenBBB-5} 
\end{align}
\end{subequations} }
for some $\mu_0$ symbol of order $0$. As before, we only put the boundary term \eqref{estapriorimauvaisbLgenBBB-5} in $g_{\beta *}$ and the rest $h_{\beta *}$. 

We can estimate \eqref{estapriorimauvaisbLgenBBB-3}, \eqref{estapriorimauvaisbLgenBBB-4}, and \eqref{estapriorimauvaisbLgenBBB-5} like in the case of the cone, like \eqref{estapriorimauvaisbCgenBBB-dersymb}, \eqref{estapriorimauvaisbCgenBBB-s-3} and \eqref{estapriorimauvaisbCgenBBB-s-5} respectively (since they already had a better decay than wanted). 

For the last one, 
\begin{align*}
\Vert \partial_t \eqref{estapriorimauvaisbLgenBBB-2} \Vert_{L^2} &\lesssim \Vert m_{A_1}(D) \nabla u(t) \Vert_{L^{\infty}} \Vert m_{A_2}(D) (x, y) f(t) \Vert_{L^2} \\
&\lesssim t^{-\frac{5}{6}} \langle t \rangle^{-\frac{1}{4}+100 \delta} ~ \langle t \rangle^{\frac{5}{48}+202\delta} \Vert u \Vert_X^2 \\
\Vert |\overline{\xi}|^{\frac{1}{2}} \langle \overline{\xi} \rangle^{\frac{1}{2}} \partial_t \eqref{estapriorimauvaisbLgenBBB-2} \Vert_{L^2} &\lesssim \Vert m_{A_1}(D) \nabla u(t) \Vert_{L^{\infty}} \Vert m_{A_2}(D) |\nabla|^{\frac{1}{2}} \langle \nabla \rangle^{\frac{1}{2}} (x, y) f(t) \Vert_{L^2} \\
&\lesssim t^{-\frac{5}{6}} \langle t \rangle^{-\frac{1}{4}+100\delta} ~ \langle t \rangle^{202\delta} \Vert u \Vert_X^2 
\end{align*}

\subsubsection{Interaction \texorpdfstring{$(BBB, \widehat{\mathcal{P}}\widehat{\mathcal{R}})$}{(BBB, PR)}} 

Let us localise by $\mu_{BBB}(\overline{\xi}, \overline{\eta}) m_{\widehat{\mathcal{P}}}(\overline{\eta}) m_{\widehat{\mathcal{R}}}(\overline{\sigma})$. 

We can then develop and proceed just like the generic $BBB$ case, except for terms of the form:  
\begin{align}
\int_0^t \int e^{i s \varphi} \mu_0 \xi_0 \widehat{f}(s, \overline{\eta}) \nabla_{\overline{\xi}} \widehat{f}(s, \overline{\sigma}) ~ d\overline{\eta} ds \label{estapriorimauvaisbLBBBPR} 
\end{align} 
for $\mu_0$ a symbol of order $0$ containing all localisations. Then 
\begin{align*}
\Vert \langle \overline{\xi} \rangle \partial_t \eqref{estapriorimauvaisbLBBBPR} \Vert_{L^2} &\lesssim \Vert \nabla u(t) \Vert_{L^{\infty}} \Vert \langle \nabla \rangle m_{\widehat{\mathcal{R}}}(D) (x, y) f(t) \Vert_{L^2} \lesssim t^{-\frac{5}{6}} \langle t \rangle^{-\frac{1}{6}+100\delta} \Vert u \Vert_X^2 \\
&\lesssim \left( t^{-\frac{5}{6}} \langle t \rangle^{-\frac{1}{6}-\delta} \right) \langle t \rangle^{101\delta} \Vert u \Vert_X^2
\end{align*}
as wanted. 

\subsubsection{Interaction \texorpdfstring{$(BBB, \widehat{\mathcal{C}}\widehat{\mathcal{R}})$}{(BBB, CR)}}

Let us localise by $\mu_{BBB}(\overline{\xi}, \overline{\eta}) m_{\widehat{\mathcal{C}}}(\overline{\eta}) m_{\widehat{\mathcal{R}}}(\overline{\sigma})$. 

Developing as in the generic $BBB$ case, the only term we need to estimate differently is 
\begin{align}
\int_0^t \int e^{i s \varphi} i s \widehat{X}_b(\overline{\xi}) \cdot \nabla_{\overline{\xi}} \varphi ~ \mu_{BBB}(\overline{\xi}, \overline{\eta}) m_{\widehat{\mathcal{L}}}(\overline{\xi}) m_{\widehat{\mathcal{C}}}(\overline{\eta}) m_{\widehat{\mathcal{R}}}(\overline{\sigma}) i \xi_0 \widehat{f}(s, \overline{\eta}) \widehat{f}(s, \overline{\sigma}) ~ d\overline{\eta} ds \label{estapriorimauvaisbLBBBCR-s} 
\end{align}

By Lemma \ref{lem-non-res-loin0Cone}, we have 
\begin{align*}
1 = O\left( \nabla_{\overline{\eta}} \varphi \right) + O(\varphi) 
\end{align*}
But by Lemma \ref{lemcalculsconecoordonneesconiquesvarphi}, here, 
\begin{align*}
\xi_t^{\overline{\eta} \overline{\sigma}} &= O\left( \widehat{X}_c(\overline{\eta}) \cdot \nabla_{\overline{\eta}} \varphi \right) \\
\overline{\sigma}_a^{\overline{\eta}} - \overline{\eta}_a^{\overline{\eta}} &= O\left( \widehat{X}_a(\overline{\eta}) \cdot \nabla_{\overline{\eta}} \varphi \right) + O\left( \xi_t^{\overline{\eta} \overline{\sigma}} \right) + O\left( \overline{\eta}_b^{\overline{\eta}} \right) \\
6 \sqrt{3} \frac{\eta_0}{|\eta_0|} \varphi &= \left( \overline{\xi}_a^{\overline{\eta}} \right)^3 + \left( \overline{\xi}_b^{\overline{\eta}} \right)^3 - \left( \overline{\eta}_a^{\overline{\eta}} \right)^3 - \left( \overline{\eta}_b^{\overline{\eta}} \right)^3 - \left( \overline{\sigma}_a^{\overline{\eta}} \right)^3 - \left( \overline{\sigma}_b^{\overline{\eta}} \right)^3 + O\left( \xi_t^{\overline{\eta} \overline{\sigma}} \right) \\
&= 3 \overline{\xi}_a^{\overline{\eta}} \overline{\eta}_a^{\overline{\eta}} \overline{\sigma}_a^{\overline{\eta}} + O\left( \xi_t^{\overline{\eta} \overline{\sigma}} \right) + O\left( \overline{\eta}_b^{\overline{\eta}} \right)
\end{align*}
where we used $\xi_t^{\overline{\eta} \overline{\sigma}} \simeq \xi_t^{\overline{\xi} \overline{\eta}} \simeq \xi_t^{\overline{\xi} \overline{\sigma}}$ due to the localisation $\mu_{BBB}$. But here $\overline{\xi}_a^{\overline{\eta}} \simeq |\overline{\xi}| \simeq \overline{\eta}_a^{\overline{\eta}} \simeq \overline{\sigma}_a^{\overline{\eta}}$, so we deduce that 
\begin{align*}
1 = O(\varphi) + \sum_{\alpha = a, b, c} O\left( m_{\alpha}(\overline{\eta}) \widehat{X}_{\alpha}(\overline{\eta}) \cdot \nabla_{\overline{\eta}} \varphi \right) 
\end{align*}
In particular, we can decompose without singularities like in the decomposition lemma, and the estimates follow. 

\subsubsection{Interaction \texorpdfstring{$(BBB, \widehat{\mathcal{P}}\widehat{\mathcal{L}})$}{(BBB, PL)} or \texorpdfstring{$(BBB, \widehat{\mathcal{C}}\widehat{\mathcal{L}})$}{(BBB, CL)}} 

Let us localise by $\mu_{BBB}(\overline{\xi}, \overline{\eta}) m_{A_1}(\overline{\eta}) m_{\widehat{\mathcal{L}}}(\overline{\sigma})$ pour $A_1 \in \{ \widehat{\mathcal{P}}, \widehat{\mathcal{C}} \}$. 

Then 
\begin{align*}
|\overline{\xi}| + |\overline{\sigma}| \simeq |\overline{\eta}| \simeq |\eta| \lesssim |\xi| + |\sigma| 
\end{align*}
In particular, we deduce that, either $|\xi| \gtrsim |\overline{\xi}|$ and we can obtain for free $m_{\beta}(\overline{\xi})$ and apply the decomposition lemma, or $|\sigma| \gtrsim |\overline{\sigma}|$ and we can recover the case $(BBB, \widehat{\mathcal{P}}\widehat{\mathcal{R}})$ or $(BBB, \widehat{\mathcal{C}}\widehat{\mathcal{R}})$ (up to choosing localisation functions with larger support). 

\subsubsection{Interaction \texorpdfstring{$(BBB, \widehat{\mathcal{P}}\widehat{\mathcal{P}})$}{(BBB, PP)}} 

Let us localise by $\mu_{BBB}(\overline{\xi}, \overline{\eta}) m_{\widehat{\mathcal{P}}}(\overline{\eta}) m_{\widehat{\mathcal{P}}}(\overline{\sigma})$. 

Again, 
\begin{align*}
|\overline{\eta}| + |\overline{\sigma}| \simeq |\overline{\xi}| \simeq |\xi_0| \lesssim |\eta_0| + |\sigma_0| 
\end{align*}
We deduce that $|\eta_0| \gtrsim |\overline{\eta}|$ or $|\sigma_0| \gtrsim |\overline{\sigma}|$, which allows to recover the case $(BBB, \widehat{\mathcal{P}} \widehat{\mathcal{R}})$. 

\subsubsection{Interaction \texorpdfstring{$(BBB, \widehat{\mathcal{C}}\widehat{\mathcal{P}})$}{(BBB, CP)}}

Let us localise by $\mu_{BBB}(\overline{\xi}, \overline{\eta}) m_{\widehat{\mathcal{C}}}(\overline{\eta}) m_{\widehat{\mathcal{P}}}(\overline{\sigma})$. 

Developing as in the generic $BBB$ case, the only term that needs to be estimated differently is 
\begin{align*}
\int_0^t \int e^{i s \varphi} i s \widehat{X}_{\beta}(\overline{\eta}) \cdot \nabla_{\overline{\xi}} \varphi ~ \mu_{BBB}(\overline{\xi}, \overline{\eta}) m_{\widehat{\mathcal{L}}}(\overline{\xi}) m_{\widehat{\mathcal{C}}}(\overline{\eta}) m_{\widehat{\mathcal{P}}}(\overline{\sigma}) i \xi_0 \widehat{f}(s, \overline{\eta}) \widehat{f}(s, \overline{\sigma}) ~ d\overline{\eta} ds
\end{align*}
But here, 
\begin{align*}
\varphi &= \xi_0^3 + \xi_0 |\xi|^2 - \eta_0^3 - \eta_0 |\eta|^2 - \sigma_0^3 - \sigma_0 |\sigma|^2 = O(\xi) + O(\sigma_0) - \eta_0 |\eta|^2 
\end{align*}
and $\eta_0 |\eta|^2 \simeq |\overline{\eta}|^2 \gg |\xi| + |\sigma_0|$ by the localisation, so $1 = O(\varphi)$. We can thus apply only an integration by parts in time and estimate like in the generic $BBB$ case. 

\subsubsection{Interaction \texorpdfstring{$(BBB, \widehat{\mathcal{C}}\widehat{\mathcal{C}})$}{(BBB, CC)}}

Let us localise by $\mu_{BBB}(\overline{\xi}, \overline{\eta}) m_{\widehat{\mathcal{C}}}(\overline{\eta}) m_{\widehat{\mathcal{C}}}(\overline{\sigma})$. 

We now develop: {\footnotesize 
\begin{subequations}
\begin{align}
&m_{\widehat{\mathcal{L}}}(\overline{\xi}) \widehat{X}_{\beta}(\overline{\xi}) \cdot \nabla_{\overline{\xi}} \widehat{I}_{\mu_{BBB} m_{\widehat{\mathcal{C}}} m_{\widehat{\mathcal{C}}}}[f, f](t, \overline{\xi}) \notag \\
&= \int_0^t \int e^{i s \varphi} i s \widehat{X}_{\beta}(\overline{\eta}) \cdot \nabla_{\overline{\xi}} \varphi ~ \mu_{BBB}(\overline{\xi}, \overline{\eta}) m_{\widehat{\mathcal{L}}}(\overline{\xi}) m_{\widehat{\mathcal{C}}}(\overline{\eta}) m_{\widehat{\mathcal{C}}}(\overline{\sigma}) i \xi_0 \widehat{f}(s, \overline{\eta}) \widehat{f}(s, \overline{\sigma}) ~ d\overline{\eta} ds \label{estapriorimauvaisbLBBBCC-s} \\
&\quad + \int_0^t \int e^{i s \varphi} \mu_{BBB}(\overline{\xi}, \overline{\eta}) m_{\widehat{\mathcal{L}}}(\overline{\xi}) m_{\widehat{\mathcal{C}}}(\overline{\eta}) m_{\widehat{\mathcal{C}}}(\overline{\sigma}) i \xi_0 \widehat{f}(s, \overline{\eta}) \widehat{X}_{\beta}(\overline{\eta}) \cdot \nabla_{\overline{\xi}} \widehat{f}(s, \overline{\sigma}) ~ d\overline{\eta} ds \label{estapriorimauvaisbLBBBCC-eta} \\
&\quad + \int_0^t \int e^{i s \varphi} m_{\widehat{\mathcal{L}}}(\overline{\xi}) \widehat{X}_{\beta}(\overline{\xi}) \cdot \nabla_{\overline{\xi}} \left( \mu_{BBB}(\overline{\xi}, \overline{\eta}) m_{\widehat{\mathcal{C}}}(\overline{\eta}) m_{\widehat{\mathcal{C}}}(\overline{\sigma}) i \xi_0 \right) \widehat{f}(s, \overline{\eta}) \widehat{f}(s, \overline{\sigma}) ~ d\overline{\eta} ds \label{estapriorimauvaisbLBBBCC-dersymb} 
\end{align}
\end{subequations} }
\eqref{estapriorimauvaisbLBBBCC-dersymb} can be estimated like in the generic $BBB$ case. 

Note that, in this case, we have $\sigma$ close to $-\eta$ so $\theta^{\overline{\eta} \overline{\sigma}}$ close to $-1$. 

Assume first that $\epsilon^{\overline{\eta} \overline{\sigma}} = -1$. Then 
\begin{align*}
\left| \overline{\eta}_b^{\overline{\eta}} \right| + \left| \overline{\sigma}_b^{\overline{\sigma}} \right| + |\xi| &\geq \left| \sqrt{3} |\eta_0| - \sqrt{3} |\sigma_0| - |\eta| + |\sigma| \right| + |\xi| \\
&= \left| \sqrt{3} \frac{\eta_0}{|\eta_0|} \xi_0 - \frac{\eta}{|\eta|} \cdot \xi + |\sigma| \left( 1 + \theta^{\overline{\eta} \overline{\sigma}} \right) \right| + |\xi| \\
&\gtrsim |\xi_0| \simeq |\overline{\eta}| + |\overline{\sigma}| 
\end{align*}
so $\left| \overline{\eta}_b^{\overline{\eta}} \right| \gtrsim |\overline{\eta}|$ or $\left| \overline{\sigma}_b^{\overline{\sigma}} \right| \gtrsim |\overline{\sigma}|$, and, up to choosing localisations with larger support, we recover the case $(BBB, \widehat{\mathcal{C}} \widehat{\mathcal{R}})$. 

We now consider the case $\epsilon^{\overline{\eta} \overline{\sigma}} = 1$, as shown on Figure \ref{figurecasdegeneremaxLbBBBCC}. In particular, the vector field $\widehat{X}_{b-\widehat{\mathcal{C}}}(\overline{\eta}, \overline{\sigma})$ is not singular and, in a similar way as for the proof of the decomposition lemma, we can apply integration by parts to rewrite $\eqref{estapriorimauvaisbLBBBCC-s} + \eqref{estapriorimauvaisbLBBBCC-eta}$ as {\footnotesize 
\begin{align}
\int_0^t \int e^{i s \varphi} i s \left( \widehat{X}_{\beta}(\overline{\xi}) \cdot \nabla_{\overline{\xi}} + P_{\beta}^b(\overline{\xi}, \overline{\sigma}) \frac{\sigma_0}{|\sigma_0|} \widehat{X}_{b-\widehat{\mathcal{C}}}(\overline{\eta}, \overline{\sigma}) \cdot \nabla_{\overline{\eta}} \right) \varphi ~ \mu_{BBB}(\overline{\xi}, \overline{\eta}) m_{\widehat{\mathcal{L}}}(\overline{\xi}) m_{\widehat{\mathcal{C}}}(\overline{\eta}) m_{\widehat{\mathcal{C}}}(\overline{\sigma}) i \xi_0 \widehat{f}(s, \overline{\eta}) \widehat{f}(s, \overline{\sigma}) ~ d\overline{\eta} ds \label{estapriorimauvaisbLBBBCC-stot} 
\end{align} }
up to terms obtained in the decomposition lemma, and therefore simpler to estimate. 

Reusing the computations done in the proof of the decomposition lemma, we have 
\begin{align*}
\xi_t^{\overline{\eta} \overline{\sigma}} \left( \overline{\eta}_a^{\overline{\eta}} - \overline{\sigma}_a^{\overline{\sigma}} \right) &= O\left( \widehat{Y}(\overline{\eta}, \overline{\sigma}) \cdot \nabla_{\overline{\eta}} \varphi \right) \\
\overline{\eta}_b^{\overline{\eta}} &= O\left( m_b(\overline{\eta}) \widehat{X}_a(\overline{\sigma}) \cdot \nabla_{\overline{\eta}} \varphi \right) + O\left( m_b(\overline{\eta}) \widehat{X}_c(\overline{\sigma}) \cdot \nabla_{\overline{\eta}} \varphi \right) \\
\overline{\sigma}_b^{\overline{\sigma}} &= O\left( m_b(\overline{\sigma}) \widehat{X}_a(\overline{\eta}) \cdot \nabla_{\overline{\eta}} \varphi \right) + O\left( m_b(\overline{\sigma}) \widehat{X}_c(\overline{\eta}) \cdot \nabla_{\overline{\eta}} \varphi \right) \\
1 + \theta^{\overline{\eta} \overline{\sigma}} &= O\left( \varphi \right) + O\left( \overline{\eta}_b^{\overline{\eta}} \right) + O\left( \overline{\sigma}_b^{\overline{\sigma}} \right) 
\end{align*}
In particular, 
\begin{align*}
|\overline{\eta}| &= 2 |\eta_0| + O\left( \overline{\eta}_b^{\overline{\eta}} \right), \quad |\overline{\sigma}| = 2 |\sigma_0| + O\left( \overline{\sigma}_b^{\overline{\sigma}} \right) \\
|\xi|^2 &= |\eta|^2 + |\sigma|^2 + 2 |\eta| |\sigma| \theta^{\overline{\eta} \overline{\sigma}} \\
&= 3 (\eta_0 - \sigma_0)^2 + O\left( \overline{\eta}_b^{\overline{\eta}} \right) + O\left( \overline{\sigma}_b^{\overline{\sigma}} \right) + O\left( 1 + \theta^{\overline{\eta} \overline{\sigma}} \right) 
\end{align*}
We can then compute that {\footnotesize 
\begin{align*}
&\left( \widehat{X}_b(\overline{\xi}) \cdot \nabla_{\overline{\xi}} + P_b^b(\overline{\xi}, \overline{\sigma}) \frac{\sigma_0}{|\sigma_0|} \widehat{X}_{b-\widehat{\mathcal{C}}}(\overline{\eta}, \overline{\sigma}) \cdot \nabla_{\overline{\eta}} \right) \varphi \\
&\quad = \frac{|\xi|}{|\overline{\xi}|} \left( 3 \xi_0^2 + |\xi|^2 - 3 \sigma_0^2 - |\sigma|^2 \right) - \frac{\xi_0 \xi}{|\overline{\xi}| |\xi|} \cdot \left( 2 \xi_0 \xi - 2 \sigma_0 \sigma \right) \\
&\quad \quad + \frac{|\xi|^2 |\sigma|^2 + \xi_0 \sigma_0 \xi \cdot \sigma}{|\overline{\xi}| |\overline{\sigma}| |\xi| |\sigma|} \frac{2 \sqrt{3}}{7 + \theta^{\overline{\eta} \overline{\sigma}}} \left( \left( \overline{\sigma}_b^{\overline{\sigma}} \right)^2 - \left( \overline{\eta}_a^{\overline{\eta}} \right)^2 - 2 \sqrt{3} \frac{1 + \theta^{\overline{\eta} \overline{\sigma}}}{3} |\sigma_0| |\sigma| \right) \\
&\quad = \frac{|\xi|}{|\overline{\xi}|} \left( \xi_0^2 + |\xi|^2 - 6 \sigma_0^2 \right) + 2 \frac{\xi_0 \sigma_0 \xi \cdot \sigma}{|\overline{\xi}| |\xi|} - 4 \sqrt{3} \frac{|\xi|^2 |\sigma|^2 + \xi_0 \sigma_0 \xi \cdot \sigma}{|\overline{\xi}| |\overline{\sigma}| |\xi| |\sigma|} \eta_0^2 
+ O\left( \overline{\eta}_b^{\overline{\eta}} \right) + O\left( \overline{\sigma}_b^{\overline{\sigma}} \right) + O\left( 1 + \theta^{\overline{\eta} \overline{\sigma}} \right) \\
&\quad = \frac{|\xi|}{|\overline{\xi}|} \left( (\eta_0 + \sigma_0)^2 + 3 (\eta_0 - \sigma_0)^2 - 6 \sigma_0^2 - 6 \eta_0^2 \right) 
+ 4 \sqrt{3} \frac{\xi_0 \sigma_0 \xi \cdot \sigma}{|\overline{\xi}| |\overline{\sigma}| |\xi| |\sigma|} \left( \sigma_0^2 - \eta_0^2 \right) 
+ O\left( \overline{\eta}_b^{\overline{\eta}} \right) + O\left( \overline{\sigma}_b^{\overline{\sigma}} \right) + O\left( 1 + \theta^{\overline{\eta} \overline{\sigma}} \right) \\
&\quad = - 2 \frac{|\xi|}{|\overline{\xi}|} \left( \eta_0^2 + \sigma_0^2 + \eta_0 \sigma_0 \right) 
+ 2 \frac{\xi_0^2 \xi \cdot \sigma}{|\overline{\xi}| |\xi| |\sigma|} (|\eta| - |\sigma|)
+ O\left( \overline{\eta}_b^{\overline{\eta}} \right) + O\left( \overline{\sigma}_b^{\overline{\sigma}} \right) + O\left( 1 + \theta^{\overline{\eta} \overline{\sigma}} \right) \\
&\quad = 2 \frac{\xi_0^2 (|\eta| - |\sigma|)^2 - \xi_0^2 |\xi|^2}{|\overline{\xi}| |\xi|}
+ O\left( \overline{\eta}_b^{\overline{\eta}} \right) + O\left( \overline{\sigma}_b^{\overline{\sigma}} \right) + O\left( 1 + \theta^{\overline{\eta} \overline{\sigma}} \right) \\
&\quad = - 4 \frac{\xi_0^2 |\eta| |\sigma| ( 1 + \theta^{\overline{\eta} \overline{\sigma}})}{|\overline{\xi}| |\xi|}
+ O\left( \overline{\eta}_b^{\overline{\eta}} \right) + O\left( \overline{\sigma}_b^{\overline{\sigma}} \right) + O\left( 1 + \theta^{\overline{\eta} \overline{\sigma}} \right) \\
&\left( \widehat{X}_c(\overline{\xi}) \cdot \nabla_{\overline{\xi}} + P_c^b(\overline{\xi}, \overline{\sigma}) \frac{\sigma_0}{|\sigma_0|} \widehat{X}_{b-\widehat{\mathcal{C}}}(\overline{\eta}, \overline{\sigma}) \cdot \nabla_{\overline{\eta}} \right) \varphi \\
&\quad = \frac{J \xi}{|\xi|} \cdot \left( 2 \xi_0 \xi - 2 \sigma_0 \sigma \right) - \frac{\sigma_0 J \xi \cdot \sigma}{|\overline{\sigma}| |\xi| |\sigma|} \frac{2 \sqrt{3}}{7 + \theta^{\overline{\eta} \overline{\sigma}}} \left( \left( \overline{\sigma}_b^{\overline{\sigma}} \right)^2 - \left( \overline{\eta}_a^{\overline{\eta}} \right)^2 - 2 \sqrt{3} \frac{1 + \theta^{\overline{\eta} \overline{\sigma}}}{3} |\sigma_0| |\sigma| \right) \\
&\quad = - 2 \sqrt{3} \sigma_0^2 \xi_t^{\overline{\xi} \overline{\sigma}} + 2 \sqrt{3} \frac{\sigma_0}{|\sigma_0|} \eta_0^2 \xi_t^{\overline{\xi} \overline{\sigma}} 
+ O\left( \overline{\eta}_b^{\overline{\eta}} \right) + O\left( \overline{\sigma}_b^{\overline{\sigma}} \right) + O\left( 1 + \theta^{\overline{\eta} \overline{\sigma}} \right) \\
&\quad = \frac{1}{2 \sqrt{3}} \left( \left( \overline{\eta}_a^{\overline{\eta}} \right)^2 - \left( \overline{\sigma}_a^{\overline{\sigma}} \right)^2 \right) \xi_t^{\overline{\xi} \overline{\sigma}} 
+ O\left( \overline{\eta}_b^{\overline{\eta}} \right) + O\left( \overline{\sigma}_b^{\overline{\sigma}} \right) + O\left( 1 + \theta^{\overline{\eta} \overline{\sigma}} \right)
\end{align*} }
In both cases, it only remains a term where we miss a factor $|\xi|$ in order to apply a non-singular decomposition; furthermore, we have 
\begin{align*}
\frac{1 + \theta^{\overline{\eta} \overline{\sigma}}}{|\xi|} \lesssim |\xi|, \quad \quad \left( \left( \overline{\eta}_a^{\overline{\eta}} \right)^2 - \left( \overline{\sigma}_a^{\overline{\sigma}} \right)^2 \right) \lesssim |\xi| 
\end{align*}
The other terms can be decomposed as in the proof of the decomposition lemma and then easily estimated. We therefore only need to consider {\footnotesize 
\begin{subequations}
\begin{align}
&\int_0^t \int e^{i s \varphi} i s \mu_0 \frac{1 + \theta^{\overline{\eta} \overline{\sigma}}}{|\xi|} |\overline{\eta}|^4 ~ \mu_{BBB}(\overline{\xi}, \overline{\eta}) m_{\widehat{\mathcal{L}}}(\overline{\xi}) m_{\widehat{\mathcal{C}}}(\overline{\eta}) m_{\widehat{\mathcal{C}}}(\overline{\sigma}) \widehat{f}(s, \overline{\eta}) \widehat{f}(s, \overline{\sigma}) ~ d\overline{\eta} ds \label{estapriorimauvaisbLBBBCC-svarphi} \\
&\quad + \int_0^t \int e^{i s \varphi} i s \mu_0' \left( \left( \overline{\eta}_a^{\overline{\eta}} \right)^2 - \left( \overline{\sigma}_a^{\overline{\sigma}} \right)^2 \right) \xi_t^{\overline{\xi} \overline{\sigma}} ~ \mu_{BBB}(\overline{\xi}, \overline{\eta}) m_{\widehat{\mathcal{L}}}(\overline{\xi}) m_{\widehat{\mathcal{C}}}(\overline{\eta}) m_{\widehat{\mathcal{C}}}(\overline{\sigma}) \widehat{f}(s, \overline{\eta}) \widehat{f}(s, \overline{\sigma}) ~ d\overline{\eta} ds \label{estapriorimauvaisbLBBBCC-sY} 
\end{align}
\end{subequations} }
for some symbols $\mu_0, \mu_0'$ of order $0$. 

We now decompose 
\begin{align*}
1 &= \sum_{n \in \mathbb{N}} \sum_{\substack{j, k \in \mathbb{Z}, \\ k \leq -10}} \psi_n\left( \langle s \rangle \right) \psi_{j, k}^{\widehat{\mathcal{L}}}(\overline{\xi}) 
\end{align*}
and get {\footnotesize 
\begin{subequations}
\begin{align}
&\int_0^t \psi_n\left( \langle s \rangle \right) \psi_{j, k}^{\widehat{\mathcal{L}}}(\overline{\xi}) \int e^{i s \varphi} i s \mu_0 \frac{1 + \theta^{\overline{\eta} \overline{\sigma}}}{|\xi|} |\overline{\eta}|^4 ~ \mu_{BBB}(\overline{\xi}, \overline{\eta}) m_{\widehat{\mathcal{L}}}(\overline{\xi}) m_{\widehat{\mathcal{C}}}(\overline{\eta}) m_{\widehat{\mathcal{C}}}(\overline{\sigma}) \widehat{f}(s, \overline{\eta}) \widehat{f}(s, \overline{\sigma}) ~ d\overline{\eta} ds \label{estapriorimauvaisbLBBBCC-svarphidec} \\
&\quad + \int_0^t \psi_n\left( \langle s \rangle \right) \psi_{j, k}^{\widehat{\mathcal{L}}}(\overline{\xi}) \int e^{i s \varphi} i s \mu_0'\left( \left( \overline{\eta}_a^{\overline{\eta}} \right)^2 - \left( \overline{\sigma}_a^{\overline{\sigma}} \right)^2 \right) \xi_t^{\overline{\xi} \overline{\sigma}}  |\overline{\eta}| ~ \mu_{BBB}(\overline{\xi}, \overline{\eta}) m_{\widehat{\mathcal{L}}}(\overline{\xi}) m_{\widehat{\mathcal{C}}}(\overline{\eta}) m_{\widehat{\mathcal{C}}}(\overline{\sigma}) \widehat{f}(s, \overline{\eta}) \widehat{f}(s, \overline{\sigma}) ~ d\overline{\eta} ds \label{estapriorimauvaisbLBBBCC-sYdec} 
\end{align}
\end{subequations} }

\paragraph{1.} Let us consider the sum restricted to 
\begin{align*}
j \leq - \frac{5n}{24}, \quad k \leq -\frac{3j}{2} - \frac{n}{2} 
\end{align*}
In this case, we can estimate directly: 
\begin{align*}
&\sum_{n \in \mathbb{N}} \sum_{\substack{j, k \in \mathbb{Z}, \\ k \leq -10}} 1_{j \leq -\frac{5n}{24}} 1_{k \leq -\frac{3j}{2}-\frac{n}{2}} \Vert \partial_t (\eqref{estapriorimauvaisbLBBBCC-svarphidec} + \eqref{estapriorimauvaisbLBBBCC-sYdec}) \Vert_{L^2} \\
&\quad \lesssim \sum_{n \in \mathbb{N}} \sum_{\substack{j, k \in \mathbb{Z}, \\ k \leq -10}} 1_{j \leq -\frac{5n}{24}} 1_{k \leq -\frac{3j}{2}-\frac{n}{2}} \psi_n\left( \langle t \rangle \right) t 2^{3j+k} \Vert u(t) \Vert_{L^4}^2 \\
&\quad \lesssim \sum_{n \in \mathbb{N}} \sum_{\substack{j \in \mathbb{Z}}} 1_{j \leq -\frac{5n}{24}} \psi_n\left( \langle t \rangle \right) 2^{\frac{3j}{2}} t^{\frac{1}{6}} \langle t \rangle^{-\frac{3}{4}+100\delta} \Vert u \Vert_X^2 \\
&\quad \lesssim t^{\frac{1}{6}} \langle t \rangle^{-\frac{17}{16}+100\delta} \Vert u \Vert_X^2 \lesssim \left( t^{\frac{1}{6}} \langle t \rangle^{-\frac{7}{6}-\delta} \right) \langle t \rangle^{\frac{5}{48}+101\delta} \Vert u \Vert_X^2 \\
&\sum_{n \in \mathbb{N}} \sum_{\substack{j, k \in \mathbb{Z}, \\ k \leq -10}} 1_{j \leq -\frac{5n}{24}} 1_{k \leq -\frac{3j}{2}-\frac{n}{2}} \Vert |\overline{\xi}|^{\frac{1}{2}} \langle \overline{\xi} \rangle^{\frac{1}{2}} \partial_t (\eqref{estapriorimauvaisbLBBBCC-svarphidec} + \eqref{estapriorimauvaisbLBBBCC-sYdec}) \Vert_{L^2} \\
&\quad \lesssim \sum_{n \in \mathbb{N}} \sum_{\substack{j, k \in \mathbb{Z}, \\ k \leq -10}} 1_{j \leq -\frac{5n}{24}} 1_{k \leq -\frac{3j}{2}-\frac{n}{2}} \psi_n\left( \langle t \rangle \right) t 2^{\frac{7j}{2}+k} \Vert u(t) \Vert_{L^4}^2 \\
&\quad \lesssim t^{\frac{1}{6}} \langle t \rangle^{-\frac{7}{6}+100\delta} \Vert u \Vert_X^2 \lesssim \left( t^{\frac{1}{6}} \langle t \rangle^{-\frac{7}{6}-\delta} \right) \langle t \rangle^{101\delta} \Vert u \Vert_X^2 
\end{align*}
which is enough. 

\paragraph{2.} Let us consider the sum restricted to 
\begin{align*}
&j \leq -\frac{5n}{24}, \quad k \geq -\frac{3j}{2} - \frac{n}{2} \\
&\mbox{or } j \geq -\frac{5n}{24}, \quad k \geq \frac{j}{2} - \frac{n}{12} 
\end{align*}
In this case, we apply an integration by parts on every term, either in time or frequency: {\footnotesize 
\begin{subequations} 
\begin{align}
&\eqref{estapriorimauvaisbLBBBCC-svarphidec} = \int_0^t \psi_n\left( \langle s \rangle \right) \psi_{j, k}^{\widehat{\mathcal{L}}}(\overline{\xi}) |\xi|^{-1} \int e^{i s \varphi} i s \mu_0 (1 + \theta^{\overline{\eta} \overline{\sigma}}) |\overline{\eta}|^4 ~ \mu_{BBB}(\overline{\xi}, \overline{\eta}) m_{\widehat{\mathcal{L}}}(\overline{\xi}) m_{\widehat{\mathcal{C}}}(\overline{\eta}) m_{\widehat{\mathcal{C}}}(\overline{\sigma}) \widehat{f}(s, \overline{\eta}) \widehat{f}(s, \overline{\sigma}) ~ d\overline{\eta} ds \notag \\
&= \int_0^t \psi_n\left( \langle s \rangle \right) \psi_{j, k}^{\widehat{\mathcal{L}}}(\overline{\xi}) |\xi|^{-1} \int e^{i s \varphi} i s \mu_0 |\overline{\eta}| ~ \mu_{BBB}(\overline{\xi}, \overline{\eta}) m_{\widehat{\mathcal{L}}}(\overline{\xi}) m_{\widehat{\mathcal{C}}}(\overline{\eta}) m_{\widehat{\mathcal{C}}}(\overline{\sigma}) \partial_s \left( \widehat{f}(s, \overline{\eta}) \widehat{f}(s, \overline{\sigma}) \right) ~ d\overline{\eta} ds \label{estapriorimauvaisbLBBBCC-svarphidec-1} \\
&\quad + \int_0^t \partial_s \left( s \psi_n\left( \langle s \rangle \right) \right) \psi_{j, k}^{\widehat{\mathcal{L}}}(\overline{\xi}) |\xi|^{-1} \int e^{i s \varphi} i \mu_0 |\overline{\eta}| ~ \mu_{BBB}(\overline{\xi}, \overline{\eta}) m_{\widehat{\mathcal{L}}}(\overline{\xi}) m_{\widehat{\mathcal{C}}}(\overline{\eta}) m_{\widehat{\mathcal{C}}}(\overline{\sigma}) \widehat{f}(s, \overline{\eta}) \widehat{f}(s, \overline{\sigma}) ~ d\overline{\eta} ds \label{estapriorimauvaisbLBBBCC-svarphidec-2} \\
&\quad + \psi_n\left( \langle t \rangle \right) \psi_{j, k}^{\widehat{\mathcal{L}}}(\overline{\xi}) |\xi|^{-1} \int e^{i t \varphi} i t \mu_0 |\overline{\eta}| ~ \mu_{BBB}(\overline{\xi}, \overline{\eta}) m_{\widehat{\mathcal{L}}}(\overline{\xi}) m_{\widehat{\mathcal{C}}}(\overline{\eta}) m_{\widehat{\mathcal{C}}}(\overline{\sigma}) \widehat{f}(t, \overline{\eta}) \widehat{f}(t, \overline{\sigma}) ~ d\overline{\eta} ds \label{estapriorimauvaisbLBBBCC-svarphidec-3} \\
&\eqref{estapriorimauvaisbLBBBCC-sYdec} = \int_0^t \psi_n\left( \langle s \rangle \right) \psi_{j, k}^{\widehat{\mathcal{L}}}(\overline{\xi}) |\xi|^{-1} \int e^{i s \varphi} \mu_0  |\overline{\eta}| ~ \mu_{BBB}(\overline{\xi}, \overline{\eta}) m_{\widehat{\mathcal{L}}}(\overline{\xi}) m_{\widehat{\mathcal{C}}}(\overline{\eta}) m_{\widehat{\mathcal{C}}}(\overline{\sigma}) \widehat{Y}(\overline{\eta}, \overline{\sigma}) \cdot \nabla_{\overline{\eta}} \left( \widehat{f}(s, \overline{\eta}) \widehat{f}(s, \overline{\sigma}) ~ d\overline{\eta} \right) ds \label{estapriorimauvaisbLBBBCC-sYdec-1} \\
&\quad + \int_0^t \psi_n\left( \langle s \rangle \right) \psi_{j, k}^{\widehat{\mathcal{L}}}(\overline{\xi}) |\xi|^{-1} \int e^{i s \varphi} \nabla_{\overline{\eta}} \cdot \left( \widehat{Y}(\overline{\eta}, \overline{\sigma}) \mu_0  |\overline{\eta}| ~ \mu_{BBB}(\overline{\xi}, \overline{\eta}) m_{\widehat{\mathcal{L}}}(\overline{\xi}) m_{\widehat{\mathcal{C}}}(\overline{\eta}) m_{\widehat{\mathcal{C}}}(\overline{\sigma}) \right) \widehat{f}(s, \overline{\eta}) \widehat{f}(s, \overline{\sigma}) ~ d\overline{\eta} ds \label{estapriorimauvaisbLBBBCC-sYdec-2} 
\end{align}
\end{subequations} }
where $\mu_0$ is a symbol of order $0$ that can change from line to line. We only put the boundary term \eqref{estapriorimauvaisbLBBBCC-svarphidec-3} in $g_{b*}$ and the rest in $h_{b*}$ or $h_{c*}$. We then estimate: {\footnotesize 
\begin{align*}
&\sum_{n \in \mathbb{N}} \sum_{\substack{j, k \in \mathbb{Z}, \\ k \leq -10}} \left( 1_{j \leq -\frac{5n}{24}} 1_{k \geq -\frac{3j}{2}-\frac{n}{2}} + 1_{j \geq -\frac{5n}{24}} 1_{k \geq \frac{j}{2}-\frac{n}{12}} \right) \Vert \partial_t \eqref{estapriorimauvaisbLBBBCC-svarphidec-1} \Vert_{L^2} \\
&\quad \lesssim \sum_{n \in \mathbb{N}} \sum_{\substack{j, k \in \mathbb{Z}, \\ k \leq -10}} \left( 1_{j \leq -\frac{5n}{24}} 1_{k \geq -\frac{3j}{2}-\frac{n}{2}} + 1_{j \geq -\frac{5n}{24}} 1_{k \geq \frac{j}{2}-\frac{n}{12}} \right) \psi_n\left( \langle t \rangle \right) t 2^{-k} \Vert |\nabla|^{-1} \partial_t f(t) \Vert_{L^2} \Vert \psi_j(D) m_{\widehat{\mathcal{C}}}(D) \nabla u(t) \Vert_{L^{\infty}} \\
&\quad \lesssim \sum_{n \in \mathbb{N}} \sum_{\substack{j \in \mathbb{Z}}} \left( 1_{j \leq -\frac{5n}{24}} 2^{\frac{3j}{2}+\frac{n}{2}} + 1_{j \geq -\frac{5n}{24}} 2^{-\frac{j}{2}+\frac{n}{12}} \right) \psi_n\left( \langle t \rangle \right) t^{-\frac{2}{3}} \langle t \rangle^{-\frac{1}{2}+200\delta} \Vert u \Vert_X^3 \\
&\quad \lesssim \langle t \rangle^{\frac{3}{16}} t^{-\frac{2}{3}} \langle t \rangle^{-\frac{1}{2}+200\delta} \Vert u \Vert_X^3 \\
&\quad \lesssim \left( t^{-\frac{2}{3}} \langle t \rangle^{-\frac{1}{3}-\delta} \right) \langle t \rangle^{\frac{1}{48}+201\delta} \Vert u \Vert_X^3 \\
&\sum_{n \in \mathbb{N}} \sum_{\substack{j, k \in \mathbb{Z}, \\ k \leq -10}} \left( 1_{j \leq -\frac{5n}{24}} 1_{k \geq -\frac{3j}{2}-\frac{n}{2}} + 1_{j \geq -\frac{5n}{24}} 1_{k \geq \frac{j}{2}-\frac{n}{12}} \right) \Vert |\overline{\xi}|^{\frac{1}{2}} \langle \overline{\xi} \rangle^{\frac{1}{2}} \partial_t \eqref{estapriorimauvaisbLBBBCC-svarphidec-1} \Vert_{L^2} \\
&\quad \lesssim \sum_{n \in \mathbb{N}} \sum_{\substack{j, k \in \mathbb{Z}, \\ k \leq -10}} \left( 1_{j \leq -\frac{5n}{24}} 1_{k \geq -\frac{3j}{2}-\frac{n}{2}} + 1_{j \geq -\frac{5n}{24}} 1_{k \geq \frac{j}{2}-\frac{n}{12}} \right) \psi_n\left( \langle t \rangle \right) t 2^{-k+\frac{j}{2}-\frac{j_{+}}{2}} \Vert |\nabla|^{-1} \langle \nabla \rangle \partial_t f(t) \Vert_{L^2} \Vert \psi_j(D) m_{\widehat{\mathcal{C}}}(D) \nabla u(t) \Vert_{L^{\infty}} \\
&\quad \lesssim \sum_{n \in \mathbb{N}} \left( 2^{\frac{n}{12}} + |n| 2^{\frac{n}{12}} \right) \psi_n\left( \langle t \rangle \right) t^{-\frac{2}{3}} \langle t \rangle^{-\frac{1}{2}+200\delta} \Vert u \Vert_X^3 \\
&\quad \lesssim \left( t^{-\frac{2}{3}} \langle t \rangle^{-\frac{1}{3}-\delta} \right) \langle t \rangle^{-\frac{1}{12}+202\delta} \Vert u \Vert_X^3 \\
&\sum_{n \in \mathbb{N}} \sum_{\substack{j, k \in \mathbb{Z}, \\ k \leq -10}} \left( 1_{j \leq -\frac{5n}{24}} 1_{k \geq -\frac{3j}{2}-\frac{n}{2}} + 1_{j \geq -\frac{5n}{24}} 1_{k \geq \frac{j}{2}-\frac{n}{12}} \right) \Vert \partial_t \eqref{estapriorimauvaisbLBBBCC-svarphidec-2} \Vert_{L^2} \\
&\quad \lesssim \sum_{n \in \mathbb{N}} \sum_{\substack{j, k \in \mathbb{Z}, \\ k \leq -10}} \left( 1_{j \leq -\frac{5n}{24}} 1_{k \geq -\frac{3j}{2}-\frac{n}{2}} + 1_{j \geq -\frac{5n}{24}} 1_{k \geq \frac{j}{2}-\frac{n}{12}} \right) \widetilde{\psi}_n\left( \langle t \rangle \right) 2^{-k} \Vert u(t) \Vert_{L^4}^2 \\
&\quad \lesssim \sum_{n \in \mathbb{N}} \sum_{\substack{j \in \mathbb{Z}}} \left( 1_{j \leq -\frac{5n}{24}} 2^{\frac{3j}{2}+\frac{n}{2}} + 1_{j \geq -\frac{5n}{24}} 2^{-\frac{j}{2}+\frac{n}{12}} \right) \widetilde{\psi}_n\left( \langle t \rangle \right) t^{-\frac{5}{6}} \langle t \rangle^{-\frac{1}{4}+100\delta} \Vert u \Vert_X^2 \\
&\quad \lesssim \langle t \rangle^{\frac{3}{16}} t^{-\frac{5}{6}} \langle t \rangle^{-\frac{1}{4}+100\delta} \Vert u \Vert_X^2 \\
&\quad \lesssim \left( t^{-\frac{5}{6}} \langle t \rangle^{-\frac{1}{6}-\delta} \right) \langle t \rangle^{\frac{5}{48}+101\delta} \Vert u \Vert_X^2 \\
&\sum_{n \in \mathbb{N}} \sum_{\substack{j, k \in \mathbb{Z}, \\ k \leq -10}} \left( 1_{j \leq -\frac{5n}{24}} 1_{k \geq -\frac{3j}{2}-\frac{n}{2}} + 1_{j \geq -\frac{5n}{24}} 1_{k \geq \frac{j}{2}-\frac{n}{12}} \right) \Vert |\overline{\xi}|^{\frac{1}{2}} \langle \overline{\xi} \rangle^{\frac{1}{2}} \partial_t \eqref{estapriorimauvaisbLBBBCC-svarphidec-2} \Vert_{L^2} \\
&\quad \lesssim \sum_{n \in \mathbb{N}} \sum_{\substack{j, k \in \mathbb{Z}, \\ k \leq -10}} \left( 1_{j \leq -\frac{5n}{24}} 1_{k \geq -\frac{3j}{2}-\frac{n}{2}} + 1_{j \geq -\frac{5n}{24}} 1_{k \geq \frac{j}{2}-\frac{n}{12}} \right) \widetilde{\psi}_n\left( \langle t \rangle \right) 2^{-k+\frac{j}{2}-\frac{3j_{+}}{2}} \Vert \langle \nabla \rangle m_{\widehat{\mathcal{C}}}(D) u(t) \Vert_{L^4}^2 \\
&\quad \lesssim \sum_{n \in \mathbb{N}} \sum_{\substack{j \in \mathbb{Z}}} \left( 1_{j \leq -\frac{5n}{24}} 2^{2j+\frac{n}{2}} + 1_{j \geq -\frac{5n}{24}} 2^{-\frac{3j_{+}}{2}+\frac{n}{12}} \right) \widetilde{\psi}_n\left( \langle t \rangle \right) t^{-\frac{5}{6}} \langle t \rangle^{-\frac{1}{4}+100\delta} \Vert u \Vert_X^2 \\
&\quad \lesssim (1 + \log(\langle t \rangle)) \langle t \rangle^{\frac{1}{12}} t^{-\frac{5}{6}} \langle t \rangle^{-\frac{1}{4}+100\delta} \Vert u \Vert_X^2 \\
&\quad \lesssim \left( (1 + \log(\langle t \rangle)) t^{-\frac{5}{6}} \langle t \rangle^{-\frac{1}{6}-\delta} \right) \langle t \rangle^{101\delta} \Vert u \Vert_X^2 \\
&\sum_{n \in \mathbb{N}} \sum_{\substack{j, k \in \mathbb{Z}, \\ k \leq -10}} \left( 1_{j \leq -\frac{5n}{24}} 1_{k \geq -\frac{3j}{2}-\frac{n}{2}} + 1_{j \geq -\frac{5n}{24}} 1_{k \geq \frac{j}{2}-\frac{n}{12}} \right) \Vert \eqref{estapriorimauvaisbLBBBCC-svarphidec-3} \Vert_{L^2} \\
&\quad \lesssim \sum_{n \in \mathbb{N}} \sum_{\substack{j, k \in \mathbb{Z}, \\ k \leq -10}} \left( 1_{j \leq -\frac{5n}{24}} 1_{k \geq -\frac{3j}{2}-\frac{n}{2}} + 1_{j \geq -\frac{5n}{24}} 1_{k \geq \frac{j}{2}-\frac{n}{12}} \right) \psi_n\left( \langle t \rangle \right) t 2^{-k} \Vert u(t) \Vert_{L^4}^2 \\
&\quad \lesssim \left( t^{\frac{1}{6}} \langle t \rangle^{-\frac{1}{6}-\delta} \right) \langle t \rangle^{\frac{5}{48}+101\delta} \Vert u \Vert_X^2 \\
&\sum_{n \in \mathbb{N}} \sum_{\substack{j, k \in \mathbb{Z}, \\ k \leq -10}} \left( 1_{j \leq -\frac{5n}{24}} 1_{k \geq -\frac{3j}{2}-\frac{n}{2}} + 1_{j \geq -\frac{5n}{24}} 1_{k \geq \frac{j}{2}-\frac{n}{12}} \right) \Vert |\overline{\xi}|^{\frac{1}{2}} \langle \overline{\xi} \rangle^{\frac{1}{2}} \eqref{estapriorimauvaisbLBBBCC-svarphidec-3} \Vert_{L^2} \\
&\quad \lesssim \sum_{n \in \mathbb{N}} \sum_{\substack{j, k \in \mathbb{Z}, \\ k \leq -10}} \left( 1_{j \leq -\frac{5n}{24}} 1_{k \geq -\frac{3j}{2}-\frac{n}{2}} + 1_{j \geq -\frac{5n}{24}} 1_{k \geq \frac{j}{2}-\frac{n}{12}} \right) \psi_n\left( \langle t \rangle \right) t 2^{-k+\frac{j}{2}-\frac{3j_{+}}{2}} \Vert \langle \nabla \rangle m_{\widehat{\mathcal{C}}}(D) u(t) \Vert_{L^4}^2 \\
&\quad \lesssim \left( t^{\frac{1}{6}} \langle t \rangle^{-\frac{1}{6}-\delta} \right) \langle t \rangle^{101\delta} \Vert u \Vert_X^2 \\
&\sum_{n \in \mathbb{N}} \sum_{\substack{j, k \in \mathbb{Z}, \\ k \leq -10}} \left( 1_{j \leq -\frac{5n}{24}} 1_{k \geq -\frac{3j}{2}-\frac{n}{2}} + 1_{j \geq -\frac{5n}{24}} 1_{k \geq \frac{j}{2}-\frac{n}{12}} \right) \Vert \nabla \mathcal{F}^{-1} \eqref{estapriorimauvaisbLBBBCC-svarphidec-3} \Vert_{L^4} \\ 
&\quad \lesssim \sum_{n \in \mathbb{N}} \sum_{\substack{j, k \in \mathbb{Z}, \\ k \leq -10}} \left( 1_{j \leq -\frac{5n}{24}} 1_{k \geq -\frac{3j}{2}-\frac{n}{2}} + 1_{j \geq -\frac{5n}{24}} 1_{k \geq \frac{j}{2}-\frac{n}{12}} \right) 2^{-k} \Vert u(t) \Vert_{L^4} \Vert \nabla m_{\widehat{\mathcal{C}}}(D) u(t) \Vert_{L^{\infty}} \\
&\quad \lesssim \left( t^{\frac{1}{6}} \langle t \rangle^{-\frac{1}{6}-\delta} \right) \left( t^{-\frac{5}{12}} \langle t \rangle^{-\frac{1}{8}+50\delta} \right) \langle t \rangle^{\frac{5}{48}+101\delta} \Vert u \Vert_X^2 \\
&\sum_{n \in \mathbb{N}} \sum_{\substack{j, k \in \mathbb{Z}, \\ k \leq -10}} \left( 1_{j \leq -\frac{5n}{24}} 1_{k \geq -\frac{3j}{2}-\frac{n}{2}} + 1_{j \geq -\frac{5n}{24}} 1_{k \geq \frac{j}{2}-\frac{n}{12}} \right) \Vert |\nabla|^{\frac{3}{2}} \langle \nabla \rangle^{\frac{1}{2}} \mathcal{F}^{-1} \eqref{estapriorimauvaisbLBBBCC-svarphidec-3} \Vert_{L^4} \\ 
&\quad \lesssim \sum_{n \in \mathbb{N}} \sum_{\substack{j, k \in \mathbb{Z}, \\ k \leq -10}} \left( 1_{j \leq -\frac{5n}{24}} 1_{k \geq -\frac{3j}{2}-\frac{n}{2}} + 1_{j \geq -\frac{5n}{24}} 1_{k \geq \frac{j}{2}-\frac{n}{12}} \right) 2^{-k+\frac{j}{2}-\frac{j_{+}}{2}} \Vert \langle \nabla \rangle m_{\widehat{\mathcal{C}}}(D) u(t) \Vert_{L^4} \Vert \nabla m_{\widehat{\mathcal{C}}}(D) u(t) \Vert_{L^{\infty}} \\
&\quad \lesssim \left( t^{\frac{1}{6}} \langle t \rangle^{-\frac{1}{6}-\delta} \right) \left( t^{-\frac{5}{12}} \langle t \rangle^{-\frac{1}{8}+50\delta} \right) \langle t \rangle^{101\delta} \Vert u \Vert_X^2 \\
&\sum_{n \in \mathbb{N}} \sum_{\substack{j, k \in \mathbb{Z}, \\ k \leq -10}} \left( 1_{j \leq -\frac{5n}{24}} 1_{k \geq -\frac{3j}{2}-\frac{n}{2}} + 1_{j \geq -\frac{5n}{24}} 1_{k \geq \frac{j}{2}-\frac{n}{12}} \right) \Vert \partial_t \eqref{estapriorimauvaisbLBBBCC-sYdec-1} \Vert_{L^2} \\
&\quad \lesssim \sum_{n \in \mathbb{N}} \sum_{\substack{j, k \in \mathbb{Z}, \\ k \leq -10}} \left( 1_{j \leq -\frac{5n}{24}} 1_{k \geq -\frac{3j}{2}-\frac{n}{2}} + 1_{j \geq -\frac{5n}{24}} 1_{k \geq \frac{j}{2}-\frac{n}{12}} \right) \psi_n\left( \langle t \rangle \right) 2^{-k} \sum_{\alpha = a, b, c} \Vert m_{\alpha}(D) X_{\alpha} f(t) \Vert_{L^2} \Vert m_{\widehat{\mathcal{C}}}(D) \nabla u(t) \Vert_{L^{\infty}} \\
&\quad \lesssim \left( t^{-\frac{5}{6}} \langle t \rangle^{-\frac{1}{6}-\delta} \right) \langle t \rangle^{\frac{5}{48}+101\delta} \Vert u \Vert_X^2 \\
&\sum_{n \in \mathbb{N}} \sum_{\substack{j, k \in \mathbb{Z}, \\ k \leq -10}} \left( 1_{j \leq -\frac{5n}{24}} 1_{k \geq -\frac{3j}{2}-\frac{n}{2}} + 1_{j \geq -\frac{5n}{24}} 1_{k \geq \frac{j}{2}-\frac{n}{12}} \right) \Vert |\overline{\xi}|^{\frac{1}{2}} \langle \overline{\xi} \rangle^{\frac{1}{2}} \partial_t \eqref{estapriorimauvaisbLBBBCC-sYdec-1} \Vert_{L^2} \\
&\quad \lesssim \sum_{n \in \mathbb{N}} \sum_{\substack{j, k \in \mathbb{Z}, \\ k \leq -10}} \left( 1_{j \leq -\frac{5n}{24}} 1_{k \geq -\frac{3j}{2}-\frac{n}{2}} + 1_{j \geq -\frac{5n}{24}} 1_{k \geq \frac{j}{2}-\frac{n}{12}} \right) \psi_n\left( \langle t \rangle \right) 2^{-k+\frac{j}{2}-\frac{j_{+}}{2}} \\
&\quad \quad \quad \sum_{\alpha = a, b, c} \Vert \langle \nabla \rangle m_{\alpha}(D) X_{\alpha} f(t) \Vert_{L^2} \Vert m_{\widehat{\mathcal{C}}}(D) \nabla u(t) \Vert_{L^{\infty}} \\
&\quad \lesssim \left( t^{-\frac{5}{6}} \langle t \rangle^{-\frac{1}{6}-\delta} \right) \langle t \rangle^{101\delta} \Vert u \Vert_X^2 
\end{align*} }
\eqref{estapriorimauvaisbLBBBCC-sYdec-2} is similar to \eqref{estapriorimauvaisbLBBBCC-svarphidec-2}. 

\paragraph{3.} Let us finally consider the sum restricted to 
\begin{align*}
j \geq -\frac{5n}{24}, \quad k \leq \frac{j}{2} - \frac{n}{12}
\end{align*}
In this case, since $1 = O\left( \nabla_{\overline{\eta}} \varphi \right)$ here, we can estimate \eqref{estapriorimauvaisbLBBBCC-svarphidec} and \eqref{estapriorimauvaisbLBBBCC-sYdec} together by applying an integration by parts in $\overline{\eta}$, up to having singular $b$-weights: {\footnotesize 
\begin{subequations}
\begin{align}
\eqref{estapriorimauvaisbLBBBCC-svarphidec} + \eqref{estapriorimauvaisbLBBBCC-sYdec}
&= \int_0^t \psi_n\left( \langle s \rangle \right) \psi_{j, k}^{\widehat{\mathcal{L}}}(\overline{\xi}) |\xi| \int e^{i s \varphi} i s \mu_0 \nabla_{\overline{\eta}} \varphi ~ \mu_{BBB}(\overline{\xi}, \overline{\eta}) m_{\widehat{\mathcal{L}}}(\overline{\xi}) m_{\widehat{\mathcal{C}}}(\overline{\eta}) m_{\widehat{\mathcal{C}}}(\overline{\sigma}) \widehat{f}(s, \overline{\eta}) \widehat{f}(s, \overline{\sigma}) ~ d\overline{\eta} ds \notag \\
&= \int_0^t \psi_n\left( \langle s \rangle \right) \psi_{j, k}^{\widehat{\mathcal{L}}}(\overline{\xi}) |\xi| \int e^{i s \varphi} \mu_0 ~ \mu_{BBB}(\overline{\xi}, \overline{\eta}) m_{\widehat{\mathcal{L}}}(\overline{\xi}) m_{\widehat{\mathcal{C}}}(\overline{\eta}) m_{\widehat{\mathcal{C}}}(\overline{\sigma}) \nabla_{\overline{\eta}} \left( \widehat{f}(s, \overline{\eta}) \widehat{f}(s, \overline{\sigma}) \right) ~ d\overline{\eta} ds \label{estapriorimauvaisbLBBBCC-ssing1} \\
&\quad + \int_0^t \psi_n\left( \langle s \rangle \right) \psi_{j, k}^{\widehat{\mathcal{L}}}(\overline{\xi}) |\xi| \int e^{i s \varphi} \nabla_{\overline{\eta}} \left( \mu_0 ~ \mu_{BBB}(\overline{\xi}, \overline{\eta}) m_{\widehat{\mathcal{L}}}(\overline{\xi}) m_{\widehat{\mathcal{C}}}(\overline{\eta}) m_{\widehat{\mathcal{C}}}(\overline{\sigma}) \right) \widehat{f}(s, \overline{\eta}) \widehat{f}(s, \overline{\sigma}) ~ d\overline{\eta} ds \label{estapriorimauvaisbLBBBCC-ssing2}
\end{align}
\end{subequations} }
for some symbol $\mu_0$ of order $0$ that can change from line to line. We can then estimate: {\footnotesize 
\begin{align*}
&\sum_{n \in \mathbb{N}} \sum_{\substack{j, k \in \mathbb{Z}, \\ k \leq -10}} 1_{j \geq -\frac{5n}{24}} 1_{k \leq \frac{j}{2}-\frac{n}{12}} \Vert \partial_t \eqref{estapriorimauvaisbLBBBCC-ssing1} \Vert_{L^2} \\
&\quad \lesssim \sum_{n \in \mathbb{N}} \sum_{\substack{j, k \in \mathbb{Z}, \\ k \leq -10}} 1_{j \geq -\frac{5n}{24}} 1_{k \leq \frac{j}{2}-\frac{n}{12}} \psi_n\left( \langle t \rangle \right) 2^{k-j} \Vert m_{\widehat{\mathcal{C}}}(D) \nabla u(t) \Vert_{L^{\infty}} \Vert m_{\widehat{\mathcal{C}}}(D) \nabla (x, y) f(t) \Vert_{L^2} \\
&\quad \lesssim \sum_{n \in \mathbb{N}} \sum_{\substack{j \in \mathbb{Z}}} 1_{j \geq -\frac{5n}{24}} \psi_n\left( \langle t \rangle \right) t^{-\frac{5}{6}} \langle t \rangle^{-\frac{1}{6}+100\delta+101\delta} 2^{-\frac{j}{2}} \Vert u \Vert_X^2 \\
&\quad \lesssim \left( t^{-\frac{5}{6}} \langle t \rangle^{-\frac{1}{6}-\delta} \right) \langle t \rangle^{\frac{5}{48}+202\delta} \Vert u \Vert_X^2 \\
&\sum_{n \in \mathbb{N}} \sum_{\substack{j, k \in \mathbb{Z}, \\ k \leq -10}} 1_{j \geq -\frac{5n}{24}} 1_{k \leq \frac{j}{2}-\frac{n}{12}} \Vert |\overline{\xi}|^{\frac{1}{2}} \langle \overline{\xi} \rangle^{\frac{1}{2}} \partial_t \eqref{estapriorimauvaisbLBBBCC-ssing1} \Vert_{L^2} \\
&\quad \lesssim \sum_{n \in \mathbb{N}} \sum_{\substack{j, k \in \mathbb{Z}, \\ k \leq -10}} 1_{j \geq -\frac{5n}{24}} 1_{k \leq \frac{j}{2}-\frac{n}{12}} \psi_n\left( \langle t \rangle \right) 2^{k-\frac{j}{2}-\frac{j_{+}}{8}} \Vert m_{\widehat{\mathcal{C}}}(D) |\nabla| \langle \nabla \rangle^{\frac{1}{2}+\frac{1}{8}} u(t) \Vert_{L^{\infty}} \Vert m_{\widehat{\mathcal{C}}}(D) \nabla (x, y) f(t) \Vert_{L^2} \\
&\quad \lesssim \sum_{n \in \mathbb{N}} \sum_{\substack{j \in \mathbb{Z}}} 1_{j \geq -\frac{5n}{24}} \psi_n\left( \langle t \rangle \right) t^{-\frac{5}{6}} \langle t \rangle^{-\frac{1}{6}+100\delta+101\delta} 2^{-\frac{j_{+}}{8}} \Vert u \Vert_X^2 \\
&\quad \lesssim \left( (1 + \log(\langle t \rangle)) t^{-\frac{5}{6}} \langle t \rangle^{-\frac{1}{6}-\delta} \right) \langle t \rangle^{\frac{5}{48}+202\delta} \Vert u \Vert_X^2 
\end{align*} }
\eqref{estapriorimauvaisbLBBBCC-ssing2} is simpler. 

\subsubsection{Generic \texorpdfstring{$HBH$}{HBH} interaction}

Let us localise by $\mu_{HBH}(\overline{\xi}, \overline{\eta}) m_{A}(\overline{\sigma})$, $A \in \{ \widehat{\mathcal{R}}, \widehat{\mathcal{C}}, \widehat{\mathcal{P}} \}$. 

In order for the integral not to be identically zero, we need $A = \widehat{\mathcal{R}}$. But then 
\begin{align*}
|\xi| + |\eta| \gtrsim |\sigma| \simeq |\overline{\sigma}| \simeq |\overline{\xi}| 
\end{align*}
In particular, either $|\overline{\eta}| \gtrsim |\overline{\xi}|$ and we recover a $BBB$ case (up to changing the symbol), or $|\xi| \gtrsim |\overline{\xi}|$ and we get for free $m_{\beta}(\overline{\xi})$ allowing to apply the decomposition lemma. 

\subsubsection{Interaction \texorpdfstring{$(HBH, \widehat{\mathcal{L}})$}{(HBH, L)}}

Let us localise by $\mu_{BBB}(\overline{\xi}, \overline{\eta}) m_{\widehat{\mathcal{L}}}(\overline{\sigma})$. 

In this case, we also localise $\overline{\eta}$ depending on its direction. 

\paragraph{1.} Let us localise $\overline{\eta}$ by $1 - m_{\widehat{\mathcal{C}}}$. We develop: {\footnotesize 
\begin{subequations}
\begin{align}
&m_{\widehat{\mathcal{L}}}(\overline{\xi}) \widehat{X}_{\beta}(\overline{\xi}) \cdot \nabla_{\overline{\xi}} \widehat{I}_{\mu_{HBH} (1 - m_{\widehat{\mathcal{C}}}) m_{\widehat{\mathcal{L}}}}[f, f](t, \overline{\xi}) \notag\\
&\quad = \int_0^t \int e^{i s \varphi} i s \widehat{X}_{\beta}(\overline{\eta}) \cdot \nabla_{\overline{\xi}} \varphi ~ \mu_{HBH}(\overline{\xi}, \overline{\eta}) m_{\widehat{\mathcal{L}}}(\overline{\xi}) \left( 1 - m_{\widehat{\mathcal{C}}}(\overline{\eta}) \right) m_{\widehat{\mathcal{L}}}(\overline{\sigma}) i \xi_0 \widehat{f}(s, \overline{\eta}) \widehat{f}(s, \overline{\sigma}) ~ d\overline{\eta} ds \label{estapriorimauvaisbLHBHL-s} \\
&\quad \quad + \int_0^t \int e^{i s \varphi} \mu_{HBH}(\overline{\xi}, \overline{\eta}) m_{\widehat{\mathcal{L}}}(\overline{\xi}) \left( 1 - m_{\widehat{\mathcal{C}}}(\overline{\eta}) \right) m_{\widehat{\mathcal{L}}}(\overline{\sigma}) i \xi_0 \widehat{f}(s, \overline{\eta}) \widehat{X}_{\beta}(\overline{\eta}) \cdot \nabla_{\overline{\xi}} \widehat{f}(s, \overline{\sigma}) ~ d\overline{\eta} ds \label{estapriorimauvaisbLHBHL-eta} \\
&\quad \quad + \int_0^t \int e^{i s \varphi} m_{\widehat{\mathcal{L}}}(\overline{\xi}) \widehat{X}_{\beta}(\overline{\xi}) \cdot \nabla_{\overline{\xi}} \left( \mu_{HBH}(\overline{\xi}, \overline{\eta}) \left( 1 - m_{\widehat{\mathcal{C}}}(\overline{\eta}) \right) m_{\widehat{\mathcal{L}}}(\overline{\sigma}) i \xi_0 \right) \widehat{f}(s, \overline{\eta}) \widehat{f}(s, \overline{\sigma}) ~ d\overline{\eta} ds \label{estapriorimauvaisbLHBHL-dersymb} 
\end{align}
\end{subequations} }
\eqref{estapriorimauvaisbLHBHL-dersymb} can be estimated as in the generic $BBB$ case. 

In \eqref{estapriorimauvaisbLHBHL-eta}, we separate into: {\footnotesize 
\begin{subequations}
\begin{align}
&\begin{aligned}
\eqref{estapriorimauvaisbLHBHL-eta} = \int_0^t \int e^{i s \varphi} \mu_{HBH}(\overline{\xi}, \overline{\eta}) \left( 1 - m_{\widehat{\mathcal{C}}}(\overline{\eta}) \right) m_{\widehat{\mathcal{L}}}(\overline{\sigma}) i \widehat{f}(s, \overline{\eta}) \frac{\xi_0 + \sigma_0}{2} \frac{m_{\widehat{\mathcal{L}}}(\overline{\xi}) + m_{\widehat{\mathcal{L}}}(\overline{\sigma})}{2} \\
\widehat{h}_{\beta *}(s, \overline{\sigma}) ~ d\overline{\eta} ds 
\end{aligned} \label{estapriorimauvaisbLHBHL-eta-symh} \\
&\quad + \int_0^t \int e^{i s \varphi} \mu_{HBH}(\overline{\xi}, \overline{\eta}) \left( 1 - m_{\widehat{\mathcal{C}}}(\overline{\eta}) \right) m_{\widehat{\mathcal{L}}}(\overline{\sigma}) i \widehat{f}(s, \overline{\eta}) \frac{\xi_0 + \sigma_0}{2} \frac{m_{\widehat{\mathcal{L}}}(\overline{\xi}) + m_{\widehat{\mathcal{L}}}(\overline{\sigma})}{2} \widehat{g}_{\beta *}(s, \overline{\sigma}) ~ d\overline{\eta} ds \label{estapriorimauvaisbLHBHL-eta-symg} \\
&\quad \begin{aligned}
+ \int_0^t \int e^{i s \varphi} \mu_{HBH}(\overline{\xi}, \overline{\eta}) m_{\widehat{\mathcal{L}}}(\overline{\sigma}) i \left( \xi_0 m_{\widehat{\mathcal{L}}}(\overline{\xi}) \widehat{X}_{\beta}(\overline{\eta}) - \frac{\xi_0 + \sigma_0}{2} \frac{m_{\widehat{\mathcal{L}}}(\overline{\xi}) + m_{\widehat{\mathcal{L}}}(\overline{\sigma})}{2} \widehat{X}_{\beta}(\overline{\sigma}) \right) \cdot \nabla_{\overline{\xi}} \widehat{f}(s, \overline{\sigma}) \\
\left( 1 - m_{\widehat{\mathcal{C}}}(\overline{\eta}) \right) \widehat{f}(s, \overline{\eta}) ~ d\overline{\eta} ds \end{aligned}  \label{estapriorimauvaisbLHBHL-eta-2} 
\end{align}
\end{subequations} }
Like \eqref{estapriorimauvaisbCBBBPC-eta-symh} in the case of the cone for the $(BBB, \widehat{\mathcal{P}}\widehat{\mathcal{C}})$ interaction, \eqref{estapriorimauvaisbLHBHL-eta-symh}'s contribution in \eqref{estapriorimauvaisbDuhamelhcarre} is identically zero. We can also estimate \eqref{estapriorimauvaisbLHBHL-eta-symg}: 
\begin{align*}
\Vert \partial_t \eqref{estapriorimauvaisbLHBHL-eta-symg} \Vert_{L^2} &\lesssim \Vert u \Vert_{L^4} \Vert m_{\widehat{\mathcal{L}}}(D) \nabla g_{\beta *}(t) \Vert_{L^4} \lesssim \left( t^{-\frac{5}{6}} \langle t \rangle^{-\frac{1}{4}+100\delta} \right) \langle t \rangle^{\frac{5}{48}+202\delta} \Vert u \Vert_X^2 \\
\Vert |\overline{\xi}|^{\frac{1}{2}} \langle \overline{\xi} \rangle^{\frac{1}{2}} \partial_t \eqref{estapriorimauvaisbLHBHL-eta-symg} \Vert_{L^2} &\lesssim \Vert u \Vert_{L^4} \Vert m_{\widehat{\mathcal{L}}}(D) |\nabla|^{\frac{3}{2}} \langle \nabla \rangle^{\frac{1}{2}} g_{\beta *}(t) \Vert_{L^4} \lesssim \left( t^{-\frac{5}{6}} \langle t \rangle^{-\frac{1}{4}+100\delta} \right) \langle t \rangle^{202 \delta} \Vert u \Vert_X^2 
\end{align*}

Then, we apply an integration by parts on \eqref{estapriorimauvaisbLHBHL-eta-2}: {\footnotesize 
\begin{subequations}
\begin{align}
&\begin{aligned}
\eqref{estapriorimauvaisbLHBHL-eta-2} = 
\int_0^t \int e^{i s \varphi} i s \left( \xi_0 m_{\widehat{\mathcal{L}}}(\overline{\xi}) \widehat{X}_{\beta}(\overline{\eta}) - \frac{\xi_0 + \sigma_0}{2} \frac{m_{\widehat{\mathcal{L}}}(\overline{\xi}) + m_{\widehat{\mathcal{L}}}(\overline{\sigma})}{2} \widehat{X}_{\beta}(\overline{\sigma}) \right) \cdot \nabla_{\overline{\eta}} \varphi ~ \\
\mu_{HBH}(\overline{\xi}, \overline{\eta}) \left( 1 - m_{\widehat{\mathcal{C}}}(\overline{\eta}) \right) m_{\widehat{\mathcal{L}}}(\overline{\sigma}) i \widehat{f}(s, \overline{\eta}) \widehat{f}(s, \overline{\sigma}) ~ d\overline{\eta} ds 
\end{aligned} \label{estapriorimauvaisbLHBHL-eta-2-s} \\
&\quad \begin{aligned}
+ \int_0^t \int e^{i s \varphi} \left( \xi_0 m_{\widehat{\mathcal{L}}}(\overline{\xi}) \widehat{X}_{\beta}(\overline{\eta}) - \frac{\xi_0 + \sigma_0}{2} \frac{m_{\widehat{\mathcal{L}}}(\overline{\xi}) + m_{\widehat{\mathcal{L}}}(\overline{\sigma})}{2} \widehat{X}_{\beta}(\overline{\sigma}) \right) \cdot \nabla_{\overline{\eta}} \widehat{f}(s, \overline{\eta}) \\
\mu_{HBH}(\overline{\xi}, \overline{\eta}) \left( 1 - m_{\widehat{\mathcal{C}}}(\overline{\eta}) \right) m_{\widehat{\mathcal{L}}}(\overline{\sigma}) i \widehat{f}(s, \overline{\sigma}) ~ d\overline{\eta} ds 
\end{aligned} \label{estapriorimauvaisbLHBHL-eta-2-eta} \\
&\quad
+ \int_0^t \int \nabla_{\overline{\eta}} \cdot \left( \left( \xi_0 m_{\widehat{\mathcal{L}}}(\overline{\xi}) \widehat{X}_{\beta}(\overline{\eta}) - \frac{\xi_0 + \sigma_0}{2} \frac{m_{\widehat{\mathcal{L}}}(\overline{\xi}) + m_{\widehat{\mathcal{L}}}(\overline{\sigma})}{2} \widehat{X}_{\beta}(\overline{\sigma}) \right) \mu_{HBH}(\overline{\xi}, \overline{\eta}) \left( 1 - m_{\widehat{\mathcal{C}}}(\overline{\eta}) \right) m_{\widehat{\mathcal{L}}}(\overline{\sigma}) \right) \notag \\
&\quad \quad \quad \quad \quad \quad \quad \quad \quad \quad \quad \quad \quad \quad \quad \quad \quad \quad e^{i s \varphi} i \widehat{f}(s, \overline{\eta}) \widehat{f}(s, \overline{\sigma}) ~ d\overline{\eta} ds
\label{estapriorimauvaisbLHBHL-eta-2-dersymb}
\end{align}
\end{subequations} }
\eqref{estapriorimauvaisbLHBHL-eta-2-dersymb} is similar to \eqref{estapriorimauvaisbLHBHL-dersymb}. For \eqref{estapriorimauvaisbLHBHL-eta-2-eta}, we notice that 
\begin{align}
\xi_0 m_{\widehat{\mathcal{L}}}(\overline{\xi}) \widehat{X}_{\beta}(\overline{\eta}) - \frac{\xi_0 + \sigma_0}{2} \frac{m_{\widehat{\mathcal{L}}}(\overline{\xi}) + m_{\widehat{\mathcal{L}}}(\overline{\sigma})}{2} \widehat{X}_{\beta}(\overline{\sigma}) &= O\left( \overline{\eta} \right) \label{estapriorimauvaisbLHBHL-gaineta} 
\end{align}
so we can estimate it like \eqref{estapriorimauvaisbLgenBBB-2} of the generic $BBB$ case. 

Finally, we group \eqref{estapriorimauvaisbLHBHL-s} with \eqref{estapriorimauvaisbLHBHL-eta-2-s}: {\footnotesize 
\begin{align}
&\eqref{estapriorimauvaisbLHBHL-s} + \eqref{estapriorimauvaisbLHBHL-eta-2-s} \notag \\
&\quad \begin{aligned}
&= \int_0^t \int e^{i s \varphi} i s \left( \xi_0 m_{\widehat{\mathcal{L}}}(\overline{\xi}) \widehat{X}_{\beta}(\overline{\eta}) \cdot \left( \nabla_{\overline{\xi}} + \nabla_{\overline{\eta}} \right) - \frac{\xi_0 + \sigma_0}{2} \frac{m_{\widehat{\mathcal{L}}}(\overline{\xi}) + m_{\widehat{\mathcal{L}}}(\overline{\sigma})}{2} \widehat{X}_{\beta}(\overline{\sigma}) \cdot \nabla_{\overline{\eta}} \right) \varphi ~ \\
&\pushright{ \mu_{HBH}(\overline{\xi}, \overline{\eta}) \left( 1 - m_{\widehat{\mathcal{C}}}(\overline{\eta}) \right) m_{\widehat{\mathcal{L}}}(\overline{\sigma}) i \widehat{f}(s, \overline{\eta}) \widehat{f}(s, \overline{\sigma}) ~ d\overline{\eta} ds }
\end{aligned} \label{estapriorimauvaisbLHBHL-s-final} 
\end{align} }

But here, 
\begin{align*}
\partial_{\eta_0} \varphi &= 3 \sigma_0^2 + |\sigma|^2 - 3 \eta_0^2 - |\eta|^2 \simeq |\overline{\sigma}|^2 
\end{align*}
so $1 = O(\partial_{\eta_0} \varphi)$. Moreover, by \eqref{estapriorimauvaisbLHBHL-gaineta}, we can apply an integration by parts and get {\footnotesize 
\begin{subequations}
\begin{align}
\eqref{estapriorimauvaisbLHBHL-s-final} &= \int_0^t \int e^{i s \varphi} \overline{\eta} \mu_0(\overline{\xi}, \overline{\eta}) ~ \mu_{HBH}(\overline{\xi}, \overline{\eta}) \left( 1 - m_{\widehat{\mathcal{C}}}(\overline{\eta}) \right) m_{\widehat{\mathcal{L}}}(\overline{\sigma}) \widehat{f}(s, \overline{\eta}) \partial_{\eta_0} \widehat{f}(s, \overline{\sigma}) ~ d\overline{\eta} ds \label{estapriorimauvaisbLHBHL-s-final-1} \\
&\quad + \int_0^t \int e^{i s \varphi} \overline{\eta} \mu_0(\overline{\xi}, \overline{\eta}) ~ \mu_{HBH}(\overline{\xi}, \overline{\eta}) \left( 1 - m_{\widehat{\mathcal{C}}}(\overline{\eta}) \right) m_{\widehat{\mathcal{L}}}(\overline{\sigma}) \partial_{\eta_0} \widehat{f}(s, \overline{\eta}) \widehat{f}(s, \overline{\sigma}) ~ d\overline{\eta} ds \label{estapriorimauvaisbLHBHL-s-final-2} \\
&\quad + \int_0^t \int e^{i s \varphi} \partial_{\eta_0} \left( \overline{\eta} \mu_0(\overline{\xi}, \overline{\eta}) ~ \mu_{HBH}(\overline{\xi}, \overline{\eta}) \left( 1 - m_{\widehat{\mathcal{C}}}(\overline{\eta}) \right) m_{\widehat{\mathcal{L}}}(\overline{\sigma}) \right) \widehat{f}(s, \overline{\eta}) \widehat{f}(s, \overline{\sigma}) ~ d\overline{\eta} ds \label{estapriorimauvaisbLHBHL-s-final-3} 
\end{align}
\end{subequations} }
for some symbol $\mu_0$ of order $0$. We can then estimate \eqref{estapriorimauvaisbLHBHL-s-final-3} like \eqref{estapriorimauvaisbLHBHL-dersymb}, and then since 
\begin{align*}
\partial_{\eta_0} &= \frac{\sigma_0}{|\overline{\sigma}|} \widehat{X}_a(\overline{\sigma}) \cdot \nabla_{\overline{\eta}} + \frac{|\sigma|}{|\overline{\sigma}|} \widehat{X}_b(\overline{\sigma}) \cdot \nabla_{\overline{\eta}}
\end{align*}
we deduce that 
\begin{align*}
\Vert \langle \overline{\xi} \rangle \partial_t \eqref{estapriorimauvaisbLHBHL-s-final-1} \Vert_{L^2} &\lesssim \Vert \nabla u(t) \Vert_{L^{\infty}} \left( \Vert \langle \nabla \rangle X_a f(t) \Vert_{L^2} + \Vert \langle \nabla \rangle m_b(D) X_b f(t) \Vert_{L^2} \right) \\
&\lesssim t^{-\frac{5}{6}} \langle t \rangle^{-\frac{1}{6}+100\delta} \Vert u \Vert_X^2 \\
\Vert \partial_t \eqref{estapriorimauvaisbLHBHL-s-final-2} \Vert_{L^2} &\lesssim \Vert (1 - m_{\widehat{\mathcal{C}}}(D)) (x, y) f(t) \Vert_{L^2} \Vert m_{\widehat{\mathcal{L}}}(D) \nabla u(t) \Vert_{L^{\infty}} \\
&\lesssim \left( t^{-\frac{5}{6}} \langle t \rangle^{-\frac{1}{4}+100\delta} \right) \langle t \rangle^{\frac{5}{48}+202\delta} \Vert u \Vert_X^2 \\
\Vert |\overline{\xi}|^{\frac{1}{2}} \langle \overline{\xi} \rangle^{\frac{1}{2}} \partial_t \eqref{estapriorimauvaisbLHBHL-s-final-2} \Vert_{L^2} &\lesssim \Vert |\nabla|^{\frac{1}{2}} \langle \nabla \rangle^{\frac{1}{2}} (1 - m_{\widehat{\mathcal{C}}}(D)) (x, y) f(t) \Vert_{L^2} \Vert m_{\widehat{\mathcal{L}}}(D) \nabla u(t) \Vert_{L^{\infty}} \\
&\lesssim \left( t^{-\frac{5}{6}} \langle t \rangle^{-\frac{1}{4}+100\delta} \right) \langle t \rangle^{202\delta} \Vert u \Vert_X^2
\end{align*}
as wanted. 

\paragraph{2.} Let us now localise by $m_{\widehat{\mathcal{C}}}(\overline{\eta})$. 

The beginning of the decomposition is exactly the same as in the case 1. above: but when we get to \eqref{estapriorimauvaisbLHBHL-eta-2}, we do not apply an integration by parts and, using again \eqref{estapriorimauvaisbLHBHL-gaineta}, we can estimate directly in a similar way as \eqref{estapriorimauvaisbLgenBBB-2} from the generic $BBB$ case. 

It then only remains \eqref{estapriorimauvaisbLHBHL-s}, but we then have that 
\begin{align*}
\widehat{X}_a(\overline{\eta}) \cdot \nabla_{\overline{\eta}} \varphi &= \frac{\eta_0}{|\overline{\eta}|} \left( 3 \sigma_0^2 + |\sigma|^2 - 3 \eta_0^2 - |\eta|^2 \right) + \frac{\eta}{|\overline{\eta}|} \cdot \left( 2 \sigma_0 \sigma - 2 \eta_0 \eta \right) \\
&\simeq |\overline{\sigma}|^2 
\end{align*}
so $1 = O\left( \widehat{X}_a(\overline{\eta}) \cdot \nabla_{\overline{\eta}} \varphi \right)$. Since we also have
\begin{align*}
\nabla_{\overline{\xi}} \varphi &= O\left( \overline{\eta} \right) 
\end{align*}
we can apply an integration by parts and estimate everything in a similar way, we skip the details.  

\subsubsection{Generic \texorpdfstring{$BHH$}{BHH} interaction}

Let us localise by $\mu_{BHH}(\overline{\xi}, \overline{\eta}) m_{A_1}(\overline{\eta}) m_{A_2}(\overline{\sigma})$ with $A_1 , A_2 \in \{ \widehat{\mathcal{R}}, \widehat{\mathcal{P}}, \widehat{\mathcal{L}} \}$. 

In this case, we have by Lemma \ref{lem-non-res-loin0Cone}
\begin{align*}
\xi_0 &= O(\varphi) + O(\nabla_{\overline{\eta}} \varphi)
\end{align*}
so we can estimate like in the generic $BBB$ case noting that we cannot have simultaneously $A_1 = \widehat{\mathcal{L}}$ and $A_2 = \widehat{\mathcal{P}}$ (ou symmetrically), which justifies that we can indeed apply exactly the same estimates (if no geometric area is equal to $\widehat{\mathcal{P}}$), or a similar estimate losing on the $L^{\infty}$ bound but gaining on the weighted $L^2$ one (if one is $\widehat{\mathcal{P}}$), in a similar way as in the $(BBB, \widehat{\mathcal{P}} \widehat{\mathcal{R}})$ case. We skip the details. 

\subsubsection{Interaction \texorpdfstring{$(BHH, \widehat{\mathcal{C}})$}{(BHH, C)}} 

Let us finally localise by $\mu_{BHH}(\overline{\xi}, \overline{\eta}) m_{A_1}(\overline{\eta}) m_{A_2}(\overline{\sigma})$ with $A_1 = \widehat{\mathcal{C}}$ or $A_2 = \widehat{\mathcal{C}}$. By symmetry we can assume $A_1 = \widehat{\mathcal{C}}$, then for the integral not to be identically zero we need $A_2 \in \{ \widehat{\mathcal{R}}, \widehat{\mathcal{C}} \}$. 

Note that, in this case, we can exactly apply the corresponding part of the proof of the decomposition lemma for the $b$ vector field: indeed, in this case, we did not use the presence of a factor $\frac{|\xi|}{|\overline{\xi}|}$ (so of $m_b$ or $m_c$). 

More precisely, we have that 
\begin{align*}
\xi_0 \simeq \sqrt{3} |\xi_0| - |\xi| 
\end{align*}
and that 
\begin{align*}
\overline{\eta}_a^{\overline{\eta}} - \overline{\sigma}_a^{\overline{\sigma}} &= \sqrt{3} (|\eta_0| - |\sigma_0|) + |\eta| - |\sigma| \\
&= \sqrt{3} \frac{\eta_0}{|\eta_0|} \xi_0 + O(\xi) \simeq \xi_0 
\end{align*}
so there is no difficulty to apply a non-singular decomposition and estimate.

\section{Continuity of the X-norm} \label{section-continuite}

In this last section, we justify the continuity of the $\Vert \cdot \Vert_X$-norm. 

\paragraph{Boundedness implies continuity.} By the estimates from previous sections, continuity of this norm follows from its boundedness. Consider $u_0$ such that $\Vert (1 + x^2 + |y|^2) u_0 \Vert_{H^1} \leq \varepsilon$ for some $\varepsilon \ll 1$, and $u$ a solution defined in $H^1$ on some time interval $[0, T]$, $T > 0$. Assume also that, on $[0, T]$, $\Vert u \Vert_{X, T} < \infty$. 

Take $0 \leq t_1 \leq t_2 \leq T$. Then the estimates from previous sections imply 
\begin{align*}
\Vert u \Vert_{X, t_2} - \Vert u \Vert_{X, t_1} \underset{t_2 - t_1 \to 0}{\longrightarrow} 0 
\end{align*}
Indeed, for all the norms not involving $g$ entering in the definition of $\Vert u \Vert_{X, t}$, for instance $\Vert h_a \Vert_{L^2}$, we can write
\begin{align*}
\Vert h_a(t_2) \Vert_{L^2}^2 - \Vert h_a(t_1) \Vert_{L^2}^2 \leq \int_{t_1}^{t_2} A(s) ~ ds 
\end{align*}
for some expression $A(s)$. The point is that $A(s)$ is at least locally in time in $s$ if $\Vert u \Vert_{X, T}$ is bounded. 

This is immediate for most of the estimates: there is only something to do differently for the terms having a factor $\varphi$ on which we applied an integration by parts in time in order to obtain a good enough estimate. Such terms are of two kinds: either they are the interaction of $u$ with some $v$ in the estimate of $\Vert v \Vert_{L^2}$, or they are interactions of lower order terms than $v$, for instance of $h_b$ with $h_b$ in the estimate of $\Vert \partial_x m_b(D) X_b h_b(t) \Vert_{L^2}$. In the first case, we can simply use $\varphi$ to distribute derivatives and, up to losing some time decay by estimating $\Vert \nabla u(t) \Vert_{L^{\infty}}$ instead of $\Vert \partial_x u(t) \Vert_{L^{\infty}}$, we still obtain that $A \in L^1_{loc}$. In the other case, we can apply the integration by parts in time nonetheless since the terms are of lower order. We may then obtain boundary terms, but those can be controlled by the $X$-norm, and the difference of the two boundary terms goes to $0$ as $t_2 - t_1 \to 0$ by continuity of these lower-order terms. 

In particular, when we do not apply another integration by parts in time, it is not necessary to assume that $\Vert u \Vert_{X, T}$ is small to obtain some (suboptimal) a priori estimate. 

On the other hand, for the norms involving $g$, since $g$ is a quadratic expression in $u$, which is typically estimated as a $L^2 \times L^{\infty}$ product, the $L^2$-continuity follows from what we already explained, and the $L^{\infty}$-continuity, by running again the $L^{\infty}$ dispersive estimate, can be decomposed itself between continuity on norms of $h_{\alpha}$ or $u$ without any weight, which we already have, and continuity of a $g$-depending norm. For this last one, we may use the smallness assumption: indeed, putting everything together, we obtain an estimate of the form 
\begin{align*}
\Vert g_{\alpha}(t_1) - g_{\alpha}(t_2) \Vert \lesssim \Theta(t_1, t_2) + \Vert u \Vert_{X, T} \left( \Vert g_{\alpha}(t_1) - g_{\alpha}(t_2) \Vert \right) 
\end{align*}
where the same argument can be made for all the desired norms. Here above, $\Theta$ is some continuous function such that $\Theta(t, t) = 0$ for every $t$. But since $\Vert u \Vert_{X, T}$ is assumed to be very small, it can absorb the implicit universal constant and the second term on the right can be absorbed on the left. 

From this we deduce that, as soon as $\Vert u \Vert_{X, T}$ is bounded, then $t \mapsto \Vert u \Vert_{X, t}$ is upper semi-continuous, and even continuous by reversibility of the equation. Then the a priori estimate applies and $\Vert u \Vert_{X, T}$ remains smaller than $2 \varepsilon$ for $\varepsilon > 0$ small enough with respect to the universal constants of the a priori estimate. 

\paragraph{High regularity case.} Let us first consider the case where $u_0 \in H^N$ for some $N \gg 1$, with always $\Vert (1 + x^2 + |y|^2) u_0 \Vert_{H^1} \leq \varepsilon$. 

From the local well-posedness in $H^N$, we can build a solution $u \in C^0([0, T], H^N)$ for some $T > 0$. 

Furthermore, to propagate weighted estimates, we consider the equation 
\begin{align*}
\left\{ \begin{array}{l}
\partial_t v + \partial_x \Delta v + \frac{1}{2} \partial_x \left( u v \right) = \left( \Delta u + 2 \partial_x^2 u + \frac{1}{2} u^2, 2 \partial_x \nabla u \right) \\
v(t = 0) = (x, y) u_0 
\end{array} \right. 
\end{align*}
where $v = (v_0, v_1, v_2)$ is vector-valued. The right-hand side is a source term independant of $v$, and very regular since we assumed $u \in C^0_t H^N_{x, y}$. On the other hand, on the left-hand side, the equation is linear in $v$, and similar estimates as for the local well-posedness of the Zakharov-Kuznetsov equation ensure that $v$ can be defined locally in time, valued in $H^1_{x, y}$. In fact, if $u_0 \in H^N$ and $(1 + x^2 + |y|^2) u_0 \in H^1$, we have automatically by standard interpolation that $(x, y) u_0 \in H^k$ for some $N \gg k \gg 1$, so we may construct $v \in C^0_t H^k_{x, y}$ locally. 

Applying the same approach for the quadratic weight, we obtain 
\begin{align*}
\left\{ \begin{array}{l}
\partial_t w + \partial_x \Delta w + \frac{1}{2} \partial_x \left( u w \right) = 6 \partial_x^2 v_0 + 6 \partial_x u + 4 \partial_x \nabla \cdot (v_1, v_2) + 2 \partial_x u + \partial_x (u v_0) \\
w(t = 0) = (x^2 + |y|^2) u_0
\end{array} \right. 
\end{align*}
Again, the right-hand side is regular and has been already constructed independantly of $w$, while the equation on $w$ is linear in $w$ and hence we may construct $w \in C^0_t H^1_{x, y}$ locally in time. 

In particular, locally in time, $\Vert (1 + x^2 + |y|^2) u(t) \Vert_{H^1} + \Vert u(t) \Vert_{H^N}$ is bounded and even continuous. But taking the supremum of this norm in time, it is strictly stronger than $\Vert u \Vert_{X, t}$, and hence $\Vert u \Vert_{X, t}$ is locally bounded in time, and therefore continuous. Here we use that $\Vert (x, y) f \Vert_{L^2} \lesssim \Vert (x, y) u \Vert_{L^2} + \Vert u \Vert_{H^3}$ and similarly for higher weights and higher regularity. 

We now propagate the estimate and the continuity. From the a priori estimate, $\Vert u \Vert_{X, t}$ remains small. Furthermore, from the estimates made on the $H^N$ norm, we have that 
\begin{align*}
\Vert u(t) \Vert_{H^N}^2 \lesssim \Vert u_0 \Vert_{H^N}^2 + \Vert u \Vert_{X, t} \Vert u \Vert_{L^{\infty}([0, t], H^N)}^2 \left( 1 + \Vert u \Vert_{X, t} \right) 
\end{align*}
In particular, since we know that $\Vert u \Vert_{X, t} \leq 2 \varepsilon$ as long as $u$ is well-defined in $H^N$, $v$ in $H^k$ and $w$ in $H^1$, if $\varepsilon$ is small enough (with respect to the universal constant) then $\Vert u(t) \Vert_{H^N}$ remains uniformly bounded in time, even if $\Vert u_0 \Vert_{H^N}^2$ is not assumed to be small. 

Then, since the equation for $v$ is linear in $v$, standard Gronwall estimate ensure that $\Vert v \Vert_{H^k}$ cannot blow up in finite time, and similarly for $w$. 

In particular, we deduce that the solution can be globally defined and scatters in $H^1$ (or in $H^s$ for any $s < N$ here). 

\paragraph{$H^1$ regularity.} Let us finally consider the case of an initial data $u_0$ such that $\Vert (1 + x^2 + |y|^2) u_0 \Vert_{H^1} \leq \varepsilon$ without any additional regularity. 

In this case, let us consider $u_0^{\tau} = u_0 \ast \rho_{\tau}$ for $\rho$ a regularisation kernel and 
\begin{align*}
\rho_{\tau}(x, y) &= \tau^{-d} \rho\left( \tau^{-1} (x, y) \right) 
\end{align*}
In particular, $\widehat{\rho_{\tau}}(\overline{\xi}) = \widehat{\rho}(\tau \overline{\xi})$. 

Note now that $u_0^{\tau} \in H^N$ for every $\tau > 0$ and that, for $\tau \in (0, 1)$, 
\begin{align*}
\Vert (1 + x^2+|y|^2) u_0^{\tau} \Vert_{H^1} &\lesssim \Vert u_0^{\tau} \Vert_{H^1} + \Vert |\overline{\xi}| \nabla_{\overline{\xi}}^2 \left( \widehat{u_0} \widehat{\rho_{\tau}} \right) \Vert_{L^2} \\
&\lesssim \Vert u_0 \Vert_{H^1} + \Vert |\overline{\xi}| \nabla_{\overline{\xi}}^2 \widehat{u_0} \Vert_{L^2} + \tau \Vert |\overline{\xi}| \nabla_{\overline{\xi}} \widehat{u_0} \Vert_{L^2} + \tau^2 \Vert |\overline{\xi}| \widehat{u_0} \Vert_{L^2} \\
&\lesssim \Vert (1 + x^2 + |y|^2) u_0 \Vert_{H^1} 
\end{align*}
In particular, uniformly in $\tau$, $u_0^{\tau}$ satisfies $\Vert (1 + x^2 + |y|^2) u_0^{\tau} \Vert_{H^1} \lesssim \varepsilon$. 

We may therefore apply the estimates of the regular case for the solution $u^{\tau}$ obtained from $u_0^{\tau}$, which is defined globally with $\Vert u^{\tau} \Vert_{X, T = \infty} \lesssim \varepsilon$ uniformly in $\tau$, and then pass to the limit $\tau \to 0$. From theglobal well-posedness in $H^1$ of the equation, we deduce that $u^{\tau}(t)$ converges in $H^1$ to $u(t)$ for every $t \geq 0$. On the other hand, for all the weights involved in the definition of $\Vert u^{\tau} \Vert_X$, the boundedness uniformly in $\tau$ ensures that there is at least a weak limit, and this weak limit has to be $u$ with the same weight (by unicity in the space of distributions). In particular, we deduce that $\Vert u \Vert_{X, T}$ is bounded for any $T > 0$. 

But as we have seen, boundedness implies continuity, and therefore global stability and scattering. This concludes.

\section*{Acknowledgements}

I am particularly grateful to Frédéric Rousset for his advices and support during the preparation of this paper. I would also like to thank several people for exchanges and discussions: Yann Brenier, Lucas Ertzbischoff, Patrick Gérard, Pierre Germain, Didier Pilod, Jean-Claude Saut, and Frédéric Valet.

\end{document}